\documentclass[letterpaper]{amsbook}

\usepackage[letterpaper]{geometry}
\usepackage{amsmath, amsthm, amsfonts, amsbsy, thmtools, amssymb, bm,textcmds,mathscinet}
\usepackage[sans]{dsfont}

\usepackage{thm-restate}
\usepackage{hyperref}
\usepackage[mathscr]{euscript}
\usepackage{tikz}
\usetikzlibrary{arrows}
\usepackage{cleveref}
\usepackage[utf8]{inputenc}
\usepackage[T1]{fontenc}
\usepackage{tikz}

\numberwithin{figure}{chapter}

\usepackage{imakeidx} 

\newtheorem{thm}{Theorem}[section]
\newtheorem{prop}[thm]{Proposition}
\newtheorem{lem}[thm]{Lemma}
\newtheorem{cor}[thm]{Corollary}

\theoremstyle{remark}
\newtheorem{rem}[thm]{Remark}
\theoremstyle{definition}
\newtheorem{definition}[thm]{Definition}
\newtheorem{que}[thm]{Question}
\newtheorem{example}[thm]{Example}

\usepackage{gensymb,epic,eepic,float}
\usepackage{caption}
\usepackage{siunitx}
\usepackage{empheq}
\usepackage[T1]{fontenc}
\usepackage{mathtools}
\usepackage{comment}
\usepackage{enumerate}
\usepackage[shortlabels]{enumitem}
\usetikzlibrary{graphs,patterns,decorations.markings,,matrix}
\usetikzlibrary{calc,decorations.pathmorphing,decorations.pathreplacing,shapes}
\allowdisplaybreaks
\numberwithin{equation}{section}

\makeatletter
\newsavebox{\@brx}
\newcommand{\llangle}[1][]{\savebox{\@brx}{\(\m@th{#1\langle}\)}%
	\mathopen{\copy\@brx\kern-0.5\wd\@brx\usebox{\@brx}}}
\newcommand{\rrangle}[1][]{\savebox{\@brx}{\(\m@th{#1\rangle}\)}%
	\mathclose{\copy\@brx\kern-0.5\wd\@brx\usebox{\@brx}}}
\makeatother

\DeclareMathOperator{\Sign}{Sign}

\DeclareMathOperator{\SeqSign}{SeqSign}

\DeclareMathOperator{\sgn}{sgn}
\DeclareMathOperator{\End}{End} 
\DeclareMathOperator{\inv}{inv}

\DeclareMathOperator{\height}{ht}
\DeclareMathOperator{\spin}{sp}

\def\b{\bar}

\def\leq{\leqslant}
\def\geq{\geqslant}

\def\l{\lambda}
\def\m{\mu}

\def\s{\mathfrak{s}}

\def\ie{{\it i.e.}\/,}
\def\cf{{\it cf.}\/}

\colorlet{lgray}{white!85!black}
\colorlet{lred}{white!85!red}
\colorlet{lgreen}{white!80!green}
\colorlet{dgreen}{black!30!green}
\definecolor{green}{rgb}{0.1,0.8,0.1}
\definecolor{yellow}{rgb}{1.0,0.85,0.25}

\newcommand{\bra}[1]{\left\langle #1\right|}
\newcommand{\ket}[1]{\left|#1\right\rangle}

\newcommand\fs{\footnotesize}

\renewcommand{\tikz}[2]{
	\begin{tikzpicture}[scale=#1,baseline=(current bounding box.center),>=stealth]
	#2
	\end{tikzpicture}}

\def\l{{\sf L}}
\def\m{{\sf M}}

\makeatletter
\newcommand{\subalign}[1]{%
	\vcenter{%
		\Let@ \restore@math@cr \default@tag
		\baselineskip\fontdimen10 \scriptfont\tw@
		\advance\baselineskip\fontdimen12 \scriptfont\tw@
		\lineskip\thr@@\fontdimen8 \scriptfont\thr@@
		\lineskiplimit\lineskip
		\ialign{\hfil$\m@th\scriptstyle##$&$\m@th\scriptstyle{}##$\crcr
			#1\crcr
		}%
	}
}
\makeatother

\newlength{\cellsize}
\newcommand{\tableaucell}[1]{{%
		\def \arg{#1}\def \void{}%
		\ifx \void \arg
		\vbox to \cellsize{\vfil \hrule width \cellsize height 0pt}%
		\else
		\unitlength=\cellsize
		\begin{picture}(1,1)
		\put(0,0){\makebox(1,1){$#1$}}
		\put(0,0){\line(1,0){1}}
		\put(0,1){\line(1,0){1}}
		\put(0,0){\line(0,1){1}}
		\put(1,0){\line(0,1){1}}
		\end{picture}%
		\fi}}

\newcommand\tableau[1]{
	\vcenter{
		\let\\=\cr
		\baselineskip=-16000pt
		\lineskiplimit=16000pt
		\lineskip=0pt
		\halign{&\tableaucell{##}\cr#1\crcr}}}

	\renewcommand\le{\leq}
	\renewcommand\ge{\geq}

	\def\b{\bar}
	
	\def\leq{\leqslant}
	\def\geq{\geqslant}

	\def\l{\lambda}
	\def\m{\mu}
	
	\def\s{\mathfrak{s}}

	\def\ie{{\it i.e.}\/,}
	\def\cf{{\it cf.}\/}
	
	\colorlet{llgray}{white!92!black}
	\colorlet{lgray}{white!85!black}
	\colorlet{lred}{white!85!red}
	\colorlet{lgreen}{white!80!green}
	\colorlet{dgreen}{black!30!green}
	\definecolor{green}{rgb}{0.1,0.8,0.1}
	\definecolor{yellow}{rgb}{1.0,0.85,0.25}

	\renewcommand{\tikz}[2]{
		\begin{tikzpicture}[scale=#1,baseline=(current bounding box.center),>=stealth]
			#2
	\end{tikzpicture}}

	\def\l{{\sf L}}
	\def\m{{\sf M}}

	\def\a{\mathfrak{a}}
	\def\b{\mathfrak{b}}
	\def\c{\mathfrak{c}}
	\def\d{\mathfrak{d}}

	\title{Colored Fermionic Vertex Models and Symmetric Functions}
	
	\author{Amol Aggarwal, Alexei Borodin, and Michael Wheeler}

	\makeindex
	
\begin{document}

	\maketitle
	
	\begin{abstract}
		
		In this text we introduce and analyze families of symmetric functions arising as partition functions for colored fermionic vertex models associated with the quantized affine Lie superalgebra $U_q \big( \widehat{\mathfrak{sl}} (1 | n) \big)$. We establish various combinatorial results for these vertex models and symmetric functions, which include the following. 
		
		\begin{enumerate}
			\item We apply the fusion procedure to the fundamental $R$-matrix for $U_q \big( \widehat{\mathfrak{sl}} (1 | n) \big)$ to obtain an explicit family of vertex weights satisfying the Yang--Baxter equation.
			\item We define families of symmetric functions as partition functions for colored, fermionic vertex models under these fused weights. We further establish several combinatorial properties for these symmetric functions, such as branching rules and Cauchy identities.
			\item We show that the Lascoux--Leclerc--Thibon (LLT) polynomials arise as special cases of these symmetric functions. This enables us to show both old and new properties about the LLT polynomials, including Cauchy identities, contour integral formulas, stability properties, and branching rules under a certain family of plethystic transformations.
			\item A different special case of our symmetric functions gives rise to a new family of polynomials called factorial LLT polynomials. We show they generalize the LLT polynomials, while also satisfying a vanishing condition reminiscent of that satisfied by the factorial Schur functions.
			\item By considering our vertex model on a cylinder, we obtain fermionic partition function formulas for both the symmetric and nonsymmetric Macdonald polynomials. 
			\item We prove combinatorial formulas for the coefficients of the LLT polynomials when expanded in the modified Hall--Littlewood basis, as partition functions for a $U_q \big( \widehat{\mathfrak{sl}} (2 | n) \big)$ vertex model. 
			\end{enumerate}
		
	\end{abstract} 
	
	\tableofcontents

\chapter{Introduction}

\label{Introduction}

The study of exactly solvable, or \emph{integrable}, lattice models is a vast domain that originated at the interface between statistical and quantum mechanics well over 50 years ago; a classical introduction to the subject is Baxter's book \cite{ESMSM}. Over the years it has sprouted multiple branches, some of which have since developed into full-fledged domains of research; a prominent example is the representation theory of quantum groups. A more recent, and so far less studied, branch connects integrable lattice models with the theory of symmetric functions. Earlier works in this direction include those of Kirillov--Reshetikhin \cite{CT}, Fomin--Kirillov \cite{ESFP,PTE}, Lascoux--Leclerc--Thibon \cite{FVE}, Gleizer--Postnikov \cite{CE}, Tsilevich \cite{QISMMSF}, Lascoux \cite{VMAP}, Zinn-Justin \cite{SLTM}, Brubaker--Bump--Friedberg \cite{SPE}, and Korff \cite{CVSFDA}.

 The theory of symmetric functions studies remarkable families of symmetric and associated nonsymmetric polynomials with origins in diverse areas, such as the theory of finite groups, multivariate statistics, representation theory of Lie and $p$-adic groups, harmonic analysis on Riemannian homogeneous spaces, probability theory, mathematical physics, algebraic geometry, and enumerative combinatorics. Macdonald's treatise \cite{SFP} provides a comprehensive introdution to this subject. Perhaps the most well-known families of symmetric and nonsymmetric functions also bear the name of Macdonald, and their special cases include other celebrated families of functions such as Schur, Jack, Hall--Littlewood, and $q$-Whittaker functions.
 
 Although symmetric functions come from very different backgrounds, their structural properties show remarkable uniformity. They typically include elements such as:
 
 \begin{itemize} 
 	\item branching rules, which are combinatorial recipes for decomposing a polynomial of a family into similar polynomials with fewer variables;
 	\item summation identities, including (skew) Cauchy and Pieri formulas;
 	\item combinatorial procedures for computing structure constants, namely, the coefficients arising from decomposing products of functions in a family on a basis of similar functions;
 	\item combinatorial understanding of transition matrices between different families of functions.
 \end{itemize} 

The theory of integrable vertex models has been shown to provide a convenient framework for simultaneously accessing each of the above features, for various families of functions.
 
 \begin{itemize} 
 	
 \item Branching rules for these functions typically follow directly from their representations as partition functions for a suitable lattice model under appropriate boundary data.\footnote{An exception exists in the case of the general Macdonald polynomials. Partition function formulas for them were provided by Cantini--de Gier--Wheeler \cite{MFP} and Borodin--Wheeler \cite{NPIVM}, which do not appear to directly relate to the corresponding branching rules.}
 
 \item Summation identities of Cauchy and Pieri type were proved in many instances by applying the Yang--Baxter commutation relations underlying the solvable lattice model. For earlier examples, we refer to the works of Bump--McNamara--Nakasuji \cite{FFE}, Motegi--Sakai \cite{VMP,QICP}, Borodin \cite{SRF}, Wheeler--Zinn-Justin \cite{RISVMPF}, and Borodin--Petrov \cite{HSVMSRF}. 
 
 \item A lattice model approach to find structure constants was developed by Gleizer--Postnikov \cite{CE}, Zinn-Justin \cite{CIT,CPI,PI,FP}, Wheeler--Zinn-Justin \cite{PI,CPI}, and Knutson--Zinn-Justin \cite{IITF}. It turns out to be closely related to the celebrated puzzle representation for the classical Littlewood--Richardson coefficients that goes back to Knutson--Tao \cite{MTPS,EC} and Knutson--Tao--Woodward \cite{TPDFC}.
 
 \item The Kostka--Foulkes transition matrix between the Schur and Hall--Littlewood symmetric functions was the subject of one of the earliest works \cite{CT} relating symmetric functions and integrable vertex models. Moreover, the Yang--Baxter elements in the Hecke algebra \cite{FVE} were central in the approach developed by Lascoux for transitioning between bases of (non)symmetric functions; see \cite{WP} and references therein.
 
 \end{itemize} 
 
 In addition to providing a new path to understanding previously known families of symmetric and nonsymmetric functions, integrable vertex models turned out to be useful for defining novel families that often have many of the desirable properties satisfied by the classical ones. Examples include the spin Hall--Littlewood functions and spin $q$-Whittaker polynomials. Their symmetric versions were introduced and analyzed in \cite{SRF,HSVMSRF,SP}, and further studied by Garbali--de Gier--Wheeler \cite{NGP} and Muccicioni--Petrov \cite{SPDQ}; their nonsymmetric versions were introduced and analyzed by Borodin--Wheeler in \cite{SVMST}.
 
 The development of the theory of spin Hall--Littlewood functions and $q$-Whittaker polynomials was based on a thorough investigation of partition functions on semi-infinite strips built from vertex weights associated with the quantum affine algebra $U_q \big( \widehat{\mathfrak{sl}}(n+1) \big)$. The rank $n = 1$ case corresponded to symmetric functions, and the higher rank $n > 1$ case to nonsymmetric ones. It was essential to consider not only weights coming from the fundamental $R$-matrix for $U_q \big( \widehat{\mathfrak{sl}} (n + 1) \big)$, but also those obtained by its fusion in both directions.\footnote{In representation theoretic terms, this means specializing the universal $R$-matrix to the tensor product of symmetric powers of the fundamental representation and further analytically continuing in the exponents of those powers.}
 
 The goal of the present work is to develop a general theory for symmetric functions associated with the quantum affine deformation of the general linear Lie superalgebra $U_q \big(\widehat{\mathfrak{sl}}(m | n) \big)$, with $m=1$. The choice $m=1$ was dictated by the fact that it is the case most different from the one related to $U_q \big(\widehat{\mathfrak{sl}}(n+1) \big)$. In the latter case, the vertex models were \emph{bosonic}, meaning that there is no constraint on the number of arrows that can exist along a given edge; in the former, they will be \emph{fermionic}, meaning that there is an exclusion rule preventing arrows of the same color from occupying the same edge. Another reason for restricting to the case $m = 1$ is that, in this setting, basic Cauchy summation identities involve only \emph{symmetric} functions for any choice of the rank $n$, while for $m > 1$ these identities must also incorporate nonsymmetric ones.\footnote{It should be noted, however, that some of our results involve vertex models with $m>1$; more details will be given below.} In this text, we will primarily focus on the $n > 1$ setting; the somewhat different $n=1$ story will be treated separately in \cite{P}. 
 
 Understanding the relationship between the bosonic and fermionic theories has been a key motivation for us. At first glance, they appear quite different; for example, the basic objects in the rank $1$ bosonic and fermionic situations are the Hall--Littlewood \cite{SRF,RISVMPF} and big Schur symmetric functions \cite{PI}, respectively. In fact, the nontrivial interplay between these two families of functions is observed throughout Chapter III of Macdonald's book \cite{SFP}. We were able to observe a few further hidden points of contact as well. One is a surprising fact that cylindrical partition functions of the type considered in \cite{NPIVM} yield symmetric and nonsymmetric Macdonald polynomials for \emph{both} bosonic and fermionic vertex weights. Another is a color merging result, stating that partition functions for fermionic vertex models can be obtained by partially anti-symmetrizing partition functions for certain bosonic vertex models of a higher rank. 
 
 The present work enjoys multiple connections with the recent works of Brubaker--Buciumas--Bump--Gustafsson \cite{VOSLMMF,CVMF,MFSLM}, where vertex models related to $U_q \big(\widehat{\mathfrak{sl}}(m | n) \big)$ also played a central role and some of the same symmetric functions appeared, such as the Lascoux--Leclerc--Thibon (LLT) polynomials of \cite{RT}. While most of our results appear to be different from those proved there, we are hopeful that the relationship between these two directions will become stronger with time. 
 
 The development of the theory of spin Hall--Littlewood and $q$-Whittaker polynomials, as well as various other recent works in the subject, came hand-in-hand with a development of their probabilistic applications; for example, an interested reader might consult the lecture notes of Borodin--Petrov \cite{IPSVMSF} in the rank $1$ case and the paper of Borodin--Gorin--Wheeler \cite{SVMP} in the higher rank case. In the present work we shy away from the probabilistic direction, but the forthcoming paper \cite{P} will include a substantial probabilistic component in the rank $n = 1$. 
 
 Let us now briefly describe our main results. 
 
 1. \emph{Fused weights}: We start by applying the fusion procedure to the fundamental $R$-matrix of Bazhanov--Shadrikov \cite{TSTESS} for the quantum affine superalgebra $U_q \big( \widehat{\mathfrak{sl}} (m | n) \big)$. This yields an explicit family of vertex weights satisfying the Yang--Baxter equation, which are in general given by a (sort of elaborate) sum. Under various specializations of their parameters, we show that these weights simplify considerably, factoring completely.
 
 2. \emph{Symmetric functions}: We define families of symmetric rational functions $F_{\boldsymbol{\lambda} / \boldsymbol{\mu}}$ and\footnote{We also introduce a third family $H_{\boldsymbol{\lambda} / \boldsymbol{\mu}}$, which will play a more subsidiary role in this text.} $G_{\boldsymbol{\lambda} / \boldsymbol{\mu}}$ as partition functions under certain boundary data for the lattice models associated with the $m = 1$ cases of the above (otherwise fully general) fused $U_q \big( \widehat{\mathfrak{sl}} (m | n) \big)$ vertex weights. These skew functions are indexed by $n$-tuples of signatures $\boldsymbol{\lambda} = \big( \lambda^{(1)}, \lambda^{(2)}, \ldots , \lambda^{(n)} \big)$ and $\boldsymbol{\mu} = \big( \mu^{(1)}, \mu^{(2)}, \ldots , \mu^{(n)} \big)$, and they depend on four sets of parameters $(\textbf{x}, \textbf{r}; \textbf{y}, \textbf{s})$. We further establish several combinatorial properties for these functions. These include branching rules, which are direct consequences of the definition of the underlying lattice model; Cauchy identities, which are proven through Yang--Baxter commutation relations, together with a fully factorized expression for the partition function of our vertex model under domain-wall boundary data; and contour integral formulas, which are shown through a reduction to the bosonic $U_q \big( \widehat{\mathfrak{sl}} (n + 1) \big)$ case of \cite{SVMST}, using a color merging result.
 
3. \emph{LLT polynomials}: We next analyze how these functions specialize under various degenerations of their parameter sets. In certain cases when only one (among four) of these sets remains generic, we show that our symmetric functions specialize to skew LLT polynomials $\mathcal{L}_{\boldsymbol{\lambda} / \boldsymbol{\mu}} (\textbf{x})$. 

These polynomials can be viewed as $q$-deformations of products of skew Schur functions $s_{\lambda / \mu} (\textbf{x})$, in that for $q = 1$ we have $\mathcal{L}_{\boldsymbol{\lambda} / \boldsymbol{\mu} } (\textbf{x}) = \prod_{j = 1}^n s_{\lambda^{(j)} / \mu^{(j)}} (\textbf{x})$. They were originally introduced in \cite{RT} as plethystic images of Schur functions under $q$-deformed power sum operators; by setting $q$ to roots of unity, such images are useful for analyzing how compositions of $\text{GL}_n (\mathbb{C})$-representations decompose into irreducible ones \cite{SSFSAP,MECPM}. However, since then, they have been found to be ubiquitous in algebra. Within algebraic combinatorics, they satisfy branching rules and Cauchy identities \cite{TRA}; are generalizations of modified\footnote{Here, the term ``modified'' refers to the plethystic image under the map sending the sum $X = \sum_{i = 1}^{\infty} x_i$ of formal variables (arguments for a symmetric function) to $(1 - q)^{-1} X$.} Hall--Littlewood polynomials \cite{RT} and (after plethysm) certain classes of quasi-symmetric chromatic polynomials \cite{S,PCQF}; and arise as coefficients of modified Macdonald polynomials when expanded in powers of $q$ \cite{CP}. Within enumerative algebraic geometry, they arise in the analysis of certain flag varieties \cite{RT,CSFV,GCMSCNCC,AAPP} and of the Frobenius series for the space of diagonal harmonics \cite{CFCDC,NDP,NSDH,S}. Within representation theory, they appear in the Fock space representation for $U_q \big( \widehat{\mathfrak{sl}} (n) \big)$ \cite{RT,CAP,TPTP} and exhibit relations with the Kazhdan--Lusztig theory for affine Hecke algebras \cite{CAP,AAPP}. 

Our interpretation of the LLT polynomials as partition functions for an integrable lattice model enables us to prove old and new combinatorial properties about them. These include branching rules under a family of plethystic substitutions; both standard and dual Cauchy identities (originally due to Lam \cite{TRA}); their specializations to modified Hall--Littlewood polynomials (originally due to Lascoux--Leclerc--Thibon \cite{RT}); and contour integral formulas for them in terms of nonsymmetric Hall--Littlewood polynomials (originally implicitly due to Grojnowski--Haiman \cite{AAPP}).  

Let us mention that the recent work of Corteel--Gitlin--Keating--Meza \cite{VMFP}, which was pursued independently of ours, also provides a solvable lattice model for the LLT polynomials\footnote{The same vertex model was also proposed by Curran--Yost-Wolff--Zhang--Zhang \cite{LMP}, where the Yang--Baxter equation was confirmed for $n \le 3$ but left open for $n > 3$.} that coincides with one of ours. They further use this vertex model to establish LLT Cauchy identities and degenerations to modified Hall--Littlewood polynomials. 

4. \emph{Factorial LLT polynomials}: By instead letting all but two of their parameter sets remain generic, our functions become inhomogeneous generalizations of the LLT polynomials that we call factorial LLT polynomials. Our reason for this terminology is that we show these polynomials possess a certain vanishing property, similar to those satisfied by factorial Schur and interpolation Macdonald polynomials. Vanishing properties of this type have proven to be central in the theory of (non)symmetric functions, as in many cases they fully characterize the underlying function; see the survey \cite{ISF} of Okounkov. Although the factorial LLT polynomials are not fully characterized (but might ``partially'' be; see \Cref{lambdam0sigma} below) by their vanishing property for general $n$, they are when $n = 1$, in which case we show that they in fact coincide with the factorial Schur functions. 

5. \emph{Macdonald polynomials}: By taking our $U_q \big( \widehat{\mathfrak{sl}} (1 | n) \big)$ vertex model on a cylinder and suitably specializing parameters, we obtain partition function formulas for both the nonsymmetric and symmetric Macdonald polynomials. The nonsymmetric formula differs from the one provided in \cite{NPIVM}, as the latter was bosonic (and therefore involved an infinite sum), while ours here is fermionic (and therefore only involves a finite sum). Comparing this lattice model interpretation for the symmetric Macdonald polynomials with that for the LLT polynomials, we deduce a new proof of an expression (originally due to Haglund--Haiman--Loehr \cite{CP}) for the modified Macdonald polynomials as a linear combination of LLT ones.

6. \emph{Transition coefficients}: We provide a combinatorial procedure for determining the transition coefficients of both modified Macdonald and LLT polynomials in the basis of modified Hall--Littlewood polynomials.\footnote{This also yields the transition coefficients between the (standard) Macdonald and Hall--Littlewood bases.} These coefficients are expressed through partition functions for a suitably specialized fused $U_q \big( \widehat{\mathfrak{sl}} (2 | n) \big)$ vertex model under certain boundary data; it bears close similarities with the puzzle interpretation for other families of expansion coefficients studied in \cite{MTPS,TPDFC,CE,CIT,PI,CPI,IITF,FP}. Let us mention that not all weights for this lattice model are nonnegative as polynomials in the underlying parameters $q$ and $t$, and they cannot be, since the expansion of modified Macdonald or LLT polynomials in the modified Hall--Littlewood basis can involve negative coefficients. By combining our results with known statements for the Kostka--Foulkes transition matrix between Schur and Hall--Littlewood functions, we deduce a (not manifestly nonnegative) lattice model representation for the expansion coefficients of the modified Macdonald and LLT polynomials in the Schur basis. 

Over the past two decades, an extensive literature has developed surrounding the decomposition of LLT and modified Macdonald polynomials into the basis of Schur functions, as the resulting coefficients admit interpretations as fundamental invariants from algebraic geometry \cite{SPP,GCMSCNCC} and representation theory \cite{RT,CAP,AAPP}. Examples of such works include those of Leclerc--Thibon \cite{CAP}, Haiman \cite{SPP}, Grojnowski--Haiman \cite{AAPP}, Assaf \cite{DEGCP, EOP}, Blasiak \cite{OCP}, Alexandersson--Uhlin \cite{CSSPP}, and Foster \cite{CETHP}; however, the combinatorial procedures there seem to be of a different nature from ours, and it is not transparent to us how to match them. The transition coefficients between the modified Macdonald and modified Hall--Littlewood bases have also been considered recently in the work of Mellit \cite{PCVPAF}, who interpreted them through conjugacy classes of nilpotent matrices over a finite field with specified Jordan form. However, enumerating over such classes is intricate, and so turning his interpretation into an explicit combinatorial algorithm for determining these coefficients seems to be difficult. 
  
  Having outlined our results, let us proceed to describe them in more detail. Throughout this text, we fix a complex number $q \in \mathbb{C}$.\index{Q@$q$; quantization parameter}

  \section{Fused Weights}
  
  \label{0Weights} 
  
  The vertex models we consider in this text will be ensembles of directed up-right paths on subdomains of the square lattice; each path in this ensemble will be labeled by a \emph{color}, which is an index in $\{ 1, 2, \ldots n \}$. Every vertex in the domain has some number of colored paths entering and exiting it and, depending on this local path configuration (or \emph{arrow configuration}), the vertex is assigned a weight. The total ensemble weight is then given by the product of the weights of its vertices. 
  
  The arrow configuration at a vertex $v \in \mathbb{Z}^2$ will be indexed by an ordered sequence of four elements $(\textbf{A}, \textbf{B}; \textbf{C}, \textbf{D})$ in $\mathbb{Z}_{\ge 0}^n$. Here, $\textbf{A} = (A_1, A_2, \ldots , A_n)$ counts the number of paths of colors $1, 2, \ldots , n$ vertically entering through $v$, respectively. In the same way, $\textbf{B}$, $\textbf{C}$, and $\textbf{D}$ count the colored paths horizontally entering, vertically exiting, and horizontally exiting $v$, respectively. Thus, one can view $\textbf{A}$, $\textbf{B}$, $\textbf{C}$, and $\textbf{D}$ as the \emph{states} of the south, west, north, and east edges adjacent to $v$, respectively. 
  
  In addition to depending on an arrow configuration $(\textbf{A}, \textbf{B}; \textbf{C}, \textbf{D}$), the weight of a vertex $v \in \mathbb{Z}^2$ will also be governed by several complex parameters. The first among them consist in two pairs of \emph{rapidity parameters} $(x; r)$ and $(y; s)$\index{X@$(x; r), (y; s)$; rapidity parameters}, which are associated with the row and column intersecting to form $v$, respectively; these rapidities $(x; r)$ and $(y; s)$ may vary across the domain but remain constant along rows or columns, respectively. The last is a \emph{quantization parameter} $q$\index{Q@$q$; quantization parameter}, which cannot vary and is fixed throughout the model. This produces five governing parameters, but the vertex weight will in fact only depend on $x$ and $y$ through their quotient $z = \frac{x}{y}$, which is sometimes referred to as a \emph{spectral parameter}. 
  
  We will diagrammatically depict vertex weights by 
  \begin{align*}
  	\tikz{.75}{	
  		\draw[] (-7.75, 0) circle node[scale = 1, left]{$W_{x / y; q} (\textbf{A}, \textbf{B}; \textbf{C}, \textbf{D} \boldsymbol{\mid} r, s) =$} ;
  		\draw[->, ultra thick, gray] (-6, 0) -- (-4.5, 0) -- (-3, 0);
  		\draw[->, ultra thick, gray] (-4.5, -1.5) -- (-4.5, 0) -- (-4.5, 1.5);
  		\draw[] (-6.25, 0) circle[radius = 0] node[scale = .9, left]{$(x; r)$};
  		\draw[] (-4.5, -1.75) circle[radius = 0] node[scale = .9, below]{$(y; s)$};
  		\draw[] (-4.5, -.75) circle[radius = 0] node[scale = .7, right]{$\textbf{A}$};
  		\draw[] (-5.25, 0) circle[radius = 0] node[scale = .7, above]{$\textbf{B}$};
  		\draw[] (-4.5, .75) circle[radius = 0] node[scale = .7, right]{$\textbf{C}$};
  		\draw[] (-3.75, 0) circle[radius = 0] node[scale = .7, above]{$\textbf{D}$};
  	}
  \end{align*} 
  
  \noindent sometimes omitting the labels $(\textbf{A}, \textbf{B}; \textbf{C}, \textbf{D})$, $(x; r)$, and $(y; s)$ from the diagram when convenient. We will also often omit the dependence of the quantization parameter $q$ from $W$, by writing $W_z (\textbf{A}, \textbf{B}; \textbf{C}, \textbf{D} \boldsymbol{\mid} r, s) = W_{x / y; q} (\textbf{A}, \textbf{B}; \textbf{C},\textbf{D} \boldsymbol{\mid} r, s)$\index{W@$W_z (\textbf{A}, \textbf{B}; \textbf{C}, \textbf{D} \boldsymbol{\mid} r, s)$; fused weight} (where we recall the spectral parameter $z = \frac{x}{y}$).
  
  Our first series of results concerns the derivation of an explicit family vertex weights that satisify the Yang--Baxter equation (\Cref{0wabcdproduct2} below); they are provided by the following definition. In what follows, for any vector $\textbf{X} = (X_1, X_2, \ldots , X_k) \in \mathbb{R}^k$ we set $|\textbf{X}| = \sum_{i = 1}^k X_i$; we recall the \emph{$q$-Pochhammer symbol} $(u; q)_k = \prod_{j = 1}^k (1 - q^{j - 1} u)$\index{U@$(u; q)_k$; $q$-Pochhammer symbol}; and we set 
  \begin{flalign}
  	\label{tufunction}
  	\varphi (\textbf{T}, \textbf{U}) = \displaystyle\sum_{1 \le i < j \le k} T_i U_j, \qquad \text{for $\textbf{T} = (T_1, T_2, \ldots , T_k) \in \mathbb{R}^k$ and $\textbf{U} = (U_1, U_2, \ldots , U_k) \in \mathbb{R}^k$}.
  \end{flalign}
	\index{0@$\varphi$}

  \begin{definition}[\Cref{wabcdrsxy} below]
  	
  	\label{0wabcdrsxy}
  	
  	If at least one of $\textbf{A}, \textbf{B}, \textbf{C}, \textbf{D}$ is not in $\{ 0, 1 \}^n$, then set $W_z (\textbf{A}, \textbf{B}; \textbf{C}, \textbf{D} \boldsymbol{\mid} r, s) = 0$. Otherwise define $\textbf{V} = (V_1, V_2, \ldots , V_n) \in \mathbb{Z}_{\ge 0}^n$ by setting $V_j = \min \{ A_j, B_j, C_j, D_j \}$ for each $j \in [1, n]$, and denote $|\textbf{X}| = x$ for each $X \in \{A, B, C, D, V \}$. Then set
  	\begin{flalign*}
  			W_z (\textbf{A}, \textbf{B}; \textbf{C}, \textbf{D} \boldsymbol{\mid} r, s) & = (-1)^v z^{d - b} r^{2c - 2a} s^{2d} q^{\varphi (\textbf{D} - \textbf{V}, \textbf{C}) + \varphi (\textbf{V}, \textbf{A}) - av + cv} \displaystyle\frac{(q^{1 - v} r^{-2} z; q)_v}{(q^{- v} s^2 r^{-2} z; q)_v} \displaystyle\frac{(r^2; q)_d}{(r^2; q)_b}  \\
  			& \quad \times \textbf{1}_{\textbf{A} + \textbf{B} = \textbf{C} + \textbf{D}} \displaystyle\sum_{p = 0}^{\min \{ b - v, c - v \}}  \displaystyle\frac{(q^{- v} s^2 r^{-2} z; q)_{c - p} (q^v r^2 z^{-1}; q)_p (z; q)_{b - p - v}}{(s^2 z; q)_{c + d - p - v}} \\
  			& \qquad \qquad \qquad \times (q^{-v} r^{-2} z)^p  \displaystyle\sum_{\textbf{P}} q^{\varphi (\textbf{B} - \textbf{D} - \textbf{P}, \textbf{P})},
   	\end{flalign*}
  	
  	\noindent where the last sum sum is over all $n$-tuples $\textbf{P} = (P_1, P_2, \ldots , P_n) \in \{ 0, 1 \}^n$ such that $|\textbf{P}| = p$ and $P_i \le \min \{ B_i - V_i, C_i - V_i \}$ for each $i \in [1, n]$. 
  	
  \end{definition}

	Let us mention three points concerning these vertex weights. The first is that they are nonzero only if $\textbf{A} + \textbf{B} = \textbf{C} + \textbf{D}$, a constraint we refer to by \emph{arrow conservation}. This essentially forces any colored arrow entering a vertex to also exit it, thereby enabling the interpretation of our vertex models as ensembles of colored paths. The second is that they are nonzero only if $\textbf{A}, \textbf{B}, \textbf{C}, \textbf{D} \in \{ 0, 1 \}^n$. We refer this to as \emph{fermionicity}, as it imposes an exclusion rule that prevents two distinct paths of the same color from passing along the same edge of the lattice. 
	
	The third concerns their origin. At $r = s = q^{-1 / 2}$ they arise as entries for the universal $R$-matrix of the quantum affine superalgebra $U_q \big( \widehat{\mathfrak{sl}} (1 | n) \big)$ under its defining (or \emph{fundamental}) representation, whose explicit form dates back to \cite{TSTESS}. In the more general case when $r = q^{-L / 2}$ and $s = q^{-M / 2}$ are nonnegative integer powers of $q^{-1 / 2}$, we obtain them through the \emph{fusion} procedure originating in the work \cite{ERT} of Kulish--Reshetikhin--Sklyanin. In particular, in this case our $W$ weights are entries of the above mentioned universal $R$-matrix, under the product of an $L$-fold and $M$-fold symmetric tensor power of the fundamental representation. These entries then happen to be rational functions in $r = q^{-L / 2}$ and $s = q^{-M / 2}$, enabling us to extend their definition to arbitrary $r, s \in \mathbb{C}$ by analytic continuation. For this reason, we refer to the weights $W_z (\textbf{A}, \textbf{B}; \textbf{C}, \textbf{D} \boldsymbol{\mid} r, s)$ from \Cref{0wabcdrsxy} as \emph{fused weights}. 
		
	We implement this procedure in detail in \Cref{Weights1}, \Cref{WeightsR}, \Cref{FusedW}, and \Cref{SymmetricBranching} below. However, our exposition there will be less representation theoretic than what was described above. As has been done in several recent works \cite{SHSVML,HSVMSRF,SVMST,MPI,SVMP} discussing fusion, we will instead proceed though an equivalent combinatorial (diagrammatic) framework to derive the fused weights. This fusion will in fact be applied on the quantum superalgebra $U_q \big( \widehat{\mathfrak{sl}} (m | n) \big)$ for arbitrary $m \ge 1$, with these more general fused weights given by \Cref{rxyml} below. We only specialize to the case $m = 1$ in \Cref{SymmetricBranching} when analytically continuing these weights in $L$ and $M$ (athough this can also be done for any $m \ge 1$).
	
	As a consequence of this framework, we establish the following Yang--Baxter equation for our fused weights. 

\begin{thm}[\Cref{wabcdproduct2} below]
	
	\label{0wabcdproduct2} 
	
	For any $x, y, z, r, s, t \in \mathbb{C}$ and $\textbf{\emph{I}}_1, \textbf{\emph{J}}_1, \textbf{\emph{K}}_1, \textbf{\emph{I}}_3, \textbf{\emph{J}}_3, \textbf{\emph{K}}_3 \in \{ 0, 1 \}^n$, we have 
	\begin{flalign*}
			& \displaystyle\sum_{\textbf{\emph{I}}_2, \textbf{\emph{J}}_2, \textbf{\emph{K}}_2} W_{x / y} ( \textbf{\emph{I}}_1, \textbf{\emph{J}}_1; \textbf{\emph{I}}_2, \textbf{\emph{J}}_2 \boldsymbol{\mid} r, s) W_{x / z} ( \textbf{\emph{K}}_1, \textbf{\emph{J}}_2; \textbf{\emph{K}}_2, \textbf{\emph{J}}_3 \boldsymbol{\mid} r, t) W_{y / z} ( \textbf{\emph{K}}_2, \textbf{\emph{I}}_2; \textbf{\emph{K}}_3, \textbf{\emph{I}}_3 \boldsymbol{\mid} s, t) \\
			& \quad = \displaystyle\sum_{\textbf{\emph{I}}_2, \textbf{\emph{J}}_2, \textbf{\emph{K}}_2} W_{y / z} ( \textbf{\emph{K}}_1, \textbf{\emph{I}}_1; \textbf{\emph{K}}_2, \textbf{\emph{I}}_2 \boldsymbol{\mid} s, t) W_{x / z} ( \textbf{\emph{K}}_2, \textbf{\emph{J}}_1; \textbf{\emph{K}}_3, \textbf{\emph{J}}_2 \boldsymbol{\mid} r, t) W_{x / y} ( \textbf{\emph{I}}_2, \textbf{\emph{J}}_2; \textbf{\emph{I}}_3, \textbf{\emph{J}}_3 \boldsymbol{\mid} r, s),   
	\end{flalign*}
	
	\noindent where both sums are over all $\textbf{\emph{I}}_2, \textbf{\emph{J}}_2, \textbf{\emph{K}}_2 \in \{ 0, 1 \}^n$. Diagrammatically, 
	
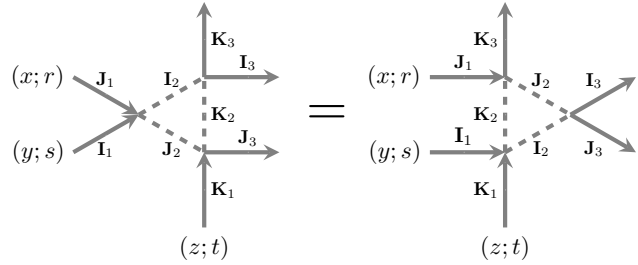
\begin{center}	
		\begin{tikzpicture}[
			>=stealth,
			auto,
			style={
				scale = 1
			}
			]			
			\draw[->, ultra thick, gray] (-.87, -.5) -- (0, 0);
			\draw[->, ultra thick, gray] (-.87, .5) -- (0, 0);
			\draw[->, ultra thick, gray] (.87, -1.5) -- (.87, -.5);
			\draw[] (-.87, .5) circle[radius = 0]  node[left, scale = .9]{$(x; r)$};
			\draw[] (-.87, -.5) circle[radius = 0]  node[left, scale = .9]{$(y; s)$};
			\draw[] (.87, -1.5) circle[radius = 0]  node[below, scale = .9]{$(z; t)$};
			\draw[gray, ultra thick, dashed] (0, 0) -- (.87, -.5); 
			\draw[gray, ultra thick, dashed] (0, 0) -- (.87, .5); 
			\draw[gray, ultra thick, dashed] (.87, -.5) -- (.87, .5); 
			\draw[->, ultra thick, gray] (.87, .5) -- (1.87, .5); 
			\draw[->, ultra thick, gray] (.87, -.5) -- (1.87, -.5); 
			\draw[->, ultra thick, gray] (.87, .5) -- (.87, 1.5);
			\draw[] (3.87, .5) circle[radius = 0]  node[left, scale = .9]{$(x; r)$};
			\draw[] (3.87, -.5) circle[radius = 0]  node[left, scale = .9]{$(y; s)$};
			\draw[] (4.87, -1.5) circle[radius = 0]  node[below, scale = .9]{$(z; t)$};
			\draw[->, ultra thick, gray] (4.87, -1.5) -- (4.87, -.5); 
			\draw[->, ultra thick, gray] (3.87, .5) -- (4.87, .5); 
			\draw[->, ultra thick, gray] (3.87, -.5) -- (4.87, -.5); 
			\draw[->, ultra thick, gray] (4.87, .5) -- (4.87, 1.5); 
			\draw[->, ultra thick, gray] (5.74, 0) -- (6.61, -.5); 
			\draw[->, ultra thick, gray] (5.74, 0) -- (6.61, .5); 
			\draw[gray, ultra thick, dashed] (4.87, -.5) -- (4.87, .5);
			\draw[-, gray, ultra thick, dashed] (4.87, -.5) -- (5.74, 0); 
			\draw[gray, ultra thick, dashed] (4.87, .5) -- (5.74, 0);
			\filldraw[fill=white, draw=black] (2.5, 0) circle [radius=0] node[scale = 2]{$=$};
			\filldraw[fill=white, draw=black] (-.44, -.275) circle [radius=0] node[below, scale = .7]{$\textbf{\emph{I}}_1$};
			\filldraw[fill=white, draw=black] (.44, .275) circle [radius=0] node[above, scale = .7]{$\textbf{\emph{I}}_2$};
			\filldraw[fill=white, draw=black] (1.45, .5) circle [radius=0] node[above, scale = .7]{$\textbf{\emph{I}}_3$};
			\filldraw[fill=white, draw=black] (-.44, .275) circle [radius=0] node[above, scale = .7]{$\textbf{\emph{J}}_1$};
			\filldraw[fill=white, draw=black] (.44, -.275) circle [radius=0] node[below, scale = .7]{$\textbf{\emph{J}}_2$};
			\filldraw[fill=white, draw=black] (1.45, -.5) circle [radius=0] node[above, scale = .7]{$\textbf{\emph{J}}_3$};
			\filldraw[fill=white, draw=black] (.87, -1) circle [radius=0] node[right, scale = .7]{$\textbf{\emph{K}}_1$};
			\filldraw[fill=white, draw=black] (.87, 0) circle [radius=0] node[right, scale = .7]{$\textbf{\emph{K}}_2$};
			\filldraw[fill=white, draw=black] (.87, 1) circle [radius=0] node[right, scale = .7]{$\textbf{\emph{K}}_3$};	
			\filldraw[fill=white, draw=black] (4.32, -.5) circle [radius=0] node[above, scale =.8]{$\textbf{\emph{I}}_1$};
			\filldraw[fill=white, draw=black] (5.35, -.26) circle [radius=0] node[below, scale = .7]{$\textbf{\emph{I}}_2$};
			\filldraw[fill=white, draw=black] (6.05, .26) circle [radius=0] node[above, scale = .7]{$\textbf{\emph{I}}_3$};
			\filldraw[fill=white, draw=black] (4.32, .5) circle [radius=0] node[above, scale = .7]{$\textbf{\emph{J}}_1$};
			\filldraw[fill=white, draw=black] (5.35, .26) circle [radius=0] node[above, scale = .7]{$\textbf{\emph{J}}_2$};
			\filldraw[fill=white, draw=black] (6.05, -.26) circle [radius=0] node[below, scale = .7]{$\textbf{\emph{J}}_3$};
			\filldraw[fill=white, draw=black] (4.87, -1) circle [radius=0] node[left, scale = .7]{$\textbf{\emph{K}}_1$};
			\filldraw[fill=white, draw=black] (4.87, 0) circle [radius=0] node[left, scale = .7]{$\textbf{\emph{K}}_2$};
			\filldraw[fill=white, draw=black] (4.87, 1) circle [radius=0] node[left, scale = .7]{$\textbf{\emph{K}}_3$};
		\end{tikzpicture}
	\end{center}
	
	\noindent where states along solid edges are fixed and those along dashed edges are summed over. 
	
	\end{thm} 
  
  	Although the fused weights given by \Cref{0wabcdrsxy} might appear a bit unpleasant, the fact that they are governed by four parameters $(q, z, r, s)$ makes them remarkably general. In particular, we show they factor completely under various specializations, including $n = 1$ (\Cref{wn1}), $r = q^{-1 / 2}$ (\Cref{rql1}), and a series of at least ten other degenerations (detailed in \Cref{DegenerationWeights} and \Cref{Limitrsz} below) that are depicted in \Cref{wlimit}. In that chart, for any $\textbf{X} = (X_1, X_2, \ldots , X_m) \in \mathbb{R}^m$ and any $1 \le j \le k \le m$ we write $X_{[j, k]} = \sum_{i = j}^k X_i$.\index{X@$X_{[j, k]}$} We also write $\textbf{X} \ge \textbf{Y}$\index{X@$\textbf{X} \ge \textbf{Y}$} for any $\textbf{Y} = (Y_1, Y_2, \ldots , Y_m) \in \mathbb{R}^m$ if $X_j \ge Y_j$ for each $j$. 
  	
  	Let us mention that directly establishing the Yang--Baxter equation for most of these degenerated weights (without realizing them as special cases of our general fused $W$ ones) does not appear to be an immediate task. For example, this was done in Appendix A of \cite{VMFP} for the specializations $x^d q^{\varphi (\textbf{D}, \textbf{C} + \textbf{D})} \textbf{1}_{v = 0}$ depicted on the bottom-right of \Cref{wlimit}, through an elaborate series of combinatorial reductions.

		\begin{figure}[t]
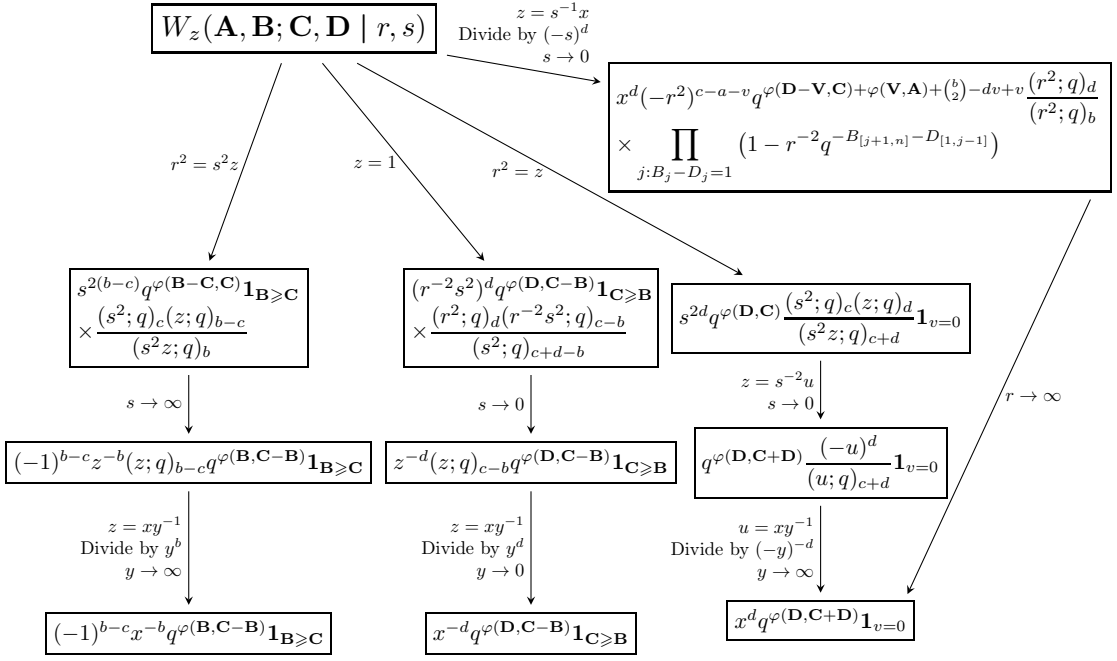

			\begin{center}
  	\tikz{.7}{	
  		\node (w) at (-6, 2) [scale = 1.1]{\framebox{$W_z (\textbf{A}, \textbf{B}; \textbf{C}, \textbf{D} \boldsymbol{\mid} r, s)$}};
  		
  		\node (wrs) at (-8, -3.5) [scale = .85, align = center] {\framebox{\parbox{3.5cm}{$s^{2 (b - c)} q^{\varphi (\textbf{B} - \textbf{C}, \textbf{C})} \textbf{1}_{\textbf{B} \ge \textbf{C}}$ \\ $\times \displaystyle\frac{(s^2; q)_c (z; q)_{b - c}}{(s^2 z; q)_b}$}}};
  		\node (wz) at (-1.5, -3.5) [scale = .85, align = center] {\framebox{\parbox{3.75cm}{$(r^{-2} s^2)^d q^{\varphi (\textbf{D}, \textbf{C} - \textbf{B})} \textbf{1}_{\textbf{C} \ge \textbf{B}}$ \\ $\times \displaystyle\frac{(r^2; q)_d (r^{-2} s^2; q)_{c - b}}{(s^2; q)_{c + d - b}}$}}};
  		\node (wrz) at (4, -3.5) [scale = .85] {\framebox{$s^{2d} q^{\varphi (\textbf{D}, \textbf{C})} \displaystyle\frac{(s^2; q)_c (z; q)_d}{(s^2 z; q)_{c + d}} \textbf{1}_{v = 0}$}};
  		
  		\node (wszx) at (4.75, .125) [scale = .85, align = center] {\framebox{\parbox{7.625cm}{$x^d (-r^2)^{c - a - v} q^{\varphi (\textbf{D} - \textbf{V}, \textbf{C}) + \varphi (\textbf{V}, \textbf{A}) + \binom{b}{2} - dv + v} \displaystyle\frac{(r^2; q)_d}{(r^2; q)_b}$ \\ $\times \displaystyle\prod_{j: B_j - D_j = 1} \big( 1 - r^{-2} q^{-B_{[j + 1, n]} - D_{[1, j - 1]}} \big)$}}};
  		
  		\node (wrss) at (-8, -6.25) [scale = .85] {\framebox{$(-1)^{b - c} z^{-b} (z; q)_{b - c} q^{\varphi (\textbf{B}, \textbf{C} - \textbf{B})} \textbf{1}_{\textbf{B} \ge \textbf{C}}$}};
  		\node (wrsz) at (-8, -9.5) [scale = .85] {\framebox{$(-1)^{b - c} x^{-b} q^{\varphi (\textbf{B}, \textbf{C} - \textbf{B})} \textbf{1}_{\textbf{B} \ge \textbf{C}}$}};
  		
  		\node (wzs) at (-1.5, -6.25) [scale = .85] {\framebox{$z^{-d} (z; q)_{c - b} q^{\varphi (\textbf{D}, \textbf{C} - \textbf{B})} \textbf{1}_{\textbf{C} \ge \textbf{B}}$}};
  		\node (wzz) at (-1.5, -9.5) [scale = .85] {\framebox{$x^{-d} q^{\varphi (\textbf{D}, \textbf{C} - \textbf{B})}\textbf{1}_{\textbf{C} \ge \textbf{B}}$}};
  		
  		\node (wrzs) at (4, -6.25) [scale = .85] {\framebox{$q^{\varphi (\textbf{D}, \textbf{C} + \textbf{D})} \displaystyle\frac{(-u)^d}{(u; q)_{c + d}} \textbf{1}_{v = 0}$}};
  		\node (wrzsz) at (4, -9.25) [scale = .85] {\framebox{$x^d q^{\varphi (\textbf{D}, \textbf{C} + \textbf{D})} \textbf{1}_{v = 0}$}};
  	
  		\draw[->] (w) --  node[left, scale = .7]{$r^2 = s^2 z$} (wrs);
  		\draw[->] (w) --  node[left, scale = .7]{$z = 1$} (wz);
  		\draw[->] (w) --  node[left, scale = .7]{$r^2 = z$} (wrz);
  		\draw[->] (w) --  node[above, scale = .7, align = right]{$z = s^{-1} x$ \\ Divide by $(-s)^d$ \\ $s \rightarrow 0$} (wszx);
  		
  		\draw[->] (wrs) -- node[left, scale = .7]{$s \rightarrow \infty$} (wrss);
  		\draw[->] (wrss) -- node[left, scale = .7, align = right]{$z = x y^{-1}$ \\ Divide by $y^b$ \\ $y \rightarrow \infty$}(wrsz);
  		
  		\draw[->] (wz) -- node[left, scale = .7]{$s \rightarrow 0$} (wzs);
  		\draw[->] (wzs) -- node[left, scale = .7, align = right]{$z = xy^{-1}$ \\ Divide by $y^d$ \\ $y \rightarrow 0$} (wzz);
  		
  		\draw[->] (wrz) -- node[left, scale = .7, align = right]{$z = s^{-2} u$ \\ $s \rightarrow 0$} (wrzs);
  		\draw[->] (wrzs) -- node[left, scale = .7, align = right]{$u = x y^{-1}$ \\ Divide by $(-y)^{-d}$ \\ $y \rightarrow \infty$} (wrzsz);
  		
  		\draw[->] (9.125, -1.25) -- node [right, scale = .7]{$r \rightarrow \infty$} (5.65, -8.75);
  	} 
  \caption{\label{wlimit} Depicted above are some fully factored cases of the fused weights $W_z (\textbf{A}, \textbf{B}; \textbf{C}, \textbf{D} \boldsymbol{\mid} r, s)$.}
  \end{center}
\end{figure}

\section{Symmetric Functions}
  
  \label{0Functions}
  
  The symmetric functions we consider in this text will be obtained from vertex models under the fused weights from \Cref{0wabcdrsxy}, on a semi-infinite strip with certain boundary data. To explain this further, we require some additional notation. In what follows, we fix $n \ge 1$ and denote the $n$-tuples $\textbf{e}_0 = (0, 0, \ldots , 0) \in \{ 0, 1 \}^n$\index{E@$\textbf{e}_0$} and $\textbf{e}_{[1, n]} = (1, 1, \ldots , 1) \in \{ 0, 1 \}^n$\index{E@$\textbf{e}_{[1, n]}$}. We also define the normalization $\widehat{W}_z$ of the fused weights $W_z$ from \Cref{0wabcdrsxy} by 
  \begin{flalign}
  	\label{wzwz}  
  \widehat{W}_z (\textbf{A}, \textbf{B}; \textbf{C}, \textbf{D} \boldsymbol{\mid} r, s) = \displaystyle\frac{(s^2 z; q)_n}{s^{2n} (z; q)_n} W_z (\textbf{A}, \textbf{B}; \textbf{C}, \textbf{D} \boldsymbol{\mid} r, s).
  \end{flalign}
	\index{W@$\widehat{W}_z (\textbf{A}, \textbf{B}; \textbf{C}, \textbf{D} \boldsymbol{\mid} r, s)$; normalized fused weight}

	\noindent Then the $W_z$ and $\widehat{W}_z$ weights always satisfy (see \Cref{wabcd01n} and \eqref{wz1} below)
  \begin{flalign}
  	\label{wz1wz1} 
  W_z (\textbf{e}_0, \textbf{e}_0; \textbf{e}_0, \textbf{e}_0 \boldsymbol{\mid} r, s) = 1; \qquad \widehat{W}_z \big( \textbf{e}_0, \textbf{e}_{[1, n]}; \textbf{e}_0, \textbf{e}_{[1, n]} \boldsymbol{\mid} r, s \big) = 1.
  \end{flalign}

  A \emph{signature} $\lambda = (\lambda_1, \lambda_2, \ldots , \lambda_{\ell}) \in \mathbb{Z}_{\ge 0}^{\ell}$ is a finite (possibly empty) non-increasing sequence of nonnegative integers. Its \emph{size} is $|\lambda| = \sum_{j = 1}^{\ell} \lambda_j$\index{0@$\lambda, \mu$; typical signatures or partitions!$\mid$$\lambda / \mu$$\mid$; size of $\lambda / \mu$!$\mid$$\lambda$$\mid$; size of $\lambda$}\index{0@$\lambda, \mu$; typical signatures or partitions} and its \emph{length} is $\ell = \ell (\lambda)$\index{0@$\lambda, \mu$; typical signatures or partitions!$\ell (\lambda)$; length of $\lambda$}. Let $\Sign_{\ell}$ denote the set of all signatures of length $\ell$; let $0^{\ell} \in \Sign_{\ell}$ denote the signature of length $\ell$\index{0@$\lambda, \mu$; typical signatures or partitions!$0^{\ell}$; zero signature} whose entries are all equal to $0$; and let $\varnothing \in \Sign_0$ denote the empty signature\index{0@$\lambda, \mu$; typical signatures or partitions!$\varnothing$; empty signature}, that is, the unique one of length $0$. Further define $\Sign = \bigcup_{\ell = 0}^{\infty} \Sign_{\ell}$\index{S@$\Sign$; set of all signatures}\index{S@$\Sign$; set of all signatures!$\Sign_{\ell}$; set of length $\ell$ signatures}, and for any $\lambda \in \Sign_{\ell}$ set 
  \begin{flalign}
  	\label{t}
  	\mathfrak{T} (\lambda) = (\lambda_1 + \ell, \lambda_2 + \ell - 1, \ldots , \lambda_{\ell} + 1) \in \mathbb{Z}_{> 0}^{\ell},
  \end{flalign}
	\index{T@$\mathfrak{T}$}
  
  \noindent whose entries are all distinct since $\lambda$ is non-increasing. For example, $\mathfrak{T} (0^M) = (M, M - 1, \ldots , 1)$.
  
  Let $\SeqSign_n$\index{S@$\SeqSign_n$; set of sequences of $n$ signatures} denote the set of sequences of $n$ signatures $\boldsymbol{\lambda} = \big( \lambda^{(1)}, \lambda^{(2)}, \ldots , \lambda^{(n)} \big)$\index{0@$\boldsymbol{\lambda}, \boldsymbol{\mu}$; typical signature sequences} and, for any $M \ge 0$ let $\SeqSign_{n; M} \subset \SeqSign_n$\index{S@$\SeqSign_n$; set of sequences of $n$ signatures!$\SeqSign_{n; M}$; set of sequences of $n$ signatures of length $M$} denote the set of those sequences such that $\ell \big( \lambda^{(j)} \big) = M$ for each $j \in [1, n]$. Further define the size of the sequence $\boldsymbol{\lambda}$ by $|\boldsymbol{\lambda}| = \sum_{j = 1}^n \big| \lambda^{(j)} \big|$.\index{0@$\boldsymbol{\lambda}, \boldsymbol{\mu}$; typical signature sequences!$\mid$$\boldsymbol{\lambda} / \boldsymbol{\mu}$$\mid$; size of $\boldsymbol{\lambda} / \boldsymbol{\mu}$!$\mid$$\boldsymbol{\lambda}$$\mid$; size of $\boldsymbol{\lambda}$} For example, let $\boldsymbol{0}^M = \boldsymbol{0}^{(M; n)} \in \SeqSign_{n; M}$ denote the (size $0$) sequence whose $n$ signatures are all equal to $0^M \in \Sign_M$,\index{0@$\boldsymbol{\lambda}, \boldsymbol{\mu}$; typical signature sequences!$\boldsymbol{0}^M$; sequence of zero signatures} and let $\boldsymbol{\varnothing} = \boldsymbol{\varnothing}_n \in \SeqSign_{n; 0}$ denote the (size $0$) sequence of $n$ empty signatures.\index{0@$\boldsymbol{\lambda}, \boldsymbol{\mu}$; typical signature sequences!$\boldsymbol{\varnothing}$; sequence of empty signatures}
  
  For any $\boldsymbol{\lambda} \in \SeqSign_n$, we introduce the following infinite sequence $\mathscr{S} (\boldsymbol{\lambda}) = \big( \textbf{S}_1 (\boldsymbol{\lambda}), \textbf{S}_2 (\boldsymbol{\lambda}), \ldots \big)$\index{S@$\mathscr{S} (\boldsymbol{\lambda})$} of elements in $\{ 0, 1 \}^n$. For each $j \ge 1$, define $\textbf{S}_j = \textbf{S}_j (\boldsymbol{\lambda}) = (S_{1, j}, S_{2, j}, \ldots , S_{n, j}) \in \{ 0, 1 \}^n$ by setting $S_{i, j} = \textbf{1}_{j \in \mathfrak{T} (\lambda^{(i)})}$ for every $i \in [1, n]$. For example, $\mathscr{S} (\boldsymbol{\varnothing}) = (\textbf{e}_0, \textbf{e}_0, \ldots )$; moreover, the first $M$ entries of $\mathscr{S} (\boldsymbol{0}^M)$ are all $\textbf{e}_{[1, n]}$, and the remaining ones are all $\textbf{e}_0$. In general, we will use these sequences $\mathscr{S} (\boldsymbol{\lambda})$ to index boundary data for the vertex models associated with our symmetric functions. 
  
  Now fix an integer $N \ge 1$; finite sequences of complex numbers $\textbf{x} = (x_1, x_2, \ldots , x_N)$ and $\textbf{r} = (r_1, r_2, \ldots , r_N)$; and infinite sequences of complex numbers $\textbf{y} = (y_1, y_2, \ldots )$ and $\textbf{s} = (s_1, s_2, \ldots )$. We will define functions $G_{\boldsymbol{\lambda} / \boldsymbol{\mu}} (\textbf{x}; \textbf{r} \boldsymbol{\mid} \textbf{y}; \textbf{s})$ and $F_{\boldsymbol{\lambda} / \boldsymbol{\mu}} (\textbf{x}; \textbf{r} \boldsymbol{\mid} \textbf{y}; \textbf{s})$ diagrammatically as \emph{partition functions} for certain vertex models, that is, the sum of the total weights of all path ensembles with given boundary data; see \Cref{fgdefinition} below for an equivalent algebraic definition. In what follows, we fix an integer $M \ge 0$. 
  
  For any signature sequences $\boldsymbol{\lambda}, \boldsymbol{\mu} \in \SeqSign_{n; M}$, define $G_{\boldsymbol{\lambda} / \boldsymbol{\mu}} (\textbf{x}; \textbf{r} \boldsymbol{\mid} \textbf{y}; \textbf{s})$\index{G@$G_{\boldsymbol{\lambda} / \boldsymbol{\mu}} (\textbf{x}; \textbf{r} \boldsymbol{\mid} \textbf{y}; \textbf{s})$} as the partition function for the vertex model
  	\begin{align}
  	\label{gvertex} 
  		\tikz{.75}{
  			\draw[->, thick, red] (.95, 0) -- (.95, 1);
  			\draw[->, thick, red] (2.95, 3) -- (2.95, 4); 
  			\draw[->, thick, red] (2.95, 0) -- (2.95, 1);
  			\draw[->, thick, red] (5.95, 3) -- (5.95,4);
  			\draw[->, thick, blue] (1.95, 0) -- (1.95, 1);
  			\draw[->, thick, blue] (5, 3) -- (5, 4);
  			\draw[->, thick, blue] (4, 0) -- (4, 1);
  			\draw[->, thick, blue] (6.05, 3) -- (6.05, 4);
  			\draw[->, thick, green] (1.05, 0) -- (1.05, 1); 
  			\draw[->, thick, green] (3.05, 3) -- (3.05, 4);
  			\draw[->, thick, green] (2.05, 0) -- (2.05, 1);
  			\draw[->, thick, green] (7, 3) -- (7, 4);
  			\draw[ultra thick, gray, dashed] (1, 1) -- (1, 3);
  			\draw[ultra thick, gray, dashed] (2, 1) -- (2, 3);
  			\draw[ultra thick, gray, dashed] (3, 1) -- (3, 3);
  			\draw[ultra thick, gray, dashed] (4, 1) -- (4, 3);
  			\draw[ultra thick, gray, dashed] (5, 1) -- (5, 3);
  			\draw[ultra thick, gray, dashed] (6, 1) -- (6, 3);
  			\draw[ultra thick, gray, dashed] (7, 1) -- (7, 3);
  			\draw[ultra thick, gray, dashed] (1, 1)node[black, left = 4, scale = .65]{$\textbf{e}_0$} -- (7, 1) node[black, right, scale = .85]{$\cdots$};
  			\draw[ultra thick, gray, dashed] (1, 2)node[black, left = 4, scale = .65]{$\textbf{e}_0$} -- (7, 2) node[black, right, scale = .85]{$\cdots$};
  			\draw[ultra thick, gray, dashed] (1, 3)node[black, left = 4, scale = .65]{$\textbf{e}_0$} -- (7, 3) node[black, right, scale = .85]{$\cdots$};
  			\draw[->, very thick] (0, 0) -- (0, 4.5);
  			\draw[->, very thick] (0, 0) -- (8, 0);
  			\draw[]  (4, 3.15) circle [radius = 0] node[above, scale = .7]{$(y_4; s_4)$};
  			\draw[]  (8.875, 2) circle [radius = 0] node[scale = .7]{$(x_2; r_2)$};
  			\draw[] (1, -.15) circle [radius = 0] node[below, scale = .65]{$\textbf{S}_1 (\boldsymbol{\mu})$};
  			\draw[] (2, -.15) circle [radius = 0] node[below, scale = .65]{$\textbf{S}_2 (\boldsymbol{\mu})$};
  			\draw[] (3, -.15) circle [radius = 0] node[below, scale = .65]{$\textbf{S}_3 (\boldsymbol{\mu})$};
  			\draw[] (4, -.15) circle [radius = 0] node[below, scale = .85]{$\cdots$};
  			\draw[] (1, 4.15) circle [radius = 0] node[above, scale = .65]{$\textbf{S}_1 (\boldsymbol{\lambda})$};
  			\draw[] (2, 4.15) circle [radius = 0] node[above, scale = .65]{$\textbf{S}_2 (\boldsymbol{\lambda})$};
  			\draw[] (3, 4.15) circle [radius = 0] node[above, scale = .65]{$\textbf{S}_3 (\boldsymbol{\lambda})$};
  			\draw[] (4, 4.15) circle [radius = 0] node[above, scale = .85]{$\cdots$};
  			\draw[]  (8, 1) circle [radius = 0] node[scale = .65]{$\textbf{e}_0$};
  			\draw[]  (8, 2) circle [radius = 0] node[scale = .65]{$\textbf{e}_0$};
  			\draw[]  (8, 3) circle [radius = 0] node[scale = .65]{$\textbf{e}_0$};
  			\draw[] (-3, 1.5) circle[radius = 0] node[above, scale = 1.1]{$G_{\boldsymbol{\lambda} / \boldsymbol{\mu}} (\textbf{x}; \textbf{r} \boldsymbol{\mid} \textbf{y}; \textbf{s}) = $};
  			\draw[] (9.75, 1) -- (10, 1) -- (10, 3) -- (9.75, 3);
  			\draw[] (10, 2) circle [radius = 0] node[right]{$N$};
  		}
  	\end{align}
  
  \noindent where the vertex weights are given by the $W$ of \Cref{0wabcdrsxy}. The model here resides on a semi-infinite strip consisting of $N$ rows. The rapidity pair in the $i$-th column is given by $(y_i; s_i)$, and that in the $j$-th row by $(x_j; r_j)$; the quantization parameter is $q$. No paths enter or exit through the left or right boundaries of the strip, but they enter along the bottom boundary as indexed by $\mathscr{S} (\boldsymbol{\mu})$ and exit along the top boundary as indexed by $\mathscr{S} (\boldsymbol{\lambda})$. All but finitely many vertices in this model have arrow configuration $(\textbf{e}_0, \textbf{e}_0; \textbf{e}_0, \textbf{e}_0)$, which has weight $1$ under $W$ by \eqref{wz1wz1}; therefore, the partition function for this model is well-defined. We abbreviate $G_{\boldsymbol{\lambda}} (\textbf{x}; \textbf{r} \boldsymbol{\mid} \textbf{y}; \textbf{s}) = G_{\boldsymbol{\lambda} / \boldsymbol{0}^M} (\textbf{x}; \textbf{r} \boldsymbol{\mid} \textbf{y}; \textbf{s})$.\index{G@$G_{\boldsymbol{\lambda} / \boldsymbol{\mu}} (\textbf{x}; \textbf{r} \boldsymbol{\mid} \textbf{y}; \textbf{s})$!$G_{\boldsymbol{\lambda}} (\textbf{x}; \textbf{r} \boldsymbol{\mid} \textbf{y}; \textbf{s})$}
 
  Next, for any signature sequences $\boldsymbol{\lambda} \in \SeqSign_{n; M + N}$ and $\boldsymbol{\mu} \in \SeqSign_{n; M}$, define $F_{\boldsymbol{\lambda} / \boldsymbol{\mu}} (\textbf{x}; \textbf{r} \boldsymbol{\mid} \textbf{y}; \textbf{s})$\index{F@$F_{\boldsymbol{\lambda} / \boldsymbol{\mu}} (\textbf{x}; \textbf{r} \boldsymbol{\mid} \textbf{y}; \textbf{s})$} as the partition function for the vertex model 
  	\begin{align}
  	\label{fvertex}
  		\tikz{.75}{
  			\draw[] (1, -6.15) circle [radius = 0] node[below, scale = .65]{$\textbf{S}_1 (\boldsymbol{\lambda})$};
  			\draw[] (2, -6.15) circle [radius = 0] node[below, scale = .65]{$\textbf{S}_2 (\boldsymbol{\lambda})$};
  			\draw[] (3, -6.15) circle [radius = 0] node[below, scale = .65]{$\textbf{S}_3 (\boldsymbol{\lambda})$};
  			\draw[] (4, -6.15) circle [radius = 0] node[below, scale = .85]{$\cdots$};
  			\draw[] (1, -1.85) circle [radius = 0] node[above, scale = .65]{$\textbf{S}_1 (\boldsymbol{\mu})$};
  			\draw[] (2, -1.85) circle [radius = 0] node[above, scale = .65]{$\textbf{S}_2 (\boldsymbol{\mu})$};
  			\draw[] (3, -1.85) circle [radius = 0] node[above, scale = .65]{$\textbf{S}_3 (\boldsymbol{\mu})$};
  			\draw[] (4, -1.85) circle [radius = 0] node[above, scale = .85]{$\cdots$}; 
  			\draw[->, thick, red] (6, -5.1) -- (7, -5.1);
  			\draw[->, thick, blue] (6, -5) -- (7, -5);
  			\draw[->, thick, green] (6, -4.9) -- (7, -4.9);
  			\draw[->, thick, red] (6, -4.1) -- (7, -4.1);
  			\draw[->, thick, blue] (6, -4) -- (7, -4);
  			\draw[->, thick, green] (6, -3.9) -- (7, -3.9);
  			\draw[->, thick, red] (6, -3.1) -- (7, -3.1);
  			\draw[->, thick, blue] (6, -3) -- (7, -3);
  			\draw[->, thick, green] (6, -2.9) -- (7, -2.9);	
  			\draw[->, thick, red] (7, -5.1) -- (8, -5.1);
  			\draw[->, thick, blue] (7, -5) -- (8, -5) node[right, black, scale = .85]{$\cdots$};
  			\draw[->, thick, green] (7, -4.9) -- (8, -4.9);
  			\draw[->, thick, red] (7, -4.1) -- (8, -4.1);
  			\draw[->, thick, blue] (7, -4) -- (8, -4) node[right, black, scale = .85]{$\cdots$};
  			\draw[->, thick, green] (7, -3.9) -- (8, -3.9);
  			\draw[->, thick, red] (7, -3.1) -- (8, -3.1);
  			\draw[->, thick, blue] (7, -3) -- (8, -3) node[right, black, scale = .85]{$\cdots$};
  			\draw[->, thick, green] (7, -2.9) -- (8, -2.9);		
  			\draw[->, thick, red] (.95, -6) -- (.95, -5);
  			\draw[->, thick, red] (5.95, -6) -- (5.95, -5);
  			\draw[->, thick, red] (4.95, -6) -- (4.95, -5);
  			\draw[->, thick, red] (2.95, -3) -- (2.95, -2); 
  			\draw[->, thick, red] (2.9, -6) -- (2.9, -5);
  			\draw[->, thick, blue] (1.95, -6) -- (1.95, -5);
  			\draw[->, thick, blue] (5, -3) -- (5, -2);
  			\draw[->, thick, blue] (3, -6) -- (3, -5);
  			\draw[->, thick, blue] (3.95, -6) -- (3.95, -5);
  			\draw[->, thick, blue] (6.05, -6) -- (6.05, -5);
  			\draw[->, thick, green] (1.05, -6) -- (1.05, -5); 
  			\draw[->, thick, green] (3.05, -3) -- (3.05, -2);
  			\draw[->, thick, green] (2.05, -6) -- (2.05, -5);
  			\draw[->, thick, green] (3.1, -6) -- (3.1, -5);
  			\draw[->, thick, green] (5.05, -6) -- (5.05, -5);
  			\draw[ultra thick, gray, dashed] (1, -5) -- (1, -3);
  			\draw[ultra thick, gray, dashed] (2, -5) -- (2, -3);
  			\draw[ultra thick, gray, dashed] (3, -5) -- (3, -3);
  			\draw[ultra thick, gray, dashed] (4, -5) -- (4, -3);
  			\draw[ultra thick, gray, dashed] (5, -5) -- (5, -3);
  			\draw[ultra thick, gray, dashed] (6, -5) -- (6, -3);
  			\draw[ultra thick, gray, dashed] (1, -5) node[black, left= 2, scale = .65]{$\textbf{e}_0$} -- (6, -5);
  			\draw[ultra thick, gray, dashed] (1, -4) node[black, left= 2, scale = .65]{$\textbf{e}_0$} -- (6, -4);
  			\draw[ultra thick, gray, dashed] (1, -3) node[black, left= 2, scale = .65]{$\textbf{e}_0$} -- (6, -3);		
  			\draw[]  (4, -2.85) circle [radius = 0] node[above, scale = .7]{$(y_4; s_4)$};
  			\draw[]  (9.125, -3) circle [radius = 0] node[scale = .65]{$\textbf{e}_{[1, n]}$};
  			\draw[]  (9.125, -4) circle [radius = 0] node[scale = .65]{$\textbf{e}_{[1, n]}$};
  			\draw[]  (9.125, -5) circle [radius = 0] node[scale = .65]{$\textbf{e}_{[1, n]}$};
  			\draw[]  (10.25, -4) circle [radius = 0] node[scale = .7]{$(x_2; r_2)$};
  			\draw[] (-3, -4.5) circle[radius = 0] node[above, scale = 1.1]{$F_{\boldsymbol{\lambda} / \boldsymbol{\mu}} (\textbf{x}; \textbf{r} \boldsymbol{\mid} \textbf{y}; \textbf{s}) = $};
  			\draw[->, very thick] (0, -6) -- (0, -1.5);
  			\draw[->, very thick] (0, -6) -- (9, -6);
  			\draw[] (11.25, -3) -- (11.5, -3) -- (11.5, -5) -- (11.25, -5);
  			\draw[] (11.5, -4) circle [radius = 0] node[right]{$N$};
	}
  	\end{align}
  	
  	\noindent whose vertex weights are given by the $\widehat{W}$ of \eqref{wzwz}. The parameters and boundary data for the vertex model here are the same as that for $G_{\boldsymbol{\lambda} / \boldsymbol{\mu}}$ above, except for two differences. First, the boundary data along the top and bottom boundaries are interchanged, that is, arrows enter along the bottom boundary as indexed by $\mathscr{S} (\boldsymbol{\lambda})$ and they exit along the top one as indexed by $\mathscr{S} (\boldsymbol{\mu})$. The second is that all horizontal edges in the strip sufficiently far to the right\footnote{In various previous works \cite{SRF,HSVMSRF,SP,SVMST}, analogous functions were defined through vertex models in which all paths enter at the left boundary (instead of exit at the right one). We will also consider functions, denoted by $H_{\boldsymbol{\lambda} / \boldsymbol{\mu}}$ in \Cref{fgdefinition} below, defined by such boundary data in this text, but they will play less of a prominent role. Our reason for focusing on the $F_{\boldsymbol{\lambda} / \boldsymbol{\mu}}$ functions here is that they seem to pair better with the $G_{\boldsymbol{\lambda} / \boldsymbol{\mu}}$ ones for proving Cauchy identities. Still, under certain limits, we show that these $F$ and $H$ functions essentially coincide; see \Cref{hf} below.} are occupied by paths of all $n$ colors. In this way, all but finitely many vertices in this model have arrow configuration $\big( \textbf{e}_0, \textbf{e}_{[1, n]}; \textbf{e}_0, \textbf{e}_{[1, n]} \big)$, which has weight $1$ under $\widehat{W}$ by \eqref{wz1wz1}; therefore, the partition function for this model is also well-defined. We abbreviate $F_{\boldsymbol{\lambda}} (\textbf{x}; \textbf{r} \boldsymbol{\mid} \textbf{y}; \textbf{s}) = F_{\boldsymbol{\lambda} / \boldsymbol{\varnothing}} (\textbf{x}; \textbf{r} \boldsymbol{\mid} \textbf{y}; \textbf{s})$.\index{F@$F_{\boldsymbol{\lambda} / \boldsymbol{\mu}} (\textbf{x}; \textbf{r} \boldsymbol{\mid} \textbf{y}; \textbf{s})$!$F_{\boldsymbol{\lambda}} (\textbf{x}; \textbf{r} \boldsymbol{\mid} \textbf{y}; \textbf{s})$}

	Let us begin with an example for one of these functions. When $(\boldsymbol{\lambda}, \boldsymbol{\mu}) = (\textbf{0}^N, \boldsymbol{\varnothing})$, $F_{\boldsymbol{\lambda} / \boldsymbol{\mu}} (\textbf{x}; \textbf{r} \boldsymbol{\mid} \textbf{y}; \textbf{s})$ becomes a partition function with \emph{domain-wall boundary data}, given by 
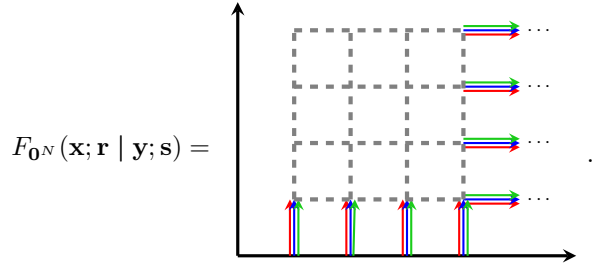
\begin{align*}
	\tikz{.75}{	 			
		\draw[->, very thick] (0, -7) -- (0, -2.5);
		\draw[->, very thick] (0, -7) -- (6, -7);
		\draw[->, thick, red] (4, -6.075) -- (5, -6.075);
		\draw[->, thick, red] (4, -5.075) -- (5, -5.075);
		\draw[->, thick, red] (4, -4.075) -- (5, -4.075);
		\draw[->, thick, red] (4, -3.075) -- (5, -3.075);
		\draw[->, thick, blue] (4, -6) -- (5, -6) node[right, black, scale = .75]{$\cdots$};
		\draw[->, thick, blue] (4, -5) -- (5, -5) node[right, black, scale = .75]{$\cdots$};
		\draw[->, thick, blue] (4, -4) -- (5, -4) node[right, black, scale = .75]{$\cdots$};
		\draw[->, thick, blue] (4, -3) -- (5, -3) node[right, black, scale = .75]{$\cdots$};
		\draw[->, thick, green] (4, -5.925) -- (5, -5.925);
		\draw[->, thick, green] (4, -4.925) -- (5, -4.925);
		\draw[->, thick, green] (4, -3.925) -- (5, -3.925);
		\draw[->, thick, green] (4, -2.925) -- (5, -2.925);
		\draw[->, thick, red] (.925, -7) -- (.925, -6);
		\draw[->, thick, red] (1.925, -7) -- (1.925, -6);
		\draw[->, thick, red] (2.925, -7) -- (2.925, -6);
		\draw[->, thick, red] (3.925, -7) -- (3.925, -6);
		\draw[->, thick, blue] (1, -7) -- (1, -6);
		\draw[->, thick, blue] (2, -7) -- (2, -6);
		\draw[->, thick, blue] (3, -7) -- (3, -6);
		\draw[->, thick, blue] (4, -7) -- (4, -6);
		\draw[->, thick, green] (1.075, -7) -- (1.075, -6);
		\draw[->, thick, green] (2.05, -7) -- (2.075, -6);
		\draw[->, thick, green] (3.075, -7) -- (3.075, -6);
		\draw[->, thick, green] (4.075, -7) -- (4.075, -6);
		\draw[ultra thick, gray, dashed] (1, -6) -- (1, -3);
		\draw[ultra thick, gray, dashed] (2, -6) -- (2, -3);
		\draw[ultra thick, gray, dashed] (3, -6) -- (3, -3);
		\draw[ultra thick, gray, dashed] (4, -6) -- (4, -3);
		\draw[ultra thick, gray, dashed] (1, -3) -- (4, -3);
		\draw[ultra thick, gray, dashed] (1, -6) -- (4, -6);
		\draw[ultra thick, gray, dashed] (1, -5) -- (4, -5);
		\draw[ultra thick, gray, dashed] (1, -4) -- (4, -4);
		\draw[] (-2.25, -5.5) circle[radius = 0] node[above, scale = 1]{$F_{\boldsymbol{0}^N} (\textbf{x}; \textbf{r} \boldsymbol{\mid} \textbf{y}; \textbf{s}) = $};
		\draw[] (6.25, -5.5) circle[radius = 0] node[above, scale = 1]{.};
	}
\end{align*}

Domain-wall partition functions of this general qualitative type have been studied extensively in the mathematical physics literature since the works of Korepin \cite{NWF}, Izergin \cite{PFSVMFV}, and Izergin--Coker--Korepin \cite{DSV}. Such analyses have gained popularity over the past several decades, especially due to their connections with algebraic combinatorics, initially observed by Kuperberg \cite{PPASM} in the context of alternating sign matrices. In many cases, domain-wall partition functions are given explicitly by a determinant. We show the specific domain-wall partition function $F_{\boldsymbol{0}^N} (\textbf{x}; \textbf{r} \boldsymbol{\mid} \textbf{y}; \textbf{s})$ depicted above in fact admits the following fully factored form, which will be useful for establishing the Cauchy identity between $F$ and $G$, given by \Cref{0fgsum2} below. 

\begin{prop}[\Cref{f0} below]
	
	\label{0f0} 
	
	We have
	\begin{flalign}
		\label{0f0n}
		\begin{aligned}
		F_{\boldsymbol{0}^N} (\textbf{\emph{x}}; \textbf{\emph{r}} \boldsymbol{\mid} \textbf{\emph{y}}; \textbf{\emph{s}}) & = \displaystyle\prod_{j = 1}^n s_j^{2n (j - N)} r_j^{2n (j - N - 1)} x_j^{n (N - j + 1)} y_j^{-jn} (r_j^2; q)_n \\
		& \qquad \times \displaystyle\prod_{1 \le i < j \le N} (r_i^2 x_i^{-1} x_j; q)_n (s_i^2 y_i^{-1} y_j; q)_n \displaystyle\prod_{i = 1}^N \displaystyle\prod_{j = 1}^N (x_j y_i^{-1}; q)_n^{-1}.
		\end{aligned} 
	\end{flalign}
	
\end{prop}

Let us mention that the above result was established as equation (42) in the work of Kulish--Ryasichenko \cite{SCRQS} in the case $n = 1$, in which setting the above vertex model for $F_{\boldsymbol{0}^N}$ becomes a free-fermionic six-vertex model. For $n > 1$, \Cref{0f0} is, to the best of our knowledge, new. 

We now proceed to more general properties of the $F$ and $G$ functions; the first concerns their symmetry. The following proposition shows that $G_{\boldsymbol{\lambda} / \boldsymbol{\mu}} (\textbf{x}; \textbf{r} \boldsymbol{\mid} \textbf{y}; \textbf{s})$ is symmetric under joint permutations of $\textbf{x}$ and $\textbf{r}$. It also shows that $F_{\boldsymbol{\lambda} / \boldsymbol{\mu}} (\textbf{x}; \textbf{r} \boldsymbol{\mid} \textbf{y}; \textbf{s})$ is ``almost symmetric,'' in that it multiplies by an explicit factor under any such joint permutation. In the below, $\mathfrak{S}_N$ denotes the symmetric group on $N$ elements,\index{S@$\mathfrak{S}_N$; symmetric group} and we set $\sigma (\mathcal{I}) = \big( i_{\sigma (1)}, i_{\sigma (2)}, \ldots , i_{\sigma (N)} \big)$ for any permutation $\sigma \in \mathfrak{S}_N$ and sequence $\mathcal{I} = (i_1, i_2, \ldots , i_N)$.\index{S@$\mathfrak{S}_N$; symmetric group!$\sigma (\mathcal{I})$}

\begin{prop}[\Cref{gxfxsigma} below] 
	
	\label{0gxfxsigma} 
	
	For any $\sigma \in \mathfrak{S}_N$,
	\begin{flalign*}
			& G_{\boldsymbol{\lambda} / \boldsymbol{\mu}} \big( \sigma (\textbf{\emph{x}}); \sigma (\textbf{\emph{r}}) \boldsymbol{\mid} \textbf{\emph{y}}; \textbf{\emph{s}} \big) = G_{\boldsymbol{\lambda} / \boldsymbol{\mu}} (\textbf{\emph{x}}; \textbf{\emph{r}} \boldsymbol{\mid} \textbf{\emph{y}}; \textbf{\emph{s}}); \\
			& F_{\boldsymbol{\lambda} / \boldsymbol{\mu}} \big( \sigma (\textbf{\emph{x}}); \sigma (\textbf{\emph{r}}) \boldsymbol{\mid} \textbf{\emph{y}}; \textbf{\emph{s}}\big) = F_{\boldsymbol{\lambda} / \boldsymbol{\mu}} (\textbf{\emph{x}}; \textbf{\emph{r}} \boldsymbol{\mid} \textbf{\emph{y}}; \textbf{\emph{s}}) \displaystyle\prod_{\substack{1 \le i < j \le N \\ \sigma (i) > \sigma (j) }} \displaystyle\frac{(r_j^2 x_i x_j^{-1}; q)_n}{(r_i^2 x_i^{-1} x_j; q)_n} \left( \displaystyle\frac{r_i^2 x_j}{r_j^2 x_i} \right)^n.
	\end{flalign*} 
	
\end{prop} 

We next have the following branching rules for $F$ and $G$.

\begin{prop}[\Cref{0fghbranching} below]
	
	\label{0fghbranching}
	
	Suppose $N = K + L$, and denote $\textbf{\emph{x}}' = (x_1, \ldots , x_K)$, $\textbf{\emph{x}}'' = (x_{K + 1}, \ldots , x_{K + L})$, $\textbf{\emph{r}}' = (r_1, \ldots,  r_K)$, and $\textbf{\emph{r}}'' = (r_{K + 1}, \ldots , r_{K + L})$. For any $\boldsymbol{\lambda}, \boldsymbol{\mu} \in \Sign_{n; M}$,
	\begin{flalign*}
		& \displaystyle\sum_{\boldsymbol{\nu} \in \SeqSign_{n; M}} G_{\boldsymbol{\lambda} / \boldsymbol{\nu}} (\textbf{\emph{x}}''; \textbf{\emph{r}}'' \boldsymbol{\mid} \textbf{\emph{y}}; \textbf{\emph{s}}) G_{\boldsymbol{\nu} / \boldsymbol{\mu}} (\textbf{\emph{x}}'; \textbf{\emph{r}}' \boldsymbol{\mid} \textbf{\emph{y}}; \textbf{\emph{s}}) = G_{\boldsymbol{\lambda} / \boldsymbol{\mu}} (\textbf{\emph{x}}; \textbf{\emph{r}} \boldsymbol{\mid} \textbf{\emph{y}}; \textbf{\emph{s}}),
	\end{flalign*} 
	
	\noindent Moreover, for any $\boldsymbol{\lambda} \in \Sign_{n; M + N}$ and $\boldsymbol{\mu} \in \Sign_{n; M}$,
	\begin{flalign*}
			& \displaystyle\sum_{\boldsymbol{\nu} \in \SeqSign_{n; M + L}} F_{\boldsymbol{\lambda} / \boldsymbol{\nu}} (\textbf{\emph{x}}'; \textbf{\emph{r}}' \boldsymbol{\mid} \textbf{\emph{y}}; \textbf{\emph{s}}) F_{\boldsymbol{\nu} / \boldsymbol{\mu}} (\textbf{\emph{x}}''; \textbf{\emph{r}}'' \boldsymbol{\mid} \textbf{\emph{y}}; \textbf{\emph{s}}) = F_{\boldsymbol{\lambda} / \boldsymbol{\mu}} (\textbf{\emph{x}}; \textbf{\emph{r}} \boldsymbol{\mid} \textbf{\emph{y}}; \textbf{\emph{s}}).
	\end{flalign*}
	
\end{prop}

We further show the following Cauchy type identity between $F_{\boldsymbol{\lambda}}$ and $G_{\boldsymbol{\lambda}}$. It is in fact a special case of a more general skew Cauchy identity between the functions $F_{\boldsymbol{\lambda} / \boldsymbol{\mu}}$ and $G_{\boldsymbol{\lambda} / \boldsymbol{\mu}}$, but we will not state the latter here and instead refer to \Cref{fgidentity} below for its precise formulation.

\begin{thm}[\Cref{fgsum2} below]
	
	\label{0fgsum2}
	
	Assume that there exists an integer $K > 1$ such that
	\begin{flalign*}
		\displaystyle\sup_{k > K} \displaystyle\max_{\substack{1 \le i \le M \\ 1 \le j \le N}}	\displaystyle\max_{\substack{a, b \in [0, n] \\ (a, b) \ne (n, 0)}} \Bigg| s_k^{2a + 2b - 2n} \displaystyle\frac{(s_k^2 u_j y_k^{-1}; q)_n (u_j y_k^{-1}; q)_a}{(u_j y_k^{-1}; q)_n (s_k^2 u_j y_k^{-1}; q)_a} \displaystyle\frac{(w_i y_k^{-1}; q)_b}{(s_k^2 w_i y_k^{-1}; q)_b} \Bigg| < 1.
	\end{flalign*}
	
	\noindent Then,
	\begin{flalign*}
			\displaystyle\sum_{\boldsymbol{\lambda} \in \SeqSign_{n; N}} & F_{\boldsymbol{\lambda}} (\textbf{\emph{u}}; \textbf{\emph{r}} \boldsymbol{\mid} \textbf{\emph{y}}; \textbf{\emph{s}}) G_{\boldsymbol{\lambda}} (\textbf{\emph{w}}; \textbf{\emph{t}} \boldsymbol{\mid} \textbf{\emph{y}}; \textbf{\emph{s}}) = F_{\boldsymbol{0}^N} (\textbf{\emph{u}}; \textbf{\emph{r}} \boldsymbol{\mid} \textbf{\emph{y}}; \textbf{\emph{s}}) \displaystyle\prod_{i = 1}^M \displaystyle\prod_{j = 1}^N \displaystyle\frac{(t_i^2 u_j w_i^{-1}; q)_n}{t_i^{2n} (u_j w_i^{-1}; q)_n},
	\end{flalign*}

	\noindent with $F_{\boldsymbol{0}^N} (\textbf{\emph{u}}; \textbf{\emph{r}} \boldsymbol{\mid} \textbf{\emph{y}}; \textbf{\emph{s}})$ given by \eqref{0f0n}.
	
\end{thm} 

Let us conclude by mentioning that we also establish a contour integral representation for $G_{\boldsymbol{\lambda}}$ as \Cref{gintegralf2} (and also for the more general skew function $G_{\boldsymbol{\lambda} / \boldsymbol{\mu}}$ as \Cref{glambdamuidentity} and \Cref{gintegralf}) below. Although we will not precisely state the fully general version of this result in this introduction, we will provide its degeneration to the LLT case as \Cref{0lintegral} below.

\section{LLT Polynomials}

\label{0Polynomials}

We next consider specializations for our symmetric functions $F_{\boldsymbol{\lambda} / \boldsymbol{\mu}} (\textbf{x}; \textbf{r} \boldsymbol{\mid} \textbf{y}; \textbf{s})$ and $G_{\boldsymbol{\lambda} / \boldsymbol{\mu}} (\textbf{x}; \textbf{r} \boldsymbol{\mid} \textbf{y}; \textbf{s})$ under some of the limit degenerations depicted in \Cref{wlimit}. In particular, here we will be mainly focused on the bottom-right degeneration (with weight $x^d q^{\varphi (\textbf{D}, \textbf{C} + \textbf{D})} \textbf{1}_{v = 0}$) shown there. In what follows, we fix a set of complex numbers $\textbf{x} = (x_1, x_2, \ldots , x_N)$.

Define\footnote{Our reason for the notation $\mathcal{G}_{\boldsymbol{\lambda} / \boldsymbol{\mu}} (\textbf{x}; \infty \boldsymbol{\mid} 0; 0)$ (and $\mathcal{F}_{\boldsymbol{\lambda} / \boldsymbol{\mu}} (\textbf{x}; \infty \boldsymbol{\mid} 0; 0)$) is that this function is obtained from the original $G_{\boldsymbol{\lambda} / \boldsymbol{\mu}} (\textbf{x}; \textbf{r} \boldsymbol{\mid} \textbf{y}; \textbf{s})$ (and $F_{\boldsymbol{\lambda} / \boldsymbol{\mu}} (\textbf{x}; \textbf{r} \boldsymbol{\mid} \textbf{y}; \textbf{s})$, respectively) by, after suitably normalizing, first sending each parameter in $\textbf{s}$ to $0$, and then sending each one in $\textbf{r}$ and $\textbf{y}$ to $\infty$ and $0$, respectively. See \eqref{limitg} (and \eqref{limithf}) below.} $\mathcal{G}_{\boldsymbol{\lambda} / \boldsymbol{\mu}} (\textbf{x}; \infty \boldsymbol{\mid} 0; 0)$ and $\mathcal{F}_{\boldsymbol{\lambda} / \boldsymbol{\mu}} (\textbf{x}; \infty \boldsymbol{\mid} 0; 0)$ as\index{G@$G_{\boldsymbol{\lambda} / \boldsymbol{\mu}} (\textbf{x}; \textbf{r} \boldsymbol{\mid} \textbf{y}; \textbf{s})$!$\mathcal{G}_{\boldsymbol{\lambda} / \boldsymbol{\mu}} (\textbf{x}; \infty \boldsymbol{\mid} 0; 0)$}\index{F@$F_{\boldsymbol{\lambda} / \boldsymbol{\mu}} (\textbf{x}; \textbf{r} \boldsymbol{\mid} \textbf{y}; \textbf{s})$!$\mathcal{F}_{\boldsymbol{\lambda} / \boldsymbol{\mu}} (\textbf{x}; \infty \boldsymbol{\mid} 0; 0)$} partition functions for the vertex models shown in \eqref{gvertex} and \eqref{fvertex}, respectively, but where for each $j \in [1, N]$ we replace the original $W_{x_j / y_i}$ and $\widehat{W}_{x_j / y_i}$ weights in the $j$-th row of the model with the degenerated ones defined by
\begin{flalign*}
W_z (\textbf{A}, \textbf{B}; \textbf{C}, \textbf{D} \boldsymbol{\mid} r, s) & \mapsto x_j^d q^{\varphi (\textbf{D}, \textbf{C} + \textbf{D})} \textbf{1}_{v = 0} \textbf{1}_{\textbf{A} + \textbf{B} = \textbf{C} + \textbf{D}}; \\ 
\widehat{W}_z (\textbf{A}, \textbf{B}; \textbf{C}, \textbf{D} \boldsymbol{\mid} r, s) & \mapsto x_j^{d - n} q^{\varphi (\textbf{D}, \textbf{C} + \textbf{D}) - \binom{n}{2}} \textbf{1}_{v = 0} \textbf{1}_{\textbf{A} + \textbf{B} = \textbf{C} + \textbf{D}}, 
\end{flalign*}

\noindent respectively. The below theorem then indicates that these specialized functions are given (up to global multiplicative factors) by the Lascoux--Leclerc--Thibon (LLT) polynomials $\mathcal{L}_{\boldsymbol{\lambda} / \boldsymbol{\mu}} (\textbf{x}; q)$\index{L@$\mathcal{L}_{\boldsymbol{\lambda} / \boldsymbol{\mu}} (\textbf{x})$; LLT polynomial} originally introduced in \cite{RT}; we refer to \Cref{Partitionsn} and \Cref{FunctionGPartitions} below for the detailed definition of these polynomials. In what follows, for any sequence of signatures $\boldsymbol{\lambda} = \big( \lambda^{(1)}, \lambda^{(2)}, \ldots , \lambda^{(n)} \big)\in \SeqSign_n$, set
\begin{flalign}
	\label{0psi} 
	\psi (\boldsymbol{\lambda}) = \displaystyle\frac{1}{2} \displaystyle\sum_{1 \le i < j \le n} \displaystyle\sum_{a \in \mathfrak{T}_i} \displaystyle\sum_{b \in \mathfrak{T}_j} \textbf{1}_{a > b},
\end{flalign}
\index{0@$\psi$}

\noindent where we have abbreviated $\mathfrak{T}_k = \mathfrak{T} \big( \lambda^{(k)} \big)$ for each index $k \in [1, n]$ (recall $\mathfrak{T}$ from \eqref{t}).

\begin{thm}[Parts \ref{lg1} and \ref{lf1} of \Cref{limitg0} below]
	
	\label{0limitg0}
	
	The following statements hold.
	
	\begin{enumerate} 
		
		\item \label{0lg1} For any sequences of signatures $\boldsymbol{\lambda}, \boldsymbol{\mu} \in \SeqSign_{n; M}$,
		\begin{flalign*}
			\mathcal{G}_{\boldsymbol{\lambda} / \boldsymbol{\mu}} (\textbf{\emph{x}}; \infty \boldsymbol{\mid} 0; 0)  = q^{\psi (\boldsymbol{\lambda}) - \psi (\boldsymbol{\mu})} \mathcal{L}_{\boldsymbol{\lambda} / \boldsymbol{\mu}} (\textbf{\emph{x}}; q). 
		\end{flalign*}
		
		\item \label{0lf1} For any sequences of signatures $\boldsymbol{\lambda} \in \Sign_{n; M + N}$ and $\boldsymbol{\mu} \in \SeqSign_{n; M}$,
		\begin{flalign*} 
			\mathcal{F}_{\boldsymbol{\lambda} / \boldsymbol{\mu}} (\textbf{\emph{x}}; \infty \boldsymbol{\mid} 0; 0)  = q^{\psi (\boldsymbol{\mu}) - \psi (\boldsymbol{\lambda}) + \binom{M}{2} \binom{n}{2} / 2 - \binom{M + N}{2} \binom{n}{2} / 2} \mathcal{L}_{\boldsymbol{\lambda} / \boldsymbol{\mu}} (q^{1 - n} \textbf{\emph{x}}^{-1}; q) \displaystyle\prod_{j = 1}^N x_j^{n (j - M - N)}.
		\end{flalign*}

	\end{enumerate}
	
\end{thm} 

In addition to using the bottom-right limiting weights shown in \Cref{wlimit} to degenerate our symmetric functions to the LLT polynomials, it is also possible to use the bottom-middle and bottom-left ones depicted there. We refer to parts \ref{lg2} and \ref{lh1} of \Cref{limitg0} below for the specific statements of these results. 

The vertex model interpretation given by \Cref{0limitg0} is well-suited for establishing various properties about the LLT polynomials. We provide four here; the first, second, and fourth are known from \cite{TRA,RT,AAPP}, respectively, and the third appears to be new. 

The first is a pair of Cauchy identities, originally shown as Theorem 35 and Proposition 36 of \cite{TRA}. In what follows, for any signature sequence $\boldsymbol{\lambda} = \big( \lambda^{(1)}, \lambda^{(2)}, \ldots , \lambda^{(n)} \big) \in \Sign_n$, we define $\boldsymbol{\lambda}' = \big( \lambda'^{(1)}, \lambda'^{(2)}, \ldots , \lambda'^{(n)} \big) \in \Sign_n$ by setting each $\lambda'^{(j)}$ to be the dual\index{0@$\boldsymbol{\lambda}, \boldsymbol{\mu}$; typical signature sequences!$\boldsymbol{\lambda}' / \boldsymbol{\mu}'$; dual!$\boldsymbol{\lambda}'$}\index{0@$\lambda, \mu$; typical signatures or partitions!$\lambda' / \mu'$; dual!$\lambda'$} of $\lambda^{(n - j)}$ (that is, obtained from the latter by transposing its Young diagram, as shown on the top of \Cref{lambdamulambdamu1} below). 

\begin{cor}[\Cref{suml1} and \Cref{suml2} below]
	
	\label{0suml1} 
	
	Fix sequences of complex numbers $\textbf{\emph{x}} = (x_1, x_2, \ldots , x_N)$ and $\textbf{\emph{y}} = (y_1, y_2, \ldots , y_M)$. If $|q| < 1$ and $|x_j|, |y_i| < 1$ for each $i \in [1, M]$ and $j \in [1, N]$, then 
	\begin{flalign}
		\label{0lsum1}
		\displaystyle\sum_{\boldsymbol{\lambda} \in \SeqSign_{n; N}} \mathcal{L}_{\boldsymbol{\lambda}} (\textbf{\emph{x}}; q) \mathcal{L}_{\boldsymbol{\lambda}} (\textbf{\emph{y}}; q) = \displaystyle\prod_{i = 1}^M \displaystyle\prod_{j = 1}^N (x_j y_i; q)_n^{-1},
	\end{flalign}
	
	\noindent and
	\begin{flalign}
		\label{0lsum2}
		\displaystyle\sum_{\boldsymbol{\lambda} \in \SeqSign_{n; N}} \mathcal{L}_{\boldsymbol{\lambda}} (\textbf{\emph{x}}; q) \mathcal{L}_{\boldsymbol{\lambda}'} (\textbf{\emph{y}}; q^{-1}) = \displaystyle\prod_{i = 1}^M \displaystyle\prod_{j = 1}^N (- q^{(n - 1) / 2} x_j y_i; q)_n.
	\end{flalign}
	
\end{cor}

The standard Cauchy identity \eqref{0lsum1} essentially follows directly from \Cref{0fgsum2} and \Cref{0limitg0}. The dual Cauchy identity \eqref{0lsum2} is shown\footnote{As discussed in \cite{TRA}, it is alternatively possible to derive the dual Cauchy identity directly from the standard one by applying to it the involutory automorphism $\omega$ on the ring of symmetric functions that interchanges complete and elementary symmetric functions.} also using \Cref{0fgsum2} and the second part of \Cref{0limitg0} expressing an LLT polynomial as a degeneration of a $F$ function; however, we further require part \ref{lg2} of \Cref{limitg0} below that expresses a dual LLT polynomial $\mathcal{L}_{\boldsymbol{\lambda}'} (\textbf{x}; q^{-1})$ as a degeneration of a $G$ function (under the bottom-middle limit depicted in \Cref{wlimit}).

The second result is a \emph{stability} property, originally shown as Theorem 6.6 of \cite{RT}, stating that if each signature in $\boldsymbol{\lambda}$ only contains one part, then the LLT polynomial $\mathcal{L}_{\boldsymbol{\lambda}} (\textbf{x}; q)$ is a modified Hall--Littlewood polynomial $Q_{\lambda}' (\textbf{x})$;\index{Q@$Q_{\lambda} (\textbf{x})$; Hall--Littlewood polynomial!$Q_{\lambda}' (\textbf{x})$; modified Hall--Littlewood polynomial} we refer to the beginning of \Cref{StabilityPolynomials} for a definition of the latter. 

\begin{prop}[\Cref{qlambdaqlambda} below]
	
	\label{0qlambdaqlambda} 
	
	Fix a signature $\lambda = (\lambda_1, \lambda_2, \ldots,  \lambda_n) \in \Sign_n$, and define $\lambda^{(j)} = (\lambda_{n - j + 1}) \in \Sign_1$, for each $j \in [1, n]$. Letting $\boldsymbol{\lambda} = \big( \lambda^{(1)}, \lambda^{(2)}, \ldots , \lambda^{(n)} \big) \in \SeqSign_{n; 1}$, we have $\mathcal{L}_{\boldsymbol{\lambda}} (\textbf{\emph{x}}; q) = Q_{\lambda}' (\textbf{\emph{x}})$. 
	
\end{prop} 

The use of the vertex model interpretation for the LLT polynomials (\Cref{0limitg0}) in establishing this stability property is that there also exists a vertex model, due to Garbali--Wheeler \cite{MPI}, for the modified Hall--Littlewood polynomials. The vertex weights for the latter model are quite similar to those of the former, enabling a comparison between the two that verifies \Cref{0qlambdaqlambda}. 

Let us mention that Theorem 3.4 of the recent work \cite{VMFP} (pursued independently from ours) also establishes the first part of \Cref{0limitg0}. Using this, Corollary 6.13 of \cite{VMFP} proves the standard Cauchy identity \eqref{0lsum1}, and Proposition 5.1 of \cite{VMFP} proves the stability property \Cref{0qlambdaqlambda}. 

Our third result provides a vertex model for the LLT polynomials under a certain family of plethysms. We will provide more detailed definitions in \Cref{Polynomialsr} below but, briefly, one views the skew LLT polynomial $\mathcal{L}_{\boldsymbol{\lambda} / \boldsymbol{\mu}} (\textbf{x}; q)$ as an element in the ring of symmetric functions in $\textbf{x}$. Defining the formal sum $X = \sum_{x \in \textbf{x}} x$, for any variable $u$ one lets $\mathcal{L}_{\boldsymbol{\lambda} / \boldsymbol{\mu}} \big[ (1 - u) X \big]$ denote the image of $\mathcal{L}_{\boldsymbol{\lambda} / \boldsymbol{\mu}} (\textbf{x}; q)$ under the plethystic substitution $X \mapsto (1 - u) X$. 

Now let $\mathcal{G}_{\boldsymbol{\lambda} / \boldsymbol{\mu}} ( \textbf{x}; r \boldsymbol{\mid} 0; 0)$\index{G@$G_{\boldsymbol{\lambda} / \boldsymbol{\mu}} (\textbf{x}; \textbf{r} \boldsymbol{\mid} \textbf{y}; \textbf{s})$!$\mathcal{G}_{\boldsymbol{\lambda} / \boldsymbol{\mu}} (\textbf{x}; r \boldsymbol{\mid} 0; 0)$} denote the partition function for the vertex model shown in \eqref{gvertex}, but where for each $j \in [1, N]$ we replace the original weight $W_{x_j / y_i}$ in the $j$-th row of the model by the degenerated ones depicted in the top-right of \Cref{wlimit}, namely, 
\begin{flalign*}
W_z (\textbf{A}, \textbf{B}; \textbf{C}, \textbf{D} \boldsymbol{\mid} r, s) & \mapsto x_j^d (-r^2)^{c - a - v} q^{\varphi (\textbf{D} - \textbf{V}, \textbf{C}) + \varphi (\textbf{V}, \textbf{A}) + \binom{b}{2} - dv + v} \textbf{1}_{\textbf{A} + \textbf{B} = \textbf{C} + \textbf{D}} \\ 
& \qquad \times \displaystyle\frac{(r^2; q)_d}{(r^2; q)_b} \displaystyle\prod_{j: B_j - D_j = 1} \big( 1 - r^{-2} q^{-B_{[j + 1, n]} - D_{[1, j - 1]}} \big).
\end{flalign*}

\noindent Then, we have the following result. 

\begin{prop}[\Cref{gl} below]
	
	\label{0gl} 
	
	For any variable $r$, we have the plethystic identity 
	\begin{flalign*}
		\mathcal{L}_{\boldsymbol{\lambda} / \boldsymbol{\mu}} \big[ (1 - r^{-2}) X \big] = q^{\psi(\boldsymbol{\mu}) - \psi(\boldsymbol{\lambda})} \mathcal{G}_{\boldsymbol{\lambda} / \boldsymbol{\mu}} (\textbf{\emph{x}}; r \boldsymbol{\mid} 0; 0).
	\end{flalign*} 
	
\end{prop}

By combining \Cref{0gl} with \Cref{0fghbranching}, one can obtain explicit branching rules for the LLT polynomials under the plethystic substitutions $X \mapsto (1 - u) X$. Let us mention that the skew LLT polynomials under such plethysms have been considered before in the algebraic combinatorics literature. For example, in certain cases when $u = q$ (that is, $r = q^{-1 / 2}$), it was shown by Carlsson--Mellit \cite{S} that they coincide with chromatic quasisymmetric functions, as introduced by Shareshian--Wachs \cite{CSF,CSFV}, associated with the incomparability graph of a unit interval order.

The fourth result provides a contour integral formula for skew LLT polynomials in terms of nonsymmetric Hall--Littlewood polynomials. To state it, we must introduce some additional notation. For any sequence of positive integers $\mu = (\mu_1, \mu_2, \ldots , \mu_n) \in \mathbb{Z}_{> 0}^n$ and sequence of complex numbers $\textbf{x} = (x_1, x_2, \ldots , x_n)$, let $\mathtt{f}_{\mu}^{(q)} (\textbf{x} \boldsymbol{\mid} 0)$ and $\mathtt{g}_{\mu}^{(q)} (\textbf{x} \boldsymbol{\mid} 0)$ denote certain normalizations of the nonsymmetric Hall--Littlewood polynomials,\index{F@$\mathtt{f}_{\mu}^{(q)} (\textbf{x} \boldsymbol{\mid} s)$}\index{G@$\mathtt{g}_{\mu}^{(q)} (\textbf{x} \boldsymbol{\mid} s)$} given more precisely in \Cref{lambdafg} below. Moreover, for any $\boldsymbol{\lambda} \in \SeqSign_{n; M}$, let $\Upsilon (\boldsymbol{\lambda})$\index{0@$\Upsilon (\boldsymbol{\lambda})$} denote the set of sequences $\kappa = (\kappa_1, \kappa_2, \ldots , \kappa_{nM}) \in \mathbb{Z}_{> 0}^{nM}$ such that $\kappa^{(i)} = \big( \kappa_{(i - 1) M + 1}, \kappa_{(i - 1) M + 2}, \ldots , \kappa_{iM} \big)$ is a permutation of $\mathfrak{T} \big( \lambda^{(i)} \big)$ (where we recall $\mathfrak{T}$ from \eqref{t}), for each $i \in [1, n]$. For any $\kappa \in \Upsilon (\boldsymbol{\lambda})$, we also set
	\begin{flalign*}
		\inv_{\boldsymbol{\lambda}} (\kappa) = \displaystyle\sum_{i = 1}^n \inv \big( \kappa^{(i)} \big), \qquad \text{where} \qquad \inv \big( \kappa^{(i)} \big) = \displaystyle\sum_{1 \le j < k \le M} \textbf{1}_{\kappa_{(i - 1) M + j} > \kappa_{(i - 1) M + k}}.
	\end{flalign*}
	\index{I@$\inv$}
	\index{I@$\inv$!$\inv_{\boldsymbol{\lambda}}$}

Then we have the following theorem, which provides a relationship between the skew LLT polynomials $\mathcal{L}_{\boldsymbol{\lambda} / \boldsymbol{\mu}}$ and the nonsymmetric Hall--Littlewood polynomials $\mathtt{f}_{\mu}^{(q)}$ and $\mathtt{g}_{\mu}^{(q)}$. It was originally implicitly shown as equation (34) of \cite{AAPP}.

\begin{thm}[\Cref{llambdamuidentity} below]
	
	\label{0llambdamuidentity} 
	
	Fix $\boldsymbol{\lambda}, \boldsymbol{\mu} \in \SeqSign_{n; M}$ and $\textbf{\emph{x}} = (x_1, x_2, \ldots , x_N) \subset \mathbb{C}$. Letting $\kappa = \Upsilon (\boldsymbol{\lambda})$ denote the unique element of $\Upsilon (\boldsymbol{\lambda})$ such that $\inv_{\boldsymbol{\lambda}} (\kappa) = 0$, we have
	\begin{flalign}
		\label{lsumintegral} 
		\begin{aligned} 
		\mathcal{L}_{\boldsymbol{\lambda} / \boldsymbol{\mu}} (\textbf{\emph{x}}) = \displaystyle\frac{1}{(2 \pi \textbf{\emph{i}})^{nM}} \displaystyle\frac{q^{\psi (\boldsymbol{\mu}) - \psi (\boldsymbol{\lambda})}}{(1 - q^{-1})^{nM}} & \displaystyle\oint \cdots \displaystyle\oint \mathtt{f}_{\kappa}^{(q)} (\textbf{\emph{u}}^{-1} \boldsymbol{\mid} 0) \displaystyle\sum_{\nu \in \Upsilon (\boldsymbol{\mu})} (-1)^{\inv_{\boldsymbol{\mu}} (\nu)} \mathtt{g}_{\nu}^{(q)} (\textbf{\emph{u}} \boldsymbol{\mid} 0)  \\
		& \qquad \times \displaystyle\prod_{1 \le i < j \le nM} \displaystyle\frac{u_j - u_i}{q^{-1} u_j - u_i} \displaystyle\prod_{i = 1}^{nM} \displaystyle\prod_{j = 1}^N \displaystyle\frac{1}{1 - u_i x_j} \displaystyle\prod_{i = 1}^{nM} \displaystyle\frac{du_i}{u_i},  
		\end{aligned} 
	\end{flalign}
	
	\noindent where $\textbf{\emph{u}} = (u_1, u_2, \ldots , u_{nM}) \subset \mathbb{C}$, and each $u_i$ is integrated along a positively oriented, closed contour $\Gamma_i$ satisfying the following two properties. First, each $\Gamma_i$ contains $0$ and does not contain $q^{-k} x_j^{-1}$ for all integers $k \in [1, nM - 1]$ and $j \in [1, M]$. Second, the $\{ \Gamma_i \}$ are mutually non-intersecting, and $\Gamma_{i - 1}$ is contained in both $\Gamma_i$ and $q^{-1} \Gamma_i$ for each $i \in [2, nM]$. 
	
\end{thm} 

Let us briefly indicate the relation between \Cref{0llambdamuidentity} and the results of \cite{AAPP}; we refer to \Cref{integrall} below for a more detailed explanation. By \eqref{gcoefficientl} below, \eqref{lsumintegral} can be reformulated as 
\begin{flalign}
	\label{0gcoefficientl}
	\text{Coeff} \big[ \mathcal{L}_{\boldsymbol{\lambda} / \boldsymbol{\mu}}, s_{\theta} \big] = q^{\psi (\boldsymbol{\mu}) - \psi (\boldsymbol{\lambda})} \text{Coeff} \Bigg[ s_{\theta} \displaystyle\sum_{\nu \in \Upsilon (\boldsymbol{\mu})} (-1)^{\inv_{\boldsymbol{\mu}} (\nu)} \mathtt{g}_{\nu}; \mathtt{g}_{\kappa} \Bigg],
\end{flalign} 

\noindent for any signature $\theta$. Here, we recall the Schur functions $s_{\lambda}$\index{S@$s_{\lambda}$; Schur function} and have abbreviated $\mathcal{L}_{\boldsymbol{\lambda} / \boldsymbol{\mu}} = \mathcal{L}_{\boldsymbol{\lambda} / \boldsymbol{\mu}} (\textbf{x}; q)$ and $\mathtt{g}_{\mu} = \mathtt{g}_{\mu}^{(q)} (\textbf{x} \boldsymbol{\mid} 0)$. Moreover, for any (non)symmetric function $F$ and basis $\{ h_{\mu} \}$ for the space of (non)symmetric functions, we have also let $\text{Coeff} [F; h_{\lambda}]$ denote the coefficient of $h_{\lambda}$ in the expansion of $F$ over $\{ h_{\mu} \}$. The equality \eqref{0gcoefficientl} was originally established as equation (34) of \cite{AAPP}, which was central in their proof of nonnegativity for the expansion of the LLT polynomials in the Schur basis. Indeed, the coefficients in this expansion are given by the left side of \eqref{0gcoefficientl}, and those on the right side can be equated with matrix entries for the action of a certain operator on a submodule of the affine Hecke algebra. It is shown in \cite{AAPP} that this operator admits an algebro-geometric interpretation that can be used to establish nonnegativity for its matrix entries.

In the case when $\boldsymbol{\mu} = \boldsymbol{0}^M$, the anti-symmetrizing sum of $\mathtt{g}_{\nu}^{(q)}$ on the right side of \eqref{lsumintegral} can be evaluated explicitly, giving rise to the following corollary.

\begin{cor}[\Cref{lintegralf2} below]
	
	\label{0lintegral} 

	Adopting the notation of \Cref{0llambdamuidentity}, we have
	\begin{flalign*}
		\mathcal{L}_{\boldsymbol{\lambda}} (\textbf{\emph{x}}; q) = \displaystyle\frac{q^{\binom{n}{2} \binom{N}{2} / 2 - \psi (\boldsymbol{\lambda})}}{(2 \pi \textbf{\emph{i}})^{nM}} \displaystyle\oint \cdots \displaystyle\oint & \displaystyle\prod_{1 \le i < j \le nM} \displaystyle\frac{u_j - u_i}{u_j - q u_i}  \displaystyle\prod_{k = 0}^{n - 1} \displaystyle\prod_{1 \le i < j \le M} (q u_{Mk + i} - u_{Mk + j}) \\
		& \times \mathtt{f}_{\kappa}^{(q)} (\textbf{\emph{u}}^{-1} \boldsymbol{\mid} 0)  \displaystyle\prod_{i = 1}^{nM} \displaystyle\prod_{j = 1}^N \displaystyle\frac{1}{1 - u_i x_j} \displaystyle\prod_{i = 1}^{nM} \displaystyle\frac{du_i}{u_i},
	\end{flalign*}
	
	\noindent where each $u_i$ is integrated along $\Gamma_i$.
	
\end{cor} 

Under the condtion when all parts of any $\lambda^{(i)}$ are less by at least $M - 1$ than all parts of any $\lambda^{(j)}$ for $i < j$, the nonsymmetric Hall--Littlewood function appearing in the integral from \Cref{0lintegral} factors completely. This gives rise to the following simplified formulas for certain LLT polynomials.

\begin{cor}[\Cref{lintegralf3} below]
	
	\label{0lintegralf3}
	
	Adopting the notation of \Cref{0llambdamuidentity}; assume that $\lambda_k^{(i)} + M - 1 \le \lambda_{k'}^{(j)}$ whenever $1 \le i < j \le n$ and $k, k' \in [1, M]$. Then,
	\begin{flalign*} 
		\mathcal{L}_{\boldsymbol{\lambda}} (\textbf{\emph{x}}; q) = \displaystyle\frac{q^{\binom{n}{2} \binom{N}{2} / 2} - \psi (\boldsymbol{\lambda})}{(2 \pi \textbf{\emph{i}})^{nM}} \displaystyle\oint \cdots \displaystyle\oint & \displaystyle\prod_{1 \le i < j \le nM} \displaystyle\frac{u_j - u_i}{u_j - q u_i}  \displaystyle\prod_{k = 0}^{n - 1} \displaystyle\prod_{1 \le i < j \le M} (q u_{Mk + i} - u_{Mk + j})  \\
		& \times \displaystyle\prod_{i = 1}^n \displaystyle\prod_{j = 1}^M u_{iM - j + 1}^{j - \lambda_j^{(i)} - M - 1} \displaystyle\prod_{i = 1}^{nM} \displaystyle\prod_{j = 1}^N \displaystyle\frac{1}{1 - u_i x_j}  \displaystyle\prod_{i = 1}^{nM} du_i,
	\end{flalign*}	
	\noindent where each $u_i$ is integrated along $\Gamma_i$.
	
\end{cor}

By applying a residue expansion of integral in \Cref{0lintegralf3}, it might be possible to obtain a symmetrization formula for LLT polynomials satisfying the ordering constraint described there. However, we will not pursue this here.

\section{Factorial LLT Polynomials}

\label{0Polynomials0}

In \Cref{0Polynomials} we explained how considering the degenerations of the $F$ and $G$ functions under the bottom-right limit depicted in \Cref{wlimit} gave rise to the LLT polynomials. Here we consider the degenerations of these functions under the limit depicted directly above it, namely, with weight $q^{\varphi (\textbf{D}, \textbf{C} + \textbf{D})} (-u)^d (u; q)_{c + d}^{-1} \textbf{1}_{v = 0}$. We will in fact be concerned with the $F$ function, which requires a normalization of this weight as in \eqref{wzwz}. In what follows, we fix a finite set of complex numbers $\textbf{x} = (x_1, x_2, \ldots , x_N)$ and an infinite set of complex numbers $\textbf{y} = (y_1, y_2, \ldots )$. 

Define $\mathcal{F}_{\boldsymbol{\lambda} / \boldsymbol{\mu}} (\textbf{x}; \infty \boldsymbol{\mid} \textbf{y}; 0)$ as the partition functions for the vertex model shown in \eqref{fvertex}, but where for each $i \ge 1$ and $j \in [1, N]$ we replace the original $\widehat{W}_{x_j / y_i}$ weight at the vertex $(i, j) \in \mathbb{Z}_{> 0}^2$ with the degenerated one given by
\begin{flalign*}
	\widehat{W}_z (\textbf{A}, \textbf{B}; \textbf{C}, \textbf{D} \boldsymbol{\mid} r, s) & \mapsto (-x_j y_i^{-1})^{d - n} q^{\varphi (\textbf{D}, \textbf{C} + \textbf{D}) - \binom{n}{2}} \displaystyle\frac{(x_j y_i^{-1}; q)_n}{(x_j y_i^{-1}; q)_{c + d}} \textbf{1}_{v = 0} \textbf{1}_{\textbf{A} + \textbf{B} = \textbf{C} + \textbf{D}}.
\end{flalign*} 

\noindent Observe from the second part of \Cref{0limitg0} that, upon letting each of the $y_i$ tend to $\infty$, this function $\mathcal{F}_{\boldsymbol{\lambda}} (\textbf{x}; \infty \boldsymbol{\mid} \textbf{y}; 0)$ degenerates to an LLT polynomial. Although $\mathcal{F}_{\boldsymbol{\lambda} / \boldsymbol{\mu}} (\textbf{x}; \infty \boldsymbol{\mid} \textbf{y}; 0)$ is not a polynomial in $\textbf{x}$ and $\textbf{y}$, a mild modification of it is. More specifically, denoting $\textbf{z}^e = (z_1^e, z_2^e, \ldots )$ for any real number $e$ and (possibly infinite) sequence $\textbf{z} = (z_1, z_2, \ldots )$ (and recalling the sequence $\mathscr{S} (\boldsymbol{\lambda}) = \big( \textbf{S}_1 (\boldsymbol{\lambda}), \textbf{S}_2 (\boldsymbol{\lambda}), \ldots \big)$ from \Cref{0Functions}), the function 
\begin{flalign*}
\check{\mathcal{F}}_{\boldsymbol{\lambda}} (\textbf{x} \boldsymbol{\mid} \textbf{y}) = \mathcal{F}_{\boldsymbol{\lambda}} (\textbf{x}^{-1}; \infty \boldsymbol{\mid} \textbf{y}^{-1}; 0) \displaystyle\prod_{j = 1}^N x_j^{n (j - N)} \displaystyle\prod_{i = 1}^{\infty} y_i^{nN - \sum_{k = 1}^i |\textbf{S}_k (\boldsymbol{\lambda})|},
\end{flalign*}

\noindent is quickly verified to be a polynomial in $(\textbf{x}, \textbf{y})$ of total degree $|\boldsymbol{\lambda}|$. We refer to $\check{\mathcal{F}}_{\boldsymbol{\lambda}} (\textbf{x} \boldsymbol{\mid} \textbf{y})$ as a \emph{factorial LLT polynomial}, since we will show that it satisfies a vanishing property reminiscent of the one satisfied by factorial Schur (and interpolation Macdonald) polynomials. 

These vanishing points for this LLT factorial polynomial will take the form $\{ x_j = q^{n - \kappa_j - 1} y_{\mathfrak{m}_j} \}$ for some integers $\mathfrak{m}_j \ge 1$ and $\kappa_j \in [0, n - 1]$. It will be convenient to express such specializations through \emph{marked sequences}, which are pairs $(\mathfrak{m}, \kappa)$ of integer sets of the same length, such the entries $\mathfrak{m}_1 > \mathfrak{m}_2 > \cdots > \mathfrak{m}_{\ell}$ of $\mathfrak{m}$ are decreasing and positive, and such that each entry of $\kappa = (\kappa_1, \kappa_2, \ldots , \kappa_{\ell})$ is in $[0, n - 1]$. In what follows, we will often only refer to $\mathfrak{m}$ as the marked sequence and view $\kappa$ as its \emph{marking}, in that each entry $\mathfrak{m}_j \in \mathfrak{m}$ is \emph{marked} by the corresponding entry $\kappa_j \in \kappa$. 

To state the vanishing property, we require the notion of a splitting for a marked sequence. 

\begin{definition}[\Cref{sequencem} below]
	
	\label{0sequencem} 
	
	Let $\mathfrak{m} = (\mathfrak{m}_1, \mathfrak{m}_2, \ldots , \mathfrak{m}_{\ell})$ denote a marked sequence with marking $\kappa = (\kappa_1, \kappa_2, \ldots , \kappa_{\ell})$. A \emph{splitting} of $(\mathfrak{m}, \kappa)$\index{M@$(\mathfrak{m}, \kappa)$; marked sequence} is a family $\mathfrak{M} = \big( \mathfrak{m}^{(1)}, \mathfrak{m}^{(2)}, \ldots , \mathfrak{m}^{(n)} \big)$ of decreasing subsequences of $\mathbb{Z}_{> 0}$, such that the following two properties holds. 
	
	\begin{enumerate}
		\item For any $i \in [1, n]$, every entry $m \in \mathfrak{m}^{(i)}$ is equal to $\mathfrak{m}_j$, for some $j = j (m) \in [1, \ell]$.
		\item For any $j \in [1, \ell]$, there exist at least $n - \kappa_j$ distinct indices $i \in [1, n]$ such that $\mathfrak{m}_j \in \mathfrak{m}^{(i)}$.
	\end{enumerate}
	
\end{definition}

Moreover, given two signatures $\mathfrak{m} = (\mathfrak{m}_1, \mathfrak{m}_2, \ldots , \mathfrak{m}_{\ell}) \in \Sign_{\ell}$ and $\mathfrak{n} = (\mathfrak{n}_1, \mathfrak{n}_2, \ldots , \mathfrak{n}_k) \in \Sign_k$, we write $\mathfrak{m} \nprec \mathfrak{n}$ if there exists an integer $j \ge 0$ such that $\mathfrak{m}_{\ell - j} > \mathfrak{n}_{k - j}$, where we set $\mathfrak{m}_i = \infty = \mathfrak{n}_i$ if $i \le 0$. The vanishing property is then provided by the following theorem.

\begin{thm}[\Cref{lambdam0} and \Cref{lambdam0f} below]
	
	\label{0lambdam0}
	
	Fix a marked sequence $(\mathfrak{m}, \kappa)$ of length $\ell \le N$ and some $\boldsymbol{\lambda} = \big( \lambda^{(1)}, \lambda^{(2)}, \ldots , \lambda^{(n)} \big)\in \SeqSign_{n; N}$. Denoting $\mathfrak{l}^{(i)} = \mathfrak{T} \big( \lambda^{(i)} \big)$ for each $i$, assume for any splitting $\mathfrak{M} = \big( \mathfrak{m}^{(1)}, \mathfrak{m}^{(2)}, \ldots , \mathfrak{m}^{(n)} \big)$ of $(\mathfrak{m}, \kappa)$ that there exists an index $h = h (\mathfrak{M}) \in [1, n]$ such that $\mathfrak{l}^{(h)} \nprec \mathfrak{m}^{(h)}$. If $x_j = q^{n -\kappa_j - 1} y_{\mathfrak{m}_j}$ holds for each $j \in [1, \ell]$, then $\check{\mathcal{F}}_{\boldsymbol{\lambda}} (\textbf{\emph{x}} \boldsymbol{\mid} \textbf{\emph{y}}) = 0$. 
	
\end{thm}
 
The proof of this theorem in \Cref{HFunction0} below makes use of the vertex model interpretation for $\mathcal{F}_{\boldsymbol{\lambda}} (\textbf{x}; \infty \boldsymbol{\mid} \textbf{y}; 0)$ by showing that, under the vanishing condition described in \Cref{0lambdam0}, there are no path ensembles with nonzero weight that contribute to this partition function. Based on our computer data, we were unable to find (for generic $q$ and $\boldsymbol{\lambda}$) any additional vanishing points for $\check{\mathcal{F}}_{\boldsymbol{\lambda}} (\textbf{x} \boldsymbol{\mid} \textbf{y})$ of the form $\{ x_j = q^{k_j} y_{\mathfrak{m}_j} \}$ that are not listed in \Cref{0lambdam0}. 

In general, it seems unlikely that the vanishing condition of \Cref{0lambdam0} admits an efficient description in terms of the individual coordinates $\{ \mathfrak{m}_j \}$ and $\{ \mathfrak{l}_j \}$. Already in the case $N = n = \ell = 2$, such a description is quite intricate; see \Cref{hlambda02} below for the precise statement. However, it is still possible to explicitly list a subset of these coordinates, as indicated by the example below.

\begin{example}[\Cref{h0nkappa0} below]
	
	\label{0h0nkappa0} 
	
	For $\kappa = 0^{\ell}$, \Cref{0lambdam0} implies $\check{\mathcal{F}}_{\boldsymbol{\lambda}} (\textbf{x} \boldsymbol{\mid} \textbf{y}) = 0$ when $x_j = q^{n - 1}  y_{\mathfrak{m}_j}$ for each $j \in [1, \ell]$, if there exists some $k \in [1, N]$ with $\max_{h \in [1, n]} \mathfrak{l}_{N - k}^{(h)} > \mathfrak{m}_{\ell - k}$.
	
\end{example} 

In the case $n = 1$, we can in fact explicitly list all vanishing points prescribed by \Cref{0lambdam0}.

\begin{example}[\Cref{h0n1} below]
	
	\label{0h0n1} 
	
	Assume $n = 1$, and abbreviate $\lambda = \lambda^{(1)}$ and $\mathfrak{l} = \mathfrak{l}^{(1)} = (\mathfrak{l}_1, \mathfrak{l}_2, \ldots , \mathfrak{l}_N)$. Then \Cref{0lambdam0} implies $\check{\mathcal{F}}_{\lambda} (\textbf{x} \boldsymbol{\mid} \textbf{y}) = 0$ when $x_j = y_{\mathfrak{m}_j}$ for each $j \in [1, \ell]$, if there exists some $k \in [1, n]$ with $\mathfrak{l}_{N - k} > \mathfrak{m}_{\ell - k}$. 
	
\end{example}

One of the remarkable uses of vanishing properties is that in many cases they fully characterize the underlying (non)symmetric function, a fact originally observed by Sahi \cite{SIDOSS,IIP} in the context of interpolation Macdonald polynomials. This phenomenon is useful for us as well. Indeed, the $n = 1$ case given by \Cref{0h0n1} of the vanishing property for the factorial LLT polynomials $\check{\mathcal{F}}_{\lambda} (\textbf{x} \boldsymbol{\mid} \textbf{y})$ happens to coincide with that satisfied for the factorial Schur polynomials $s_{\lambda} (\textbf{x} \boldsymbol{\mid} \textbf{y})$ (whose definition we recall as \eqref{sxy} below). The latter vanishing condition is known \cite{IIP} to fully characterize the factorial Schur polynomials (among those of total degree $|\lambda|$ that are symmetric in $\textbf{x}$), and so we deduce the following corollary equating them with the factorial LLT polynomials at $n = 1$.

\begin{cor}[\Cref{hfunctionn1} and \Cref{hf} below]
	
	\label{0hfunctionn1} 
	
	If $n = 1$, $\check{\mathcal{F}}_{\lambda} (\textbf{\emph{x}} \boldsymbol{\mid} \textbf{\emph{y}}) = s_{\lambda} (\textbf{\emph{x}} \boldsymbol{\mid} \textbf{\emph{y}})$. 
	
\end{cor} 

Let us mention that \Cref{0hfunctionn1} was also established as the $t = 0$ case of Theorem 1 in \cite{FFE}. 

The vanishing condition of \Cref{0lambdam0} unfortunately does not fully characterize the factorial LLT polynomials $\check{\mathcal{F}}_{\boldsymbol{\lambda}}$ for $n \ge 2$. Indeed, denoting $\sigma (\boldsymbol{\lambda}) = \big( \lambda^{(\sigma (1))}, \lambda^{(\sigma (2))}, \ldots , \lambda^{(\sigma (n))} \big)$ for any $\sigma \in \mathfrak{S}_n$, a marked sequence $(\mathfrak{m}, \kappa)$ satisfies the vanishing condition described in \Cref{lambdam0} with respect to $\boldsymbol{\lambda}$ if and only if it does with respect to any $\sigma (\boldsymbol{\lambda})$. Thus, the $\big\{ \check{\mathcal{F}}_{\sigma(\boldsymbol{\lambda})} (\textbf{x} \boldsymbol{\mid} \textbf{y}) \big\}$ over $\sigma \in \mathfrak{S}_n$ all satisfy the same vanishing condition but typically do not coincide up to a global factor.\footnote{A similar extensive, but non-characterizing, vanishing condition also holds for a different family of polynomials studied by Sahi--Salmasian--Serganova in \cite{OSSI}.} 

Still, one can ask whether this failure of characterization is exhaustive, in the following sense. 

\begin{que}
	
	\label{lambdam0sigma}
	
	Let $P (\textbf{x}; \textbf{y})$ denote a polynomial of total degree $|\boldsymbol{\lambda}|$, which is symmetric in $\textbf{x}$. If $P (\textbf{x}; \textbf{y}) = 0$ for any $(\textbf{x}; \textbf{y})$ satisfying the conditions of \Cref{0lambdam0}, then do there exist constants $\{ a_{\sigma} \}_{\sigma \in \mathfrak{S}_n} \subset \mathbb{C}$ such that $P (\textbf{x}; \textbf{y}) = \sum_{\sigma \in \mathfrak{S}_n} a_{\sigma} \check{\mathcal{F}}_{\sigma (\boldsymbol{\lambda})} (\textbf{x} \boldsymbol{\mid} \textbf{y})$?
	
\end{que} 	

We have not yet found a proof or counterexample to this statement.

\section{Macdonald Polynomials}

\label{0PolynomialsP}

We next provide fermionic vertex model representations for the nonsymmetric and integral Macdonald polynomials, to which end we first require some notation. A \emph{partition} is a finite, non-increasing sequence $\lambda = (\lambda_1, \lambda_2, \ldots , \lambda_{\ell}) \in \mathbb{Z}_{> 0}^{\ell}$ of positive integers. A \emph{(positive) composition} is a sequence $\mu = (\mu_1, \mu_2, \ldots , \mu_{\ell}) \in \mathbb{Z}_{> 0}^{\ell}$ of positive integers (without any constraint on its ordering), and a \emph{nonnegative composition} is a composition in which some of its entries may equal $0$. For any (positive or nonnegative) composition $\mu$, there exists a unique signature, denoted by $\mu^+$,\index{0@$\mu^+$; dominant ordering of $\mu$} obtained from $\mu$ by permuting its entries into non-increasing order; $\mu^+$ is a partition if $\mu$ is a positive composition. Moreover, for any $i \in [1, n]$ let $\textbf{e}_i \in \{ 0, 1 \}^n$\index{E@$\textbf{e}_i$} denote the $n$-tuple whose $i$-th coordinate is $1$ and whose remaining coordinates are $0$. Further let $\textbf{e}_{[1, i]} \in \{ 0, 1 \}^n$ denote the $n$-tuple whose first $i$ elements are $1$ and whose remaining ones are $0$.\index{E@$\textbf{e}_{[1, i]}$}

For any sequence $\textbf{x} = (x_1, x_2, \ldots , x_n)$, nonnegative composition $\mu$, and partition $\lambda$, let $f_{\mu} (\textbf{x})$ denote the nonsymmetric Macdonald polynomial\index{F@$f_{\mu} (\textbf{x})$; nonsymmetric Macdonald polynomial} and $J_{\lambda} (\textbf{x})$ denote the integral Macdonald polynomial,\index{J@$J_{\lambda} (\textbf{x})$; integral Macdonald polynomial} recalled in \Cref{Polynomialf} and \Cref{SumJ} below, respectively.  Both of these functions depend on two parameters $q, t \in \mathbb{C}$, which are fixed throughout this section. Let us clarify that it will be the Macdonald parameter $t$ (and not $q$) that will correspond to the quantization parameter for our vertex models. 

Next, we require the vertex weights for our model, which are obtained as the $r = q^{-1 / 2}$ specializations of the top-right limiting weights depicted in \Cref{wlimit}. In this case, it is quickly verified (see \Cref{rql1} below) that arrow configurations admitting more than one arrow along a horizontal edge are assigned weight $0$. Thus, any arrow configuration with nonzero weight must be of the form $(\textbf{A}, \textbf{e}_b; \textbf{C}, \textbf{e}_d)$ for some $b, d \in [0, n]$, and so we will abbreviate them by $(\textbf{A}, b; \textbf{C}, d)$. 

These limiting vertex weights, which we denote by $L_x (\textbf{A}, b; \textbf{C}, d) = L_{x, t} (\textbf{A}, b; \textbf{C}, d)$ are then given diagrammatically by
\begin{center}
	
	\begin{tikzpicture}[
		>=stealth,
		scale = .75
		]	
		\draw[-, black] (-7.5, 3.1) -- (10, 3.1);
		\draw[-, black] (-7.5, -2.1) -- (10, -2.1);
		\draw[-, black] (-7.5, -1.1) -- (10, -1.1);
		\draw[-, black] (-7.5, -.4) -- (10, -.4);
		\draw[-, black] (-7.5, 2.4) -- (10, 2.4);
		\draw[-, black] (-7.5, -2.1) -- (-7.5, 3.1);
		\draw[-, black] (10, -2.1) -- (10, 3.1);
		\draw[-, black] (7.5, -2.1) -- (7.5, 3.1);
		\draw[-, black] (-5, -2.1) -- (-5, 2.4);
		\draw[-, black] (5, -2.1) -- (5, 3.1);
		\draw[-, black] (-2.5, -2.1) -- (-2.5, 2.4);
		\draw[-, black] (2.5, -2.1) -- (2.5, 2.4);
		\draw[-, black] (0, -2.1) -- (0, 3.1);
		\draw[->, thick, blue] (-6.3, .1) -- (-6.3, 1.9); 
		\draw[->, thick, green] (-6.2, .1) -- (-6.2, 1.9); 
		\draw[->, thick, blue] (-3.8, .1) -- (-3.8, 1) -- (-2.85, 1);
		\draw[->, thick, green] (-3.7, .1) -- (-3.7, 1.9);
		\draw[->, thick, blue] (-1.35, .1) -- (-1.35, 1.9);
		\draw[->, thick, green] (-1.25, .1) -- (-1.25, 1.9);
		\draw[->, thick,  orange] (-2.15, 1.1) -- (-1.15, 1.1) -- (-1.15, 1.9);
		
		\draw[->, thick, red] (.35, 1) -- (1.15, 1) -- (1.15, 1.9);
		\draw[->, thick, blue] (1.25, .1) -- (1.25, 1.9);
		\draw[->, thick, green] (1.35, .1) -- (1.35, 1.1) -- (2.15, 1.1);
		\draw[->, thick, blue] (3.65, .1) -- (3.65, 1) -- (4.65, 1);
		\draw[->, thick, green] (3.75, .1) -- (3.75, 1.9);
		\draw[->, thick, orange] (2.85, 1.1) -- (3.85, 1.1) -- (3.85, 1.9); 
		\draw[->, thick, red] (5.35, 1) -- (7.15, 1); 
		\draw[->, thick, blue] (6.2, .1) -- (6.2, 1.9);
		\draw[->, thick, green] (6.3, .1) -- (6.3, 1.9); 
		\draw[->, thick, blue] (7.85, 1) -- (9.65, 1); 
		\draw[->, thick, blue] (8.7, .1) -- (8.7, 1.9);
		\draw[->, thick, green] (8.8, .1) -- (8.8, 1.9); 
		\filldraw[fill=gray!50!white, draw=black] (-2.85, 1) circle [radius=0] node [black, right = -1, scale = .7] {$i$};
		\filldraw[fill=gray!50!white, draw=black] (2.15, 1) circle [radius=0] node [black, right = -1, scale = .7] {$j$};
		\filldraw[fill=gray!50!white, draw=black] (4.65, 1) circle [radius=0] node [black, right = -1, scale = .7] {$i$};
		\filldraw[fill=gray!50!white, draw=black] (7.15, 1) circle [radius=0] node [black, right = -1, scale = .7] {$i$};
		\filldraw[fill=gray!50!white, draw=black] (9.65, 1) circle [radius=0] node [black, right = -1, scale = .7] {$i$};
		\filldraw[fill=gray!50!white, draw=black] (7.85, 1) circle [radius=0] node [black, left = -1, scale = .7] {$i$};
		\filldraw[fill=gray!50!white, draw=black] (5.35, 1) circle [radius=0] node [black, left = -1, scale = .7] {$i$};
		\filldraw[fill=gray!50!white, draw=black] (2.85, 1) circle [radius=0] node [black, left = -1, scale = .7] {$j$};
		\filldraw[fill=gray!50!white, draw=black] (.35, 1) circle [radius=0] node [black, left = -1, scale = .7] {$i$};
		\filldraw[fill=gray!50!white, draw=black] (-2.15, 1) circle [radius=0] node [black, left = -1, scale = .7] {$i$};
		\filldraw[fill=gray!50!white, draw=black] (-6.25, 1.9) circle [radius=0] node [black, above = -1, scale = .7] {$\textbf{A}$};
		\filldraw[fill=gray!50!white, draw=black] (-3.75, 1.9) circle [radius=0] node [black, above = -1, scale = .7] {$\textbf{A}_i^-$};
		\filldraw[fill=gray!50!white, draw=black] (-1.25, 1.9) circle [radius=0] node [black, above = -1, scale = .7] {$\textbf{A}_i^+$};
		\filldraw[fill=gray!50!white, draw=black] (1.25, 1.9) circle [radius=0] node [black, above = -1, scale = .65] {$\textbf{A}_{ij}^{+-}$};
		\filldraw[fill=gray!50!white, draw=black] (3.75, 1.9) circle [radius=0] node [black, above = -1, scale = .65] {$\textbf{A}_{ji}^{+-}$};
		\filldraw[fill=gray!50!white, draw=black] (6.25, 1.9) circle [radius=0] node [black, above = -1, scale = .7] {$\textbf{A}$};
		\filldraw[fill=gray!50!white, draw=black] (8.75, 1.9) circle [radius=0] node [black, above = -1, scale = .7] {$\textbf{A}$};
		\filldraw[fill=gray!50!white, draw=black] (-6.25, .1) circle [radius=0] node [black, below = -1, scale = .7] {$\textbf{A}$};
		\filldraw[fill=gray!50!white, draw=black] (-3.75, .1) circle [radius=0] node [black, below = -1, scale = .7] {$\textbf{A}$};
		\filldraw[fill=gray!50!white, draw=black] (-1.25, .1) circle [radius=0] node [black, below = -1, scale = .7] {$\textbf{A}$};
		\filldraw[fill=gray!50!white, draw=black] (1.25, .1) circle [radius=0] node [black, below = -1, scale = .7] {$\textbf{A}$};
		\filldraw[fill=gray!50!white, draw=black] (3.75, .1) circle [radius=0] node [black, below = -1, scale = .7] {$\textbf{A}$};
		\filldraw[fill=gray!50!white, draw=black] (6.25, .1) circle [radius=0] node [black, below = -1, scale = .7] {$\textbf{A}$};
		\filldraw[fill=gray!50!white, draw=black] (8.75, .1) circle [radius=0] node [black, below = -1, scale = .7] {$\textbf{A}$};
		\filldraw[fill=gray!50!white, draw=black] (-3.75, 2.75) circle [radius=0] node [black] {$1 \le i \le n$};
		\filldraw[fill=gray!50!white, draw=black] (2.5, 2.75) circle [radius=0] node [black] {$1 \le i < j \le n$}; 
		\filldraw[fill=gray!50!white, draw=black] (6.25, 2.75) circle [radius=0] node [black] {$A_i = 0$};
		\filldraw[fill=gray!50!white, draw=black] (8.75, 2.75) circle [radius=0] node [black] {$A_i = 1$};
		\filldraw[fill=gray!50!white, draw=black] (-6.25, -.75) circle [radius=0] node [black, scale = .8] {$(\textbf{A}, 0; \textbf{A}, 0)$};
		\filldraw[fill=gray!50!white, draw=black] (-3.75, -.75) circle [radius=0] node [black, scale = .8] {$\big( \textbf{A}, 0; \textbf{A}_i^-, i \big)$};
		\filldraw[fill=gray!50!white, draw=black] (-1.25, -.75) circle [radius=0] node [black, scale = .8] {$\big( \textbf{A}, i; \textbf{A}_i^+, 0 \big)$};
		\filldraw[fill=gray!50!white, draw=black] (1.25, -.75) circle [radius=0] node [black, scale = .8] {$\big( \textbf{A}, i; \textbf{A}_{ij}^{+-}, j \big)$};
		\filldraw[fill=gray!50!white, draw=black] (3.75, -.75) circle [radius=0] node [black, scale = .8] {$\big( \textbf{A}, j; \textbf{A}_{ji}^{+-}, i \big)$};
		\filldraw[fill=gray!50!white, draw=black] (6.25, -.75) circle [radius=0] node [black, scale = .8] {$(\textbf{A}, i; \textbf{A}, i)$};
		\filldraw[fill=gray!50!white, draw=black] (8.75, -.75) circle [radius=0] node [black, scale = .8] {$(\textbf{A}, i; \textbf{A}, i)$};
		\filldraw[fill=gray!50!white, draw=black] (-6.25, -1.6) circle [radius=0] node [black, scale = .8] {$1$};
		\filldraw[fill=gray!50!white, draw=black] (-3.75, -1.6) circle [radius=0] node [black, scale = .75] {$t^{A_{[i + 1, n]}} x (1 - t)$};
		\filldraw[fill=gray!50!white, draw=black] (-1.25, -1.6) circle [radius=0] node [black, scale = .8] {$1$};
		\filldraw[fill=gray!50!white, draw=black] (1.25, -1.6) circle [radius=0] node [black, scale = .75] {$t^{A_{[j + 1, n]}} x (1 - t)$};
		\filldraw[fill=gray!50!white, draw=black] (3.75, -1.6) circle [radius=0] node [black, scale = .8] {$0$};
		\filldraw[fill=gray!50!white, draw=black] (6.25, -1.6) circle [radius=0] node [black, scale = .8] {$t^{A_{[i, n]}} x$};
		\filldraw[fill=gray!50!white, draw=black] (8.75, -1.6) circle [radius=0] node [black, scale = .8] {$- t^{A_{[i, n]}} x$};
	\end{tikzpicture}
\end{center}
\index{L@$L_x (\textbf{A}, b; \textbf{C}, d)$}

\noindent where arrow configurations not of the above form are assigned weight $0$, and we have denoted
\begin{flalign}
	\label{0aij} 
	\textbf{A}_i^+ = \textbf{A} + \textbf{e}_i; \qquad \textbf{A}_j^- = \textbf{A} - \textbf{e}_j; \qquad \textbf{A}_{ij}^{+-} = \textbf{A} + \textbf{e}_i - \textbf{e}_j.
\end{flalign} 
\index{A@$\textbf{A}_-^+, \textbf{A}_j^-, \textbf{A}_{ij}^{+-}$}

\noindent Let us mention that almost all of the $L_x$ vertex weights shown above coincide with those given by equation (3.7) of \cite{NPIVM}, which were used to provide bosonic partition function formulas for the nonsymmetric Macdonald polynomials. The one exceptional weight is the rightmost one (namely, $L_x (\textbf{A}, i; \textbf{A}, i)$ for $A_i = 1$), which differs from its counterpart in \cite{NPIVM}. This contrasting vertex weight encapsulates the difference between the bosonic model used in \cite{NPIVM} and the fermionic one we will use here to access the Macdonald polynomials; it will account for the finiteness of our partition function formulas here, and it will also enable us to compare these formulas with the LLT ones  given by \Cref{0limitg0}.   

Next, for any nonnegative composition $\mu$ of length $n$, we introduce the following infinite sequence\footnote{At certain points in this text (particularly when we consider only positive compositions), we might start the indexing for $\mathscr{I} (\mu)$ at $1$ instead of $0$, but we will make a point about such changes in convention whenever they arise.} $\mathscr{I} (\mu) = \big( \textbf{I}_0 (\mu), \textbf{I}_1 (\mu), \ldots \big)$ of elements in $\{ 0, 1 \}^n$. For each $j \ge 0$, define $\textbf{I}_j = \textbf{I}_j (\mu) = (I_{1, j}, I_{2, j}, \ldots , I_{n, j})$ so that $I_{k, j} = \textbf{1}_{\mu_k = j}$\index{I@$\mathscr{I} (\mu)$} for every $k \in [1, n]$. For any $i \in [1, n]$ and $j \ge 0$, we further define the \emph{twist parameters} 
\begin{flalign*}
v_{i, j} (\mu) = q^{\mu_i - j} t^{\gamma_{i, j} (\mu)} \textbf{1}_{\mu_i > j},
\end{flalign*}
\index{V@$v_{i, j} (\mu)$; twist parameters}

\noindent where the exponents $\gamma_{i, j} (\mu)$ are given by 
\begin{align*}
	\gamma_{i,j}(\mu)
	=
	-\#\{k < i : \mu_k > \mu_i\} + \#\{k > i : j \le \mu_k < \mu_i\}.
\end{align*}
\index{0@$\gamma_{i, j} (\mu)$}

\noindent We also define the normalization constant 
	\begin{align*}
	\Omega_{\mu}(q,t)
	=
	\prod_{i=1}^{n}
	\prod_{j=0}^{\mu_i-1}
	\frac{1}{1-q^{\mu_i-j} t^{\alpha_{i,j}(\mu)+1}},
\end{align*}
\index{0@$\Omega_{\mu} (q, t)$}

\noindent where the exponents $\alpha_{i, j} (\mu)$ are given by
\begin{align*}
	\alpha_{i,j}(\mu)
	=
	\#\{k < i: \mu_k = \mu_i\}
	+
	\#\{k \not= i: j < \mu_k < \mu_i\}
	+
	\#\{k > i: j = \mu_k\}.
\end{align*}
\index{0@$\alpha_{i, j} (\mu)$}

Under this notation, we have the following partition function formula for the nonsymmetric Macdonald polynomials, expressed diagrammatically.

\begin{thm}[\Cref{thm:f-formula} below]
	
	\label{0thm:f-formula}

	Denoting $N = \max_{i \in [1, n]} \mu_i$ and $v_{i, j} = v_{i, j} (\mu)$, we have
\begin{align}
	\label{0vertexfmu}
	\tikz{1}{
		\draw[->, thick, red] (3.95, 4) -- (3.95, 5); 
		\draw[->, thick, blue] (5, 4) -- (5, 5);
		\draw[->, thick, orange] (1, 4) -- (1, 5);
		\draw[->, thick, blue] (2, 4) -- (2, 5);
		\draw[->, thick, blue] (2, 0) -- (2, .95);
		\draw[->, thick, green] (4.05, 4) -- (4.05, 5);
		\draw[ultra thick, gray, dashed] (1, 1) -- (1, 4);
		\draw[ultra thick, gray, dashed] (2, 1) -- (2, 4);
		\draw[ultra thick, gray, dashed] (3, 1) -- (3, 4);
		\draw[ultra thick, gray, dashed] (4, 1) -- (4, 4);
		\draw[ultra thick, gray, dashed] (5, 1) -- (5, 4);
		\draw[->, thick, red] (0, 1) node[black, left, scale = .9]{$x_1$} -- node[black, above, scale = .75]{$1$} (.95, 1);
		\draw[->, thick, blue] (0, 2) node[black, left, scale = .9]{$x_2$} -- node[black, above, scale = .75]{$2$} (.95, 2);
		\draw[->, thick, green] (0, 3) node[black, left, scale = .9]{$\vdots$} -- node[black, above, scale = .75]{$\vdots$} (.95, 3);
		\draw[->, thick, orange] (0, 4) node[black, left, scale = .9]{$x_n$} -- node[black, above, scale = .75]{$n$} (.95, 4);
		\draw[ultra thick, gray, dashed] (1, 1) -- (5, 1) node[black, right, scale = .75]{$0$};
		\draw[ultra thick, gray, dashed] (1, 2) -- (5, 2) node[black, right, scale = .75]{$0$};
		\draw[ultra thick, gray, dashed] (1, 3) -- (5, 3) node[black, right, scale = .75]{$\vdots$};
		\draw[ultra thick, gray, dashed] (1, 4) -- (5, 4) node[black, right, scale = .75]{$0$};
		\draw[->, very thick] (0, 0) -- (0, 5.5);
		\draw[->, very thick] (0, 0) -- (6, 0);
		\draw[] (1, -.15) circle [radius = 0] node[below, scale = .75]{$\textbf{\emph{M}}_0$};
		\draw[] (2, -.15) circle [radius = 0] node[below, scale = .75]{$\textbf{\emph{M}}_1$};
		\draw[] (3, -.15) circle [radius = 0] node[below, scale = .75]{$\textbf{\emph{M}}_2$};
		\draw[] (4, -.15) circle [radius = 0] node[below, scale = .85]{$\cdots$};
		\draw[] (5, -.15) circle [radius = 0] node[below, scale = .75]{$\textbf{\emph{M}}_N$};
		\draw[] (1, 5.15) circle [radius = 0] node[above, scale = .65]{$\textbf{\emph{M}}_0 + \textbf{\emph{I}}_0$};
		\draw[] (2, 5.15) circle [radius = 0] node[above, scale = .65]{$\textbf{\emph{M}}_1 + \textbf{\emph{I}}_1$};
		\draw[] (3, 5.15) circle [radius = 0] node[above, scale = .65]{$\textbf{\emph{M}}_2 + \textbf{\emph{I}}_2$};
		\draw[] (4, 5.15) circle [radius = 0] node[above, scale = .85]{$\cdots$};
		\draw[] (5, 5.15) circle[radius = 0] node[above, scale = .65]{$\textbf{\emph{M}}_N + \textbf{\emph{I}}_N$};
		\draw[] (-3.65, 1.5) circle[radius = 0] node[above, scale = 1.1]{$f_{\mu} (\textbf{\emph{x}}) = \Omega_{\mu} (q, t) \displaystyle\sum_{\mathscr{M}} \displaystyle\prod_{i = 1}^n \displaystyle\prod_{j = 1}^N v_{i, j}^{M_{i, j}} \times $};
	}
\end{align}

\noindent under the vertex weights $L_{x_j, t} (\textbf{\emph{A}}, b; \textbf{\emph{C}}, d)$ in row $j$. Here the sum is over all length $N$ sequences $\mathscr{M} = (\textbf{\emph{M}}_0, \textbf{\emph{M}}_1, \ldots , \textbf{\emph{M}}_N)$ of elements in $\{ 0, 1 \}^n$, with $\textbf{\emph{M}}_j = (M_{1, j}, M_{2, j}, \ldots , M_{n, j})$, such that $\textbf{\emph{I}}_k (\mu) + \textbf{\emph{M}}_k \in \{ 0, 1 \}^n$ for each $k \ge 0$. Here, we have abbreviated $\textbf{\emph{I}}_k = \textbf{\emph{I}}_k (\mu)$ for each $k \ge 0$. 
\end{thm} 

Up to the overall factor of $\Omega_{\mu} (q, t)$, one may interpret the right side of \eqref{0vertexfmu} as the partition function for our vertex model (under the $L_x$ weights) on the cylinder obtained by identifying the top and bottom boundaries of the strip $\mathbb{Z}_{> 0} \times [1, n]$. The sum over $\mathscr{M}$ allows our paths the option of ``wrapping'' around the cylinder at most once in each column; whenever a path of color $i$ does this in column $j$, it contributes the twist parameter $v_{i, j}$ to the partition function. 

The next theorem states that the integral, symmetric Macdonald polynomials admit similar partition function formulas, but with slightly different boundary data. Here, a composition $\nu = (\nu_1, \nu_2, \ldots , \nu_p)$ is \emph{anti-dominant} if $\nu_1 \le \nu_2 \le \cdots \le \nu_p$.

\begin{thm}[\Cref{jpd} below]
	
	\label{0jpd} 
	
	Fix an integer $p \in [1, n]$ and an anti-dominant, positive composition $\nu = (\nu_1, \nu_2, \ldots , \nu_p)$ of length $p$. Denoting $N = \max_{i \in [1, n]} \nu_i$ and\footnote{Since $\nu$ is anti-dominant, the exponents $\gamma_{i, j}$ under this specialization are all equal to $0$.} $w_{i, j} = v_{i, j} (\nu) = q^{\nu_i - j} \textbf{\emph{1}}_{\nu_i > j}$, we have
\begin{align}
	\label{0vertexjnu} 
	\tikz{1}{
		\draw[->, thick, red] (3.95, 4) -- (3.95, 5); 
		\draw[->, thick, green] (5, 4) -- (5, 5);
		\draw[->, thick, blue] (3, 4) -- (3, 5);
		\draw[->, thick, blue] (3, 0) -- (3, .95);
		\draw[->, thick, blue] (4.05, 4) -- (4.05, 5);
		\draw[ultra thick, gray, dashed] (1, 1) -- (1, 4);
		\draw[ultra thick, gray, dashed] (2, 1) -- (2, 4);
		\draw[ultra thick, gray, dashed] (3, 1) -- (3, 4);
		\draw[ultra thick, gray, dashed] (4, 1) -- (4, 4);
		\draw[ultra thick, gray, dashed] (5, 1) -- (5, 4);
		\draw[] (0, 1) circle [radius = 0] node[left, scale = .9]{$x_1$};
		\draw[] (0, 2) circle [radius = 0] node[left, scale = .9]{$x_2$};
		\draw[] (0, 3) circle [radius = 0] node[left, scale = .9]{$\vdots$};
		\draw[] (0, 4) circle [radius = 0] node[left, scale = .9]{$x_n$};
		\draw[->, red, thick] (.925, 0) -- (.925, .95);
		\draw[->, blue, thick] (1, 0) -- (1, .95);
		\draw[->, green, thick] (1.075, 0) -- (1.075, .95);
		\draw[ultra thick, gray, dashed] (1, 1) node[black, left, scale = .75]{$0$} -- (5, 1) node[black, right, scale = .75]{$0$};
		\draw[ultra thick, gray, dashed] (1, 2) node[black, left, scale = .75]{$0$} -- (5, 2) node[black, right, scale = .75]{$0$};
		\draw[ultra thick, gray, dashed] (1, 3) node[black, left, scale = .75]{$\vdots$} -- (5, 3) node[black, right, scale = .75]{$\vdots$};
		\draw[ultra thick, gray, dashed] (1, 4) node[black, left, scale = .75]{$0$} -- (5, 4) node[black, right, scale = .75]{$0$};
		\draw[->, very thick] (0, 0) -- (0, 5.5);
		\draw[->, very thick] (0, 0) -- (6, 0);
		\draw[] (1, -.15) circle [radius = 0] node[below, scale = .75]{$\textbf{\emph{e}}_{[1, p]}$};
		\draw[] (2, -.15) circle [radius = 0] node[below, scale = .75]{$\textbf{\emph{M}}_1$};
		\draw[] (3, -.15) circle [radius = 0] node[below, scale = .75]{$\textbf{\emph{M}}_2$};
		\draw[] (4, -.15) circle [radius = 0] node[below, scale = .85]{$\cdots$};
		\draw[] (5, -.15) circle [radius = 0] node[below, scale = .75]{$\textbf{\emph{M}}_N$};
		\draw[] (1, 5.15) circle [radius = 0] node[above, scale = .75]{$\textbf{\emph{e}}_0$};
		\draw[] (2, 5.15) circle [radius = 0] node[above, scale = .65]{$\textbf{\emph{M}}_1 + \textbf{\emph{I}}_1$};
		\draw[] (3, 5.15) circle [radius = 0] node[above, scale = .65]{$\textbf{\emph{M}}_2 + \textbf{\emph{I}}_2$};
		\draw[] (4, 5.15) circle [radius = 0] node[above, scale = .85]{$\cdots$};
		\draw[] (5, 5.15) circle[radius = 0] node[above, scale = .65]{$\textbf{\emph{M}}_N + \textbf{\emph{I}}_N$};
		\draw[] (-3.15, 1.5) circle[radius = 0] node[above, scale = 1.1]{$J_{\nu^+} (\textbf{\emph{x}}) = \displaystyle\sum_{\mathscr{M}} \displaystyle\prod_{i = 1}^n \displaystyle\prod_{j = 1}^N w_{i, j}^{M_{i, j}} \times $};
	}
\end{align}

\noindent under the vertex weights $L_{x_j, t} (\textbf{\emph{A}}, b; \textbf{\emph{C}}, d)$ in row $j$. Here the sum is over all length $N$ sequences $\mathscr{M} = (\textbf{\emph{M}}_1, \textbf{\emph{M}}_2, \ldots , \textbf{\emph{M}}_N)$ of elements in $\{ 0, 1 \}^n$, with $\textbf{\emph{M}}_j = (M_{1, j}, M_{2, j}, \ldots , M_{n, j})$, such that $\textbf{\emph{I}}_k (\nu) + \textbf{\emph{M}}_k \in \{ 0, 1 \}^n$ for each $k \ge 1$. Here, we have abbreviated $\textbf{\emph{I}} = \textbf{\emph{I}}_k (\nu)$ for each $k \ge 1$. 

\end{thm}

As before, to the right of the first column, one may interpret the vertex model \eqref{0vertexjnu} as existing on a cylinder. However, the first column behaves differently; paths of colors $\{ 1, 2, \ldots , p \}$ vertically enter through it, and none vertically exit through it. In particular, no paths of colors in $\{ p + 1, p + 2, \ldots , n \}$ exist in this model. 

By comparing the vertex model representations \Cref{0limitg0} and \Cref{0jpd} for the LLT and Macdonald polynomials, respectively, we deduce the following expression for the modified Macdonald polynomials $\widetilde{J}_{\lambda} (\textbf{x})$ (recalled in \Cref{SumJ} below) as linear combinations of skew LLT ones.\index{J@$J_{\lambda} (\textbf{x})$; integral Macdonald polynomial!$\widetilde{J}_{\lambda} (\textbf{x})$; modified Macdonald polynomial} It was originally established as Theorem 2.2, equation (23), and Proposition 3.4 of \cite{CP}. In the below, we recall $\mathscr{S}$ from \Cref{0Functions}. 
	
\begin{cor}[\Cref{psuml} below]

\label{0psuml}

Let $\nu$ denote an anti-dominant composition of length $p \le n$. For each $i \in [1, n]$ and $j \ge 1$, define
\begin{flalign*}
	u_{i, j} = u_{i, j} (\nu) = q^{\nu_i - j} t^{\beta_{i, j} (\nu)} \textbf{\emph{1}}_{\nu_i > j},
\end{flalign*}

\noindent where $\nu_k = 0$ for $k > n$, and the exponents $\beta_{i, j} (\nu)$ are given by 
\begin{align*}
	\beta_{i,j} (\nu) = \displaystyle\frac{1}{2} \Big( \#\{ k < i: \nu_k > j \} - p + i \Big).
\end{align*}

\noindent For any infinite sequence $\mathscr{K} = (\textbf{\emph{K}}_1, \textbf{\emph{K}}_2, \ldots )$ of elements in $\{ 0, 1 \}^n$, define $\boldsymbol{\mu} (\mathscr{K}) \in \SeqSign_n$ so that $\mathscr{S} \big( \boldsymbol{\mu} (\mathscr{K}) \big) = \mathscr{K}$. Then, 
\begin{flalign}
	\label{0hsuml} 
\widetilde{J}_{\nu^+} (\textbf{\emph{x}}) = \displaystyle\sum_{\mathscr{M}} \mathcal{L}_{\boldsymbol{\mu} (\mathscr{M}_0 + \mathscr{I}_0 (\nu)) / \boldsymbol{\mu} (\mathscr{M}_p)} (\textbf{\emph{x}}; q) \displaystyle\prod_{i = 1}^n \displaystyle\prod_{j = 1}^{\infty} u_{i, j}^{M_{i, j}},
\end{flalign}

\noindent where the sum is over all infinite sequences $\mathscr{M} = (\textbf{\emph{M}}_1, \textbf{\emph{M}}_2, \ldots )$ of elements in $\{ 0, 1 \}^n$, with $\textbf{\emph{M}}_j = (M_{1, j}, M_{2, j}, \ldots , M_{n, j})$, such that $\textbf{\emph{I}}_k (\nu) + \textbf{\emph{M}}_k \in \{ 0, 1 \}^n$ for each $k \ge 1$. Here, we have denoted $\mathscr{M}_0 = (\textbf{\emph{e}}_0, \textbf{\emph{M}}_1, \textbf{\emph{M}}_2, \ldots )$; $\mathscr{I}_0 (\nu) = \big( \textbf{\emph{e}}_0, \textbf{\emph{I}}_1 (\nu), \textbf{\emph{I}}_2 (\nu), \ldots \big)$; and $\mathscr{M}_p = \big( \textbf{\emph{e}}_{[1, p]}, \textbf{\emph{M}}_1, \textbf{\emph{M}}_2, \ldots \big)$.

\end{cor} 

Here, the reason for introducing the signature sequence $\boldsymbol{\mu} (\mathscr{K})$ is to account for the different indexing of the boundary data for the LLT and Macdonald vertex models; the former implements the shift $\mathfrak{T}$ from \eqref{t} (using $\mathscr{S}$), and the latter does not (using $\mathscr{I}$). Let us further mention that the sum on the right side of \eqref{0hsuml} is in fact finite, since $u_{i, j} = 0$ for $j \ge \nu_p$.

\section{Transition Coefficients}

\label{0Coefficients}

Here we provide vertex model interpretations for the transition coefficients from both the LLT and modified Macdonald polynomials to the modified Hall--Littlewood basis of symmetric functions. In what follows, we fix Macdonald parameters $q, t \in \mathbb{C}$, and we recall that the LLT, modified Macdonald, and modified Hall--Littlewood polynomials are denoted by $\mathcal{L}_{\boldsymbol{\lambda} / \boldsymbol{\mu}} (\textbf{x}; q)$, $\widetilde{J}_{\lambda / \mu} (\textbf{x})$, and $Q_{\lambda}' (\textbf{x})$, respectively. Then, for a signature $\lambda \in \Sign$ and signature sequences $\boldsymbol{\lambda}, \boldsymbol{\mu} \in \SeqSign_{n; M}$, the above mentioned transition coefficients are implicitly defined by 
\begin{flalign}
	\label{0qf} 
q^{\psi (\boldsymbol{\lambda}) - \psi (\boldsymbol{\mu})} \mathcal{L}_{\boldsymbol{\lambda} / \boldsymbol{\mu}} (\textbf{x}; q) = \displaystyle\sum_{\nu} f_{\boldsymbol{\lambda} / \boldsymbol{\mu}}^{\nu} (q) Q_{\nu}' (\textbf{x}); \qquad 
\widetilde{J}_{\lambda} (\textbf{x}) = \displaystyle\sum_{\nu} g_{\lambda}^{\nu} (q, t) Q_{\nu}' (\textbf{x}), 
\end{flalign}
\index{F@$f_{\boldsymbol{\lambda} / \boldsymbol{\mu}}^{\nu} (q)$}\index{G@$g_{\lambda}^{\nu} (q, t)$}

\noindent where we recall $\psi$ from \eqref{0psi}, and both sums are over all partitions $\nu$. Our reason for introducing the powers of $q$ on the left side of \eqref{0qf} is to ensure that the transition coefficients $f_{\boldsymbol{\lambda} / \boldsymbol{\mu}}^{\nu} (q)$ are polynomials in $q$ (instead of in only $q^{1 / 2}$). 

In this section we provide combinatorial formulas for these coefficients $f_{\boldsymbol{\lambda} / \boldsymbol{\mu}}^{\nu}$ and $g_{\lambda}^{\nu}$ as partition functions for vertex models under certain weights. In order to state this result, we introduce two new types of vertex weights. Here, we depict these vertices as tiles instead of as the intersection between two lines; this will enable us to more visibly distinguish the two families of weights by shading the tiles under the second family and not shading them under the first. 

\begin{definition}[\Cref{ssec:weights1} below] 
	
\label{0weights1}
 
Fix nonnegative integers $\mathfrak{a}, \mathfrak{b}, \mathfrak{c}, \mathfrak{d} \ge 0$ and $n$-tuples $\textbf{A}, \textbf{B}, \textbf{C}, \textbf{D} \in \{ 0, 1 \}^n$. Define $\textbf{V} = (V_1,\dots,V_n) \in \{ 0, 1 \}^n$ with $V_i = \min\{A_i,B_i,C_i,D_i\}$ for each $i \in [1, n]$. We introduce the lattice weights
\begin{align}
	\begin{aligned}
		\label{0weights-light}
		\tikz{0.6}{
			\draw[lgray,line width=1.5pt] (-1,-1) -- (1,-1) -- (1,1) -- (-1,1) -- (-1,-1);
			\node[left] at (-1,0) {\tiny $(\mathfrak{b},\textbf{B})$};\node[right] at (1,0) {\tiny $(\mathfrak{d},\textbf{D})$};
			\node[below] at (0,-1) {\tiny $(\mathfrak{a},\textbf{A})$};\node[above] at (0,1) {\tiny $(\mathfrak{c},\textbf{C})$};
		}
		& =
		(-1)^{\mathfrak{c}+|\textbf{V}|} q^{\chi}
		\frac{(q;q)_{\mathfrak{b}-|\textbf{B}|}}{(q;q)_{\mathfrak{d}-|\textbf{D}|}}
		\frac{(q^{\mathfrak{b}-|\textbf{B}|-\mathfrak{a}+1};q)_\mathfrak{c}}{(q;q)_\mathfrak{c}}
		{\bm 1}_{|\textbf{B}| \le \mathfrak{b}}
		{\bm 1}_{|\textbf{D}| \le \mathfrak{d}}
		{\bm 1}_{\mathfrak{a}  +\mathfrak{d} = \mathfrak{b} + \mathfrak{c}}
		\\
		& \qquad \times
		{\bm 1}_{\textbf{A}+\textbf{B}=\textbf{C}+\textbf{D}}
		\prod_{j:B_j - D_j=1}
		\big( 1-q^{\mathfrak{d}-B_{[j+1,n]}-D_{[1, j - 1]}} \big),
	\end{aligned} 
\end{align}
where the exponent $\chi \equiv \chi(\mathfrak{a}, \mathfrak{b}, \mathfrak{c}, \mathfrak{d};\textbf{A},\textbf{B},\textbf{C},\textbf{D})$ is given by
\begin{align*}
	\chi
	=
	\binom{\mathfrak{d}-|\textbf{D}|}{2}
	+
	\binom{\mathfrak{c}+1}{2}
	-
	\big( \mathfrak{c}+|\textbf{C}| \big)\mathfrak{d}
	+
	|\textbf{V}| \big( \mathfrak{d}-|\textbf{D}|+1 \big)
	+
	\varphi(\textbf{D},\textbf{C})+\varphi(\textbf{V},\textbf{D}-\textbf{B}).
\end{align*}

\end{definition} 

Observe that there are two quadruples we consider here. The first is the quadruple $(\textbf{A}, \textbf{B}; \textbf{C}, \textbf{D})$ of elements in $\{ 0, 1 \}^n$, which satisfies $\textbf{A} + \textbf{B} = \textbf{C} + \textbf{D}$. We once again view $\textbf{A}$, $\textbf{B}$, $\textbf{C}$, and $\textbf{D}$ as indexing the up-right directed fermionic paths passing through the south, west, north, and east boundaries of a tile, respectively. The second is the quadruple $(\mathfrak{a}, \mathfrak{b}; \mathfrak{c}, \mathfrak{d})$ of integers, which instead satisfies $\mathfrak{a} + \mathfrak{d} = \mathfrak{b} + \mathfrak{c}$. In this way, we may view $\mathfrak{a}$, $\mathfrak{b}$, $\mathfrak{c}$, and $\mathfrak{d}$ as counting the numbers of directed down-right paths of the same bosonic color that pass through the bottom, left, top, and right boundaries of a tile, respectively. We provide depictions of these vertices in \Cref{0lambdamun2} below.

We will require a further set of weights, which are distinguished from those in \eqref{0weights-light} by means of shading. 

\begin{definition}[\Cref{ssec:weights2} below]
	
	\label{0weights2} 
	
 	Fix nonnegative integers $\mathfrak{a}, \mathfrak{b}, \mathfrak{c}, \mathfrak{d} \ge 0$ and $n$-tuples $\textbf{A}, \textbf{B}, \textbf{C}, \textbf{D} \in \{ 0, 1 \}^n$ with coordinates indexed by $[1, n]$. We define the weights
\begin{align}
	\label{0weights-dark}
	\tikz{0.6}{
		\filldraw[lgray,line width=1.5pt,fill=llgray] (-1,-1) -- (1,-1) -- (1,1) -- (-1,1) -- (-1,-1);
		\node[left] at (-1,0) {\tiny $(\mathfrak{b},\textbf{B})$};\node[right] at (1,0) {\tiny $(\mathfrak{d},\textbf{D})$};
		\node[below] at (0,-1) {\tiny $(\mathfrak{a},\textbf{A})$};\node[above] at (0,1) {\tiny $(\mathfrak{c},\textbf{C})$};
	}
	=
	{\bm 1}_{|\textbf{B}| = \mathfrak{b} \le 1}
	{\bm 1}_{|\textbf{D}| \le \mathfrak{d} \le 1}
	{\bm 1}_{\mathfrak{a} + \mathfrak{d} = \mathfrak{b} + \mathfrak{c}}
	{\bm 1}_{\textbf{A}+\textbf{B}=\textbf{C}+\textbf{D}}
	\cdot
	W \big(\mathfrak{a},\mathfrak{b},\mathfrak{c},\mathfrak{d};\textbf{A},\textbf{B},\textbf{C},\textbf{D}\big),
\end{align}
where (recalling $\textbf{A}_i^+, \textbf{A}_j^-, \textbf{A}_{ij}^{+-}$ from \eqref{0aij}) the final function appearing in \eqref{0weights-dark} is given by
\begin{flalign}
	\begin{aligned}
		\label{0W-def}
		W & \big( \mathfrak{a}, \mathfrak{b}, \mathfrak{c}, \mathfrak{d};\textbf{A},\textbf{B},\textbf{C},\textbf{D} \big) =
		\\
		& \left\{
		\begin{array}{ll}
			q^{|\textbf{A}|},
			&
			\quad {\text{if}} \quad
			(\mathfrak{b},\textbf{B}) = (0,\textbf{e}_0),
			\quad
			(\mathfrak{c},\textbf{C}) = (\mathfrak{a},\textbf{A}),
			\quad
			(\mathfrak{d},\textbf{D}) = (0,\textbf{e}_0),
			\\
			\\
			1,
			&
			\quad {\text{if}} \quad
			(\mathfrak{b},\textbf{B}) = (0,\textbf{e}_0),
			\quad
			(\mathfrak{c},\textbf{C}) = (\mathfrak{a}+1,\textbf{A}),
			\quad
			(\mathfrak{d},\textbf{D}) = (1,\textbf{e}_0),
			\\
			\\
			(1-q^{A_i})q^{A_{[i+1,n]}},
			&
			\quad {\text{if}} \quad
			(\mathfrak{b},\textbf{B}) = (0,\textbf{e}_0),
			\quad
			(\mathfrak{c},\textbf{C}) = (\mathfrak{a}+1,\textbf{A}^{-}_i),
			\quad
			(\mathfrak{d},\textbf{D}) = (1,\textbf{e}_i),
			\\
			\\
			(1 - q^{A_j}) q^{A_{[j + 1, n]}} \textbf{1}_{A_i = 0},
			&
			\quad {\text{if}} \quad
			(\mathfrak{b},\textbf{B}) = (1, \textbf{e}_i),
			\quad
			(\mathfrak{c},\textbf{C}) = (\mathfrak{a},\textbf{A}^{+-}_{ij}),
			\quad
			(\mathfrak{d},\textbf{D}) = (1,\textbf{e}_j),
			\\
			\\
			\bm{1}_{A_i=0},
			&
			\quad {\text{if}} \quad
			(\mathfrak{b},\textbf{B}) = (1,\textbf{e}_i),
			\quad
			(\mathfrak{c},\textbf{C}) = (\mathfrak{a},\textbf{A}^{+}_i),
			\quad
			\quad
			(\mathfrak{d},\textbf{D}) = (1,\textbf{e}_0),
			\\
			\\
			(-1)^{A_i}
			q^{A_{[i, n]}},
			&
			\quad {\text{if}} \quad
			(\mathfrak{b},\textbf{B}) = (1,\textbf{e}_i),
			\quad
			(\mathfrak{c},\textbf{C}) = (\mathfrak{a},\textbf{A}),
			\quad
			(\mathfrak{d},\textbf{D}) = (1,\textbf{e}_i),
		\end{array}
		\right.
	\end{aligned} 
\end{flalign}

\noindent where we assume that $i < j$ whenever $i$ and $j$ both appear. In all cases of $(\mathfrak{a},\mathfrak{b},\mathfrak{c},\mathfrak{d};\textbf{A},\textbf{B},\textbf{C},\textbf{D})$ not listed in \eqref{0W-def}, we set $W \big(\mathfrak{a},\mathfrak{b},\mathfrak{c},\mathfrak{d};\textbf{A},\textbf{B},\textbf{C},\textbf{D}\big) = 0$.

\end{definition}

Now we can state the following partition function formula for the coefficients $f_{\boldsymbol{\lambda} / \boldsymbol{\mu}}^{\nu}$ from \eqref{0qf}.

\begin{thm}[\Cref{thm:comb} below]
	\label{0thm:comb}
	
	Fix integers $M \ge 0$ and $m \in [1, n]$; signature sequences $\boldsymbol{\lambda}, \boldsymbol{\mu} \in \SeqSign_{n; M}$; and a partition $\nu \in \Sign_{m}$, such that $|\boldsymbol{\lambda}|-|\boldsymbol{\mu}| = |\nu|$. Fix any integer 
	\begin{align*}
		K \ge \max \Big\{ \max_{i \in [1, n]} \mathfrak{T} \big(\lambda^{(i)} \big),\nu_1 + 1 \Big\},
	\end{align*}
	and write $\bar{\nu}=(\bar{\nu}_1, \bar{\nu}_2, \cdots , \bar{\nu}_{m})$ for the signature obtained by complementing $\nu$ in a $(K - 1) \times m$ box; its parts are given by $\bar{\nu}_j = K - \nu_{m - j +1} - 1$, for each $j \in [1, m]$. Denote $\mathfrak{T} (\bar{\nu}) = (\mathfrak{n}_1, \mathfrak{n}_2, \ldots , \mathfrak{n}_{m})$ so, for each $j \in [1, m]$, 
	\begin{align}
		\label{0nj}
		\mathfrak{n}_j = \bar{\nu}_j + m - j + 1 = K - \nu_{m - j + 1} + m - j.
	\end{align}
	The coefficients \eqref{0qf} are given by the partition function
	\begin{align}
		\label{0main-formula}
		f^{\nu}_{\boldsymbol{\lambda}/\boldsymbol{\mu}}(q)
		=
		\frac{(-1)^{|\bar{\nu}|}}{b_{\nu}(q)(q;q)_{m}}
		\times
		\tikz{1.1}{
			\filldraw[lgray,line width=1.5pt,fill=llgray] (1.5,3.5) -- (7.5,3.5) -- (7.5,4.5) -- (1.5,4.5) -- (1.5,3.5);
			\filldraw[lgray,line width=1.5pt,fill=llgray] (1.5,0.5) -- (7.5,0.5) -- (7.5,1.5) -- (1.5,1.5) -- (1.5,0.5);
			\foreach\y in {0,...,6}{
				\draw[lgray,line width=1.5pt] (1.5,0.5+\y) -- (7.5,0.5+\y);
			}
			\foreach\x in {0,...,6}{
				\draw[lgray,line width=1.5pt] (1.5+\x,0.5) -- (1.5+\x,6.5);
			}
			%bottom labels
			\node[below] at (2,0.5) {\footnotesize $\big( 0, \textbf{\emph{S}}_1 (\boldsymbol{\mu}) \big)$};
			\node[above,text centered] at (4,0) {$\cdots$};
			\node[above,text centered] at (5,0) {$\cdots$};
			\node[below] at (7,0.5) {\footnotesize $\big(0, \textbf{\emph{S}}_K (\boldsymbol{\mu}) \big)$};
			%top labels
			\node[above] at (2,6.5) {\footnotesize $\big(0,\textbf{\emph{S}}_1 (\boldsymbol{\lambda}) \big)$};
			\node[above,text centered] at (4,6.5) {$\cdots$};
			\node[above,text centered] at (5,6.5) {$\cdots$};
			\node[above] at (7,6.5) {\footnotesize $\big( 0, \textbf{\emph{S}}_K (\boldsymbol{\lambda}) \big)$};
			%right labels
			\node[right] at (7.5,1) {\footnotesize $(1,\textbf{\emph{e}}_0)$};
			\node[right] at (7.5,2) {\footnotesize $(0,\textbf{\emph{e}}_0)$};
			\node[right] at (7.5,3.1) {$\vdots$};
			\node[right] at (7.5,4) {\footnotesize $(1,\textbf{\emph{e}}_0)$};
			\node[right] at (7.5,5) {$\vdots$};
			\node[right] at (7.5,6) {\footnotesize $(0,\textbf{\emph{e}}_0)$};
			%left labels
			\node[left] at (1.5,1) {\footnotesize $(0,\textbf{\emph{e}}_0)$};
			\node[left] at (1.5,2) {\footnotesize $(0,\textbf{\emph{e}}_0)$};
			\node[left] at (1.5,3.1) {$\vdots$};
			\node[left] at (1.5,4) {\footnotesize $(0,\textbf{\emph{e}}_0)$};
			\node[left] at (1.5,5) {$\vdots$};
			\node[left] at (1.5,6) {\footnotesize $(m,\textbf{\emph{e}}_0)$};
		}
	\end{align}

	\noindent consisting of $\mathfrak{n}_1 + 1$ rows, where the $i$-th row of the lattice (counted from the top to bottom, starting at $i=0$) takes the form
	\begin{align}
		\label{0maya-seq2}
		\left\{
		\begin{array}{ll}
			\tikz{0.4}{
				\filldraw[lgray,line width=1.5pt,fill=llgray] (-1,-1) -- (7,-1) -- (7,1) -- (-1,1) -- (-1,-1);
				\draw[lgray,line width=1.5pt] (1,-1) -- (1,1);
				\draw[lgray,line width=1.5pt] (3,-1) -- (3,1);
				\draw[lgray,line width=1.5pt] (5,-1) -- (5,1);
				\node[left] at (-1,0) {\tiny $(0,\textbf{\emph{e}}_0)$};\node[right] at (7,0) {\tiny $(1,\textbf{\emph{e}}_0)$};
			}
			\quad\quad & i \in \mathfrak{T} (\bar{\nu}),
			\\ \\
			\tikz{0.4}{
				\draw[lgray,line width=1.5pt] (-1,-1) -- (7,-1) -- (7,1) -- (-1,1) -- (-1,-1);
				\draw[lgray,line width=1.5pt] (1,-1) -- (1,1);
				\draw[lgray,line width=1.5pt] (3,-1) -- (3,1);
				\draw[lgray,line width=1.5pt] (5,-1) -- (5,1);
				\node[left] at (-1,0) {\tiny $(0,\textbf{\emph{e}}_0)$};\node[right] at (7,0) {\tiny $(0,\textbf{\emph{e}}_0)$};
			}
			\quad\quad & i \not\in \mathfrak{T} (\bar{\nu}),\ i\not=0,
			\\ \\ 
			\tikz{0.4}{
				\draw[lgray,line width=1.5pt] (-1,-1) -- (7,-1) -- (7,1) -- (-1,1) -- (-1,-1);
				\draw[lgray,line width=1.5pt] (1,-1) -- (1,1);
				\draw[lgray,line width=1.5pt] (3,-1) -- (3,1);
				\draw[lgray,line width=1.5pt] (5,-1) -- (5,1);
				\node[left] at (-1,0) {\tiny $(m,\textbf{\emph{e}}_0)$};\node[right] at (7,0) {\tiny $(0,\textbf{\emph{e}}_0)$};
			}
			\quad\quad & i=0.
		\end{array} 
		\right.
	\end{align}
	The constant $b_{\nu}(q)$ appearing in \eqref{0main-formula} is given by $b_{\nu}(q) = \prod_{j = 1}^{\infty} (q;q)_{m_j (\nu)}$.\index{B@$b_{\nu} (q)$}
\end{thm}

 Since \Cref{0psuml} provides an expansion for $\widetilde{J}_{\lambda}$ over the skew LLT polynomials, we deduce as a quick consequence of \Cref{0thm:comb} the following combinatorial formula for the coefficients $g_{\lambda}^{\nu}$ from \eqref{0qf}. 

\begin{cor}[\Cref{jsumg} below]
	
	\label{0jsumg} 
	
	Fix integers $p, m \in [1, n]$; an anti-dominant (positive) composition $\lambda$ of length $p$; a partition $\nu$ of length $m$; and an integer $K \ge \max \{ \lambda_p, \nu_1 + 1 \}$. For any integers $i \in [1, p]$ and $j \ge 1$, set $w_{i,j} = q^{\lambda_i-j} \bm{1}_{\lambda_i > j}$. Then, the coefficient $g_{\lambda^+}^{\nu}(q,t)$ is given by
	\begin{align}
		\label{0gsumfunction} 
		g_{\lambda^+}^{\nu}(q,t)
		=
		\frac{(-1)^{|\bar{\nu}|}}{b_{\nu}(t)(t;t)_{m}}
		\sum_{\mathscr{M}}
		\prod_{i=1}^p
		\prod_{j = 1}^K 
		%\left(
		%q^{\lambda_i-j}
		%\cdot
		%\bm{1}_{\lambda_i > j}
		%\right)^{N_{i,j}}
		w_{i,j}^{M_{i,j}} \times
		\tikz{1}{
			\filldraw[lgray,line width=1.5pt,fill=llgray] (1.5,3.5) -- (7.5,3.5) -- (7.5,4.5) -- (1.5,4.5) -- (1.5,3.5);
			\filldraw[lgray,line width=1.5pt,fill=llgray] (1.5,0.5) -- (7.5,0.5) -- (7.5,1.5) -- (1.5,1.5) -- (1.5,0.5);
			\foreach\y in {0,...,6}{
				\draw[lgray,line width=1.5pt] (1.5,0.5+\y) -- (7.5,0.5+\y);
			}
			\foreach\x in {0,...,6}{
				\draw[lgray,line width=1.5pt] (1.5+\x,0.5) -- (1.5+\x,6.5);
			}
			%bottom labels
			\node[above] at (1.875,-.0625) {\footnotesize $(0,\textbf{\emph{e}}_{[1, p]})$};
			\node[above] at (3.125,0) {\footnotesize $(0,\textbf{\emph{M}}_1)$};
			\node[above,text centered] at (5,0) {$\cdots$};
			\node[above] at (7,0) {\footnotesize $(0,\textbf{\emph{M}}_K)$};
			%top labels
			\node[above] at (1.875,6.5) {\footnotesize $(0,\textbf{\emph{e}}_0)$};
			\node[above] at (3.25,6.5) {\footnotesize $(0,\textbf{\emph{M}}_1+\textbf{\emph{I}}_1)$};
			\node[above,text centered] at (5,6.5) {$\cdots$};
			\node[above] at (7,6.5) {\footnotesize $(0,\textbf{\emph{M}}_K + \textbf{\emph{I}}_K)$};
			%right labels
			\node[right] at (7.5,1) {\footnotesize $(1,\textbf{\emph{e}}_0)$};
			\node[right] at (7.5,2) {\footnotesize $(0,\textbf{\emph{e}}_0)$};
			\node[right] at (7.5,3.1) {$\vdots$};
			\node[right] at (7.5,4) {\footnotesize $(1,\textbf{\emph{e}}_0)$};
			\node[right] at (7.5,5) {$\vdots$};
			\node[right] at (7.5,6) {\footnotesize $(0,\textbf{\emph{e}}_0)$};
			%left labels
			\node[left] at (1.5,1) {\footnotesize $(0,\textbf{\emph{e}}_0)$};
			\node[left] at (1.5,2) {\footnotesize $(0,\textbf{\emph{e}}_0)$};
			\node[left] at (1.5,3.1) {$\vdots$};
			\node[left] at (1.5,4) {\footnotesize $(0,\textbf{\emph{e}}_0)$};
			\node[left] at (1.5,5) {$\vdots$};
			\node[left] at (1.5,6) {\footnotesize $(m,\textbf{\emph{e}}_0)$};
		}
	\end{align}

	\noindent where the rows of the partition function \eqref{0gsumfunction} are specified by \eqref{0nj} and \eqref{0maya-seq2}, using the same weights as in \eqref{0weights-light} and \eqref{0weights-dark} but with the $q$ there replaced by $t$ here. The sum is over all sequences $\mathscr{M} = (\textbf{\emph{M}}_1, \textbf{\emph{M}}_2, \ldots, \textbf{\emph{M}}_K)$ of elements in $\{ 0, 1 \}^n$, with $\textbf{\emph{M}}_j = (M_{1,j}, M_{2, j}, \ldots,M_{n,j})$, such that $M_{i, k} = 0$ for $i > p$ and $\textbf{\emph{I}}_k (\lambda) + \textbf{\emph{M}}_k \in \{ 0, 1 \}^n$ for each $k \in [1, K]$. Here, we abbreviated $\textbf{\emph{I}}_k = \textbf{\emph{I}}_j(\lambda)$ for each $k$, and we recall $b_{\nu} (t)$ from \Cref{0thm:comb}.
	
\end{cor}

\begin{rem} 
	
	\label{jq} 
	
	Letting $Q_{\lambda} (\textbf{x}; t)$ and $J_{\lambda} (\textbf{x}; q, t)$ denote the Hall--Littlewood and integral Macdonald polynomials,\index{Q@$Q_{\lambda} (\textbf{x})$; Hall--Littlewood polynomial} respectively, we have 
	\begin{flalign}
		\label{jsums}
	J_{\lambda} (\textbf{x}; q, t) = \displaystyle\sum_{\nu} g_{\lambda}^{\nu} (q, t) Q_{\nu} (\textbf{x}; t),
	\end{flalign}

	\noindent which follows from applying the plethystic substitution $X \mapsto (1 - t) X$ in the second statement of \eqref{0qf}, where $X = \sum_{x \in \textbf{x}} x$. Thus, \Cref{0jsumg} equivalently yields a combinatorial formula for the expansion coefficients of the (standard) Macdonald polynomials in the Hall--Littlewood basis. 
	
\end{rem} 

	Let us briefly explain how \Cref{thm:comb} and \Cref{0jsumg} can be used to obtain (not manifestly nonnegative) combinatorial formulas for the expansion coefficients of the LLT and modified Macdonald polynomials in the Schur basis $s_{\lambda} (\textbf{x})$, defined implicitly by 
	\begin{flalign*}
	\mathcal{L}_{\boldsymbol{\lambda} / \boldsymbol{\mu}} (\textbf{x}; q) = \displaystyle\sum_{\nu} K_{\boldsymbol{\lambda} / \boldsymbol{\mu}}^{\nu} (q) s_{\nu} (\textbf{x}); \qquad \widetilde{J}_{\lambda} (\textbf{x}; q, t) = \displaystyle\sum_{\nu} K_{\lambda}^{\nu} (q, t) s_{\lambda} (\textbf{x}).
	\end{flalign*}

	\noindent We have by (6.8.4(i)) of \cite{SFP} that 
	\begin{flalign*}
	J_{\lambda} (\textbf{x}; q, q) = c_{\lambda} (q) s_{\lambda} (\textbf{x}), \quad \text{where} \quad c_{\lambda} (q) = \displaystyle\prod_{(i, j) \in \lambda} (1 - q^{h_{\lambda} (i, j)}),
	\end{flalign*}

	\noindent and $h_{\lambda} (i, j)$ denotes the hook length of cell $(i, j)$ in $\lambda$. Thus, setting $q = t$ in \eqref{jsums} gives
	\begin{flalign*}
	s_{\lambda} (\textbf{x}) = c_{\lambda} (q)^{-1} \displaystyle\sum_{\nu} g_{\lambda}^{\nu} (q, q) Q_{\nu} (\textbf{x}; q).
	\end{flalign*}
	
	\noindent By Proposition 1.2 of \cite{NCFMP}, it follows that 
	\begin{flalign*}
	Q_{\lambda}' (\textbf{x}) = \displaystyle\sum_{\nu} b_{\nu} (q) c_{\lambda} (q)^{-1} g_{\lambda}^{\nu} (q, q) s_{\nu} (x).
	\end{flalign*}

	\noindent Inserting this into \eqref{0qf} then yields 
\begin{align}
	\label{ksumk} 
	K_{\boldsymbol{\lambda}/\boldsymbol{\mu}}^{\nu} (q) = \sum_{\theta} b_{\nu} (q) c_{\theta} (q)^{-1} g_{\theta}^{\nu} (q, q) f_{\boldsymbol{\lambda}/\boldsymbol{\mu}}^{\theta}(q); \qquad K_{\lambda}^{\nu} (q, t) = \displaystyle\sum_{\theta} b_{\nu} (q) c_{\theta} (q)^{-1} g_{\theta}^{\nu} (q, q) g_{\lambda}^{\theta} (q, t).
\end{align}

\noindent Since $f_{\boldsymbol{\lambda}/\boldsymbol{\mu}}^{\theta} (q)$ and $g_{\lambda}^{\theta} (q, t)$ are given by \Cref{0thm:comb} and \Cref{0jsumg}, respectively, this provides explicit combinatorial formulas for the coefficients $K_{\boldsymbol{\lambda}/\boldsymbol{\mu}}^{\nu}(q)$ and $K_{\lambda}^{\nu} (q, t)$. It was conjectured in \cite{RT} that $K_{\boldsymbol{\lambda}/\boldsymbol{\mu}}^{\nu} (q) \in \mathbb{Z}_{\ge 0} [q]$; this was later proved in \cite{AAPP}, although the proof there did not yield any combinatorial expression for the coefficients.  As stated, we emphasize that the formulas \eqref{ksumk} are not manifestly positive, since $f_{\boldsymbol{\lambda}/\boldsymbol{\mu}}^{\theta} (q)$ and $g_{\lambda}^{\theta}$ can have negative coefficients. It is natural to hope that, with $f_{\boldsymbol{\lambda}/\boldsymbol{\mu}}^{\theta}(q)$ and $g_{\nu}^{\theta}(q)$ in \eqref{ksumk} both known explicitly, one might by suitable manipulation obtain a positive rule for $K_{\boldsymbol{\lambda}/\boldsymbol{\mu}}^{\nu}(q)$ and $K_{\lambda}^{\nu} (q, t)$, but this is outside the scope of the present work.

We conclude with an example for \Cref{0thm:comb}. In what follows, we draw bosonic paths in the color blue. Unlike in previous diagrams, to simplify the figures, we depict paths as sometimes making diagonal steps, which are equivalent to making one step up or down and one step right. Let us mention that additional examples are provided in \Cref{Models2n} below, which were generated by a \texttt{Mathematica} implementation of \eqref{0main-formula}; the code is available from the authors upon request. 

\begin{example}[\Cref{lambdamun2} below]
	
	\label{0lambdamun2} 
	
	Let $n=2$, $M = 2$, $\boldsymbol{\lambda}=\left(\lambda^{(1)}, \lambda^{(2)}\right)$, and $\boldsymbol{\mu} = \left(\mu^{(1)}, \mu^{(2)}\right)$ with $\lambda^{(1)} = (3,1)$, $\lambda^{(2)}=(2,2)$, $\mu^{(1)} = (2, 0)$, and $\mu^{(2)}=(1,1)$. Let $\nu = (2,1,1)$ and choose $K = 5$, so that $\bar{\nu} = (3, 3, 2)$ and $\mathfrak{T} (\bar{\nu}) = (6,5,3)$. Given that $|\bar{\nu}| = 8$, $m = 3$, and $b_{\nu}(q) = (1-q)^2(1-q^2)$, we have
	\begin{align}
		\label{0factor3}
		\frac{(-1)^{|\bar{\nu}|}}{b_{\nu}(q)(q;q)_{m}}
		=
		\frac{1}{(1-q)^3(1-q^2)^2(1-q^3)}.
	\end{align}
	The formula \eqref{0main-formula} provides five non-vanishing lattice configurations, indicated below with their weights, where we have multiplied by the overall factor \eqref{0factor3}. Here, red is color $1$ and green is color $2$.
	\begin{align*}
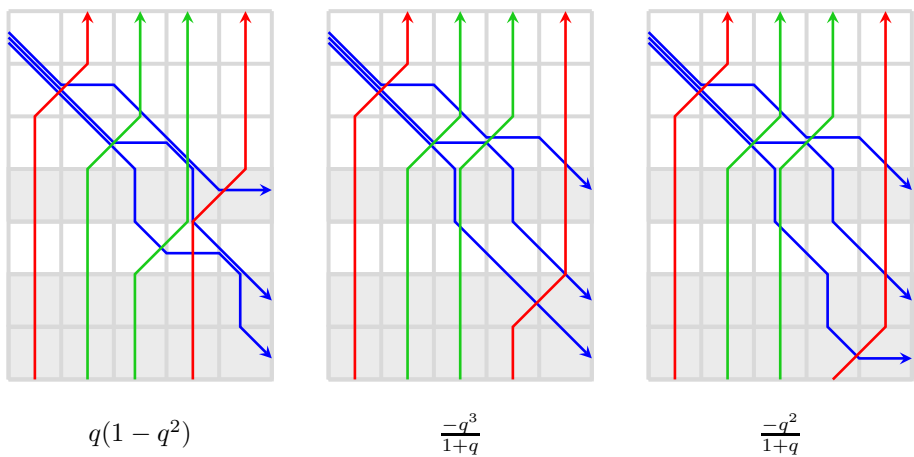
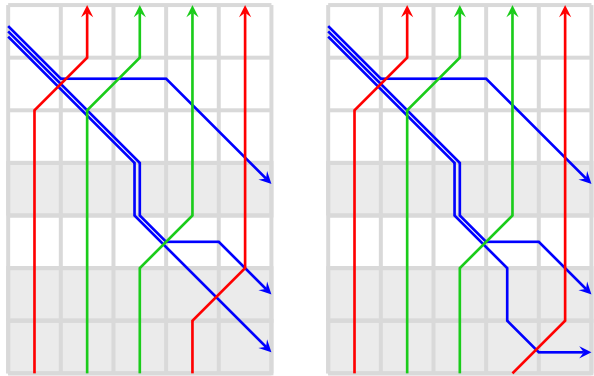
		\begin{array}{ccc}
			\tikz{0.7}{
				\filldraw[lgray,line width=1.5pt,fill=llgray] (2.5,2.5) -- (7.5,2.5) -- (7.5,3.5) -- (2.5,3.5) -- (2.5,2.5);
				\filldraw[lgray,line width=1.5pt,fill=llgray] (2.5,0.5) -- (7.5,0.5) -- (7.5,1.5) -- (2.5,1.5) -- (2.5,0.5);
				\filldraw[lgray,line width=1.5pt,fill=llgray] (2.5,-0.5) -- (7.5,-0.5) -- (7.5,0.5) -- (2.5,0.5) -- (2.5,-0.5);
				\foreach\y in {-1,...,6}{
					\draw[lgray,line width=1.5pt] (2.5,0.5+\y) -- (7.5,0.5+\y);
				}
				\foreach\x in {1,...,6}{
					\draw[lgray,line width=1.5pt] (1.5+\x,-0.5) -- (1.5+\x,6.5);
				}
				%%right labels
				%\node[right] at (7.5,0) {\footnotesize $(1,\textbf{e}_0)$};
				%\node[right] at (7.5,1) {\footnotesize $(1,\textbf{e}_0)$};
				%\node[right] at (7.5,2) {\footnotesize $(0,\textbf{e}_0)$};
				%\node[right] at (7.5,3) {\footnotesize $(1,\textbf{e}_0)$};
				%\node[right] at (7.5,4) {\footnotesize $(0,\textbf{e}_0)$};
				%\node[right] at (7.5,5) {\footnotesize $(0,\textbf{e}_0)$};
				%\node[right] at (7.5,6) {\footnotesize $(0,\textbf{e}_0)$};
				%%left labels
				%\node[left] at (2.5,0) {\footnotesize $(0,\textbf{e}_0)$};
				%\node[left] at (2.5,1) {\footnotesize $(0,\textbf{e}_0)$};
				%\node[left] at (2.5,2) {\footnotesize $(0,\textbf{e}_0)$};
				%\node[left] at (2.5,3) {\footnotesize $(0,\textbf{e}_0)$};
				%\node[left] at (2.5,4) {\footnotesize $(0,\textbf{e}_0)$};
				%\node[left] at (2.5,5) {\footnotesize $(0,\textbf{e}_0)$};
				%\node[left] at (2.5,6) {\footnotesize $(3,\textbf{e}_0)$};
				%%boson paths
				\draw[blue,line width=1pt,<-] (7.5,0-0.1) -- (7-0.1,0.5) -- (7-0.1,1.5) -- (6.5,2-0.1) -- (5.5,2-0.1) -- (5-0.1,2.5) -- (5-0.1,3.5) -- (2.5,6-0.1);
				\draw[blue,line width=1pt,<-] (7.5,1) -- (6,2.5) -- (6,3.5) -- (5.5,4) -- (4.5,4) -- (2.5,6);
				\draw[blue,line width=1pt,<-] (7.5,3+0.1) -- (6.5,3+0.1) -- (4.5,5+0.1) -- (3.5,5+0.1) -- (2.5,6+0.1);
				%red fermion paths
				\draw[red,line width=1pt,->] (3,-0.5) -- (3,4.5) -- (4,5.5) -- (4,6.5);
				\draw[red,line width=1pt,->] (6,-0.5) -- (6,2.5) -- (7,3.5) -- (7,6.5);
				%green fermion paths
				\draw[green,line width=1pt,->] (4,-0.5) -- (4,3.5) -- (5,4.5) -- (5,6.5);
				\draw[green,line width=1pt,->] (5-0.1,-0.5) -- (5-0.1,1.5) -- (6-0.1,2.5) -- (6-0.1,6.5);
			}
			&
			\quad
			\tikz{0.7}{
				\filldraw[lgray,line width=1.5pt,fill=llgray] (2.5,2.5) -- (7.5,2.5) -- (7.5,3.5) -- (2.5,3.5) -- (2.5,2.5);
				\filldraw[lgray,line width=1.5pt,fill=llgray] (2.5,0.5) -- (7.5,0.5) -- (7.5,1.5) -- (2.5,1.5) -- (2.5,0.5);
				\filldraw[lgray,line width=1.5pt,fill=llgray] (2.5,-0.5) -- (7.5,-0.5) -- (7.5,0.5) -- (2.5,0.5) -- (2.5,-0.5);
				\foreach\y in {-1,...,6}{
					\draw[lgray,line width=1.5pt] (2.5,0.5+\y) -- (7.5,0.5+\y);
				}
				\foreach\x in {1,...,6}{
					\draw[lgray,line width=1.5pt] (1.5+\x,-0.5) -- (1.5+\x,6.5);
				}
				%%right labels
				%\node[right] at (7.5,0) {\footnotesize $(1,\textbf{e}_0)$};
				%\node[right] at (7.5,1) {\footnotesize $(1,\textbf{e}_0)$};
				%\node[right] at (7.5,2) {\footnotesize $(0,\textbf{e}_0)$};
				%\node[right] at (7.5,3) {\footnotesize $(1,\textbf{e}_0)$};
				%\node[right] at (7.5,4) {\footnotesize $(0,\textbf{e}_0)$};
				%\node[right] at (7.5,5) {\footnotesize $(0,\textbf{e}_0)$};
				%\node[right] at (7.5,6) {\footnotesize $(0,\textbf{e}_0)$};
				%%left labels
				%\node[left] at (2.5,0) {\footnotesize $(0,\textbf{e}_0)$};
				%\node[left] at (2.5,1) {\footnotesize $(0,\textbf{e}_0)$};
				%\node[left] at (2.5,2) {\footnotesize $(0,\textbf{e}_0)$};
				%\node[left] at (2.5,3) {\footnotesize $(0,\textbf{e}_0)$};
				%\node[left] at (2.5,4) {\footnotesize $(0,\textbf{e}_0)$};
				%\node[left] at (2.5,5) {\footnotesize $(0,\textbf{e}_0)$};
				%\node[left] at (2.5,6) {\footnotesize $(3,\textbf{e}_0)$};
				%%boson paths
				\draw[blue,line width=1pt,<-] (7.5,0-0.1) -- (5.5,2-0.1) -- (5-0.1,2.5) -- (5-0.1,3.5) -- (2.5,6-0.1);
				\draw[blue,line width=1pt,<-] (7.5,1) -- (6,2.5) -- (6,3.5) -- (5.5,4) -- (4.5,4) -- (2.5,6);
				\draw[blue,line width=1pt,<-] (7.5,3+0.1) -- (6.5,4+0.1) -- (5.5,4+0.1) -- (4.5,5+0.1) -- (3.5,5+0.1) -- (2.5,6+0.1);
				%red fermion paths
				\draw[red,line width=1pt,->] (3,-0.5) -- (3,4.5) -- (4,5.5) -- (4,6.5);
				\draw[red,line width=1pt,->] (6,-0.5) -- (6,0.5) -- (7,1.5) -- (7,6.5);
				%green fermion paths
				\draw[green,line width=1pt,->] (4,-0.5) -- (4,3.5) -- (5,4.5) -- (5,6.5);
				\draw[green,line width=1pt,->] (5,-0.5) -- (5,3.5) -- (6,4.5) -- (6,6.5);
			}
			&
			\quad
			\tikz{0.7}{
				\filldraw[lgray,line width=1.5pt,fill=llgray] (2.5,2.5) -- (7.5,2.5) -- (7.5,3.5) -- (2.5,3.5) -- (2.5,2.5);
				\filldraw[lgray,line width=1.5pt,fill=llgray] (2.5,0.5) -- (7.5,0.5) -- (7.5,1.5) -- (2.5,1.5) -- (2.5,0.5);
				\filldraw[lgray,line width=1.5pt,fill=llgray] (2.5,-0.5) -- (7.5,-0.5) -- (7.5,0.5) -- (2.5,0.5) -- (2.5,-0.5);
				\foreach\y in {-1,...,6}{
					\draw[lgray,line width=1.5pt] (2.5,0.5+\y) -- (7.5,0.5+\y);
				}
				\foreach\x in {1,...,6}{
					\draw[lgray,line width=1.5pt] (1.5+\x,-0.5) -- (1.5+\x,6.5);
				}
				%%right labels
				%\node[right] at (7.5,0) {\footnotesize $(1,\textbf{e}_0)$};
				%\node[right] at (7.5,1) {\footnotesize $(1,\textbf{e}_0)$};
				%\node[right] at (7.5,2) {\footnotesize $(0,\textbf{e}_0)$};
				%\node[right] at (7.5,3) {\footnotesize $(1,\textbf{e}_0)$};
				%\node[right] at (7.5,4) {\footnotesize $(0,\textbf{e}_0)$};
				%\node[right] at (7.5,5) {\footnotesize $(0,\textbf{e}_0)$};
				%\node[right] at (7.5,6) {\footnotesize $(0,\textbf{e}_0)$};
				%%left labels
				%\node[left] at (2.5,0) {\footnotesize $(0,\textbf{e}_0)$};
				%\node[left] at (2.5,1) {\footnotesize $(0,\textbf{e}_0)$};
				%\node[left] at (2.5,2) {\footnotesize $(0,\textbf{e}_0)$};
				%\node[left] at (2.5,3) {\footnotesize $(0,\textbf{e}_0)$};
				%\node[left] at (2.5,4) {\footnotesize $(0,\textbf{e}_0)$};
				%\node[left] at (2.5,5) {\footnotesize $(0,\textbf{e}_0)$};
				%\node[left] at (2.5,6) {\footnotesize $(3,\textbf{e}_0)$};
				%%boson paths
				\draw[blue,line width=1pt,<-] (7.5,0-0.1) -- (6.5,0-0.1) -- (6-0.1,0.5) -- (6-0.1,1.5) -- (5.5,2-0.1) -- (5-0.1,2.5) -- (5-0.1,3.5) -- (2.5,6-0.1);
				\draw[blue,line width=1pt,<-] (7.5,1) -- (6,2.5) -- (6,3.5) -- (5.5,4) -- (4.5,4) -- (2.5,6);
				\draw[blue,line width=1pt,<-] (7.5,3+0.1) -- (6.5,4+0.1) -- (5.5,4+0.1) -- (4.5,5+0.1) -- (3.5,5+0.1) -- (2.5,6+0.1);
				%red fermion paths
				\draw[red,line width=1pt,->] (3,-0.5) -- (3,4.5) -- (4,5.5) -- (4,6.5);
				\draw[red,line width=1pt,->] (6,-0.5) -- (7,0.5) -- (7,6.5);
				%green fermion paths
				\draw[green,line width=1pt,->] (4,-0.5) -- (4,3.5) -- (5,4.5) -- (5,6.5);
				\draw[green,line width=1pt,->] (5,-0.5) -- (5,3.5) -- (6,4.5) -- (6,6.5);
			}
			\\
			\\
			q (1-q^2) &\quad \frac{-q^3}{1+q} &\quad \frac{-q^2}{1+q}
		\end{array}
	\end{align*}
	
	\begin{align*}
		\begin{array}{cc}
			\tikz{0.7}{
				\filldraw[lgray,line width=1.5pt,fill=llgray] (2.5,2.5) -- (7.5,2.5) -- (7.5,3.5) -- (2.5,3.5) -- (2.5,2.5);
				\filldraw[lgray,line width=1.5pt,fill=llgray] (2.5,0.5) -- (7.5,0.5) -- (7.5,1.5) -- (2.5,1.5) -- (2.5,0.5);
				\filldraw[lgray,line width=1.5pt,fill=llgray] (2.5,-0.5) -- (7.5,-0.5) -- (7.5,0.5) -- (2.5,0.5) -- (2.5,-0.5);
				\foreach\y in {-1,...,6}{
					\draw[lgray,line width=1.5pt] (2.5,0.5+\y) -- (7.5,0.5+\y);
				}
				\foreach\x in {1,...,6}{
					\draw[lgray,line width=1.5pt] (1.5+\x,-0.5) -- (1.5+\x,6.5);
				}
				%%right labels
				%\node[right] at (7.5,0) {\footnotesize $(1,\textbf{e}_0)$};
				%\node[right] at (7.5,1) {\footnotesize $(1,\textbf{e}_0)$};
				%\node[right] at (7.5,2) {\footnotesize $(0,\textbf{e}_0)$};
				%\node[right] at (7.5,3) {\footnotesize $(1,\textbf{e}_0)$};
				%\node[right] at (7.5,4) {\footnotesize $(0,\textbf{e}_0)$};
				%\node[right] at (7.5,5) {\footnotesize $(0,\textbf{e}_0)$};
				%\node[right] at (7.5,6) {\footnotesize $(0,\textbf{e}_0)$};
				%%left labels
				%\node[left] at (2.5,0) {\footnotesize $(0,\textbf{e}_0)$};
				%\node[left] at (2.5,1) {\footnotesize $(0,\textbf{e}_0)$};
				%\node[left] at (2.5,2) {\footnotesize $(0,\textbf{e}_0)$};
				%\node[left] at (2.5,3) {\footnotesize $(0,\textbf{e}_0)$};
				%\node[left] at (2.5,4) {\footnotesize $(0,\textbf{e}_0)$};
				%\node[left] at (2.5,5) {\footnotesize $(0,\textbf{e}_0)$};
				%\node[left] at (2.5,6) {\footnotesize $(3,\textbf{e}_0)$};
				%%boson paths
				\draw[blue,line width=1pt,<-] (7.5,0-0.1) -- (5.5,2-0.1) -- (5-0.1,2.5) -- (5-0.1,3.5) -- (2.5,6-0.1);
				\draw[blue,line width=1pt,<-] (7.5,1) -- (6.5,2) -- (5.5,2) -- (5,2.5) -- (5,3.5) -- (2.5,6);
				\draw[blue,line width=1pt,<-] (7.5,3+0.1) --(5.5,5+0.1) -- (3.5,5+0.1) -- (2.5,6+0.1);
				%red fermion paths
				\draw[red,line width=1pt,->] (3,-0.5) -- (3,4.5) -- (4,5.5) -- (4,6.5);
				\draw[red,line width=1pt,->] (6,-0.5) -- (6,0.5) -- (7,1.5) -- (7,6.5);
				%green fermion paths
				\draw[green,line width=1pt,->] (4,-0.5) -- (4,4.5) -- (5,5.5) -- (5,6.5);
				\draw[green,line width=1pt,->] (5,-0.5) -- (5,1.5) -- (6,2.5) -- (6,6.5);
			}
			& \quad
			\tikz{0.7}{
				\filldraw[lgray,line width=1.5pt,fill=llgray] (2.5,2.5) -- (7.5,2.5) -- (7.5,3.5) -- (2.5,3.5) -- (2.5,2.5);
				\filldraw[lgray,line width=1.5pt,fill=llgray] (2.5,0.5) -- (7.5,0.5) -- (7.5,1.5) -- (2.5,1.5) -- (2.5,0.5);
				\filldraw[lgray,line width=1.5pt,fill=llgray] (2.5,-0.5) -- (7.5,-0.5) -- (7.5,0.5) -- (2.5,0.5) -- (2.5,-0.5);
				\foreach\y in {-1,...,6}{
					\draw[lgray,line width=1.5pt] (2.5,0.5+\y) -- (7.5,0.5+\y);
				}
				\foreach\x in {1,...,6}{
					\draw[lgray,line width=1.5pt] (1.5+\x,-0.5) -- (1.5+\x,6.5);
				}
				%%right labels
				%\node[right] at (7.5,0) {\footnotesize $(1,\textbf{e}_0)$};
				%\node[right] at (7.5,1) {\footnotesize $(1,\textbf{e}_0)$};
				%\node[right] at (7.5,2) {\footnotesize $(0,\textbf{e}_0)$};
				%\node[right] at (7.5,3) {\footnotesize $(1,\textbf{e}_0)$};
				%\node[right] at (7.5,4) {\footnotesize $(0,\textbf{e}_0)$};
				%\node[right] at (7.5,5) {\footnotesize $(0,\textbf{e}_0)$};
				%\node[right] at (7.5,6) {\footnotesize $(0,\textbf{e}_0)$};
				%%left labels
				%\node[left] at (2.5,0) {\footnotesize $(0,\textbf{e}_0)$};
				%\node[left] at (2.5,1) {\footnotesize $(0,\textbf{e}_0)$};
				%\node[left] at (2.5,2) {\footnotesize $(0,\textbf{e}_0)$};
				%\node[left] at (2.5,3) {\footnotesize $(0,\textbf{e}_0)$};
				%\node[left] at (2.5,4) {\footnotesize $(0,\textbf{e}_0)$};
				%\node[left] at (2.5,5) {\footnotesize $(0,\textbf{e}_0)$};
				%\node[left] at (2.5,6) {\footnotesize $(3,\textbf{e}_0)$};
				%%boson paths
				\draw[blue,line width=1pt,<-] (7.5,0-0.1) -- (6.5,0-0.1) -- (6-0.1,0.5) -- (6-0.1,1.5) -- (5.5,2-0.1) -- (5-0.1,2.5) -- (5-0.1,3.5) -- (2.5,6-0.1);
				\draw[blue,line width=1pt,<-] (7.5,1) -- (6.5,2) -- (5.5,2) -- (5,2.5) -- (5,3.5) -- (2.5,6);
				\draw[blue,line width=1pt,<-] (7.5,3+0.1) --(5.5,5+0.1) -- (3.5,5+0.1) -- (2.5,6+0.1);
				%red fermion paths
				\draw[red,line width=1pt,->] (3,-0.5) -- (3,4.5) -- (4,5.5) -- (4,6.5);
				\draw[red,line width=1pt,->] (6,-0.5) -- (7,0.5) -- (7,6.5);
				%green fermion paths
				\draw[green,line width=1pt,->] (4,-0.5) -- (4,4.5) -- (5,5.5) -- (5,6.5);
				\draw[green,line width=1pt,->] (5,-0.5) -- (5,1.5) -- (6,2.5) -- (6,6.5);
			}
			\\
			\\
			\frac{q(1-q)}{1+q} &\quad \frac{1-q}{1+q}
		\end{array}
	\end{align*}
	We therefore find that $f^{\nu}_{\boldsymbol{\lambda}/\boldsymbol{\mu}}(q)
	=
	q(1-q^2)-\frac{q^3}{1+q}-\frac{q^2}{1+q}+\frac{q(1-q)}{1+q}+\frac{1-q}{1+q}
	=
	1-q^2 -q^3$.
\end{example}

\subsection*{Acknowledgements}

The authors thank Johnny Fonseca and Jason Saied for pointing out an inconsistency in the original definition of our fused weights. The authors are also grateful to Daniel Bump for bringing \cite{MFSLM} to our attention; Anton Mellit for explaining his work to us; and Siddhartha Sahi for bringing \cite{OSSI} to our attention. Amol Aggarwal was partially supported by an NSF Graduate Research Fellowship, under grant number DGE-1144152, and a Clay Research Fellowship. Alexei Borodin was partially supported by the NSF grants DMS-1664619, DMS-1853981 and the Simons Investigator program. Michael Wheeler was supported by an Australian Research Council Future Fellowship, grant FT200100981.

\chapter{\texorpdfstring{$U_q \big(\widehat{\mathfrak{sl}} (m|n) \big)$}{} Vertex Model} 

\label{Weights1}

In this chapter we introduce and provide notation for the $U_q \big( \widehat{\mathfrak{sl}} (m | n) \big)$ vertex model that will serve as the basis for the fused model we will derive in this text. We further establish a ``color merging'' property for it, indicating how a $U_q \big( \widehat{\mathfrak{sl}} (m' | n') \big)$ model can be obtained from a $U_q \big( \widehat{\mathfrak{sl}} (m | n) \big)$ one, if $m \ge m'$ and $n \ge n'$.

\section{Fundamental \texorpdfstring{$R$}{}-Matrix}

\label{Weightszq}

The vertex weights for the model analyzed in this text will arise from applying the fusion procedure to the $R$-matrix associated with the fundamental representation of $U_q \big( \widehat{\mathfrak{sl}} (m | n) \big)$. In this section we recall this $R$-matrix and the Yang--Baxter equation it satisfies. Throughout, we fix two integers $m \ge 1$ and $n \ge 0$ such that $m + n \ge 2$. 

Let $a, b$ denote two indices and let $V_a \simeq \mathbb{C}^{m + n} \simeq V_b$ denote two $(m + n)$-dimensional complex vector spaces, spanned by basis vectors $| 0 \rangle, | 1 \rangle, \ldots , | m + n - 1 \rangle$. For any complex number $z \in \mathbb{C}$, define the \emph{fundamental $R$-matrix} $R_{ab} (z) = R_{a, b} (z) \in \End (V_a \otimes V_b)$\index{R@$R_{ab} (z)$; fundamental $R$-matrix} by setting 
\begin{flalign}
\label{rzdefinition} 
\begin{aligned}
R_{ab} (z) & = \displaystyle\sum_{i = 0}^{m - 1} R_z (i, i; i, i) E_a^{(ii)} \otimes E_b^{(ii)} + \displaystyle\sum_{j = m}^{m + n - 1} R_z (j, j; j, j) E_a^{(jj)} \otimes E_b^{(jj)}  \\
& \quad + \displaystyle\sum_{0 \le i < j \le m + n - 1} R_z (i, j; i, j) E_a^{(jj)} \otimes E_b^{(ii)} + \displaystyle\sum_{0 \le i < j \le m + n - 1} R_z (j, i; j, i) E_a^{(ii)} \otimes E_b^{(jj)} \\
& \quad + \displaystyle\sum_{0 \le i < j \le m + n - 1} R_z (i, j; j, i) E_a^{(ji)} \otimes E_b^{(ij)} + \displaystyle\sum_{0 \le i < j \le m + n - 1} R_z (j, i; i, j) E_a^{(ij)} \otimes E_b^{(ji)}.
\end{aligned} 
\end{flalign}

\noindent Here, for any integers $i, j \in [0, m + n - 1]$ and index $k \in \{ a, b \}$, $E_k^{(ij)} \in \End (V_k)$ denotes the $(m + n) \times (m + n)$ matrix whose entries are all equal to $0$, except for its $(i, j)$-entry, which is equal to $1$; stated alternatively, it satisfies 
\begin{flalign*} 	
E_k^{(ij)} | j \rangle = | i \rangle, \quad \text{and} \quad E_k^{(ij)} | i' \rangle = 0, \quad \text{for each $i' \ne i$}.
\end{flalign*} 

\noindent The $R$-matrix entries $R_z (i_1, j_1; i_2, j_2) = R_z^{(m; n)} (i_1, j_1; i_2, j_2)$ in \eqref{rzdefinition} are defined by setting
\begin{flalign}
\label{rzij}
\begin{aligned} 
& R_z (j, i; j, i) = \displaystyle\frac{q (1 - z)}{1 - qz}; \qquad R_z (i, j; i, j) = \displaystyle\frac{1 - z}{1 - qz}; \\
& R_z (j, i; i, j) = \displaystyle\frac{1 - q}{1 - qz}; \qquad R_z (i, j; j, i) = \displaystyle\frac{z (1 - q)}{1 - qz},
\end{aligned}
\end{flalign} 

\noindent for $0 \le i < j \le m + n - 1$; 
\begin{flalign}
\label{rzi}
R_z (i, i; i, i) = 1, \quad \text{for $i \in [0, m - 1]$;} \qquad R_z (j, j; j, j) = \displaystyle\frac{z - q}{1 - qz}, \quad \text{for $j \in [m, m + n - 1]$}.    
\end{flalign}

\noindent We also set $R_z (i_1, j_1; i_2, j_2) = 0$ for any $(i_1, j_1; i_2, j_2)$ not of the above form. These $R$-matrix entries were originally introduced under a different parameterization as equation (14) of \cite{TSTESS}, but they appear essentially as above in Section 5.1 of \cite{PI}. 

\begin{rem}
	
\label{m1m2rr}

The $R_z (i_1, j_1; i_2, j_2)$ do not depend on $n$ unless $i_1 = j_1 = i_2 = j_2 \ge m$, namely, 
\begin{flalign*}
R_z^{(m; n)} (i_1, j_1; i_2, j_2) = R_z^{(m + n; 0)} (i_1, j_1; i_2, j_2), \qquad \text{unless $i_1 = j_1 = i_2 = j_2 \in [m, m + n - 1]$}.
\end{flalign*}

\end{rem}

We interpret these entries $R_z (i_1, j_1; i_2, j_2)$ as \emph{vertex weights} as follows. A \emph{vertex} $v$ is a transverse intersection between two oriented curves, typically (after a rotation, if necessary) a horizontal line $\ell_1$ directed east and a vertical line $\ell_2$ directed north. Associated with each of these two curves is a \emph{rapidity parameter}, which is a complex number. Denoting that for $\ell_1$ by $x \in \mathbb{C}$ and that for $\ell_2$ by $y \in \mathbb{C}$, the \emph{spectral parameter} of $v$ is the ratio $\frac{y}{x}$.

\begin{figure}[t]
	
	\begin{center}
		
		\begin{tikzpicture}[
		>=stealth,
		scale = .75
		]

		\draw[-, black] (-7.5, 3.1) -- (7.5, 3.1);
		\draw[-, black] (-7.5, -2.5) -- (7.5, -2.5);
		\draw[-, black] (-7.5, -1.1) -- (7.5, -1.1);
		\draw[-, black] (-7.5, -.4) -- (7.5, -.4);
		\draw[-, black] (-7.5, 2.4) -- (7.5, 2.4);
		\draw[-, black] (-7.5, -2.5) -- (-7.5, 3.1);
		\draw[-, black] (7.5, -2.5) -- (7.5, 3.1);
		\draw[-, black] (-5, -2.5) -- (-5, 3.1);
		\draw[-, black] (5, -2.5) -- (5, 2.4);
		\draw[-, black] (-2.5, -2.5) -- (-2.5, 3.1);
		\draw[-, black] (2.5, -2.5) -- (2.5, 2.4);
		\draw[-, black] (0, -2.5) -- (0, 2.4);

		\draw[->, thick, red] (-7.15, 1) -- (-5.35, 1); 
		\draw[->, thick, red] (-6.25, .1) -- (-6.25, 1.9); 
		
		\draw[->, thick, blue] (3.75, .1) -- (3.75, 1) -- (4.65, 1);
		\draw[->, thick, red] (2.85, 1) -- (3.75, 1) -- (3.75, 1.9); 
		
		\draw[->, thick, blue] (-1.25, .1) -- (-1.25, 1.9);
		\draw[->, thick,  red] (-2.15, 1) -- (-.35, 1);
		
		\draw[->, thick, blue] (.35, 1) -- (2.15, 1);
		\draw[->, thick, red] (1.25, .1) -- (1.25, 1.9);
		
		\draw[->, thick, blue] (5.35, 1) -- (6.25, 1) -- (6.25, 1.9);
		\draw[->, thick, red] (6.25, .1) -- (6.25, 1) -- (7.15, 1); 
		
		\draw[->, thick, blue] (-3.75, .1) -- (-3.75, 1.9);
		\draw[->, thick, blue] (-4.65, 1) -- (-2.85, 1);

		\filldraw[fill=gray!50!white, draw=black] (-5.35, 1) circle [radius=0] node [black, right = -1, scale = .75] {$i$};
		\filldraw[fill=gray!50!white, draw=black] (-2.85, 1) circle [radius=0] node [black, right = -1, scale = .75] {$j$};
		\filldraw[fill=gray!50!white, draw=black] (-.35, 1) circle [radius=0] node [black, right = -1, scale = .75] {$i$};
		\filldraw[fill=gray!50!white, draw=black] (2.15, 1) circle [radius=0] node [black, right = -1, scale = .75] {$j$};
		\filldraw[fill=gray!50!white, draw=black] (4.65, 1) circle [radius=0] node [black, right = -1, scale = .75] {$j$};
		\filldraw[fill=gray!50!white, draw=black] (7.15, 1) circle [radius=0] node [black, right = -1, scale = .75] {$i$};
		
		\filldraw[fill=gray!50!white, draw=black] (5.35, 1) circle [radius=0] node [black, left = -1, scale = .75] {$j$};
		\filldraw[fill=gray!50!white, draw=black] (2.85, 1) circle [radius=0] node [black, left = -1, scale = .75] {$i$};
		\filldraw[fill=gray!50!white, draw=black] (.35, 1) circle [radius=0] node [black, left = -1, scale = .75] {$j$};
		\filldraw[fill=gray!50!white, draw=black] (-2.15, 1) circle [radius=0] node [black, left = -1, scale = .75] {$i$};
		\filldraw[fill=gray!50!white, draw=black] (-4.65, 1) circle [radius=0] node [black, left = -1, scale = .75] {$j$};
		\filldraw[fill=gray!50!white, draw=black] (-7.15, 1) circle [radius=0] node [black, left = -1, scale = .75] {$i$};
		
		\filldraw[fill=gray!50!white, draw=black] (-6.25, 1.9) circle [radius=0] node [black, above = -1, scale = .75] {$i$};
		\filldraw[fill=gray!50!white, draw=black] (-3.75, 1.9) circle [radius=0] node [black, above = -1, scale = .75] {$j$};
		\filldraw[fill=gray!50!white, draw=black] (-1.25, 1.9) circle [radius=0] node [black, above = -1, scale = .75] {$j$};
		\filldraw[fill=gray!50!white, draw=black] (1.25, 1.9) circle [radius=0] node [black, above = -1, scale = .75] {$i$};
		\filldraw[fill=gray!50!white, draw=black] (3.75, 1.9) circle [radius=0] node [black, above = -1, scale = .75] {$i$};
		\filldraw[fill=gray!50!white, draw=black] (6.25, 1.9) circle [radius=0] node [black, above = -1, scale = .75] {$j$};

		\filldraw[fill=gray!50!white, draw=black] (-6.25, .1) circle [radius=0] node [black, below = -1, scale = .75] {$i$};
		\filldraw[fill=gray!50!white, draw=black] (-3.75, .1) circle [radius=0] node [black, below = -1, scale = .75] {$j$};
		\filldraw[fill=gray!50!white, draw=black] (-1.25, .1) circle [radius=0] node [black, below = -1, scale = .75] {$j$};
		\filldraw[fill=gray!50!white, draw=black] (1.25, .1) circle [radius=0] node [black, below = -1, scale = .75] {$i$};
		\filldraw[fill=gray!50!white, draw=black] (3.75, .1) circle [radius=0] node [black, below = -1, scale = .75] {$j$};
		\filldraw[fill=gray!50!white, draw=black] (6.25, .1) circle [radius=0] node [black, below = -1, scale = .75] {$i$};

		\filldraw[fill=gray!50!white, draw=black] (-6.25, .1) circle [radius=0] node [black, below = -1, scale = .75] {$i$};
		\filldraw[fill=gray!50!white, draw=black] (-3.75, .1) circle [radius=0] node [black, below = -1, scale = .75] {$j$};
		\filldraw[fill=gray!50!white, draw=black] (-1.25, .1) circle [radius=0] node [black, below = -1, scale = .75] {$j$};
		\filldraw[fill=gray!50!white, draw=black] (1.25, .1) circle [radius=0] node [black, below = -1, scale = .75] {$i$};
		\filldraw[fill=gray!50!white, draw=black] (3.75, .1) circle [radius=0] node [black, below = -1, scale = .75] {$j$};
		\filldraw[fill=gray!50!white, draw=black] (6.25, .1) circle [radius=0] node [black, below = -1, scale = .75] {$i$};
		
		\filldraw[fill=gray!50!white, draw=black] (-6.25, 2.75) circle [radius=0] node [black, scale = .75] {$i \in [0, m)$};
		\filldraw[fill=gray!50!white, draw=black] (-3.75, 2.75) circle [radius=0] node [black, scale = .75] {$j \in [m, m + n)$};
		\filldraw[fill=gray!50!white, draw=black] (2.5, 2.75) circle [radius=0] node [black] {$0 \le i < j < m + n$}; 
		
		\filldraw[fill=gray!50!white, draw=black] (-6.25, -.75) circle [radius=0] node [black] {$(i, i; i, i)$};
		\filldraw[fill=gray!50!white, draw=black] (-3.75, -.75) circle [radius=0] node [black] {$(j, j; j, j)$};
		\filldraw[fill=gray!50!white, draw=black] (-1.25, -.75) circle [radius=0] node [black] {$(j, i; j, i)$};
		\filldraw[fill=gray!50!white, draw=black] (1.25, -.75) circle [radius=0] node [black] {$(i, j; i, j)$};
		\filldraw[fill=gray!50!white, draw=black] (3.75, -.75) circle [radius=0] node [black] {$(j, i; i, j)$};
		\filldraw[fill=gray!50!white, draw=black] (6.25, -.75) circle [radius=0] node [black] {$(i, j; j, i)$};
		
		\filldraw[fill=gray!50!white, draw=black] (-6.25, -1.8) circle [radius=0] node [black] {$1$};
		\filldraw[fill=gray!50!white, draw=black] (-3.75, -1.8) circle [radius=0] node [black] {$\displaystyle\frac{z - q}{1 - qz}$};
		\filldraw[fill=gray!50!white, draw=black] (-1.25, -1.8) circle [radius=0] node [black] {$\displaystyle\frac{q (1 - z)}{1 - qz}$};
		\filldraw[fill=gray!50!white, draw=black] (1.25, -1.8) circle [radius=0] node [black] {$\displaystyle\frac{1 - z}{1 - qz}$};
		\filldraw[fill=gray!50!white, draw=black] (3.75, -1.8) circle [radius=0] node [black] {$\displaystyle\frac{1 - q}{1 - qz}$};
		\filldraw[fill=gray!50!white, draw=black] (6.25, -1.8) circle [radius=0] node [black] {$\displaystyle\frac{z (1 - q)}{1 - qz}$};

		\end{tikzpicture}
		
	\end{center}
	
	\caption{\label{sixvertexfigureclass} The possible colored arrow configurations, and their vertex weights, are depicted above.}
\end{figure}

Each of the four segments, also called \emph{arrows}, of $\ell_1$ and $\ell_2$ adjacent to $v$ is assigned a \emph{color}, which is a label in $\{ 0, 1, \ldots, m + n - 1 \}$. Let $i_1$ and $j_1$ denote the colors of the vertical and horizontal segments  entering $v$, respectively; similarly, let $i_2$ and $j_2$ denote the colors of the vertical and horizontal arrows exiting $v$, respectively. Then, $i_1, j_1, i_2, j_2 \in \{ 0, 1, \ldots , m + n - 1 \}$. We further assume that $\{ i_1, j_1 \} = \{ i_2, j_2 \}$ as multi-sets, so that the same number of arrows of any given color enters $v$ as exits $v$; this is known as \emph{arrow conservation}. We refer to the quadruple $(i_1, j_1; i_2, j_2)$ as the \emph{arrow configuration} at $v$, and we view $R_z (i_1, j_1; i_2, j_2)$ as the weight of a vertex with arrow configuration $(i_1, j_1; i_2, j_2)$ and spectral parameter $z$. The possible arrow configurations are depicted in \Cref{sixvertexfigureclass}. Under this convention, at most one arrow can occupy any edge; we will later remove this condition through fusion in \Cref{WeightsR} below. 

Observe from \eqref{rzi} that vertices adjacent to four arrows of the same color $i$ have different weights depending on whether $i \in [0, m - 1]$ or $i \in [m, m + n - 1]$. Throughout this text, we refer to colors in the interval $[0, m - 1]$ as \emph{bosonic} and to the remaining colors (in the interval $[m, m + n - 1]$) as \emph{fermionic}. This terminology will in a sense be later justified by \Cref{rb1d2} below, which implies that fermionic colors satisfy a certain exclusion property (which does not hold for bosonic colors). 

The following proposition indicates that the $R$-matrix $R_{ab} (z)$ from \eqref{rzdefinition} satisfies the Yang--Baxter equation. It was originally due to \cite{TSTESS}, but it can also be verified directly from the explicit forms \eqref{rzij} and \eqref{rzi} for its entries. 

\begin{prop}[{\cite[Section 3]{TSTESS}}]
	
	\label{rrr1} 
	
	Let $V_1$, $V_2$, and $V_3$ denote three $(m + n)$-dimensional complex vector spaces, and let $x, y, z \in \mathbb{C}$ denote complex numbers. As operators on $V_1 \otimes V_2 \otimes V_3$, we have
	\begin{flalign}
	\label{matrixrrr}
	R_{12} \Big( \frac{y}{x} \Big) R_{13} \Big( \frac{z}{x} \Big) R_{23} \Big( \frac{z}{y} \Big) = R_{23} \Big( \frac{z}{y} \Big) R_{13} \Big( \frac{z}{x} \Big) R_{12} \Big( \frac{y}{x} \Big),
	\end{flalign}
	
	\noindent where $R_{ij}$ acts as the identity on $V_k$, for $k \notin \{ i, j \}$.  
	
\end{prop}

For fixed indices $i_1, j_1, k_1, i_3, j_3, k_3 \in [0, m + n - 1]$, comparing the entries in \eqref{matrixrrr} in the $(i_1, j_1, k_1)$-th row and $(i_3, j_3, k_3)$-column, we deduce that 
\begin{flalign}
\label{rrrijk}
\begin{aligned}
\displaystyle\sum_{i_2, j_2, k_2} & R_{y / x} (i_1, j_1; i_2, j_2) R_{z / x} (k_1, j_2; k_2, j_3) R_{z / y} (k_2, i_2; k_3, i_3) \\
& = \displaystyle\sum_{i_2, j_2, k_2} R_{z / y} (k_1, i_1; k_2, i_2) R_{z / x} (k_2, j_1; k_3, j_2) R_{y / x} (i_2, j_2; i_3, j_3),
\end{aligned}
\end{flalign}

\noindent where on both sides $i_2, j_2, k_2$ are each ranged over all indices in $[0, m + n - 1]$. The diagrammatic interpretation of this formulation of the Yang--Baxter equation is then given by

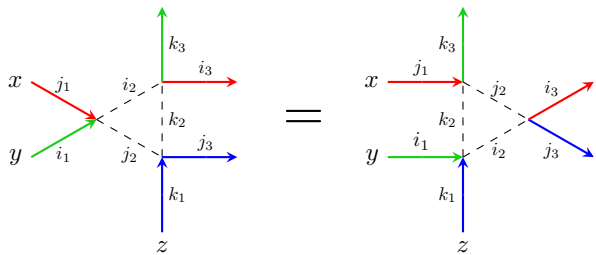
\begin{center}

	\begin{tikzpicture}[
		>=stealth,
		auto,
		style={
			scale = 1
		}
		]
		
		\draw[->, thick, green] (-.87, -.5) -- (0, 0);
		\draw[->, thick, red] (-.87, .5) -- (0, 0);
		\draw[->, thick, blue] (.87, -1.5) -- (.87, -.5); 
		
		\draw[] (-.87, .5) circle[radius = 0]  node[left, scale = .95]{$x$};
		\draw[] (-.87, -.5) circle[radius = 0]  node[left, scale = .95]{$y$};
		\draw[] (.87, -1.5) circle[radius = 0]  node[below, scale = .95]{$z$};
		
		\draw[black, dashed] (0, 0) -- (.87, -.5); 
		\draw[black, dashed] (0, 0) -- (.87, .5); 
		\draw[black, dashed] (.87, -.5) -- (.87, .5);
		
		\draw[->, thick, red] (.87, .5) -- (1.87, .5); 
		\draw[->, thick, blue] (.87, -.5) -- (1.87, -.5); 
		\draw[->, thick, green] (.87, .5) -- (.87, 1.5);

		\draw[] (3.87, .5) circle[radius = 0]  node[left, scale = .95]{$x$};
		\draw[] (3.87, -.5) circle[radius = 0]  node[left, scale = .95]{$y$};
		\draw[] (4.87, -1.5) circle[radius = 0]  node[below, scale = .95]{$z$};
		
		\draw[->, thick, blue] (4.87, -1.5) -- (4.87, -.5); 
		\draw[->, thick, red] (3.87, .5) -- (4.87, .5); 
		\draw[->, thick, green] (3.87, -.5) -- (4.87, -.5); 
		\draw[->, thick, green] (4.87, .5) -- (4.87, 1.5); 
		\draw[->, thick, blue] (5.74, 0) -- (6.61, -.5); 
		\draw[->, thick, red] (5.74, 0) -- (6.61, .5); 
		
		\draw[black, dashed] (4.87, -.5) -- (4.87, .5);
		\draw[-, black, dashed] (4.87, -.5) -- (5.74, 0); 
		\draw[black, dashed] (4.87, .5) -- (5.74, 0);

		\filldraw[fill=white, draw=black] (2.75, 0) circle [radius=0] node[scale = 2]{$=$};
		
		\filldraw[fill=white, draw=black] (-.44, -.275) circle [radius=0] node[below, scale = .7]{$i_1$};
		\filldraw[fill=white, draw=black] (.44, .275) circle [radius=0] node[above, scale = .7]{$i_2$};
		\filldraw[fill=white, draw=black] (1.45, .5) circle [radius=0] node[above, scale = .7]{$i_3$};
		\filldraw[fill=white, draw=black] (-.44, .275) circle [radius=0] node[above, scale = .7]{$j_1$};
		\filldraw[fill=white, draw=black] (.44, -.275) circle [radius=0] node[below, scale = .7]{$j_2$};
		\filldraw[fill=white, draw=black] (1.45, -.5) circle [radius=0] node[above, scale = .7]{$j_3$};
		\filldraw[fill=white, draw=black] (.87, -1) circle [radius=0] node[right, scale = .7]{$k_1$};
		\filldraw[fill=white, draw=black] (.87, 0) circle [radius=0] node[right, scale = .7]{$k_2$};
		\filldraw[fill=white, draw=black] (.87, 1) circle [radius=0] node[right, scale = .7]{$k_3$};
		
		\filldraw[fill=white, draw=black] (4.32, -.5) circle [radius=0] node[above, scale =.8]{$i_1$};
		\filldraw[fill=white, draw=black] (5.35, -.225) circle [radius=0] node[below, scale = .7]{$i_2$};
		\filldraw[fill=white, draw=black] (6.05, .225) circle [radius=0] node[above, scale = .7]{$i_3$};
		\filldraw[fill=white, draw=black] (4.32, .5) circle [radius=0] node[above, scale = .7]{$j_1$};
		\filldraw[fill=white, draw=black] (5.35, .225) circle [radius=0] node[above, scale = .7]{$j_2$};
		\filldraw[fill=white, draw=black] (6.05, -.225) circle [radius=0] node[below, scale = .7]{$j_3$};
		\filldraw[fill=white, draw=black] (4.87, -1) circle [radius=0] node[left, scale = .7]{$k_1$};
		\filldraw[fill=white, draw=black] (4.87, 0) circle [radius=0] node[left, scale = .7]{$k_2$};
		\filldraw[fill=white, draw=black] (4.87, 1) circle [radius=0] node[left, scale = .7]{$k_3$};

	\end{tikzpicture}
	
\end{center}

\noindent where on either side of the equation is a family of vertices, and we view the weight of each family as the product of the weights of its constituent vertices. Along the solid edges the colors are fixed, and along the dashed ones they are summed over.

\section{Domains, Boundary Data, and Partition Functions}

\label{DomainBoundary} 

In this section we introduce notation for vertex models on the general domains we consider. In what follows, we view any vertex $v \in \mathbb{Z}^2$ as a two-dimensional vector, so that vertices may be added or subtracted to form elements of $\mathbb{Z}^2$. 

\begin{figure}[t]
	
	\begin{center}		
		
		\begin{tikzpicture}[
		>=stealth,
		auto,
		style={
			scale = .85
		}
		]
		\draw[dashed] (1, 8) -- (1, 7) -- (3, 7);
		\draw[dashed] (2, 8) -- (2, 6) -- (3 ,6);
		\draw[dashed] (3, 4) -- (4, 4) -- (4, 3);
		\draw[dashed] (4, 5) -- (4, 4) -- (5, 4);
		
		\draw[] (-.8, 8) circle[radius = 0] node[left, black, scale = .7]{$e_1$};
		\draw[] (-.8, 7) circle[radius = 0] node[left, black, scale = .7]{$e_2$};
		\draw[] (0, 6.5) circle[radius = 0] node[left, black, scale = .7]{$e_3$};
		\draw[] (.5, 6) circle[radius = 0] node[above, black, scale = .7]{$e_4$};
		\draw[] (1, 5.5) circle[radius = 0] node[left, black, scale = .7]{$e_5$};
		\draw[] (1.5, 5) circle[radius = 0] node[above, black, scale = .7]{$e_6$};
		\draw[] (2, 4.5) circle[radius = 0] node[left, black, scale = .7]{$e_7$};
		\draw[] (2.5, 4) circle[radius = 0] node[above, black, scale = .7]{$e_8$};
		\draw[] (2.5, 3) circle[radius = 0] node[above, black, scale = .7]{$e_9$};
		\draw[] (3, 2.2) circle[radius = 0] node[below, black, scale = .7]{$e_{10}$};
		\draw[] (4, 2.2) circle[radius = 0] node[below, black, scale = .7]{$e_{11}$};
		\draw[] (5, 2.2) circle[radius = 0] node[below, black, scale = .7]{$e_{12}$};
		
		\draw[] (0, 8.8) circle[radius = 0] node[above, black, scale = .7]{$f_1$};
		\draw[] (1, 8.8) circle[radius = 0] node[above, black, scale = .7]{$f_2$};
		\draw[] (2, 8.8) circle[radius = 0] node[above, black, scale = .7]{$f_3$};
		\draw[] (3, 8.8) circle[radius = 0] node[above, black, scale = .7]{$f_4$};
		\draw[] (3.4, 8) circle[radius = 0] node[above, black, scale = .7]{$f_5$};
		\draw[] (3.4, 7) circle[radius = 0] node[above, black, scale = .7]{$f_6$};
		\draw[] (3.4, 6) circle[radius = 0] node[above, black, scale = .7]{$f_7$};
		\draw[] (4, 5.4) circle[radius = 0] node[left, black, scale = .7]{$f_8$};
		\draw[] (5, 5.4) circle[radius = 0] node[left, black, scale = .7]{$f_9$};
		\draw[] (5.8, 5) circle[radius = 0] node[right, black, scale = .7]{$f_{10}$};
		\draw[] (5.8, 4) circle[radius = 0] node[right, black, scale = .7]{$f_{11}$};
		\draw[] (5.8, 3) circle[radius = 0] node[right, black, scale = .7]{$f_{12}$};
		
		\draw[thick, blue, ->] (-.8, 8) -- (0, 8);
		\draw[thick, red, ->] (-.8, 7) -- (0, 7);
		\draw[thick, red, ->] (.2, 6) -- (1, 6);
		\draw[thick, green, ->] (1.2, 5) -- (2, 5);
		\draw[thick, green, ->] (2.2, 4) -- (3, 4);
		\draw[thick, blue, ->] (2.2, 3) -- (3, 3);
		
		\draw[thick, red, ->] (0, 6.2) -- (0, 7);
		\draw[thick, green, ->] (1, 5.2) -- (1, 6);
		\draw[thick, green, ->] (2, 4.2) -- (2, 5);
		\draw[thick, orange, ->] (3, 2.2) -- (3, 3);
		\draw[thick, blue, ->] (4, 2.2) -- (4, 3);
		\draw[thick, blue, ->] (5, 2.2) -- (5, 3);
		
		\draw[thick, red, ->] (3, 8) -- (3.8, 8);
		\draw[thick, red, ->] (3, 7) -- (3.8, 7);
		\draw[thick, blue, ->] (3, 6) -- (3.8, 6);
		\draw[thick, blue, ->] (5, 5) -- (5.8, 5);
		\draw[thick, blue, ->] (5, 4) -- (5.8, 4);
		\draw[thick, orange, ->] (5, 3) -- (5.8, 3);
		
		\draw[thick, red, ->] (0, 8) -- (0, 8.8);
		\draw[thick, green, ->] (1, 8) -- (1, 8.8);
		\draw[thick, blue, ->] (2, 8) -- (2, 8.8);
		\draw[thick, green, ->] (3, 8) -- (3, 8.8);
		\draw[thick, green, ->] (4, 5) -- (4, 5.8);
		\draw[thick, green, ->] (5, 5) -- (5, 5.8);
		
		\draw[ultra thick] (0, 8) -- (0, 7) -- (1, 7) -- (1, 6) -- (2, 6) -- (2, 5) -- (3, 5) -- (3, 3) -- (5, 3);
		\draw[ultra thick] (0, 8) -- (3, 8) -- (3, 7) -- (3, 5) -- (4, 5) -- (5, 5) -- (5, 3);

		\draw[thick, red, ->] (9.2, 7) -- (10, 7) -- (10, 8.8);
		\draw[thick, red, ->] (10, 6.2) -- (10, 7) -- (12, 7) -- (12, 8) -- (13.8, 8);
		\draw[thick, red, ->] (10.2, 6) -- (13, 6) -- (13, 7) -- (13.8, 7);
		\draw[thick, blue, ->] (12.2, 3) -- (13, 3) -- (13, 6) -- (13.8, 6);
		\draw[thick, blue, ->] (9.2, 8) -- (12, 8) -- (12, 8.8);
		\draw[thick, blue, ->] (14, 2.2) -- (14, 5) -- (15.8, 5);
		\draw[thick, blue, ->] (15, 2.2) -- (15, 4) -- (15.8, 4);
		\draw[thick, green, ->] (11.2, 5) -- (12, 5) -- (12, 7) -- (13, 7) -- (13, 8.8);
		\draw[thick, green, ->] (12.2, 4) -- (15, 4) -- (15, 5.8);
		\draw[thick, green, ->] (11, 5.2) -- (11, 8.8);
		\draw[thick, green, ->] (12, 4.2) -- (12, 5) -- (14, 5) -- (14, 5.8);
		\draw[thick, orange, ->] (13, 2.2) -- (13, 3) -- (15.8, 3);
		
		 \draw[fill=black] (14, 4) circle[radius = .075];
		 \draw[] (14.15, 4.15) circle[radius = 0] node[black, scale = .6]{$v$};
		 \draw[] (13.4, 4) circle[radius = 0] node[above, black, scale = .7]{$a(v)$};
		 \draw[] (14, 4.6) circle[radius = 0] node[right, black, scale = .7]{$b(v)$};
		 \draw[] (14.6, 4) circle[radius = 0] node[below, black, scale = .7]{$c(v)$};
		 \draw[] (14, 3.4) circle[radius = 0] node[left, black, scale = .7]{$d(v)$};
		
		\end{tikzpicture}
		
	\end{center}
	
	\caption{\label{domainpaths} Shown to the left is an east-south domain with some boundary data, and shown to the right is a path ensemble on it. } 	
\end{figure}
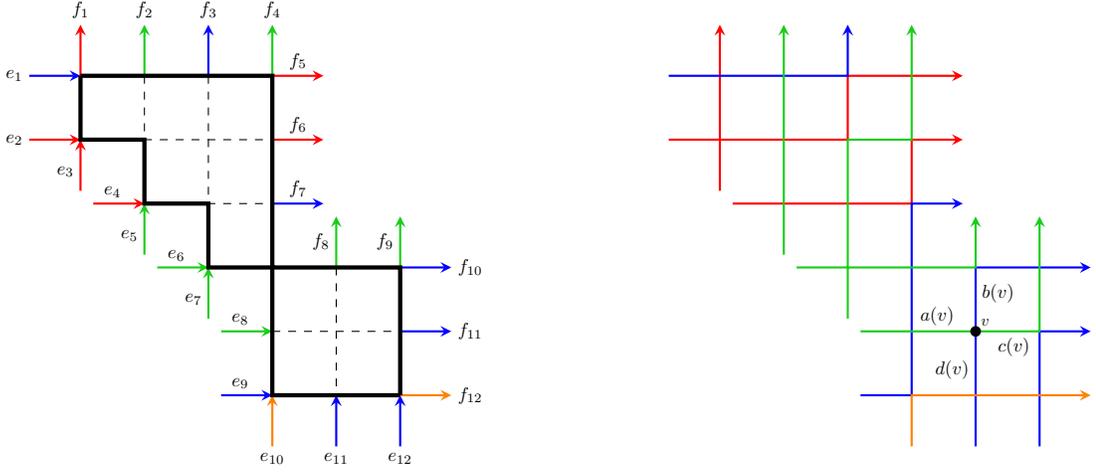

An \emph{east-south path} is an ordered sequence of vertices $\textbf{p} = (v_1, v_2, \ldots,  v_k) \subset \mathbb{Z}^2$ such that $v_{i + 1} - v_i \in \big\{ (1, 0), (0, -1) \big\}$ for each $i \in [1, k - 1]$. We call $k$ the \emph{length} of $\textbf{p}$ and refer to any pair of consecutive vertices $(v_i, v_{i + 1})$ as an \emph{edge} of $\textbf{p}$; this edge is horizontal if $v_{i + 1} - v_i = (1, 0)$ and vertical if $v_{i + 1} - v_i = (0, -1)$. Given two east-south paths  $\textbf{p} = (v_1, v_2, \ldots , v_k)$ and $\textbf{p}' = (v_1', v_2', \ldots , v_k')$ of the same length, we write $\textbf{p}' \ge \textbf{p}$ (or equivalently $\textbf{p} \le \textbf{p}'$) if $v_i' - v_i \in \mathbb{Z}_{\ge 0}^2$ for each $i \in [0, k]$, that is, if each $v_i'$ is northeast of $v_i$. In this case, if we further have $v_0 = v_0'$ and $v_k = v_k'$ (namely, if $\textbf{p}$ and $\textbf{p}'$ share the same starting and ending points), then let $\mathcal{D} (\textbf{p}, \textbf{p}') \subset \mathbb{Z}^2$ denote the domain bounded by $\textbf{p}$ and $\textbf{p}'$. More specifically, it consists of those vertices $u \in \mathbb{Z}^2$ for which there exists some $i \in [1, k]$ such that $u - v_i, v_i' - u \in \mathbb{Z}_{\ge 0}^2$. A domain of the form $\mathcal{D} (\textbf{p}, \textbf{p}')$, for some east-south paths $\textbf{p} \le \textbf{p}'$, is an \emph{east-south domain}. We refer to the left side of \Cref{domainpaths} for an example.

Given an east-south domain $\mathcal{D} \subset \mathbb{Z}^2$, we call any $v \in \mathbb{Z}^2 \setminus \mathcal{D}$ a \emph{boundary vertex} of $\mathcal{D}$ if there exists some $u \in \mathcal{D}$ that is adjacent to $v$ through an edge of $\mathbb{Z}^2$; in this case, the edge connecting $(v, u)$ is a \emph{boundary edge}. This edge is \emph{incoming} if $u - v \in \big\{ (1, 0), (0, 1) \big\}$ and otherwise \emph{outgoing} if $u - v \in \big\{ (-1, 0), (0, -1) \big\}$. Similarly, a boundary vertex $v \in \mathbb{Z}^2 \setminus \mathcal{D}$ is \emph{incoming} or \emph{outgoing} if it belongs to an incoming or outgoing boundary edge, respectively.

Any east-south domain $\mathcal{D} (\textbf{p}, \textbf{p}')$, whose paths $\textbf{p}$ and $\textbf{p}'$ are both of length $k$, has $k + 1$ incoming and outgoing boundary vertices, and thus $k + 1$ incoming and outgoing boundary edges. We index these incoming (and outgoing) boundary edges by $\{ 1, 2, \ldots , k + 1 \}$ from northwest to southeast, that is, so that the incoming boundary vertex in the $i$-th incoming edge is weakly northwest of that in the $j$-th one if and only if  $i \le j$ (and similarly for outgoing edges). 

We next consider vertex models on an east-south domain $\mathcal{D} = \mathcal{D} (\textbf{p}, \textbf{p}') \subset \mathbb{Z}^2$. A \emph{path ensemble} on $\mathcal{D}$ is a \emph{consistent} assignment of an arrow configution $\big( a(v), b(v); c(v), d(v) \big)$ to each vertex $v \in \mathcal{D}$; the consistency here means that for any $u, v \in \mathcal{D}$ we have $b(u) = d(v)$ if $u - v = (1, 0)$ and $a(u) = c(v)$ if $u - v = (0, 1)$. The arrow conservation condition $\big\{ a(v), b(v) \big\} = \big\{ c(v), d(v) \big\}$ implies that any path ensemble may be viewed as a collection of (possibly crossing) colored paths. Each such path emanates from and exits through a boundary vertex of $\mathcal{D}$, and no two paths share an edge; this last condition will later be removed in \Cref{Weights1n} below through fusion. We refer to the right side of \Cref{domainpaths} for a depiction.

Boundary data for a path ensemble prescribes the colors of its incoming and outgoing paths. More specifically, for any sequences of indices $\mathfrak{E} = (e_1, e_2, \ldots , e_{k + 1})$ and $\mathfrak{F} = (f_1, f_2, \ldots , f_{k + 1})$ in $[0, m + n - 1]$, a path ensemble has \emph{boundary data} $(\mathfrak{E}; \mathfrak{F})$ if the following holds. For each $i \in [1, k + 1]$, an arrow of color $e_i$ enters through the $i$-th incoming edge of $\mathcal{D}$ and an arrow of color $f_i$ exits through the $i$-th outgoing edge of $\mathcal{D}$; see the left side of \Cref{domainpaths}. We sometimes refer to $\mathfrak{E}$ as \emph{entrance data} on $\mathcal{D}$ and to $\mathfrak{F}$ as \emph{exit data}. 

Given a set of spectral parameters\footnote{Recall from \Cref{Weightszq} that the spectral parameter $z = z(i, j)$ at any vertex $(i, j) \in \mathcal{D}$ is set to be the ratio $x_j^{-1} y_i$, where $x_j$ and $y_i$ are rapidity parameters associated with row $j$ and column $i$ of the model, respectively. Although this constraint is central to the integrability of the underlying model, we will omit it in \Cref{zef1} in order to ease notation (but later restore it in contexts where we directly require it, such as in \Cref{zxy1} below).} $\textbf{z} = \big( z(v) \big)_{v \in \mathcal{D}}$ for each vertex of $\mathcal{D}$, the following definition provides notation for the partition function (sum of weights of all path ensembles) of a vertex model on $\mathcal{D}$ with boundary data $(\mathfrak{E}; \mathfrak{F})$, under the weights $R_z (a, b; c, d)$ from \eqref{rzij} and \eqref{rzi}. 

\begin{definition} 
	
	\label{zef1}
	
	Let $\mathcal{D} \subset \mathbb{Z}^2$ denote an east-south domain. Fix boundary data $(\mathfrak{E}; \mathfrak{F})$ on $\mathcal{D}$, and a set of complex numbers $\textbf{z} = \big( z(v) \big)_{v \in \mathcal{D}}$. Define the \emph{partition function} 
	\begin{flalign*}
	Z_{\mathcal{D}}^{(m; n)} (\mathfrak{E}; \mathfrak{F} \boldsymbol{\mid} \textbf{z})  = \displaystyle\sum \displaystyle\prod_{v \in \mathcal{D}} R_{z(v)} \big( a(v), b(v); c(v), d(v) \big),
	\end{flalign*} \index{Z@$Z_{\mathcal{D}}^{(m; n)} (\mathfrak{E}; \mathfrak{F} \boldsymbol{\mid} \textbf{z})$}
	
	\noindent where the sum is over all path ensembles on $\mathcal{D}$ with boundary data $(\mathfrak{E}; \mathfrak{F})$. 
	
\end{definition}

\section{Color Merging} 

\label{Colors} 

In this section we describe a way of merging colors of a $U_q \big(\widehat{\mathfrak{sl}} (m | n) \big)$ vertex model, analogous to the ``colorblind projection'' results from Section 2.3 of \cite{SVMST}. In particular, we show that certain signed sums of partition functions of such a vertex model is equal to the partition function for this model with several of its colors identified (as in \Cref{vertexsum1} below), thereby reducing its \emph{rank} $m + n - 1$.

To explain this in more detail, we first require some notation. An \emph{interval partition} of an integer set $I = \{ a, a + 1, \ldots , b \} \subset \mathbb{Z}_{\ge 0}$ is a collection of mutually disjoint integer intervals $\mathbb{J} = (J_0, J_1, \ldots, J_{\ell})$ such that $\bigcup_{k = 0}^{\ell} J_k = I$; here, the intervals in $\mathbb{J}$ are ordered so that any $j_i \in J_i$ is less than any $j_k \in J_k$ whenever $i < k$. Given an interval partition $\mathbb{J} = (J_0, J_1, \ldots , J_{\ell})$ of $I$, we define the function $\theta_{\mathbb{J}}: I \rightarrow [0, \ell]$\index{0@$\theta_{\mathbb{J}}$} by setting $\theta_{\mathbb{J}} (j) = k$ for each $j \in J_k$ and $k \in [0, \ell]$. In what follows, we will typically view $I = \{ 0, 1, \ldots , m + n - 1 \}$ as the set of colors in our model and $\theta_{\mathbb{J}}$ as a prescription for ``merging'' them, so that all colors in any $J_k$ are identified and renamed to color $k$.

We further require a certain inversion count. To define it, for any sequence $\mathfrak{I} = (i_1, i_2, \ldots , i_k)$ of indices in $[0, m + n - 1]$ and integer interval $J \subseteq [0, m + n - 1]$, we set
\begin{flalign}
\label{ijsum} 
\inv (\mathfrak{I}; J) = \displaystyle\sum_{1 \le h < j \le k} \textbf{1}_{i_h \in J} \textbf{1}_{i_j \in J} \textbf{1}_{i_h > i_j}.
\end{flalign}

\noindent Observe for example that if $J = \{ 0, 1, \ldots , m + n - 1 \}$ then $\inv (\mathfrak{I}; J) = \inv (\mathfrak{I})$, which denotes the number of pairs of indices $(h, j) \in [1, k]^2$ such that $h < j$ and $i_h > i_j$.

 Now we have the below proposition, which essentially states the following for any interval partition $\mathbb{J}$ of $[0, m + n - 1]$ separating $[0, m - 1]$ (bosonic colors) from $[m, m + n - 1]$ (fermionic ones). Consider a $U_q \big(\widehat{\mathfrak{sl}} (m | n) \big)$ vertex model partition function, ``symmetrize'' it over all choices consistent with $\mathbb{J}$ of colors for its bosonic exiting arrows, and ``anti-symmetrize'' it over all choices consistent with $\mathbb{J}$ of colors for its fermionic entering arrows. This yields the partition function for the vertex model obtained from the original one by merging its colors as prescribed by $\theta_{\mathbb{J}}$. 
 
 In what follows, we recall the partition function $Z_{\mathcal{D}}^{(m; n)} (\mathfrak{E}; \mathfrak{F} \boldsymbol{\mid} \textbf{z})$ from \Cref{zef1}, and for any function $f: \mathbb{Z} \rightarrow \mathbb{Z}$ and sequence $\mathfrak{I} = (i_1, i_2, \ldots , i_k)$ we denote $f(\mathfrak{I}) = \big( f(i_1), f(i_2), \ldots , f(i_k) \big)$. We will establish this proposition in \Cref{ProofSumZ} below. 

\begin{prop} 
	
\label{zsumef} 

Fix integers $m \ge m' \ge 1$, $n \ge n' \ge 0$, and $k \ge 1$; an east-south domain $\mathcal{D} = \mathcal{D} (\textbf{\emph{p}}, \textbf{\emph{p}}')$ with boundary paths $\textbf{\emph{p}}$ and $\textbf{\emph{p}}'$ of length $k$; and a set of complex numbers $\textbf{\emph{z}} = \big( z(v) \big)_{v \in \mathcal{D}}$. Let $\mathbb{J} = (J_0, J_1, \ldots , J_{m' + n' - 1})$ denote an interval partition of $\{ 0, 1, \ldots , m + n - 1\}$ such that
\begin{flalign}
\label{jmn} 
\bigcup_{k = 0}^{m' - 1} J_k = \{ 0, 1, \ldots , m - 1\}; \qquad \bigcup_{k = m'}^{m' + n' - 1} J_k = \{ m, m + 1, \ldots , m + n - 1 \}.
\end{flalign} 

\noindent For any fixed sequences of indices $\mathfrak{E} = (e_1, e_2, \ldots , e_{k + 1})$ and $\mathfrak{F} = (f_1, f_2, \ldots , f_{k + 1})$, with entries in $[0, m + n - 1]$, constituting entrance and exit data on $\mathcal{D}$, respectively, we have 
\begin{flalign}
\label{sumefzsum}
\displaystyle\sum_{\breve{\mathfrak{E}}, \breve{\mathfrak{F}}} Z_{\mathcal{D}}^{(m; n)} (\breve{\mathfrak{E}}; \breve{\mathfrak{F}} \boldsymbol{\mid} \textbf{\emph{z}}) \displaystyle\prod_{i = m'}^{m' + n' - 1} (-1)^{\inv (\breve{\mathfrak{E}}; J_i) - \inv (\mathfrak{F}; J_i)} = Z_{\mathcal{D}}^{(m'; n')} \big( \theta_{\mathbb{J}} (\mathfrak{E}); \theta_{\mathbb{J}} (\mathfrak{F}) \boldsymbol{\mid} \textbf{\emph{z}}).
\end{flalign} 

\noindent Here, the sum is over all sequences of indices $\breve{\mathfrak{E}} = (\breve{e}_1, \breve{e}_2, \ldots , \breve{e}_{k + 1})$ and $\breve{\mathfrak{F}} = (\breve{f}_1, \breve{f}_2, \ldots , \breve{f}_{k + 1})$ with entries in $[0, m + n - 1]$ such that $\theta_{\mathbb{J}} (\breve{\mathfrak{E}}) = \theta_{\mathbb{J}} (\mathfrak{E})$; $\theta_{\mathbb{J}} (\breve{\mathfrak{F}}) = \theta_{\mathbb{J}} (\mathfrak{F})$; $\breve{e}_i = e_i$ whenever $e_i \in [0, m - 1]$; and $\breve{f}_i = f_i$ whenever $f_i \in [m, m + n - 1]$. 

\end{prop} 

For $n = n' = 0$, \Cref{zsumef} essentially states that symmetrizing a $U_q \big( \widehat{\mathfrak{sl}} (m) \big)$ vertex model partition function over its exit data (in a way consistent with $\mathbb{J}$) in effect merges colors in the model, as prescribed by $\theta_{\mathbb{J}}$; this statement appeared as Proposition 4.11 of \cite{SVMP} (when all but one interval in $\mathbb{J}$ had length $1$). The necessity of anti-symmetrizing the entering fermionic colors for this model appears to be a new effect present for $n > n' > 0$. 

Now let us consider certain special cases of \Cref{zsumef}, first when $m' + n' = m + n$ and second (in a certain scenario of) when $m' + n' = m + n - 1$ and $|\mathcal{D}| = 1$. They will facilitate the proof of \Cref{zsumef} when $|\mathcal{D}| = 1$, before proceeding to fully general domains. 

\begin{example}
	
	\label{mnmn} 
	
	Suppose that $m' + n' = m + n$. Then, $m = m'$ and $n = n'$, which implies that each interval in $\mathbb{J}$ is a singleton. Thus $\theta_{\mathbb{J}} (i) = i$ for each $i \in [0, m + n - 1]$, and so the left side of \eqref{sumefzsum} is supported on the term $(\breve{\mathfrak{E}}, \breve{\mathfrak{F}}) = (\mathfrak{E}, \mathfrak{F})$. Since $\inv (\mathfrak{E}; J_i) = 0 = \inv (\mathfrak{F}; J_i)$ for each $i \in [0, m' + n' - 1]$ (by \eqref{ijsum}, since each $J_i$ is a singleton), this implies that both sides of \eqref{sumefzsum} are equal $Z_{\mathcal{D}}^{(m; n)} (\mathfrak{E}; \mathfrak{F} \boldsymbol{\mid} \textbf{z})$, thereby verifying \Cref{zsumef} in this case.

\end{example} 

\begin{lem}
	
	\label{mnmn1d1} 
	
	Suppose that $\mathcal{D} = \{ v \}$ consists of a single vertex (so $k = 1$), and that $m' + n' = m + n - 1$. Further assume that there exists an index $h \in [0, m + n - 1]$ with $e_1, e_2, f_1, f_2 \in \{ h, h + 1 \}$ and
	\begin{flalign}
	\label{jh1}
	J_i = \{ i \}, \quad \text{for $i \in [1, h - 1]$}; \qquad J_h = \{ h, h + 1\}; \qquad J_i = \{ i + 1\}, \quad \text{for $i \in [h + 1, m + n - 1]$}.
	\end{flalign}
	
	\noindent Then \Cref{zsumef} holds. 
	
\end{lem} 

\begin{figure}
	
	\begin{center}		
		
		\begin{tikzpicture}[
		>=stealth,
		auto,
		style={
			scale = .7
		}
		] 
		
		\draw[->, red, thick] (1, 3.5) -- (1, 5.5);
		\draw[->, blue, thick] (0, 4.5) -- (2, 4.5);
		
		\draw[black] (1, 3.5) circle[radius = 0] node[below, scale = .7]{$h$};
		\draw[black] (1, 5.5) circle[radius = 0] node[above, scale = .7]{$h$};
		\draw[black] (0, 4.5) circle[radius = 0] node[above, scale = .7]{$h + 1$};
		\draw[black] (2, 4.5) circle[radius = 0] node[above, scale = .7]{$h + 1$};

		\draw[black] (3, 4.5) circle[radius = 0] node[scale = 1.1]{$+$};
		
		\draw[->, red, thick] (5, 3.5) -- (5, 4.5) -- (6, 4.5); 
		\draw[->, blue, thick] (4, 4.5) -- (5, 4.5) -- (5, 5.5);

		\draw[black] (5, 3.5) circle[radius = 0] node[below, scale = .7]{$h$};
		\draw[black] (5, 5.5) circle[radius = 0] node[above, scale = .7]{$h + 1$};
		\draw[black] (4, 4.5) circle[radius = 0] node[above, scale = .7]{$h + 1$};
		\draw[black] (6, 4.5) circle[radius = 0] node[above, scale = .7]{$h$};

		\draw[black] (5, -1) circle[radius = 0] node[below, scale = .7]{$h$};
		\draw[black] (5, 1) circle[radius = 0] node[above, scale = .7]{$h$};
		\draw[black] (4, 0) circle[radius = 0] node[above, scale = .7]{$h + 1$};
		\draw[black] (6, 0) circle[radius = 0] node[above, scale = .7]{$h + 1$};
		
		\draw[->, red, thick] (9.5, 3.5) -- (9.5, 5.5); 
		\draw[->, red, thick] (8.5, 4.5) -- (10.5, 4.5);
		
		\draw[black] (9.5, 3.5) circle[radius = 0] node[below, scale = .7]{$h$};
		\draw[black] (9.5, 5.5) circle[radius = 0] node[above, scale = .7]{$h$};
		\draw[black] (8.5, 4.5) circle[radius = 0] node[above, scale = .7]{$h$};
		\draw[black] (10.5, 4.5) circle[radius = 0] node[above, scale = .7]{$h$};
		
		\draw[black] (9.5, -1) circle[radius = 0] node[below, scale = .7]{$h$};
		\draw[black] (9.5, 1) circle[radius = 0] node[above, scale = .7]{$h$};
		\draw[black] (8.5, 0) circle[radius = 0] node[above, scale = .7]{$h$};
		\draw[black] (10.5, 0) circle[radius = 0] node[above, scale = .7]{$h$};
		
		\draw[black] (7.25, 4.5) circle[radius = 0] node[scale = 1.25]{$=$};
		
		\draw[black] (11.75, 4.5) circle[radius = 0] node[scale = 1.25]{$=$};
		
		\draw[->, blue, thick] (14, 3.5) -- (14, 4.5) -- (15, 4.5); 
		\draw[->, red, thick] (13, 4.5) -- (14, 4.5) -- (14, 5.5);

		\draw[black] (16, 4.5) circle[radius = 0] node[scale = 1.1]{$+$};
		
		\draw[->, blue, thick] (18, 3.5) -- (18, 5.5); 
		\draw[->, red, thick] (17, 4.5) -- (19, 4.5);

		\draw[black] (7.25, 0) circle[radius = 0] node[scale = 1.25]{$=$};
		\draw[black] (11.75, 0) circle[radius = 0] node[scale = 1.25]{$=$};
		
		\draw[->, red, thick] (0, 0) -- (1, 0) -- (1, 1);
		\draw[->, blue, thick] (1, -1) -- (1, 0) -- (2, 0);

		\draw[black] (1, -1) circle[radius = 0] node[below, scale = .7]{$h + 1$};
		\draw[black] (1, 1) circle[radius = 0] node[above, scale = .7]{$h$};
		\draw[black] (0, 0) circle[radius = 0] node[above, scale = .7]{$h$};
		\draw[black] (2, 0) circle[radius = 0] node[above, scale = .7]{$h + 1$};
		
		\draw[black] (3, 0) circle[radius = 0] node[scale = 1.1]{$-$};
		
		\draw[->, red, thick] (5, -1) -- (5, 1); 
		\draw[->, blue, thick] (4, 0) -- (6, 0);
		
		\draw[->, red, thick] (9.5, -1) -- (9.5, 1); 
		\draw[->, red, thick] (8.5, 0) -- (10.5, 0);
		
		\draw[->, red, thick] (14, -1) -- (14, 0) -- (15, 0); 
		\draw[->, blue, thick] (13, 0) -- (14, 0) -- (14, 1);
		
		\draw[black] (14, 3.5) circle[radius = 0] node[below, scale = .7]{$h + 1$};
		\draw[black] (14, 5.5) circle[radius = 0] node[above, scale = .7]{$h$};
		\draw[black] (13, 4.5) circle[radius = 0] node[above, scale = .7]{$h$};
		\draw[black] (15, 4.5) circle[radius = 0] node[above, scale = .7]{$h + 1$};
		
		\draw[black] (14, -1) circle[radius = 0] node[below, scale = .7]{$h$};
		\draw[black] (14, 1) circle[radius = 0] node[above, scale = .7]{$h + 1$};
		\draw[black] (13, 0) circle[radius = 0] node[above, scale = .7]{$h + 1$};
		\draw[black] (15, 0) circle[radius = 0] node[above, scale = .7]{$h$};
		
		\draw[black] (16, 0) circle[radius = 0] node[scale = 1.1]{$-$};
		
		\draw[->, blue, thick] (18, -1) -- (18, 1); 
		\draw[->, red, thick] (17, 0) -- (19, 0);
		
		\draw[black] (18, -1) circle[radius = 0] node[below, scale = .7]{$h + 1$};
		\draw[black] (18, 1) circle[radius = 0] node[above, scale = .7]{$h + 1$};
		\draw[black] (17, 0) circle[radius = 0] node[above, scale = .7]{$h$};
		\draw[black] (19, 0) circle[radius = 0] node[above, scale = .7]{$h$}; 
		
		\draw[black] (18, 3.5) circle[radius = 0] node[below, scale = .7]{$h + 1$};
		\draw[black] (18, 5.5) circle[radius = 0] node[above, scale = .7]{$h + 1$};
		\draw[black] (17, 4.5) circle[radius = 0] node[above, scale = .7]{$h$};
		\draw[black] (19, 4.5) circle[radius = 0] node[above, scale = .7]{$h$};
		
		\draw[very thick] (-.5, -1.75) -- (-.5, 7.25) -- (19.5, 7.25) -- (19.5, -1.75) -- (-.5, -1.75);
		\draw[very thick] (-.5, 1.75) -- (19.5, 1.75);
		\draw[very thick] (-.5, 2.75) -- (19.5, 2.75);
		\draw[very thick] (-.5, 6.25) -- (19.5, 6.25);
		
		\draw[] (9.5, 6.75) circle[radius = 0] node[scale = 1]{$h \in [0, m -2]$};
		\draw[] (9.5, 2.25) circle[radius = 0] node[scale = 1]{$h \in [m, m + n - 1]$};

		\end{tikzpicture}
		
	\end{center}
	
	\caption{\label{vertexsum1} Shown above is the diagrammatic interpretation for \eqref{rzrzsum1} when $h < m - 1$, and shown below is that for \eqref{rzrzsum2} when $h \ge m$.} 	
\end{figure}
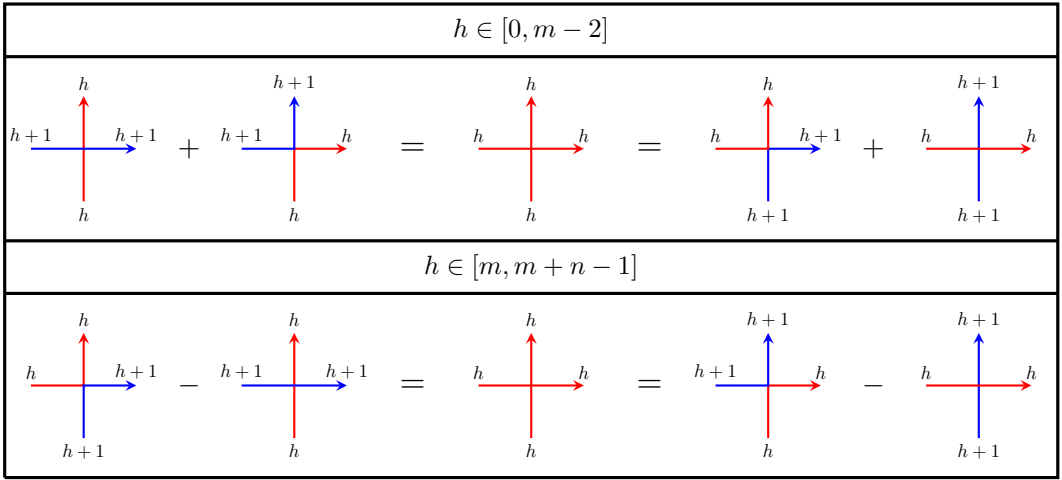 

\begin{proof}
	
	Since $|\mathcal{D}| = 1$, the set $\textbf{z} = \big( z(v) \big)$ consists of a single element, which we abbreviate $z = z(v)$. Then, for any entrance data $\breve{\mathfrak{E}} = (\breve{e}_1, \breve{e}_2)$ and exit data $\breve{\mathfrak{F}} = (\breve{f}_1, \breve{f}_2)$ on $v$, we have $Z_{\mathcal{D}}^{(m; n)} (\breve{\mathfrak{E}}; \breve{\mathfrak{F}} \boldsymbol{\mid} \textbf{z}) = R_z (\breve{e}_2, \breve{e}_1; \breve{f}_1, \breve{f}_2)$. Moreover, since $\theta_{\mathbb{J}} (\breve{e}_i) = \theta_{\mathbb{J}} (e_i)$ and $\theta_{\mathbb{J}} (\breve{f}_i) = \theta_{\mathbb{J}} (f_i)$ for each $i \in \{ 1, 2 \}$, we have $\theta_{\mathbb{J}} (e_i) = h = \theta_{\mathbb{J}} (f_i)$ for each $i$. In particular, the right side of \eqref{sumefzsum} is equal to $R_z^{(m'; n')} (h, h; h, h)$. Furthermore, by \eqref{jmn}, $\mathbb{J}$ separates $[0, m - 1]$ from $[m, m + n - 1]$, so \eqref{jh1} implies that $h \ne m - 1$. We will analyze the cases $h < m - 1$ and $h \ge m$ separately. 
	
	Let us first suppose $h < m - 1$, in which case $(m', n') = (m - 1, n)$. Then, on the left side of \eqref{sumefzsum} we have $\breve{e}_i = e_i$ for each $i$ (since $e_i \le h + 1 \le m - 1$), and we sum over each $\breve{f}_i \in \{ h, h + 1 \}$ (since $\theta_{\mathbb{J}} (\breve{f}_i) = \theta_{\mathbb{J}} (f_i) = h$). Also, as mentioned above, the right side of \eqref{sumefzsum} is $R_z^{(m'; n')} (h, h; h, h) = 1$, where the equality holds by the first statement of \eqref{rzi} (as $h < m - 1 = m'$). Additionally $\inv (\breve{\mathfrak{E}}; J_i) = 0 = \inv (\mathfrak{F}; J_i)$ for each $i \in [m', m' + n' - 1]$, since each entry of $\breve{\mathfrak{E}}$ and $\mathfrak{F}$ is $h$ or $h + 1$, and neither is contained in $J_i$ for $i \ge m'$ (as $J_i \subseteq [m, m + n - 1]$ if $i \ge m'$). 
	
	Hence, to show \eqref{sumefzsum} for $h < m - 1$ we must verify
	\begin{flalign}
	\label{rzf1f2hh12}
	\displaystyle\sum_{\breve{f}_1, \breve{f}_2 \in \{ h, h + 1 \}} R_z^{(m; n)} (e_2, e_1; \breve{f}_1, \breve{f}_2) = R_z^{(m'; n')} \big( h, h; h, h \big) = 1.
	\end{flalign}
	
	\noindent This follows by using \eqref{rzij} and \eqref{rzi} to check the four cases for $(e_1, e_2)$ (equal to either $(h, h)$, $(h + 1, h + 1)$, $(h, h + 1)$, or $(h + 1, h)$) individually. Indeed, for $h \le m - 2$ these definitions imply
	\begin{flalign}
	\label{rzrzsum1} 
	\begin{aligned}
	& R_z^{(m; n)} (h, h; h, h) = 1 = R_z^{(m; n)} (h + 1, h + 1; h + 1, h + 1); \\
	& R_z^{(m; n)} (h, h + 1; h, h + 1) + R_z^{(m; n)} (h, h + 1; h + 1, h) = 1; \\
	& R_z^{(m; n)} (h + 1, h; h, h + 1) + R_z^{(m; n)} (h + 1, h; h + 1, h) = 1;
	\end{aligned}  
	\end{flalign} 
	
	\noindent see the top of \Cref{vertexsum1} for a depiction. This confirms \eqref{rzf1f2hh12} and thus the lemma for $h < m - 1$.

	So, let us suppose instead that $h \ge m$, in which case $(m', n') = (m, n - 1)$. Then, on the left side of \eqref{sumefzsum} we have $\breve{f}_i = f_i$ for each $i$ (since $f_i \ge h \ge m$), and we sum over each $\breve{e}_i \in \{ h, h + 1 \}$. Also, as mentioned above, the right side of \eqref{sumefzsum} is $R_z^{(m; n)} (h, h; h, h) = \frac{z - q}{1 - qz}$, where the equality follows from the second statement of \eqref{rzi}. Additionally, $\inv (\breve{\mathfrak{E}}; J_i) = 0 = \inv (\mathfrak{F}; J_i)$ for $i \ne h$, but since now $h \in [m', m' + n' - 1]$ we must take into account the facts that $\inv (\breve{\mathfrak{E}}; J_h) = \textbf{1}_{\breve{e}_1 > \breve{e}_2}$ and $\inv (\mathfrak{F}; J_h) = \textbf{1}_{f_1 > f_2}$. 
	
	Hence, to show \eqref{sumefzsum} for $h \ge m$ we must verify
	\begin{flalign}
	\label{rzf1f2hh22}
	\displaystyle\sum_{\breve{e}_1, \breve{e}_2 \in \{ h, h + 1 \}} (-1)^{\textbf{1}_{\breve{e}_1 > \breve{e}_2} - \textbf{1}_{f_1 > f_2}} R_z^{(m; n)} (\breve{e}_2, \breve{e}_1; f_1, f_2) = R_z^{(m'; n')} \big( h, h; h, h \big) = \displaystyle\frac{z - q}{1 - qz}.
	\end{flalign}
	
	\noindent As in the case $h < m - 1$, this follows by using \eqref{rzij} and \eqref{rzi} to check the four cases for $(f_1, f_2)$ individually. Indeed, for $h \ge m$ these definitions imply
	\begin{flalign}
	\label{rzrzsum2} 
	\begin{aligned}
	& R_z^{(m'; n')} (h, h; h, h) = \displaystyle\frac{z - q}{1 - qz} = R_z^{(m'; n')} (h + 1, h + 1; h + 1, h + 1); \\
	& R_z^{(m'; n')} (h + 1, h; h, h + 1) - R_z^{(m'; n')} (h, h + 1; h, h + 1) = \displaystyle\frac{z - q}{1 - qz}; \\
	& R_z^{(m'; n')} (h, h + 1; h + 1, h) - R_z^{(m'; n')} (h + 1, h; h + 1, h) = \displaystyle\frac{z - q}{1 - qz};
	\end{aligned} 
	\end{flalign} 
	
	\noindent see the bottom of \Cref{vertexsum1}. This confirms \eqref{rzf1f2hh22} and thus the lemma when $h \ge m$.
\end{proof}

We can now deduce the following corollary verifying \Cref{zsumef} when $|\mathcal{D}| = 1$. 

\begin{cor}
	
	\label{sumefzd1} 
	
	If $|\mathcal{D}| = 1$, then \Cref{zsumef} holds.
\end{cor}

\begin{proof}
	
	Since $|\mathcal{D}| = 1$, we have $k = 1$. Thus, arrow conservation implies that any summand supported by the left side of \eqref{sumefzsum} must satisfy $\{ \breve{e}_1, \breve{e}_2 \} = \{ \breve{f}_1, \breve{f}_2 \}$ as unordered multi-sets. Since $\breve{e}_i = e_i$ whenever $e_i \in [0, m - 1]$ and $\breve{f}_i = f_i$ whenever $f_i \in [m, m + n - 1]$, these unordered multi-sets are determined by $(\mathfrak{E}, \mathfrak{F})$. So, let $\{ \breve{e}_1, \breve{e}_2 \} = \{ a, b \} = \{ \breve{f}_1, \breve{f}_2 \}$, where $a$ and $b$ are fixed by $(\mathfrak{E}, \mathfrak{F})$ (and, in particular, are not summed over on the left side of \eqref{sumefzsum}).
	
	Now recall from the explicit forms \eqref{rzij} and \eqref{rzi} of the weights $R_z (i, j; i', j')$ that they only depend on whether $i, j, i', j'$ are bosonic (lie in $[0, m - 1]$) or fermionic (lie in $[m, m + n -  1]$), and on the relative ordering of $(i, j; i', j')$. Thus, we may assume in what follows that $a \in \{ b - 1, b \}$, that is, $a$ and $b$ either coincide or are consecutive. Since no color other than $a$ or $b$ gives rise to a nonzero contribution to either side of \eqref{sumefzsum}, we may further assume that each $J_i$ is a singleton, except for possibly at most one, which would be given by $J_i = \{ a, b \}$. 
	
	This reduces to the case when $m' + n' = m + n$ if either $a = b$ or $a, b$ reside in different $J_i$, or when $m' + n' = m + n - 1$ if $a \ne b$ reside in the same $J_i$. The former scenario was verified as \Cref{mnmn}, and the latter was addressed by \Cref{mnmn1d1}. 
\end{proof}

\section{Proof of \Cref{zsumef}}

\label{ProofSumZ}

In this section we establish \Cref{zsumef}, which will follow from \Cref{sumefzd1}, together with induction on $|\mathcal{D}|$. 

\begin{proof}[Proof of \Cref{zsumef}]
	
	We induct on $|\mathcal{D}|$. If $|\mathcal{D}| = 1$, then the lemma is verified by \Cref{sumefzd1}. So, let us assume in what follows that it holds whenever $|\mathcal{D}| < r$ for some integer $r > 1$, and we will show it also holds whenever $|\mathcal{D}| = r$. 
	
	To that end, fix an east-south domain $\mathcal{D} = \mathcal{D} (\textbf{p}, \textbf{p}') \subset \mathbb{Z}^2$ with $|\mathcal{D}| = r$. Then, there exists an east-south domain $\widetilde{\mathcal{D}} = \mathcal{D} (\widetilde{\textbf{p}}, \widetilde{\textbf{p}}') \subset \mathcal{D}$ obtained from $\mathcal{D}$ by removing one vertex $v \in \textbf{p}$ on its southwest boundary; in particular, $|\widetilde{\mathcal{D}}| = |\mathcal{D}| - 1$.
	
	First suppose that $\textbf{p} \ne \textbf{p}'$, in which case we may choose $v \notin \textbf{p} \cap \textbf{p}'$. Then, the northeast boundary $\widetilde{\textbf{p}}'$ of $\widetilde{\mathcal{D}}$ is $\textbf{p}'$; the southwest boundary $\widetilde{\textbf{p}}$ of $\widetilde{\mathcal{D}}$ is obtained from $\textbf{p}$ by replacing its corner at $v$, which is formed by a vertical edge followed by a horizontal one, with a corner formed by a horizontal edge followed by a vertical one. Let $v$ be the $K$-th vertex in $\textbf{p}$, for some  index $K \in [1, k]$. So, under entrance data $\mathfrak{E} = (e_1, e_2, \ldots , e_{k + 1})$, arrows of colors $e_K$ and $e_{K + 1}$ horizontally and vertically enter through $v$, respectively. We refer to \Cref{domaindpathsd} for a depiction, where there $K = 9$.

\begin{figure}
	
	\begin{center}		
		
		\begin{tikzpicture}[
		>=stealth,
		auto,
		style={
			scale = .7
		}
		]
				
		\draw[dashed] (1, 8) -- (1, 7) -- (3, 7);
		\draw[dashed] (2, 8) -- (2, 6) -- (3 ,6);
		\draw[dashed] (3, 4) -- (4, 4) -- (4, 3);
		\draw[dashed] (4, 5) -- (4, 4) -- (5, 4);
		
		\draw[] (-.8, 8) circle[radius = 0] node[left, black, scale = .7]{$\breve{e}_1$};
		\draw[] (-.8, 7) circle[radius = 0] node[left, black, scale = .7]{$\breve{e}_2$};
		\draw[] (0, 6.5) circle[radius = 0] node[left, black, scale = .7]{$\breve{e}_3$};
		\draw[] (.5, 6) circle[radius = 0] node[above, black, scale = .7]{$\breve{e}_4$};
		\draw[] (1, 5.5) circle[radius = 0] node[left, black, scale = .7]{$\breve{e}_5$};
		\draw[] (1.5, 5) circle[radius = 0] node[above, black, scale = .7]{$\breve{e}_6$};
		\draw[] (2, 4.5) circle[radius = 0] node[left, black, scale = .7]{$\breve{e}_7$};
		\draw[] (2.5, 4) circle[radius = 0] node[above, black, scale = .7]{$\breve{e}_8$};
		\draw[] (2.5, 3) circle[radius = 0] node[above, black, scale = .7]{$\breve{e}_K$};
		\draw[] (3, 2.2) circle[radius = 0] node[below, black, scale = .7]{$\breve{e}_{K + 1}$};
		\draw[] (4, 2.2) circle[radius = 0] node[below, black, scale = .7]{$\breve{e}_{11}$};
		\draw[] (5, 2.2) circle[radius = 0] node[below, black, scale = .7]{$\breve{e}_{12}$};
		
		\draw[] (0, 8.8) circle[radius = 0] node[above, black, scale = .7]{$\breve{f}_1$};
		\draw[] (1, 8.8) circle[radius = 0] node[above, black, scale = .7]{$\breve{f}_2$};
		\draw[] (2, 8.8) circle[radius = 0] node[above, black, scale = .7]{$\breve{f}_3$};
		\draw[] (3, 8.8) circle[radius = 0] node[above, black, scale = .7]{$\breve{f}_4$};
		\draw[] (3.4, 8) circle[radius = 0] node[above, black, scale = .7]{$\breve{f}_5$};
		\draw[] (3.4, 7) circle[radius = 0] node[above, black, scale = .7]{$\breve{f}_6$};
		\draw[] (3.4, 6) circle[radius = 0] node[above, black, scale = .7]{$\breve{f}_7$};
		\draw[] (4, 5.4) circle[radius = 0] node[left, black, scale = .7]{$\breve{f}_8$};
		\draw[] (5, 5.4) circle[radius = 0] node[left, black, scale = .7]{$\breve{f}_9$};
		\draw[] (5.8, 5) circle[radius = 0] node[right, black, scale = .7]{$\breve{f}_{10}$};
		\draw[] (5.8, 4) circle[radius = 0] node[right, black, scale = .7]{$\breve{f}_{11}$};
		\draw[] (5.8, 3) circle[radius = 0] node[right, black, scale = .7]{$\breve{f}_{12}$};
		
		\draw[thick, blue, ->] (-.8, 8) -- (0, 8);
		\draw[thick, red, ->] (-.8, 7) -- (0, 7);
		\draw[thick, red, ->] (.2, 6) -- (1, 6);
		\draw[thick, green, ->] (1.2, 5) -- (2, 5);
		\draw[thick, green, ->] (2.2, 4) -- (3, 4);
		\draw[thick, blue, ->] (2.2, 3) -- (2.8, 3);
		
		\draw[thick, red, ->] (0, 6.2) -- (0, 7);
		\draw[thick, green, ->] (1, 5.2) -- (1, 6);
		\draw[thick, green, ->] (2, 4.2) -- (2, 5);
		\draw[thick, orange, ->] (3, 2.2) -- (3, 2.8);
		\draw[thick, blue, ->] (4, 2.2) -- (4, 3);
		\draw[thick, blue, ->] (5, 2.2) -- (5, 3);
		
		\draw[thick, red, ->] (3, 8) -- (3.8, 8);
		\draw[thick, red, ->] (3, 7) -- (3.8, 7);
		\draw[thick, blue, ->] (3, 6) -- (3.8, 6);
		\draw[thick, blue, ->] (5, 5) -- (5.8, 5);
		\draw[thick, blue, ->] (5, 4) -- (5.8, 4);
		\draw[thick, orange, ->] (5, 3) -- (5.8, 3);
		
		\draw[thick, red, ->] (0, 8) -- (0, 8.8);
		\draw[thick, green, ->] (1, 8) -- (1, 8.8);
		\draw[thick, blue, ->] (2, 8) -- (2, 8.8);
		\draw[thick, green, ->] (3, 8) -- (3, 8.8);
		\draw[thick, green, ->] (4, 5) -- (4, 5.8);
		\draw[thick, green, ->] (5, 5) -- (5, 5.8);
		
		\draw[ultra thick] (0, 8) -- (0, 7) -- (1, 7) -- (1, 6) -- (2, 6) -- (2, 5) -- (3, 5) -- (3, 3) -- (5, 3);
		\draw[ultra thick] (0, 8) -- (3, 8) -- (3, 7) -- (3, 5) -- (4, 5) -- (5, 5) -- (5, 3);
		
		\draw[fill=white!50!gray] (3, 3) circle[radius = .2];
		
		\draw[] (3.3, 3.3) circle[radius = 0] node[scale = .8]{$v$};
		
		\draw[dashed] (11, 8) -- (11, 7) -- (13, 7);
		\draw[dashed] (12, 8) -- (12, 6) -- (13 ,6);
		\draw[dashed] (14, 5) -- (14, 4) -- (15, 4);
		
		\draw[] (9.2, 8) circle[radius = 0] node[left, black, scale = .7]{$\breve{e}_1$};
		\draw[] (9.2, 7) circle[radius = 0] node[left, black, scale = .7]{$\breve{e}_2$};
		\draw[] (10, 6.5) circle[radius = 0] node[left, black, scale = .7]{$\breve{e}_3$};
		\draw[] (.5, 6) circle[radius = 0] node[above, black, scale = .7]{$\breve{e}_4$};
		\draw[] (11, 5.5) circle[radius = 0] node[left, black, scale = .7]{$\breve{e}_5$};
		\draw[] (11.5, 5) circle[radius = 0] node[above, black, scale = .7]{$\breve{e}_6$};
		\draw[] (12, 4.5) circle[radius = 0] node[left, black, scale = .7]{$\breve{e}_7$};
		\draw[] (12.5, 4) circle[radius = 0] node[above, black, scale = .7]{$\breve{e}_8$};
		\draw[] (13, 3.5) circle[radius = 0] node[left, black, scale = .7]{$\widetilde{e}_K$};
		\draw[] (13.55, 3) circle[radius = 0] node[above, black, scale = .7]{$\widetilde{e}_{K + 1}$};
		\draw[] (14, 2.2) circle[radius = 0] node[below, black, scale = .7]{$\breve{e}_{11}$};
		\draw[] (15, 2.2) circle[radius = 0] node[below, black, scale = .7]{$\breve{e}_{12}$};
		
		\draw[] (10, 8.8) circle[radius = 0] node[above, black, scale = .7]{$\breve{f}_1$};
		\draw[] (11, 8.8) circle[radius = 0] node[above, black, scale = .7]{$\breve{f}_2$};
		\draw[] (12, 8.8) circle[radius = 0] node[above, black, scale = .7]{$\breve{f}_3$};
		\draw[] (13, 8.8) circle[radius = 0] node[above, black, scale = .7]{$\breve{f}_4$};
		\draw[] (13.4, 8) circle[radius = 0] node[above, black, scale = .7]{$\breve{f}_5$};
		\draw[] (13.4, 7) circle[radius = 0] node[above, black, scale = .7]{$\breve{f}_6$};
		\draw[] (13.4, 6) circle[radius = 0] node[above, black, scale = .7]{$\breve{f}_7$};
		\draw[] (14, 5.4) circle[radius = 0] node[left, black, scale = .7]{$\breve{f}_8$};
		\draw[] (15, 5.4) circle[radius = 0] node[left, black, scale = .7]{$\breve{f}_9$};
		\draw[] (15.8, 5) circle[radius = 0] node[right, black, scale = .7]{$\breve{f}_{10}$};
		\draw[] (15.8, 4) circle[radius = 0] node[right, black, scale = .7]{$\breve{f}_{11}$};
		\draw[] (15.8, 3) circle[radius = 0] node[right, black, scale = .7]{$\breve{f}_{12}$};
		
		\draw[thick, blue, ->] (9.2, 8) -- (10, 8);
		\draw[thick, red, ->] (9.2, 7) -- (10, 7);
		\draw[thick, red, ->] (10.2, 6) -- (11, 6);
		\draw[thick, green, ->] (11.2, 5) -- (12, 5);
		\draw[thick, green, ->] (12.2, 4) -- (13, 4);
		\draw[thick, orange, ->] (13.2, 3) -- (14, 3);
		
		\draw[thick, red, ->] (10, 6.2) -- (10, 7);
		\draw[thick, green, ->] (11, 5.2) -- (11, 6);
		\draw[thick, green, ->] (12, 4.2) -- (12, 5);
		\draw[thick, blue, ->] (13, 3.2) -- (13, 4);
		\draw[thick, blue, ->] (14, 2.2) -- (14, 3);
		\draw[thick, blue, ->] (15, 2.2) -- (15, 3);
		
		\draw[thick, red, ->] (13, 8) -- (13.8, 8);
		\draw[thick, red, ->] (13, 7) -- (13.8, 7);
		\draw[thick, blue, ->] (13, 6) -- (13.8, 6);
		\draw[thick, blue, ->] (15, 5) -- (15.8, 5);
		\draw[thick, blue, ->] (15, 4) -- (15.8, 4);
		\draw[thick, orange, ->] (15, 3) -- (15.8, 3);
		
		\draw[thick, red, ->] (10, 8) -- (10, 8.8);
		\draw[thick, green, ->] (11, 8) -- (11, 8.8);
		\draw[thick, blue, ->] (12, 8) -- (12, 8.8);
		\draw[thick, green, ->] (13, 8) -- (13, 8.8);
		\draw[thick, green, ->] (14, 5) -- (14, 5.8);
		\draw[thick, green, ->] (15, 5) -- (15, 5.8);
		
		\draw[ultra thick] (10, 8) -- (10, 7) -- (11, 7) -- (11, 6) -- (12, 6) -- (12, 5) -- (13, 5) -- (13, 4) -- (14, 4) -- (14, 3) -- (15, 3);
		\draw[ultra thick] (10, 8) -- (13, 8) -- (13, 7) -- (13, 5) -- (14, 5) -- (15, 5) -- (15, 3);
		
		\draw[fill=white!50!gray] (13, 3) circle[radius = .2];
		\draw[] (12.7, 2.7) circle[radius = 0] node[scale = .8]{$v$};
		
		\end{tikzpicture}
		
	\end{center}
	
	\caption{\label{domaindpathsd} Shown to the left and right are the domains $\mathcal{D} = \mathcal{D} (\textbf{p}, \textbf{p}')$ and $\widetilde{\mathcal{D}} = \mathcal{D} (\widetilde{\textbf{p}}, \textbf{p}')$, respectively, used in the proof of \Cref{zsumef} if $\textbf{p} \ne \textbf{p}'$. } 	
\end{figure}

Letting $\widetilde{\textbf{z}} = \big( z(v) \big)_{v \in \widetilde{\mathcal{D}}} = \textbf{z} \setminus \big\{ z(v) \big\}$, we have 
\begin{flalign}
\label{zdmnefsum} 
Z_{\mathcal{D}}^{(m; n)} (\breve{\mathfrak{E}}; \breve{\mathfrak{F}} \boldsymbol{\mid} \textbf{z}) = \displaystyle\sum_{\widetilde{e}_K = 0}^{m + n - 1} \displaystyle\sum_{\widetilde{e}_{K + 1} = 0}^{m + n - 1} Z_{\{ v \}}^{(m; n)} \big( (\breve{e}_K, \breve{e}_{K + 1}); (\widetilde{e}_K, \widetilde{e}_{K + 1}) \boldsymbol{\mid} z(v) \big) Z_{\widetilde{\mathcal{D}}}^{(m; n)} (\widetilde{\mathfrak{E}}; \breve{\mathfrak{F}} \boldsymbol{\mid} \widetilde{\textbf{z}}),
\end{flalign}

\noindent for any sequences $\breve{\mathfrak{E}}$ and $\breve{\mathfrak{F}}$ of indices in $[0, m + n - 1]$; here, $\widetilde{\mathfrak{E}}$ is obtained from $\breve{\mathfrak{E}}$ by replacing $(\breve{e}_K, \breve{e}_{K + 1})$ with $(\widetilde{e}_K, \widetilde{e}_{K + 1})$, so that $\widetilde{\mathfrak{E}} = (\breve{e}_1, \ldots , \breve{e}_{K - 1}, \widetilde{e}_K, \widetilde{e}_{K + 1}, \breve{e}_{K + 2}, \ldots , \breve{e}_{k + 1})$. 

Now we will essentially perform a signed sum of the right side of \eqref{zdmnefsum} over $\breve{\mathfrak{E}}$ and $\breve{\mathfrak{F}}$. To that end, recall that $\mathfrak{E} = (e_1, e_2, \ldots , e_{k + 1})$ and $\mathfrak{F} = (f_1, f_2, \ldots , f_{k + 1})$, and set $\mathfrak{e} = (e_K, e_{K + 1})$. Further fix pairs $\widehat{\mathfrak{e}} = (\widehat{e}_K, \widehat{e}_{K + 1})$ of indices in $[0, m + n - 1]$, and define $\widehat{\mathfrak{E}} = (\widehat{e}_1, \widehat{e}_2, \ldots , \widehat{e}_{k + 1})$ from $\mathfrak{E}$ by replacing $\mathfrak{e}$ with $\widehat{\mathfrak{e}}$, so that $\widehat{\mathfrak{E}} = (e_1, \ldots , e_{k - 1}, \widehat{e}_K, \widehat{e}_{K + 1}, e_{K + 2}, \ldots , e_{k + 1})$. Then, \Cref{sumefzd1} yields 
\begin{flalign}
\label{sumzv1}
\displaystyle\sum_{\breve{\mathfrak{e}}} \displaystyle\sum_{\widetilde{\mathfrak{e}}} Z_{\{ v \}}^{(m; n)} \big( \breve{\mathfrak{e}}; \widetilde{\mathfrak{e}} \boldsymbol{\mid} z(v) \big) \displaystyle\prod_{i = m'}^{m' + n' - 1} (-1)^{\inv (\breve{\mathfrak{e}}, J_i) - \inv (\widehat{\mathfrak{e}}, J_i)} = Z_{ \{v \}}^{(m'; n')} \big( \theta_{\mathbb{J}} (\mathfrak{e}); \theta_{\mathbb{J}} (\widehat{\mathfrak{e}}) \boldsymbol{\mid} z (v) \big),
\end{flalign}

\noindent where we sum over all $\breve{\mathfrak{e}} = (\breve{e}_K, \breve{e}_{K + 1})$ and $\widetilde{\mathfrak{e}} = (\widetilde{e}_K, \widetilde{e}_{K + 1})$ such that $\theta_{\mathbb{J}} (\breve{\mathfrak{e}}) = \theta_{\mathbb{J}} (\mathfrak{e})$; $\theta_{\mathbb{J}} (\widetilde{\mathfrak{e}}) = \theta_{\mathbb{J}} (\widehat{\mathfrak{e}})$; $\breve{e}_i = e_i$ whenever $e_i \in [0, m - 1]$; and $\widetilde{e}_i = \widehat{e}_i$ whenver $\widehat{e}_i \in [m, m + n - 1]$. 

Moreover, since $|\widetilde{\mathcal{D}}| = |\mathcal{D}| - 1 = r - 1$, \eqref{sumefzsum} applies and gives 
\begin{flalign}
\label{sumzv2}
\displaystyle\sum_{\widetilde{\mathfrak{E}}} \displaystyle\sum_{\breve{\mathfrak{F}}} Z_{\widetilde{\mathcal{D}}}^{(m; n)} \big( \widetilde{\mathfrak{E}}; \breve{\mathfrak{F}} \boldsymbol{\mid} \widetilde{\textbf{z}} \big) \displaystyle\prod_{i = m'}^{m' + n' - 1} (-1)^{\inv (\widetilde{\mathfrak{E}}, J_i) - \inv (\mathfrak{F}, J_i)} = Z_{\widetilde{\mathcal{D}}}^{(m'; n')} \big( \theta_{\mathbb{J}} (\widehat{\mathfrak{E}}); \theta_{\mathbb{J}} (\mathfrak{F}) \boldsymbol{\mid} \widetilde{\textbf{z}} \big),
\end{flalign}

\noindent where we sum over all $\widetilde{\mathfrak{E}} = (\widetilde{e}_1, \widetilde{e}_2, \ldots , \widetilde{e}_{k + 1})$ and $\breve{\mathfrak{F}} = (\breve{f}_1, \breve{f}_2, \ldots, \breve{f}_{k + 1})$ such that $\theta_{\mathbb{J}} (\widetilde{\mathfrak{E}}) = \theta_{\mathbb{J}} (\widehat{\mathfrak{E}})$; $\theta_{\mathbb{J}} (\breve{\mathfrak{F}}) = \theta_{\mathbb{J}} (\mathfrak{F})$; $\widetilde{e}_i = \widehat{e}_i$ whenever $\widehat{e}_i \in [0, m - 1]$; and $\breve{f}_i = f_i$ whenever $f_i \in [m, m + n - 1]$. 

We next combine \eqref{sumzv1} and \eqref{sumzv2}. To that end, set $\breve{e}_i = \widetilde{e}_i $ for each $i \notin \{ K, K + 1 \}$, so $\widetilde{\mathfrak{E}} = (\breve{e}_1, \ldots , \breve{e}_{K - 1}, \widetilde{e}_K, \widetilde{e}_{K + 1}, \breve{e}_{K + 2}, \ldots , \breve{e}_{k + 1})$ (as above), and for each $i \in [m', m' + n' - 1]$ observe 
\begin{flalign*}
\inv (\widetilde{\mathfrak{E}}; J_i) + \inv (\breve{\mathfrak{e}}; J_i) - \inv (\widehat{\mathfrak{e}}, J_i) = \inv (\widetilde{\mathfrak{E}}; J_i) + \inv (\breve{\mathfrak{e}}; J_i) - \inv (\widetilde{\mathfrak{e}}, J_i) = \inv (\breve{\mathfrak{E}}, J_i),
\end{flalign*} 

\noindent where the first equality follows from the fact that $\widetilde{e}_i = \widehat{e}_i$ whenever $\widetilde{e}_i \in [m, m + n - 1]$ and the second follows from \eqref{ijsum}. This together with \eqref{sumzv1} and \eqref{sumzv2} yields, for a fixed choice of $\theta (\widehat{\mathfrak{e}}) = \big( \theta (\widehat{e}_K), \theta (\widehat{e}_{K + 1}) \big)$, that 
\begin{flalign}
\label{sumef12}
\begin{aligned}
\displaystyle\sum_{\breve{\mathfrak{E}}} \displaystyle\sum_{\breve{\mathfrak{F}}} \displaystyle\sum_{\widetilde{e}_K} \displaystyle\sum_{\widetilde{e}_{K + 1}} Z_{\{ v \}}^{(m; n)} \big( \breve{\mathfrak{e}}; \widetilde{\mathfrak{e}} \boldsymbol{\mid} z(v) \big) & Z_{\widetilde{\mathcal{D}}}^{(m; n)} (\widetilde{\mathfrak{E}}; \breve{\mathfrak{F}} \boldsymbol{\mid} \widetilde{\textbf{z}}) \displaystyle\prod_{i = m'}^{m' + n' - 1} (-1)^{\inv (\breve{\mathfrak{E}}, J_i) - \inv (\mathfrak{F}, J_i)} \\
& = Z_{ \{v \}}^{(m'; n')} \big( \theta_{\mathbb{J}} (\mathfrak{e}); \theta_{\mathbb{J}} (\widehat{\mathfrak{e}}) \boldsymbol{\mid} z (v) \big) Z_{\widetilde{\mathcal{D}}}^{(m'; n')} \big( \theta_{\mathbb{J}} (\widehat{\mathfrak{E}}); \theta_{\mathbb{J}} (\mathfrak{F}) \boldsymbol{\mid} \widetilde{\textbf{z}} \big).
\end{aligned} 
\end{flalign}

\noindent Here, we sum over all $\breve{\mathfrak{E}}$ and $\breve{\mathfrak{F}}$ such that $\theta_{\mathbb{J}} (\breve{\mathfrak{E}}) = \theta_{\mathbb{J}} (\mathfrak{E})$; $\theta_{\mathbb{J}} (\breve{\mathfrak{F}}) = \theta_{\mathbb{J}} (\mathfrak{F})$; $\breve{e}_i = e_i$ whenever $e_i \in [0, m - 1]$; and $\breve{f}_i = f_i$ whenever $f_i \in [m, m + n - 1]$. We further sum over all indices $\widetilde{e}_K$ and $\widetilde{e}_{K + 1}$ in $[0, m + n - 1]$ such that $\theta_{\mathbb{J}} (\widetilde{e}_K) = \theta_{\mathbb{J}} (\widehat{e}_K)$ and $\theta_{\mathbb{J}} (\widetilde{e}_{K + 1}) = \theta_{\mathbb{J}} (\widehat{e}_{K + 1})$. To combine \eqref{sumzv1} and \eqref{sumzv2} in this way, we have used the fact that the following hold for any index $i \in \{ K, K + 1 \}$. The color $\widetilde{e}_i$ is summed over in \eqref{sumzv1} if and only if it is fixed in \eqref{sumzv2} (as both hold if $\widehat{e}_i \in [0, m - 1]$), and that it is fixed in \eqref{sumzv1} if and only if it is summed over in \eqref{sumzv2} (both hold if $\widehat{e}_i \in [m, m + n - 1]$). 

Summing \eqref{sumef12} over all choices of $\theta (\widehat{\mathfrak{e}}) = \big( \theta_{\mathbb{J}} (\widehat{e}_K), \theta_{\mathbb{J}} (\widehat{e}_{K + 1}) \big)$ then yields
\begin{flalign}
\label{2sumef12} 
\begin{aligned}
\displaystyle\sum_{\breve{\mathfrak{E}}} \displaystyle\sum_{\breve{\mathfrak{F}}} \displaystyle\sum_{\widetilde{e}_K = 0}^{m + n - 1} \displaystyle\sum_{\widetilde{e}_{K + 1} = 0}^{m + n - 1} & Z_{\{ v \}}^{(m; n)} \big( \breve{\mathfrak{e}}; \widetilde{\mathfrak{e}} \boldsymbol{\mid} z(v) \big) Z_{\widetilde{\mathcal{D}}}^{(m; n)} (\widetilde{\mathfrak{E}}; \breve{\mathfrak{F}} \boldsymbol{\mid} \widetilde{\textbf{z}}) \displaystyle\prod_{i = m'}^{m' + n' - 1} (-1)^{\inv (\breve{\mathfrak{E}}, J_i) - \inv (\mathfrak{F}, J_i)} \\
& = \displaystyle\sum_{\theta_{\mathbb{J}} (\widehat{e}_K)} \displaystyle\sum_{\theta_{\mathbb{J}} (\widehat{e}_{K + 1})} Z_{ \{v \}}^{(m'; n')} \big( \theta_{\mathbb{J}} (\mathfrak{e}); \theta_{\mathbb{J}} (\widehat{\mathfrak{e}}) \boldsymbol{\mid} z (v) \big) Z_{\widetilde{\mathcal{D}}}^{(m'; n')} \big( \theta_{\mathbb{J}} (\widehat{\mathfrak{E}}); \theta_{\mathbb{J}} (\mathfrak{F}) \boldsymbol{\mid} \widetilde{\textbf{z}} \big),
\end{aligned} 
\end{flalign} 

\noindent where on the left side $\breve{\mathfrak{E}}$ and $\breve{\mathfrak{F}}$ are summed as previously, and on the right side $\theta_{\mathbb{J}} (\widehat{e}_K)$ and $\theta_{\mathbb{J}} (\widehat{e}_{K + 1})$ are both summed over $[0, m' + n' - 1]$. Then \eqref{sumefzsum} follows from applying \eqref{zdmnefsum} on both sides of \eqref{2sumef12}, thereby verifying the proposition if $\textbf{p} \ne \textbf{p}'$. 

Hence, let us instead assume that $\textbf{p} = \textbf{p}'$, in which case $\mathcal{D} = \textbf{p} = \textbf{p}'$. Let $v$ denote the common ending vertex of $\textbf{p}$ and $\textbf{p}'$; define the domain $\widetilde{\mathcal{D}} = \mathcal{D} \setminus \{ v \} = \textbf{p} \setminus \{ v \}$; and set $\widetilde{\textbf{z}} = \textbf{z} \setminus \big\{ z (v) \big\}$. Then, any entrance data $\breve{\mathfrak{E}} = (\breve{e}_1, \breve{e}_2, \ldots , \breve{e}_{k + 1})$ and exit data $\breve{\mathfrak{F}} = (\breve{f}_1, \breve{f}_2, \ldots , \breve{f}_{k + 1})$ on $\mathcal{D}$ fix the colors of three among the four arrows adjacent to $v$; the fourth is then determined by arrow conservation. This gives rise to entrance data $\mathfrak{e}' = (e_k', e_{k + 1}')$ and exit data $\mathfrak{f}' = (f_k', f_{k + 1}')$ on $v$, which satisfies $(e_{k + 1}', f_{k + 1}') = (\breve{e}_{k + 1}, \breve{f}_{k + 1})$. It also satisfies either $e_k' = \breve{e}_k$ if the last edge of $\textbf{p}$ is vertical, or $f_k' = \breve{f}_k$ if the last edge of $\textbf{p}$ is horizontal. Let us assume in what follows that the former scenario $e_k' = \breve{e}_k$ holds, for the proof in the latter is entirely analogous. Then $(\breve{\mathfrak{E}}, \breve{\mathfrak{F}})$ induces entrance data $\widetilde{\mathfrak{E}} = (\widetilde{e}_1, \widetilde{e}_2, \ldots , \widetilde{e}_k)$ and exit data $\widetilde{\mathfrak{F}} = (\widetilde{f}_1, \widetilde{f}_2, \ldots , \widetilde{f}_k)$ on $\widetilde{\mathcal{D}}$, which are in particular given by $\widetilde{\mathfrak{E}} = (\breve{e}_1, \breve{e}_2, \ldots , \breve{e}_{k - 1}, f_k')$ and $\widetilde{\mathfrak{F}} = (\breve{f}_1, \breve{f}_2, \ldots , \breve{f}_k)$. Under this notation, we have
\begin{flalign}
\label{efzproduct} 
Z_{\mathcal{D}}^{(m; n)} (\breve{\mathfrak{E}}; \breve{\mathfrak{F}} \boldsymbol{\mid} \textbf{z}) = Z_{\{ v \} }^{(m; n)} \big( \mathfrak{e}'; \mathfrak{f}' \boldsymbol{\mid} z(v) \big) Z_{\widetilde{\mathcal{D}}}^{(m; n)} \big( \widetilde{\mathfrak{E}}; \widetilde{\mathfrak{F}} \boldsymbol{\mid} \widetilde{\textbf{z}}).
\end{flalign} 

\noindent We refer to \Cref{pathd} for a depiction.  

\begin{figure}
	
	\begin{center}		
		
		\begin{tikzpicture}[
		>=stealth,
		auto,
		style={
			scale = .8
		}
		]

		\draw[] (-1.4, 2) circle[radius = 0] node[above, black, scale = .7]{$\breve{e}_1$};
		\draw[] (-1, 1.2) circle[radius = 0] node[below, black, scale = .7]{$\breve{e}_2$};
		\draw[] (0, 1.2) circle[radius = 0] node[below, black, scale = .7]{$\breve{e}_3$};
		\draw[] (1, 1.5) circle[radius = 0] node[left, black, scale = .7]{$\breve{e}_4$};
		\draw[] (1.6, 1) circle[radius = 0] node[above, black, scale = .7]{$\breve{e}_5$};
		\draw[] (2, .5) circle[radius = 0] node[left, black, scale = .7]{$\breve{e}_6$};
		\draw[] (3, -.8) circle[radius = 0] node[below, black, scale = .7]{$\breve{e}_{k + 1}$};
		\draw[] (2.2, 0) circle[radius = 0] node[below, black, scale = .7]{$\breve{e}_k$};
		
		\draw[] (-1, 2.8) circle[radius = 0] node[above, black, scale = .7]{$\breve{f}_1$};
		\draw[] (0, 2.8) circle[radius = 0] node[above, black, scale = .7]{$\breve{f}_2$};
		\draw[] (1, 2.8) circle[radius = 0] node[above, black, scale = .7]{$\breve{f}_3$};
		\draw[] (2, 2.8) circle[radius = 0] node[above, black, scale = .7]{$\breve{f}_4$};
		\draw[] (2.5, 2) circle[radius = 0] node[above, black, scale = .7]{$\breve{f}_5$};
		\draw[] (3, 1.5) circle[radius = 0] node[left, black, scale = .7]{$\breve{f}_6$};
		\draw[] (3.8, 1) circle[radius = 0] node[right, black, scale = .7]{$\breve{f}_k$};
		\draw[] (3.8, 0) circle[radius = 0] node[right, black, scale = .7]{$\breve{f}_{k + 1}$};
		
		\draw[thick, red, ->] (-1.8, 2) -- (-1, 2);
		\draw[thick, blue, ->] (1.2, 1) -- (2, 1);
		
		\draw[thick, blue, ->] (-1, 1.2) -- (-1, 2);
		\draw[thick, green, ->] (0, 1.2) -- (0, 2);
		\draw[thick, blue, ->] (1, 1.2) -- (1, 2);
		\draw[thick, red, ->] (2, .2) -- (2, 1);
		
		\draw[thick, blue, ->] (2, 2) -- (2.8, 2);
		\draw[thick, blue, ->] (3, 1) -- (3.8, 1);
		
		\draw[thick, blue, ->] (-1, 2) -- (-1, 2.8);
		\draw[thick, red, ->] (0, 2) -- (0, 2.8);
		\draw[thick, green, ->] (1, 2) -- (1, 2.8);
		\draw[thick, red, ->] (2, 2) -- (2, 2.8);
		\draw[thick, green, ->] (3, 1) -- (3, 1.8);
		
		\draw[thick, orange, ->] (2.2, 0) -- (2.8, 0);
		\draw[thick, orange, ->] (3.2, 0) -- (3.8, 0);
		\draw[thick, green, ->] (3, -.8) -- (3, -.2);

		\draw[ultra thick] (-1, 2) -- (0, 2) -- (-1, 2) -- (2, 2) -- (2, 1) -- (3, 1) -- (3, 0);
		
		\draw[fill=white!50!gray] (3, 0) circle[radius = .2];
		
		\draw[] (2.7, -.3) circle[radius = 0] node[scale = .8]{$v$};

		\draw[] (7.6, 2) circle[radius = 0] node[above, black, scale = .7]{$\breve{e}_1$};
		\draw[] (8, 1.2) circle[radius = 0] node[below, black, scale = .7]{$\breve{e}_2$};
		\draw[] (9, 1.2) circle[radius = 0] node[below, black, scale = .7]{$\breve{e}_3$};
		\draw[] (10, 1.5) circle[radius = 0] node[left, black, scale = .7]{$\breve{e}_4$};
		\draw[] (10.6, 1) circle[radius = 0] node[above, black, scale = .7]{$\breve{e}_5$};
		\draw[] (11, .5) circle[radius = 0] node[left, black, scale = .7]{$\breve{e}_6$};
		\draw[] (12, .5) circle[radius = 0] node[right, black, scale = .7]{$\widetilde{e}_k = f_k'$};
		\draw[] (12, -.8) circle[radius = 0] node[below, black, scale = .7]{$e_{k + 1}' = \breve{e}_{k + 1}$};
		\draw[] (11, 0) circle[radius = 0] node[below, black, scale = .7]{$e_k' = \breve{e}_k$};
		
		\draw[] (8, 2.8) circle[radius = 0] node[above, black, scale = .7]{$\breve{f}_1$};
		\draw[] (9, 2.8) circle[radius = 0] node[above, black, scale = .7]{$\breve{f}_2$};
		\draw[] (10, 2.8) circle[radius = 0] node[above, black, scale = .7]{$\breve{f}_3$};
		\draw[] (11, 2.8) circle[radius = 0] node[above, black, scale = .7]{$\breve{f}_4$};
		\draw[] (11.5, 2) circle[radius = 0] node[above, black, scale = .7]{$\breve{f}_5$};
		\draw[] (12, 1.5) circle[radius = 0] node[left, black, scale = .7]{$\breve{f}_6$};
		\draw[] (12.8, 1) circle[radius = 0] node[right, black, scale = .7]{$\breve{f}_k$};
		\draw[] (12.8, 0) circle[radius = 0] node[right, black, scale = .7]{$f_{k + 1}' = \breve{f}_{k + 1}$};
		
		\draw[thick, red, ->] (7.2, 2) -- (8, 2);
		\draw[thick, blue, ->] (10.2, 1) -- (11, 1);
		
		\draw[thick, blue, ->] (8, 1.2) -- (8, 2);
		\draw[thick, green, ->] (9, 1.2) -- (9, 2);
		\draw[thick, blue, ->] (10, 1.2) -- (10, 2);
		\draw[thick, red, ->] (11, .2) -- (11, 1);
		\draw[thick, green, ->] (12, .2) -- (12, 1);
		
		\draw[thick, blue, ->] (11, 2) -- (11.8, 2);
		\draw[thick, blue, ->] (12, 1) -- (12.8, 1);
		
		\draw[thick, blue, ->] (8, 2) -- (8, 2.8);
		\draw[thick, red, ->] (9, 2) -- (9, 2.8);
		\draw[thick, green, ->] (10, 2) -- (10, 2.8);
		\draw[thick, red, ->] (11, 2) -- (11, 2.8);
		\draw[thick, green, ->] (12, 1) -- (12, 1.8);
		
		\draw[thick, orange, ->] (11.2, 0) -- (11.8, 0);
		\draw[thick, orange, ->] (12.2, 0) -- (12.8, 0);
		\draw[thick, green, ->] (12, -.8) -- (12, -.2);
		
		\draw[ultra thick] (8, 2) -- (9, 2) -- (10, 2) -- (11, 2) -- (11, 1) -- (12, 1);
		
		\draw[fill=white!50!gray] (12, 0) circle[radius = .2];
		\draw[] (12.3, -.3) circle[radius = 0] node[scale = .8]{$v$};
		
		\end{tikzpicture}
		
	\end{center}
	
	\caption{\label{pathd} Shown to the left and right are the domains $\mathcal{D}$ and $\widetilde{\mathcal{D}}$, respectively, used in the proof of \Cref{zsumef} if $\textbf{p} = \textbf{p}'$. } 	
\end{figure} 

 As previously for \eqref{zdmnefsum}, we will essentially perform a signed sum of the right side of \eqref{efzproduct} over $\breve{\mathfrak{E}}$ and $\breve{\mathfrak{F}}$. To that end, observe (using the same reasoning as that implemented in the derivation of \eqref{efzproduct}) that the boundary data $(\mathfrak{E}; \mathfrak{F})$ on $\mathcal{D}$ determines entrance data $\mathfrak{e}'' = (e_k'', e_{k + 1}'')$ and exit data $\mathfrak{f}'' = (f_k'', f_{k + 1}'')$ on $v$ satisfying $e_k'' = e_k$ and $(e_{k + 1}'', f_{k + 1}'') = (e_{k + 1}, f_{k + 1})$. Once again, this fixes entrace data $\widehat{\mathfrak{E}} = (\widehat{e}_1, \widehat{e}_2, \ldots , \widehat{e}_k) = (e_1, e_2, \ldots , e_{k - 1}, f_k'')$ and exit data $\widehat{\mathfrak{F}} = (\widehat{f}_1, \widehat{f}_2, \ldots , \widehat{f}_k) = (f_1, f_2, \ldots , f_k)$ on $\widetilde{\mathcal{D}}$. Then, \Cref{sumefzd1} yields 
\begin{flalign}
\label{sumzv3}
\displaystyle\sum_{\mathfrak{e}'} \displaystyle\sum_{\mathfrak{f}'} Z_{\{ v \}}^{(m; n)} \big( \mathfrak{e}'; \mathfrak{f}' \boldsymbol{\mid} z(v) \big) \displaystyle\prod_{i = m'}^{m' + n' - 1} (-1)^{\inv (\mathfrak{e}', J_i) - \inv (\mathfrak{f}'', J_i)} = Z_{ \{v \}}^{(m'; n')} \big( \theta_{\mathbb{J}} (\mathfrak{e}''); \theta_{\mathbb{J}} (\mathfrak{f}'') \boldsymbol{\mid} z (v) \big),
\end{flalign}

\noindent where we sum over all $\mathfrak{e}' = (e_k', e_{k + 1}')$ and $\mathfrak{f}' = (f_k', f_{k + 1}')$ such that $\theta_{\mathbb{J}} (\mathfrak{e}') = \theta_{\mathbb{J}} (\mathfrak{e}'')$; $\theta_{\mathbb{J}} (\mathfrak{f}') = \theta_{\mathbb{J}} (\mathfrak{f}'')$; $e_i' = e_i''$ whenever $e_i'' \in [0, m - 1]$; and $f_i' = f_i''$ whenever $f_i'' \in [m, m + n - 1]$. 

Moreover, since $|\widetilde{\mathcal{D}}| = |\mathcal{D}| - 1 = r - 1$, \eqref{sumefzsum} applies and gives 
\begin{flalign}
\label{sumzv4}
\displaystyle\sum_{\widetilde{\mathfrak{E}}} \displaystyle\sum_{\widetilde{\mathfrak{F}}} Z_{\widetilde{\mathcal{D}}}^{(m; n)} \big( \widetilde{\mathfrak{E}}; \widetilde{\mathfrak{F}} \boldsymbol{\mid} \widetilde{\textbf{z}} \big) \displaystyle\prod_{i = m'}^{m' + n' - 1} (-1)^{\inv (\widetilde{\mathfrak{E}}, J_i) - \inv (\widehat{\mathfrak{F}}, J_i)} = Z_{\widetilde{\mathcal{D}}}^{(m'; n')} \big( \theta_{\mathbb{J}} (\widehat{\mathfrak{E}}); \theta_{\mathbb{J}} (\widehat{\mathfrak{F}}) \boldsymbol{\mid} \widetilde{\textbf{z}} \big),
\end{flalign}

\noindent where we sum over all $\widetilde{\mathfrak{E}} = (\widetilde{e}_1, \widetilde{e}_2, \ldots , \widetilde{e}_k)$ and $\widetilde{\mathfrak{F}} = (\widetilde{f}_1, \widetilde{f}_2, \ldots, \widetilde{f}_k)$ such that $\theta_{\mathbb{J}} (\widetilde{\mathfrak{E}}) = \theta_{\mathbb{J}} (\widehat{\mathfrak{E}})$; $\theta_{\mathbb{J}} (\widetilde{\mathfrak{F}}) = \theta_{\mathbb{J}} (\widehat{\mathfrak{F}})$; $\widetilde{e}_i = \widehat{e}_i$ whenever $\widehat{e}_i \in [0, m - 1]$; and $\widetilde{f}_i = \widehat{f}_i$ whenever $\widehat{f}_i \in [m, m + n - 1]$. 

We now combine \eqref{sumzv3} and \eqref{sumzv4}. To that end, recall the correspondence between $\breve{\mathfrak{E}}$ and $(\widetilde{\mathfrak{E}}, \mathfrak{e}')$, and between $\breve{\mathfrak{F}}$ and $(\widetilde{\mathfrak{F}}, \mathfrak{f}')$, described above \eqref{efzproduct}. We claim for each $i \in [0, m' + n' - 1]$ that, whenever $(\mathfrak{e}', \mathfrak{f}')$ and $(\widetilde{\mathfrak{E}}, \widetilde{\mathfrak{F}})$ satisfy arrow conservation on $v$ and on $\widetilde{\mathcal{D}}$, respectively, we have
\begin{flalign}
\label{ejifji}
\inv (\breve{\mathfrak{E}}, J_i) - \inv (\widetilde{\mathfrak{E}}; J_i) - \inv (\mathfrak{e}'; J_i) = \inv (\breve{\mathfrak{F}}, J_i) - \inv (\widetilde{\mathfrak{F}}; J_i) - \inv (\mathfrak{f}'; J_i).
\end{flalign} 

\noindent To verify \eqref{ejifji}, fix some $i \in [0, m' + n' - 1]$. We may assume (for notational convenience) that each element of $\breve{\mathfrak{E}}, \widetilde{\mathfrak{E}}, \mathfrak{e}', \breve{\mathfrak{F}}, \widetilde{\mathfrak{F}}, \mathfrak{f}'$ is in $J_i$, for otherwise any elements not in $J_i$ contained in at least one of these sets can be removed without effecting either side of \eqref{ejifji}. Then, since $\widetilde{e}_i = \breve{e}_i$ for each $i \in [1, k - 1]$ and $(e_k', e_{k + 1}') = (\breve{e}_k, \breve{e}_{k + 1})$, \eqref{ijsum} imples that the left side of \eqref{ejifji} is equal to 
\begin{flalign}
\label{ejifji1}
\inv (\breve{\mathfrak{E}}, J_i) - \inv (\widetilde{\mathfrak{E}}; J_i) - \inv (\mathfrak{e}'; J_i) = \displaystyle\sum_{i = 1}^{k - 1} ( \textbf{1}_{\breve{e}_i > \breve{e}_k} + \textbf{1}_{\breve{e}_i > \breve{e}_{k + 1}} - \textbf{1}_{\breve{e}_i > \widetilde{e}_k}).
\end{flalign}

\noindent Similarly, since $\widetilde{f}_i = \breve{f}_i$ for each $i \in [1, k]$ and $(f_k', f_{k + 1}') = (\widetilde{e}_k, \breve{f}_{k + 1})$, the right side of \eqref{ejifji} equals 
\begin{flalign}
\label{ejifji2} 
\inv (\breve{\mathfrak{F}}, J_i) - \inv (\widetilde{\mathfrak{F}}; J_i) - \inv (\mathfrak{f}'; J_i) = \displaystyle\sum_{i = 1}^k \textbf{1}_{\breve{f}_i > \breve{f}_{k + 1}} - \textbf{1}_{\widetilde{e}_k > \breve{f}_{k + 1}}.
\end{flalign}

\noindent Next, arrow conservation for $(\mathfrak{e}', \mathfrak{f}')$ implies $(e_k', e_{k + 1}') \in \big\{ (f_k', f_{k + 1}'), (f_{k + 1}', f_k') \big\}$, so $(\breve{e}_k, \breve{e}_{k + 1})$ is either $(\widetilde{e}_k, \breve{f}_{k + 1})$ or $(\breve{f}_{k + 1}, \widetilde{e}_k)$. Moreover, arrow conservation for $(\widetilde{\mathfrak{E}}, \widetilde{\mathfrak{F}})$ yields $\{ \breve{e}_1, \breve{e}_2, \ldots , \breve{e}_{k - 1}, \widetilde{e}_k \} = \{ \breve{f}_1, \breve{f}_2, \ldots , \breve{f}_k \}$ as (unordered) multi-sets. From this, it is directly verified that the right sides of \eqref{ejifji1} and \eqref{ejifji2} are equal in both cases $(\breve{e}_k, \breve{e}_{k + 1}) \in \big\{ (\widetilde{e}_k, \breve{f}_{k + 1}), (\breve{f}_{k + 1}, \widetilde{e}_k) \big\}$, establishing \eqref{ejifji}. 

Now, observe that $\inv (\breve{\mathfrak{F}}, J_i) = \inv (\mathfrak{F}, J_i)$ holds for each $i \in [m', m' + n' - 1]$, since $J_i \subseteq [m, m + n - 1]$ and each $\breve{f}_j = f_j$ whenever $f_j \in [m, m + n - 1]$. For the same reason, we also have for each $i \in [m', m' + n' - 1]$ that $\inv (\widetilde{\mathfrak{F}}; J_i) = \inv(\widehat{\mathfrak{F}}; J_i)$ and $\inv (\mathfrak{f}'; J_i) = \inv (\mathfrak{f}''; J_i)$. Combining these three equalities with \eqref{ejifji} yields for each $i \in [m', m' + n' - 1]$ that
\begin{flalign*}
\inv (\widetilde{\mathfrak{E}}; J_i) - \inv (\widehat{\mathfrak{F}}; J_i) +  \inv (\mathfrak{e}'; J_i) - \inv (\mathfrak{f}''; J_i) = \inv (\breve{\mathfrak{E}}, J_i) - \inv (\mathfrak{F}, J_i),
\end{flalign*} 

\noindent whenever $(\mathfrak{e}', \mathfrak{f}')$ and $(\widetilde{\mathfrak{E}}, \widetilde{\mathfrak{F}})$ satisfy arrow conservation. This, \eqref{sumzv3}, and \eqref{sumzv4} then together yield 
\begin{flalign}
\label{3sumef12}
\begin{aligned}
\displaystyle\sum_{\widetilde{\mathfrak{E}}} \displaystyle\sum_{\widetilde{\mathfrak{F}}} & \displaystyle\sum_{\mathfrak{e}'} \displaystyle\sum_{\mathfrak{f}'} Z_{\{ v \}}^{(m; n)} \big( \mathfrak{e}'; \mathfrak{f}' \boldsymbol{\mid} z(v) \big) Z_{\widetilde{\mathcal{D}}}^{(m; n)} (\widetilde{\mathfrak{E}}; \widetilde{\mathfrak{F}} \boldsymbol{\mid} \widetilde{\textbf{z}}) \displaystyle\prod_{i = m'}^{m' + n' - 1} (-1)^{\inv (\breve{\mathfrak{E}}, J_i) - \inv (\mathfrak{F}, J_i)} \\
& = Z_{ \{v \}}^{(m'; n')} \big( \theta_{\mathbb{J}} (\mathfrak{e}''); \theta_{\mathbb{J}} (\mathfrak{f}'') \boldsymbol{\mid} z (v) \big) Z_{\widetilde{\mathcal{D}}}^{(m'; n')} \big( \theta_{\mathbb{J}} (\widehat{\mathfrak{E}}); \theta_{\mathbb{J}} (\widehat{\mathfrak{F}}) \boldsymbol{\mid} \widetilde{\textbf{z}} \big) = Z_{\mathcal{D}}^{(m'; n')} \big( \theta_{\mathbb{J}} (\mathfrak{E}); \theta_{\mathbb{J}} (\mathfrak{F}) \boldsymbol{\mid} \textbf{z} \big),
\end{aligned} 
\end{flalign}

\noindent where to deduce the last equality we applied \eqref{efzproduct}. Here, $(\mathfrak{e}', \mathfrak{f}')$ and $(\widetilde{\mathfrak{E}}, \widetilde{\mathfrak{F}})$ are summed as in \eqref{sumzv3} and \eqref{sumzv4}, respectively, with the additional constraint that $\widetilde{e}_k = f_k'$. As in the derivation of \eqref{sumef12}, here we have used the fact that the index $\widetilde{e}_k = f_k'$ common to \eqref{sumzv3} and \eqref{sumzv4} is fixed in one if and only if it is summed over in the other. 

Due to the correpsondences between $\breve{\mathfrak{E}}$ and $(\widetilde{\mathfrak{E}}, \mathfrak{e}')$ and between $\breve{\mathfrak{F}}$ and $(\widetilde{\mathfrak{F}}, \mathfrak{f}')$, we may instead sum the left side of \eqref{3sumef12} over $(\breve{\mathfrak{E}}, \breve{\mathfrak{F}})$. This gives
\begin{flalign}
\label{4sumef12} 
\begin{aligned}
\displaystyle\sum_{\breve{\mathfrak{E}}} \displaystyle\sum_{\breve{\mathfrak{F}}}  Z_{\{ v \}}^{(m; n)} \big( \mathfrak{e}'; \mathfrak{f}' \boldsymbol{\mid} z(v) \big) Z_{\widetilde{\mathcal{D}}}^{(m; n)} (\widetilde{\mathfrak{E}}; \widetilde{\mathfrak{F}} \boldsymbol{\mid} \widetilde{\textbf{z}}) \displaystyle\prod_{i = m'}^{m' + n' - 1} & (-1)^{\inv (\breve{\mathfrak{E}}, J_i) - \inv (\mathfrak{F}, J_i)} \\
& = Z_{\mathcal{D}}^{(m'; n')} \big( \theta_{\mathbb{J}} (\mathfrak{E}); \theta_{\mathbb{J}} (\mathfrak{F}) \boldsymbol{\mid} \textbf{z} \big),
\end{aligned} 
\end{flalign} 

\noindent where we sum over all $\breve{\mathfrak{E}}$ and $\breve{\mathfrak{F}}$ such that $\theta_{\mathbb{J}} (\breve{\mathfrak{E}}) = \theta_{\mathbb{J}} (\mathfrak{E})$; $\theta_{\mathbb{J}} (\breve{\mathfrak{F}}) = \theta_{\mathbb{J}} (\mathfrak{F})$; $\breve{e}_i = e_i$ whenever $e_i \in [0, m - 1]$; and $\breve{f}_i = f_i$ whenever $f_i \in [m, m + n - 1]$. We now deduce \eqref{sumefzsum} from applying \eqref{efzproduct} in \eqref{4sumef12}; this establishes the proposition if $\textbf{p} = \textbf{p}'$. 
\end{proof}

\chapter{Fusion of Weights}

\label{WeightsR}

The vertex weights provided in \Cref{Weights1} are nonzero only if each edge incident to the vertex accommodates at most one colored arrow. We now proceed to remove this condition by applying the fusion procedure, which originated in \cite{ERT}, to these weights by suitably concatenating rows of vertices whose spectral parameters are in a certain geometric progression.

\section{Specialized Rectangular Partition Functions} 

\label{Weightss}

In this section we provide notation for, and establish properties of, certain rectangular partition functions that will be relevant for the fusion procedure. To that end, we begin with the following partition function, which is the special case of \Cref{zef1} when the domain $\mathcal{D}$ there is a $M \times L$ rectangle.

\begin{definition} 
	
	\label{zxy1}
	
	Fix integers $L, M \ge 1$ and sequences of complex numbers $\textbf{x} = (x_1, x_2, \ldots , x_L)$ and $\textbf{y} = (y_1, y_2, \ldots , y_M)$. For any sequences $\mathfrak{A} = (a_1, a_2, \ldots , a_M)$, $\mathfrak{B} = (b_1, b_2, \ldots , b_L)$, $\mathfrak{C} = (c_1, c_2, \ldots , c_M)$, and $\mathfrak{D} = (d_1, d_2, \ldots,  d_L)$ of indices in $[0, m + n - 1]$, define $Z (\mathfrak{A}, \mathfrak{B}; \mathfrak{C}, \mathfrak{D} \boldsymbol{\mid} \textbf{x}, \textbf{y}) = Z^{(m; n)} (\mathfrak{A}, \mathfrak{B}; \mathfrak{C}, \mathfrak{D} \boldsymbol{\mid} \textbf{x}, \textbf{y})$ by setting
	\begin{flalign}
	\label{zabcd}
	Z (\mathfrak{A}, \mathfrak{B}; \mathfrak{C}, \mathfrak{D} \boldsymbol{\mid} \textbf{x}, \textbf{y})  = \displaystyle\sum \displaystyle\prod_{i = 1}^M \displaystyle\prod_{j = 1}^L R_{y_i / x_j} \big( v_{i, j}, u_{i, j}; v_{i, j + 1}, u_{i + 1, j} \big),
	\end{flalign} \index{Z@$Z (\mathfrak{A}, \mathfrak{B}; \mathfrak{C}, \mathfrak{D} \boldsymbol{\mid} \textbf{x}, \textbf{y})$}
	
	\noindent where the sum is over all sequences $(u_{i, j})$ and $(v_{i, j})$ of indices in $[0, m + n - 1]$ such that $v_{k, 1} = a_k$, $u_{1, k} = b_k$, $v_{k, L + 1} = c_k$, and $u_{M + 1, k} = d_k$, for each $k$. 
	
	If there exist $x, y \in \mathbb{C}$ such that $\textbf{x} = (x, qx, \ldots , q^{L - 1} x)$ and $y = (q^{M - 1} y, q^{M - 2} y, \ldots , y)$, then we further denote $Z_{x, y} (\mathfrak{A}, \mathfrak{B}; \mathfrak{C}, \mathfrak{D}) = Z_{x, y}^{(m; n)} (\mathfrak{A}, \mathfrak{B}; \mathfrak{C}, \mathfrak{D}) = Z (\mathfrak{A}, \mathfrak{B}; \mathfrak{C}, \mathfrak{D} \boldsymbol{\mid} \textbf{x}, \textbf{y})$. \index{Z@$Z_{x, y} (\mathfrak{A}, \mathfrak{B}; \mathfrak{C}, \mathfrak{D})$}
	
\end{definition} 

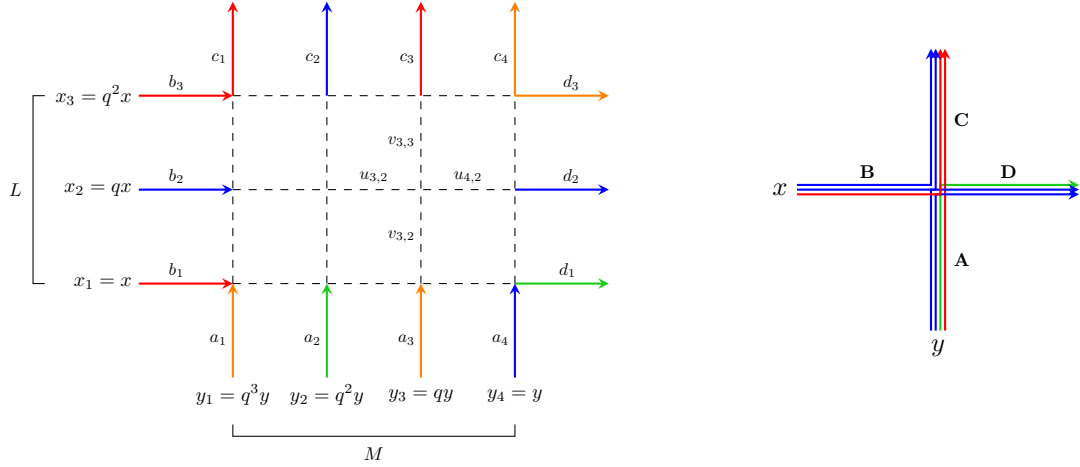
\begin{figure}

	\begin{center}

		\begin{tikzpicture}[
		>=stealth,
		auto,
		style={
			scale = 1.25
		}
		]
		
		\draw[->, thick, red] (-1, 0) -- (0, 0);
		\draw[->, thick, blue] (-1, 1) -- (0, 1);
		\draw[->, thick, red] (-1, 2) -- (0, 2);
		
		\draw[->, thick, green] (3, 0) -- (4, 0);
		\draw[->, thick, blue] (3, 1) -- (4, 1);
		\draw[->, thick, orange] (3, 2) -- (4, 2);
		
		\draw[->, thick, orange] (0, -1) -- (0, 0);
		\draw[->, thick, green] (1, -1) -- (1, 0);
		\draw[->, thick, orange] (2, -1) -- (2, 0);
		\draw[->, thick, blue] (3, -1) -- (3, 0);
		
		\draw[->, thick, red] (0, 2) -- (0, 3);
		\draw[->, thick, blue] (1, 2) -- (1, 3);
		\draw[->, thick, red] (2, 2) -- (2, 3);
		\draw[->, thick, orange] (3, 2) -- (3, 3);
		
		\draw[-, dashed] (0, 0) -- (3, 0);
		\draw[-, dashed] (0, 1) -- (3, 1);
		\draw[-, dashed] (0, 2) -- (3, 2);
		
		\draw[-, dashed] (0, 0) -- (0, 2);
		\draw[-, dashed] (1, 0) -- (1, 2);
		\draw[-, dashed] (2, 0) -- (2, 2);
		\draw[-, dashed] (3, 0) -- (3, 2);

		\draw[] (0, -1) circle[radius = 0] node[below, scale = .75]{$y_1 = q^3 y$};
		\draw[] (1, -1) circle[radius = 0] node[below, scale = .75]{$y_2 = q^2 y$};
		\draw[] (2, -1) circle[radius = 0] node[below = 2, scale = .75]{$y_3 = q y$};
		\draw[] (3, -1) circle[radius = 0] node[below = 2, scale = .75]{$y_4 = y$};
		
		\draw[] (-1, 0) circle[radius = 0] node[left, scale = .75]{$x_1 = x$};
		\draw[] (-1, 1) circle[radius = 0] node[left, scale = .75]{$x_2 = q x$};
		\draw[] (-1, 2) circle[radius = 0] node[left, scale = .75]{$x_3 = q^2 x$};
		
		\draw[] (-.6, 0) circle[radius = 0] node[above, scale = .65]{$b_1$};
		\draw[] (-.6, 1) circle[radius = 0] node[above, scale = .65]{$b_2$};
		\draw[] (-.6, 2) circle[radius = 0] node[above, scale = .65]{$b_3$};
		
		\draw[] (3.56, 0) circle[radius = 0] node[above, scale = .65]{$d_1$};
		\draw[] (3.6, 1) circle[radius = 0] node[above, scale = .65]{$d_2$};
		\draw[] (3.6, 2) circle[radius = 0] node[above, scale = .65]{$d_3$};
		
		\draw[] (0, -.6) circle[radius = 0] node[left, scale = .65]{$a_1$};
		\draw[] (1, -.6) circle[radius = 0] node[left, scale = .65]{$a_2$};
		\draw[] (2, -.6) circle[radius = 0] node[left, scale = .65]{$a_3$};
		\draw[] (3, -.6) circle[radius = 0] node[left, scale = .65]{$a_4$};
		
		\draw[] (0, 2.4) circle[radius = 0] node[left, scale = .65]{$c_1$};
		\draw[] (1, 2.4) circle[radius = 0] node[left, scale = .65]{$c_2$};
		\draw[] (2, 2.4) circle[radius = 0] node[left, scale = .65]{$c_3$};
		\draw[] (3, 2.4) circle[radius = 0] node[left, scale = .65]{$c_4$};

		\draw[] (1.5, 1) circle[radius = 0] node[above, scale = .65]{$u_{3, 2}$};
		\draw[] (2, .5) circle[radius = 0] node[left, scale = .65]{$v_{3, 2}$};
		\draw[] (2.5, 1) circle[radius = 0] node[above, scale = .65]{$u_{4, 2}$};
		\draw[] (2, 1.5) circle[radius = 0] node[left, scale = .65]{$v_{3, 3}$};			
		
		\draw[] (0, -1.5) -- (0, -1.625) -- (3, -1.625) -- (3, -1.5); 
		\draw[] (-2, 0) -- (-2.125, 0) -- (-2.125, 2) -- (-2, 2);
		
		\draw[] (1.5, -1.675) circle[radius = 0] node[below, scale = .7]{$M$};
		\draw[] (-2.175, 1) circle[radius = 0] node[left, scale = .7]{$L$};

		\draw[->, thick, green] (7.525, -.5) -- (7.525, 1.05) -- (9, 1.05);
		\draw[->, thick, red] (7.575, -.5) -- (7.575, 2.5) ;
		\draw[->, thick, blue] (7.475, -.5) -- (7.475, .95) -- (9, .95);
		\draw[->, thick, blue] (7.425, -.5) -- (7.425, 1) -- (9, 1);
		
		\draw[->, thick, red] (6, .95) -- (7.525, .95) -- (7.525, 2.5);
		\draw[->, thick, blue] (6, 1) -- (7.475, 1) -- (7.475, 2.5);
		\draw[->, thick, blue] (6, 1.05) -- (7.425, 1.05) -- (7.425, 2.5);

		\draw[] (6, 1) circle[radius = 0] node[scale = 1, left]{$x$};
		\draw[] (7.5, -.5) circle[radius = 0] node[scale = 1, below]{$y$};
		
		\draw[] (7.6, .25) circle[radius = 0] node[scale = .7, right]{$\textbf{A}$};
		\draw[] (6.75, 1.05) circle[radius = 0] node[scale = .7, above]{$\textbf{B}$};
		\draw[] (7.6, 1.75) circle[radius = 0] node[scale = .7, right]{$\textbf{C}$};
		\draw[] (8.25, 1.05) circle[radius = 0] node[scale = .7, above]{$\textbf{D}$};
		
		\end{tikzpicture}
		
	\end{center}
	
	\caption{\label{zxy} To the left is a diagrammatic interpretation for $Z_{x, y} (\mathfrak{A}, \mathfrak{B}; \mathfrak{C}, \mathfrak{D})$ and to the right is a depiction of a vertex for $\mathcal{R}_{x, y} (\textbf{A}, \textbf{B}; \textbf{C}, \textbf{D})$.}	
\end{figure}

For any ordered set $\mathfrak{Z} = (z_1, z_2, \ldots , z_k)$, let $\overleftarrow{\mathfrak{Z}} = (z_k, z_{k - 1}, \ldots , z_1)$ denote the reverse ordering of $\mathfrak{Z}$.\index{X@$\overleftarrow{\mathscr{X}}$; reverse ordering of $\mathscr{X}$} Then, \eqref{zabcd} denotes the partition function as in \Cref{zef1} for the vertex model on the rectangular domain $[1, M] \times [1, L]$, whose entrance and exit data are given by $\overleftarrow{\mathfrak{B}} \cup \mathfrak{A}$ and $\mathfrak{C} \cup \overleftarrow{\mathfrak{D}}$, respectively. Here, we view $x_k$ and $y_k$ as the rapidities in the $k$-th row (from the bottom) and the $k$-th column (from the left) of this model, respectively, so that the spectral parameter $z(i, j)$ at $(i, j) \in [1, M] \times [1, L]$ is given by $x_j^{-1} y_i$. We refer to the left side of \Cref{zxy} for a depiction. 

\begin{rem}
	
	\label{tzisum} 

For any ordered sequence $\mathcal{T} = (t_1, t_2, \ldots, t_K)$ and indices $1 \le i \le j \le K$, define the restriction $\mathcal{T}_{[i, j]} = (t_i, t_{i + 1}, \ldots , t_j)$. Then observe from the definition \eqref{zabcd} for $Z$ that, for any integer $h \ge 0$, we have
\begin{flalign}
\label{zabcdzabcd}
\begin{aligned}
Z (\mathfrak{A}, \mathfrak{B}; \mathfrak{C}, \mathfrak{D} \boldsymbol{\mid} \textbf{x}, \textbf{y}) & = \displaystyle\sum_{\mathfrak{I}} Z \big( \mathfrak{A}, \mathfrak{B}_{[1, h]}; \mathfrak{I}, \mathfrak{D}_{[1, h]} \boldsymbol{\mid} \textbf{x}_{[1, h]}, \textbf{y} \big) Z \big( \mathfrak{I}, \mathfrak{B}_{[h + 1, L]}; \mathfrak{C}, \mathfrak{D}_{[h + 1, L]} \boldsymbol{\mid} \textbf{x}_{[h + 1, L]}, \textbf{y} \big); \\
Z (\mathfrak{A}, \mathfrak{B}; \mathfrak{C}, \mathfrak{D} \boldsymbol{\mid} \textbf{x}, \textbf{y}) & = \displaystyle\sum_{\mathfrak{J}} Z \big( \mathfrak{A}_{[1, h]}, \mathfrak{B}; \mathfrak{C}_{[1, h]}, \mathfrak{J} \boldsymbol{\mid} \textbf{x}, \textbf{y}_{[1, h]} \big) Z \big( \mathfrak{A}_{[h + 1, M]}, \mathfrak{J}; \mathfrak{C}, \mathfrak{D}_{[h + 1, M]} \boldsymbol{\mid} \textbf{x}, \textbf{y}_{[h + 1, M]} \big),
\end{aligned}
\end{flalign}

\noindent where $\mathfrak{I} = (i_1, i_2, \ldots , i_M)$ and $\mathfrak{J} = (j_1, j_2, \ldots , j_L)$ range over all sequences of indices in $[0, m + n - 1]$. Indeed, in the former equality of \eqref{zabcdzabcd}, $\mathfrak{I}$ is interpreted as the ordered set of colors $(v_{1, h + 1}, v_{2, h + 1}, \ldots , v_{M, h + 1})$ (from \eqref{zabcd}) of vertical arrows intersecting the line $y = h + \frac{1}{2}$; a similar interpretation holds for $\mathfrak{J}$ in the latter statement of \eqref{zabcdzabcd}.

\end{rem} 

Next, for any index set $\mathcal{I} = (i_1, i_2, \ldots , i_{\ell})$ and real number $k \in \mathbb{R}$, let $m_k (\mathcal{I})$ denote the number of indices $j \in [1, \ell]$ such that $i_j = k$.\index{M@$m_k (\mathcal{I})$; multiplicity of $k$ in $\mathcal{I}$} Furthermore, for any sequence $\textbf{I} = (I_0, I_1, \ldots , I_k) \in \mathbb{Z}_{\ge 0}^{k + 1}$, let $\mathcal{M} (\textbf{I})$ denote the family of sequences $\mathcal{I} = (i_1, i_2, \ldots , i_{|\textbf{I}|})$ with elements in $\{ 0, 1, \ldots , k \}$ such that $I_j = m_j (\mathcal{I})$ for each $j \in [0, k]$ (where we recall $|\textbf{I}| = \sum_{j = 0}^k I_j$).\index{X@$\mid$$\textbf{X}$$\mid$}\index{M@$\mathcal{M} (\textbf{I})$; set of sequences with fixed multiplicities}

The fused weights to be considered below will be closely related to the following linear combinations of rectangular partition functions. Here we adopt the convention that symbols of the form $\mathfrak{X} = (x_1, x_2, \ldots , x_K)$ denote ordered sequences of indices (colors) in $[0, m + n - 1]$, as above, and those of the form $\textbf{X} = (X_0, X_1, \ldots , X_{m + n - 1}) \in \mathbb{Z}_{\ge 0}^{m + n}$ (in \Cref{SymmetricBranching} and afterwards we will start the indexing at $1$ instead of at $0$) denote nonnegative compositions. The two will typically be related by stipulating $X_j = m_j (\mathfrak{X})$ for each index $j \in [0, m + n - 1]$, that is, $\mathfrak{X} \in \mathcal{M} (\textbf{X})$.

In the below, we recall for any sequence $\mathscr{X} = (x_1, x_2, \ldots , x_{\ell}) \in \mathbb{R}^{\ell}$ that $\inv (\mathfrak{X})$ denotes the number of index pairs $(i, j) \in [1, \ell]^2$ such that $i < j$ and $x_i > x_j$.\index{I@$\inv$}

\begin{definition}
	
\label{zabcd2} 

Adopt the notation of \Cref{zxy1}, and let $\textbf{A} = (A_0, A_1, \ldots , A_{m + n - 1}) \in \mathbb{Z}_{\ge 0}^{m + n}$ and $\textbf{B} = (B_0, B_1, \ldots , B_{m + n - 1}) \in \mathbb{Z}_{\ge 0}^{m + n}$ denote sequences of nonnegative integers such that $|\textbf{A}| = M$ and $|\textbf{B}| = L$. Define 
\begin{flalign*}
\mathcal{Z}_{x, y} (\textbf{A}, \textbf{B}; \mathfrak{C}, \mathfrak{D}) = \displaystyle\sum_{\substack{\mathfrak{A}\in \mathcal{M}(\textbf{A}) \\ \mathfrak{B} \in \mathcal{M}(\textbf{B})}} q^{\inv (\mathfrak{A}) + \inv (\overleftarrow{\mathfrak{B}})} Z_{x, y} (\mathfrak{A}, \mathfrak{B}; \mathfrak{C}, \mathfrak{D}). 
\end{flalign*} \index{Z@$\mathcal{Z}_{x, y} (\textbf{A}, \textbf{B}; \mathfrak{C}, \mathfrak{D})$}

\end{definition} 

The following lemma provides a $q$-exchangeability property for these weights $\mathcal{Z}_{x, y} (\textbf{A}, \textbf{B}; \mathfrak{C}, \mathfrak{D})$, stating that they multiply by (explicit) powers of $q$ upon permuting $\mathfrak{C}$ and $\mathfrak{D}$, assuming that $\max \{ A_h, B_h \} \le 1$ for each $h \in [m, m + n - 1]$. 
	
\begin{lem} 
	
\label{zabcdcd}

Adopt the notation of \Cref{zabcd2}, and assume that $\max \{ A_h, B_h \} \le 1$ for each index $h \in [m, m + n - 1]$. For any permutations $\mathfrak{C}' = (c_1', c_2', \ldots , c_M')$ and $\mathfrak{D}' = (d_1', d_2', \ldots , d_L')$ of $\mathfrak{C}$ and $\mathfrak{D}$, respectively, we have 
\begin{flalign}
\label{qcd1}
q^{ - \inv (\mathfrak{C}) - \inv (\overleftarrow{\mathfrak{D}})} \mathcal{Z}_{x, y} (\textbf{\emph{A}}, \textbf{\emph{B}}; \mathfrak{C}, \mathfrak{D}) = q^{- \inv (\mathfrak{C}') - \inv (\overleftarrow{\mathfrak{D}}')} \mathcal{Z}_{x, y} (\textbf{\emph{A}}, \textbf{\emph{B}}; \mathfrak{C}', \mathfrak{D}').
\end{flalign}
	
\end{lem} 

\begin{proof} 
	
	Since the symmetric group $\mathfrak{S}_n$ is generated by the transpositions $\{ \mathfrak{s}_i \}_{i \in [1, n - 1]}$ interchanging $(i, i + 1)$,\index{S@$\mathfrak{s}_i$; transposition $(i, i + 1)$} it suffices to establish \eqref{qcd1} assuming that either $(\mathfrak{C}', \mathfrak{D}') = \big( \mathfrak{s}_k (\mathfrak{C}), \mathfrak{D} \big)$ or $(\mathfrak{C}', \mathfrak{D}') = \big( \mathfrak{C}, \mathfrak{s}_k (\mathfrak{D}) \big)$ for some $k$. Let us assume that the former holds, as the proof in the latter case is entirely analogous (by replacing our use of the top diagram in \Cref{rvertical} below with the bottom one there). Thus, we must show that 
	\begin{flalign}
	\label{qcd2}
	\mathcal{Z}_{x, y} (\textbf{A}, \textbf{B}; \mathfrak{C}, \mathfrak{D}) = q^{\sgn (c_k - c_{k + 1})} \mathcal{Z}_{x, y} \big(\textbf{A}, \textbf{B}; \mathfrak{C}', \mathfrak{D} \big).
	\end{flalign}
	
	\noindent We will in fact establish \eqref{qcd2} under the more general hypothesis that only imposes $A_h \le 1$ for each $h \in [m, m + n - 1]$ (while allowing $B_h \ge 2$ for $h \in [m, m + n - 1]$). 
	
	Let us reduce \eqref{qcd2} to the case $M = 2$. Assuming \eqref{qcd2} holds for $M = 2$ (and arbitrary $L \ge 1$), \Cref{zabcd2} and repeated application of the second identity in \eqref{zabcdzabcd} gives for arbitrary $M \ge 2$ that 
	\begin{flalign*}
	\mathcal{Z}_{x, y} (\textbf{A}, \textbf{B}; \mathfrak{C}, \mathfrak{D}) & = \displaystyle\sum_{\substack{\mathfrak{A} \in \mathcal{M} (\textbf{A}) \\ \mathfrak{B} \in \mathcal{M} (\textbf{B})}} q^{\inv (\mathfrak{A}) + \inv (\overleftarrow{\mathfrak{B}})} Z_{x, y} (\mathfrak{A}, \mathfrak{B}; \mathfrak{C}, \mathfrak{D}) \\
	& = \displaystyle\sum_{\substack{\mathfrak{A} \in \mathcal{M} (\textbf{A}) \\ \mathfrak{B} \in \mathcal{M}(\textbf{B})}} \displaystyle\sum_{\mathfrak{J}_1, \mathfrak{J}_2} q^{\inv (\mathfrak{A}) + \inv (\overleftarrow{\mathfrak{B}})} Z_{x, q^{M - k + 1} y} \big( \mathfrak{A}_{[1, k - 1]}, \mathfrak{B}; \mathfrak{C}_{[1, k - 1]}, \mathfrak{J}_1 \big) \\
	& \qquad \times Z_{x, q^{M - k - 1} y} \big( \mathfrak{A}_{[k, k + 1]}, \mathfrak{J}_1; \mathfrak{C}_{[k, k + 1]}, \mathfrak{J}_2 \big) Z_{x, y} \big( \mathfrak{A}_{[k + 2, M]}, \mathfrak{J}_2; \mathfrak{C}_{[k + 2, M]}, \mathfrak{D} \big) \\
	& = \displaystyle\sum_{\substack{\mathfrak{A} \in \mathcal{M} (\textbf{A}) \\ \mathfrak{B}\in \mathcal{M}(\textbf{B})}} \displaystyle\sum_{\mathfrak{J}_1, \mathfrak{J}_2} q^{\inv (\mathfrak{A}) + \inv (\overleftarrow{\mathfrak{B}}) + \sgn (c_k - c_{k + 1})} Z_{x, q^{M - k + 1} y} \big( \mathfrak{A}_{[1, k - 1]}, \mathfrak{B}; \mathfrak{C}_{[1, k - 1]}, \mathfrak{J}_1 \big) \\
	& \qquad \times Z_{x, q^{M - k - 1} y} \Big( \mathfrak{A}_{[k, k + 1]}, \mathfrak{J}_1; \mathfrak{s}_1 \big( \mathfrak{C}_{[k, k + 1]} \big), \mathfrak{J}_2 \Big) Z_{x, y} \big( \mathfrak{A}_{[k + 2, M]}, \mathfrak{J}_2; \mathfrak{C}_{[k + 2, M]}, \mathfrak{D} \big) \\
	& = q^{\sgn (c_k - c_{k + 1})} Z_{x, y} ( \textbf{A}, \textbf{B}; \mathfrak{s}_k (\mathfrak{C}), \mathfrak{D}), 
	\end{flalign*}
	
	\noindent where $\mathfrak{J}_1, \mathfrak{J}_2$ are summed over $\{ 0, 1, \ldots , m + n - 1 \}^L$. This verifies the $M \ge 2$ case of \eqref{qcd2} assuming that the $M = 2$ case holds. So, we may assume in what follows that $M = 2$ and thus that $k = 1$, in which case $\mathfrak{A} = (a_1, a_2)$ and $\mathfrak{C} = (c_1, c_2)$. 
	
	To that end, recalling the definition of $\mathcal{Z}$ from \Cref{zabcd2}, it suffices to show that 
	\begin{flalign}
	\label{qcd3}
	\begin{aligned} 
	q^{\textbf{1}_{a_1 > a_2}} & Z_{x, y} \big( (a_1, a_2), \mathfrak{B}; (c_1, c_2), \mathfrak{D} \big) + q^{\textbf{1}_{a_2 > a_1}} Z_{x, y} \big( (a_2, a_1), \mathfrak{B}; (c_1, c_2), \mathfrak{D} \big) \\
	& = q^{\sgn (c_1 - c_2)} \Big( q^{\textbf{1}_{a_1 > a_2}} Z_{x, y} \big( (a_1, a_2), \mathfrak{B}; (c_2, c_1), \mathfrak{D} \big) + q^{\textbf{1}_{a_2 > a_1}} Z_{x, y} \big( (a_2, a_1), \mathfrak{B}; (c_2, c_1), \mathfrak{D} \big) \Big), 
	\end{aligned} 
	\end{flalign}
	
	\noindent which would impy \eqref{qcd2} upon multiplying both sides by $q^{\inv (\overleftarrow{\mathfrak{B}})}$ and summing over $\mathfrak{B} \in \mathcal{M} (\textbf{B})$. 
	
	Since \eqref{qcd3} holds if $c_1 = c_2$, we will assume in what follows that $c_1 \ne c_2$, in which case \eqref{qcd3} will follow from a suitable application of the Yang--Baxter equation. Indeed, applying \eqref{rrrijk} $L$ times, we deduce for any fixed ordered pair of indices $(i_1, i_2)$ constituting a permutation of $\mathfrak{A}$ that 
	\begin{flalign} 
	\label{rqzxyk1k2} 
	\begin{aligned} 
	\displaystyle\sum_{k_1, k_2} R_q (i_1, i_2; k_1, k_2) & Z_{x, y} \big( (k_1, k_2), \mathfrak{B}; (c_1, c_2), \mathfrak{D} \big) \\
	& = \displaystyle\sum_{k_1, k_2} Z_{x, y} \big( (i_2, i_1), \mathfrak{B}; (k_2, k_1), \mathfrak{D} \big) R_q (k_1, k_2; c_1, c_2),
	\end{aligned} 
	\end{flalign}
	
	\noindent where both sums range over all $k_1, k_2 \in [0, m + n - 1]$; we refer to the top of \Cref{rvertical} for a depiction. By arrow conservation, any summand on either the left or right side is nonzero only if $(k_1, k_2)$ is a permutation of $\mathfrak{A}$ or $\mathfrak{C}$, respectively. So, we may instead sum over $(k_1, k_2) \in \big\{ (a_1, a_2), (a_2, a_1) \big\}$ on the left side of \eqref{rqzxyk1k2} and $(k_1, k_2) \in \big\{ (c_1, c_2), (c_2, c_1) \big\}$ on the right side of \eqref{rqzxyk1k2}.

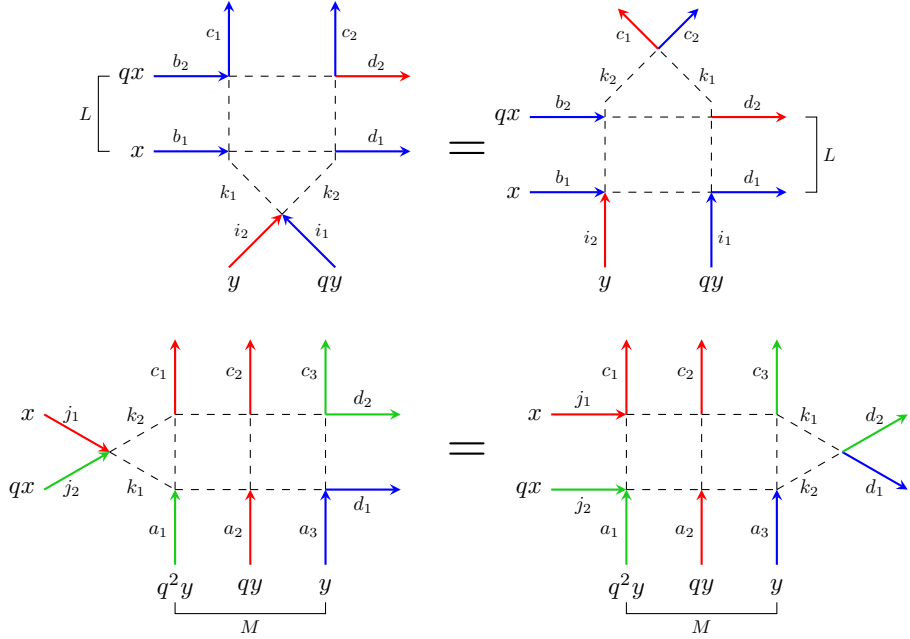
\begin{figure}[t]

	\begin{center}

		\begin{tikzpicture}[
		>=stealth,
		auto,
		style={
			scale = 1
		}
		]

		\draw[-, dashed] (-1.4142, 0) -- (0, 0);
		\draw[->, thick, blue] (0, 0) -- (1, 0);
		\draw[->, thick, blue] (-2.4142, 0) -- (-1.4142, 0);
		
		\draw[-, dashed] (-1.4142, 1) -- (0, 1); 
		\draw[->, thick, red] (0, 1) -- (1, 1);
		\draw[->, thick, blue] (-2.4142, 1) -- (-1.4142, 1); 
		
		\draw[->, thick, red] (-1.4142, -1.5435) -- (-.7071, -.8364);
		\draw[-, dashed] (-.7071, -.8364) -- (0, -.1193);
		
		\draw[->, thick, blue] (0, -1.5435) -- (-.7071, -.8364);
		\draw[-, dashed] (-.7071, -.8364) -- (-1.4142, -.1193);
		
		\draw[-, dashed] (-1.4142, 0) -- (-1.4142, 1);
		\draw[-, dashed] (0, 0) -- (0, 1);
		
		\draw[-, dashed] (0, -.1193) -- (0, 0);
		\draw[-, dashed] (-1.4142, -.1193) -- (-1.4142, 0);
		
		\draw[->, thick, blue] (0, 1) -- (0, 2);
		\draw[->, thick, blue] (-1.4142, 1) -- (-1.4142, 2);

		\draw[] (1.75, 0) circle [radius = 0] node[scale = 2]{$=$};
		
		\draw[] (-2.4142, 0) circle[radius = 0] node[scale = .95, left]{$x$};
		\draw[] (-2.4142, 1) circle[radius = 0] node[scale = .95, left]{$qx$};
		
		\draw[] (-.0707, -1.5435) circle[radius = 0] node[scale = .95, below]{$q y$};
		\draw[] (-1.3435, -1.5435) circle[radius = 0] node[scale = .95, below]{$y$};
		
		\draw[] (-1.0606, -1.06) circle[radius = 0] node[scale = .7, left]{$i_2$};
		\draw[] (-.3535, -1.06) circle[radius = 0] node[scale = .7, right]{$i_1$};
		\draw[] (-1.2106, -.55) circle[radius = 0] node[scale = .7, left]{$k_1$};
		\draw[] (-.25	35, -.55) circle[radius = 0] node[scale = .7, right]{$k_2$};
		\draw[] (-1.4142, 1.55) circle[radius = 0] node[scale = .7, left]{$c_1$};
		\draw[] (0, 1.55) circle[radius = 0] node[scale = .7, right]{$c_2$};
		\draw[] (-2.05, 0) circle[radius = 0] node[scale = .7, above]{$b_1$};
		\draw[] (.55, 0) circle[radius = 0] node[scale = .7, above]{$d_1$};
		\draw[] (-2.05, 1) circle[radius = 0] node[scale = .7, above]{$b_2$};
		\draw[] (.55, 1) circle[radius = 0] node[scale = .7, above]{$d_2$};

		\draw[-, dashed] (3.6858, -.5435) -- (5, -.5435); 
		\draw[->, blue, thick] (5, -.5435) -- (6, -.5435);
		\draw[->, blue, thick] (2.5858, -.5435) -- (3.5858, -.5435); 
		
		\draw[-, dashed] (3.6858, .4565) -- (5, .4565); 
		\draw[->, red, thick] (5, .4565) -- (6, .4565);
		\draw[->, blue, thick] (2.5858, .4565) -- (3.5858, .4565); 
	
		\draw[-, dashed] (3.5858, .6565) -- (4.2929, 1.3636);
		\draw[->, blue, thick] (4.2929, 1.3636) -- (4.8293, 1.9);
		
		\draw[-, dashed] (5, .6565) -- (4.2929, 1.3636);
		\draw[->, red, thick] (4.2929, 1.3636) -- (3.7565, 1.9);
		
		\draw[->, red, thick] (3.5858, -1.5435) -- (3.5858, -.5435);
		\draw[->, blue, thick] (5, -1.5435) -- (5, -.5435);
		
		\draw[-, dashed] (3.5858, -.5435) -- (3.5858, .4565);
		\draw[-, dashed] (5, -.5435) -- (5, .4565);
		
		\draw[-, dashed] (3.5858, .4565) -- (3.5858, .6565);
		\draw[-, dashed] (5, .4565) -- (5, .6565);
		
		\draw[] (2.5858, -.5435) circle[radius = 0] node[scale = .95, left]{$x$};
		\draw[] (2.5858, .4565) circle[radius = 0] node[scale = .95, left]{$qx$};
		\draw[] (5, -1.5435) circle[radius = 0] node[scale = .95, below]{$qy$};
		\draw[] (3.5858, -1.5435) circle[radius = 0] node[scale = .95, below]{$y$};
		
		\draw[] (3.5858, -1.15) circle[radius = 0] node[scale = .7, left]{$i_2$};
		\draw[] (5, -1.15) circle[radius = 0] node[scale = .7, right]{$i_1$};
		\draw[] (3.8393, 1) circle[radius = 0] node[scale = .7, left]{$k_2$};
		\draw[] (4.7464, 1) circle[radius = 0] node[scale = .7, right]{$k_1$};
		\draw[] (4.0393, 1.5272) circle[radius = 0] node[scale = .7, left]{$c_1$};
		\draw[] (4.5464, 1.5272) circle[radius = 0] node[scale = .7, right]{$c_2$};
		\draw[] (3.0358, -.5435) circle[radius = 0] node[scale = .7, above]{$b_1$};
		\draw[] (5.55, -.5435) circle[radius = 0] node[scale = .7, above]{$d_1$};
		\draw[] (3.0358, .4565) circle[radius = 0] node[scale = .7, above]{$b_2$};
		\draw[] (5.55, .4565) circle[radius = 0] node[scale = .7, above]{$d_2$};
		
		\draw[] (-3, 0) -- (-3.15, 0) -- (-3.15, 1) -- (-3, 1);
		\draw[] (6.25, -.5435) -- (6.4, -.5435) -- (6.4, .4565) -- (6.25, .4565); 
		
		\draw[] (-3.15, .5) circle[radius = 0] node[scale = .7, left] {$L$};
		\draw[] (6.4, -.0435) circle[radius = 0] node[scale = .7, right]{$L$};

		\draw[->, thick, green] (-3.87, -4.5) -- (-3, -4);
		\draw[->, thick, red] (-3.87, -3.5) -- (-3, -4);
		\draw[->, thick, blue] (-.13, -5.5) -- (-.13, -4.5); 
		
		\draw[] (-3.87, -3.5) circle[radius = 0]  node[left, scale = .95]{$x$};
		\draw[] (-3.87, -4.5) circle[radius = 0]  node[left, scale = .95]{$qx$};
		\draw[] (-.13, -5.5) circle[radius = 0]  node[below = 2, scale = .95]{$y$};
		\draw[] (-1.13, -5.5) circle[radius = 0]  node[below = 2, scale = .95]{$q y$};
		\draw[] (-2.13, -5.5) circle[radius = 0]  node[below, scale = .95]{$q^2 y$};
		\draw[] (-1.13, -6.15) circle[radius = 0]  node[below, scale = .7]{$M$};
		
		\draw[] (-2.13, -5.05) circle[radius = 0]  node[left, scale = .75]{$a_1$};
		\draw[] (-1.13, -5.05) circle[radius = 0]  node[left, scale = .75]{$a_2$};
		\draw[] (-.13, -5.05) circle[radius = 0]  node[left, scale = .75]{$a_3$};	
		\draw[] (-2.13, -2.95) circle[radius = 0]  node[left, scale = .75]{$c_1$};
		\draw[] (-1.13, -2.95) circle[radius = 0]  node[left, scale = .75]{$c_2$};
		\draw[] (-.13, -2.95) circle[radius = 0]  node[left, scale = .75]{$c_3$};
		
		\draw[] (.37, -4.5) circle[radius = 0]  node[below, scale = .75]{$d_1$};
		\draw[] (.37, -3.5) circle[radius = 0]  node[above, scale = .75]{$d_2$};
		
		\draw[] (-3.5, -4.275) circle[radius = 0]  node[below, scale = .75]{$j_2$};
		\draw[] (-3.5, -3.7) circle[radius = 0]  node[above, scale = .75]{$j_1$};
		
		\draw[] (-2.65, -4.275) circle[radius = 0]  node[below, scale = .75]{$k_1$};
		\draw[] (-2.65, -3.725) circle[radius = 0]  node[above, scale = .75]{$k_2$};
		
		\draw[] (-2.13, -6) -- (-2.13, -6.15) -- (-.13, -6.15) -- (-.13, -6);
		
		\draw[-, black, dashed] (-3, -4) -- (-2.13, -4.5); 
		\draw[-, black, dashed] (-3, -4) -- (-2.13, -3.5); 
		\draw[-, black, dashed] (-.13, -4.5) -- (-.13, -3.5);
		\draw[-, black, dashed] (-1.13, -4.5) -- (-1.13, -3.5);
		\draw[-, black, dashed] (-2.13, -4.5) -- (-2.13, -3.5);
		\draw[-, black, dashed] (-2.13, -4.5) -- (-1.13, -4.5);
		\draw[-, black, dashed] (-2.13, -3.5) -- (-1.13, -3.5);
		\draw[-, black, dashed] (-1.13, -4.5) -- (-.13, -4.5);
		\draw[-, black, dashed] (-1.13, -3.5) -- (-.13, -3.5);
		
		\draw[->, green, thick] (-2.13, -5.5) -- (-2.13, -4.5);
		\draw[->, red, thick] (-2.13, -3.5) -- (-2.13, -2.5);
		\draw[->, red, thick] (-1.13, -5.5) -- (-1.13, -4.5);
		\draw[->, red, thick] (-1.13, -3.5) -- (-1.13, -2.5);
		\draw[->, green, thick] (-.13, -3.5) -- (.87, -3.5); 
		\draw[->, blue, thick] (-.13, -4.5) -- (.87, -4.5); 
		\draw[->, green, thick] (-.13, -3.5) -- (-.13, -2.5);
		
		\draw[] (2.87, -3.5) circle[radius = 0]  node[left, scale = .95]{$x$};
		\draw[] (2.87, -4.5) circle[radius = 0]  node[left, scale = .95]{$qx$};
		\draw[] (5.87, -5.5) circle[radius = 0]  node[below = 2, scale = .95]{$y$};
		\draw[] (4.87, -5.5) circle[radius = 0]  node[below = 2, scale = .95]{$q y$};
		\draw[] (3.87, -5.5) circle[radius = 0]  node[below, scale = .95]{$q^2 y$};
	
		\draw[] (3.87, -5.05) circle[radius = 0]  node[left, scale = .75]{$a_1$};
		\draw[] (4.87, -5.05) circle[radius = 0]  node[left, scale = .75]{$a_2$};
		\draw[] (5.87, -5.05) circle[radius = 0]  node[left, scale = .75]{$a_3$}; 
		\draw[] (3.87, -2.95) circle[radius = 0]  node[left, scale = .75]{$c_1$};
		\draw[] (4.87, -2.95) circle[radius = 0]  node[left, scale = .75]{$c_2$};
		\draw[] (5.87, -2.95) circle[radius = 0]  node[left, scale = .75]{$c_3$};
		
		\draw[] (3.3, -3.5) circle[radius = 0]  node[above, scale = .75]{$j_1$};
		\draw[] (3.3, -4.5) circle[radius = 0]  node[below, scale = .75]{$j_2$};
		
		\draw[] (6.305, -3.725) circle[radius = 0]  node[above, scale = .75]{$k_1$};
		\draw[] (6.305, -4.275) circle[radius = 0]  node[below, scale = .75]{$k_2$};
		
		\draw[] (7.177, -3.725) circle[radius = 0]  node[above, scale = .75]{$d_2$};
		\draw[] (7.177, -4.275) circle[radius = 0]  node[below, scale = .75]{$d_1$};
		
		\draw[] (4.87, -6.15) circle[radius = 0]  node[below, scale = .7]{$M$};
		
		\draw[] (3.87, -6) -- (3.87, -6.15) -- (5.87, -6.15) -- (5.87, -6);
		
		\draw[->, red, thick] (4.87, -5.5) -- (4.87, -4.5); 
		\draw[->, blue, thick] (5.87, -5.5) -- (5.87, -4.5); 
		\draw[->, red, thick] (4.87, -3.5) -- (4.87, -2.5); 
		\draw[->, green, thick] (5.87, -3.5) -- (5.87, -2.5); 
		\draw[->, green, thick] (3.87, -5.5) -- (3.87, -4.5); 
		\draw[->, red, thick] (2.87, -3.5) -- (3.87, -3.5); 
		\draw[->, green, thick] (2.87, -4.5) -- (3.87, -4.5); 
		\draw[->, red, thick] (3.87, -3.5) -- (3.87, -2.5); 
		\draw[->, blue, thick] (6.74, -4) -- (7.61, -4.5); 
		\draw[->, green, thick] (6.74, -4) -- (7.61, -3.5); 
		
		\draw[-, black, dashed] (3.87, -4.5) -- (3.87, -3.5);
		\draw[-, black, dashed] (4.87, -4.5) -- (4.87, -3.5);
		\draw[-, black, dashed] (5.87, -4.5) -- (5.87, -3.5);
		\draw[-, black, dashed] (5.87, -4.5) -- (6.74, -4); 
		\draw[-, black, dashed] (5.87, -3.5) -- (6.74, -4);
		\draw[-, black, dashed] (3.87, -4.5) -- (4.87, -4.5);
		\draw[-, black, dashed] (3.87, -3.5) -- (4.87, -3.5);
		\draw[-, black, dashed] (4.87, -4.5) -- (5.87, -4.5);
		\draw[-, black, dashed] (4.87, -3.5) -- (5.87, -3.5);
		
		\filldraw[fill=white, draw=black] (1.75, -4) circle [radius=0] node[scale = 2]{$=$};
		
		\end{tikzpicture}
		
	\end{center}

	\caption{\label{rvertical} Shown on the top is a diagrammatic interpretation for \eqref{rqzxyk1k2}; shown on the bottom is its rotated analog.} 	
\end{figure}

Now, if $a_1 \ne a_2$, then by the explicit form (given by \eqref{rzdefinition}, \eqref{rzij}, and \eqref{rzi}) for $R_{ab} (z)$, we have $R_q (i, j; i', j') = q^{\textbf{1}_{i' > j'}} (q + 1)^{-1}$ if $i' \ne j'$ and $\{ i, j \} = \{ i', j' \}$. So, \eqref{rqzxyk1k2} implies upon multiplying both sides by $q + 1$ that  
\begin{flalign}
\label{ibib1} 
\begin{aligned}
q^{\textbf{1}_{a_1 > a_2}} Z_{x, y} & \big( (a_1, a_2), \mathfrak{B}; (c_1, c_2), \mathfrak{D} \big) + q^{\textbf{1}_{a_2 > a_1}} Z_{x, y} \big( (a_2, a_1), \mathfrak{B}; (c_1, c_2), \mathfrak{D} \big) \\
& = q^{\textbf{1}_{c_1 > c_2}} \Big( Z_{x, y} \big( (i_2, i_1), \mathfrak{B}; (c_1, c_2), \mathfrak{D} \big) + Z_{x, y} \big( (i_2, i_1), \mathfrak{B}; (c_2, c_1), \mathfrak{D} \big) \Big), 
\end{aligned} 
\end{flalign}

\noindent and, by similar reasoning, we find that 
\begin{flalign}
\label{ib1ib}
\begin{aligned}
q^{\textbf{1}_{a_1 > a_2}} Z_{x, y} & \big( (a_1, a_2), \mathfrak{B}; (c_2, c_1), \mathfrak{D} \big) + q^{\textbf{1}_{a_2 > a_1}} Z_{x, y} \big( (a_2, a_1), \mathfrak{B}; (c_2, c_1), \mathfrak{D} \big) \\
& = q^{\textbf{1}_{c_2 > c_1}} \Big( Z_{x, y} \big( (i_2, i_1), \mathfrak{B}; (c_1, c_2), \mathfrak{D} \big) + Z_{x, y} \big( (i_2, i_1), \mathfrak{B}; (c_2, c_1), \mathfrak{D} \big) \Big).
\end{aligned} 
\end{flalign}

\noindent Then, \eqref{qcd3} follows from comparing \eqref{ibib1} and \eqref{ib1ib}, since $\sgn (c_1 - c_2) = \textbf{1}_{c_1 > c_2} - \textbf{1}_{c_2 > c_1}$. 

If instead $a_1 = a_2$, then letting $a_1 = a = a_2$, we must have that $a < m$ (since $A_h \le 1$ holds for each $h \in [m, m + n - 1]$). Then, spin conservation implies that the left side of \eqref{rqzxyk1k2} is supported on the term $(k_1, k_2) = (a, a)$ and so, using the facts that $R_q (a, a; a, a) = 1$ for $a \in [0, m - 1]$ and $R_q (k_1, k_2; c_1, c_2) = (q + 1)^{-1} q^{\textbf{1}_{c_1 > c_2}}$, \eqref{rqzxyk1k2} yields 
\begin{flalign}
\label{ibib3}
\begin{aligned}
Z_{x, y} \big( (a, a), & \mathfrak{B}; (c_1, c_2), \mathfrak{D} \big) \\
& = (q + 1)^{-1} q^{\textbf{1}_{c_1 > c_2}} \Big( Z_{x, y} \big( (a, a), \mathfrak{B}; (c_1, c_2), \mathfrak{D} \big) + Z_{x, y} \big( (a, a), \mathfrak{B}; (c_2, c_1), \mathfrak{D} \big) \Big).
\end{aligned} 
\end{flalign}

\noindent By similar reasoning, 
\begin{flalign}
\label{ibib4}
\begin{aligned}
Z_{x, y} \big( (a, a), & \mathfrak{B}; (c_2, c_1), \mathfrak{D} \big) \\
& = (q + 1)^{-1} q^{\textbf{1}_{c_2 > c_1}} \Big( Z_{x, y} \big( (a, a), \mathfrak{B}; (c_1, c_2), \mathfrak{D} \big) + Z_{x, y} \big( (a, a), \mathfrak{B}; (c_2, c_1), \mathfrak{D} \big) \Big).
\end{aligned} 
\end{flalign}

\noindent Then, \eqref{ibib3} and \eqref{ibib4} together imply \eqref{qcd3} again since $\sgn (c_1 - c_2) = \textbf{1}_{c_1 > c_2} - \textbf{1}_{c_2 > c_1}$. 

This addresses the case $(\mathfrak{C}', \mathfrak{D}') = \big( \mathfrak{s}_k (\mathfrak{C}), \mathfrak{D} \big)$; as mentioned above, we omit the proof in the alternative case $(\mathfrak{C}', \mathfrak{D}') = \big( \mathfrak{C}, \mathfrak{s}_k (\mathfrak{D}) \big)$. 
\end{proof}

We next have the following lemma indicating a sort of exclusion principle for arrows of color $i \in [m, m + n - 1]$, in that if the west and south boundaries of any vertex model corresponding to $\mathcal{Z}_{x, y}$ each admit at most one arrow of any color $i \in [m, m + n - 1]$, then the same holds for its east and north boundaries. 

\begin{lem}
	
	\label{zb1d2} 
	
	Assume $n \ge 1$; fix $x, y \in \mathbb{C}$; and let $\textbf{\emph{A}}, \textbf{\emph{B}}, \textbf{\emph{C}}, \textbf{\emph{D}} \in \mathbb{Z}_{\ge 0}^{m + n}$ denote nonnegative integer sequences with coordinates indexed by $[0, m + n - 1]$, such that $|\textbf{\emph{A}}| = M = |\textbf{\emph{C}}|$ and $|\textbf{\emph{B}}| = L = |\textbf{\emph{D}}|$. Assume that $\max \{ A_i, B_i \} \le 1$ for each $i \in [m, m + n - 1]$ but that $\max \{C_i, D_i \} \ge 2$ for some $i \in [m, m + n - 1]$. Then, for any sequences of indices $\mathfrak{C} = (c_1, c_2, \ldots, c_M) \in \mathcal{M}(\textbf{\emph{C}})$ and $\mathfrak{D} = (d_1, d_2, \ldots,  d_L) \in \mathcal{M} (\textbf{\emph{D}})$, we have that $\mathcal{Z}_{x, y} (\textbf{\emph{A}}, \textbf{\emph{B}}; \mathfrak{C}, \mathfrak{D}) = 0$.	
\end{lem}

\begin{proof} 
	
	It suffices to show that $\mathcal{Z}_{x, y} (\textbf{A}, \textbf{B}; \mathfrak{C}, \mathfrak{D}) = 0$ if there exists some $h \in [m, m + n - 1]$ such that either $C_h \ge 2$ or $D_h \ge 2$. Since the two cases are entirely analogous, let us assume in what follows that $C_h \ge 2$. By \Cref{zabcdcd}, we may further assume that $c_1 = h = c_2$ and, by following the first part of the proof of \Cref{zabcdcd} (in particular, by using the second relation in \eqref{zabcdzabcd}), we may assume that $M = 2$. To that end, letting $\textbf{A} = \textbf{e}_{a_1} + \textbf{e}_{a_2} \in \mathbb{Z}_{\ge 0}^{m + n}$ for some indices $a_1, a_2 \in [0, m + n - 1]$, it suffices to show that 
	\begin{flalign}
	\label{qa1a2hh}
	q^{\textbf{1}_{a_1 > a_2}} Z_{x, y} & \big( (a_1, a_2), \mathfrak{B}; (h, h), \mathfrak{D} \big) + q^{\textbf{1}_{a_2 > a_1}} Z_{x, y} \big( (a_2, a_1), \mathfrak{B}; (h, h), \mathfrak{D} \big) = 0. 
	\end{flalign} 
	
	This will follow from an application of the Yang--Baxter equation similar to the one used to show \Cref{zabcdcd}. To that end, we may first assume that $a_1 \ne a_2$ since otherwise we would either have $a_1 = h = a_2$, meaning $A_h = 2$ or $B_h \ge 2$ (by arrow conservation); both contradict the fact that $\max \{ A_h, B_h \} \le 1$, and so $a_1 \ne a_2$. Next, as in \eqref{rqzxyk1k2}, $L$ applications of \eqref{rrrijk} (again, see the top diagram in \Cref{rvertical}) implies for any fixed permutation $(i_1, i_2)$ of $(a_1, a_2)$ that
	\begin{flalign} 
	\label{rqzxyii}
	\displaystyle\sum_{k_1, k_2} R_q (i_1, i_2; k_1, k_2) Z_{x, y} \big( (k_1, k_2), \mathfrak{B}; (h, h), \mathfrak{D} \big) = Z_{x, y} \big( (i_2, i_1), \mathfrak{B}; (h, h), \mathfrak{D} \big) R_q (h, h; h, h),
	\end{flalign}
	
	\noindent where the sum ranges over all pairs of indices $(k_1, k_2) \in \big\{ (a_1, a_2), (a_2, a_1) \big\}$. Here, on the right side, we have used the fact that $R_q (k_1, k_2; h, h)$ is nonzero only if $k_1 = h = k_2$, by arrow conservation. 
	
	By the explicit form (given by \eqref{rzdefinition}, \eqref{rzij}, and \eqref{rzi}) for $R_{ab} (z)$, we have $R_q (i, j; i', j') = q^{\textbf{1}_{i' > j'}} (q + 1)^{-1}$ if $i' \ne j'$ and $\{ i, j \} = \{ i', j' \}$, and also that $R_q (h, h; h, h) = 0$ for $h \in [m, m + n - 1]$. So, \eqref{qa1a2hh} follows from multiplying both sides of \eqref{rqzxyii} by $q + 1$, thereby implying the lemma. 
\end{proof}

\section{Fused Weights and the Yang--Baxter Equation}

\label{WeightsZR}

In this section we define the fused weights for the $U_q \big( \widehat{\mathfrak{sl}} (m | n) \big)$ vertex model and show they satisfy the Yang--Baxter equation. To that end, we begin with the following definition for these weights.

\begin{definition}
	
	\label{rijkh}
	
	Fix $L, M \in \mathbb{Z}_{\ge 1}$ and, for each index $X \in \{ A, B, C, D \}$, fix an $(m + n)$-tuple $\textbf{X} = (X_0, X_1, \ldots , X_{m + n - 1}) \in \mathbb{Z}_{\ge 0}^{m + n}$ such that $|\textbf{A}| = M = |\textbf{C}|$ and $|\textbf{B}| = L = |\textbf{D}|$. Letting $\mathfrak{C} = (c_1, c_2, \ldots , c_M) \in \mathcal{M} (\textbf{C})$ and $\mathfrak{D} = (d_1, d_2, \ldots , d_L) \in \mathcal{M} (\textbf{D})$ denote the unique sequences such that $c_1 \le c_2 \le \cdots \le c_M$ and $d_1 \ge d_2 \ge \cdots \ge d_L$, define the \emph{fused weight} 
	\begin{flalign*} 
	\mathcal{R}_{x, y} (\textbf{A}, \textbf{B}; \textbf{C}, \textbf{D}) = \mathcal{R}_{x, y}^{(m; n)} (\textbf{A}, \textbf{B}; \textbf{C}, \textbf{D}) = \mathcal{Z}_{x, y} (\textbf{A}, \textbf{B}; \mathfrak{C}, \mathfrak{D}) \cdot \displaystyle\prod_{i = 0}^{m + n - 1} \displaystyle\frac{(q; q)_{A_i} (q; q)_{B_i}}{(q; q)_{C_i} (q; q)_{D_i}},
	\end{flalign*} 

	\noindent where we recall that the right side is given by \Cref{zabcd2}.\index{R@$\mathcal{R}_{x, y}^{(m; n)} (\textbf{A}, \textbf{B}; \textbf{C}, \textbf{D})$}\index{R@$\mathcal{R}_{x, y}^{(m; n)} (\textbf{A}, \textbf{B}; \textbf{C}, \textbf{D})$!$\mathcal{R}_{x, y} (\textbf{A}, \textbf{B}; \textbf{C}, \textbf{D})$}
	
\end{definition}

\begin{rem}

\label{rabcd1} 

Fix nonnegative integer sequences $\textbf{A}, \textbf{B}, \textbf{C}, \textbf{D} \in \mathbb{Z}_{\ge 0}^{m + n - 1}$ such that $|\textbf{A}| = M = |\textbf{C}|$; $|\textbf{B}| = L = |\textbf{D}|$; and $\max \{ A_h, B_h \} \le 1$ for each $h \in [m, m + n - 1]$. Further fix index sequences $\mathfrak{C} \in \mathcal{M} (\textbf{C})$ and $\mathfrak{D} \in \mathcal{M} (\textbf{D})$. By \Cref{zabcdcd}, we have 
\begin{flalign}
	\label{rabcdrabcd} 
	\mathcal{Z}_{x, y} (\textbf{A}, \textbf{B}; \mathfrak{C}, \mathfrak{D}) = q^{\inv (\mathfrak{C}) + \inv(\overleftarrow{\mathfrak{D}})} \mathcal{R}_{x, y} (\textbf{A}, \textbf{B}; \textbf{C}, \textbf{D}) \cdot \displaystyle\prod_{i = 0}^{m + n - 1} \displaystyle\frac{(q; q)_{C_i} (q; q)_{D_i}}{(q; q)_{A_i} (q; q)_{B_i}}.
\end{flalign}

\noindent Thus, summing over all $\mathfrak{C} \in \mathcal{M} (\textbf{C})$ and $\mathfrak{D} \in \mathcal{M} (\textbf{D})$, and further using the fact that for any $\textbf{X} = (X_0, X_1, \ldots , X_{m + n - 1}) \in \mathbb{Z}_{\ge 0}^{m + n - 1}$ we have
\begin{flalign*}
	\displaystyle\sum_{\mathscr{X} \in \mathcal{M} (\textbf{X})} q^{\inv (\mathscr{X})} = \displaystyle\frac{(q; q)_{|\textbf{X}|}}{\prod_{i = 0}^{m + n - 1} (q; q)_{X_i}},
\end{flalign*}

\noindent we deduce 
\begin{flalign*}
\mathcal{R}_{x, y} (\textbf{A}, \textbf{B}; \textbf{C}, \textbf{D}) = \displaystyle\frac{\prod_{i = 0}^{m + n - 1} (q; q)_{A_i}}{(q; q)_M} \cdot \displaystyle\frac{\prod_{i = 0}^{m + n - 1} (q; q)_{B_i}}{(q; q)_L} \cdot \displaystyle\sum_{\substack{\mathfrak{C} \in \mathcal{M} (\textbf{C}) \\ \mathfrak{D} \in \mathcal{M} (\textbf{D})}} \mathcal{Z}_{x, y} (\textbf{A}, \textbf{B}; \textbf{C}, \textbf{D}).
\end{flalign*} 

\noindent Thus, in the case $n = 0$, the fused weights described in \Cref{rijkh} coincide with those given by Definition 8.3 (see also equation (8.3)) of \cite{SVMP}.

\end{rem}

Similar to in \Cref{Weightszq} for the case $L = 1 = M$, we view the quantity $\mathcal{R}_{x, y} (\textbf{A}, \textbf{B}; \textbf{C}, \textbf{D})$ as the weight of a vertex $v$ whose row and column rapidities are given by $x$ and $y$, respectively, with arrow configuration $(\textbf{A}, \textbf{B}; \textbf{C}, \textbf{D})$. The latter point now means that, for each $i \in [0, m + n - 1]$, $A_i$ and $B_i$ arrows of color $i$ vertically and horizontally enter $v$, respectively; similarly, $C_i$ and $D_i$ arrows of color $i$ vertically and horizontally exit $v$, respectively. Unlike in \Cref{Weightszq}, the case $L, M > 1$ allows multiple arrows to exist along vertical and horizontal edges adjacent to $v$; we refer to the right side of \Cref{zxy} for a depiction. 

The following corollary, which follows directly from \Cref{zb1d2}, shows that the weights $\mathcal{R}_{x, y}$ impose an exclusion restriction on colors $i \in [m, m + n - 1]$ along any edge.

\begin{cor}
	
	\label{rb1d2} 
	
	Assume $n \ge 1$; fix $x, y \in \mathbb{C}$; and let $\textbf{\emph{A}}, \textbf{\emph{B}}, \textbf{\emph{C}}, \textbf{\emph{D}} \in \mathbb{Z}_{\ge 0}^{m + n}$ denote nonnegative integer sequences with coordinates indexed by $[0, m + n - 1]$. Assume that $\max \{ A_i, B_i \} \le 1$ for each $i \in [m, m + n - 1]$. Then, $\mathcal{R}_{x, y} (\textbf{\emph{A}}, \textbf{\emph{B}}; \textbf{\emph{C}}, \textbf{\emph{D}}) = 0$ unless $\max \{ C_i, D_i \} \le 1$ for each $i \in [m, m + n - 1]$.	
\end{cor}

\begin{proof} 
	
	This follows from \Cref{zb1d2} and \Cref{rijkh}.
\end{proof} 

In what follows we will typically consider vertex models whose boundary data admits at most one arrow of any color $i \in [m, m + n - 1]$ along any edge entering the domain. Then, \Cref{rb1d2} shows that this exclusion property is retained along interior or exiting boundary edges of the domain. 

The following proposition now states that the fused $\mathcal{R}_{x, y}$ weights satisfy the Yang--Baxter equation.

\begin{prop}
	
	\label{wabcdproduct} 
	
	Fix $L, M, N \in \mathbb{Z}_{\ge 1}$ and $x, y, z \in \mathbb{C}$. Fix $\textbf{\emph{I}}_1, \textbf{\emph{J}}_1, \textbf{\emph{K}}_1, \textbf{\emph{I}}_3, \textbf{\emph{J}}_3, \textbf{\emph{K}}_3 \in \mathbb{Z}_{\ge 0}^{m + n}$ such that $|\textbf{\emph{I}}_1| = L = |\textbf{\emph{I}}_3|$, $|\textbf{\emph{J}}_1| = M = |\textbf{\emph{J}}_3|$, and $|\textbf{\emph{K}}_1| = N = |\textbf{\emph{K}}_3|$. If $X_h \le 1$ for any indices $h \in [m, m + n - 1]$ and $\textbf{\emph{X}} = (X_0, X_1, \ldots , X_{m + n - 1}) \in \{ \textbf{\emph{I}}_1, \textbf{\emph{J}}_1, \textbf{\emph{K}}_1 \}$, then
	\begin{flalign}
	\label{requation}
	\begin{aligned} 
	\displaystyle\sum_{\textbf{\emph{I}}_2, \textbf{\emph{J}}_2, \textbf{\emph{K}}_2} \mathcal{R}_{x, y} & ( \textbf{\emph{I}}_1, \textbf{\emph{J}}_1; \textbf{\emph{I}}_2, \textbf{\emph{J}}_2) \mathcal{R}_{x, z} ( \textbf{\emph{K}}_1, \textbf{\emph{J}}_2; \textbf{\emph{K}}_2, \textbf{\emph{J}}_3) \mathcal{R}_{y, z} ( \textbf{\emph{K}}_2, \textbf{\emph{I}}_2; \textbf{\emph{K}}_3, \textbf{\emph{I}}_3) \\
	& = \displaystyle\sum_{\textbf{\emph{I}}_2, \textbf{\emph{J}}_2, \textbf{\emph{K}}_2} \mathcal{R}_{y, z} ( \textbf{\emph{K}}_1, \textbf{\emph{I}}_1; \textbf{\emph{K}}_2, \textbf{\emph{I}}_2) \mathcal{R}_{x, z} ( \textbf{\emph{K}}_2, \textbf{\emph{J}}_1; \textbf{\emph{K}}_3, \textbf{\emph{J}}_2) \mathcal{R}_{x, y} ( \textbf{\emph{I}}_2, \textbf{\emph{J}}_2; \textbf{\emph{I}}_3, \textbf{\emph{J}}_3), 
	\end{aligned} 
	\end{flalign}
	
	\noindent where both sums are over all $\textbf{\emph{I}}_2, \textbf{\emph{J}}_2, \textbf{\emph{K}}_2 \in \mathbb{Z}_{\ge 0}^{m + n}$ with $|\textbf{\emph{I}}_2| = L$, $|\textbf{\emph{J}}_2| = M$, and $|\textbf{\emph{K}}_2| = N$. Diagrammatically, 
	
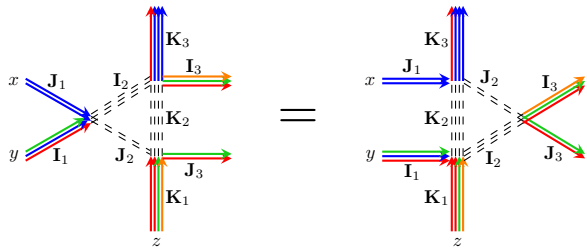
\begin{center}			
		\begin{tikzpicture}[
			>=stealth,
			auto,
			style={
				scale = 1
			}
			]
			\draw[->, thick, blue] (-.87, .47) -- (0, -.03);
			\draw[->, thick, blue] (-.87, .53) -- (0, .03);		
			\draw[->, thick, red] (-.87, -.56) -- (0, -.06);
			\draw[->, thick, blue] (-.87, -.5) -- (-.0435, -.025);
			\draw[->, thick, green] (-.87, -.44) -- (-.087, .01);			
			\draw[->, thick, red] (.845, -1.5) -- (.845, -.5);
			\draw[->, thick, red] (.795, -1.5) -- (.795, -.5);
			\draw[->, thick, green] (.895, -1.5) -- (.895, -.5);
			\draw[->, thick, orange] (.945, -1.5) -- (.945, -.5);		
			\draw[-, dashed] (.895, -.5) -- (.895, .5);
			\draw[-, dashed] (.845, -.5) -- (.845, .5);
			\draw[-, dashed] (.945, -.5) -- (.945, .5);
			\draw[-, dashed] (.795, -.5) -- (.795, .5);		
			\draw[->, thick, blue] (.895, .5) -- (.895, 1.5);
			\draw[->, thick, blue] (.845, .5) -- (.845, 1.5);
			\draw[->, thick, blue] (.945, .5) -- (.945, 1.5);
			\draw[->, thick, red] (.795, .5) -- (.795, 1.5);
			\draw[] (-.87, .5) circle[radius = 0]  node[left, scale = .7]{$x$};
			\draw[] (-.87, -.5) circle[radius = 0]  node[left, scale = .7]{$y$};
			\draw[] (.87, -1.5) circle[radius = 0]  node[below, scale = .7]{$z$};
			\draw[-, dashed] (0, .03) -- (.783, -.425); 
			\draw[-, dashed] (0, -.03) -- (.783, -.475); 
			\draw[-, dashed] (0, -.06) -- (.783, .44); 
			\draw[-, dashed] (0, 0) -- (.783, .5); 
			\draw[-, dashed] (0, .06) -- (.783, .56); 	
			\draw[->, thick, red] (.95, -.53) -- (1.87, -.53); 
			\draw[->, thick, green] (.95, -.47) -- (1.87, -.47); 			
			\draw[->, thick, red] (.95, .44) -- (1.87, .44);
			\draw[->, thick, green] (.95, .5) -- (1.87, .5); 
			\draw[->, thick, orange] (.95, .56) -- (1.87, .56); 			
			\draw[] (3.87, .5) circle[radius = 0]  node[left, scale = .7]{$x$};
			\draw[] (3.87, -.5) circle[radius = 0]  node[left, scale = .7]{$y$};
			\draw[] (4.87, -1.5) circle[radius = 0]  node[below, scale = .7]{$z$}; 			
			\draw[->, thick, red] (4.795, -1.5) -- (4.795, -.5); 
			\draw[->, thick, red] (4.845, -1.5) -- (4.845, -.5); 
			\draw[->, thick, green] (4.895, -1.5) -- (4.895, -.5); 
			\draw[->, thick, orange] (4.945, -1.5) -- (4.945, -.5); 		
			\draw[->, thick, blue] (3.87, .47) -- (4.77, .47); 
			\draw[->, thick, blue] (3.87, .53) -- (4.77, .53);			
			\draw[->, thick, red] (3.87, -.56) -- (4.77, -.56);
			\draw[->, thick, blue] (3.87, -.5) -- (4.77, -.5);
			\draw[->, thick, green] (3.87, -.44) -- (4.77, -.44);		
			\draw[->, thick, red] (4.795, .5) -- (4.795, 1.5); 
			\draw[->, thick, blue] (4.845, .5) -- (4.845, 1.5); 
			\draw[->, thick, blue] (4.895, .5) -- (4.895, 1.5); 
			\draw[->, thick, blue] (4.945, .5) -- (4.945, 1.5); 		
			\draw[->, thick, red] (5.74, -.03) -- (6.567, -.53);
			\draw[->, thick, green] (5.74, .03) -- (6.567, -.47);		
			\draw[->, thick, red] (5.74, -.06) -- (6.567, .44); 
			\draw[->, thick, green] (5.74, 0) -- (6.567, .5); 
			\draw[->, thick, orange] (5.74, .06) -- (6.567, .56); 		
			\draw[-, dashed] (4.895, -.5) -- (4.895, .5);
			\draw[-, dashed] (4.795, -.5) -- (4.795, .5); 
			\draw[-, dashed] (4.845, -.5) -- (4.845, .5);
			\draw[-, dashed] (4.945, -.5) -- (4.945, .5);	
			\draw[-, dashed] (4.957, -.44) -- (5.74, .06); 
			\draw[-, dashed] (4.957, -.5) -- (5.74, 0); 
			\draw[-, dashed] (4.957, -.56) -- (5.74, -.06); 	
			\draw[-, dashed] (4.957, .47) -- (5.74, -.03);
			\draw[-, dashed] (4.957, .53) -- (5.74, .03);
			\filldraw[fill=white, draw=black] (2.75, 0) circle [radius=0] node[scale = 2]{$=$};
			\draw[] (.9, -1.05) circle[radius = 0]  node[right, scale = .7]{$\textbf{\emph{K}}_1$};
			\draw[] (.9, 1.05) circle[radius = 0]  node[right, scale = .7]{$\textbf{\emph{K}}_3$};
			\draw[] (.9, 0) circle[radius = 0]  node[right, scale = .7]{$\textbf{\emph{K}}_2$};
			\draw[] (-.4, -.31) circle[radius = 0]  node[below, scale = .7]{$\textbf{\emph{I}}_1$};
			\draw[] (-.45, .3) circle[radius = 0]  node[above, scale = .7]{$\textbf{\emph{J}}_1$};
			\draw[] (.4, .31) circle[radius = 0]  node[above, scale = .7]{$\textbf{\emph{I}}_2$};
			\draw[] (.45, -.3) circle[radius = 0]  node[below, scale = .7]{$\textbf{\emph{J}}_2$};
			\draw[] (1.37, .5) circle[radius = 0]  node[above, scale = .7]{$\textbf{\emph{I}}_3$};
			\draw[] (1.37, -.5) circle[radius = 0]  node[below, scale = .7]{$\textbf{\emph{J}}_3$};	
			\draw[] (4.27, -.52) circle[radius = 0]  node[below, scale = .7]{$\textbf{\emph{I}}_1$};
			\draw[] (4.27, .5) circle[radius = 0]  node[above, scale = .7]{$\textbf{\emph{J}}_1$};
			\draw[] (4.85, -1.05) circle[radius = 0]  node[left, scale = .7]{$\textbf{\emph{K}}_1$};
			\draw[] (4.85, 0) circle[radius = 0]  node[left, scale = .7]{$\textbf{\emph{K}}_2$};
			\draw[] (4.85, 1.05) circle[radius = 0]  node[left, scale = .7]{$\textbf{\emph{K}}_3$};
			\draw[] (5.3, .275) circle[radius = 0]  node[above = 2, scale = .7]{$\textbf{\emph{J}}_2$};
			\draw[] (5.35, -.275) circle[radius = 0]  node[below = 2, scale = .7]{$\textbf{\emph{I}}_2$};
			\draw[] (6.15, -.275) circle[radius = 0]  node[below, scale = .7]{$\textbf{\emph{J}}_3$};
			\draw[] (6.1, .275) circle[radius = 0]  node[above, scale = .7]{$\textbf{\emph{I}}_3$};
		\end{tikzpicture}
	\end{center}
	
	\noindent where states along the solid edges are fixed and those along dashed edges are summed over.

\end{prop}

\begin{proof}
	
	It will be useful here to introduce a gauge transformation $\widetilde{\mathcal{R}}_{x, y} (\textbf{A}, \textbf{B}; \textbf{C}, \textbf{D})$ of the $\mathcal{R}_{x, y} (\textbf{A}, \textbf{B}; \textbf{C}, \textbf{D})$ fused weights, given as follows. Recalling the notation of \Cref{rijkh}, we set 
	\begin{flalign}
		\label{r1xyabcd}
		\widetilde{\mathcal{R}}_{x, y} (\textbf{A}, \textbf{B}; \textbf{C}, \textbf{D}) = \mathcal{R}_{x, y} (\textbf{A}, \textbf{B}; \textbf{C}, \textbf{D}) \cdot \displaystyle\prod_{i = 0}^{m + n - 1} \displaystyle\frac{(q; q)_{C_i} (q; q)_{D_i}}{(q; q)_{A_i} (q; q)_{B_i}}.
	\end{flalign} 
	
	\noindent Then, to establish \eqref{requation}, it suffices to show 
	\begin{flalign}
		\label{requation1}
		\begin{aligned} 
			\displaystyle\sum_{\textbf{I}_2, \textbf{J}_2, \textbf{K}_2} \widetilde{\mathcal{R}}_{x, y} & ( \textbf{I}_1, \textbf{J}_1; \textbf{I}_2, \textbf{J}_2) \widetilde{\mathcal{R}}_{x, z} ( \textbf{K}_1, \textbf{J}_2; \textbf{K}_2, \textbf{J}_3) \widetilde{\mathcal{R}}_{y, z} ( \textbf{K}_2, \textbf{I}_2; \textbf{K}_3, \textbf{I}_3) \\
			& = \displaystyle\sum_{\textbf{I}_2, \textbf{J}_2, \textbf{K}_2} \widetilde{\mathcal{R}}_{y, z} ( \textbf{K}_1, \textbf{I}_1; \textbf{K}_2, \textbf{I}_2) \widetilde{\mathcal{R}}_{x, z} ( \textbf{K}_2, \textbf{J}_1; \textbf{K}_3, \textbf{J}_2) \widetilde{\mathcal{R}}_{x, y} ( \textbf{I}_2, \textbf{J}_2; \textbf{I}_3, \textbf{J}_3), 
		\end{aligned} 
	\end{flalign}

	\noindent where the sum is as in \eqref{requation}.
	
	To that end, by $LMN$ applications of \eqref{rrrijk}, we obtain for any fixed $\mathfrak{I}_h \in \mathcal{M}(\textbf{I}_h)$, $\mathfrak{J}_h \in \mathcal{M}(\textbf{J}_h)$, and $\mathfrak{K}_h \in \mathcal{M}(\textbf{K}_h)$, for $h \in \{ 1, 3 \}$, that
	\begin{flalign} 
	\label{lmzrrrrz}
	\begin{aligned} 
	& \displaystyle\sum_{\mathfrak{I}_2, \mathfrak{J}_2, \mathfrak{K}_2} Z_{x, y} \big( \overleftarrow{\mathfrak{I}_1}, \mathfrak{J}_1; \overleftarrow{\mathfrak{I}_2}, \mathfrak{J}_2 \big) Z_{x, z} \big( \mathfrak{K}_1, \mathfrak{J}_2; \mathfrak{K}_2, \mathfrak{J}_3 \big) Z_{y, z} (\mathfrak{K}_2, \mathfrak{I}_2; \mathfrak{K}_3, \mathfrak{I}_3) \\
	& \qquad = \displaystyle\sum_{\mathfrak{I}_2, \mathfrak{J}_2, \mathfrak{K}_2}  Z_{y, z} \big( \mathfrak{K}_1, \mathfrak{I}_1; \mathfrak{K}_2, \mathfrak{I}_2 \big) Z_{x, z} (\mathfrak{K}_2, \mathfrak{J}_1; \mathfrak{K}_3, \mathfrak{J}_2) Z_{x, y} \big( \overleftarrow{\mathfrak{I}_2}, \mathfrak{J}_2; \overleftarrow{\mathfrak{I}_3}, \mathfrak{J}_3 \big),
	\end{aligned} 
	\end{flalign} 
	
	\noindent where the sum is over all sequences $\mathfrak{I}_2, \mathfrak{J}_2, \mathfrak{K}_2$ of indices in $[0, m + n - 1]$; diagrammatically,
	\begin{center}

		\begin{tikzpicture}[
		>=stealth,
		auto,
		style={
			scale = 1
		}
		]
		
		\draw[->, thick, red] (-.5, -.5) -- (0, 0);
		\draw[->, thick, blue] (0, -1) -- (.5, -.5);
		\draw[->, thick, blue] (.5, -1.5) -- (1, -1);
		
		\draw[->, thick, blue] (-.5, .5) -- (0, 0);
		\draw[->, thick, green] (0, 1) -- (.5, .5);

		\draw[dashed] (.5, .5) -- (1, 0);
		\draw[dashed] (1, 0) -- (1.5, -.5);
		
		\draw[dashed] (0, 0) -- (.5, -.5);
		\draw[dashed] (.5, -.5) -- (1, -1);
		
		\draw[dashed] (0, 0) -- (.5, .5);
		\draw[dashed] (.5, -.5) -- (1, 0);
		\draw[dashed] (1, -1) -- (1.5, -.5);

		\draw[-, dashed] (.5, .5) -- (1, 1);
		\draw[-, dashed] (1, 0) -- (1.5, .5);
		\draw[-, dashed] (1.5, -.5) -- (2, 0);
		\draw[-, dashed] (1.5, -.5) -- (1.75, -.75);
		\draw[-, dashed] (1, -1) -- (1.25, -1.25);

		\draw[dashed] (2, 0) -- (2.5, 0);
		\draw[dashed] (1.5, .5) -- (2.5, .5);
		\draw[dashed] (1, 1) -- (2.5, 1);
		\draw[dashed] (1.75, -.75) -- (2.5, -.75);
		\draw[dashed] (1.25, -1.25) -- (2.5, -1.25);
		
		\draw[-, dashed] (2.5, -1.25) -- (3, -1.25);
		\draw[-, dashed] (2.5, -.75) -- (3, -.75);
		\draw[-, dashed] (2.5, 0) -- (3, 0);
		\draw[-, dashed] (2.5, .5) -- (3, .5);
		\draw[-, dashed] (2.5, 1) -- (3, 1);
		
		\draw[->, thick, blue] (3, -1.25) -- (3.5, -1.25);
		\draw[->, thick, red] (3, -.75) -- (3.5, -.75);
		\draw[->, thick, green] (3, 0) -- (3.5, 0);
		\draw[->, thick, blue] (3, .5) -- (3.5, .5);
		\draw[->, thick, blue] (3, 1) -- (3.5, 1);
		
		\draw[->, thick, red] (3, 1) -- (3, 1.6); 
		\draw[->, thick, red] (3, -2) -- (3, -1.25);
		
		\draw[->, thick, green] (2.5, 1) -- (2.5, 1.6);
		\draw[->, thick, green] (2.5, -2) -- (2.5, -1.25);
		
		\draw[] (4.3, -.125) circle [radius = 0] node[scale = 2]{$=$};
	
		\draw[dashed] (2.5, -1.25) -- (2.5, -.75);
		\draw[dashed] (2.5, -.75) -- (2.5, 0);
		\draw[dashed] (2.5, .5) -- (2.5, 1);
		\draw[dashed] (2.5, 0) -- (2.5, .5);
		
		\draw[dashed] (3, -1.25) -- (3, -.75);
		\draw[dashed] (3, -.75) -- (3, 0);
		\draw[dashed] (3, .5) -- (3, 1);
		\draw[dashed] (3, 0) -- (3, .5);		
		
		\draw[] (2.57, -1.65) circle[radius = 0] node[left, scale = .7]{$k_1^{(1)}$};
		\draw[] (2.57, 1.6) circle[radius = 0] node[above, scale = .7]{$k_1^{(3)}$};
		
		\draw[] (3.07, -1.65) circle[radius = 0] node[left, scale = .7]{$k_2^{(1)}$};
		\draw[] (3.07, 1.6) circle[radius = 0] node[above, scale = .7]{$k_2^{(3)}$};
		
		\draw[] (.9, -1.1) circle[radius = 0] node[left, scale = .7]{$i_1^{(1)}$};
		\draw[] (.4, -.6) circle[radius = 0] node[left, scale = .7]{$i_2^{(1)}$};
		\draw[] (-.1, -.1) circle[radius = 0] node[left, scale = .7]{$i_3^{(1)}$};
		
		\draw[] (-.1, .25) circle[radius = 0] node[above, scale = .7]{$j_1^{(1)}$};
		\draw[] (.4, .75) circle[radius = 0] node[above, scale = .7]{$j_2^{(1)}$};
		
		\draw[] (3.3, 0) circle[radius = 0] node[above, scale = .7]{$i_1^{(3)}$};
		\draw[] (3.3, .5) circle[radius = 0] node[above, scale = .7]{$i_2^{(3)}$};
		\draw[] (3.3, 1) circle[radius = 0] node[above, scale = .7]{$i_3^{(3)}$};
		
		\draw[] (3.3, -.75) circle[radius = 0] node[above, scale = .7]{$j_2^{(3)}$};
		\draw[] (3.3, -1.25) circle[radius = 0] node[above, scale = .7]{$j_1^{(3)}$}; 
		
		\draw[] (2.15, 0) circle[radius = 0] node[above, scale = .7]{$i_1^{(2)}$};
		\draw[] (1.9, .5) circle[radius = 0] node[above, scale = .7]{$i_2^{(2)}$};
		\draw[] (1.65, 1) circle[radius = 0] node[above, scale = .7]{$i_3^{(2)}$};
		
		\draw[] (2.15, -.75) circle[radius = 0] node[above, scale = .7]{$j_2^{(2)}$};
		\draw[] (1.85, -1.25) circle[radius = 0] node[above, scale = .7]{$j_1^{(2)}$};

		\draw[] (7.55, -1.25) circle[radius = 0] node[above, scale = .7]{$i_1^{(2)}$};
		\draw[] (7.25, -.75) circle[radius = 0] node[above, scale = .7]{$i_2^{(2)}$};
		\draw[] (6.95, -.25) circle[radius = 0] node[above, scale = .7]{$i_3^{(2)}$};
		
		\draw[] (7.35, 1) circle[radius = 0] node[above, scale = .7]{$j_2^{(2)}$};
		\draw[] (6.95, .5) circle[radius = 0] node[above, scale = .7]{$j_1^{(2)}$}; 
		
		\draw[dashed] (6.5, -1.25) -- (6.5, -.75);
		\draw[dashed] (6.5, -.75) -- (6.5, -.25); 
		\draw[dashed] (6.5, .5) -- (6.5, 1);
		\draw[dashed] (6.5, -.25) -- (6.5, .5);
		
		\draw[-, dashed] (6, -.25) -- (6.5, -.25);
		\draw[-, dashed] (6, .5) -- (6.5, .5);
		\draw[-, dashed] (6, 1) -- (6.5, 1);
		\draw[-, dashed] (6, -.75) -- (6.5, -.75);
		\draw[-, dashed] (6, -1.25) -- (6.5, -1.25);
		
		\draw[->, thick, red] (6.5, -2) -- (6.5, -1.25);
		\draw[->, thick, red] (6.5, 1) -- (6.5, 1.6);
		
		\draw[dashed] (6, -1.25) -- (6, -.75);
		\draw[dashed] (6, -.75) -- (6, -.25); 
		\draw[dashed] (6, .5) -- (6, 1);
		\draw[dashed] (6, -.25) -- (6, .5);
		
		\draw[->, thick, red] (5.5, -.25) -- (6, -.25);
		\draw[->, thick, blue] (5.5, .5) -- (6, .5);
		\draw[->, thick, green] (5.5, 1) -- (6, 1);
		\draw[->, thick, blue] (5.5, -.75) -- (6, -.75);
		\draw[->, thick, blue] (5.5, -1.25) -- (6, -1.25);
		
		\draw[->, thick, green] (6, -2) -- (6, -1.25);
		\draw[->, thick, green] (6, 1) -- (6, 1.6);
		
		\draw[dashed] (7.75, .25) -- (8.25, .75); 
		\draw[dashed] (8.25, -.25) -- (8.75, .25);
		\draw[dashed] (8.75, -.75) -- (9.25, -.25);
		
		\draw[dashed] (8, 1) -- (8.25, .75); 
		\draw[dashed] (7.5, .5) -- (7.75, .25);
		\draw[dashed] (7.25, -.25) -- (7.75, .25);
		\draw[dashed] (7.75, -.75) -- (8.25, -.25);
		\draw[dashed] (8.25, -1.25) -- (8.75, -.75);
		
		\draw[dashed] (8.25, .75) -- (8.75, .25);
		\draw[dashed] (8.75, .25) -- (9.25, -.25);
		\draw[dashed] (8.25, -.25) -- (8.75, -.75);
		\draw[dashed] (7.75, .25) -- (8.25, -.25);
		
		\draw[-, dashed] (6.5, -.25) -- (7.25, -.25);
		\draw[-, dashed] (6.5, .5) -- (7.5, .5);
		\draw[-, dashed] (6.5, 1) -- (8, 1);
		\draw[-, dashed] (6.5, -.75) --(7.75, -.75);
		\draw[-, dashed] (6.5, -1.25) -- (8.25, -1.25);
		
		\draw[->, thick, red] (9.25, -.25) -- (9.75, -.75);
		\draw[->, thick, blue] (8.75, -.75) -- (9.25, -1.25);
		\draw[->, thick, green] (9.25, -.25) -- (9.75, .25);
		\draw[->, thick, blue] (8.75, .25) -- (9.25, .75);
		\draw[->, thick, blue] (8.25, .75) -- (8.75, 1.25);

		\draw[] (-.5, .5) circle[radius = 0] node[left, scale = .7]{$x$};
		\draw[] (0, 1) circle[radius = 0] node[left, scale = .7]{$q x$};
		
		\draw[] (.5, -1.5) circle[radius = 0] node[left, scale = .7]{$y$};
		\draw[] (0, -1) circle[radius = 0] node[left, scale = .7]{$qy$};
		\draw[] (-.5, -.5) circle[radius = 0] node[left, scale = .7]{$q^2 y$};
		
		\draw[] (2.5, -1.97) circle[radius = 0] node[below, scale = .7]{$qz$};
		\draw[] (3, -1.97) circle[radius = 0] node[below, scale = .7]{$z$};
		
		\draw[] (5.4, .5) circle[radius = 0] node[left, scale = .7]{$x$};
		\draw[] (5.4, 1) circle[radius = 0] node[left, scale = .7]{$q x$};
		
		\draw[] (5.4, -1.25) circle[radius = 0] node[left, scale = .7]{$y$};
		\draw[] (5.4, -.75) circle[radius = 0] node[left, scale = .7]{$qy$};
		\draw[] (5.4, -.25) circle[radius = 0] node[left, scale = .7]{$q^2 y$};
		
		\draw[] (6, -1.97) circle[radius = 0] node[below, scale = .7]{$qz$};
		\draw[] (6.5, -1.97) circle[radius = 0] node[below, scale = .7]{$z$};
		
		\draw[] (5.97, -1.65) circle[radius = 0] node[right, scale = .7]{$k_1^{(1)}$};
		\draw[] (5.97, 1.6) circle[radius = 0] node[above, scale = .7]{$k_1^{(3)}$};
		
		\draw[] (6.47, -1.65) circle[radius = 0] node[right, scale = .7]{$k_2^{(1)}$};
		\draw[] (6.47, 1.6) circle[radius = 0] node[above, scale = .7]{$k_2^{(3)}$};
		
		\draw[] (8.65, 1.2) circle[radius = 0] node[left, scale = .7]{$i_3^{(3)}$};
		\draw[] (9.15, .7) circle[radius = 0] node[left, scale = .7]{$i_2^{(3)}$};
		\draw[] (9.65, .2) circle[radius = 0] node[left, scale = .7]{$i_1^{(3)}$};
		
		\draw[] (9.75, -.6) circle[radius = 0] node[above, scale = .7]{$j_2^{(3)}$};
		\draw[] (9.25, -1.1) circle[radius = 0] node[above, scale = .7]{$j_1^{(3)}$};
		
		\draw[] (5.7, -1.25) circle[radius = 0] node[above, scale = .7]{$i_1^{(1)}$};
		\draw[] (5.7, -.75) circle[radius = 0] node[above, scale = .7]{$i_2^{(1)}$};
		\draw[] (5.7, -.25) circle[radius = 0] node[above, scale = .7]{$i_3^{(1)}$};
		
		\draw[] (5.7, .5) circle[radius = 0] node[above, scale = .7]{$j_1^{(1)}$};
		\draw[] (5.7, 1) circle[radius = 0] node[above, scale = .7]{$j_2^{(1)}$};
		
		\draw[-] (-.875, -.875) -- (-1, -1)-- (0, -2) -- (.125, -1.875);
		\draw[-] (-.875, .875) -- (-1, 1) -- (-.5, 1.5) -- (-.375, 1.375); 
		
		\draw[-] (9, 1.5) -- (9.125, 1.625) -- (10.125, .625) -- (10, .5);
		\draw[] (10, -1) -- (10.125, -1.125) -- (9.625, -1.625) -- (9.5, -1.5);
		
		\draw[] (2.5, -2.3) -- (2.5, -2.45)  -- (3, -2.45) -- (3, -2.3);
		\draw[] (6, -2.3) -- (6, -2.45) -- (6.5, -2.45) -- (6.5, -2.3);
						
		\draw[] (-.5, -1.6) circle[radius = 0] node[left, scale = .7] {$L$};
		\draw[] (-.75, 1.35) circle[radius = 0] node[left, scale = .7] {$M$};
		
		\draw[] (9.625, 1.225) circle[radius = 0] node[right, scale = .7] {$L$};
		\draw[] (9.875, -1.475) circle[radius = 0] node[right, scale = .7] {$M$};
		\draw[] (2.75, -2.45) circle[radius = 0] node[below, scale = .7]{$N$};
		\draw[] (6.25, -2.45) circle[radius = 0] node[below, scale = .7]{$N$};	
		\draw[] (11, 0) circle [radius = 0] node[]{.};			
		\end{tikzpicture}
		
	\end{center}

	\noindent To deduce \eqref{requation1} from \eqref{lmzrrrrz}, we will repeatedly first $q$-symmetrize the inputs of the $Z$-weights, thereby turning them into the $\mathcal{Z}$ weights of \Cref{zabcd2}. Then, we will use the $q$-exchangeability property given by \Cref{zabcdcd} (together with \Cref{rijkh} and \eqref{r1xyabcd}) to turn the $\mathcal{Z}$-weights into the $\widetilde{\mathcal{R}}$ weights, which will give rise to \eqref{requation1}. 
	
	To that end, multiplying both sides by $q^{\inv (\overleftarrow{\mathfrak{I}_1}) + \inv (\overleftarrow{\mathfrak{J}}_1)}$, summing over $\mathfrak{I}_1 \in \mathcal{M}(\textbf{I}_1)$ and $\mathfrak{J}_1 \in \mathcal{M}(\textbf{J}_1)$, and using \Cref{zabcd2} then gives	
	\begin{flalign} 
	\label{lmzrrrr2}
	\begin{aligned} 	
	& \displaystyle\sum_{\mathfrak{I}_2, \mathfrak{J}_2, \mathfrak{K}_2} \mathcal{Z}_{x, y} \big( \textbf{I}_1, \textbf{J}_1; \overleftarrow{\mathfrak{I}_2}, \mathfrak{J}_2 \big) Z_{x, z} \big( \mathfrak{K}_1, \mathfrak{J}_2; \mathfrak{K}_2, \mathfrak{J}_3 \big) Z_{y, z} (\mathfrak{K}_2, \mathfrak{I}_2; \mathfrak{K}_3, \mathfrak{I}_3) \\
	& \quad = \displaystyle\sum_{\mathfrak{I}_1, \mathfrak{J}_1, \mathfrak{I}_2, \mathfrak{J}_2, \mathfrak{K}_2} q^{\inv (\overleftarrow{\mathfrak{I}_1}) + \inv(\overleftarrow{\mathfrak{J}}_1)} Z_{x, y} \big( \overleftarrow{\mathfrak{I}_1}, \mathfrak{J}_1; \overleftarrow{\mathfrak{I}_2}, \mathfrak{J}_2 \big) Z_{x, z} \big( \mathfrak{K}_1, \mathfrak{J}_2; \mathfrak{K}_2, \mathfrak{J}_3 \big) Z_{y, z} (\mathfrak{K}_2, \mathfrak{I}_2; \mathfrak{K}_3, \mathfrak{I}_3) \\
	& \quad = \displaystyle\sum_{\mathfrak{I}_1, \mathfrak{J}_1, \mathfrak{I}_2, \mathfrak{J}_2, \mathfrak{K}_2} q^{\inv(\overleftarrow{\mathfrak{I}_1}) + \inv(\overleftarrow{\mathfrak{J}}_1)} Z_{y, z} \big( \mathfrak{K}_1, \mathfrak{I}_1; \mathfrak{K}_2, \mathfrak{I}_2 \big) Z_{x, z} (\mathfrak{K}_2, \mathfrak{J}_1; \mathfrak{K}_3, \mathfrak{J}_2) Z_{x, y} \big( \overleftarrow{\mathfrak{I}_2}, \mathfrak{J}_2; \overleftarrow{\mathfrak{I}_3}, \mathfrak{J}_3 \big), 
	\end{aligned} 
	\end{flalign} 
	
	\noindent where we always sum over $\mathfrak{X}_h \in \mathcal{M}(\textbf{X}_h)$, for any indices $X \in \{ I, J, K \}$ and $h \in \{ 1, 3 \}$. Multiplying both sides by $q^{\inv (\mathfrak{K}_1)}$; summing over $\mathfrak{K}_1 \in \mathcal{M} (\textbf{K}_1)$; and using \Cref{zabcd2} again then gives
	\begin{flalign} 
	\label{lmzrrrr3}
	\begin{aligned} 	
	\displaystyle\sum_{\mathfrak{K}_1, \mathfrak{I}_2, \mathfrak{J}_2, \mathfrak{K}_2} & q^{\inv (\mathfrak{K}_1)} \mathcal{Z}_{x, y} \big( \textbf{I}_1, \textbf{J}_1; \overleftarrow{\mathfrak{I}_2}, \mathfrak{J}_2 \big) Z_{x, z} \big( \mathfrak{K}_1, \mathfrak{J}_2; \mathfrak{K}_2, \mathfrak{J}_3 \big) Z_{y, z} (\mathfrak{K}_2, \mathfrak{I}_2; \mathfrak{K}_3, \mathfrak{I}_3) \\
	& = \displaystyle\sum_{\mathfrak{J}_1, \mathfrak{I}_2; \mathfrak{J}_2, \mathfrak{K}_2} q^{\inv (\overleftarrow{\mathfrak{J}}_1)} \mathcal{Z}_{y, z} \big( \textbf{K}_1, \textbf{I}_1; \mathfrak{K}_2, \mathfrak{I}_2 \big) Z_{x, z} (\mathfrak{K}_2, \mathfrak{J}_1; \mathfrak{K}_3, \mathfrak{J}_2) Z_{x, y} \big( \overleftarrow{\mathfrak{I}_2}, \mathfrak{J}_2; \overleftarrow{\mathfrak{I}_3}, \mathfrak{J}_3 \big).
	\end{aligned} 
	\end{flalign} 
	
	\noindent By \Cref{zb1d2}, we may restrict the sum on the left side of \eqref{lmzrrrr3} over $(\mathfrak{I}_2, \mathfrak{J}_2)$ such that, for each index $h \in [m, m + n - 1]$, we have that $\max \big\{ m_h (\mathfrak{I}_2), m_h (\mathfrak{J}_2) \big\} \le 1$. Similarly, we may restrict the sum on the right side of \eqref{lmzrrrr3} over $(\mathfrak{J}_2, \mathfrak{K}_2)$ such that $\max \big\{ m_h (\mathfrak{J}_2), m_h (\mathfrak{K}_2) \big\} \le 1$ for each $h \in [m, m + n - 1]$. 
	
	Now, fix nonnegative integer sequences $\textbf{A}, \textbf{B}, \textbf{C}, \textbf{D} \in \mathbb{Z}_{\ge 0}^{m + n - 1}$ such that $|\textbf{A}| = M = |\textbf{C}|$; $|\textbf{B}| = L = |\textbf{D}|$; and $\max \{ A_h, B_h \} \le 1$ for each $h \in [m, m + n - 1]$. Further fix index sequences $\mathfrak{C} \in \mathcal{M} (\textbf{C})$ and $\mathfrak{D} \in \mathcal{M} (\textbf{D})$. By \Cref{rabcdrabcd} and \eqref{r1xyabcd}, we have 
	\begin{flalign}
		\label{rabcdrabcd1} 
		\mathcal{Z}_{x, y} (\textbf{A}, \textbf{B}; \mathfrak{C}, \mathfrak{D}) = q^{\inv (\mathfrak{C}) + \inv(\overleftarrow{\mathfrak{D}})} \widetilde{\mathcal{R}}_{x, y} (\textbf{A}, \textbf{B}; \textbf{C}, \textbf{D}) \cdot \displaystyle\prod_{i = 0}^{m + n - 1} \displaystyle\frac{(q; q)_{C_i} (q; q)_{D_i}}{(q; q)_{A_i} (q; q)_{B_i}}.
	\end{flalign}
	
	Repeated application of \eqref{rabcdrabcd1} in \eqref{lmzrrrr3}, and also decomposing each sum over an index set $\mathfrak{X}_2$ in \eqref{lmzrrrr3} as a sum over the associated nonnegative composition $\textbf{X}_2 \in \mathbb{Z}_{\ge 0}^{m + n}$ and one over $\mathfrak{X}_2 \in \mathcal{M} (\textbf{X}_2)$, then gives
	\begin{flalign} 
	\label{lmzrrrr4}
	\begin{aligned} 	
	& \displaystyle\sum_{\textbf{I}_2, \textbf{J}_2, \textbf{K}_2} \widetilde{\mathcal{R}}_{x, y} \big( \textbf{I}_1, \textbf{J}_1; \textbf{I}_2, \textbf{J}_2 \big) \\ 
	& \qquad \qquad \times \displaystyle\sum_{\mathfrak{K}_1, \mathfrak{I}_2, \mathfrak{J}_2, \mathfrak{K}_2} q^{\inv (\mathfrak{K}_1) + \inv (\overleftarrow{\mathfrak{I}_2}) + \inv(\overleftarrow{\mathfrak{J}}_2)}  Z_{x, z} \big( \mathfrak{K}_1, \mathfrak{J}_2; \mathfrak{K}_2, \mathfrak{J}_3 \big) Z_{y, z} (\mathfrak{K}_2, \mathfrak{I}_2; \mathfrak{K}_3, \mathfrak{I}_3) \\
	& \qquad = \displaystyle\sum_{\textbf{I}_2, \textbf{J}_2, \textbf{K}_2} \widetilde{\mathcal{R}}_{y, z} \big( \textbf{K}_1, \textbf{I}_1; \textbf{K}_2, \textbf{I}_2 \big)  \\
	& \qquad \qquad \qquad \times \displaystyle\sum_{\mathfrak{J}_1, \mathfrak{I}_2, \mathfrak{J}_2, \mathfrak{K}_2} q^{ \inv (\mathfrak{K}_2) + \inv (\overleftarrow{\mathfrak{I}}_2) + \inv(\overleftarrow{\mathfrak{J}}_1)} Z_{x, z} (\mathfrak{K}_2, \mathfrak{J}_1; \mathfrak{K}_3, \mathfrak{J}_2) Z_{x, y} \big( \overleftarrow{\mathfrak{I}_2}, \mathfrak{J}_2; \overleftarrow{\mathfrak{I}_3}, \mathfrak{J}_3 \big),
	\end{aligned} 
	\end{flalign} 
	
	\noindent where the sums over $\textbf{X}_2$ range in $\mathbb{Z}_{\ge 0}^{m + n}$ and those over $\mathfrak{X}_2$ range in $\mathcal{M}(\textbf{X}_2)$, for each $X \in \{ I, J, K \}$; here, we have also used \Cref{rb1d2} to sum over all $\textbf{X}_2$ (instead of only those satisfying whose $h$-th components are bounded above by $1$ for each $h \in [m, m + n - 1]$). 
	
	For fixed $\textbf{X}_h$ for $X \in \{ I, J, K \}$ and $h \in \{ 1, 2, 3 \}$ (whose $i$-th component is bounded above by $1$ for each $i \in [m, m + n - 1]$), \Cref{zb1d2}, \Cref{rb1d2}, \Cref{rijkh}, and \eqref{rabcdrabcd1} together yield for fixed $\textbf{I}_2, \textbf{J}_2, \textbf{K}_2 \in \{ 0, 1 \}^n$ that 
	\begin{flalign}
	\label{ryzj2} 
	\begin{aligned} 
	\displaystyle\sum_{\mathfrak{K}_1, \mathfrak{I}_2, \mathfrak{J}_2, \mathfrak{K}_2} & q^{\inv (\mathfrak{K}_1) + \inv (\overleftarrow{\mathfrak{I}_2}) + \inv(\overleftarrow{\mathfrak{J}}_2)}  Z_{x, z} \big( \mathfrak{K}_1, \mathfrak{J}_2; \mathfrak{K}_2, \mathfrak{J}_3 \big) Z_{y, z} (\mathfrak{K}_2, \mathfrak{I}_2; \mathfrak{K}_3, \mathfrak{I}_3) \\
	& = \displaystyle\sum_{\mathfrak{J}_2, \mathfrak{K}_2} q^{\inv(\overleftarrow{\mathfrak{I}}_2)}  \mathcal{Z}_{x, z} \big( \textbf{K}_1, \textbf{J}_2; \mathfrak{K}_2, \mathfrak{J}_3 \big) Z_{y, z} (\mathfrak{K}_2, \mathfrak{I}_2; \mathfrak{K}_3, \mathfrak{I}_3) \\
	& = \displaystyle\sum_{\mathfrak{J}_2, \mathfrak{K}_2} q^{\inv (\mathfrak{K}_2) + \inv (\overleftarrow{\mathfrak{I}_2}) + \inv(\overleftarrow{\mathfrak{J}}_3)}  \widetilde{\mathcal{R}}_{x, z} \big( \textbf{K}_1, \textbf{J}_2; \textbf{K}_2, \textbf{J}_3 \big) Z_{y, z} (\mathfrak{K}_2, \mathfrak{I}_2; \mathfrak{K}_3, \mathfrak{I}_3) \\
	& = q^{\inv (\mathfrak{K}_3) + \inv (\overleftarrow{\mathfrak{I}}_3) + \inv (\overleftarrow{\mathfrak{J}}_3)}  \widetilde{\mathcal{R}}_{x, z} \big( \textbf{K}_1, \textbf{J}_2; \textbf{K}_2, \textbf{J}_3 \big) \widetilde{\mathcal{R}}_{y, z} (\textbf{K}_2, \textbf{I}_2; \textbf{K}_3, \textbf{I}_3).
	\end{aligned}
	\end{flalign}

\noindent By similar reasoning,
\begin{flalign}
\label{rxzi2}
\begin{aligned} 
\displaystyle\sum_{\mathfrak{J}_1, \mathfrak{I}_2, \mathfrak{J}_2, \mathfrak{K}_2} & q^{\inv (\mathfrak{K}_2) + \inv (\overleftarrow{\mathfrak{I}}_2) + \inv(\overleftarrow{\mathfrak{J}}_1)} Z_{x, z} (\mathfrak{K}_2, \mathfrak{J}_1; \mathfrak{K}_3, \mathfrak{J}_2) Z_{x, y} \big( \overleftarrow{\mathfrak{I}_2}, \mathfrak{J}_2; \overleftarrow{\mathfrak{I}_3}, \mathfrak{J}_3 \big) \\
& = q^{\inv (\mathfrak{K}_3) + \inv (\overleftarrow{\mathfrak{I}}_3) + \inv (\overleftarrow{\mathfrak{J}}_3)} \widetilde{\mathcal{R}}_{x, z} \big( \textbf{K}_2, \textbf{J}_1; \textbf{K}_3, \textbf{J}_2 \big) \widetilde{\mathcal{R}}_{x, y} ( \textbf{I}_2, \textbf{J}_2; \textbf{I}_3, \textbf{J}_3).
\end{aligned}
\end{flalign}

\noindent Now \eqref{requation1}, and thus the proposition, follows from inserting \eqref{ryzj2} and \eqref{rxzi2} into \eqref{lmzrrrr4}. 
\end{proof}

\chapter{Evaluation of the Fused Weights} 

\label{FusedW}

Observe that the fused weights $\mathcal{R}_{x, y}$ were defined in \Cref{rijkh} as (linear combinations of) certain rectangular partition functions with specialized rapidity parameters. In this chapter we provide a closed form for these quantities, given as \Cref{rxyml} below.

\section{The \texorpdfstring{$U_q \big( \widehat{\mathfrak{sl}} (m) \big)$ Fused Weights}{}}

\label{WeightFusedm0}

In this section we recall from \cite{STR} explicit forms for the fused weights of the $U_q \big( \widehat{\mathfrak{sl}} (m) \big)$-vertex model, thereby providing exact expressions for the weights $\mathcal{R}_{x, y}^{(m; n)} (\textbf{A}, \textbf{B}; \textbf{C}, \textbf{D})$ from \Cref{rijkh} if $n = 0$. To that end, following Theorem C.1.1 of \cite{SVMST}, for any $x, y \in \mathbb{C}$ and $m$-tuples $\lambda = (\lambda_1, \lambda_2, \ldots , \lambda_m) \in \mathbb{Z}_{\ge 0}^m$ and $\mu = (\mu_1, \mu_2, \ldots , \mu_m) \in \mathbb{Z}_{\ge 0}^m$ of nonnegative integers, let 
\begin{flalign}
	\label{lambdamuxyfunction} 
\Phi (\lambda, \mu; x, y) = \displaystyle\frac{(x; q)_{|\lambda|} (x^{-1} y; q)_{|\mu| - |\lambda|}}{(y; q)_{|\mu|}} \left( \displaystyle\frac{y}{x} \right)^{|\lambda|} q^{\varphi (\mu - \lambda, \lambda)} \displaystyle\prod_{i = 1}^m \displaystyle\frac{(q; q)_{\mu_i}}{(q; q)_{\lambda_i} (q; q)_{\mu_i - \lambda_i}}, 
\end{flalign}
\index{0@$\Phi (\lambda, \mu; x, y)$}

\noindent where we have recalled $\varphi$ from \eqref{tufunction}. 

The following proposition that was essentially originally established as equation (7.8) of \cite{STR}, but appears in its below form as Theorem 8.5 of \cite{SVMP}, provides an explicit form for the fused weights $\mathcal{R}_{x, y}^{(m; n)} (\textbf{A}, \textbf{B}; \textbf{C}, \textbf{D})$ if $n = 0$. 

\begin{prop}[{\cite[Theorem 8.5]{SVMP}}]
	
\label{unmweights}

Assume $n = 0$. Fix integers $L, M \ge 1$ and $x, y \in \mathbb{C}$; set $z = \frac{x}{y}$. Let $\textbf{\emph{A}}, \textbf{\emph{B}}, \textbf{\emph{C}}, \textbf{\emph{D}} \in \mathbb{Z}_{\ge 0}^m$ be such that $|\textbf{\emph{A}}| = M = |\textbf{\emph{C}}|$ and $|\textbf{\emph{B}}| = L = |\textbf{\emph{D}}|$; set $\textbf{\emph{X}} = (X_0, X_1, \ldots , X_{m - 1})$ and $\check{\textbf{\emph{X}}} = (X_1, X_2, \ldots,  X_{m - 1})$, for each $X \in \{ A, B, C, D \}$. Then, $\mathcal{R}_{x, y}^{(m; 0)} (\textbf{\emph{A}}, \textbf{\emph{B}}; \textbf{\emph{C}}, \textbf{\emph{D}}) = \omega_{x, y} (\textbf{\emph{A}}, \textbf{\emph{B}}; \textbf{\emph{C}}, \textbf{\emph{D}})$, where $\omega_{x, y} (\textbf{\emph{A}}, \textbf{\emph{B}}; \textbf{\emph{C}}, \textbf{\emph{D}}) = \omega_{x, y}^{(L; M)} (\textbf{\emph{A}}, \textbf{\emph{B}}; \textbf{\emph{C}}, \textbf{\emph{D}})$ is given by 
\begin{flalign}
\label{omegaxy}
\begin{aligned} 
\omega_{x, y} (\textbf{\emph{A}}, \textbf{\emph{B}}; \textbf{\emph{C}}, \textbf{\emph{D}}) & = z^{|\check{\textbf{\emph{D}}}| - |\check{\textbf{\emph{B}}}|} q^{|\check{\textbf{\emph{A}}}| L - |\check{\textbf{\emph{D}}} | M} \textbf{\emph{1}}_{\textbf{\emph{A}} + \textbf{\emph{B}} = \textbf{\emph{C}} + \textbf{\emph{D}}} \\
& \qquad \times  \displaystyle\sum_{\check{\textbf{\emph{P}}}} \Phi (\check{\textbf{\emph{C}}} - \check{\textbf{\emph{P}}}, \check{\textbf{\emph{C}}} + \check{\textbf{\emph{D}}} - \check{\textbf{\emph{P}}}; q^{L - M} z, q^{-M} z)  \Phi ( \check{\textbf{\emph{P}}}, \check{\textbf{\emph{B}}}; q^{-L} z^{-1}, q^{-L}),
\end{aligned}
\end{flalign}

\noindent where $\check{\textbf{\emph{P}}} = (P_1, \ldots , P_{m - 1}) \in \mathbb{Z}_{\ge 0}^{m - 1}$ is summed over $(m - 1)$-tuples of nonnegative integers such that $P_i \le \min \{ B_i, C_i \}$, for each $i \in [1, m - 1]$. 

\end{prop}
\index{0@$\omega_{x, y} (\textbf{A}, \textbf{B}; \textbf{C}, \textbf{D})$}

Let us provide an example for these weights in the case $L = 1$, which will later be useful for us. In what follows, we recall for any set of real numbers $(X_0, X_1, \ldots , X_{\ell}) \in \mathbb{R}^{\ell + 1}$ and indices $0 \le i \le j \le \ell$ that $X_{[i, j]} = \sum_{k = i}^j X_k$. \index{X@$X_{[j, k]}$}

\begin{example}
	
	\label{l1rxyn0}
	
	Assume $L = 1$, and abbreviate $\mathcal{R}_{x, y} (\textbf{A}, b; \textbf{C}, d) = \mathcal{R}_{x, y} \big( \textbf{A}, \textbf{e}_b; \textbf{C}, \textbf{e}_d \big)$,\index{R@$\mathcal{R}_{x, y}^{(m; n)} (\textbf{A}, \textbf{B}; \textbf{C}, \textbf{D})$!$\mathcal{R}_{x, y} (\textbf{A}, \textbf{B}; \textbf{C}, \textbf{D})$!$\mathcal{R}_{x, y} (\textbf{A}, b; \textbf{C}, d)$} for any $\textbf{A}, \textbf{C} \in \mathbb{Z}_{\ge 0}^m$ and $b, d \in \{ 0, 1, \ldots , m - 1 \}$. Then, \Cref{unmweights} (see also equation (B.4.1) of \cite{SVMST}) becomes 
	\begin{flalign*}
	 \mathcal{R}_{x, y} (\textbf{A}, j; \textbf{A}, j) = & \displaystyle\frac{(x - q^{A_j} y) q^{A_{[j + 1, m - 1]}}}{x - q^M y}; \qquad  \mathcal{R}_{x, y} (\textbf{A}, j; \textbf{A} + \textbf{e}_j - \textbf{e}_k, k) = \displaystyle\frac{(1 - q^{A_k}) q^{A_{[k + 1, m - 1]}} x}{x - q^M y}; \\
	 & \mathcal{R}_{x, y} (\textbf{A}, k; \textbf{A} + \textbf{e}_k - \textbf{e}_j, j) = \displaystyle\frac{(1 - q^{A_j}) q^{A_{[j + 1, m - 1]}} y}{x - q^M y},
	\end{flalign*} 
	
	\noindent for any $0 \le j < k \le M - 1$, and $\mathcal{R}_{x, y} (\textbf{A}, b; \textbf{C}, d) = 0$ for any $(\textbf{A}, b; \textbf{C}, d)$ not of the above form.
	
\end{example}

\section{Special Cases of the Fused Weights} 

\label{Horizontal} 

In this section we evaluate the fused weights $\mathcal{R}_{x, y}^{(m; n)} (\textbf{A}, \textbf{B}; \textbf{C}, \textbf{D})$ for general $n \ge 0$ under two special cases. The first assumes that $\textbf{A}, \textbf{B}, \textbf{C}, \textbf{D}$ do not share any fermionic colors, in which case the weight is given by \eqref{omegaxy}, as in the $n = 0$ situation; the second assumes that $L = 1$, in which case the weight is given by a modification of those in \Cref{l1rxyn0}.

\begin{lem}
	
	\label{zrmn}
	
	Fix $x, y \in \mathbb{C}$ and integers $L, M \ge 1$. Let $\textbf{\emph{A}}, \textbf{\emph{B}}, \textbf{\emph{C}}, \textbf{\emph{D}} \in \mathbb{Z}_{\ge 0}^{m + n}$, whose coordinates are indexed by $[0, m + n - 1]$; assume they satisfy $|\textbf{\emph{A}}| = M = |\textbf{\emph{C}}|$ and $|\textbf{\emph{B}}| = L = |\textbf{\emph{D}}|$. If there does not exist any index $h \in [m, m + n - 1]$ such that $\max \{ A_h, B_h, C_h, D_h \} \ge 1$, then $\mathcal{R}_{x, y} (\textbf{\emph{A}}, \textbf{\emph{B}}, \textbf{\emph{C}}, \textbf{\emph{D}}) = \omega_{x, y} (\textbf{\emph{A}}, \textbf{\emph{B}}, \textbf{\emph{C}}, \textbf{\emph{D}})$, where we recall the weight $\omega_{x, y}$ from \eqref{omegaxy}. 
	 
\end{lem} 

\begin{proof}
	
	Recall from \Cref{zxy1} the partition function $Z_{x, y}^{(m; n)}$. By \Cref{zabcd2}, \Cref{rijkh}, and \Cref{unmweights}, it suffices to show $Z_{x, y}^{(m; n)} (\mathfrak{A}, \mathfrak{B}; \mathfrak{C}, \mathfrak{D}) = Z_{x, y}^{(m + n; 0)} (\mathfrak{A}, \mathfrak{B}; \mathfrak{C}, \mathfrak{D})$, for any $\mathfrak{A} \in \mathcal{M} (\textbf{A})$, $\mathfrak{B} \in \mathcal{M} (\textbf{B})$, $\mathfrak{C} \in \mathcal{M} (\textbf{C})$, and $\mathfrak{D} \in \mathcal{M} (\textbf{D})$. To that end, observe since there does not exist any $h \in [m, m + n - 1]$ such that $A_h, B_h, C_h, D_h \ge 1$, there does not exist any arrow configuration $(v_{i, j}, u_{i, j}; v_{i, j + 1}, u_{i + 1, j})$ on the right side of \eqref{zabcd} in the expansion of $Z_{x, y}^{(m; n)}$ that is of the form $(h, h; h, h)$ for some $h \in [m, m + n - 1]$. Thus, by \Cref{m1m2rr}, we have that $R_z^{(m; n)} (v_{i, j}, u_{i, j}; v_{i, j + 1}, u_{i + 1, j}) = R_z^{(m + n; 0)} (v_{i, j}, u_{i, j}; v_{i, j + 1}, u_{i + 1, j})$ for any such vertex weight appearing there. Inserting this fact into \eqref{zabcd} then yields $Z_{x, y}^{(m; n)} (\mathfrak{A}, \mathfrak{B}; \mathfrak{C}, \mathfrak{D}) = Z_{x, y}^{(m + n; 0)} (\mathfrak{A}, \mathfrak{B}; \mathfrak{C}, \mathfrak{D})$, which as mentioned above implies the lemma. 
\end{proof}

Now let us evaluate the fused weights $\mathcal{R}_{x, y} (\textbf{A}, \textbf{B}; \textbf{C}, \textbf{D})$ in the case $L = 1$. Letting $\textbf{B} = \textbf{e}_b$ and $\textbf{D} = \textbf{e}_d$ for some indices $b, d \in [0, m + n - 1]$, this is done by \Cref{zrmn} unless $b = d \in [m, m + n - 1]$ and $A_b = C_b = 1$. This remaining situation is addressed through the following lemma.

\begin{lem}
	
	\label{rz1vmn}
	
	Fix $x, y \in \mathbb{C}$ and an integer $M \ge 1$. Let $\textbf{\emph{A}} \in \mathbb{Z}_{\ge 0}^{m + n}$ satisfy $|\textbf{\emph{A}}| = M$, and denote $\textbf{\emph{A}} = (A_0, A_1, \ldots , A_{m + n - 1})$; assume that $A_i \le 1$ for each index $i \in [m, m + n - 1]$. If $h \in [m, m + n - 1]$ is such that $A_h = 1$, then 
	\begin{flalign*}
	\mathcal{R}_{x, y} (\textbf{\emph{A}}, \textbf{\emph{e}}_h; \textbf{\emph{A}}, \textbf{\emph{e}}_h) = q^{A_{[h + 1, m + n - 1]}} \displaystyle\frac{y - qx}{x - q^M y}. 
	\end{flalign*}
	
\end{lem}

\begin{figure}

	\begin{center}

		\begin{tikzpicture}[
		>=stealth,
		auto,
		style={
			scale = 1
		}
		]
		
		\draw[->, red, thick] (-1, 0) -- (0, 0);
		\draw[->, red, thick] (3, 0) -- (4, 0);
		
		\draw[->, blue, thick] (0, -1) -- (0, 0);
		\draw[->, green, thick] (1, -1) -- (1, 0);
		\draw[->, blue, thick] (2, -1) -- (2, 0);
		\draw[->, red, thick] (3, -1) -- (3, 0);
		
		\draw[-, red, thick, dashed] (0, 0) -- (1, 0);
		\draw[-, red, thick, dashed] (1, 0) -- (2, 0);
		\draw[-, red, thick, dashed] (2, 0) -- (3, 0);
		\draw[-, red, thick, dashed] (3, 0) -- (4, 0);

		\draw[->, blue, thick] (0, 0) -- (0, 1);
		\draw[->, green, thick] (1, 0) -- (1, 1);
		\draw[->, blue, thick] (2, 0) -- (2, 1);
		\draw[->, red, thick] (3, 0) -- (3, 1);

		\draw[] (0, -1) circle[radius = 0] node[below, scale = .7]{$q^3 y$};
		\draw[] (1, -1) circle[radius = 0] node[below, scale = .7]{$q^2 y$};
		\draw[] (2, -1) circle[radius = 0] node[below = 2, scale = .7]{$q y$};
		\draw[] (3, -1) circle[radius = 0] node[below = 2, scale = .7]{$y$};
		
		\draw[] (-1, 0) circle[radius = 0] node[left, scale = .7]{$x$};
		
		\draw[] (-.6, 0) circle[radius = 0] node[above, scale = .7]{$h$};
		
		\draw[] (2.5, 0) circle[radius = 0] node[above, scale = .7]{$h$};
		\draw[] (3.5, 0) circle[radius = 0] node[above, scale = .7]{$h$};
		\draw[] (1.5, 0) circle[radius = 0] node[above, scale = .7]{$h$};
		\draw[] (.5, 0) circle[radius = 0] node[above, scale = .7]{$h$};
		
		\draw[] (0, -.6) circle[radius = 0] node[right, scale = .6]{$a_1' = a_1$};
		\draw[] (1, -.6) circle[radius = 0] node[right, scale = .6]{$a_2' = a_2$};
		\draw[] (2, -.6) circle[radius = 0] node[right, scale = .6]{$a_3' = a_3$};
		\draw[] (3, -.6) circle[radius = 0] node[right, scale = .6]{$a_4' = h$};
		
		\draw[] (0, .5) circle[radius = 0] node[right, scale = .7]{$a_1$};
		\draw[] (1, .5) circle[radius = 0] node[right, scale = .7]{$a_2$};
		\draw[] (2, .5) circle[radius = 0] node[right, scale = .7]{$a_3$};
		\draw[] (3, .52) circle[radius = 0] node[right, scale = .7]{$h$};

		\draw[] (0, -1.625) -- (0, -1.75) -- (3, -1.75) -- (3, -1.625); 
		
		\draw[] (1.5, -1.75) circle[radius = 0] node[below, scale = .7]{$M$};

		\end{tikzpicture}
		
	\end{center}
	
	\caption{\label{rxyz} Shown above is a diagrammatic interpretation for \eqref{zxyaa}. } 	
\end{figure}

\begin{proof}
	
	Let $\mathfrak{A}= (a_1, a_2, \ldots , a_M) \in \mathcal{M} (\textbf{A})$ denote the sequence of indices in $[0, m + n -1]$ such that $a_M = h$ and $a_1 \le a_2 \le \cdots \le a_{M - 1}$. Then \Cref{rijkh}, \Cref{zabcd2}, and \Cref{zabcdcd} together imply that 
	\begin{flalign}
	\label{rxyaa}
	\mathcal{R}_{x, y} (\textbf{A}, \textbf{e}_h; \textbf{A}, \textbf{e}_h) = q^{-\inv (\mathfrak{A})} \displaystyle\sum_{\mathfrak{A}' \in \mathcal{M}(\textbf{A})} q^{\inv (\mathfrak{A}')} Z_{x, y} (\mathfrak{A}', h; \mathfrak{A}, h). 
	\end{flalign}  
	
	\noindent Expanding $Z_{x, y} (\mathfrak{A}', h; \mathfrak{A}, h)$ as on the right side of \eqref{zabcd}, and then using arrow conservation, the fact that $a_j \ne h$ for any $j < M$, and induction on $M - i$, we deduce that each nonzero vertex weight $R_z (v_{i, 1}, u_{i, 1}; v_{i, 2}, u_{i + 1, 1})$ there is associated with the arrow configuration $(a_i, h; a_i, h)$; we refer to \Cref{rxyz} for a depiction. Thus, only the $\mathfrak{A}' = \mathfrak{A}$ summand on the right side of \eqref{rxyaa} is nonzero, which yields
	\begin{flalign}
	\label{rxyaa1}
	\mathcal{R}_{x, y} (\textbf{A}, \textbf{e}_h; \textbf{A}, \textbf{e}_h) = Z_{x, y} (\mathfrak{A}, h; \mathfrak{A}, h).
	\end{flalign}
	
	\noindent Further let $\widetilde{\mathfrak{A}} = (a_1, a_2, \ldots , a_{M - 1})$ denote the $(M - 1)$-tuple obtained by removing the last entry of $\mathfrak{A}$. Using the fact that the rightmost arrow configuration $(v_{M, 1}, u_{M, 1}; v_{M, 2}, u_{M + 1, 1})$ in the expansion of $Z_{x, y} (\mathfrak{A}, h; \mathfrak{A}, h)$ given by the right side of \eqref{zabcd} is $(h, h; h, h)$ (since $a_M = h$), we find that 
	\begin{flalign}
	\label{zxyaa}
	Z_{x, y} (\mathfrak{A}, h; \mathfrak{A}, h) = Z_{x, qy} \big( \widetilde{\mathfrak{A}}, h; \widetilde{\mathfrak{A}}, h \big) R_{y / x} (h, h; h, h).
	\end{flalign}
	
	Letting $\widetilde{\textbf{A}} \in \mathbb{Z}_{\ge 0}^{m + n}$ be such that $\widetilde{\mathfrak{A}} \in \mathcal{M} \big( \widetilde{\textbf{A}} \big)$, we have that $h \notin \widetilde{\mathfrak{A}}$ since $A_h = 1$ and $a_M = h$. Thus, arrow conservation again implies that, in the expansion of $Z_{x, y} \big( \widetilde{\mathfrak{A}}, h; \widetilde{\mathfrak{A}}, h \big)$ as on the right side of \eqref{zabcd}, each arrow configuration $(v_{i, 1}, u_{i, 1}; v_{i, 2}, u_{i + 1, 1})$ is of the form $(a_i, h; a_i, h)$. Thus, similarly to in \eqref{rxyaa1}, we find by \Cref{rijkh}, \Cref{zabcd2}, and \Cref{zabcdcd} (and also the fact that $\inv \big( \widetilde{\mathfrak{A}} \big) = 0$) that 
	\begin{flalign*}
	Z_{x, qy} \big( \widetilde{\mathfrak{A}}, h; \widetilde{\mathfrak{A}}, h \big) = \displaystyle\sum_{\widetilde{\mathfrak{A}}' \in \mathcal{M}(\widetilde{\textbf{A}})} q^{\inv (\widetilde{\mathfrak{A}'})} Z_{x, qy} \big( \widetilde{\mathfrak{A}}', h; \widetilde{\mathfrak{A}}, h \big) = \mathcal{R}_{x, qy} \big( \widetilde{\textbf{A}}, \textbf{e}_h; \widetilde{\textbf{A}}, \textbf{e}_h \big).
	\end{flalign*}
	
	\noindent Inserting this into \eqref{zxyaa} and using the identity (from \eqref{rzi}) that $R_{y / x} (h, h; h, h) = \frac{y - qx}{x - qy}$ yields
	\begin{flalign}
	\label{rxyah} 
	\mathcal{R}_{x, y} (\textbf{A}, h; \textbf{A}, h) = \mathcal{R}_{x, qy} \big( \widetilde{\textbf{A}}, \textbf{e}_h; \widetilde{\textbf{A}}, \textbf{e}_h \big) \displaystyle\frac{y - qx}{x - qy}.
	\end{flalign}  
	
	Now observe, upon denoting $\widetilde{\textbf{A}} = \big( \widetilde{A}_0, \widetilde{A}_1, \ldots , \widetilde{A}_{m + n - 1} \big)$, we have that $\widetilde{A}_h = 0$. Thus, \Cref{zrmn}, \Cref{unmweights}, and \Cref{l1rxyn0} together imply that
	\begin{flalign*}
	\mathcal{R}_{x, qy} \big( \widetilde{\textbf{A}}, h; \widetilde{\textbf{A}}, h \big) = q^{\widetilde{A}_{[h + 1, m + n - 1]}} \displaystyle\frac{x - qy}{x - q^M y} =  q^{A_{[h + 1, m + n - 1]}} \displaystyle\frac{x - qy}{x - q^M y}.
	\end{flalign*}
	
	\noindent Inserting this into \eqref{rxyah} then yields the lemma.	
\end{proof}

\section{General Fused Weights} 

\label{Vertical} 

In this section we evaluate the general fused weights $\mathcal{R}_{x, y} (\textbf{A}, \textbf{B}; \textbf{C}, \textbf{D})$ through \Cref{rxyml} below. This was done in \Cref{zrmn} assuming that $\textbf{A}, \textbf{B}, \textbf{C}, \textbf{D}$ do not share any fermionic colors. If instead these four sequences share a fermionic color $h \in [m, m + n - 1]$, then the following lemma provides a recursion for the weight $\mathcal{R}_{x, y} (\textbf{A}, \textbf{B}; \textbf{C}, \textbf{D})$ that proceeds by removing a shared fermionic color $h$ from $\textbf{B}$ and $\textbf{D}$. 

\begin{prop}

\label{rxy1l}

Fix $x, y \in \mathbb{C}$ and $L, M \in \mathbb{Z}_{\ge 1}$. Let $\textbf{\emph{A}}, \textbf{\emph{B}}, \textbf{\emph{C}}, \textbf{\emph{D}} \in \mathbb{Z}_{\ge 0}^{m + n}$, whose coordinates are indexed by $[0, m + n - 1]$; suppose they satisfy $|\textbf{\emph{A}}| = M = |\textbf{\emph{C}}|$ and $|\textbf{\emph{B}}| = L = |\textbf{\emph{D}}|$. Assume that $\max \{ A_i, B_i \} \le 1$ for each index $i \in [m, m + n - 1]$ and that there exists an index $h \in [m, m + n - 1]$ such that $A_h = B_h = C_h = D_h = 1$. Letting $\widetilde{\textbf{\emph{X}}} = \textbf{\emph{X}} - \textbf{\emph{e}}_h$ for each $X \in \{ B, D \}$, we have 
\begin{flalign*}
\mathcal{R}_{x, y} (\textbf{\emph{A}}, \textbf{\emph{B}}; \textbf{\emph{C}}, \textbf{\emph{D}}) = q^{A_{[h + 1, m + n - 1]}} \displaystyle\frac{y - q^L x}{q^{L - 1} x - q^M y} \mathcal{R}_{x, y} \big( \textbf{\emph{A}}, \widetilde{\textbf{\emph{B}}}; \textbf{\emph{C}}, \widetilde{\textbf{\emph{D}}} \big).
\end{flalign*} 

\end{prop} 

\begin{proof} 
	
	By \Cref{zabcd2}, \Cref{zabcdcd}, and \Cref{rijkh}, we have for any fixed $\mathfrak{C} \in \mathcal{M} (\textbf{C})$ and $\mathfrak{D} \in \mathcal{M}( \textbf{D})$ that 
	\begin{flalign}
	\label{rxysumqz} 
	\mathcal{R}_{x, y} (\textbf{A}, \textbf{B}; \textbf{C}, \textbf{D}) = q^{- \inv (\overleftarrow{\mathfrak{D}}) - \inv (\mathfrak{C})}  \displaystyle\prod_{i = 0}^{m + n - 1} \displaystyle\frac{(q; q)_{A_i} (q; q)_{B_i}}{(q; q)_{C_i} (q; q)_{D_i}} \cdot \displaystyle\sum_{\substack{\mathfrak{A} \in \mathcal{M} (\textbf{A}) \\ \mathfrak{B} \in \mathcal{M} (\textbf{B})}} q^{\inv (\mathfrak{A}) + \inv (\overleftarrow{\mathfrak{B}})} Z_{x, y} (\mathfrak{A}, \mathfrak{B}; \mathfrak{C}, \mathfrak{D}).
	\end{flalign}
	
	\noindent For each $X \in \{B, D \}$, let $\mathfrak{X} = (x_1, x_2, \ldots , x_L)$ and $\widetilde{\mathfrak{X}}' = (x_1, x_2, \ldots , x_{L - 1})$. Then, we have by the definition of $Z_{x, y}$ (in particular, the first statement of \eqref{zabcdzabcd}) that 
	\begin{flalign*} 
	Z_{x, y} (\mathfrak{A}, \mathfrak{B}; \mathfrak{C}, \mathfrak{D}) = \displaystyle\sum_{\mathfrak{I}} Z_{x, y} \big( \mathfrak{A}, \widetilde{\mathfrak{B}}'; \mathfrak{I}, \widetilde{\mathfrak{D}}' \big) Z_{q^{L - 1} x, y} \big( \mathfrak{I}, b_L; \mathfrak{C}, d_L),
	\end{flalign*}
	
	\noindent where the sum is over all sequences $\mathfrak{I} = (i_1, i_2, \ldots , i_M)$ of indices in $[0, m + n - 1]$. Inserting this into \eqref{rxysumqz} yields
	\begin{flalign}
	\label{rxysumz2}
	\begin{aligned}
	\mathcal{R}_{x, y} (\textbf{A}, \textbf{B}; \textbf{C}, \textbf{D}) & \cdot \displaystyle\prod_{i = 0}^{m + n - 1} \displaystyle\frac{(q; q)_{C_i} (q; q)_{D_i}}{(q; q)_{A_i} (q; q)_{B_i}} \\
	 & = q^{- \inv (\overleftarrow{\mathfrak{D}}) - \inv (\mathfrak{C})}  \displaystyle\sum_{\mathfrak{A} \in \mathcal{M} (\textbf{A})} \displaystyle\sum_{b_L = 0}^{m + n - 1} \displaystyle\sum_{\widetilde{\mathfrak{B}}' \in \mathcal{M} (\textbf{B} - \textbf{e}_{b_L})} q^{\inv (\mathfrak{A}) + \inv (\overleftarrow{\mathfrak{B}})} \\
	& \qquad \qquad \qquad \qquad \qquad \times \displaystyle\sum_{\mathfrak{I}} Z_{x, y} \big( \mathfrak{A}, \widetilde{\mathfrak{B}}'; \mathfrak{I}, \widetilde{\mathfrak{D}}' \big) Z_{q^{L - 1} x, y} ( \mathfrak{I}, b_L; \mathfrak{C}, d_L) \\
	& = q^{- \inv (\overleftarrow{\mathfrak{D}}) - \inv (\mathfrak{C})} \displaystyle\sum_{b_L = 0}^{m + n - 1} q^{B_{[0, b_L - 1]}} \displaystyle\sum_{\mathfrak{I}} Z_{q^{L - 1} x, y} ( \mathfrak{I}, b_L; \mathfrak{C}, d_L) \\
	& \qquad \qquad \qquad \qquad \times \displaystyle\sum_{\mathfrak{A} \in \mathcal{M} (\textbf{A})}  \displaystyle\sum_{\widetilde{\mathfrak{B}}' \in \mathcal{M} (\textbf{B} - \textbf{e}_{b_L})} q^{\inv (\mathfrak{A}) + \inv (\overleftarrow{\widetilde{\mathfrak{B}}'})} Z_{x, y} \big( \mathfrak{A}, \widetilde{\mathfrak{B}}'; \mathfrak{I}, \widetilde{\mathfrak{D}}' \big) \\
	& = q^{- \inv (\overleftarrow{\mathfrak{D}}) - \inv (\mathfrak{C})} \displaystyle\sum_{b_L = 0}^{m + n - 1} q^{B_{[0, b_L - 1]}} \\
	& \qquad \qquad \qquad \qquad \times \displaystyle\sum_{\mathfrak{I}} Z_{q^{L - 1} x, y} ( \mathfrak{I}, b_L; \mathfrak{C}, d_L) \mathcal{Z}_{x, y} \big( \textbf{A}, \textbf{B} - \textbf{e}_{b_L}; \mathfrak{I}, \widetilde{\mathfrak{D}}' \big),
	\end{aligned} 
	\end{flalign}
	
	\noindent where in the last equalities we used the identity $\inv \big(\overleftarrow{\widetilde{\mathfrak{B}}'} \big) = \inv (\overleftarrow{\mathfrak{B}}) - B_{[0, b_L - 1]}$ and applied \Cref{zabcd2}.
	
	Now, let us suppose that $\mathfrak{D}\in \mathcal{M}(\textbf{D})$ is such that $d_L = h$. Since $\max \{ A_j, B_j \} \le 1$ for each $j \in [m, m + n - 1]$, we have by \Cref{zb1d2} that $\mathcal{Z}_{x, y} \big( \textbf{A}, \textbf{B} - \textbf{e}_{b_L}; \mathfrak{I}, \widetilde{\mathfrak{D}}' \big) = 0$ unless $I_h \le 1$, where $\textbf{I} = (I_0, I_1, \ldots , I_{m + n - 1}) \in \mathbb{Z}_{\ge 0}^{m + n}$ is such that $\mathfrak{I} \in \mathcal{M}(\textbf{I})$. Furthermore, since $d_L = h$ and $h \in \mathfrak{C}$ (as $C_h = 1$), we have by arrow conservation that $Z_{q^{L - 1} x, y} (\mathfrak{I}, b_L; \mathfrak{C}, d_L) = 0$ unless $I_h = 2 - \textbf{1}_{b_L = h}$. Thus, for both of these quantities to be nonzero, we must have that $b_L = h$, and so $\textbf{B} - \textbf{e}_{b_L} = \widetilde{\textbf{B}}$. Again by arrow conservation, this implies that $\textbf{I} = \textbf{C}$ for $Z_{x, y} (\mathfrak{I}, b_L; \mathfrak{C}, d_L) = Z_{x, y} (\mathfrak{I}, h; \mathfrak{C}, h)$ to be nonzero. Inserting these facts into \eqref{rxysumz2}, we deduce that 
	\begin{flalign}
	\label{rxyabcdzsum}
	\begin{aligned}
	\mathcal{R}_{x, y} & (\textbf{A}, \textbf{B}; \textbf{C}, \textbf{D}) \cdot \displaystyle\prod_{i = 0}^{m + n - 1} \displaystyle\frac{(q; q)_{C_i} (q; q)_{D_i}}{(q; q)_{A_i} (q; q)_{B_i}} \\
	& = q^{B_{[0, h - 1]} - \inv (\overleftarrow{\mathfrak{D}}) - \inv (\mathfrak{C})} \displaystyle\sum_{\mathfrak{I} \in \mathcal{M} (\textbf{C})} Z_{q^{L - 1} x, y} ( \mathfrak{I}, h; \mathfrak{C}, h) \mathcal{Z}_{x, y} \big( \textbf{A}, \widetilde{\textbf{B}}; \mathfrak{I}, \widetilde{\mathfrak{D}}' \big) \\
	& = q^{B_{[0, h - 1]} - \inv (\overleftarrow{\mathfrak{D}}) + \inv (\overleftarrow{\widetilde{\mathfrak{D}}'})} \mathcal{R}_{x, y} \big( \textbf{A}, \widetilde{\textbf{B}}; \textbf{C}, \widetilde{\textbf{D}} \big) \cdot \displaystyle\prod_{i = 0}^{m + n - 1} \displaystyle\frac{(q; q)_{C_i} (q; q)_{D_i}}{(q; q)_{A_i} (q; q)_{B_i}} \\
	& \qquad \times \displaystyle\sum_{\mathfrak{I} \in \mathcal{M} (\textbf{C})} q^{\inv (\mathfrak{I}) - \inv (\mathfrak{C})} Z_{q^{L - 1} x, y} ( \mathfrak{I}, h; \mathfrak{C}, h) \\
	& = q^{B_{[0, h - 1]} - D_{[0, h - 1]}} \mathcal{R}_{x, y} \big( \textbf{A}, \widetilde{\textbf{B}}; \textbf{C}, \widetilde{\textbf{D}} \big) \mathcal{R}_{q^{L - 1} x, y} ( \textbf{C}, \textbf{e}_h; \textbf{C}, \textbf{e}_h) \cdot \displaystyle\prod_{i = 0}^{m + n - 1} \displaystyle\frac{(q; q)_{C_i} (q; q)_{D_i}}{(q; q)_{A_i} (q; q)_{B_i}}.
	\end{aligned} 
	\end{flalign}
	
	\noindent Here, in the second equality we applied \Cref{zabcdcd} and \Cref{rijkh}, together with the facts that $B_h = 1 = D_h$ and that the $h$-th coordinates of $\widetilde{\textbf{B}}$ and $\widetilde{\textbf{D}}$ are equal to $0$. In the last, we applied those two statements, \Cref{zabcd2}, and the identity $\inv (\overleftarrow{\mathfrak{D}}) = \inv \big( \overleftarrow{\widetilde{\mathfrak{D}}}' \big) + D_{[0, h - 1]}$ (as $d_L = h$). 
	
	By \Cref{rz1vmn}, we have that 
	\begin{flalign}
	\label{rql1xy}
	\mathcal{R}_{q^{L - 1} x, y} (\textbf{C}, \textbf{e}_h; \textbf{C}, \textbf{e}_h) = q^{C_{[h + 1, m + n - 1]}} \displaystyle\frac{y - q^L x}{q^{L - 1} x - q^M y}.
	\end{flalign}
	
	\noindent Since arrow conservation and the fact that $B_h = 1 = D_h$ together imply that $\mathcal{R}_{x, y} (\textbf{A}, \textbf{B}; \textbf{C}, \textbf{D}) = 0 = \mathcal{R}_{x, y} \big( \textbf{A}, \widetilde{\textbf{B}}; \textbf{C}, \widetilde{\textbf{D}} \big)$ unless 
	\begin{flalign*}
	A_{[h + 1, m + n - 1]} - B_{[0, h - 1]} & = A_{[h + 1, m + n - 1]} + B_{[h + 1, m + n - 1]} - L + 1 \\
	& = C_{[h + 1, m + n - 1]} + D_{[h + 1, m + n - 1]} - L + 1 = C_{[h + 1, m + n - 1]} - D_{[0, h - 1]},
	\end{flalign*} 
	
	\noindent the proposition follows from inserting \eqref{rql1xy} into \eqref{rxyabcdzsum}. 
	\end{proof}
	
	Repeated application of \Cref{rxy1l} reduces a general fused weight $\mathcal{R}_{x, y} (\textbf{A}, \textbf{B}; \textbf{C}, \textbf{D})$ to one satisfying the conditions of \Cref{zrmn}; this gives rise to the following theorem.
	
	\begin{thm}
		
		\label{rxyml}

		Fix $x, y \in \mathbb{C}$ and integers $L, M \ge 1$. Let $\textbf{\emph{A}}, \textbf{\emph{B}}, \textbf{\emph{C}}, \textbf{\emph{D}} \in \mathbb{Z}_{\ge 0}^{m + n}$, with coordinates indexed by $[0, m + n - 1]$, satisfy $|\textbf{\emph{A}}| = M = |\textbf{\emph{C}}|$, $|\textbf{\emph{B}}| = L = |\textbf{\emph{D}}|$; suppose that $\max \{ A_j, B_j \} \le 1$ for each index $j \in [m, m + n - 1]$. 
		
		\begin{enumerate}
			\item If $\max \{ C_j, D_j \} \ge 2$ for some index $j \in [m, m + n - 1]$, then $\mathcal{R}_{x, y} (\textbf{\emph{A}}, \textbf{\emph{B}}; \textbf{\emph{C}}, \textbf{\emph{D}}) = 0$. 
			
			\item Otherwise, define $\textbf{\emph{V}} = (V_0, V_1, \ldots,  V_{m + n - 1}) \in \{ 0, 1 \}^{m + n}$ by setting $V_i = 0$ for $i \in [0, m - 1]$ and $V_j = \min \{ A_j, B_j, C_j, D_j \}$ for $j \in [m, m + n - 1]$; let $v = |\textbf{\emph{V}}|$; and recall $\omega_{x, y}^{(L; M)}$ from \eqref{omegaxy}. Then, denoting $z = \frac{x}{y}$, we have
		\begin{flalign*}
		\mathcal{R}_{x, y} (\textbf{\emph{A}}, \textbf{\emph{B}}; \textbf{\emph{C}}, \textbf{\emph{D}}) & = (-1)^v q^{\sum_{h: V_h = 1} A_{[h + 1, m + n - 1]} - Mv} \displaystyle\frac{(q^{L - v + 1} z; q)_v}{(q^{L - M - v} z; q)_v} \\
		& \qquad \times \omega_{x, y}^{(L - v; M)} (\textbf{\emph{A}}, \textbf{\emph{B}} - \textbf{\emph{V}}; \textbf{\emph{C}}, \textbf{\emph{D}} - \textbf{\emph{V}}).
		\end{flalign*} 
		
		\end{enumerate}
		
	\end{thm} 

\begin{proof}
	
	The first statement of the theorem follows from \Cref{rb1d2}, so it remains to establish the latter. To that end, observe by applying \Cref{rxy1l} $v$ times (once for each $h \in \textbf{V}$, in increasing order) that 
	\begin{flalign}
	\label{rxyabcdv} 
	\begin{aligned} 
	\mathcal{R}_{x, y} (\textbf{A}, \textbf{B}; \textbf{C}, \textbf{D}) & = \mathcal{R}_{x, y} (\textbf{A}, \textbf{B} - \textbf{V}; \textbf{C}, \textbf{D} - \textbf{V}) q^{\sum_{h: V_h = 1} A_{[h + 1, m + n - 1]}} \displaystyle\prod_{j = 0}^{v - 1} \displaystyle\frac{y - q^{L - j} x}{q^{L - j - 1} x - q^M y} \\
	& = (-1)^v q^{\sum_{h: V_h = 1} A_{[h + 1, m + n - 1]} - Mv} \mathcal{R}_{x, y} (\textbf{A}, \textbf{B} - \textbf{V}; \textbf{C}, \textbf{D} - \textbf{V}) \displaystyle\prod_{j = 0}^{v - 1} \displaystyle\frac{1 - q^{L - j} z}{1 - q^{L - M - j - 1} z}.
	\end{aligned}
	\end{flalign}
	
	\noindent By \Cref{zrmn}, we have that $\mathcal{R}_{x, y} (\textbf{A}, \textbf{B} - \textbf{V}; \textbf{C}, \textbf{D} - \textbf{V}) = \omega_{x, y}^{(L - v; M)} (\textbf{A}, \textbf{B} - \textbf{V}; \textbf{C}, \textbf{D} - \textbf{V})$, since $\min \{ A_j, B_j - V_j, C_j, D_j - V_j \} = 0$ for each $j \in [m, m + n - 1]$. Inserting this into \eqref{rxyabcdv} yields the theorem. 
\end{proof}

\chapter{The \texorpdfstring{$U_q \big( \widehat{\mathfrak{sl}} (1 | n) \big)$}{} Specialization}

\label{SymmetricBranching} 

In this chapter we give analytic continuations for the $m = 1$ specializations of the fused weights $\mathcal{R}_{x, y}^{(m; n)}$ from \Cref{rijkh} and provide a color merging result for them.

\section{The \texorpdfstring{$U_q \big( \widehat{\mathfrak{sl}} (1 | n) \big)$}{} Weights and Analytic Continuation}

\label{Weights1n} 

In this section we consider the $m = 1$ specializations of the fused weights $\mathcal{R}_{x, y}$ from \Cref{rijkh}; this corresponds to the case when one arrow color is bosonic and the remaining ones are all fermionic. The following definition provides a family of vertex weights $W_z$, which we will show as \Cref{wabcdsxz} below are analytic continuations (in $q^{-L / 2}$ and $q^{-M / 2}$) of these specialized fused $\mathcal{R}_{x, y} (\textbf{A}, \textbf{B}; \textbf{C}, \textbf{D})$ weights if each fermionic color appears at most once in $\textbf{A}$ and $\textbf{B}$.

In what follows, we recall the function $\varphi$ from \eqref{tufunction}.

\begin{definition}
	
	\label{wabcdrsxy}
	
	Fix an integer $n \ge 1$; complex numbers $r, s, z \in \mathbb{C}$; and $n$-tuples $\textbf{A}, \textbf{B}, \textbf{C}, \textbf{D} \in \mathbb{Z}_{\ge 0}^n$. For each index $X \in \{ A, B, C, D \}$, denote $\textbf{X} = (X_1, X_2, \ldots , X_n)$ and set $|\textbf{X}| = x$. Further define $\textbf{V} = (V_1, V_2, \ldots , V_n) \in \mathbb{Z}_{\ge 0}^n$ by setting $V_j = \min \{ A_j, B_j, C_j, D_j \}$ for each $j \in [1, n]$, and let $|\textbf{V}| = v$. If $\textbf{A} + \textbf{B} = \textbf{C} + \textbf{D}$ and $\textbf{A}, \textbf{B}, \textbf{C}, \textbf{D} \in \{ 0, 1 \}^n$, then define\index{W@$W_z (\textbf{A}, \textbf{B}; \textbf{C}, \textbf{D} \boldsymbol{\mid} r, s)$; fused weight} $W_z (\textbf{A}, \textbf{B}; \textbf{C}, \textbf{D} \boldsymbol{\mid} r, s) = W_{z; q}^{(1; n)} (\textbf{A}, \textbf{B}; \textbf{C}, \textbf{D} \boldsymbol{\mid} r, s)$
	\begin{flalign}
	\label{wabcdp}
	\begin{aligned}  
	W_z (\textbf{A}, \textbf{B}; \textbf{C}, \textbf{D} \boldsymbol{\mid} r, s) & = (-1)^v z^{d - b} r^{2c - 2a} s^{2d} q^{\varphi (\textbf{D} - \textbf{V}, \textbf{C}) + \varphi (\textbf{V}, \textbf{A}) - av + cv} \displaystyle\frac{(q^{1 - v} r^{-2} z; q)_v}{(q^{- v} s^2 r^{-2} z; q)_v} \displaystyle\frac{(r^2; q)_d}{(r^2; q)_b}  \\
	& \quad \times  \displaystyle\sum_{p = 0}^{\min \{ b - v, c - v \}}  \displaystyle\frac{(q^{- v} s^2 r^{-2} z; q)_{c - p} (q^v r^2 z^{-1}; q)_p (z; q)_{b - p - v}}{(s^2 z; q)_{c + d - p - v}} (q^{-v} r^{-2} z)^p \\
	& \qquad \qquad \qquad \times  \displaystyle\sum_{\textbf{P}} q^{\varphi (\textbf{B} - \textbf{D} - \textbf{P}, \textbf{P})},
	\end{aligned}
	\end{flalign}
	
	\noindent where the last sum sum is over all $n$-tuples $\textbf{P} = (P_1, P_2, \ldots , P_n) \in \{ 0, 1 \}^n$ such that $|\textbf{P}| = p$ and $P_i \le \min \{ B_i - V_i, C_i - V_i \}$ for each $i \in [1, n]$. Otherwise, set $W_z (\textbf{A}, \textbf{B}; \textbf{C}, \textbf{D} \boldsymbol{\mid} r, s) = 0$. 
	
\end{definition}

Let us provide several examples of these $W_z$ weights that will be useful for us later.

\begin{example} 
	
	\label{wabcd01n} 
	
	Recall the notation $\textbf{e}_0 = (0, 0, \ldots , 0) \in \mathbb{Z}^n$ and $\textbf{e}_{[1, n]} = (1, 1, \ldots , 1) \in \mathbb{Z}^n$.\index{E@$\textbf{e}_0$}\index{E@$\textbf{e}_{[1, n]}$} For any $r, s, z \in \mathbb{C}$, we have that 
	\begin{flalign}
	\label{wabcd01nequation} 
	\begin{aligned} 
	W_z (\textbf{e}_0, & \textbf{e}_0; \textbf{e}_0, \textbf{e}_0 \boldsymbol{\mid} r, s) = 1; \qquad W_z ( \textbf{e}_0, \textbf{e}_{[1, n]}; \textbf{e}_0, \textbf{e}_{[1, n]} \boldsymbol{\mid} r, s) = \displaystyle\frac{s^{2n} (z; q)_n}{(s^2 z; q)_n}; \\
	& W_z ( \textbf{e}_{[1, n]}, \textbf{e}_{[1, n]}; \textbf{e}_{[1, n]}, \textbf{e}_{[1, n]} \boldsymbol{\mid} r, s) = (s^2 r^{-2} z)^n \displaystyle\frac{(r^2 z^{-1}; q)_n}{(s^2 z; q)_n}.
	\end{aligned} 
	\end{flalign} 
	
	\noindent Indeed, each of the three weights in \eqref{wabcd01nequation} are of the form $W_z (\textbf{A}, \textbf{B}; \textbf{C}, \textbf{D} \boldsymbol{\mid} r, s)$, with $\textbf{C} = \textbf{V}$ (where, as in \Cref{wabcdrsxy}, $\textbf{V}$ denotes the $n$-tuple obtained by taking the entrywise minimum of $\textbf{A}$, $\textbf{B}$, $\textbf{C}$, and $\textbf{D}$). Thus, the sum over $\textbf{P}$ on the right side of \eqref{wabcdp} is supported on the single term $\textbf{P} = \textbf{e}_0$, from which one quickly deduces \eqref{wabcd01nequation}. Through similar reasoning, generalizing the first two statements of \eqref{wabcd01nequation}, we have for any $\textbf{B} \in \{ 0, 1 \}^n$ with $|\textbf{B}| = b$ that
	\begin{flalign}
	\label{wze0b}
	W_z (\textbf{e}_0, \textbf{B}; \textbf{e}_0, \textbf{B} \boldsymbol{\mid} r, s) = \displaystyle\frac{s^{2b} (z; q)_b}{(s^2 z; q)_b}.
	\end{flalign}
\end{example}

For any $r, s, z \in \mathbb{C}$ and $n$-tuples $\textbf{A}, \textbf{B}, \textbf{C}, \textbf{D} \in \mathbb{Z}_{\ge 0}^n$, we define a normalization of the weight $W_z (\textbf{A}, \textbf{B}; \textbf{C}, \textbf{D} \boldsymbol{\mid} r, s)$\index{W@$\widehat{W}_z (\textbf{A}, \textbf{B}; \textbf{C}, \textbf{D} \boldsymbol{\mid} r, s)$; normalized fused weight} from \Cref{wabcdrsxy}, given by 
\begin{flalign}
\label{wabcd2} 
\widehat{W}_z (\textbf{A}, \textbf{B}; \textbf{C}, \textbf{D} \boldsymbol{\mid} r, s) = \displaystyle\frac{(s^2 z; q)_n}{s^{2n} (z; q)_n} W_z (\textbf{A}, \textbf{B}; \textbf{C}, \textbf{D} \boldsymbol{\mid} r, s).
\end{flalign}

\noindent Observe in particular by the second statement of \eqref{wabcd01nequation} that 
\begin{flalign}
\label{wz1} 
\widehat{W}_z \big( \textbf{e}_0, \textbf{e}_{[1, n]}; \textbf{e}_0, \textbf{e}_{[1, n]} \boldsymbol{\mid} r, s \big) = 1.
\end{flalign}

 Similarly to in \Cref{Weightszq}, we interpret $W_z (\textbf{A}, \textbf{B}; \textbf{C}, \textbf{D} \boldsymbol{\mid} r, s)$ (or $\widehat{W}_z (\textbf{A}, \textbf{B}; \textbf{C}, \textbf{D} \boldsymbol{\mid} r, s)$, depending on the context) as a vertex weight in the following way. As there, a vertex $v$ is the intersection between two directed transverse curves $\ell_1$ (typically oriented east) and $\ell_2$ (typically oriented north), but now associated with each curve is an ordered pair of rapidity parameters. Let the one associated with $\ell_1$ be $(x; r)$ and that associated with $\ell_2$ be $(y; s)$.\index{X@$(x; r), (y; s)$; rapidity parameters} The spectral parameter associated with the vertex $v$ is then given by the ratio\footnote{Observe that this is slightly different from in \Cref{Weightszq}, where the spectral parameter was instead given by $\frac{y}{x}$.} $\frac{x}{y}$. 
 
 Each of the four edges (segments of the curves $\ell_1$ and $\ell_2$) adjacent to $v$ may accommodate arrows of colors in $\{ 1, 2, \ldots , n \}$. Let $A_i$, $B_i$, $C_i$, and $D_i$ denote the numbers of arrows of any color $i \in [1, n]$ that vertically enter, horizontally enter, vertically exit, and horizontally exit $v$, respectively. As in \Cref{wabcdrsxy}, define $\textbf{X} = (X_1, X_2, \ldots , X_n)$ for each index $X \in \{A, B, C, D \}$. We will assume in what follows that $\textbf{A}, \textbf{B}, \textbf{C}, \textbf{D} \in \{ 0, 1 \}^n$, so that each edge accommodates at most one arrow of any given color (that is, all colors are fermionic). We will additionally impose $\textbf{A} + \textbf{B} = \textbf{C} + \textbf{D}$, which is a form of arrow conservation. The quadruple $(\textbf{A}, \textbf{B}; \textbf{C}, \textbf{D})$ is the \emph{(fused) arrow configuration} at $v$, and we interpret $W_{x / y} (\textbf{A}, \textbf{B}; \textbf{C}, \textbf{D} \boldsymbol{\mid} r, s)$ (or $\widehat{W}_{x / y} (\textbf{A}, \textbf{B}; \textbf{C}, \textbf{D} \boldsymbol{\mid} r, s)$) as the weight of this vertex $v$. We diagrammatically depict such vertices as on the left side of \Cref{equationpathsw}.

 \begin{figure}

 	\begin{center}

 		\begin{tikzpicture}[
 		>=stealth,
 		auto,
 		style={
 			scale = 1
 		}
 		]

 		\draw[->, ultra thick, gray] (-6, 0) -- (-4.5, 0) -- (-3, 0);
 		\draw[->, ultra thick, gray] (-4.5, -1.5) -- (-4.5, 0) -- (-4.5, 1.5);
 		
 		\draw[] (-6.25, 0) circle[radius = 0] node[scale = .9, left]{$(x; r)$};
 		\draw[] (-4.5, -1.75) circle[radius = 0] node[scale = .9, below]{$(y; s)$};
 		
 		\draw[] (-4.5, -.75) circle[radius = 0] node[scale = .7, right]{$\textbf{A}$};
 		\draw[] (-5.25, 0) circle[radius = 0] node[scale = .7, above]{$\textbf{B}$};
 		\draw[] (-4.5, .75) circle[radius = 0] node[scale = .7, right]{$\textbf{C}$};
 		\draw[] (-3.75, 0) circle[radius = 0] node[scale = .7, above]{$\textbf{D}$};

 		\draw[->, ultra thick, gray] (-.87, -.5) -- (0, 0);
 		\draw[->, ultra thick, gray] (-.87, .5) -- (0, 0);
 		\draw[->, ultra thick, gray] (.87, -1.5) -- (.87, -.5); 
 		
 		\draw[] (-.87, .5) circle[radius = 0]  node[left, scale = .9]{$(x; r)$};
 		\draw[] (-.87, -.5) circle[radius = 0]  node[left, scale = .9]{$(y; s)$};
 		\draw[] (.87, -1.5) circle[radius = 0]  node[below, scale = .9]{$(z; t)$};
 		
 		\draw[gray, ultra thick, dashed] (0, 0) -- (.87, -.5); 
 		\draw[gray, ultra thick, dashed] (0, 0) -- (.87, .5); 
 		\draw[gray, ultra thick, dashed] (.87, -.5) -- (.87, .5); 
 		
 		\draw[->, ultra thick, gray] (.87, .5) -- (1.87, .5); 
 		\draw[->, ultra thick, gray] (.87, -.5) -- (1.87, -.5); 
 		\draw[->, ultra thick, gray] (.87, .5) -- (.87, 1.5);

 		\draw[] (3.87, .5) circle[radius = 0]  node[left, scale = .9]{$(x; r)$};
 		\draw[] (3.87, -.5) circle[radius = 0]  node[left, scale = .9]{$(y; s)$};
 		\draw[] (4.87, -1.5) circle[radius = 0]  node[below, scale = .9]{$(z; t)$};
 		
 		\draw[->, ultra thick, gray] (4.87, -1.5) -- (4.87, -.5); 
 		\draw[->, ultra thick, gray] (3.87, .5) -- (4.87, .5); 
 		\draw[->, ultra thick, gray] (3.87, -.5) -- (4.87, -.5); 
 		\draw[->, ultra thick, gray] (4.87, .5) -- (4.87, 1.5); 
 		\draw[->, ultra thick, gray] (5.74, 0) -- (6.61, -.5); 
 		\draw[->, ultra thick, gray] (5.74, 0) -- (6.61, .5); 
 		
 		\draw[gray, ultra thick, dashed] (4.87, -.5) -- (4.87, .5);
 		\draw[-, gray, ultra thick, dashed] (4.87, -.5) -- (5.74, 0); 
 		\draw[gray, ultra thick, dashed] (4.87, .5) -- (5.74, 0);

 		\filldraw[fill=white, draw=black] (2.5, 0) circle [radius=0] node[scale = 2]{$=$};
 		
 		\filldraw[fill=white, draw=black] (-.44, -.275) circle [radius=0] node[below, scale = .7]{$\textbf{I}_1$};
 		\filldraw[fill=white, draw=black] (.44, .275) circle [radius=0] node[above, scale = .7]{$\textbf{I}_2$};
 		\filldraw[fill=white, draw=black] (1.45, .5) circle [radius=0] node[above, scale = .7]{$\textbf{I}_3$};
 		\filldraw[fill=white, draw=black] (-.44, .275) circle [radius=0] node[above, scale = .7]{$\textbf{J}_1$};
 		\filldraw[fill=white, draw=black] (.44, -.275) circle [radius=0] node[below, scale = .7]{$\textbf{J}_2$};
 		\filldraw[fill=white, draw=black] (1.45, -.5) circle [radius=0] node[above, scale = .7]{$\textbf{J}_3$};
 		\filldraw[fill=white, draw=black] (.87, -1) circle [radius=0] node[right, scale = .7]{$\textbf{K}_1$};
 		\filldraw[fill=white, draw=black] (.87, 0) circle [radius=0] node[right, scale = .7]{$\textbf{K}_2$};
 		\filldraw[fill=white, draw=black] (.87, 1) circle [radius=0] node[right, scale = .7]{$\textbf{K}_3$};
 		
 		\filldraw[fill=white, draw=black] (4.32, -.5) circle [radius=0] node[above, scale =.8]{$\textbf{I}_1$};
 		\filldraw[fill=white, draw=black] (5.35, -.26) circle [radius=0] node[below, scale = .7]{$\textbf{I}_2$};
 		\filldraw[fill=white, draw=black] (6.05, .26) circle [radius=0] node[above, scale = .7]{$\textbf{I}_3$};
 		\filldraw[fill=white, draw=black] (4.32, .5) circle [radius=0] node[above, scale = .7]{$\textbf{J}_1$};
 		\filldraw[fill=white, draw=black] (5.35, .26) circle [radius=0] node[above, scale = .7]{$\textbf{J}_2$};
 		\filldraw[fill=white, draw=black] (6.05, -.26) circle [radius=0] node[below, scale = .7]{$\textbf{J}_3$};
 		\filldraw[fill=white, draw=black] (4.87, -1) circle [radius=0] node[left, scale = .7]{$\textbf{K}_1$};
 		\filldraw[fill=white, draw=black] (4.87, 0) circle [radius=0] node[left, scale = .7]{$\textbf{K}_2$};
 		\filldraw[fill=white, draw=black] (4.87, 1) circle [radius=0] node[left, scale = .7]{$\textbf{K}_3$};

 		\end{tikzpicture}
 		
 	\end{center}

 	\caption{\label{equationpathsw} A fused vertex is shown to the left. The diagrammatic interpretation of the Yang--Baxter equation for the $W_z$ weights is depicted to the right.  }
 \end{figure}

  Observe that, unlike in \Cref{Weights1}, in \Cref{wabcdrsxy} we do not impose either $|\textbf{A}| = |\textbf{C}|$ or $|\textbf{B}| = |\textbf{D}|$, and we do not count arrows of the originally bosonic color $0$ (which we might view as ``empty''). Still, as the following proposition indicates, if $r$ and $s$ are negative half-integer powers of $q$, then the $W_z$ weights can be matched with the $\mathcal{R}_z$ weights from \Cref{rijkh} under a change of variables. 

\begin{prop} 
	
	\label{wabcdsxz} 
	
	Fix integers $n \ge 1$ and $L, M \ge 1$. Let $\textbf{\emph{A}}, \textbf{\emph{B}}, \textbf{\emph{C}}, \textbf{\emph{D}} \in \mathbb{Z}_{\ge 0}^{n + 1}$ and, for each index $X \in \{ A, B, C, D \}$, set $\textbf{\emph{X}} = (X_0, X_1, \ldots , X_n)$ and $\check{\textbf{\emph{X}}} = (X_1, X_2, \ldots , X_n)$. If $\check{\textbf{\emph{A}}}, \check{\textbf{\emph{B}}} \in \{ 0, 1 \}^n$, $|\textbf{\emph{A}}| = M = |\textbf{\emph{C}}|$, and $|\textbf{\emph{B}}| = L = |\textbf{\emph{D}}|$, then
	\begin{flalign*}
	\mathcal{R}_{x, y}^{(1; n)} (\textbf{\emph{A}}, \textbf{\emph{B}}; \textbf{\emph{C}}, \textbf{\emph{D}}) = W_{x / y} ( \check{\textbf{\emph{A}}}, \check{\textbf{\emph{B}}}, \check{\textbf{\emph{C}}}, \check{\textbf{\emph{D}}} \boldsymbol{\mid} q^{-L / 2}, q^{-M / 2}),
	\end{flalign*} 

	\noindent where the $W_{x / y}$ weights are given by \Cref{wabcdrsxy}.
	
\end{prop} 

\begin{proof}
	
	If either $\textbf{A} + \textbf{B} \ne \textbf{C} + \textbf{D}$, $\check{\textbf{C}} \notin \{ 0, 1 \}^n$, or $\check{\textbf{D}} \notin \{ 0, 1 \}^n$, then \Cref{rxyml} implies that $\mathcal{R}_{x, y} (\textbf{A}, \textbf{B}; \textbf{C}, \textbf{D}) = 0$, and so \Cref{wabcdrsxy} implies $\mathcal{R}_{x, y} (\textbf{A}, \textbf{B}; \textbf{C}, \textbf{D}) = W_{x / y} (\check{\textbf{A}}, \check{\textbf{B}}; \check{\textbf{C}}, \check{\textbf{D}})$. Hence, let us assume in what follows that $\textbf{A} + \textbf{B} = \textbf{C} + \textbf{D}$ and that $\check{\textbf{C}}, \check{\textbf{D}} \in \{ 0, 1 \}^n$. Define $\textbf{V} = (V_0, V_1, \ldots , V_n) \in \{ 0, 1 \}^{n + 1}$ by setting $V_0 = 0$ and $V_j = \min \{ A_j, B_j, C_j, D_j \}$ for each $j \in [1, n]$; we also let $\check{\textbf{V}} = (V_1, V_2, \ldots , V_n) \in \{ 0, 1 \}^n$ and, for each index $X \in \{ A, B, C, D, V \}$, we set $|\check{\textbf{X}}| = x$. 

	Then, recalling the definition of $\omega_{x, y}$ from \eqref{omegaxy} and setting $z = \frac{x}{y}$, we find that
	\begin{flalign}
	\label{omegaxyr}
	\begin{aligned}
	& \omega_{x, y}^{(L - v; M)} (\textbf{A}, \textbf{B} - \textbf{V}, \textbf{C}, \textbf{D} - \textbf{V}) \\
	& \quad = z^{d - b} q^{a (L - v) - (d - v) M} \displaystyle\sum_{p = 0}^{\min \{ b - v, c - v \}} \displaystyle\frac{(q^{L - M - v} z; q)_{c - p} (q^{v - L}; q)_{d - v}}{(q^{-M} z; q)_{c + d - v - p}} \displaystyle\frac{(q^{v - L} z^{-1}; q)_p (z; q)_{b - v - p}}{(q^{v - L}; q)_{b - v}} \\
	& \qquad \qquad \qquad \qquad \qquad \qquad \qquad \qquad \times q^{(L - v) (p - c)} z^p \displaystyle\sum_{\check{\textbf{P}}} q^{\varphi (\check{\textbf{D}} - \check{\textbf{V}}, \check{\textbf{C}} - \check{\textbf{P}}) + \varphi (\check{\textbf{B}} - \check{\textbf{V}} - \check{\textbf{P}}, \check{\textbf{P}})},
	\end{aligned}
	\end{flalign}

\noindent where $\check{\textbf{P}} = (P_1, P_2, \ldots , P_n) \in \{ 0, 1 \}^n$ is summed over all $n$-tuples such that $|\check{\textbf{P}}| = p$ and $P_i \le \min \{ B_i - V_i, C_i \} = \min \{ B_i - V_i, C_i - V_i \}$ for each $i \in [1, n]$ (and the last equality holds since $V_i \ne 0$ only if $A_i = B_i = C_i = D_i = 1$). Here, we have used the fact that all factors of the form $(q; q)_{B_i - V_i} (q; q)_{P_i}^{-1} (q; q)_{B_i - V_i - P_i}^{-1}$ and $(q; q)_{C_i + D_i - V_i - P_i} (q; q)_{D_i - V_i}^{-1} (q; q)_{C_i - P_i}^{-1}$ are equal to $1$. Indeed, these follow from the facts that $B_i - V_i, C_i + D_i - V_i \in \{ 0, 1 \}$; the first holds since $A_i, B_i, C_i, D_i \in \{ 0, 1 \}$ and the latter holds since $C_i + D_i \ge 2$ implies $A_i = B_i = C_i = D_i = V_i = 1$ (as $A_i + B_i = C_i + D_i$ and $A_i, B_i, C_i, D_i \in \{ 0, 1 \}$).  

Inserting \eqref{omegaxyr} into \Cref{rxyml} yields
\begin{flalign*}
\mathcal{R}_{x, y} (\textbf{A}, \textbf{B}; \textbf{C}, \textbf{D}) & = z^{d - b} q^{a (L - v) - (d - v) M} (-1)^v q^{\sum_{h: V_h = 1} A_{[h + 1, n]} - Mv} \displaystyle\frac{(q^{L - v + 1} z; q)_v}{(q^{L - M - v} z; q)_v} \\
& \qquad \times  \displaystyle\sum_{p = 0}^{\min \{ b - v, c - v \}}  \displaystyle\frac{(q^{L - M - v} z; q)_{c - p} (q^{v - L}; q)_{d - v}}{(q^{-M} z; q)_{c + d - v - p}} \displaystyle\frac{(q^{v - L} z^{-1}; q)_p (z; q)_{b - v - p}}{(q^{v - L}; q)_{b - v}} \\
& \qquad \qquad \qquad \qquad \times q^{(L - v) (p - c)} z^p \displaystyle\sum_{\check{\textbf{P}}} q^{\varphi (\check{{\textbf{D}}} - \check{\textbf{V}}, \check{\textbf{C}}) + \varphi (\check{{\textbf{B}}} - \check{\textbf{D}} - \check{\textbf{P}}, \check{\textbf{P}})},
\end{flalign*}

\noindent where we have also used the bilinearity of the function $\varphi$ from \eqref{tufunction}. Further using the facts that $\sum_{h: V_h = 1} A_{[h + 1, n]} = \varphi (\check{\textbf{V}}, \check{\textbf{A}})$ and that $q^L = r^{-2}$ and $q^M = s^{-2}$, we deduce that 
\begin{flalign*}
 \mathcal{R}_{x, y} (\textbf{A}, \textbf{B}; \textbf{C}, \textbf{D}) & = z^{d - b} r^{-2a} s^{2d} (-1)^v q^{\varphi (\check{\textbf{V}}, \check{\textbf{A}}) - av} \displaystyle\frac{(q^{1 - v} r^{-2} z; q)_v}{(q^{- v} s^2 r^{-2} z; q)_v}  \\
& \qquad \times  \displaystyle\sum_{p = 0}^{\min \{ b - v, c \}}  \displaystyle\frac{(q^{- v} s^2 r^{-2} z; q)_{c - p} (q^v r^2; q)_{d - v}}{(s^2 z; q)_{c + d - v - p}} \displaystyle\frac{(q^v r^2 z^{-1}; q)_p (z; q)_{b - p - v}}{(q^v r^2; q)_{b - v}} \\ 
& \qquad \qquad \qquad \times r^{2c - 2p} q^{v (c - p)} z^p \displaystyle\sum_{\check{\textbf{P}}} q^{\varphi (\check{{\textbf{D}}} - \check{\textbf{V}}, \check{\textbf{C}}) + \varphi (\check{{\textbf{B}}} - \check{\textbf{D}} - \check{\textbf{P}}, \check{\textbf{P}})},
\end{flalign*}

\noindent which yields the proposition in view of \eqref{wabcdp}.
\end{proof}

The following proposition states that these $W_z$ weights satisfy the Yang--Baxter equation; we refer to the right side of \Cref{equationpathsw} for a depiction. 

\begin{prop}
	
	\label{wabcdproduct2} 
	
	Fix an integer $n \ge 1$ and $x, y, z, r, s, t \in \mathbb{C}$. For any $\textbf{\emph{I}}_1, \textbf{\emph{J}}_1, \textbf{\emph{K}}_1, \textbf{\emph{I}}_3, \textbf{\emph{J}}_3, \textbf{\emph{K}}_3 \in \{ 0, 1 \}^n$, we have that 
	\begin{flalign}
	\label{ijkw} 
	\begin{aligned}
	& \displaystyle\sum_{\textbf{\emph{I}}_2, \textbf{\emph{J}}_2, \textbf{\emph{K}}_2} W_{x / y} ( \textbf{\emph{I}}_1, \textbf{\emph{J}}_1; \textbf{\emph{I}}_2, \textbf{\emph{J}}_2 \boldsymbol{\mid} r, s) W_{x / z} ( \textbf{\emph{K}}_1, \textbf{\emph{J}}_2; \textbf{\emph{K}}_2, \textbf{\emph{J}}_3 \boldsymbol{\mid} r, t) W_{y / z} ( \textbf{\emph{K}}_2, \textbf{\emph{I}}_2; \textbf{\emph{K}}_3, \textbf{\emph{I}}_3 \boldsymbol{\mid} s, t) \\
	& \qquad = \displaystyle\sum_{\textbf{\emph{I}}_2, \textbf{\emph{J}}_2, \textbf{\emph{K}}_2} W_{y / z} ( \textbf{\emph{K}}_1, \textbf{\emph{I}}_1; \textbf{\emph{K}}_2, \textbf{\emph{I}}_2 \boldsymbol{\mid} s, t) W_{x / z} ( \textbf{\emph{K}}_2, \textbf{\emph{J}}_1; \textbf{\emph{K}}_3, \textbf{\emph{J}}_2 \boldsymbol{\mid} r, t) W_{x / y} ( \textbf{\emph{I}}_2, \textbf{\emph{J}}_2; \textbf{\emph{I}}_3, \textbf{\emph{J}}_3 \boldsymbol{\mid} r, s), 
	\end{aligned}  
	\end{flalign}
	
	\noindent where both sums are over all $\textbf{\emph{I}}_2, \textbf{\emph{J}}_2, \textbf{\emph{K}}_2 \in \{ 0, 1 \}^n$. 
\end{prop} 

\begin{proof}
	
	Observe by the explicit form \eqref{wabcdp} for $W$ that, for fixed $\textbf{I}_1, \textbf{J}_1, \textbf{K}_1, \textbf{I}_3, \textbf{J}_3, \textbf{K}_3 \in \{ 0, 1 \}^n$ and nonzero $x, y, z \in \mathbb{C}$, both sides of \eqref{ijkw} are rational functions in $r$, $s$, and $t$. Thus, it suffices to verify \eqref{ijkw} assuming that there exist integers $L, M, N \ge |\textbf{I}|_3 + |\textbf{J}_3| + |\textbf{K}_3|$ such that $r = q^{-L / 2}$, $s = q^{-M / 2}$, and $t = q^{-N / 2}$. We may further assume that $\textbf{I}_1 + \textbf{J}_1 + \textbf{K}_1 = \textbf{I}_3 + \textbf{J}_3 + \textbf{K}_3$, for otherwise both sides of \eqref{ijkw} are equal to $0$, by arrow conservation. 
	
	Define $\textbf{I}_1', \textbf{J}_1', \textbf{K}_1', \textbf{I}_3', \textbf{J}_3', \textbf{K}_3' \in \mathbb{Z}_{\ge 0}^{n + 1}$ by setting $\textbf{I}_h' = \big( L - |\textbf{I}_h|, \textbf{I}_h \big)$, $\textbf{J}_h' = \big( M - |\textbf{J}_h|, \textbf{J}_h \big)$, and $\textbf{K}_h' = \big( N - |\textbf{K}_h|, \textbf{K}_h \big)$, for each $h \in \{ 1, 2, 3 \}$. Then, \Cref{wabcdproduct} implies that
	\begin{flalign}
	\label{ijkw2}
	\begin{aligned}
	\displaystyle\sum_{\textbf{I}_2', \textbf{J}_2', \textbf{K}_2'} & \mathcal{R}_{x, y} ( \textbf{I}_1', \textbf{J}_1'; \textbf{I}_2', \textbf{J}_2') \mathcal{R}_{x, z} ( \textbf{K}_1', \textbf{J}_2'; \textbf{K}_2', \textbf{J}_3' ) \mathcal{R}_{y, z} ( \textbf{K}_2', \textbf{I}_2'; \textbf{K}_3', \textbf{I}_3') \\
	& = \displaystyle\sum_{\textbf{I}_2', \textbf{J}_2', \textbf{K}_2'} \mathcal{R}_{y, z} ( \textbf{K}_1', \textbf{I}_1'; \textbf{K}_2', \textbf{I}_2') \mathcal{R}_{x, z} ( \textbf{K}_2', \textbf{J}_1'; \textbf{K}_3', \textbf{J}_2') \mathcal{R}_{x, y} ( \textbf{I}_2, \textbf{J}_2; \textbf{I}_3, \textbf{J}_3), 
	\end{aligned}  
	\end{flalign}
	
	\noindent where both sums are over all $\textbf{I}_2', \textbf{J}_2', \textbf{K}_2' \in \mathbb{Z}_{\ge 0}^{n + 1}$ such that $|\textbf{I}_2'| = L$, $|\textbf{J}_2'| = M$, and $|\textbf{K}_2'| = N$. Since $\textbf{I}_1, \textbf{J}_1, \textbf{K}_1 \in \{ 0, 1 \}^n$, the first statement of \Cref{rxyml} implies that we may instead sum both sides of \eqref{ijkw2} over those $\textbf{I}_2', \textbf{J}_2', \textbf{K}_2'$ satisfying the above conditions and also that $\check{\textbf{I}}_2', \check{\textbf{J}}_2', \check{\textbf{K}}_2' \in \{ 0, 1 \}^n$, where $\check{\textbf{X}} = (X_1, X_2, \ldots, X_n)$ for any $\textbf{X} = (X_0, X_1, \ldots , X_n) \in \mathbb{Z}_{\ge 0}^{n + 1}$. Then, since $\check{\textbf{I}}_h' = \textbf{I}_h$, $\check{\textbf{J}}_h' = \textbf{J}_h$, and $\check{\textbf{K}}_h' = \textbf{K}_h$ for each index $h \in \{ 1, 3 \}$, \Cref{wabcdsxz} and \eqref{ijkw2} together imply \eqref{ijkw}.
\end{proof}

\section{Fused Color Merging} 

\label{ColorsFused}

In this section we establish a generalization to the fused setting of the color merging procedure, given by \Cref{zsumef} in the case of fundamental weights. Throughout, we recall from \Cref{DomainBoundary} the notions of east-south paths and domains. We will once again consider vertex models on these domains, but with the difference that they will now be fused, in that paths may share edges. 

To describe this in more detail, we first introduce the fused versions of the notation from \Cref{DomainBoundary}. To that end, let $\mathcal{D} = \mathcal{D} (\textbf{p}, \textbf{p}') \subset \mathbb{Z}^2$ denote an east-south domain, whose boundary paths $\textbf{p}$ and $\textbf{p}'$ are of length $k$. A \emph{(fused) path ensemble} on $\mathcal{D}$ is a consistent assignment of a (fused) arrow configuration\footnote{Throughout the remainder of this text, nearly all arrow configurations and path ensembles will be fused, so we will typically not mention this explicitly in what follows.} $\big( \textbf{A} (v), \textbf{B} (v); \textbf{C} (v), \textbf{D} (v) \big)$ to each vertex $v \in \mathcal{D}$; the consistency here means $\textbf{B} (u) = \textbf{D} (v)$ if $u - v = (1, 0)$ and $\textbf{A} (u) = \textbf{C} (v)$ if $u - v = (0, 1)$. Arrow conservation $\textbf{A} + \textbf{B} = \textbf{C} + \textbf{D}$ implies that any path ensemble may be viewed as a collection of colored paths, which may cross or share edges; the condition $\textbf{A}, \textbf{B}, \textbf{C}, \textbf{D} \in \{ 0, 1 \}^n$ indicates that two paths of the same color cannot share edges (but paths of different colors can).

\begin{figure}
	
	\begin{center}		
		
		\begin{tikzpicture}[
		>=stealth,
		auto,
		style={
			scale = .75
		}
		]
		
		\draw[dashed] (1, 8) -- (1, 7) -- (3, 7);
		\draw[dashed] (2, 8) -- (2, 6) -- (3 ,6);
		\draw[dashed] (3, 4) -- (4, 4) -- (4, 3);
		\draw[dashed] (4, 5) -- (4, 4) -- (5, 4);
		
		\draw[] (-.8, 8) circle[radius = 0] node[left, black, scale = .7]{$\textbf{E}_1$};
		\draw[] (-.8, 7) circle[radius = 0] node[left, black, scale = .7]{$\textbf{E}_2$};
		\draw[] (0, 6.5) circle[radius = 0] node[left, black, scale = .7]{$\textbf{E}_3$};
		\draw[] (.5, 6) circle[radius = 0] node[above, black, scale = .7]{$\textbf{E}_4$};
		\draw[] (1, 5.5) circle[radius = 0] node[left, black, scale = .7]{$\textbf{E}_5$};
		\draw[] (1.5, 5) circle[radius = 0] node[above, black, scale = .7]{$\textbf{E}_6$};
		\draw[] (2, 4.5) circle[radius = 0] node[left, black, scale = .7]{$\textbf{E}_7$};
		\draw[] (2.5, 4.05) circle[radius = 0] node[above, black, scale = .7]{$\textbf{E}_8$};
		\draw[] (2.5, 3) circle[radius = 0] node[above, black, scale = .7]{$\textbf{E}_9$};
		\draw[] (3, 2.2) circle[radius = 0] node[below, black, scale = .7]{$\textbf{E}_{10}$};
		\draw[] (4, 2.2) circle[radius = 0] node[below, black, scale = .7]{$\textbf{E}_{11}$};
		\draw[] (5, 2.2) circle[radius = 0] node[below, black, scale = .7]{$\textbf{E}_{12}$};
		
		\draw[] (0, 8.8) circle[radius = 0] node[above, black, scale = .7]{$\textbf{F}_1$};
		\draw[] (1, 8.8) circle[radius = 0] node[above, black, scale = .7]{$\textbf{F}_2$};
		\draw[] (2, 8.8) circle[radius = 0] node[above, black, scale = .7]{$\textbf{F}_3$};
		\draw[] (3, 8.8) circle[radius = 0] node[above, black, scale = .7]{$\textbf{F}_4$};
		\draw[] (3.4, 8.05) circle[radius = 0] node[above, black, scale = .7]{$\textbf{F}_5$};
		\draw[] (3.4, 7) circle[radius = 0] node[above, black, scale = .7]{$\textbf{F}_6$};
		\draw[] (3.4, 6) circle[radius = 0] node[above, black, scale = .7]{$\textbf{F}_7$};
		\draw[] (4, 5.4) circle[radius = 0] node[left, black, scale = .7]{$\textbf{F}_8$};
		\draw[] (5, 5.4) circle[radius = 0] node[left, black, scale = .7]{$\textbf{F}_9$};
		\draw[] (5.8, 5) circle[radius = 0] node[right, black, scale = .7]{$\textbf{F}_{10}$};
		\draw[] (5.8, 4) circle[radius = 0] node[right, black, scale = .7]{$\textbf{F}_{11}$};
		\draw[] (5.8, 3) circle[radius = 0] node[right, black, scale = .7]{$\textbf{F}_{12}$};
		
		\draw[thick, blue, ->] (-.8, 8.05) -- (0, 8.05);
		\draw[thick, red, ->] (-.8, 7.95) -- (0, 7.95);
		
		\draw[dotted] (-.8, 7) -- (0, 7);
		
		\draw[thick, red, ->] (.2, 5.95) -- (1, 5.95);
		\draw[thick, green, ->] (.2, 6.05) -- (1, 6.05);
		
		\draw[thick, orange, ->] (1.2, 5.05) -- (2, 5.05);
		\draw[thick, green, ->] (1.2, 4.95) -- (2, 4.95);
		
		\draw[thick, orange, ->] (2.2, 4.1) -- (3, 4.1);
		\draw[thick, green, ->] (2.2, 4) -- (3, 4);
		\draw[thick, blue, ->] (2.2, 3.9) -- (3, 3.9);
		
		\draw[thick, blue, ->] (2.2, 3.05) -- (3, 3.05);
		\draw[thick, red, ->] (2.2, 2.95) -- (3, 2.95);

		\draw[thick, blue, ->] (0, 6.2) -- (0, 7);
		\draw[dotted] (1, 5.2) -- (1, 6);
		\draw[thick, green, ->] (2, 4.2) -- (2, 5);
		
		\draw[thick, green, ->] (3.05, 2.2) -- (3.05, 3);
		\draw[thick, blue, ->] (2.95, 2.2) -- (2.95, 3);
		
		\draw[thick, orange, ->] (4, 2.2) -- (4, 3);
		
		\draw[dotted] (5, 2.2) -- (5, 3);
		
		\draw[thick, blue, ->] (3, 7.9) -- (3.8, 7.9);
		\draw[thick, green, ->] (3, 8) -- (3.8, 8);
		\draw[thick, orange, ->] (3, 8.1) -- (3.8, 8.1);

		\draw[thick, red, ->] (3, 7) -- (3.8, 7);
		\draw[thick, blue, ->] (3, 6) -- (3.8, 6);
		
		\draw[thick, green, ->] (5, 4.95) -- (5.8, 4.95);
		\draw[thick, blue, ->] (5, 5.05) -- (5.8, 5.05);
		
		\draw[thick, blue, ->] (5, 4.05) -- (5.8, 4.05);
		\draw[thick, red, ->] (5, 3.95) -- (5.8, 3.95);
		
		\draw[thick, orange, ->] (5, 3.05) -- (5.8, 3.05);
		\draw[thick, green, ->] (5, 2.95) -- (5.8, 2.95);
		
		\draw[dotted] (0, 8) -- (0, 8.8);
		
		\draw[thick, green, ->] (1.05, 8) -- (1.05, 8.8);
		\draw[thick, blue, ->] (.95, 8) -- (.95, 8.8);
		
		\draw[thick, red, ->] (2, 8) -- (2, 8.8);
		
		\draw[thick, orange, ->] (3, 8) -- (3, 8.8);
		
		\draw[thick, green, ->] (4, 5) -- (4, 5.8);
		
		\draw[dotted] (5, 5) -- (5, 5.8);
		
		\draw[ultra thick] (0, 8) -- (0, 7) -- (1, 7) -- (1, 6) -- (2, 6) -- (2, 5) -- (3, 5) -- (3, 3) -- (5, 3);
		\draw[ultra thick] (0, 8) -- (3, 8) -- (3, 7) -- (3, 5) -- (4, 5) -- (5, 5) -- (5, 3);

		\end{tikzpicture}
		
	\end{center}
	
	\caption{\label{domainpathsfused} Shown above is an east-south domain with boundary data. } 	
\end{figure} 

For any $(k + 1)$-tuples $\mathscr{E} = (\textbf{E}_1, \textbf{E}_2, \ldots , \textbf{E}_{k + 1})$ and $\mathscr{F} = (\textbf{F}_1, \textbf{F}_2, \ldots , \textbf{F}_{k + 1})$ of elements in $\{ 0, 1 \}^n$, a path ensemble has boundary data $(\mathscr{E}; \mathscr{F})$ if the following holds for each $i \in [1, k + 1]$ and $h \in [1, n]$. An arrow of color $h$ enters through the $i$-th incoming edge in $\mathcal{D}$ if and only if the $h$-th coordinate of $\textbf{E}_i$ is $1$, and one exits through its $i$-th outgoing edge if and only if the $h$-th coordinate of $\textbf{F}_i$ is $1$. We refer to $\mathscr{E}$ as entrance data on $\mathcal{D}$ and $\mathscr{F}$ as exit data. For example, if red, blue, green, and orange are colors $1$, $2$, $3$, and $4$, respectively, then $\textbf{E}_9 = (1, 1, 0, 0)$ and $\textbf{F}_5 = (0, 1, 1, 1)$ in \Cref{domainpathsfused}. 

The following definition provides notation for the partition functions of a vertex model with given boundary data, under the fused weights $W_z$ from \Cref{wabcdrsxy}. 

\begin{definition} 
	
	\label{zef1rsz}
	
	Let $\mathcal{D} \subset \mathbb{Z}^2$ denote an east-south domain. For any boundary data $(\mathscr{E}; \mathscr{F})$ on $\mathcal{D}$ and any sets of complex numbers $\textbf{z} = \big( z(v) \big)_{v \in \mathcal{D}}$, $\textbf{r} = \big( r(v) \big)_{v \in \mathcal{D}}$, and $\textbf{s} = \big( s(v) \big)_{v \in \mathcal{D}}$, define 
	\begin{flalign*}
	& W_{\mathcal{D}}^{(m; n)} (\mathscr{E}; \mathscr{F} \boldsymbol{\mid} \textbf{z} \boldsymbol{\mid} \textbf{r}, \textbf{s})  = \displaystyle\sum \displaystyle\prod_{v \in \mathcal{D}} W_{z(v)} \big( \textbf{A} (v), \textbf{B} (v); \textbf{C} (v), \textbf{D} (v) \boldsymbol{\mid} r(v), s(v) \big),
	\end{flalign*} \index{W@$W_{\mathcal{D}}^{(m; n)} (\mathscr{E}; \mathscr{F} \boldsymbol{\mid} \textbf{z} \boldsymbol{\mid} \textbf{r}, \textbf{s})$}
	
	\noindent where the sum is over all (fused) path ensembles on $\mathcal{D}$ with boundary data $(\mathscr{E}; \mathscr{F})$. We further set $W_{\mathcal{D}}^{(m; n)} (\mathscr{E}; \mathscr{F} \boldsymbol{\mid} \textbf{z} \boldsymbol{\mid} \textbf{r}, \textbf{s}) = 0$ if there exists some $\textbf{E}_i \in \mathscr{E}$ or $\textbf{F}_i \in \mathscr{F}$ not in $\{ 0, 1 \}^n$. 	
	
\end{definition}

To state the color merging result, we further require an inversion count, similar to \eqref{ijsum} in the fundamental case. For any $k$-tuple $\mathscr{I} = (\textbf{I}_1, \textbf{I}_2, \ldots , \textbf{I}_k)$ of elements in $\{ 0, 1 \}^n$, with $\textbf{I}_j = (I_{1, j}, I_{2, j}, \ldots , I_{n, j}) \in \{ 0, 1 \}^n$ for each $j \in [1, k]$, and any integer interval $J \subseteq [1, n]$, define
\begin{flalign}
\label{sumij} 
\inv (\mathscr{I}; J) = \displaystyle\sum_{1 \le h < j \le k} \displaystyle\sum_{1 \le a < b \le n} \textbf{1}_{a \in J} \textbf{1}_{b \in J} I_{b, h} I_{a, j}. 
\end{flalign}

\noindent Observe if $J = \{ 1, 2, \ldots , n \}$, then $\inv (\mathscr{I}; J) = \sum_{1 \le h < j \le k} \varphi (\textbf{I}_j, \textbf{I}_h)$, where $\varphi$ is given by \eqref{tufunction}.\index{I@$\inv$}

Next, recall the notion of interval partitions from \Cref{Colors}, and fix one $\mathbb{J} = (J_1, J_2, \ldots , J_{\ell})$ of $\{ 1, 2, \ldots , n \}$. Define the function $\vartheta_{\mathbb{J}}: \mathbb{Z}_{\ge 0}^n \rightarrow \mathbb{Z}_{\ge 0}^{\ell}$\index{0@$\vartheta_{\mathbb{J}}$} by setting 
\begin{flalign*} 
\vartheta_{\mathbb{J}} (\textbf{I}) = \Bigg( \displaystyle\sum_{h \in J_1} I_h, \displaystyle\sum_{h \in J_2} I_h, \ldots , \displaystyle\sum_{h \in J_{\ell}} I_h \Bigg), \qquad \text{for any $\textbf{I} = (I_1, I_2, \ldots , I_n) \in \{ 0, 1 \}^n$}. 
\end{flalign*}

\noindent Similar to the function $\theta_{\mathbb{J}}$ defined in \Cref{Colors}, $\vartheta_{\mathbb{J}}$ identifies all colors in any $J_i$ and renames them to color $i$. For any sequence $\mathscr{I} = (\textbf{I}_1, \textbf{I}_2, \ldots , \textbf{I}_k)$ of elements in $\mathbb{Z}_{\ge 0}^n$, let $\vartheta (\mathscr{I}) = \big( \vartheta (\textbf{I}_1), \vartheta (\textbf{I}_2), \ldots , \vartheta (\textbf{I}_k) \big)$. 

We now have the following theorem providing color merging for the fused $U_q \big( \widehat{\mathfrak{sl}} (1 | n) \big)$ model. Observe here that since the number of bosonic colors is $m = 1$, we do not symmetrize over the exit data for the model, as we did in \Cref{zsumef} for general $m$.  

\begin{thm} 
	
	\label{zsumfused} 
	
	Fix integers $n \ge n' \ge 1$, and $k \ge 1$; an east-south domain $\mathcal{D} = \mathcal{D} (\textbf{\emph{p}}, \textbf{\emph{p}}')$ with boundary paths $\textbf{\emph{p}}$ and $\textbf{\emph{p}}'$ of length $k$; and sets of complex numbers $\textbf{\emph{z}} = \big( z(v) \big)_{v \in \mathcal{D}}$, $\textbf{\emph{r}} = \big( r(v) \big)_{v \in \mathcal{D}}$, and $\textbf{\emph{s}} = \big( s(v) \big)_{v \in \mathcal{D}}$. Let $\mathbb{J} = (J_1, J_2, \ldots , J_{n'})$ denote an interval partition of $\{ 1, \ldots , n \}$, and let $\mathscr{E} = (\textbf{\emph{E}}_1, \textbf{\emph{E}}_2, \ldots , \textbf{\emph{E}}_{k + 1})$ and $\mathscr{F} = (\textbf{\emph{F}}_1, \textbf{\emph{F}}_2, \ldots , \textbf{\emph{F}}_{k + 1})$ denote sequences of elements in $\{ 0, 1 \}^n$ constituting entrance and exit data on $\mathcal{D}$, respectively; assume that $\vartheta_{\mathbb{J}} (\textbf{\emph{E}}_i) \in \{ 0, 1 \}^{n'}$ for each $i$. Then, 
	\begin{flalign}
	\label{sumwfused}
	\displaystyle\sum_{\breve{\mathscr{E}}} W_{\mathcal{D}}^{(1; n)} (\breve{\mathscr{E}}; \mathscr{F} \boldsymbol{\mid} \textbf{\emph{z}} \boldsymbol{\mid} \textbf{\emph{r}}, \textbf{\emph{s}}) \displaystyle\prod_{i = 1}^{n'} (-1)^{\inv (\breve{\mathscr{E}}; J_i) - \inv (\mathscr{F}; J_i)} = W_{\mathcal{D}}^{(1; n')} \big( \vartheta_{\mathbb{J}} (\mathscr{E}); \vartheta_{\mathbb{J}} (\mathscr{F}) \boldsymbol{\mid} \textbf{\emph{z}} \boldsymbol{\mid} \textbf{\emph{r}}, \textbf{\emph{s}}),
	\end{flalign}
	
	\noindent where the sum is over all sequences $\breve{\mathscr{E}} = (\breve{\textbf{\emph{E}}}_1, \breve{\textbf{\emph{E}}}_2, \ldots , \breve{\textbf{\emph{E}}}_{k + 1})$ of elements in $\{ 0, 1 \}^n$ such that $\vartheta_{\mathbb{J}} (\breve{\mathscr{E}}) = \vartheta_{\mathbb{J}} (\mathscr{E})$. 
	
\end{thm} 

\begin{proof}
	
	Throughout this proof, we will assume that $|\mathcal{D}| = 1$ since, given this, the proof in the general case follows from an induction very similar to the one implemented in \Cref{ProofSumZ} (alternatively, one can use the same method as applied in the proof of the Yang-Baxter equation in \Cref{wabcdproduct}). 
	
	Then letting $\mathcal{D} = \{ v \}$, the sets $\textbf{z}$, $\textbf{r}$, and $\textbf{s}$ each consist of one element, abbreviated by $z = z(v)$, $r = r(v)$, and $s = s(v)$, respectively. Hence, recalling $\mathscr{E} = (\textbf{E}_1, \textbf{E}_2)$ and $\mathscr{F} = (\textbf{F}_1, \textbf{F}_2)$, we have for any entrance data $\breve{\mathscr{E}} = (\breve{\textbf{E}}_1, \breve{\textbf{E}}_2)$ on $\mathcal{D}$ that $W_{\mathcal{D}}^{(1; n)} (\breve{\mathscr{E}}; \mathscr{F} \boldsymbol{\mid} \textbf{z} \boldsymbol{\mid} \textbf{r}, \textbf{s}) = W_z^{(1; n)} (\breve{\textbf{E}}_2, \breve{\textbf{E}}_1; \textbf{F}_1; \textbf{F}_2 \boldsymbol{\mid} r, s)$ and $W_{\mathcal{D}}^{(1; n')} \big( \vartheta_{\mathbb{J}} (\mathscr{E}); \vartheta_{\mathbb{J}} (\mathscr{F}) \boldsymbol{\mid} \textbf{z} \boldsymbol{\mid} \textbf{r}, \textbf{s}) = W_z^{(1; n')} \big( \vartheta_{\mathbb{J}} (\textbf{E}_2), \vartheta_{\mathbb{J}} (\textbf{E}_1);  \vartheta_{\mathbb{J}} (\textbf{F}_1),  \vartheta_{\mathbb{J}} (\textbf{F}_2) \boldsymbol{\mid} r, s \big)$. The explicit form \eqref{wabcdrsxy} for $W$ implies, for any fixed boundary data $(\mathscr{E}; \mathscr{F})$ on $v$, that both sides of \eqref{sumwfused} are rational in $r$ and $s$. Thus it suffices to verify \eqref{sumwfused} assuming $r = q^{-L / 2}$ and $s = q^{-M / 2}$, for some integers $L, M > n$.

	Next, we use \Cref{wabcdsxz} to express both sides of \eqref{sumwfused} in terms of the weights $\mathcal{R}_{z, 1}$ from \Cref{rijkh}. To do this, let $\mathbb{J}' = (J_0, J_1, \ldots , J_{n'})$ denote the interval partition of $\{ 0, 1, \ldots , n \}$ obtained by appending the singleton $J_0 = \{ 0 \}$ to $\mathbb{J}$. As in the proof of \Cref{wabcdproduct2}, define $\textbf{E}_1', \textbf{E}_2', \breve{\textbf{E}}_1', \breve{\textbf{E}}_2', \textbf{F}_1', \textbf{F}_2' \in \mathbb{Z}_{\ge 0}^{n + 1}$, with entries indexed by $[0, n + 1]$, by setting $\textbf{X}' = \big( L - |\textbf{X}|, \textbf{X} \big)$ for each index $X \in \{ E_1, \breve{E}_1, F_2 \}$ and $\textbf{X}' = \big( M - |\textbf{X}|, \textbf{X} \big)$ for each $X \in \{ E_2, \breve{E}_2, F_1 \}$. Then, \Cref{wabcdsxz} implies 
	\begin{flalign*}
	W_z^{(1; n)} (\breve{\textbf{E}}_2, \breve{\textbf{E}}_1; \textbf{F}_1, \textbf{F}_2 \boldsymbol{\mid} r, s) & = \mathcal{R}_{z, 1}^{(1; n)} (\breve{\textbf{E}}_2', \breve{\textbf{E}}_1'; \textbf{F}_1', \textbf{F}_2'); \\	
	W_z^{(m'; n')} \big( \vartheta_{\mathbb{J}} (\textbf{E}_2), \vartheta_{\mathbb{J}} (\textbf{E}_1);  \vartheta_{\mathbb{J}} (\textbf{F}_1),  \vartheta_{\mathbb{J}} (\textbf{F}_2) \boldsymbol{\mid} r, s \big) & = \mathcal{R}_{z, 1}^{(1; n')} \big( \vartheta_{\mathbb{J}'} (\textbf{E}_2'), \vartheta_{\mathbb{J}'} (\textbf{E}_1'); \vartheta_{\mathbb{J}'} (\textbf{F}_1'), \vartheta_{\mathbb{J}} (\textbf{F}_2') \big).
	\end{flalign*} 
	
	\noindent Moreover, recall the rectangular partition function $Z_{x, y}^{(m, n)} (\mathfrak{A}, \mathfrak{B}; \mathfrak{C}, \mathfrak{D})$ from \Cref{zxy1}, and observe by \Cref{zabcd2}, \Cref{zabcdcd}, and \Cref{rijkh} that 
	\begin{flalign}
	\label{rxy1z}
	\begin{aligned}
	\mathcal{R}_{x, y}^{(1; n)} (\textbf{A}, \textbf{B}; \textbf{C}, \textbf{D}) & \cdot \displaystyle\prod_{i = 0}^{m + n - 1} \displaystyle\frac{(q; q)_{C_i} (q; q)_{D_i}}{(q; q)_{A_i} (q; q)_{B_i}} \\
	& = q^{- \inv (\overleftarrow{\mathfrak{D}}) - \inv(\mathfrak{C})} \displaystyle\sum_{\mathfrak{A} \in \mathcal{M} (\textbf{A})} \displaystyle\sum_{\mathfrak{B} \in \mathcal{M} (\textbf{B})} q^{\inv (\mathfrak{A}) + \inv (\overleftarrow{\mathfrak{B}})} Z_{x, y}^{(1; n)} (\mathfrak{A}, \mathfrak{B}; \mathfrak{C}, \mathfrak{D}),
	\end{aligned} 
	\end{flalign}
	
	\noindent for any fixed complex numbers $x, y \in \mathbb{C}$; integer sequences $\textbf{A}, \textbf{B} \in \{ 0, 1\}^{n + 1}$ and $\textbf{C}, \textbf{D} \in \mathbb{Z}_{\ge 0}^{n + 1}$, with coordinates indexed by $[0, n]$; $\mathfrak{C} \in \mathcal{M} (\textbf{C})$; and $\mathfrak{D} \in \mathcal{M} (\textbf{D})$. Hence, to establish \eqref{sumwfused} we must show that\footnote{Fix an index $\textbf{X} \in \{ \textbf{E}_1', \textbf{E}_2', \textbf{F}_1', \textbf{F}_2' \}$; let $\breve{\textbf{X}} = (\breve{X}_0, \breve{X}_1, \ldots , \breve{X}_n)$; and let $\textbf{Y} = (Y_0, Y_1, \ldots , Y_{n'}) = \vartheta_{\mathbb{J}'} (\textbf{X})$. Then, it is quickly verified since $\vartheta_{\mathbb{J}'} \big( \breve{\textbf{X}} \big) = \vartheta_{\mathbb{J}'} (\textbf{X}) = \textbf{Y} \in \{ 0, 1 \}^{n'}$ that $\prod_{i = 0}^n (q; q)_{\breve{X}_i} = \prod_{i = 0}^{n'} (q; q)_{Y_i}$. In passing from \eqref{rxy1z} to \eqref{sumfusedz}, factors of this type were canceled on both sides.}
	\begin{flalign}
	\label{sumfusedz}
	\begin{aligned}
	 \displaystyle\sum_{\breve{\mathscr{E}}} \displaystyle\sum_{\breve{\mathfrak{E}}_1'} \displaystyle\sum_{\breve{\mathfrak{E}}_2'} q^{\inv(\breve{\mathfrak{E}}_2') + \inv (\overleftarrow{\breve{\mathfrak{E}}_1'})} & Z_{z, 1}^{(1; n)} (\breve{\mathfrak{E}}_2', \breve{\mathfrak{E}}_1'; \mathfrak{F}_1', \mathfrak{F}_2') \displaystyle\prod_{i = 1}^{n'} (-1)^{\inv (\breve{\mathscr{E}}; J_i) - \inv (\mathscr{F}; J_i)} \\
	 & = \displaystyle\sum_{\mathfrak{E}_1''} \displaystyle\sum_{\mathfrak{E}_2''} q^{\inv (\mathfrak{E}_2'') + \inv (\overleftarrow{\mathfrak{E}_1''})}  Z_{z, 1}^{(1; n')} \big( \mathfrak{E}_2'', \mathfrak{E}_1''; \mathfrak{F}_1'', \mathfrak{F}_2'' \big),
	 \end{aligned} 
	\end{flalign}
	
	\noindent where we sum $\breve{\mathscr{E}}$ as in \eqref{sumwfused}; $\breve{\mathfrak{E}}_1'$ over $\mathcal{M} (\breve{\textbf{E}}_1')$; $\breve{\mathfrak{E}}_2'$ over $\mathcal{M} (\breve{\textbf{E}}_2')$; $\mathfrak{E}_1''$ over $\mathcal{M} \big( \vartheta_{\mathbb{J}'} (\textbf{E}_1') \big)$; and $\mathfrak{E}_2''$ over $\mathcal{M} \big( \vartheta_{\mathbb{J}'} (\textbf{E}_2') \big)$. Here, $\mathfrak{F}_1' \in \mathcal{M} (\textbf{F}_1')$, $\mathfrak{F}_2' \in \mathcal{M} (\textbf{F}_2')$, $\mathfrak{F}_1'' \in \mathcal{M} \big( \vartheta_{\mathbb{J}'} (\textbf{F}_1') \big)$, and $\mathfrak{F}_2'' \in \mathcal{M} \big( \vartheta_{\mathbb{J}'} (\textbf{F}_2') \big)$ are fixed by the condition $\inv (\mathfrak{X}) = 0$ for each $\mathfrak{X} \in \big\{ \mathfrak{F}_1', \overleftarrow{\mathfrak{F}}_2', \mathfrak{F}_1'', \overleftarrow{\mathfrak{F}}_2'' \big\}$. Recalling the function $\theta_{\mathbb{J}'}$ from \Cref{Colors}, this ensures that $\theta_{\mathbb{J}'} (\mathfrak{F}_i') = \mathfrak{F}_i''$ for each $i \in \{ 1, 2 \}$, since $\mathfrak{X}' \in \mathcal{M} (\textbf{X}')$ implies $\theta_{\mathbb{J}'} (\mathfrak{X}') \in \mathcal{M} \big( \vartheta_{\mathbb{J}'} (\textbf{X}') \big)$ for any sequence $\mathfrak{X}'$ of indices in $[0, n]$ and any $\textbf{X} \in \mathbb{Z}_{\ge 0}^{n + 1}$. 
	
	To establish \eqref{sumfusedz}, we will apply \Cref{zsumef}, using the fact that the $Z_{z, 1}$ on both sides of \eqref{sumfusedz} are partition functions (as in \Cref{zef1}) for a fundamental $U_q \big( \widehat{\mathfrak{sl}} (1 | n) \big)$ six-vertex model on the rectangular domain $\mathcal{D}' = [1, M] \times [1, L] \subset \mathbb{Z}^2$. To implement this, let us fix $\mathfrak{E}_1'' \in \mathcal{M} \big( \vartheta_{\mathbb{J}'} (\textbf{E}_1') \big)$ and $\mathfrak{E}_2'' \in \mathcal{M} \big( \vartheta_{\mathbb{J}'} (\textbf{E}_2') \big)$. Then the $m = 1$ case of \Cref{zsumef} yields
	\begin{flalign}
	\label{sumzr}
	\displaystyle\sum_{\breve{\mathfrak{E}}'} Z_{z, 1}^{(1; n)} (\breve{\mathfrak{E}}_2', \breve{\mathfrak{E}_1}; \mathfrak{F}_1', \mathfrak{F}_2') \displaystyle\prod_{i = 1}^{n'} (-1)^{\inv (\breve{\mathfrak{E}}'; J_i) - \inv (\mathfrak{F}'; J_i)} = Z_{z, 1}^{(1; n')} \big( \mathfrak{E}_2'', \mathfrak{E}_1''; \mathfrak{F}_1'', \mathfrak{F}_2'' \big),
	\end{flalign}
	
	\noindent where $\breve{\mathfrak{E}}' = (\breve{\mathfrak{E}}_2', \breve{\mathfrak{E}}_1')$ is summed over all sequences of indices in $[0, n]$ such that $\theta_{\mathbb{J}'} (\breve{\mathfrak{E}}_i') = \mathfrak{E}_i''$ for each $i \in \{ 1, 2 \}$, and we have set $\mathfrak{F}' = (\mathfrak{F}_1', \mathfrak{F}_2')$ (which satisfies $\theta_{\mathbb{J}'} (\mathfrak{F}') = \mathfrak{F}''$, as mentioned above). 
	
	We next match the powers of $-1$ and $q$ in \eqref{sumfusedz} and \eqref{sumzr}, to which end we claim for any $i \in [1, n']$ that 
	\begin{flalign}
	\label{efe} 
	\inv (\breve{\mathfrak{E}}_2') = \inv (\mathfrak{E}_2''); \quad \inv \big( \overleftarrow{\breve{\mathfrak{E}}_1'} \big) = \inv \big( \overleftarrow{\mathfrak{E}_1''} \big); \quad \inv (\breve{\mathscr{E}}; J_i) = \inv (\breve{\mathfrak{E}}'; J_i); \quad \inv (\mathscr{F}; J_i) = \inv (\mathfrak{F}'; J_i).
	\end{flalign}
	
	\noindent Given \eqref{efe}, \eqref{sumfusedz} follows by multiplying \eqref{sumzr} by $q^{\inv (\breve{\mathfrak{E}}_2') + \inv (\overleftarrow{\breve{\mathfrak{E}}_1'})} = q^{\inv (\mathfrak{E}_2'') + \inv (\overleftarrow{\mathfrak{E}_1''})}$ and summing both sides over each $\mathfrak{E}_i'' \in \mathcal{M} \big( \vartheta_{\mathbb{J}'} (\textbf{E}_i') \big)$. Indeed, the left side of \eqref{sumfusedz} is summed over $\breve{\mathscr{E}} = (\breve{\textbf{E}}_2, \breve{\textbf{E}}_1)$ and $\breve{\mathfrak{E}}_i'$ with each $\vartheta_{\mathbb{J}} (\breve{\textbf{E}}_i) = \vartheta_{\mathbb{J}} (\textbf{E}_i)$ (equivalently, $\vartheta_{\mathbb{J}'} (\breve{\textbf{E}}_i') = \vartheta_{\mathbb{J}'} (\textbf{E}_i')$) and $\breve{\mathfrak{E}}_i' \in \mathcal{M} (\breve{\textbf{E}}_i')$. Since $\mathfrak{X}' \in \mathcal{M} (\textbf{X}')$ implies $\theta_{\mathbb{J}'} (\mathfrak{X}') \in \mathcal{M} \big( \vartheta_{\mathbb{J}'} (\textbf{X}') \big)$ for any $\mathfrak{X}'$ and $\textbf{X}'$, this is equivalent to first summing over each $\breve{\mathfrak{E}}_i'$ with $\theta_{\mathbb{J}'} (\breve{\mathfrak{E}}_i') = \mathfrak{E}_i''$ for fixed $\mathfrak{E}_i'' \in \mathcal{M} \big( \vartheta_{\mathbb{J}'} (\textbf{E}_i') \big)$ and then over each such $\mathfrak{E}_i''$. 
	
	Now let us establish \eqref{efe}. To that end, observe since each $\breve{\mathfrak{E}}_2' \in \mathcal{M} (\breve{\textbf{E}}_2')$ and $\vartheta (\breve{\textbf{E}}_2) = \vartheta (\textbf{E}_2) \in \{ 0, 1 \}^{n'}$ that $|\breve{\mathfrak{E}}_2' \cap J_i| \le 1$ for each $i \in [1, n']$. Thus, applying $\theta_{\mathbb{J}'}$ to $\breve{\mathfrak{E}}_2'$ does not change the relative ordering of its nonzero elements, which  yields $\inv (\breve{\mathfrak{E}}_2') = \inv \big( \theta_{\mathbb{J}'} (\breve{\mathfrak{E}}_2') \big) = \inv (\mathfrak{E}_2'')$. This verifies the first statement of \eqref{efe}; the proof of the second is entirely analogous and therefore omitted.
	
	To confirm the third, observe for any $i \in [1, n']$ that
	\begin{flalign*}
	\inv (\breve{\mathfrak{E}}'; J_i) & = \inv (\breve{\mathfrak{E}}_1'; J_i) + \inv (\breve{\mathfrak{E}}_2'; J_i) + \displaystyle\sum_{a \in \breve{\mathfrak{E}}_1'} \displaystyle\sum_{b \in \breve{\mathfrak{E}}_2'} \textbf{1}_{a > b} \textbf{1}_{a \in J_i} \textbf{1}_{b \in J_i} \\
	 & = \displaystyle\sum_{a \in \breve{\mathfrak{E}}_1'} \displaystyle\sum_{b \in \breve{\mathfrak{E}}_2'} \textbf{1}_{a > b} \textbf{1}_{a \in J_i} \textbf{1}_{b \in J_i} =  \displaystyle\sum_{a \in \breve{\mathfrak{E}}_1} \displaystyle\sum_{b \in \breve{\mathfrak{E}}_2} \textbf{1}_{a > b} \textbf{1}_{a \in J_i} \textbf{1}_{b \in J_i} = \inv (\breve{\mathscr{E}}; J_i).
	\end{flalign*} 
	
	\noindent Here, to deduce first equality, we used \eqref{ijsum}; for the second, we used the fact that each $\inv (\breve{\mathfrak{E}}_h'; J_i) = 0$ (again since $|\breve{\mathfrak{E}}_h' \cap J_i| \le 1$); for the third, we used the fact that each $\breve{\mathfrak{E}}_h' \cap J_i = \breve{\mathfrak{E}}_h \cap J_i$; and for the fourth we used \eqref{sumij} and the fact that each $\breve{\mathfrak{E}}_h \in \mathcal{M} (\breve{\textbf{E}}_h)$. This confirms the third equality in \eqref{efe}; the proof of the fourth is omitted since it is very similar. This establishes \eqref{efe}, thus implying \eqref{sumfusedz} and hence the theorem.
\end{proof}

\chapter{Transfer Operators}

\label{Operators}

In this chapter we use the $W_z$ (and $\widehat{W}_z$) weights from \Cref{wabcdrsxy} (and \eqref{wabcd2}) to define the actions of certain operators on a vector space; this will be used later in \Cref{Functions1} to define and analyze properties of certain families of (symmetric) functions. 

\section{Single-Row and Double-Row Partition Functions} 

\label{Row12Functions} 

 In this section we provide properties for partition functions of vertex models under the $W_z$ or $\widehat{W}_z$ weights (from \Cref{wabcdrsxy} and \eqref{wabcd2}, respectively) on domains consisting of one or two rows. This will later be useful in \Cref{OperatorRow} to define certain operator actions on vector spaces. 

To that end, we begin with the following definition for the single-row partition functions. 

\begin{definition}
	
	\label{mwabcd}
	
	Fix an integer $M \ge 1$; complex numbers $x, r \in \mathbb{C}$; and finite sequences $\textbf{y} = (y_1, y_2, \ldots , y_M)$ and $\textbf{s} = (s_1, s_2, \ldots , s_M)$ of complex numbers. For any elements $\textbf{B} \in \{ 0, 1 \}^n$ and $\textbf{D} \in \{ 0, 1 \}^n$, and sequences $\mathscr{A} = (\textbf{A}_1, \textbf{A}_2, \ldots , \textbf{A}_M) \subset \{ 0, 1 \}^n$ and $\mathscr{C} = (\textbf{C}_1, \textbf{C}_2, \ldots , \textbf{C}_M) \in \{ 0, 1 \}^n$, define the \emph{single-row partition functions}
	\begin{flalign}
		\label{wxyac}
		\begin{aligned} 
			& W_{x; \textbf{y}} \big( \mathscr{A}, \textbf{B}; \mathscr{C}, \textbf{D} \boldsymbol{\mid} r, \textbf{s} \big) = \displaystyle\sum_{\mathscr{J}} \displaystyle\prod_{i = 1}^M W_{x / y_i} (\textbf{A}_i, \textbf{J}_i; \textbf{C}_i, \textbf{J}_{i + 1} \boldsymbol{\mid} r, s );\\
			& \widehat{W}_{x; \textbf{y}} \big( \mathscr{A}, \textbf{B}; \mathscr{C}, \textbf{D} \boldsymbol{\mid} r, \textbf{s} \big) = \displaystyle\sum_{\mathscr{J}} \displaystyle\prod_{i = 1}^M \widehat{W}_{x / y_i} (\textbf{A}_i, \textbf{J}_i; \textbf{C}_i, \textbf{J}_{i + 1} \boldsymbol{\mid} r, s),
		\end{aligned} 
	\end{flalign}
	\index{W@$W_{x; \textbf{y}} (\mathscr{A}, \textbf{B}; \mathscr{C}, \textbf{D})$; single-row partition function}
	
	\noindent where both sums are over all sequences $\mathscr{J} = (\textbf{J}_1, \textbf{J}_2, \ldots , \textbf{J}_{M + 1}) \subset \{ 0, 1 \}^n$ such that $\textbf{J}_1 = \textbf{B}$ and $\textbf{J}_{M + 1} = \textbf{D}$; by arrow conservation, both of these sums are supported on at most one term. 

	These quantities can also be defined when $M = \infty$, in which case we set
	\begin{flalign}
	\label{wxyaclimit} 
	\begin{aligned} 
	& W_{x; \textbf{y}} \big( \mathscr{A}, \textbf{B}; \mathscr{C}, \textbf{D} \boldsymbol{\mid} r, \textbf{s} \big) = \displaystyle\sum_{\mathscr{J}} \displaystyle\prod_{i = 1}^{\infty} W_{x / y_i} (\textbf{A}_i, \textbf{J}_i; \textbf{C}_i, \textbf{J}_{i + 1} \boldsymbol{\mid} r, s_i);\\
	& \widehat{W}_{x; \textbf{y}} \big( \mathscr{A}, \textbf{B}; \mathscr{C}, \textbf{D} \boldsymbol{\mid} r, \textbf{s} \big) = \displaystyle\sum_{\mathscr{J}} \displaystyle\prod_{i = 1}^{\infty} \widehat{W}_{x / y_i} (\textbf{A}_i, \textbf{J}_i; \textbf{C}_i, \textbf{J}_{i + 1} \boldsymbol{\mid} r, s_i),
	\end{aligned} 
	\end{flalign}
	\index{W@$\widehat{W}_{x; \textbf{y}} (\mathscr{A}, \textbf{B}; \mathscr{C}, \textbf{D})$; normalized single-row partition function}
	
	\noindent where, in both equations $\mathscr{J} = (\textbf{J}_1, \textbf{J}_2, \ldots )$ is summed over all infinite sequences in $\{ 0, 1 \}^n$ such that $\textbf{J}_1 = \textbf{B}$ and $\textbf{J}_i = \textbf{D}$ for sufficiently large $i$; as above, arrow conservation implies that both sums on the right sides of \eqref{wxyaclimit} are supported on at most one term. In this $M = \infty$ case, we must further assume that parameters $x$, $r$, $\textbf{y}$, and $\textbf{s}$ are chosen so the infinite products on the right sides of \eqref{wxyaclimit} converge. 
	
\end{definition}

Observe that $W_{x; \textbf{y}}$ and $\widehat{W}_{x; \textbf{y}}$ are partition functions under the weights $W_z$ and $\widehat{W}_z$, respectively, for the vertex model with $M$ vertices in a row depicted on the left side of \Cref{wxyabcd}. More explicitly, the boundary data for this vertex model is given as follows. The vertical entrance data (meaning the colors of the incoming vertical arrows) from west to east along this row is given by $(\textbf{A}_1, \textbf{A}_2, \ldots , \textbf{A}_M)$; the horizontal entrance data is given by $\textbf{B}$; the vertical exit data by $(\textbf{C}_1, \textbf{C}_2, \ldots , \textbf{C}_M)$; and the horizontal exit data by $\textbf{D}$. 

\begin{figure}

	\begin{center}

		\begin{tikzpicture}[
		>=stealth,
		auto,
		style={
			scale = 1.2
		}
		]
		
		\draw[->, gray, ultra thick] (-1, 0) -- (0, 0);
		\draw[->, gray, ultra thick] (3, 0) -- (4, 0);
		
		\draw[->, gray, ultra thick] (0, -1) -- (0, 0);
		\draw[->, gray, ultra thick] (1, -1) -- (1, 0);
		\draw[->, gray, ultra thick] (2, -1) -- (2, 0);
		\draw[->, gray, ultra thick] (3, -1) -- (3, 0);
		
		\draw[-, gray, ultra thick, dashed] (0, 0) -- (1, 0);
		\draw[-, gray, ultra thick, dashed] (1, 0) -- (2, 0);
		\draw[-, gray, ultra thick, dashed] (2, 0) -- (3, 0);
		\draw[-, gray, ultra thick, dashed] (3, 0) -- (4, 0);

		\draw[->, gray, ultra thick] (0, 0) -- (0, 1);
		\draw[->, gray, ultra thick] (1, 0) -- (1, 1);
		\draw[->, gray, ultra thick] (2, 0) -- (2, 1);
		\draw[->, gray, ultra thick] (3, 0) -- (3, 1);

		\draw[] (0, -1) circle[radius = 0] node[below, scale = .7]{$(y_1; s_1)$};
		\draw[] (1, -1) circle[radius = 0] node[below, scale = .7]{$(y_2; s_2)$};
		\draw[] (2, -1) circle[radius = 0] node[below, scale = .7]{$(y_3; s_3)$};
		\draw[] (3, -1) circle[radius = 0] node[below, scale = .7]{$(y_4; s_4)$};
		
		\draw[] (-1, 0) circle[radius = 0] node[left, scale = .7]{$(x; r)$};
		
		\draw[] (-.6, 0) circle[radius = 0] node[above, scale = .7]{$\textbf{B}$};
		
		\draw[] (2.5, 0) circle[radius = 0] node[above, scale = .7]{$\textbf{J}_4$};
		\draw[] (3.5, 0) circle[radius = 0] node[above, scale = .7]{$\textbf{D}$};
		\draw[] (1.5, 0) circle[radius = 0] node[above, scale = .7]{$\textbf{J}_3$};
		\draw[] (.5, 0) circle[radius = 0] node[above, scale = .7]{$\textbf{J}_2$};
		
		\draw[] (0, -.6) circle[radius = 0] node[right, scale = .6]{$\textbf{A}_1$};
		\draw[] (1, -.6) circle[radius = 0] node[right, scale = .6]{$\textbf{A}_2$};
		\draw[] (2, -.6) circle[radius = 0] node[right, scale = .6]{$\textbf{A}_3$};
		\draw[] (3, -.6) circle[radius = 0] node[right, scale = .6]{$\textbf{A}_4$};
		
		\draw[] (0, .5) circle[radius = 0] node[right, scale = .7]{$\textbf{C}_1$};
		\draw[] (1, .5) circle[radius = 0] node[right, scale = .7]{$\textbf{C}_2$};
		\draw[] (2, .5) circle[radius = 0] node[right, scale = .7]{$\textbf{C}_3$};
		\draw[] (3, .52) circle[radius = 0] node[right, scale = .7]{$\textbf{C}_4$};

		\draw[] (0, -1.625) -- (0, -1.75) -- (3, -1.75) -- (3, -1.625); 
		
		\draw[] (1.5, -1.75) circle[radius = 0] node[below, scale = .7]{$M$};

		\draw[->, gray, ultra thick] (5, -.5) -- (6, -.5);
		\draw[->, gray, ultra thick] (9, -.5) -- (10, -.5);
		
		\draw[->, gray, ultra thick] (6, -1.5) -- (6, -.5);
		\draw[->, gray, ultra thick] (7, -1.5) -- (7, -.5);
		\draw[->, gray, ultra thick] (8, -1.5) -- (8, -.5);
		\draw[->, gray, ultra thick] (9, -1.5) -- (9, -.5);
		
		\draw[-, gray, ultra thick, dashed] (6, -.5) -- (7, -.5);
		\draw[-, gray, ultra thick, dashed] (7, -.5) -- (8, -.5);
		\draw[-, gray, ultra thick, dashed] (8, -.5) -- (9, -.5);
		\draw[-, gray, ultra thick, dashed] (9, -.5) -- (10, -.5);
		
		\draw[-, dashed, gray, ultra thick] (6, -.5) -- (6, .5);
		\draw[-, dashed, gray, ultra thick] (7, -.5) -- (7, .5);
		\draw[-, dashed, gray, ultra thick] (8, -.5) -- (8, .5);
		\draw[-, dashed, gray, ultra thick] (9, -.5) -- (9, .5);
		
		\draw[] (6, -1.5) circle[radius = 0] node[below, scale = .7]{$(y_1; s_1)$};
		\draw[] (7, -1.5) circle[radius = 0] node[below, scale = .7]{$(y_2; s_2)$};
		\draw[] (8, -1.5) circle[radius = 0] node[below, scale = .7]{$(y_3; s_3)$};
		\draw[] (9, -1.5) circle[radius = 0] node[below, scale = .7]{$(y_4; s_4)$};
		
		\draw[] (5, -.5) circle[radius = 0] node[left, scale = .7]{$(x_1; r_1)$};
		
		\draw[] (5.4, -.5) circle[radius = 0] node[above, scale = .7]{$\textbf{B}_1$};
		\draw[] (9.5, -.5) circle[radius = 0] node[above, scale = .7]{$\textbf{D}_1$};
		
		\draw[] (6, -1.1) circle[radius = 0] node[right, scale = .6]{$\textbf{A}_1$};
		\draw[] (7, -1.1) circle[radius = 0] node[right, scale = .6]{$\textbf{A}_2$};
		\draw[] (8, -1.1) circle[radius = 0] node[right, scale = .6]{$\textbf{A}_3$};
		\draw[] (9, -1.1) circle[radius = 0] node[right, scale = .6]{$\textbf{A}_4$};
		
		\draw[] (6, 0) circle[radius = 0] node[right, scale = .7]{$\textbf{K}_1$};
		\draw[] (7, 0) circle[radius = 0] node[right, scale = .7]{$\textbf{K}_2$};
		\draw[] (8, 0) circle[radius = 0] node[right, scale = .7]{$\textbf{K}_3$};
		\draw[] (9, .02) circle[radius = 0] node[right, scale = .7]{$\textbf{K}_4$};
		
		\draw[] (6, -2.125) -- (6, -2.25) -- (9,  -2.25) -- (9, -2.125); 
		
		\draw[] (7.5, -2.25) circle[radius = 0] node[below, scale = .7]{$M$};
		
		\draw[->, gray, ultra thick] (5, .5) -- (6, .5);
		\draw[->, gray, ultra thick] (9, .5) -- (10, .5);
		
		\draw[-, gray, ultra thick, dashed] (6, .5) -- (7, .5);
		\draw[-, gray, ultra thick, dashed] (7, .5) -- (8, .5);
		\draw[-, gray, ultra thick, dashed] (8, .5) -- (9, .5);
		\draw[-, gray, ultra thick, dashed] (9, .5) -- (10, .5);
		
		\draw[->, gray, ultra thick] (6, .5) -- (6, 1.5);
		\draw[->, gray, ultra thick] (7, .5) -- (7, 1.5);
		\draw[->, gray, ultra thick] (8, .5) -- (8, 1.5);
		\draw[->, gray, ultra thick] (9, .5) -- (9, 1.5);
		
		\draw[] (5, .5) circle[radius = 0] node[left, scale = .7]{$(x_2; r_2)$};
		
		\draw[] (5.4, .5) circle[radius = 0] node[above, scale = .7]{$\textbf{B}_2$};
		\draw[] (9.5, .5) circle[radius = 0] node[above, scale = .7]{$\textbf{D}_2$};

		\draw[] (6, 1) circle[radius = 0] node[right, scale = .7]{$\textbf{C}_1$};
		\draw[] (7, 1) circle[radius = 0] node[right, scale = .7]{$\textbf{C}_2$};
		\draw[] (8, 1) circle[radius = 0] node[right, scale = .7]{$\textbf{C}_3$};
		\draw[] (9, 1.02) circle[radius = 0] node[right, scale = .7]{$\textbf{C}_4$};

		\end{tikzpicture}
		
	\end{center}
	
	\caption{\label{wxyabcd} Shown to the left is a diagrammatic interpretation for \eqref{wxyac}; shown to the right is one for \eqref{w2}. } 	
\end{figure}

It will also be useful to define partition functions for two rows, one on top the other. This is given by the following definition. 

\begin{definition}
	
	\label{mwabcd2}
	
	Let $M$, $\textbf{y}$, $\textbf{s}$, $\mathscr{A}$, and $\mathscr{C}$ be as in \Cref{mwabcd} (with $M$ possibly infinite). Further let $\textbf{x} = (x_1, x_2)$ and $\textbf{r} = (r_1, r_2)$ denote pairs of complex numbers. Then, for any sequences $\mathscr{B} = (\textbf{B}_1, \textbf{B}_2) \subset \{ 0, 1 \}^n$ and $\mathscr{D} = (\textbf{D}_1, \textbf{D}_2) \subset \{ 0, 1 \}^n$, define the \emph{double-row partition functions}\footnote{Observe that the orderings $\mathscr{B} = (\textbf{B}_1, \textbf{B}_2)$ and $\mathscr{D} = (\textbf{D}_1, \textbf{D}_2)$ are the reverse of what they were in \Cref{ColorsFused} (see \Cref{domainpathsfused}). The reason for this discrepancy is that our rapidity parameters (which were not relevant in \Cref{ColorsFused}) will be indexed from bottom to top, so we would like to make the ordering of $\mathscr{B}$ and $\mathscr{D}$ consistent with this.}
	\begin{flalign}
	\label{w2}
	\begin{aligned} 	
	W_{\textbf{x}; \textbf{y}} ( \mathscr{A}, \mathscr{B}; \mathscr{C}, \mathscr{D} \boldsymbol{\mid} \textbf{r}, \textbf{s}) & = \displaystyle\sum_{\mathscr{K}} W_{x_1; \textbf{y}} ( \mathscr{A}, \textbf{B}_1; \mathscr{K}; \textbf{D}_1 \boldsymbol{\mid} r_1, \textbf{s}) W_{x_2, \textbf{y}} ( \mathscr{K}, \textbf{B}_2; \mathscr{C}, \textbf{D}_2 \boldsymbol{\mid} r_2, \textbf{s}); \\
	\widehat{W}_{\textbf{x}; \textbf{y}} ( \mathscr{A}, \mathscr{B}; \mathscr{C}, \mathscr{D} \boldsymbol{\mid} \textbf{r}, \textbf{s}) & = \displaystyle\sum_{\mathscr{K}} \widehat{W}_{x_1; \textbf{y}} ( \mathscr{A}, \textbf{B}_1; \mathscr{K}; \textbf{D}_1 \boldsymbol{\mid} r_1, \textbf{s}) \widehat{W}_{x_2, \textbf{y}} ( \mathscr{K}, \textbf{B}_2; \mathscr{C}, \textbf{D}_2 \boldsymbol{\mid} r_2, \textbf{s} ).
	\end{aligned} 
	\end{flalign}
	
	\noindent Here $\mathscr{K} = (\textbf{K}_1, \textbf{K}_2, \ldots , \textbf{K}_M)$ is summed over all sequences of elements in $\{ 0, 1 \}^n$, which we assume to satisfy $\textbf{K}_i = \textbf{e}_0$ for sufficiently large $i$ if $M = \infty$, and we recall the single-row partition functions $W_{x_i; \textbf{y}}$ from \eqref{wxyac}. 
	
	\end{definition}

	Observe that $W_{\textbf{x}; \textbf{y}}$ and $\widehat{W}_{\textbf{x}; \textbf{y}}$ are partition functions under the weights $W_z$ and $\widehat{W}_z$, respectively, for the vertex model with $2M$ vertices, arranged in two rows with $M$ vertices each, as depicted on the right side of \Cref{wxyabcd}. More precisely, the vertical entrance data for this model is given by $(\textbf{A}_1, \textbf{A}_2, \ldots , \textbf{A}_M)$; the horizontal entrance data (from south to north) is given by $(\textbf{B}_1, \textbf{B}_2)$; the vertical exit data is given by $(\textbf{C}_1, \textbf{C}_2, \ldots , \textbf{C}_M)$; and the horizontal exit data is given by $(\textbf{D}_1, \textbf{D}_2)$. 

	The following lemma is a quick consequence of the Yang--Baxter equation.

	\begin{lem} 
		
	\label{w2w}
	
	Fix a finite integer $M \ge 1$ and sequences of complex numbers $\textbf{\emph{x}} = (x_1, x_2)$; $\textbf{\emph{r}} = (r_1, r_2)$; $\textbf{\emph{y}} = (y_1, y_2, \ldots , y_M)$; and $\textbf{\emph{s}} = (s_1, s_2, \ldots , s_M)$. For any $\{ 0, 1 \}^n$-sequences $\mathscr{A} = (\textbf{\emph{A}}_1, \textbf{\emph{A}}_2, \ldots , \textbf{\emph{A}}_M)$; $\mathscr{B} = (\textbf{\emph{B}}_1, \textbf{\emph{B}}_2)$; $\mathscr{C} = (\textbf{\emph{C}}_1, \textbf{\emph{C}}_2, \ldots , \textbf{\emph{C}}_M)$; and $\mathscr{D} = (\textbf{\emph{D}}_1, \textbf{\emph{D}}_2)$, we have
	\begin{flalign}
	\label{w2widentity}
	\begin{aligned} 
	& \displaystyle\sum_{\textbf{\emph{I}}, \textbf{\emph{J}} \in \{ 0, 1 \}^n} W_{x_2 / x_1} (\textbf{\emph{B}}_1, \textbf{\emph{B}}_2; \textbf{\emph{I}}, \textbf{\emph{J}} \boldsymbol{\mid} r_2, r_1)  W_{(x_2, x_1); \textbf{\emph{y}}} \big( \mathscr{A}, (\textbf{\emph{J}}, \textbf{\emph{I}}); \mathscr{C}, (\textbf{\emph{D}}_2, \textbf{\emph{D}}_1) \boldsymbol{\mid} (r_2, r_1), \textbf{\emph{s}} \big) \\
	& \quad = \displaystyle\sum_{\textbf{\emph{I}}, \textbf{\emph{J}} \in \{ 0, 1 \}^n} W_{(x_1, x_2); \textbf{\emph{y}}} \big( \mathscr{A}, (\textbf{\emph{B}}_1, \textbf{\emph{B}}_2); \mathscr{C}, (\textbf{\emph{I}}, \textbf{\emph{J}}) \boldsymbol{\mid} (r_1, r_2), \textbf{\emph{s}} \big) W_{x_2 / x_1} ( \textbf{\emph{I}}, \textbf{\emph{J}}; \textbf{\emph{D}}_1, \textbf{\emph{D}}_2 \boldsymbol{\mid} r_2, r_1).
	\end{aligned} 
	\end{flalign}
	
	\noindent Diagrammatically,
		
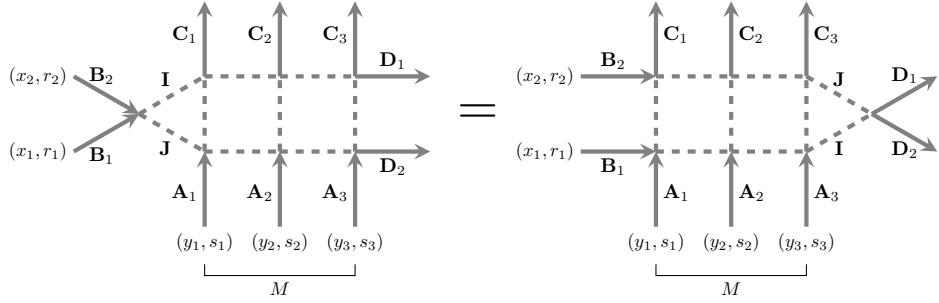
\begin{center}		
			\begin{tikzpicture}[
			>=stealth,
			auto,
			style={
				scale = 1
			}
			]
			\draw[->, ultra thick, gray] (-3.87, -4.5) -- (-3, -4);
			\draw[->, ultra thick, gray] (-3.87, -3.5) -- (-3, -4);
			\draw[->, ultra thick, gray] (-.13, -5.5) -- (-.13, -4.5);
			\draw[] (-3.87, -3.5) circle[radius = 0]  node[left, scale = .7]{$(x_2, r_2)$};
			\draw[] (-3.87, -4.5) circle[radius = 0]  node[left, scale = .7]{$(x_1, r_1)$};
			\draw[] (-.13, -5.5) circle[radius = 0]  node[below, scale = .7]{$(y_3, s_3)$};
			\draw[] (-1.13, -5.5) circle[radius = 0]  node[below, scale = .7]{$(y_2, s_2)$};
			\draw[] (-2.13, -5.5) circle[radius = 0]  node[below, scale = .7]{$(y_1, s_1)$};
			\draw[] (-1.13, -6.15) circle[radius = 0]  node[below, scale = .7]{$M$};
			\draw[] (-2.13, -5.05) circle[radius = 0]  node[left, scale = .75]{$\textbf{\emph{A}}_1$};
			\draw[] (-1.13, -5.05) circle[radius = 0]  node[left, scale = .75]{$\textbf{\emph{A}}_2$};
			\draw[] (-.13, -5.05) circle[radius = 0]  node[left, scale = .75]{$\textbf{\emph{A}}_3$}; 	
			\draw[] (-2.13, -2.95) circle[radius = 0]  node[left, scale = .75]{$\textbf{\emph{C}}_1$};
			\draw[] (-1.13, -2.95) circle[radius = 0]  node[left, scale = .75]{$\textbf{\emph{C}}_2$};
			\draw[] (-.13, -2.95) circle[radius = 0]  node[left, scale = .75]{$\textbf{\emph{C}}_3$};
			\draw[] (.37, -4.5) circle[radius = 0]  node[below, scale = .75]{$\textbf{\emph{D}}_2$};
			\draw[] (.37, -3.5) circle[radius = 0]  node[above, scale = .75]{$\textbf{\emph{D}}_1$};
			\draw[] (-3.5, -4.275) circle[radius = 0]  node[below = 2, scale = .75]{$\textbf{\emph{B}}_1$};
			\draw[] (-3.5, -3.7) circle[radius = 0]  node[above, scale = .75]{$\textbf{\emph{B}}_2$};
			\draw[] (-2.65, -4.275) circle[radius = 0]  node[below, scale = .75]{$\textbf{\emph{J}}$};
			\draw[] (-2.65, -3.725) circle[radius = 0]  node[above, scale = .75]{$\textbf{\emph{I}}$};
			\draw[] (-2.13, -6) -- (-2.13, -6.15) -- (-.13, -6.15) -- (-.13, -6);
			\draw[-, gray,ultra thick, dashed] (-3, -4) -- (-2.13, -4.5); 
			\draw[-, gray,ultra thick, dashed] (-3, -4) -- (-2.13, -3.5); 
			\draw[-, gray,ultra thick, dashed] (-.13, -4.5) -- (-.13, -3.5);
			\draw[-, gray,ultra thick, dashed] (-1.13, -4.5) -- (-1.13, -3.5);
			\draw[-, gray,ultra thick, dashed] (-2.13, -4.5) -- (-2.13, -3.5);
			\draw[-, gray,ultra thick, dashed] (-2.13, -4.5) -- (-1.13, -4.5);
			\draw[-, gray,ultra thick, dashed] (-2.13, -3.5) -- (-1.13, -3.5);
			\draw[-, gray,ultra thick, dashed] (-1.13, -4.5) -- (-.13, -4.5);
			\draw[-, gray,ultra thick, dashed] (-1.13, -3.5) -- (-.13, -3.5);
			\draw[->, gray, ultra thick] (-2.13, -5.5) -- (-2.13, -4.5);
			\draw[->, gray, ultra thick] (-2.13, -3.5) -- (-2.13, -2.5);
			\draw[->, gray, ultra thick] (-1.13, -5.5) -- (-1.13, -4.5);
			\draw[->, gray, ultra thick] (-1.13, -3.5) -- (-1.13, -2.5);
			\draw[->, gray, ultra thick] (-.13, -3.5) -- (.87, -3.5); 
			\draw[->, gray, ultra thick] (-.13, -4.5) -- (.87, -4.5); 
			\draw[->, gray, ultra thick] (-.13, -3.5) -- (-.13, -2.5);	
			\draw[] (2.87, -3.5) circle[radius = 0]  node[left, scale = .7]{$(x_2, r_2)$};
			\draw[] (2.87, -4.5) circle[radius = 0]  node[left, scale = .7]{$(x_1, r_1)$};
			\draw[] (5.87, -5.5) circle[radius = 0]  node[below, scale = .7]{$(y_3, s_3)$};
			\draw[] (4.87, -5.5) circle[radius = 0]  node[below, scale = .7]{$(y_2, s_2)$};
			\draw[] (3.87, -5.5) circle[radius = 0]  node[below, scale = .7]{$(y_1, s_1)$};
			\draw[] (3.87, -5.05) circle[radius = 0]  node[right, scale = .75]{$\textbf{\emph{A}}_1$};
			\draw[] (4.87, -5.05) circle[radius = 0]  node[right, scale = .75]{$\textbf{\emph{A}}_2$};
			\draw[] (5.87, -5.05) circle[radius = 0]  node[right, scale = .75]{$\textbf{\emph{A}}_3$}; 
			\draw[] (3.87, -2.95) circle[radius = 0]  node[right, scale = .75]{$\textbf{\emph{C}}_1$};
			\draw[] (4.87, -2.95) circle[radius = 0]  node[right, scale = .75]{$\textbf{\emph{C}}_2$};
			\draw[] (5.87, -2.95) circle[radius = 0]  node[right, scale = .75]{$\textbf{\emph{C}}_3$};			
			\draw[] (3.3, -3.5) circle[radius = 0]  node[above, scale = .75]{$\textbf{\emph{B}}_2$};
			\draw[] (3.3, -4.5) circle[radius = 0]  node[below, scale = .75]{$\textbf{\emph{B}}_1$};			
			\draw[] (6.305, -3.725) circle[radius = 0]  node[above, scale = .75]{$\textbf{\emph{J}}$};
			\draw[] (6.305, -4.275) circle[radius = 0]  node[below, scale = .75]{$\textbf{\emph{I}}$};			
			\draw[] (7.177, -3.725) circle[radius = 0]  node[above = 1, scale = .75]{$\textbf{\emph{D}}_1$};
			\draw[] (7.177, -4.275) circle[radius = 0]  node[below, scale = .75]{$\textbf{\emph{D}}_2$};			
			\draw[] (4.87, -6.15) circle[radius = 0]  node[below, scale = .7]{$M$};			
			\draw[] (3.87, -6) -- (3.87, -6.15) -- (5.87, -6.15) -- (5.87, -6);		
			\draw[->, gray, ultra thick] (4.87, -5.5) -- (4.87, -4.5); 
			\draw[->, gray, ultra thick] (5.87, -5.5) -- (5.87, -4.5); 
			\draw[->, gray, ultra thick] (4.87, -3.5) -- (4.87, -2.5); 
			\draw[->, gray, ultra thick] (5.87, -3.5) -- (5.87, -2.5); 
			\draw[->, gray, ultra thick] (3.87, -5.5) -- (3.87, -4.5); 
			\draw[->, gray, ultra thick] (2.87, -3.5) -- (3.87, -3.5); 
			\draw[->, gray, ultra thick] (2.87, -4.5) -- (3.87, -4.5); 
			\draw[->, gray, ultra thick] (3.87, -3.5) -- (3.87, -2.5); 
			\draw[->, gray, ultra thick] (6.74, -4) -- (7.61, -4.5); 
			\draw[->, gray, ultra thick] (6.74, -4) -- (7.61, -3.5); 		
			\draw[-, gray, ultra thick, dashed] (3.87, -4.5) -- (3.87, -3.5);
			\draw[-, gray, ultra thick, dashed] (4.87, -4.5) -- (4.87, -3.5);
			\draw[-, gray, ultra thick, dashed] (5.87, -4.5) -- (5.87, -3.5);
			\draw[-, gray, ultra thick, dashed] (5.87, -4.5) -- (6.74, -4); 
			\draw[-, gray, ultra thick, dashed] (5.87, -3.5) -- (6.74, -4);
			\draw[-, gray, ultra thick, dashed] (3.87, -4.5) -- (4.87, -4.5);
			\draw[-, gray, ultra thick, dashed] (3.87, -3.5) -- (4.87, -3.5);
			\draw[-, gray, ultra thick, dashed] (4.87, -4.5) -- (5.87, -4.5);
			\draw[-, gray, ultra thick, dashed] (4.87, -3.5) -- (5.87, -3.5);		
			\filldraw[fill=white, draw=black] (1.5, -4) circle [radius=0] node[scale = 2]{$=$};
			\draw[] (8, -4) circle [radius = 0] node[]{.};
			\end{tikzpicture}
		\end{center}

\end{lem} 

\begin{proof}
This follows from $M$ applications of the Yang--Baxter equation \Cref{wabcdproduct2}.
\end{proof}

\section{Transfer Operators and Commutation Relations}

\label{OperatorRow}

The functions we will consider in this text are expressible in terms of actions of certain operators on vector spaces. In this section we define these operators and provide some relations they satisfy. Throughout this section, we fix an integer $n \ge 1$ and infinite sequences of complex numbers $\textbf{y} = (y_1, y_2, \ldots )$ and $\textbf{s} = (s_1, s_2, \ldots )$. 

Let $\mathbb{V} = \mathbb{V}_{\textbf{y}; \textbf{s}}$ denote the infinte-dimensional vector space spanned by basis vectors of the form $| \mathscr{A} \rangle$, where $\mathscr{A} = (\textbf{A}_1, \textbf{A}_2, \ldots )$ ranges over all \emph{finitary} sequences of elements in $\{ 0, 1 \}^n$, namely, those that satisfy $\textbf{A}_j = \textbf{e}_0 = (0, 0, \ldots , 0)$ for all but finitely many $j$. Similarly, let $\mathbb{V}^* = \mathbb{V}_{\textbf{y}; \textbf{s}}^*$ denote the space spanned by finitary dual vectors, namely, those of the form $\langle \mathscr{C} |$ over all finitary sequences of elements $\mathscr{C} = (\textbf{C}_1, \textbf{C}_2, \ldots )$ in $\{ 0, 1 \}^n$.  We impose an inner product on $\mathbb{V}^* \times \mathbb{V}$ by first setting $\langle \mathscr{C}| \mathscr{A} \rangle = \textbf{1}_{\mathscr{A} = \mathscr{C}}$, for any finitary $\mathscr{A}$ and $\mathscr{C}$, and then extending to all of $\mathbb{V}^* \times \mathbb{V}$ by bilinearity. 
\index{A@$\mid$$\mathscr{A} \rangle$}
\index{C@$\langle \mathscr{C}$$\mid$}
\index{V@$\mathbb{V}, \mathbb{V}^*$; vector spaces}

Now, for any complex numbers $x, r \in \mathbb{C}$ and sequences $\textbf{B}, \textbf{D} \in \{ 0, 1 \}^n$ we define \emph{transfer operators} $\mathds{T}_{\textbf{B}; \textbf{D}} = \mathds{T}_{\textbf{B}; \textbf{D}} (x; r): \mathbb{V} \rightarrow \mathbb{V}$ and $\widehat{\mathds{T}}_{\textbf{B}; \textbf{D}} = \widehat{\mathds{T}}_{\textbf{B}; \textbf{D}} (x; r): \mathbb{V} \rightarrow \mathbb{V}$ by first setting, for any finitary sequence $\mathscr{A}$, 
\begin{flalign}
\label{wta} 
\begin{aligned} 
& \mathds{T}_{\textbf{B}; \textbf{D}} | \mathscr{A} \rangle = \displaystyle\sum_{\mathscr{C}} W_{x; \textbf{y}} (\mathscr{A}, \textbf{B}; \mathscr{C}, \textbf{D} \boldsymbol{\mid} r, \textbf{s} ) | \mathscr{C} \rangle; \qquad \widehat{\mathds{T}}_{\textbf{B}; \textbf{D}} | \mathscr{A} \rangle = \displaystyle\sum_{\mathscr{C}}  \widehat{W}_{x; \textbf{y}} (\mathscr{A}, \textbf{B}; \mathscr{C}, \textbf{D} \boldsymbol{\mid} r, \textbf{s}) | \mathscr{C} \rangle,
\end{aligned} 
\end{flalign} \index{T@$\mathds{T}_{\textbf{B}; \textbf{D}} (x; r)$; transfer operator}\index{T@$\widehat{\mathds{T}}_{\textbf{B}; \textbf{D}} (x; r)$; normalized operator}

\noindent where both sums are over all finitary sequences $\mathscr{C}$, and then extending their actions to all of $\mathbb{V}$ by linearity. Here, we recall the $M = \infty$ single-row partition functions $W_{x; \textbf{y}}$ and $\widehat{W}_{x; \textbf{y}}$ from \eqref{wxyaclimit} and assume the parameters $x, r, \textbf{y}, \textbf{s}$ are chosen such that the infinite products there converge. 

 By \eqref{wta}, we find for any finitary sequences $\mathscr{A}$ and $\mathscr{C}$ that
\begin{flalign}
\label{twxyinner}
\begin{aligned}
& W_{x; \textbf{y}} (\mathscr{A}, \textbf{B}; \mathscr{C}, \textbf{D} \boldsymbol{\mid} r, \textbf{s} ) = \big\langle \mathscr{C} \big| \mathds{T}_{\textbf{B}; \textbf{D}} (x; r) \big| \mathscr{A} \big\rangle; \qquad \widehat{W}_{x; \textbf{y}} (\mathscr{A}, \textbf{B}; \mathscr{C}, \textbf{D} \boldsymbol{\mid} r, \textbf{s} ) = \big\langle \mathscr{C} \big| \widehat{\mathds{T}}_{\textbf{B}; \textbf{D}} (x; r) \big| \mathscr{A} \big\rangle,
\end{aligned} 
\end{flalign} 

\noindent and so $\mathds{T}_{\textbf{B}; \textbf{D}}$ and $\widehat{\mathds{T}}_{\textbf{B}; \textbf{D}}$ admit dual actions on $\mathbb{V}^*$ given by 
\begin{flalign}
\label{wta2}
\begin{aligned} 
& \langle \mathscr{C} | \mathds{T}_{\textbf{B}; \textbf{D}} = \displaystyle\sum_{\mathscr{A}} W_{x; \textbf{y}} ( \mathscr{A}, \textbf{B}; \mathscr{C}, \textbf{D} \boldsymbol{\mid} r, \textbf{s}) \langle \mathscr{A} |; \qquad  \langle \mathscr{C} | \widehat{\mathds{T}}_{\textbf{B}; \textbf{D}} = \displaystyle\sum_{\mathscr{A}}  \widehat{W}_{x; \textbf{y}} ( \mathscr{A}, \textbf{B}; \mathscr{C} , \textbf{D} \boldsymbol{\mid} r, \textbf{s} ) \langle \mathscr{A}|.
\end{aligned} 
\end{flalign}

Now, \Cref{w2w} will imply certain commutation relations between the operators $\mathds{T}_{\textbf{B}; \textbf{D}}$ and $\widehat{\mathds{T}}_{\textbf{B}; \textbf{D}}$ that will be useful for us. To explain them, for any $x, r \in \mathbb{C}$, we define the operators 
\begin{flalign}
\label{bd}
\mathds{B} (x; r) = \widehat{\mathds{T}}_{\textbf{e}_0; \textbf{e}_{[1, n]}} (x; r); \qquad \mathds{C} (x; r) = \mathds{T}_{\textbf{e}_{[1, n]}; \textbf{e}_0}; \qquad \mathds{D} (x; r) = \mathds{T}_{\textbf{e}_0; \textbf{e}_0} (x; r),
\end{flalign}
\index{B@$\mathds{B} (x; r)$}
\index{C@$\mathds{C} (x; r)$}
\index{D@$\mathds{D} (x; r)$}

\noindent where we recall that $\textbf{e}_0 = (0, 0, \ldots , 0)$ and $\textbf{e}_{[1, n]} = (1, 1, \ldots , 1)$.

Observe that these operators are well-defined for any choices of parameters $(x; r; \textbf{y}; \textbf{s})$. Indeed, if $\textbf{D} = \textbf{e}_0$ then, since $\mathscr{A}$ and $\mathscr{C}$ are finitary, all but finitely many terms in the product appearing on the left side of the first equation of \eqref{wxyaclimit} are of the form $W_z (\textbf{e}_0, \textbf{e}_0; \textbf{e}_0, \textbf{e}_0 \boldsymbol{\mid} r, s) = 1$; so the product defining the $W$ weight there converges, implying that $\mathds{C}$ and $\mathds{D}$ are well-defined. A similar statement holds for the $\widehat{W}$ weights when $\textbf{D} = \textbf{e}_{[1, n]}$, implying that $\mathds{B}$ is also well-defined.  

The following two results, which are both consequences of (limits of) \Cref{w2w}, provide commutation relations for these $\mathds{B}$, $\mathds{C}$, and $\mathds{D}$ operators.

\begin{lem}
	
	\label{bd1limit}
	
	Fix  $x_1, x_2, r_1, r_2 \in \mathbb{C}$ and $\textbf{\emph{y}} = (y_1, y_2, \ldots ) \subset \mathbb{C}$ and $\textbf{\emph{s}} = (s_1, s_2, \ldots ) \subset \mathbb{C}$. Then, as operators on $\mathbb{V}$, we have that 
	\begin{flalign}
	\label{bbdd12limit}
	\begin{aligned}
	& \mathds{B} (x_1; r_1) \mathds{B} (x_2; r_2) = \left( \displaystyle\frac{r_1^2 x_2}{r_2^2 x_1} \right)^n \displaystyle\frac{(r_2^2 x_1 x_2^{-1}; q)_n}{(r_1^2 x_1^{-1} x_2; q)_n} \cdot \mathds{B} (x_2; r_2) \mathds{B} (x_1; r_1); \\
	& \mathds{C} (x_1; r_1) \mathds{C} (x_2; r_2) = \left( \displaystyle\frac{r_2^2 x_1}{r_1^2 x_2} \right)^n \displaystyle\frac{(r_1^2 x_1^{-1} x_2; q)_n}{(r_2^2 x_1 x_2^{-1}; q)_n} \cdot \mathds{C} (x_2; r_2) \mathds{C} (x_1; r_1); \\
	& \mathds{D} (x_1; r_1) \mathds{D} (x_2; r_2) = \mathds{D} (x_2; r_2) \mathds{D} (x_1; r_1).
	\end{aligned}
	\end{flalign}
	
\end{lem}

\begin{proof}
	
	First observe for any $\textbf{B}_1, \textbf{B}_2, \textbf{D}_1, \textbf{D}_2 \in \{ 0, 1 \}^n$ and finitary $\mathscr{A}$ and $\mathscr{C}$ that
	\begin{flalign}
	\label{ct1t2a} 
	\begin{aligned} 
	\big\langle \mathscr{C} \big| \mathds{T}_{\textbf{B}_2; \textbf{D}_2} (x_2; r_2) \mathds{T}_{\textbf{B}_1; \textbf{D}_1} (x_1; r_1) \big| \mathscr{A} \big\rangle & = \displaystyle\sum_{\mathscr{K}} \big\langle \mathscr{C} \big| \mathds{T}_{\textbf{B}_2; \textbf{D}_2} (x_2; r_2) \big| \mathscr{K} \big\rangle \big\langle \mathscr{K} \big| \mathds{T}_{\textbf{B}_1; \textbf{D}_1} (x_1; r_1) \big| \mathscr{A} \big\rangle \\
	& =  \displaystyle\sum_{\mathscr{K}} W_{x_2; \textbf{y}} (\mathscr{K}, \textbf{B}_2; \mathscr{C}, \textbf{D}_2 \boldsymbol{\mid} r_2, \textbf{s}) W_{x_1; \textbf{y}} (\mathscr{A}, \textbf{B}_1; \mathscr{K}, \textbf{D}_1 \boldsymbol{\mid} r_1, \textbf{s})  \\
	& = W_{(x_1, x_2); \textbf{y}} \big( \mathscr{A}, (\textbf{B}_1, \textbf{B}_2); \mathscr{C}, (\textbf{D}_1, \textbf{D}_2) \boldsymbol{\mid} (r_1, r_2), \textbf{s} \big),
	\end{aligned}
	\end{flalign}
	
	\noindent where in the second statement we applied \eqref{twxyinner} and in the last we applied \eqref{w2}. Similarly, we have that 
	\begin{flalign}
	\label{t1t2innerac} 
	\begin{aligned}
	\big\langle \mathscr{C} \big| \widehat{\mathds{T}}_{\textbf{B}_2; \textbf{D}_2} (x_2; r_2) \widehat{\mathds{T}}_{\textbf{B}_1; \textbf{D}_1} (x_1; r_1) \big| \mathscr{A} \big\rangle = \widehat{W}_{(x_1, x_2); \textbf{y}} \big( \mathscr{A}, (\textbf{B}_1, \textbf{B}_2); \mathscr{C}, (\textbf{D}_1, \textbf{D}_2) \boldsymbol{\mid} (r_1, r_2), \textbf{s} \big). 
	\end{aligned}
	\end{flalign}
	
	Given these facts, \eqref{bbdd12limit} will follow from specializations of \eqref{w2widentity}. In particular, let us apply \Cref{w2w} with the $(\textbf{B}_1, \textbf{B}_2)$ and $(\textbf{D}_1, \textbf{D}_2)$ there both equal to $(\textbf{e}_0, \textbf{e}_0)$. By arrow conservation, the sums on both sides of \eqref{w2widentity} are supported on the $(\textbf{I}, \textbf{J}) = (\textbf{e}_0, \textbf{e}_0)$ term. Since $W_z ( \textbf{e}_0, \textbf{e}_0; \textbf{e}_0, \textbf{e}_0 \boldsymbol{\mid} r, s) = 1$ by \eqref{wabcd01nequation}, this yields for any integer $M \ge 1$ that
	\begin{flalign}
	\label{wx2x1mw} 
	\begin{aligned} 
	 W_{(x_2, x_1); \textbf{y}_{[1, M]}} & \big( \mathscr{A}_{[1, M]}, (\textbf{e}_0, \textbf{e}_0); \mathscr{C}_{[1, M]}, (\textbf{e}_0, \textbf{e}_0) \boldsymbol{\mid} (r_2, r_1), \textbf{s}_{[1, M]} \big) \\
	& =  W_{(x_1, x_2); \textbf{y}_{[1, M]}} \big( \mathscr{A}_{[1, M]}, (\textbf{e}_0, \textbf{e}_0); \mathscr{C}_{[1, M]}, (\textbf{e}_0, \textbf{e}_0) \boldsymbol{\mid} (r_1, r_2), \textbf{s}_{[1, M]} \big),
	\end{aligned} 
	\end{flalign}
	
	\noindent for any finitary $\mathscr{A}$ and $\mathscr{C}$, where we have set $\mathscr{X}_{[1, M]} = (\textbf{X}_1, \textbf{X}_2, \ldots , \textbf{X}_M)$ for any sequence $\mathscr{X} = (\textbf{X}_1, \textbf{X}_2, \ldots )$. Taking the limit as $M$ tends to $\infty$ on both sides of \eqref{wx2x1mw} yields 
	\begin{flalign*}
	 W_{(x_2, x_1); \textbf{y}} & \big( \mathscr{A}, (\textbf{e}_0, \textbf{e}_0); \mathscr{C}, (\textbf{e}_0, \textbf{e}_0) \boldsymbol{\mid} (r_2, r_1), \textbf{s} \big)  =  W_{(x_1, x_2); \textbf{y}} \big( \mathscr{A}, (\textbf{e}_0, \textbf{e}_0); \mathscr{C}, (\textbf{e}_0, \textbf{e}_0) \boldsymbol{\mid} (r_1, r_2), \textbf{s} \big), 
	\end{flalign*}
	
	\noindent which by \eqref{ct1t2a} yields the third statement of \eqref{bbdd12limit}. 
	
	To establish the first statement of \eqref{bbdd12limit}, we apply the $(\textbf{B}_1, \textbf{B}_2) = (\textbf{e}_0, \textbf{e}_0)$ and $(\textbf{D}_1, \textbf{D}_2) = \big( \textbf{e}_{[1, n]}, \textbf{e}_{[1, n]} \big)$ case of \Cref{w2w}; multiply both sides of \eqref{w2widentity} by
	\begin{flalign*}
	\displaystyle\prod_{j = 1}^M \displaystyle\frac{(s_j^2 x_1 y_j^{-1}; q)_n}{s_j^{2n} (x_1 y_j^{-1}; q)_n} \displaystyle\frac{(s_j^2 x_2 y_j^{-1}; q)_n}{s_j^{2n} (x_2 y_j^{-1}; q)_n};
	\end{flalign*}
	
	\noindent and use \eqref{wabcd2} to deduce that 
	\begin{flalign}
	\label{w2widentity2}
	\begin{aligned} 
	& \displaystyle\sum_{\textbf{I}, \textbf{J} \in \{ 0, 1 \}^n} \widehat{W}_{(x_2, x_1); \textbf{y}_{[1, M]}} \Big( \mathscr{A}, (\textbf{J}, \textbf{I}); \mathscr{C}, \big(\textbf{e}_{[1, n]}, \textbf{e}_{[1, n]} \big) \boldsymbol{\mid} (r_2, r_1), \textbf{s}_{[1, M]} \Big) W_{x_2 / x_1} (\textbf{e}_0, \textbf{e}_0; \textbf{I}, \textbf{J} \boldsymbol{\mid} r_2, r_1) \\
	& \quad = \displaystyle\sum_{\textbf{I}, \textbf{J} \in \{ 0, 1 \}^n} \widehat{W}_{(x_1, x_2); \textbf{y}_{[1, M]}} \big( \mathscr{A}, ( \textbf{e}_0, \textbf{e}_0 ); \mathscr{C}, (\textbf{I}, \textbf{J}) \boldsymbol{\mid} (r_1, r_2), \textbf{s}_{[1, M]} \big) W_{x_2 / x_1} ( \textbf{I}, \textbf{J}; \textbf{e}_{[1, n]}, \textbf{e}_{[1, n]} \boldsymbol{\mid} r_2, r_1).
	\end{aligned} 
	\end{flalign}
	
	\noindent By arrow conservation, the left side of \eqref{w2widentity2} is supported on the term $(\textbf{I}, \textbf{J}) = (\textbf{e}_0, \textbf{e}_0)$ and the right side is supported on the term $(\textbf{I}, \textbf{J}) = \big( \textbf{e}_{[1, n]}, \textbf{e}_{[1, n]})$. Using the expressions from \eqref{wabcd01nequation} for the weights $W_z ( \textbf{e}_0, \textbf{e}_0; \textbf{e}_0, \textbf{e}_0 \boldsymbol{\mid} r, s)$ and $W_z \big( \textbf{e}_{[1, n]}, \textbf{e}_{[1, n]}; \textbf{e}_{[1, n]}, \textbf{e}_{[1, n]} \boldsymbol{\mid} r, s \big)$, it follows from \eqref{w2widentity2} that 
	\begin{flalign}
	\label{w2widentity3}
	\begin{aligned} 
	& \widehat{W}_{(x_2, x_1); \textbf{y}_{[1, M]}} \Big( \mathscr{A}, (\textbf{e}_0, \textbf{e}_0); \mathscr{C}, \big( \textbf{e}_{[1, n]}, \textbf{e}_{[1, n]} \big) \boldsymbol{\mid} (r_2, r_1), \textbf{s}_{[1, M]} \Big) \\
	& \quad =  \widehat{W}_{(x_1, x_2); \textbf{y}_{[1, M]}} \Big( \mathscr{A}, ( \textbf{e}_0, \textbf{e}_0); \mathscr{C}, \big( \textbf{e}_{[1, n]}, \textbf{e}_{[1, n]} \big) \boldsymbol{\mid} (r_1, r_2), \textbf{s}_{[1, M]} \Big) \bigg( \displaystyle\frac{r_1^2 x_2}{r_2^2 x_1} \bigg)^n \displaystyle\frac{(r_2^2 x_1 x_2^{-1}; q)_n}{(r_1^2 x_1^{-1} x_2; q)_n}.
	\end{aligned} 
	\end{flalign}
	
	\noindent Now the first statement of \eqref{bbdd12limit} follows from taking the limit of \eqref{w2widentity3} as $M$ tends to $\infty$ and applying \eqref{t1t2innerac}. The proof of the second statement in \eqref{bbdd12limit} is entirely analogous (obtained from the $(\textbf{B}_1, \textbf{B}_2) = \big( \textbf{e}_{[1, n]}, \textbf{e}_{[1, n]} \big)$ and $(\textbf{D}_1, \textbf{D}_2) = (\textbf{e}_0, \textbf{e}_0)$ case of \Cref{w2w}) and is therefore omitted.
	\end{proof}

	\begin{prop}
		
		\label{bd2limit}
		
		Fix complex numbers $u, w, r, t \in \mathbb{C}$ and sequences $\textbf{\emph{y}} = (y_1, y_2, \ldots )$ and $\textbf{\emph{s}} = (s_1, s_2, \ldots )$ of complex numbers. If there exists an integer $K > 0$ such that 
		\begin{flalign}
		\label{x1x2bdestimate}
		\displaystyle\sup_{j > K}	\displaystyle\max_{\substack{a, b \in [0, n] \\ (a, b) \ne (n, 0)}} \Bigg| s_j^{2a + 2b - 2n} \displaystyle\frac{(s_j^2 u y_j^{-1}; q)_n (u y_j^{-1}; q)_a}{(u y_j^{-1}; q)_n (s_j^2 u y_j^{-1}; q)_a} \displaystyle\frac{(w y_j^{-1}; q)_b}{(s_j^2 w y_j^{-1}; q)_b} \Bigg| < 1,
		\end{flalign}
		
		\noindent then, as operators on $\mathbb{V}$, we have 
		\begin{flalign}
		\label{bdlimit}
		\mathds{B} (u; r) \mathds{D} (w; t) = \displaystyle\frac{(t^2 u w^{-1}; q)_n}{t^{2n} (u w^{-1}; q)_n} \cdot \mathds{D} (w; t) \mathds{B} (u; r).
		\end{flalign}
		
	\end{prop}
	
	\begin{proof} 
	
	For any integer $M \ge 1$ (including $M = \infty$); pairs $\textbf{x} = (x_1, x_2)$ and $\textbf{r} = (r_1, r_2)$ of complex numbers; and (finitary, if $M = \infty$) sequences $\mathscr{A} = (\textbf{A}_1, \textbf{A}_2, \ldots , \textbf{A}_M)$ and $\mathscr{C} = (\textbf{C}_1, \textbf{C}_2, \ldots , \textbf{C}_M)$ of elements in $\{ 0, 1 \}^n$, denote 
	\begin{flalign}
	\label{wabbcddw} 
	\begin{aligned}
	& \widetilde{W}_{\textbf{x}; \textbf{y}} \big( \mathscr{A}, (\textbf{B}_1, \textbf{B}_2); \mathscr{C}, (\textbf{D}_1, \textbf{D}_2) \boldsymbol{\mid} \textbf{r}, \textbf{s} \big) = \displaystyle\sum_{\mathscr{K}} \widehat{W}_{x_1; \textbf{y}} ( \mathscr{A}, \textbf{B}_1; \mathscr{K}, \textbf{D}_1 \boldsymbol{\mid} r_1, \textbf{s}) W_{x_2, \textbf{y}} (\mathscr{K}, \textbf{B}_2; \mathscr{C}, \textbf{D}_2 \boldsymbol{\mid} r_2, \textbf{s}); \\
	& \widetilde{W}_{\textbf{x}; \textbf{y}}' (\mathscr{A}, (\textbf{B}_1, \textbf{B}_2); \mathscr{C}, (\textbf{D}_1, \textbf{D}_2) \boldsymbol{\mid} \textbf{r}, \textbf{s}) = \displaystyle\sum_{\mathscr{K}} W_{x_1; \textbf{y}} ( \mathscr{A}, \textbf{B}_1; \mathscr{K}, \textbf{D}_1 \boldsymbol{\mid} r_1, \textbf{s} )\widehat{W}_{x_2, \textbf{y}} (\mathscr{K}, \textbf{B}_2; \mathscr{C}, \textbf{D}_2 \boldsymbol{\mid} r_2, \textbf{s}),
	\end{aligned} 
	\end{flalign}
	
	\noindent where $\mathscr{K} = (\textbf{K}_1, \textbf{K}_2, \ldots , \textbf{K}_M)$ is summed over all (finitary, if $M = \infty$) sequences of elements in $\{ 0, 1 \}^n$. Then, analogous reasoning as applied in \eqref{ct1t2a} yields for any finitary $\mathscr{A}$ and $\mathscr{C}$ that
	\begin{flalign*}
	& \big\langle \mathscr{C} \big|  \mathds{T}_{\textbf{B}_2; \textbf{D}_2} (x_2; r_2) \widehat{\mathds{T}}_{\textbf{B}_1; \textbf{D}_1} (x_1; r_1) \big| \mathscr{A} \big\rangle = \widetilde{W}_{\textbf{x}; \textbf{y}} \big( \mathscr{A}, (\textbf{B}_1, \textbf{B}_2); \mathscr{C}, (\textbf{D}_1, \textbf{D}_2) \boldsymbol{\mid} \textbf{r}, \textbf{s} \big); \\
	& \big\langle \mathscr{C} \big| \widehat{\mathds{T}}_{\textbf{B}_2; \textbf{D}_2} (x_2; r_2) \mathds{T}_{\textbf{B}_1; \textbf{D}_1} (x_1; r_1) \big| \mathscr{A} \big\rangle = \widetilde{W}_{\textbf{x}; \textbf{y}}' \big( \mathscr{A}, (\textbf{B}_1, \textbf{B}_2); \mathscr{C}, (\textbf{D}_1, \textbf{D}_2) \boldsymbol{\mid} \textbf{r}, \textbf{s} \big). 
	\end{flalign*}
	
	\noindent Therefore, to establish \eqref{bdlimit} it suffices to show under \eqref{x1x2bdestimate} that for any finitary $\mathscr{A}$ and $\mathscr{C}$ we have 
	\begin{flalign}
	\label{wbwd}
	\begin{aligned} 
	 \widetilde{W}_{(w, u); \textbf{y}}' \Big( & \mathscr{A}, (\textbf{e}_0, \textbf{e}_0); \mathscr{C}, \big(\textbf{e}_0, \textbf{e}_{[1, n]} \big) \Big| (t, r), \textbf{s} \Big) \\
	 & = \displaystyle\frac{(t^2 u w^{-1}; q)_n}{t^{2n} (u w^{-1}; q)_n}  \widetilde{W}_{(u, w); \textbf{y}} \Big( \mathscr{A}, (\textbf{e}_0, \textbf{e}_0); \mathscr{C}, \big( \textbf{e}_{[1, n]}, \textbf{e}_0 \big) \Big| (r, t); \textbf{s} \Big). 
	 \end{aligned}  
	\end{flalign}
	
	 To that end, we apply the $(\textbf{B}_1, \textbf{B}_2) = (\textbf{e}_0, \textbf{e}_0)$ and $(\textbf{D}_1, \textbf{D}_2) = \big( \textbf{e}_0, \textbf{e}_{[1, n]} \big)$ case of \Cref{w2w}; multiply both sides of \eqref{w2widentity} by
	\begin{flalign*}
	\displaystyle\prod_{j = 1}^M \displaystyle\frac{(s_j^2 u y_j^{-1}; q)_n}{s_j^{2n} (u y_j^{-1}; q)_n};
	\end{flalign*}
	
	\noindent use \eqref{wabcd2}; and recall by arrow conservation that $W_{u / w} (\textbf{e}_0, \textbf{e}_0; \textbf{I}, \textbf{J}) \boldsymbol{\mid} t, r) = 0$ unless $\textbf{I} = \textbf{e}_0 = \textbf{J}$, in which case it is equal to $1$ by \eqref{wabcd01nequation}, to deduce for any integer $M \ge 1$ that 
	\begin{flalign}
	\label{w2widentity4}
	\begin{aligned} 
	& \widetilde{W}_{(u, w); \textbf{y}_{[1, M]}} \Big( \mathscr{A}_{[1, M]}, (\textbf{e}_0, \textbf{e}_0); \mathscr{C}_{[1, M]}, \big( \textbf{e}_{[1, n]}, \textbf{e}_0 \big) \Big| (r, t), \textbf{s}_{[1, M]} \Big) \\
	& \quad = \displaystyle\sum_{\textbf{I}, \textbf{J} \in \{ 0, 1 \}^n} \widetilde{W}_{(w, u); \textbf{y}_{[1, M]}}' \big( \mathscr{A}_{[1, M]}, (\textbf{e}_0, \textbf{e}_0); \mathscr{C}_{[1, M]}, (\textbf{I}, \textbf{J}) \boldsymbol{\mid} (t, r), \textbf{s}_{[1, M]} \big) W_{u / w} ( \textbf{I}, \textbf{J}; \textbf{e}_0, \textbf{e}_{[1, n]} \boldsymbol{\mid} r, t),
	\end{aligned} 
	\end{flalign}
	
	\noindent where we have set $\mathscr{X}_{[1, M]} = (\textbf{X}_1, \textbf{X}_2, \ldots , \textbf{X}_M)$ for any sequence $\mathscr{X} = (\textbf{X}_1, \textbf{X}_2, \ldots )$. 
	
	We will show under \eqref{x1x2bdestimate} that the right side of \eqref{w2widentity4} is asymptotically (in the limit as $M$ tends to $\infty$) supported on the term $(\textbf{I}, \textbf{J}) = \big( \textbf{e}_0, \textbf{e}_{[1, n]} \big)$, in the sense that 
	\begin{flalign}
	\label{limitmw} 
	\displaystyle\lim_{M \rightarrow \infty} \widetilde{W}_{(w, u); \textbf{y}_{[1, M]}}' \big( \mathscr{A}_{[1, M]}, (\textbf{e}_0, \textbf{e}_0); \mathscr{C}_{[1, M]}, (\textbf{I}, \textbf{J}) \boldsymbol{\mid} (t, r), \textbf{s}_{[1, M]} \big) = 0, \quad \text{unless $(\textbf{I}, \textbf{J}) = \big( \textbf{e}_0, \textbf{e}_{[1, n]} \big)$}.
	\end{flalign}
	
	Given \eqref{limitmw}, \eqref{wbwd} (and therefore the lemma) follows from first inserting \eqref{limitmw} into limit of \eqref{w2widentity4} as $M$ tends to $\infty$, and then using the explicit form (given by the second statement of \eqref{wabcd01nequation}) for $W_{u / w} \big( \textbf{e}_0, \textbf{e}_{[1, n]}; \textbf{e}_0, \textbf{e}_{[1, n]} \boldsymbol{\mid} r, t \big)$. Thus, it remains to establish \eqref{limitmw}. 
	
	To that end, first observe from \eqref{wze0b} and \eqref{wabcd2} that, for any integer $i \ge 1$ and sequence $\textbf{B} \in \{ 0, 1 \}^n$ with $|\textbf{B}| = b$, 
	\begin{flalign*}
	& \widehat{W}_{u / y_j} ( \textbf{e}_0, \textbf{B}; \textbf{e}_0, \textbf{B} \boldsymbol{\mid} r, s_j \big) = s_j^{2b - 2n} \displaystyle\frac{(s_1^2 u y_j^{-1}; q)_n (u y_j^{-1}; q)_b}{(u y_j^{-1}; q)_n (s_j^2 u y_j^{-1}; q)_b}; \\
	& W_{w / y_j} ( \textbf{e}_0, \textbf{B}; \textbf{e}_0, \textbf{B} \boldsymbol{\mid} t, s_j) = s_j^{2b} \displaystyle\frac{(w y_j^{-1}; q)_b}{(s_j^2 w y_j^{-1}; q)_b}. 
	\end{flalign*}
	
	\noindent So, \eqref{x1x2bdestimate} implies the existence of some real number $\varepsilon > 0$ such that, for sufficiently large $K$,
	\begin{flalign}
	\label{wkj}
	\displaystyle\sup_{j > K} \displaystyle\max_{\substack{\textbf{I}, \textbf{J} \in \{ 0, 1 \}^n \\ (\textbf{I}, \textbf{J}) \ne (\textbf{e}_0, \textbf{e}_{[1, n]})}} \big|  \widehat{W}_{u / y_j} ( \textbf{e}_0, \textbf{J}; \textbf{e}_0, \textbf{J} \boldsymbol{\mid} r, s_j ) W_{w / y_j} ( \textbf{e}_0, \textbf{I}; \textbf{e}_0, \textbf{I} \boldsymbol{\mid} t, s_j) \big| < 1 -\varepsilon.
	\end{flalign}

	\noindent We may also assume that $K$ is sufficiently large so that $\textbf{A}_j = \textbf{e}_0 = \textbf{C}_j$ for $j > K$ and so that $\sum_{j = 1}^{\infty} | \textbf{A}_j| < K$. 
	
	Inserting the definition \eqref{wxyac} of the single-row partition functions $W$ and $\widehat{W}$ into the second equation of \eqref{wabbcddw} yields
	\begin{flalign}
	\label{wwuwwv}
	\begin{aligned}
	& \Big| \widetilde{W}_{(w, u); \textbf{y}_{[1, M]}}' \big( \mathscr{A}_{[1, M]}, (\textbf{e}_0, \textbf{e}_0); \mathscr{C}_{[1, M]}, ( \textbf{I}, \textbf{J}) \boldsymbol{\mid} (t, r); \textbf{s}_{[1, M]} \big) \Big| \\
	& \quad \le  \displaystyle\sum_{\mathscr{K}} \displaystyle\sum_{\mathscr{J}} \displaystyle\sum_{\mathscr{I}} \displaystyle\prod_{m = 1}^M \Big| W_{w / y_m} (\textbf{A}_m, \textbf{I}_m; \textbf{K}_m , \textbf{I}_{m + 1} \boldsymbol{\mid} t, s_m) \widehat{W}_{u / y_m} (\textbf{K}_m, \textbf{J}_m; \textbf{C}_m, \textbf{J}_{m + 1}  \boldsymbol{\mid} r, s_m) \Big|,
	\end{aligned} 
	\end{flalign}
	
	\noindent where $\mathscr{I} = (\textbf{I}_1, \textbf{}_2, \ldots , \textbf{I}_{M + 1})$, $\mathscr{J} = (\textbf{J}_1, \textbf{J}_2, \ldots , \textbf{J}_{M + 1})$, and $\mathscr{K} = (\textbf{K}_1, \textbf{K}_2, \ldots , \textbf{K}_M)$ are each summed over all sequences of elements in $\{ 0, 1 \}^n$ such that $\textbf{I}_{M + 1} = \textbf{I}$, $\textbf{J}_{M + 1} = \textbf{J}$, and $\textbf{I}_1 = \textbf{e}_0 = \textbf{J}_1$. Since $\textbf{A}_m = \textbf{e}_0 = \textbf{C}_m$ for $m > K$, by arrow conservation any nonzero summand on the right side of \eqref{wwuwwv} must satisfy $|\textbf{I}_{m + 1}| \le |\textbf{I}_m|$ and $|\textbf{J}_{m + 1}| \ge |\textbf{J}_m|$. Thus, if $(\textbf{I}, \textbf{J}) \ne \big( \textbf{e}_0, \textbf{e}_{[1, n]} \big)$, then \eqref{wkj} implies that 
	\begin{flalign}
	\label{wwymwuym}
	\displaystyle\prod_{m = K + 1}^M \Big| W_{w / y_m} (\textbf{A}_m, \textbf{I}_m; \textbf{K}_m , \textbf{I}_{m + 1} \boldsymbol{\mid} t, s_m) \widehat{W}_{u / y_m} (\textbf{K}_m, \textbf{J}_m; \textbf{C}_m, \textbf{J}_{m + 1}  \boldsymbol{\mid} r, s_m) \Big| < (1 - \varepsilon)^{M - K},
	\end{flalign}
	
	\noindent for any $\mathscr{I}$, $\mathscr{J}$, and $\mathscr{K}$ supported by the sum on the right side of \eqref{wwuwwv}.
	
	Furthermore, arrow conservation implies that $\sum_{j = 1}^M |\textbf{K}_j| \le \sum_{j = 1}^M |\textbf{A}_j| < nK$ for any $\mathscr{K}$ supported on the right side of \eqref{wwuwwv}, from which it follows that there are at most $2^n \binom{M}{nK} \le (2M)^{nK}$ choices of $\mathscr{K}$ for which the corresponding summand is nonzero. Given any fixed such $\mathscr{K}$, arrow conservation implies that there is at most one choice of $\mathscr{I}$ and $\mathscr{J}$ on which the right side of \eqref{wwuwwv} is supported. Inserting this and \eqref{wwymwuym} into \eqref{wwuwwv}, and further denoting 
	\begin{flalign*}
	\Xi = \displaystyle\max_{m \in [1, K]} \displaystyle\max \Big| W_{w / y_m} ( \textbf{I}_1, \textbf{J}_1; \textbf{I}_2, \textbf{J}_2 \boldsymbol{\mid} t, s_m) \widehat{W}_{u / y_m} ( \textbf{I}_1', \textbf{J}_1'; \textbf{I}_2', \textbf{J}_2' \boldsymbol{\mid} r, s_m) \Big|,
	\end{flalign*}
	
	\noindent where the second maximum is taken over $\textbf{I}_1, \textbf{I}_2, \textbf{I}_1', \textbf{I}_2', \textbf{J}_1, \textbf{J}_2, \textbf{J}_1', \textbf{J}_2' \in \{ 0, 1 \}^n$ yields if $(\textbf{I}, \textbf{J}) \ne \big( \textbf{e}_0, \textbf{e}_{[1, n]} \big)$ that
	\begin{flalign*}
	\Big| \widetilde{W}_{(w, u); \textbf{y}_{[1, M]}}' \big( \mathscr{A}_{[1, M]}, (\textbf{e}_0, \textbf{e}_0); \mathscr{C}_{[1, M]}, (\textbf{I}, \textbf{J}) \boldsymbol{\mid} (t, r); \textbf{s}_{[1, M]} \big) \Big| \le \Xi^K (2M)^{nK} (1 - \varepsilon)^{M - K}.
	\end{flalign*}
		
	\noindent This gives \eqref{limitmw} by letting $M$ tend to $\infty$, which as mentioned above implies the lemma.
\end{proof} 

Commutation relations between the operators $(\mathds{B}, \mathds{C})$ and $(\mathds{C}, \mathds{D})$ can also be derived as a result of \Cref{w2w}, but we will not state them here since we will not require them in full generality.

\section{Composition of \texorpdfstring{$\mathds{D}$}{} Operators}

\label{ProductD}

In this section we establish the following result for composing the operators $\mathds{D} (x; r)$ from \eqref{bd}, which will be useful for analyzing principal specializations of symmetric functions in \Cref{IdentitiesLambdaMu} below. 

\begin{prop}
	
	\label{dxrdxr}
	
	Fix complex numbers $x, r_1, r_2 \in \mathbb{C}$ and infinite sequences $\textbf{\emph{y}} = (y_1, y_2, \ldots )$ and $\textbf{\emph{s}} = (s_1, s_2, \ldots )$ of complex numbers. As operators on $\mathbb{V}$, we have 
	\begin{flalign}
	\label{xrd} 
	\mathds{D} (r_1^{-2} x; r_2) \mathds{D} (x; r_1) = \mathds{D} (x; r_1 r_2).
	\end{flalign} 
	
\end{prop} 

Composition results similar to \Cref{dxrdxr} also hold for the operators $\mathds{B}$ and $\mathds{C}$ from \eqref{bd}, but we will not pursue this here. 

We will deduce \Cref{dxrdxr} essentially by interpreting both sides of \eqref{xrd} as special cases of the rectangular partition functions given by \Cref{zxy1}, and then applying the first statement of the branching-type result \eqref{zabcdzabcd} for the latter quantities. To that end, we require the following partition function, which is for the vertex model consisting of horizontally adjacent rectangles whose rapidity parameters are in a geometric progression (as in \Cref{zxy1}).

\begin{definition}

\label{zxymk} 

Fix integers $K, L \ge 1$; a set of positive integers $\textbf{M} = (M_1, M_2, \ldots , M_K)$; a complex number $x \in \mathbb{C}$; and a set of complex numbers $\textbf{y}_{[1, K]} = (y_1, y_2, \ldots , y_K)$. For each $i \in [1, K]$, define the sequence $\textbf{w}^{(i)} = (q^{M_i - 1} y_i, q^{M_i - 2} y_i, \ldots , y_i)$, and further let
\begin{flalign*}
\textbf{x} = (x, qx, \ldots , q^{L - 1} x); \qquad \textbf{w} = \textbf{w}^{(1)} \cup \textbf{w}^{(2)} \cup \cdots \cup \textbf{w}^{(K)},
\end{flalign*}

\noindent where the union is ordered, so $\textbf{y} = (q^{M_1 - 1} y_1, q^{M_1 - 2} y_1, \ldots , y_1, \ldots , q^{M_K - 1} y_K, q^{M_K - 2} y_K, \ldots , y_K)$. 

Additionally, let $\mathfrak{B} = (b_1, b_2, \ldots , b_L)$ and $\mathfrak{D} = (d_1, d_2, \ldots , d_L)$ denote sequences of indices in $[0, n]$. Also let $\mathcal{A} = \big( \mathfrak{A}^{(1)}, \mathfrak{A}^{(2)}, \ldots , \mathfrak{A}^{(K)} \big)$ and $\mathcal{C} = \big( \mathfrak{C}^{(1)}, \mathfrak{C}^{(2)}, \ldots , \mathfrak{C}^{(K)} \big)$, where for each $i \in [1, K]$ the $\mathfrak{A}^{(i)} = \big( a_1^{(i)}, a_2^{(i)}, \ldots , a_{M_i}^{(i)} \big)$ and $\mathfrak{C}^{(i)} = \big( c_1^{(i)}, c_2^{(i)}, \ldots , c_{M_i}^{(i)} \big)$ are sequences of indices in $[0, n]$. Moreover define (where the unions below are again ordered)
\begin{flalign*}
\mathfrak{A} = \mathfrak{A}^{(1)} \cup \mathfrak{A}^{(2)} \cup \cdots \mathfrak{A}^{(K)}; \qquad \mathfrak{C} = \mathfrak{C}^{(1)} \cup \mathfrak{C}^{(2)} \cup \cdots \mathfrak{C}^{(K)}. 
\end{flalign*}

\noindent Then, recalling the notation of \Cref{zxy1}, denote the partition function
\begin{flalign*}
Z_{L; \textbf{M}} \big( \mathcal{A}, \mathfrak{B}; \mathcal{C}, \mathfrak{D} \boldsymbol{\mid} x; \textbf{y}_{[1, K]} \big) = Z (\mathfrak{A}, \mathfrak{B}; \mathfrak{C}, \mathfrak{D} \boldsymbol{\mid} \textbf{x}, \textbf{w}).
\end{flalign*}

\end{definition} 

Observe under that $Z_{L; \textbf{M}} \big( \mathcal{A}, \mathfrak{B}; \mathcal{C}, \mathfrak{D} \boldsymbol{\mid} x, \textbf{y}_{[1, K]} \big)$ is the partition function for the vertex model obtained by horizontally juxtaposing the $K$ rectangles corresponding to $Z_{x, y_i} \big( \mathfrak{A}^{(i)}, \mathfrak{B}; \mathfrak{C}^{(i)}, \mathfrak{D} \big)$ from \Cref{zxy1} (see the left side of \Cref{zxy}). We refer to \Cref{zxy2} for a depiction. 

\begin{figure}

	\begin{center}

		\begin{tikzpicture}[
		>=stealth,
		auto,
		style={
			scale = 1.2
		}
		]
		
		\draw[->, thick, red] (-1, 0) -- (0, 0);
		\draw[->, thick, blue] (-1, 1) -- (0, 1);
		\draw[->, thick, red] (-1, 2) -- (0, 2);
		
		\draw[->, thick, green] (8, 0) -- (9, 0);
		\draw[->, thick, orange] (8, 1) -- (9, 1);
		\draw[->, thick, blue] (8, 2) -- (9, 2);
		
		\draw[->, thick, red] (0, -1) -- (0, 0);
		\draw[->, thick, green] (1, -1) -- (1, 0);
		\draw[->, thick, blue] (2, -1) -- (2, 0);
		\draw[->, thick, red] (3, -1) -- (3, 0);
		\draw[->, thick, red] (4, -1) -- (4, 0);
		\draw[->, thick, green] (5, -1) -- (5, 0);
		\draw[->, thick, blue] (6, -1) -- (6, 0);
		\draw[->, thick, blue] (7, -1) -- (7, 0);
		\draw[->, thick, orange] (8, -1) -- (8, 0);
		
		\draw[->, thick, red] (0, 2) -- (0, 3);
		\draw[->, thick, blue] (1, 2) -- (1, 3);
		\draw[->, thick, red] (2, 2) -- (2, 3);
		\draw[->, thick, red] (3, 2) -- (3, 3);
		\draw[->, thick, blue] (4, 2) -- (4, 3);
		\draw[->, thick, green] (5, 2) -- (5, 3);
		\draw[->, thick, red] (6, 2) -- (6, 3);
		\draw[->, thick, blue] (7, 2) -- (7, 3);
		\draw[->, thick, red] (8, 2) -- (8, 3);
		
		\draw[-, dashed] (0, 0) -- (8, 0);
		\draw[-, dashed] (0, 1) -- (8, 1);
		\draw[-, dashed] (0, 2) -- (8, 2);
		
		\draw[-, dashed] (0, 0) -- (0, 2);
		\draw[-, dashed] (1, 0) -- (1, 2);
		\draw[-, dashed] (2, 0) -- (2, 2);
		\draw[-, dashed] (3, 0) -- (3, 2);
		\draw[-, dashed] (4, 0) -- (4, 2);
		\draw[-, dashed] (5, 0) -- (5, 2);
		\draw[-, dashed] (6, 0) -- (6, 2);
		\draw[-, dashed] (7, 0) -- (7, 2);
		\draw[-, dashed] (8, 0) -- (8, 2);

		\draw[] (0, -1) circle[radius = 0] node[below, scale = .75]{$q^3 y_1$};
		\draw[] (1, -1) circle[radius = 0] node[below, scale = .75]{$q^2 y_1$};
		\draw[] (2, -1) circle[radius = 0] node[below = 2, scale = .75]{$q y_1$};
		\draw[] (3, -1) circle[radius = 0] node[below = 2, scale = .75]{$y_1$};
		\draw[] (4, -1) circle[radius = 0] node[below = 2, scale = .75]{$q y_2$};
		\draw[] (5, -1) circle[radius = 0] node[below = 2, scale = .75]{$y_2$};
		\draw[] (6, -1) circle[radius = 0] node[below = 2, scale = .75]{$y_3$};
		\draw[] (7, -1) circle[radius = 0] node[below = 2, scale = .75]{$q y_3$};
		\draw[] (8, -1) circle[radius = 0] node[below, scale = .75]{$q^2 y_3$};
		
		\draw[] (-1, 0) circle[radius = 0] node[left, scale = .75]{$x$};
		\draw[] (-1, 1) circle[radius = 0] node[left, scale = .75]{$q x$};
		\draw[] (-1, 2) circle[radius = 0] node[left, scale = .75]{$q^2 x$};
		
		\draw[] (-.6, 0) circle[radius = 0] node[above, scale = .65]{$b_1$};
		\draw[] (-.6, 1) circle[radius = 0] node[above, scale = .65]{$b_2$};
		\draw[] (-.6, 2) circle[radius = 0] node[above, scale = .65]{$b_3$};
		
		\draw[] (8.56, 0) circle[radius = 0] node[above, scale = .65]{$d_1$};
		\draw[] (8.6, 1) circle[radius = 0] node[above, scale = .65]{$d_2$};
		\draw[] (8.6, 2) circle[radius = 0] node[above, scale = .65]{$d_3$};
		
		\draw[] (0, -.6) circle[radius = 0] node[left, scale = .65]{$a_1^{(1)}$};
		\draw[] (1, -.6) circle[radius = 0] node[left, scale = .65]{$a_2^{(1)}$};
		\draw[] (2, -.6) circle[radius = 0] node[left, scale = .65]{$a_3^{(1)}$};
		\draw[] (3, -.6) circle[radius = 0] node[left, scale = .65]{$a_4^{(1)}$};
		\draw[] (4, -.6) circle[radius = 0] node[left, scale = .65]{$a_1^{(2)}$};
		\draw[] (5, -.6) circle[radius = 0] node[left, scale = .65]{$a_2^{(2)}$};
		\draw[] (6, -.6) circle[radius = 0] node[left, scale = .65]{$a_1^{(3)}$};
		\draw[] (7, -.6) circle[radius = 0] node[left, scale = .65]{$a_2^{(3)}$};
		\draw[] (8, -.6) circle[radius = 0] node[left, scale = .65]{$a_3^{(3)}$};
		
		\draw[] (0, 2.5) circle[radius = 0] node[left, scale = .65]{$c_1^{(1)}$};
		\draw[] (1, 2.5) circle[radius = 0] node[left, scale = .65]{$c_2^{(1)}$};
		\draw[] (2, 2.5) circle[radius = 0] node[left, scale = .65]{$c_3^{(1)}$};
		\draw[] (3, 2.5) circle[radius = 0] node[left, scale = .65]{$c_4^{(1)}$};
		\draw[] (4, 2.5) circle[radius = 0] node[left, scale = .65]{$c_1^{(2)}$};
		\draw[] (5, 2.5) circle[radius = 0] node[left, scale = .65]{$c_2^{(2)}$};
		\draw[] (6, 2.5) circle[radius = 0] node[left, scale = .65]{$c_1^{(3)}$};
		\draw[] (7, 2.5) circle[radius = 0] node[left, scale = .65]{$c_2^{(3)}$};
		\draw[] (8, 2.5) circle[radius = 0] node[left, scale = .65]{$c_3^{(3)}$};
		
		\draw[] (0, -1.5) -- (0, -1.625) -- (3, -1.625) -- (3, -1.5); 
		\draw[] (4, -1.5) -- (4, -1.625) -- (5, -1.625) -- (5, -1.5); 
		\draw[] (6, -1.5) -- (6, -1.625) -- (8, -1.625) -- (8, -1.5); 
		\draw[] (-1.625, 0) -- (-1.75, 0) -- (-1.75, 2) -- (-1.625, 2);
		
		\draw[thick, gray!50!white, dotted] (3.5, -1.25) -- (3.5, 3.25);
		\draw[thick, gray!50!white, dotted] (5.5, -1.25) -- (5.5, 3.25);
		
		\draw[] (1.5, -1.675) circle[radius = 0] node[below, scale = .7]{$M_1$};
		\draw[] (4.5, -1.675) circle[radius = 0] node[below, scale = .7]{$M_2$};
		\draw[] (7, -1.675) circle[radius = 0] node[below, scale = .7]{$M_3$};
		\draw[] (-1.75, 1) circle[radius = 0] node[left, scale = .7]{$L$};
		
		\end{tikzpicture}
		
	\end{center}
	
	\caption{\label{zxy2} Shown above is a diagrammatic interpretation for $Z_{L; \textbf{M}} (\mathcal{A}, \mathfrak{B}; \mathcal{C}, \mathfrak{D} \boldsymbol{\mid} x; \textbf{y}_{[1, K]} )$.}
\end{figure}

We now have the following lemma that expresses $\big\langle \mathscr{C} \big| \mathds{D} (x; r) \big| \mathscr{A} \big\rangle$ as linear combinations of the partition functions from \Cref{zxymk}, if $r = q^{-L / 2}$ and $s_i = q^{-M_i / 2}$ for each $i \in [1, K]$ and some integers $L, M_1, M_2, \ldots , M_K \ge 0$. 

\begin{lem} 

\label{dsac}

Let $K, L, \textbf{\emph{M}}, x, \textbf{\emph{y}}_{[1, K]}, \textbf{\emph{x}}, \textbf{\emph{w}}$ be as in \Cref{zxymk}, and define the sequences of indices $\mathfrak{B} = (0, 0, \ldots , 0) = \mathfrak{D}$ (where $0$ appears with multiplicity $L$). Let $\mathscr{A} = (\textbf{\emph{A}}_1, \textbf{\emph{A}}_2, \ldots )$ and $\mathscr{C} = (\textbf{\emph{C}}_1, \textbf{\emph{C}}_2, \ldots )$ denote finitary sequences of elements in $\{ 0, 1 \}^n$; assume that $\textbf{\emph{A}}_i = \textbf{\emph{e}}_0 = \textbf{\emph{C}}_i$ for each $i > K$, and that $|\textbf{\emph{A}}_i|, |\textbf{\emph{C}}_i| \le M_i$ for each $i \in [1, K]$. Set $\textbf{\emph{A}}_i' = \big( M_i - |\textbf{\emph{A}}_i|, \textbf{\emph{A}}_i \big)$ and $\textbf{\emph{C}}_i' = \big( M_i - |\textbf{\emph{C}}_i|, \textbf{\emph{C}}_i \big)$ for each $i \in [1, K]$; fix sequences of indices $\mathfrak{C}^{(i)} \in \mathcal{M} \big( \textbf{\emph{C}}_i' \big)$; and let $\mathcal{C} = \big( \mathfrak{C}^{(1)}, \mathfrak{C}^{(2)}, \ldots , \mathfrak{C}^{(K)} \big)$. If $r = q^{-L / 2}$ and $s_i = q^{-M_i / 2}$ for each $i \in [1, K]$, then we have 
\begin{flalign*}
\big\langle \mathscr{A} \big| \mathds{D} (x; r) \big| \mathscr{C} \big\rangle = \displaystyle\sum_{\mathfrak{A}^{(i)} \in \mathcal{M} (\textbf{\emph{A}}_i')} Z_{L; \textbf{\emph{M}}} \big( \mathcal{A}, \mathfrak{B}; \mathcal{C}, \mathfrak{D} \boldsymbol{\mid} x; \textbf{\emph{y}}_{[1, K]} \big) \displaystyle\prod_{i = 1}^K \Bigg( q^{\inv (\mathfrak{A}^{(i)}) - \inv(\mathfrak{C}^{(i)})} \displaystyle\prod_{j = 0}^n \displaystyle\frac{(q; q)_{m_j (\mathfrak{A}^{(i)})}}{(q; q)_{m_j (\mathfrak{C}^{(i)})}} \Bigg),
\end{flalign*}

\noindent where $\mathcal{A} = \big( \mathfrak{A}^{(1)}, \mathfrak{A}^{(2)}, \ldots , \mathfrak{A}^{(K)} \big)$, and each $\mathfrak{A}^{(i)}$ is summed over $\mathcal{M} \big( \textbf{\emph{A}}_i' \big)$.

\end{lem} 

\begin{proof}
	
	By \eqref{twxyinner} and \eqref{wxyaclimit}, we have that 
	\begin{flalign}
	\label{acdw}
	\big\langle \mathscr{A} \big| \mathds{D} (x; r) \big| \mathscr{C} \big\rangle = \displaystyle\sum_{\mathscr{J}} \displaystyle\prod_{i = 1}^{\infty} W_{x / y_i} \big( \textbf{A}_i, \textbf{J}_i; \textbf{C}_i, \textbf{J}_{i + 1} \boldsymbol{\mid} r, s_i) = \displaystyle\sum_{\mathscr{J}} \displaystyle\prod_{i = 1}^K W_{x / y_i} \big( \textbf{A}_i, \textbf{J}_i; \textbf{C}_i, \textbf{J}_{i + 1} \boldsymbol{\mid} r, s_i),
	\end{flalign}
	
	\noindent where in the second term $\mathscr{J} = (\textbf{J}_0, \textbf{J}_1, \ldots)$ is summed over all infinite sequences of elements in $\{ 0, 1 \}^n$ such that $\textbf{J}_0 = \textbf{e}_0 = \textbf{J}_i$ for sufficiently large $i$, and in the third $\mathscr{J} = (\textbf{J}_0, \textbf{J}_1, \ldots \textbf{J}_{K + 1})$ is summed over all $K$-term sequences of elements in $\{ 0, 1\}^n$ such that $\textbf{J}_0 = \textbf{e}_0 = \textbf{J}_{K + 1}$. Here, the second equality follows from the fact that $W_{x / y_i} (\textbf{e}_0, \textbf{e}_0; \textbf{e}_0, \textbf{e}_0 \boldsymbol{\mid} r, s_i) = 1$ and $\textbf{A}_i = \textbf{C}_i = \textbf{J}_i = \textbf{e}_0$ for $i > K$ in the second term of \eqref{acdw} (where the latter holds by arrow conservation, since $\textbf{A}_i = \textbf{e}_0 = \textbf{C}_i$ for $i > K$ and $\textbf{J}_i = \textbf{e}_0$ for sufficiently large $i$). By arrow conservation, the right side of \eqref{acdw} is supported on a single sequence, which we also call $\mathscr{J} = (\textbf{J}_0, \textbf{J}_1, \ldots , \textbf{J}_{K + 1})$. 

	Then, setting $\textbf{J}_i' = \big( L - |\textbf{J}_i|, \textbf{J}_i \big) \in \mathbb{Z}_{\ge 0}^{n + 1}$ and applying \Cref{wabcdsxz} gives
	\begin{flalign}
	\label{adc2} 
	\big\langle \mathscr{A} \big| \mathds{D} (x; r) \big| \mathscr{C} \big\rangle = \displaystyle\prod_{i = 1}^K \mathcal{R}_{x, y_i}^{(1; n)} \big( \textbf{A}_i', \textbf{J}_i'; \textbf{C}_i', \textbf{J}_{i + 1}' \boldsymbol{\mid} r, s_i).
	\end{flalign} 
	
	\noindent By \Cref{zabcd2}, \Cref{zabcdcd}, and \Cref{rijkh} (see also \eqref{rabcdrabcd}), for any $\mathfrak{J}^{(i + 1)} \in \mathcal{M} (\textbf{J}_{i + 1}')$ we have that  
	\begin{flalign*}
	\mathcal{R}_{x, y_i} (\textbf{A}_i', \textbf{J}_i'; \textbf{C}_i', \textbf{J}_{i + 1}' \boldsymbol{\mid} r, s_i) = \displaystyle\sum_{\mathfrak{A}^{(i)} \in \mathcal{M} (\textbf{A}_i')} & \displaystyle\sum_{\mathfrak{J}^{(i)} \in \mathcal{M} (\textbf{J}_i')} q^{\inv (\mathfrak{A}^{(i)}) - \inv (\mathfrak{C}^{(i)}) + \inv (\overleftarrow{\mathfrak{J}}^{(i)}) - \inv(\overleftarrow{\mathfrak{J}}^{(i + 1)})} \\
	& \times Z_{x, y_i} \big( \mathfrak{A}^{(i)}, \mathfrak{J}^{(i)}; \mathfrak{C}^{(i)}, \mathfrak{J}^{(i + 1)} \big) \cdot \displaystyle\prod_{j = 0}^n \displaystyle\frac{(q; q)_{m_j (\mathfrak{A}^{(i)})}}{(q; q)_{m_j (\mathfrak{C}^{(i)})}}.
	\end{flalign*}
	
	\noindent Inserting this into \eqref{adc2} and using the fact that $\inv \big( \mathfrak{J}^{(0)} \big) = 0 = \inv \big( \mathfrak{J}^{(K + 1)} \big)$ (as $\textbf{J}_0 = \textbf{e}_0 = \textbf{J}_{K + 1}$), we deduce 
	\begin{flalign}
	\label{adc3}
	\begin{aligned} 
	\big\langle \mathscr{A} \big| \mathds{D} (x; r) \big| \mathscr{C} \big\rangle = \displaystyle\sum_{\mathfrak{A}^{(i)} \in \mathcal{M} (\textbf{A}_i')} \displaystyle\sum_{\mathfrak{J}^{(i)} \in \mathcal{M} (\textbf{J}_i')} & Z_{x, y_i} \big( \mathfrak{A}^{(i)}, \mathfrak{J}^{(i)}; \mathfrak{C}^{(i)}, \mathfrak{J}^{(i + 1)}) \\
	& \times \displaystyle\prod_{i = 1}^K \Bigg( q^{\inv (\mathfrak{A}^{(i)}) - \inv (\mathfrak{C}^{(i)})} \displaystyle\prod_{j = 0}^n \displaystyle\frac{(q; q)_{m_j (\mathfrak{A}^{(i)})}}{(q; q)_{m_j (\mathfrak{C}^{(i)})}} \Bigg).
	\end{aligned}
	\end{flalign} 
	
	Now, applying the second statement of \eqref{zabcdzabcd} $K - 1$ times, for the $h$ there in $\big\{ M_1, M_1 + M_2, \ldots , \sum_{i = 1}^{K - 1} M_i \big\}$ (see \Cref{zxy2}) gives
	\begin{flalign*}
	\displaystyle\sum_{\mathfrak{J}^{(i)} \in \mathcal{M} (\textbf{J}_i')} \displaystyle\prod_{i = 1}^K Z_{x, y_i} \big( \mathfrak{A}^{(i)}, \mathfrak{J}^{(i)}; \mathfrak{C}^{(i)}, \mathfrak{J}^{(i + 1)} \boldsymbol{\mid} r, s_i) = Z_{L; \textbf{M}} \big( \mathcal{A}, \mathfrak{B}; \mathcal{C}, \mathfrak{D} \boldsymbol{\mid} x; \textbf{y}_{[1, K]} \big),
	\end{flalign*}
	
	\noindent which upon insertion into \eqref{adc3} yields the lemma. 
\end{proof}

Now we can establish \Cref{dxrdxr}. 

\begin{proof}[Proof of \Cref{dxrdxr}]
	
	For any finitary $\mathscr{A} = (\textbf{A}_1, \textbf{A}_2, \ldots )$ and $\mathscr{C} = (\textbf{C}_1, \textbf{C}_2, \ldots)$, it suffices to show 
	\begin{flalign*}
	 \big\langle \mathscr{C} \big| \mathds{D} (x; r_1 r_2) \big| \mathscr{A} \big\rangle = \big\langle \mathscr{C} \big| \mathds{D} (r_1^{-2} x; r_2) \mathds{D} (x; r_1) \big| \mathscr{A} \big\rangle.
	\end{flalign*}
	
	\noindent Thus, since 
	\begin{flalign*}
	\big\langle \mathscr{C} \big| \mathds{D} (r_1^{-2} x; r_2) \mathds{D} (x; r_1) \big| \mathscr{A} \big\rangle = \displaystyle\sum_{\mathscr{H}} \big\langle \mathscr{C} \big| \mathds{D} (r_1^{-2} x; r_2) \big| \mathscr{H} \big\rangle \big\langle \mathscr{H} \big| \mathds{D} (x; r_1) \big| \mathscr{A} \big\rangle, 
	\end{flalign*}
	
	\noindent where we sum over all finitary $\mathscr{H} = (\textbf{H}_1, \textbf{H}_2, \ldots )$, it suffices to show 
	\begin{flalign}
	\label{cdda} 
	\big\langle \mathscr{C} \big| \mathds{D} (x; r_1 r_2) \big| \mathscr{A} \big\rangle = \displaystyle\sum_{\mathscr{H}} \big\langle \mathscr{C} \big| \mathds{D} (r_1^{-2} x; r_2) \big| \mathscr{H} \big\rangle \big\langle \mathscr{H} \big| \mathds{D} (x; r_1) \big| \mathscr{A} \big\rangle, 
	\end{flalign}
	
	Now let $K > 0$ be such that $\textbf{A}_j = \textbf{e}_0 = \textbf{C}_j$ for $j > K$. It then follows by arrow conservation that any finitary sequence $\mathscr{H}$ supported by the sum on the right side of \eqref{cdda} must also satisfy $\textbf{H}_j = \textbf{e}_0$ for $j > K$. Thus, from \eqref{bd}, \eqref{twxyinner}, and \Cref{mwabcd}, together with the explicit form \Cref{wabcdrsxy} of the $W$ weights, both sides of \eqref{cdda} are rational functions in the variables $r_1, r_2, s_1, s_2, \ldots , s_K$. So, to establish \eqref{cdda}, we may assume in what follows there exist integers $L_1, L_2, M_1, M_2, \ldots , M_K \ge 1$ such that $r_j = q^{-L_j / 2}$ for each $j \in \{ 1, 2 \}$; such that $s_j = q^{-M_j / 2}$ for each $j \in \{ 1, 2, \ldots , K \}$; and such that $M_j \ge |\textbf{A}_j|, |\textbf{C}_j|$ holds for each $j \in \{ 1, 2, \ldots , K \}$. 
	
	In this scenario, \Cref{dsac} applies. In particular, fix some finitary sequence $\mathscr{H} = (\textbf{H}_1, \textbf{H}_2, \ldots )$ of elements in $\{ 0, 1 \}^n$ as above; set $\textbf{X}_j' = \big( M_j - |\textbf{X}_j|, \textbf{X}_j \big)$ for each index $X \in \{ A, C, H \}$; and, letting $0^j = (0, 0, \ldots , 0) \in \mathbb{Z}^j$ for any integer $j \ge 1$, set
	\begin{flalign*}
	\mathfrak{B}^{(1)} = 0^{L_1} = \mathfrak{D}^{(1)}; \qquad \mathfrak{B}^{(2)} = 0^{L_2} = \mathfrak{D}^{(2)}; \qquad \mathfrak{B} = 0^{L_1 + L_2} = \mathfrak{D}.
	\end{flalign*}
	
	\noindent Then, since $r_1 = q^{-L_1 / 2}$ and $r_2 = q^{-L_2 / 2}$, \Cref{dsac} gives 
	\begin{flalign}
	\label{cdacdhhda} 
	\begin{aligned}
	\big\langle \mathscr{C} \big| \mathds{D} (x; r_1 r_2) \big| \mathscr{A} \big\rangle & = \displaystyle\sum_{\mathfrak{A}^{(i)} \in \mathcal{M} (\textbf{A}_i')} Z_{L_1 + L_2; \textbf{M}} \big( \mathcal{A}, \mathfrak{B}; \mathcal{C}, \mathfrak{D} \boldsymbol{\mid} x; \textbf{y}_{[1, K]} \big) \\
	& \qquad \qquad \times \displaystyle\prod_{i = 1}^K \Bigg( q^{\inv (\mathfrak{A}^{(i)}) - \inv(\mathfrak{C}^{(i)})} \displaystyle\prod_{j = 0}^n \displaystyle\frac{(q; q)_{m_j (\mathfrak{A}^{(i)})}}{(q; q)_{m_j (\mathfrak{C}^{(i)})}} \Bigg); \\
	\big\langle \mathscr{C} \big| \mathds{D} (r_1^{-2} x; r_2) \big| \mathscr{H} \big\rangle & = \displaystyle\sum_{\mathfrak{H}^{(i)} \in \mathcal{M} (\textbf{H}_i')} Z_{L_2; \textbf{M}} \big( \mathcal{H}, \mathfrak{B}; \mathcal{C}, \mathfrak{D} \boldsymbol{\mid} q^{L_1} x; \textbf{y}_{[1, K]} \big) \\
	& \qquad \qquad \times \displaystyle\prod_{i = 1}^K \Bigg( q^{\inv (\mathfrak{H}^{(i)}) - \inv(\mathfrak{C}^{(i)})} \displaystyle\prod_{j = 0}^n \displaystyle\frac{(q; q)_{m_j (\mathfrak{H}^{(i)})}}{(q; q)_{m_j (\mathfrak{C}^{(i)})}} \Bigg); \\
	\big\langle \mathscr{H} \big| \mathds{D} (x; r_1) \big| \mathscr{A} \big\rangle & = \displaystyle\sum_{\mathfrak{A}^{(i)} \in \mathcal{M} (\textbf{A}_i')} Z_{L_1; \textbf{M}} \big( \mathcal{A}, \mathfrak{B}; \mathcal{H}, \mathfrak{D} \boldsymbol{\mid} x; \textbf{y}_{[1, K]} \big) \\
	& \qquad \qquad \times \displaystyle\prod_{i = 1}^K \Bigg( q^{\inv (\mathfrak{A}^{(i)}) - \inv(\mathfrak{H}^{(i)})} \displaystyle\prod_{j = 0}^n \displaystyle\frac{(q; q)_{m_j (\mathfrak{A}^{(i)})}}{(q; q)_{m_j (\mathfrak{H}^{(i)})}} \Bigg).
	\end{aligned} 
	\end{flalign} 
	
	\noindent In first and third statements of \eqref{cdacdhhda}, we sum over all sets $\mathcal{A} = \big( \mathfrak{A}^{(1)}, \mathfrak{A}^{(2)}, \ldots , \mathfrak{A}^{(K)} \big)$ with $\mathfrak{A}^{(j)} \in \mathcal{M} (\textbf{A}_j')$ for each $j$, and in the second we sum over all sets $\mathcal{H} = \big( \mathfrak{H}^{(1)}, \mathfrak{H}^{(2)}, \ldots , \mathfrak{H}^{(K)} \big)$ with $\mathfrak{H}^{(j)} \in \mathcal{M} (\textbf{H}_j')$ for each $j$. Moreover, in the first and second statements of \eqref{cdacdhhda}, we have fixed a set $\mathcal{C} = \big( \mathfrak{C}^{(1)}, \mathfrak{C}^{(2)}, \ldots , \mathfrak{C}^{(K)} \big)$ with $\mathfrak{C}^{(j)} \in \mathcal{M} (\textbf{C}_j')$ for each $j$. Similarly, in the third statement of \eqref{cdacdhhda}, we have fixed $\mathcal{H} = \big( \mathfrak{H}^{(1)}, \mathfrak{H}^{(2)}, \ldots , \mathfrak{H}^{(K)} \big)$ with $\mathfrak{H}^{(j)} \in \mathcal{M} (\textbf{H}_j')$ for each $j$. 

	Next, the first statement of \eqref{zabcdzabcd} (applied with the $h$ there equal to $L_1$ here) gives
	\begin{flalign*}
	Z_{L_1 + L_2; \textbf{M}} \big( \mathcal{A}, \mathfrak{B}; \mathcal{C}, \mathfrak{D} \boldsymbol{\mid} x; \textbf{y}_{[1, K]} \big) = \displaystyle\sum_{\mathcal{H}} & Z_{L_1; \textbf{M}} \big( \mathcal{A}, \mathfrak{B}; \mathcal{H}, \mathfrak{D} \boldsymbol{\mid} x; \textbf{y}_{[1, K]} \big) \\
	& \times Z_{L_2; \textbf{M}} \big( \mathcal{H}, \mathfrak{B}; \mathcal{C}, \mathfrak{D} \boldsymbol{\mid} q^{L_1} x; \textbf{y}_{[1, K]} \big),
	\end{flalign*}  
	
	\noindent which together with \eqref{cdacdhhda} implies \eqref{cdda} and thus the proposition.
\end{proof}

\chapter{Functions and Identities}

\label{Functions1}

In this chapter we define the functions we will study in this text and establish symmetry, branching, and skew Cauchy identities for them.

\section{Symmetric Functions}

\label{Symmetric} 

In this section we define the symmetric functions we will study in this text as partition functions for vertex models with certain boundary conditions. Throughout this section, we fix infinite sequences $\textbf{y} = (y_1, y_2, \ldots )$ and $\textbf{s} = (s_1, s_2, \ldots )$ of complex numbers. 

Recall the vector spaces $\mathbb{V} = \mathbb{V}_{\textbf{y}; \textbf{s}}$ and $\mathbb{V}^* = \mathbb{V}_{\textbf{y}; \textbf{s}}^*$ from \Cref{OperatorRow}, which are spanned by basis vectors of the form $\big\{ | \mathscr{A} \rangle \big\}$ and $\big\{ \langle \mathscr{C} | \big\}$, respectively, where $\mathscr{A}$ and $\mathscr{C}$ range over all infinite, finitary sequences of elements in $\{ 0, 1 \}^n$. Let us introduce relabelings for these bases of $\mathbb{V}$ and $\mathbb{V}^*$ that will be more convenient for defining the symmetric functions below. 

To that end, recall the notions of signatures and signature sequences from \Cref{0Functions}, as well as the shift operator $\mathfrak{T}$ from \eqref{t}. In what follows, we will denote any signature $\lambda \in \Sign_{\ell}$ by $\lambda = (\lambda_1, \lambda_2, \ldots , \lambda_{\ell})$\index{0@$\lambda, \mu$; typical signatures or partitions} and any signature sequence $\boldsymbol{\lambda} \in \SeqSign_n$ by $\boldsymbol{\lambda} = \big( \lambda^{(1)}, \lambda^{(2)}, \ldots , \lambda^{(n)} \big)$,\index{0@$\boldsymbol{\lambda}, \boldsymbol{\mu}$; typical signature sequences} even when not explicitly mentioned. For any $\boldsymbol{\lambda} \in \SeqSign_n$, further recall from \Cref{0Functions} the infinite, finitary sequence $\mathscr{S} (\boldsymbol{\lambda}) = \big( \textbf{S}_1 (\boldsymbol{\lambda}), \textbf{S}_2 (\boldsymbol{\lambda}), \ldots \big)$\index{S@$\mathscr{S} (\boldsymbol{\lambda})$} of elements in $\{ 0, 1 \}^n$ defined as follows. For each $j \ge 1$, let $\textbf{S}_j = \textbf{S}_j (\boldsymbol{\lambda}) = (S_{1, j}, S_{2, j}, \ldots , S_{n, j}) \in \{ 0, 1 \}^n$, where $S_{i, j} = \textbf{1}_{j \in \mathfrak{T} (\lambda^{(i)})}$ for every $i \in [1, n]$. Then $\mathscr{S}$ induces a bijection between $\SeqSign_n$ and finitary sequences of elements in $\{ 0, 1 \}^n$.

Denoting $| \boldsymbol{\lambda} \rangle = \big| \mathscr{S} (\boldsymbol{\lambda}) \big\rangle$ and $\langle \boldsymbol{\lambda} | = \big\langle \mathscr{S} (\boldsymbol{\lambda}) \big|$\index{0@$\mid$$\boldsymbol{\lambda} \rangle$,$\langle \boldsymbol{\lambda}$$\mid$} for any $\boldsymbol{\lambda} \in \SeqSign_n$, these vectors constitute bases of $\mathbb{V}$ and $\mathbb{V}^*$, respectively. Given this relabeling, we can define the following functions; in the below, we recall the operators $\mathds{B}$, $\mathds{C}$, and $\mathds{D}$ from \eqref{bd}. 

\begin{definition} 
	
	\label{fgdefinition} 
	
	Fix integers $N \ge 1$ and $M \ge 0$; sequences of complex numbers $\textbf{x} = (x_1, x_2, \ldots , x_N)$ and $\textbf{r} = (r_1, r_2, \ldots , r_N)$; and infinite sequences of complex numbers $\textbf{y} = (y_1, y_2, \ldots )$ and $\textbf{s} = (s_1, s_2, \ldots )$. For any $\boldsymbol{\lambda}, \boldsymbol{\mu} \in \SeqSign_{n; M}$, define
	\begin{flalign*}
	G_{\boldsymbol{\lambda} / \boldsymbol{\mu}} (\textbf{x}; \textbf{r} \boldsymbol{\mid} \textbf{y}; \textbf{s} ) = \big\langle \boldsymbol{\lambda} \big| \mathds{D} (x_N; r_N) \cdots \mathds{D} (x_1; r_1) \big| \boldsymbol{\mu} \big\rangle. 
	\end{flalign*}
	
	\noindent We further abbreviate $G_{\boldsymbol{\lambda}} (\textbf{x}; \textbf{r} \boldsymbol{\mid} \textbf{y}; \textbf{s}) = G_{\boldsymbol{\lambda} / \boldsymbol{0}^M} (\textbf{x}; \textbf{r} \boldsymbol{\mid} \textbf{y}; \textbf{s})$.\index{G@$G_{\boldsymbol{\lambda} / \boldsymbol{\mu}} (\textbf{x}; \textbf{r} \boldsymbol{\mid} \textbf{y}; \textbf{s})$}\index{G@$G_{\boldsymbol{\lambda} / \boldsymbol{\mu}} (\textbf{x}; \textbf{r} \boldsymbol{\mid} \textbf{y}; \textbf{s})$!$G_{\boldsymbol{\lambda}} (\textbf{x}; \textbf{r} \boldsymbol{\mid} \textbf{y}; \textbf{s})$} 
	
	 Moreover, for any $\boldsymbol{\lambda} \in \Sign_{n; M + N}$ and $\boldsymbol{\mu} \in \SeqSign_{n; M}$, define 
	\begin{flalign*}
	& F_{\boldsymbol{\lambda} / \boldsymbol{\mu}} (\textbf{x}; \textbf{r} \boldsymbol{\mid} \textbf{y}; \textbf{s}) = \big\langle \boldsymbol{\mu} \big| \mathds{B} (x_N; r_N) \cdots \mathds{B} (x_1; r_1) \big| \boldsymbol{\lambda} \big\rangle; \\
	& H_{\boldsymbol{\lambda} / \boldsymbol{\mu}} (\textbf{x}; \textbf{r} \boldsymbol{\mid} \textbf{y}; \textbf{s}) = \big\langle \boldsymbol{\lambda} \big| \mathds{C} (x_N; r_N) \cdots \mathds{C} (x_1; r_1) \big| \boldsymbol{\mu} \big\rangle,
	\end{flalign*}
	
	\noindent We further abbreviate $F_{\boldsymbol{\lambda}} (\textbf{x}; \textbf{r} \boldsymbol{\mid} \textbf{y}; \textbf{s}) = F_{\boldsymbol{\lambda} / \boldsymbol{\varnothing}} (\textbf{x}; \textbf{r} \boldsymbol{\mid} \textbf{y}; \textbf{s})$ and $H_{\boldsymbol{\lambda}} (\textbf{x}; \textbf{r} \boldsymbol{\mid} \textbf{y}, \textbf{s}) = H_{\boldsymbol{\lambda} / \boldsymbol{\varnothing}} (\textbf{x}; \textbf{r} \boldsymbol{\mid} \textbf{y}; \textbf{s})$.\index{F@$F_{\boldsymbol{\lambda} / \boldsymbol{\mu}} (\textbf{x}; \textbf{r} \boldsymbol{\mid} \textbf{y}; \textbf{s})$}\index{H@$H_{\boldsymbol{\lambda} / \boldsymbol{\mu}} (\textbf{x}; \textbf{r} \boldsymbol{\mid} \textbf{y}; \textbf{s})$}\index{F@$F_{\boldsymbol{\lambda} / \boldsymbol{\mu}} (\textbf{x}; \textbf{r} \boldsymbol{\mid} \textbf{y}; \textbf{s})$!$F_{\boldsymbol{\lambda}} (\textbf{x}; \textbf{r} \boldsymbol{\mid} \textbf{y}; \textbf{s})$}\index{H@$H_{\boldsymbol{\lambda} / \boldsymbol{\mu}} (\textbf{x}; \textbf{r} \boldsymbol{\mid} \textbf{y}; \textbf{s})$!$H_{\boldsymbol{\lambda}} (\textbf{x}; \textbf{r} \boldsymbol{\mid} \textbf{y}; \textbf{s})$}

\end{definition}

\begin{rem}
	
	\label{fgmn}
	
	Although we will not pursue this here, one can define $G_{\boldsymbol{\lambda} / \boldsymbol{\mu}}, F_{\boldsymbol{\lambda} / \boldsymbol{\mu}}, H_{\boldsymbol{\lambda} / \boldsymbol{\mu}}$ for arbitrary sequences $\boldsymbol{\lambda}, \boldsymbol{\mu} \in \Sign_n$ (without stipulations on the lengths of their signatures) similarly to in \Cref{fgdefinition}. All identities to be shown in \Cref{IdentitiesLambdaMu} and \Cref{IdentitiesSum} below will still hold in this more general setting. Under this convention, it is quickly verified that $G_{\boldsymbol{\lambda} / \boldsymbol{\mu}} (\textbf{x}; \textbf{r} \boldsymbol{\mid} \textbf{y}; \textbf{s}) = 0$ for any $(\textbf{x}, \textbf{y}, \textbf{r}, \textbf{s})$ unless $\ell \big( \lambda^{(j)} \big) = \ell \big( \mu^{(j)} \big)$ for each $j \in [1, n]$. Similarly, if $\textbf{x} = (x_1, x_2, \ldots , x_N)$, then $F_{\boldsymbol{\lambda} / \boldsymbol{\mu}} (\textbf{x}; \textbf{r} \boldsymbol{\mid} \textbf{y}; \textbf{s}) = 0 = H_{\boldsymbol{\lambda} / \boldsymbol{\mu}} (\textbf{x}; \textbf{r} \boldsymbol{\mid} \textbf{y}; \textbf{s})$ unless $\ell \big( \lambda^{(j)} \big) = \ell \big( \mu^{(j)} \big) + N$ for each $j \in [1, n]$.  

\end{rem}

Let us describe diagrammatic interpretations of $G_{\boldsymbol{\lambda} / \boldsymbol{\mu}}$, $F_{\boldsymbol{\lambda} / \boldsymbol{\mu}}$, and $H_{\boldsymbol{\lambda} / \boldsymbol{\mu}}$ as partition functions for vertex models on the domain 
\begin{flalign}
\label{dn} 
\mathcal{D} = \mathcal{D}_N = \mathbb{Z}_{> 0} \times \{ 1, 2, \ldots , N \} \subset \mathbb{Z}_{> 0}^2.
\end{flalign}

The following definition provides notation for sets of path ensembles on $\mathcal{D}_N$ with certain types of boundary data.

\begin{definition}
	
	\label{pgpfph} 
	
	Fix $M, N \in \mathbb{Z}_{\ge 0}$ with $N \ge 1$, and $\boldsymbol{\lambda}, \boldsymbol{\mu} \in \SeqSign_n$. 
	
	If $\boldsymbol{\lambda}, \boldsymbol{\mu} \in \SeqSign_{n; M}$, then let $\mathfrak{P}_G (\boldsymbol{\lambda} / \boldsymbol{\mu}; N)$ denote the set of path ensembles on $\mathcal{D}$ satisfying the following two properties; see the left side of \Cref{fgpaths}.\index{P@$\mathfrak{P}_G (\boldsymbol{\lambda} / \boldsymbol{\mu}; N)$}
	
	\begin{enumerate}
		\item For every $c \in [1, n]$, one color $c$ arrow vertically enters $\mathcal{D}$ through $( \mathfrak{m}, 1)$ for each $\mathfrak{m} \in \mathfrak{T} \big( \mu^{(c)} \big)$.  
		\item For every $c \in [1, n]$, one color $c$ arrow vertically exits $\mathcal{D}$ through $(\mathfrak{l}, N)$ for each $\mathfrak{l} \in \mathfrak{T} \big( \lambda^{(c)} \big)$. 
	\end{enumerate}

	If $\boldsymbol{\lambda} \in \SeqSign_{n; M + N}$ and $\boldsymbol{\mu} \in \SeqSign_{n; M}$, then let $\mathfrak{P}_F (\boldsymbol{\lambda} / \boldsymbol{\mu})$ denote\footnote{Observe here that $N$ is fixed by $\boldsymbol{\lambda}$ and $\boldsymbol{\mu}$ and is therefore not incorporated in the notation for $\mathfrak{P}_F$.} the set of path ensembles on $\mathcal{D}$ satisfying the following three properties for each $c \in [1, n]$; see the middle of \Cref{fgpaths}. \index{P@$\mathfrak{P}_F (\boldsymbol{\lambda} / \boldsymbol{\mu})$}
	
	\begin{enumerate}
		\item For every $c \in [1, n]$, one color $c$ arrow vertically enters $\mathcal{D}$ through $(\mathfrak{l}, 1)$ for each $\mathfrak{l} \in \mathfrak{T} \big( \lambda^{(c)} \big)$. 
		\item For every $c \in [1, n]$, one color $c$ arrow vertically exits $\mathcal{D}$ through $(\mathfrak{m}, N)$ for each $\mathfrak{m} \in \mathfrak{T} \big( \mu^{(c)} \big)$. 
		\item For every $c \in [1, n]$, one color $c$ arrow horizontally exits $\mathcal{D}$ through $(\infty, j)$ for each $j \in [1, N]$.\footnote{This means that every edge connecting $(i, j)$ to $(i + 1, j)$ for sufficiently large $i$ contains an arrow of color $c$.}
	\end{enumerate}

	If $\boldsymbol{\lambda} \in \SeqSign_{n; M + N}$ and $\boldsymbol{\mu} \in \SeqSign_{n; M}$, then let $\mathfrak{P}_H (\boldsymbol{\lambda} / \boldsymbol{\mu})$ denote the set of path ensembles on $\mathcal{D}$ satisfying the following three properties; see the right side of \Cref{fgpaths}. \index{P@$\mathfrak{P}_H (\boldsymbol{\lambda} / \boldsymbol{\mu})$}
	
	\begin{enumerate}
		\item For every $c \in [1, n]$, one color $c$ arrow vertically enters $\mathcal{D}$ through $(\mathfrak{m}, 1)$ for each $\mathfrak{m} \in \mathfrak{T} \big( \mu^{(c)} \big)$. 
		\item For every $c \in [1, n]$, one color $c$ arrow vertically exits $\mathcal{D}$ through $( \mathfrak{l}, N)$ for each $\mathfrak{l} \in \mathfrak{T} \big( \lambda^{(c)} \big)$. 
		\item For every $c \in [1, n]$, one color $c$ arrow horizontally enters $\mathcal{D}$ through $(1, j)$ for each $j \in [1, N]$.
	\end{enumerate}
	
\end{definition}

\begin{figure}
	
	\begin{center}

		\begin{tikzpicture}[
		>=stealth,
		scale = .5
		]
		
		\draw[->, thick, red] (.95, 0) -- (.95, 1);
		\draw[->, thick, red] (2.95, 3) -- (2.95, 4); 
		\draw[->, thick, red] (2.95, 0) -- (2.95, 1);
		\draw[->, thick, red] (5.95, 3) -- (5.95,4);
		\draw[->, thick, blue] (1.95, 0) -- (1.95, 1);
		\draw[->, thick, blue] (5, 3) -- (5, 4);
		\draw[->, thick, blue] (4, 0) -- (4, 1);
		\draw[->, thick, blue] (6.05, 3) -- (6.05, 4);
		\draw[->, thick, green] (1.05, 0) -- (1.05, 1); 
		\draw[->, thick, green] (3.05, 3) -- (3.05, 4);
		\draw[->, thick, green] (2.05, 0) -- (2.05, 1);
		\draw[->, thick, green] (7, 3) -- (7, 4);
		
		\draw[ultra thick, gray, dashed] (1, 1) -- (1, 3);
		\draw[ultra thick, gray, dashed] (2, 1) -- (2, 3);
		\draw[ultra thick, gray, dashed] (3, 1) -- (3, 3);
		\draw[ultra thick, gray, dashed] (4, 1) -- (4, 3);
		\draw[ultra thick, gray, dashed] (5, 1) -- (5, 3);
		\draw[ultra thick, gray, dashed] (6, 1) -- (6, 3);
		\draw[ultra thick, gray, dashed] (7, 1) -- (7, 3);
		
		\draw[ultra thick, gray, dashed] (1, 1) -- (7, 1);
		\draw[ultra thick, gray, dashed] (1, 2) -- (7, 2);
		\draw[ultra thick, gray, dashed] (1, 3) -- (7, 3);
		
		\draw[->, very thick] (0, 0) -- (0, 4.5);
		\draw[->, very thick] (0, 0) -- (7.5, 0);

		\draw[]  (4, -.15) circle [radius = 0] node[below, scale = .55]{$(y_4; s_4)$};
		\draw[]  (-.15, 2) circle [radius = 0] node[left, scale = .55]{$(x_2; r_2)$};
		
		\draw[-] (1, -.5) -- (1, -1) -- (7, -1) -- (7, -.5);
		\draw[-] (1, 4.5) -- (1, 5) -- (7, 5) -- (7, 4.5);
		
		\draw[] (4, -1) circle[radius = 0] node[below, scale = .7]{$\mathscr{S} (\boldsymbol{\mu})$};
		\draw[] (4, 5) circle[radius = 0] node[above, scale = .7]{$\mathscr{S} (\boldsymbol{\lambda})$};
		\draw[] (4, 6) circle[radius = 0] node[above, scale = .8]{$G_{\boldsymbol{\lambda} / \boldsymbol{\mu}}$};

		\draw[->, very thick] (10, 0) -- (10, 4.5);
		\draw[->, very thick] (10, 0) -- (17.5, 0);
		
		\draw[->, thick, red] (10.95, 0) -- (10.95, 1);
		\draw[->, thick, red] (15.95, 0) -- (15.95, 1);
		\draw[->, thick, red] (14.95, 0) -- (14.95, 1);
		\draw[->, thick, red] (12.95, 3) -- (12.95, 4); 
		\draw[->, thick, red] (12.9, 0) -- (12.9, 1);
		
		\draw[->, thick, blue] (11.95, 0) -- (11.95, 1);
		\draw[->, thick, blue] (15, 3) -- (15, 4);
		\draw[->, thick, blue] (13, 0) -- (13, 1);
		\draw[->, thick, blue] (13.95, 0) -- (13.95, 1);
		\draw[->, thick, blue] (16.05, 0) -- (16.05, 1);
		
		\draw[->, thick, green] (11.05, 0) -- (11.05, 1); 
		\draw[->, thick, green] (13.05, 3) -- (13.05, 4);
		\draw[->, thick, green] (12.05, 0) -- (12.05, 1);
		\draw[->, thick, green] (13.1, 0) -- (13.1, 1);
		\draw[->, thick, green] (15.05, 0) -- (15.05, 1);
		
		\draw[->, thick, red] (16.1, .9) -- (17, .9);
		\draw[->, thick, blue] (16.1, 1) -- (17, 1);
		\draw[->, thick, green] (16.1, 1.1) -- (17, 1.1);
		
		\draw[->, thick, red] (16.1, 1.9) -- (17, 1.9);
		\draw[->, thick, blue] (16.1, 2) -- (17, 2);
		\draw[->, thick, green] (16.1, 2.1) -- (17, 2.1);
		
		\draw[->, thick, red] (16.1, 2.9) -- (17, 2.9);
		\draw[->, thick, blue] (16.1, 3) -- (17, 3);
		\draw[->, thick, green] (16.1, 3.1) -- (17, 3.1);
		
		\draw[ultra thick, gray, dashed] (11, 1) -- (11, 3);
		\draw[ultra thick, gray, dashed] (12, 1) -- (12, 3);
		\draw[ultra thick, gray, dashed] (13, 1) -- (13, 3);
		\draw[ultra thick, gray, dashed] (14, 1) -- (14, 3);
		\draw[ultra thick, gray, dashed] (15, 1) -- (15, 3);
		\draw[ultra thick, gray, dashed] (16, 1) -- (16, 3);
		
		\draw[ultra thick, gray, dashed] (11, 1) -- (16, 1);
		\draw[ultra thick, gray, dashed] (11, 2) -- (16, 2);
		\draw[ultra thick, gray, dashed] (11, 3) -- (16, 3);
		
		\draw[]  (14, -.15) circle [radius = 0] node[below, scale = .55]{$(y_4; s_4)$};
		\draw[]  (9.85, 2) circle [radius = 0] node[left, scale = .55]{$(x_2; r_2)$};
		
		\draw[-] (11, 4.5) -- (11, 5) -- (16, 5) -- (16, 4.5);
		\draw[-] (11, -.5) -- (11, -1) -- (16, -1) -- (16, -.5);
		
		\draw[] (13.5, -1) circle[radius = 0] node[below, scale = .7]{$\mathscr{S} (\boldsymbol{\lambda})$};
		\draw[] (13.5, 5) circle[radius = 0] node[above, scale = .7]{$\mathscr{S} (\boldsymbol{\mu})$};
		\draw[] (13.5, 6) circle[radius = 0] node[above, scale = .8]{$F_{\boldsymbol{\lambda} / \boldsymbol{\mu}}$};
		
		\draw[->, very thick] (20, 0) -- (20, 4.5);
		\draw[->, very thick] (20, 0) -- (27.5, 0);
		
		\draw[->, thick, red] (20.95, 3) -- (20.95, 4);
		\draw[->, thick, red] (25.95, 3) -- (25.95, 4);
		\draw[->, thick, red] (24.95, 3) -- (24.95, 4);
		\draw[->, thick, red] (22.95, 0) -- (22.95, 1); 
		\draw[->, thick, red] (22.9, 3) -- (22.9, 4);
		
		\draw[->, thick, blue] (21.95, 3) -- (21.95, 4);
		\draw[->, thick, blue] (25, 0) -- (25, 1);
		\draw[->, thick, blue] (23, 3) -- (23, 4);
		\draw[->, thick, blue] (23.95, 3) -- (23.95, 4);
		\draw[->, thick, blue] (26.05, 3) -- (26.05, 4);
		
		\draw[->, thick, green] (21.05, 3) -- (21.05, 4); 
		\draw[->, thick, green] (23.05, 0) -- (23.05, 1);
		\draw[->, thick, green] (22.05, 3) -- (22.05, 4);
		\draw[->, thick, green] (23.1, 3) -- (23.1, 4);
		\draw[->, thick, green] (25.05, 3) -- (25.05, 4);
		
		\draw[->, thick, red] (20.1, .9) -- (21, .9);
		\draw[->, thick, blue] (20.1, 1) -- (21, 1);
		\draw[->, thick, green] (20.1, 1.1) -- (21, 1.1);
		
		\draw[->, thick, red] (20.1, 1.9) -- (21, 1.9);
		\draw[->, thick, blue] (20.1, 2) -- (21, 2);
		\draw[->, thick, green] (20.1, 2.1) -- (21, 2.1);
		
		\draw[->, thick, red] (20.1, 2.9) -- (21, 2.9);
		\draw[->, thick, blue] (20.1, 3) -- (21, 3);
		\draw[->, thick, green] (20.1, 3.1) -- (21, 3.1);
		
		\draw[ultra thick, gray, dashed] (21, 1) -- (21, 3);
		\draw[ultra thick, gray, dashed] (22, 1) -- (22, 3);
		\draw[ultra thick, gray, dashed] (23, 1) -- (23, 3);
		\draw[ultra thick, gray, dashed] (24, 1) -- (24, 3);
		\draw[ultra thick, gray, dashed] (25, 1) -- (25, 3);
		\draw[ultra thick, gray, dashed] (26, 1) -- (26, 3);
		
		\draw[ultra thick, gray, dashed] (21, 1) -- (26, 1);
		\draw[ultra thick, gray, dashed] (21, 2) -- (26, 2);
		\draw[ultra thick, gray, dashed] (21, 3) -- (26, 3);
		
		\draw[]  (24, -.15) circle [radius = 0] node[below, scale = .55]{$(y_4; s_4)$};
		\draw[]  (19.85, 2) circle [radius = 0] node[left, scale = .55]{$(x_2; r_2)$};
		
		\draw[-] (21, 4.5) -- (21, 5) -- (26, 5) -- (26, 4.5);
		\draw[-] (21, -.5) -- (21, -1) -- (26, -1) -- (26, -.5);
		
		\draw[] (23.5, -1) circle[radius = 0] node[below, scale = .7]{$\mathscr{S} (\boldsymbol{\mu})$};
		\draw[] (23.5, 5) circle[radius = 0] node[above, scale = .7]{$\mathscr{S} (\boldsymbol{\lambda})$};
		\draw[] (23.5, 6) circle[radius = 0] node[above, scale = .8]{$H_{\boldsymbol{\lambda} / \boldsymbol{\mu}}$};
				
		\end{tikzpicture}
		
	\end{center}
	
	\caption{\label{fgpaths} To the left, middle, and right are the vertex models $\mathfrak{P}_G (\boldsymbol{\lambda} / \boldsymbol{\mu}; 3)$, $\mathfrak{P}_F (\boldsymbol{\lambda} / \boldsymbol{\mu})$, and $\mathfrak{P}_H (\boldsymbol{\lambda} / \boldsymbol{\mu})$, respectively. On the left, we have $\boldsymbol{\lambda} = \big( (4, 2), (4, 4), (5, 2) \big)$ and $\boldsymbol{\mu} = \big( (1, 0), (2, 1), (0, 0) \big)$, and on the middle and right we have $\boldsymbol{\lambda} = \big( (2, 2, 1, 0), (2, 1, 1, 1), (1, 0, 0, 0) \big)$ and $\boldsymbol{\mu} = \big( (2), (4), (2) \big)$. Here, red is color $1$, blue is color $2$, and green is color $3$. } 
	
\end{figure}

Now, fix finite sequences of complex numbers $\textbf{x} = (x_1, x_2, \ldots , x_N)$ and $\textbf{r} = (r_1, r_2, \ldots , r_N)$. We define the weight (with respect to either the vertex weights $W_z$ from \Cref{wabcdrsxy} or $\widehat{W}_z$ from \eqref{wabcd2}) of any fused ensemble $\mathcal{E}$ on $\mathcal{D}$ to be
\begin{flalign}
\label{weighte} 
\begin{aligned}
& W (\mathcal{E} \boldsymbol{\mid} \textbf{x}; \textbf{r} \boldsymbol{\mid} \textbf{y}; \textbf{s}) = \displaystyle\prod_{(i, j) \in \mathcal{D}} W_{x_j / y_i} \big( \textbf{A} (i, j), \textbf{B} (i, j); \textbf{C} (i, j), \textbf{D} (i, j) \boldsymbol{\mid} r_j, s_i \big); \\
& \widehat{W} (\mathcal{E} \boldsymbol{\mid} \textbf{x}; \textbf{r} \boldsymbol{\mid} \textbf{y}; \textbf{s}) = \displaystyle\prod_{(i, j) \in \mathcal{D}} \widehat{W}_{x_j / y_i} \big( \textbf{A} (i, j), \textbf{B} (i, j); \textbf{C} (i, j), \textbf{D} (i, j) \boldsymbol{\mid} r_j, s_i \big), 
\end{aligned} 
\end{flalign}
\index{W@$W (\mathcal{E} \boldsymbol{\mid} \textbf{x}; \textbf{r} \boldsymbol{\mid} \textbf{y}; \textbf{s})$; ensemble weight}
\index{W@$\widehat{W} (\mathcal{E} \boldsymbol{\mid} \textbf{x}; \textbf{r} \boldsymbol{\mid} \textbf{y}; \textbf{s})$; normalized ensemble weight}

\noindent where $\big( \textbf{A} (v), \textbf{B} (v); \textbf{C} (v), \textbf{D} (v) \big)$ denotes the arrow configuration under $\mathcal{E}$ at any vertex $v \in \mathcal{D}$, and we must assume that the above infinite products converge. Observe in particular that if $\mathcal{E} \in \mathfrak{P}_G (\boldsymbol{\lambda} / \boldsymbol{\mu}; N)$ or $\mathcal{E} \in \mathfrak{P}_H (\boldsymbol{\lambda} / \boldsymbol{\mu})$ for some $\boldsymbol{\lambda}, \boldsymbol{\mu} \in \SeqSign_n$, then the product defining $W (\mathcal{E})$ converges since $\big( \textbf{A} (v), \textbf{B} (v); \textbf{C} (v), \textbf{D} (v) \big) = (\textbf{e}_0, \textbf{e}_0; \textbf{e}_0, \textbf{e}_0)$ for all but finitely many $v \in \mathcal{D}$, and this arrow configuration has weight $1$ under $W$ by \eqref{wabcd01nequation}. Similarly, if $\mathcal{E} \in \mathfrak{P}_F (\boldsymbol{\lambda} / \boldsymbol{\mu})$ for some $\boldsymbol{\lambda}, \boldsymbol{\mu} \in \SeqSign_n$, then the product defining $\widehat{W} (\mathcal{E})$ converges since $\big( \textbf{A} (v), \textbf{B} (v); \textbf{C} (v), \textbf{D} (v) \big) = \big( \textbf{e}_0, \textbf{e}_{[1, n]}; \textbf{e}_0, \textbf{e}_{[1, n]} \big)$ for all but finitely many $v \in \mathcal{D}$, and this arrow configuration has weight $1$ under $\widehat{W}$ by \eqref{wz1}.

Then, under the above notation, we have that 
\begin{flalign}
\label{gfhe}
\begin{aligned} 
G_{\boldsymbol{\lambda} / \boldsymbol{\mu}} (\textbf{x}; \textbf{r} \boldsymbol{\mid} \textbf{y}; \textbf{s}) = \displaystyle\sum_{\mathcal{E} \in \mathfrak{P}_G (\boldsymbol{\lambda} / \boldsymbol{\mu}; N)} & W (\mathcal{E} \boldsymbol{\mid} \textbf{x}; \textbf{r} \boldsymbol{\mid} \textbf{y}; \textbf{s}); \quad H_{\boldsymbol{\lambda} / \boldsymbol{\mu}} (\textbf{x}; \textbf{r} \boldsymbol{\mid} \textbf{y}; \textbf{s}) = \displaystyle\sum_{\mathcal{E} \in \mathfrak{P}_H (\boldsymbol{\lambda} / \boldsymbol{\mu})} W (\mathcal{E} \boldsymbol{\mid} \textbf{x}; \textbf{r} \boldsymbol{\mid} \textbf{y}; \textbf{s}); \\
& F_{\boldsymbol{\lambda} / \boldsymbol{\mu}} (\textbf{x}; \textbf{r} \boldsymbol{\mid} \textbf{y}; \textbf{s}) = \displaystyle\sum_{\mathcal{E} \in \mathfrak{P}_F (\boldsymbol{\lambda} / \boldsymbol{\mu})} \widehat{W} (\mathcal{E} \boldsymbol{\mid} \textbf{x}; \textbf{r} \boldsymbol{\mid} \textbf{y}; \textbf{s}).
\end{aligned} 
\end{flalign} 

This indicates that $G_{\boldsymbol{\lambda} / \boldsymbol{\mu}}$, $F_{\boldsymbol{\lambda} / \boldsymbol{\mu}}$, and $H_{\boldsymbol{\lambda} / \boldsymbol{\mu}}$ can be interpreted as partition functions for the vertex models $\mathfrak{P}_G (\boldsymbol{\lambda} / \boldsymbol{\mu}; N)$ (under the weight $W_z$), $\mathfrak{P}_F (\boldsymbol{\lambda} / \boldsymbol{\mu})$ (under $\widehat{W}_z$), and $\mathfrak{P}_H (\boldsymbol{\lambda} / \boldsymbol{\mu})$ (under $W_z$), respectively. In each of these vertex models, the rapidity pair for the $i$-th column from the left is $(y_i; s_i)$ and that for the $j$-th row from the bottom is $(x_j; r_j)$; we refer to \Cref{fgpaths} for depictions in all three cases.

\section{Symmetry, Branching, and Principal Specializations}

\label{IdentitiesLambdaMu} 

In this section we establish branching identities, symmetry relations, and results concerning principal specializations for the $F, G, H$ functions from \Cref{fgdefinition}. We begin with the former. 

\begin{prop}
	
	\label{fghbranching}
	
	Fix integers $n, K, L, M \ge 1$; finite sequences of complex numbers $\textbf{\emph{x}}' = (x_1, \ldots , x_K)$, $\textbf{\emph{x}}'' = (x_{K + 1}, \ldots , x_{K + L})$, $\textbf{\emph{r}}' = (r_1, \ldots,  r_K)$, and $\textbf{\emph{r}}'' = (r_{K + 1}, \ldots , r_{K + L})$; and infinite sequences of complex numbers $\textbf{\emph{y}} = (y_1, y_2, \ldots )$ and $\textbf{\emph{s}} = (s_1, s_2, \ldots )$. Define $\textbf{\emph{x}} = \textbf{\emph{x}}' \cup \textbf{\emph{x}}'' = (x_1, x_2, \ldots , x_{K + L})$ and $\textbf{\emph{r}} = \textbf{\emph{r}}' \cup \textbf{\emph{r}}'' = (r_1, r_2, \ldots , r_{K + L})$. For any $n$-tuples of signatures $\boldsymbol{\lambda}, \boldsymbol{\mu} \in \Sign_{n; M}$, we have
	\begin{flalign}
	\label{sumg}
	& \displaystyle\sum_{\boldsymbol{\nu} \in \SeqSign_{n; M}} G_{\boldsymbol{\lambda} / \boldsymbol{\nu}} (\textbf{\emph{x}}''; \textbf{\emph{r}}'' \boldsymbol{\mid} \textbf{\emph{y}}; \textbf{\emph{s}}) G_{\boldsymbol{\nu} / \boldsymbol{\mu}} (\textbf{\emph{x}}'; \textbf{\emph{r}}' \boldsymbol{\mid} \textbf{\emph{y}}; \textbf{\emph{s}}) = G_{\boldsymbol{\lambda} / \boldsymbol{\mu}} (\textbf{\emph{x}}; \textbf{\emph{r}} \boldsymbol{\mid} \textbf{\emph{y}}; \textbf{\emph{s}}),
	\end{flalign} 
	
	\noindent Moreover, for any $n$-tuples of signatures $\boldsymbol{\lambda} \in \Sign_{n; M + K + L}$ and $\boldsymbol{\mu} \in \Sign_{n; M}$, we have that
	\begin{flalign}
	\label{sumf} 
	\begin{aligned}
	& \displaystyle\sum_{\boldsymbol{\nu} \in \SeqSign_{n; M + L}} F_{\boldsymbol{\lambda} / \boldsymbol{\nu}} (\textbf{\emph{x}}'; \textbf{\emph{r}}' \boldsymbol{\mid} \textbf{\emph{y}}; \textbf{\emph{s}}) F_{\boldsymbol{\nu} / \boldsymbol{\mu}} (\textbf{\emph{x}}''; \textbf{\emph{r}}'' \boldsymbol{\mid} \textbf{\emph{y}}; \textbf{\emph{s}}) = F_{\boldsymbol{\lambda} / \boldsymbol{\mu}} (\textbf{\emph{x}}; \textbf{\emph{r}} \boldsymbol{\mid} \textbf{\emph{y}}; \textbf{\emph{s}}); \\
	& \displaystyle\sum_{\boldsymbol{\nu} \in \SeqSign_{n; M + K}} H_{\boldsymbol{\lambda} / \boldsymbol{\nu}} (\textbf{\emph{x}}''; \textbf{\emph{r}}'' \boldsymbol{\mid} \textbf{\emph{y}}; \textbf{\emph{s}}) H_{\boldsymbol{\nu} / \boldsymbol{\mu}} (\textbf{\emph{x}}'; \textbf{\emph{r}}' \boldsymbol{\mid} \textbf{\emph{y}}; \textbf{\emph{s}}) = H_{\boldsymbol{\lambda} / \boldsymbol{\mu}} (\textbf{\emph{x}}; \textbf{\emph{r}} \boldsymbol{\mid} \textbf{\emph{y}}; \textbf{\emph{s}}).
	\end{aligned}
	\end{flalign}
	
\end{prop}

\begin{proof}

	By \Cref{fgdefinition}, we have that
	\begin{flalign*}
	G_{\boldsymbol{\lambda} / \boldsymbol{\mu}} & (\textbf{x}; \textbf{r} \boldsymbol{\mid} \textbf{y}; \textbf{s}) \\
	& = \big\langle \boldsymbol{\lambda} \boldsymbol{\mid} \mathds{D} (x_{K + L}; r_{K + L}) \cdots \mathds{D} (x_1; r_1) \big| \boldsymbol{\mu} \big\rangle \\
	& = \displaystyle\sum_{\boldsymbol{\nu} \in \SeqSign_n} \big\langle \boldsymbol{\lambda} \big| \mathds{D} (x_{K + L}; r_{K + L}) \cdots \mathds{D} (x_{K + 1}; r_{K + 1}) \big| \boldsymbol{\nu} \big\rangle \big\langle \boldsymbol{\nu} \big| \mathds{D} (x_K; r_K) \cdots \mathds{D} (x_1; r_1) \big| \boldsymbol{\mu} \big\rangle \\
	& = \displaystyle\sum_{\boldsymbol{\nu} \in \SeqSign_n} G_{\boldsymbol{\lambda} / \boldsymbol{\nu}} (\textbf{x}''; \textbf{r}'' \boldsymbol{\mid} \textbf{y}; \textbf{s}) G_{\boldsymbol{\nu} / \boldsymbol{\mu}} (\textbf{x}'; \textbf{r}' \boldsymbol{\mid} \textbf{y}; \textbf{s}),
	\end{flalign*}
	
	\noindent which establishes \eqref{sumg} (by using the last statement of \Cref{fgmn} to restrict the sum over $\nu \in \SeqSign_{n; M}$). The proof of \eqref{sumf} is entirely analogous and is therefore omitted. 
\end{proof}

We next have the following proposition, which shows that $G$ is symmetric in $\textbf{x}$ and $\textbf{r}$. It also shows that $F$ and $H$ are not exactly symmetric in those variables, but are up to explicit factors; this lack of exact symmetry is a consequence of the fact that $W_z \big( \textbf{e}_{[1, n]}, \textbf{e}_{[1, n]}; \textbf{e}_{[1, n]}, \textbf{e}_{[1, n]} \boldsymbol{\mid} r, s \big) \ne 1$.

\begin{prop} 

\label{gxfxsigma} 

Fix integers $n, N \ge 1$; finite sequences of complex numbers $\textbf{\emph{x}} = (x_1, x_2, \ldots , x_N)$ and $\textbf{\emph{r}} = (r_1, r_2, \ldots,  r_N)$; and infinite sequences of complex numbers $\textbf{\emph{y}} = (y_1, y_2, \ldots )$ and $\textbf{\emph{s}} = (s_1, s_2, \ldots )$. For any $n$-tuples $\boldsymbol{\lambda}, \boldsymbol{\mu} \in \SeqSign_n$ of signatures and permutation $\sigma \in \mathfrak{S}_N$, we have that 
\begin{flalign}
\label{sigmafg} 
\begin{aligned} 
& G_{\boldsymbol{\lambda} / \boldsymbol{\mu}} \big( \sigma (\textbf{\emph{x}}); \sigma (\textbf{\emph{r}}) \boldsymbol{\mid} \textbf{\emph{y}}; \textbf{\emph{s}} \big) = G_{\boldsymbol{\lambda} / \boldsymbol{\mu}} (\textbf{\emph{x}}; \textbf{\emph{r}} \boldsymbol{\mid} \textbf{\emph{y}}; \textbf{\emph{s}}); \\
& F_{\boldsymbol{\lambda} / \boldsymbol{\mu}} \big( \sigma (\textbf{\emph{x}}); \sigma (\textbf{\emph{r}}) \boldsymbol{\mid} \textbf{\emph{y}}; \textbf{\emph{s}}\big) = F_{\boldsymbol{\lambda} / \boldsymbol{\mu}} (\textbf{\emph{x}}; \textbf{\emph{r}} \boldsymbol{\mid} \textbf{\emph{y}}; \textbf{\emph{s}}) \displaystyle\prod_{\substack{1 \le i < j \le N \\ \sigma (i) > \sigma (j) }} \displaystyle\frac{(r_j^2 x_i x_j^{-1}; q)_n}{(r_i^2 x_i^{-1} x_j; q)_n} \left( \displaystyle\frac{r_i^2 x_j}{r_j^2 x_i} \right)^n; \\
& H_{\boldsymbol{\lambda} / \boldsymbol{\mu}} \big( \sigma (\textbf{\emph{x}}); \sigma (\textbf{\emph{r}}) \boldsymbol{\mid} \textbf{\emph{y}}; \textbf{\emph{s}}\big) = H_{\boldsymbol{\lambda} / \boldsymbol{\mu}} (\textbf{\emph{x}}; \textbf{\emph{r}} \boldsymbol{\mid} \textbf{\emph{y}}; \textbf{\emph{s}}) \displaystyle\prod_{\substack{1 \le i < j \le N \\ \sigma (i) > \sigma (j) }} \displaystyle\frac{(r_i^2 x_i^{-1} x_j; q)_n}{(r_j^2 x_i x_j^{-1}; q)_n} \left( \displaystyle\frac{r_j^2 x_i}{r_i^2 x_j} \right)^n.
\end{aligned}
\end{flalign} 

\end{prop} 

\begin{proof} 
	
	Since $G_{\boldsymbol{\lambda} / \boldsymbol{\mu}} (\textbf{x}; \textbf{r} \boldsymbol{\mid} \textbf{y}; \textbf{s}) = \big\langle \boldsymbol{\lambda} \big| \mathds{D} (x_N; r_N) \cdots \mathds{D} (x_1; r_1) \big| \boldsymbol{\mu} \big\rangle$, the first statement of \eqref{sigmafg} follows from the third statement of \eqref{bbdd12limit}. So, it remains to establish the second and third ones there; since their proofs are very similar, we only detail the former.
	
	To that end, since $\mathfrak{S}_N$ is generated by the transpositions $\{ \mathfrak{s}_i \}$ for $i \in [1, N - 1]$ interchanging $(i, i + 1)$, we may assume that $\sigma = \mathfrak{s}_k$ for some $k \in [1, N - 1]$. In this case, \Cref{fgdefinition} and the first statement of \eqref{bbdd12limit} together yield 
	\begin{flalign*}
	F_{\boldsymbol{\lambda} / \boldsymbol{\mu}} \big( & \mathfrak{s}_k (\textbf{x}); \mathfrak{s}_k (\textbf{r}) \boldsymbol{\mid} \textbf{y}; \textbf{s}\big) \\
	& = \big\langle \boldsymbol{\mu} \big| \mathds{B} (x_N; r_N) \cdots \mathds{B} (x_{k + 2}; r_{k + 2}) \mathds{B} (x_k; r_k) \mathds{B} (x_{k + 1}; r_{k + 1}) \mathds{B} (x_{k - 1}; r_{k - 1}) \cdots \mathds{B} (x_1; r_1) \big| \boldsymbol{\lambda} \big\rangle \\
	& = \left( \displaystyle\frac{r_k^2 x_{k + 1}}{r_{k + 1}^2 x_k} \right)^n \displaystyle\frac{(r_{k + 1}^2 x_k x_{k + 1}^{-1}; q)_n}{(r_k^2 x_k^{-1} x_{k + 1}; q)_n} \big\langle \boldsymbol{\mu} \big| \mathds{B} (x_N; r_N) \cdots \mathds{B} (x_1; r_1) \big| \boldsymbol{\lambda} \big\rangle \\
	& = \left( \displaystyle\frac{r_k^2 x_{k + 1}}{r_{k + 1}^2 x_k} \right)^n \displaystyle\frac{(r_{k + 1}^2 x_k x_{k + 1}^{-1}; q)_n}{(r_k^2 x_k^{-1} x_{k + 1}; q)_n} F_{\boldsymbol{\lambda} / \boldsymbol{\mu}} ( \textbf{x}; \textbf{r} \boldsymbol{\mid} \textbf{y}; \textbf{s}),
	\end{flalign*} 
	
	\noindent which gives the second statement of \eqref{sigmafg}. 
\end{proof} 

The following proposition indicates that, if each entry of $\textbf{r}$ is a negative half-integer power of $q$, then $G_{\boldsymbol{\lambda} / \boldsymbol{\mu}} (\textbf{x}; \textbf{r} \boldsymbol{\mid} \textbf{y}; \textbf{s})$ may be recovered from a function of the form $G_{\boldsymbol{\lambda} / \boldsymbol{\mu}} (\textbf{w}; \textbf{q}^{-1 / 2} \boldsymbol{\mid} \textbf{y}; \textbf{s})$, where $\textbf{q}^{-1 / 2} = (q^{-1 / 2}, q^{-1 / 2}, \ldots , q^{-1 / 2})$. Here, $\textbf{w}$ will be a union of \emph{principal specializations}, that is, geometric progressions of the form $(x_i, qx_i, \ldots , q^{L_i - 1} x_i)$. 

\begin{prop}
	
	\label{gq}
	
	Fix integers $M \ge 0$ and $N \ge 1$; finite sequences of complex numbers $\textbf{\emph{x}} = (x_1, x_2, \ldots , x_N)$ and $\textbf{\emph{r}} = (r_1, r_2, \ldots , r_N)$; and infinite sequences of complex numbers $\textbf{\emph{s}} = (s_1, s_2, \ldots )$ and $\textbf{\emph{y}} = (y_1, y_2, \ldots)$. Assume for each $i \in [1, N]$ that there exists an integer $L_i \in \mathbb{Z}_{\ge 1}$ such that $r_i = q^{-L_i / 2}$, and define $\textbf{\emph{w}} = \textbf{\emph{w}}^{(1)} \cup \textbf{\emph{w}}^{(2)} \cup \cdots \cup \textbf{\emph{w}}^{(N)}$, where $\textbf{\emph{w}}^{(i)} = (x_i, qx_i, \ldots , q^{L_i - 1} x_i)$. Further set $\textbf{\emph{q}}^{-1 / 2} = (q^{-1 / 2}, q^{-1 / 2}, \ldots , q^{-1 / 2})$, where $q^{-1 / 2}$ appears with multiplicity $N$. Then, for any signature sequences $\boldsymbol{\lambda}, \boldsymbol{\mu} \in \SeqSign_{n; M}$, we have that $G_{\boldsymbol{\lambda} / \boldsymbol{\mu}} (\textbf{\emph{x}}; \textbf{\emph{r}} \boldsymbol{\mid} \textbf{\emph{y}}; \textbf{\emph{s}}) = G_{\boldsymbol{\lambda} / \boldsymbol{\mu}} (\textbf{\emph{w}}; \textbf{\emph{q}}^{-1 / 2} \boldsymbol{\mid} \textbf{\emph{y}}; \textbf{\emph{s}})$. 
	
\end{prop}

\begin{proof}
	
	It suffices to establish the proposition for $N = 1$, as the result for general $N$ follows from this together with the branching relation \eqref{sumg}. Thus, we may assume in what follows that 
	\begin{flalign*} 
	\textbf{x} = (x); \qquad \textbf{r} = (r); \qquad r = q^{-L / 2}; \qquad \textbf{w} = (x, qx, \ldots , q^{L - 1} x),
	\end{flalign*} 
	
	\noindent for some complex numbers $x, r \in \mathbb{C}$ and integer $L \in \mathbb{Z}_{\ge 1}$. Then, \Cref{fgdefinition} gives
	\begin{flalign}
	\label{glambdamuxr}
	\begin{aligned}
	G_{\boldsymbol{\lambda} / \boldsymbol{\mu}} (\textbf{x}; \textbf{r} \boldsymbol{\mid} \textbf{y}; \textbf{s}) & = \big\langle \boldsymbol{\lambda} \big| \mathds{D} (x; r) \big| \boldsymbol{\mu} \big\rangle; \\
	 G_{\boldsymbol{\lambda} / \boldsymbol{\mu}} (\textbf{w}; \textbf{q}^{-1 / 2} \boldsymbol{\mid} \textbf{y}; \textbf{s}) & = \big\langle \boldsymbol{\lambda} \big| \mathds{D} (q^{L - 1} x; q^{-1 / 2}) \mathds{D} (q^{L - 2} x; q^{-1 / 2}) \cdots \mathds{D} (x; q^{-1 / 2}) \big| \boldsymbol{\mu} \big\rangle.
	\end{aligned} 
	\end{flalign}
	 
	 \noindent Moreover, \Cref{dxrdxr} yields by induction on $L$ that  
	 \begin{flalign*} 
	 \mathds{D} (q^{L - 1} x; q^{-1 / 2}) \mathds{D} (q^{L - 2} x; q^{-1 / 2}) \cdots \mathds{D} (x; q^{-1 / 2}) = \mathds{D} (x; q^{-L / 2}) = \mathds{D} (x; r),
	\end{flalign*} 
	
	\noindent which together with with \eqref{glambdamuxr} implies the proposition.
\end{proof}

Results similar to \Cref{gq} can also be derived for the functions $F_{\boldsymbol{\lambda} / \boldsymbol{\mu}}$ and $H_{\boldsymbol{\lambda} / \boldsymbol{\mu}}$, but we will not pursue this here.

\section{Cauchy Identities}

\label{IdentitiesSum}

In this section we establish various Cauchy identities for the $F$ and $G$ functions. We begin with the following proposition, which provides the most general skew-Cauchy identity. 

\begin{thm} 
	
\label{fgidentity} 

Fix integers $L, M \ge 0$ and $N \ge 1$; finite sequences $\textbf{\emph{u}} = (u_1, u_2, \ldots , u_N)$, $\textbf{\emph{r}} = (r_1, r_2, \ldots , r_N)$, $\textbf{\emph{w}} = (w_1, w_2, \ldots , w_M)$, and $\textbf{\emph{t}} = (t_1, t_2, \ldots , t_M)$ of complex numbers; and infinite sequences $\textbf{\emph{y}} = (y_1, y_2, \ldots )$ and $\textbf{\emph{s}} = (s_1, s_2, \ldots )$ of complex numbers. Assume that there exists an integer $K > 1$ such that
\begin{flalign}
\label{x1x2bdestimate2}
\displaystyle\sup_{k > K} \displaystyle\max_{\substack{1 \le i \le M \\ 1 \le j \le N}}	\displaystyle\max_{\substack{a, b \in [0, n] \\ (a, b) \ne (n, 0)}} \Bigg| s_k^{2a + 2b - 2n} \displaystyle\frac{(s_k^2 u_j y_k^{-1}; q)_n (u_j y_k^{-1}; q)_a}{(u_j y_k^{-1}; q)_n (s_k^2 u_j y_k^{-1}; q)_a} \displaystyle\frac{(w_i y_k^{-1}; q)_b}{(s_k^2 w_i y_k^{-1}; q)_b} \Bigg| < 1.
\end{flalign}

\noindent Then for any $\boldsymbol{\mu} \in \Sign_{n; L + N}$ and $\boldsymbol{\nu} \in \SeqSign_{n; L}$, we have that 
\begin{flalign}
\label{fgsum1}
\begin{aligned} 
\displaystyle\sum_{\boldsymbol{\lambda} \in \SeqSign_{n; L + N}} & F_{\boldsymbol{\lambda} / \boldsymbol{\nu}} (\textbf{\emph{u}}; \textbf{\emph{r}} \boldsymbol{\mid} \textbf{\emph{y}}; \textbf{\emph{s}}) G_{\boldsymbol{\lambda} / \boldsymbol{\mu}} (\textbf{\emph{w}}; \textbf{\emph{t}} \boldsymbol{\mid} \textbf{\emph{y}}; \textbf{\emph{s}}) \\
& \qquad = \displaystyle\prod_{i = 1}^M \displaystyle\prod_{j = 1}^N \displaystyle\frac{(t_i^2 u_j w_i^{-1}; q)_n}{t_i^{2n} (u_j w_i^{-1}; q)_n} \displaystyle\sum_{\boldsymbol{\kappa} \in \SeqSign_{n; L}} G_{\boldsymbol{\nu} / \boldsymbol{\kappa}} (\textbf{\emph{w}}; \textbf{\emph{t}} \boldsymbol{\mid} \textbf{\emph{y}}; \textbf{\emph{s}}) F_{\boldsymbol{\mu} / \boldsymbol{\kappa}} (\textbf{\emph{u}}; \textbf{\emph{r}} \boldsymbol{\mid} \textbf{\emph{y}}; \textbf{\emph{s}}).
\end{aligned} 
\end{flalign}

\end{thm}

\begin{proof}
	
	By \Cref{fgdefinition}, we have that 
	\begin{flalign*}
	\displaystyle\sum_{\boldsymbol{\lambda} \in \SeqSign_{n; L + N}} F_{\boldsymbol{\lambda} / \boldsymbol{\nu}} & (\textbf{u}; \textbf{r} \boldsymbol{\mid} \textbf{y}; \textbf{s}) G_{\boldsymbol{\lambda} / \boldsymbol{\mu}} (\textbf{w}; \textbf{t} \boldsymbol{\mid} \textbf{y}; \textbf{s}) \\
	 & = \displaystyle\sum_{\boldsymbol{\lambda} \in \SeqSign_{n; L + N}} \big\langle \boldsymbol{\nu} \big| \mathds{B} (u_N; r_N) \cdots \mathds{B} (u_1; r_1) \big| \boldsymbol{\lambda} \big\rangle \big\langle \boldsymbol{\lambda} \big| \mathds{D} (w_M; t_M) \cdots \mathds{D} (w_1; t_1) \big| \boldsymbol{\mu} \big\rangle \\
	 & = \big\langle \boldsymbol{\nu} \big| \mathds{B} (u_N; r_N) \cdots \mathds{B} (u_1; r_1) \mathds{D} (w_M; t_M) \cdots \mathds{D} (w_1; t_1) \big| \boldsymbol{\mu} \big\rangle,
	\end{flalign*}
	
	\noindent and by similar reasoning 
	\begin{flalign*}
	 \displaystyle\sum_{\boldsymbol{\kappa} \in \SeqSign_{n; L}} G_{\boldsymbol{\nu} / \boldsymbol{\kappa}} & (\textbf{w}; \textbf{t} \boldsymbol{\mid} \textbf{y}; \textbf{\emph{s}}) F_{\boldsymbol{\mu} / \boldsymbol{\kappa}} (\textbf{u}; \textbf{r} \boldsymbol{\mid} \textbf{y}; \textbf{s}) \\
	 & = \big\langle \boldsymbol{\nu} \big| \mathds{D} (w_M; t_M) \cdots \mathds{D} (w_1; t_1) \mathds{B} (u_N; r_N) \cdots \mathds{B} (u_1; r_1) \big| \boldsymbol{\mu} \big\rangle. 
	 \end{flalign*} 
	 
	 \noindent These identities, together with $MN$ applications of \Cref{bd2limit} (where the condition \eqref{x1x2bdestimate} there is verified by \eqref{x1x2bdestimate2}) yield the theorem. 
\end{proof}

Let us mention that $H_{\boldsymbol{\lambda} / \boldsymbol{\mu}}$ does not appear to satisfy a Cauchy identity when paired with either $G_{\boldsymbol{\lambda} / \boldsymbol{\mu}}$ or $F_{\boldsymbol{\lambda} / \boldsymbol{\mu}}$; it instead would satisfy one when paired with a different function (which did not appear in \Cref{fgdefinition} above) of the form $\big\langle \boldsymbol{\lambda} \big| \mathds{A} (x_N; r_N) \cdots \mathds{A} (x_1; r_1) \big| \boldsymbol{\mu} \big\rangle$, where $\mathds{A} (x; r) = \widehat{\mathds{T}}_{\textbf{e}_{[1, n]}; \textbf{e}_{[1, n]}} (x; r)$ (and we recall the latter from \eqref{wta}). However, we will not describe this type of Cauchy identity in detail here, since we will not need it.

We will next derive from the general skew-Cauchy identity \Cref{fgidentity} a Cauchy identity, to which end we first require the following proposition that provides a factored form for the function $F_{\boldsymbol{0}^N} (\textbf{x}; \textbf{r} \boldsymbol{\mid} \textbf{y}; \textbf{s})$. Its proof will appear in \Cref{ProofF0} below. 

\begin{prop}
	
	\label{f0} 
	
	Fix an integer $N \ge 1$; finite sequences of complex numbers $\textbf{\emph{x}} = (x_1, x_2, \ldots , x_N)$ and $\textbf{\emph{r}} = (r_1, r_2, \ldots , r_N)$; and infinite sequences of complex numbers $\textbf{\emph{y}} = (y_1, y_2, \ldots )$ and $\textbf{\emph{s}} = (s_1, s_2, \ldots )$. Recalling the signature sequence $\boldsymbol{0}^N = (0^N, 0^N, \ldots , 0^N) \in \Sign_{n; N}$, we have
	\begin{flalign*}
	F_{\boldsymbol{0}^N} (\textbf{\emph{x}}; \textbf{\emph{r}} \boldsymbol{\mid} \textbf{\emph{y}}; \textbf{\emph{s}}) & = \displaystyle\prod_{j = 1}^n s_j^{2n (j - N)} r_j^{2n (j - N - 1)} x_j^{n (N - j + 1)} y_j^{-jn} (r_j^2; q)_n \\
	& \qquad \times \displaystyle\prod_{1 \le i < j \le N} (r_i^2 x_i^{-1} x_j; q)_n (s_i^2 y_i^{-1} y_j; q)_n \displaystyle\prod_{i = 1}^N \displaystyle\prod_{j = 1}^N (x_j y_i^{-1}; q)_n^{-1}.
	\end{flalign*}
	
\end{prop}

We next have the following consequence of \Cref{fgidentity}.

\begin{thm}
	
	\label{fgsum2}
	
	Adopt the notation and assumptions of \Cref{fgidentity} with $L = 0$. Then,
	\begin{flalign}
	\label{sumfg2}
	\begin{aligned} 
	\displaystyle\sum_{\boldsymbol{\lambda} \in \SeqSign_{n; N}} & F_{\boldsymbol{\lambda}} (\textbf{\emph{u}}; \textbf{\emph{r}} \boldsymbol{\mid} \textbf{\emph{y}}; \textbf{\emph{s}}) G_{\boldsymbol{\lambda}} (\textbf{\emph{w}}; \textbf{\emph{t}} \boldsymbol{\mid} \textbf{\emph{y}}; \textbf{\emph{s}}) \\
	& = \displaystyle\prod_{j = 1}^n s_j^{2n (j - N)} r_j^{2n (j - N - 1)} u_j^{n (N - j + 1)} y_j^{-jn} (r_j^2; q)_n \displaystyle\prod_{i = 1}^N \displaystyle\prod_{j = 1}^N (u_j y_i^{-1}; q)_n^{-1} \\ 
	& \qquad \times \displaystyle\prod_{1 \le i < j \le N} (r_i^2 u_i^{-1} u_j; q)_n (s_i^2 y_i^{-1} y_j; q)_n \displaystyle\prod_{i = 1}^M \displaystyle\prod_{j = 1}^N \displaystyle\frac{(t_i^2 u_j w_i^{-1}; q)_n}{t_i^{2n} (u_j w_i^{-1}; q)_n}.
	\end{aligned}
	\end{flalign}
	
\end{thm}

\begin{proof}
	
	Applying \Cref{fgidentity} in the case when $\boldsymbol{\mu} = \boldsymbol{0}^N$ and $\boldsymbol{\nu} = \boldsymbol{\varnothing}$, it is quickly verified by the first identity in \eqref{gfhe} (for $G_{\boldsymbol{\lambda} / \boldsymbol{\mu}}$) that the right side of \eqref{fgsum1} is supported on $\boldsymbol{\kappa} = \boldsymbol{\varnothing}$. This, together with the fact that $G_{\boldsymbol{\varnothing} / \boldsymbol{\varnothing}} (\textbf{w}; \textbf{t} \boldsymbol{\mid} \textbf{y}; \textbf{s}) = 1$ (since all signatures in $\boldsymbol{\varnothing}$ are empty) yields 
	\begin{flalign*}
	\displaystyle\sum_{\boldsymbol{\lambda} \in \SeqSign_{n; N}} & F_{\boldsymbol{\lambda}} (\textbf{u}; \textbf{r} \boldsymbol{\mid} \textbf{y}; \textbf{s}) G_{\boldsymbol{\lambda}} (\textbf{w}; \textbf{t} \boldsymbol{\mid} \textbf{y}; \textbf{s}) = \displaystyle\prod_{i = 1}^M \displaystyle\prod_{j = 1}^N \displaystyle\frac{(t_i^2 u_j w_i^{-1}; q)_n}{t_i^{2n} (u_j w_i^{-1}; q)_n}  F_{\boldsymbol{0}^N} (\textbf{u}; \textbf{r} \boldsymbol{\mid} \textbf{y}; \textbf{s}).
	\end{flalign*}
	
	\noindent This, together with \Cref{f0}, implies the theorem.
\end{proof}

\section{Proof of \Cref{f0}}

\label{ProofF0}

In this section we establish \Cref{f0}. To that end, we begin with the following lemma indicating that, upon suitable normalization, $F_{\boldsymbol{0}^N} (\textbf{x}; \textbf{r} \boldsymbol{\mid} \textbf{y}; \textbf{s})$ is a polynomial of specified degree in $\textbf{x} \cup \textbf{y}$.

\begin{lem}
	
	\label{f0np1}

	For any integers $n, N \ge 1$ and sequences of complex numbers $\textbf{\emph{r}} = (r_1 ,r_2, \ldots , r_N)$ and $\textbf{\emph{s}} = (s_1, s_2, \ldots )$, there exists a constant $C = C_{n; N} (\textbf{\emph{r}}, \textbf{\emph{s}})$ such that the following holds. For any sequences of complex numbers $\textbf{\emph{x}} = (x_1, x_2, \ldots , x_N)$ and $\textbf{\emph{y}} = (y_1, y_2, \ldots )$,  there exists a polynomial $P(\textbf{\emph{x}}, \textbf{\emph{y}})$ of total degree $n N (N - 1)$ in the variables $\textbf{\emph{x}}$ and $\textbf{\emph{y}}$, with coefficients in $\mathbb{C} (q, \textbf{\emph{r}}, \textbf{\emph{s}})$, such that
	\begin{flalign}
	\label{f0np} 
	F_{\boldsymbol{0}^N} (\textbf{\emph{x}}; \textbf{\emph{r}} \boldsymbol{\mid} \textbf{\emph{y}}; \textbf{\emph{s}}) \displaystyle\prod_{i = 1}^N \displaystyle\prod_{j = 1}^N \displaystyle\prod_{k = 0}^{n - 1} (y_i - q^k x_j) = P (\textbf{\emph{x}}, \textbf{\emph{y}}) \displaystyle\prod_{j = 1}^N x_j^n.
	\end{flalign}
	
\end{lem}

\begin{proof}
	
	Since by \eqref{gfhe} $F_{\boldsymbol{0}^N} (\textbf{x}; \textbf{r} \boldsymbol{\mid}  \textbf{y}; \textbf{s})$ is the partition function for the vertex model $\mathfrak{P}_F (\boldsymbol{0}^N / \boldsymbol{\varnothing})$ (from \Cref{pgpfph}; see also the middle of \Cref{fgpaths}) on the domain $\mathcal{D}_N = \mathbb{Z}_{> 0} \times [1, N]$, under the weights $\widehat{W}_z (\textbf{A}, \textbf{B}; \textbf{C}, \textbf{D})$ from \Cref{wabcdrsxy} and \eqref{wabcd2}, it will be useful to analyze properties of these weights. In particular, let us show for any $\textbf{A}, \textbf{B}, \textbf{C}, \textbf{D} \in \{ 0, 1 \}^n$ that $(z; q)_n \widehat{W}_z (\textbf{A}, \textbf{B}; \textbf{C}, \textbf{D} \boldsymbol{\mid} r, s)$ is a polynomial of degree at most $n$ in $z$, with coefficients in $\mathbb{C} (q, r, s)$. By \eqref{wabcd2}, it suffices to show that $(s^2 z; q)_n W_z (\textbf{A}, \textbf{B}; \textbf{C}, \textbf{D} \boldsymbol{\mid} r, s)$ is a polynomial of degree at most $n$ in $z$. 
	
	The polynomiality of this normalized weight follows from the definition \eqref{wabcdp} of $W_z$, since
	\begin{flalign*}
	\displaystyle\frac{(s^2 z; q)_n}{(s^2 z; q)_{c + d - p - v}}  \in \mathbb{C} (q, s) [z]; \qquad \displaystyle\frac{(q^{-v} s^2 r^{-2} z; q)_{c - p}}{(q^{-v} s^2 r^{-2} z; q)_v} = (s^2 r^{-2} z; q)_{c - p - v} \in \mathbb{C} (q, r, s) [z].
	\end{flalign*}
	
	\noindent Here, the first inclusion follows from the fact that $c + d - p - v \le c + d - v \le n$ (since if $\textbf{C} = (C_1, C_2, \ldots , C_n) \in \{ 0, 1 \}^n$ and $\textbf{D} = (D_1, D_2, \ldots , D_n) \in \{ 0, 1 \}^n$, then $c + d - v$ counts the number of indices $i \in [1, n]$ for which $\max \{ C_i, D_i \} = 1$), and the second follows from the fact that $p \le c - v$ (in the sum on the right side of \eqref{wabcdp}). The fact that $(s^2 z; q)_n W_z (\textbf{A}, \textbf{B}; \textbf{C}, \textbf{D} \boldsymbol{\mid} r, s)$ is of degree at most $n$ in $z$ then follows from the fact that the degree (in $z$) of the numerator of any summand on the right side of \eqref{wabcdp} is at most that of the denominator. This verifies that $(z; q)_n \widehat{W}_z (\textbf{A}, \textbf{B}; \textbf{C}, \textbf{D} \boldsymbol{\mid} r, s)$ is a polynomial of degree at most $n$ in $z$. 
	
	Now, observe that $F_{\boldsymbol{0}^N} (\textbf{x}; \textbf{r} \boldsymbol{\mid}  \textbf{y}; \textbf{s})$ is the partition function for the vertex model $\mathfrak{P}_F (\boldsymbol{0}^N / \boldsymbol{\varnothing})$ on the domain $[1, N] \times [1, N]$, whose weight at any vertex $(i, j) \in \mathcal{D}_N = \mathbb{Z}_{> 0} \times [1, N]$ is given by $\widehat{W}_{x_j / y_i} \big( \textbf{A}, \textbf{B}; \textbf{C}, \textbf{D} \boldsymbol{\mid} r_j, s_i \big)$; here, we replaced the original domain $\mathcal{D}_N$ with $[1, N] \times [1, N]$, which is permitted as the arrow configuration at each vertex $(i, j) \in \mathcal{D}_N$ with $i > N$ under any path ensemble in $\mathfrak{P}_F (\boldsymbol{0}^N / \boldsymbol{\varnothing})$ is given by $\big( \textbf{e}_0, \textbf{e}_{[1, n]}; \textbf{e}_0, \textbf{e}_{[1, n]} \big)$, and $\widehat{W}_z \big( \textbf{e}_0, \textbf{e}_{[1, n]}; \textbf{e}_0, \textbf{e}_{[1, n]} \boldsymbol{\mid} r_j, s_i \big) = 1$ by \eqref{wz1}. Thus, the polynomiality of $(z; q)_n \widehat{W}_z (\textbf{A}, \textbf{B}; \textbf{C}, \textbf{D} \boldsymbol{\mid} r, s)$ implies that 
	\begin{flalign}
	\label{1f0n} 
	F_{\boldsymbol{0}^N} (\textbf{x}; \textbf{r} \boldsymbol{\mid} \textbf{y}; \textbf{s}) \displaystyle\prod_{i = 1}^N \displaystyle\prod_{j = 1}^N (x_j y_i^{-1}; q)_n
	\end{flalign}
	
	\noindent is a polynomial in the variables $\{ x_j y_i^{-1} \}$, whose degree in any individual $x_j y_i^{-1}$ is at most $n$. Multiplying by $\prod_{j = 1}^N y_j^{nN}$ yields that the left side of \eqref{f0np} is a polynomial in $\textbf{x}$ and $\textbf{y}$ of total degree at most $nN^2$. 
	
	It therefore remains to show that the left side of \eqref{f0np} is a multiple of $x_j^n$ for each $j \in [1, N]$. To that end, observe from \eqref{wabcdp} and \eqref{wabcd2} that $z^{|\textbf{D}| - |\textbf{B}|}$ divides $(z; q)_n \widehat{W}_z (\textbf{A}, \textbf{B}; \textbf{C}, \textbf{D} \boldsymbol{\mid} r, s)$. Now observe under each path ensemble $\mathcal{E} \in \mathfrak{P}_F (\boldsymbol{0}^N)$, with arrow configuration $\big( \textbf{A} (v), \textbf{B} (v); \textbf{C} (v), \textbf{D} (v) \big)$ at any vertex $v \in \mathcal{D}_N$, that for any $j \in [1, N]$ we have $\sum_{i = 1}^N \big( |\textbf{D}| (i, j) - |\textbf{B}| (i, j) \big) = \big| \textbf{D} (N, j) \big| - \big| \textbf{B} (1, j) \big| = n$, since $\textbf{B} (1, j) = \textbf{e}_0$ and $\textbf{D} (N, j) = \textbf{e}_{[1, n]}$. Hence, $x_j^n$ divides the numerator of \eqref{1f0n}, which upon multiplication by $\prod_{j = 1}^N y_j^{nN}$ implies that it also divides the left side of \eqref{f0np}.
	
	Hence, this left side is a polynomial of degree $n N^2$ in $\textbf{x} \cup \textbf{y}$ and is a multiple of $\prod_{j = 1}^N x_j^n$. This implies that there exists a polynomial $P \in \mathbb{C} (q, \textbf{r}, \textbf{s}) [\textbf{x}, \textbf{y}]$ of degree $n (N^2 - N)$ such that \eqref{f0np} holds.
\end{proof}

We must next determine the polynomial $P(\textbf{x}, \textbf{y})$ in \Cref{f0np1}, which we do by specifying sufficiently many of its zeroes. To that end, the following symmetry relation for $F_{\boldsymbol{0}^N} (\textbf{x}; \textbf{r} \boldsymbol{\mid} \textbf{y}; \textbf{s})$ in $\textbf{y}$ and $\textbf{s}$ (similar to \Cref{gxfxsigma} for the $(\textbf{x}, \textbf{r})$ variables) will be useful. 

\begin{lem}
	
	\label{f0nxyy}
	
	Fix integers $n, N \ge 1$; finite sequences of complex numbers $\textbf{\emph{x}} = (x_1, x_2, \ldots , x_N)$ and $\textbf{\emph{r}} = (r_1, r_2, \ldots , r_N)$; and infinite sequences of complex numbers $\textbf{\emph{y}} = (y_1, y_2, \ldots )$ and $\textbf{\emph{s}} = (s_1, s_2, \ldots)$. For any permutation $\sigma \in \mathfrak{S}_N$, we have 
	\begin{flalign*}
	F_{\boldsymbol{0}^N} \big( \textbf{\emph{x}}; \textbf{\emph{r}} \boldsymbol{\mid} \sigma (\textbf{\emph{y}}); \sigma (\textbf{\emph{s}}) \big) = F_{\boldsymbol{0}^N} (\textbf{\emph{x}}; \textbf{\emph{r}} \boldsymbol{\mid} \textbf{\emph{y}}; \textbf{\emph{s}} \big) \displaystyle\prod_{\substack{1 \le i < j \le N \\ \sigma (i) > \sigma (j)}} \displaystyle\frac{(s_j^2 y_i y_j^{-1}; q)_n}{(s_i^2 y_i^{-1} y_j; q)_n} \bigg( \displaystyle\frac{s_i^2 y_j}{s_j^2 y_i} \bigg)^n.
	\end{flalign*}
	
\end{lem}

\begin{proof}
	
	Since $\mathfrak{S}_N$ is generated by the transpositions $\{ \mathfrak{s}_i \}$ for $i \in [1, N]$ interchanging $(i, i + 1)$, we may assume that $\sigma = \mathfrak{s}_k$ for some $k \in [1, N - 1]$. Then, by $N$ applications of the Yang--Baxter equation \Cref{wabcdproduct2}, we obtain
	\begin{flalign}
	\label{wyy} 
	\begin{aligned}
	\widehat{W}_{y_k / y_{k + 1}} \big( \textbf{e}_{[1, n]}, \textbf{e}_{[1, n]};  \textbf{e}_{[1, n]},  \textbf{e}_{[1, n]} & \boldsymbol{\mid} s_k, s_{k + 1} \big) F_{\boldsymbol{0}^N} \big( \textbf{x}; \textbf{r} \boldsymbol{\mid} \sigma (\textbf{y}); \sigma (\textbf{s}) \big) \\
	& = F_{\boldsymbol{0}^N} (\textbf{x}; \textbf{r} \boldsymbol{\mid} \textbf{y}; \textbf{s}) \widehat{W}_{y_k / y_{k + 1}} \big( \textbf{e}_0,  \textbf{e}_0;  \textbf{e}_0,  \textbf{e}_0 \boldsymbol{\mid} s_k, s_{k + 1} \big).
	\end{aligned} 
	\end{flalign} 

	\noindent Diagrammatically, 
		
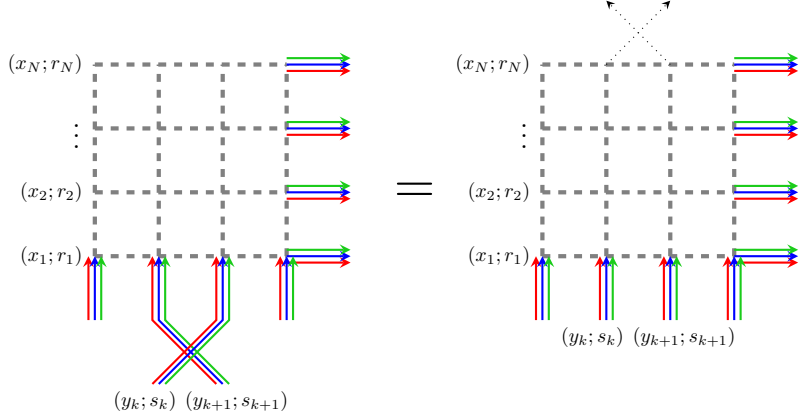
\begin{center}
		\begin{tikzpicture}[
			>=stealth,
			scale = .85
			]
			\draw[gray, ultra thick, dashed] (1, 1) -- (4, 1);
			\draw[gray, ultra thick, dashed] (1, 2) -- (4, 2);
			\draw[gray, ultra thick, dashed] (1, 3) -- (4, 3);
			\draw[gray, ultra thick, dashed] (1, 4) -- (4, 4);	
			\draw[gray, ultra thick, dashed] (1, 1) -- (1, 4);
			\draw[gray, ultra thick, dashed] (2, 1) -- (2, 4);
			\draw[gray, ultra thick, dashed] (3, 1) -- (3, 4);
			\draw[gray, ultra thick, dashed] (4, 1) -- (4, 4);
			\draw[->, thick, red] (.9, 0) -- (.9, 1);
			\draw[->, thick, red] (4, 3.9) -- (5, 3.9);
			\draw[->, thick, blue] (1, 0) -- (1, 1);
			\draw[->, thick, blue] (4, 4) -- (5, 4);
			\draw[->, thick, green] (1.1, 0) -- (1.1, 1); 
			\draw[->, thick, green] (4, 4.1) -- (5, 4.1);
			\draw[->, thick, red] (2.9, -1) -- (1.9, 0) -- (1.9, 1);
			\draw[->, thick, red] (4, 2.9) -- (5, 2.9);
			\draw[->, thick, blue] (3, -1) node[black, scale = .75, right = 7, below]{$(y_{k + 1}; s_{k + 1})$} -- (2, 0) -- (2, 1); 
			\draw[->, thick, blue] (4, 3) -- (5, 3);
			\draw[->, thick, green] (3.1, -1) -- (2.1, 0) -- (2.1, 1);
			\draw[->, thick, green] (4, 3.1) -- (5, 3.1);
			\draw[->, thick, red] (1.9, -1) -- (2.9, 0) -- (2.9, 1);
			\draw[->, thick, red] (4, 1.9) -- (5, 1.9);
			\draw[->, thick, blue] (2, -1) node[black, scale = .75, left = 7, below]{$(y_k; s_k)$}-- (3, 0) -- (3, 1);
			\draw[->, thick, blue] (4, 2) -- (5, 2);
			\draw[->, thick, green] (2.1, -1) -- (3.1, 0) -- (3.1, 1);
			\draw[->, thick, green] (4, 2.1) -- (5, 2.1);
			\draw[->, thick, red] (3.9, 0) -- (3.9, 1);
			\draw[->, thick, red] (4, .9) -- (5, .9);
			\draw[->, thick, blue] (4, 0) -- (4, 1);
			\draw[->, thick, blue] (4, 1) -- (5, 1);
			\draw[->, thick, green] (4.1, 0) -- (4.1, 1);
			\draw[->, thick, green] (4, 1.1) -- (5, 1.1);
			\draw[] (6, 2) circle[radius = 0] node[scale = 2]{$=$}; 
			\draw[gray, ultra thick, dashed] (8, 1) -- (11, 1);
			\draw[gray, ultra thick, dashed] (8, 2) -- (11, 2);
			\draw[gray, ultra thick, dashed] (8, 3) -- (11, 3);
			\draw[gray, ultra thick, dashed] (8, 4) -- (11, 4);
			\draw[gray, ultra thick, dashed] (8, 1) -- (8, 4);
			\draw[gray, ultra thick, dashed] (9, 1) -- (9, 4);
			\draw[gray, ultra thick, dashed] (10, 1) -- (10, 4);
			\draw[gray, ultra thick, dashed] (11, 1) -- (11, 4);	
			\draw[dotted, ->] (9, 4) -- (10, 5);
			\draw[dotted, ->] (10, 4) -- (9, 5);	
			\draw[->, thick, red] (7.9, 0) -- (7.9, 1);
			\draw[->, thick, red] (11, 3.9) -- (12, 3.9);
			\draw[->, thick, blue] (8, 0) -- (8, 1);
			\draw[->, thick, blue] (11, 4) -- (12, 4);
			\draw[->, thick, green] (8.1, 0) -- (8.1, 1); 
			\draw[->, thick, green] (11, 4.1) -- (12, 4.1);	
			\draw[->, thick, red] (8.9, 0) -- (8.9, 1);
			\draw[->, thick, red] (11, 2.9) -- (12, 2.9);
			\draw[->, thick, blue] (9, 0) node[scale = .75, black, left = 7, below]{$(y_k; s_k)$}-- (9, 1); 
			\draw[->, thick, blue] (11, 3) -- (12, 3);
			\draw[->, thick, green] (9.1, 0) -- (9.1, 1);
			\draw[->, thick, green] (11, 3.1) -- (12, 3.1);	
			\draw[->, thick, red] (9.9, 0) -- (9.9, 1);
			\draw[->, thick, red] (11, 1.9) -- (12, 1.9);
			\draw[->, thick, blue] (10, 0) node[scale = .75, black, right = 7, below]{$(y_{k + 1}; s_{k + 1})$}-- (10, 1);
			\draw[->, thick, blue] (11, 2) -- (12, 2);
			\draw[->, thick, green] (10.1, 0) -- (10.1, 1);
			\draw[->, thick, green] (11, 2.1) -- (12, 2.1);			
			\draw[->, thick, red] (10.9, 0) -- (10.9, 1);
			\draw[->, thick, red] (11, .9) -- (12, .9);
			\draw[->, thick, blue] (11, 0) -- (11, 1);
			\draw[->, thick, blue] (11, 1) -- (12, 1);
			\draw[->, thick, green] (11.1, 0) -- (11.1, 1);
			\draw[->, thick, green] (11, 1.1) -- (12, 1.1);	
			\draw[] (1, 1) circle [radius = 0] node[left = 2, scale = .75]{$(x_1; r_1)$};
			\draw[] (1, 2) circle [radius = 0] node[left = 2, scale = .75]{$(x_2; r_2)$};
			\draw[] (1, 3) circle [radius = 0] node[left = 2, scale = 1]{$\vdots$};
			\draw[] (1, 4) circle [radius = 0] node[left = 2, scale = .75]{$(x_N; r_N)$};
			
			\draw[] (8, 1) circle [radius = 0] node[left = 2, scale = .75]{$(x_1; r_1)$};
			\draw[] (8, 2) circle [radius = 0] node[left = 2, scale = .75]{$(x_2; r_2)$};
			\draw[] (8, 3) circle [radius = 0] node[left = 2, scale = 1]{$\vdots$};
			\draw[] (8, 4) circle [radius = 0] node[left = 2, scale = .75]{$(x_N; r_N)$};
		\end{tikzpicture}
	\end{center}
	
	\noindent where the partition functions on both sides are with respect to the $\widehat{W}_z$ weights. 
	
	Next, by \eqref{wabcd01nequation} and \eqref{wabcd2}, we have
	\begin{flalign*}
	\widehat{W}_{y_k / y_{k + 1}} \big( \textbf{e}_{[1, n]},  \textbf{e}_{[1, n]};  \textbf{e}_{[1, n]},  \textbf{e}_{[1, n]} \boldsymbol{\mid} s_k, s_{k + 1} \big) & = \bigg( \displaystyle\frac{y_k}{s_k^2 y_{k + 1}} \bigg)^n \displaystyle\frac{(s_k^2 y_k^{-1} y_{k + 1}; q)_n}{(y_k y_{k + 1}^{-1}; q)_n}; \\
	\widehat{W}_{y_k / y_{k + 1}} \big( \textbf{e}_0, \textbf{e}_0; \textbf{e}_0, \textbf{e}_0 \boldsymbol{\mid} s_k, s_{k + 1} \big) & = \displaystyle\frac{(s_{k + 1}^2 y_k y_{k + 1}^{-1}; q)_n}{s_{k + 1}^{2n} (y_k y_{k + 1}^{-1}; q)_n},
	\end{flalign*} 
	
	\noindent which upon insertion into \eqref{wyy} yields
	\begin{flalign*}
	F_{\boldsymbol{0}^N} \big( \textbf{x}; \textbf{r} \boldsymbol{\mid} \sigma (\textbf{y}); \sigma (\textbf{s}) \big) = \bigg( \displaystyle\frac{s_k^2 y_{k + 1}}{s_{k + 1}^2 y_k} \bigg)^n \displaystyle\frac{(s_{k + 1}^2 y_k y_{k + 1}^{-1}; q)_n}{(s_k^2 y_k^{-1} y_{k + 1}; q)_n} F_{\boldsymbol{0}^N} (\textbf{x}; \textbf{r} \boldsymbol{\mid} \textbf{y}; \textbf{s}).
	\end{flalign*} 
	
	\noindent This yields the lemma for $\sigma = \mathfrak{s}_k$, which as mentioned above implies it in general.
\end{proof}

The following corollary then determines $F_{\boldsymbol{0}^N} (\textbf{x}; \textbf{r} \boldsymbol{\mid} \textbf{y}; \textbf{s})$, up to an overall constant independent of $\textbf{x}$ and $\textbf{y}$. 

\begin{cor}
	
	\label{f0nxy} 
	
	For any integers $n, N \ge 1$ and sequences of complex numbers $\textbf{\emph{r}} = (r_1 ,r_2, \ldots , r_N)$ and $\textbf{\emph{s}} = (s_1, s_2, \ldots )$, there exists a constant $C = C_{n; N} (\textbf{\emph{r}}, \textbf{\emph{s}})$ such that the following holds. For any sequences of complex numbers $\textbf{\emph{x}} = (x_1, x_2, \ldots , x_N)$ and $\textbf{\emph{y}} = (y_1, y_2, \ldots )$, we have that 
	\begin{flalign*}
	F_{\boldsymbol{0}^N} (\textbf{\emph{x}}; \textbf{\emph{r}} \boldsymbol{\mid} \textbf{\emph{y}}; \textbf{\emph{s}}) = C \displaystyle\prod_{j = 1}^n x_j^n  \displaystyle\prod_{1 \le i < j \le N} \displaystyle\prod_{k = 0}^{n - 1} (x_i - q^k r_i^2 x_j) (y_i - q^k s_i^2 y_j) \displaystyle\prod_{i = 1}^N \displaystyle\prod_{j = 1}^N \displaystyle\prod_{k = 0}^{n - 1} (y_i - q^k x_j)^{-1}.
	\end{flalign*}
	
\end{cor}

\begin{proof}
	
	Let $\sigma_0 \in \mathfrak{S}_N$ denote the longest element, that is, it is defined so that $\sigma_0 (j) = N - j + 1$ for each $j \in [1, N]$. Then, \Cref{gxfxsigma} and \Cref{f0nxyy} give
	\begin{flalign}
	\label{sigma0f0} 
	\begin{aligned}
	F_{\boldsymbol{0}^N} ( \textbf{x}; \textbf{r} \boldsymbol{\mid} \textbf{y}; \textbf{s}) & = F_{\boldsymbol{0}^N} \big( \sigma_0 (\textbf{x}); \sigma_0 (\textbf{r}) \boldsymbol{\mid} \textbf{y}; \textbf{s} \big) \displaystyle\prod_{1 \le i < j \le N} \displaystyle\frac{(r_i^2 x_i^{-1} x_j; q)_n}{(r_j^2 x_i x_j^{-1}; q)_n} \bigg( \displaystyle\frac{r_j^2 x_i}{r_i^2 x_j} \bigg)^n \\
	F_{\boldsymbol{0}^N} ( \textbf{x}; \textbf{r} \boldsymbol{\mid} \textbf{y}; \textbf{s} \big) & = F_{\boldsymbol{0}^N} \big( \textbf{x}; \textbf{r} \boldsymbol{\mid} \sigma_0 (\textbf{y}); \sigma_0 (\textbf{s}) \big)  \displaystyle\prod_{1 \le i < j \le N} \displaystyle\frac{(s_i^2 y_i^{-1} y_j; q)_n}{(s_j^2 y_i y_j^{-1}; q)_n} \bigg( \displaystyle\frac{s_j^2 y_i}{s_i^2 y_j} \bigg)^n. 
	\end{aligned} 
	\end{flalign}
	
	\noindent Hence, since \Cref{f0np1} implies that 
	\begin{flalign}
	\label{f0xjp} 
	F_{\boldsymbol{0}^N} (\textbf{x}; \textbf{r} \boldsymbol{\mid} \textbf{y}; \textbf{s}) \displaystyle\prod_{j = 1}^N x_j^{-n} \displaystyle\prod_{i = 1}^N \displaystyle\prod_{j = 1}^N \displaystyle\prod_{k = 0}^{n - 1} (y_i - q^k x_j)
	\end{flalign} 
	
	\noindent is a polynomial in $\textbf{x}$ and $\textbf{y}$ of total degree $n N (N - 1)$, it follows from \eqref{sigma0f0} that \eqref{f0xjp} is equal to $0$ if there exist $1 \le i < j \le N$ and $k \in [0, n - 1]$ such that either $x_i = q^k r_i^2 x_j$ or $y_i = q^k s_i^2 y_j$. This implies that 
	\begin{flalign}
	\label{productxrys} 
	\displaystyle\prod_{1 \le i < j \le N} \displaystyle\prod_{k = 0}^{n - 1} (x_i - q^k r_i^2 x_j) (y_i - q^k s_i^2 y_j)
	\end{flalign}
	
	\noindent divides \eqref{f0xjp}. The corollary then follows since both \eqref{f0xjp} and \eqref{productxrys} are of total degree $n N (N - 1)$ in $\textbf{x}$ and $\textbf{y}$. 
\end{proof}

We next state the following lemma, which will be established in \Cref{Proofc} below, that provides an explicit form for the constant $C_{n; N} (\textbf{r}, \textbf{s})$ from \Cref{f0nxy}.

\begin{lem}
	
	\label{f0c}
	
	Adopting the notation of \Cref{f0nxy}, we have 
	\begin{flalign*}
	C_{n; N} (\textbf{\emph{r}}, \textbf{\emph{s}}) = \displaystyle\prod_{j = 1}^N s_j^{2n (j - N)} r_j^{2n (j - N - 1)} (r_j^2; q)_n.
	\end{flalign*}
	
\end{lem}

Given the above, we can now quickly establish \Cref{f0}. 

\begin{proof}[Proof of \Cref{f0}]

By \Cref{f0nxy} and \Cref{f0c}, we have 
\begin{flalign*}
F_{\boldsymbol{0}^N} (\textbf{x}; \textbf{r} \boldsymbol{\mid} \textbf{y}; \textbf{s}) & = \displaystyle\prod_{1 \le i < j \le N} \displaystyle\prod_{k = 0}^{n - 1} (x_i - q^k r_i^2 x_j) (y_i - q^k s_i^2 y_j) \displaystyle\prod_{j = 1}^N s_j^{2n (j - N)} x_j^n r_j^{2n (j - N - 1)} (r_j^2; q)_n \\
& \qquad \quad \times \displaystyle\prod_{i = 1}^N \displaystyle\prod_{j = 1}^N \displaystyle\prod_{k = 0}^{n - 1} (y_i - q^k x_j)^{-1}.
\end{flalign*}

\noindent This, together with the facts that 
\begin{flalign}
\label{productxyqk} 
\begin{aligned} 
\displaystyle\prod_{1 \le i < j \le N} \displaystyle\prod_{k = 0}^{n - 1} (x_i - q^k r_i^2 x_j) (y_i - q^k s_i^2 y_j) & = \displaystyle\prod_{1 \le i < j \le N} (r_i^2 x_i^{-1} x_j; q)_n (s_i^2 y_i^{-1} y_j; q)_n \displaystyle\prod_{j = 1}^N (x_j y_j)^{n (N - j)}; \\
\displaystyle\prod_{i = 1}^N \displaystyle\prod_{j = 1}^N \displaystyle\prod_{k = 0}^{n - 1} (y_i - q^k x_j) & = \displaystyle\prod_{j = 1}^N y_j^{nN} \displaystyle\prod_{i = 1}^N \displaystyle\prod_{j = 1}^N (x_j y_i^{-1}; q)_n,
\end{aligned} 
\end{flalign}

\noindent imply the proposition.
\end{proof}

\section{Proof of \Cref{f0c}}

\label{Proofc} 

In this section we establish \Cref{f0c}, which explicitly determines the constant $C_{n; N} (\textbf{r}, \textbf{s})$ from \Cref{f0nxy}. To do this, it essentially suffices to evaluate $F_{\boldsymbol{0}^N} (\textbf{x}; \textbf{r} \boldsymbol{\mid} \textbf{y}; \textbf{s})$ for any choice of parameters $(\textbf{x}, \textbf{y})$. Recalling the interpretation \eqref{gfhe} of $F_{\boldsymbol{0}^N} (\textbf{x}; \textbf{r} \boldsymbol{\mid} \textbf{y}; \textbf{s})$ as a partition function under the weights $\widehat{W}_{x_j/ y_i} (\textbf{A}, \textbf{B}; \textbf{C}, \textbf{D} \boldsymbol{\mid} r_j, s_i)$ (from \eqref{wabcd2}) for the vertex model $\mathfrak{P}_F (\boldsymbol{0}^N / \boldsymbol{\varnothing})$ (from \Cref{pgpfph}; see the middle of \Cref{fgpaths}), we will choose these parameters $(\textbf{x}, \textbf{y})$ in such a way that this partition function \emph{freezes}, that is, it admits at most one path ensemble with nonzero weight. The latter task will be facilitated through the following proposition, which provides a vanishing condition for the weights $\widehat{W}_z (\textbf{A}, \textbf{B}; \textbf{C}, \textbf{D} \boldsymbol{\mid} r, s)$ when $r^2 = s^2 z$.

\begin{lem}
	
	\label{szr0} 
	
	Fix complex numbers $s, r, z \in \mathbb{C}$ such that $r^2 = s^2 z$. For any $\textbf{\emph{A}}, \textbf{\emph{B}},\textbf{\emph{C}}, \textbf{\emph{D}}, \in \{ 0, 1 \}^n$, we have that $W_z (\textbf{\emph{A}}, \textbf{\emph{B}}; \textbf{\emph{C}}, \textbf{\emph{D}} \boldsymbol{\mid} r, s) = 0 = \widehat{W}_z (\textbf{\emph{A}}, \textbf{\emph{B}}; \textbf{\emph{C}}, \textbf{\emph{D}} \boldsymbol{\mid} r, s)$ unless $\textbf{\emph{A}} + \textbf{\emph{B}} = \textbf{\emph{C}} + \textbf{\emph{D}}$ and $\textbf{\emph{B}} \ge \textbf{\emph{C}}$. 
	
\end{lem}  

\begin{proof} 

	It suffices to show that the lemma holds for the weight $W_z (\textbf{A}, \textbf{B}; \textbf{C}, \textbf{D} \boldsymbol{\mid} r, s)$, for then \eqref{wabcd2} would imply that it also holds for $\widehat{W}_z (\textbf{A}, \textbf{B}; \textbf{C}, \textbf{D} \boldsymbol{\mid} r, s)$. We may assume in what follows that $\textbf{A} + \textbf{B} = \textbf{C} + \textbf{D}$, for otherwise arrow conservation implies $W_z (\textbf{A}, \textbf{B}; \textbf{C}, \textbf{D} \boldsymbol{\mid} r, s) = 0$ for any $r, s, z \in \mathbb{C}$. Then, setting $r^2 = s^2 z$ in \Cref{wabcdrsxy} gives 
	\begin{flalign}
	\label{wzrszsum}
	\begin{aligned} 
	W_z (\textbf{A}, \textbf{B}; \textbf{C}, \textbf{D} \boldsymbol{\mid} r, s) &= (-1)^v z^{c + d - a - b} s^{2c + 2d - 2a} q^{\varphi (\textbf{D} - \textbf{V}, \textbf{C})+ \varphi (\textbf{V}, \textbf{A}) - av + cv} \\
	& \quad \times \displaystyle\frac{(q^{1 - v} s^{-2}; q)_v}{(q^{-v}; q)_v} \displaystyle\frac{(s^2 z; q)_d}{(s^2 z; q)_b} \displaystyle\sum_{p = 0}^{\min \{ b - v, c - v \}} \displaystyle\frac{(q^{-v}; q)_{c - p} (q^v s^2; q)_p (z; q)_{b - p - v}}{(s^2 z; q)_{c + d - p - v}} \\
	& \qquad \qquad \qquad \times q^{-pv} s^{-2p} \displaystyle\sum_{\textbf{\emph{P}}} q^{\varphi (\textbf{B} - \textbf{D} - \textbf{P}, \textbf{P})}, 
	\end{aligned}
	\end{flalign}
	
	\noindent where the sum is over all $n$-tuples $\textbf{P} = (P_1, P_2, \ldots , P_n) \in \{ 0, 1 \}^n$ such that $|\textbf{P}| = p$ and $P_i \le \min \{ B_i - V_i, C_i -V_i \}$, for each $i \in [1, n]$. 
	
	Since $c - p \ge v$, we have $(q^{-v}; q)_{c - p} = \textbf{1}_{v = c - p} (q^{-v}; q)_v$. Thus, the sum over $p$ on the right side of \eqref{wzrszsum} is supported on the term $p = c - v$; since $P_i \le C_i - V_i$ for each $i \in [1, n]$, this implies that the sum over $\textbf{P}$ there is supported on $\textbf{P} = \textbf{C} - \textbf{V}$. Since we must also have that $P_i \le B_i - V_i$ for each $i \in [1, n]$, this implies that $W_z ( \textbf{A}, \textbf{B}; \textbf{C}, \textbf{D} \boldsymbol{\mid} r, s) = 0$ unless $C_i \le B_i$ for each $i \in [1, n]$ (that is, unless $\textbf{B} \ge \textbf{C}$). 
\end{proof} 

\begin{rem}

\label{wsrz} 

As \Cref{sxz1weights} below, we will in fact derive a fully factored form for this weight $W_z (\textbf{A}, \textbf{B}; \textbf{C}, \textbf{D} \boldsymbol{\mid} r, s)$ under the specialization $r^2 = s^2 z$, which will follow from the fact that the sum on the right side of \eqref{wzrszsum} is supported on a single term.

\end{rem}

The following corollary then evaluates $F_{\boldsymbol{0}^N} (\textbf{x}; \textbf{r} \boldsymbol{\mid} \textbf{y}; \textbf{s})$ when $s_j^2 x_j = r_j^2 y_j$ for each $j$.

\begin{cor}
	
	\label{f0rxsy}
	
	Adopt the notation and assumptions of \Cref{f0nxy}, and further assume for each $j \in [1, N]$ that $s_j^2 x_j = r_j^2 y_j$. Then,
	\begin{flalign*}
	F_{\boldsymbol{0}^N} (\textbf{\emph{x}}; \textbf{\emph{r}} \boldsymbol{\mid} \textbf{\emph{y}}; \textbf{\emph{s}}) & = \displaystyle\prod_{j = 1}^N s_j^{2n (j - N)} x_j^{n (N - j + 1)} r_j^{2n (j - N - 1)} y_j^{-jn} (r_j^2; q)_n \\
	& \qquad \times \displaystyle\prod_{1 \le i < j \le N} (r_i^2 x_i^{-1} x_j; q)_n (s_i^2 y_i^{-1} y_j; q)_n \displaystyle\prod_{i = 1}^N \displaystyle\prod_{j = 1}^N (x_j y_i^{-1}; q)_n^{-1}.
	\end{flalign*}
\end{cor}

\begin{proof}

	Recall from \eqref{gfhe} the interpretation of $F_{\boldsymbol{\lambda} / \boldsymbol{\mu}}$ as a partition function under the $\widehat{W}_z$ weights (from \eqref{wabcd2}) for the vertex model $\mathfrak{P}_F (\boldsymbol{0}^N / \boldsymbol{\varnothing})$ (defined in \Cref{pgpfph} and depicted in the middle of \Cref{fgpaths}). Thus, for any indices $c \in [1, n]$ and $j \in [1, N]$, one path of color $c$ vertically enters the model through $(j, 1)$ and horizontally exits it through $(N, j)$; we refer to \Cref{f0paths} for an example.

	\begin{figure}
		
		\begin{center}

			\begin{tikzpicture}[
			>=stealth,
			scale = .75
			]
			
			\draw[->, thick, red] (.9, 0) -- (.9, .9) -- (1.9, .9) -- (1.9, 1.9) -- (2.9, 1.9) -- (2.9, 2.9) -- (3.9, 2.9) -- (3.9, 3.9) -- (6, 3.9);
			\draw[->, thick, blue] (1, 0) node[below, black, scale = .8]{$1$} -- (1, 1) -- (2, 1) -- (2, 2) -- (3, 2) -- (3, 3) -- (4, 3) -- (4, 4) -- (6, 4);
			\draw[->, thick, green] (1.1, 0) -- (1.1, 1.1) -- (2.1, 1.1) -- (2.1, 2.1) -- (3.1, 2.1) -- (3.1, 3.1) -- (4.1, 3.1) -- (4.1, 4.1) -- (6, 4.1);
			
			\draw[->, thick, red] (1.9, 0) -- (1.9, .9) -- (2.9, .9) -- (2.9, 1.9) -- (3.9, 1.9) -- (3.9, 2.9) -- (6, 2.9);
			\draw[->, thick, blue] (2, 0) node[below, black, scale = .8]{$2$} -- (2, 1) -- (3, 1) -- (3, 2) -- (4, 2) -- (4, 3) -- (6, 3);
			\draw[->, thick, green] (2.1, 0) -- (2.1, 1.1) -- (3.1, 1.1) -- (3.1, 2.1) -- (4.1, 2.1) -- (4.1, 3.1) -- (6, 3.1);
			
			\draw[->, thick, red] (2.9, 0) -- (2.9, .9) -- (3.9, .9) -- (3.9, 1.9) -- (6, 1.9);
			\draw[->, thick, blue] (3, 0) node[below, black, scale = .8]{$3$} -- (3, 1) -- (4, 1) -- (4, 2) -- (6, 2);
			\draw[->, thick, green] (3.1, 0) -- (3.1, 1.1) -- (4.1, 1.1) -- (4.1, 2.1) -- (6, 2.1);
			
			\draw[->, thick, red] (3.9, 0) -- (3.9, .9)-- (6, .9);
			\draw[->, thick, blue] (4, 0) node[below, black, scale = .8]{$4$} -- (4, 1) -- (6, 1);
			\draw[->, thick, green] (4.1, 0) -- (4.1, 1.1) -- (6, 1.1);
			
			\draw[] (0, 1) circle[radius = 0] node[left = 2, scale = .8]{$1$};
			\draw[] (0, 2) circle[radius = 0] node[left = 2, scale = .8]{$2$};
			\draw[] (0, 3) circle[radius = 0] node[left = 2, scale = .8]{$3$};
			\draw[] (0, 4) circle[radius = 0] node[left = 2, scale = .8]{$4$};
			
			\draw[->, very thick] (0, 0) -- (0, 4.75);
			\draw[->, very thick] (0, 0) -- (6.75, 0);
			
			\end{tikzpicture}
			
		\end{center}
		
		\caption{\label{f0paths} Shown above is the frozen vertex model correpsonding to the $F_{\boldsymbol{0}^N}$ when $s_j^2 x_j = r_j^2 y_j$ from the proof of \Cref{f0rxsy}. } 
		
	\end{figure}

	By \Cref{szr0}, $\widehat{W}_{x_j / y_j} (\textbf{A}, \textbf{B}; \textbf{C}, \textbf{D} \boldsymbol{\mid} r_j, s_j) = 0$ unless $\textbf{B} \ge \textbf{C}$. From this fact, it is quickly verified by induction on $N$ that this vertex model is frozen, that is, it admits only one path ensemble with a nonzero weight; it is given by the one depicted in \Cref{f0paths}, in which to the left of the $(N + 1)$-st column each path alternates between moving one space north and one space east (and after the $(N + 1)$-st column, all paths proceed east). More specifically, in this ensemble, vertices $(i, j)$ have arrow configuration $(\textbf{e}_0, \textbf{e}_0; \textbf{e}_0, \textbf{e}_0)$ if $1 \le i < j \le N$; arrow configuration $\big( \textbf{e}_{[1, n]}, \textbf{e}_0; \textbf{e}_0, \textbf{e}_{[1, n]} \big)$ if $1 \le i = j \le N$; arrow configuration $\big( \textbf{e}_{[1, n]}, \textbf{e}_{[1, n]}; \textbf{e}_{[1, n]}, \textbf{e}_{[1, n]} \big)$ if $ 1 \le j < i \le N$; and arrow configuration $\big( \textbf{e}_0, \textbf{e}_{[1, n]}; \textbf{e}_0, \textbf{e}_{[1, n]} \big)$ if $1 \le j \le N < i$. Hence,
\begin{flalign}
\label{fproductw} 
\begin{aligned} 
F_{\boldsymbol{0}^N} (\textbf{x}; \textbf{r} \boldsymbol{\mid} \textbf{y}; \textbf{s}) & = \displaystyle\prod_{1 \le i < j \le N}  \widehat{W}_{x_j / y_i} \big( \textbf{e}_0, \textbf{e}_0; \textbf{e}_0, \textbf{e}_0 \boldsymbol{\mid} r_j, s_i\big) \displaystyle\prod_{j = 1}^N \widehat{W}_{x_j / y_j} \big( \textbf{e}_{[1, n]}, \textbf{e}_0; \textbf{e}_0, \textbf{e}_{[1, n]} \boldsymbol{\mid} r_j, s_j \big) \\
& \qquad \quad \times \displaystyle\prod_{1 \le j < i \le N} \widehat{W}_{x_j / y_i} \big( \textbf{e}_{[1, n]}, \textbf{e}_{[1, n]}; \textbf{e}_{[1, n]}, \textbf{e}_{[1, n]} \boldsymbol{\mid} r_j, s_i \big),
\end{aligned} 
\end{flalign}

\noindent where we have used the fact that $\widehat{W}_z \big( \textbf{e}_0, \textbf{e}_{[1, n]}; \textbf{e}_0, \textbf{e}_{[1, n]} \boldsymbol{\mid} r, s \big) = 1$ (by \eqref{wz1}). By \eqref{wabcd01nequation} and \eqref{wabcd2}, we have 
\begin{flalign}
\label{w1w}
\begin{aligned} 
\widehat{W}_{x_j / y_i} (\textbf{e}_0, \textbf{e}_0; \textbf{e}_0, \textbf{e}_0 \boldsymbol{\mid} r_j, s_i) & = \displaystyle\frac{(s_i^2 x_j y_i^{-1}; q)_n}{s_i^{2n} (x_j y_i^{-1}; q)_n}; \\
\widehat{W}_{x_j / y_i} \big( \textbf{e}_{[1, n]}, \textbf{e}_{[1, n]}; \textbf{e}_{[1, n]}, \textbf{e}_{[1, n]} \boldsymbol{\mid} r_j, s_i \big) & = \bigg( \displaystyle\frac{x_j}{r_j^2 y_i} \bigg)^n \displaystyle\frac{(r_j^2 x_j^{-1} y_i; q)_n}{(x_j y_i^{-1}; q)_n}. 
\end{aligned} 
\end{flalign}

\noindent We further have by \eqref{wabcdp} that 
\begin{flalign}
\label{wze}
W_z \big(\textbf{e}_{[1, n]}, \textbf{e}_0; \textbf{e}_0, \textbf{e}_{[1, n]} \boldsymbol{\mid} r, s \big) = \bigg( \displaystyle\frac{s^2 z}{r^2} \bigg)^n \displaystyle\frac{(r^2; q)_n}{(s^2 z; q)_n}.
\end{flalign}

\noindent Indeed, under this weight, the $n$-tuple $\textbf{C} \in \{ 0, 1 \}^n$ from \Cref{wabcdrsxy} is equal to $\textbf{e}_0$. So, the sum over $\textbf{P}$ on the right side of \eqref{wabcdp} is supported on $\textbf{P} = \textbf{e}_0$, from which \eqref{wze} quickly follows. Applying \eqref{wze}, together with \eqref{wabcd2}, yields
\begin{flalign*}
\widehat{W}_{x_j / y_i} \big( \textbf{e}_{[1, n]}, \textbf{e}_0; \textbf{e}_0, \textbf{e}_{[1, n]} \boldsymbol{\mid} r_j, s_i \big) = \bigg( \displaystyle\frac{x_j}{r_j^2 y_i} \bigg)^n \displaystyle\frac{(r_j^2; q)_n}{(x_j y_i^{-1}; q)_n}.
\end{flalign*}

\noindent Inserting this, together with \eqref{fproductw}, into \eqref{w1w} yields 
\begin{flalign*}
F_{\boldsymbol{0}^N} (\textbf{x}; \textbf{r} \boldsymbol{\mid} \textbf{y}; \textbf{s}) & = \displaystyle\prod_{j = 1}^N s_j^{2n (j - N)} x_j^{n (N - j + 1)} r_j^{2n (j - N - 1)} y_j^{-jn} (r_j^2; q)_n \\
& \qquad \times \displaystyle\prod_{1 \le i < j \le N} (s_i^2 x_j y_i^{-1}; q)_n (r_i^2 x_i^{-1} y_j; q)_n \displaystyle\prod_{i = 1}^N \displaystyle\prod_{j = 1}^N (x_j y_i^{-1}; q)_n^{-1} \\
& = \displaystyle\prod_{j = 1}^N s_j^{2n (j - N)} x_j^{n (N - j + 1)} r_j^{2n (j - N - 1)} y_j^{-jn} (r_j^2; q)_n \\
& \qquad \times \displaystyle\prod_{1 \le i < j \le N} (r_i^2 x_i^{-1} x_j; q)_n (s_i^2 y_i^{-1} y_j; q)_n \displaystyle\prod_{i = 1}^N \displaystyle\prod_{j = 1}^N (x_j y_i^{-1}; q)_n^{-1},
\end{flalign*}

\noindent where to deduce the second equality we used the fact that $s_i^2 x_j y_i^{-1} = r_i^2 x_i^{-1} x_j$ and $r_i^2 x_i^{-1} y_j = s_i^2 y_i^{-1} y_j$ (as $s_i^2 x_i = r_i^2 y_i$). This implies the corollary. 
\end{proof}

Now we can quickly establish \Cref{f0c}.

\begin{proof}[Proof of \Cref{f0c}]
	
	This follows from \Cref{f0rxsy} and the fact that 
	\begin{flalign*}  
	\displaystyle\prod_{j = 1}^N & s_j^{2n (j - N)} x_j^{n (N - j + 1)} r_j^{2n (j - N - 1)} y_j^{-jn} (r_j^2; q)_n \\
	& \qquad \times \displaystyle\prod_{1 \le i < j \le N} (r_i^2 x_i^{-1} x_j; q)_n (s_i^2 y_i^{-1} y_j; q)_n \displaystyle\prod_{i = 1}^N \displaystyle\prod_{j = 1}^N (x_j y_i^{-1}; q)_n^{-1} \\
	& =  \displaystyle\prod_{j = 1}^N s_j^{2n (j - N)} r_j^{2n (j - N - 1)} (r_j^2; q)_n \\
	& \qquad \times \displaystyle\prod_{j = 1}^n x_j^n \displaystyle\prod_{1 \le i < j \le N} \displaystyle\prod_{k = 0}^{n - 1} (x_i - q^k r_i^2 x_j) (y_i - q^k s_i^2 y_j) \displaystyle\prod_{i = 1}^N \displaystyle\prod_{j = 1}^N \displaystyle\prod_{k = 0}^{n - 1} (y_i - q^k x_j)^{-1},
	\end{flalign*} 
	
	\noindent where we applied \eqref{productxyqk}.
\end{proof}

\chapter{Degenerations} 

\label{Degeneration} 

Although the weights $W_z$ from \Cref{wabcdrsxy} might at first glance appear a bit unpleasant, the fact that they are governed by four parameters $(r, s, z, q)$ makes them quite general. We will see in this chapter how imposing certain relations among these parameters drastically simplifies the $W_z$ weights, enabling them to factor completely. We will furthermore analyze the symmetric functions resulting from such simplifications and the Cauchy identities they satisfy.

\section{Restricting Horizontal Arrows} 

\label{Horizontal1} 

In this section we explain two degenerations of the $W_z$ weights from \Cref{wabcdrsxy} that are obtained by restricting the parameters so that at most one arrow is permitted along any horizontal edge. Throughout this section, we adopt the notation of \Cref{wabcdrsxy}.

The first way of doing this is to set $n = 1$, as in the following example.

\begin{example}
	
	\label{wn1} 
	
	Suppose $n = 1$, so that $W_z (\textbf{A}, \textbf{B}; \textbf{C}, \textbf{D} \boldsymbol{\mid} r, s) \ne 0$ only if $\textbf{A}, \textbf{B}, \textbf{C}, \textbf{D} \in \{ \textbf{e}_0, \textbf{e}_1 \}$. Abbreviating $W (i, j; i', j') = W_z \big( \textbf{e}_i, \textbf{e}_j; \textbf{e}_{i'}, \textbf{e}_{j'} \boldsymbol{\mid} r, s \big)$ for any $i, j, i', j' \in \{ 0, 1 \}$\, we have by \eqref{wabcdp} that 
	\begin{flalign*}
		&   W (0, 0; 0, 0) = 1; \qquad  \qquad \qquad \quad W (1, 0; 1, 0) = \displaystyle\frac{1 - s^2 r^{-2} z}{1 - s^2 z}; \qquad W (1, 0; 0, 1) = \displaystyle\frac{s^2 z (r^{-2} - 1)}{1 - s^2 z}; \\
		& W (1, 1; 1, 1) = \displaystyle\frac{s^2 (r^{-2} z - 1)}{1 - s^2 z}; \qquad W (0, 1; 0, 1) = \displaystyle\frac{s^2 (1 - z)}{1 - s^2 z}; \quad \qquad W (0, 1; 1, 0) = \displaystyle\frac{1 - s^2}{1 - s^2 z},
	\end{flalign*}	
	
	\noindent and $W_z (\textbf{A}, \textbf{B}; \textbf{C}, \textbf{D} \boldsymbol{\mid} r, s) = 0$ otherwise. This recovers the vertex weights for a generalized free-fermionic six-vertex model studied in \cite{AWM}. We refer to \Cref{n1weights} for a depiction.

\end{example}

\begin{figure}[t]
	
	\begin{center}
		
		\begin{tikzpicture}[
		>=stealth,
		scale = .75
		]

		\draw[-, black] (-7.5, -2.5) -- (7.5, -2.5);
		\draw[-, black] (-7.5, -1.1) -- (7.5, -1.1);
		\draw[-, black] (-7.5, -.4) -- (7.5, -.4);
		\draw[-, black] (-7.5, 2.4) -- (7.5, 2.4);
		\draw[-, black] (-7.5, -2.5) -- (-7.5, 2.4);
		\draw[-, black] (7.5, -2.5) -- (7.5, 2.4);
		\draw[-, black] (-5, -2.5) -- (-5, 2.4);
		\draw[-, black] (5, -2.5) -- (5, 2.4);
		\draw[-, black] (-2.5, -2.5) -- (-2.5, 2.4);
		\draw[-, black] (2.5, -2.5) -- (2.5, 2.4);
		\draw[-, black] (0, -2.5) -- (0, 2.4);

		\draw[-, dotted] (-7.15, 1) -- (-5.35, 1); 
		\draw[-, dotted] (-6.25, .1) -- (-6.25, 1.9); 
		
		\draw[->, thick] (3.75, .1) -- (3.75, 1) -- (4.65, 1);
		\draw[-, dotted] (2.85, 1) -- (3.75, 1) -- (3.75, 1.9); 
		
		\draw[->, thick] (-1.25, .1) -- (-1.25, 1.9);
		\draw[-,  dotted] (-2.15, 1) -- (-.35, 1);
		
		\draw[->, thick] (.35, 1) -- (2.15, 1);
		\draw[-, dotted] (1.25, .1) -- (1.25, 1.9);
		
		\draw[->, thick] (5.35, 1) -- (6.25, 1) -- (6.25, 1.9);
		\draw[-, dotted] (6.25, .1) -- (6.25, 1) -- (7.15, 1); 
		
		\draw[->, thick] (-3.75, .1) -- (-3.75, 1.9);
		\draw[->, thick] (-4.65, 1) -- (-2.85, 1);
		
		\filldraw[fill=gray!50!white, draw=black] (-6.25, -.75) circle [radius=0] node [black, scale = .9] {$(0, 0; 0, 0)$};
		\filldraw[fill=gray!50!white, draw=black] (-3.75, -.75) circle [radius=0] node [black, scale = .9] {$(1, 1; 1,1)$};
		\filldraw[fill=gray!50!white, draw=black] (-1.25, -.75) circle [radius=0] node [black, scale = .9] {$(1, 0; 1, 0)$};
		\filldraw[fill=gray!50!white, draw=black] (1.25, -.75) circle [radius=0] node [black, scale = .9] {$(0, 1; 0, 1)$};
		\filldraw[fill=gray!50!white, draw=black] (3.75, -.75) circle [radius=0] node [black, scale = .9] {$(1, 0; 0, 1)$};
		\filldraw[fill=gray!50!white, draw=black] (6.25, -.75) circle [radius=0] node [black, scale = .9] {$(0, 1; 1, 0)$};
		
		\filldraw[fill=gray!50!white, draw=black] (-6.25, -1.8) circle [radius=0] node [scale = .85, black] {$1$};
		\filldraw[fill=gray!50!white, draw=black] (-3.75, -1.8) circle [radius=0] node [scale = .85, black] {$\displaystyle\frac{s^2 (r^{-2} z - 1)}{1 - s^2 z}$};
		\filldraw[fill=gray!50!white, draw=black] (-1.25, -1.8) circle [radius=0] node [scale = .85, black] {$\displaystyle\frac{1 - s^2 r^{-2} z}{1 - s^2 z}$};
		\filldraw[fill=gray!50!white, draw=black] (1.25, -1.8) circle [radius=0] node [scale = .85, black] {$\displaystyle\frac{s^2 (1 - z)}{1 - s^2 z}$};
		\filldraw[fill=gray!50!white, draw=black] (3.75, -1.8) circle [radius=0] node [scale = .85, black] {$\displaystyle\frac{s^2 z (r^{-2} - 1)}{1 - s^2 z}$};
		\filldraw[fill=gray!50!white, draw=black] (6.25, -1.8) circle [radius=0] node [scale = .85, black] {$\displaystyle\frac{1 - s^2}{1 - s^2 z}$};

		\end{tikzpicture}
		
	\end{center}
	
	\caption{\label{n1weights} The  $n = 1$ cases of the $W_z$ weights, and their arrow configurations, are depicted above.}
\end{figure}

The second way of doing this is to set $r = q^{-1 / 2}$; this corresponds to the situation when $L = 1$ in \Cref{Weightss}, which again implies that any horizontal edge can be occupied by at most one arrow.

\begin{example}
	
	\label{rql1}

	Suppose $r = q^{-1 / 2}$. In this case, $W_z (\textbf{A}, \textbf{B}; \textbf{C}, \textbf{D} \boldsymbol{\mid} r, s) \ne 0$ only if $\textbf{A}, \textbf{C} \in \{ 0, 1 \}^n$ and $\textbf{B}, \textbf{D} \in \{ \textbf{e}_0, \textbf{e}_1, \ldots , \textbf{e}_n \}$, so let us assume this to be the case below. For any $i, j \in [0, n]$, abbreviate $W (\textbf{A}, i; \textbf{C}, j) = \textbf{W}_z \big( \textbf{A}, \textbf{e}_i; \textbf{C}, \textbf{e}_j \boldsymbol{\mid} r, s \big)$\index{W@$W_z (\textbf{A}, \textbf{B}; \textbf{C}, \textbf{D} \boldsymbol{\mid} r, s)$; fused weight!$W (\textbf{A}, i; \textbf{C}, j)$} and, as in (1.2.2) of \cite{SVMST}, set
	\begin{flalign}
	\label{aij} 
	\textbf{A}_i^+ = \textbf{A} + \textbf{e}_i; \qquad \textbf{A}_j^- = \textbf{A} - \textbf{e}_j; \qquad \textbf{A}_{ij}^{+-} = \textbf{A} + \textbf{e}_i - \textbf{e}_j.
	\end{flalign} 
	\index{A@$\textbf{A}_-^+, \textbf{A}_j^-, \textbf{A}_{ij}^{+-}$}

	\noindent Then, it follows as a direct consequence of \eqref{wabcdp} that for each $\textbf{A} \in \{ 0, 1 \}^n$ and $i \in [1, n]$, we have (letting $|\textbf{A}| = a$)
	\begin{flalign*}
	W (\textbf{A}, 0; \textbf{A}, 0) = & \displaystyle\frac{1 - q^a s^2 z}{1 - s^2 z}; \quad W (\textbf{A}, 0; \textbf{A}_i^-, i) = q^{A_{[i + 1, n]}} \displaystyle\frac{s^2 z (q - 1)}{1 - s^2 z}; \quad  W (\textbf{A}, i; \textbf{A}_i^+, 0) = \displaystyle\frac{1 - q^a s^2}{1 - s^2 z}.
	\end{flalign*} 
	
	\noindent Moreover, for any $1 \le i < j \le n$, we have
	\begin{flalign*} 
	W (\textbf{A}, i; \textbf{A}_{ij}^{+-}, j) = q^{A_{[j + 1, n]}} \displaystyle\frac{s^2 z (q - 1)}{1 - s^2 z}; \qquad W (\textbf{A}, j; \textbf{A}_{ji}^{+-}, i) = q^{A_{[i + 1, n]}} \displaystyle\frac{s^2 (q - 1)}{1 - s^2 z}.
	\end{flalign*}
	
	\noindent Additionally, 
	\begin{flalign*}
	& W (\textbf{A}, i; \textbf{A}, i) = q^{A_{[i + 1, n]}} \displaystyle\frac{s^2 (1 - z)}{1 - s^2 z}, \quad \text{if $A_i = 0$}; \qquad W (\textbf{A}, i; \textbf{A}; i) = q^{A_{[i + 1, n]}} \displaystyle\frac{s^2 (qz - 1)}{1 - s^2 z}, \quad \text{if $A_i = 1$}.
	\end{flalign*}	
	
	\noindent We further have $W_z (\textbf{A}, b; \textbf{C}, d \boldsymbol{\mid} r, s) = 0$ if $(\textbf{A}, b; \textbf{C}, d)$ is not of the above form, with $\textbf{A}, \textbf{C} \in \{ 0, 1 \}^n$. We refer to \Cref{vertexfigurerq} for a depiction.
	
\end{example}

\begin{rem}
	
	\label{wmnwmn} 
	
	All vertex weights $ W(\textbf{A}, i; \textbf{B}, j)$ in \Cref{vertexfigurerq}, except for the rightmost one there (correpsonding to when $i = j$ and $A_i = 1$), coincide with the $\tilde{L}_{z / s} (\textbf{A}, i; \textbf{B}, j)$ weights from equation (2.5.1) of \cite{SVMST} (see also equation (B.4.7) there) for the $U_q \big( \widehat{\mathfrak{sl}} (n + 1) \big)$-vertex model. This can be viewed as a consequence of \Cref{zrmn}, stating that the $U_q \big( \widehat{\mathfrak{sl}} (m | n) \big)$ weights $\mathcal{R}_{x, y}^{(m; n)} (\textbf{A}, \textbf{B}; \textbf{C}, \textbf{D})$ coincide with the $U_q \big( \widehat{\mathfrak{sl}} (m + n) \big)$ ones if $(\textbf{A}, \textbf{B}, \textbf{C}, \textbf{D})$ share no fermionic colors.
	
\end{rem}

\begin{figure}[t]
	
	\begin{center}
		
		\begin{tikzpicture}[
		>=stealth,
		scale = .85
		]

		\draw[-, black] (-7.5, 3.1) -- (10, 3.1);
		\draw[-, black] (-7.5, -2.5) -- (10, -2.5);
		\draw[-, black] (-7.5, -1.1) -- (10, -1.1);
		\draw[-, black] (-7.5, -.4) -- (10, -.4);
		\draw[-, black] (-7.5, 2.4) -- (10, 2.4);
		\draw[-, black] (-7.5, -2.5) -- (-7.5, 3.1);
		\draw[-, black] (10, -2.5) -- (10, 3.1);
		\draw[-, black] (7.5, -2.5) -- (7.5, 3.1);
		\draw[-, black] (-5, -2.5) -- (-5, 2.4);
		\draw[-, black] (5, -2.5) -- (5, 3.1);
		\draw[-, black] (-2.5, -2.5) -- (-2.5, 2.4);
		\draw[-, black] (2.5, -2.5) -- (2.5, 2.4);
		\draw[-, black] (0, -2.5) -- (0, 3.1);

		\draw[->, thick, blue] (-6.3, .1) -- (-6.3, 1.9); 
		\draw[->, thick, green] (-6.2, .1) -- (-6.2, 1.9); 
		
		\draw[->, thick, blue] (-3.8, .1) -- (-3.8, 1) -- (-2.85, 1);
		\draw[->, thick, green] (-3.7, .1) -- (-3.7, 1.9);
		
		\draw[->, thick, blue] (-1.35, .1) -- (-1.35, 1.9);
		\draw[->, thick, green] (-1.25, .1) -- (-1.25, 1.9);
		\draw[->, thick,  orange] (-2.15, 1.1) -- (-1.15, 1.1) -- (-1.15, 1.9);
		
		\draw[->, thick, red] (.35, 1) -- (1.15, 1) -- (1.15, 1.9);
		\draw[->, thick, blue] (1.25, .1) -- (1.25, 1.9);
		\draw[->, thick, green] (1.35, .1) -- (1.35, 1.1) -- (2.15, 1.1);
		
		\draw[->, thick, blue] (3.65, .1) -- (3.65, 1) -- (4.65, 1);
		\draw[->, thick, green] (3.75, .1) -- (3.75, 1.9);
		\draw[->, thick, orange] (2.85, 1.1) -- (3.85, 1.1) -- (3.85, 1.9); 
		
		\draw[->, thick, red] (5.35, 1) -- (7.15, 1); 
		\draw[->, thick, blue] (6.2, .1) -- (6.2, 1.9);
		\draw[->, thick, green] (6.3, .1) -- (6.3, 1.9); 
		
		\draw[->, thick, blue] (7.85, 1) -- (9.65, 1); 
		\draw[->, thick, blue] (8.7, .1) -- (8.7, 1.9);
		\draw[->, thick, green] (8.8, .1) -- (8.8, 1.9); 
	
		\filldraw[fill=gray!50!white, draw=black] (-2.85, 1) circle [radius=0] node [black, right = -1, scale = .75] {$i$};
		\filldraw[fill=gray!50!white, draw=black] (2.15, 1) circle [radius=0] node [black, right = -1, scale = .75] {$j$};
		\filldraw[fill=gray!50!white, draw=black] (4.65, 1) circle [radius=0] node [black, right = -1, scale = .75] {$i$};
		\filldraw[fill=gray!50!white, draw=black] (7.15, 1) circle [radius=0] node [black, right = -1, scale = .75] {$i$};
		\filldraw[fill=gray!50!white, draw=black] (9.65, 1) circle [radius=0] node [black, right = -1, scale = .75] {$i$};
		
		\filldraw[fill=gray!50!white, draw=black] (7.85, 1) circle [radius=0] node [black, left = -1, scale = .75] {$i$};
		\filldraw[fill=gray!50!white, draw=black] (5.35, 1) circle [radius=0] node [black, left = -1, scale = .75] {$i$};
		\filldraw[fill=gray!50!white, draw=black] (2.85, 1) circle [radius=0] node [black, left = -1, scale = .75] {$j$};
		\filldraw[fill=gray!50!white, draw=black] (.35, 1) circle [radius=0] node [black, left = -1, scale = .75] {$i$};
		\filldraw[fill=gray!50!white, draw=black] (-2.15, 1) circle [radius=0] node [black, left = -1, scale = .75] {$i$};
		
		\filldraw[fill=gray!50!white, draw=black] (-6.25, 1.9) circle [radius=0] node [black, above = -1, scale = .75] {$\textbf{A}$};
		\filldraw[fill=gray!50!white, draw=black] (-3.75, 1.9) circle [radius=0] node [black, above = -1, scale = .75] {$\textbf{A}_i^-$};
		\filldraw[fill=gray!50!white, draw=black] (-1.25, 1.9) circle [radius=0] node [black, above = -1, scale = .75] {$\textbf{A}_i^+$};
		\filldraw[fill=gray!50!white, draw=black] (1.25, 1.9) circle [radius=0] node [black, above = -1, scale = .75] {$\textbf{A}_{ij}^{+-}$};
		\filldraw[fill=gray!50!white, draw=black] (3.75, 1.9) circle [radius=0] node [black, above = -1, scale = .75] {$\textbf{A}_{ji}^{+-}$};
		\filldraw[fill=gray!50!white, draw=black] (6.25, 1.9) circle [radius=0] node [black, above = -1, scale = .75] {$\textbf{A}$};
		\filldraw[fill=gray!50!white, draw=black] (8.75, 1.9) circle [radius=0] node [black, above = -1, scale = .75] {$\textbf{A}$};

		\filldraw[fill=gray!50!white, draw=black] (-6.25, .1) circle [radius=0] node [black, below = -1, scale = .7] {$\textbf{A}$};
		\filldraw[fill=gray!50!white, draw=black] (-3.75, .1) circle [radius=0] node [black, below = -1, scale = .7] {$\textbf{A}$};
		\filldraw[fill=gray!50!white, draw=black] (-1.25, .1) circle [radius=0] node [black, below = -1, scale = .7] {$\textbf{A}$};
		\filldraw[fill=gray!50!white, draw=black] (1.25, .1) circle [radius=0] node [black, below = -1, scale = .7] {$\textbf{A}$};
		\filldraw[fill=gray!50!white, draw=black] (3.75, .1) circle [radius=0] node [black, below = -1, scale = .7] {$\textbf{A}$};
		\filldraw[fill=gray!50!white, draw=black] (6.25, .1) circle [radius=0] node [black, below = -1, scale = .7] {$\textbf{A}$};
		\filldraw[fill=gray!50!white, draw=black] (8.75, .1) circle [radius=0] node [black, below = -1, scale = .75] {$\textbf{A}$};

		\filldraw[fill=gray!50!white, draw=black] (-3.75, 2.75) circle [radius=0] node [black] {$1 \le i \le n$};
		\filldraw[fill=gray!50!white, draw=black] (2.5, 2.75) circle [radius=0] node [black] {$1 \le i < j \le n$}; 
		\filldraw[fill=gray!50!white, draw=black] (6.25, 2.75) circle [radius=0] node [black] {$A_i = 0$};
		\filldraw[fill=gray!50!white, draw=black] (8.75, 2.75) circle [radius=0] node [black] {$A_i = 1$};
		
		\filldraw[fill=gray!50!white, draw=black] (-6.25, -.75) circle [radius=0] node [black, scale = .9] {$(\textbf{A}, 0; \textbf{A}, 0)$};
		\filldraw[fill=gray!50!white, draw=black] (-3.75, -.75) circle [radius=0] node [black, scale = .9] {$\big( \textbf{A}, 0; \textbf{A}_i^-, i \big)$};
		\filldraw[fill=gray!50!white, draw=black] (-1.25, -.75) circle [radius=0] node [black, scale = .9] {$\big( \textbf{A}, i; \textbf{A}_i^+, 0 \big)$};
		\filldraw[fill=gray!50!white, draw=black] (1.25, -.75) circle [radius=0] node [black, scale = .9] {$\big( \textbf{A}, i; \textbf{A}_{ij}^{+-}, j \big)$};
		\filldraw[fill=gray!50!white, draw=black] (3.75, -.75) circle [radius=0] node [black, scale = .9] {$\big( \textbf{A}, j; \textbf{A}_{ji}^{+-}, i \big)$};
		\filldraw[fill=gray!50!white, draw=black] (6.25, -.75) circle [radius=0] node [black, scale = .9] {$(\textbf{A}, i; \textbf{A}, i)$};
		\filldraw[fill=gray!50!white, draw=black] (8.75, -.75) circle [radius=0] node [black, scale = .9] {$(\textbf{A}, i; \textbf{A}, i)$};
		
		\filldraw[fill=gray!50!white, draw=black] (-6.25, -1.8) circle [radius=0] node [black, scale = .85] {$\displaystyle\frac{1 - q^a s^2 z}{1 - s^2 z}$};
		\filldraw[fill=gray!50!white, draw=black] (-3.75, -1.8) circle [radius=0] node [black, scale = .7] {$q^{A_{[i + 1, n]}} \displaystyle\frac{s^2 z (q - 1)}{1 - s^2 z}$};
		\filldraw[fill=gray!50!white, draw=black] (-1.25, -1.8) circle [radius=0] node [black, scale = .85] {$\displaystyle\frac{1 - q^a s^2}{1 - s^2 z}$};
		\filldraw[fill=gray!50!white, draw=black] (1.25, -1.8) circle [radius=0] node [black, scale = .7] {$q^{A_{[j + 1, n]}} \displaystyle\frac{s^2 z (q - 1)}{1 - s^2 z}$};
		\filldraw[fill=gray!50!white, draw=black] (3.75, -1.8) circle [radius=0] node [black, scale = .7] {$q^{A_{[i + 1, n]}} \displaystyle\frac{s^2 (q - 1)}{1 - s^2 z}$};
		\filldraw[fill=gray!50!white, draw=black] (6.25, -1.8) circle [radius=0] node [black, scale = .7] {$q^{A_{[i + 1, n]}} \displaystyle\frac{s^2 (1 - z)}{1 - s^2 z}$};
		\filldraw[fill=gray!50!white, draw=black] (8.75, -1.8) circle [radius=0] node [black, scale = .7] {$q^{A_{[i + 1, n]}} \displaystyle\frac{s^2 (qz - 1)}{1 - s^2 z}$};

		\end{tikzpicture}
		
	\end{center}
	
	\caption{\label{vertexfigurerq} The  $r= q^{-1 / 2}$ cases of the $W_z$ weights, and their arrow configurations, are depicted above.}
\end{figure}

\section{Additional Simplifications} 

\label{DegenerationWeights}

In this section we provide four additional special choices of the parameters $(r, s, z)$ under which the sums from \eqref{wabcdp} describing the weights $W_z (\textbf{A}, \textbf{B}; \textbf{C}, \textbf{D} \boldsymbol{\mid} r, s)$ again simplify into factored products, but which now permit multiple arrows to exist along horizontal edges. Throughout this section, we adopt the notation of \Cref{wabcdrsxy}.

Before detailing these degenerations, we begin with the following lemma that will be useful in their derivations. In what follows, we recall the function $\varphi$ from \eqref{tufunction}. 

\begin{lem}
	
	\label{xyfunction} 
	
	Fix $\textbf{\emph{X}}, \textbf{\emph{Y}} \in \{ 0, 1 \}^n$, whose coordinates are both indexed by $[1, n]$, and define $\textbf{\emph{Z}} = (Z_1, Z_2, \ldots , Z_n) \in \{ 0, 1 \}^n$ by setting $Z_i = \min \{ X_i, Y_i \}$ for each $i \in [1, n]$. Then, 
	\begin{flalign}
	\label{xyz} 
	\varphi (\textbf{\emph{X}}, \textbf{\emph{Y}}) + \varphi (\textbf{\emph{Y}}, \textbf{\emph{X}}) = |\textbf{\emph{X}}| |\textbf{\emph{Y}}| - |\textbf{\emph{Z}}|.
	\end{flalign} 
	
	\noindent In particular, recalling that $\textbf{\emph{e}}_{[1, n]} = (1, 1, \ldots , 1) \in \{ 0, 1 \}^n$, we have
	\begin{flalign}
	\label{xn1}
	\varphi \big( \textbf{\emph{X}}, \textbf{\emph{e}}_{[1, n]} \big) + \varphi \big( \textbf{\emph{e}}_{[1, n]}, \textbf{\emph{X}} \big) = (n - 1) |\textbf{\emph{X}}|.
	\end{flalign}
\end{lem}

\begin{proof}
	
	Observe that \eqref{xn1} follows from the $\textbf{Y} = \textbf{e}_{[1, n]}$ case of \eqref{xyz} (since then $\textbf{Z} = \textbf{X}$ and $|\textbf{Y}| = n$), so it suffices to establish the latter. To that end, observe from the definition \eqref{tufunction} of $\varphi$ that 
	\begin{flalign*}
	\varphi (\textbf{X}, \textbf{Y}) + \varphi (\textbf{Y}, \textbf{X}) & = \displaystyle\sum_{1 \le i < j \le n} X_i Y_j + \displaystyle\sum_{1 \le j < i \le n} X_i Y_j \\
	& = \displaystyle\sum_{j = 1}^n Y_j \Bigg( \displaystyle\sum_{i \ne j} X_i \Bigg) = \displaystyle\sum_{j = 1}^n Y_j \big( |\textbf{X}| - X_j \big) = |\textbf{X}| |\textbf{Y}| - \displaystyle\sum_{j = 1}^n X_j Y_j = |\textbf{X}| |\textbf{Y}| -  |\textbf{Z}|,
	\end{flalign*} 
	
	\noindent where in the last equality we used the fact that $Z_j = \min \{ X_j, Y_j \} = X_j Y_j$ for each $j \in [1, n]$ (as $X_j, Y_j \in \{ 0, 1 \}$). 
\end{proof}

Now, the following three propositions address simplifications of $W_z$ under certain relations between the three parameters $(r, s, z) \in \mathbb{C}$. The source of these simplifications comes from the three factors $(q^{-v} s^2 r^{-2} z; q)_{c - p}$, $(z; q)_{b - p - v}$, and $(q^{1 - v} r^{-2} z; q)_v (q^v r^2 z^{-1}; q)_p$ appearing on the right side of \eqref{wabcdp}. If $s^2 r^{-2} z = 1$, $z = 1$, or $r^2 z^{-1} = 1$, then the first, second, or third of these quantities will be supported on a single value of $p$, respectively,\footnote{For the first quantity $(q^{-v} s^2 r^{-2}; q)_{c - p}$, this follows from the fact that we additionally have $c - v \ge p$.} thereby giving rise to a factored form for $W_z$.

We begin with the $s^2 r^{-2} z = 1$ specialization, which is similar to the $q$-Hahn type specialization originally explained in Proposition 6.7 of \cite{SRF} and Proposition 6.11 of \cite{HSVMSRF} for the $U_q \big( \widehat{\mathfrak{sl}} (2) \big)$ vertex model, and later in Proposition 7 of \cite{SRM} for the more general $U_q \big( \widehat{\mathfrak{sl}} (n) \big)$ vertex model.

\begin{prop}
	
	\label{sxz1weights}
	
	For any $s, z \in \mathbb{C}$, we have\footnote{In the below, and throughout this text, quantities such as $z^{1 / 2}$ will not depend on the choice for the root of $z$.}
	\begin{flalign}
	\label{wszs} 
	W_z ( \textbf{\emph{A}}, \textbf{\emph{B}}; \textbf{\emph{C}}, \textbf{\emph{D}} \boldsymbol{\mid} s z^{1/2}, s) = s^{2 (b - c)} q^{\varphi (\textbf{\emph{B}} - \textbf{\emph{C}}, \textbf{\emph{C}})} \displaystyle\frac{(s^2; q)_c (z; q)_{b - c}}{(s^2 z; q)_b} \textbf{\emph{1}}_{\textbf{\emph{A}} + \textbf{\emph{B}} = \textbf{\emph{C}} + \textbf{\emph{D}}} \textbf{\emph{1}}_{\textbf{\emph{B}} \ge \textbf{\emph{C}}}.
	\end{flalign}
\end{prop}

\begin{proof}
	
	We may assume in what follows that $\textbf{A} + \textbf{B} = \textbf{C} + \textbf{D}$, for otherwise both sides of \eqref{wszs} are equal to $0$ by arrow conservation. We may further assume by \Cref{szr0} that $\textbf{B} \ge \textbf{C}$. Inserting $r = s z^{1 / 2}$ into \eqref{wabcdp}, we have that
	\begin{flalign}
	\label{wzrszsum2}
	\begin{aligned} 
	W_z (\textbf{A}, \textbf{B}; \textbf{C}, \textbf{D} \boldsymbol{\mid} s z^{1/2}; s) &= (-1)^v z^{c + d - a - b} s^{2c + 2d - 2a} q^{\varphi (\textbf{D} - \textbf{V}, \textbf{C})+ \varphi (\textbf{V}, \textbf{A}) - av + cv} \\
	& \quad \times \displaystyle\frac{(q^{1 - v} s^{-2}; q)_v}{(q^{-v}; q)_v} \displaystyle\frac{(s^2 z; q)_d}{(s^2 z; q)_b} \displaystyle\sum_{p = 0}^{\min \{ b - v, c - v \}} \displaystyle\frac{(q^{-v}; q)_{c - p} (q^v s^2; q)_p (z; q)_{b - p - v}}{(s^2 z; q)_{c + d - p - v}} \\
	& \qquad \qquad \qquad \times q^{-pv} s^{-2p} \displaystyle\sum_{\textbf{\emph{P}}} q^{\varphi (\textbf{B} - \textbf{D} - \textbf{P}, \textbf{P})}, 
	\end{aligned}
	\end{flalign}
	
	\noindent where the sum is over all $n$-tuples $\textbf{P} = (P_1, P_2, \ldots , P_n) \in \{ 0, 1 \}^n$ such that $|\textbf{P}| = p$ and $P_i \le \min \{ B_i - V_i, C_i -V_i \}$, for each $i \in [1, n]$. As indicated at the end of the proof of \Cref{szr0}, the sum on the right side of \eqref{wzrszsum2} is supported on the term $\textbf{P} = \textbf{C} - \textbf{V}$. Inserting this into \eqref{wzrszsum2}, and using the fact that $a + b = c + d$, yields 
	\begin{flalign*}
	W_z (\textbf{A}, \textbf{B}; \textbf{C}, \textbf{D} \boldsymbol{\mid} s z^{1 / 2}; s) & = (-1)^v s^{2b - 2c + 2v} q^{\varphi (\textbf{D} - \textbf{V}, \textbf{C}) + \varphi (\textbf{V}, \textbf{A}) - av + v^2} \displaystyle\frac{(q^{1 - v} s^{-2}; q)_v}{(s^2 z; q)_b} \\
	& \qquad \times (q^v s^2; q)_{c - v} (z; q)_{b - c} q^{\varphi (\textbf{B} - \textbf{D} - \textbf{C} + \textbf{V}, \textbf{C} - \textbf{V})}.
	\end{flalign*}
	
	\noindent Using the identities $\textbf{A} + \textbf{B} = \textbf{C} + \textbf{D}$ and
	\begin{flalign*}
	(q^{1 - v} s^{-2}; q)_v = (-1)^v q^{-\binom{v}{2}} s^{-2v} (s^2; q)_v,
	\end{flalign*}
	
	\noindent it follows that
	\begin{flalign}
	\label{wabcdzsz}
	\begin{aligned} 
	W_z (\textbf{A}, \textbf{B}; \textbf{C}, \textbf{D} \boldsymbol{\mid} s z^{1 / 2}; s) & = s^{2b - 2c} q^{\varphi (\textbf{D} - \textbf{V}, \textbf{C}) + \varphi (\textbf{V}, \textbf{A}) - av + v^2} \displaystyle\frac{(z; q)_{b - c}}{(s^2 z; q)_b} \\
	& \qquad \times (s^2; q)_v (q^v s^2; q)_{c - v} q^{\varphi (\textbf{V} - \textbf{A}, \textbf{C} - \textbf{V}) - \binom{v}{2}}.
	\end{aligned} 
	\end{flalign}
	
	\noindent Next, by the bilinearity of $\varphi$, we have that 
	\begin{flalign}
	\label{sumdvc}
	\begin{aligned}
	\varphi (\textbf{D} - \textbf{V}, \textbf{C}) + \varphi (\textbf{V}, \textbf{A}) + \varphi (\textbf{V} - \textbf{A}, \textbf{C} - \textbf{V}) & = \varphi (\textbf{D} - \textbf{A}, \textbf{C}) + \varphi (\textbf{V}, \textbf{A}) - \varphi (\textbf{V}, \textbf{V}) + \varphi (\textbf{A}, \textbf{V}) \\
	& = \varphi ( \textbf{B} - \textbf{C}, \textbf{C}) + av - v - \binom{v}{2},
	\end{aligned}
	\end{flalign}
	
	\noindent where in the second statement we have used the facts that $\textbf{A} + \textbf{B} = \textbf{C} + \textbf{D}$, that $\varphi (\textbf{V}, \textbf{V}) = \binom{v}{2}$, and that $\varphi (\textbf{A}, \textbf{V}) + \varphi (\textbf{V}, \textbf{A}) = av - v$ (where the last statement holds by \Cref{xyfunction}, since $\textbf{V} \le \textbf{A}$). 

	Now the proposition follows from inserting \eqref{sumdvc} into \eqref{wabcdzsz} and using the identities $v^2 = 2 \binom{v}{2} + v$ and $(s^2; q)_v (q^v s^2; q)_{c - v} = (s^2; q)_c$.
\end{proof}

Next we consider the case $z = 1$.

\begin{prop}
	
	\label{z1weights} 
	
	For any $r, s \in \mathbb{C}$, we have
	\begin{flalign}
	\label{w1abcdrs}
	W_1 ( \textbf{\emph{A}}, \textbf{\emph{B}}; \textbf{\emph{C}}, \textbf{\emph{D}} \boldsymbol{\mid} r, s) = (r^{-2} s^2)^d (r^2; q)_d \displaystyle\frac{(r^{-2} s^2; q)_{c - b}}{(s^2; q)_{c + d - b}} q^{\varphi (\textbf{\emph{D}}, \textbf{\emph{C}} - \textbf{\emph{B}})} \textbf{\emph{1}}_{\textbf{\emph{A}} + \textbf{\emph{B}} = \textbf{\emph{C}} + \textbf{\emph{D}}} \textbf{\emph{1}}_{\textbf{\emph{C}} \ge \textbf{\emph{B}}}.
	\end{flalign}
\end{prop}

\begin{proof}
	
	We may again assume in what follows that $\textbf{A} + \textbf{B} = \textbf{C} + \textbf{D}$, for otherwise arrow conservation implies $W_z (\textbf{A}, \textbf{B}; \textbf{C}, \textbf{D} \boldsymbol{\mid} r, s) = 0$ for any $r, s, z \in \mathbb{C}$. Then setting $z = 1$ in \eqref{wabcdp}, we find that the factor $(z; q)_{b - p - v} = \textbf{1}_{p = b - v}$ on the right side there implies that the sum is supported on the term $\textbf{P} = \textbf{B} - \textbf{V}$. In particular, since we must have $P_i \le \min \{ B_i - V_i, C_i - V_i \}$, this implies that $W_1 (\textbf{A}, \textbf{B}; \textbf{C}, \textbf{D} \boldsymbol{\mid} r, s) = 0$ unless $\textbf{C} \ge \textbf{B}$. Thus, we will assume in in what follows that $\textbf{C} \ge \textbf{B}$, for otherwise both sides of \eqref{w1abcdrs} are equal to $0$.
	
	Then inserting $z = 1$ into \eqref{wabcdp} (and using that it is supported on $\textbf{P} = \textbf{B} - \textbf{V}$) yields
	\begin{flalign*}
	W_1 (\textbf{A}, \textbf{B}; \textbf{C}, \textbf{D} \boldsymbol{\mid} r, s) &= (-1)^v r^{2c - 2a} s^{2d} q^{\varphi (\textbf{D} - \textbf{V}, \textbf{C})+ \varphi (\textbf{V}, \textbf{A}) - av + cv} \displaystyle\frac{(q^{1 - v} r^{-2}; q)_v}{(q^{-v} s^2 r^{-2}; q)_v} \displaystyle\frac{(r^2; q)_d}{(r^2; q)_b} \\
	& \qquad \times \displaystyle\frac{(q^{-v} s^2 r^{-2}; q)_{c + v - b} (q^v r^2; q)_{b - v}}{(s^2; q)_{c + d - b}} (q^{-v} r^{-2})^{b - v} q^{\varphi (\textbf{V} - \textbf{D}, \textbf{B} - \textbf{V})}.
	\end{flalign*}
	
	\noindent Using the identities
	\begin{flalign*}
	& \displaystyle\frac{(q^{1 -  v} r^{-2}; q)_v (q^v r^2; q)_{b - v}}{(r^2; q)_b} = \displaystyle\frac{(q^{1 - v} r^{-2}; q)_v}{(r^2; q)_v} = (-1)^v r^{-2v} q^{-\binom{v}{2}}; \\
	& \displaystyle\frac{(q^{-v} s^2 r^{-2}; q)_{c + v - b}}{(q^{-v} s^2  r^{-2}; q)_v} = (r^{-2} s^2; q)_{c - b},
	\end{flalign*}
	
	\noindent it follows that 
	\begin{flalign}
	\label{wabcd1rs}
	\begin{aligned}
	W_1 (\textbf{A}, \textbf{B}; \textbf{C}, \textbf{D} \boldsymbol{\mid} r, s) &= r^{2c - 2a - 2b} s^{2d} q^{\varphi (\textbf{D} - \textbf{V}, \textbf{C})+ \varphi (\textbf{V}, \textbf{A}) + \varphi (\textbf{V} - \textbf{D}, \textbf{B} - \textbf{V}) - av - bv + cv + v^2 - \binom{v}{2}} \\
	& \qquad \times (r^2; q)_d \displaystyle\frac{(r^{-2} s^2; q)_{c - b}}{(s^2; q)_{c + d - b}}.
	\end{aligned}
	\end{flalign}
	
	\noindent Next, by the bilinearity of $\varphi$, we find that 
	\begin{flalign}
	\label{dvcsum} 
	\begin{aligned}
	\varphi (\textbf{D} - \textbf{V}, \textbf{C}) & + \varphi (\textbf{V}, \textbf{A}) + \varphi (\textbf{V} - \textbf{D}, \textbf{B} - \textbf{V}) \\
	& = \varphi (\textbf{D}, \textbf{C} - \textbf{B}) + \varphi (\textbf{V}, \textbf{A} + \textbf{B} - \textbf{C}) + \varphi (\textbf{D}, \textbf{V}) - \varphi (\textbf{V}, \textbf{V}) \\
	& = \varphi (\textbf{D}, \textbf{C} - \textbf{B}) + \varphi (\textbf{V}, \textbf{D}) + \varphi (\textbf{D}, \textbf{V}) - \varphi (\textbf{V}, \textbf{V}),
	\end{aligned}
	\end{flalign} 
	
	\noindent where in the last equality we used the fact that $\textbf{A} + \textbf{B} = \textbf{C} + \textbf{D}$. Since $\textbf{V} \le \textbf{D}$, \Cref{xyfunction} implies that $\varphi (\textbf{V}, \textbf{D}) + \varphi (\textbf{D}, \textbf{V}) = |\textbf{D}| |\textbf{V}| - |\textbf{V}| = dv - v$. Inserting this and the fact that $\varphi (\textbf{V}, \textbf{V}) = \binom{v}{2}$ into \eqref{dvcsum} yields
	\begin{flalign*}
	\varphi (\textbf{D} - \textbf{V}, \textbf{C}) & + \varphi (\textbf{V}, \textbf{A}) + \varphi (\textbf{V} - \textbf{D}, \textbf{B} - \textbf{V}) = \varphi (\textbf{D}, \textbf{C} - \textbf{B}) + dv - v - \binom{v}{2}.
	\end{flalign*}
	
	\noindent The proposition then follows from inserting this into \eqref{wabcd1rs} and applying the identities $a + b = c + d$ and $v^2 = 2 \binom{v}{2} + v$.
\end{proof}

Now we address the specialization $r^2 z^{-1} = 1$.

\begin{prop}
	
	\label{zrw}
	
	For any $s, z \in \mathbb{C}$, we have 
	\begin{flalign*}
	W_z (\textbf{\emph{A}}, \textbf{\emph{B}}; \textbf{\emph{C}}, \textbf{\emph{D}} \boldsymbol{\mid} z^{1 / 2}, s) = s^{2d} q^{\varphi (\textbf{\emph{D}}, \textbf{\emph{C}})} \displaystyle\frac{(s^2; q)_c (z; q)_d}{(s^2 z; q)_{c + d}} \textbf{\emph{1}}_{\textbf{\emph{A}} + \textbf{\emph{B}} = \textbf{\emph{C}} + \textbf{\emph{D}}} \textbf{\emph{1}}_{v = 0}.
	\end{flalign*} 
\end{prop}

\begin{proof} 
	
	As in the proofs of \Cref{sxz1weights} and \Cref{z1weights}, we may assume that $\textbf{A} + \textbf{B} = \textbf{C} + \textbf{D}$. Then, setting $r = z^{1/2}$ in \eqref{wabcdp} yields
	\begin{flalign*}
	W_z (\textbf{A}, \textbf{B}; \textbf{C}, \textbf{D} \boldsymbol{\mid} z^{1 / 2}, s) & = (-1)^v z^{c + d - a - b} s^{2d} q^{\varphi (\textbf{D} - \textbf{V}, \textbf{C}) + \varphi (\textbf{V}, \textbf{A}) - av + cv} \displaystyle\frac{(q^{1 - v}; q)_v}{(q^{- v} s^2; q)_v} \displaystyle\frac{(z; q)_d}{(z; q)_b}  \\
	& \quad \times  \displaystyle\sum_{p = 0}^{\min \{ b - v, c - v \}}  \displaystyle\frac{(q^{- v} s^2; q)_{c - p} (q^v; q)_p (z; q)_{b - p - v}}{(s^2 z; q)_{c + d - p - v}} q^{-vp} \displaystyle\sum_{\textbf{P}} q^{\varphi (\textbf{B} - \textbf{D} - \textbf{P}, \textbf{P})},
	\end{flalign*}
	
	\noindent where the last sum sum is over all $n$-tuples $\textbf{P} = (P_1, P_2, \ldots , P_n) \in \{ 0, 1 \}^n$ such that $|\textbf{P}| = p$ and $P_i \le \min \{ B_i - V_i, C_i - V_i \}$ for each $i \in [1, n]$. Since $(q^{1 - v}; q)_v = \textbf{1}_{v = 0}$ and since $a + b = c + d$ (as $\textbf{A} + \textbf{B} = \textbf{C} + \textbf{D})$, it follows that 
	\begin{flalign*}
	W_z (\textbf{A}, \textbf{B}; \textbf{C}, \textbf{D} \boldsymbol{\mid} z^{1 / 2}, s) = s^{2d} q^{\varphi (\textbf{D}, \textbf{C})} \displaystyle\frac{(z; q)_d}{(z; q)_b} \textbf{1}_{v = 0} \displaystyle\sum_{p = 0}^{\min \{ b, c \}} & \displaystyle\frac{(s^2; q)_{c - p} (q^v; q)_p (z; q)_{b - p}}{(s^2 z; q)_{c + d - p}} \\
	& \times \displaystyle\sum_{\textbf{P}} q^{\varphi (\textbf{B} - \textbf{D} - \textbf{P}, \textbf{P})}.
	\end{flalign*}
	
	\noindent As $\textbf{1}_{v = 0} (q^v; q)_p = \textbf{1}_{v = 0} \textbf{1}_{p = 0}$, it follows that the above sum is supported on the term $\textbf{P} = \textbf{e}_0$, from which we deduce the proposition. 
\end{proof}  

A fourth simplification arises in the limiting regime when $(s, z)$ tends to $(0, \infty)$ in such a way that $x = sz$ is fixed. Ensuring the existence of this limit will require normalizing the weight by a factor of $(-s)^{-d}$.

\begin{prop}
	
	\label{ws0}
	
	For any fixed $x \in \mathbb{C}$, we have 
	\begin{flalign} 
		\label{wxsr} 
		\begin{aligned}
			\displaystyle\lim_{s \rightarrow 0} (-s)^{-d} W_{x / s} (\textbf{\emph{A}}, \textbf{\emph{B}}; \textbf{\emph{C}}, \textbf{\emph{D}} \boldsymbol{\mid} r, s) & = x^d (-r^2)^{c - a - v} q^{\varphi (\textbf{\emph{D}} - \textbf{\emph{V}}, \textbf{\emph{C}}) + \varphi (\textbf{\emph{V}}, \textbf{\emph{A}}) + \binom{b}{2} - dv + v} \displaystyle\frac{(r^2; q)_d}{(r^2; q)_b} \\
			& \quad \times \textbf{\emph{1}}_{\textbf{\emph{A}} + \textbf{\emph{B}} = \textbf{\emph{C}} + \textbf{\emph{D}}} \displaystyle\prod_{j : B_j - D_j = 1} \big(1 - r^{-2} q^{-B_{[j + 1, n]} - D_{[1, j - 1]}} \big).
		\end{aligned} 
	\end{flalign}

\end{prop}

\begin{proof}
	
	As in the proof of \Cref{sxz1weights}, we may assume in what follows that $\textbf{A} + \textbf{B} = \textbf{C} + \textbf{D}$. Then, replacing the $z$ in \Cref{wabcdrsxy} with $s^{-1} x$ yields 
	\begin{flalign}
		\label{wzs1}
		\begin{aligned}
			(-s)^d W_{x / s} (\textbf{A}, \textbf{B}; \textbf{C}, \textbf{D} \boldsymbol{\mid} r, s) & = (-1)^{v + d} s^b x^{d - b} r^{2c - 2a} q^{\varphi (\textbf{D} - \textbf{V}, \textbf{C}) + \varphi (\textbf{V}, \textbf{A}) - av + cv} \displaystyle\frac{(q^{1 - v} r^{-2} s^{-1} x; q)_v}{(q^{- v} s r^{-2} x; q)_v} \\
			& \quad \times \displaystyle\frac{(r^2; q)_d}{(r^2; q)_b} \displaystyle\sum_{p = 0}^{\min \{ b - v, c - v \}}  \displaystyle\frac{(q^{- v} s r^{-2} x; q)_{c - p} (q^v r^2 s x^{-1}; q)_p}{(s x; q)_{c + d - p - v}} \\
			& \qquad \qquad \times (q^{-v} r^{-2} s^{-1} x)^p (s^{-1} x; q)_{b - p - v} \displaystyle\sum_{\textbf{P}} q^{\varphi (\textbf{B} - \textbf{D} - \textbf{P}, \textbf{P})},
		\end{aligned} 
	\end{flalign} 
	
	\noindent where the sum is over all $n$-tuples $\textbf{P} = (P_1, P_2, \ldots , P_n) \in \{ 0, 1 \}^n$ such that $|\textbf{P}| = p$ and $P_i \le \min \{ B_i - V_i, C_i - V_i \}$ for each $i \in [1, n]$. Next, observe that 
	\begin{flalign*}
		\displaystyle\lim_{s \rightarrow 0} (-s)^v (q^{1 - v} r^{-2} s^{-1} x; q)_v & = q^{-\binom{v}{2}} r^{-2v} x^v; \\
		\displaystyle\lim_{s \rightarrow 0} s^{b - p - v} (s^{-1} x; q)_{b - p - v} & = q^{\binom{b - p - v}{2}} (-1)^{b - p - v} x^{b - p - v},
	\end{flalign*}
	
	\noindent which upon insertion into \eqref{wzs1} gives 
	\begin{flalign}
		\label{wzs2} 
		\begin{aligned}
			\displaystyle\lim_{s \rightarrow 0} (-s)^d W_{x / s} (\textbf{A}, \textbf{B}; \textbf{C}, \textbf{D} \boldsymbol{\mid} r, s) & = (-1)^{v + b + d} x^d r^{2c - 2a - 2v} q^{\varphi (\textbf{D} - \textbf{V}, \textbf{C}) + \varphi (\textbf{V}, \textbf{A}) - av + cv -\binom{v}{2}} \\
			& \quad \times \displaystyle\frac{(r^2; q)_d}{(r^2; q)_b} \displaystyle\sum_{p = 0}^{\min \{ b - v, c - v \}} (-r^{-2})^p q^{\binom{b - p - v}{2} - pv} \displaystyle\sum_{\textbf{P}} q^{\varphi (\textbf{B} - \textbf{D} - \textbf{P}, \textbf{P})}.
		\end{aligned} 
	\end{flalign} 
	
	\noindent By the bilinearity of $\varphi$ and the fact that $\varphi (\textbf{X}, \textbf{X}) = \binom{x}{2}$ whenever $|\textbf{X}| = x$, we have 
	\begin{flalign*}
		\binom{b - p - v}{2} & - \binom{v}{2} + \varphi (\textbf{B} - \textbf{D} - \textbf{P}, \textbf{P}) \\
		& = \varphi (\textbf{B} - \textbf{P} - \textbf{V}, \textbf{B} - \textbf{P} - \textbf{V}) - \varphi (\textbf{V}, \textbf{V}) +  \varphi (\textbf{B} - \textbf{D} - \textbf{P}, \textbf{P}) \\
		& =  \varphi (\textbf{B} - \textbf{P}, \textbf{B} - \textbf{P}) - \varphi (\textbf{B} - \textbf{P}, \textbf{V}) - \varphi (\textbf{V}, \textbf{B} - \textbf{P}) + \varphi (\textbf{B} - \textbf{P}, \textbf{P}) - \varphi (\textbf{D}, \textbf{P}) \\
		& = \binom{b}{2} - \varphi (\textbf{P}, \textbf{B}) - (b - p - 1) v - \varphi (\textbf{D}, \textbf{P}),
	\end{flalign*}
	
	\noindent where in the last equality we also applied \Cref{xyfunction} (with the fact that $\textbf{P} \le \textbf{B} - \textbf{V}$). Upon insertion into \eqref{wzs2} and using the fact that $a + b = c + d$, this yields 
	\begin{flalign}
		\label{sw0} 
		\begin{aligned} 
			\displaystyle\lim_{s \rightarrow 0} (-s)^{-d} W_{x / s} (\textbf{A}, \textbf{B}; \textbf{C}, \textbf{D} \boldsymbol{\mid} r, s) & = x^d (-r^2)^{c - a - v} q^{\varphi (\textbf{D} - \textbf{V}, \textbf{C}) + \varphi (\textbf{V}, \textbf{A}) + \binom{b}{2} - dv + v} \displaystyle\frac{(r^2; q)_d}{(r^2; q)_b} \\
			& \qquad \times \displaystyle\sum_{p = 0}^{\min \{ b - v, c - v \}} (-r^{-2})^p \displaystyle\sum_{\textbf{P}} q^{- \varphi (\textbf{P}, \textbf{B}) - \varphi (\textbf{D}, \textbf{P})}.
		\end{aligned} 
	\end{flalign} 
	
	Next, observe again by the bilinearity of $\varphi$ and \Cref{xyfunction}, applied with the $(\textbf{X}, \textbf{Y})$ there equal to $(\textbf{P}, \textbf{V})$ here, gives
	\begin{flalign*}
		\varphi (\textbf{P}, \textbf{B}) + \varphi (\textbf{D}, \textbf{P}) & = \varphi (\textbf{P}, \textbf{B} - \textbf{V}) + \varphi (\textbf{D} - \textbf{V}, \textbf{P}) + \varphi (\textbf{P}, \textbf{V}) + \varphi (\textbf{V}, \textbf{P}) \\
		& = \varphi (\textbf{P}, \textbf{B} - \textbf{V}) + \varphi (\textbf{D} - \textbf{V}, \textbf{P}) + pv,
	\end{flalign*}
	
	\noindent where we have further used the fact that the $\textbf{Z}$ in \Cref{xyfunction} is equal to $\textbf{e}_0$ here since $\textbf{P} \le \textbf{B} - \textbf{V}$. Inserting this into \eqref{sw0} then gives
	\begin{flalign}
		\label{1sw0} 
		\begin{aligned}
			\displaystyle\lim_{s \rightarrow 0} (-s)^{-d} W_{x / s} (\textbf{A}, \textbf{B}; \textbf{C}, \textbf{D} \boldsymbol{\mid} r, s) & = x^d (-r^2)^{c - a - v} q^{\varphi (\textbf{D} - \textbf{V}, \textbf{C}) + \varphi (\textbf{V}, \textbf{A}) + \binom{b}{2} - dv + v} \displaystyle\frac{(r^2; q)_d}{(r^2; q)_b} \\
			& \qquad \displaystyle\sum_{p = 0}^{\min \{ b - v, c - v \}} (-q^{-v} r^{-2})^p \displaystyle\sum_{\textbf{P}} q^{\varphi (\textbf{P}, \textbf{V} - \textbf{B}) + \varphi (\textbf{V} - \textbf{D}, \textbf{P})}.
		\end{aligned} 
	\end{flalign} 
	
	\noindent Now observe for any $w \in \mathbb{C}$ and $\textbf{X}, \textbf{Y}, \textbf{Z} \in \{ 0, 1 \}^n$ (whose coordinates are indexed by $[1, n]$) that
	\begin{flalign}
		\label{sumuq}
		\displaystyle\sum_{\substack{\textbf{P} \in \{ 0, 1 \}^n \\ \textbf{P} \le \textbf{Z}}} w^{|\textbf{P}|} q^{-\varphi (\textbf{P}, \textbf{X}) - \varphi (\textbf{Y}, \textbf{P})} = \displaystyle\prod_{j : Z_j = 1} \big( 1 + q^{-X_{[j + 1, n]} - Y_{[1, j - 1]}} w \big),
	\end{flalign}
	
	\noindent which holds since \eqref{tufunction} implies $\varphi (\textbf{P}, \textbf{X}) = \sum_{j: P_j = 1}^n X_{[j + 1, n]}$ and $\varphi (\textbf{Y}, \textbf{P}) = \sum_{j : P_j = 1} Y_{[1, j - 1]}$. Applying \eqref{1sw0}, together with the
	\begin{flalign*} 
		(w, \textbf{X}, \textbf{Y}, \textbf{Z}) = \big( - q^{-v} r^{-2}, \textbf{B} - \textbf{V}, \textbf{D} - \textbf{V}, \min \{ \textbf{B} - \textbf{V}, \textbf{C} - \textbf{V} \} \big)
	\end{flalign*}
	
	\noindent case of \eqref{sumuq} (where $\min \{ \textbf{B} - \textbf{V}, \textbf{C} - \textbf{V} \} \in \{ 0, 1 \}^n$ denotes the entrywise minimum of $\textbf{B} - \textbf{V}$ and $\textbf{C} - \textbf{V}$) yields 
	\begin{flalign*}  
		\displaystyle\lim_{s \rightarrow 0} (-s)^{-d} W_{x / s} (\textbf{A}, \textbf{B}; \textbf{C}, \textbf{D} \boldsymbol{\mid} r, s) & = x^d (-r^2)^{c - a - v} q^{\varphi (\textbf{D} - \textbf{V}, \textbf{C}) + \varphi (\textbf{V}, \textbf{A}) + \binom{b}{2} - dv + v} \displaystyle\frac{(r^2; q)_d}{(r^2; q)_b} \\
		& \qquad \times \displaystyle\prod_{j : \min \{ B_j - V_j, C_j - V_j \} = 1} \big(1 - r^{-2} q^{-B_{[j + 1, n]} - D_{[1, j - 1]}} \big),
	\end{flalign*} 
	
	\noindent where we have used the fact that $V_{[1, j - 1]} + V_{[j + 1, n]} = |\textbf{V}| - V_j = v$ holds for any $j \in [0, n]$ with $B_j - V_j = 1$ (as $B_j, V_j \in \{ 0, 1 \}$). Now the lemma follows from this, together with the fact that $\min \{ B_j - V_j, C_j - V_j \} = 1$ holds if and only if $(B_j, C_j, V_j) = (1, 1, 0)$, which in turn holds if and only if $B_j - D_j = 1$ (as $A_j + B_j = C_j + D_j$ and $A_j, B_j, C_j, D_j \in \{ 0, 1 \}$). 
\end{proof}

\begin{rem} 
	
	\label{wr} 
	
	It is quickly verified that the right side of \eqref{wxsr} is a polynomial in $r^{-2}$, as the potential denominator in the term $(r^2; q)_b^{-1} (r^2; q)_d$ there can be shown to be canceled by (part of) the product $\prod_{j: B_j - D_j = 1} \big( 1 - r^{-2} q^{-B_{[j + 1, n]} - D_{[1, j - 1]}} \big)$.

\end{rem}

\section{Limit Degenerations}

\label{Limitrsz}

In this section we consider the special cases of the $W_z$ weights given by \Cref{sxz1weights}, \Cref{z1weights}, and \Cref{zrw}, and consider their limits as some of the parameters $(r, s, z)$ tend to $0$ or $\infty$. Although there are in principle many possible choices for these limit degenerations, we only consider two in each case; the first (given by \Cref{weightsxz}, \Cref{s0z1weights}, and \Cref{szrzw} below) applies the limit as one parameter tends to $0$ or $\infty$, and the second (given by \Cref{weightz}, \Cref{weights1s0z}, and \Cref{sxs} below) applies it as two do. We will later see as \Cref{limitg0} below that, under the latter three double limits, the $F, G, H$ functions from \Cref{fgdefinition} degenerate to the Lascoux--Leclerc--Thibon (LLT) polynomials. Throughout this section we again adopt the notation of \Cref{wabcdrsxy}. 

We begin by degenerating \Cref{sxz1weights}, first when $s$ tends to $\infty$ and then when $z$ tends to $0$. 

\begin{cor} 
	
	\label{weightsxz} 
	
	For any fixed $z \in \mathbb{C}$, we have 
	\begin{flalign*}
	\displaystyle\lim_{s \rightarrow \infty} W_z ( \textbf{\emph{A}}, \textbf{\emph{B}}, \textbf{\emph{C}}, \textbf{\emph{D}} \boldsymbol{\mid} s z^{1/2}; s)  = (-1)^{b - c} q^{\varphi (\textbf{\emph{B}}, \textbf{\emph{C}} - \textbf{\emph{B}})} z^{-b} (z; q)_{b - c} \textbf{\emph{1}}_{\textbf{\emph{A}} + \textbf{\emph{B}} = \textbf{\emph{C}} + \textbf{\emph{D}}} \textbf{\emph{1}}_{\textbf{\emph{B}} \ge \textbf{\emph{C}}}.
	\end{flalign*} 
	 
\end{cor}

\begin{proof}
	
	Since for any fixed $z \in \mathbb{C}$ we have that
	\begin{flalign*}
	\displaystyle\lim_{s \rightarrow \infty} s^{2b - 2c} \displaystyle\frac{(s^2; q)_c}{(s^2 z; q)_b} = (-1)^{b - c} z^{-b} q^{\binom{c}{2} - \binom{b}{2}},
	\end{flalign*}
	
	\noindent and, by the bilinearity of $\varphi$, 
	\begin{flalign*}
	\varphi (\textbf{B} - \textbf{C}, \textbf{C}) + \binom{c}{2} - \binom{b}{2} = \varphi (\textbf{B} - \textbf{C}, \textbf{C}) + \varphi (\textbf{C}, \textbf{C}) - \varphi (\textbf{B}, \textbf{B}) = \varphi (\textbf{B}, \textbf{C} - \textbf{B}),
	\end{flalign*} 
	
	\noindent the corollary follows from \Cref{sxz1weights}. 
\end{proof}

\begin{cor}
	
	\label{weightz} 
	
	For any fixed $x \in \mathbb{C}$, we have  
	\begin{flalign}
	\label{wrslimit} 
	\displaystyle\lim_{y \rightarrow \infty} y^{-b} \bigg( \displaystyle\lim_{s \rightarrow \infty} W_{x / y} \big(\textbf{\emph{A}}, \textbf{\emph{B}}; \textbf{\emph{C}}, \textbf{\emph{D}} \boldsymbol{\mid} s x^{1 / 2}; s y^{1 / 2} \big) \bigg) = (-1)^{b - c} x^{-b} q^{\varphi (\textbf{\emph{B}}, \textbf{\emph{C}} - \textbf{\emph{B}})} \textbf{\emph{1}}_{\textbf{\emph{A}} + \textbf{\emph{B}} = \textbf{\emph{C}} + \textbf{\emph{D}}} \textbf{\emph{1}}_{\textbf{\emph{B}} \ge \textbf{\emph{C}}}.
	\end{flalign} 
	
\end{cor}

\begin{proof}
	
	This follows by setting $(s, z)$ to $\big( s y^{1 / 2}, \frac{x}{y} \big)$ in \Cref{weightsxz} and letting $y$ tend to $\infty$.
\end{proof}

Next we degenerate \Cref{z1weights}, first when $(r, s) = (sx^{1 / 2}, sy^{1 / 2})$ and $s$ tends to $0$, and then when $y$ tends to $0$. 

\begin{cor}
	
	\label{s0z1weights}
	
	For any fixed $r \in \mathbb{C}$, we have 
	\begin{flalign*}
	\displaystyle\lim_{s \rightarrow 0} W_1 \big( \textbf{\emph{A}}, \textbf{\emph{B}}; \textbf{\emph{C}}, \textbf{\emph{D}} \boldsymbol{\mid} s x^{1/2}, s y^{1/2} \big) = x^{-d} y^d (x^{-1} y; q)_{c - b} q^{\varphi (\textbf{\emph{D}}, \textbf{\emph{C}} - \textbf{\emph{B}})} \textbf{\emph{1}}_{\textbf{\emph{A}} + \textbf{\emph{B}} = \textbf{\emph{C}} + \textbf{\emph{D}}} \textbf{\emph{1}}_{\textbf{\emph{C}} \ge \textbf{\emph{B}}}. 
	\end{flalign*} 
\end{cor}

\begin{proof}
	This follows by setting $(r, s)$ to $(sx^{1 / 2}, sy^{1 / 2})$ in \Cref{z1weights} and letting $s$ tend to $0$.
\end{proof}

\begin{cor}
	
	\label{weights1s0z}

	For any fixed $z \in \mathbb{C}$, we have 
	\begin{flalign}
	\label{limitwz1} 
	\displaystyle\lim_{y \rightarrow 0} y^{-d} \bigg( \displaystyle\lim_{s \rightarrow 0} W_1 (\textbf{\emph{A}}, \textbf{\emph{B}}; \textbf{\emph{C}}, \textbf{\emph{D}} \boldsymbol{\mid} s x^{1/2}, s y^{1 / 2}) \bigg) = x^{-d} q^{\varphi (\textbf{\emph{D}}, \textbf{\emph{C}} - \textbf{\emph{B}})} \textbf{\emph{1}}_{\textbf{\emph{A}} + \textbf{\emph{B}} = \textbf{\emph{C}} + \textbf{\emph{D}}} \textbf{\emph{1}}_{\textbf{\emph{C}} \ge \textbf{\emph{B}}}. 
	\end{flalign}
	
\end{cor}

\begin{proof}
	This follows by letting $y$ tend to $0$ in \Cref{s0z1weights}.
\end{proof}

Now we degenerate \Cref{zrw}, first when $s$ tends to $0$ and then when $y$ tends to $0$.

\begin{cor}
	
	\label{szrzw} 
	
	For any fixed $z \in \mathbb{C}$, we have  
	\begin{flalign*}
	\displaystyle\lim_{s \rightarrow 0} s^{-2d} W_z (\textbf{\emph{A}}, \textbf{\emph{B}}; \textbf{\emph{C}}, \textbf{\emph{D}} \boldsymbol{\mid} z^{1 / 2}; s) = q^{\varphi (\textbf{\emph{D}}, \textbf{\emph{C}})} (z; q)_d \textbf{\emph{1}}_{\textbf{\emph{A}} + \textbf{\emph{B}} = \textbf{\emph{C}} + \textbf{\emph{D}}} \textbf{\emph{1}}_{v = 0}.
	\end{flalign*}
	
\end{cor}
	
\begin{proof}
	
	This follows by letting $s$ tend to $0$ in \Cref{zrw}.
\end{proof}

\begin{cor}
	
	\label{sxs}

	For any fixed $x \in \mathbb{C}$, we have 
\begin{flalign}
\label{sxslimit} 
\displaystyle\lim_{y \rightarrow 0} (-y)^d \bigg( \displaystyle\lim_{s \rightarrow 0} s^{-2d} W_{x / y} \big( \textbf{\emph{A}}, \textbf{\emph{B}}; \textbf{\emph{C}}, \textbf{\emph{D}} \boldsymbol{\mid} (xy^{-1})^{1/2}, s \big) \bigg) = x^d q^{\varphi (\textbf{\emph{D}}, \textbf{\emph{C}} + \textbf{\emph{D}})} \textbf{\emph{1}}_{\textbf{\emph{A}} + \textbf{\emph{B}} = \textbf{\emph{C}} + \textbf{\emph{D}}} \textbf{\emph{1}}_{v = 0}.
\end{flalign}

\end{cor} 

\begin{proof}
	
	This follows by setting $z = \frac{x}{y}$ in \Cref{szrzw}, letting $y$ tend to $0$, the fact that 
	\begin{flalign*}
	\displaystyle\lim_{y \rightarrow 0} (-y)^d (xy^{-1}; q)_d = x^d q^{\binom{d}{2}} = x^d q^{\varphi (\textbf{D}, \textbf{D})},
	\end{flalign*}
	
	\noindent and the bilinearity of $\varphi$.
\end{proof}

\begin{rem}
	
	\label{limitsy}
	
	It is quickly verified that one obtains the same results as in \Cref{weightz}, \Cref{weights1s0z}, and \Cref{sxs} by taking the double limits of \eqref{wabcdp} in $y$ and $s$ simultaneously, assuming $s^2 y^{-1}$ tends to $\infty$ in \Cref{weightz} or that it tends to $0$ in \Cref{weights1s0z} and \Cref{sxs}.
	
\end{rem}

Before proceeding, it will be useful to provide one further degeneration of the fused weights $W_z (\textbf{A}, \textbf{B}; \textbf{C}, \textbf{D} \boldsymbol{\mid} r, s)$ from \Cref{wabcdrsxy}. It corresponds to the limit of the weight from from \Cref{zrw} as $s^2$ tends to $0$ and $z$ tends $\infty$ simultaneously, so that $s^2 z$ remains fixed; this weight will be useful in \Cref{Function3LG} and \Cref{HFunction0} below. 

\begin{cor}
	
	\label{szrzw2}
	
	For any fixed $z \in \mathbb{C}$, we have  
	\begin{flalign*}
	\displaystyle\lim_{s \rightarrow 0} W_{z / s^2} (\textbf{\emph{A}}, \textbf{\emph{B}}; \textbf{\emph{C}}, \textbf{\emph{D}} \boldsymbol{\mid} s^{-1} z^{1 / 2}; s) =  q^{\varphi (\textbf{\emph{D}}, \textbf{\emph{C}} + \textbf{\emph{D}})} \displaystyle\frac{(-z)^d}{(z; q)_{c + d}} \textbf{\emph{1}}_{\textbf{\emph{A}} + \textbf{\emph{B}} = \textbf{\emph{C}} + \textbf{\emph{D}}} \textbf{\emph{1}}_{v = 0}.
	\end{flalign*}
	
\end{cor}

\begin{proof}
	
	By \Cref{zrw}, we have that 
	\begin{flalign*}
	W_{z / s^2} ( \textbf{A}, \textbf{B}; \textbf{C}, \textbf{D} \boldsymbol{\mid} s^{-1} z^{1 / 2}; s) = s^{2d} q^{\varphi (\textbf{D}, \textbf{C})} \displaystyle\frac{(s^2; q)_c (s^{-2} z; q)_d}{(z; q)_{c + d}} \textbf{1}_{\textbf{A} + \textbf{B} = \textbf{C} + \textbf{D}} \textbf{1}_{v = 0}.
	\end{flalign*}
	
	\noindent Letting $s$ tend to $0$, using the fact that 
	\begin{flalign*}
	\displaystyle\lim_{s \rightarrow 0} s^{2d} (s^{-2} z; q)_d = (-z)^d q^{\binom{d}{2}} = (-z)^d q^{\varphi (\textbf{D}, \textbf{D})},
	\end{flalign*}
	
	\noindent and applying the bilinearity of $\varphi$, we deduce the corollary. 	
\end{proof}

\section{Limiting Weights and Degenerated Functions}

\label{WeightsLimits}

In this section we introduce notation for the weights and functions degenerated according to the limits considered in \Cref{Limitrsz}. In view of \Cref{weightsxz}, \Cref{weightz}, \Cref{weights1s0z}, \Cref{sxs}, \Cref{szrzw2}, and \Cref{ws0}, we define the following limiting weights.

\begin{definition}
	
	\label{wlimits} 
	
	If $\textbf{A} + \textbf{B} = \textbf{C} + \textbf{D}$, then set
	\begin{flalign}
	\label{limitw}  
	\begin{aligned}
	 \mathcal{W}_{x; y} (\textbf{A}, \textbf{B}; \textbf{C}, \textbf{D} \boldsymbol{\mid} \infty, \infty) & = (-1)^{b - c} x^{-b} y^b (xy^{-1}; q)_{b - c} q^{\varphi (\textbf{B}, \textbf{C} - \textbf{B})}  \textbf{1}_{\textbf{B} \ge \textbf{C}}; \\
	 \mathcal{W}_x (\textbf{A}, \textbf{B}; \textbf{C}, \textbf{D} \boldsymbol{\mid} \infty, \infty) & = (-1)^{b - c} x^{-b} q^{\varphi (\textbf{B}, \textbf{C} - \textbf{B})}  \textbf{1}_{\textbf{B} \ge \textbf{C}}; \\
	 \mathcal{W}_x (\textbf{A}, \textbf{B}; \textbf{C}, \textbf{D} \boldsymbol{\mid} 0, 0) & = x^{-d} q^{\varphi (\textbf{D}, \textbf{C} - \textbf{B})}  \textbf{1}_{\textbf{C} \ge \textbf{B}}; \\
	\mathcal{W}_x (\textbf{A}, \textbf{B}; \textbf{C}, \textbf{D} \boldsymbol{\mid} \infty, 0) & = x^d q^{\varphi (\textbf{D}, \textbf{C} + \textbf{D})}  \textbf{1}_{v = 0}; \\
	\mathcal{W}_{x; y} (\textbf{A}, \textbf{B}; \textbf{C}, \textbf{D} \boldsymbol{\mid} \infty, 0) & = q^{\varphi(\textbf{D}, \textbf{C} + \textbf{D})} \displaystyle\frac{(-x)^d y^{-d}}{(x y^{-1}; q)_{c + d}} \textbf{1}_{\textbf{A} + \textbf{B} = \textbf{C} + \textbf{D}} \textbf{1}_{v = 0}; \\
	\mathcal{W}_x (\textbf{A}, \textbf{B}; \textbf{C}, \textbf{D} \boldsymbol{\mid} r) & = x^d (-r^2)^{c - a - v} q^{\varphi (\textbf{D} - \textbf{V}, \textbf{C}) + \varphi (\textbf{V}, \textbf{A}) + \binom{b}{2} - dv + v} \displaystyle\frac{(r^2; q)_d}{(r^2; q)_b} \\
	& \qquad \times \textbf{1}_{\textbf{A} + \textbf{B} = \textbf{C} + \textbf{D}} \displaystyle\prod_{j : B_j - D_j = 1} \big(1 - r^{-2} q^{-B_{[j + 1, n]} - D_{[1, j - 1]}} \big),
	\end{aligned} 
	\end{flalign}
\index{W@$W_z (\textbf{A}, \textbf{B}; \textbf{C}, \textbf{D} \boldsymbol{\mid} r, s)$; fused weight!$\mathcal{W}_{x; y} (\textbf{A}, \textbf{B}; \textbf{C}, \textbf{D} \boldsymbol{\mid} \infty, \infty)$}
\index{W@$W_z (\textbf{A}, \textbf{B}; \textbf{C}, \textbf{D} \boldsymbol{\mid} r, s)$; fused weight!$\mathcal{W}_x (\textbf{A}, \textbf{B}; \textbf{C}, \textbf{D} \boldsymbol{\mid} \infty, \infty)$}
\index{W@$W_z (\textbf{A}, \textbf{B}; \textbf{C}, \textbf{D} \boldsymbol{\mid} r, s)$; fused weight!$\mathcal{W}_x (\textbf{A}, \textbf{B}; \textbf{C}, \textbf{D} \boldsymbol{\mid} 0, 0)$}
\index{W@$W_z (\textbf{A}, \textbf{B}; \textbf{C}, \textbf{D} \boldsymbol{\mid} r, s)$; fused weight!$\mathcal{W}_x (\textbf{A}, \textbf{B}; \textbf{C}, \textbf{D} \boldsymbol{\mid} \infty, 0)$}
\index{W@$W_z (\textbf{A}, \textbf{B}; \textbf{C}, \textbf{D} \boldsymbol{\mid} r, s)$; fused weight!$\mathcal{W}_{x; y} (\textbf{A}, \textbf{B}; \textbf{C}, \textbf{D} \boldsymbol{\mid} \infty, 0)$}
\index{W@$W_z (\textbf{A}, \textbf{B}; \textbf{C}, \textbf{D} \boldsymbol{\mid} r, s)$; fused weight!$\mathcal{W}_x (\textbf{A}, \textbf{B}; \textbf{C}, \textbf{D} \boldsymbol{\mid} r)$}
	
	\noindent and
	\begin{flalign}
	\label{limitw2} 
	\begin{aligned}
	\widehat{\mathcal{W}}_{x; y} (\textbf{A}, \textbf{B}; \textbf{C}, \textbf{D} \boldsymbol{\mid} \infty, \infty) & = (-1)^{b - c - n} x^{n - b} y^{b - n} \displaystyle\frac{(x y^{-1}; q)_{b - c}}{(x y^{-1}; q)_n} q^{\varphi (\textbf{B}, \textbf{C} - \textbf{B}) + \binom{n}{2}} \textbf{1}_{\textbf{B} \ge \textbf{C}}; \\
	\widehat{\mathcal{W}}_x (\textbf{A}, \textbf{B}; \textbf{C}, \textbf{D} \boldsymbol{\mid} \infty, \infty) & = (-1)^{b - c - n} x^{n - b} q^{\varphi (\textbf{B}, \textbf{C} - \textbf{B}) + \binom{n}{2}}  \textbf{1}_{\textbf{B} \ge \textbf{C}}; \\
	\widehat{\mathcal{W}}_x (\textbf{A}, \textbf{B}; \textbf{C}, \textbf{D} \boldsymbol{\mid} \infty, 0) & = x^{d - n} q^{\varphi (\textbf{D}, \textbf{C} + \textbf{D}) - \binom{n}{2}}  \textbf{1}_{v = 0}, \\
	\widehat{\mathcal{W}}_{x; y} (\textbf{A}, \textbf{B}; \textbf{C}, \textbf{D} \boldsymbol{\mid} \infty, 0) & = q^{\varphi(\textbf{D}, \textbf{C} + \textbf{D}) - \binom{n}{2}} \displaystyle\frac{(-x)^{d - n} y^{n - d} (x y^{-1}; q)_n}{(x y^{-1}; q)_{c + d}} \textbf{1}_{\textbf{A} + \textbf{B} = \textbf{C} + \textbf{D}} \textbf{1}_{v = 0},
	\end{aligned} 
	\end{flalign}
\index{W@$\widehat{W}_z (\textbf{A}, \textbf{B}; \textbf{C}, \textbf{D} \boldsymbol{\mid} r, s)$; normalized fused weight!$\widehat{\mathcal{W}}_{x; y} (\textbf{A}, \textbf{B}; \textbf{C}, \textbf{D} \boldsymbol{\mid} \infty, \infty)$}
\index{W@$\widehat{W}_z (\textbf{A}, \textbf{B}; \textbf{C}, \textbf{D} \boldsymbol{\mid} r, s)$; normalized fused weight!$\widehat{\mathcal{W}}_x (\textbf{A}, \textbf{B}; \textbf{C}, \textbf{D} \boldsymbol{\mid} \infty, \infty)$}
\index{W@$\widehat{W}_z (\textbf{A}, \textbf{B}; \textbf{C}, \textbf{D} \boldsymbol{\mid} r, s)$; normalized fused weight!$\widehat{\mathcal{W}}_x (\textbf{A}, \textbf{B}; \textbf{C}, \textbf{D} \boldsymbol{\mid} \infty, 0)$}
\index{W@$\widehat{W}_z (\textbf{A}, \textbf{B}; \textbf{C}, \textbf{D} \boldsymbol{\mid} r, s)$; normalized fused weight!$\widehat{\mathcal{W}}_{x; y} (\textbf{A}, \textbf{B}; \textbf{C}, \textbf{D} \boldsymbol{\mid} \infty, 0)$}
	
	\noindent Otherwise, set all of these quantities to $0$.
		
\end{definition}

Observe that the definitions \eqref{limitw} and \eqref{limitw2} are consistent with \eqref{wabcd2}. We did not introduce $\widehat{\mathcal{W}}_x (\textbf{A}, \textbf{B}; \textbf{C}, \textbf{D} \boldsymbol{\mid} 0, 0)$, since it is ill-defined as $\mathcal{W}_x \big( \textbf{e}_0, \textbf{e}_{[1, n]}; \textbf{e}_0, \textbf{e}_{[1, n]} \boldsymbol{\mid} 0, 0 \big) = 0$ (we also did not introduce $\widehat{\mathcal{W}}_x (\textbf{A}, \textbf{B}; \textbf{C}, \textbf{D} \boldsymbol{\mid} r)$ since we will not use it). 

The next definition provides the corresponding limit degenerations for the functions $G_{\boldsymbol{\lambda} / \boldsymbol{\mu}}$, $F_{\boldsymbol{\lambda} / \boldsymbol{\mu}}$, and $H_{\boldsymbol{\lambda} / \boldsymbol{\mu}}$ from \Cref{fgdefinition}. There are various additional degenerations (using all possible weights from \Cref{wlimits}) that could be defined and considered, but here we only introduce the ones we will use in this text. In what follows, for any complex number $t \in \mathbb{C}$ and (possibly infinite) sequence $\textbf{z} = (z_1, z_2, \ldots , z_k) \subset \mathbb{C}$, we define the sequence $t \textbf{z} = (tz_1, tz_2, \ldots , tz_k) \subset \mathbb{C}$.

\begin{definition}
	
	\label{sgf0} 
	
	Fix integers $M \ge 0$ and $N \ge 1$; a complex number $r \in \mathbb{C}$; a finite sequence $\textbf{x} = (x_1, x_2, \ldots , x_N)$ of complex numbers; and an infinite sequence $\textbf{y} = (y_1, y_2, \ldots )$ of complex numbers. For any sequences of $n$ signatures $\boldsymbol{\lambda}, \boldsymbol{\mu} \in \SeqSign_{n; M}$, set
	\begin{flalign}
	\label{limitg}
	\begin{aligned}
	\mathcal{G}_{\boldsymbol{\lambda} / \boldsymbol{\mu}} (\textbf{x} ; \infty \boldsymbol{\mid} \textbf{y}; \infty) & = \displaystyle\lim_{s \rightarrow \infty} G_{\boldsymbol{\lambda} / \boldsymbol{\mu}} (\textbf{x}; s \textbf{x}^{1/2} \boldsymbol{\mid} \textbf{y}; s \textbf{y}^{1/2}); \\
	\mathcal{G}_{\boldsymbol{\lambda} / \boldsymbol{\mu}} (0; \textbf{x} \boldsymbol{\mid} 0; 0) & = \displaystyle\lim_{y \rightarrow 0} y^{|\boldsymbol{\mu}| - |\boldsymbol{\lambda}|} \bigg( \displaystyle\lim_{s \rightarrow 0} G_{\boldsymbol{\lambda} / \boldsymbol{\mu}} \big( (y, y, \ldots ); s \textbf{x}^{1 / 2} \boldsymbol{\mid} (y, y, \ldots ); s (y^{1 / 2}, y^{1 / 2}, \ldots ) \big) \bigg); \\
	\mathcal{G}_{\boldsymbol{\lambda} / \boldsymbol{\mu}} (\textbf{x} ; \infty \boldsymbol{\mid} 0; 0) & = \displaystyle\lim_{y \rightarrow 0}  (-y)^{|\boldsymbol{\lambda}| - |\boldsymbol{\mu}|} \bigg( \displaystyle\lim_{s \rightarrow 0} s^{2 |\boldsymbol{\mu}| - 2 |\boldsymbol{\lambda}|} G_{\boldsymbol{\lambda} / \boldsymbol{\mu}} \big( \textbf{x}; y^{-1/2} \textbf{x}^{1/2} \boldsymbol{\mid} (y, y, \ldots ); (s, s, \ldots ) \big) \bigg); \\
	\mathcal{G}_{\boldsymbol{\lambda} / \boldsymbol{\mu}} (\textbf{x}; r \boldsymbol{\mid} 0; 0) & = \displaystyle\lim_{s \rightarrow 0} (-s)^{|\boldsymbol{\mu}| - |\boldsymbol{\lambda}|} G_{\boldsymbol{\lambda} / \boldsymbol{\mu}} \big( \textbf{x}; (r, r, \ldots ) \boldsymbol{\mid} (s, s, \ldots ); (s, s, \ldots ) \big).
	\end{aligned}
	\end{flalign}
	\index{G@$G_{\boldsymbol{\lambda} / \boldsymbol{\mu}} (\textbf{x}; \textbf{r} \boldsymbol{\mid} \textbf{y}; \textbf{s})$!$\mathcal{G}_{\boldsymbol{\lambda} / \boldsymbol{\mu}} (\textbf{x}; \infty \boldsymbol{\mid} \textbf{y}; \infty)$}\index{G@$G_{\boldsymbol{\lambda} / \boldsymbol{\mu}} (\textbf{x}; \textbf{r} \boldsymbol{\mid} \textbf{y}; \textbf{s})$!$\mathcal{G}_{\boldsymbol{\lambda} / \boldsymbol{\mu}} (0; \textbf{x} \boldsymbol{\mid} 0; 0)$}\index{G@$G_{\boldsymbol{\lambda} / \boldsymbol{\mu}} (\textbf{x}; \textbf{r} \boldsymbol{\mid} \textbf{y}; \textbf{s})$!$\mathcal{G}_{\boldsymbol{\lambda} / \boldsymbol{\mu}} (\textbf{x}; \infty \boldsymbol{\mid} 0; 0)$}\index{G@$G_{\boldsymbol{\lambda} / \boldsymbol{\mu}} (\textbf{x}; \textbf{r} \boldsymbol{\mid} \textbf{y}; \textbf{s})$!$\mathcal{G}_{\boldsymbol{\lambda} / \boldsymbol{\mu}} (\textbf{x}; r \boldsymbol{\mid} 0; 0)$}
	
	\noindent Moreover, for any sequences of $n$ signatures $\boldsymbol{\lambda} \in \SeqSign_{n; M + N}$ and $\boldsymbol{\mu} \in \SeqSign_{n; M}$, set
	\begin{flalign}
	\label{limithf}
	\begin{aligned} 
	\mathcal{H}_{\boldsymbol{\lambda} / \boldsymbol{\mu}} (\textbf{x} ; \infty \boldsymbol{\mid} \textbf{y}; \infty) & = \displaystyle\lim_{s \rightarrow \infty} H_{\boldsymbol{\lambda} / \boldsymbol{\mu}} (\textbf{x}; s \textbf{x}^{1/2} \boldsymbol{\mid} \textbf{y}; s \textbf{y}^{1/2}); \\
	\mathcal{H}_{\boldsymbol{\lambda} / \boldsymbol{\mu}} (\textbf{x} ; \infty \boldsymbol{\mid} \infty; \infty) & = \displaystyle\lim_{y \rightarrow \infty} y^{|\boldsymbol{\mu}| - |\boldsymbol{\lambda}| + n \binom{M + 1}{2} - n \binom{M + N + 1}{2}} \\
	& \qquad \times \bigg( \displaystyle\lim_{s \rightarrow \infty} H_{\boldsymbol{\lambda} / \boldsymbol{\mu}} \big( \textbf{x}; s^{1/2} \textbf{x} \boldsymbol{\mid} (y, y, \ldots ); s^{1 / 2} (y, y, \ldots ) \big) \bigg); \\
	\mathcal{F}_{\boldsymbol{\lambda} / \boldsymbol{\mu}} (\textbf{x} ; \infty \boldsymbol{\mid} 0; 0) & = \displaystyle\lim_{y \rightarrow 0} (-y)^{|\boldsymbol{\mu}| - |\boldsymbol{\lambda}| + n \binom{M}{2} - n \binom{M + N}{2}} \\
	& \quad \times  \bigg( \displaystyle\lim_{s \rightarrow 0} s^{2(|\boldsymbol{\lambda}| - |\boldsymbol{\mu}| - n \binom{M}{2} + n \binom{M + N}{2})} F_{\boldsymbol{\lambda} / \boldsymbol{\mu}} \big( \textbf{x}; y^{-1 / 2} \textbf{x}^{1/2} \boldsymbol{\mid} (y, y, \ldots ); (s, s, \ldots ) \big) \bigg); \\
	\mathcal{F}_{\boldsymbol{\lambda} / \boldsymbol{\mu}} (\textbf{x}; \infty \boldsymbol{\mid} \textbf{y}; 0) & = \displaystyle\lim_{s \rightarrow 0} F_{\boldsymbol{\lambda} / \boldsymbol{\mu}} \Big( \textbf{x}; \big( s^{-1} (x_1 y_1^{-1})^{1 / 2}, s^{-1} (x_2 y_2^{-1})^{1 / 2}, \ldots \big) \boldsymbol{\mid} s^2 \textbf{y}; (s, s, \ldots) \Big).
	\end{aligned} 
	\end{flalign}
	\index{F@$F_{\boldsymbol{\lambda} / \boldsymbol{\mu}} (\textbf{x}; \textbf{r} \boldsymbol{\mid} \textbf{y}; \textbf{s})$!$\mathcal{F}_{\boldsymbol{\lambda} / \boldsymbol{\mu}} (\textbf{x}; \infty \boldsymbol{\mid} \textbf{y}; 0)$} \index{F@$F_{\boldsymbol{\lambda} / \boldsymbol{\mu}} (\textbf{x}; \textbf{r} \boldsymbol{\mid} \textbf{y}; \textbf{s})$!$\mathcal{F}_{\boldsymbol{\lambda} / \boldsymbol{\mu}} (\textbf{x}; \infty \boldsymbol{\mid} 0; 0)$} \index{H@$H_{\boldsymbol{\lambda} / \boldsymbol{\mu}} (\textbf{x}; \textbf{r} \boldsymbol{\mid} \textbf{y}; \textbf{s})$!$\mathcal{H}_{\boldsymbol{\lambda} / \boldsymbol{\mu}} (\textbf{x}; \infty \boldsymbol{\mid} \infty; \infty)$} \index{H@$H_{\boldsymbol{\lambda} / \boldsymbol{\mu}} (\textbf{x}; \textbf{r} \boldsymbol{\mid} \textbf{y}; \textbf{s})$!$\mathcal{H}_{\boldsymbol{\lambda} / \boldsymbol{\mu}} (\textbf{x}; \infty \boldsymbol{\mid} \textbf{y}; \infty)$} 
	
\end{definition}

\begin{rem}
	
	\label{sumbw} 
	
	Let us briefly account for the powers of $y$ (and $s$) appearing in the limits on the right sides of \eqref{limitg} and \eqref{limithf}, taking $\mathcal{H}_{\boldsymbol{\lambda} / \boldsymbol{\mu}} (\textbf{x}; \infty \boldsymbol{\mid} \infty; \infty)$ as an example. As indicated in \eqref{gfhe}, $H_{\boldsymbol{\lambda} / \boldsymbol{\mu}} (\textbf{x}; \textbf{r} \boldsymbol{\mid} \textbf{y}; \textbf{s})$ is the partition function under the weights $W_{x / y} (\textbf{A}, \textbf{B}; \textbf{C}, \textbf{D} \boldsymbol{\mid} r, s)$ for the vertex model $\mathfrak{P}_H (\boldsymbol{\lambda} / \boldsymbol{\mu})$ from \Cref{pgpfph}. By \Cref{weightz}, the weights $\mathcal{W}_x (\textbf{A}, \textbf{B}; \textbf{C}, \textbf{D} \boldsymbol{\mid} \infty, \infty)$ are obtained as the limit of $y^{-b} W_{x / y} (\textbf{A}, \textbf{B}; \textbf{C}, \textbf{D} \boldsymbol{\mid} sx^{1 / 2}, sy^{1 / 2})$ as $s$ and $y$ tend to $\infty$; in particular, a factor of $y^{-b}$ is required as a normalization in these weights. Thus, in view of \eqref{gfhe} (and recalling the notation around there), a normalization of $y^{-\xi (\boldsymbol{\lambda} / \boldsymbol{\mu})}$ is required in the definition of $\mathcal{H}_{\boldsymbol{\lambda} / \boldsymbol{\mu}} (\textbf{x}; \infty \boldsymbol{\mid} \infty; \infty)$, where $\xi (\boldsymbol{\lambda} / \boldsymbol{\mu})$ is the common value over all $\mathcal{E} \in \mathfrak{P}_H (\boldsymbol{\lambda} / \boldsymbol{\mu})$ of
	\begin{flalign*}
	\xi (\boldsymbol{\lambda} / \boldsymbol{\mu}) = \displaystyle\sum_{(i, j) \in \mathcal{D}} \big| \textbf{B} (i, j) \big| & = \displaystyle\sum_{c = 1}^n \Bigg( \displaystyle\sum_{\mathfrak{l} \in \mathfrak{T} (\lambda^{(c)})} \mathfrak{l} - \displaystyle\sum_{\mathfrak{m} \in \mathfrak{T} (\mu^{(c)})} \mathfrak{m} \Bigg) \\
	& = \displaystyle\sum_{c = 1}^n \Bigg( \displaystyle\sum_{i = 1}^{M + N} \big( \lambda_i^{(c)} + M + N - i + 1 \big) - \displaystyle\sum_{i = 1}^M \big( \mu_i^{(c)} + M - i + 1 \big) \Bigg) \\
	& = |\boldsymbol{\lambda}| + \binom{M + N + 1}{2} n - |\boldsymbol{\mu}| - \binom{M + 1}{2} n,
	\end{flalign*} 
	
	\noindent where $\big( \textbf{A} (v), \textbf{B} (v); \textbf{C} (v), \textbf{D} (v) \big)$ denotes the arrow configuration at any vertex $v \in \mathcal{D}$ under $\mathcal{E}$. This indeed recovers the normalizing factor of $y$ on the right side of the second equation of \eqref{limithf}; similar statements hold for the remaining limits in \Cref{sgf0}.

\end{rem}

\begin{rem} 
	
	\label{weightesum}
	
Observe that the degenerated functions from \Cref{sgf0} are partition functions for vertex models under the limiting weights from \Cref{wlimits}. Indeed, let us extend the definitions of the weights $W (\mathcal{E} \boldsymbol{\mid} \textbf{x}; \textbf{r} \boldsymbol{\mid} \textbf{y}; \textbf{s})$ and $\widehat{W} (\mathcal{E} \boldsymbol{\mid} \textbf{x}; \textbf{r} \boldsymbol{\mid} \textbf{y}; \textbf{s})$ from \eqref{weighte} to when $\textbf{x}, \textbf{y}, \textbf{r}, \textbf{s} \in \{ 0, \infty \}$, so that for example 
\begin{flalign*}
W (\mathcal{E} \boldsymbol{\mid} \textbf{x}; \infty \boldsymbol{\mid} \infty; \infty) = \displaystyle\prod_{(i, j) \in \mathcal{D}} \mathcal{W}_{x_i} \big( \textbf{A} (i, j), \textbf{B} (i, j); \textbf{C} (i, j), \textbf{D} (i, j) \boldsymbol{\mid} \infty, \infty \big),
\end{flalign*}

\noindent for any path ensemble $\mathcal{E}$ on $\mathcal{D}$ (recall \eqref{dn}), where $\big( \textbf{A} (v), \textbf{B} (v); \textbf{C} (v), \textbf{D} (v) \big)$ denotes the arrow configuration at any vertex $v \in \mathcal{D}$ under $\mathcal{E}$. Then, \eqref{gfhe} continues to hold for $\textbf{x}, \textbf{r}, \textbf{s}, \textbf{y} \in \{ 0, \infty \}$. 

\end{rem}

\section{Degenerations of Cauchy Identities} 

\label{IdentitiesDegenerations}

In this section we analyze the special cases of the Cauchy identity \Cref{fgsum2} under the degenerations considered in \Cref{DegenerationWeights} and \Cref{Limitrsz}. Throughout this section, we fix integers $N, M \ge 1$; finite sequences of complex numbers $\textbf{u} = (u_1, u_2, \ldots , u_N)$, $\textbf{r} = (r_1, r_2, \ldots , r_N)$, $\textbf{w} = (w_1, w_2, \ldots , w_M)$, and $\textbf{t} = (t_1, t_2, \ldots , t_M)$; and infinite sequences of complex numbers $\textbf{y} = (y_1, y_2, \ldots )$ and $\textbf{s} = (s_1, s_2, \ldots )$. 

There are in principle numerous such degenerations that can be analyzed, but let us only restrict our attention to two, corresponding to the sums over $\boldsymbol{\lambda}$ of $\mathcal{F}_{\boldsymbol{\lambda}} (\textbf{u}; \infty \boldsymbol{\mid} 0; 0) \mathcal{G}_{\boldsymbol{\lambda}} (\textbf{w}; \infty \boldsymbol{\mid} 0; 0)$ and $\mathcal{F}_{\boldsymbol{\lambda}} (\textbf{u}; \infty \boldsymbol{\mid} 0; 0) \mathcal{G}_{\boldsymbol{\lambda}} (0; \textbf{w} \boldsymbol{\mid} 0; 0)$. We will later see that these correspond to the Cauchy identity and the dual Cauchy identity for the LLT polynomials, respectively. 

\begin{cor}
	
	\label{sxzlimitidentity}
	
	If 
	\begin{flalign}
	\label{x1x2bdestimate6} 
	\displaystyle\max_{\substack{1 \le i \le M \\ 1 \le j \le N}} \displaystyle\max_{\substack{a, b \in [0, n] \\ (a, b) \ne (n, 0)}} \big| q^{\binom{a}{2} + \binom{b}{2} - \binom{n}{2}} u_j^{a - n} w_i^b \big| < 1,
	\end{flalign}
	
	\noindent then
	\begin{flalign}
	\label{fgsums1}
	\begin{aligned}
	\displaystyle\sum_{\boldsymbol{\lambda} \in \SeqSign_{n; N}} \mathcal{F}_{\boldsymbol{\lambda}} (\textbf{\emph{u}}; \infty \boldsymbol{\mid} 0; 0) \mathcal{G}_{\boldsymbol{\lambda}} (\textbf{\emph{w}}; \infty \boldsymbol{\mid} 0; 0) & = q^{-\binom{n}{2} \binom{N}{2}} \displaystyle\prod_{i = 1}^M \displaystyle\prod_{j = 1}^N (u_j^{-1} w_i; q^{-1})_n^{-1}  \displaystyle\prod_{j = 1}^N u_j^{n (j - N)}.
	\end{aligned} 
	\end{flalign}
	
\end{cor}

\begin{proof}
	
	Fix $s, y \in \mathbb{C}$, and denote the infinite sequences $\textbf{s} = (s, s, \ldots )$ and $\textbf{y} = (y, y, \ldots )$. Applying \Cref{fgsum2}, we obtain
	\begin{flalign}
	\label{f0g0sum1}
	\begin{aligned}
	& \displaystyle\sum_{\boldsymbol{\lambda} \in \SeqSign_{n; N}} (- y^{-1} s^2)^{|\boldsymbol{\lambda}| + n \binom{N}{2}}  F_{\boldsymbol{\lambda}} \big( \textbf{u}; y^{-1 / 2} \textbf{u}^{1 / 2} \boldsymbol{\mid} \textbf{y}; \textbf{s} \big) (- y^{-1} s^2)^{- |\boldsymbol{\lambda}|} G_{\boldsymbol{\lambda}} \big( \textbf{w}; y^{-1 / 2} \textbf{w}^{1 / 2} \boldsymbol{\mid} \textbf{y}; \textbf{s} \big) \\
	& \qquad = (- y^{-1} s^2)^{n \binom{N}{2}} \displaystyle\prod_{k = 1}^N s^{2n (k -N)} (s^2; q)_n^{N - k} \displaystyle\prod_{j = 1}^N (u_j y^{-1}; q)_n^{j - N} \displaystyle\prod_{i = 1}^M \displaystyle\prod_{j = 1}^N \displaystyle\frac{y^n (u_j y^{-1}; q)_n}{w_i^n (u_j w_i^{-1}; q)_n} \\
	& \qquad = \displaystyle\prod_{k = 1}^N (s^2; q)_n^{N - k} \displaystyle\prod_{j = 1}^N (-y)^{n (j - N)} (u_j y^{-1}; q)_n^{j - N} \displaystyle\prod_{i = 1}^M \displaystyle\prod_{j = 1}^N \displaystyle\frac{y^n (u_j y^{-1}; q)_n}{w_i^n (u_j w_i^{-1}; q)_n},
	\end{aligned} 
	\end{flalign}
	
	\noindent assuming there exists an integer $K > 1$ and a real number $\varepsilon > 0$ such that
	\begin{flalign}
	\label{x1x2bdestimate5}
	\begin{aligned}
	\displaystyle\sup_{k > K} & \displaystyle\max_{\substack{1 \le i \le M \\ 1 \le j \le N}} \displaystyle\max_{\substack{a, b \in [0, n] \\ (a, b) \ne (n, 0)}} |s|^{2a + 2b - 2n} \displaystyle\max_{|\textbf{B}| = a} \big| (-y^{-1} s^2)^{n - a} \widehat{W}_{u_j / y} (\textbf{e}_0, \textbf{B}; \textbf{e}_0, \textbf{B} \boldsymbol{\mid} y^{-1 / 2} u_j^{1 / 2}; s) \big|\\
	& \qquad \qquad \qquad \qquad \qquad \quad \times \displaystyle\max_{|\textbf{B}| = b} \big| (-y^{-1} s^2)^{-b} W_{w_i / y} (\textbf{e}_0, \textbf{B}; \textbf{e}_0, \textbf{B} \boldsymbol{\mid} y^{-1 / 2} w_i^{1 / 2}; s) \big|  \\
	& \quad = \displaystyle\sup_{k > K} \displaystyle\max_{\substack{1 \le i \le M \\ 1 \le j \le N}}	\displaystyle\max_{\substack{a, b \in [0, n] \\ (a, b) \ne (n, 0)}} \Bigg| y^{a + b - n} \displaystyle\frac{(s^2 u_j y^{-1}; q)_n (u_j y^{-1}; q)_a}{(u_j y^{-1}; q)_n (s^2 u_j y^{-1}; q)_a} \displaystyle\frac{(w_i y^{-1}; q)_b}{(s^2 w_i y^{-1}; q)_b} \Bigg| < 1 - \varepsilon,
	\end{aligned}
	\end{flalign}
	
	\noindent where the additional factors of $(- y^{-1} s^2)^{n - a}$ and $(- y^{-1} s^2)^{-b}$ arise due to the correpsonding gauge factors on the left side of \eqref{f0g0sum1}, and we have used \eqref{wze0b} and \eqref{wabcd2} in the first equality in \eqref{x1x2bdestimate5}. In particular,  since $\lim_{y \rightarrow 0} (-y)^k (xy^{-1}; q)_k = x^k q^{\binom{k}{2}}$ for any $x \in \mathbb{C}$ and $k \in \mathbb{Z}_{\ge 0}$, \eqref{x1x2bdestimate6} implies the existence of a constant $\delta = \delta (\varepsilon) > 0$ such that \eqref{x1x2bdestimate5} holds whenever $|s|, |y|, |s^2 y^{-1}| < \delta$. 
	
	Now let us show that the left side of \eqref{f0g0sum1} converges uniformly over $s, y \in \mathbb{C}$ with $|s| \le |y| \le \delta$ (so that we may take the limit there as $s$ and $y$ to $0$). To that end, observe under \eqref{x1x2bdestimate5} that 
	\begin{flalign}
	\label{fglambda} 
	\begin{aligned}
	\Big|(- y^{-1} s^2)^{|\boldsymbol{\lambda}| + n \binom{N}{2}}  F_{\boldsymbol{\lambda}} \big( \textbf{u}; y^{-1 / 2} \textbf{u}^{1 / 2} \boldsymbol{\mid} \textbf{y}; \textbf{s} \big) & (- y^{-1} s^2)^{- |\boldsymbol{\lambda}|}  G_{\boldsymbol{\lambda}} \big( \textbf{w}; y^{-1 / 2} \textbf{w}^{1 / 2} \boldsymbol{\mid} \textbf{y}; \textbf{s} \big) \Big| \\
	& \quad < C \binom{|\boldsymbol{\lambda}| + N}{N}^{2nN} (1 - \varepsilon)^{2 |\boldsymbol{\lambda}| - 2 nN^2},
	\end{aligned} 
	\end{flalign}
	
	\noindent for any $\boldsymbol{\lambda} \in \SeqSign_{n; N}$ and some constant $C > 0$ independent of $\boldsymbol{\lambda}$. Indeed, any path ensemble $\mathcal{E} \in \mathfrak{P}_G (\boldsymbol{\lambda} / \boldsymbol{\varnothing}; N)$ (from \Cref{pgpfph}; see the left side of \Cref{fgpaths}), has at least $|\boldsymbol{\lambda}| - nN^2$ vertices with arrow configuration of the form $(\textbf{e}_0, \textbf{B}; \textbf{e}_0; \textbf{B})$ for some $\textbf{B} \in \{ 0, 1 \}^n$ with $|\textbf{B}| > 0$, as there exist at most $nN^2$ vertices at which any of the $nN$ paths in the ensemble can go up. Similarly, any path ensemble $\mathcal{E} \in \mathfrak{P}_F (\boldsymbol{\lambda} / \boldsymbol{\varnothing})$ (from \Cref{pgpfph}; see the middle of \Cref{fgpaths}) has at least $|\boldsymbol{\lambda}| - nN^2$ vertices with arrow configuration of the form $(\textbf{e}_0, \textbf{B}; \textbf{e}_0, \textbf{B})$ with $|\textbf{B}| < n$. Moreover, there are at most $\binom{|\boldsymbol{\lambda}| + N}{N}$ choices for any of the $nN$ paths in $\mathcal{E} \in \mathfrak{P}_G (\boldsymbol{\lambda} / \boldsymbol{\varnothing}; N) \cup \mathfrak{P}_F (\boldsymbol{\lambda} / \boldsymbol{\varnothing})$, which gives $\big| \mathfrak{P}_G (\boldsymbol{\lambda} / \boldsymbol{\varnothing}; N) \mathfrak{P}_F (\boldsymbol{\lambda} / \boldsymbol{\varnothing}) \big| \le \binom{|\boldsymbol{\lambda}| + N}{N}^{2nN}$. Thus, \eqref{fglambda} follows from first applying \eqref{x1x2bdestimate5} (which holds for $|s| \le |y| \le \delta$) to bound the contribution of weights to the right of the $K$-th column in any such vertex model (where $K$ is from \eqref{x1x2bdestimate5}), and using the constant $C$ to account for the products of weights in the leftmost $K$ columns. 
	
	This implies that the sum on the left side of \eqref{f0g0sum1} converges uniformly for $|s| \le |y| \le \delta$ (and fixed $n, N \ge 1$). So, first letting $s$ tend to $0$ and then letting $y$ tend to $0$ in \eqref{f0g0sum1}, it follows that
	\begin{flalign}
	\label{sumf0g0sum2} 
	\begin{aligned} 
	\displaystyle\sum_{\boldsymbol{\lambda} \in \SeqSign_{n; N}} \mathcal{F}_{\boldsymbol{\lambda}} & \big( \textbf{u}; \infty \boldsymbol{\mid} 0; 0 \big) \mathcal{G}_{\boldsymbol{\lambda}} \big( \textbf{w}; \infty \boldsymbol{\mid} 0; 0 \big) \\
	& = \displaystyle\lim_{y \rightarrow 0} \displaystyle\prod_{j = 1}^N (-y)^{n (j - N)} (u_j y^{-1}; q)_n^{j - N} \displaystyle\prod_{i = 1}^M \displaystyle\prod_{j = 1}^N \displaystyle\frac{y^n (u_j y^{-1}; q)_n}{w_i^n (u_j w_i^{-1}; q)_n}.
	\end{aligned}
	\end{flalign}
	
	\noindent Thus, since 
	\begin{flalign}
	\label{sylimit}
	& \displaystyle\lim_{y \rightarrow 0} (-y)^n (u_j y^{-1}; q)_n = u_j^n q^{\binom{n}{2}}; \qquad (u_j w_i^{-1}; q)_n = (-u_j w_i^{-1})^n q^{\binom{n}{2}} (u_j^{-1} w_i; q^{-1})_n,
	\end{flalign}
	
	\noindent we deduce \eqref{fgsums1} from \eqref{sumf0g0sum2}.
\end{proof}

\begin{cor}
	
	\label{sxzlimitidentity3}
	
	We have that
	\begin{flalign*}
	\displaystyle\sum_{\boldsymbol{\lambda} \in \SeqSign_{n; N}} (-1)^{|\boldsymbol{\lambda}|} \mathcal{F}_{\boldsymbol{\lambda}} (\textbf{\emph{u}}; \infty \boldsymbol{\mid} 0; 0) \mathcal{G}_{\boldsymbol{\lambda}} (0; \textbf{\emph{t}} \boldsymbol{\mid} 0; 0) & = q^{- \binom{n}{2} \binom{N}{2}} \displaystyle\prod_{j = 1}^N u_j^{n (j - N)} \displaystyle\prod_{i = 1}^M \displaystyle\prod_{j = 1}^N (t_i^{-1} u_j^{-1}; q^{-1})_n.
	\end{flalign*}
	
\end{cor}

\begin{proof}
	
	Fix a complex number $y \in \mathbb{C}$, and denote the infinite sequence $\textbf{s} = (s, s, \ldots )$. Applying \Cref{fgsum2}, we obtain
	\begin{flalign}
	\label{sumf0g01}
	\begin{aligned} 
	\displaystyle\sum_{\boldsymbol{\lambda} \in \SeqSign_{n; N}} & (-1)^{|\boldsymbol{\lambda}|} (-s^{-2})^{- |\boldsymbol{\lambda}| - n \binom{N}{2}} F_{\boldsymbol{\lambda}} (\textbf{u}; s^{-1} \textbf{u}^{1 / 2} \boldsymbol{\mid} \textbf{s}^2; \textbf{s}^2) (s^2)^{-|\boldsymbol{\lambda}|} G_{\boldsymbol{\lambda}} (\textbf{s}^2; s \textbf{t}^{1 / 2} \boldsymbol{\mid} \textbf{s}^2; \textbf{s}^2) \\
	& \qquad \quad = (-s^2)^{-n \binom{N}{2}} (s^4; q)_n^{\binom{N}{2}} \displaystyle\prod_{j = 1}^N (s^{-2} u_j; q)_n^{j - N} \displaystyle\prod_{i = 1}^M \displaystyle\prod_{j = 1}^N \displaystyle\frac{(t_i u_j; q)_n}{s^{2n} t_i^n (s^{-2} u_j; q)_n},
	\end{aligned} 
	\end{flalign}
	
	\noindent where here we do not require the convergence condition since the sum on the left side of \eqref{sumf0g01} is finite, as there are only finitely many $\boldsymbol{\lambda} \in \Sign_{n; N}$ so that $G_{\boldsymbol{\lambda}} (\textbf{s}^2; s \textbf{t}^{1 / 2} \boldsymbol{\mid} \textbf{s}^2; \textbf{s}^2) \ne 0$. Indeed, to verify the latter point, observe by induction on $N$ that if there exists some index pair $(i, j) \in [1, n] \times [1, N]$ such that $\lambda_j^{(i)} > n (N + 1)$, then any path ensemble in the vertex model $\mathfrak{P}_G (\boldsymbol{\lambda} / \boldsymbol{\mu}; N)$ (from \Cref{pgpfph}; see also the left side of \Cref{fgpaths}) describing $G_{\boldsymbol{\lambda}}$ must contain a vertex of arrow configuration $(\textbf{A}, \textbf{B}; \textbf{C}, \textbf{D})$ with $B_i > C_i$ (where $\textbf{B} = (B_1, B_2, \ldots , B_n) \in \{ 0, 1 \}^n$ and $\textbf{C} = (C_1, C_2, \ldots , C_n) \in \{ 0, 1 \}^n$). By \Cref{z1weights}, the weight under $W_1$ of such a vertex is $0$, and so $G_{\boldsymbol{\lambda}} (\textbf{s}^2; s \textbf{t}^{1 / 2} \boldsymbol{\mid} \textbf{s}^2; \textbf{s}^2) = 0$ for all but finitely many $\boldsymbol{\lambda}$. 
	
	Now by \Cref{weights1s0z}, \Cref{sxs}, \Cref{limitsy}, \Cref{wlimits}, and \Cref{sgf0}, we deduce 
	\begin{flalign*}
	& \displaystyle\lim_{s \rightarrow 0} (-s^{-2})^{-|\boldsymbol{\lambda}| - n \binom{N}{2}} F_{\boldsymbol{\lambda}} (\textbf{u}; s^{-1} \textbf{u}^{1 / 2} \boldsymbol{\mid} \textbf{s}^2; \textbf{s}^2) = \mathcal{F}_{\lambda} (\textbf{u}; \infty \boldsymbol{\mid} 0; 0); \\
	& \displaystyle\lim_{s \rightarrow 0} (s^2)^{-|\boldsymbol{\lambda}|} G_{\boldsymbol{\lambda}} (\textbf{s}^2; s \textbf{t}^{1 / 2} \boldsymbol{\mid} \textbf{s}^2; \textbf{s}^2) = \mathcal{G}_{\boldsymbol{\lambda}} (0; \textbf{t} \boldsymbol{\mid} 0; 0).
	\end{flalign*}
	
	\noindent This, letting $s$ tend to $0$ in \eqref{sumf0g01}, and the facts that 
	\begin{flalign*}
	\displaystyle\lim_{s \rightarrow 0} s^{2n} (s^{-2} u_j; & q)_n = (-u_j)^n q^{\binom{n}{2}}; \qquad (t_i u_j; q)_n = q^{\binom{n}{2}} (-t_i u_j)^n (t_i^{-1} u_j^{-1}; q^{-1})_n,
	\end{flalign*}
	
	\noindent then together imply the corollary.
\end{proof}

\chapter{Degeneration to the LLT Polynomials}

\label{Polynomialsq}

In this chapter we explain through \Cref{limitg0} below (whose proof will appear in \Cref{ProofFGL}) how some of the degenerations of the $G$, $F$, and $H$ functions from \Cref{sgf0} can be matched with the Lascoux--Leclerc--Thibon polynomials originally introduced in \cite{RT}. 

\section{The LLT Polynomials}

\label{Partitionsn}

In this section we recall the definitions of the Lascoux--Leclerc--Thibon (LLT) polynomials introduced in \cite{RT}. Although these polynomials are typically indexed by partitions, here we will index them by signatures $\lambda = (\lambda_1, \lambda_2, \ldots , \lambda_{nN}) \in \Sign_{n N}$, where $n, N \ge 1$ are integers that will be fixed throughout this section. Any such signature $\lambda$ is associated with the partition obtained by removing the zero entries of $\lambda$, under which correpsondence the notions below will coincide with the ones from \cite{RT}. In particular, it is quickly verified that the definitions in this section will only depend on the nonzero entries of $\lambda$ (and not on the number of its entries equal to $0$).

Let us begin with some notation. Given two signatures $\lambda \in \Sign_{\ell}$ and $\mu \in \Sign_m$, we say $\mu \subseteq \lambda$\index{0@$\lambda, \mu$; typical signatures or partitions!$\mu \subseteq \lambda$} if $m \le \ell$ and $\mu_i \le \lambda_i$ for each $i \in [1, m]$. Moreover, for any sequences of $n$ signatures $\boldsymbol{\lambda}, \boldsymbol{\mu} \in \SeqSign_n$,\index{0@$\boldsymbol{\lambda}, \boldsymbol{\mu}$; typical signature sequences!$\boldsymbol{\mu} \subseteq \boldsymbol{\lambda}$} we say $\boldsymbol{\mu} \subseteq \boldsymbol{\lambda}$ if $\mu^{(i)} \subseteq \lambda^{(i)}$ for each $i \in [1, n]$. If $\mu \subseteq \lambda$, we call $\lambda / \mu$ a \emph{skew-shape}; the \emph{size} of this skew-shape is $|\lambda / \mu| = |\lambda| - |\mu| \ge 0$.\index{0@$\lambda, \mu$; typical signatures or partitions!$\mid$$\lambda / \mu$$\mid$; size of $\lambda / \mu$} Similarly, if $\boldsymbol{\mu} \subseteq \boldsymbol{\lambda}$, we call $\boldsymbol{\lambda} / \boldsymbol{\mu}$ a \emph{sequence of skew-shapes}; its size is $|\boldsymbol{\lambda} / \boldsymbol{\mu}| = |\boldsymbol{\lambda}| - |\boldsymbol{\mu}| \ge 0$.\index{0@$\boldsymbol{\lambda}, \boldsymbol{\mu}$; typical signature sequences!$\mid$$\boldsymbol{\lambda} / \boldsymbol{\mu}$$\mid$; size of $\boldsymbol{\lambda} / \boldsymbol{\mu}$} Similarly to our convention for sequence of signatures, we will denote any sequence $\boldsymbol{\lambda} / \boldsymbol{\mu}$ of $n$ skew-shapes by $\boldsymbol{\lambda} / \boldsymbol{\mu} = \big( \lambda^{(1)} / \mu^{(1)}, \lambda^{(2)} / \mu^{(2)}, \ldots , \lambda^{(n)} / \mu^{(n)} \big)$, unless mentioned otherwise. Any signature $\lambda$ is associated with the skew-shape $\lambda / \varnothing$, and any sequence $\boldsymbol{\lambda}$ of signatures is associated with the sequence of skew-shapes $\boldsymbol{\lambda} / \boldsymbol{\varnothing}$. 

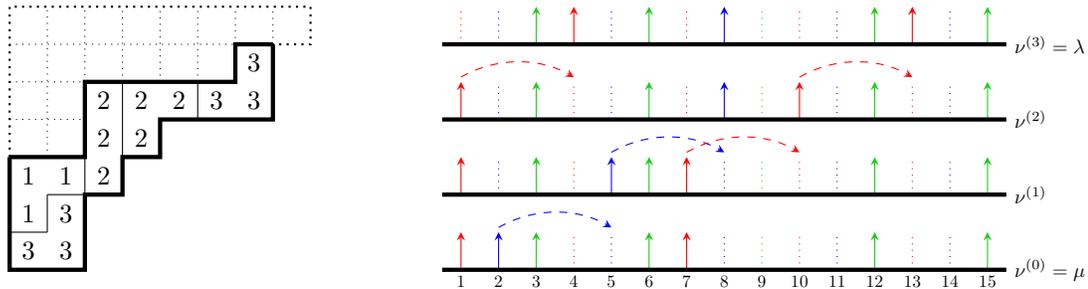
\begin{figure}

	\begin{center}

		\begin{tikzpicture}[
		>=stealth,
		scale = .5
		]
		
		\draw[ultra thick] (0, 4) -- (0, 1) -- (2, 1) -- (2, 3) -- (3, 3) -- (3, 4) -- (4, 4) -- (4, 5) -- (7, 5) -- (7, 7);
		\draw[ultra thick] (0, 4) -- (2, 4) -- (2, 6) -- (4, 6) -- (6, 6) -- (6, 7) -- (7, 7);
		\draw[thick, dotted] (0, 4) -- (0, 8) -- (8, 8) -- (8, 7) -- (7, 7);
		
		\draw (1, 2) -- (1, 3);
		\draw (2, 2) -- (2, 4);
		\draw (3, 4) -- (3, 6);
		\draw (5, 5) -- (5, 6);

		\draw (0, 2) -- (1, 2);
		\draw (1, 3) -- (2, 3);
		
		\draw[dotted] (1, 4) -- (1, 8);
		\draw[dotted] (2, 5) -- (2, 8);
		\draw[dotted] (3, 6) -- (3, 8);
		\draw[dotted] (4, 6) -- (4, 8);
		\draw[dotted] (5, 6) -- (5, 8);
		\draw[dotted] (6, 7) -- (6, 8);
		\draw[dotted] (7, 7) -- (7, 8);
		\draw[dotted] (0, 5) -- (2, 5);
		\draw[dotted] (0, 6) -- (4, 6);
		\draw[dotted]  (0, 7) -- (6, 7);
		
		\draw[] (5.5, 5.5) circle [radius = 0] node[]{$3$};
		\draw[] (6.5, 5.5) circle [radius = 0] node[]{$3$};
		\draw[] (6.5, 6.5) circle [radius = 0] node[]{$3$};
		
		\draw[] (3.5, 4.5) circle [radius = 0] node[]{$2$};
		\draw[] (3.5, 5.5) circle [radius = 0] node[]{$2$};
		\draw[] (4.5, 5.5) circle [radius = 0] node[]{$2$};
		
		\draw[] (2.5, 4.5) circle [radius = 0] node[]{$2$};
		\draw[] (2.5, 3.5) circle [radius = 0] node[]{$2$};
		\draw[] (2.5, 5.5) circle [radius = 0] node[]{$2$};
		
		\draw[] (.5, 1.5) circle [radius = 0] node[]{$3$};
		\draw[] (1.5, 2.5) circle [radius = 0] node[]{$3$};
		\draw[] (1.5, 1.5) circle [radius = 0] node[]{$3$};
		
		\draw[] (.5, 2.5) circle [radius = 0] node[]{$1$};
		\draw[] (.5, 3.5) circle [radius = 0] node[]{$1$};
		\draw[] (1.5, 3.5) circle [radius = 0] node[]{$1$}; 
		
		\draw[dotted, red] (12, 7) -- (12, 8);
		\draw[dotted, blue] (13, 7) -- (13, 8);
		\draw[->, green] (14, 7) -- (14, 8);
		\draw[->, red] (15, 7) -- (15, 8);
		\draw[dotted, blue] (16, 7) -- (16, 8);
		\draw[->, green] (17, 7) -- (17, 8);
		\draw[dotted, red] (18, 7) -- (18, 8);
		\draw[->, blue] (19, 7) -- (19, 8);
		\draw[dotted, green] (20, 7) -- (20, 8);
		\draw[dotted, red] (21, 7) -- (21, 8);
		\draw[dotted, blue] (22, 7) -- (22, 8);
		\draw[->, green] (23, 7) -- (23, 8);
		\draw[->, red] (24, 7) -- (24, 8);	
		\draw[dotted, blue] (25, 7) -- (25, 8);
		\draw[->, green] (26, 7) -- (26, 8);	
		
		\draw[->, red] (12, 5) -- (12, 6);
		\draw[dotted, blue] (13, 5) -- (13, 6);
		\draw[->, green] (14, 5) -- (14, 6);
		\draw[dotted, red] (15, 5) -- (15, 6);
		\draw[dotted, blue] (16, 5) -- (16, 6);
		\draw[->, green] (17, 5) -- (17, 6);
		\draw[dotted, red] (18, 5) -- (18, 6);
		\draw[->, blue] (19, 5) -- (19, 6);
		\draw[dotted, green] (20, 5) -- (20, 6);
		\draw[->, red] (21, 5) -- (21, 6);
		\draw[dotted, blue] (22, 5) -- (22, 6);
		\draw[->, green] (23, 5) -- (23, 6);
		\draw[dotted, red] (24, 5) -- (24, 6);	
		\draw[dotted, blue] (25, 5) -- (25, 6);
		\draw[->, green] (26, 5) -- (26, 6);	
		
		\draw[->, red] (12, 3) -- (12, 4);
		\draw[dotted, blue] (13, 3) -- (13, 4);
		\draw[->, green] (14, 3) -- (14, 4);
		\draw[dotted, red] (15, 3) -- (15, 4);
		\draw[->, blue] (16, 3) -- (16, 4);
		\draw[->, green] (17, 3) -- (17, 4);
		\draw[->, red] (18, 3) -- (18, 4);
		\draw[dotted, blue] (19, 3) -- (19, 4);
		\draw[dotted, green] (20, 3) -- (20, 4);
		\draw[dotted, red] (21, 3) -- (21, 4);
		\draw[dotted, blue] (22, 3) -- (22, 4);
		\draw[->, green] (23, 3) -- (23, 4);
		\draw[dotted, red] (24, 3) -- (24, 4);	
		\draw[dotted, blue] (25, 3) -- (25, 4);
		\draw[->, green] (26, 3) -- (26, 4);	
		
		\draw[->, red] (12, 1) node[below, black, scale = .65]{$1$} -- (12, 2);
		\draw[->, blue] (13, 1) node[below, black, scale = .65]{$2$} -- (13, 2);
		\draw[->, green] (14, 1) node[below, black, scale = .65]{$3$} -- (14, 2);
		\draw[dotted, red] (15, 1) node[below, black, scale = .65]{$4$} -- (15, 2);
		\draw[dotted, blue] (16, 1) node[below, black, scale = .65]{$5$} -- (16, 2);
		\draw[->, green] (17, 1) node[below, black, scale = .65]{$6$} -- (17, 2);
		\draw[->, red] (18, 1) node[below, black, scale = .65]{$7$} -- (18, 2);
		\draw[dotted, blue] (19, 1) node[below, black, scale = .65]{$8$} -- (19, 2);
		\draw[dotted, green] (20, 1) node[below, black, scale = .65]{$9$} -- (20, 2);
		\draw[dotted, red] (21, 1) node[below, black, scale = .65]{$10$} -- (21, 2);
		\draw[dotted, blue] (22, 1) node[below, black, scale = .65]{$11$} -- (22, 2);
		\draw[->, green] (23, 1) node[below, black, scale = .65]{$12$} -- (23, 2);
		\draw[dotted, red] (24, 1) node[below, black, scale = .65]{$13$} -- (24, 2);	
		\draw[dotted, blue] (25, 1) node[below, black, scale = .65]{$14$} -- (25, 2);
		\draw[->, green] (26, 1) node[below, black, scale = .65]{$15$} -- (26, 2); 
		
		\draw[dashed, blue, ->] (13, 2.125) arc (120:60:3); 
		\draw[dashed, blue, ->] (16, 4.125) arc (120:60:3); 
		\draw[dashed, red, ->] (18, 4.125) arc (120:60:3); 
		\draw[dashed, red, ->] (12, 6.125) arc (120:60:3); 
		\draw[dashed, red, ->] (21, 6.125) arc (120:60:3); 
		
		\draw[ultra thick] (11.5, 7) -- (26.5, 7) node[right, scale = .75]{$\nu^{(3)} = \lambda$};
		\draw[ultra thick] (11.5, 5) -- (26.5, 5) node[right, scale = .75]{$\nu^{(2)}$};
		\draw[ultra thick] (11.5, 3) -- (26.5, 3) node[right, scale = .75]{$\nu^{(1)}$};
		\draw[ultra thick] (11.5, 1) -- (26.5, 1) node[right, scale = .75]{$\nu^{(0)} = \mu$};

		\end{tikzpicture}
		
	\end{center}
	
	\caption{\label{skew} To the left is a $3$-ribbon tableau for the skew shape $\lambda / \mu$, where $\lambda = (8, 7, 7, 4, 3, 2, 2)$ and $\mu = (8, 6, 2, 2)$; to the right is the associated sequence of colored Maya diagrams. Here, red is color $1$, blue is color $2$, and green is color $3$.} 
	
\end{figure}

We depict any signature $\lambda = (\lambda_1, \lambda_2, \ldots , \lambda_{\ell})$ through its \emph{Young diagram} $\mathscr{Y} (\lambda)$,\index{Y@$\mathscr{Y} (\lambda / \mu)$; Young diagram of $\lambda / \mu$!$\mathscr{Y} (\lambda)$; Young diagram of $\lambda$} which is a left-justified collection of $|\lambda|$ boxes arranged in $\ell (\lambda) = \ell$ rows, whose $i$-th row from the top contains $\lambda_i$ boxes. The condition $\mu \subseteq \lambda$ is then equivalent to imposing that $\mathscr{Y} (\mu)$ is contained inside $\mathscr{Y} (\lambda)$, where these Young diagrams are translated so that their top-left corners coincide.  The Young diagram $\mathscr{Y} (\lambda / \mu)$ for the skew-shape $\lambda / \mu$ is then given by $\mathscr{Y} (\lambda) \setminus \mathscr{Y} (\mu)$, obtained by removing the Young diagram of $\mu$ from that of $\lambda$; we refer to the left side of \Cref{skew} for a depiction, where there dotted part is $\mathscr{Y} (\mu)$ and the solid part is $\mathscr{Y} (\lambda / \mu)$.\index{Y@$\mathscr{Y} (\lambda / \mu)$; Young diagram of $\lambda / \mu$}

A skew-shape $\lambda / \mu$ is \emph{connected} if its Young diagram is edge-connected, that is, for any two boxes $b, b' \in \mathscr{Y} (\lambda / \mu)$ there exists a sequence of boxes $b_0 = b, b_1, \ldots , b_k = b' \in \mathscr{Y} (\lambda / \mu)$ such that $b_i$ and $b_{i + 1}$ share an edge for each $i \in [0, k - 1]$. A \emph{ribbon} is then a connected skew-shape that does not contain any $2 \times 2$ block of boxes, and an \emph{$n$-ribbon} is a ribbon of size $n$; we refer to the left side of \Cref{skew} for a tiling of a skew-shape by $3$-ribbons. The \emph{height} $\text{ht} (\lambda / \mu)$\index{0@$\lambda, \mu$; typical signatures or partitions!$\height(\lambda / \mu)$; height of a ribbon} of an $n$-ribbon $\lambda / \mu$ is the number of rows it spans; following equation (15) of \cite{RT}, the \emph{spin} $\text{sp} (\lambda / \mu)$ of $\lambda / \mu$\index{0@$\lambda, \mu$; typical signatures or partitions!$\spin (\lambda / \mu)$; spin statistic} is defined by
\begin{flalign}
\label{lambdamuheight}
\text{sp} (\lambda / \mu) = \displaystyle\frac{\text{ht} (\lambda / \mu) - 1}{2}.
\end{flalign}

A \emph{horizontal $n$-ribbon strip} is a (not necessarily connected) skew-shape $\lambda / \mu$ that can be tiled by $n$-ribbons, each of whose top-right boxes is adjacent to the upper boundary of $\lambda / \mu$. The spin $\text{sp} (\lambda / \mu)$ of this horizontal $n$-ribbon strip is given by the sum of the spins of the $n$-ribbons composing it. 

An \emph{$n$-ribbon tableau} for a skew-shape $\lambda / \mu$ is a sequence of signatures $T = \big( \nu^{(0)}, \nu^{(1)}, \ldots \nu^{(k)} \big)$ such that $\mu = \nu^{(0)} \subseteq \nu^{(1)} \subseteq \cdots \subseteq \nu^{(k)} = \lambda$ and $\nu^{(i)} / \nu^{(i - 1)}$ is a horizontal $n$-ribbon strip for each index $i \in [1, k]$. We associate such a tableau $T$ with the function $T: \mathscr{Y} (\lambda / \mu) \rightarrow [1, k]$ defined by, for any box $b \in \mathscr{Y} (\lambda / \mu)$, setting $T(b) = i$ where $i$ is such that $b \in \mathscr{Y} \big( \nu^{(i)} / \nu^{(i - 1)} \big)$; see the left side of \Cref{skew} for a depiction. The spin $\text{sp} (T)$ of $T$ is then given by the sum of the spins of the horizontal $n$-ribbon strips composing it, namely,
\begin{flalign*}
\text{sp} (T) = \displaystyle\sum_{i = 1}^k \text{sp} \big( \nu^{(i)} / \nu^{(i - 1)} \big).
\end{flalign*} 

\noindent For example, the spin of the $3$-ribbon tableaux shown on the left side of \Cref{skew} is $3$. 

For any (possibly infinite) sequence of complex numbers $\textbf{x} = (x_1, x_2, \ldots , x_N)$; skew-shape $\lambda / \mu$; and $n$-ribbon tableau $T$ for $\lambda / \mu$, define the monomial  
\begin{flalign*}
\textbf{x}^T = \displaystyle\prod_{b \in \mathscr{Y} (\lambda / \mu)} x_{T (b)}^{1 / n},
\end{flalign*}

\noindent where we set $\textbf{x}^T = 0$ if $N$ is less than the number of horizontal $n$-ribbon strips composing $T$. Observe that $\textbf{x}^T$ indeed defines a monomial since, for any index $i$, the number of boxes $b \in \mathscr{Y} (\lambda / \mu)$ such that $T(b) = i$ is $\big| \nu^{(i)} / \nu^{(i - 1)} \big|$, which is a multiple of $n$ (where this multiple is given by the number of $n$-ribbons composing $\nu^{(i)} / \nu^{(i - 1)}$). 

Now, for any complex number $q \in \mathbb{C}$ and finite sequence $\textbf{x} = (x_1, x_2, \ldots , x_N)$ of complex numbers, we can define the following polynomials, which originated (up to an overall factor of $q$) in equation (26) of \cite{RT}. For any signatures $\lambda, \mu \in \Sign$, the \emph{Lascoux--Leclerc--Thibon (LLT) polynomial} $\mathcal{L}_{\lambda / \mu}^{(n)} (\textbf{x}) = \mathcal{L}_{\lambda / \mu}^{(n)} (\textbf{x}; q)$\index{L@$\mathcal{L}_{\boldsymbol{\lambda} / \boldsymbol{\mu}} (\textbf{x})$; LLT polynomial!$\mathcal{L}_{\lambda / \mu}^{(n)} (\textbf{x})$} is defined by setting
\begin{flalign}
\label{lambdamug} 
\mathcal{L}_{\lambda / \mu}^{(n)} (\textbf{x}) = \displaystyle\sum_T q^{\text{sp} (T)} \textbf{x}^T,
\end{flalign}

\noindent where $T$ is summed over all $n$-ribbon tableaux for $\lambda / \mu$. Observe in particular that these functions are homogeneous of degree $n^{-1} \big( |\lambda| - |\mu| \big)$ in $\textbf{x}$. The fact that they are symmetric in $\textbf{x}$ is not immediate from the definition \eqref{lambdamug} but is known from Theorem 6.1 of \cite{RT}, and it will also follow as a consequence of \eqref{1gl} below (and \Cref{gxfxsigma}). 

\begin{rem}
	
	\label{lbranching} 
	
	The LLT polynomials $\mathcal{L}_{\lambda / \mu}^{(n)}$ satisfy a branching relation (similar to \Cref{fghbranching}) in the following sense. If $\textbf{x} = (x_1, x_2, \ldots , x_{K + M}) \subset \mathbb{C}$, and we denote $\textbf{x}' = (x_1, x_2, \ldots , x_K) \subset \mathbb{C}$ and $\textbf{x}'' = (x_{K + 1}, x_{K + 2}, \ldots , x_{K + M}) \subset \mathbb{C}$, then
	\begin{flalign}
	\label{branchingl} 
	\mathcal{L}_{\lambda / \mu}^{(n)} (\textbf{x}) = \displaystyle\sum_{\nu \in \Sign_{\ell}} \mathcal{L}_{\lambda / \nu}^{(n)} (\textbf{x}'') \mathcal{L}_{\nu / \mu}^{(n)} (\textbf{x}'),
	\end{flalign}
	
	\noindent for any sufficiently large integer $\ell$. Indeed, this follows from decomposing the ribbon tableau $T$ for $\lambda / \mu$ on the right side of \eqref{lambdamug} into a union of ribbon tableaux $(T', T'')$, where $T'$ and $T''$ consist of those boxes in $T$ with entries in the intervals $[1, K]$ and $[K + 1, K + M]$, respectively. Then, the $\nu \in \mathbb{Y}$ on the right side of \eqref{branchingl} is the partition such that $T'$ and $T''$ are ribbon tableaux for $\nu / \mu$ and $\lambda / \nu$, respectively. 
	
\end{rem}

\section{$n$-Quotients and Colored Maya Diagrams}

\label{FunctionGPartitions}

The LLT polynomials defined in \eqref{lambdamug} are indexed by skew-shapes $\lambda / \mu$, while our functions from \Cref{fgdefinition} are indexed by sequences of $n$ skew-shapes $\boldsymbol{\lambda} / \boldsymbol{\mu}$. The correspondence between these two different indexings is given by setting $\boldsymbol{\lambda} / \boldsymbol{\mu}$ to be the $n$-quotient of $\lambda / \mu$. In this section we recall this procedure in detail, as this will in any case be useful for degenerating our functions to the LLT polynomials. Throughout this section, we again fix an integer $n \ge 1$. 

Recall from \eqref{t} that, for any $\lambda \in \Sign_{\ell}$, we set $\mathfrak{T} (\lambda) = (\lambda_1 + \ell, \lambda_2 + \ell - 1, \ldots , \lambda_{\ell} + 1) \in \mathbb{Z}_{> 0}^{\ell}$, which is sometimes referred to as the \emph{Maya diagram} for $\lambda$. We visually represent this Maya diagram $\mathfrak{T} (\lambda)$ through a particle configuration on the positive integer lattice $\mathbb{Z}_{> 0}$ where a site $i \in \mathbb{Z}_{> 0}$ is occupied by a particle, depicted by an arrow\footnote{Typically, these are depicted by beads, but here we use arrows to emphasize the similarity with vertex models.}, if and only if $i \in \mathfrak{T} (\lambda)$. We further \emph{color} this Maya diagram by associating with each arrow $i \in \mathfrak{T} (\lambda)$ the (unique) index $c (i) \in \{ 1, 2, \ldots , n \}$ such that $n$ divides $i - c(i)$; this yields a \emph{colored Maya diagram} (more specifically, an \emph{$n$-colored Maya diagram}). The right side of \Cref{skew} depicts $3$-colored Maya diagrams associated with all signatures in the $3$-ribbon tableau $(8, 6, 2, 2, 0, 0, 0) \subseteq (8, 6, 2, 2, 2, 1, 0) \subseteq (8, 6, 5, 4, 3, 1, 0) \subseteq (8, 7, 7, 4, 3, 2, 2)$.

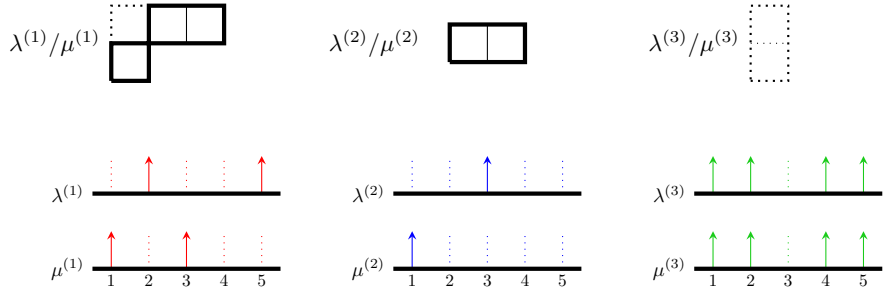
\begin{figure}

	\begin{center}

		\begin{tikzpicture}[
		>=stealth,
		scale = .5
		]
		
		\draw[thick, dotted] (1, 6) -- (1, 7) -- (2, 7);
		\draw[ultra thick] (1, 5) -- (2, 5) -- (2, 6) -- (4, 6) -- (4, 7) -- (2, 7) -- (2, 6) -- (1, 6) -- (1, 5);
		\draw (3, 6) -- (3, 7);
		
		\draw[ultra thick] (10, 5.5) -- (12, 5.5) -- (12, 6.5) -- (10, 6.5) -- (10, 5.5);
		\draw (11, 5.5) -- (11, 6.5);
		
		\draw[thick, dotted] (18, 5) -- (19, 5) -- (19, 7) -- (18, 7) -- (18, 5);
		\draw[dotted] (18, 6) -- (19, 6);
			
		\draw[] (-.5, 6) circle[radius = 0] node[scale = .9]{$\lambda^{(1)} / \mu^{(1)}$};	
		\draw[] (8, 6) circle[radius = 0] node[scale = .9]{$\lambda^{(2)} / \mu^{(2)}$};	
		\draw[] (16.5, 6) circle[radius = 0] node[scale = .9]{$\lambda^{(3)} / \mu^{(3)}$};		
				
		\draw[red, ->] (1, 0) node[below, black, scale = .65]{$1$} -- (1, 1);
		\draw[red, dotted] (2, 0) node[below, black, scale = .65]{$2$} -- (2, 1);
		\draw[red, ->] (3, 0) node[below, black, scale = .65]{$3$}  -- (3, 1);
		\draw[red, dotted] (4, 0) node[below, black, scale = .65]{$4$}  -- (4, 1);
		\draw[red, dotted] (5, 0) node[below, black, scale = .65]{$5$}  -- (5, 1);
		
		\draw[red, dotted]  (1, 2) -- (1, 3);
		\draw[red, ->] (2, 2) -- (2, 3);
		\draw[red, dotted] (3, 2) -- (3, 3);
		\draw[red, dotted] (4, 2) -- (4, 3);
		\draw[red, ->] (5, 2) -- (5, 3);
		
		\draw[blue, ->] (9, 0) node[below, black, scale = .65]{$1$}  -- (9, 1); 
		\draw[blue, dotted] (10, 0) node[below, black, scale = .65]{$2$}  -- (10, 1); 
		\draw[blue, dotted] (11, 0) node[below, black, scale = .65]{$3$} -- (11, 1); 
		\draw[blue, dotted] (12, 0) node[below, black, scale = .65]{$4$} -- (12, 1); 
		\draw[blue, dotted] (13, 0) node[below, black, scale = .65]{$5$} -- (13, 1); 
		
		\draw[blue, dotted] (9, 2) -- (9, 3); 
		\draw[blue, dotted] (10, 2) -- (10, 3); 
		\draw[blue, ->] (11, 2) -- (11, 3); 
		\draw[blue, dotted] (12, 2) -- (12, 3); 
		\draw[blue, dotted] (13, 2) -- (13, 3); 
		
		\draw[green, ->] (17, 0) node[below, black, scale = .65]{$1$} -- (17, 1);
		\draw[green, ->] (18, 0) node[below, black, scale = .65]{$2$} -- (18, 1);
		\draw[green, dotted] (19, 0) node[below, black, scale = .65]{$3$} -- (19, 1);
		\draw[green, ->] (20, 0) node[below, black, scale = .65]{$4$} -- (20, 1);
		\draw[green, ->] (21, 0) node[below, black, scale = .65]{$5$} -- (21, 1);
		
		\draw[green, ->] (17, 2) -- (17, 3);
		\draw[green, ->] (18, 2) -- (18, 3);
		\draw[green, dotted] (19, 2) -- (19, 3);
		\draw[green, ->] (20, 2) -- (20, 3);
		\draw[green, ->] (21, 2) -- (21, 3);
		
		\draw[ultra thick] (.5, 2) node[left, scale  = .75]{$\lambda^{(1)}$} -- (5.5, 2) ;
		\draw[ultra thick] (.5, 0) node[left, scale = .75]{$\mu^{(1)}$} -- (5.5, 0);
		
		\draw[ultra thick] (8.5, 2) node[left, scale = .75]{$\lambda^{(2)}$} -- (13.5, 2) ;
		\draw[ultra thick] (8.5, 0) node[left, scale = .75]{$\mu^{(2)}$} -- (13.5, 0);
		
		\draw[ultra thick] (16.5, 2) node[left, scale = .75]{$\lambda^{(3)}$} -- (21.5, 2);
		\draw[ultra thick] (16.5, 0) node[left, scale = .75]{$\mu^{(3)}$} -- (21.5, 0);

		\end{tikzpicture}
		
	\end{center}
	
	\caption{\label{skewquotient} To the top is the $3$-quotient for $\lambda / \mu$, where $\lambda = (8, 7, 7, 4, 3, 2, 2)$ and $\mu = (8, 6, 2, 2)$. To the bottom are the Maya diagrams for the $3$-quotients of $\lambda$ and $\mu$.} 
	
\end{figure}

Any Maya diagram $\mathfrak{T} (\lambda)$ gives rise to a sequence of $n$ Maya diagrams $\big( \mathcal{T}^{(1)}, \mathcal{T}^{(2)}, \ldots , \mathcal{T}^{(n)} \big)$ by setting $j \in \mathcal{T}^{(i)}$ if and only if $n (j - 1) + i \in \mathcal{T}$, for any indices $i \in [1, n]$ and $j \in \mathbb{Z}_{> 0}$. In this way, one might view $\mathcal{T}^{(i)}$ as the restriction to color $i$ of the Maya diagram $\mathfrak{T} (\lambda)$. Since $\mathfrak{T}$ induces a bijection between signatures and subsets of $\mathbb{Z}_{> 0}$, there exists a sequence of $n$ signatures $\boldsymbol{\lambda} = \big( \lambda^{(1)}, \lambda^{(2)}, \ldots , \lambda^{(n)} \big) \in \SeqSign_n$ such that $\mathcal{T}^{(i)} = \mathfrak{T} \big( \lambda^{(i)} \big)$ for each $i \in [1, n]$; the sequence $\boldsymbol{\lambda}$ is called the \emph{$n$-quotient} of $\lambda$. The $n$-quotient of a skew shape $\lambda / \mu$ is then given by the sequence of skew shapes $\boldsymbol{\lambda} / \boldsymbol{\mu} = \big( \lambda^{(1)} / \mu^{(1)}, \lambda^{(2)} / \mu^{(2)}, \ldots , \lambda^{(n)} / \mu^{(n)} \big)$ where $\boldsymbol{\lambda}$ and $\boldsymbol{\mu}$ are the $n$-quotients of $\lambda$ and $\mu$, respectively. The $n$-quotient for the skew shape $\lambda / \mu$ from \Cref{skew}, which is given by $\big( (3, 1) / (1, 0), (2) / (0), (1, 1, 0, 0) / (1, 1, 0, 0) \big)$, is depicted in \Cref{skewquotient}. Under this notation, we alternatively write the LLT polynomials from \eqref{lambdamug} by 
\begin{flalign}
\label{llambdamu} 
\mathcal{L}_{\boldsymbol{\lambda} / \boldsymbol{\mu}} (\textbf{x}) = \mathcal{L}_{\boldsymbol{\lambda} / \boldsymbol{\mu}} (\textbf{x}; q) = \mathcal{L}_{\lambda / \mu}^{(n)} (\textbf{x}). 
\end{flalign}
\index{L@$\mathcal{L}_{\boldsymbol{\lambda} / \boldsymbol{\mu}} (\textbf{x})$; LLT polynomial}

Now let us introduce some (nonstandard) terminology on how the notions from \Cref{Partitionsn} (such as ribbons, horizontal ribbon strips, and ribbon tableaux) can be interpreted in terms of Maya diagrams. To that end, let $\mathcal{U} = (u_1, u_2, \ldots , u_k) \in \mathbb{Z}_{> 0}^k$ denote a Maya diagram. An \emph{$n$-jump} of $\mathcal{U}$ is a Maya diagram $\mathcal{T} = (t_1, t_2, \ldots , t_k) \in \mathbb{Z}_{> 0}^k$ for which there exist indices $i, j \in [1, k]$ such that $u_i + n = t_j$ and $\mathcal{U} \setminus \{ u_i \} = \mathcal{T} \setminus \{ t_j \}$. Stated alternatively, $\mathcal{T}$ is obtained from $\mathcal{U}$ by having the particle at site $u_i$ ``jump'' to the right by $n$ spaces, so we call $u_i$ the \emph{jumper} of this $n$-jump. The quantity $\big| [u_i + 1, u_i + n - 1] \cap \mathcal{U} \big|$ then denotes the number of particles that the jumper jumps over. 

For signatures $\lambda, \mu \in \Sign_N$, it is quickly verified that $\mathfrak{T} (\lambda)$ is an $n$-jump of $\mathfrak{T} (\mu)$ if and only if $\lambda / \mu$ is an $n$-ribbon. In this case, the number of particles jumped over is $\height (\lambda / \mu) - 1 = 2 \spin (\lambda / \mu)$, where we recall from \Cref{Partitionsn} that $\height (\lambda / \mu)$ denotes the height of the ribbon $\lambda / \mu$ and $\spin (\lambda / \mu)$ denotes its spin given by \eqref{lambdamuheight}. We refer to the bottom of the right side of \Cref{skew} for an example, where $\mathfrak{T} \big( \nu^{(1)} \big) = (1, 3, 5, 6, 7, 12, 15)$ is a $3$-jump of $\mathfrak{T} \big( \nu^{(0)} \big) = (1, 2, 3, 6, 7, 12, 15)$; the jumper there is $2$, and it jumps over one particle at site $3$.

Next, we say that $\mathcal{T}$ is obtained by a \emph{series of $n$-jumps} from $\mathcal{U}$ if there exists a sequence $\mathcal{U} = \mathcal{T}_0, \mathcal{T}_1, \ldots , \mathcal{T}_m = \mathcal{T}$ of Maya diagrams such that $\mathcal{T}_i$ is an $n$-jump of $\mathcal{T}_{i - 1}$ for each $i \in [1, m]$. Such a series is called \emph{nonconflicting} if, denoting the jumper from $\mathcal{T}_{i - 1}$ to $\mathcal{T}_i$ by $j_{i - 1} \in \mathcal{T}_{i - 1}$ for each $i \in [1, m]$, we have $j_0 < j_1 < \cdots < j_{m - 1}$. It is quickly verified that $\lambda / \mu$ is a horizontal $n$-ribbon strip if and only if $\mathfrak{T} (\lambda)$ can be obtained by a nonconflicting series of $n$-jumps from $\mathfrak{T} (\mu)$. In this case, such a nonconflicting series must be unique, and $2 \spin (\lambda / \mu)$ denotes the total number of particles jumped over in this series. 

Thus, $n$-ribbon tableaux $T$ of $\lambda / \mu$ are in bijection with sequences $(\mathcal{T}_0, \mathcal{T}_1, \ldots , \mathcal{T}_m)$ of Maya diagrams such that $\mathcal{T}_0 = \mathfrak{T} (\mu)$, $\mathcal{T}_m = \mathfrak{T} (\lambda)$, and $\mathcal{T}_i$ can be obtained by a nonconflicting series of $n$-jumps from $\mathcal{T}_{i - 1}$ for each $i \in [1, m]$.  In this case, $2 \spin (T)$ gives the total number of particles jumped over in this sequence of nonconflicting series of $n$-jumps. For example, the left side of \Cref{skew} depicts a $3$-ribbon tableau, with each $\mathfrak{T} \big( \nu^{(i)} \big)$ obtained by a non-conflicting series of $n$-jumps from $\mathfrak{T} \big( \nu^{(i - 1)} \big)$; the total number of particles jumped over there is $6$. 

For a more detailed exposition on some of these notions (under slightly different notation), we refer to Section 2 and Section 3 of \cite{ESRTAAP}.

\section{Correspondences With LLT Polynomials}

\label{FunctionsFGL}

In this section we state \Cref{limitg0} below, which indicates how some of the various specializations of the $F$, $G$, and $H$ functions introduced in \Cref{sgf0} degenerate to the LLT polynomials given by \eqref{lambdamug} (or, equivalently, \eqref{llambdamu}). In what follows, for any sequence of signatures $\boldsymbol{\lambda} \in \SeqSign_n$, set
\begin{flalign}
\label{lambdamupsi} 
\psi (\boldsymbol{\lambda}) = \displaystyle\frac{1}{2} \displaystyle\sum_{1 \le i < j \le n} \displaystyle\sum_{a \in \mathfrak{T}_i} \displaystyle\sum_{b \in \mathfrak{T}_j} \textbf{1}_{a > b},
\end{flalign}
\index{0@$\psi$}

\noindent where we have abbreviated $\mathfrak{T}_k = \mathfrak{T} \big( \lambda^{(k)} \big)$ for each index $k \in [1, n]$. 

\begin{example} 
	
\label{psi0n}

For any integer $N \ge 0$, we have 	
\begin{flalign}
\label{psi0}
\psi (\textbf{0}^N) = \frac{1}{2} \displaystyle\sum_{i = 1}^n (n - i) \displaystyle\sum_{a = 1}^N (a - 1) = \displaystyle\frac{1}{2} \binom{n}{2} \binom{N}{2}.
\end{flalign} 

\end{example}

Moreover, for any skew-shape $\lambda / \mu$, let $\lambda' / \mu'$ denote its \emph{dual},\index{0@$\lambda, \mu$; typical signatures or partitions!$\lambda' / \mu'$; dual} whose Young diagram $\mathscr{Y} (\lambda' / \mu')$ is obtained by reflecting $\mathscr{Y} (\lambda / \mu)$ across the diagonal line $x + y = 0$ (where the top-left corner of the top-left box in $\mathscr{Y} (\lambda/ \mu)$ is located at $(0, 0)$); see the top of \Cref{lambdamulambdamu1} below for a depiction.\footnote{In principle, there is an ambiguity in the numbers of entries equal to zero in $\lambda'$ and $\mu'$. However, we will ignore this point, since (as mentioned in the beginning of \Cref{Partitionsn}) $\mathcal{L}_{\lambda / \mu}^{(n)} (\textbf{x}; q)$ does not depend on these numbers.} In addition, if $\boldsymbol{\lambda} / \boldsymbol{\mu}$ denotes the $n$-quotient of $\lambda / \mu$, then we let $\boldsymbol{\lambda}' / \boldsymbol{\mu}'$ denote\index{0@$\boldsymbol{\lambda}, \boldsymbol{\mu}$; typical signature sequences!$\boldsymbol{\lambda}' / \boldsymbol{\mu}'$; dual} that of its dual\footnote{Let us clarify that, under this notation, the $i$-th component $\lambda'^{(i)}$ of $\boldsymbol{\lambda}'$ happens to not be the dual of $\lambda^{(i)}$. Instead, it is that of $\lambda^{(n - i)}$ (and similarly for $\boldsymbol{\mu}$).} $\lambda' / \mu'$. 

Further recall for any signature sequence $\boldsymbol{\lambda} = \big( \lambda^{(1)}, \lambda^{(2)}, \ldots , \lambda^{(n)} \big) \in \SeqSign_n$ that its reverse ordering is denoted by $\overleftarrow{\boldsymbol{\lambda}} = \big( \lambda^{(n)}, \lambda^{(n - 1)}, \ldots , \lambda^{(1)} \big) \in \SeqSign_n$.\index{X@$\overleftarrow{\mathscr{X}}$; reverse ordering of $\mathscr{X}$} We may view the latter as obtained from the former by reversing its colors (that is, interchanging colors $i$ and $n - i + 1$ for each $i \in [1, n]$).

We can now state the following theorem, which describes four ways of obtaining the LLT polynomials from the degenerated functions $\mathcal{F}$, $\mathcal{G}$, and $\mathcal{H}$ from \Cref{sgf0}. 

\begin{thm}
	
	\label{limitg0}
	
	For any fixed integers $N \ge 1$ and $M \ge 0$, and set $\textbf{\emph{x}} = (x_1, x_2, \ldots , x_N)$ of complex numbers, the following statements hold.
	
	\begin{enumerate} 
	
	\item \label{lg1} For any sequences of signatures $\boldsymbol{\lambda}, \boldsymbol{\mu} \in \SeqSign_{n; M}$, we have
	\begin{flalign}
	\label{1gl} 
	\mathcal{G}_{\boldsymbol{\lambda} / \boldsymbol{\mu}} (\textbf{\emph{x}}; \infty \boldsymbol{\mid} 0; 0)  = q^{\psi (\boldsymbol{\lambda}) - \psi (\boldsymbol{\mu})} \mathcal{L}_{\boldsymbol{\lambda} / \boldsymbol{\mu}} (\textbf{\emph{x}}; q). 
	\end{flalign}
	
	\item \label{lf1} For any sequences of signatures $\boldsymbol{\lambda} \in \Sign_{n; M + N}$ and $\boldsymbol{\mu} \in \SeqSign_{n; M}$, we have
	\begin{flalign}
	\label{1fl} 
	\mathcal{F}_{\boldsymbol{\lambda} / \boldsymbol{\mu}} (\textbf{\emph{x}}; \infty \boldsymbol{\mid} 0; 0)  = q^{\psi (\boldsymbol{\mu}) - \psi (\boldsymbol{\lambda}) + \binom{M}{2} \binom{n}{2} / 2 - \binom{M + N}{2} \binom{n}{2} / 2} \mathcal{L}_{\boldsymbol{\lambda} / \boldsymbol{\mu}} (q^{1 - n} \textbf{\emph{x}}^{-1}; q) \displaystyle\prod_{j = 1}^N x_j^{n (j - M - N)}.
	\end{flalign}
	
	\item \label{lg2} For any sequences of signatures $\boldsymbol{\lambda}, \boldsymbol{\mu} \in \SeqSign_{n; M}$, we have
	\begin{flalign}
	\label{2gl} 
	\mathcal{G}_{\boldsymbol{\lambda} / \boldsymbol{\mu}} (0; \textbf{\emph{x}} \boldsymbol{\mid} 0; 0) = q^{\psi (\boldsymbol{\lambda}) - \psi (\boldsymbol{\mu})} \mathcal{L}_{\boldsymbol{\lambda}' / \boldsymbol{\mu}'} (q^{(n - 1) / 2} \textbf{\emph{x}}^{-1}; q^{-1}).
	\end{flalign}
	
	\item \label{lh1} For any sequences of signatures $\boldsymbol{\lambda} \in \SeqSign_{n; M + N}$ and $\boldsymbol{\mu} \in \SeqSign_{n; M}$, we have
		\begin{flalign}
		\label{1hl} 
		\mathcal{H}_{\boldsymbol{\lambda} / \boldsymbol{\mu}} (\textbf{\emph{x}}; \infty \boldsymbol{\mid} \infty; \infty) = q^{\psi (\overleftarrow{\boldsymbol{\mu}}) - \psi (\overleftarrow{\boldsymbol{\lambda}}) + \binom{M + N}{2} \binom{n}{2} / 2 - \binom{M}{2} \binom{n}{2} / 2} \mathcal{L}_{\overleftarrow{\boldsymbol{\lambda}} / \overleftarrow{\boldsymbol{\mu}}} (-\textbf{\emph{x}}^{-1}; q^{-1}) \displaystyle\prod_{j = 1}^N x_j^{n (j - M - N - 1)}.
		\end{flalign}

	 \end{enumerate}
	 
\end{thm} 

\begin{rem} 

\label{lambdallambdag} 

The skew-shapes $\lambda / \mu$ whose $n$-quotients are of the form $\boldsymbol{\lambda} / \boldsymbol{\mu} \in \SeqSign_{n; M}$, for some $M \ge 1$, are those for which $\lambda$ and $\mu$ both have an empty $n$-core.\footnote{The \emph{$n$-core} of a signature $\lambda \in \Sign$ is a signature $\kappa \in \Sign$ of minimal size such that $\lambda / \kappa$ admits an $n$-ribbon tableau. This signature $\kappa$ is unique, up to its number of entries equal to $0$.} To allow for more general $\lambda / \mu$, one must permit the signatures in $\boldsymbol{\lambda}$ (and $\boldsymbol{\mu}$) to have different lengths. The framework introduced above allows this with very minor modifications (see \Cref{fgmn}), under which \eqref{1gl} and \eqref{2gl} continue to hold as written, and the exponents of $q$ and the $x_j$ on the right side of \eqref{1fl} become slightly different. However, to ease notation, we will not pursue this generalization here in detail.

\end{rem}

In addition to the curiosity that the LLT polynomials can be recovered from the $F$, $G$, and $H$ functions in numerous ways,\footnote{The list of degenerations in \Cref{limitg0} from these functions to the LLT polynomials is likely not at all exhaustive, but we will not pursue further ones here.} the four parts of \Cref{limitg0} will have their own uses and interpretations. More specifically, we will see as \Cref{suml1} below that \eqref{1gl} and \eqref{1fl} together imply that pairing $\mathcal{G}_{\boldsymbol{\lambda}} (\textbf{y}; \infty \boldsymbol{\mid} 0; 0)$ and $\mathcal{F}_{\boldsymbol{\lambda}} (\textbf{x}; \infty \boldsymbol{\mid} 0; 0)$ recovers the Cauchy identity for the LLT polynomials. Similarly, by \Cref{suml2} below, \eqref{2gl} and \eqref{1fl} together imply pairing $\mathcal{G}_{\boldsymbol{\lambda}} (0; \textbf{y} \boldsymbol{\mid} 0; 0)$ and $\mathcal{F}_{\boldsymbol{\lambda}} (\textbf{x}; \infty \boldsymbol{\mid} 0; 0)$ recovers the dual LLT Cauchy identity. We will further show in \Cref{HFunction0} below that the degeneration of $\mathcal{H}_{\boldsymbol{\lambda} / \boldsymbol{\mu}}$ given by the fourth part \eqref{1hl} of \Cref{limitg0} admits a deformation that satisfies a vanishing property similar to those defining various families of interpolation polynomials.

Using the branching property \Cref{fghbranching}, we can quickly reduce \Cref{limitg0} to the following proposition, which verifies it in the case when $N = 1$; we will establish the latter result in \Cref{ProofFGL} below. 

\begin{prop}
	
	\label{limitg0horizontal} 
	
	If $N = 1$, then \Cref{limitg0} holds.
	
\end{prop}

\begin{proof}[Proof of \Cref{limitg0} Assuming \Cref{limitg0horizontal}]
	
	All four parts of \Cref{limitg0} will follow from the facts that the $\mathcal{G}_{\boldsymbol{\lambda} / \boldsymbol{\mu}}$, $\mathcal{F}_{\boldsymbol{\lambda} / \boldsymbol{\mu}}$, or $\mathcal{H}_{\boldsymbol{\lambda} / \boldsymbol{\mu}}$ functions satisfy the same branching relation as do the $\mathcal{L}_{\boldsymbol{\lambda} / \boldsymbol{\mu}}$ ones. We will only describe this in detail for the first two statements of \Cref{limitg0} (assuming their counterparts in \Cref{limitg0horizontal}); the proofs of the remaining two statements there are entirely analogous and are therefore omitted. 
	
	We begin by verifying \eqref{1gl}, assuming it holds when $N = 1$.  To that end, observe by \Cref{lbranching} and \eqref{llambdamu} that if $\textbf{x} = \textbf{x}' \cup \textbf{x}''$ then we have the branching relation
	\begin{flalign}
	\label{identitylbranching} 
	\mathcal{L}_{\boldsymbol{\lambda} / \boldsymbol{\mu}} (\textbf{x}) = \displaystyle\sum_{\boldsymbol{\nu} \in \SeqSign_{n; M}} \mathcal{L}_{\boldsymbol{\lambda} / \boldsymbol{\nu}} (\textbf{x}'') \mathcal{L}_{\boldsymbol{\nu} / \boldsymbol{\mu}} (\textbf{x}').
	\end{flalign}
	
	\noindent Thus, since the same branching relation holds for $\mathcal{G}_{\boldsymbol{\lambda} / \boldsymbol{\mu}} (\textbf{x}; \infty \boldsymbol{\mid} 0; 0)$ by \Cref{fghbranching} and \eqref{limitg}, and since \Cref{limitg0horizontal} implies that $\mathcal{G}_{\boldsymbol{\lambda} / \boldsymbol{\mu}} (x_i; \infty \boldsymbol{\mid} 0; 0) = \mathcal{L}_{\boldsymbol{\lambda} / \boldsymbol{\mu}} (x_i)$ for each $i \in [1, N]$, we deduce \Cref{limitg0} in general by induction on $N$.
	
	The verification of \eqref{1fl} is largely similar, except that we must first introduce a normalized variant of $\mathcal{F}_{\boldsymbol{\lambda} / \boldsymbol{\mu}} (x; \infty \boldsymbol{\mid} 0; 0)$ as suggested by \eqref{1fl} and then verify that its branching coincides with the one indicated by \eqref{identitylbranching}. More specifically, for any sequences of $n$ signatures $\boldsymbol{\lambda} \in \SeqSign_{n; M + N}$ and $\boldsymbol{\mu} \in \SeqSign_{n; M}$, define 
	\begin{flalign}
	\label{f2x1} 
	\widetilde{\mathcal{F}}_{\boldsymbol{\lambda} / \boldsymbol{\mu}} (\textbf{x}; \infty \boldsymbol{\mid} 0; 0) = q^{\psi (\boldsymbol{\lambda}) - \psi (\boldsymbol{\mu}) + \binom{M + N}{2} \binom{n}{2} / 2 - \binom{M}{2} \binom{n}{2} / 2} \mathcal{F}_{\boldsymbol{\lambda} / \boldsymbol{\mu}} (\textbf{x}; \infty \boldsymbol{\mid} 0; 0) \displaystyle\prod_{j = 1}^N x_j^{n (M + N - j)}.
	\end{flalign}
	
	\noindent We must show $\widetilde{\mathcal{F}}_{\boldsymbol{\lambda} / \boldsymbol{\mu}} (\textbf{x}; \infty \boldsymbol{\mid} 0; 0) = \mathcal{L}_{\boldsymbol{\lambda} / \boldsymbol{\mu}} (\textbf{x})$, assuming this holds when $N = 1$. As in the proof of \eqref{1gl}, to do this it will suffice by induction to show that they satisfy the same branching property. 
	
	To that end, fix integers $K, L \ge 1$ such that $N = K + L$. Let $\textbf{x}' = (x_1, x_2, \ldots , x_K)$ and $\textbf{x}'' = (x_{K + 1}, x_{K + 2}, \ldots , x_{K + L})$, so that $\textbf{x} = \textbf{x}' \cup \textbf{x}''$. Then, we find that 
	\begin{flalign*}
	\widetilde{\mathcal{F}}_{\boldsymbol{\lambda} / \boldsymbol{\mu}} (\textbf{x}; \infty \boldsymbol{\mid} 0; 0) & = q^{\psi (\boldsymbol{\lambda}) - \psi (\boldsymbol{\mu}) + \binom{M + N}{2} \binom{n}{2} / 2 - \binom{M}{2} \binom{n}{2} / 2} \displaystyle\prod_{j = 1}^N x_j^{n (M + N - j)} \\
	& \qquad \qquad \times \displaystyle\sum_{\boldsymbol{\nu} \in \SeqSign_{n; M + L}} \mathcal{F}_{\boldsymbol{\lambda} / \boldsymbol{\nu}} (\textbf{x}'; \infty \boldsymbol{\mid} 0; 0) \mathcal{F}_{\boldsymbol{\nu} / \boldsymbol{\mu}} (\textbf{x}''; \infty \boldsymbol{\mid} 0; 0) \\
	& = \displaystyle\prod_{j = 1}^N x_j^{n (M + N - j)} \displaystyle\prod_{j = 1}^K x_j^{n (j - M - N)} \displaystyle\prod_{j = 1}^L x_{j + K}^{n (j - M - L)} \\
	& \qquad \qquad \times \displaystyle\sum_{\boldsymbol{\nu} \in \SeqSign_{n; M + L}} \widetilde{\mathcal{F}}_{\boldsymbol{\lambda} / \boldsymbol{\nu}} (\textbf{x}'; \infty \boldsymbol{\mid} 0; 0) \widetilde{\mathcal{F}}_{\boldsymbol{\nu} / \boldsymbol{\mu}} (\textbf{x}''; \infty \boldsymbol{\mid} 0; 0) \\
	& = \displaystyle\sum_{\boldsymbol{\nu} \in \SeqSign_n} \widetilde{\mathcal{F}}_{\boldsymbol{\lambda} / \boldsymbol{\nu}} (\textbf{x}'; \infty \boldsymbol{\mid} 0; 0) \widetilde{\mathcal{F}}_{\boldsymbol{\nu} / \boldsymbol{\mu}} (\textbf{x}''; \infty \boldsymbol{\mid} 0; 0),
	\end{flalign*}
	
	\noindent where in the first equality we used the definition \eqref{f2x1} for $\widetilde{\mathcal{F}}$ and its branching (which follows from \Cref{fghbranching} and \eqref{limithf}); in the second we used \eqref{f2x1} again; and in the third we used the fact that $N = K + L$.
	
	By \Cref{lbranching}, $\mathcal{L}_{\boldsymbol{\lambda} / \boldsymbol{\mu}} (\textbf{x})$ satisfies the same branching identity, given by \eqref{identitylbranching}. Since by \Cref{limitg0horizontal} we have $\widetilde{\mathcal{F}}_{\boldsymbol{\lambda} / \boldsymbol{\mu}} (\textbf{x}; \infty \boldsymbol{\mid} 0; 0) = \mathcal{L}_{\boldsymbol{\lambda} / \boldsymbol{\mu}} (\textbf{x})$ when $N = 1$, it follows by induction on $N$ that this equality holds in general. 
\end{proof}

\section{Cauchy Identities for LLT Polynomials} 

\label{SumL}

In this section we establish, as consequences of \Cref{limitg0} and the results from \Cref{IdentitiesDegenerations}, two identities for the LLT polynomials that were first derived in \cite{TRA}. The first is a Cauchy identity for the LLT polynomials and the second is a dual Cauchy identity. 

\begin{cor}[{\cite[Theorem 35]{TRA}}]
	
	\label{suml1} 
	
	Fix integers $N, M \ge 1$ and sequences of complex numbers $\textbf{\emph{x}} = (x_1, x_2, \ldots , x_N)$ and $\textbf{\emph{y}} = (y_1, y_2, \ldots , y_M)$. If $|q| < 1$ and $|x_j|, |y_i| < 1$ for each $i \in [1, M]$ and $j \in [1, N]$, then 
	\begin{flalign*}
	\displaystyle\sum_{\boldsymbol{\lambda} \in \SeqSign_{n; N}} \mathcal{L}_{\boldsymbol{\lambda}} (\textbf{\emph{x}}) \mathcal{L}_{\boldsymbol{\lambda}} (\textbf{\emph{y}}) = \displaystyle\prod_{i = 1}^M \displaystyle\prod_{j = 1}^N (x_j y_i; q)_n^{-1}. 
	\end{flalign*}
	
\end{cor}

\begin{proof} 
	
	Applying \Cref{sxzlimitidentity}, with the $(\textbf{u}, \textbf{w})$ there equal to $(q^{1 - n} \textbf{x}^{-1}, \textbf{y})$ here, we deduce that 
	\begin{flalign}
	\label{sumfgl}
	\displaystyle\sum_{\boldsymbol{\lambda} \in \SeqSign_{n; N}} \mathcal{F}_{\boldsymbol{\lambda}} (q^{1 - n} \textbf{x}^{-1}; \infty \boldsymbol{\mid} 0; 0) \mathcal{G}_{\boldsymbol{\lambda}} (\textbf{y}; \infty \boldsymbol{\mid} 0; 0) = q^{\binom{n}{2} \binom{N}{2}} \displaystyle\prod_{j = 1}^N x_j^{n (N - j)} \displaystyle\prod_{i = 1}^M \displaystyle\prod_{j = 1}^N (q^{n - 1} x_j y_i; q^{-1})_n^{-1},
	\end{flalign} 
	
	\noindent where \eqref{x1x2bdestimate6} holds since the bounds $|q| < 1$ and $\binom{a}{2} - \binom{n}{2} + (n - 1) (n - a) \ge 0$ for $a \le n$ give
	\begin{flalign*}
	\displaystyle\max_{\substack{1 \le i \le M \\ 1 \le j \le N}} \displaystyle\max_{\substack{a, b \in [0, n] \\ (a, b) \ne (n, 0)}} \big| q^{\binom{a}{2} + \binom{b}{2} - \binom{n}{2}} q^{(n - 1) (n - a)} x_j^{n - a} y_i^b \big| \le \displaystyle\max_{\substack{1 \le i \le M \\ 1 \le j \le N}} \big\{ |x_j|, |y_i| \big\} < 1.
	\end{flalign*}
	
	\noindent Since \eqref{1gl} and \eqref{1fl} imply that $\mathcal{G}_{\boldsymbol{\lambda}} (\textbf{y}; \infty \boldsymbol{\mid} 0; 0) = q^{\psi (\boldsymbol{\lambda}) - \psi (\boldsymbol{0}^N)}  \mathcal{L}_{\boldsymbol{\lambda}} (\textbf{y})$ and 
	\begin{flalign}
	\label{flambda0}
	\mathcal{F}_{\boldsymbol{\lambda}} (q^{1 - n} \textbf{x}^{-1}; \infty \boldsymbol{\mid} 0; 0) = q^{- \psi (\boldsymbol{\lambda}) - \binom{n}{2} \binom{N}{2} / 2} \mathcal{L} (\textbf{x}) \displaystyle\prod_{j = 1}^N q^{n (n - 1) (N - j)} x_j^{n (N - j)},
	\end{flalign}
	
	\noindent the corollary follows from \eqref{sumfgl}, \eqref{psi0}, and the fact that $(q^{n - 1} z; q^{-1})_n = (z; q)_n$, for any $z \in \mathbb{C}$.
\end{proof}

\begin{cor}[{\cite[Proposition 36]{TRA}}]
	
	\label{suml2}
	
	Fix integers $N, M \ge 1$ and sequences of complex numbers $\textbf{\emph{x}} = (x_1, x_2, \ldots , x_N)$ and $\textbf{\emph{y}} = (y_1, y_2, \ldots , y_M)$. Then,
	\begin{flalign*}
	\displaystyle\sum_{\boldsymbol{\lambda} \in \SeqSign_{n; N}} \mathcal{L}_{\boldsymbol{\lambda}} (\textbf{\emph{x}}; q) \mathcal{L}_{\boldsymbol{\lambda}'} (\textbf{\emph{y}}; q^{-1}) = \displaystyle\prod_{i = 1}^M \displaystyle\prod_{j = 1}^N (- q^{(n - 1) / 2} x_j y_i; q)_n.
	\end{flalign*}
	
\end{cor}

\begin{proof} 
	
	Applying \Cref{sxzlimitidentity3}, with the $(\textbf{u}, \textbf{t})$ there equal to $(q^{1 - n} \textbf{x}^{-1}, -q^{(n - 1) / 2} \textbf{y}^{-1})$ here, we deduce that 
	\begin{flalign}
	\label{sumfgl2}
	\begin{aligned}
	\displaystyle\sum_{\boldsymbol{\lambda} \in \SeqSign_{n; N}} (-1)^{|\boldsymbol{\lambda}|} & \mathcal{F}_{\boldsymbol{\lambda}} (q^{1 - n} \textbf{x}^{-1}; \infty \boldsymbol{\mid} 0; 0) \mathcal{G}_{\boldsymbol{\lambda}} (0; -q^{(n - 1) / 2} \textbf{y}^{-1} \boldsymbol{\mid} 0; 0) \\
	& = q^{\binom{n}{2} \binom{N}{2}} \displaystyle\prod_{j = 1}^N x_j^{n (N - j)} \displaystyle\prod_{i = 1}^M \displaystyle\prod_{j = 1}^N (-q^{(n - 1) / 2} x_j y_i; q^{-1})_n^{-1}.
	\end{aligned}
	\end{flalign} 
	
	\noindent Since \eqref{2gl} and the homogeneity of the LLT polynomials $\mathcal{L}_{\boldsymbol{\lambda}}$ together imply that
	\begin{flalign*}
	(-1)^{|\boldsymbol{\lambda}|} \mathcal{G}_{\boldsymbol{\lambda}} (0; -q^{(n - 1) / 2} \textbf{y}^{-1} \boldsymbol{\mid} 0; 0) & = (-1)^{|\boldsymbol{\lambda}|} q^{\psi  (\boldsymbol{\lambda}) - \psi (\boldsymbol{0}^N)} \mathcal{L}_{\boldsymbol{\lambda}'} (-\textbf{y}; q^{-1}) = q^{\psi (\boldsymbol{\lambda}) - \psi (\boldsymbol{0}^N)} \mathcal{L}_{\boldsymbol{\lambda}'} (\textbf{y}; q^{-1}),
	\end{flalign*}
	
	\noindent the corollary follows from \eqref{sumfgl2}, \eqref{flambda0}, and \eqref{psi0}.	
\end{proof}

\section{Plethystic Substitution of LLT Polynomials}

\label{Polynomialsr}

In this section we explain how the $\mathcal{G}_{\boldsymbol{\lambda} / \boldsymbol{\mu}} \big( (x_1, x_2, \ldots , x_N), r \boldsymbol{\mid} 0; 0)$ functions from \eqref{limitg} can be interpreted as images of the LLT polynomials under a certain plethystic substitution. In what follows, we fix an infinite set of variables $\textbf{x} = (x_1, x_2, \ldots )$; let $\Lambda (\textbf{x})$\index{0@$\Lambda (\textbf{x})$; ring of symmetric functions in $\textbf{x}$} denote the ring of symmetric functions in $\textbf{x}$; and set $\textbf{x}_{[1, N]} = (x_1, x_2, \ldots , x_N)$ for any integer $N \ge 1$. 

As defined in \eqref{limitg}, the $\mathcal{G}_{\boldsymbol{\lambda} / \boldsymbol{\mu}}$ are functions of finitely many variables and are therefore not directly elements of $\Lambda (\textbf{x})$. However, the following lemma shows that the sequence of functions $\big\{ \mathcal{G}_{\boldsymbol{\lambda} / \boldsymbol{\mu}} (\textbf{x}_{[1, N]}, r; \infty, 0) \big\}_{N \ge 1}$ satisfies the compatibility condition that makes it an inverse system in $\Lambda (\textbf{x})$. 

\begin{lem} 
	
\label{g0}	

Fix integers $M \ge 0$ and $N \ge 1$; signature sequences $\boldsymbol{\lambda}, \boldsymbol{\mu} \in \SeqSign_{n; M}$; and a complex number $r \in \mathbb{C}$. For any integer $N \ge 1$, the function $\mathcal{G}_{\boldsymbol{\lambda} / \boldsymbol{\mu}} (\textbf{\emph{x}}_{[1, N]}; r \boldsymbol{\mid} 0; 0)$ from \eqref{limitg} is a symmetric, homogeneous polynomial in $\textbf{\emph{x}}_{[1, N]}$ of degree $|\boldsymbol{\lambda}| - |\boldsymbol{\mu}|$. Moreover, if $x_{N + 1} = 0$, then $\mathcal{G}_{\boldsymbol{\lambda} / \boldsymbol{\mu}} \big( \textbf{\emph{x}}_{[1, N + 1]}; r \boldsymbol{\mid} 0; 0 \big) = \mathcal{G}_{\boldsymbol{\lambda} / \boldsymbol{\mu}} \big( \textbf{\emph{x}}_{[1, N]}; r \boldsymbol{\mid} 0, 0 \big)$.

\end{lem} 

\begin{proof}
	
	The symmetry of $\mathcal{G}_{\boldsymbol{\lambda} / \boldsymbol{\mu}} (\textbf{x}_{[1, N]}; r \boldsymbol{\mid} 0, 0)$ follows from \Cref{gxfxsigma}. To verify that it is homogenous and of degree $|\boldsymbol{\lambda}| - |\boldsymbol{\mu}|$, recall from \eqref{gfhe} and \Cref{weightesum} that $\mathcal{G}_{\boldsymbol{\lambda} / \boldsymbol{\mu}} \big( \textbf{x}_{[1, N]}; r \boldsymbol{\mid} 0; 0 \big)$ is the partition function under the weights $\mathcal{W}_x (\textbf{A}, \textbf{B}; \textbf{C}, \textbf{D} \boldsymbol{\mid} r)$ (from \eqref{limitw}) for the vertex model $\mathfrak{P}_G (\boldsymbol{\lambda} / \boldsymbol{\mu}; N)$ (from \Cref{pgpfph}; see also the left side of \Cref{fgpaths}) on the domain $\mathcal{D}_N = \mathbb{Z}_{> 0} \times \{ 1, 2, \ldots , N \}$. The explicit form \eqref{limitw} for these weights implies that they only depend on $x$ through a factor of $x^{|\textbf{D}|}$. Since for each path ensemble $\mathcal{E} \in \mathfrak{P}_G (\boldsymbol{\lambda} / \boldsymbol{\mu}; N)$, under which the arrow configuration at any vertex $v \in \mathcal{D}_N$ is denoted by $\big( \textbf{A} (v), \textbf{B} (v); \textbf{C} (v), \textbf{D} (v) \big)$, we have $\sum_{v \in \mathcal{D}_N} \big| \textbf{D} (v) \big| = |\boldsymbol{\lambda}| - |\boldsymbol{\mu}|$, it follows that $\mathcal{G}_{\boldsymbol{\lambda} / \boldsymbol{\mu}} (\textbf{x}_{[1, N]}; r \boldsymbol{\mid} 0; 0)$ is a homogeneous polynomial in $\textbf{x}_{[1, N]}$ of degree $|\boldsymbol{\lambda}| - |\boldsymbol{\mu}|$. This verifies the first statement of the lemma.
	
	Next let us verify the second statement of the lemma. Once again, $\mathcal{G}_{\boldsymbol{\lambda} / \boldsymbol{\mu}} \big( \textbf{x}_{[1, N + 1]}; r \boldsymbol{\mid} 0; 0 \big)$ is the partition function under the weights $\mathcal{W}_x (\textbf{A}, \textbf{B}; \textbf{C}, \textbf{D} \boldsymbol{\mid} r)$ for $\mathfrak{P}_G (\boldsymbol{\lambda} / \boldsymbol{\mu}; N + 1)$ on the domain $\mathcal{D}_{N + 1} = \mathbb{Z}_{> 0} \times \{ 1, 2, \ldots , N + 1 \}$. Since spectral parameter $x_{N + 1}$ for this top row in this model is equal to $0$, and since \eqref{limitw} implies that $\mathcal{W}_0 (\textbf{A}, \textbf{B}; \textbf{C}, \textbf{D} \boldsymbol{\mid} r) = 0$ unless $|\textbf{D}| = 0$, the arrow configurations $\big( \textbf{A} (v), \textbf{B} (v); \textbf{C} (v), \textbf{D} (v) \big)$ under any path ensemble supported by this vertex model must satisfy $\textbf{D} (v) = \textbf{e}_0$ whenever $v \in \mathbb{Z}_{> 0} \times \{ N + 1 \}$. In particular, paths in any such ensemble cannot move to the right in this top row, meaning that for any $v \in \mathbb{Z}_{> 0} \times \{ N + 1 \}$ we must have that $\big( \textbf{A} (v), \textbf{B} (v); \textbf{C} (v), \textbf{D} (v) \big)$ is of the form $(\textbf{A}, \textbf{e}_0; \textbf{A}, \textbf{e}_0)$, for some $\textbf{A} \in \{ 0, 1 \}^n$. By the last statement of \eqref{limitw}, $\mathcal{W}_x (\textbf{A}, \textbf{e}_0; \textbf{A}, \textbf{e}_0 \boldsymbol{\mid} r) = 1$ for any $\textbf{A}$, meaning that the weight of the topmost row of any ensemble $\mathcal{E}$ supported by $\mathcal{G}_{\boldsymbol{\lambda} / \boldsymbol{\mu}} \big( \textbf{x}_{[1, N + 1]}; r \boldsymbol{\mid} 0; 0 \big)$ is equal to $1$. 
	
	Hence, we can omit this row from model, deducing from \Cref{weightesum} that
	\begin{flalign*}
	\mathcal{G}_{\boldsymbol{\lambda} / \boldsymbol{\mu}} \big( \textbf{x}_{[1, N + 1]}; r \boldsymbol{\mid} 0; 0 \big) & = \displaystyle\sum_{\mathcal{E} \in \mathfrak{P}_G (\boldsymbol{\lambda} / \boldsymbol{\mu}; N + 1)} \mathcal{W} (\mathcal{E} \boldsymbol{\mid} \textbf{x}; r \boldsymbol{\mid} \infty, 0) \\
	& = \displaystyle\sum_{\mathcal{E} \in \mathfrak{P}_G (\boldsymbol{\lambda} / \boldsymbol{\mu}; N)} \mathcal{W} (\mathcal{E} \boldsymbol{\mid} \textbf{x}; r \boldsymbol{\mid} \infty, 0) = \mathcal{G}_{\boldsymbol{\lambda} / \boldsymbol{\mu}} \big( \textbf{x}_{[1, N]}; r \boldsymbol{\mid} 0; 0 \big),
	\end{flalign*}
	
	\noindent and thereby establishing the lemma.
\end{proof}

\begin{definition}
	
	\label{r0g} 
	
	For any $r \in \mathbb{C}$ and $\boldsymbol{\lambda}, \boldsymbol{\mu} \in \SeqSign_n$, let $\mathcal{G}_{\boldsymbol{\lambda} / \boldsymbol{\mu}} (\textbf{x}; r \boldsymbol{\mid} 0; 0) \in \Lambda (\textbf{x})$ denote the element of $\Lambda (\textbf{x})$ associated with the inverse system of symmetric, degree $|\boldsymbol{\lambda}| - |\boldsymbol{\mu}|$ polynomials $\big\{ \mathcal{G}_{\boldsymbol{\lambda} / \boldsymbol{\mu}} (\textbf{x}; r \boldsymbol{\mid} 0; 0) \big\}_{N \ge 1}$. 
\end{definition}

\begin{rem}

\label{lg0} 

Observe from the explicit forms \eqref{limitw} of the $\mathcal{W}$ weights (or, alternatively, as a consequence of \Cref{limitsy}) that, for any $\textbf{A}, \textbf{B}, \textbf{C}, \textbf{D} \in \{ 0, 1 \}^n$, we have
\begin{flalign*}
\displaystyle\lim_{r \rightarrow \infty} \mathcal{W}_x (\textbf{A}, \textbf{B}; \textbf{C}, \textbf{D} \boldsymbol{\mid} r) = \mathcal{W}_x (\textbf{A}, \textbf{B}; \textbf{C}, \textbf{D} \boldsymbol{\mid} \infty, 0). 
\end{flalign*}

\noindent Thus, \eqref{gfhe} and \Cref{weightesum} together imply for any $\boldsymbol{\lambda}, \boldsymbol{\mu} \in \SeqSign_n$ that
\begin{flalign}
\label{glrlimit} 
\displaystyle\lim_{r \rightarrow \infty} \mathcal{G}_{\boldsymbol{\lambda} / \boldsymbol{\mu}} (\textbf{x}; r \boldsymbol{\mid} 0; 0) = \mathcal{G}_{\boldsymbol{\lambda} / \boldsymbol{\mu}} (\textbf{x}; \infty \boldsymbol{\mid} 0; 0).
\end{flalign}

\noindent This, together with \eqref{1gl}, implies that we may view LLT polynomial $\mathcal{L}_{\boldsymbol{\lambda} / \boldsymbol{\mu}} (\textbf{x})$ as an element of $\Lambda (\textbf{x})$.  

\end{rem} 

Next, we briefly recall the notion of plethystic substituion. Letting $p_k (\textbf{x}) \in \Lambda (\textbf{x})$ denote the power sum $p_k (\textbf{x}) = \sum_{i = 1}^{\infty} x_i^k$ for each integer $k \ge 1$, the ring $\Lambda (\textbf{x})$ is freely generated as an algebra by the $\big\{ p_k (\textbf{x}) \big\}$. For any symmetric function $f \in \Lambda (\textbf{x})$ and formal sum of monomials $A = \sum_{i = 1}^{\infty} a_i$, define the \emph{plethysm} $f[A]$ to be the image of $f (\textbf{x})$ under the substitution $p_k (\textbf{x}) = \sum_{i = 1}^{\infty} a_i^k$. Let us mention that, under this notation, we will view additional parameters such as $q$ and $r$ below as variables. In particular, if for instance $B = (1 - r^{-2}) A$, then 
\begin{flalign*} 
p_k [B] = \displaystyle\sum_{i = 1}^{\infty} a_i^k - \displaystyle\sum_{i = 1}^{\infty} r^{-2k} a_i^k = \displaystyle\sum_{i = 1}^{\infty} (1 - r^{-2k}) a_i^k.
\end{flalign*}

Now we can state the following consequence of \Cref{limitg0} that expresses the function $\mathcal{G}_{\boldsymbol{\lambda} / \boldsymbol{\mu}} (\textbf{x}; r \boldsymbol{\mid} 0; 0) \in \Lambda (\textbf{x})$ from \Cref{r0g} as a certain plethystic substitution of the LLT polynomial $\mathcal{L}_{\boldsymbol{\lambda} / \boldsymbol{\mu}} (\textbf{x}) \in \Lambda (\textbf{x})$ from \eqref{lambdamug} (and \eqref{llambdamu}). In what follows, we recall the function $\psi$ from \eqref{lambdamupsi}. 

\begin{prop}
	
	\label{gl} 
	
	Fix an integer $M \ge 0$ and signature sequences $\boldsymbol{\lambda}, \boldsymbol{\mu} \in \SeqSign_{n; M}$; further let $r$ and $q$ denote variables. Defining the formal sum $X = \sum_{i = 1}^{\infty} x_i$, we have the plethystic identity
	\begin{flalign}
	\label{lgr} 
	\mathcal{L}_{\boldsymbol{\lambda} / \boldsymbol{\mu}} \big[ (1 - r^{-2}) X \big] = q^{\psi(\boldsymbol{\mu}) - \psi(\boldsymbol{\lambda})} \mathcal{G}_{\boldsymbol{\lambda} / \boldsymbol{\mu}} (\textbf{\emph{x}}; r \boldsymbol{\mid} 0; 0).
	\end{flalign} 
	
\end{prop}

\begin{proof}
	
	For the purposes of this proof, define the symmetric function $Y_{\boldsymbol{\lambda} / \boldsymbol{\mu}} (\textbf{x}) \in \Lambda (\textbf{x})$ by setting
	\begin{flalign}
	\label{y1} 
	Y_{\boldsymbol{\lambda} / \boldsymbol{\mu}} (\textbf{x}) = \mathcal{G}_{\boldsymbol{\lambda} / \boldsymbol{\mu}} (\textbf{x}; q^{-1 / 2} \boldsymbol{\mid} 0; 0). 
	\end{flalign}
	
	\noindent We will establish the plethystic identity
	\begin{flalign}
	\label{yg} 
	Y_{\boldsymbol{\lambda} / \boldsymbol{\mu}} \big[ (1 - r^{-2}) (1 - q)^{-1} X \big] = \mathcal{G}_{\boldsymbol{\lambda} / \boldsymbol{\mu}} (\textbf{x}; r \boldsymbol{\mid} 0; 0). 
	\end{flalign}
	
	Let us first show the proposition assuming \eqref{yg}. To that end, observe by specializing $r = \infty$ in \eqref{yg}, and applying \eqref{glrlimit} and \eqref{1gl}, yields
	\begin{flalign*}
	Y_{\boldsymbol{\lambda} / \boldsymbol{\mu}} \big[ (1 - q)^{-1} X \big] = \mathcal{G}_{\boldsymbol{\lambda} / \boldsymbol{\mu}} (\textbf{x}; \infty \boldsymbol{\mid} 0; 0) = q^{\psi (\boldsymbol{\lambda}) - \psi(\boldsymbol{\mu})} \mathcal{L}_{\boldsymbol{\lambda} / \boldsymbol{\mu}} (X).
	\end{flalign*}	
	
	\noindent Applying this with $X$ replaced by $(1 - r^{-2}) X$, together with \eqref{yg}, implies \eqref{lgr}. 
	
	So, it remains to verify \eqref{yg}. To that end, observe that from the explicit form \eqref{limitw} of the weights $\mathcal{W}_x (\textbf{A}, \textbf{B}; \textbf{C}, \textbf{D} \boldsymbol{\mid} r)$ that they are rational functions of $r$. Hence, by \eqref{gfhe} and \Cref{weightesum}, the coefficients of $\mathcal{G}_{\boldsymbol{\lambda} / \boldsymbol{\mu}} (\textbf{x}; r \boldsymbol{\mid} 0; 0)$ when expanded in the monomial basis of $\Lambda (\textbf{x})$ are also rational functions of $r$. Hence, it suffices to show that \eqref{yg} holds for infinitely many values of $r$. In particular, we may assume that $r = q^{-L / 2}$ for some integer $L \ge 1$. 
	
	In this case, $(1 - r^{-2}) (1 - q)^{-1} = (1 - q^L) (1 - q)^{-1} = \sum_{j = 0}^{L - 1} q^j$. In particular, the plethysm mapping $X$ to $(1 - r^{-2}) (1 - q)^{-1} X$ is realized by replacing the variable sequence $\textbf{x}$ with $\textbf{w} = \textbf{w}^{(1)} \cup \textbf{w}^{(2)} \cup \cdots$, where we have denoted the principal specialization $\textbf{w}^{(j)} = (x_j, q x_j, \ldots , q^{L - 1} x_j)$ for each $j \ge 1$. Thus, in view of the definition \eqref{y1} of $Y_{\boldsymbol{\lambda} / \boldsymbol{\mu}}$, it suffices to show that 
	\begin{flalign*}
	\mathcal{G}_{\boldsymbol{\lambda} / \boldsymbol{\mu}} (\textbf{w}; q^{-1 / 2} \boldsymbol{\mid} 0; 0) = \mathcal{G}_{\boldsymbol{\lambda} / \boldsymbol{\mu}} (\textbf{x}; q^{-L / 2} \boldsymbol{\mid} 0; 0).
	\end{flalign*}
	
	\noindent If there exists some integer $N \ge 0$ such that $x_j = 0$ for $j > N$, then this follows from \Cref{gq}, from which we deduce the general case by letting $N$ tend to $\infty$.
\end{proof}

\begin{rem} 
	
	\label{yfunction} 

	Plethysms of the LLT polynomials have been considered and studied previously in the algebraic combinatorics literature. In particular, in the case when each skew-shape $\lambda^{(i)} / \mu^{(i)}$ in the sequence $\boldsymbol{\lambda} / \boldsymbol{\mu}$ has size $1$, the functions $Y_{\boldsymbol{\lambda} / \boldsymbol{\mu}} (\textbf{x})$ from \eqref{y1} coincide with chromatic quasisymmetric functions \cite{CSFV,CSF} associated with the incomparability graph of a unit interval order; see Remark 3.6 of \cite{S} and Section 5.4 of \cite{PCQF}. 
		
\end{rem}

\section{Stability}

\label{StabilityPolynomials}

In this section we explain how \Cref{limitg0} can be used to recover a \emph{stability} property for the LLT polynomials, which was originally established in \cite{RT}. This property essentially states that the LLT polynomial $\mathcal{L}_{n \lambda}^{(n)} (\textbf{x})$ from \eqref{branchingl} becomes independent of $n$ once $n \ge \ell (\lambda)$, at which point it coincides with a modified Hall--Littlewood polynomial. To define the latter, for any partition $\lambda = (\lambda_1, \lambda_2, \ldots , \lambda_M)$ and sequence of complex variables $\textbf{x} = (x_1, x_2, \ldots,  x_N)$ (possibly with $N = \infty$), recall  the Hall--Littlewood function $Q_{\lambda} (\textbf{x}) = Q_{\lambda} (\textbf{x}; q)$ from Section 3.2 of \cite{SFP}.\index{Q@$Q_{\lambda} (\textbf{x})$; Hall--Littlewood polynomial} Then, (following, for instance, Definition 1.1 of \cite{NCFMP}) set the \emph{modified Hall--Littlewood polynomial} $Q_{\lambda}' (\textbf{x})$\index{Q@$Q_{\lambda} (\textbf{x})$; Hall--Littlewood polynomial!$Q_{\lambda}' (\textbf{x})$; modified Hall--Littlewood polynomial} to be 
\begin{flalign*}
Q_{\lambda}' (\textbf{x}) = Q_{\lambda} \Bigg( \bigcup_{j = 0}^{\infty} q^j \textbf{x} \Bigg).
\end{flalign*}

\noindent Equivalently, it is given by the plethystic substitution $Q_{\lambda}' [X] = Q_{\lambda} \big[ (1 - q)^{-1} X \big]$, where $X$ denotes the formal sum $X = \sum_{i = 1}^{\infty} x_i$. 

Under this notation, we can state the following stability result, first shown in \cite{RT}.

\begin{prop}[{\cite[Theorem 6.6]{RT}}]

\label{qlambdaqlambda} 

Fix integers $N \ge 1$ and $n \ge M \ge 1$; a set of complex numbers $\textbf{\emph{x}} = (x_1, x_2, \ldots , x_N)$; and a partition $\lambda = (\lambda_1, \lambda_2, \ldots , \lambda_M)$ of length $M$. Define the signature $\widetilde{\lambda} = (\widetilde{\lambda}_1, \widetilde{\lambda}_2, \ldots , \widetilde{\lambda}_n) \in \Sign_n$ such that $\widetilde{\lambda}_i = n \lambda_i$ for each $i \in [1, M]$ and $\widetilde{\lambda}_i = 0$ for $i \in [M + 1, n]$. Then, $\mathcal{L}_{\widetilde{\lambda}}^{(n)} (\textbf{\emph{x}}) = Q_{\lambda}' (\textbf{\emph{x}})$. 

\end{prop}

\begin{rem}
	
	\label{fnf1} 
	
	Through similar methods, one can also establish variants of \Cref{qlambdaqlambda} for the more general functions $G_{\boldsymbol{\lambda} / \boldsymbol{\mu}} (\textbf{x}; \textbf{r} \boldsymbol{\mid} \textbf{y}; \textbf{s})$ from \Cref{fgdefinition} (and not only their degenerations as LLT polnomials). Here, the stable limits will instead be fully fused analogs of the $U_q \big( \widehat{\mathfrak{sl}} (2) \big)$ vertex model partition functions $\mathsf{G}_{\lambda / \mu}$ introduced in Definition 4.3 of \cite{HSVMSRF}. For the sake of brevity, we will not pursue this further here. 
	
\end{rem}

To establish \Cref{qlambdaqlambda}, we require notation for partition functions similar to those given in \Cref{zef1rsz} (see \Cref{domainpathsfused}), under the weights $\mathcal{W}_z (\textbf{A}, \textbf{B}; \textbf{C}, \textbf{D} \boldsymbol{\mid} \infty, 0)$ from \eqref{limitw}; in what follows, we recall the notions of fused path ensembles, east-south domains,\footnote{Here, we will allow our east-south domains to be infinite, in which scenario all results and notation from \Cref{ColorsFused} continue to apply.} and boundary data from \Cref{ColorsFused}. Let $\mathcal{D} \subset \mathbb{Z}^2$ denote some east-south domain, and fix some boundary data $(\mathscr{E}; \mathscr{F})$ on $\mathcal{D}$. Letting $\textbf{z} = \big( z(v) \big)_{v \in \mathcal{D}}$ denote a set of complex numbers, define the partition function
\begin{flalign*}
\mathcal{W}_{\mathcal{D}} (\mathscr{E}; \mathscr{F} \boldsymbol{\mid} \textbf{z} \boldsymbol{\mid} \infty, 0) = \displaystyle\sum \displaystyle\prod_{v \in \mathcal{D}} \mathcal{W}_{z(v)} \big( \textbf{A} (v), \textbf{B} (v); \textbf{C} (v), \textbf{D} (v) \boldsymbol{\mid} \infty, 0 \big),
\end{flalign*}

\noindent where the sum is over all ensembles on $\mathcal{D}$ with boundary data $(\mathscr{E}; \mathscr{F})$. 

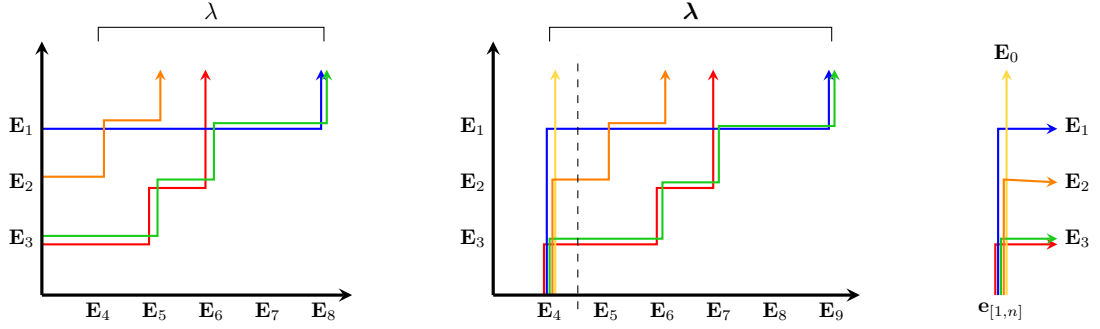
\begin{figure}
	
	\begin{center}

		\begin{tikzpicture}[
		>=stealth,
		scale = .75
		]
		
		\draw[->, thick, red]  (1, .9) -- (2.9, .9) -- (2.9, 1.9) -- (3.9, 1.9) -- (3.9, 4);
		\draw[->, thick, blue]  (1, 2.95) -- (1.95, 2.95) -- (5.95, 2.95) -- (5.95, 4);
		\draw[->, thick, green]  (1, 1.05) -- (3.05, 1.05) -- (3.05, 2.05) -- (4.05, 2.05) -- (4.05, 3.05) -- (6.05, 3.05) -- (6.05, 4);
		\draw[->, thick, orange] (1, 2.1) -- (2.1, 2.1) -- (2.1, 3.1) -- (3.1, 3.1) -- (3.1, 4);
	
		\draw[-] (2, 4.5) -- (2, 4.75) -- (6, 4.75) -- (6, 4.5);
		
		\draw[] (4, 4.75) circle[radius = 0] node[above, scale = .9]{$\lambda$};
		
		\draw[->, very thick] (1, 0) -- (1, 4.5);
		\draw[->, very thick] (1, 0) -- (6.5, 0);

		\draw[->, thick, red] (9.9, 0) -- (9.9, .9) -- (11.9, .9) -- (11.9, 1.9) -- (12.9, 1.9) -- (12.9, 4);
		\draw[->, thick, blue] (9.95, 0) -- (9.95, 2.95) -- (10.95, 2.95) -- (14.95, 2.95) -- (14.95, 4);
		\draw[->, thick, yellow] (10.1, 0) -- (10.1, 4);
		\draw[->, thick, green] (10, 0) -- (10, 1) -- (12, 1) -- (12, 2) -- (13, 2) -- (13, 3) -- (15.05, 3) -- (15.05, 4);
		\draw[->, thick, orange] (10.05, 0) -- (10.05, 2.05) -- (11.05, 2.05) -- (11.05, 3.05) -- (12.05, 3.05) -- (12.05, 4);

		\draw[-] (10, 4.5) -- (10, 4.75) -- (15, 4.75) -- (15, 4.5);
		
		\draw[] (12.5, 4.75) circle[radius = 0] node[above, scale = .9]{$\boldsymbol{\lambda}$};
		
		\draw[->, very thick] (9, 0) -- (9, 4.5);
		\draw[->, very thick] (9, 0) -- (15.5, 0);
		
		\draw[] (1, 1) circle[radius = 0] node[left = 0, scale = .8]{$\textbf{E}_3$};
		\draw[] (1, 2) circle[radius = 0] node[left = 0, scale = .8]{$\textbf{E}_2$};
		\draw[] (1, 3) circle[radius = 0] node[left = 0, scale = .8]{$\textbf{E}_1$};
		
		\draw[] (2, 0) circle[radius = 0] node[below = 0, scale = .8]{$\textbf{E}_4$};
		\draw[] (3, 0) circle[radius = 0] node[below = 0, scale = .8]{$\textbf{E}_5$};
		\draw[] (4, 0) circle[radius = 0] node[below = 0, scale = .8]{$\textbf{E}_6$};
		\draw[] (5, 0) circle[radius = 0] node[below = 0, scale = .8]{$\textbf{E}_7$};
		\draw[] (6, 0) circle[radius = 0] node[below = 0, scale = .8]{$\textbf{E}_8$};

		\draw[] (9, 1) circle[radius = 0] node[left = 0, scale = .8]{$\textbf{E}_3$};
		\draw[] (9, 2) circle[radius = 0] node[left = 0, scale = .8]{$\textbf{E}_2$};
		\draw[] (9, 3) circle[radius = 0] node[left = 0, scale = .8]{$\textbf{E}_1$};
		
		\draw[] (10, 0) circle[radius = 0] node[below = 0, scale = .8]{$\textbf{E}_4$};
		\draw[] (11, 0) circle[radius = 0] node[below = 0, scale = .8]{$\textbf{E}_5$};
		\draw[] (12, 0) circle[radius = 0] node[below = 0, scale = .8]{$\textbf{E}_6$};
		\draw[] (13, 0) circle[radius = 0] node[below = 0, scale = .8]{$\textbf{E}_7$};
		\draw[] (14, 0) circle[radius = 0] node[below = 0, scale = .8]{$\textbf{E}_8$};
		\draw[] (15, 0) circle[radius = 0] node[below = 0, scale = .8]{$\textbf{E}_9$};
		
		\draw[dashed] (10.5, -.25) -- (10.5, 4.25);

		\draw[->, thick, red] (17.9, 0) -- (17.9, .9) -- (19, .9);
		\draw[->, thick, blue] (17.95, 0) -- (17.95, 2.95) -- (19, 2.95);
		\draw[->, thick, yellow] (18.1, 0) -- (18.1, 4) node[above = 0, scale = .8, black]{$\textbf{E}_0$};
		\draw[->, thick, green] (18, 0) -- (18, 1) -- (19, 1);
		\draw[->, thick, orange] (18.05, 0) -- (18.05, 2.05) -- (19, 2.0);
		
		\draw[] (19, 1) circle [radius = 0] node[right = 0, scale = .8]{$\textbf{E}_3$};
		\draw[] (19, 2) circle [radius = 0] node[right = 0, scale = .8]{$\textbf{E}_2$};
		\draw[] (19, 3) circle [radius = 0] node[right = 0, scale = .8]{$\textbf{E}_1$};
		
		\draw[] (18, 0) circle [radius = 0] node[below = 0, scale = .8]{$\textbf{e}_{[1, n]}$};
		
		\end{tikzpicture}
		
	\end{center}
	
	\caption{\label{sumw} To the left is a path ensemble on $\mathcal{D}_3$ with boundary data $\big( \mathscr{E}; \mathscr{I} (\lambda) \big)$, where $\lambda = (3, 5, 5, 2)$ and $\textbf{E}_i = \textbf{e}_0$ for $i > 3$.  In the middle is a path ensemble counted by $\mathcal{G}_{\boldsymbol{\lambda}}$, for $\boldsymbol{\lambda} = \big( (3), (5), (5), (2), (0) \big)$, whose leftmost column is depicted to the right. Here red, blue, green, orange, and yellow are colors $1$, $2$, $3$, $4$, and $5$, respectively.} 
	
\end{figure}

Certain entrance and exit data on the east-south domain $\mathcal{D}_N = \mathbb{Z}_{> 0} \times \{ 1, 2, \ldots , N \} \subset \mathbb{Z}_{\ge 0}^2$ will be of particular use to us. First, any entrance data $\mathscr{E} = (\textbf{E}_1, \textbf{E}_2, \ldots )$ on $\mathcal{D}_N$ such that $\textbf{E}_k = \textbf{e}_0$ for $k > N$ corresponds to prohibiting paths from entering $\mathcal{D}_N$ along its bottom boundary $\{ y = 1 \}$. Second, for any (not necessarily non-increasing) sequence $\lambda = (\lambda_1, \lambda_2, \ldots , \lambda_{\ell}) \in \mathbb{Z}_{> 0}^{\ell}$ of positive integers of length $\ell \le n$, define the sequence\footnote{The indexing here, which starts at $1$, is slightly different from that in \Cref{0PolynomialsP}, which started at $0$.} $\mathscr{I} (\lambda) = \big( \textbf{I}_1 (\lambda), \textbf{I}_2 (\lambda), \ldots \big)$\index{I@$\mathscr{I} (\mu)$} of elements in $\{ 0, 1 \}^n$, such that the $j$-th entry of $\textbf{I}_k (\lambda) \in \{ 0, 1 \}^n$ equals $\textbf{1}_{\lambda_j = k}$, for each $j \in [1, \ell]$ and $k > 0$ (and this entry equals $0$ for each $j \in [\ell + 1, n]$). In this way, $\mathscr{I} (\lambda)$ prescribes exit data on $\mathcal{D}_N$, under which one arrow of color $j$ vertically exits through $(\lambda_j, N)$, for each $j \in [1, \ell]$. In particular, any path ensemble with such exit data has at most one path of any given color; observe that this notation does not specify the order in which the colored paths enter through the left boundary of $\mathcal{D}_N$. We refer to the left side of \Cref{sumw} for a depiction of an ensemble under such entrance and exit data.

Now we have the following lemma expressing the LLT polynomial $\mathcal{L}_{\widetilde{\lambda}}^{(n)} (\textbf{x})$ as linear combinations of the partition functions $\mathcal{W}_{\mathcal{D}} (\mathscr{E}; \mathscr{F} \boldsymbol{\mid} \textbf{z} \boldsymbol{\mid} \infty, 0)$, for $(\mathscr{E}; \mathscr{F})$ of the above form. In the below, we set $\inv (\mathscr{E}) = \inv \big( \mathscr{E}; [1, n] \big)$, where we recall the latter from \eqref{sumij}.\index{I@$\inv$}

\begin{lem}
	
	\label{llambdasum}  
	
	Adopt the notation of \Cref{qlambdaqlambda}, and define the variable set $\textbf{\emph{z}} = \big( z(v) \big)_{v \in \mathcal{D}_N}$ so that $z(v) = z (i, j) = x_j$ for each $v = (i, j) \in \mathcal{D}_N$. Then, 
	\begin{flalign*}
	\mathcal{L}_{\widetilde{\lambda}}^{(n)} (\textbf{\emph{x}}) = \displaystyle\sum_{\mathscr{E}} q^{\inv (\mathscr{E})} \mathcal{W}_{\mathcal{D}_N} \big( \mathscr{E}; \mathscr{I} (\lambda) \boldsymbol{\mid} \textbf{\emph{z}} \boldsymbol{\mid} \infty, 0 \big) \displaystyle\prod_{j = 1}^N q^{\binom{|\textbf{\emph{E}}_j|}{2}} x_{N - j + 1}^{|\textbf{\emph{E}}_j|},
	\end{flalign*}
	
	\noindent where we sum over all entrance data $\mathscr{E} = (\textbf{\emph{E}}_1, \textbf{\emph{E}}_2, \ldots )$ on $\mathcal{D}_N$, with $\textbf{\emph{E}}_k = \textbf{\emph{e}}_0$ for $k > N$. 
\end{lem}

\begin{proof}
	
	Observe that the $n$-quotient $\widetilde{\boldsymbol{\lambda}} = \big( \widetilde{\lambda}^{(1)}, \widetilde{\lambda}^{(2)}, \ldots , \widetilde{\lambda}^{(n)} \big)$ of $\widetilde{\lambda}$ is given by setting $\widetilde{\lambda}^{(i)} = (\lambda_{n - i + 1}) \in \Sign_1$ for each $i \in [n - M + 1, n]$ and $\widetilde{\lambda}^{(i)} = (0)$ for $i \in [1, n - M]$. Indeed, recalling $\mathfrak{T}$ from \eqref{t}, we have 
	\begin{flalign*} 
	\mathfrak{T} (\widetilde{\lambda}) = (n \lambda_1 + n, n \lambda_2 + n - 1, \ldots , n \lambda_M + n - M + 1, n - M, n - M - 1, \ldots , 1).
	\end{flalign*}  
	
	\noindent Hence the content in \Cref{FunctionGPartitions} implies for each $i \in [n - M + 1, n]$ that $\mathfrak{T} \big( \widetilde{\lambda}^{(i)} \big) = (\lambda_{n - i + 1} + 1)$, meaning $\widetilde{\lambda}^{(i)} = (\lambda_{n - i + 1})$, and for each $i \in [1, n - M]$ that $\mathfrak{T} \big( \widetilde{\lambda}^{(i)} \big) = (1)$, meaning $\widetilde{\lambda}^{(i)} = (0)$.
	
	Since $\lambda_1 \ge \lambda_2 \ge \cdots \ge \lambda_M$, it follows from \eqref{lambdamupsi} that $\psi (\widetilde{\boldsymbol{\lambda}}) = \psi (\boldsymbol{0}^1) = 0$, meaning by \eqref{1gl} that $\mathcal{G}_{\widetilde{\boldsymbol{\lambda}}} (\textbf{x}; \infty \boldsymbol{\mid} 0; 0) = \mathcal{L}_{\widetilde{\boldsymbol{\lambda}}} (\textbf{x}) = \mathcal{L}_{\widetilde{\lambda}}^{(n)} (\textbf{x})$. Therefore, it suffices to show that 
	\begin{flalign}
	\label{gsumw1}
	\mathcal{G}_{\widetilde{\boldsymbol{\lambda}}} (\textbf{x}; \infty \boldsymbol{\mid} 0; 0) = \displaystyle\sum_{\mathscr{E}} q^{\inv (\mathscr{E})} \mathcal{W}_{\mathcal{D}_N} \big( \mathscr{E}; \mathscr{I} (\lambda) \boldsymbol{\mid} \textbf{z} \boldsymbol{\mid} \infty, 0 \big) \displaystyle\prod_{j = 1}^N q^{\binom{|\textbf{E}_j|}{2}} x_{N - j + 1}^{|\textbf{E}_j|},
	\end{flalign}
	
	\noindent where we sum over all entrance data $\mathscr{E} = (\textbf{E}_1, \textbf{E}_2, \ldots )$ on $\mathcal{D}_N$ such that $\textbf{E}_k = \textbf{e}_0$ for $k > N$.
	
	To that end, observe from \eqref{gfhe} and \Cref{weightesum} that $\mathcal{G}_{\widetilde{\boldsymbol{\lambda}}} (\textbf{x}; \infty \boldsymbol{\mid} 0; 0)$ is the partition function under the $\mathcal{W}_x (\textbf{A}, \textbf{B}; \textbf{C}, \textbf{D} \boldsymbol{\mid} \infty, 0)$ weights for the vertex model $\mathfrak{P}_G (\boldsymbol{\lambda} / \boldsymbol{\varnothing}; N)$ from \Cref{pgpfph} (see also the middle of \Cref{sumw}). Thus the second statement of \eqref{zabcdzabcd} (or, diagrammatically, by separating out the leftmost column of the model, as on the middle and right of \Cref{sumw}) gives 
	\begin{flalign}
	\label{gsumw2} 
	\mathcal{G}_{\widetilde{\boldsymbol{\lambda}}} (\textbf{x}; \infty \boldsymbol{\mid} 0; 0) = \displaystyle\sum_{\mathscr{E}} \mathcal{W}_{\mathcal{D}_N} \big( \mathscr{E}; \mathscr{I} (\lambda) \boldsymbol{\mid} \textbf{z} \boldsymbol{\mid} \infty, 0 \big) \zeta_N (\mathscr{E}; \textbf{x})
	\end{flalign}
	
	\noindent where we sum over all entrance data $\mathscr{E} = (\textbf{E}_1, \textbf{E}_2, \ldots )$ on $\mathcal{D}_N$ such that $\textbf{E}_k = \textbf{e}_0$ for $k > N$. Here, $\zeta_N (\mathscr{E}; \textbf{x})$ is defined to be the partition function on the column $\mathcal{D}_{1, N} = \{ 1 \} \times \{ 1, 2, \ldots , N \} \subset \mathbb{Z}^2$ under the weight $\mathcal{W}_{x_j} \big( \textbf{A} (v), \textbf{B} (v); \textbf{C} (v), \textbf{D} (v) \big)$ at any $v = (1, j) \in \mathcal{D}_{1, N}$, with the following boundary conditions. Its entrance data is $\big( \textbf{e}_0, \textbf{e}_0, \ldots , \textbf{e}_0, \textbf{e}_{[1, n]} \big)$ (where $\textbf{e}_0$ appears with multiplicity $N$), and its exit data is $(\textbf{E}_0, \textbf{E}_1, \textbf{E}_2, \ldots , \textbf{E}_N)$, where $\textbf{E}_0 \in \{ 0, 1 \}^n$ is defined so that $\sum_{i = 0}^N \textbf{E}_i = \textbf{e}_{[1, n]}$. We refer to the right side of \Cref{sumw} for a depiction. 
	
	Thus, by \eqref{gsumw2}, to establish \eqref{gsumw1} we must verify that
	\begin{flalign}
	\label{zetaw}
	\zeta_N (\mathscr{E}; \textbf{x}) = q^{\inv (\mathscr{E})} \displaystyle\prod_{j = 1}^N q^{\binom{|\textbf{E}_j|}{2}} x_{N - j + 1}^{|\textbf{E}_j|}.
	\end{flalign}
	
	\noindent To that end, observe by arrow conservation that there is a single path ensemble supported by $\zeta_N (\mathscr{E}; \textbf{x})$, for any fixed sequence $\mathscr{E} = (\textbf{E}_1, \textbf{E}_2, \ldots )$ of elements in $\{ 0, 1 \}^n$. Denote its arrow configuration at any vertex $v \in \mathcal{D}_{1, N}$ by $\big( \textbf{A} (v), \textbf{B} (v); \textbf{C} (v), \textbf{D} (v) \big)$, which must satisfy $\textbf{B} (1, j) = \textbf{e}_0$ and $\textbf{D} (1, j) = \textbf{E}_{N - j + 1}$ for each $j \in [1, N]$ (see the right side of \Cref{sumw}). Then, by the explicit form \eqref{limitw} for the $\mathcal{W}_z (\textbf{A}, \textbf{B}; \textbf{C}, \textbf{D} \boldsymbol{\mid} \infty, 0)$ weights and the equality $\varphi (\textbf{E}_{N - j + 1}, \textbf{E}_{N - j + 1}) = \binom{|\textbf{E}_{N - j + 1}|}{2}$ (recall $\varphi$ from \eqref{tufunction}), we obtain
	\begin{flalign*}
	\zeta_N (\mathscr{E}; \textbf{x}) & = \displaystyle\prod_{j = 1}^N \mathcal{W}_{x_j} \big( \textbf{A} (1, j), \textbf{e}_0; \textbf{C} (1, j), \textbf{E}_{N - j + 1} \boldsymbol{\mid} \infty, 0  \big) \\
	& = \displaystyle\prod_{j = 1}^N x_j^{|\textbf{E}_{N - j + 1}|} q^{\binom{|\textbf{E}_{N - j + 1}|}{2}} q^{\varphi (\textbf{E}_{N - j + 1}, \textbf{C} (1, j))} = q^{\inv (\mathscr{E})} \displaystyle\prod_{j = 1}^N q^{\binom{|\textbf{E}_{N - j + 1}|}{2}} x_j^{|\textbf{E}_{N - j + 1}|}.
	\end{flalign*}
	
	\noindent Here, we have used the fact that 
	\begin{flalign*}
	\displaystyle\sum_{j = 1}^N \varphi \big(  \textbf{E}_{N - j + 1}, \textbf{C} (1, j) \big) & = \displaystyle\sum_{j = 1}^N \varphi \Bigg( \textbf{E}_{N - j + 1}, \displaystyle\sum_{i = 0}^{N - j} \textbf{E}_i \Bigg) = \displaystyle\sum_{1 \le i < j \le N} \varphi (\textbf{E}_j, \textbf{E}_i) = \inv (\mathscr{E}),
	\end{flalign*}
	
	\noindent where the first equality follows from arrow conservation; the second from the bilinearity of $\varphi$ and the fact that $\varphi (\textbf{E}_i, \textbf{E}_0) = 0$ for $i > 0$, which holds since each color vertically exiting the leftmost column is less than each color horizontally exiting it (as $\widetilde{\lambda}^{(i)} = (0)$ if and only if $i \in [1, n - M]$); and the third from the definition of $\inv$ (see \eqref{sumij}). This yields \eqref{zetaw} and thus the lemma. 
\end{proof}

By \Cref{llambdasum}, the LLT polynomial $\mathcal{L}_{\widetilde{\lambda}}^{(n)} (\textbf{x})$ is given by a linear combination of partition functions for (fermionic) colored vertex models. To proceed, we first recall from \cite{MPI} an expression for the modified Hall--Littlewood polynomial $Q_{\lambda}' (\textbf{x})$ through (bosonic) uncolored vertex models (see \Cref{qlambdasum} below). Then we will provide a color merging result, similar to (but slightly different from) \Cref{zsumfused}, that equates linear combinations of colored vertex model partition functions with uncolored vertex model partition functions (see \Cref{psump} below). 

To that end, we require some additional notation. For any nonnegative integers $a, b, c, d \in \mathbb{Z}_{\ge 0}$ and complex number $z \in \mathbb{C}$, define (similarly to in equation (59) of \cite{MPI})
\begin{flalign}
\label{uweight} 
\mathcal{U}_z (a, b; c, d) & = x^d q^{\binom{d}{2}} \displaystyle\frac{(q; q)_{a + b}}{(q; q)_a (q; q)_b} \textbf{1}_{a + b = c + d}; \qquad \mathcal{U}_z' (a, b; c, d) = x^b q^{\binom{b}{2}} \displaystyle\frac{(q; q)_{a + b}}{(q; q)_a (q; q)_b} \textbf{1}_{a + b = c + d}.
\end{flalign} 

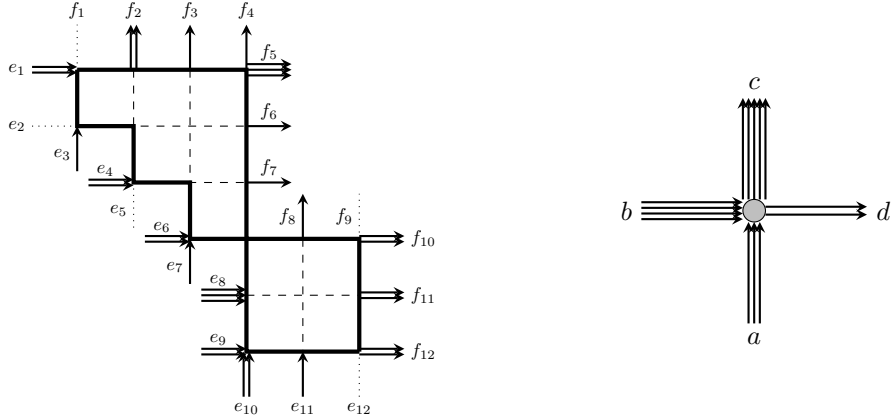
\begin{figure}
	
	\begin{center}		
		
		\begin{tikzpicture}[
		>=stealth,
		auto,
		style={
			scale = .75
		}
		]
		
		\draw[dashed] (1, 8) -- (1, 7) -- (3, 7);
		\draw[dashed] (2, 8) -- (2, 6) -- (3 ,6);
		\draw[dashed] (3, 4) -- (4, 4) -- (4, 3);
		\draw[dashed] (4, 5) -- (4, 4) -- (5, 4);
		
		\draw[] (-.8, 8) circle[radius = 0] node[left, black, scale = .7]{$e_1$};
		\draw[] (-.8, 7) circle[radius = 0] node[left, black, scale = .7]{$e_2$};
		\draw[] (0, 6.5) circle[radius = 0] node[left, black, scale = .7]{$e_3$};
		\draw[] (.5, 6) circle[radius = 0] node[above, black, scale = .7]{$e_4$};
		\draw[] (1, 5.5) circle[radius = 0] node[left, black, scale = .7]{$e_5$};
		\draw[] (1.5, 5) circle[radius = 0] node[above, black, scale = .7]{$e_6$};
		\draw[] (2, 4.5) circle[radius = 0] node[left, black, scale = .7]{$e_7$};
		\draw[] (2.5, 4.05) circle[radius = 0] node[above, black, scale = .7]{$e_8$};
		\draw[] (2.5, 3) circle[radius = 0] node[above, black, scale = .7]{$e_9$};
		\draw[] (3, 2.2) circle[radius = 0] node[below, black, scale = .7]{$e_{10}$};
		\draw[] (4, 2.2) circle[radius = 0] node[below, black, scale = .7]{$e_{11}$};
		\draw[] (5, 2.2) circle[radius = 0] node[below, black, scale = .7]{$e_{12}$};
		
		\draw[] (0, 8.8) circle[radius = 0] node[above, black, scale = .7]{$f_1$};
		\draw[] (1, 8.8) circle[radius = 0] node[above, black, scale = .7]{$f_2$};
		\draw[] (2, 8.8) circle[radius = 0] node[above, black, scale = .7]{$f_3$};
		\draw[] (3, 8.8) circle[radius = 0] node[above, black, scale = .7]{$f_4$};
		\draw[] (3.4, 8.05) circle[radius = 0] node[above, black, scale = .7]{$f_5$};
		\draw[] (3.4, 7) circle[radius = 0] node[above, black, scale = .7]{$f_6$};
		\draw[] (3.4, 6) circle[radius = 0] node[above, black, scale = .7]{$f_7$};
		\draw[] (4, 5.4) circle[radius = 0] node[left, black, scale = .7]{$f_8$};
		\draw[] (5, 5.4) circle[radius = 0] node[left, black, scale = .7]{$f_9$};
		\draw[] (5.8, 5) circle[radius = 0] node[right, black, scale = .7]{$f_{10}$};
		\draw[] (5.8, 4) circle[radius = 0] node[right, black, scale = .7]{$f_{11}$};
		\draw[] (5.8, 3) circle[radius = 0] node[right, black, scale = .7]{$f_{12}$};
		
		\draw[thick, ->] (-.8, 8.05) -- (0, 8.05);
		\draw[thick, ->] (-.8, 7.95) -- (0, 7.95);
		
		\draw[dotted] (-.8, 7) -- (0, 7);
		
		\draw[thick, ->] (.2, 5.95) -- (1, 5.95);
		\draw[thick, ->] (.2, 6.05) -- (1, 6.05);
		
		\draw[thick, ->] (1.2, 5.05) -- (2, 5.05);
		\draw[thick, ->] (1.2, 4.95) -- (2, 4.95);
		
		\draw[thick, ->] (2.2, 4.1) -- (3, 4.1);
		\draw[thick, ->] (2.2, 4) -- (3, 4);
		\draw[thick, ->] (2.2, 3.9) -- (3, 3.9);
		
		\draw[thick, ->] (2.2, 3.05) -- (3, 3.05);
		\draw[thick, ->] (2.2, 2.95) -- (3, 2.95);

		\draw[thick, ->] (0, 6.2) -- (0, 7);
		\draw[dotted] (1, 5.2) -- (1, 6);
		\draw[thick, ->] (2, 4.2) -- (2, 5);
		
		\draw[thick, ->] (3.05, 2.2) -- (3.05, 3);
		\draw[thick, ->] (2.95, 2.2) -- (2.95, 3);
		
		\draw[thick, ->] (4, 2.2) -- (4, 3);
		
		\draw[dotted] (5, 2.2) -- (5, 3);
		
		\draw[thick, ->] (3, 7.9) -- (3.8, 7.9);
		\draw[thick, ->] (3, 8) -- (3.8, 8);
		\draw[thick, ->] (3, 8.1) -- (3.8, 8.1);

		\draw[thick, ->] (3, 7) -- (3.8, 7);
		\draw[thick, ->] (3, 6) -- (3.8, 6);
		
		\draw[thick, ->] (5, 4.95) -- (5.8, 4.95);
		\draw[thick, ->] (5, 5.05) -- (5.8, 5.05);
		
		\draw[thick, ->] (5, 4.05) -- (5.8, 4.05);
		\draw[thick, ->] (5, 3.95) -- (5.8, 3.95);
		
		\draw[thick, ->] (5, 3.05) -- (5.8, 3.05);
		\draw[thick, ->] (5, 2.95) -- (5.8, 2.95);
		
		\draw[dotted] (0, 8) -- (0, 8.8);
		
		\draw[thick, ->] (1.05, 8) -- (1.05, 8.8);
		\draw[thick, ->] (.95, 8) -- (.95, 8.8);
		
		\draw[thick, ->] (2, 8) -- (2, 8.8);
		
		\draw[thick, ->] (3, 8) -- (3, 8.8);
		
		\draw[thick, ->] (4, 5) -- (4, 5.8);
		
		\draw[dotted] (5, 5) -- (5, 5.8);
		
		\draw[ultra thick] (0, 8) -- (0, 7) -- (1, 7) -- (1, 6) -- (2, 6) -- (2, 5) -- (3, 5) -- (3, 3) -- (5, 3);
		\draw[ultra thick] (0, 8) -- (3, 8) -- (3, 7) -- (3, 5) -- (4, 5) -- (5, 5) -- (5, 3);

		\draw[fill = gray!50!white] (12, 5.5) circle[radius = .2];
		
		\draw[thick, ->] (12, 3.5) -- (12, 5.3);
		\draw[thick, ->] (11.9, 3.5) -- (11.9, 5.3);
		\draw[thick, ->] (12.1, 3.5) -- (12.1, 5.3);
		
		\draw[thick, ->] (10, 5.35) -- (11.8, 5.35);
		\draw[thick, ->] (10, 5.45) -- (11.8, 5.45);
		\draw[thick, ->] (10, 5.55) -- (11.8, 5.55);
		\draw[thick, ->] (10, 5.65) -- (11.8, 5.65);
		
		\draw[thick, ->] (12.2, 5.43) -- (14, 5.43);
		\draw[thick, ->] (12.2, 5.57) -- (14, 5.57);
		
		\draw[thick, ->] (11.8, 5.7) -- (11.8, 7.5);
		\draw[thick, ->] (11.9, 5.7) -- (11.9, 7.5);
		\draw[thick, ->] (12, 5.7) -- (12, 7.5);
		\draw[thick, ->] (12.1, 5.7) -- (12.1, 7.5);
		\draw[thick, ->] (12.2, 5.7) -- (12.2, 7.5);
		
		\draw[] (12, 3.5) circle[radius = 0] node[below]{$a$};
		\draw[] (10, 5.5) circle[radius = 0] node[left]{$b$};
		\draw[] (12, 7.5) circle[radius = 0] node[above]{$c$};
		\draw[] (14, 5.5) circle[radius = 0] node[right]{$d$};

		\end{tikzpicture}
		
	\end{center}
	
	\caption{\label{uncoloredpaths} Shown to the left is an east-south domain with uncolored boundary data, and shown to the right is an uncolored vertex with $(a, b; c, d) = (3, 4; 5, 2)$. } 	
\end{figure}

Once again, $\mathcal{U}_z$ and $\mathcal{U}_z'$ can be interpreted diagrammatically, now as uncolored vertex weights. As previously, a vertex $v$ is the intersection between two transverse curves $\ell_1$ and $\ell_2$, and associated with each vertex is a spectral parameter $z$. Each of the four edges adjacent to $v$ may accommodate any nonnegative number of \emph{uncolored} arrows, that is, arrows are not assigned any color (equivalently, we may view all arrows to be of the same, unlabeled color). Let $a$, $b$, $c$, and $d$ denote the numbers of arrows that vertically enter, horizontally enter, vertically exit, and horizontally exit $v$ respectively.\footnote{Unlike as for the colored arrow configurations in \Cref{Weightszq}, the entries of an uncolored arrow configuration $(a, b; c, d)$ denote the numbers of arrows entering or exiting the vertex (and do not index the colors of the arrows).} We will assume that $a + b = c + d$, which is a form of arrow conservation. The quadruple $(a, b; c, d)$ is the \emph{uncolored arrow configuration} at $v$, and we interpret $\mathcal{U}_z (a, b; c, d)$ (or $\mathcal{U}_z' (a, b; c, d)$) as the associated vertex weight. Such vertices are bosonic, since we allow multiple arrows (of the same, unlabeled color) to exist along edges; we refer to the right side of \Cref{uncoloredpaths} for a depiction. 

A \emph{(bosonic) uncolored path ensemble} on a domain $\mathcal{D} \subset \mathbb{Z}^2$ is a consistent assignment of a  arrow configuration $\big( a (v), b (v); c (v), d (v) \big)$ to each vertex $v \in \mathcal{D}$, meaning that $b (u) = d(v)$ if $u - v = (1, 0)$ and $a (u) = c (v)$ if $u - v = (0, 1)$. Arrow conservation $a + b = c + d$ implies that any path ensemble may be viewed as a collection of colored paths, which may share edges but do not cross. 

Now fix an east-south domain $\mathcal{D} \subset \mathbb{Z}^2$ whose boundary paths are both of length $k$ (possibly infinite). For any $(k + 1)$-tuples $\mathsf{E} = (e_1, e_2, \ldots , e_{k + 1})$ and $\mathsf{F} = (f_1, f_2, \ldots , f_{k + 1})$ of nonnnegative integers, an uncolored path ensemble has \emph{uncolored boundary data} $(\mathsf{E}; \mathsf{F})$ if, for each $i \in [1, k + 1]$, $e_i$ arrows enter through the $i$-th incoming edge in $\mathcal{D}$ and $f_i$ arrows exit through its $i$-th outgoing edge.  We refer to $\mathsf{E}$ as uncolored entrance data on $\mathcal{D}$ and $\mathsf{F}$ as exit data. For example, $\mathsf{E} = (2, 0, 1, 2, 0, 2, 1, 3, 2, 2, 1, 0)$ and $\mathsf{F} = (0, 2, 1, 1, 3, 1, 1, 1, 0, 2, 2, 2)$ on the left side of \Cref{uncoloredpaths}.

Next let us consider partition functions for uncolored path ensembles. Fix some east-south domain $\mathcal{D} \subset \mathbb{Z}^2$; some uncolored boundary data $(\mathsf{E}; \mathsf{F})$ on $\mathcal{D}$; and a set of complex numbers $\textbf{z} = \big( z(v) \big)_{v \in \mathcal{D}}$. For each index $\mathcal{X} \in \{ \mathcal{U}, \mathcal{U}' \}$, define the \emph{uncolored partition function}
\begin{flalign*}
\mathcal{X}_{\mathcal{D}} (\mathsf{E}; \mathsf{F} \boldsymbol{\mid} \textbf{z}) & = \displaystyle\sum \displaystyle\prod_{v \in \mathcal{D}} \mathcal{X}_{z(v)} \big( a(v), b(v); c(v), d (v)\big),
\end{flalign*} 

\noindent where the sum is over all uncolored path ensembles with boundary data $(\mathsf{E}; \mathsf{F})$. 

Again, certain uncolored entrance and exit data on $\mathcal{D}_N \subset \mathbb{Z}_{\ge 0}^2$ will be of particular use to us. First, any uncolored entrance data $\mathsf{E} = (e_1, e_2, \ldots )$ on $\mathcal{D}_N$ with $e_k = 0$ for $k > N$ corresponds to prohibiting paths from entering $\mathcal{D}_N$ along its bottom boundary. Second, for any signature $\lambda \in \Sign$, define the sequence $\mathsf{I} (\lambda) = \big( m_1 (\lambda), m_2 (\lambda), \ldots \big)$ of nonnegative integers, which is obtained from the colored exit data $\mathscr{I} (\lambda)$ by ignoring the colors of all arrows (viewing them as uncolored). In particular, $m_k (\lambda)$ denotes the multiplicity of $k$ in $\lambda$, for each integer $k \ge 1$.  

Under this notation, we can recall from \cite{MPI} the expression for modified Hall--Littlewood polynomials through such uncolored partition functions under the above boundary data.

\begin{prop}[{\cite[Equation (66)]{MPI}}]

\label{qlambdasum} 

Fix integers $N, M \ge 1$; a partition $\lambda \in \Sign_M$; and a sequence of complex numbers $\textbf{\emph{x}} = (x_1, x_2, \ldots , x_N)$. Defining the variable set $\textbf{\emph{z}} = \big( z(v) \big)_{\mathcal{D}_N}$ so that $z(v) = z (i, j) = x_j$ for each $v = (i, j) \in \mathcal{D}_N$, we have 
\begin{flalign}
\label{qu} 
Q_{\lambda}' (\textbf{\emph{x}}) = \displaystyle\sum_{\mathsf{E}} \mathcal{U}_{\mathcal{D}_N}' \big( \mathsf{E}; \mathsf{I} (\lambda) \boldsymbol{\mid} \textbf{\emph{x}} \big) = \displaystyle\sum_{\mathsf{E}} \mathcal{U}_{\mathcal{D}_N} \big( \mathsf{E}; \mathsf{I} (\lambda) \boldsymbol{\mid} \textbf{\emph{x}} \big)  \displaystyle\prod_{i = 1}^N q^{\binom{e_i}{2}} x_{N - i + 1}^{e_i},
\end{flalign}

\noindent where both sums are over all uncolored entrance data $\mathsf{E} = (e_1, e_2, \ldots )$ on $\mathcal{D}_N$, with $e_k = 0$ for $k > N$. 

\end{prop}

\begin{rem}
	
	\label{qlambdau}
	
	The first statement of \eqref{qu} is what was established in and below equation (66) of \cite{MPI} (after rotating the vertex model there by $180$ degrees). The second follows from the first, together with the gauge transformation $\mathcal{U}_z' (a, b; c, d) = x^{b - d} q^{\binom{b}{2} - \binom{d}{2}} \mathcal{U}_z (a, b; c, d)$ relating the weights $\mathcal{U}_z$ and $\mathcal{U}_z'$ from \eqref{uweight}.
\end{rem}

We next state the below color merging result, which essentially states the following. Consider a colored vertex model partition function, and $q$-symmetrize its entrance data $\mathscr{E}$ over all ``colorings'' of some fixed uncolored entrance data $\mathsf{E}$. This yields the partition function for the uncolored vertex model obtained from the original one by ignoring all of its arrow colors.

\begin{lem}
	
	\label{psump}
	
	Fix an integer $k \ge 1$; an east-south domain $\mathcal{D} \subset \mathbb{Z}^2$ whose boundary paths are both length $k$; and a set of complex numbers $\textbf{\emph{z}} = \big( z(v) \big)_{v \in \mathcal{D}}$. Further fix a sequence of elements $\mathscr{F} = \big( \textbf{\emph{F}}_1, \textbf{\emph{F}}_2, \ldots , \textbf{\emph{F}}_k \big)$ in $\{ 0, 1 \}^n$ constituting (colored) exit data on $\mathcal{D}$, and sequences $\mathsf{E} = (e_1, e_2, \ldots , e_k) \in \mathbb{Z}_{\ge 0}^k$ and $\mathsf{F} = (f_1, f_2, \ldots , f_k) \in \mathbb{Z}_{\ge 0}^k$ constituting uncolored entrance and exit data on $\mathcal{D}$, respectively; assume that $|\textbf{\emph{F}}_i| = f_i$ for each $i \in [1, k]$ and that $\sum_{i = 1}^k \textbf{\emph{F}}_i \le \textbf{\emph{e}}_{[1, n]}$. Then,
	\begin{flalign}
	\label{sumqw} 
	\displaystyle\sum_{\mathscr{E}} q^{\inv (\mathscr{E}) - \inv (\mathscr{F})} \mathcal{W}_{\mathcal{D}} (\mathscr{E}; \mathscr{F} \boldsymbol{\mid} \textbf{\emph{z}} \boldsymbol{\mid} \infty, 0) = \mathcal{U}_{\mathcal{D}} (\mathsf{E}; \mathsf{F} \boldsymbol{\mid} \textbf{\emph{z}}),
	\end{flalign}
	
	\noindent where the sum is over all sequences $\mathscr{E} = (\textbf{\emph{E}}_1, \textbf{\emph{E}}_2, \ldots , \textbf{\emph{E}}_k)$ of elements in $\{ 0, 1 \}^n$ such that $|\textbf{\emph{E}}_i| = e_i$ for each $i \in [1, k]$. 
\end{lem}

\begin{proof}[Proof (Outline)]
	
	We only establish the lemma when $\mathcal{D} = \{ v \}$ consists of a single vertex, as the general case follows from this one similarly to as in the proof of \Cref{zsumef} in \Cref{ProofSumZ}. So, let us abbreviate $z = z(v)$; we may also assume that $e_1 + e_2 = f_1 + f_2$, as otherwise both sides of \eqref{sumqw} are equal to $0$ by arrow conservation.  
	
	In this case, $k = 2$, and so the sum in \eqref{sumqw} is over all pairs $\mathscr{E} = (\textbf{E}_1, \textbf{E}_2)$ of elements in $\{ 0, 1 \}^n$ such that $|\textbf{E}_i| = e_i$, for each $i \in \{ 1, 2 \}$. Moreover, by \eqref{sumij}, we have $\inv (\mathscr{E}) = \varphi (\textbf{E}_2, \textbf{E}_1)$ and $\inv (\mathscr{F}) = \varphi (\textbf{F}_2, \textbf{F}_1)$. So, since $\mathcal{W}_{\mathcal{D}} (\mathscr{E}; \mathscr{F} \boldsymbol{\mid} \textbf{z} \boldsymbol{\mid} \infty, 0) = \mathcal{W}_z (\textbf{E}_2, \textbf{E}_1; \textbf{F}_1, \textbf{F}_2 \boldsymbol{\mid} \infty, 0)$ and $\mathcal{U}_{\mathcal{D}} (\mathsf{E}; \mathsf{F} \boldsymbol{\mid} \textbf{z}) = \mathcal{U}_z (e_2, e_1; f_1, f_2)$ it suffices by \eqref{uweight} to show that 
	\begin{flalign*}
	\displaystyle\sum_{\mathscr{E}} q^{\varphi (\textbf{E}_2, \textbf{E}_1) - \varphi (\textbf{F}_2, \textbf{F}_1)} \mathcal{W}_z (\textbf{E}_2, \textbf{E}_1; \textbf{F}_1, \textbf{F}_2 \boldsymbol{\mid} \infty, 0) = \mathcal{U}_z (e_2, e_1; f_1, f_2) =  q^{\binom{f_2}{2}} z^{f_2} \displaystyle\frac{(q; q)_{e_1 + e_2}}{(q; q)_{e_1} (q; q)_{e_2}},
	\end{flalign*}
	
	\noindent where we sum over all pairs $\mathscr{E} = (\textbf{E}_1, \textbf{E}_2)$ with $|\textbf{E}_i| = e_i$, for each $i \in \{ 1, 2 \}$. Recalling from \eqref{limitw} that $\mathcal{W}_z (\textbf{E}_1, \textbf{E}_1; \textbf{F}_1, \textbf{F}_2 \boldsymbol{\mid} \infty, 0) = q^{\binom{f_2}{2} + \varphi (\textbf{F}_2, \textbf{F}_1)} z^{f_2}$, we must therefore verify
	\begin{flalign}
	\label{sumef} 
	\displaystyle\sum_{\mathscr{E}} q^{\varphi (\textbf{E}_2, \textbf{E}_1)} = \displaystyle\frac{(q; q)_{e_1 + e_2}}{(q; q)_{e_1} (q; q)_{e_2}},
	\end{flalign}
	
	\noindent where we sum over all pairs $\mathscr{E} = (\textbf{E}_1, \textbf{E}_2)$ with $|\textbf{E}_i| = e_i$ for each $i \in \{ 1, 2 \}$ and $\textbf{E}_1 + \textbf{E}_2 = \textbf{F}_1 + \textbf{F}_2$ (which may be imposed by arrow conservation). 
	
	To that end, let $\textbf{E}_i = (E_{i, 1}, E_{i, 2}, \ldots , E_{I, n})$ for each $i \in \{ 1, 2 \}$; since $\mathscr{F} = (\textbf{F}_1, \textbf{F}_2)$ is fixed, the condition $\textbf{E}_1 + \textbf{E}_2 = \textbf{F}_1 + \textbf{F}_2 \le \textbf{e}_{[1, n]}$ (which was imposed by the statement of the lemma) fixes the indices $j \in [1, n]$ for which we have $E_{1, j} + E_{2, j} = 0$ or $E_{1, j} + E_{2, j} = 1$. Indices $j$ for which $E_{1, j} + E_{2, j} = 0$ impose $E_{1, j} = E_{2, j} = 0$ and do not contribute to the exponent $q^{\varphi (\textbf{E}_2, \textbf{E}_1)}$ on the left side of \eqref{sumef}. Thus, we may remove such indices from consideration, thereby assuming that $n = e_1 + e_2$. Then, \eqref{sumef} is equivalent to
	\begin{flalign}
	\label{sumef2}
	\displaystyle\sum_{\substack{\textbf{E}_1 + \textbf{E}_2 = \textbf{e}_{[1, n]} \\ |\textbf{E}_1| = e_1 \\ |\textbf{E}_2| = e_2}} \displaystyle\prod_{1 \le j < k \le n} q^{E_{1, k} E_{2, j}} = \displaystyle\frac{(q; q)_{e_1 + e_2}}{(q; q)_{e_1} (q; q)_{e_2}}.
	\end{flalign} 
	
	\noindent Observe in this sum that $\textbf{E}_1$ is fixed by $\textbf{E}_2$ and that $E_{1, k} E_{2, j} = 1$ holds if and only if $(E_{2, j}, E_{2, k}) = (1 ,0)$ or, equivalently, $E_{2, j} > E_{2, k}$. Thus, \eqref{sumef2} follows from the interpretation of the $q$-binomial coefficient $(q; q)_{e_1 + e_2} (q; q)_{e_1}^{-1} (q; q)_{e_2}^{-1}$ as the generating series for inversions of $\{ 0, 1 \}^{e_1 + e_2}$ strings with $e_1$ zeroes and $e_2$ ones. 
\end{proof}

Now we can establish \Cref{qlambdaqlambda}.

\begin{proof}[Proof of \Cref{qlambdaqlambda}]
	
	In what follows, we adopt the notation of \Cref{llambdasum}. That lemma and the fact that $\inv \big( \mathscr{I} (\lambda) \big)= 0$ for any partition $\lambda$, together imply
	\begin{flalign}
	\label{llambdasum1}
	\begin{aligned} 
	\mathcal{L}_{\widetilde{\lambda}}^{(n)} (\textbf{x}) & = \displaystyle\sum_{\mathscr{E}} q^{\inv (\mathscr{E})} \mathcal{W}_{\mathcal{D}_N} \big( \mathscr{E}; \mathscr{I} (\lambda) \boldsymbol{\mid} \textbf{z} \boldsymbol{\mid} \infty, 0 \big) \displaystyle\prod_{j = 1}^N q^{\binom{|\textbf{E}_j|}{2}} x_{N - j + 1}^{|\textbf{E}_j|} \\
	& = \displaystyle\sum_{\mathsf{E}} \displaystyle\prod_{j = 1}^N q^{\binom{e_j}{2}} x_{N - j + 1}^{e_j} \displaystyle\sum_{\mathscr{E}} q^{\inv (\mathscr{E}) - \inv (\mathscr{I} (\lambda))} \mathcal{W}_{\mathcal{D}_N} \big( \mathscr{E}; \mathscr{I} (\lambda) \boldsymbol{\mid} \textbf{z} \boldsymbol{\mid} \infty, 0 \big).
	\end{aligned} 
	\end{flalign}
	
	\noindent Here, the sum in the second term is over all sequences $\mathscr{E} = (\textbf{E}_1, \textbf{E}_2, \ldots)$ of elements in $\{ 0, 1 \}^n$ such that $\textbf{E}_k = \textbf{e}_0$ for $k > N$; the first sum in the third term is over all integer sequences $\mathsf{E} = (e_1, e_2, \ldots )$ such that $e_k = 0$ for $k > N$; and the second sum in the third term is over all sequences $\mathscr{E} = (\textbf{E}_1, \textbf{E}_2, \ldots )$ of elements in $\{ 0, 1 \}^n$ such that $|\textbf{E}_i| = e_i$ for each $i \ge 1$. 
	
	Applying \eqref{sumqw} in \eqref{llambdasum1}, we deduce 
	\begin{flalign*}
	\mathcal{L}_{\widetilde{\lambda}}^{(n)} (\textbf{\emph{x}}) & = \displaystyle\sum_{\mathsf{E}} \mathcal{U}_{\mathcal{D}_N} \big( \mathsf{E}; \mathsf{I} (\lambda) \boldsymbol{\mid} \textbf{z} \big) \displaystyle\prod_{j = 1}^N q^{\binom{e_j}{2}} x_{N - j + 1}^{e_j},
	\end{flalign*}

	\noindent which together with \Cref{qlambdasum} yields the proposition.
\end{proof}

\chapter{Proof of \Cref{limitg0horizontal}}

\label{ProofFGL}

In this chapter we establish \Cref{limitg0horizontal} (that is, the $N = 1$ case of \Cref{limitg0}). This consists in four parts, given by the corresponding parts of \Cref{limitg0}. We establish the first, second, third, and fourth in \Cref{GPolynomials}, \Cref{Fx1FunctionL}, \Cref{Functions2GL}, and \Cref{Function3LG}, respectively. Throughout this section, we adopt the notation of \Cref{limitg0} and set $x_1 = x$.

\section{LLT Polynomials From \texorpdfstring{$\mathcal{G}_{\boldsymbol{\lambda} / \boldsymbol{\mu}} (\textbf{x}; \infty \boldsymbol{\mid} 0; 0)$}{}}

\label{GPolynomials}

\begin{figure}

	\begin{center}

		\begin{tikzpicture}[
		>=stealth,
		scale = .5
		]

		\draw[red, thick, ->] (-.1, .5) -- (-.1, 2.4) -- (1.9, 2.4) -- (1.9, 4.5);
		\draw[blue, thick, ->] (0, .5) -- (0, 2.5) -- (4, 2.5) -- (4, 4.5);
		\draw[green, thick, ->] (.1, .5) -- (.1, 4.5);
		\draw[red, thick, ->] (3.9, .5) -- (3.9, 2.4) -- (7.9, 2.4) -- (7.9, 4.5);
		\draw[green, thick, ->] (6.1, .5) -- (6.1, 2.6) -- (8.1, 2.6) -- (8.1, 4.5);
		\draw[blue, thick, ->] (6, .5) -- (6, 4.5);
	
		\draw[] (0, .5) circle[radius = 0] node[below, scale = .65]{$\textbf{A}(1)$};
		\draw[] (2, .5) circle[radius = 0] node[below, scale = .65]{$\textbf{A}(2)$};
		\draw[] (4, .5) circle[radius = 0] node[below, scale = .65]{$\textbf{A}(3)$};
		\draw[] (6, .5) circle[radius = 0] node[below, scale = .65]{$\textbf{A}(4)$};
		\draw[] (8, .5) circle[radius = 0] node[below, scale = .65]{$\textbf{A}(5)$};
		
		\draw[] (0, 4.5) circle[radius = 0] node[above, scale = .65]{$\textbf{C}(1)$};
		\draw[] (2, 4.5) circle[radius = 0] node[above, scale = .65]{$\textbf{C}(2)$};
		\draw[] (4, 4.5) circle[radius = 0] node[above, scale = .65]{$\textbf{C}(3)$};
		\draw[] (6, 4.5) circle[radius = 0] node[above, scale = .65]{$\textbf{C}(4)$};
		\draw[] (8, 4.5) circle[radius = 0] node[above, scale = .65]{$\textbf{C}(5)$};
		
		\draw[] (-1.5, 4.5) circle[radius = 0] node[above, scale = .85]{$\boldsymbol{\lambda}$};
		\draw[] (-1.5, .5) circle[radius = 0] node[below, scale = .85]{$\boldsymbol{\mu}$};

		\draw[dotted, red] (12, 3) -- (12, 4);
		\draw[dotted, blue] (13, 3) -- (13, 4);
		\draw[->, green] (14, 3) -- (14, 4);
		\draw[->, red] (15, 3) -- (15, 4);
		\draw[dotted, blue] (16, 3) -- (16, 4);
		\draw[dotted, green] (17, 3) -- (17, 4);
		\draw[dotted, red] (18, 3) -- (18, 4);
		\draw[->, blue] (19, 3) -- (19, 4);
		\draw[dotted, green] (20, 3) -- (20, 4);
		\draw[dotted, red] (21, 3) -- (21, 4);
		\draw[->, blue] (22, 3) -- (22, 4);
		\draw[dotted, green] (23, 3) -- (23, 4);
		\draw[->, red] (24, 3) -- (24, 4);	
		\draw[dotted, blue] (25, 3) -- (25, 4);
		\draw[->, green] (26, 3) -- (26, 4);	
		
		\draw[->, red] (12, 1) -- (12, 2);
		\draw[->, blue] (13, 1) -- (13, 2);
		\draw[->, green] (14, 1) -- (14, 2);
		\draw[dotted, red] (15, 1) -- (15, 2);
		\draw[dotted, blue] (16, 1) -- (16, 2);
		\draw[dotted, green] (17, 1) -- (17, 2);
		\draw[->, red] (18, 1) -- (18, 2);
		\draw[dotted, blue] (19, 1) -- (19, 2);
		\draw[dotted, green] (20, 1) -- (20, 2);
		\draw[dotted, red] (21, 1) -- (21, 2);
		\draw[->, blue] (22, 1) -- (22, 2);
		\draw[->, green] (23, 1) -- (23, 2);
		\draw[dotted, red] (24, 1) -- (24, 2);	
		\draw[dotted, blue] (25, 1) -- (25, 2);
		\draw[dotted, green] (26, 1) -- (26, 2); 
		
		\draw[dashed, blue, ->] (13, 2.125) arc (120:60:3); 
		\draw[dashed, blue, ->] (16, 2.125) arc (120:60:3); 
		\draw[dashed, red, ->] (18, 2.125) arc (120:60:3); 
		\draw[dashed, red, ->] (12, 2.125) arc (120:60:3); 
		\draw[dashed, red, ->] (21, 2.125) arc (120:60:3); 
		\draw[dashed, green, ->] (23, 2.125) arc (120:60:3); 
		
		\draw[ultra thick] (11.5, 3) node[left, scale = .75]{$\lambda$} -- (26.5, 3);
		\draw[ultra thick] (11.5, 1) node[left, scale = .75]{$\mu$} -- (26.5, 1);
		
		\draw[] (12, .75) -- (12, .5) -- (14, .5) -- (14, .75); 
		\draw[] (15, .75) -- (15, .5) -- (17, .5) -- (17, .75); 
		\draw[] (18, .75) -- (18, .5) -- (20, .5) -- (20, .75); 
		\draw[] (21, .75) -- (21, .5) -- (23, .5) -- (23, .75); 
		\draw[] (24, .75) -- (24, .5) -- (26, .5) -- (26, .75); 
		
		\draw[] (12, 4.25) -- (12, 4.5) -- (14, 4.5) -- (14, 4.25); 
		\draw[] (15, 4.25) -- (15, 4.5) -- (17, 4.5) -- (17, 4.25); 
		\draw[] (18, 4.25) -- (18, 4.5) -- (20, 4.5) -- (20, 4.25); 
		\draw[] (21, 4.25) -- (21, 4.5) -- (23, 4.5) -- (23, 4.25); 
		\draw[] (24, 4.25) -- (24, 4.5) -- (26, 4.5) -- (26, 4.25); 
		
		\draw[] (13, .5) circle[radius = 0] node[scale = .65, below]{$\textbf{A} (1)$};
		\draw[] (16, .5) circle[radius = 0] node[scale = .65, below]{$\textbf{A} (2)$};
		\draw[] (19, .5) circle[radius = 0] node[scale = .65, below]{$\textbf{A} (3)$};
		\draw[] (22, .5) circle[radius = 0] node[scale = .65, below]{$\textbf{A} (4)$};
		\draw[] (25, .5) circle[radius = 0] node[scale = .65, below]{$\textbf{A} (5)$};
		
		\draw[] (13, 4.5) circle[radius = 0] node[scale = .65, above]{$\textbf{C} (1)$};
		\draw[] (16, 4.5) circle[radius = 0] node[scale = .65, above]{$\textbf{C} (2)$};
		\draw[] (19, 4.5) circle[radius = 0] node[scale = .65, above]{$\textbf{C} (3)$};
		\draw[] (22, 4.5) circle[radius = 0] node[scale = .65, above]{$\textbf{C} (4)$};
		\draw[] (25, 4.5) circle[radius = 0] node[scale = .65, above]{$\textbf{C} (5)$};
		
		\end{tikzpicture}
		
	\end{center}
	
	\caption{\label{sequencearrows} To the left is a vertex model for a $3$-tuple of skew-shapes. To the right is its decompression into a pair of colored Maya diagrams.}
\end{figure}
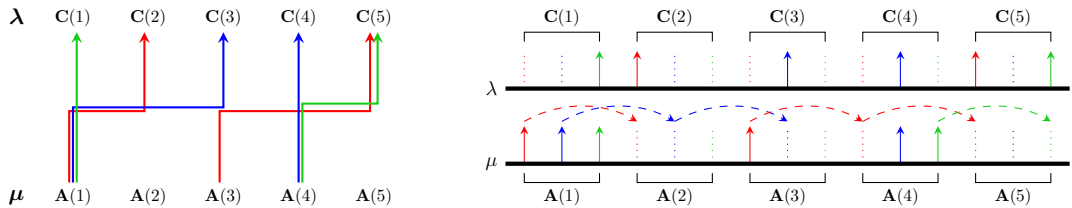

In this section we establish the first part of \Cref{limitg0horizontal}, given by the $N = 1$ case of \eqref{1gl}. To that end, recall from \eqref{gfhe} and \Cref{weightesum} that $\mathcal{G}_{\boldsymbol{\lambda} / \boldsymbol{\mu}} (\textbf{x}; \infty \boldsymbol{\mid} 0; 0)$ is the partition function under the $\mathcal{W}_x (\textbf{A}, \textbf{B}; \textbf{C}, \textbf{D} \boldsymbol{\mid} \infty, 0)$ weights (from \eqref{limitw}) for the vertex model $\mathfrak{P}_G (\boldsymbol{\lambda} / \boldsymbol{\mu}; 1)$ (from \Cref{pgpfph}; see also the left side of \Cref{fgpaths}) on the domain $\mathcal{D} = \mathcal{D}_1 = \mathbb{Z}_{> 0} \times \{ 1 \}$. Since this domain has one row, $\mathfrak{P}_G (\boldsymbol{\lambda} / \boldsymbol{\mu}; 1)$ only consists in a single path ensemble, which we denote by $\mathcal{E} (\boldsymbol{\lambda} / \boldsymbol{\mu})$.\index{E@$\mathcal{E} (\boldsymbol{\lambda} / \boldsymbol{\mu})$} For each $j \ge 1$, let $\big(\textbf{A} (j), \textbf{B} (j); \textbf{C} (j), \textbf{D} (j) \big)$ denote the arrow configuration at the vertex $(j, 1) \in \mathcal{D}$ under $\mathcal{E} (\boldsymbol{\lambda} / \boldsymbol{\mu})$; we refer to the left side of \Cref{sequencearrows} for an example in the case when $\boldsymbol{\lambda} = \big( (3, 1), (2, 2), (3, 0) \big)$ and $\boldsymbol{\mu} = \big( (1, 0), (2, 0), (2, 0) \big)$.

Let us now ``decompress'' $\mathcal{E} (\boldsymbol{\lambda} / \boldsymbol{\mu})$ to form a pair of colored Maya diagrams $\big( \mathfrak{T} (\lambda), \mathfrak{T} (\mu) \big)$, as follows. We first split each vertex $(j, 1)$ of $\mathcal{E} (\boldsymbol{\lambda} / \boldsymbol{\mu})$ into the interval $\{ nj - n + 1, nj - n + 2, \ldots , nj \}$ consisting of $n$ integers. Then, if an arrow of color $c \in [1, n]$ enters through $(j, 1)$, we set $(j - 1) n + c \in \mathfrak{T} (\mu)$; similarly, if an arrow of color $c$ exits through $(j, 1)$, we set $(j - 1) n + c \in \mathfrak{T} (\lambda)$. Observe from the content in \Cref{FunctionGPartitions} that $\lambda, \mu \in \Sign_{nM}$ are characterized by the property that $\boldsymbol{\lambda} / \boldsymbol{\mu}$ is the $n$-quotient of $\lambda / \mu$. We refer to \Cref{sequencearrows} for an example, where $\lambda / \mu = (9, 8, 7, 5, 2, 2) / (6, 6, 3, 0, 0, 0)$. 

Next we establish the following two lemmas expressing $\spin (\lambda / \mu)$ and $\psi (\boldsymbol{\lambda} / \boldsymbol{\mu})$ through $\mathcal{E} (\boldsymbol{\lambda} / \boldsymbol{\mu})$. 

\begin{lem}
	
	\label{spinlambdamu}
	
	If $\lambda / \mu$ is a horizontal $n$-ribbon strip, then (recalling $\varphi$ from \eqref{tufunction}) we have that
	\begin{flalign*}
	\spin (\lambda / \mu) = \displaystyle\sum_{j = 1}^{\infty} \varphi \big( \textbf{\emph{D}} (j), \textbf{\emph{D}} (j) \big) + \displaystyle\frac{1}{2} \displaystyle\sum_{j = 1}^{\infty} \Big( \varphi \big( \textbf{\emph{D}} (j), \textbf{\emph{C}} (j) \big) + \varphi \big( \textbf{\emph{A}} (j), \textbf{\emph{B}} (j) \big) \Big).
	\end{flalign*}
	
\end{lem}

\begin{proof} 
	
	Recall from \Cref{FunctionGPartitions} that, if $\lambda / \mu$ is a horizontal $n$-ribbon strip, then $\spin (\lambda / \mu)$ equals half of the total number of particles jumped over in the nonconflicting series of $n$-jumps $\mathcal{T}_0, \mathcal{T}_1, \ldots , \mathcal{T}_m$ that transforms the Maya diagram $\mathfrak{T} (\mu) = \mathcal{T}_0$ into $\mathfrak{T} (\lambda) = \mathcal{T}_m$. Observe that the paths in $\mathcal{E} (\boldsymbol{\lambda} / \boldsymbol{\mu})$ correpsond to the jump trajectories for particles in this nonconflicting series, in the sense that a horizontal arrow of color $c \in [1, n]$ vertically enters and horizontally exits through a vertex $(j, 1)$ in this model if and only if a particle at site $n (j - 1) + c$ eventually jumps in this nonconflicting series. 
	
	Now, let us classify the ways in which a particle at some site $J \in \mathcal{T}_i$ can jump over one at another site $K \in \mathcal{T}_i$ into three types; in what follows we let $J = n (j - 1) + a$ and $K = n (k - 1) + b$, so that $k \in \{ j, j + 1 \}$ (since all jumps are of size $n$, that is, $J < K < J + n$) with $a < b$ if $k = j$ and $a > b$ if $k = j + 1$. This jump is of type $1$ if it satisfies $K \in \mathcal{T}_m = \mathfrak{T} (\lambda)$ and $k = j$; of type $2$ if it satisfies $K \in \mathcal{T}_0 = \mathfrak{T} (\mu)$ and $k = j + 1$; and of type $3$ if it satisfies neither of these two conditions. For instance, on the right side of \Cref{sequencearrows}, a jump of the first type comes from the pair $(J, K) = (1, 3)$; a jump of the second type comes from $(J, K) = (5, 7)$; and a jump of the third type comes from $(J, K) = (2, 4)$. 
	
	For fixed $j,  k \ge 1$, jumps of type $1$ are in bijection with pairs of arrows at $(j, 1)$ in $\mathcal{E} (\boldsymbol{\lambda} / \boldsymbol{\mu})$, of distinct colors $a$ and $b$, satisfying the following three conditions. First, $a < b$ (since $j = k$); second, the arrow of color $a$ exits $(j, 1)$ horizontally (since it corresponds to the jumper at site $J = n (j - 1) + a$); and, third, the arrow of color $b$ exits $(j, 1)$ vertically (since $K = n (k - 1) + b \in \mathfrak{T} (\lambda)$). Recalling the definition of $\varphi$ from \eqref{tufunction}, the number of such pairs is given by $\varphi \big( \textbf{D} (j), \textbf{C} (j) \big)$. By similar reasoning, jumps of type $2$ are in bijection with pairs of arrows at $(k, 1)$, of distinct colors $a$ and $b$, such that $a > b$ (since $k = j + 1$); the arrow of color $a$ enters $(k, 1)$ horizontally (since the destination of the jumper at site $J = n (j - 1) + a$ is $n (k - 1) + a$); and the arrow of color $b$ enters $(k, 1)$ vertically (since $K = n (k - 1) + b \in \mathfrak{T} (\mu)$). Thus, there are $\varphi \big( \textbf{A} (k), \textbf{B} (k) \big)$ such pairs. 
	
	Next, we claim that the number of jumps of the third type is twice the number of unordered pairs (that is, subsets of size $2$) of horizontal arrows in $\mathcal{E} (\boldsymbol{\lambda} / \boldsymbol{\mu})$ passing from $(j, 1)$ to $(j + 1, 1)$ of distinct colors $a$ and $b$. Indeed, recalling that $J = n (j - 1) + a$ and $K = n (k - 1) + b$, let us consider the cases $k \in \{ j, j + 1 \}$ individually. If $k = j$, then we must have $K \notin \mathfrak{T} (\lambda)$ for the jump to not be of type $1$. Thus, the particle at site $K$ must eventually jump horizontally in the nonconflicting series $\mathcal{T}_0, \mathcal{T}_1, \ldots , \mathcal{T}_m$, meaning that an arrow of color $b$ must horizontally exit $(j, 1)$ in $\mathcal{E} (\boldsymbol{\lambda} / \boldsymbol{\mu})$. Then, we associate this jump with the unordered pair $\{ a, b \}$ of arrows horizontally exiting $(j, 1)$ (the arrow of color $a$ must horizontally exit this vertex during the jump). If instead $k = j + 1$, then we must have $K \notin \mathfrak{T} (\mu)$ for the jump to not be of type $2$. Thus, the particle at site $K$ must have jumped from site $K - n$ at some previous point in the nonconflicting series $\mathcal{T}_0, \mathcal{T}_1, \ldots , \mathcal{T}_m$, meaning that an arrow of color $b$ must horizontally enter $(k, 1) = (j + 1, 1)$ in $\mathcal{E} (\boldsymbol{\lambda} / \boldsymbol{\mu})$. Then, we again associate this jump with the unordered pair $\{ a, b \}$ of arrows horizontally entering $(j + 1, 1)$ (the arrow of color $a$ must horizontally enter this vertex during the jump). This establishes the claim, so the number of jumps of the third type is given by $2 \binom{|\textbf{D} (j)|}{2}  = 2 \varphi \big( \textbf{D} (j) \big)$. 
	
	Summing over these three types, we deduce that the total number of jumps in the nonconflicting series $\mathcal{T}_1, \mathcal{T}_2, \ldots , \mathcal{T}_m$ is given by
	\begin{flalign*}
	\displaystyle\sum_{j = 1}^{\infty} \Big( 2 \varphi \big( \textbf{D} (j), \textbf{D} (j) \big) + \varphi \big( \textbf{D} (j), \textbf{C} (j) \big) + \varphi \big( \textbf{A} (j), \textbf{B} (j) \big) \Big).
	\end{flalign*}
	
	\noindent This implies the lemma, since $2 \spin (\lambda / \mu)$ is the total number of these jumps.
\end{proof} 

\begin{lem}
	
	\label{psilambdamu}
	
	Recalling $\psi$ from \eqref{lambdamupsi}, we have that 
	\begin{flalign}
	\label{psilambdamuequation} 
	\psi (\boldsymbol{\lambda}) - \psi (\boldsymbol{\mu}) = \displaystyle\frac{1}{2} \displaystyle\sum_{j = 1}^{\infty} \Big( \varphi \big( \textbf{\emph{D}} (j), \textbf{\emph{C}} (j) \big) - \varphi \big( \textbf{\emph{A}} (j), \textbf{\emph{B}} (j) \big) \Big).
	\end{flalign}
\end{lem}

\begin{figure}

	\begin{center}

		\begin{tikzpicture}[
		>=stealth,
		scale = .5
		]

		\draw[red, thick, ->] (-.1, .5) -- (-.1, 2.4) -- (1.9, 2.4) -- (1.9, 4.5);
		\draw[blue, thick, ->] (0, .5) -- (0, 2.5) -- (4, 2.5) -- (4, 4.5);
		\draw[green, thick, ->] (.1, .5) -- (.1, 4.5);
		\draw[red, thick, ->] (3.9, .5) -- (3.9, 2.4) -- (7.9, 2.4) -- (7.9, 4.5);
		\draw[green, thick, ->] (6.1, .5) -- (6.1, 2.6) -- (8.1, 2.6) -- (8.1, 4.5);
		\draw[blue, thick, ->] (6, .5) -- (6, 4.5);
		
		\draw[] (0, .5) circle[radius = 0] node[below, scale = .65]{$\textbf{A}(1)$};
		\draw[] (2, .5) circle[radius = 0] node[below, scale = .65]{$\textbf{A}(2)$};
		\draw[] (4, .5) circle[radius = 0] node[below, scale = .65]{$\textbf{A}(3)$};
		\draw[] (6, .5) circle[radius = 0] node[below, scale = .65]{$\textbf{A}(4)$};
		\draw[] (8, .5) circle[radius = 0] node[below, scale = .65]{$\textbf{A}(5)$};
		
		\draw[] (0, 4.5) circle[radius = 0] node[above, scale = .65]{$\textbf{C}(1)$};
		\draw[] (2, 4.5) circle[radius = 0] node[above, scale = .65]{$\textbf{C}(2)$};
		\draw[] (4, 4.5) circle[radius = 0] node[above, scale = .65]{$\textbf{C}(3)$};
		\draw[] (6, 4.5) circle[radius = 0] node[above, scale = .65]{$\textbf{C}(4)$};
		\draw[] (8, 4.5) circle[radius = 0] node[above, scale = .65]{$\textbf{C}(5)$};
		
		\draw[] (-1.5, 4.5) circle[radius = 0] node[above, scale = .85]{$\boldsymbol{\lambda}$};
		\draw[] (-1.5, .5) circle[radius = 0] node[below, scale = .85]{$\boldsymbol{\mu}$};

		\draw[red, thick, ->] (15.9, .5) -- (15.9, 2.4) -- (17.9, 2.4) -- (17.9, 4.5);
		\draw[blue, thick, ->] (16, .5) -- (16, 2.5) -- (18, 2.5) -- (18, 4.5);
		\draw[green, thick, ->] (16.1, .5) -- (16.1, 4.5);
		\draw[red, thick, ->] (19.9, .5) -- (19.9, 2.4) -- (23.9, 2.4) -- (23.9, 4.5);
		\draw[green, thick, ->] (22.1, .5) -- (22.1, 2.6) -- (24.1, 2.6) -- (24.1, 4.5);
		\draw[blue, thick, ->] (22, .5) -- (22, 4.5);
		
		\draw[] (16, .5) circle[radius = 0] node[below, scale = .65]{$\textbf{A}'(1)$};
		\draw[] (18, .5) circle[radius = 0] node[below, scale = .65]{$\textbf{A}'(2)$};
		\draw[] (20, .5) circle[radius = 0] node[below, scale = .65]{$\textbf{A}'(3)$};
		\draw[] (22, .5) circle[radius = 0] node[below, scale = .65]{$\textbf{A}'(4)$};
		\draw[] (24, .5) circle[radius = 0] node[below, scale = .65]{$\textbf{A}'(5)$};
		
		\draw[] (2, -.5) circle[radius = 0] node[below, scale = .65]{$K_0$};
		\draw[] (4, -.5) circle[radius = 0] node[below, scale = .65]{$K_0 + 1$};
		
		\draw[] (18, -.5) circle[radius = 0] node[below, scale = .65]{$K_0$};
		\draw[] (20, -.5) circle[radius = 0] node[below, scale = .65]{$K_0 + 1$};
		
		\draw[] (16, 4.5) circle[radius = 0] node[above, scale = .65]{$\textbf{C}'(1)$};
		\draw[] (18, 4.5) circle[radius = 0] node[above, scale = .65]{$\textbf{C}'(2)$};
		\draw[] (20, 4.5) circle[radius = 0] node[above, scale = .65]{$\textbf{C}'(3)$};
		\draw[] (22, 4.5) circle[radius = 0] node[above, scale = .65]{$\textbf{C}'(4)$};
		\draw[] (24, 4.5) circle[radius = 0] node[above, scale = .65]{$\textbf{C}'(5)$};
		
		\draw[] (14.5, 4.5) circle[radius = 0] node[above, scale = .85]{$\boldsymbol{\nu}$};
		\draw[] (14.5, .5) circle[radius = 0] node[below, scale = .85]{$\boldsymbol{\mu}$};

		\draw[dashed, blue, ->] (4, 5.5) arc (80:100:6); 
		
		\end{tikzpicture}
		
	\end{center}
	
	\caption{\label{lambdamunumuvertex} Shown above are the single-row path ensembles $\mathcal{E} (\boldsymbol{\lambda} / \boldsymbol{\mu})$ and $\mathcal{E} (\boldsymbol{\nu} / \boldsymbol{\mu})$ used in the proof of \Cref{psilambdamu}.}
\end{figure}
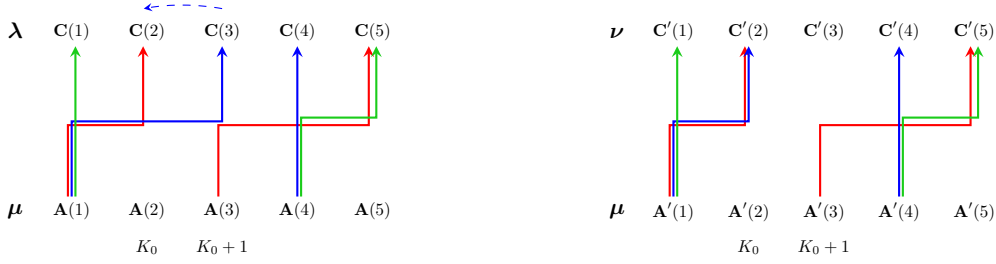

\begin{proof}
	
	We induct on $|\boldsymbol{\lambda}| - |\boldsymbol{\mu}|$. If $|\boldsymbol{\lambda}| = |\boldsymbol{\mu}|$ then $\boldsymbol{\lambda} = \boldsymbol{\mu}$, so $\textbf{B} (j) = \textbf{e}_0 = \textbf{D} (j)$ for each $j \ge 1$, which implies that both sides of \eqref{psilambdamuequation} equal $0$. Thus, let us assume in what follows that \eqref{psilambdamuequation} holds whenever $|\boldsymbol{\lambda}| - |\boldsymbol{\mu}| < m$ for some integer $m > 0$, and we will show it also holds whenever $|\boldsymbol{\lambda}| - |\boldsymbol{\mu}| = m$.
	
	To that end, fix $\boldsymbol{\lambda}, \boldsymbol{\mu} \in \SeqSign_n$ such that $|\boldsymbol{\lambda}| - |\boldsymbol{\mu}| = m$ and $\boldsymbol{\mu} \subseteq \boldsymbol{\lambda}$. Then, there exists a sequence $\boldsymbol{\nu} \in \SeqSign_n$ such that $|\boldsymbol{\nu}| - |\boldsymbol{\mu}| = m - 1$ and $\boldsymbol{\mu} \subseteq \boldsymbol{\nu} \subseteq \boldsymbol{\lambda}$. In particular, there exist indices $h \in [1, n]$ and $k \in [1, M]$ such that $\nu_j^{(i)} = \lambda_j^{(i)}$ for all $(i, j) \ne (h, k)$ and $\nu_k^{(h)} = \lambda_k^{(h)} - 1$. In what follows, we set $K_0 = \nu_k^{(h)} + M - k + 1 \in \mathfrak{T} \big( \nu^{(h)} \big)$, so that $K_0 + 1 = \lambda_k^{(h)} + M - k + 1 \in \mathfrak{T} \big( \lambda^{(h)} \big)$. 
		
	Now let us compare the single-row path ensembles $\mathcal{E} (\boldsymbol{\lambda} / \boldsymbol{\mu})$ and $\mathcal{E} (\boldsymbol{\nu} / \boldsymbol{\mu})$ on $\mathbb{Z}_{> 0} \times \{ 1 \}$. We may interpret $\mathcal{E} (\boldsymbol{\nu} / \boldsymbol{\mu})$ as the ensemble obtained by moving the output of the color $h$ arrow vertically exiting through $(K_0 + 1, 1)$ in $\mathcal{E} (\boldsymbol{\lambda} / \boldsymbol{\mu})$ to the left by one space to $(K_0, 1)$; see \Cref{lambdamunumuvertex} for a depiction, where there $K_0 = 2$ and the output $(K_0 + 1, 1)$ for the color $2$ (blue) arrow is shifted to $(K_0, 1)$. In this way, $\mathcal{E} (\boldsymbol{\lambda} / \boldsymbol{\mu})$ and $\mathcal{E} (\boldsymbol{\nu} / \boldsymbol{\mu})$ coincide everywhere in all colors except for the $h$-th color, and and they also coincide in the $h$-th color at all vertices except for $(K_0, 1)$ and $(K_0 + 1, 1)$. 
	
	To make the latter point more precise, recall for each $j \ge 1$ that $\big( \textbf{A} (j), \textbf{B} (j); \textbf{C} (j), \textbf{D} (j) \big)$ denotes the arrow configuration at $(j, 1)$ under $\mathcal{E} (\boldsymbol{\lambda} / \boldsymbol{\mu})$. Similarly let $\big( \textbf{A}' (j), \textbf{B}' (j); \textbf{C}' (j), \textbf{D}' (j) \big)$ denote the arrow configuration at $(j, 1)$ under $\mathcal{E} (\boldsymbol{\nu} / \boldsymbol{\mu})$. For each index $X \in \{ A, B, C, D, A', B', C', D' \}$, set $\textbf{X} (j) = \big( X_1 (j), X_2 (j), \ldots , X_n (j) \big) \in \{ 0, 1 \}^n$. Then, it is quickly verified that $X_i (j) = X_i' (j)$ for any integers $i \in [1, n]$ and $j \ge 1$ and index $X \in \{ A, B, C, D \}$, unless $i = h$ and $(X, j) \in \big\{ (C, K_0), (D, K_0), (B, K_0 + 1), (C, K_0 + 1) \big\}$. In those exceptional cases, we have that 
	\begin{flalign}
	\label{chdhk0k01}
	\begin{aligned} 
	& C_h (K_0 + 1) = 1 = C_h' (K_0); \qquad C_h (K_0) = 0 = C_h' (K_0 + 1); \\
	& D_h (K_0) = 1 = B_h (K_0 + 1); \qquad D_h' (K_0) = 0 = B_h' (K_0 + 1).
	\end{aligned} 
	\end{flalign}
	
	Now, observe that 
	\begin{flalign}
	\label{psilambdanuidentity} 
	\psi (\boldsymbol{\lambda}) - \psi (\boldsymbol{\nu}) = \displaystyle\frac{1}{2} \big( C_{[h + 1, n]} (K_0) - C_{[1, h - 1]} (K_0 + 1) \big),
	\end{flalign}
	
	\noindent since any contribution to $\psi (\boldsymbol{\lambda})$ not in $\psi (\boldsymbol{\nu})$ in the right side of \eqref{lambdamupsi} comes from some pair $(K_0 + 1, K_0) \in \mathfrak{T} \big( \lambda^{(h)} \big) \times \mathfrak{T} \big( \lambda^{(j)} \big)$ for some $j \in [h + 1, n]$, and any contribution to $\psi (\boldsymbol{\nu})$ not in $\psi (\boldsymbol{\lambda})$ comes from some pair $(K_0 + 1, K_0) \in \mathfrak{T} \big( \nu^{(i)} \big) \times \mathfrak{T} \big( \nu^{(h)} \big) = \mathfrak{T} \big( \lambda^{(i)} \big) \times \mathfrak{T} \big( \nu^{(h)} \big)$ for some $i \in [1, h - 1]$.
	
	Moreover, applying \eqref{psilambdamuequation} for the $(\boldsymbol{\lambda}, \boldsymbol{\mu})$ there equal to $(\boldsymbol{\nu}, \boldsymbol{\mu})$ here, we deduce that
	\begin{flalign*}
	\psi (\boldsymbol{\nu}) - \psi (\boldsymbol{\mu}) & = \displaystyle\frac{1}{2} \displaystyle\sum_{j = 1}^{\infty} \Big( \varphi \big( \textbf{D}' (j), \textbf{C}' (j) \big) - \varphi \big( \textbf{A}' (j), \textbf{B}' (j) \big) \Big) \\
	& = \displaystyle\frac{1}{2} \displaystyle\sum_{j = 1}^{\infty} \Big( \varphi \big( \textbf{D} (j), \textbf{C} (j) \big) - \varphi \big( \textbf{A} (j), \textbf{B} (j) \big) \Big) \\
	& \qquad + \displaystyle\frac{1}{2} \big( D_{[1, h - 1]} (K_0) - C_{[h + 1, n]} (K_0) - D_{[1, h - 1]} (K_0 + 1) + A_{[1, h - 1]} (K_0 + 1) \big), 
	\end{flalign*}
	
	\noindent where in the last equality we used the matching between $\mathcal{E} (\boldsymbol{\lambda} / \boldsymbol{\mu})$ and $\mathcal{E} (\boldsymbol{\nu} / \boldsymbol{\mu})$ at all vertices and indices except those indicated in \eqref{chdhk0k01}, together with the facts that
	\begin{flalign*}
	& \varphi \big( \textbf{D}' (K_0), \textbf{C}' (K_0) \big) = \varphi \big( \textbf{D} (K_0), \textbf{C} (K_0) \big) + D_{[1, h - 1]} (K_0) - C_{[h + 1, n]} (K_0); \\
	& \varphi \big( \textbf{D}' (K_0 + 1), \textbf{C}' (K_0 + 1) \big) = \varphi \big( \textbf{D} (K_0 + 1), \textbf{C} (K_0 + 1) \big) - D_{[1, h - 1]} (K_0 + 1); \\
	& \varphi \big( \textbf{A}' (K_0), \textbf{B}' (K_0) \big) = \varphi \big( \textbf{A} (K_0), \textbf{B} (K_0) \big); \\
	& \varphi \big( \textbf{A}' (K_0 + 1), \textbf{B}' (K_0 + 1) \big) = \varphi \big( \textbf{A} (K_0 + 1), \textbf{B} (K_0 + 1) \big) - A_{[1, h - 1]} (K_0 + 1),
	\end{flalign*}
	
	\noindent which follow from \eqref{chdhk0k01} and the explicit form \eqref{tufunction} for $\varphi$. Combining this with the identities $\textbf{D} (K_0) = \textbf{B} (K_0 + 1)$ and $\textbf{A} (K_0 + 1) + \textbf{B} (K_0 + 1) = \textbf{C} (K_0 + 1) + \textbf{D} (K_0 + 1)$ gives
	\begin{flalign*}
	\psi (\boldsymbol{\nu}) - \psi (\boldsymbol{\mu}) & = \displaystyle\frac{1}{2} \displaystyle\sum_{j = 1}^{\infty} \Big( \varphi \big( \textbf{D} (j), \textbf{C} (j) \big) - \varphi \big( \textbf{A} (j), \textbf{B} (j) \big) \Big) + \displaystyle\frac{1}{2} \big( C_{[1, h - 1]} (K_0 + 1) - C_{[h + 1, n]} (K_0) \big),
	\end{flalign*}
	\noindent which, together with \eqref{psilambdanuidentity}, yields \eqref{psilambdamuequation} and thus the lemma.
\end{proof}

We can now establish the first part of \Cref{limitg0horizontal}. 

\begin{proof}[Proof of Part \ref{lg1} of \Cref{limitg0horizontal}]
	
	Let us first assume that $\lambda / \mu$ is not a horizontal $n$-ribbon strip. Then, any series of $n$-jumps transforming $\mathfrak{T} (\mu)$ into $\mathfrak{T} (\lambda)$ must be conflicting, meaning that in any such series there exists a jumper at some site $J \ge 1$ such that its destination $J + r n \in \mathfrak{T} (\mu)$ (where $r \ge 1$ is some integer); stated alternatively, the destination of this jumper is occupied in $\mathfrak{T} (\mu)$. Letting $J + r n = n (k - 1) + h$ for some integers $k \ge 2$ and $h \in [1, n]$, this means that horizontal and vertical arrows of color $h$ both enter and exit through the vertex $(k, 1)$ in the path ensemble $\mathcal{E} (\boldsymbol{\lambda} / \boldsymbol{\mu})$, that is, $A_h (k) = B_h (k) = C_h (k) = D_h (k) = 1$. 
	
	Due to the factor of $\textbf{1}_{v = 0}$ in the expression \eqref{limitw} for $\mathcal{W}_z (\textbf{A}, \textbf{B}; \textbf{C}, \textbf{D} \boldsymbol{\mid} \infty, 0)$, it follows that the weight of $\mathcal{E} (\boldsymbol{\lambda} / \boldsymbol{\mu})$ is $0$, and so $\mathcal{G}_{\boldsymbol{\lambda} / \boldsymbol{\mu}} (x, \infty \boldsymbol{\mid} 0; 0) = 0$. Similarly, from the definition \eqref{lambdamug} (and \eqref{llambdamu}) for $\mathcal{L}_{\boldsymbol{\lambda} / \boldsymbol{\mu}}$, we also find that $\mathcal{L}_{\boldsymbol{\lambda} / \boldsymbol{\mu}} (x) = 0$ if $\lambda / \mu$ is not a ribbon strip, thereby verifying the proposition in this case. 
	
	So, let us assume that $\lambda / \mu$ is a horizontal ribbon strip. Then, \eqref{lambdamug} implies that
	\begin{flalign}
	\label{xllambdamu} 
	\mathcal{L}_{\boldsymbol{\lambda} / \boldsymbol{\mu}} (x) = \mathcal{L}_{\lambda / \mu}^{(n)} (x) = q^{\spin (\lambda / \mu)} x^{(|\lambda| - |\mu|) / n} = q^{\spin (\lambda / \mu)} x^{|\boldsymbol{\lambda}| - |\boldsymbol{\mu}|}.
	\end{flalign} 
		
	\noindent Additionally, since $\mathcal{G}_{\boldsymbol{\lambda} / \boldsymbol{\mu}} (x, \infty; 0; 0)$ is the weight of $\mathcal{E} (\boldsymbol{\lambda} / \boldsymbol{\mu})$ under $\mathcal{W}_z (\textbf{A}, \textbf{B}; \textbf{C}, \textbf{D} \boldsymbol{\mid} \infty, 0)$, the explicit form from (the fourth statement of) \eqref{limitw} for this weight yields
	\begin{flalign*}
	\mathcal{G}_{\boldsymbol{\lambda} / \boldsymbol{\mu}} (x; \infty \boldsymbol{\mid} 0; 0) = \displaystyle\prod_{j = 1}^{\infty} \mathcal{W}_x \big( \textbf{A} (j), \textbf{B} (j); \textbf{C} (j), \textbf{D} (j) \boldsymbol{\mid} \infty, 0 \big) = \displaystyle\prod_{j = 1}^{\infty} x^{|\textbf{D} (j)|} q^{\varphi (\textbf{D} (j), \textbf{C} (j) + \textbf{D} (j))}.
	\end{flalign*}
	
	\noindent Using \Cref{spinlambdamu}, \Cref{psilambdamu}, and the fact that $\sum_{j = 1}^{\infty} \big| \textbf{D} (j) \big| = |\boldsymbol{\lambda}| - |\boldsymbol{\mu}|$, it follows that
	\begin{flalign*}
	\mathcal{G}_{\boldsymbol{\lambda} / \boldsymbol{\mu}} (x; \infty \boldsymbol{\mid} 0; 0)  = q^{\spin (\lambda / \mu) + \psi (\boldsymbol{\lambda}) - \psi (\boldsymbol{\mu})} x^{|\boldsymbol{\lambda}| - |\boldsymbol{\mu}|},
	\end{flalign*} 
	
	\noindent which together with \eqref{xllambdamu} implies the first part of the proposition.
\end{proof}

\section{LLT Polynomials From \texorpdfstring{$\mathcal{F}_{\boldsymbol{\lambda} / \boldsymbol{\mu}} (\textbf{x}; \infty \boldsymbol{\mid} 0; 0)$}{}}

\label{Fx1FunctionL}

In this section, we establish the second part of \Cref{limitg0horizontal}, namely, the $N = 1$ case of \eqref{1fl}. To that end, as in \Cref{GPolynomials}, recall from \eqref{gfhe} and \Cref{weightesum} that $\mathcal{F}_{\boldsymbol{\lambda} / \boldsymbol{\mu}} (\textbf{x}; \infty \boldsymbol{\mid} 0; 0)$ is the partition function under the $\widehat{\mathcal{W}}_x (\textbf{A}, \textbf{B}; \textbf{C}, \textbf{D} \boldsymbol{\mid} \infty, 0)$ weights from \eqref{limitw2} for the vertex model $\mathfrak{P}_F (\boldsymbol{\lambda} / \boldsymbol{\mu})$ (from \Cref{pgpfph} and depicted in the middle of \Cref{fgpaths}) on the domain $\mathcal{D} = \mathcal{D}_1 = \mathbb{Z}_{> 0} \times \{ 1 \}$. Since this domain has one row, $\mathfrak{P}_F (\boldsymbol{\lambda} / \boldsymbol{\mu})$ only consists in a single path ensemble, which we denote by $\widehat{\mathcal{E}} (\boldsymbol{\lambda} / \boldsymbol{\mu})$. For each $j \ge 1$, let $\big( \widehat{\textbf{A}} (j), \widehat{\textbf{B}} (j); \widehat{\textbf{C}} (j), \widehat{\textbf{D}} (j) \big)$ denote the arrow configuration at the vertex $(j, 1) \in \mathcal{D}$ under $\widehat{\mathcal{E}} (\boldsymbol{\lambda} / \boldsymbol{\mu})$. In particular, there exists a minimal integer $L_0 > 1$ such that $\big( \widehat{\textbf{A}} (K), \widehat{\textbf{B}} (K); \widehat{\textbf{C}} (K), \widehat{\textbf{D}} (K) \big) = \big( \textbf{e}_0, \textbf{e}_{[1, n]}; \textbf{e}_0, \textbf{e}_{[1, n]})$ for all $K \ge L_0$; it is given by 
\begin{flalign}
\label{l0} 
L_0 = \displaystyle\max_{i \in [1, n]} \Big( \max \mathfrak{T} \big( \lambda^{(i)} \big) \Big) + 1.
\end{flalign} 

\noindent See the left side of \Cref{lambdamu1mu2mu3vertex} for an example when $\boldsymbol{\lambda} = \big( (1, 0), (2, 0), (2, 0) \big)$ and $\boldsymbol{\mu} = \big( (1), (2), (0) \big)$. 

Recall that $M$ denotes the length of each signature in $\boldsymbol{\mu}$ and $M + 1$ denotes that of each one in $\boldsymbol{\lambda}$. To show the $N = 1$ case of \eqref{1fl}, we will compare $\widehat{\mathcal{E}} (\boldsymbol{\lambda} / \boldsymbol{\mu})$ with $\mathcal{E} \big( \widetilde{\boldsymbol{\mu}} / \mathring{\boldsymbol{\lambda}})$, for some sequences $\widetilde{\boldsymbol{\mu}}, \mathring{\boldsymbol{\lambda}} \in \SeqSign_{n; M + 1}$, where we recall the path ensemble $\mathcal{E}$ from \Cref{GPolynomials}. In particular, to define these sequences, set
\begin{flalign}
\label{lambda2mu}
\begin{aligned} 
\mathring{\lambda}_j^{(i)} = \lambda_j^{(i)} + 1, \quad & \text{for all $j \in [1, M +1]$}; \\
\widetilde{\mu}_1^{(i)} = L_0 - M - 1, \quad & \text{and} \quad  \widetilde{\mu}_j^{(i)} = \mu_{j - 1}^{(i)} + 1, \quad \text{for all $j \in [2, M + 1]$},
\end{aligned}
\end{flalign}

\noindent for each $i \in [1, n]$. Stated alternatively, $\mathring{\boldsymbol{\lambda}} \in \SeqSign_{n; M + 1}$ is defined by forming each $\mathfrak{T} \big( \mathring{\lambda}^{(i)} \big)$ through shifting every entry in $\mathfrak{T} \big( \lambda^{(i)} \big)$ to the right by one space, and $\widetilde{\boldsymbol{\mu}} \in \SeqSign_{n; M + 1}$ is defined by forming each $\mathfrak{T} \big(\widetilde{\mu}^{(i)} \big)$ through first shifting every entry in $\mathfrak{T} \big( \mu^{(i)} \big)$ to the right by one space, and then appending $L_0$. Observe that defining $\mu^{(i)}$ in this way indeed gives rise to a valid signature (with non-increasing entries), since $\widetilde{\mu}_1^{(i)} = L_0 - M - 1 \ge \lambda_1^{(i)} + 1 \ge \mu_1^{(i)} + 1 = \widetilde{\mu}_2^{(i)}$. For each $j \ge 1$, let $\big( \widetilde{\textbf{A}} (j), \widetilde{\textbf{B}} (j); \widetilde{\textbf{C}} (j), \widetilde{\textbf{D}} (j) \big)$ denote the arrow configuration at $(j, 1) \in \mathcal{D}$ under $\mathcal{E} \big( \widetilde{\boldsymbol{\mu}} / \mathring{\boldsymbol{\lambda}} \big)$. We refer to the middle of \Cref{lambdamu1mu2mu3vertex} for a depiction. 

Observe in particular that the path ensembles $\widehat{\mathcal{E}} (\boldsymbol{\lambda} / \boldsymbol{\mu})$ with $\mathcal{E} \big( \widetilde{\boldsymbol{\mu}} / \mathring{\boldsymbol{\lambda}})$ nearly coincide, except for two differences. First, all vertices in the latter are shifted to the right by one with respect to the former. Second, arrows of all colors exit horizontally through $(L_0 - 1, 1)$ (and through all vertices east of it) in $\widehat{\mathcal{E}} (\boldsymbol{\lambda} / \boldsymbol{\mu})$, while arrows of all colors exit vertically through $(L_0, 1)$ in $\mathcal{E} \big( \widetilde{\boldsymbol{\mu}} / \mathring{\boldsymbol{\lambda}})$.

\begin{figure}

	\begin{center}

		\begin{tikzpicture}[
		>=stealth,
		scale = .42
		]

		\draw[red, thick, ->] (-.1, .5) -- (-.1, 2.4) -- (1.9, 2.4) -- (1.9, 4.5);
		\draw[blue, thick, ->] (0, .5) -- (0, 2.5) -- (4, 2.5) -- (4, 4.5);
		\draw[green, thick, ->] (.1, .5) -- (.1, 4.5);
		\draw[red, thick, ->] (3.9, .5) -- (3.9, 2.4) -- (7.9, 2.4) -- (8, 2.4);
		\draw[blue, thick, ->] (6, .5) -- (6, 2.5) -- (8, 2.5);
		\draw[green, thick, ->] (6.1, .5) -- (6.1, 2.6) -- (8, 2.6);
		
		\draw[] (0, .5) circle[radius = 0] node[below, scale = .65]{$\widehat{\textbf{A}}(1)$};
		\draw[] (2, .5) circle[radius = 0] node[below, scale = .65]{$\widehat{\textbf{A}}(2)$};
		\draw[] (4, .5) circle[radius = 0] node[below, scale = .65]{$\widehat{\textbf{A}}(3)$};
		\draw[] (6, .5) circle[radius = 0] node[below, scale = .65]{$\widehat{\textbf{A}}(4)$};
		
		\draw[] (0, 4.5) circle[radius = 0] node[above, scale = .65]{$\widehat{\textbf{C}}(1)$};
		\draw[] (2, 4.5) circle[radius = 0] node[above, scale = .65]{$\widehat{\textbf{C}}(2)$};
		\draw[] (4, 4.5) circle[radius = 0] node[above, scale = .65]{$\widehat{\textbf{C}}(3)$};
		\draw[] (6, 4.5) circle[radius = 0] node[above, scale = .65]{$\widehat{\textbf{C}}(4)$};
		
		\draw[] (-1.5, 4.5) circle[radius = 0] node[above, scale = .85]{$\boldsymbol{\mu}$};
		\draw[] (-1.5, .5) circle[radius = 0] node[below, scale = .85]{$\boldsymbol{\lambda}$};
		
		\draw[] (0, -.5) circle[radius = 0] node[below, scale = .65]{$1$};
		\draw[] (2, -.5) circle[radius = 0] node[below, scale = .65]{$2$};
		\draw[] (4, -.5) circle[radius = 0] node[below, scale = .65]{$3$};
		\draw[] (6, -.5) circle[radius = 0] node[below, scale = .55]{$L_0 - 1 = 4$};

		\draw[red, thick, ->] (13.9, .5) -- (13.9, 2.4) -- (15.9, 2.4) -- (15.9, 4.5);
		\draw[blue, thick, ->] (14, .5) -- (14, 2.5) -- (18, 2.5) -- (18, 4.5);
		\draw[green, thick, ->] (14.1, .5) -- (14.1, 4.5);
		\draw[red, thick, ->] (17.9, .5) -- (17.9, 2.4) -- (19.9, 2.4) -- (19.9, 4.5);
		\draw[blue, thick, ->] (20, .5) -- (20, 4.5);
		\draw[green, thick, ->] (20.1, .5) -- (20.1, 4.5);
		
		\draw[] (12, .5) circle[radius = 0] node[below, scale = .65]{$\widetilde{\textbf{A}}(1)$};
		\draw[] (14, .5) circle[radius = 0] node[below, scale = .65]{$\widetilde{\textbf{A}}(2)$};
		\draw[] (16, .5) circle[radius = 0] node[below, scale = .65]{$\widetilde{\textbf{A}}(3)$};
		\draw[] (18, .5) circle[radius = 0] node[below, scale = .65]{$\widetilde{\textbf{A}}(4)$};
		\draw[] (20, .5) circle[radius = 0] node[below, scale = .65]{$\widetilde{\textbf{A}}(5)$};
		
		\draw[] (12, 4.5) circle[radius = 0] node[above, scale = .65]{$\widetilde{\textbf{C}}(1)$};
		\draw[] (14, 4.5) circle[radius = 0] node[above, scale = .65]{$\widetilde{\textbf{C}}(2)$};
		\draw[] (16, 4.5) circle[radius = 0] node[above, scale = .65]{$\widetilde{\textbf{C}}(3)$};
		\draw[] (18, 4.5) circle[radius = 0] node[above, scale = .65]{$\widetilde{\textbf{C}}(4)$};
		\draw[] (20, 4.5) circle[radius = 0] node[above, scale = .65]{$\widetilde{\textbf{C}}(5)$};
		
		\draw[] (10.5, 4.5) circle[radius = 0] node[above, scale = .85]{$\widetilde{\boldsymbol{\mu}}$};
		\draw[] (10.5, .5) circle[radius = 0] node[below, scale = .85]{$\mathring{\boldsymbol{\lambda}}$};
		
		\draw[] (12, -.5) circle[radius = 0] node[below, scale = .65]{$1$};
		\draw[] (14, -.5) circle[radius = 0] node[below, scale = .65]{$2$};
		\draw[] (16, -.5) circle[radius = 0] node[below, scale = .65]{$3$};
		\draw[] (18, -.5) circle[radius = 0] node[below, scale = .65]{$4$};
		\draw[] (20, -.5) circle[radius = 0] node[below, scale = .55]{$L_0 = 5$};

		\draw[red, thick, ->] (24.9, .5) -- (24.9, 2.4) -- (26.9, 2.4) -- (26.9, 4.5);
		\draw[blue, thick, ->] (25, .5) -- (25, 2.5) -- (27, 2.5) -- (27, 4.5);
		\draw[green, thick, ->] (25.1, .5) -- (25.1, 2.6) -- (27.1, 2.6) -- (27.1, 4.5);
		\draw[red, thick, ->] (30.9, .5) -- (30.9, 4.5);
		\draw[blue, thick, ->] (33, .5) -- (33, 4.5);
		\draw[green, thick, ->] (29.1, .5) -- (29.1, 2.6) -- (33.1, 2.6) -- (33.1, 4.5);
		
		\draw[] (25, .5) circle[radius = 0] node[below, scale = .65]{$\mathring{\textbf{A}}(1)$};
		\draw[] (27, .5) circle[radius = 0] node[below, scale = .65]{$\mathring{\textbf{A}}(2)$};
		\draw[] (29, .5) circle[radius = 0] node[below, scale = .65]{$\mathring{\textbf{A}}(3)$};
		\draw[] (31, .5) circle[radius = 0] node[below, scale = .65]{$\mathring{\textbf{A}}(4)$};
		\draw[] (33, .5) circle[radius = 0] node[below, scale = .65]{$\mathring{\textbf{A}}(5)$};
		
		\draw[] (25, 4.5) circle[radius = 0] node[above, scale = .65]{$\mathring{\textbf{C}}(1)$};
		\draw[] (27, 4.5) circle[radius = 0] node[above, scale = .65]{$\mathring{\textbf{C}}(2)$};
		\draw[] (29, 4.5) circle[radius = 0] node[above, scale = .65]{$\mathring{\textbf{C}}(3)$};
		\draw[] (31, 4.5) circle[radius = 0] node[above, scale = .65]{$\mathring{\textbf{C}}(4)$};
		\draw[] (33, 4.5) circle[radius = 0] node[above, scale = .65]{$\mathring{\textbf{C}}(5)$};
		
		\draw[] (23.5, .5) circle[radius = 0] node[below, scale = .85]{$\mathring{\boldsymbol{\mu}}$};
		\draw[] (23.5, 4.5) circle[radius = 0] node[above, scale = .85]{$\mathring{\boldsymbol{\lambda}}$};
		
		\draw[] (25, -.5) circle[radius = 0] node[below, scale = .65]{$1$};
		\draw[] (27, -.5) circle[radius = 0] node[below, scale = .65]{$2$};
		\draw[] (29, -.5) circle[radius = 0] node[below, scale = .65]{$3$};
		\draw[] (31, -.5) circle[radius = 0] node[below, scale = .65]{$4$};
		\draw[] (33, -.5) circle[radius = 0] node[below, scale = .55]{$L_0 = 5$};

		\end{tikzpicture}
		
	\end{center}
	
	\caption{\label{lambdamu1mu2mu3vertex} Shown above are the single-row path ensembles $\widehat{\mathcal{E}} (\boldsymbol{\lambda} / \boldsymbol{\mu})$, $\mathcal{E} \big( \widetilde{\boldsymbol{\mu}} / \mathring{\boldsymbol{\lambda}} \big)$, and $\mathcal{E}\big( \mathring{\boldsymbol{\lambda}} / \mathring{\boldsymbol{\mu}} \big)$.}
\end{figure}

For this reason, the partition function of $\widehat{\mathcal{E}} (\boldsymbol{\lambda} / \boldsymbol{\mu})$, under the weights $\widehat{\mathcal{W}}_x (\textbf{A}, \textbf{B}; \textbf{C}, \textbf{D} \boldsymbol{\mid} \infty, 0)$, will coincide with that of $\mathcal{E} \big( \widetilde{\boldsymbol{\mu}} / \mathring{\boldsymbol{\lambda}} \big)$, under the weights $\mathcal{W}_x (\textbf{A}, \textbf{B}; \textbf{C}, \textbf{D} \boldsymbol{\mid} \infty, 0)$, up to a factor coming from the relation \eqref{wabcd2} between $\mathcal{W}_x$ and $\widehat{\mathcal{W}}_x$. As shown in \Cref{GPolynomials}, partition functions of the latter type give rise to a factor of $q^{\spin (\widetilde{\boldsymbol{\mu}} / \mathring{\boldsymbol{\lambda}})}$. Since the LLT polynomial $\mathcal{L}_{\boldsymbol{\lambda} / \boldsymbol{\mu}}$ involves a power of $q^{\spin (\boldsymbol{\lambda} / \boldsymbol{\mu})}$, we must express these spins (in the exponents) through one another other. 

To that end, we have the following lemma, where in the below $\lambda \in \Sign_{n (M + 1)}$ and $\mu \in \Sign_{nM}$ are defined so that their $n$-quotients are $\boldsymbol{\lambda}$ and $\boldsymbol{\mu}$, respectively, and similarly $\mathring{\lambda}, \widetilde{\mu} \in \Sign_{n (M + 1)}$ are defined so that their $n$-quotients are $\mathring{\boldsymbol{\lambda}}$ and $\widetilde{\boldsymbol{\mu}}$, respectively. 

\begin{lem}
	
	\label{spinmulambdamu} 
	
	If $\lambda / \mu$ is a horizontal $n$-ribbon strip, then 
	\begin{flalign*}
	\spin \big( \widetilde{\mu} / \mathring{\lambda} \big) = \spin (\lambda / \mu) + (n - 1) \big( |\boldsymbol{\mu}| - |\boldsymbol{\lambda}| - n \big) + (L_0 - M) \binom{n}{2}.
	\end{flalign*}

\end{lem}

\begin{proof}
	
	Recall from \eqref{lambda2mu} that the sequence $\mathring{\boldsymbol{\lambda}} \in \SeqSign_{n; M + 1}$ was defined from $\boldsymbol{\lambda}$ by increasing every entry of each of its signatures by one. Let us define $\mathring{\boldsymbol{\mu}} \in \SeqSign_{n; M + 1}$ from $\boldsymbol{\mu}$ in a similar way, but where we also append a zero to each of its signatures so that they all have length $M + 1$ (as in $\mathring{\boldsymbol{\lambda}}$). More specifically, for each $i \in [1, n]$, we set
	\begin{flalign}
	\label{2mu} 
	\mathring{\mu}_j^{(i)} = \mu_j^{(i)} + 1, \qquad \text{for each $j \in [1, M]$, and} \quad \mathring{\mu}_{M + 1}^{(i)} = 0.
	\end{flalign}
	
	\noindent Equivalently, it is defined by forming each $\mathfrak{T} \big( \mathring{\mu}^{(j)} \big)$ through first shifting every entry in $\mathfrak{T} \big( \mu^{(j)} \big)$ to the right two spaces, and then appending $1$. We refer to the right side of \Cref{lambdamu1mu2mu3vertex} for a depiction. Let $\mathring{\mu} \in \Sign_{n (M + 1)}$ denote the signature whose $n$-quotient is $\mathring{\boldsymbol{\mu}}$. 
	
	Since $\boldsymbol{\mu} \subseteq \boldsymbol{\lambda}$, we have that $\mathring{\boldsymbol{\mu}} \subseteq \mathring{\boldsymbol{\lambda}}$. So, as $\mathcal{D} = \mathbb{Z}_{> 0} \times \{ 1 \}$ has one only row, there is a unique path ensemble $\mathcal{E} (\mathring{\boldsymbol{\lambda}} / \mathring{\boldsymbol{\mu}})$ on $\mathcal{D}$ in the set $\mathfrak{P}_G (\mathring{\boldsymbol{\lambda}} / \mathring{\boldsymbol{\mu}}; 1)$ from \Cref{pgpfph}. For each $j \ge 1$, let $\big( \mathring{\textbf{A}} (j), \mathring{\textbf{B}} (j); \mathring{\textbf{C}} (j), \mathring{\textbf{D}} (j) \big)$ denote the arrow configuration at the vertex $(j, 1) \in \mathcal{D}$ under $\mathcal{E} (\mathring{\boldsymbol{\lambda}} / \mathring{\boldsymbol{\mu}})$; see the right side of \Cref{lambdamu1mu2mu3vertex} for a depiction. In what follows, the coordinates of any element $\textbf{X} \in \{ 0, 1 \}^n$ will be indexed by $[1, n]$.
		
	Let us show that
	\begin{flalign}
	\label{mu1mu2lambdaspin}
	\spin (\lambda / \mu) + \binom{n}{2} = \spin \big( \mathring{\lambda} / \mathring{\mu}) = \spin \big( \widetilde{\mu} / \mathring{\lambda} \big) + (n - 1) \big( |\mathring{\boldsymbol{\lambda}}| - |\mathring{\boldsymbol{\mu}}| \big) + (M - L_0 + 1) \binom{n}{2},
	\end{flalign}
	
	\noindent from which the lemma would follow since \eqref{lambda2mu} and \eqref{2mu} imply
	\begin{flalign*} 
	|\mathring{\boldsymbol{\lambda}}| = |\boldsymbol{\lambda}| + (M + 1)n; \qquad |\mathring{\boldsymbol{\mu}}| = |\boldsymbol{\mu}| + Mn.
	\end{flalign*}
	
	To establish the first relation in \eqref{mu1mu2lambdaspin}, observe by \eqref{lambda2mu} and \eqref{2mu} that $\mathring{\lambda}$ is obtained from $\lambda$ by increasing each entry by $n$, and $\mathring{\mu}$ is obtained from $\mu$ by first increasing each entry by $n$ and then appending $n$ entries equal to zero; see the left and middle parts of \Cref{lambdamulambdamu} for depictions. Thus, the skew-shape $\mathring{\lambda} / \mathring{\mu}$ is obtained from $\lambda / \mu$ by adding an $n \times n$ block of boxes. Since this block is tiled by $n$ ribbons of shape $1 \times n$, each of which has height $n - 1$, we deduce $\spin (\mathring{\lambda} / \mathring{\mu}) - \spin (\lambda / \mu) = \frac{n (n- 1)}{2}$. This verifies the first statement of \eqref{mu1mu2lambdaspin}, so it suffices to establish the latter.

	\begin{figure}[t]

		\begin{center}

			\begin{tikzpicture}[
			>=stealth,
			scale = .4
			]
					
			\draw[] (0, 6) circle[radius = 0] node[above]{$\lambda / \mu$};
			
			\draw[thick, dotted] (2, 3) -- (0, 3) -- (0, 6) -- (5, 6);
			\draw[ultra thick] (2, 3) -- (3, 3) -- (3, 4) -- (6, 4) -- (6, 6) -- (5, 6) -- (5, 5) -- (2, 5) -- (2, 3);
			
			\draw[dotted] (0, 4) -- (2, 4);
			\draw[dotted] (0, 5) -- (2, 5);
			\draw[dotted] (1, 3) -- (1, 6);
			\draw[dotted] (2, 5) -- (2, 6);
			\draw[dotted] (3, 5) -- (3, 6);
			\draw[dotted] (4, 5) -- (4, 6);
			\draw[] (2, 4) -- (3, 4) -- (3, 5);
			\draw[] (4, 4) -- (4, 5);
			\draw[] (5, 4) -- (5, 5) -- (6, 5);
			
			\draw[] (10, 6) circle[radius = 0] node[above]{$\mathring{\lambda} / \mathring{\mu}$};	
			
			\draw[thick, dotted] (10, 3) -- (10, 6) -- (18, 6);
			\draw[ultra thick] (10, 0) -- (10, 3) -- (13, 3) -- (13, 0) -- (10, 0); 
			\draw[thick, dotted] (13, 3) -- (15, 3); 
			\draw[ultra thick] (15, 3) -- (16, 3) -- (16, 4) -- (19, 4) -- (19, 6) -- (18, 6) -- (18, 5) -- (15, 5) -- (15, 3);
			
			\draw[] (15, 4) -- (16, 4) -- (16, 5);
			\draw[] (17, 4) -- (17, 5);
			\draw[] (18, 4) -- (18, 5) -- (19, 5);
			\draw[] (11, 0) -- (11, 3);
			\draw[] (12, 0) -- (12, 3);
			\draw[] (10, 1) -- (13, 1);
			\draw[] (10, 2) -- (13, 2);

			\draw[dotted] (11, 3) -- (11, 6);
			\draw[dotted] (12, 3) -- (12, 6);
			\draw[dotted] (13, 3) -- (13, 6);
			\draw[dotted] (14, 3) -- (14, 6);
			\draw[dotted] (15, 5) -- (15, 6);
			\draw[dotted] (16, 5) -- (16, 6);
			\draw[dotted] (17, 5) -- (17, 6);
			
			\draw[dotted] (10, 4) -- (15, 4);
			\draw[dotted] (10, 5) -- (15, 5);
			
			\draw[] (23, 6) circle[radius = 0] node[above]{$\widetilde{\mu} / \mathring{\mu}$};
			
			\draw[thick, dotted] (23, 3) -- (23, 6) -- (31, 6);
			\draw[ultra thick] (23, 0) -- (23, 3) -- (28, 3) -- (28, 5) -- (31, 5) -- (31, 6) -- (34, 6) -- (34, 3) -- (31, 3) -- (31, 2) -- (28, 2) -- (28, 0) -- (23, 0);
			
			\draw[] (24, 0) -- (24, 3);
			\draw[] (25, 0) -- (25, 3);
			\draw[] (26, 0) -- (26, 3);
			\draw[] (27, 0) -- (27, 3);
			\draw[] (29, 2) -- (29, 5);
			\draw[] (30, 2) -- (30, 5);
			\draw[] (32, 3) -- (32, 6);
			\draw[] (33, 3) -- (33, 6);
			\draw[] (23, 1) -- (28, 1);
			\draw[] (23, 2) -- (28, 2);
			\draw[] (28, 3) -- (31, 3);
			\draw[] (28, 4) -- (34, 4);
			\draw[] (31, 5) -- (34, 5);
			\draw[] (28, 3) -- (28, 2);
			\draw[] (31, 3) -- (31, 5);
			
			\draw[dotted] (24, 3) -- (24, 6);
			\draw[dotted] (25, 3) -- (25, 6);
			\draw[dotted] (26, 3) -- (26, 6);
			\draw[dotted] (27, 3) -- (27, 6);
			\draw[dotted] (28, 5) -- (28, 6);
			\draw[dotted] (29, 5) -- (29, 6);
			\draw[dotted] (30, 5) -- (30, 6);
			
			\draw[dotted] (23, 4) -- (28, 4);
			\draw[dotted] (23, 5) -- (28, 5);
			
			\end{tikzpicture}
			
		\end{center}
		
		\caption{\label{lambdamulambdamu} Shown to the left, middle, and right are the Young diagrams for skew-shapes $\lambda / \mu$, $\mathring{\lambda} / \mathring{\mu}$, and $\widetilde{\mu} / \mathring{\mu}$ from \Cref{lambdamu1mu2mu3vertex}, respectively, where here $n = 3$.}
	\end{figure}

	To that end, we induct\footnote{Here, we will ignore the relationship \eqref{lambda2mu} between $\boldsymbol{\lambda}$ and $\mathring{\boldsymbol{\lambda}}$, as the former is not present in the second equality of \eqref{mu1mu2lambdaspin}. In doing so, we will establish a more general statement by allowing some entries of signatures in $\mathring{\boldsymbol{\lambda}}$ to be zero, which is in principle not permitted by \eqref{lambda2mu}.} on $|\mathring{\boldsymbol{\lambda}}| - |\mathring{\boldsymbol{\mu}}|$, as in the proof of \Cref{psilambdamu}. If $|\mathring{\boldsymbol{\lambda}}| = |\mathring{\boldsymbol{\mu}}|$, then $\mathring{\lambda} = \mathring{\mu}$, and it suffices to show that 
	\begin{flalign}
	\label{muspinmu}
	\spin \big( \widetilde{\mu} / \mathring{\mu} \big) = \big( L_0 - M - 1 \big) \binom{n}{2}.
	\end{flalign}
	
	To do this, we claim that $\widetilde{\mu} \in \Sign_{(M + 1) n}$ is obtained from $\mathring{\mu}$ by adding a $1 \times n$ ribbon in each of the leftmost $n \big( L_0 - M - 1 \big)$ columns of its Young diagram; see the right side of \Cref{lambdamulambdamu} for a depiction. Since each such ribbon has height $n - 1$, we would then deduce $\spin (\widetilde{\mu} / \mathring{\mu}) = \frac{(n - 1) }{2} n \big(  L_0 - M - 1 \big)$, which implies \eqref{muspinmu}. 
	
	To confirm the claim, we must show that
	\begin{flalign}
	\label{mu1mu2} 
	\widetilde{\mu}_j = n \big( L_0 - M - 1 \big), \quad \text{for each $j \in [1, n]$}; \qquad  \widetilde{\mu}_i = \mathring{\mu}_{i - n}, \quad \text{for each $i \in [n + 1, Mn + n]$.}
	\end{flalign}
	
	\noindent The fact that $\widetilde{\boldsymbol{\mu}}$ is the $n$-quotient of $\widetilde{\mu}$ yields for any $j \in [1, n]$ that the $j$-th largest element in $\mathfrak{T} \big( \widetilde{\mu} \big) = \bigcup_{i = 1}^n \big( n \mathfrak{T} (\widetilde{\mu}^{(i)}) - n + i \big)$ is
	\begin{flalign*} 
	\widetilde{\mu}_j + (M + 1) n - j + 1 = n \max \mathfrak{T} \big( \widetilde{\mu}^{(j)} \big) - j + 1 = n L_0 - j + 1,
	\end{flalign*} 
	
	\noindent where we have also used the fact that \eqref{lambda2mu} implies $\max \mathfrak{T} \big( \widetilde{\mu}^{(i)} \big) = \widetilde{\mu}_1^{(i)} + M + 1 = L_0$ for each $i$. This establishes the first statement of \eqref{mu1mu2}. Further observe by \eqref{lambda2mu} and \eqref{2mu} that $\mathring{\mu}_j^{(i)} = \mu_j^{(i)} + 1 = \widetilde{\mu}_{j - 1}^{(i)}$ for each $j \in [2, M + 1]$, which quickly implies that each entry in $\widetilde{\mu}$ (after removing the $n$ largest ones equal to $L_0 - M - 1$) appears in $\mathring{\mu}$ with the same multiplicity. This implies the second statement of \eqref{mu1mu2}, thereby establishing \eqref{muspinmu} and hence the $|\mathring{\boldsymbol{\lambda}}| = |\mathring{\boldsymbol{\mu}}|$ case of \eqref{mu1mu2lambdaspin}. 
	
	Thus, let us assume in what follows that \eqref{mu1mu2lambdaspin} holds whenever $|\mathring{\boldsymbol{\lambda}}| - |\mathring{\boldsymbol{\mu}}| < m$ for some integer $m > 0$, and we will show it also holds whenever $|\mathring{\boldsymbol{\lambda}}| - |\mathring{\boldsymbol{\mu}}| = m$. So, fix $\mathring{\boldsymbol{\lambda}} \in \SeqSign_{n; M + 1}$ such that $|\mathring{\boldsymbol{\lambda}}| - |\mathring{\boldsymbol{\mu}}| = m$ and $\mathring{\boldsymbol{\mu}} \subseteq \mathring{\boldsymbol{\lambda}} \subseteq \widetilde{\boldsymbol{\mu}}$. Then, as in the proof of \Cref{psilambdamu}, there exists a sequence $\mathring{\boldsymbol{\nu}} \in \SeqSign_{n; M + 1}$ such that $|\mathring{\boldsymbol{\nu}}| - |\mathring{\boldsymbol{\mu}}| = m - 1$ and $\mathring{\boldsymbol{\mu}} \subseteq \mathring{\boldsymbol{\nu}} \subseteq \mathring{\boldsymbol{\lambda}}$. In particular, there exist indices $h \in [1, n]$ and $k \in [1, M]$ such that $\mathring{\nu}_j^{(i)} = \mathring{\lambda}_j^{(i)}$ for all $(i, j) \ne (h, k)$ and $\mathring{\nu}_k^{(h)} = \mathring{\lambda}_k^{(h)} - 1$. In what follows, we set $K_0 = \mathring{\nu}_k^{(h)} + M - k + 2 \in \mathfrak{T} \big( \mathring{\nu}^{(h)} \big)$, so that $K_0 + 1 = \mathring{\lambda}_k^{(h)} + M - k + 2 \in \mathfrak{T} \big( \lambda^{(h)} \big)$. 
	
	Now let us compare the single-row path ensembles in $\big( \mathcal{E} (\mathring{\boldsymbol{\lambda}} / \mathring{\boldsymbol{\mu}}), \mathcal{E} (\mathring{\boldsymbol{\nu}} / \mathring{\boldsymbol{\mu}}) \big)$ and in $\big( \mathcal{E} (\widetilde{\boldsymbol{\mu}} / \mathring{\boldsymbol{\lambda}}), \mathcal{E} (\widetilde{\boldsymbol{\mu}} / \mathring{\boldsymbol{\nu}}) \big)$ on $\mathcal{D}_1$. We may interpret $\mathcal{E} (\mathring{\boldsymbol{\nu}} / \mathring{\boldsymbol{\mu}})$ as the ensemble obtained by moving the output of the color $h$ arrow vertically exiting through $(K_0 + 1, 1)$ in $\mathcal{E} (\mathring{\boldsymbol{\lambda}} / \mathring{\boldsymbol{\mu}})$ to the left by one space to $(K_0, 1)$. Similarly, we may interpret $\mathcal{E} \big( \widetilde{\boldsymbol{\mu}} / \mathring{\boldsymbol{\nu}} \big)$ as the ensemble obtained by moving the input of the color $h$ arrow vertically entering through $(K_0 + 1, 1)$ in $\mathcal{E} \big(\widetilde{\boldsymbol{\mu}} / \mathring{\boldsymbol{\lambda}} \big)$ to the left by one space to $(K_0, 1)$. We refer to \Cref{lambdamue} for depictions in both cases, where there $K_0 = 1$. In this way, the ensembles in $\big( \mathcal{E} (\mathring{\boldsymbol{\lambda}} / \mathring{\boldsymbol{\mu}}), \mathcal{E} (\mathring{\boldsymbol{\nu}} / \mathring{\boldsymbol{\mu}}) \big)$ and in $\big( \mathcal{E} (\widetilde{\boldsymbol{\mu}} / \mathring{\boldsymbol{\lambda}}), \mathcal{E} (\widetilde{\boldsymbol{\mu}} / \mathring{\boldsymbol{\nu}}) \big)$ coincide everywhere in all colors except for the $h$-th color, and and they also coincide in the $h$-th color at all vertices except for $(K_0, 1)$ and $(K_0 + 1, 1)$. 
	
		\begin{figure}

		\begin{center}

			\begin{tikzpicture}[
			>=stealth,
			scale = .42
			]

			\draw[red, thick, ->] (13.9, .5) -- (13.9, 2.4) -- (15.9, 2.4) -- (15.9, 4.5);
			\draw[blue, thick, ->] (14, .5) -- (14, 2.5) -- (18, 2.5) -- (18, 4.5);
			\draw[green, thick, ->] (14.1, .5) -- (14.1, 4.5);
			\draw[red, thick, ->] (17.9, .5) -- (17.9, 2.4) -- (19.9, 2.4) -- (19.9, 4.5);
			\draw[blue, thick, ->] (20, .5) -- (20, 4.5);
			\draw[green, thick, ->] (20.1, .5) -- (20.1, 4.5);
			
			\draw[] (12, .5) circle[radius = 0] node[below, scale = .65]{$\widetilde{\textbf{A}}(1)$};
			\draw[] (14, .5) circle[radius = 0] node[below, scale = .65]{$\widetilde{\textbf{A}}(2)$};
			\draw[] (16, .5) circle[radius = 0] node[below, scale = .65]{$\widetilde{\textbf{A}}(3)$};
			\draw[] (18, .5) circle[radius = 0] node[below, scale = .65]{$\widetilde{\textbf{A}}(4)$};
			\draw[] (20, .5) circle[radius = 0] node[below, scale = .65]{$\widetilde{\textbf{A}}(5)$};
			
			\draw[] (12, 4.5) circle[radius = 0] node[above, scale = .65]{$\widetilde{\textbf{C}}(1)$};
			\draw[] (14, 4.5) circle[radius = 0] node[above, scale = .65]{$\widetilde{\textbf{C}}(2)$};
			\draw[] (16, 4.5) circle[radius = 0] node[above, scale = .65]{$\widetilde{\textbf{C}}(3)$};
			\draw[] (18, 4.5) circle[radius = 0] node[above, scale = .65]{$\widetilde{\textbf{C}}(4)$};
			\draw[] (20, 4.5) circle[radius = 0] node[above, scale = .65]{$\widetilde{\textbf{C}}(5)$};
			
			\draw[] (10.5, 4.5) circle[radius = 0] node[above, scale = .85]{$\widetilde{\boldsymbol{\mu}}$};
			\draw[] (10.5, .5) circle[radius = 0] node[below, scale = .85]{$\mathring{\boldsymbol{\lambda}}$};
			
			\draw[] (12, -.5) circle[radius = 0] node[below, scale = .65]{$K_0$};
			\draw[] (14, -.5) circle[radius = 0] node[below, scale = .65]{$K_0 + 1$};
			\draw[] (20, -.5) circle[radius = 0] node[below, scale = .55]{$L_0$};

			\draw[red, thick, ->] (26.9, .5) -- (26.9, 2.4) -- (30.9, 2.4) -- (30.9, 4.5);
			\draw[blue, thick, ->] (29, .5) -- (29, 2.5) -- (33, 2.5) -- (33, 4.5);
			\draw[green, thick, ->] (29.1, .5) -- (29.1, 4.5);
			\draw[red, thick, ->] (32.9, .5) -- (32.9, 2.4) -- (34.9, 2.4) -- (34.9, 4.5);
			\draw[blue, thick, ->] (35, .5) -- (35, 4.5);
			\draw[green, thick, ->] (35.1, .5) -- (35.1, 4.5);
			
			\draw[] (27, .5) circle[radius = 0] node[below, scale = .65]{$\widetilde{\textbf{A}}' (1)$};
			\draw[] (29, .5) circle[radius = 0] node[below, scale = .65]{$\widetilde{\textbf{A}}' (2)$};
			\draw[] (31, .5) circle[radius = 0] node[below, scale = .65]{$\widetilde{\textbf{A}}' (3)$};
			\draw[] (33, .5) circle[radius = 0] node[below, scale = .65]{$\widetilde{\textbf{A}}' (4)$};
			\draw[] (35, .5) circle[radius = 0] node[below, scale = .65]{$\widetilde{\textbf{A}}' (5)$};
			
			\draw[] (27, 4.5) circle[radius = 0] node[above, scale = .65]{$\widetilde{\textbf{C}}' (1)$};
			\draw[] (29, 4.5) circle[radius = 0] node[above, scale = .65]{$\widetilde{\textbf{C}}' (2)$};
			\draw[] (31, 4.5) circle[radius = 0] node[above, scale = .65]{$\widetilde{\textbf{C}}' (3)$};
			\draw[] (33, 4.5) circle[radius = 0] node[above, scale = .65]{$\widetilde{\textbf{C}}' (4)$};
			\draw[] (35, 4.5) circle[radius = 0] node[above, scale = .65]{$\widetilde{\textbf{C}}' (5)$};
			
			\draw[] (25.5, 4.5) circle[radius = 0] node[above, scale = .85]{$\widetilde{\boldsymbol{\mu}}$};
			\draw[] (25.5, .5) circle[radius = 0] node[below, scale = .85]{$\mathring{\boldsymbol{\nu}}$};
			
			\draw[] (27, -.5) circle[radius = 0] node[below, scale = .65]{$K_0$};
			\draw[] (29, -.5) circle[radius = 0] node[below, scale = .65]{$K_0 + 1$};
			\draw[] (35, -.5) circle[radius = 0] node[below, scale = .55]{$L_0$};

			\draw[dashed, red, ->] (14, 13.75) arc (60:120:2); 
			\draw[dashed, red, ->] (14, -1.5) arc (300:240:2); 
			
			\draw[red, thick, ->] (26.9, 8.5) -- (26.9, 12.5);
			\draw[blue, thick, ->] (27, 8.5) -- (27, 10.5) -- (29, 10.5) -- (29, 12.5);
			\draw[green, thick, ->] (27.1, 8.5) -- (27.1, 10.6) -- (29.1, 10.6) -- (29.1, 12.5);
			\draw[red, thick, ->] (32.9, 8.5) -- (32.9, 12.5);
			\draw[blue, thick, ->] (35, 8.5) -- (35, 12.5);
			\draw[green, thick, ->] (31.1, 8.5) -- (31.1, 10.6) -- (35.1, 10.6) -- (35.1, 12.5);
			
			\draw[] (27, 8.5) circle[radius = 0] node[below, scale = .65]{$\mathring{\textbf{A}}' (1)$};
			\draw[] (29, 8.5) circle[radius = 0] node[below, scale = .65]{$\mathring{\textbf{A}}' (2)$};
			\draw[] (31, 8.5) circle[radius = 0] node[below, scale = .65]{$\mathring{\textbf{A}}' (3)$};
			\draw[] (33, 8.5) circle[radius = 0] node[below, scale = .65]{$\mathring{\textbf{A}}' (4)$};
			\draw[] (35, 8.5) circle[radius = 0] node[below, scale = .65]{$\mathring{\textbf{A}}' (5)$};
			
			\draw[] (27, 12.5) circle[radius = 0] node[above, scale = .65]{$\mathring{\textbf{C}}' (1)$};
			\draw[] (29, 12.5) circle[radius = 0] node[above, scale = .65]{$\mathring{\textbf{C}}' (2)$};
			\draw[] (31, 12.5) circle[radius = 0] node[above, scale = .65]{$\mathring{\textbf{C}}' (3)$};
			\draw[] (33, 12.5) circle[radius = 0] node[above, scale = .65]{$\mathring{\textbf{C}}' (4)$};
			\draw[] (35, 12.5) circle[radius = 0] node[above, scale = .65]{$\mathring{\textbf{C}}' (5)$};
			
			\draw[] (25.5, 8.5) circle[radius = 0] node[below, scale = .85]{$\mathring{\boldsymbol{\mu}}$};
			\draw[] (25.5, 12.5) circle[radius = 0] node[above, scale = .85]{$\mathring{\boldsymbol{\nu}}$};
			
			\draw[] (27, 7.5) circle[radius = 0] node[below, scale = .65]{$K_0$};
			\draw[] (29, 7.5) circle[radius = 0] node[below, scale = .65]{$K_0 + 1$};
			\draw[] (35, 7.5) circle[radius = 0] node[below, scale = .55]{$L_0$};

			\draw[red, thick, ->] (11.9, 8.5) -- (11.9, 10.4) -- (13.9, 10.4) -- (13.9, 12.5);
			\draw[blue, thick, ->] (12, 8.5) -- (12, 10.5) -- (14, 10.5) -- (14, 12.5);
			\draw[green, thick, ->] (12.1, 8.5) -- (12.1, 10.6) -- (14.1, 10.6) -- (14.1, 12.5);
			\draw[red, thick, ->] (17.9, 8.5) -- (17.9, 12.5);
			\draw[blue, thick, ->] (20, 8.5) -- (20, 12.5);
			\draw[green, thick, ->] (16.1, 8.5) -- (16.1, 10.6) -- (20.1, 10.6) -- (20.1, 12.5);
			
			\draw[] (12, 8.5) circle[radius = 0] node[below, scale = .65]{$\mathring{\textbf{A}}(1)$};
			\draw[] (14, 8.5) circle[radius = 0] node[below, scale = .65]{$\mathring{\textbf{A}}(2)$};
			\draw[] (16, 8.5) circle[radius = 0] node[below, scale = .65]{$\mathring{\textbf{A}}(3)$};
			\draw[] (18, 8.5) circle[radius = 0] node[below, scale = .65]{$\mathring{\textbf{A}}(4)$};
			\draw[] (20, 8.5) circle[radius = 0] node[below, scale = .65]{$\mathring{\textbf{A}}(5)$};
			
			\draw[] (12, 12.5) circle[radius = 0] node[above, scale = .65]{$\mathring{\textbf{C}}(1)$};
			\draw[] (14, 12.5) circle[radius = 0] node[above, scale = .65]{$\mathring{\textbf{C}}(2)$};
			\draw[] (16, 12.5) circle[radius = 0] node[above, scale = .65]{$\mathring{\textbf{C}}(3)$};
			\draw[] (18, 12.5) circle[radius = 0] node[above, scale = .65]{$\mathring{\textbf{C}}(4)$};
			\draw[] (20, 12.5) circle[radius = 0] node[above, scale = .65]{$\mathring{\textbf{C}}(5)$};
			
			\draw[] (10.5, 8.5) circle[radius = 0] node[below, scale = .85]{$\mathring{\boldsymbol{\mu}}$};
			\draw[] (10.5, 12.5) circle[radius = 0] node[above, scale = .85]{$\mathring{\boldsymbol{\lambda}}$};
			
			\draw[] (12, 7.5) circle[radius = 0] node[below, scale = .65]{$K_0$};
			\draw[] (14, 7.5) circle[radius = 0] node[below, scale = .65]{$K_0 + 1$};
			\draw[] (20, 7.5) circle[radius = 0] node[below, scale = .55]{$L_0$};

			\end{tikzpicture}
			
		\end{center}
		
		\caption{\label{lambdamue} Shown above are the path ensembles $\mathcal{E} (\mathring{\boldsymbol{\lambda}} / \mathring{\boldsymbol{\mu}})$, $\mathcal{E} (\mathring{\boldsymbol{\nu}} / \mathring{\boldsymbol{\mu}})$, $\mathcal{E}\big( \widetilde{\boldsymbol{\mu}} / \mathring{\boldsymbol{\lambda}} \big)$, and $\mathcal{E}\big( \widetilde{\boldsymbol{\mu}} / \mathring{\boldsymbol{\nu}} \big)$.}
	\end{figure}

	To make this more precise, recall that $\big( \mathring{\textbf{A}} (j), \mathring{\textbf{B}} (j); \mathring{\textbf{C}} (j), \mathring{\textbf{D}} (j) \big)$ and $\big( \widetilde{\textbf{A}} (j), \widetilde{\textbf{B}} (j); \widetilde{\textbf{C}} (j), \widetilde{\textbf{D}} (j) \big)$ are the arrow configurations at any $(j, 1) \in \mathcal{D}_1$ under $\mathcal{E} ( \mathring{\boldsymbol{\lambda}} / \mathring{\boldsymbol{\mu}})$ and $\mathcal{E} \big( \widetilde{\boldsymbol{\mu}} / \mathring{\boldsymbol{\lambda}} \big)$, respectively. For each $j \ge 1$, further let $\big( \mathring{\textbf{A}}' (j), \mathring{\textbf{B}}' (j); \mathring{\textbf{C}}' (j), \mathring{\textbf{D}}' (j) \big)$ and $\big( \widetilde{\textbf{A}}' (j), \widetilde{\textbf{B}}' (j); \widetilde{\textbf{C}}' (j), \widetilde{\textbf{D}}' (j) \big)$ denote the arrow configurations at $(j, 1) \in \mathcal{D}_1$ under $\mathcal{E} ( \mathring{\boldsymbol{\nu}} / \mathring{\boldsymbol{\mu}})$ and $\mathcal{E} \big( \widetilde{\boldsymbol{\mu}} / \mathring{\boldsymbol{\nu}} \big)$, respectively. See \Cref{lambdamue} for depictions.
	
	Then, it is quickly verified that $X_i (j) = X_i' (j)$ for any integers $i \in [1, n]$ and $j \ge 1$ and index $X \in \{ \mathring{A}, \mathring{B}, \mathring{C}, \mathring{D} \}$, unless $i = h$ and $(X, j) \in \big\{ (\mathring{C}, K_0), (\mathring{D}, K_0), (\mathring{B}, K_0 + 1), (\mathring{C}, K_0 + 1) \big\}$. It is further quickly verified that $X_i (j) = X_i' (j)$ for any $i \in [1, n]$, $j \ge 1$, and $X \in \big\{ \widetilde{A}, \widetilde{B}, \widetilde{C}, \widetilde{D} \big\}$, unless $i = h$ and $(X, j) \in \big\{ (\widetilde{A}, K_0), (\widetilde{D}, K_0), (\widetilde{A}, K_0 + 1), (\widetilde{B}, K_0 + 1) \big\}$. In those exceptional cases, we have
	\begin{flalign}
	\label{abcdabcd1}
	\begin{aligned}  
	& \mathring{C}_h (K_0 + 1) = 1 = \mathring{C}_h' (K_0); \qquad \mathring{C}_h (K_0) = 0 = \mathring{C}_h' (K_0 + 1); \\
	& \mathring{D}_h (K_0) = 1 = \mathring{B}_h (K_0 + 1); \qquad \mathring{D}_h' (K_0) = 0 = \mathring{B}_h' (K_0 + 1),
	\end{aligned} 
	\end{flalign} 
	
	\noindent and 
	\begin{flalign}
	\label{abcdabcd2}
	\begin{aligned} 
	& \widetilde{A}_h (K_0 + 1) = 1 = \widetilde{A}_h' (K_0); \qquad \widetilde{A}_h (K_0) = 0 = \widetilde{A}_h' (K_0 + 1); \\
	& \widetilde{D}_h (K_0) = 0 = \widetilde{B}_h (K_0 + 1); \qquad \widetilde{D}_h' (K_0) = 1 = \widetilde{B}_h' (K_0 + 1).
	\end{aligned} 
	\end{flalign} 
	
	\noindent We again refer to \Cref{lambdamue} for a depiction. 
	
	In particular, \eqref{abcdabcd1} and the bilinearity (and definition \eqref{tufunction}) of $\varphi$ together imply that 
	\begin{flalign*} 
	\varphi \big( \mathring{\textbf{D}} (K_0), \mathring{\textbf{D}} (K_0) \big) & = \varphi \big( \mathring{\textbf{D}}' (K_0), \mathring{\textbf{D}}' (K_0) \big) + \varphi \big(\mathring{\textbf{D}} (K_0), \textbf{e}_h \big) + \varphi \big( \textbf{e}_h, \mathring{\textbf{D}} (K_0) \big); \\
	\varphi \big( \mathring{\textbf{A}} (K_0 + 1), \mathring{\textbf{B}} (K_0 + 1) \big) & = \varphi \big( \mathring{\textbf{A}}' (K_0 + 1), \mathring{\textbf{B}}' (K_0 + 1) \big) + \varphi \big(\mathring{\textbf{A}} (K_0 + 1), \textbf{e}_h \big); \\
	\varphi \big( \mathring{\textbf{D}} (K_0), \mathring{\textbf{C}} (K_0) \big) & = \varphi \big( \mathring{\textbf{D}}' (K_0), \mathring{\textbf{C}}' (K_0) \big) - \varphi \big(\mathring{\textbf{D}} (K_0), \textbf{e}_h \big) + \varphi \big( \textbf{e}_h, \mathring{\textbf{C}} (K_0) \big); \\
	\varphi \big( \mathring{\textbf{D}} (K_0 + 1), \mathring{\textbf{C}} (K_0 + 1) \big) & = \varphi \big( \mathring{\textbf{D}}' (K_0 + 1), \mathring{\textbf{C}}' (K_0 + 1) \big) + \varphi \big(\mathring{\textbf{D}} (K_0 + 1), \textbf{e}_h \big).
	\end{flalign*} 
	
	\noindent So, by \Cref{spinlambdamu} we have 
	\begin{flalign}
	\label{spin1}
	\begin{aligned} 
	\spin \big( \mathring{\lambda} / \mathring{\mu} \big) &= \displaystyle\sum_{j = 1}^{L_0 - 1}  \varphi \big( \mathring{\textbf{D}} (j), \mathring{\textbf{D}} (j) \big) + \displaystyle\frac{1}{2} \displaystyle\sum_{j = 2}^{L_0} \varphi \big( \mathring{\textbf{A}} (j), \mathring{\textbf{B}} (j) \big) + \displaystyle\frac{1}{2} \displaystyle\sum_{j = 1}^{L_0 - 1} \varphi \big( \mathring{\textbf{D}} (j), \mathring{\textbf{C}} (j) \big) \\
	& = \displaystyle\sum_{j = 1}^{L_0 - 1}  \varphi \big( \mathring{\textbf{D}}' (j), \mathring{\textbf{D}}' (j) \big) + \displaystyle\frac{1}{2} \displaystyle\sum_{j = 2}^{L_0} \varphi \big( \mathring{\textbf{A}}' (j), \mathring{\textbf{B}}' (j) \big) + \displaystyle\frac{1}{2} \displaystyle\sum_{j = 1}^{L_0 - 1} \varphi \big( \mathring{\textbf{D}}' (j), \mathring{\textbf{C}} (j) \big)  \\
	& \qquad + \displaystyle\frac{1}{2} \Big( \varphi \big( \mathring{\textbf{A}} (K_0 + 1) + \mathring{\textbf{D}} (K_0) + \mathring{\textbf{D}} (K_0 + 1), \textbf{e}_h \big) + \varphi \big( \textbf{e}_h, \mathring{\textbf{C}} (K_0) + 2 \mathring{\textbf{D}} (K_0) \big) \Big) \\
	& = \spin (\mathring{\nu} / \mathring{\mu}) + \displaystyle\frac{1}{2} \varphi \big( \mathring{\textbf{A}} (K_0 + 1) + \mathring{\textbf{D}} (K_0) + \mathring{\textbf{D}} (K_0 + 1), \textbf{e}_h \big) + \displaystyle\frac{1}{2} \varphi \big( \textbf{e}_h, \mathring{\textbf{C}} (K_0) + 2 \mathring{\textbf{D}} (K_0) \big).
	\end{aligned} 
	\end{flalign}
	
	\noindent By similar reasoning, again using \Cref{spinlambdamu}, \eqref{abcdabcd2}, and the blinearity of $\varphi$, we deduce
	\begin{flalign}
	\label{spin2}
	\begin{aligned}
	\spin \big( \widetilde{\mu} / \mathring{\lambda} \big)= \spin \big( \widetilde{\mu} / \mathring{\nu} \big) & + \displaystyle\frac{1}{2} \varphi \big( \textbf{e}_h, \widetilde{\textbf{B}} (K_0 + 1) - \widetilde{\textbf{B}} (K_0) - 2 \widetilde{\textbf{D}} (K_0) - \widetilde{\textbf{C}} (K_0)\big) \\
	& - \displaystyle\frac{1}{2} \varphi \big( \widetilde{\textbf{A}} (K_0 + 1) + 2 \widetilde{\textbf{D}} (K_0), \textbf{e}_h \big).
	\end{aligned} 
	\end{flalign}
	
	\noindent Subtracting \eqref{spin1} from \eqref{spin2} then yields
	\begin{flalign*}	
	\spin \big( \widetilde{\mu} / \mathring{\lambda} \big) - \spin (\mathring{\lambda} / \mathring{\mu}) & = \spin \big( \widetilde{\mu} / \mathring{\nu} \big) - \spin \big( \mathring{\nu} / \mathring{\mu} \big) - \displaystyle\frac{1}{2} \varphi ( \textbf{Y}, \textbf{e}_h) - \displaystyle\frac{1}{2} \varphi (\textbf{e}_h, \textbf{Z}),
	\end{flalign*}
	
	\noindent where 
	\begin{flalign}
	\label{yz} 
	\begin{aligned}
	& \textbf{Y} = \widetilde{\textbf{A}} (K_0 + 1) + 2 \widetilde{\textbf{D}} (K_0) + \mathring{\textbf{A}} (K_0 + 1) + \mathring{\textbf{D}} (K_0) + \mathring{\textbf{D}} (K_0 + 1); \\
	& \textbf{Z} = \widetilde{\textbf{B}} (K_0) + 2 \widetilde{\textbf{D}} (K_0) + \widetilde{\textbf{C}} (K_0) + \mathring{\textbf{C}} (K_0) + 2 \mathring{\textbf{D}} (K_0) - \widetilde{\textbf{B}} (K_0 + 1).
	\end{aligned}
	\end{flalign}
	
	\noindent Now, let us simplify $\textbf{Y}$ and $\textbf{Z}$. To do so, observe from \eqref{lambda2mu} and \eqref{2mu} (see also \Cref{lambdamu1mu2mu3vertex}) that 
	\begin{flalign}
	\label{acac} 
	\widetilde{\textbf{A}} (j) = \mathring{\textbf{C}} (j), \quad \text{for each $j \ge 1$}; \qquad \widetilde{\textbf{C}} (j) = \mathring{\textbf{A}} (j + 1), \quad \text{for each $1 \le j \le L_0 - 1$}.
	\end{flalign} 
	 
	\noindent By spin conservation we also have that
	\begin{flalign}
	\label{bdbd} 
	\begin{aligned}
	& \mathring{\textbf{B}} (j + 1) = \mathring{\textbf{D}} (j) = \displaystyle\sum_{i = 1}^j \big( \mathring{\textbf{A}} (i) - \mathring{\textbf{C}} (i) \big), \quad \text{for each $j \ge 1$}; \\
	& \widetilde{\textbf{B}} (j + 1) = \widetilde{\textbf{D}} (j) = \displaystyle\sum_{i = 1}^j \big( \widetilde{\textbf{A}} (i) - \widetilde{\textbf{C}} (i) \big) =  \displaystyle\sum_{i = 1}^j \big( \mathring{\textbf{C}} (i) - \mathring{\textbf{A}} (i + 1) \big), \quad \text{for each $1 \le j \le L_0 - 1$},
	\end{aligned} 
	\end{flalign}
	
	\noindent where in the last equality we used \eqref{acac}. In particular, \eqref{acac} and \eqref{bdbd} together imply
	\begin{flalign*}
	\mathring{\textbf{D}} (K_0) + \widetilde{\textbf{D}} (K_0) = \mathring{\textbf{A}} (1) - \mathring{\textbf{A}} (K_0 & + 1); \qquad \mathring{\textbf{D}} (K_0 + 1) + \widetilde{\textbf{D}} (K_0) = \mathring{\textbf{A}} (1) - \mathring{\textbf{C}} (K_0 + 1); \\
	& \widetilde{\textbf{A}} (K_0 + 1) = \mathring{\textbf{C}} (K_0 + 1),
	\end{flalign*}
	
	\noindent which upon insertion into \eqref{yz} yields $\textbf{Y} = 2 \mathring{\textbf{A}} (1) = 2 \textbf{e}_{[1, n]}$ (where the last equality follows from the fact that $\mathring{\mu}_{M + 1}^{(i)} = 0$ for each $i \in [1, n]$, by \eqref{2mu}). Similarly, \eqref{acac} and \eqref{bdbd} together imply
	\begin{flalign*}
	\widetilde{\textbf{B}} (K_0 + 1) = \widetilde{\textbf{D}} (K_0); \quad \widetilde{\textbf{B}} (K_0) + \mathring{\textbf{D}} (K_0) = \mathring{\textbf{A}} (1) - \mathring{\textbf{C}} (K_0); \quad \widetilde{\textbf{D}} (K_0) + \mathring{\textbf{D}} (K_0) = \mathring{\textbf{A}} (1) - \mathring{\textbf{C}} (K_0),
	\end{flalign*}
	
	\noindent which upon insertion into \eqref{yz} yields $\textbf{Z} = 2 \mathring{\textbf{A}} (1) = 2 \textbf{e}_{[1, n]}$. 
	
	Thus, $\textbf{Y} = 2 \textbf{e}_{[1, n]} = \textbf{Z}$, and so by the $\textbf{X} = \textbf{e}_h$ case of \eqref{xn1} we have that $\varphi (\textbf{Y}, \textbf{e}_h) + \varphi (\textbf{e}_h, \textbf{Z}) = 2n - 2$. Inserting this into \eqref{spin2} gives
	\begin{flalign*}
	\spin \big( \widetilde{\mu} / \mathring{\lambda} \big) - \spin ( \mathring{\lambda} / \mathring{\mu} ) & = \spin \big( \widetilde{\mu} / \mathring{\nu} \big) - \spin ( \mathring{\nu} / \mathring{\mu}) - n + 1 \\
	& = (n - 1) \big( |\mathring{\boldsymbol{\mu}}| - |\mathring{\boldsymbol{\nu}}| \big) + \big( M - L_0 + 1 \big) \binom{n}{2} - n + 1,
	\end{flalign*}
	
	\noindent where in the last equality we applied the ($\mathring{\boldsymbol{\lambda}} = \mathring{\boldsymbol{\nu}}$ case of) the second statement in \eqref{mu1mu2lambdaspin}. Recalling that $|\mathring{\boldsymbol{\lambda}}| = |\mathring{\boldsymbol{\nu}}| + 1$, it follows that 
	\begin{flalign*}
	\spin \big( \widetilde{\mu} / \mathring{\lambda} \big) - \spin ( \mathring{\lambda} / \mathring{\mu} ) & = (n - 1) \big( |\mathring{\boldsymbol{\mu}}| - |\mathring{\boldsymbol{\lambda}}| \big) + \big( M - L_0 + 1 \big) \binom{n}{2},
	\end{flalign*}
	
	\noindent which establishes the second statement of \eqref{mu1mu2lambdaspin} and therefore the lemma.
\end{proof}

Now we can establish the second part of \Cref{limitg0horizontal}. 

\begin{proof}[Proof of Part \ref{lf1} of \Cref{limitg0horizontal}]

	Let us first assume that $\lambda / \mu$ is not a horizontal $n$-ribbon strip. In this case, recall the signature sequence $\mathring{\boldsymbol{\mu}} \in \SeqSign_{n; M + 1}$ defined in \eqref{2mu}, which is the $n$-quotient of some $\mathring{\mu} \in \Sign_{(M + 1) n}$. Since $\mathring{\boldsymbol{\lambda}}$ is obtained from $\boldsymbol{\lambda}$ by increasing every entry in each of its signatures by one, and since $\mathring{\boldsymbol{\mu}}$ is obtained from $\boldsymbol{\mu}$ in the same way (and by then appending a zero to each signature), it follows that $\mathring{\lambda} / \mathring{\mu}$ is also not a horizontal $n$-strip. 
	
	Now, recall from the proof of \eqref{muspinmu} that $\widetilde{\mu} \in \Sign_{(M + 1) n}$ is obtained by appending a $1 \times n$ ribbon in each of the leftmost $n \big( L_0 - M - 1 \big)$ columns in the Young diagram of $\mathring{\mu}$. In particular, $\widetilde{\mu}$ is the signature with maximal size such that $\widetilde{\mu}_1 \le n \big(L_0 - M - 1 \big)$ and $\widetilde{\mu} / \mathring{\mu}$ is a horizontal $n$-ribbon strip. Consequently, its Maya diagram $\mathfrak{T} (\widetilde{\mu})$ is obtained from $\mathfrak{T} (\mathring{\mu})$ by having each particle in $\mathfrak{T} (\mathring{\mu})$ perform $n$-jumps until colliding with another particle in $\mathfrak{T} (\mathring{\mu})$; stated alternatively, each particle $i \in \mathfrak{T} (\mathring{\mu})$ jumps to site $i + j_i n$, where $j_i \in \mathbb{Z}_{\ge 0}$ is such that $i + n, i + 2n, \ldots , i + j_i n \notin \mathfrak{T} (\mathring{\mu})$ but either $i + (j_i + 1) n \in \mathfrak{T} (\mathring{\mu})$ or $i + (j_i + 1) n > n L_0$. This characterization quickly implies that $\widetilde{\mu} / \nu$ is a horizontal $n$-ribbon strip for some signature $\nu \in \Sign_{n (M + 1)}$ with $\mathring{\mu} \subseteq \nu \subseteq \widetilde{\mu}$ only if $\nu / \mathring{\mu}$ is. Hence, since $\mathring{\lambda} / \mathring{\mu}$ is not a horizontal $n$-ribbon strip, it follows that $\widetilde{\mu} / \mathring{\lambda}$ is not one as well. 
	
	Therefore, as explained in the proof of the first part of \Cref{limitg0horizontal} in \Cref{GPolynomials}, there exists a vertex $(j, 1) \in \mathcal{D}_1$ and a color $h \in [1, n]$ such that horizontal and vertical arrows of color $h$ both enter and exit through $(j, 1)$ in the path ensemble $\mathcal{E} \big( \widetilde{\boldsymbol{\mu}} / \mathring{\boldsymbol{\lambda}} \big)$. It follows that the same statement at $(j - 1, 1) \in \mathcal{D}$ holds in $\widehat{\mathcal{E}} \big( \boldsymbol{\lambda} / \boldsymbol{\mu} \big)$, since the latter is obtained from $\mathcal{E} \big( \widetilde{\boldsymbol{\mu}} / \mathring{\boldsymbol{\lambda}} \big)$ by first shifting each vertex to the left by one space and then modifying the horizontally and vertically exiting arrows at only the vertex $\big( L_0 - 1, 1 \big) \in \mathcal{D}_1$. Hence, $\widehat{A}_h (j - 1) = \widehat{B}_h (j - 1) = \widehat{C}_h (j - 1) = \widehat{D}_h (j - 1) = 1$. 
	
	Due to the factor of $\textbf{1}_{v = 0}$ in the expression \eqref{limitw2} for $\widehat{\mathcal{W}}_z (\textbf{A}, \textbf{B}; \textbf{C}, \textbf{D} \boldsymbol{\mid} \infty, 0)$, it follows that the weight of $\widehat{\mathcal{E}} (\boldsymbol{\lambda} / \boldsymbol{\mu})$ is $0$, and so $\mathcal{F}_{\boldsymbol{\lambda} / \boldsymbol{\mu}} (x; \infty \boldsymbol{\mid} 0; 0) = 0$. Similarly, from the definition \eqref{lambdamug} (and \eqref{llambdamu}) for $\mathcal{L}_{\boldsymbol{\lambda} / \boldsymbol{\mu}}$, we also find that $\mathcal{L}_{\boldsymbol{\lambda} / \boldsymbol{\mu}} (x) = 0$ if $\lambda / \mu$ is not a ribbon strip, thereby verifying the proposition in this case. 
	
	So, let us assume that $\lambda / \mu$ is a horizontal $n$-strip. Then the definition of $\mathcal{F}_{\boldsymbol{\lambda} / \boldsymbol{\mu}} (x; \infty \boldsymbol{\mid} 0; 0)$ as the partition function, under the $\widehat{\mathcal{W}}_z (\textbf{A}, \textbf{B}; \textbf{C}, \textbf{D} \boldsymbol{\mid} \infty, 0)$ weights, for the path ensemble $\widehat{\mathcal{E}} (\boldsymbol{\lambda} / \boldsymbol{\mu})$ (together with the explicit forms \eqref{limitw2} for these weights) gives 
	\begin{flalign}
	\label{fx1} 
	\mathcal{F}_{\boldsymbol{\lambda} / \boldsymbol{\mu}} (x; \infty \boldsymbol{\mid} 0; 0) = \displaystyle\prod_{j = 1}^{L_0 - 1} \widehat{\mathcal{W}}_x \big( \widehat{\textbf{A}} (j), \widehat{\textbf{B}} (j); \widehat{\textbf{C}} (j), \widehat{\textbf{D}} (j) \big) = \displaystyle\prod_{j = 1}^{L_0 - 1} x^{|\widehat{\textbf{D}} (j)| - n} q^{\varphi (\widehat{\textbf{D}} (j), \widehat{\textbf{C}} (j) + \widehat{\textbf{D}} (j)) - \binom{n}{2}},
	\end{flalign}
	
	\noindent where here we have used the fact that $\big( \widehat{\textbf{A}} (j), \widehat{\textbf{B}} (j); \widehat{\textbf{C}} (j), \widehat{\textbf{D}} (j) \big) = \big( \textbf{e}_0, \textbf{e}_{[1, n]}; \textbf{e}_0, \textbf{e}_{[1, n]} \big)$ for $j \ge L_0 - 1$ and the weight under $\widehat{\mathcal{W}}_x$ of this arrow configuration is equal to $1$ by \eqref{wz1} and \eqref{limitw2}.
	
	Now, observe since $\widehat{\textbf{B}} (1) = \textbf{e}_0$ and $\widehat{\textbf{B}} (j) = \textbf{e}_{[1, n]}$ for $j > L_0$ that
	\begin{flalign*}
	\displaystyle\sum_{j = 2}^{L_0} \Big( n - \big| \widehat{\textbf{B}} (j) \big| \Big) + n & = \displaystyle\sum_{j = 1}^{\infty} \Big( n - \big| \widehat{\textbf{B}} (j) \big| \Big) \\
	& = \displaystyle\sum_{i = 1}^n \Bigg( \displaystyle\sum_{\mathfrak{l} \in \mathfrak{T} (\lambda^{(i)})} \mathfrak{l} - \displaystyle\sum_{\mathfrak{m} \in \mathfrak{T} (\mu^{(i)})} \mathfrak{m} \Bigg) \\
	& = \displaystyle\sum_{i = 1}^n \Bigg( \displaystyle\sum_{j = 1}^{M + 1} \big( \lambda_j^{(i)} + M - j + 2 \big)- \displaystyle\sum_{j = 1}^M \big( \mu_j^{(i)} + M - j + 1  \big)\Bigg) \\
	& = |\boldsymbol{\lambda}| - |\boldsymbol{\mu}| + n (M + 1),
	\end{flalign*}
	
	\noindent and so, since $\widehat{\textbf{D}} (j) = \widehat{\textbf{B}} (j + 1)$ for each $j \ge 1$, we have
	\begin{flalign*} 
	\displaystyle\sum_{j = 1}^{L_0 - 1} \Big( n - \big| \widehat{\textbf{D}} (j) \big|  \Big) = |\boldsymbol{\lambda}| - |\boldsymbol{\mu}| + nM.
	\end{flalign*} 
	
	\noindent Inserting this into \eqref{fx1} then yields
	\begin{flalign}
	\label{f1xw}
	\begin{aligned} 
	\mathcal{F}_{\boldsymbol{\lambda} / \boldsymbol{\mu}} (x; \infty \boldsymbol{\mid} 0; 0) & = x^{|\boldsymbol{\mu}| - |\boldsymbol{\lambda}| - nM}  \displaystyle\prod_{j = 2}^{L_0} q^{\varphi (\widehat{\textbf{D}} (j - 1), \widehat{\textbf{C}} (j - 1) + \widehat{\textbf{D}} (j - 1)) - \binom{n}{2}} \\
	& = x^{|\boldsymbol{\mu}| - |\boldsymbol{\lambda}| - nM}  \displaystyle\prod_{j = 2}^{L_0 - 1} q^{\varphi (\widehat{\textbf{D}} (j - 1), \widehat{\textbf{C}} (j - 1) + \widehat{\textbf{D}} (j - 1)) - \binom{n}{2}},
	\end{aligned}
	\end{flalign}
	
	\noindent where in the second equality we used the identity $\varphi \big( \widehat{\textbf{D}} (L_0 - 1), \widehat{\textbf{C}} (L_0 - 1) + \widehat{\textbf{D}} (L_0 - 1) \big) = \binom{n}{2}$, which holds since $\widehat{\textbf{D}} (L_0 - 1) = \textbf{e}_{[1, n]}$ and $\widehat{\textbf{C}} (L_0 - 1) = \textbf{e}_0$. 
	
	 Now, since \eqref{lambda2mu} implies for any indices $X \in \{ A, B, C, D \}$ and $j \in [2, L_0 - 1]$ that $\widehat{\textbf{X}} (j - 1) = \widetilde{\textbf{X}} (j)$, we deduce
	\begin{flalign}
	\label{sum2l0}
	\displaystyle\sum_{j = 2}^{L_0 - 1} \Bigg( \varphi \big( \widehat{\textbf{D}} (j - 1), \widehat{\textbf{C}} (j - 1) + \widehat{\textbf{D}} (j - 1) \big) - \binom{n}{2} \Bigg) = \displaystyle\sum_{j = 1}^{L_0} \varphi \big( \widetilde{\textbf{D}} (j), \widetilde{\textbf{C}} (j) + \widetilde{\textbf{D}} (j) \big) + \big( 2 - L_0 \big) \binom{n}{2},
	\end{flalign} 
	
	\noindent where we have additionally used the facts that $\widetilde{\textbf{D}} (1) = \textbf{e}_0 = \widetilde{\textbf{D}} \big( L_0 \big)$. Moreover, \Cref{spinlambdamu}, \Cref{psilambdamu}, and \Cref{spinmulambdamu} together yield
	\begin{flalign}
	\label{sumlambdamul0}
	\begin{aligned}
	\displaystyle\sum_{j = 1}^{L_0} \varphi \big( \widetilde{\textbf{D}} (j), \widetilde{\textbf{C}} (j) + \widetilde{\textbf{D}} (j) \big) & =  \spin (\widetilde{\mu} / \mathring{\lambda}) + \psi (\widetilde{\boldsymbol{\mu}}) - \psi (\mathring{\boldsymbol{\lambda}}) \\
	& =  \spin (\lambda / \mu) + (n - 1) \big( |\boldsymbol{\mu}| - |\boldsymbol{\lambda}| - n \big) + \big( L_0 - M \big) \binom{n}{2} + \psi (\widetilde{\boldsymbol{\mu}}) - \psi (\mathring{\boldsymbol{\lambda}}).
	\end{aligned}
	\end{flalign} 
	
	\noindent Inserting \eqref{sum2l0} and \eqref{sumlambdamul0} into \eqref{f1xw}, we obtain
	\begin{flalign}
	\label{flambdamux1} 
	\mathcal{F}_{\boldsymbol{\lambda} / \boldsymbol{\mu}} (x; \infty \boldsymbol{\mid} 0; 0) = x^{-nM} (q^{1 - n} x^{-1})^{|\boldsymbol{\lambda}| - |\boldsymbol{\mu}|} q^{\spin (\lambda / \mu) + \psi (\widetilde{\boldsymbol{\mu}}) - \psi (\mathring{\boldsymbol{\lambda}}) - M \binom{n}{2}}.
	\end{flalign}
	
	\noindent Next, the definition \eqref{lambdamupsi} of $\psi$ implies by \eqref{lambda2mu} that 
	\begin{flalign}
	\label{psilambdapsimu} 
	\psi (\mathring{\boldsymbol{\lambda}}) = \psi (\boldsymbol{\lambda}); \qquad \psi \big( \widetilde{\boldsymbol{\mu}} \big) - \psi (\boldsymbol{\mu}) = \frac{M}{2} \binom{n}{2},
	\end{flalign} 
	
	\noindent where the latter holds since the summands $(a, b)$ on the right side of \eqref{lambdamupsi} contributing to $\psi \big( \widetilde{\boldsymbol{\mu}} \big)$ but not $\psi (\boldsymbol{\mu})$ arise when $a = L_0$, and $1 \le i < j \le n$ and $b \in \mathfrak{T} \big( \mu^{(j)} \big)$ are arbitrary. Inserting these into \eqref{flambdamux1} yields
	\begin{flalign*}
	\mathcal{F}_{\boldsymbol{\lambda} / \boldsymbol{\mu}} (x; \infty \boldsymbol{\mid} 0; 0) & = q^{\spin (\lambda / \mu) + \psi (\boldsymbol{\mu}) - M \binom{n}{2} / 2 - \psi (\boldsymbol{\lambda})} x^{-nM} (q^{1 - n} x^{-1})^{|\boldsymbol{\lambda}| - |\boldsymbol{\mu}|} \\
	& = q^{\psi (\boldsymbol{\mu}) - \psi(\boldsymbol{\lambda}) - M \binom{n}{2} / 2} x^{-nM} \mathcal{L}_{\boldsymbol{\lambda} / \boldsymbol{\mu}} (q^{1 - n} x^{-1}),
	\end{flalign*}
	
	\noindent where the last equality follows from the definition \eqref{lambdamug} (and \eqref{llambdamu}) of $\mathcal{L}_{\boldsymbol{\lambda} / \boldsymbol{\mu}}$. This establishes the second part of the proposition.
\end{proof}

\section{LLT Polynomials from \texorpdfstring{$\mathcal{G}_{\boldsymbol{\lambda} / \boldsymbol{\mu}} (0; \textbf{x} \boldsymbol{\mid} 0; 0)$}{}}

\label{Functions2GL}

In this section we establish the third part of \Cref{limitg0horizontal}, namely, that \eqref{2gl} holds for $N = 1$. This will follow from the first part of \Cref{limitg0horizontal} (given by the $N = 1$ case of \eqref{1gl}), after applying a reversal and complementation procedure for the Maya diagrams of $\lambda$ and $\mu$. 

To that end, we recall the unique path ensemble $\mathcal{E} (\boldsymbol{\lambda} / \boldsymbol{\mu}) \in \mathfrak{P}_G (\boldsymbol{\lambda} / \boldsymbol{\mu}; 1)$ from \Cref{GPolynomials}, and the associated arrow configuration $\big( \textbf{A} (j), \textbf{B} (j); \textbf{C} (j), \textbf{D} (j) \big)$ at any vertex $(j, 1) \in \mathcal{D}_1 = \mathbb{Z}_{> 0} \times \{ 1 \}$. Since $\big( \textbf{A} (K), \textbf{B} (K); \textbf{C} (K), \textbf{D} (K) \big) = (\textbf{e}_0, \textbf{e}_0; \textbf{e}_0, \textbf{e}_0)$ for sufficiently large $K$, there exists a maximal integer $L_0 > 0$ for which $\big( \textbf{A} (L_0), \textbf{B} (L_0); \textbf{C} (L_0), \textbf{D} (L_0) \big) \ne \big( \textbf{e}_0, \textbf{e}_0; \textbf{e}_0, \textbf{e}_0 \big)$. 

For any $j \in [1, L_0]$, define the elements $\textbf{A}' (j), \textbf{B}' (j), \textbf{C}' (j), \textbf{D}' (j) \in \{ 0, 1 \}^n$ by
\begin{flalign}
\label{aabbccdd}
\begin{aligned}
& \textbf{A}' (j) = \textbf{e}_{[1, n]} - \overleftarrow{\textbf{A}} (L_0 - j + 1); \qquad \textbf{B}' (j) = \overleftarrow{\textbf{D}} (L_0 - j + 1); \\
& \textbf{C}' (j) = \textbf{e}_{[1, n]} - \overleftarrow{\textbf{C}} (L_0 - j + 1); \qquad \textbf{D}' (j) = \overleftarrow{\textbf{B}} (L_0 - j + 1),
\end{aligned}
\end{flalign}

\noindent and define $\textbf{A}' (j) = \textbf{B}' (j) = \textbf{C}' (j) = \textbf{D}' (j) = \textbf{e}_0$ for each $j > L_0$. Here, we have set $\overleftarrow{\textbf{X}} = (X_n, X_{n - 1}, \ldots , X_1)$\index{X@$\overleftarrow{\mathscr{X}}$; reverse ordering of $\mathscr{X}$} to be the order reversal of any $\textbf{X} = (X_1, X_2, \ldots , X_n) \in \{ 0, 1 \}^n$. Additionally define $\textbf{V}' (j) \in \{ 0, 1 \}^n$ by setting $V_i' (j) = \min \big\{ A_i' (j), B_i' (j), C_i' (j), D_i' (j) \big\}$ for each $i \in [1, n]$. Let $x' (j) = \big| \textbf{X}' (j) \big|$ for each index $X \in \{ A, B, C, D, V \}$. 

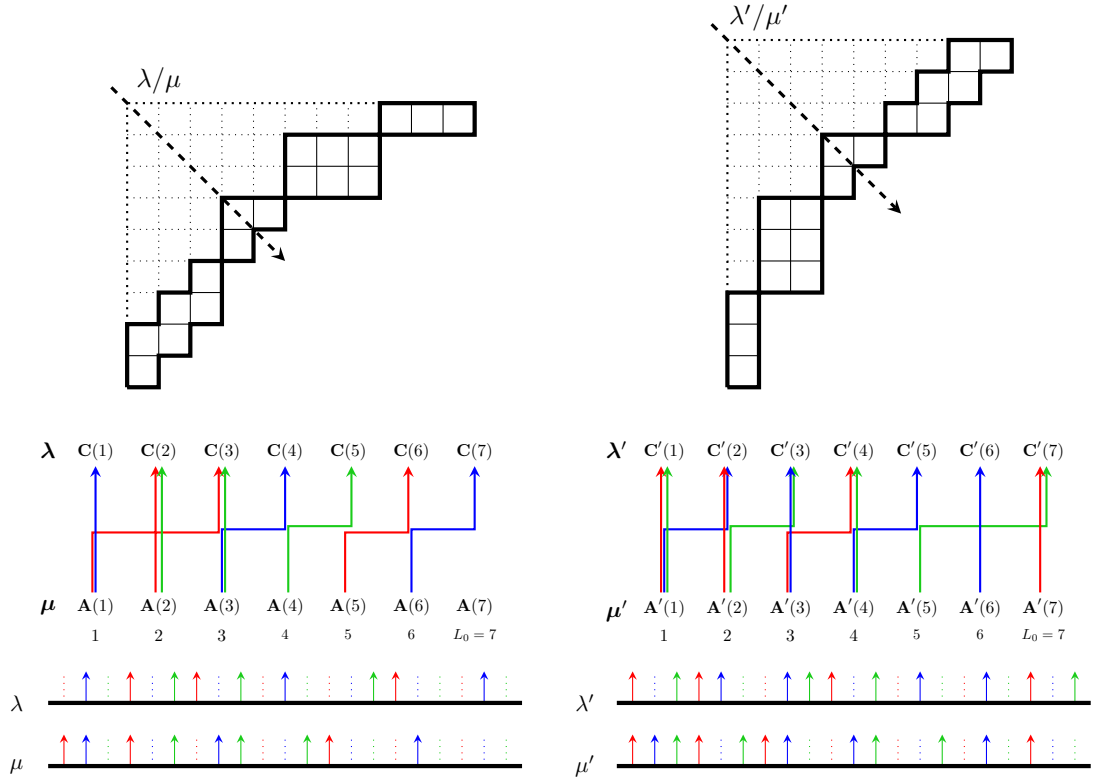
\begin{figure}[t]

	\begin{center}

		\begin{tikzpicture}[
		>=stealth,
		scale = .42
		]

		\draw[red, thick, ->] (-.1, .5) -- (-.1, 2.4) -- (1.9, 2.4) -- (1.9, 4.5);
		\draw[blue, thick, ->] (0, .5) -- (0, 4.5);
		
		\draw[red, thick, ->] (1.9, .5) -- (1.9, 2.4) -- (3.9, 2.4) -- (3.9, 4.5);
		\draw[green, thick, ->] (2.1, .5) -- (2.1, 4.5);
		
		\draw[blue, thick, ->] (4, .5) -- (4, 2.5) -- (6, 2.5) -- (6, 4.5);
		\draw[green, thick, ->] (4.1, .5) -- (4.1, 4.5);
		
		\draw[green, thick, ->] (6.1, .5) -- (6.1, 2.6) -- (8.1, 2.6) -- (8.1, 4.5);
		
		\draw[red, thick, ->] (7.9, .5) -- (7.9, 2.4) -- (9.9, 2.4) -- (9.9, 4.5);
		
		\draw[blue, thick, ->] (10, .5) -- (10, 2.5) -- (12, 2.5) -- (12, 4.5);
		
		\draw[] (0, .5) circle[radius = 0] node[below, scale = .65]{$\textbf{A}(1)$};
		\draw[] (2, .5) circle[radius = 0] node[below, scale = .65]{$\textbf{A}(2)$};
		\draw[] (4, .5) circle[radius = 0] node[below, scale = .65]{$\textbf{A}(3)$};
		\draw[] (6, .5) circle[radius = 0] node[below, scale = .65]{$\textbf{A}(4)$};
		\draw[] (8, .5) circle[radius = 0] node[below, scale = .65]{$\textbf{A}(5)$};
		\draw[] (10, .5) circle[radius = 0] node[below, scale = .65]{$\textbf{A}(6)$};
		\draw[] (12, .5) circle[radius = 0] node[below, scale = .65]{$\textbf{A}(7)$};
		
		\draw[] (0, 4.5) circle[radius = 0] node[above, scale = .65]{$\textbf{C}(1)$};
		\draw[] (2, 4.5) circle[radius = 0] node[above, scale = .65]{$\textbf{C}(2)$};
		\draw[] (4, 4.5) circle[radius = 0] node[above, scale = .65]{$\textbf{C}(3)$};
		\draw[] (6, 4.5) circle[radius = 0] node[above, scale = .65]{$\textbf{C}(4)$};
		\draw[] (8, 4.5) circle[radius = 0] node[above, scale = .65]{$\textbf{C}(5)$};
		\draw[] (10, 4.5) circle[radius = 0] node[above, scale = .65]{$\textbf{C}(6)$};
		\draw[] (12, 4.5) circle[radius = 0] node[above, scale = .65]{$\textbf{C}(7)$};
		
		\draw[] (-1.5, 4.5) circle[radius = 0] node[above, scale = .85]{$\boldsymbol{\lambda}$};
		\draw[] (-1.5, .5) circle[radius = 0] node[below, scale = .85]{$\boldsymbol{\mu}$};
		
		\draw[] (0, -.5) circle[radius = 0] node[below, scale = .65]{$1$};
		\draw[] (2, -.5) circle[radius = 0] node[below, scale = .65]{$2$};
		\draw[] (4, -.5) circle[radius = 0] node[below, scale = .65]{$3$};
		\draw[] (6, -.5) circle[radius = 0] node[below, scale = .55]{$4$};
		\draw[] (8, -.5) circle[radius = 0] node[below, scale = .55]{$5$};
		\draw[] (10, -.5) circle[radius = 0] node[below, scale = .55]{$6$};
		\draw[] (12, -.5) circle[radius = 0] node[below, scale = .55]{$L_0 = 7$};
		
		\draw[] (-2.5, -3) circle[radius = 0] node[scale = .85]{$\lambda$};
		\draw[] (-2.5, -5) circle[radius = 0] node[scale = .85]{$\mu$};

		\draw[] (15.5, -3) circle[radius = 0] node[scale = .85]{$\lambda'$};
		\draw[] (15.5, -5) circle[radius = 0] node[scale = .85]{$\mu'$};
		
		\draw[red, dotted] (-1, -3) -- (-1, -2);
		\draw[blue, ->] (-.3, -3) -- (-.3, -2);
		\draw[green, dotted] (.4, -3) -- (.4, -2);
		\draw[red, ->] (1.1, -3) -- (1.1, -2);
		\draw[blue, dotted] (1.8, -3) -- (1.8, -2);
		\draw[green, ->] (2.5, -3) -- (2.5, -2);
		\draw[red, ->] (3.2, -3) -- (3.2, -2);
		\draw[blue, dotted] (3.9, -3) -- (3.9, -2);
		\draw[green, ->] (4.6, -3) -- (4.6, -2);
		\draw[red, dotted] (5.3, -3) -- (5.3, -2);
		\draw[blue, ->] (6, -3) -- (6, -2);
		\draw[green, dotted] (6.7, -3) -- (6.7, -2);
		\draw[red, dotted] (7.4, -3) -- (7.4, -2);
		\draw[blue, dotted] (8.1, -3) -- (8.1, -2);
		\draw[green, ->] (8.8, -3) -- (8.8, -2);
		\draw[red, ->] (9.5, -3) -- (9.5, -2);
		\draw[blue, dotted] (10.2, -3) -- (10.2, -2);
		\draw[green, dotted] (10.9, -3) -- (10.9, -2);
		\draw[red, dotted] (11.6, -3) -- (11.6, -2);
		\draw[blue, ->] (12.3, -3) -- (12.3, -2);
		\draw[green, dotted] (13, -3) -- (13, -2);

		\draw[red, ->] (-1, -5) -- (-1, -4);
		\draw[blue, ->] (-.3, -5) -- (-.3, -4);
		\draw[green, dotted] (.4, -5) -- (.4, -4);
		\draw[red, ->] (1.1, -5) -- (1.1, -4);
		\draw[blue, dotted] (1.8, -5) -- (1.8, -4);
		\draw[green, ->] (2.5, -5) -- (2.5, -4);
		\draw[red, dotted] (3.2, -5) -- (3.2, -4);
		\draw[blue, ->] (3.9, -5) -- (3.9, -4);
		\draw[green, ->] (4.6, -5) -- (4.6, -4);
		\draw[red, dotted] (5.3, -5) -- (5.3, -4);
		\draw[blue, dotted] (6, -5) -- (6, -4);
		\draw[green, ->] (6.7, -5) -- (6.7, -4);
		\draw[red, ->] (7.4, -5) -- (7.4, -4);
		\draw[blue, dotted] (8.1, -5) -- (8.1, -4);
		\draw[green, dotted] (8.8, -5) -- (8.8, -4);
		\draw[red, dotted] (9.5, -5) -- (9.5, -4);
		\draw[blue, ->] (10.2, -5) -- (10.2, -4);
		\draw[green, dotted] (10.9, -5) -- (10.9, -4);
		\draw[red, dotted] (11.6, -5) -- (11.6, -4);
		\draw[blue, dotted] (12.3, -5) -- (12.3, -4);
		\draw[green, dotted] (13, -5) -- (13, -4);
		
		\draw[ultra thick, black] (-1.5, -3) -- (13.5, -3);
		\draw[ultra thick, black] (-1.5, -5) -- (13.5, -5);

		\draw[red, ->] (17, -3) -- (17, -2);
		\draw[blue, dotted] (17.7, -3) -- (17.7, -2);
		\draw[green, ->] (18.4, -3) -- (18.4, -2);
		\draw[red, ->] (19.1, -3) -- (19.1, -2);
		\draw[blue, ->] (19.8, -3) -- (19.8, -2);
		\draw[green, dotted] (20.5, -3) -- (20.5, -2);
		\draw[red, dotted] (21.2, -3) -- (21.2, -2);
		\draw[blue, ->] (21.9, -3) -- (21.9, -2);
		\draw[green, ->] (22.6, -3) -- (22.6, -2);
		\draw[red, ->] (23.3, -3) -- (23.3, -2);
		\draw[blue, dotted] (24, -3) -- (24, -2);
		\draw[green, ->] (24.7, -3) -- (24.7, -2);
		\draw[red, dotted] (25.4, -3) -- (25.4, -2);
		\draw[blue, ->] (26.1, -3) -- (26.1, -2);
		\draw[green, dotted] (26.8, -3) -- (26.8, -2);
		\draw[red, dotted] (27.5, -3) -- (27.5, -2);
		\draw[blue, ->] (28.2, -3) -- (28.2, -2);
		\draw[green, dotted] (28.9, -3) -- (28.9, -2);
		\draw[red, ->] (29.6, -3) -- (29.6, -2);
		\draw[blue, dotted] (30.3, -3) -- (30.3, -2);
		\draw[green, ->] (31, -3) -- (31, -2);
		
		\draw[red, ->] (17, -5) -- (17, -4);
		\draw[blue, ->] (17.7, -5) -- (17.7, -4);
		\draw[green, ->] (18.4, -5) -- (18.4, -4);
		\draw[red, ->] (19.1, -5) -- (19.1, -4);
		\draw[blue, dotted] (19.8, -5) -- (19.8, -4);
		\draw[green, ->] (20.5, -5) -- (20.5, -4);
		\draw[red, ->] (21.2, -5) -- (21.2, -4);
		\draw[blue, ->] (21.9, -5) -- (21.9, -4);
		\draw[green, dotted] (22.6, -5) -- (22.6, -4);
		\draw[red, dotted] (23.3, -5) -- (23.3, -4);
		\draw[blue, ->] (24, -5) -- (24, -4);
		\draw[green, ->] (24.7, -5) -- (24.7, -4);
		\draw[red, dotted] (25.4, -5) -- (25.4, -4);
		\draw[blue, dotted] (26.1, -5) -- (26.1, -4);
		\draw[green, ->] (26.8, -5) -- (26.8, -4);
		\draw[red, dotted] (27.5, -5) -- (27.5, -4);
		\draw[blue, ->] (28.2, -5) -- (28.2, -4);
		\draw[green, dotted] (28.9, -5) -- (28.9, -4);
		\draw[red, ->] (29.6, -5) -- (29.6, -4);
		\draw[blue, dotted] (30.3, -5) -- (30.3, -4);
		\draw[green, dotted] (31, -5) -- (31, -4);
		
		\draw[] (16.5, .5) circle[radius = 0] node[below, scale = .85]{$\boldsymbol{\mu}'$};
		\draw[] (16.5, 4.5) circle[radius = 0] node[above, scale = .85]{$\boldsymbol{\lambda}'$};

		\draw[red, thick, ->] (17.9, .5) -- (17.9, 4.5);
		\draw[blue, thick, ->] (18, .5) -- (18, 2.5) -- (20, 2.5) -- (20, 4.5);
		\draw[green, thick, ->] (18.1, .5) -- (18.1, 4.5);
		
		\draw[red, thick, ->] (19.9, .5) -- (19.9, 4.5);
		\draw[green, thick, ->] (20.1, .5) -- (20.1, 2.6) -- (22.1, 2.6) -- (22.1, 4.5); 
		
		\draw[red, thick, ->] (21.9, .5) -- (21.9, 2.4) -- (23.9, 2.4) -- (23.9, 4.5);
		\draw[blue, thick, ->] (22, .5) -- (22, 4.5);
		
		\draw[blue, thick, ->] (24, .5) -- (24, 2.5) -- (26, 2.5) -- (26, 4.5);
		\draw[green, thick, ->] (24.1, .5) -- (24.1, 4.5);
		
		\draw[green, thick, ->] (26.1, .5) -- (26.1, 2.6) -- (30.1, 2.6) -- (30.1, 4.5);
		
		\draw[blue, thick, ->] (28, .5) -- (28, 4.5);
		
		\draw[red, thick, ->] (29.9, .5) -- (29.9, 4.5);

		\draw[] (18, .5) circle[radius = 0] node[below, scale = .65]{$\textbf{A}'(1)$};
		\draw[] (20, .5) circle[radius = 0] node[below, scale = .65]{$\textbf{A}'(2)$};
		\draw[] (22, .5) circle[radius = 0] node[below, scale = .65]{$\textbf{A}'(3)$};
		\draw[] (24, .5) circle[radius = 0] node[below, scale = .65]{$\textbf{A}'(4)$};
		\draw[] (26, .5) circle[radius = 0] node[below, scale = .65]{$\textbf{A}'(5)$};
		\draw[] (28, .5) circle[radius = 0] node[below, scale = .65]{$\textbf{A}'(6)$};
		\draw[] (30, .5) circle[radius = 0] node[below, scale = .65]{$\textbf{A}'(7)$};
		
		\draw[] (18, 4.5) circle[radius = 0] node[above, scale = .65]{$\textbf{C}'(1)$};
		\draw[] (20, 4.5) circle[radius = 0] node[above, scale = .65]{$\textbf{C}'(2)$};
		\draw[] (22, 4.5) circle[radius = 0] node[above, scale = .65]{$\textbf{C}'(3)$};
		\draw[] (24, 4.5) circle[radius = 0] node[above, scale = .65]{$\textbf{C}'(4)$};
		\draw[] (26, 4.5) circle[radius = 0] node[above, scale = .65]{$\textbf{C}'(5)$};
		\draw[] (28, 4.5) circle[radius = 0] node[above, scale = .65]{$\textbf{C}'(6)$};
		\draw[] (30, 4.5) circle[radius = 0] node[above, scale = .65]{$\textbf{C}'(7)$};
		
		\draw[] (18, -.5) circle[radius = 0] node[below, scale = .65]{$1$};
		\draw[] (20, -.5) circle[radius = 0] node[below, scale = .65]{$2$};
		\draw[] (22, -.5) circle[radius = 0] node[below, scale = .65]{$3$};
		\draw[] (24, -.5) circle[radius = 0] node[below, scale = .65]{$4$};
		\draw[] (26, -.5) circle[radius = 0] node[below, scale = .55]{$5$};
		\draw[] (28, -.5) circle[radius = 0] node[below, scale = .55]{$6$};
		\draw[] (30, -.5) circle[radius = 0] node[below, scale = .55]{$L_0 = 7$};
		
		\draw[ultra thick, black] (16.5, -3) -- (31.5, -3);
		\draw[ultra thick, black] (16.5, -5) -- (31.5, -5);

		\draw[thick, dotted] (20, 10) -- (20, 18) -- (27, 18);
		\draw[ultra thick] (20, 7) -- (21, 7) -- (21, 10) -- (23, 10) --(23, 13) -- (24, 13) -- (24, 14) -- (25, 14) -- (25, 15) -- (27, 15) -- (27, 16) -- (28, 16) -- (28, 17) -- (29, 17) -- (29, 18) -- (27, 18) -- (27, 17) -- (26, 17) -- (26, 16) -- (25, 16) -- (25, 15) -- (23, 15) -- (23, 13) -- (21, 13) -- (21, 10) -- (20, 10) -- (20, 7);
		
		\draw[] (20, 8) -- (21, 8);
		\draw[] (20, 9) -- (21, 9);
		\draw[] (21, 11) -- (23, 11);
		\draw[] (21, 12) -- (23, 12);
		\draw[] (22, 10) -- (22, 13);
		\draw[] (23, 14) -- (24, 14) -- (24, 15);
		\draw[] (26, 15) -- (26, 16) -- (27, 16) -- (27, 17) -- (28, 17) -- (28, 18);
		
		\draw[dotted] (20, 11) -- (21, 11);
		\draw[dotted] (20, 12) -- (21, 12);
		\draw[dotted] (20, 13) -- (21, 13) -- (21, 18);
		\draw[dotted] (22, 13) -- (22, 18);
		\draw[dotted] (20, 14) -- (23, 14);
		\draw[dotted] (20, 15) -- (23, 15) -- (23, 18);
		\draw[dotted] (24, 15) -- (24, 18);
		\draw[dotted] (20, 16) -- (25, 16) -- (25, 18);
		\draw[dotted] (20, 17) -- (26, 17) -- (26, 18);

		\draw[] (21, 18) circle[radius = 0] node[above]{$\lambda' / \mu'$};

		\draw[thick, dotted] (1, 9) -- (1, 16) -- (9, 16);
		\draw[ultra thick] (1, 7) -- (1, 9) -- (2, 9) -- (2, 10) -- (3, 10) -- (3, 11) -- (4, 11) -- (4, 13) -- (6, 13) -- (6, 15) -- (9, 15) -- (9, 16) -- (12, 16) -- (12, 15) -- (9, 15) -- (9, 13) -- (6, 13) -- (6, 12) -- (5, 12) -- (5, 11) -- (4, 11) -- (4, 9) -- (3, 9) -- (3, 8) -- (2, 8) -- (2, 7) -- (1, 7);
		\draw[] (1, 8) -- (2, 8) -- (2, 9) -- (3, 9) -- (3, 10) -- (4, 10); 
		\draw[] (4, 12) -- (5, 12) -- (5, 13);
		\draw[] (6, 14) -- (9, 14);
		\draw[] (7, 13) -- (7, 15);
		\draw[] (8, 13) -- (8, 15);
		\draw[] (10, 15) -- (10, 16);
		\draw[] (11, 15) -- (11, 16);
		
		\draw[dotted] (2, 10) -- (2, 16);
		\draw[dotted] (3, 11) -- (3, 16);
		\draw[dotted] (4, 13) -- (4, 16);
		\draw[dotted] (5, 13) -- (5, 16);
		\draw[dotted] (6, 15) -- (6, 16);
		\draw[dotted] (7, 15) -- (7, 16);
		\draw[dotted] (8, 15) -- (8, 16);
		
		\draw[dotted] (1, 10) -- (2, 10);
		\draw[dotted] (1, 11) -- (3, 11);
		\draw[dotted] (1, 12) -- (4, 12);
		\draw[dotted] (1, 13) -- (4, 13);
		\draw[dotted] (1, 14) -- (6, 14);
		\draw[dotted] (1, 15) -- (6, 15);
		
		\draw[very thick,->, dashed] (.5, 16.5) -- (6, 11);
		
		\draw[very thick, ->, dashed] (19.5, 18.5) -- (25.5, 12.5);

		\draw[] (2, 16) circle[radius = 0] node[above]{$\lambda/ \mu$};	
		
		\end{tikzpicture}
		
	\end{center}
	
	\caption{\label{lambdamulambdamu1} Shown above are the Young diagrams for skew-shapes $\lambda / \mu$ and $\lambda' / \mu'$; shown to the bottom are their colored Maya diagrams; and shown to the middle are the single-row path ensembles $\mathcal{E} (\boldsymbol{\lambda} / \boldsymbol{\mu})$ and $\mathcal{E}' (\boldsymbol{\lambda} / \boldsymbol{\mu}) = \mathcal{E} (\boldsymbol{\lambda}' / \boldsymbol{\mu}')$. Here, red is color $1$, blue is color $2$, and green is color $3$.}
\end{figure}

Then, since $\textbf{A} (j) + \textbf{B} (j) = \textbf{C} (j) + \textbf{D} (j)$, \eqref{aabbccdd} implies $\textbf{A}' (j) +  \textbf{B}' (j) = \textbf{C}' (j) + \textbf{D}' (j)$ for each $j > 0$. In particular, this arrow conservation implies the existence of a single-row path ensemble $\mathcal{E}' (\boldsymbol{\lambda} / \boldsymbol{\mu})$ on $\mathbb{Z}_{> 0} \times \{ 1 \}$ whose arrow configuration at $(j, 1)$ is $\big( \textbf{A}' (j), \textbf{B}' (j); \textbf{C}' (j), \textbf{D}' (j) \big)$. In this way, $\mathcal{E}' (\boldsymbol{\lambda} / \boldsymbol{\mu})$ is obtained by first reversing the vertices and colors in $\mathcal{E} (\boldsymbol{\lambda} / \boldsymbol{\mu})$ and then \emph{complementing} all vertical arrows, that is, interchanging a particle of any color with its absence. We refer to the middle of \Cref{lambdamulambdamu1} for a depiction.

As explained in \Cref{GPolynomials}, one may decompress the top and bottom boundaries of $\mathcal{E}' (\boldsymbol{\lambda} / \boldsymbol{\mu})$ to form two colored Maya diagrams. By \eqref{aabbccdd}, the decompression of the top boundary is obtained by first reversing $\mathfrak{T} (\lambda)$ and then complementing it (interchanging particles with empty sites); applying the same reversal and complementation procedure to $\mathfrak{T} (\mu)$ yields the bottom boundary of this decompression. We refer to the bottom part of \Cref{lambdamulambdamu1} for a depiction. As such, it is quickly verified (see, for example, (1.7) of \cite{SFP}) that the top and bottom boundaries of this decompressed vertex model are given by $\mathfrak{T} (\lambda')$ and $\mathfrak{T} (\mu')$ (namely, the colored Maya diagrams for the duals of $\lambda$ and $\mu$), respectively. Hence $\mathcal{E}' (\boldsymbol{\lambda} / \boldsymbol{\mu}) = \mathcal{E} (\boldsymbol{\lambda}' / \boldsymbol{\mu}')$, the unique path ensemble in $\mathfrak{P}_G (\boldsymbol{\lambda}' / \boldsymbol{\mu}'; 1)$. 

Now, by \eqref{gfhe} and \Cref{weightesum}, $\mathcal{G}_{\boldsymbol{\lambda} / \boldsymbol{\mu}} (0; x \boldsymbol{\mid} 0; 0)$ is the partition function for $\mathcal{E} (\boldsymbol{\lambda} / \boldsymbol{\mu})$ under the weights $\mathcal{W}_x (\textbf{A}, \textbf{B}; \textbf{C}, \textbf{D} \boldsymbol{\mid} 0, 0)$. Hence, using the explicit form \eqref{limitw} for these weights, we obtain
\begin{flalign*}
\mathcal{G}_{\boldsymbol{\lambda} / \boldsymbol{\mu}} (0; x \boldsymbol{\mid} 0; 0) & = \displaystyle\prod_{j = 1}^{L_0} \mathcal{W}_x \big( \textbf{A} (j), \textbf{B} (j); \textbf{C} (j), \textbf{D} (j) \boldsymbol{\mid} 0, 0 \big) \\
& = \displaystyle\prod_{j = 1}^{L_0} x^{-d (j)} q^{\varphi (\textbf{D} (j), \textbf{C} (j) - \textbf{B} (j))} \textbf{1}_{\textbf{C} (j) \ge \textbf{B} (j)}.
\end{flalign*} 

\noindent Then \eqref{aabbccdd} and the fact that $\sum_{j = 1}^{L_0} d(j) = \sum_{j = 1}^{\infty} d (j) = |\boldsymbol{\lambda}| - |\boldsymbol{\mu}|$ together yield
\begin{flalign*}
\mathcal{G}_{\boldsymbol{\lambda} / \boldsymbol{\mu}} (0; x \boldsymbol{\mid} 0; 0) = x^{|\boldsymbol{\mu}| - |\boldsymbol{\lambda}|} \displaystyle\prod_{j = 1}^{L_0} q^{\varphi (\overleftarrow{\textbf{B}}' (j), \textbf{e}_{[1, n]} - \overleftarrow{\textbf{C}'} (j) - \overleftarrow{\textbf{D}}' (j))} \textbf{1}_{\textbf{e}_{{[1, n]}} \ge \overleftarrow{\textbf{C}}' (j) + \overleftarrow{\textbf{D}}' (j)}. 
\end{flalign*}

\noindent  Thus, since the definition \eqref{tufunction} of $\varphi$ implies
\begin{flalign} 
\label{xyxy} 
\varphi (\overleftarrow{\textbf{X}}, \overleftarrow{\textbf{Y}}) = \varphi (\textbf{Y}, \textbf{X}), \qquad \text{for any $\textbf{X}, \textbf{Y} \in \{ 0, 1 \}^n$},
\end{flalign} 

\noindent we find that
\begin{flalign} 
\label{gproductq} 
\begin{aligned}
\mathcal{G}_{\boldsymbol{\lambda} / \boldsymbol{\mu}} (0; \textbf{x} \boldsymbol{\mid} 0; 0) = x^{|\boldsymbol{\mu}| - |\boldsymbol{\lambda}|} \displaystyle\prod_{j = 1}^{L_0} q^{\varphi (\textbf{e}_{[1, n]} - \textbf{C}' (j) - \textbf{D}' (j), \textbf{B}' (j))} \textbf{1}_{v' (j) = 0},
\end{aligned} 
\end{flalign}

\noindent where we have also used the fact that $\textbf{e}_{[1, n]} \ge \overleftarrow{\textbf{C}} (j) + \overleftarrow{\textbf{D}} (j)$ holds if and only if $v' (j) = 0$, for any $j \ge 1$. Indeed, since $A_i' (j), B_i' (j), C_i' (j), D_i' (j) \in \{ 0, 1 \}$, this follows from the fact that $C_i' (j) + D_i' (j) \le 1$ holds for some $i \in [1, n]$ if and only if $V_i' (j) = 0$ does (using arrow conservation). 

We now require the following lemma to analyze the power of $q$ on the right side of \eqref{gproductq}. 

\begin{lem}
	
\label{sumeb} 

If $\lambda' / \mu'$ is a horizontal $n$-ribbon strip, then 
\begin{flalign*}
\psi (\boldsymbol{\lambda}) - \psi (\boldsymbol{\lambda'}) - \psi (\boldsymbol{\mu}) + \psi (\boldsymbol{\mu}') - \displaystyle\sum_{j = 1}^{L_0} \varphi \big( \textbf{\emph{e}}_{[1, n]}, \textbf{\emph{B}}' (j) \big) = \displaystyle\frac{n - 1}{2} \big( |\boldsymbol{\mu}| - |\boldsymbol{\lambda}| \big).
\end{flalign*}
	
\end{lem} 

\begin{proof}
	
	First observe using \eqref{aabbccdd}, \eqref{xyxy}, and the definition \eqref{tufunction} of $\varphi$ that 
	\begin{flalign}
	\label{sume1nb}
	\begin{aligned} 
	\displaystyle\sum_{j = 1}^{L_0} \varphi \big( \textbf{e}_{[1, n]}, \textbf{B}' (j) \big) = \displaystyle\sum_{j = 1}^{L_0} \varphi \big( \textbf{D} (j), \textbf{e}_{[1, n]} \big) & = \displaystyle\sum_{i = 1}^n (n - i) \displaystyle\sum_{j = 1}^{L_0} D_i (j) = \displaystyle\sum_{i = 1}^n (n - i) \Big( \big| \lambda^{(i)} \big| - \big| \mu^{(i)} \big| \Big),
	\end{aligned} 
	\end{flalign}
	
	\noindent where in the last equality we used the identity $\sum_{j = 1}^{L_0} D_i (j) = \sum_{j = 1}^{\infty} D_i (j) = \big| \lambda^{(i)} \big| - \big| \mu^{(i)} \big|$ for any $i \in [1, n]$. Now, we claim that
	\begin{flalign}
	\label{psilambdalambdasum}
	\begin{aligned} 
	& \psi \big( \boldsymbol{\lambda} \big) - \psi (\boldsymbol{\lambda}') + \displaystyle\sum_{i = 1}^n (i - 1) \big| \lambda^{(i)} \big|  = \displaystyle\frac{(n - 1) |\boldsymbol{\lambda}|}{2} + \displaystyle\frac{1}{2} \binom{n}{2} \Bigg( \binom{M}{2} - \binom{L_0 - M}{2} \Bigg); \\
	& \psi \big( \boldsymbol{\mu} \big) - \psi (\boldsymbol{\mu}') + \displaystyle\sum_{i = 1}^n (i - 1) \big| \mu^{(i)} \big|  = \displaystyle\frac{(n - 1) |\boldsymbol{\mu}|}{2} + \displaystyle\frac{1}{2} \binom{n}{2} \Bigg( \binom{M}{2} - \binom{L_0 - M}{2} \Bigg).
	\end{aligned}
	\end{flalign}
	
	\noindent Let us show the lemma assuming \eqref{psilambdalambdasum}. To that end, subtracting the two equalities in \eqref{psilambdalambdasum} yields
	\begin{flalign*}
	& \psi \big( \boldsymbol{\lambda} \big) - \psi (\boldsymbol{\lambda}') - \psi (\boldsymbol{\mu}) + \psi (\boldsymbol{\mu}') + \displaystyle\sum_{i = 1}^n (i - 1) \Big( \big| \lambda^{(i)} \big| - \big| \mu^{(i)} \big| \Big) = \displaystyle\frac{(n - 1)}{2} \big( |\boldsymbol{\lambda}| - | \boldsymbol{\mu}| \big),
	\end{flalign*}
	
	\noindent and so subtracting $\sum_{i = 1}^n (n - 1) \big( \big| \lambda^{(i)} \big| - \big| \mu^{(i)} \big| \big) = (n - 1) \big( |\boldsymbol{\lambda}| - |\boldsymbol{\mu}| \big)$ yields
	\begin{flalign*}
	& \psi \big( \boldsymbol{\lambda} \big) - \psi (\boldsymbol{\lambda}') - \psi (\boldsymbol{\mu}) + \psi (\boldsymbol{\mu}') + \displaystyle\sum_{i = 1}^n (i - n) \Big( \big| \lambda^{(i)} \big| - \big| \mu^{(i)} \big| \Big) = \displaystyle\frac{(n - 1)}{2} \big( |\boldsymbol{\mu}| - | \boldsymbol{\lambda}| \big),
	\end{flalign*}
	
	\noindent which implies the lemma by \eqref{sume1nb}. Thus, it remains to establish \eqref{psilambdalambdasum}; we will only verify the first statement there, as the proof of the latter is entirely analogous. 
	
	To that end, we induct on $|\boldsymbol{\lambda}|$, as in the proof of \Cref{psilambdamu}. Observe that if $\boldsymbol{\lambda} = \boldsymbol{0}^M$ then \eqref{psi0} imples that $\psi (\boldsymbol{\lambda}) = \frac{1}{2} \binom{n}{2} \binom{M}{2}$ and $\psi (\boldsymbol{\lambda}') = \frac{1}{2} \binom{n}{2} \binom{L_0 - M}{2}$, which verifies \eqref{psilambdalambdasum}. Thus, let us assume \eqref{psilambdalambdasum} holds whenever $|\boldsymbol{\lambda}| < m$ for some integer $m > 0$, and we will show it holds for $|\boldsymbol{\lambda}| = m$.
	
	So, fix $\boldsymbol{\lambda} \in \SeqSign_{n; M}$ with $|\boldsymbol{\lambda}| = m$. Then, there exists a sequence of $n$ signatures $\boldsymbol{\nu} \in \SeqSign_{n; M}$ such that $|\boldsymbol{\nu}| = m - 1$ and the following holds. There exists integers $h \in [1, n]$ and $k \in [1, M]$ such that $\nu_j^{(i)} = \lambda_j^{(i)}$ whenever $(i, j) \ne (h, k)$ and $\nu_k^{(h)} = \lambda_k^{(h)} - 1$. Stated alternatively, $\boldsymbol{\lambda}$ and $\boldsymbol{\nu}$ coincide in every entry of every component, except for in one entry which is smaller by one in $\boldsymbol{\nu}$ than in $\boldsymbol{\lambda}$. 
	
	Letting $K_0 = \nu_k^{(h)} + M - k + 1 \in \mathfrak{T} \big( \nu^{(h)} \big)$, so that $K_0 + 1 = \lambda_k^{(h)} + M - k + 1 \in \mathfrak{T} \big( \lambda^{(h)} \big)$, we have
	\begin{flalign}
	\label{lambdanu}
	\psi (\boldsymbol{\lambda}) = \psi (\boldsymbol{\nu}) + \displaystyle\frac{1}{2} \displaystyle\sum_{i = h + 1}^n \textbf{1}_{K_0 \in \mathfrak{T} (\lambda^{(i)})} - \displaystyle\frac{1}{2} \displaystyle\sum_{i = 1}^{h - 1} \textbf{1}_{K_0 + 1 \in \mathfrak{T} (\lambda^{(i)})}.
	\end{flalign} 
	
	\noindent Here, the first sum on the right side arises since each appearance of $K_0$ in some $\mathfrak{T} \big( \lambda^{(i)} \big)$ with $i > h$ gives rise to the additional summand $(a, b) = (K_0 + 1, K_0) \in \mathfrak{T} \big( \lambda^{(h)} \big) \times \mathfrak{T} \big( \lambda^{(i)} \big)$ on the right side of the definition \eqref{lambdamupsi} of $\psi (\boldsymbol{\lambda})$ not present in the corresponding sum for $\psi (\boldsymbol{\nu})$. Similarly, the second sum arises since each appearance of $K_0 + 1$ in some $\mathfrak{T} \big( \lambda^{(i)} \big)$ for $i < h$ gives rise to the additional summand $(a, b) = (K_0 + 1, K_0) \in \mathfrak{T} \big( \nu^{(i)} \big) \times \mathfrak{T} \big( \nu^{(h)} \big) = \mathfrak{T} \big( \lambda^{(i)} \big) \times \mathfrak{T} \big( \nu^{(h)} \big)$ in the definition of $\psi  (\boldsymbol{\nu})$ not present in the correpsonding sum for $\psi (\boldsymbol{\lambda})$. By similar reasoning and also using the facts that particles are complemented, colors are reversed, and positions are reversed in $\boldsymbol{\lambda}'$ and $\boldsymbol{\nu}'$ with respect to $\boldsymbol{\lambda}'$ and $\boldsymbol{\nu}'$, respectively, we deduce that 
	\begin{flalign}
	\label{nulambda}
	\begin{aligned} 
	\psi (\boldsymbol{\lambda}') & = \psi (\boldsymbol{\nu}') + \displaystyle\frac{1}{2} \displaystyle\sum_{i = n - h + 2}^n \textbf{1}_{L_0 - K_0 \in \mathfrak{T} (\lambda'^{(i)})} - \displaystyle\frac{1}{2} \displaystyle\sum_{i = 1}^{n - h} \textbf{1}_{L_0 - K_0 + 1 \in \mathfrak{T} (\lambda'^{(i)})} \\
	& = \psi (\boldsymbol{\nu}') + \displaystyle\frac{1}{2} \displaystyle\sum_{i = 1}^{h - 1} \textbf{1}_{K_0 + 1 \not\in \mathfrak{T} (\lambda^{(i)})} - \displaystyle\frac{1}{2} \displaystyle\sum_{i = h + 1}^n \textbf{1}_{K_0 \notin \mathfrak{T} (\lambda^{(i)})}.
	\end{aligned} 
	\end{flalign}
	
	Subtracting \eqref{nulambda} from \eqref{lambdanu} yields
	\begin{flalign*}
	\psi (\boldsymbol{\lambda}) - \psi (\boldsymbol{\lambda}') - \psi (\boldsymbol{\nu}) + \psi (\boldsymbol{\nu}') & = \displaystyle\frac{1}{2} \displaystyle\sum_{i = h + 1}^n \big( \textbf{1}_{K_0 \in \mathfrak{T} (\lambda^{(i)})} + \textbf{1}_{K_0 \notin \mathfrak{T} (\lambda^{(i)})} \big) \\
	& \qquad - \displaystyle\frac{1}{2} \displaystyle\sum_{i = 1}^{h - 1} \big( \textbf{1}_{K_0 + 1 \in \mathfrak{T} (\lambda^{(i)})} + \textbf{1}_{K_0 + 1 \not\in \mathfrak{T} (\lambda^{(i)})} \big) = \displaystyle\frac{n - 1}{2} - (h - 1),
	\end{flalign*}
	
	\noindent which upon adding  to \eqref{psilambdalambdasum} (with the $\boldsymbol{\lambda}$ or $\boldsymbol{\mu}$ there equal to $\boldsymbol{\nu}$ here) gives
	\begin{flalign*}
	\psi \big( \boldsymbol{\lambda} \big) - \psi (\boldsymbol{\lambda}') + \displaystyle\sum_{i = 1}^n (i - 1) \big| \nu^{(i)} \big| + h - 1  = \displaystyle\frac{(n - 1) |\boldsymbol{\nu}|}{2} + \displaystyle\frac{n - 1}{2} + \displaystyle\frac{1}{2} \binom{n}{2} \Bigg( \binom{M}{2} - \binom{L_0 - M}{2} \Bigg).
	\end{flalign*}
	
	\noindent This, together with the facts that $\big| \lambda^{(i)} \big| = \big| \nu^{(i)} \big| + \textbf{1}_{i = h}$ and $|\boldsymbol{\lambda}| = |\boldsymbol{\nu}| + 1$, imples the first statement of \eqref{psilambdalambdasum} and thus the lemma. 
\end{proof}

Now we can establish the third part of \Cref{limitg0horizontal}. 

\begin{proof}[Proof of Part \ref{lg2} of \Cref{limitg0horizontal}]
	
As indicated in the proof of the first part of \Cref{limitg0horizontal} from \Cref{GPolynomials}, if $\lambda' / \mu'$ is not a horizontal $n$-ribbon strip then there exists some integer $j \in [1, L_0]$ such that $v' (j) \ne 0$. By \eqref{gproductq}, this implies that $\mathcal{G}_{\boldsymbol{\lambda} / \boldsymbol{\mu}} (0; x \boldsymbol{\mid} 0; 0) = 0$; since the same holds for $\mathcal{L}_{\boldsymbol{\lambda}' / \boldsymbol{\mu}'} (q^{(1 - n) / 2} x^{-1}; q^{-1})$, this verifies \eqref{2gl}. So, let us assume below that $\lambda' / \mu'$ is a horizontal $n$-ribbon strip. 

 Applying arrow conservation and the fact that $\textbf{B}' (j + 1) = \textbf{D}' (j)$ for each $j \ge 1$ (and that $\textbf{B}' (1) = \textbf{D} (L_0) = \textbf{B} (L_0 + 1) = \textbf{e}_0$), we deduce
\begin{flalign*}
\displaystyle\sum_{j = 1}^{\infty} \varphi \big( \textbf{C}' (j) + \textbf{D}' (j), \textbf{B}' (j) \big) & = \displaystyle\sum_{j = 1}^{\infty} \varphi \big( \textbf{A}' (j) + \textbf{B}' (j), \textbf{B}' (j) \big) \\
& = \displaystyle\sum_{j = 1}^{\infty} \Big( \varphi \big( \textbf{B}' (j), \textbf{B}' (j) \big) + \varphi \big( \textbf{A}' (j), \textbf{B}' (j) \big) \Big) \\
& = \displaystyle\sum_{j = 1}^{\infty} \Big( \varphi \big( \textbf{D}' (j), \textbf{D}' (j) \big) + \varphi \big( \textbf{A}' (j), \textbf{B}' (j) \big) \Big) \\
& = \spin \big( \lambda' / \mu' \big) - \psi \big( \boldsymbol{\lambda}' \big) + \psi \big( \boldsymbol{\mu}' \big),
\end{flalign*}

\noindent where in the last equality we applied \Cref{spinlambdamu}, and \Cref{psilambdamu}. Inserting this into \eqref{gproductq} (and using the bilinearity of $\varphi$) yields 
\begin{flalign*}
\mathcal{G}_{\boldsymbol{\lambda}} (0; x \boldsymbol{\mid} 0; 0) & = q^{\psi (\boldsymbol{\lambda}') - \psi (\boldsymbol{\mu}') - \spin (\lambda' / \mu')} x^{|\boldsymbol{\mu}| - |\boldsymbol{\lambda}|} \displaystyle\prod_{j = 1}^{\infty} q^{\varphi (\textbf{e}_{[1, n]}, \textbf{B}' (j))} \\
& = q^{\psi (\boldsymbol{\lambda}) - \psi (\boldsymbol{\mu}) - \spin (\lambda' / \mu')} (q^{(n - 1) / 2} x^{-1})^{|\boldsymbol{\lambda}| - |\boldsymbol{\mu}|} = q^{\psi (\boldsymbol{\lambda}) - \psi (\boldsymbol{\mu})} \mathcal{L}_{\boldsymbol{\lambda}' / \boldsymbol{\mu}'} (q^{(n - 1) / 2} x^{-1}; q^{-1}),
\end{flalign*}

\noindent where in the second equality we applied \Cref{sumeb} and in the last we used the definition \eqref{lambdamug} (and \eqref{llambdamu}) of the LLT polynomial $\mathcal{L}_{\boldsymbol{\lambda} / \boldsymbol{\mu}}$. This establishes the third part of \Cref{limitg0horizontal}. 
\end{proof}

\section{LLT Polynomials from \texorpdfstring{$\mathcal{H}_{\boldsymbol{\lambda} / \boldsymbol{\mu}} (\textbf{x}; \infty \boldsymbol{\mid} \infty; \infty)$}{}}

\label{Function3LG}

In this section we establish the fourth part of \Cref{limitg0horizontal}, namely, that \eqref{1hl} holds for $N = 1$. Similarly to in \Cref{Functions2GL}, this will follow from the second part of \Cref{limitg0}, by expressing $\mathcal{H}_{\boldsymbol{\lambda} / \boldsymbol{\mu}}^{(q)} (x; \infty \boldsymbol{\mid} \infty; \infty)$ in terms of $\mathcal{F}_{\boldsymbol{\lambda} / \boldsymbol{\mu}}^{(1 / q)} (x; \infty \boldsymbol{\mid} 0; 0)$; here, we write $\mathcal{H}^{(q)}$ to emphasize the dependence of $\mathcal{H}$ on $q$ and $\mathcal{F}^{(1 / q)}$ to indicate that the parameter $q$ involved in the definition of $\mathcal{F}$ is replaced by $q^{-1}$. To implement this, we establish the following statement comparing the more general functions $\mathcal{H}_{\boldsymbol{\lambda} / \boldsymbol{\mu}}^{(q)} (\textbf{x}; \infty \boldsymbol{\mid} \textbf{y}; \infty)$ and $\mathcal{F}_{\boldsymbol{\lambda} / \boldsymbol{\mu}}^{(1 / q)} (q^{n - 1} \textbf{x}; \infty \boldsymbol{\mid} \textbf{y}; 0)$ (from \eqref{limithf}). In the below, we recall $\psi$ from \eqref{lambdamupsi}.

\begin{prop}
	
	\label{hf}
	
	Fix integers $N \ge 1$ and $M \ge 0$; sequences of complex numbers $\textbf{\emph{x}} = (x_1, x_2, \ldots , x_N)$ and $\textbf{\emph{y}} = (y_1, y_2, \ldots )$; and a nonzero complex number $q \in \mathbb{C}$. For any signature sequences $\boldsymbol{\lambda} \in \SeqSign_{n; M + N}$ and $\boldsymbol{\mu} \in \SeqSign_{n; M}$, we have
	\begin{flalign*}
	\mathcal{H}_{\boldsymbol{\lambda} / \boldsymbol{\mu}}^{(q)} (\textbf{\emph{x}}; \infty \boldsymbol{\mid} \textbf{\emph{y}}; \infty) & = (-1)^{n \binom{M + N}{2} - n \binom{M}{2}} q^{2 \psi (\overleftarrow{\boldsymbol{\mu}}) - 2 \psi (\overleftarrow{\boldsymbol{\lambda}}) + 2 \binom{M + N}{2} \binom{n}{2} - 2 \binom{M}{2} \binom{n}{2}} \mathcal{F}_{\overleftarrow{\boldsymbol{\lambda}} / \overleftarrow{\boldsymbol{\mu}}}^{(1 / q)} (q^{n - 1} \textbf{\emph{x}}; \infty \boldsymbol{\mid} \textbf{\emph{y}}; 0) \\
	& \qquad \times \displaystyle\prod_{j = 1}^N x_j^{-n} \displaystyle\prod_{i = 1}^n \Bigg( \displaystyle\prod_{j= 1}^M y_{\mu_j^{(i)} + M - j}^{-1} \displaystyle\prod_{k = 1}^{M + N} y_{\lambda_k^{(i)} + M + N - k} \Bigg),
	\end{flalign*}
\end{prop}

Assuming \Cref{hf}, we can quickly establish the fourth part of \Cref{limitg0horizontal}.

\begin{proof}[Proof of \ref{lh1} of \Cref{limitg0horizontal} Assuming \Cref{hf}]

	For any $x \in \mathbb{C}$, \Cref{wlimits} implies
	\begin{flalign}
	\label{wxwxylimit}
	\begin{aligned}
	\widehat{\mathcal{W}}_x (\textbf{A}, \textbf{B}; \textbf{C}, \textbf{D} \boldsymbol{\mid} \infty, 0) & = \displaystyle\lim_{y \rightarrow \infty} (-y)^{d - n} \widehat{\mathcal{W}}_{x; y} (\textbf{A}, \textbf{B}; \textbf{C}, \textbf{D} \boldsymbol{\mid} \infty, 0); \\
	\mathcal{W}_x (\textbf{A}, \textbf{B}; \textbf{C}, \textbf{D} \boldsymbol{\mid} \infty, \infty) & = \displaystyle\lim_{y \rightarrow \infty} y^{-b} \mathcal{W}_{x; y} (\textbf{A}, \textbf{B}; \textbf{C}, \textbf{D} \boldsymbol{\mid} \infty, \infty).
	\end{aligned}
	\end{flalign}

	Now, by \Cref{weightesum} and \eqref{gfhe}, $\mathcal{F}_{\boldsymbol{\lambda} / \boldsymbol{\mu}} (\textbf{x}; \infty \boldsymbol{\mid} \textbf{y}; 0)$ and $\mathcal{F}_{\boldsymbol{\lambda} / \boldsymbol{\mu}} (\textbf{x}; \infty \boldsymbol{\mid} 0; 0)$ are partition functions for the vertex model $\mathfrak{P}_F (\boldsymbol{\lambda} / \boldsymbol{\mu})$ (recall \Cref{pgpfph}; see also the middle of \Cref{fgpaths}) under the vertex weights $\widehat{\mathcal{W}}_{x; y} (\textbf{A}, \textbf{B}; \textbf{C}, \textbf{D} \boldsymbol{\mid} \infty, 0)$ and $\widehat{\mathcal{W}}_x (\textbf{A}, \textbf{B}; \textbf{C}, \textbf{D} \boldsymbol{\mid} \infty, 0)$, respectively. Thus, defining $\textbf{y} = (y, y, \ldots)$, the first statement of \eqref{wxwxylimit} quickly implies (similarly to in \Cref{sumbw}) that
	\begin{flalign}
	\label{limitf}
	\mathcal{F}_{\boldsymbol{\lambda} / \boldsymbol{\mu}} (\textbf{x}; \infty \boldsymbol{\mid} 0; 0) & = \displaystyle\lim_{y \rightarrow \infty} (-y)^{|\boldsymbol{\mu}| - |\boldsymbol{\lambda}| + n \binom{M}{2} - n \binom{M + N}{2}} \mathcal{F}_{\boldsymbol{\lambda} / \boldsymbol{\mu}} (\textbf{x}; \infty \boldsymbol{\mid} \textbf{y}; 0).
	\end{flalign}
	
	\noindent Similarly, \Cref{weightesum} and \eqref{gfhe} together indicate that $\mathcal{H}_{\boldsymbol{\lambda} / \boldsymbol{\mu}} (\textbf{x}; \infty \boldsymbol{\mid} \textbf{y}; \infty)$ and $\mathcal{H}_{\boldsymbol{\lambda} / \boldsymbol{\mu}} (\textbf{x}; \infty \boldsymbol{\mid} \infty; \infty)$ are partition functions for the vertex model $\mathfrak{P}_H (\boldsymbol{\lambda} / \boldsymbol{\mu})$ (recall \Cref{pgpfph}; see also the right side of \Cref{fgpaths}) under the vertex weights $\mathcal{W}_{x; y} (\textbf{A}, \textbf{B}; \textbf{C}, \textbf{D} \boldsymbol{\mid} \infty, \infty)$ and $\mathcal{W}_x (\textbf{A}, \textbf{B}; \textbf{C}, \textbf{D} \boldsymbol{\mid} \infty, \infty)$, respectively. Thus, again letting $\textbf{y} = (y, y, \ldots)$, the second statement of \eqref{wxwxylimit} quickly gives
	\begin{flalign}
	\label{limith}
	\mathcal{H}_{\boldsymbol{\lambda} / \boldsymbol{\mu}} (\textbf{x}; \infty \boldsymbol{\mid} \infty; \infty) & = \displaystyle\lim_{y \rightarrow \infty} y^{|\boldsymbol{\mu}| - |\boldsymbol{\lambda}| + n \binom{M + 1}{2} - n \binom{M + N + 1}{2}} \mathcal{H}_{\boldsymbol{\lambda} / \boldsymbol{\mu}} (\textbf{x}; \infty \boldsymbol{\mid} \textbf{y}; \infty).
	\end{flalign}
	
	Applying the $\textbf{y} = (y, y, \ldots )$ case of \Cref{hf}, and letting $y$ tend to $\infty$, it follows from \eqref{limitf} and \eqref{limith} that 
	\begin{flalign}
	\label{hlimitf} 
	\mathcal{H}_{\boldsymbol{\lambda} / \boldsymbol{\mu}} (\textbf{x}; \infty \boldsymbol{\mid} \infty; \infty) = (-1)^{|\boldsymbol{\lambda}| - |\boldsymbol{\mu}|} q^{2 \psi(\overleftarrow{\boldsymbol{\mu}}) - 2 \psi(\overleftarrow{\boldsymbol{\lambda}}) + 2 \binom{M + N}{2} \binom{n}{2} - 2 \binom{M}{2} \binom{n}{2}} \mathcal{F}_{\overleftarrow{\boldsymbol{\lambda}} / \overleftarrow{\boldsymbol{\mu}}}^{(1 / q)} (q^{n - 1} \textbf{x}; \infty \boldsymbol{\mid} 0; 0) \displaystyle\prod_{j = 1}^N x_j^{-n}.
	\end{flalign}
	
	\noindent By \eqref{1fl}, we have that
	\begin{flalign*}
	\mathcal{F}_{\overleftarrow{\boldsymbol{\lambda}} / \overleftarrow{\boldsymbol{\mu}}}^{(1 / q)}  (q^{n - 1} \textbf{x}; \infty \boldsymbol{\mid} 0; 0) = q^{\psi (\overleftarrow{\boldsymbol{\lambda}}) - \psi (\overleftarrow{\boldsymbol{\mu}}) + \binom{M + N}{2} \binom{n}{2} / 2 - \binom{M}{2} \binom{n}{2} / 2} \mathcal{L}_{\overleftarrow{\boldsymbol{\lambda}} / \overleftarrow{\boldsymbol{\mu}}} (\textbf{x}^{-1}; q^{-1}) \displaystyle\prod_{j = 1}^N (q^{n - 1} x_j)^{n (j - M - N)}.
	\end{flalign*}
	
	\noindent Together with \eqref{hlimitf} and the fact that $(-1)^{|\boldsymbol{\lambda}| - |\boldsymbol{\mu}|} \mathcal{L}_{\overleftarrow{\boldsymbol{\lambda}} / \overleftarrow{\boldsymbol{\mu}}} (\textbf{x}; q^{-1}) = \mathcal{L}_{\overleftarrow{\boldsymbol{\lambda}} / \overleftarrow{\boldsymbol{\mu}}} (-\textbf{x}^{-1}; q^{-1})$ (which holds by the homogeneity of $\mathcal{L}$), this implies \eqref{1hl}.
\end{proof}

To establish \Cref{hf}, we require the following lemma stating a relation between the weights $\mathcal{W}_{x; y} (\textbf{A}, \textbf{B}; \textbf{C}, \textbf{D} \boldsymbol{\mid} \infty, \infty)$ and $\widehat{\mathcal{W}}_{x; y} (\textbf{A}, \textbf{B}; \textbf{C}, \textbf{D} \boldsymbol{\mid} \infty, 0)$ from \Cref{wlimits}, under complementation of horizontal arrows and reversal of colors; we refer to \Cref{lambdamulambdamu2} for a depiction. In what follows, we again write $\mathcal{W}_{x; y}^{(q)}$ to emphasize the dependence of $\mathcal{W}_{x; y}$ on $q$ and $\mathcal{W}_{x; y}^{(1 / q)}$ to indicate that its parameter $q$ is replaced by $q^{-1}$.

\begin{lem}
	
	\label{wequalityw} 
	
	For any $\textbf{\emph{A}}, \textbf{\emph{B}}, \textbf{\emph{C}}, \textbf{\emph{D}} \in \{ 0, 1 \}^n$ and $x, y, q \in \mathbb{C}$, we have
	\begin{flalign}
	\label{wxywxy}
	\begin{aligned}
	& \mathcal{W}_{x; y}^{(q)} (\overleftarrow{\textbf{\emph{C}}}, \textbf{\emph{e}}_{[1, n]} - \overleftarrow{\textbf{\emph{B}}}; \overleftarrow{\textbf{\emph{A}}}, \textbf{\emph{e}}_{[1, n]} - \overleftarrow{\textbf{\emph{D}}} \boldsymbol{\mid} \infty; \infty) \\
	& \quad = (-1)^c q^{\varphi (\textbf{\emph{A}}, \textbf{\emph{e}}_{[1, n]}) - \varphi (\textbf{\emph{A}}, \textbf{\emph{B}}) + \varphi (\textbf{\emph{D}}, \textbf{\emph{C}}) - \varphi (\textbf{\emph{B}}, \textbf{\emph{B}}) + \varphi (\textbf{\emph{D}}, \textbf{\emph{D}})} (q^{n - 1} x y^{-1})^{b - d} \widehat{\mathcal{W}}_{q^{n - 1} x; y}^{(1 / q)} (\textbf{\emph{A}}, \textbf{\emph{B}}; \textbf{\emph{C}}, \textbf{\emph{D}} \boldsymbol{\mid} \infty; 0).
	\end{aligned} 
	\end{flalign}
\end{lem}

\begin{proof}
	
Observe that if $\textbf{A} + \textbf{B} \ne \textbf{C} + \textbf{D}$, then both sides of \eqref{wxywxy} are equal to $0$. Moreover, defining $\textbf{V} = (V_1, V_2, \ldots , V_n) \in \{ 0, 1 \}^n$ by setting $V_i = \min \{A_i, B_i, C_i, D_i \}$ for each $i \in [1, n]$ and letting $v = |\textbf{V}|$ (as in \eqref{wabcdrsxy}), \eqref{wxywxy} also holds if $v \ne 0$. Indeed, in this case, the right side of \eqref{wxywxy} is equal to $0$, by the $\textbf{1}_{v = 0}$ weight in the definition \eqref{limitw} of $\widehat{\mathcal{W}}_{x; y} (\textbf{A}, \textbf{B}; \textbf{C}, \textbf{D} \boldsymbol{\mid} \infty, 0)$. Furthermore, if $v \ne 0$, then there exists some $i \in [1, n]$ for which $V_i = 1$, in which case $A_i = 1 = B_i$, and so $A_i + B_i > 1$. Thus, $\textbf{e}_{[1, n]} - \overleftarrow{\textbf{B}} \ge \overleftarrow{\textbf{A}}$ does not hold, and so the left side of \eqref{wxywxy} also equals $0$, due to the $\textbf{1}_{\textbf{B} \ge \textbf{C}}$ factor in the definition \eqref{limitw} of $\mathcal{W}_{x; y} (\textbf{A}, \textbf{B}; \textbf{C}, \textbf{D} \boldsymbol{\mid} \infty, \infty)$.  

So, let us assume in what follows that $\textbf{A} + \textbf{B} = \textbf{C} + \textbf{D}$ and $v = 0$. Then, \eqref{limitw2} gives
\begin{flalign}
\label{qn1xy}
\begin{aligned}
\widehat{\mathcal{W}}_{q^{n - 1} x; y}^{(1 / q)} (\textbf{A}, \textbf{B}; \textbf{C}, \textbf{D} \boldsymbol{\mid} \infty, 0) & = (-q^{n - 1} xy^{-1})^{d - n} q^{\binom{n}{2} - \varphi (\textbf{D}, \textbf{C} + \textbf{D})} \displaystyle\frac{(q^{n - 1} xy^{-1}; q^{-1})_n}{(q^{n - 1} xy^{-1}; q^{-1})_{a + b}} \\
& = (-q^{n - 1} xy^{-1})^{d - n} q^{-\varphi (\textbf{D}, \textbf{C}) - \varphi (\textbf{D}, \textbf{D}) - \binom{n}{2}} (x y^{-1}; q)_{n - a - b},
\end{aligned} 
\end{flalign}

\noindent where in the last equality we used the fact that $(q^{n - 1} z; q^{-1})_n (q^{n - 1} z; q^{-1})_k^{-1} = (z; q)_{n - k}$ for any $z \in \mathbb{C}$ and $k \in [0, n]$. Moreover, \eqref{limitw} gives
\begin{flalign}
\label{wxy1} 
\begin{aligned}
\mathcal{W}_{x; y}^{(q)} & (\overleftarrow{\textbf{C}}, \textbf{e}_{[1, n]} - \overleftarrow{\textbf{B}}; \overleftarrow{\textbf{A}}, \textbf{e}_{[1, n]} - \overleftarrow{\textbf{D}} \boldsymbol{\mid} \infty; \infty) \\
& = (-1)^{n - a - b} x^{b - n} y^{n - b} q^{\varphi (\textbf{e}_{[1, n]} - \overleftarrow{\textbf{B}}, \overleftarrow{\textbf{A}} + \overleftarrow{\textbf{B}} - \textbf{e}_{[1, n]})} (xy^{-1}; q)_{n - a - b}.
\end{aligned} 
\end{flalign}

\noindent By \eqref{xyxy}, the bilinearity of $\varphi$, the fact that $\varphi \big( \textbf{e}_{[1, n]}, \textbf{e}_{[1, n]} \big) = \binom{n}{2}$, and the $\textbf{X} = \textbf{D}$ case of \eqref{xn1}, we obtain
\begin{flalign*}
\varphi \big( \textbf{e}_{[1, n]} - & \overleftarrow{\textbf{B}}, \overleftarrow{\textbf{A}} + \overleftarrow{\textbf{B}} - \textbf{e}_{[1, n]} \big) \\
& = \varphi \big( \textbf{A} + \textbf{B} - \textbf{e}_{[1, n]}, \textbf{e}_{[1, n]} - \textbf{B} \big) \\
& = \varphi \big( \textbf{B}, \textbf{e}_{[1, n]} \big) + \varphi \big( \textbf{e}_{[1, n]}, \textbf{B} \big) - \varphi \big( \textbf{A} + \textbf{B}, \textbf{B} \big) + \varphi \big( \textbf{A}, \textbf{e}_{[1, n]} \big) - \varphi \big( \textbf{e}_{[1, n]}, \textbf{e}_{[1, n]} \big) \\
& = (n - 1) b - \binom{n}{2} + \varphi \big( \textbf{A}, \textbf{e}_{[1, n]} \big) - \varphi \big( \textbf{A}, \textbf{B} \big) - \varphi \big( \textbf{B}, \textbf{B} \big).
\end{flalign*}

\noindent Together with \eqref{wxy1}, this gives
\begin{flalign*}
\mathcal{W}_{x; y}^{(q)} & (\overleftarrow{\textbf{C}}, \textbf{e}_{[1, n]} - \overleftarrow{\textbf{B}}; \overleftarrow{\textbf{A}}, \textbf{e}_{[1, n]} - \overleftarrow{\textbf{D}} \boldsymbol{\mid} \infty; \infty) \\
& = (-1)^{n - a - b} x^{b - n} y^{n - b} q^{(n - 1) b + \varphi (\textbf{A}, \textbf{e}_{[1, n]}) - \varphi (\textbf{A}, \textbf{B}) - \varphi (\textbf{B}, \textbf{B}) - \binom{n}{2}} (xy^{-1}; q)_{n - a - b}.
\end{flalign*}

\noindent This, \eqref{qn1xy}, and the fact that $a + b = c + d$ together imply the lemma.
\end{proof}

\begin{rem}
	
	\label{wweightw}
	
	The type of ``horizontal complementation symmetry'' described by \Cref{wequalityw} does not appear to admit a transparent generalization to the most general fused weights given by \Cref{wabcdrsxy}. Indeed, if $\textbf{B} \le \textbf{A}, \textbf{C}, \textbf{D}$ then the weight $W_z (\textbf{A}, \textbf{B}; \textbf{C}, \textbf{D} \boldsymbol{\mid} r, s)$ factors completely, since the sum on the right side of \eqref{wabcdp} is supported on the $p = 0$ term (as then $b - v = 0$, since $\textbf{V} = \textbf{B}$ there). However, for generic choices of parameters, the weight $W_z (\textbf{C}, \textbf{e}_{[1, n]} - \textbf{B}; \textbf{A}, \textbf{e}_{[1, n]} - \textbf{D} \boldsymbol{\mid} r, s)$ does not factor under the condition $\textbf{B} \le \textbf{A}, \textbf{C}, \textbf{D}$.

\end{rem}

Now we can establish \Cref{hf}.

\begin{proof}[Proof of \Cref{hf}]

Throughout this proof, let us assume for notational simplicity that $N = 1$, as the proof for general $N$ is entirely analogous (alternatively, given the $N = 1$ case, the proof for general $N$ can be established using the branching identities \eqref{sumf}, as in the proof of \Cref{limitg0} assuming \Cref{limitg0horizontal} in \Cref{FunctionsFGL}). 

Then, recall from \eqref{gfhe} and \Cref{weightesum} that $\mathcal{H}_{\boldsymbol{\lambda} / \boldsymbol{\mu}} (\textbf{x}; \infty \boldsymbol{\mid} \textbf{y}; \infty)$ is the partition function under the $\mathcal{W}_{x; y} (\textbf{A}, \textbf{B}; \textbf{C}, \textbf{D} \boldsymbol{\mid} \infty; \infty)$ weights (from \eqref{limitw}) for the vertex model $\mathfrak{P}_H (\boldsymbol{\lambda} / \boldsymbol{\mu})$ (from \Cref{pgpfph}) on the domain $\mathcal{D} = \mathcal{D}_1 = \mathbb{Z}_{> 0} \times \{ 1 \}$. Since this domain has one row, $\mathfrak{P}_H (\boldsymbol{\lambda} / \boldsymbol{\mu})$ only consists in a single path ensemble, which we denote by $\mathcal{E}^* (\boldsymbol{\lambda} / \boldsymbol{\mu})$. For each $j \ge 1$, let $\big( \textbf{A}^* (j), \textbf{B}^* (j); \textbf{C}^* (j), \textbf{D}^* (j) \big)$ denote the arrow configuration at the vertex $(j, 1) \in \mathcal{D}$ under $\mathcal{E}^* (\boldsymbol{\lambda} / \boldsymbol{\mu})$. In particular, there exists a minimal integer $L_0 > 1$ (given explicitly by \eqref{l0}) such that we have $\big( \textbf{A}^* (K), \textbf{B}^* (K); \textbf{C}^* (K), \textbf{D}^* (K) \big) = (\textbf{e}_0, \textbf{e}_0; \textbf{e}_0, \textbf{e}_0)$ for $K \ge L_0$; see the left side of \Cref{lambdamulambdamu2}. 

For any $j \ge 1$, define $\widehat{\textbf{A}} (j), \widehat{\textbf{B}} (j), \widehat{\textbf{C}} (j), \widehat{\textbf{D}} (j) \in \{ 0, 1 \}^n$ by
\begin{flalign}
	\label{aabbccdd2} 
	\begin{aligned} 
		& \widehat{\textbf{A}} (j) = \overleftarrow{\textbf{C}}^* (j); \qquad \widehat{\textbf{B}} (j) = \textbf{e}_{[1, n]} - \overleftarrow{\textbf{B}}^* (j); \qquad \widehat{\textbf{C}} (j) = \overleftarrow{\textbf{A}}^* (j); \qquad \widehat{\textbf{D}} (j) = \textbf{e}_{[1, n]} - \overleftarrow{\textbf{D}}^* (j).
	\end{aligned} 
\end{flalign}

\noindent For each index $X \in \{A, B, C, D \}$, let $x^* (j) = \big| \textbf{X}^* (j) \big|$ and $\widehat{x} = \big| \widehat{\textbf{X}} \big|$.

Then, since $\textbf{A}^* (j) + \textbf{B}^* (j) = \textbf{C}^* (j) + \textbf{D}^* (j)$, \eqref{aabbccdd2} implies each $\widehat{\textbf{A}} (j) + \widehat{\textbf{B}} (j) = \widehat{\textbf{C}} (j) + \widehat{\textbf{D}} (j)$. In particular, this yields the existence of a path ensemble $\widehat{\mathcal{E}} \big( \overleftarrow{\boldsymbol{\lambda}} / \overleftarrow{\boldsymbol{\mu}} \big)$ whose arrow configuration at $(j, 1)$ is $\big( \widehat{\textbf{A}} (j), \widehat{\textbf{B}} (j); \widehat{\textbf{C}} (j), \widehat{\textbf{D}} (j) \big)$. It is quickly verified that this definition is consistent with the notation from \Cref{Fx1FunctionL}, namely, that $\widehat{\mathcal{E}} \big( \overleftarrow{\boldsymbol{\lambda}} / \overleftarrow{\boldsymbol{\mu}} \big)$ is the unique element of $\mathfrak{P}_F \big( \overleftarrow{\boldsymbol{\lambda}} / \overleftarrow{\boldsymbol{\mu}} \big)$; see \Cref{lambdamulambdamu2}.

\begin{figure}

	\begin{center}

		\begin{tikzpicture}[
		>=stealth,
		scale = .42
		]
		
		\draw[red, thick, ->] (-2, 2.4) -- (1.9, 2.4) -- (1.9, 4.5);
		\draw[blue, thick, ->] (-2, 2.5) -- (0, 2.5) -- (0, 4.5);
		\draw[green, thick, ->] (-2, 2.6) -- (2.1, 2.6) -- (2.1, 4.5);
		
		\draw[blue, thick, ->] (0, .5) -- (0, 2.5) -- (2, 2.5) -- (2, 4.5);
		
		\draw[red, thick, ->] (1.9, .5) -- (1.9, 2.4) -- (3.9, 2.4) -- (3.9, 4.5);
		\draw[green, thick, ->] (2.1, .5) -- (2.1, 2.6) -- (4.1, 2.6) -- (4.1, 4.5);
		
		\draw[red, thick, ->] (3.9, .5) -- (3.9, 2.4) -- (7.9, 2.4) -- (7.9, 4.5);
		\draw[blue, thick, ->] (4, .5) -- (4, 2.5) -- (6, 2.5) -- (6, 4.5);
		
		\draw[green, thick, ->] (6.1, .5) -- (6.1, 2.6) -- (8.1, 2.6) -- (8.1, 4.5);
		
		\draw[red, thick, ->] (7.9, .5) -- (7.9, 2.4) -- (9.9, 2.4) -- (9.9, 4.5);
		\draw[blue, thick, ->] (8, .5) -- (8, 2.5) -- (12, 2.5) -- (12, 4.5);
		\draw[green, thick, ->] (8.1, .5) -- (8.1, 2.6) -- (10.1, 2.6) -- (10.1, 4.5);
		
		\draw[] (0, .5) circle[radius = 0] node[below, scale = .65]{$\textbf{A}^* (1)$};
		\draw[] (2, .5) circle[radius = 0] node[below, scale = .65]{$\textbf{A}(2)$};
		\draw[] (4, .5) circle[radius = 0] node[below, scale = .65]{$\textbf{A}^*(3)$};
		\draw[] (6, .5) circle[radius = 0] node[below, scale = .65]{$\textbf{A}^*(4)$};
		\draw[] (8, .5) circle[radius = 0] node[below, scale = .65]{$\textbf{A}^*(5)$};
		\draw[] (10, .5) circle[radius = 0] node[below, scale = .65]{$\textbf{A}^*(6)$};
		\draw[] (12, .5) circle[radius = 0] node[below, scale = .65]{$\textbf{A}^*(7)$};
		
		\draw[] (0, 4.5) circle[radius = 0] node[above, scale = .65]{$\textbf{C}^*(1)$};
		\draw[] (2, 4.5) circle[radius = 0] node[above, scale = .65]{$\textbf{C}^*(2)$};
		\draw[] (4, 4.5) circle[radius = 0] node[above, scale = .65]{$\textbf{C}^*(3)$};
		\draw[] (6, 4.5) circle[radius = 0] node[above, scale = .65]{$\textbf{C}^*(4)$};
		\draw[] (8, 4.5) circle[radius = 0] node[above, scale = .65]{$\textbf{C}^*(5)$};
		\draw[] (10, 4.5) circle[radius = 0] node[above, scale = .65]{$\textbf{C}^*(6)$};
		\draw[] (12, 4.5) circle[radius = 0] node[above, scale = .65]{$\textbf{C}^*(7)$};
		
		\draw[] (-1.5, .5) circle[radius = 0] node[below, scale = .85]{$\boldsymbol{\mu}$};
		\draw[] (-1.5, 4.5) circle[radius = 0] node[above, scale = .85]{$\boldsymbol{\lambda}$};
		
		\draw[] (0, -.5) circle[radius = 0] node[below, scale = .65]{$1$};
		\draw[] (2, -.5) circle[radius = 0] node[below, scale = .65]{$2$};
		\draw[] (4, -.5) circle[radius = 0] node[below, scale = .65]{$3$};
		\draw[] (6, -.5) circle[radius = 0] node[below, scale = .55]{$4$};
		\draw[] (8, -.5) circle[radius = 0] node[below, scale = .55]{$5$};
		\draw[] (10, -.5) circle[radius = 0] node[below, scale = .55]{$6$};
		\draw[] (12, -.5) circle[radius = 0] node[below, scale = .55]{$L_0 - 1= 7$};

		\draw[] (16.5, .5) circle[radius = 0] node[below, scale = .85]{$\overleftarrow{\boldsymbol{\lambda}}$};
		\draw[] (16.5, 4.5) circle[radius = 0] node[above, scale = .85]{$\overleftarrow{\boldsymbol{\mu}}$};

		\draw[blue, thick, ->] (18, .5) -- (18, 4.5);
		
		\draw[red, thick, ->] (19.9, .5) -- (19.9, 4.5);
		\draw[blue, thick, ->] (20, .5) -- (20, 2.5) -- (22, 2.5) -- (22, 4.5);
		\draw[green, thick, ->] (20.1, .5) -- (20.1, 4.5);
		
		\draw[red, thick, ->] (21.9, .5) -- (21.9, 2.4) -- (23.9, 2.4) -- (23.9, 4.5);
		\draw[green, thick, ->] (22.1, .5) -- (22.1, 4.5); 
		
		\draw[blue, thick, ->] (24, .5) -- (24, 2.5) -- (26, 2.5) -- (26, 4.5);
		
		\draw[red, thick, ->] (25.9, .5) -- (25.9, 4.5);
		\draw[green, thick, ->] (26.1, .5) -- (26.1, 4.5);
		
		\draw[red, thick, ->] (27.9, .5) -- (27.9, 2.4) -- (32, 2.4);
		\draw[green, thick, ->] (28.1, .5) -- (28.1, 2.6) -- (32, 2.6);
		
		\draw[blue, thick, ->] (30, .5) -- (30, 2.5) -- (32, 2.5);

		\draw[] (18, .5) circle[radius = 0] node[below, scale = .65]{$\widehat{\textbf{A}} (1)$};
		\draw[] (20, .5) circle[radius = 0] node[below, scale = .65]{$\widehat{\textbf{A}} (2)$};
		\draw[] (22, .5) circle[radius = 0] node[below, scale = .65]{$\widehat{\textbf{A}} (3)$};
		\draw[] (24, .5) circle[radius = 0] node[below, scale = .65]{$\widehat{\textbf{A}} (4)$};
		\draw[] (26, .5) circle[radius = 0] node[below, scale = .65]{$\widehat{\textbf{A}} (5)$};
		\draw[] (28, .5) circle[radius = 0] node[below, scale = .65]{$\widehat{\textbf{A}} (6)$};
		\draw[] (30, .5) circle[radius = 0] node[below, scale = .65]{$\widehat{\textbf{A}} (7)$};
		
		\draw[] (18, 4.5) circle[radius = 0] node[above, scale = .65]{$\widehat{\textbf{C}} (1)$};
		\draw[] (20, 4.5) circle[radius = 0] node[above, scale = .65]{$\widehat{\textbf{C}} (2)$};
		\draw[] (22, 4.5) circle[radius = 0] node[above, scale = .65]{$\widehat{\textbf{C}} (3)$};
		\draw[] (24, 4.5) circle[radius = 0] node[above, scale = .65]{$\widehat{\textbf{C}} (4)$};
		\draw[] (26, 4.5) circle[radius = 0] node[above, scale = .65]{$\widehat{\textbf{C}} (5)$};
		\draw[] (28, 4.5) circle[radius = 0] node[above, scale = .65]{$\widehat{\textbf{C}} (6)$};
		\draw[] (30, 4.5) circle[radius = 0] node[above, scale = .65]{$\widehat{\textbf{C}} (7)$};
		
		\draw[] (18, -.5) circle[radius = 0] node[below, scale = .65]{$1$};
		\draw[] (20, -.5) circle[radius = 0] node[below, scale = .65]{$2$};
		\draw[] (22, -.5) circle[radius = 0] node[below, scale = .65]{$3$};
		\draw[] (24, -.5) circle[radius = 0] node[below, scale = .65]{$4$};
		\draw[] (26, -.5) circle[radius = 0] node[below, scale = .55]{$5$};
		\draw[] (28, -.5) circle[radius = 0] node[below, scale = .55]{$6$};
		\draw[] (30, -.5) circle[radius = 0] node[below, scale = .55]{$L_0 - 1= 7$};

		\end{tikzpicture}
		
	\end{center}
	
	\caption{\label{lambdamulambdamu2} Shown to the left and right are examples of the single-row vertex models $\mathcal{E}^* (\boldsymbol{\lambda} / \boldsymbol{\mu})$ and $\widehat{\mathcal{E}} \big( \protect\overleftarrow{\boldsymbol{\lambda}} / \protect\overleftarrow{\boldsymbol{\mu}} \big)$, respectively. }
\end{figure}

Thus, $\mathcal{H}_{\boldsymbol{\lambda} / \boldsymbol{\mu}}^{(q)} (x; \infty \boldsymbol{\mid} \textbf{y}; \infty)$ and $\mathcal{F}_{\overleftarrow{\boldsymbol{\lambda}} / \overleftarrow{\boldsymbol{\mu}}}^{(1 / q)} (q^{n - 1} x; \infty \boldsymbol{\mid} \textbf{y}; 0)$ are partition functions, under the weights $\mathcal{W}_{x; y}^{(q)} (\textbf{A}, \textbf{B}; \textbf{C}, \textbf{D} \boldsymbol{\mid} \infty, \infty)$ and $\widehat{\mathcal{W}}_{q^{n - 1} x; y}^{(1 / q)} (\textbf{A}, \textbf{B}; \textbf{C}, \textbf{D} \boldsymbol{\mid} \infty, 0)$, for the path ensembles $\mathcal{E}^* (\boldsymbol{\lambda} / \boldsymbol{\mu})$ and $\widehat{\mathcal{E}} \big( \overleftarrow{\boldsymbol{\lambda}} / \overleftarrow{\boldsymbol{\mu}} \big)$, respectively. Therefore, \Cref{wequalityw} yields 
\begin{flalign}
\label{hfproduct} 
\begin{aligned} 
\mathcal{H}_{\boldsymbol{\lambda} / \boldsymbol{\mu}}^{(q)} (x; \infty \boldsymbol{\mid} \textbf{y}; \infty) & = \displaystyle\prod_{j = 1}^{L_0 - 1} \mathcal{W}_{x; y}^{(q)} \big( \textbf{A}^* (j), \textbf{B}^* (j); \textbf{C}^* (j), \textbf{D}^* (j) \boldsymbol{\mid} \infty, \infty \big) \\
& = \displaystyle\prod_{j = 1}^{L_0 - 1} (-1)^{\widehat{c} (j)} q^{\varphi (\widehat{\textbf{A}} (j), \textbf{e}_{[1, n]}) + \varphi (\widehat{\textbf{D}} (j), \widehat{\textbf{C}} (j)) - \varphi (\widehat{\textbf{A}} (j), \widehat{\textbf{B}} (j)) + \varphi (\widehat{\textbf{D}} (j), \widehat{\textbf{D}} (j)) - \varphi (\widehat{\textbf{B}} (j), \widehat{\textbf{B}} (j))} \\
& \qquad \qquad \times (q^{n - 1} x y_j^{-1})^{\widehat{b} (j) - \widehat{d} (j)} \widehat{\mathcal{W}}_{q^{n - 1} x; y}^{(1 / q)} \big( \widehat{\textbf{A}} (j), \widehat{\textbf{B}} (j); \widehat{\textbf{C}} (j), \widehat{\textbf{D}} (j) \boldsymbol{\mid} \infty, 0 \big) \\
& = \mathcal{F}_{\overleftarrow{\boldsymbol{\lambda}} / \overleftarrow{\boldsymbol{\mu}}}^{(1 / q)} (q^{n - 1} x; \infty \boldsymbol{\mid} \textbf{y}; 0) \displaystyle\prod_{j = 1}^{L_0 - 1} (-1)^{\widehat{c} (j)} q^{\varphi (\widehat{\textbf{A}} (j), \textbf{e}_{[1, n]}) + \varphi (\widehat{\textbf{D}} (j), \widehat{\textbf{C}} (j)) - \varphi (\widehat{\textbf{A}} (j), \widehat{\textbf{B}} (j))} \\
& \qquad \qquad \qquad \qquad \qquad \qquad \qquad \times q^{\varphi (\widehat{\textbf{D}} (j), \widehat{\textbf{D}} (j)) - \varphi (\widehat{\textbf{B}} (j), \widehat{\textbf{B}} (j))} (q^{n - 1} x y_j^{-1})^{\widehat{b} (j) - \widehat{d} (j)}.
\end{aligned} 
\end{flalign}

\noindent Now let us analyze the terms in the product appearing on the right side of \eqref{hfproduct}. To that end, observe that 
\begin{flalign}
\label{sumcae}
\displaystyle\sum_{j = 1}^{L_0 - 1} \widehat{c} (j) = nM; \qquad \displaystyle\sum_{j = 1}^{L_0 - 1} \varphi \big( \widehat{\textbf{A}} (j), \textbf{e}_{[1, n]} \big) = (M + 1) \varphi \big( \textbf{e}_{[1, n]}, \textbf{e}_{[1, n]} \big) = (M + 1) \binom{n}{2},
\end{flalign} 

\noindent where the first equality holds since $\boldsymbol{\mu} \in \SeqSign_{n; M}$, and the second holds since $\sum_{j = 1}^{L_0 - 1} \widehat{\textbf{A}} (j) = (M + 1) \textbf{e}_{[1, n]}$ (as $\boldsymbol{\lambda} \in \SeqSign_{n; M + 1}$). We also have that
\begin{flalign}
\label{dbsum}
 \displaystyle\sum_{j = 1}^{L_0 - 1} \Big( \varphi \big( \widehat{\textbf{D}} (j), \widehat{\textbf{D}} (j) \big) - \varphi \big( \widehat{\textbf{B}} (j), \widehat{\textbf{B}} (j) \big) \Big) = \varphi \big( \widehat{\textbf{B}} (L_0), \widehat{\textbf{B}} (L_0) \big) - \varphi \big( \widehat{\textbf{B}} (1), \widehat{\textbf{B}} (1) \big) = \binom{n}{2},
\end{flalign}

\noindent where the first equality holds since each $\widehat{\textbf{D}} (j) = \widehat{\textbf{B}} (j + 1)$, and the second holds as $\varphi \big( \widehat{\textbf{B}} (1), \widehat{\textbf{B}} (1) \big) = \varphi (\textbf{e}_0, \textbf{e}_0) = 0$ and $\varphi \big( \widehat{\textbf{B}} (L_0), \widehat{\textbf{B}} (L_0) \big) = \varphi \big( \textbf{e}_{[1, n]}, \textbf{e}_{[1, n]} \big) = \binom{n}{2}$. Moreover, 
\begin{flalign}
\label{productxq}
\begin{aligned}
& \displaystyle\prod_{j = 1}^{L_0 - 1} (q^{n - 1} x)^{\widehat{b} (j) - \widehat{d} (j)} = (q^{n - 1} x)^{\widehat{b} (1) - \widehat{b} (L_0)} = q^{- 2 \binom{n}{2}} x^{-n}; \\
& \displaystyle\prod_{j = 1}^{L_0 - 1} y_j^{\widehat{d} (j) - \widehat{b} (j)} = \displaystyle\prod_{i = 1}^n \displaystyle\prod_{j = 1}^{L_0 - 1} y_j^{\widehat{D}_i (j) - \widehat{B}_i (j)} = \displaystyle\prod_{i = 1}^n \Bigg( \displaystyle\prod_{j = 1}^M y_{\mu_j^{(i)} + M - j}^{-1} \displaystyle\prod_{k = 1}^{M + 1} y_{\lambda_k^{(i)} + M - k + 1} \Bigg),
\end{aligned}
\end{flalign}

\noindent where we have denoted $\widehat{\textbf{X}} (j) = \big( \widehat{X}_1 (j), \widehat{X}_2 (j), \ldots , \widehat{X}_n (j) \big)$, for each index $X \in \{ B, D \}$. The first equality of the first statement in \eqref{productxq} holds since each $\widehat{d} (j) = \widehat{b} (j + 1)$; the second equality of the first statement holds since $\widehat{b} (1) = 0$ and $\widehat{b} (L_0) = n$; the first equality of the second statement holds since $\widehat{d} (j) - \widehat{b} (j) = \sum_{j = 1}^{L_0 - 1} \big( \widehat{D}_i (j) - \widehat{B}_i (j) \big)$; and the second equality of the second statement holds since $\widehat{D}_i (j) - \widehat{B}_i (j) = 1$ if and only if $j \in \mathfrak{T} \big( \lambda^{(n - i + 1)} \big)$ and $j \notin \mathfrak{T} \big( \mu^{(n - i + 1)} \big)$, and $\widehat{D}_i (j) - \widehat{B}_i (j) = -1$ holds if and only if $j \in \mathfrak{T} \big( \mu^{(n - i + 1)} \big)$ and $j \notin \mathfrak{T} \big( \lambda^{(n - i + 1)} \big)$. 

Together, \eqref{hfproduct}, \eqref{sumcae}, \eqref{dbsum}, and \eqref{productxq} imply
\begin{flalign*}
\mathcal{H}_{\boldsymbol{\lambda} / \boldsymbol{\mu}}^{(q)} (x; \infty \boldsymbol{\mid} \textbf{y}; \infty) & = (-1)^{nM} q^{M \binom{n}{2}} x^{-n} \mathcal{F}_{\boldsymbol{\lambda} / \boldsymbol{\mu}}^{(1 / q)} (q^{n - 1} x; \infty \boldsymbol{\mid} \textbf{y}; \infty) \\
& \qquad \times \displaystyle\prod_{j = 1}^{L_0 - 1} q^{\varphi (\widehat{\textbf{D}} (j), \widehat{\textbf{C}} (j)) - \varphi (\widehat{\textbf{A}} (j), \widehat{\textbf{B}} (j))} \displaystyle\prod_{i = 1}^n \Bigg( \displaystyle\prod_{j = 1}^M y_{\mu_j^{(i)} + M - j}^{-1} \displaystyle\prod_{k = 1}^{M + 1} y_{\lambda_k^{(i)} + M - k + 1} \Bigg). 
\end{flalign*} 

\noindent Thus, to establish the proposition it suffices to show that 
\begin{flalign}
\label{sumpsi} 
\displaystyle\sum_{j = 1}^{L_0 - 1} \Big( \varphi \big( \widehat{\textbf{D}} (j), \widehat{\textbf{C}} (j)  \big) - \varphi \big( \widehat{\textbf{A}} (j), \widehat{\textbf{B}} (j) \big) \Big) = 2 \psi (\overleftarrow{\boldsymbol{\mu}}) - 2 \psi (\overleftarrow{\boldsymbol{\lambda}}) + M \binom{n}{2}.
\end{flalign} 

To that end, it will be useful to produce an analog of the path ensemble $\mathcal{E} \big( \widetilde{\boldsymbol{\mu}} / \mathring{\boldsymbol{\lambda}} \big)$ from \Cref{Fx1FunctionL}, following \eqref{lambda2mu}. So recalling the notation $\overleftarrow{\boldsymbol{\lambda}} = \big( \overleftarrow{\lambda}^{(1)}, \overleftarrow{\lambda}^{(2)}, \ldots , \overleftarrow{\lambda}^{(n)} \big)$ and $\overleftarrow{\lambda}^{(i)} = \big( \overleftarrow{\lambda}_1^{(i)}, \overleftarrow{\lambda}_2^{(i)}, \ldots , \overleftarrow{\lambda}_{M + 1}^{(i)} \big)$ (and similarly for $\overleftarrow{\boldsymbol{\mu}}$ and the $\overleftarrow{\mu}^{(i)}$), following \eqref{lambda2mu} we define the sequences $\mathring{\boldsymbol{\lambda}}, \widetilde{\boldsymbol{\mu}} \in \SeqSign_{n; M + 1}$ by setting 
\begin{flalign*}
\mathring{\lambda}_j^{(i)} = \overleftarrow{\lambda}_j^{(i)} + 1, \quad & \text{for all $j \in [1, M +1]$}; \\
\widetilde{\mu}_1^{(i)} = L_0 - M - 1, \quad & \text{and} \quad  \widetilde{\mu}_j^{(i)} = \overleftarrow{\mu}_{j - 1}^{(i)} + 1, \quad \text{for all $j \in [2, M + 1]$}.
\end{flalign*}

\noindent For any $j \ge 1$, let $\big( \widetilde{\textbf{A}} (j), \widetilde{\textbf{B}} (j); \widetilde{\textbf{C}} (j), \widetilde{\textbf{D}} (j) \big)$ denote the arrow configuration at the vertex $(j, 1)$ under the path ensemble $\mathcal{E} \big( \widetilde{\boldsymbol{\mu}} / \mathring{\boldsymbol{\lambda}} \big) \in \mathfrak{P}_G \big( \widetilde{\boldsymbol{\mu}} / \mathring{\boldsymbol{\lambda}}; 1 \big)$. We refer to \Cref{lambdamu1mu2mu3vertex} for a depiction (where, in the present context, the left path ensemble there is $\widehat{\mathcal{E}} \big( \overleftarrow{\boldsymbol{\lambda}} / \overleftarrow{\boldsymbol{\mu}} \big)$, instead of $\widehat{\mathcal{E}} ( \boldsymbol{\lambda} / \boldsymbol{\mu})$).

In particular, we have that $\widetilde{\textbf{A}} (j) = \widehat{\textbf{A}} (j - 1)$ for each $j \ge 2$; that $\widetilde{\textbf{B}} (j + 1) = \widetilde{\textbf{D}} (j) = \widehat{\textbf{B}} (j) = \widehat{\textbf{D}} (j - 1)$ for each $j < L_0$; and that $\widetilde{\textbf{A}} (K + 1) = \widetilde{\textbf{B}} (K + 1) = \widetilde{\textbf{D}} (K) = \textbf{e}_0$ for each $K \ge L_0$. Thus,
\begin{flalign*}
\displaystyle\sum_{j = 1}^{L_0 - 1} \Big( \varphi \big( \widehat{\textbf{D}} (j), \widehat{\textbf{C}} (j) \big) - \varphi \big( \widehat{\textbf{A}} (j), \widehat{\textbf{B}} (j) \big) \Big) & = \displaystyle\sum_{j = 1}^{\infty} \Big( \varphi \big( \widetilde{\textbf{D}} (j), \widetilde{\textbf{C}} (j) \big) - \varphi \big( \widetilde{\textbf{A}} (j), \widetilde{\textbf{B}} (j) \big) \Big)  \\
& = 2 \psi \big( \widetilde{\boldsymbol{\mu}} \big) - 2 \psi \big( \mathring{\boldsymbol{\lambda}} \big) = 2 \psi \big( \overleftarrow{\boldsymbol{\mu}} \big) - 2 \psi \big(\overleftarrow{\boldsymbol{\lambda}} \big) + M \binom{n}{2},
\end{flalign*}

\noindent where in the second equality we applied \Cref{psilambdamu} (with the $\boldsymbol{\lambda} / \boldsymbol{\mu}$ there equal to $\widetilde{\boldsymbol{\mu}} / \mathring{\boldsymbol{\lambda}}$ here) and in the third we applied \eqref{psilambdapsimu}. This establishes \eqref{sumpsi} and thus the proposition. 
\end{proof}

\chapter{Contour Integral Formulas for \texorpdfstring{$G_{\boldsymbol{\lambda}}$}{}}

\label{IntegralG}

In this chapter we establish contour integral formulas for the functions $G_{\boldsymbol{\lambda}} (\textbf{x}; \textbf{r} \boldsymbol{\mid} \textbf{y}; \textbf{s})$. For simplicity, we only implement this in the homogeneous regime when all $y_i = 1$ and $s_i = s$, although more general results allowing for distinct $\{ y_i \}$ and $\{ s_i \}$ can be obtained similarly.

\section{Nonsymmetric Functions and Integral Formulas} 

\label{FunctionsNonsymmetric} 

We will establish our contour integral formulas for $G_{\boldsymbol{\lambda}} (\textbf{x}; \textbf{r} \boldsymbol{\mid} \textbf{y}; \textbf{s})$, with $\boldsymbol{\lambda} \in \SeqSign_{n; M}$, by using the color merging result \Cref{zsumfused} to compare the $U_q \big( \widehat{\mathfrak{sl}} (1 | n) \big)$ symmetric partition function $G_{\boldsymbol{\lambda} / \boldsymbol{\mu}}$ to certain $U_q \big( \widehat{\mathfrak{sl}} (1 | nM) \big)$ nonsymmetric ones (which coincide with the $U_q \big( \widehat{\mathfrak{sl}} (nM + 1) \big)$ nonsymmetric functions from \cite{SVMST}). Thus, we begin by recalling from \cite{SVMST} the definitions and properties of these nonsymmetric functions. Throughout this section, we fix an integer $n \ge 1$ (which in later sections will be set to $nM$) and a complex number $s \in \mathbb{C}$. 

There will be two nonsymmetric functions of interest to us here, which will be denoted by $\mathtt{f}_{\lambda}$ and $\mathtt{g}_{\lambda}$. The first function $\mathtt{f}_{\lambda}$ will be defined as a certain transformation of the second one $\mathtt{g}_{\lambda}$, and the latter will be defined as the partition function for a particular $U_q \big(\widehat{\mathfrak{sl}} (1 | n) \big)$ vertex model. 

The weights for this vertex model, denoted by $M_{z, q, s}^{(1; n)} (\textbf{A}, \textbf{B}; \textbf{C}, \textbf{D})$, will be given by a transformation of the $W_{z; q}^{(1; n)} (\textbf{A}, \textbf{B}; \textbf{C}, \textbf{D} \boldsymbol{\mid} r, s)$ from \Cref{wabcdrsxy}. More specifically, for any complex numbers $z, q \in \mathbb{C}$ and elements $\textbf{A}, \textbf{B}, \textbf{C}, \textbf{D} \in \{ 0, 1 \}^n$, set $M (\textbf{A}, \textbf{B}; \textbf{C}, \textbf{D}) = M_z (\textbf{A}, \textbf{B}; \textbf{C}, \textbf{D}) = M_{z, q, s}^{(1; n)} (\textbf{A}, \textbf{B}; \textbf{C}, \textbf{D})$ to be\footnote{These $M$ weights are closely related to those given by equation (2.2.6) of \cite{SVMST}; see \Cref{fg1} below.}
\begin{flalign}
	\label{m}
	M (\textbf{A}, \textbf{B}; \textbf{C}, \textbf{D}) = (-s)^{\textbf{1}_{|\textbf{B}| > 0}} W_{s / z; 1 / q}^{(1; n)} (\textbf{A}, \textbf{B}; \textbf{C}, \textbf{D} \boldsymbol{\mid} q^{1 / 2}, s^{-1}). 
\end{flalign}
\index{M@$M_x (\textbf{A}, \textbf{B}; \textbf{C}, \textbf{D})$}

\noindent In particular, these $M$ weights are (up to a change of variables and an overall multiplicative factor) given by the $r = q^{-1 / 2}$ cases of the $W (\textbf{A}, \textbf{B}; \textbf{C}, \textbf{D} \boldsymbol{\mid} r, s)$, so \Cref{rql1} implies that $M_z (\textbf{A}, \textbf{B}; \textbf{C}, \textbf{D}) = 0$ unless $\textbf{B}, \textbf{D} \in \{ \textbf{e}_0, \textbf{e}_1, \ldots , \textbf{e}_n \}$. Hence, in what follows, let us abbreviate $M_z = (\textbf{A}, b; \textbf{C}, d) = M_z (\textbf{A}, \textbf{e}_b; \textbf{C}, \textbf{e}_d)$ for any $b, d \in \{ 0, 1, \ldots , n\}$. These weights are given more explicitly through the following definition.

\begin{definition} 
	
	\label{mqzs}	
	
	For any $\textbf{A}, \textbf{C} \in \mathbb{Z}_{\ge 0}^n$ and $b, d \in \{ 0, 1, \ldots ,n \}$, set $M (\textbf{A}, b; \textbf{C}, d) = 0$ unless $\textbf{A}, \textbf{C} \in \{ 0, 1 \}^n$. Letting $|\textbf{A}| = a$ and $\textbf{A}_i^+, \textbf{A}_j^-, \textbf{A}_{ij}^{+-}$ be as in \eqref{aij} for each $\textbf{A} \in \{ 0, 1 \}^n$ and $i \in [1, n]$, set
	\begin{flalign}
		\label{m1} 
		M (\textbf{A}, 0; \textbf{A}, 0) = & \displaystyle\frac{q^{-a} - s z}{1 - s z}; \quad M (\textbf{A}, 0; \textbf{A}_i^-, i) = q^{-A_{[i + 1, n]}} \displaystyle\frac{1 - q^{-1}}{1 - s z}; \quad  M (\textbf{A}, i; \textbf{A}_i^+, 0) = \displaystyle\frac{z (s^2 - q^{-a})}{1 - s z}.
	\end{flalign} 
	
	\noindent Moreover, for any $1 \le i < j \le n$, set
	\begin{flalign*} 
		M (\textbf{A}, i; \textbf{A}_{ij}^{+-}, j) = q^{-A_{[j + 1, n]}} \displaystyle\frac{s (q^{-1} - 1)}{1 - s z}; \qquad M (\textbf{A}, j; \textbf{A}_{ji}^{+-}, i) = q^{-A_{[i + 1, n]}} \displaystyle\frac{z (q^{-1} - 1)}{1 - s z}.
	\end{flalign*}
	
	\noindent Additionally set
	\begin{flalign}
		\label{m3} 
		& M (\textbf{A}, i; \textbf{A}, i) = q^{-A_{[i + 1, n]}} \displaystyle\frac{z - s}{1 - s z}, \quad \text{if $A_i = 0$}; \qquad M (\textbf{A}, i; \textbf{A}; i) = q^{-A_{[i + 1, n]}} \displaystyle\frac{q^{-1} s - z}{1 - s z}, \quad \text{if $A_i = 1$}.
	\end{flalign}
	
	\noindent We further set $M (\textbf{A}, b; \textbf{C}, d) = 0$ if $(\textbf{A}, b; \textbf{C}, d)$ is not of the above form. 
	
\end{definition}

\begin{figure}[t]

\begin{center}
	
	\begin{tikzpicture}[
		>=stealth,
		scale = .85
		]

		\draw[-, black] (-7.5, 3.1) -- (10, 3.1);
		\draw[-, black] (-7.5, -2.5) -- (10, -2.5);
		\draw[-, black] (-7.5, -1.1) -- (10, -1.1);
		\draw[-, black] (-7.5, -.4) -- (10, -.4);
		\draw[-, black] (-7.5, 2.4) -- (10, 2.4);
		\draw[-, black] (-7.5, -2.5) -- (-7.5, 3.1);
		\draw[-, black] (10, -2.5) -- (10, 3.1);
		\draw[-, black] (7.5, -2.5) -- (7.5, 3.1);
		\draw[-, black] (-5, -2.5) -- (-5, 2.4);
		\draw[-, black] (5, -2.5) -- (5, 3.1);
		\draw[-, black] (-2.5, -2.5) -- (-2.5, 2.4);
		\draw[-, black] (2.5, -2.5) -- (2.5, 2.4);
		\draw[-, black] (0, -2.5) -- (0, 3.1);

		\draw[->, thick, blue] (-6.3, .1) -- (-6.3, 1.9); 
		\draw[->, thick, green] (-6.2, .1) -- (-6.2, 1.9); 
		
		\draw[->, thick, blue] (-3.8, .1) -- (-3.8, 1) -- (-2.85, 1);
		\draw[->, thick, green] (-3.7, .1) -- (-3.7, 1.9);
		
		\draw[->, thick, blue] (-1.35, .1) -- (-1.35, 1.9);
		\draw[->, thick, green] (-1.25, .1) -- (-1.25, 1.9);
		\draw[->, thick,  orange] (-2.15, 1.1) -- (-1.15, 1.1) -- (-1.15, 1.9);
		
		\draw[->, thick, red] (.35, 1) -- (1.15, 1) -- (1.15, 1.9);
		\draw[->, thick, blue] (1.25, .1) -- (1.25, 1.9);
		\draw[->, thick, green] (1.35, .1) -- (1.35, 1.1) -- (2.15, 1.1);
		
		\draw[->, thick, blue] (3.65, .1) -- (3.65, 1) -- (4.65, 1);
		\draw[->, thick, green] (3.75, .1) -- (3.75, 1.9);
		\draw[->, thick, orange] (2.85, 1.1) -- (3.85, 1.1) -- (3.85, 1.9); 
		
		\draw[->, thick, red] (5.35, 1) -- (7.15, 1); 
		\draw[->, thick, blue] (6.2, .1) -- (6.2, 1.9);
		\draw[->, thick, green] (6.3, .1) -- (6.3, 1.9); 
		
		\draw[->, thick, blue] (7.85, 1) -- (9.65, 1); 
		\draw[->, thick, blue] (8.7, .1) -- (8.7, 1.9);
		\draw[->, thick, green] (8.8, .1) -- (8.8, 1.9); 
		
		\filldraw[fill=gray!50!white, draw=black] (-2.85, 1) circle [radius=0] node [black, right = -1, scale = .75] {$i$};
		\filldraw[fill=gray!50!white, draw=black] (2.15, 1) circle [radius=0] node [black, right = -1, scale = .75] {$j$};
		\filldraw[fill=gray!50!white, draw=black] (4.65, 1) circle [radius=0] node [black, right = -1, scale = .75] {$i$};
		\filldraw[fill=gray!50!white, draw=black] (7.15, 1) circle [radius=0] node [black, right = -1, scale = .75] {$i$};
		\filldraw[fill=gray!50!white, draw=black] (9.65, 1) circle [radius=0] node [black, right = -1, scale = .75] {$i$};
		
		\filldraw[fill=gray!50!white, draw=black] (7.85, 1) circle [radius=0] node [black, left = -1, scale = .75] {$i$};
		\filldraw[fill=gray!50!white, draw=black] (5.35, 1) circle [radius=0] node [black, left = -1, scale = .75] {$i$};
		\filldraw[fill=gray!50!white, draw=black] (2.85, 1) circle [radius=0] node [black, left = -1, scale = .75] {$j$};
		\filldraw[fill=gray!50!white, draw=black] (.35, 1) circle [radius=0] node [black, left = -1, scale = .75] {$i$};
		\filldraw[fill=gray!50!white, draw=black] (-2.15, 1) circle [radius=0] node [black, left = -1, scale = .75] {$i$};
		
		\filldraw[fill=gray!50!white, draw=black] (-6.25, 1.9) circle [radius=0] node [black, above = -1, scale = .75] {$\textbf{A}$};
		\filldraw[fill=gray!50!white, draw=black] (-3.75, 1.9) circle [radius=0] node [black, above = -1, scale = .75] {$\textbf{A}_i^-$};
		\filldraw[fill=gray!50!white, draw=black] (-1.25, 1.9) circle [radius=0] node [black, above = -1, scale = .75] {$\textbf{A}_i^+$};
		\filldraw[fill=gray!50!white, draw=black] (1.25, 1.9) circle [radius=0] node [black, above = -1, scale = .75] {$\textbf{A}_{ij}^{+-}$};
		\filldraw[fill=gray!50!white, draw=black] (3.75, 1.9) circle [radius=0] node [black, above = -1, scale = .75] {$\textbf{A}_{ji}^{+-}$};
		\filldraw[fill=gray!50!white, draw=black] (6.25, 1.9) circle [radius=0] node [black, above = -1, scale = .75] {$\textbf{A}$};
		\filldraw[fill=gray!50!white, draw=black] (8.75, 1.9) circle [radius=0] node [black, above = -1, scale = .75] {$\textbf{A}$};

		\filldraw[fill=gray!50!white, draw=black] (-6.25, .1) circle [radius=0] node [black, below = -1, scale = .7] {$\textbf{A}$};
		\filldraw[fill=gray!50!white, draw=black] (-3.75, .1) circle [radius=0] node [black, below = -1, scale = .7] {$\textbf{A}$};
		\filldraw[fill=gray!50!white, draw=black] (-1.25, .1) circle [radius=0] node [black, below = -1, scale = .7] {$\textbf{A}$};
		\filldraw[fill=gray!50!white, draw=black] (1.25, .1) circle [radius=0] node [black, below = -1, scale = .7] {$\textbf{A}$};
		\filldraw[fill=gray!50!white, draw=black] (3.75, .1) circle [radius=0] node [black, below = -1, scale = .7] {$\textbf{A}$};
		\filldraw[fill=gray!50!white, draw=black] (6.25, .1) circle [radius=0] node [black, below = -1, scale = .7] {$\textbf{A}$};
		\filldraw[fill=gray!50!white, draw=black] (8.75, .1) circle [radius=0] node [black, below = -1, scale = .75] {$\textbf{A}$};

		\filldraw[fill=gray!50!white, draw=black] (-3.75, 2.75) circle [radius=0] node [black] {$1 \le i \le n$};
		\filldraw[fill=gray!50!white, draw=black] (2.5, 2.75) circle [radius=0] node [black] {$1 \le i < j \le n$}; 
		\filldraw[fill=gray!50!white, draw=black] (6.25, 2.75) circle [radius=0] node [black] {$A_i = 0$};
		\filldraw[fill=gray!50!white, draw=black] (8.75, 2.75) circle [radius=0] node [black] {$A_i = 1$};
		
		\filldraw[fill=gray!50!white, draw=black] (-6.25, -.75) circle [radius=0] node [black, scale = .9] {$(\textbf{A}, 0; \textbf{A}, 0)$};
		\filldraw[fill=gray!50!white, draw=black] (-3.75, -.75) circle [radius=0] node [black, scale = .9] {$\big( \textbf{A}, 0; \textbf{A}_i^-, i \big)$};
		\filldraw[fill=gray!50!white, draw=black] (-1.25, -.75) circle [radius=0] node [black, scale = .9] {$\big( \textbf{A}, i; \textbf{A}_i^+, 0 \big)$};
		\filldraw[fill=gray!50!white, draw=black] (1.25, -.75) circle [radius=0] node [black, scale = .9] {$\big( \textbf{A}, i; \textbf{A}_{ij}^{+-}, j \big)$};
		\filldraw[fill=gray!50!white, draw=black] (3.75, -.75) circle [radius=0] node [black, scale = .9] {$\big( \textbf{A}, j; \textbf{A}_{ji}^{+-}, i \big)$};
		\filldraw[fill=gray!50!white, draw=black] (6.25, -.75) circle [radius=0] node [black, scale = .9] {$(\textbf{A}, i; \textbf{A}, i)$};
		\filldraw[fill=gray!50!white, draw=black] (8.75, -.75) circle [radius=0] node [black, scale = .9] {$(\textbf{A}, i; \textbf{A}, i)$};
		
		\filldraw[fill=gray!50!white, draw=black] (-6.25, -1.8) circle [radius=0] node [black, scale = .85] {$\displaystyle\frac{q^{-a} -  sz}{1 - sz}$};
		\filldraw[fill=gray!50!white, draw=black] (-3.75, -1.8) circle [radius=0] node [black, scale = .7] {$q^{-A_{[i + 1, n]}} \displaystyle\frac{1 - q^{-1}}{1 - sz}$};
		\filldraw[fill=gray!50!white, draw=black] (-1.25, -1.8) circle [radius=0] node [black, scale = .85] {$\displaystyle\frac{z (s^2 - q^{-a})}{1 - sz}$};
		\filldraw[fill=gray!50!white, draw=black] (1.25, -1.8) circle [radius=0] node [black, scale = .65] {$q^{-A_{[j + 1, n]}} \displaystyle\frac{s (q^{-1} - 1)}{1 - s z}$};
		\filldraw[fill=gray!50!white, draw=black] (3.75, -1.8) circle [radius=0] node [black, scale = .65] {$q^{-A_{[i + 1, n]}} \displaystyle\frac{z (q^{-1} - 1)}{1 - s z}$};
		\filldraw[fill=gray!50!white, draw=black] (6.25, -1.8) circle [radius=0] node [black, scale = .7] {$q^{-A_{[i + 1, n]}} \displaystyle\frac{z - s}{1 - s z}$};
		\filldraw[fill=gray!50!white, draw=black] (8.75, -1.8) circle [radius=0] node [black, scale = .7] {$q^{-A_{[i + 1, n]}} \displaystyle\frac{q^{-1} s - z}{1 - s z}$};

	\end{tikzpicture}
	
\end{center}

\caption{\label{vertexfigurerqm} The vertex weights $M_z$, and their arrow configurations, are depicted above. Here red, blue, green, and orange denote the colors $1$, $2$, $3$, and $4$, respectively.}

\end{figure}

As before, we interpret $M_z (\textbf{A}, b; \textbf{C}, d)$ as the weight associated with a vertex $v$ whose arrow configuration is $(\textbf{A}, \textbf{e}_b; \textbf{C}, \textbf{e}_d)$ and whose spectral parameter is $z$. We refer to \Cref{vertexfigurerqm} for a depiction of these weights; it is quickly verified that \eqref{m} holds for them. 

It will additionally be useful to set notation for ensemble weights (as in \eqref{weighte}). So, let $\textbf{x} = (x_1, x_2, \ldots )$ and $\textbf{y} = (y_1, y_2, \ldots )$ denote (possibly infinite) sequences of complex numbers. For any path ensemble $\mathcal{E}$ on some finite domain $\mathcal{D} \subset \mathbb{Z}_{> 0}^2$, define the ensemble weight of $\mathcal{E}$ with respect to $M$ by 
\begin{flalign}
	\label{eweightm}
M (\mathcal{E} \boldsymbol{\mid} \textbf{x} \boldsymbol{\mid} \textbf{y}) = \displaystyle\prod_{(i, j) \in \mathcal{D}} M_{x_j / y_i} \big( \textbf{A} (i, j), \textbf{B} (i, j); \textbf{C} (i, j), \textbf{D} (i, j) \big),
\end{flalign}

\noindent where $\big( \textbf{A} (v), \textbf{B} (v); \textbf{C} (v), \textbf{D} (v) \big)$ denotes the arrow configuration under $\mathcal{E}$ at any vertex $v \in \mathcal{D}$.

Next, for $K, N \in \mathbb{Z}_{\ge 1}$, we introduce a family of path ensembles on the rectangular domain $\mathcal{D}_{K, N} = \{ 1, 2, \ldots , K \} \times \{ 1, 2, \ldots , N \}$.

\begin{figure}
	
	\begin{center}

		\begin{tikzpicture}[
			>=stealth,
			scale = .65
			]
			
			\draw[ultra thick, gray, dashed] (10, 1) -- (10, 4);
			\draw[ultra thick, gray, dashed] (11, 1) -- (11, 4);
			\draw[ultra thick, gray, dashed] (12, 1) -- (12, 4);
			\draw[ultra thick, gray, dashed] (13, 1) -- (13, 4);
			\draw[ultra thick, gray, dashed] (14, 1) -- (14, 4);
			\draw[ultra thick, gray, dashed] (15, 1) -- (15, 4);
			
			\draw[ultra thick, gray, dashed] (10, 1) -- (15, 1);
			\draw[ultra thick, gray, dashed] (10, 2) -- (15, 2);
			\draw[ultra thick, gray, dashed] (10, 3) -- (15, 3);
			\draw[ultra thick, gray, dashed] (10, 4) -- (15, 4);
			
			\draw[->, thick, red] (15, 0) -- (15, 1);
			\draw[->, thick, blue] (11.95, 0) -- (11.95, 1);
			\draw[->, thick, green] (10, 0) -- (10, 1); 
			\draw[->, thick, orange] (12.05, 0) -- (12.05, 1);
			
			\draw[->, thick, red] (15, 1) -- (16, 1);
			\draw[->, thick, blue] (15, 2) -- (16, 2);
			\draw[->, thick, green] (15, 3) -- (16, 3);
			\draw[->, thick, orange] (15, 4) -- (16, 4);
			
			\draw[->, very thick] (9, 0) -- (9, 4.75);
			\draw[->, very thick] (9, 0) -- (16.5, 0);
			
			\draw[]  (10, -.15) circle [radius = 0] node[below, scale = .55]{$(3)$};
			\draw[]  (12, -.15) circle [radius = 0] node[below, scale = .55]{$(2, 4)$};
			\draw[]  (15, -.15) circle [radius = 0] node[below, scale = .55]{$(1)$};
			
			\draw[]  (16.15, 1) circle [radius = 0] node[right, scale = .55]{$1$};
			\draw[]  (16.15, 2) circle [radius = 0] node[right, scale = .55]{$2$};
			\draw[]  (16.15, 3) circle [radius = 0] node[right, scale = .55]{$3$};
			\draw[]  (16.15, 4) circle [radius = 0] node[right, scale = .55]{$4$};
			
			\draw[] (12.5, 4.75) circle[radius = 0] node[above, scale = .9]{$\mathtt{g}_{\lambda}$};

			\draw[ultra thick, gray, dashed] (21, 1) -- (21, 4);
			\draw[ultra thick, gray, dashed] (22, 1) -- (22, 4);
			\draw[ultra thick, gray, dashed] (23, 1) -- (23, 4);
			\draw[ultra thick, gray, dashed] (24, 1) -- (24, 4);
			\draw[ultra thick, gray, dashed] (25, 1) -- (25, 4);
			\draw[ultra thick, gray, dashed] (26, 1) -- (26, 4);
			
			\draw[ultra thick, gray, dashed] (21, 1) -- (26, 1);
			\draw[ultra thick, gray, dashed] (21, 2) -- (26, 2);
			\draw[ultra thick, gray, dashed] (21, 3) -- (26, 3);
			\draw[ultra thick, gray, dashed] (21, 4) -- (26, 4);

			\draw[->, thick, red] (22.95, 0) -- (22.95, 1);
			\draw[->, thick, blue] (23.05, 0) -- (23.05, 1);
			\draw[->, thick, blue] (21, 0) -- (21, 1);
			\draw[->, thick, red] (26, 0) -- (26, 1);
			
			\draw[->, thick, red] (26, 1) -- (27, 1);
			\draw[->, thick, red] (26, 2) -- (27, 2);
			\draw[->, thick, blue] (26, 3) -- (27, 3);
			\draw[->, thick, blue] (26, 4) -- (27, 4);
			
			\draw[->, very thick] (20, 0) -- (20, 4.75);
			\draw[->, very thick] (20, 0) -- (27.5, 0);
			
			\draw[]  (21, -.15) circle [radius = 0] node[below, scale = .55]{$(2)$};
			\draw[]  (23, -.15) circle [radius = 0] node[below, scale = .55]{$(1, 2)$};
			\draw[]  (26, -.15) circle [radius = 0] node[below, scale = .55]{$(1)$};
			
			\draw[]  (27.15, 1) circle [radius = 0] node[right, scale = .55]{$1$};
			\draw[]  (27.15, 2) circle [radius = 0] node[right, scale = .55]{$1$};
			\draw[]  (27.15, 3) circle [radius = 0] node[right, scale = .55]{$2$};
			\draw[]  (27.15, 4) circle [radius = 0] node[right, scale = .55]{$2$};
			
			\draw[] (23.5, 4.75) circle[radius = 0] node[above, scale = .9]{$\widetilde{\mathtt{F}}_{\boldsymbol{\lambda}}$};
			
		\end{tikzpicture}
		
	\end{center}
	
	\caption{\label{fgpaths3} Depicted to the left and right are the vertex models $\mathfrak{Q}_{\mathtt{g}} (\lambda)$ and $\mathfrak{Q}_{\mathtt{F}} (\boldsymbol{\lambda})$, respectively, for $\lambda = (1, 4, 6, 4)$ we have $\boldsymbol{\lambda} = \big( (2, 0), (4, 3) \big)$. Here red, blue, green, and orange are colors $1$, $2$, $3$, and $4$, respectively.} 
	
\end{figure}
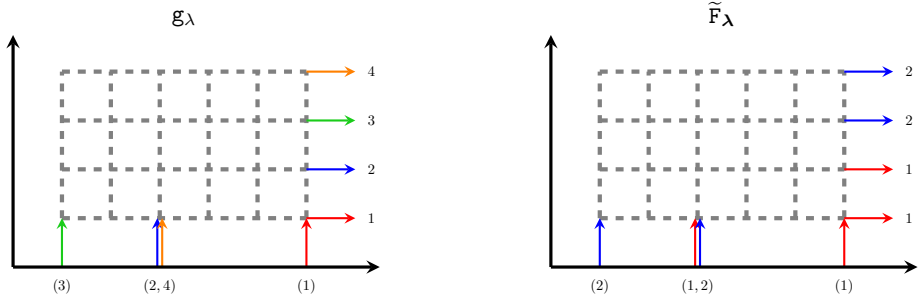

\begin{definition}
	
	\label{qg}

Fix an integer $n \ge 1$; let $\lambda = (\lambda_1, \lambda_2, \ldots , \lambda_n) \in \mathbb{Z}_{> 0}^n$ denote a composition of length $n$; and set $L = L (\lambda) = \max_{j \in [1, n]} \lambda_j$. Let $\mathfrak{Q}_{\mathtt{g}} (\lambda)$\index{Q@$\mathfrak{Q}_{\mathtt{g}} (\lambda)$} denote the set of path ensembles on $\mathcal{D}_{L, n}$ with the following boundary data; we refer to the left side of \Cref{fgpaths3} for a depiction.

\begin{enumerate}
	\item For each $c \in [1, n]$, one color $c$ arrow vertically enters $\mathcal{D}_{L, n}$ at $(L - \lambda_c + 1, 1)$.
	\item For each $c \in [1, n]$, one color $c$ arrow horizontally exits $\mathcal{D}_{L, n}$ at $(L, c)$.
\end{enumerate}

\end{definition}

Now we can define the functions $\mathtt{f}_{\lambda}$ and $\mathtt{g}_{\lambda}$. In what follows, for any sequence $\mathscr{X} = (x_1, x_2, \ldots , x_{\ell})$ we recall $\overleftarrow{\mathscr{X}} = (x_{\ell}, x_{\ell - 1}, \ldots , x_1)$\index{X@$\overleftarrow{\mathscr{X}}$; reverse ordering of $\mathscr{X}$} denotes the reverse ordering of $\mathscr{X}$. We further recall that if $\mathscr{X} \in \mathbb{R}^{\ell}$ that $\inv (\mathscr{X})$ denotes the number of index pairs $(i, j) \in [1, \ell]^2$ such that $i < j$ and $x_i > x_j$.

\begin{definition}
	
	\label{lambdafg}
	
	Fix an integer $n \ge 1$; a sequence of complex numbers $\textbf{x} = (x_1, x_2, \ldots , x_n)$; and a (positive) composition $\lambda = (\lambda_1, \lambda_2, \ldots , \lambda_n) \in \mathbb{Z}_{> 0}^n$ of length $n$. Further let $\textbf{y} = (1, 1, \ldots )$. Define $\mathtt{g}_{\lambda} (\textbf{x}) = g_{\lambda}^{(q)} (\textbf{x} \boldsymbol{\mid} q)$ by setting 
	\begin{flalign}
	\label{xg}
	\mathtt{g}_{\lambda} (\textbf{x}) = \displaystyle\sum_{\mathcal{E} \in \mathfrak{Q}_{\mathtt{g}} (\lambda)} M (\mathcal{E} \boldsymbol{\mid} \textbf{x} \boldsymbol{\mid} \textbf{y}).
	\end{flalign}
	\index{F@$\mathtt{f}_{\mu}^{(q)} (\textbf{x} \boldsymbol{\mid} s)$}
	
	\noindent Stated alternatively, it is the partition function for $\mathfrak{Q}_{\mathtt{g}} (\lambda)$ under the weights $M_{x_j}$ in the $j$-th row. 
	
	Further define $\mathtt{f}_{\lambda} (\textbf{x}) = \mathtt{f}_{\lambda}^{(q)} (\textbf{x} \boldsymbol{\mid} s)$ from $\mathtt{g}$ by 
	\begin{flalign}
		\label{g}
		\mathtt{f}_{\lambda}^{(q)} (\textbf{x} \boldsymbol{\mid} s) = s^n q^{-\inv (\overleftarrow{\lambda})} \mathtt{g}_{\overleftarrow{\lambda}}^{(1 / q)} \big( \overleftarrow{\textbf{x}}^{-1} \boldsymbol{\mid} s^{-1})  \displaystyle\prod_{i = 1}^n \big(  x_i (q - 1) \big)^{-1} \displaystyle\prod_{j = 0}^{\infty} (s^{-2}; q^{-1})_{m_j (\lambda)}.
	\end{flalign}
	\index{G@$\mathtt{g}_{\mu}^{(q)} (\textbf{x} \boldsymbol{\mid} s)$}
	
\end{definition}

In addition to the fuctions $\mathtt{f}_{\lambda}$ and $\mathtt{g}_{\lambda}$, an additional pair of (symmetric) functions given by the below definition will be useful for us. 

\begin{definition}
	
	\label{lambdag} 
	
	Fix an integer $n \ge 1$; a sequence of complex numbers $\textbf{x} = (x_1, x_2, \ldots , x_n)$; and two (positive) compositions $\lambda = (\lambda_1, \lambda_2, \ldots , \lambda_n)$ and $\mu = (\mu_1, \mu_2, \ldots , \mu_n)$ of length $n$. Define the signature sequences 
	\begin{flalign*} 
	\boldsymbol{\lambda} & = \big( (\lambda_1 - 1), (\lambda_2 - 1), \ldots , (\lambda_n - 1) \big) \in \SeqSign_{n; 1}; \\
	 \boldsymbol{\mu} & = \big( (\mu_1 - 1), (\mu_2 - 1), \ldots , (\mu_n - 1) \big) \in \SeqSign_{n; 1},
	\end{flalign*} 

	\noindent as well as the sequences of complex numbers 
	\begin{flalign*} 
	& s^{-1} \textbf{x} = (s^{-1} x_1, s^{-1} x_2, \ldots , s^{-1} x_n); \qquad \textbf{y} = (1, 1, \ldots); \\
	& \textbf{q}^{-1 / 2} = (q^{-1 / 2}, q^{-1 / 2}, \ldots , q^{-1 / 2}); \qquad \textbf{s} = (s, s, \ldots ).
	\end{flalign*} 

	\noindent  Then, define $\widetilde{\mathtt{G}}_{\lambda / \mu} (\textbf{x}) = \widetilde{\mathtt{G}}_{\lambda / \mu}^{(q)} (\textbf{x} \boldsymbol{\mid} s)$ by
	\begin{flalign}
	\label{glambdamu} 
	\widetilde{\mathtt{G}}_{\lambda / \mu} (\textbf{x}) = G_{\boldsymbol{\lambda} / \boldsymbol{\mu}} (s^{-1} \textbf{x}; \textbf{q}^{-1 / 2} \boldsymbol{\mid} \textbf{y}; \textbf{s}),
	\end{flalign}
	\index{G@$\widetilde{\mathtt{G}}_{\lambda / \mu} (\textbf{x})$}

	\noindent and $\mathtt{G}_{\lambda / \mu} (\textbf{x}) = \mathtt{G}_{\lambda / \mu}^{(q)} (\textbf{x} \boldsymbol{\mid} s)$ by
	\begin{flalign}
		\label{g3} 
		\mathtt{G}_{\lambda / \mu} (\textbf{x}) = (-s)^{|\lambda| - |\mu|}  q^{\inv (\mu) - \inv (\lambda)} \widetilde{\mathtt{G}}_{\overleftarrow{\lambda} / \overleftarrow{\mu}}^{(1 / q)} (\overleftarrow{\textbf{x}}^{-1} \boldsymbol{\mid} s^{-1}) \displaystyle\prod_{j = 0}^{\infty} \displaystyle\frac{(s^2; q)_{m_j (\mu)}}{(s^2; q)_{m_j (\lambda)}},
	\end{flalign}
	\index{G@$\mathtt{G}_{\lambda / \mu} (\textbf{x})$}

	\noindent where we recall that $m_j (\nu)$ denote the multiplicity of $j$ in $\nu$. 
\end{definition}

\begin{rem}
	
	\label{fg1}
	
	Let us briefly indicate the relation between the functions from \Cref{lambdafg} and \Cref{lambdag} and some of those introduced in \cite{SVMST}. Here, we let $\lambda$ denote a composition of length $n$ and $\widetilde{\lambda}$ denote the composition obtained from subtracting $1$ from each entry of $\lambda$. 

	First, the $\mathtt{g}_{\lambda} (\textbf{x})$ from \Cref{lambdag} coincides with the $g_{\widetilde{\lambda}} (\textbf{x})$ given by Definition 3.4.6 of \cite{SVMST}. To see this, first observe that since the vertex model $\mathfrak{Q}_{\mathtt{g}} (\lambda)$ from \Cref{qg} contains at most one arrow of any color, the rightmost weight depicted in \Cref{vertexfigurerqm} (with arrow configuration $(\textbf{A}, i; \textbf{A}, i)$ for $A_i = 1$) is irrelevant for evaluating the partition function \eqref{xg} for $\mathtt{g}_{\lambda}$; see also \Cref{wmnwmn}. Omitting this weight from consideration, the remaining ones from \Cref{mqzs} match with those given by equation (2.2.6) of \cite{SVMST} (after reflecting them into the $y$-axis). Thus, the partition function expression for $g_{\widetilde{\lambda}} (\textbf{x})$ given by equation (3.4.10) of \cite{SVMST} matches (after reflection into the vertical line $x = \frac{L + 1}{2}$) the expression \eqref{xg} for $g_{\lambda} (\textbf{x})$. 
	
	Then, equation (1.4.1) of \cite{SVMST} implies that the $\mathtt{f}_{\lambda} (\textbf{x})$ from \Cref{lambdag} coincides with the $f_{\widetilde{\lambda}} (\textbf{x})$ from Definition 3.4.3 of \cite{SVMST} (after taking into account the fact that the function $\inv$ there is defined slightly differently from how it is here). One can further show that the $\mathtt{G}_{\lambda / \mu} (\textbf{x})$ from \Cref{lambdag} coincides with the $G_{\widetilde{\lambda} / \widetilde{\mu}} (\textbf{x})$ given by Definition 4.4.1 of \cite{SVMST}. Given the partition function representation for $G_{\widetilde{\lambda} / \widetilde{\mu}} (\textbf{x})$ from equation (4.4.2) of \cite{SVMST}, the proof of this equality very closely follows that of Proposition 5.6.1 of \cite{SVMST}, so we will not provide it here.
	
\end{rem} 

Now, let us state the following integral formula for $\mathtt{G}_{\lambda / \mu}$ from \cite{SVMST}, which will be useful for us in what follows.

\begin{prop}[{\cite[Equation (9.5.2)]{SVMST}}]
	
	\label{integralg1} 
	
	Fix an integer $N \ge 1$; compositions $\lambda$ and $\mu$ of length $n$; and a sequence of complex numbers $\textbf{\emph{x}} = (x_1, x_2, \ldots , x_N)$. We have that 
	\begin{flalign}
	\label{glambdamu1} 
	\mathtt{G}_{\lambda / \mu} (\textbf{\emph{x}}) = \displaystyle\frac{1}{(2 \pi \textbf{\emph{i}})^n}\displaystyle\frac{q^{\binom{n + 1}{2} -nN}}{(q - 1)^n} \displaystyle\oint \cdots \displaystyle\oint \mathtt{f}_{\mu} (\textbf{\emph{u}}^{-1}) \mathtt{g}_{\lambda} (\textbf{\emph{u}}) \displaystyle\prod_{1 \le i < j \le n} \displaystyle\frac{u_j - u_i}{u_j - q u_i} \displaystyle\prod_{i = 1}^n \displaystyle\prod_{j = 1}^N \displaystyle\frac{u_i - q x_j}{u_i - x_j} \displaystyle\prod_{i = 1}^n \displaystyle\frac{du_i}{u_i},
	\end{flalign}
	
	\noindent where $\textbf{\emph{u}} = (u_1, u_2, \ldots , u_n) \subset \mathbb{C}$, and each $u_i$ is integrated along a positively oriented, closed contour $\gamma_i$ satisfying the following two properties. First, each $\gamma_i$ contains $s$ and the $x_j$, and does not contain $s^{-1}$. Second, the $\{ \gamma_i \}$ are mutually non-intersecting, and $\gamma_{i + 1}$ contains both $\gamma_i$ and $q \gamma_i$ for each $i \in [1, n - 1]$. 
	
\end{prop} 

Let us mention two additional points. First, the value of the right side of \eqref{glambdamu1} does not depend on whether the contours $\gamma_i$ contains $0$ or $\infty$, as it can be verified that those points do not constitute poles of the integrand. Second, whenever interpreting contour integrals such as those on the right side of \eqref{glambdamu1}, we either assume that the underlying parameters enable the existence of such contours or view the integrals as sums of residues (that can be analytically continued to regimes of parameters for which the contours no longer exist).

Although \Cref{integralg1} provides an integral formula for $\mathtt{G}_{\lambda / \mu}$, it will be a bit more useful for us to have one for $\widetilde{\mathtt{G}}_{\lambda / \mu}$. This is done through the following corollary, which will be stated with the parameter $n$ above replaced by $m$ here (since we will eventually apply it for $m = nN$, where $N$ is the number of $\textbf{x}$ variables in the function $G_{\boldsymbol{\lambda} / \boldsymbol{\mu}} (\textbf{x}; \textbf{r} \boldsymbol{\mid} \textbf{y}; \textbf{s})$).

\begin{cor}
	
	\label{integralg}
	
	Fix integers $m, N \ge 1$; compositions $\kappa$ and $\nu$ of length $m$; a complex number $s \in \mathbb{C}$; and a sequence of complex numbers $\textbf{\emph{x}} = (x_1, x_2, \ldots , x_N)$. We have that 
	\begin{flalign}
		\label{integral1} 
		\begin{aligned}
			\widetilde{\mathtt{G}}_{\kappa / \nu} (\textbf{\emph{x}}) = \displaystyle\frac{1}{(2 \pi \textbf{\emph{i}})^m} \displaystyle\frac{(-s)^{|\kappa| - |\nu|} q^{\binom{m + 1}{2}}}{(q - 1)^m} \displaystyle\oint \cdots \displaystyle\oint & \mathtt{g}_{\nu} (\textbf{\emph{u}}) \mathtt{f}_{\kappa} (\textbf{\emph{u}}^{-1})  \displaystyle\prod_{1 \le i < j \le m} \displaystyle\frac{u_j - u_i}{u_j - q u_i} \\
			& \times \displaystyle\prod_{i = 1}^m \displaystyle\prod_{j = 1}^N \displaystyle\frac{1 - q u_i x_j}{1 - u_i x_j} \displaystyle\prod_{i = 1}^m \displaystyle\frac{du_i}{u_i},
		\end{aligned} 
	\end{flalign}
	
	\noindent where $\textbf{\emph{u}} = (u_1, u_2, \ldots , u_m) \subset \mathbb{C}$, and each $u_i$ is integrated along a negatively oriented, closed contour $\gamma_i$ satisfying the following properties. First, each $\gamma_i$ contains $s^{-1}$ and the $x_j^{-1}$, and does not contain $s$. Second, the $\{ \gamma_i \}$ are mutually non-intersecting, and $\gamma_{i - 1}$ contains both $\gamma_i$ and $q^{-1} \gamma_i$ for each $i \in [2, m]$.

\end{cor}

\begin{proof} 
	
	Applying \Cref{integralg1} with the integration variables $\textbf{u}$ there replaced by $\overleftarrow{\textbf{u}}$ here, and the parameters $(x_i, q, s, \lambda, \mu, n)$ there equal to $\big( x_{N - i + 1}^{-1}, q^{-1}, s^{-1}, \overleftarrow{\kappa}, \overleftarrow{\nu}, m \big)$ here, and also using the facts that 
	\begin{flalign*}
		\displaystyle\prod_{1 \le i < j \le m} \displaystyle\frac{u_i - u_j}{u_i - q^{-1} u_j}  = q^{\binom{m}{2}} \displaystyle\prod_{1 \le i < j \le m} \displaystyle\frac{u_j - u_i}{u_j - q u_i} ; \qquad \displaystyle\prod_{i = 1}^m \displaystyle\prod_{j = 1}^N \displaystyle\frac{u_i - q^{-1} x_j^{-1}}{u_i - x_j^{-1}} = q^{-mN} \displaystyle\prod_{i = 1}^m \displaystyle\prod_{j = 1}^N \displaystyle\frac{1 - q u_i x_j}{1 - u_i x_j},
	\end{flalign*}
	
	\noindent we deduce that 
	\begin{flalign*}
		\mathtt{G}_{\overleftarrow{\kappa} / \overleftarrow{\nu}}^{(1 / q)} (\overleftarrow{\textbf{x}}^{-1} \boldsymbol{\mid} s^{-1}) & = \displaystyle\frac{1}{(2 \pi \textbf{i})^m}\displaystyle\frac{1}{(q - 1)^m} \displaystyle\oint \cdots \displaystyle\oint \mathtt{f}_{\overleftarrow{\nu}}^{(1 / q)} \big( \overleftarrow{\textbf{u}}^{-1} \boldsymbol{\mid} s^{-1} \big) \mathtt{g}_{\overleftarrow{\kappa}}^{(1 / q)} \big( \overleftarrow{\textbf{u}} \boldsymbol{\mid} s^{-1} \big) \\
		& \qquad \qquad \qquad \qquad \qquad \quad \times \displaystyle\prod_{1 \le i < j \le m} \displaystyle\frac{u_j - u_i}{u_j - q u_i} \displaystyle\prod_{i = 1}^m \displaystyle\prod_{j = 1}^N \displaystyle\frac{1 - q u_i x_j}{1 - u_i x_j} \displaystyle\prod_{i = 1}^m \displaystyle\frac{du_i}{u_i},
	\end{flalign*}
	
	\noindent where each $u_i$ is integrated along $\gamma_i$. By \eqref{g}, and using the facts that
	\begin{flalign}
		\label{sqsum}
		\begin{aligned}
			& \qquad (s^{-2}; q^{-1})_{m_j (\kappa)} = (-s^{-2})^{m_j (\kappa)} q^{-\binom{m_j (\kappa)}{2}} (s^2; q)_{m_j (\kappa)}; \\
			& \displaystyle\sum_{j = 0}^{\infty} m_j (\kappa) = m; \qquad \inv (\kappa) + \inv \big( \overleftarrow{\kappa} \big) + \displaystyle\sum_{j = 0}^{\infty} \binom{m_j (\kappa)}{2} = \binom{m}{2},
		\end{aligned} 
	\end{flalign} 
	
	\noindent we therefore obtain 
	\begin{flalign*}
		\mathtt{G}_{\overleftarrow{\kappa} / \overleftarrow{\nu}}^{(1 / q)} (\overleftarrow{\textbf{x}}^{-1} \boldsymbol{\mid} s^{-1}) & = \displaystyle\frac{1}{(2 \pi \textbf{i})^m} (-s)^m q^{\binom{m}{2} - \inv (\kappa)} \displaystyle\oint \cdots \displaystyle\oint \mathtt{f}_{\overleftarrow{\nu}}^{(1 / q)} \big( \overleftarrow{\textbf{u}}^{-1} \boldsymbol{\mid} s^{-1} \big) \mathtt{f}_{\kappa}^{(q)} ( \textbf{u}^{-1} \boldsymbol{\mid} s) \\
		& \qquad \qquad \quad \times \displaystyle\prod_{j = 0}^{\infty} \displaystyle\frac{1}{(s^2; q)_{m_j (\kappa)}} \displaystyle\prod_{1 \le i < j \le m} \displaystyle\frac{u_j - u_i}{u_j - q u_i} \displaystyle\prod_{i = 1}^m \displaystyle\prod_{j = 1}^N \displaystyle\frac{1 - q u_i x_j}{1 - u_i x_j} \displaystyle\prod_{i = 1}^m \displaystyle\frac{du_i}{u_i^2}.
	\end{flalign*}

	\noindent Again applying \eqref{g}, it follows that 
	\begin{flalign*}
	\mathtt{G}_{\overleftarrow{\kappa} / \overleftarrow{\nu}}^{(1 / q)} (\overleftarrow{\textbf{x}}^{-1} \boldsymbol{\mid} s^{-1}) & = \displaystyle\frac{1}{(2 \pi \textbf{i})^m} (1 - q^{-1} )^{-m} q^{\binom{m}{2} + \inv (\nu) - \inv (\kappa)} \displaystyle\oint \cdots \displaystyle\oint \mathtt{g}_{\nu}^{(q)} (\textbf{u} \boldsymbol{\mid} s) \mathtt{f}_{\kappa}^{(q)} ( \textbf{u}^{-1} \boldsymbol{\mid} s) \\
	& \qquad \qquad \quad \times \displaystyle\prod_{j = 0}^{\infty} \displaystyle\frac{(s^2; q)_{m_j (\nu)}}{(s^2; q)_{m_j (\kappa)}} \displaystyle\prod_{1 \le i < j \le m} \displaystyle\frac{u_j - u_i}{u_j - q u_i} \displaystyle\prod_{i = 1}^m \displaystyle\prod_{j = 1}^N \displaystyle\frac{1 - q u_i x_j}{1 - u_i x_j} \displaystyle\prod_{i = 1}^m \displaystyle\frac{du_i}{u_i}.	
	\end{flalign*}
	
	\noindent Then we deduce the corollary by first multiplying both sides by 	
	\begin{flalign*} 
		(-s)^{|\kappa| - |\nu|} q^{\inv (\overleftarrow{\nu}) - \inv (\overleftarrow{\kappa})} \displaystyle\prod_{j = 0}^{\infty} \displaystyle\frac{(s^{-2}; q^{-1})_{m_j (\kappa)}}{(s^{-2}; q^{-1})_{m_j (\nu)}} = (-s)^{|\kappa| - |\nu|} q^{\inv (\kappa) - \inv (\nu)} \displaystyle\prod_{j = 0}^{\infty} \displaystyle\frac{(s^2; q)_{m_j (\kappa)}}{(s^2; q)_{m_j (\nu)}},
	\end{flalign*} 
	
	\noindent where to deduce the equality we used \eqref{sqsum}, and then applying \eqref{g3}. 
\end{proof}

\section{Applications of Color Merging} 

\label{Colorfg} 

In this section we provide two applications of the color merging result \Cref{zsumfused}, which will be useful for deriving contour integral representations for $G_{\boldsymbol{\lambda}}$. The first expresses $G_{\boldsymbol{\lambda} / \boldsymbol{\mu}}$ as an anti-symmetrization over $\nu$ of the functions $\widetilde{\mathtt{G}}_{\kappa / \nu}$; together with \Cref{integralg}, this will imply a contour integral formula for $G_{\boldsymbol{\lambda} / \boldsymbol{\mu}}$ involving an anti-symmetrization over $\nu$ of the functions $\mathtt{g}_{\nu}$. The second expresses this anti-symmetrization of $\mathtt{g}_{\nu}$ as another partition function $\widetilde{\mathtt{F}}_{\boldsymbol{\mu}}$. 

To define the latter, we begin with the following vertex model $\mathfrak{Q}_{\mathtt{F}} (\boldsymbol{\lambda})$. It is similar to $\mathfrak{Q}_{\mathtt{g}} (\lambda)$ from \Cref{qg}, but with three differences. First, instead of consisting of $n$ paths all of different colors, the model $\mathfrak{Q}_{\mathtt{F}} (\boldsymbol{\lambda})$ consists of $nM$ paths, with $M$ paths of each color $1, 2, \ldots , n$. Second, along the east edge of the domain, the colors of the horizontally exiting paths are (from bottom to top) given by $1, \ldots , 1, 2, \ldots , 2, \ldots , n, \ldots n$, where each color appears with multiplicity $M$. Third, along the bottom edge of the domain, the boundary data is indexed slightly differently, using the shifts $\mathfrak{T} \big( \lambda^{(i)} \big)$ from \eqref{t}. In this way, the model $\mathfrak{Q}_{\mathtt{F}}$ will be essentially obtained from $\mathfrak{Q}_{\mathtt{g}}$ by merging all colors in intervals of the form $\big\{ (k - 1) M + 1, (k - 1) M + 2, \ldots , kM \big\}$ to a single color $k \in [1, n]$. We refer to \Cref{fgpaths3} for a depiction.

\begin{definition} 
	
	\label{qgqf}

	Fix an integer $n \ge 1$; let $\boldsymbol{\lambda} \in \SeqSign_{n; M}$ denote a signature sequence; and set $K = K (\boldsymbol{\lambda}) = \max_{i \in [1, n]} \max \mathfrak{T} \big( \lambda^{(i)} \big)$, where we recall $\mathfrak{T}$ from \eqref{t}. Let $\mathfrak{Q}_{\mathtt{F}} (\boldsymbol{\lambda})$\index{Q@$\mathfrak{Q}_{\mathtt{F}} (\lambda)$} denote the set of path ensembles on $\mathcal{D}_{K, nM} = [1, K] \times [1, nM]$ with the following boundary data; we refer to the right side of \Cref{fgpaths3} for a depiction. 
	
	\begin{enumerate}
		\item For each $c \in [1, n]$, one color $c$ arrow vertically enters $\mathcal{D}_{K, nM}$ at $(K - \mathfrak{l} + 1, 1)$, for all $\mathfrak{l} \in \mathfrak{T} \big( \lambda^{(c)} \big)$.
		\item For each $c \in [1, n]$ and $j \in [1, M]$, one color $c$ arrow vertically exits $\mathcal{D}_{K, nM}$ at $\big(K, (c - 1) M + j \big)$.
	\end{enumerate}

\end{definition}

Now let us define the partition function $\widetilde{\mathtt{F}}_{\boldsymbol{\lambda}}$. In what follows, we recall the $M$ weights from \Cref{mqzs} and \eqref{eweightm}.

\begin{definition}
	
	\label{qff} 
	
	Fix integers $n, M \ge 1$; a complex number $s \in \mathbb{C}$; two sequences of complex numbers $\textbf{x} = (x_1, x_2, \ldots,  x_{nM})$ and $\textbf{y} = (y_1, y_2, \ldots )$; and a signature sequence $\boldsymbol{\lambda} \in \SeqSign_{n; M}$. Then, define $\widetilde{\mathtt{F}}_{\boldsymbol{\lambda}} (\textbf{x} \boldsymbol{\mid} \textbf{y})$ by setting
	\begin{flalign*}
		\widetilde{\mathtt{F}}_{\boldsymbol{\lambda}} (\textbf{x} \boldsymbol{\mid} \textbf{y}) = \displaystyle\sum_{\mathcal{E} \in \mathfrak{Q}_{\mathtt{F}} (\boldsymbol{\lambda})} M (\mathcal{E} \boldsymbol{\mid} \textbf{x} \boldsymbol{\mid} \textbf{y}).
	\end{flalign*}
	\index{F@$\widetilde{\mathtt{F}}_{\boldsymbol{\lambda}} (\textbf{x} \boldsymbol{\mid} \textbf{y})$}
	
	\noindent Stated alternatively, it is the partition function for the vertex model $\mathfrak{Q}_{\mathtt{F}} (\boldsymbol{\lambda})$ under the weights $M_{x_j / y_i}$ at any vertex $(i, j)$. If $\textbf{y} = (1, 1, \ldots )$, we abbreviate $\widetilde{\mathtt{F}} (\textbf{x}) = \widetilde{\mathtt{F}} (\textbf{x} \boldsymbol{\mid} \textbf{y})$. \index{F@$\widetilde{\mathtt{F}}_{\boldsymbol{\lambda}} (\textbf{x} \boldsymbol{\mid} \textbf{y})$!$\widetilde{\mathtt{F}}_{\boldsymbol{\lambda}} (\textbf{x})$}
	
\end{definition} 

To proceed, we require a certain set of compositions $\Upsilon (\boldsymbol{\lambda})$, which can be interpreted as those compositions $\kappa$ satisfying the following property. The boundary data along the bottom edge for the vertex model $\mathfrak{Q}_{\mathtt{g}} (\kappa)$ reduces to that of $\mathfrak{Q}_{\mathtt{F}} (\boldsymbol{\lambda})$ upon merging all colors in intervals of the form $\big\{ (k - 1) M + 1, (k - 1) M + 2, \ldots , kM \big\}$ and identifying them as color $k$. 

\begin{definition} 
	
	\label{omegalambda} 
	
	For any integers $n, M \ge 1$ and signature sequence $\boldsymbol{\lambda} \in \SeqSign_{n; M}$, let $\Upsilon (\boldsymbol{\lambda})$\index{0@$\Upsilon (\boldsymbol{\lambda})$} denote the set of length $nM$ compositions $\kappa = (\kappa_1, \kappa_2, \ldots , \kappa_{nM})$ such that, for each $i \in [1, n]$, the length $M$ composition $\kappa^{(i)} = \big( \kappa_{(i - 1) M + 1}, \kappa_{(i - 1) M + 2}, \ldots , \kappa_{iM} \big)$ is a permutation of $\mathfrak{T} \big( \lambda^{(i)} \big)$. Moreover, for any $\kappa \in \Upsilon (\boldsymbol{\lambda})$, set
	\begin{flalign*}
		\inv_{\boldsymbol{\lambda}} (\kappa) = \displaystyle\sum_{i = 1}^n \inv \big( \kappa^{(i)} \big) = \displaystyle\sum_{i = 1}^n \displaystyle\sum_{1 \le h < j \le M} \textbf{1}_{\kappa_h^{(i)} > \kappa_j^{(i)}}.
	\end{flalign*}
	\index{I@$\inv$}
	\index{I@$\inv$!$\inv_{\boldsymbol{\lambda}}$}
	
\end{definition}

Now we can state the following anti-symmetrization identities.

\begin{lem} 
	
	\label{gsum} 
	
	Fix integers $n, M, N \ge 1$; a complex number $s$; a sequence of complex numbers $\textbf{\emph{x}} = (x_1, x_2, \ldots , x_{nN})$; two sequences of signatures $\boldsymbol{\lambda}, \boldsymbol{\mu} \in \SeqSign_{n; M}$; and any composition $\kappa \in \Upsilon (\boldsymbol{\lambda})$. Setting $\textbf{\emph{q}}^{-1 / 2} = (q^{-1 / 2}, q^{-1 / 2}, \ldots , q^{-1 / 2})$ (of length $nM$); $\textbf{\emph{y}} = (1, 1, \ldots )$; and $\textbf{\emph{s}} = (s, s, \ldots )$, we have
	\begin{flalign}
	\label{fsumg} 
	\begin{aligned}
		\displaystyle\sum_{\nu \in \Upsilon (\boldsymbol{\mu})} (-1)^{\inv_{\boldsymbol{\mu}} (\nu)} \widetilde{\mathtt{G}}_{\kappa / \nu} (\textbf{\emph{x}}) & = (-1)^{\inv_{\boldsymbol{\lambda}} (\kappa)}  G_{\boldsymbol{\lambda} / \boldsymbol{\mu}} (s^{-1} \textbf{\emph{x}}; \textbf{\emph{q}}^{-1 / 2} \boldsymbol{\mid} \textbf{\emph{y}}; \textbf{\emph{s}}); \\
		\displaystyle\sum_{\nu \in \Upsilon (\boldsymbol{\mu})} (-1)^{\inv_{\boldsymbol{\lambda}}(\nu)} \mathtt{g}_{\nu} (\textbf{\emph{x}}) & = (-1)^{n \binom{M}{2}} \widetilde{\mathtt{F}}_{\boldsymbol{\mu}} (\textbf{\emph{x}}),
	\end{aligned}
	\end{flalign}

	\noindent where for the latter statement we assume $N = M$.
	
\end{lem} 

\begin{proof} 
	
	The proofs of both statements of \eqref{fsumg} will follow in an entirely analogous way from \Cref{zsumfused}; so, let us only establish the first one there. 

	We will deduce it as an application of \Cref{zsumfused}, so let us first match the notation here with that there; in what follows, we let $K = \max_{i \in [1, n]} \max \mathfrak{T} \big( \lambda^{(i)} \big)$. Then, observe that $G_{\boldsymbol{\lambda} / \boldsymbol{\mu}} (s^{-1} \textbf{x}; \textbf{q}^{-1 / 2} \boldsymbol{\mid} \textbf{y}; \textbf{s})$ and $\widetilde{\mathtt{G}}_{\kappa / \nu} (\textbf{x})$ are both partition functions for vertex models on the rectangle $\mathcal{D}_{K, nN} = \{ 1, 2, \ldots , K \} \times \{ 1, 2, \ldots , nN \}$, under the weights $W_{x_j / s} \big(\textbf{A}, \textbf{B}; \textbf{C}, \textbf{D} \boldsymbol{\mid} q^{-1 / 2}, s \big)$ in the $j$-th row. However, the boundary data for these models are slightly different. 
	
	In particular, \Cref{pgpfph} (see also the left side of \Cref{fgpaths}) implies that the model describing $G_{\boldsymbol{\lambda} / \boldsymbol{\mu}}$ has, for each $c \in [1, n]$, one color $c$ arrow vertically entering through $(\mathfrak{m}, 1)$ for every $\mathfrak{m} \in \mathfrak{T} \big( \mu^{(c)} \big)$, and one color $c$ arrow vertically exiting through $(\mathfrak{l}, nN)$ for every $\mathfrak{l} \in \mathfrak{T} \big( \lambda^{(c)}\big)$. Recalling the notation introduced in \Cref{ColorsFused}, we denote this boundary data by $(\mathscr{E}; \mathscr{F}) = \big( \mathscr{E} (\boldsymbol{\mu}); \mathscr{F} (\boldsymbol{\lambda}) \big)$. Similarly, \eqref{glambdamu} and \Cref{pgpfph} together imply that the vertex model describing $\widetilde{\mathtt{G}}_{\kappa / \nu}$ has, for each $c \in [1, nN]$, one color $c$ arrow vertically entering through $(\nu_c, 1)$ and one color $c$ arrow vertically exiting through $(\kappa_c, nN)$. Denote this boundary data by $\big( \breve{\mathscr{E}} (\kappa); \breve{\mathscr{F}} (\nu) \big)$. Thus, defining the sequence of complex numbers $\textbf{z} = \big( z(v) \big)_{v \in \mathcal{D}}$ such that $z(v) = s^{-1} x_j$ for any $v = (i, j) \in \mathcal{D}$, and also recalling the partition function $W_{\mathcal{D}}^{(m; n)} (\mathscr{E}; \mathscr{F} \boldsymbol{\mid} \textbf{z} \boldsymbol{\mid} \textbf{r}, \textbf{s})$ from \Cref{zef1rsz}, we have 
	\begin{flalign}
		\label{gkappanu1} 
		\begin{aligned} 
		G_{\boldsymbol{\lambda} / \boldsymbol{\mu}} (s^{-1} \textbf{x}; \textbf{r} \boldsymbol{\mid} \textbf{y}; \textbf{s}) & = W_{\mathcal{D}}^{(1; n)} \big( \mathscr{E} (\boldsymbol{\mu}); \mathscr{F} (\boldsymbol{\lambda}) \boldsymbol{\mid} \textbf{z} \boldsymbol{\mid} \textbf{q}^{-1 / 2}, \textbf{s} \big); \\ 
		\widetilde{\mathtt{G}}_{\kappa / \nu} (\textbf{x}) & = W_{\mathcal{D}}^{(1; nN)} \big( \breve{\mathscr{E}} (\nu); \breve{\mathscr{F}} (\kappa) \boldsymbol{\mid} \textbf{z} \boldsymbol{\mid} \textbf{q}^{-1 / 2}, \textbf{s} \big).
		\end{aligned}
	\end{flalign}

	Now, recall the notions of interval partitions from \Cref{Colors} and the associated merging prescription $\vartheta_{\mathbb{J}}$ from \Cref{ColorsFused}. Then, \Cref{omegalambda} quickly implies that $\Upsilon (\boldsymbol{\mu})$ is equivalently given by the set of compositions $\nu$ such that $\vartheta_{\mathbb{J}} \big( \breve{\mathscr{E}} (\nu) \big) = \mathscr{E} (\boldsymbol{\mu})$, where $\mathbb{J} = (J_1, J_2, \ldots , J_n)$ is the interval partition of $\{ 1, 2, \ldots , nM \}$ defined by setting $J_i = \big\{ (i - 1) M + 1, (i - 1) M + 2, \ldots , iM \big\}$ for each $i \in [1, n]$. Under this notation, we further have that $\inv_{\boldsymbol{\mu}} (\nu) = \sum_{i = 1}^n \inv \big( \breve{\mathscr{E}} (\nu); J_i \big)$ (and $\inv_{\boldsymbol{\lambda}} (\kappa) = \sum_{i = 1}^n \inv \big( \breve{\mathscr{F}} (\kappa); J_i \big)$), where we recall the latter from \eqref{sumij}. 
	
	Thus, \Cref{zsumfused} gives
	\begin{flalign*} 
		\displaystyle\sum_{\nu \in \Upsilon (\boldsymbol{\mu})} (-1)^{\inv_{\boldsymbol{\mu}} (\nu) - \inv_{\boldsymbol{\lambda}} (\kappa)} W_{\mathcal{D}}^{(1; nN)} \big( \breve{\mathscr{E}} (\nu); \breve{\mathscr{F}} (\kappa) \boldsymbol{\mid} \textbf{z} \boldsymbol{\mid} \textbf{q}^{-1 / 2}, \textbf{s} \big) = W_{\mathcal{D}}^{(1; n)} \big( \mathscr{E} (\boldsymbol{\mu}); \mathscr{F} (\boldsymbol{\lambda}) \boldsymbol{\mid} \textbf{z} \boldsymbol{\mid} \textbf{q}^{-1 / 2}, \textbf{s} \big),
	\end{flalign*}
	
	\noindent which together with \eqref{gkappanu1}, implies the first statement of \eqref{fsumg}. 
	
	As mentioned above, the proof of the second is entirely analogous and is therefore omitted. However, let us briefly indicate the source of the factor $(-1)^{n \binom{M}{2}}$ there. Observe from \Cref{qg} that the boundary data along the right edge of the vertex model $\mathfrak{Q}_{\mathtt{g}} (\nu)$ is (from top to bottom, as in \Cref{ColorsFused}; see \Cref{domainpathsfused}) in the reverse order $nM, nM - 1, \ldots , 1$. Hence, its inversion count under \Cref{omegalambda} is equal to $n \binom{M}{2}$, which accounts for this exponent of $-1$. 
\end{proof}

\section{Integral Formulas for \texorpdfstring{$G_{\boldsymbol{\lambda} / \boldsymbol{\mu}}$}{}}

\label{GIntegral} 

In this section we establish two contour integral formulas for $G_{\boldsymbol{\lambda} / \boldsymbol{\mu}} (\textbf{x})$. The first is given as follows. 

\begin{prop}
	
	\label{glambdamuidentity}
	
	Fix integers $n, M, N \ge 1$; signature sequences $\boldsymbol{\lambda}, \boldsymbol{\mu} \in \SeqSign_{n; M}$; a complex number $s \in \mathbb{C}$; and sequences of complex numbers $\textbf{\emph{r}} = (r_1, r_2, \ldots , r_N)$ and $\textbf{\emph{x}} = (x_1, x_2, \ldots , x_N)$. Denote $\textbf{\emph{s}} = (s, s, \ldots )$ and $\textbf{\emph{y}} = (1, 1, \ldots )$, and let $\kappa \in \Upsilon (\boldsymbol{\lambda})$ denote the unique element of $\Upsilon (\boldsymbol{\lambda})$ such that $\inv_{\boldsymbol{\lambda}} (\kappa) = 0$. Then,
	\begin{flalign}
	\label{gintegral1} 
	\begin{aligned}
	G_{\boldsymbol{\lambda} / \boldsymbol{\mu}} (\textbf{\emph{x}}; \textbf{\emph{r}} \boldsymbol{\mid} \textbf{\emph{y}}; \textbf{\emph{s}}) & = \displaystyle\frac{(-s)^{|\boldsymbol{\lambda}| - |\boldsymbol{\mu}|} }{(2 \pi \textbf{\emph{i}})^{nM}} \displaystyle\frac{q^{\binom{nM + 1}{2}}}{(q - 1)^{nM}} \displaystyle\oint \cdots \displaystyle\oint \mathtt{f}_{\kappa} (\textbf{\emph{u}}^{-1}) \displaystyle\sum_{\nu \in \Upsilon (\boldsymbol{\mu})} (-1)^{\inv_{\boldsymbol{\mu}} (\nu)} \mathtt{g}_{\nu}(\textbf{\emph{u}})  \\
	& \qquad \qquad \qquad \qquad \times \displaystyle\prod_{1 \le i < j \le nM} \displaystyle\frac{u_j - u_i}{u_j - q u_i} \displaystyle\prod_{i = 1}^{nM} \displaystyle\prod_{j = 1}^N \displaystyle\frac{1 - s r_j^{-2} u_i x_j}{1 - s u_i x_j} \displaystyle\prod_{i = 1}^{nM} \displaystyle\frac{du_i}{u_i}, 
	\end{aligned} 
	\end{flalign}
	
	\noindent where $\textbf{\emph{u}} = (u_1, u_2, \ldots , u_{nM}) \subset \mathbb{C}$, and each $u_i$ is integrated along a negatively oriented, closed contour $\gamma_i$ satisfying the following three properties. First, each $\gamma_i$ contains $s^{-1}$ and $q^{-k} s^{-1} x_j^{-1}$ for all integers $k \in [1, nM - 1]$ and $j \in [1, M]$. Second, each $\gamma_i$ does not contain $s$. Third, the $\{ \gamma_i \}$ are mutually non-intersecting, and $\gamma_{i - 1}$ contains both $\gamma_i$ and $q^{-1} \gamma_i$ for each $i \in [2, nM]$. 
	 
\end{prop}

\begin{proof}
	
	Since both sides of \eqref{gintegral1} are rational functions in $\textbf{r}$, it suffices to establish the theorem assuming there exist integers $L_1, L_2, \ldots , L_M \ge 1$ such that $r_i = q^{- L_i / 2}$ for each $i$. 
	
	Let us first assume that each $L_i = 1$, that is, we have $\textbf{r} = \textbf{q}^{-1 / 2}$ where have denoted $\textbf{q}^{-1 / 2} = (q^{-1 / 2}, q^{-1 / 2}, \ldots , q^{-1 / 2})$ (where $q^{-1 / 2}$ appears with multiplicity $N$). Applying \Cref{integralg} with the $(\textbf{x}, \kappa, \nu, m)$ there equal to $(s \textbf{x}, \kappa, \nu, nM)$ here for some fixed $\nu \in \Upsilon (\boldsymbol{\mu})$; multiplying both sides of the resulting \eqref{integralg1} by $(-1)^{\inv_{\boldsymbol{\mu}} (\nu)}$; summing over $\nu \in \Upsilon (\boldsymbol{\mu})$; and applying \Cref{gsum} gives
	\begin{flalign}
	\label{gintegralq12} 
	\begin{aligned}
	G_{\boldsymbol{\lambda} / \boldsymbol{\mu}} \big( \textbf{x}; \textbf{q}^{-1 / 2} \boldsymbol{\mid} \textbf{y}; \textbf{s} \big) & = \displaystyle\frac{(-s)^{|\boldsymbol{\lambda}| - |\boldsymbol{\mu}|}}{(2 \pi \textbf{i})^{nM}} \displaystyle\frac{q^{\binom{nM + 1}{2}}}{(q - 1)^{nM}} \displaystyle\oint \cdots \displaystyle\oint \mathtt{f}_{\kappa} (\textbf{u}^{-1}) \displaystyle\sum_{\nu \in \Upsilon (\boldsymbol{\mu})} (-1)^{\inv_{\boldsymbol{\mu}} (\nu)} \mathtt{g}_{\nu} (\textbf{u})  \\
	& \qquad \qquad \qquad \qquad \quad \times \displaystyle\prod_{1 \le i < j \le nM} \displaystyle\frac{u_j - u_i}{u_j - q u_i} \displaystyle\prod_{i = 1}^{nM} \displaystyle\prod_{j = 1}^N \displaystyle\frac{1 - q s u_i x_j}{1 - s u_i x_j} \displaystyle\prod_{i = 1}^{nM} \displaystyle\frac{du_i}{u_i},
	\end{aligned}
\end{flalign}

\noindent where each $u_i$ is integrated along $\gamma_i$. This establishes \eqref{gintegral1} if $\textbf{r} = \textbf{q}^{-1 / 2}$.

Now assume that each $r_i = q^{-L_i / 2}$ for some integers $L_1, L_2, \ldots , L_N \ge 1$, and set $\textbf{w} = \textbf{w}^{(1)} \cup \textbf{w}^{(2)} \cup \cdots \cup \textbf{w}^{(N)}$, where $\textbf{w}^{(i)} = (x_i, qx_i, \ldots , q^{L - 1} x_i)$ for each $i \in [1, N]$. By \Cref{gq}, we have that $G_{\boldsymbol{\lambda} / \boldsymbol{\mu}} (\textbf{x}; \textbf{r} \boldsymbol{\mid} \textbf{y}; \textbf{s}) = G_{\boldsymbol{\lambda} / \boldsymbol{\mu}} (\textbf{w}; \textbf{q}^{-1 / 2} \boldsymbol{\mid} \textbf{y}; \textbf{s})$, and so applying \eqref{gintegralq12}, with the $\textbf{x}$ there equal to $\textbf{w}$ here, gives
\begin{flalign*}
G_{\boldsymbol{\lambda} / \boldsymbol{\mu}} (\textbf{w}; \textbf{r} \boldsymbol{\mid} \textbf{y}; \textbf{s}) & = \displaystyle\frac{(-s)^{|\boldsymbol{\lambda}| - |\boldsymbol{\mu}|} }{(2 \pi \textbf{i})^{nM}} \displaystyle\frac{q^{\binom{nM + 1}{2}}}{(q - 1)^{nM}} \displaystyle\oint \cdots \displaystyle\oint \mathtt{f}_{\kappa} (\textbf{u}^{-1}) \displaystyle\sum_{\nu \in \Upsilon (\boldsymbol{\mu})} (-1)^{\inv_{\boldsymbol{\mu}} (\nu)} \mathtt{g}_{\nu} (\textbf{u})  \\
& \qquad \qquad \qquad \qquad \quad \times \displaystyle\prod_{1 \le i < j \le nM} \displaystyle\frac{u_j - u_i}{u_j - q u_i} \displaystyle\prod_{i = 1}^{nM} \displaystyle\prod_{j = 1}^N \displaystyle\frac{1 - q^{L_i} s u_i x_j}{1 - s u_i x_j} \displaystyle\prod_{i = 1}^{nM} \displaystyle\frac{du_i}{u_i},
\end{flalign*}

\noindent where again each $u_i$ is integrated along $\gamma_i$. Since $q^{L_i} = r_i^{-2}$, this confirms the theorem if each $r_i = q^{-L_i / 2}$ for some $L_i \in \mathbb{Z}_{\ge 1}$, which (as mentioned previously) establishes it in general.
\end{proof} 

The nonsymmetric functions $\mathtt{f}_{\kappa}$ and $\mathtt{g}_{\nu}$ appearing in the integrand on the right side of \eqref{gintegral1} admit explicit (but elaborate) summation formulas, which are given by Theorem 1.5.5 of \cite{SVMST}. Thus \Cref{glambdamuidentity} provides an in principle fully explicit, although quite intricate, contour integral representation for $G_{\boldsymbol{\lambda} / \boldsymbol{\mu}}$. 

Part of this intricacy is contained in the sum over $\nu \in \Upsilon (\boldsymbol{\mu})$ on the right side of \eqref{gintegral1}. Using the second statement of \eqref{fsumg}, we can express this sum in terms of the partition function $\widetilde{\mathtt{F}}_{\boldsymbol{\mu}}$.

\begin{cor} 
	
\label{gintegralf} 

Adopting the notation of \Cref{glambdamuidentity}, we have that 
\begin{flalign}
\label{gintegral2} 
\begin{aligned}
G_{\boldsymbol{\lambda} / \boldsymbol{\mu}} (\textbf{\emph{x}}; \textbf{\emph{r}} \boldsymbol{\mid} \textbf{\emph{y}}; \textbf{\emph{s}}) = \displaystyle\frac{(-1)^{n \binom{M}{2}}}{(2 \pi \textbf{\emph{i}})^{nM}} \displaystyle\frac{(-s)^{|\boldsymbol{\lambda}| - |\boldsymbol{\mu}|} q^{\binom{nM + 1}{2}}}{(q - 1)^{nM}} \displaystyle\oint \cdots \displaystyle\oint & \mathtt{f}_{\kappa} (\textbf{\emph{u}}^{-1}) \widetilde{\mathtt{F}}_{\boldsymbol{\mu}} (\textbf{\emph{u}}) \displaystyle\prod_{1 \le i < j \le nM} \displaystyle\frac{u_j - u_i}{u_j - q u_i} \\
& \times  \displaystyle\prod_{i = 1}^{nM} \displaystyle\prod_{j = 1}^N \displaystyle\frac{1 - s r_j^{-2} u_i x_j}{1 - s u_i x_j} \displaystyle\prod_{i = 1}^{nM} \displaystyle\frac{du_i}{u_i}.
\end{aligned}
\end{flalign}

\noindent where each $u_i$ is integrated along $\gamma_i$. 

\end{cor} 

\begin{proof}
	
	This follows from inserting the second statement of \eqref{fsumg} into \Cref{glambdamuidentity}.
\end{proof}

\section{Integral Formulas for \texorpdfstring{$G_{\boldsymbol{\lambda}}$}{}}
	
\label{Gn1}

One benefit of \Cref{gintegralf} over \Cref{glambdamuidentity} is that, when $\boldsymbol{\mu} = \boldsymbol{0}^N$, the following proposition indicates that the function $\widetilde{\mathtt{F}}_{\boldsymbol{\mu}} (\textbf{x})$ factors completely. Its proof will appear in \Cref{F0nProduct} below. 

\begin{prop}
	
	\label{f0n1} 
	
	For any integer $N \ge 1$ and set of complex numbers $\textbf{\emph{x}} = (x_1, x_2, \ldots , x_N)$,
	\begin{flalign*}
		\widetilde{\mathtt{F}}_{\boldsymbol{0}^N} (\textbf{\emph{x}}) = q^{-\binom{n}{2} N^2} (1 - q^{-1})^{nN} (s^2; q)_n^{\binom{N}{2}} \displaystyle\prod_{k = 0}^{n - 1} \displaystyle\prod_{1 \le i < j \le N} (q^{-1} x_{Nk + j} - x_{Nk + i}) \displaystyle\prod_{j = 1}^{nN} (1 - s x_j)^{-N}.
	\end{flalign*} 
\end{prop}

Applying \Cref{f0n1} in \Cref{gintegralf}, we deduce the following, more concise, contour integral representation for $G_{\boldsymbol{\lambda}}$. 

\begin{thm}
	
\label{gintegralf2} 

Adopting the notation of \Cref{glambdamuidentity}, we have 
\begin{flalign}
	\label{kappag} 
	\begin{aligned} 
	G_{\boldsymbol{\lambda}} (\textbf{\emph{x}}; \textbf{\emph{r}} \boldsymbol{\mid} \textbf{\emph{y}}; \textbf{\emph{s}}) = \displaystyle\frac{(-s)^{|\boldsymbol{\lambda}|}}{(2 \pi \textbf{\emph{i}})^{nM}} & (s^2; q)_n^{\binom{M}{2}} \displaystyle\oint \cdots \displaystyle\oint \mathtt{f}_{\kappa} (\textbf{\emph{u}}^{-1})   \displaystyle\prod_{i = 1}^{nM} \displaystyle\prod_{j = 1}^N \displaystyle\frac{1 - s r_j^{-2} u_i x_j}{1 - s u_i x_j}  \displaystyle\prod_{j = 1}^{nM} (1 - s u_j)^{-M} \\
	& \quad \times \displaystyle\prod_{1 \le i < j \le nM} \displaystyle\frac{u_j - u_i}{u_j - q u_i}  \displaystyle\prod_{k = 0}^{n - 1} \displaystyle\prod_{1 \le i < j \le M} (q u_{Mk + i} - u_{Mk + j}) \displaystyle\prod_{i = 1}^{nM} \displaystyle\frac{du_i}{u_i},
	\end{aligned} 
\end{flalign}

\noindent where each $u_i$ is integrated along $\gamma_i$. 

\end{thm}

\begin{proof}
	
	Since 
	\begin{flalign*}
	q^{\binom{nM + 1}{2} - \binom{n}{2} M^2} (1 - q^{-1})^{nM} & = q^{n \binom{M}{2}} (q - 1)^{nM}; \\
	 (-q)^{n \binom{M}{2}} \displaystyle\prod_{k = 0}^{n - 1} \displaystyle\prod_{1 \le i < j \le M} (q^{-1} u_{Mk + j} - u_{Mk + i}) & = \displaystyle\prod_{k = 0}^{n - 1} \displaystyle\prod_{1 \le i < j \le M} (q u_{Mk + i} - u_{Mk + j}),
	\end{flalign*}

	\noindent this follows from inserting \Cref{f0n1} into the $\boldsymbol{\mu} = \boldsymbol{0}^M$ case of \Cref{gintegralf}.
\end{proof}

The function $\mathtt{f}_{\kappa} (\textbf{u}^{-1})$ appearing on the right side of \eqref{kappag} can be simplified in the special case when $\kappa$ is \emph{anti-dominant}, meaning that $\kappa_1 \le \kappa_2 \le \cdots \le \kappa_{nM}$. In this case, it is shown in \cite{SVMST} that $\mathtt{f}_{\kappa} (\textbf{w})$ admits an factored form. 

\begin{lem}[{\cite[Proposition 5.1.1]{SVMST}}]
	
	\label{flambdaw}
	
	Fix an integer $M \ge 1$; a sequence of complex numbers $\textbf{\emph{x}} = (x_1, x_2, \ldots , x_M)$; a complex number $s \in \mathbb{C}$; and an anti-dominant composition $\kappa$ of length $M$. We have
	\begin{flalign}
		\label{uflambda}
		\mathtt{f}_{\kappa} (\textbf{\emph{x}}) = \displaystyle\prod_{j = 1}^{\infty} (s^2; q)_{m_j (\kappa)} \displaystyle\prod_{i = 1}^M \displaystyle\frac{1}{1 - sx_i} \left( \displaystyle\frac{x_i - s}{1 - sx_i} \right)^{\kappa_i - 1}.
	\end{flalign}
	
\end{lem} 

This, together with \Cref{gintegralf2}, implies the following simpler expression for $G_{\boldsymbol{\lambda}}$ under a certain ``ordering constraint'' on $\boldsymbol{\lambda}$. 

\begin{thm}
	
	\label{gintegralf3} 
	
	Adopting the notation of \Cref{glambdamuidentity}, assume that $\lambda_1^{(i)} + M - 1 \le \lambda_M^{(j)}$ whenever $1 \le i < j \le n$. For each integer $k \ge 1$, let $m_k \big( \mathfrak{T} (\boldsymbol{\lambda}) \big)$ denote the number of indices $i \in [1, n]$ such that $k \in \mathfrak{T} \big( \lambda^{(i)} \big)$. Then, we have
	\begin{flalign*} 
			G_{\boldsymbol{\lambda}} (\textbf{\emph{x}}; \textbf{\emph{r}} \boldsymbol{\mid} \textbf{\emph{y}}; \textbf{\emph{s}}) & = \displaystyle\frac{(-s)^{|\boldsymbol{\lambda}|}}{(2 \pi \textbf{\emph{i}})^{nM}} (s^2; q)_n^{\binom{M}{2}} \displaystyle\prod_{j = 1}^{\infty} (s^2; q)_{m_j (\mathfrak{T} (\boldsymbol{\lambda}))} \\
			& \quad \times	 \displaystyle\oint \cdots \displaystyle\oint \displaystyle\prod_{i = 1}^n \displaystyle\prod_{j = 1}^M \big( 1 - s u_{iM - j + 1} \big)^{\lambda_j^{(i)} - j} \big( u_{iM - j + 1} - s \big)^{j - \lambda_j^{(i)} - M - 1} \\ 
			& \quad \times \displaystyle\prod_{i = 1}^{nM} \displaystyle\prod_{j = 1}^N \displaystyle\frac{1 - s r_j^{-2} u_i x_j}{1 - s u_i x_j} \displaystyle\prod_{1 \le i < j \le nM} \displaystyle\frac{u_j - u_i}{u_j - q u_i}  \displaystyle\prod_{k = 0}^{n - 1} \displaystyle\prod_{1 \le i < j \le M} (q u_{Mk + i} - u_{Mk + j}) \displaystyle\prod_{i = 1}^{nM} du_i,
	\end{flalign*}
	
	\noindent where each $u_i$ is integrated along $\gamma_i$. 
	
\end{thm}

\begin{proof} 
	
	Let us first show that $\kappa$ is anti-dominant. To that end, observe since $\inv_{\boldsymbol{\lambda}} (\kappa) = 0$ that the sequence $\kappa^{(i)} = \big( \kappa_{(i - 1) M + 1}, \kappa_{(i - 1) M + 2}, \ldots , \kappa_{iM} \big)$ is non-decreasing for each $i \in [1, n]$. Moreover, since $\kappa \in \Upsilon (\boldsymbol{\lambda})$, each $\kappa^{(i)}$ is a permutation of $\mathfrak{T} \big( \lambda^{(i)} \big)$, which implies (as $\kappa^{(i)}$ is non-decreasing) that 
	\begin{flalign}
	\label{kappalambda} 
	\kappa_{(i - 1) M + k} = \lambda_{M - k + 1}^{(i)} + k, \qquad \text{for each $i \in [1, n]$ and $k \in [1, M]$}.
	\end{flalign} 

	\noindent Then, since $\lambda_1^{(i)} + M - 1 \le \lambda_M^{(j)}$ holds if $i < j$, we have for any $1 \le i < j \le n$ and $k, k' \in [1, M]$ that
	\begin{flalign*} 
	\kappa_{(i - 1) M + k} = \lambda_{M - k + 1}^{(i)} + k \le \lambda_1^{(i)} + M \le \lambda_M^{(j)} + 1 \le \lambda_{M - k' + 1}^{(j)} + k' \le \kappa_{(j - 1) M + k'}.
	\end{flalign*} 

	\noindent This, together with the fact that each $\kappa^{(i)}$ is non-decreasing, implies that $\kappa$ is anti-dominant.
	
	Thus, \Cref{flambdaw} applies and, together with \eqref{kappalambda}, yields 
	\begin{flalign*}
	\mathtt{f}_{\kappa} (\textbf{u}^{-1}) & = \displaystyle\prod_{j = 1}^{\infty} (s^2; q)_{m_j (\kappa)} \displaystyle\prod_{j = 1}^{nM} \displaystyle\frac{u_j}{u_j - s} \bigg( \displaystyle\frac{1 - s u_j}{u_j - s} \bigg)^{\kappa_j - 1} \\
	& = \displaystyle\prod_{j = 1}^{\infty} (s^2; q)_{m_j (\mathfrak{T} (\boldsymbol{\lambda}))} \displaystyle\prod_{i = 1}^n \displaystyle\prod_{j = 1}^M \displaystyle\frac{u_{i M - j + 1}}{u_{iM - j + 1}- s} \bigg( \displaystyle\frac{1 - s u_{i M - j + 1}}{u_{i M - j + 1} - s} \bigg)^{\lambda_j^{(i)} + M - j}.
	\end{flalign*} 

	\noindent Upon insertion into \Cref{gintegralf3}, this gives the theorem.
\end{proof}

\section{Degenerations of the Integral Formulas}

\label{DegenerationsG}

In this section we provide various degenerations of the integral formulas from \Cref{Gn1}. We begin with the case when $n = 1$.

\begin{cor}
	
	\label{integraln1g} 
	
	Fix integers $N, M \ge 1$; a signature $\lambda \in \Sign_M$; a complex number $s \in \mathbb{C}$; and sequences of complex numbers $\textbf{\emph{r}} = (r_1, r_2, \ldots , r_N)$ and $\textbf{\emph{x}} = (x_1, x_2, \ldots , x_N)$. Denote $\boldsymbol{\lambda} = (\lambda) \in \SeqSign_{1; M}$; $\textbf{\emph{s}} = (s, s, \ldots)$; and $\textbf{\emph{y}} = (1, 1, \ldots)$. Then, we have that 
	\begin{flalign}
	\label{gn1lambda}
	\begin{aligned}
	G_{\boldsymbol{\lambda}} (\textbf{\emph{x}}; \textbf{\emph{r}} \boldsymbol{\mid} \textbf{\emph{y}}; \textbf{\emph{s}}) = \displaystyle\frac{(-s)^{|\boldsymbol{\lambda}|}}{(2 \pi \textbf{\emph{i}})^M} (1 - s^2)^{\binom{M + 1}{2}} & \displaystyle\oint \cdots \displaystyle\oint \displaystyle\prod_{j = 1}^M (1 - su_{M - j + 1})^{\lambda_j - j} (u_{M - j + 1} - s)^{j - \lambda_j - M - 1} \\
	& \qquad \times \displaystyle\prod_{1 \le i < j \le M} (u_i - u_j) \displaystyle\prod_{i = 1}^M \displaystyle\prod_{j = 1}^N \displaystyle\frac{1 - u_i r_j^{-2} x_j}{1 - u_i x_j} \displaystyle\prod_{i = 1}^M d u_i,
	\end{aligned} 
	\end{flalign}
	
	\noindent where each $u_i$ is integrated along the contour $\gamma_i$ from \Cref{glambdamuidentity}. 
\end{cor}

\begin{proof}
	
	This follows from the $n = 1$ case of \Cref{gintegralf3} since $\prod_{j = 1}^{\infty} (s^2; q)_{m_j (\mathfrak{T} (\lambda))} = (1 - s^2)^M$, which holds since $\lambda$ has length $M$ and since all entries in $\mathfrak{T} (\lambda)$ are mutually distinct.
\end{proof}

\begin{rem}
	
	\label{glambdan1} 
	
	One can use \eqref{gn1lambda} to identify the $n = 1$ case of $G_{\boldsymbol{\lambda}}$ with a generic supersymmetric Schur function, but we will not pursue this here since this $n = 1$ scenario will be addressed in substantially more detail in the forthcoming work \cite{P}. 
	
\end{rem}

Next, we consider the degenerations to the LLT case of \Cref{glambdamuidentity}, \Cref{gintegralf2}, and \Cref{gintegralf3}. We begin with the former, in which case the result reduces to one originally implicitly shown in \cite{AAPP}. In what follows, we recall the function $\psi$ from \eqref{lambdamupsi}. 

\begin{cor}[{\cite[Equation (34)]{AAPP}}]
	
	\label{llambdamuidentity}

	Adopting the notation of \Cref{glambdamuidentity}, we have 
	\begin{flalign}
		\label{lintegral1} 
		\begin{aligned} 
			\mathcal{L}_{\boldsymbol{\lambda} / \boldsymbol{\mu}} (\textbf{\emph{x}}) = \displaystyle\frac{1}{(2 \pi \textbf{\emph{i}})^{nM}} \displaystyle\frac{q^{\psi (\boldsymbol{\mu}) - \psi (\boldsymbol{\lambda})}}{(1 - q^{-1})^{nM}} & \displaystyle\oint \cdots \displaystyle\oint \mathtt{f}_{\kappa}^{(q)} (\textbf{\emph{u}}^{-1} \boldsymbol{\mid} 0) \displaystyle\sum_{\nu \in \Upsilon (\boldsymbol{\mu})} (-1)^{\inv_{\boldsymbol{\mu}} (\nu)} \mathtt{g}_{\nu}^{(q)} (\textbf{\emph{u}} \boldsymbol{\mid} 0)  \\
			& \qquad \times \displaystyle\prod_{1 \le i < j \le nM} \displaystyle\frac{u_j - u_i}{q^{-1} u_j - u_i} \displaystyle\prod_{i = 1}^{nM} \displaystyle\prod_{j = 1}^N \displaystyle\frac{1}{1 - u_i x_j} \displaystyle\prod_{i = 1}^{nM} \displaystyle\frac{du_i}{u_i},  
		\end{aligned} 
	\end{flalign}
	
	\noindent where $\textbf{\emph{u}} = (u_1, u_2, \ldots , u_{nM}) \subset \mathbb{C}$, and each $u_i$ is integrated along a positively oriented, closed contour $\Gamma_i$ satisfying the following two properties. First, each $\Gamma_i$ contains $0$ and does not contain $q^{-k} x_j^{-1}$ for all integers $k \in [1, nM - 1]$ and $j \in [1, M]$. Second, the $\{ \Gamma_i \}$ are mutually non-intersecting, and $\Gamma_{i - 1}$ is contained in both $\Gamma_i$ and $q^{-1} \Gamma_i$ for each $i \in [2, nM]$. 
	
\end{cor}

\begin{proof}

	Define the infinite sequence $\textbf{y} = (1, 1, \ldots )$. By \eqref{1gl}, \eqref{glrlimit}, and the last statement of \eqref{limitg}, we have
	\begin{flalign}
		\label{glintegrallimit}
		\begin{aligned}
			q^{\psi (\boldsymbol{\lambda}) - \psi (\boldsymbol{\mu})} \mathcal{L}_{\boldsymbol{\lambda} / \boldsymbol{\mu}} (\textbf{x}) & = \mathcal{G}_{\boldsymbol{\lambda} / \boldsymbol{\mu}} (\textbf{x}; \infty \boldsymbol{\mid} 0; 0) \\
			& = \displaystyle\lim_{r \rightarrow \infty} \mathcal{G}_{\boldsymbol{\lambda} / \boldsymbol{\mu}} (\textbf{x}; r \boldsymbol{\mid} 0; 0) \\
			& = \displaystyle\lim_{r \rightarrow \infty} \Big( \displaystyle\lim_{s \rightarrow 0} (-s)^{|\boldsymbol{\mu}| - |\boldsymbol{\lambda}|} G_{\boldsymbol{\lambda} / \boldsymbol{\mu}} \big( s^{-1} \textbf{x}; (r, r, \ldots) \boldsymbol{\mid} \textbf{y}; (s, s, \ldots ) \big) \Big).
		\end{aligned} 
	\end{flalign}

	\noindent We then deduce the corollary by inserting this into \eqref{gintegral1} (with the $\textbf{x}$ there replaced by $s^{-1} \textbf{x}$ here and the $s$ there equal to $0$ here), with the contours there reversed.	
\end{proof}

\begin{rem} 

\label{integrall}

Let us briefly explain how \Cref{llambdamuidentity} and equation (34) of \cite{AAPP} are equivalent, for the latter is not directly stated as a contour integral. In what follows, we recall the Schur polynomial $s_{\lambda} (\textbf{z})$\index{S@$s_{\lambda}$; Schur function} associated with any signature $\lambda$ and (possibly infinite) set of variables $\textbf{z} = (z_1, z_2, \ldots )$. Moreover, for any (non)symmetric function $F$ in $\textbf{z}$ and basis $\{ h_{\mu} \}$ for the space of (non)symmetric functions in $\textbf{z}$, we let $\text{Coeff} [F; h_{\lambda}]$ denote the coefficient of $h_{\lambda}$ in the expansion of $F$ over $\{ h_{\mu} \}$. 

First observe by the Cauchy identity for Schur polynomials that 
\begin{flalign*}
	\displaystyle\sum_{\theta \in \Sign_N} s_{\theta} (\textbf{x}) s_{\theta} (\textbf{u}) = \displaystyle\prod_{i = 1}^{nM} \displaystyle\prod_{j = 1}^N \displaystyle\frac{1}{1 - u_i x_j}.
\end{flalign*}
 
 \noindent Upon insertion into \eqref{lintegral1}, this gives 
 \begin{flalign*}
 	\mathcal{L}_{\boldsymbol{\lambda} / \boldsymbol{\mu}} (\textbf{x}) = \displaystyle\frac{q^{\psi (\boldsymbol{\mu}) - \psi (\boldsymbol{\lambda})}}{(1 - q^{-1})^{nM}} \displaystyle\sum_{\theta \in \Sign_N} s_{\theta} (\textbf{x}) \displaystyle\frac{1}{(2 \pi \textbf{i})^{nM}} & \displaystyle\oint \cdots \displaystyle\oint \mathtt{f}_{\kappa}^{(q)} (\textbf{u}^{-1} \boldsymbol{\mid} 0) \displaystyle\prod_{1 \le i < j \le nM} \displaystyle\frac{u_j - u_i}{q^{-1} u_j - u_i} \\
 	& \quad \times s_{\theta} (\textbf{u}) \displaystyle\sum_{\nu \in \Upsilon (\boldsymbol{\mu})} (-1)^{\inv_{\boldsymbol{\mu}} (\nu)} \mathtt{g}_{\nu}^{(q)} (\textbf{u} \boldsymbol{\mid} 0) \displaystyle\prod_{i = 1}^{nM} \displaystyle\frac{du_i}{u_i}.
 \end{flalign*}

\noindent Consequently, for any $\theta \in \Sign_N$ we have 
 \begin{flalign*}
 	\text{Coeff} \big[ \mathcal{L}_{\boldsymbol{\lambda} / \boldsymbol{\mu}} (\textbf{x}); s_{\theta} (\textbf{x}) \big] = \displaystyle\frac{q^{\psi (\boldsymbol{\mu}) - \psi (\boldsymbol{\lambda})}}{(1 - q^{-1})^{nM}} \displaystyle\frac{1}{(2 \pi \textbf{i})^{nM}} & \displaystyle\oint \cdots \displaystyle\oint \mathtt{f}_{\kappa}^{(q)} (\textbf{u}^{-1} \boldsymbol{\mid} 0) \displaystyle\prod_{1 \le i < j \le nM} \displaystyle\frac{u_j - u_i}{q^{-1} u_j - u_i} \\
	& \quad \times s_{\theta} (\textbf{u}) \displaystyle\sum_{\nu \in \Upsilon (\boldsymbol{\mu})} (-1)^{\inv_{\boldsymbol{\mu}} (\nu)} \mathtt{g}_{\nu}^{(q)} (\textbf{u} \boldsymbol{\mid} 0) \displaystyle\prod_{i = 1}^{nM} \displaystyle\frac{du_i}{u_i}.
\end{flalign*} 

\noindent Next, it follows from Theorem 9.4.1 of \cite{SVMST} that for any polynomial $H \in \mathbb{C}(q) [\textbf{u}]$ we have
\begin{flalign*}
\text{Coeff} \Big[ H (\textbf{u}), \mathtt{g}_{\kappa}^{(q)} (\textbf{u} \boldsymbol{\mid} 0)  \Big] = \displaystyle\frac{q^{\binom{nM}{2}}}{(1 - q^{-1})^{nM}} \displaystyle\frac{1}{(2 \pi \textbf{i})^{nM}} & \displaystyle\oint \cdots \displaystyle\oint H (\textbf{u}) \mathtt{f}_{\kappa}^{(q)} (\textbf{u}^{-1} \boldsymbol{\mid} 0) \\
& \quad \times \displaystyle\prod_{1 \le i < j \le nM} \displaystyle\frac{u_j - u_i}{u_j - q u_i} \displaystyle\prod_{i = 1}^{nM} \displaystyle\frac{du_i}{u_i},
\end{flalign*} 

\noindent from which we deduce that 
\begin{flalign}
	\label{gcoefficientl}
	\text{Coeff} \big[ \mathcal{L}_{\boldsymbol{\lambda} / \boldsymbol{\mu}} (\textbf{x}), s_{\theta} (\textbf{x}) \big] = q^{\psi (\boldsymbol{\mu}) - \psi (\boldsymbol{\lambda})} \text{Coeff} \Bigg[ s_{\theta} (\textbf{u}) \displaystyle\sum_{\nu \in \Upsilon (\boldsymbol{\mu})} (-1)^{\inv_{\boldsymbol{\mu}} (\nu)} \mathtt{g}_{\nu}^{(q)} (\textbf{u} \boldsymbol{\mid} 0); \mathtt{g}_{\kappa}^{(q)} (\textbf{u} \boldsymbol{\mid} 0) \Bigg].
\end{flalign} 

\noindent Now, Theorem 5.9.4 of \cite{SVMST} indicates that $\mathtt{g}_{\lambda}^{(q)} (\textbf{x} \boldsymbol{\mid} 0)$ is given, after a suitable normalization and change of variables, by a nonsymmetric Hall--Littlewood polynomial (namely, the $q = 0$ degenerations of the nonsymmetric Macdonald polynomials introduced in \cite{AAOP}). Under this identification, \eqref{gcoefficientl} coincides with equation (34) of \cite{AAPP}. 

\end{rem}

We next have the following corollary of \Cref{gintegralf2}, which simplifies \Cref{llambdamuidentity} in the case when $\boldsymbol{\mu} = \boldsymbol{0}^M$.

\begin{cor}
	
	\label{lintegralf2} 
	
	Adopting the notation of \Cref{glambdamuidentity}, we have 
	\begin{flalign*}
			\mathcal{L}_{\boldsymbol{\lambda}} (\textbf{\emph{x}}) = \displaystyle\frac{q^{\binom{n}{2} \binom{N}{2} / 2 - \psi (\boldsymbol{\lambda})}}{(2 \pi \textbf{\emph{i}})^{nM}} \displaystyle\oint \cdots \displaystyle\oint & \displaystyle\prod_{1 \le i < j \le nM} \displaystyle\frac{u_j - u_i}{u_j - q u_i}  \displaystyle\prod_{k = 0}^{n - 1} \displaystyle\prod_{1 \le i < j \le M} (q u_{Mk + i} - u_{Mk + j}) \\
			& \times \mathtt{f}_{\kappa}^{(q)} (\textbf{\emph{u}}^{-1} \boldsymbol{\mid} 0)  \displaystyle\prod_{i = 1}^{nM} \displaystyle\prod_{j = 1}^N \displaystyle\frac{1}{1 - u_i x_j} \displaystyle\prod_{i = 1}^{nM} \displaystyle\frac{du_i}{u_i},
	\end{flalign*}
	 
	 \noindent where each $u_i$ is integrated along the $\Gamma_i$ from \Cref{llambdamuidentity}.
	 
\end{cor}

\begin{proof}
	
	This follows from inserting the $\boldsymbol{\mu} = \boldsymbol{0}^N$ case of \eqref{glintegrallimit} into \Cref{gintegralf2} (with the $\textbf{x}$ there replaced by $s^{-1} \textbf{x}$ here, and the $s$ there set to $0$ here), with the contours there reversed, and \Cref{psi0n}. 
\end{proof} 

Under the ordering constraint described by \Cref{gintegralf3}, this nonsymmetric Hall--Littlewood function factors completely, giving rise to the following simplified integral formulas for certain LLT polynomials.

\begin{cor}
	
	\label{lintegralf3} 
	
	Adopting the notation of \Cref{glambdamuidentity}, assume that $\lambda_1^{(i)} + M - 1 \le \lambda_M^{(j)}$ whenever $1 \le i < j \le n$. Then, we have
	\begin{flalign*} 
		\mathcal{L}_{\boldsymbol{\lambda}} (\textbf{\emph{x}}) = \displaystyle\frac{q^{\binom{n}{2} \binom{N}{2} / 2 - \psi (\boldsymbol{\lambda})}}{(2 \pi \textbf{\emph{i}})^{nM}} \displaystyle\oint \cdots \displaystyle\oint & \displaystyle\prod_{1 \le i < j \le nM} \displaystyle\frac{u_j - u_i}{u_j - q u_i}  \displaystyle\prod_{k = 0}^{n - 1} \displaystyle\prod_{1 \le i < j \le M} (q u_{Mk + i} - u_{Mk + j})  \\
		& \times \displaystyle\prod_{i = 1}^n \displaystyle\prod_{j = 1}^M u_{iM - j + 1}^{j - \lambda_j^{(i)} - M - 1} \displaystyle\prod_{i = 1}^{nM} \displaystyle\prod_{j = 1}^N \displaystyle\frac{1}{1 - u_i x_j}  \displaystyle\prod_{i = 1}^{nM} du_i,
	\end{flalign*}
	
	\noindent where each $u_i$ is integrated along the $\Gamma_i$ from \Cref{llambdamuidentity}.
	
\end{cor}

\begin{proof} 
	
	This follows from inserting the $\boldsymbol{\mu} = \boldsymbol{0}^M$ case of \eqref{glintegrallimit} into \Cref{gintegralf3} (with the $\textbf{x}$ there replaced by $s^{-1} \textbf{x}$ here, and the $s$ there set to $0$ here), with the contours there reversed, and \Cref{psi0n}.
\end{proof}

\section{Proof of \Cref{f0n1}}

\label{F0nProduct}

In this section we establish \Cref{f0n1}, which provides an explicit form for $\widetilde{\mathtt{F}}_{\boldsymbol{0}^N} (\textbf{x})$. To that end, it will in fact be useful to analyze the more general quantity given by $\widetilde{\mathtt{F}}_{\boldsymbol{0}^N} (\textbf{x} \boldsymbol{\mid} \textbf{y})$ from \Cref{qff}. So, throughout, we fix integers $n, N \ge 1$ and sequences of complex numbers $\textbf{x} = (x_1, x_2, \ldots , x_N)$ and $\textbf{y} = (y_1, y_2, \ldots )$. 

We will establish the following proposition evaluating $\widetilde{\mathtt{F}}_{\boldsymbol{0}^N} (\textbf{x} \boldsymbol{\mid} \textbf{y})$, from which \Cref{f0n1} directly follows. 

\begin{prop}

\label{xyf0}

We have that
\begin{flalign*}
	\widetilde{\mathtt{F}}_{\boldsymbol{0}^N} (\textbf{\emph{x}} \boldsymbol{\mid} \textbf{\emph{y}}) & = q^{-\binom{n}{2} N^2} (1 - q^{-1})^{nN} \displaystyle\prod_{k = 0}^{n - 1} \displaystyle\prod_{1 \le i < j \le N} (q^{-1} x_{Nk + j} - x_{Nk + i}) (y_i - q^k s^2 y_j) \\
	& \qquad \times \displaystyle\prod_{i = 1}^N y_i^n \displaystyle\prod_{i = 1}^N \displaystyle\prod_{j = 1}^{nN} (y_i - s x_j)^{-1}.
\end{flalign*}

\end{prop}

\begin{proof}[Proof of \Cref{f0n1} Assuming \Cref{xyf0}]
	
	This follows from setting each $y_i$ equal to $1$ in \Cref{xyf0}.
\end{proof}

The proof of \Cref{xyf0} will be similar to that of \Cref{f0} in \Cref{ProofF0}, by (after suitable normalization) realizing $\widetilde{\mathtt{F}}_{\boldsymbol{0}^N} (\textbf{x} \boldsymbol{\mid} \textbf{y})$ as a polynomial that is characterized by a specific set of zeroes guaranteed by the Yang--Baxter equation. To that end, we have the following lemma. 

\begin{lem}
	
	\label{xyf0c} 
	
	There exists a constant $C = C_{n; N} (s; q)$ such that 
	\begin{flalign*}
	\widetilde{\mathtt{F}}_{\boldsymbol{0}^N} (\textbf{\emph{x}} \boldsymbol{\mid} \textbf{\emph{y}}) = C \displaystyle\prod_{k = 0}^{n - 1} \displaystyle\prod_{1 \le i < j \le N} (x_{Nk + j} - q x_{Nk + i}) (y_i - q^k s^2 y_j) \displaystyle\prod_{i = 1}^N y_i^n \displaystyle\prod_{i = 1}^N \displaystyle\prod_{j = 1}^{nN} (y_i - s x_j)^{-1}.
	\end{flalign*} 
\end{lem} 

\begin{proof}[Proof (Outline)]
	
	Since this proof is similar to that of \Cref{f0nxy}, we only outline it. 
	
	To that end, following \Cref{f0np1}, we first claim that  
	\begin{flalign} 
	\label{zf}
	\mathtt{Z} (\textbf{x} \boldsymbol{\mid} \textbf{y}) = \widetilde{\mathtt{F}}_{\boldsymbol{0}^N} \displaystyle\prod_{i = 1}^N \displaystyle\prod_{j = 1}^{nN} (y_i - s x_j),
	\end{flalign}

	\noindent is a polynomial in $(\textbf{x}, \textbf{y})$ of total degree at most $nN^2$. To see this, observe from \Cref{qff} that $\mathtt{Z} (\textbf{x} \boldsymbol{\mid} \textbf{y})$ is the partition function for the model $\mathfrak{Q}_{\mathtt{F}} (\boldsymbol{0}^N)$ from \Cref{qgqf}, whose weight at any vertex $(i, j)$ in the domain $\mathcal{D} = \mathcal{D}_{N, nN} = [1, N] \times [1, nN]$ is given by
	\begin{flalign*}
	\widetilde{M}_{x_j; y_i} (\textbf{A}, b; \textbf{C}, d) = (y_i - sx_j) M_{x_j / y_i} (\textbf{A}, b; \textbf{C}, d). 
	\end{flalign*} 

	\noindent By \Cref{mqzs}, these weights are always linear in $(\textbf{x}, \textbf{y})$, from which it follows that $\mathtt{Z} (\textbf{x} \boldsymbol{\mid} \textbf{y})$ is a polynomial in $(\textbf{x}, \textbf{y})$ of total degree at most $|\mathcal{D}| = nN^2$. 
	
	Next, we identify a divisibility property for the polynomial $\mathtt{Z}$, which will follow from certain exchange relations. To that end, recalling the transposition $\mathfrak{s}_i$ of $(i, i + 1)$, we have
	\begin{flalign}
		\label{fxysigma} 
		\begin{aligned} 
		\widetilde{\mathtt{F}}_{\boldsymbol{0}^N} (\textbf{x} \boldsymbol{\mid} \textbf{y} ) & = \widetilde{\mathtt{F}}_{\boldsymbol{0}^N} \big( \textbf{x} \boldsymbol{\mid} \mathfrak{s}_i (\textbf{y}) \big) \displaystyle\prod_{j = 0}^{n - 1} \displaystyle\frac{y_i - q^j s^2 y_{i + 1}}{y_{i + 1} - q^j s^2 y_i}, \qquad \text{for $i \in [1, N]$}; \\
		\widetilde{\mathtt{F}}_{\boldsymbol{0}^N} (\textbf{x} \boldsymbol{\mid} \textbf{y}) & = \widetilde{\mathtt{F}}_{\boldsymbol{0}^N} \big( \mathfrak{s}_i (\textbf{x}) \boldsymbol{\mid} \textbf{y} \big) \displaystyle\frac{x_{i + 1} - q x_i}{x_i - q x_{i + 1}}, \qquad \text{for $i \in [1, nN]$ with $i \notin \{ N, 2N, \ldots , nN \}$}.
		\end{aligned} 
	\end{flalign}

	\noindent We omit the proof of \eqref{fxysigma}, since its first and second statements follow from applications of the Yang--Baxter equation very similar to the ones used to show (the $s_j = s$ case of) \Cref{f0nxyy} and (the $r_j = q^{-1 / 2}$ case of) \Cref{bd1limit}, respectively. 
	
	Thus, by \eqref{fxysigma}, \eqref{zf}, and the polynomiality of $\mathtt{Z}$, we deduce that 
	\begin{flalign*}
		\displaystyle\prod_{k = 0}^{n - 1} \displaystyle\prod_{1 \le i < j \le N} (x_{Nk + j} - q x_{Nk + i}) \displaystyle\prod_{k = 0}^{n - 1} \displaystyle\prod_{1 \le i < j \le N} (y_i - q^k s^2 y_j) \quad \text{divides} \quad  \mathtt{Z} (\textbf{x} \boldsymbol{\mid} \textbf{y}).
	\end{flalign*}
	
	\noindent We further claim for each $i \in [1, N]$ that $y_i^n$ divides $\mathtt{Z} (\textbf{x} \boldsymbol{\mid} \textbf{y})$. To see this, observe that any path ensemble $\mathcal{E} \in \mathfrak{Q}_{\mathtt{F}} (\boldsymbol{0}^N)$ satisfies the following property. For each $i \in [1, N]$, there exist $n (i - 1)$ arrows horizontally entering the $i$-th column of $\mathcal{E}$ and $in$ arrows horizontally exiting it (see the left side of \Cref{pathsf0}). Thus, this $i$-th column contains at least $n$ distinct vertices whose arrow configurations are of the form $(\textbf{A}, 0; \textbf{A}_i^-, i)$. Since \eqref{m1} implies that $y_i$ divides the weight of this configuration under $\widetilde{M}_{x; y_i}$, it follows that $y_i^n$ divides $\mathtt{Z} (\textbf{x} \boldsymbol{\mid} \textbf{y})$. 
	
	Hence, there exists a constant $C = C_{n; N} (s; q)$ such that
	\begin{flalign*} 
		\mathtt{Z} (\textbf{x} \boldsymbol{\mid} \textbf{y}) = C \displaystyle\prod_{k = 0}^{n - 1} \displaystyle\prod_{1 \le i < j \le N} (x_{Nk + j} - q x_{Nk + i}) \displaystyle\prod_{k = 0}^{n - 1} \displaystyle\prod_{1 \le i < j \le N} (y_i - q^k s^2 y_j) \displaystyle\prod_{i = 1}^N y_i^n,
	\end{flalign*}
	
	\noindent since both sides are of degree $n N^2$, and the right side divides the left. This, together with \eqref{zf}, implies the lemma.
\end{proof}

Now we can establish \Cref{xyf0}. 

\begin{figure}

	\begin{center}
		
		\begin{tikzpicture}[
		>=stealth,
		scale = .75
		]
		
		\draw[dashed] (11, 1) -- (11, 6);
		\draw[dashed] (12, 1) -- (12, 6);
		\draw[dashed] (13, 1) -- (13, 6);
		
		\draw[dashed] (11, 1) -- (13, 1);
		\draw[dashed] (11, 2) -- (13, 2);
		\draw[dashed] (11, 3) -- (13, 3);
		\draw[dashed] (11, 4) -- (13, 4);
		\draw[dashed] (11, 5) -- (13, 5);
		\draw[dashed] (11, 6) -- (13, 6);
		
		\draw[->, thick, red] (10.95, 0) -- (10.95, 1);
		\draw[->, thick, blue] (11.05, 0) -- (11.05, 1);
		\draw[->, thick, red] (11.95, 0) -- (11.95, 1);
		\draw[->, thick, blue] (12.05, 0) -- (12.05, 1);
		\draw[->, thick, red] (12.95, 0) -- (12.95, 1);
		\draw[->, thick, blue] (13.05, 0) -- (13.05, 1);
		
		\draw[->, thick, red] (13, 1) -- (14, 1);
		\draw[->, thick, red] (13, 2) -- (14, 2);
		\draw[->, thick, red] (13, 3) -- (14, 3);
		\draw[->, thick, blue] (13, 4) -- (14, 4);
		\draw[->, thick, blue] (13, 5) -- (14, 5);
		\draw[->, thick, blue] (13, 6) -- (14, 6);
		
		\draw[->, very thick] (10, 0) -- (10, 6.75);
		\draw[->, very thick] (10, 0) -- (14.75, 0);
		
		\draw[]  (11, 0) circle [radius = 0] node[below = 2, scale = .8]{$1$};
		\draw[]  (12, 0) circle [radius = 0] node[below = 2, scale = .8]{$2$};
		\draw[]  (13, 0) circle [radius = 0] node[below = 2, scale = .8]{$3$};
		
		\draw[]  (10, 1) circle [radius = 0] node[left = 2, scale = .8]{$1$};
		\draw[]  (10, 2) circle [radius = 0] node[left = 2, scale = .8]{$2$};
		\draw[]  (10, 3) circle [radius = 0] node[left = 2, scale = .8]{$3$};
		\draw[]  (10, 4) circle [radius = 0] node[left = 2, scale = .8]{$4$};
		\draw[]  (10, 5) circle [radius = 0] node[left = 2, scale = .8]{$5$};
		\draw[]  (10, 6) circle [radius = 0] node[left = 2, scale = .8]{$6$};
		
		\draw[->, very thick] (20, 0) -- (24.75, 0);
		\draw[->, very thick] (20, 0) -- (20, 6.75);

		\draw (21, 0) circle[radius = 0] node[below = 2, black, scale = .8]{$1$};
		\draw (22, 0) circle [radius = 0] node[below = 2, black, scale = .8]{$2$};
		\draw (23, 0) circle [radius = 0] node[below = 2, black, scale = .8]{$3$};
		
		\draw[thick, red, ->] (20.95, 0) -- (20.95, 3) -- (24, 3);
		\draw[thick, red, ->] (21.95, 0) -- (21.95, 2) -- (24, 2);
		\draw[thick, red, ->] (22.95, 0) -- (22.95, 1) -- (24, 1);
		
		\draw[thick, blue, ->] (21.05, 0) -- (21.05, 6) -- (24, 6);
		\draw[thick, blue, ->] (22.05, 0) -- (22.05, 5) -- (24, 5);
		\draw[thick, blue, ->] (23.05, 0) -- (23.05, 4) -- (24, 4);
		
		\draw[] (20, 1) circle[radius = 0] node[left = 2, scale = .8]{$1$};
		\draw[] (20, 2) circle[radius = 0] node[left = 2, scale = .8]{$2$};
		\draw[] (20, 3) circle[radius = 0] node[left = 2, scale = .8]{$3$};
		\draw[] (20, 4) circle[radius = 0] node[left = 2, scale = .8]{$4$};
		\draw[] (20, 5) circle[radius = 0] node[left = 2, scale = .8]{$5$};
		\draw[] (20, 6) circle[radius = 0] node[left = 2, scale = .8]{$6$};
		
		\end{tikzpicture}
		
	\end{center}
	
	\caption{\label{pathsf0} Shown to the left is the vertex model for $\widetilde{\mathtt{F}}_{\boldsymbol{0}^M} (\textbf{x} \boldsymbol{\mid} \textbf{y})$. Shown to the right is the frozen path ensemble corresponding to this function when each $x_j = q^{-1} s$ and $y_i = 1$. } 
	
\end{figure}

\begin{proof}[Proof of \Cref{xyf0}]
	
	In view of \Cref{xyf0c}, it suffices to determine $C_{n; N} (s; q)$. To that end, as in the proof of \Cref{f0rxsy}, we will use a special choice for the parameters $(\textbf{x}, \textbf{y})$ that freezes the partition function $\widetilde{\mathtt{F}}_{\boldsymbol{0}^N} (\textbf{x} \boldsymbol{\mid} \textbf{y})$, enabling us to evaluate it. 
	
	Let us set $x_j = q^{-1} s$ for each $j \in [1, nN]$ and $y_i = 1$ for each $i \in [1, n]$. Under this specialization, \Cref{xyf0c} yields  
\begin{flalign}
	\label{fz1}
	\widetilde{\mathtt{F}}_{\boldsymbol{0}^N} (\textbf{x} \boldsymbol{\mid} \textbf{y}) = C (q^{-1} s)^{n \binom{N}{2}} (1 - q)^{n \binom{N}{2}} (s^2; q)_n^{\binom{N}{2}} (1 - q^{-1} s^2)^{-nN^2}.
\end{flalign} 

To evaluate $\widetilde{\mathtt{F}}_{\boldsymbol{0}^N}$ directly under this specialization observe, since $(x_j, y_i) = (q^{-1} s, 1)$, that \eqref{m3} gives $M_{x_j / y_i, s} (\textbf{A}, h; \textbf{A}; h) = 0$ whenever $A_h = 1$. From this, it is quickly verified that there is a unique path ensemble with nonzero weight in the vertex model $\mathfrak{Q}_{\mathtt{F}} (\boldsymbol{0}^N)$ from \Cref{qgqf}. It is the one depicted on the right side of \Cref{pathsf0}, where the path of color $k$ entering the domain $\mathcal{D} = \mathcal{D}_{N, nN} = [1, N] \times [1, nN]$ at $(j, 1)$ proceeds as north until it reaches $(j, kN - j + 1)$, and then proceeds east until it exits $\mathcal{D}$ at $(N, kN - j + 1)$. In particular, under this ensemble, for any $k \in [1, n]$ the arrow configuration at any vertex $(i, kN - j + 1) \in \mathcal{D}$ is $(\textbf{e}_{[k, n]}, 0; \textbf{e}_{[k, n]}, 0)$ if $1 \le i < j \le N$; $(\textbf{e}_{[k, n]}, 0; \textbf{e}_{[k + 1, n]}, k)$ if $1 \le i = j \le N$; and is $(\textbf{e}_{[k + 1, n]}, k; \textbf{e}_{[k + 1, n]}, k)$ if $1 \le j < i \le N$. 

This gives
\begin{flalign}
	\label{xf0}
	\begin{aligned} 
	\widetilde{\mathtt{F}}_{\boldsymbol{0}^N} (\textbf{x} \boldsymbol{\mid} \textbf{y}) & = \displaystyle\prod_{k = 1}^n \displaystyle\prod_{i = 1}^N M_{s / q} \big( \textbf{e}_{[k, n]}, 0; \textbf{e}_{[k, n]}, 0 \big)^{N - i} \displaystyle\prod_{k = 1}^n \displaystyle\prod_{j = 1}^N M_{s / q} \big( \textbf{e}_{[k, n]}, 0; \textbf{e}_{[k + 1, n]}, k \big) \\
	& \qquad \times \displaystyle\prod_{k = 1}^n \displaystyle\prod_{i = 1}^N M_{s / q} \big( \textbf{e}_{[k + 1, n]}, k; \textbf{e}_{[k + 1, n]}, k \big)^{i - 1}.
	\end{aligned}
\end{flalign}

\noindent By \eqref{m1} and \eqref{m3}, we have for any $k \in [1, n]$ that 
\begin{flalign*}
	M_{s / q} \big( \textbf{e}_{[k, n]}, 0; \textbf{e}_{[k, n]}, 0 \big) & = q^{k - n - 1} \displaystyle\frac{1 - q^{n - k} s^2}{1 - q^{-1} s^2}; \qquad M_{s / q} \big( \textbf{e}_{[k, n]}, 0; \textbf{e}_{[k + 1, n]}, k \big) = q^{k - n - 1} \displaystyle\frac{q - 1}{1 - q^{-1} s^2}; \\
	& M_{s / q} \big( \textbf{e}_{[k + 1, n]}, k; \textbf{e}_{[k + 1, n]}, k \big)  = q^{k - n - 1} \displaystyle\frac{s (1 - q)}{1 - q^{-1} s^2},
\end{flalign*}

\noindent which by \eqref{xf0} gives
\begin{flalign*}
	\widetilde{\mathtt{F}}_{\boldsymbol{0}^N} (\textbf{x} \boldsymbol{\mid} \textbf{y}) = q^{- \binom{n + 1}{2} N^2} (-s)^{n \binom{N}{2}} (q - 1)^{n \binom{N + 1}{2}} (s^2; q)_n^{\binom{N}{2}} (1 - q^{-1} s^2)^{-nN^2}.
\end{flalign*}

\noindent Comparing this with \eqref{fz1} yields 
\begin{flalign*}
	C = q^{n \binom{N}{2} - \binom{n + 1}{2} N^2} (q - 1)^{nN},
\end{flalign*}

\noindent which implies the proposition upon insertion into \Cref{xyf0c}. 
\end{proof}

\chapter{Vanishing Properties}

\label{HFunction0}

Since \eqref{1hl} implies that the polynomials $\mathcal{H}_{\boldsymbol{\lambda}} (\textbf{x}; \infty \boldsymbol{\mid} \textbf{y}; \infty)$ from \Cref{sgf0} degenerate to the LLT ones as each $y_i$ tends to $\infty$, we may view the $\mathcal{H}_{\boldsymbol{\lambda}} (\textbf{x}; \infty \boldsymbol{\mid} \textbf{y}; \infty)$ as inhomogeneous deformations of the LLT polynomials. In this chapter we will show that these inhomogeneous polynomials further satisfy a vanishing property, which appears similar to ones satisfied by the families of factorial Schur or interpolation Macdonald polynomials \cite{ISF, IIP}. Thus, we may interpret the $\mathcal{H}_{\boldsymbol{\lambda}} (\textbf{x}; \infty \boldsymbol{\mid} \textbf{y}; \infty)$ as\footnote{In this chapter, we largely consider the $\mathcal{H}_{\boldsymbol{\lambda}}$ functions, instead of the $\mathcal{F}_{\boldsymbol{\lambda}}$ ones from \Cref{0Polynomials0}. However, \Cref{lambdam0f} below indicates their equivalence.} ``factorial LLT polynomials.'' In the case $n = 1$, the vanishing condition for these polynomials coincides with that for the factorial Schur functions, which characterizes them completely.

\section{Zeroes of \texorpdfstring{$\mathcal{H}_{\boldsymbol{\lambda} / \boldsymbol{\mu}}$}{}}

\label{HLambdaMu0}

Let $N \ge 1$ denote an integer, and $\textbf{x} = (x_1, x_2, \ldots , x_N)$ and $\textbf{y} = (y_1, y_2, \ldots )$ denote sequences of complex numbers. In this section we will state a result, given by \Cref{lambdam0} below, that provides a family of vanishing points for the functions $\mathcal{H}_{\boldsymbol{\lambda} / \boldsymbol{\mu}} (\textbf{x}; \infty \boldsymbol{\mid} \textbf{y}; \infty)$ from \eqref{limithf}.  

These vanishing points will take the form $\{ x_j = q^{-\kappa_j} y_{\mathfrak{m}_j} \}$ for some integers $\mathfrak{m}_j \ge 1$ and $\kappa_j \in [0, n - 1]$. It will be convenient to express such specializations through \emph{marked sequences}, which are pairs $(\mathfrak{m}, \kappa)$\index{M@$(\mathfrak{m}, \kappa)$; marked sequence} of integer sets of the same length, such the entries $\mathfrak{m}_1 > \mathfrak{m}_2 > \cdots > \mathfrak{m}_{\ell}$ of $\mathfrak{m}$ are decreasing and positive, and such that each entry of $\kappa = (\kappa_1, \kappa_2, \ldots , \kappa_{\ell})$ is in $[0, n - 1]$. In what follows, we will often only refer to $\mathfrak{m}$ as the marked sequence and view $\kappa$ as its \emph{marking}, in that each entry $\mathfrak{m}_j \in \mathfrak{m}$ is \emph{marked} by the corresponding entry $\kappa_j \in \kappa$. 

To proceed, we require the notion of a splitting for a marked sequence. 

\begin{definition}
	
	\label{sequencem} 
	
	Let $\mathfrak{m} = (\mathfrak{m}_1, \mathfrak{m}_2, \ldots , \mathfrak{m}_{\ell})$ denote a marked sequence with marking given by $\kappa = (\kappa_1, \kappa_2, \ldots , \kappa_{\ell})$. A \emph{splitting} of $\mathfrak{m}$ (more precisely, of $(\mathfrak{m}, \kappa)$) is a sequence $\mathfrak{M} = \big( \mathfrak{m}^{(1)}, \mathfrak{m}^{(2)}, \ldots , \mathfrak{m}^{(n)} \big)$ of Maya diagrams (that is, decreasing subsets of $\mathbb{Z}_{> 0}$) such that the following two properties holds. 
	
	\begin{enumerate}
		\item Every entry $m \in \mathfrak{m}^{(i)}$ of any Maya diagram in $\mathfrak{M}$ is equal to $\mathfrak{m}_j$, for some $j = j (m) \in [1, \ell]$.
		\item For any $j \in [1, \ell]$, there exist at least $n - \kappa_j$ distinct indices $i \in [1, n]$ such that $\mathfrak{m}_j \in \mathfrak{m}^{(i)}$.
	\end{enumerate}
	
\end{definition}

Moreover, given two signatures $\mathfrak{m} = (\mathfrak{m}_1, \mathfrak{m}_2, \ldots , \mathfrak{m}_{\ell}) \in \Sign_{\ell}$ and $\mathfrak{n} = (\mathfrak{n}_1, \mathfrak{n}_2, \ldots , \mathfrak{n}_k) \in \Sign_k$, we (nonstandardly) write $\mathfrak{m} \nprec \mathfrak{n}$ if there exists an integer $j \ge 0$ such that $\mathfrak{m}_{\ell - j} > \mathfrak{n}_{k - j}$, where we set $\mathfrak{m}_i = \infty = \mathfrak{n}_i$ if $i \le 0$. This is equivalent to stipulating that the Young diagram for $\mathfrak{m}$ not be contained in that of $\mathfrak{n}$, if the two are superimposed to share the same bottom left corner.\footnote{Observe that this superimposition is different from the one used in \Cref{Polynomialsq}, where there two Young diagrams were always superimposed to share their top left corners.}

Now we can state the following vanishing result, which will be established in \Cref{H0Proof} below. Here, we recall the function $\mathfrak{T}$ from \eqref{t}. 

\begin{thm}
	
	\label{lambdam0}
	
	Fix integers $N \ge \ell \ge 1$; a marked sequence $(\mathfrak{m}, \kappa)$, with coordinates indexed by $[1, \ell]$; and a signature sequence $\boldsymbol{\lambda} \in \SeqSign_{n; N}$. Denoting $\mathfrak{l}^{(i)} = \mathfrak{T} \big( \lambda^{(i)} \big)$ for each $i \in [1, n]$, assume for any splitting $\mathfrak{M} = \big( \mathfrak{m}^{(1)}, \mathfrak{m}^{(2)}, \ldots , \mathfrak{m}^{(n)} \big)$ of $(\mathfrak{m}, \kappa)$ that there exists an index $h = h (\mathfrak{M}) \in [1, n]$ such that $\mathfrak{l}^{(h)} \nprec \mathfrak{m}^{(h)}$. If $\textbf{\emph{x}} = (x_1, x_2, \ldots , x_N)$ and $\textbf{\emph{y}} = (y_1, y_2, \ldots )$ are sequences of complex numbers such that $x_j = q^{-\kappa_j} y_{\mathfrak{m}_j}$ for each $j \in [1, \ell]$, then $\mathcal{H}_{\boldsymbol{\lambda}} (\textbf{\emph{x}}; \infty \boldsymbol{\mid} \textbf{\emph{y}}; \infty) = 0$. 
	
\end{thm}

\begin{rem} 
	
	\label{lambdam0f}

	By \Cref{lambdam0} and \Cref{hf}, the function $\mathcal{F}_{\boldsymbol{\lambda} / \boldsymbol{\mu}} (\textbf{x}; \infty \boldsymbol{\mid} \textbf{y}; 0)$ from \eqref{limithf} also satisfies a vanishing condition. In particular, adopting the notation and assumptions of \Cref{lambdam0}, but assuming instead that $x_j = q^{\kappa_j - n + 1} y_{\mathfrak{m}_j}$ for each $k \in [1, \ell]$, we have $\mathcal{F}_{\boldsymbol{\lambda}} (\textbf{x}; \infty \boldsymbol{\mid} \textbf{y}; 0) = 0$. Thus, these $\mathcal{F}_{\boldsymbol{\lambda} / \boldsymbol{\mu}} (\textbf{x}; \infty \boldsymbol{\mid} \textbf{y}; 0)$ may also be viewed as factorial LLT functions. Moreover, \Cref{fgsum2} can be used to show that they also satisfy a Cauchy identity, similarly to other families of interpolation polynomials \cite{CDFI,IPI}.  
	
\end{rem}

The vanishing condition prescribed by \Cref{lambdam0} appears similar to ones satisfied by various families of symmetric functions, such as factorial Schur polynomials \cite{TRF,FV} and interpolation Macdonald polynomials \cite{SIDOSS,IIP,ISF}. Although $\mathcal{H}_{\boldsymbol{\lambda}}$ is not exactly symmetric in $\textbf{x}$, it is quickly verified from \Cref{gxfxsigma} that (recalling $\mathscr{S} (\boldsymbol{\lambda})$ from \Cref{Symmetric}) its normalization
\begin{flalign}
	\label{h2}
	\check{\mathcal{H}}_{\boldsymbol{\lambda}} (\textbf{x} \boldsymbol{\mid} \textbf{y}) = \mathcal{H}_{\boldsymbol{\lambda}} (\textbf{x}^{-1}; \infty \boldsymbol{\mid} \textbf{y}^{-1}; \infty) \displaystyle\prod_{j = 1}^N x_j^{-j} \displaystyle\prod_{i = 1}^{\infty} y_i^{nN - \sum_{k = 1}^{i - 1} |\textbf{S}_k (\boldsymbol{\lambda})|},
\end{flalign}

\noindent is a polynomial in $\textbf{x}$ and $\textbf{y}$ of total degree $|\boldsymbol{\lambda}|$ that is symmetric in $\textbf{x}$. Due to the resemblance between the vanishing properties for $\check{\mathcal{H}}$ and the factorial Schur polynomials, and the facts that the former degenerate to the LLT polynomials (by \eqref{1hl}) and the latter to the Schur polynomials, one might view the functions $\check{\mathcal{H}}_{\boldsymbol{\lambda}} (\textbf{x} \boldsymbol{\mid} \textbf{y})$ as ``factorial variants'' of LLT polynomials. 

Before proceeding to the proof of \Cref{lambdam0}, let us analyze several consequences of it; throughout these examples, we adopt the notation of \Cref{lambdam0}. We begin with the case $n = 1$. 

\begin{example} 
	
	\label{h0n1} 
	
	Suppose that $n = 1$, and abbreviate $\lambda = \lambda^{(1)}$ and $\mathfrak{l} = \mathfrak{l}^{(1)} = (\mathfrak{l}_1, \mathfrak{l}_2, \ldots , \mathfrak{l}_N)$. Then, we must have $\kappa = 0^{\ell}$, so \Cref{sequencem} indicates that the unique splitting of $\mathfrak{m}$ is $\mathfrak{M} = (\mathfrak{m})$. Thus, recalling $\check{\mathcal{H}}$ from \eqref{h2}, \Cref{lambdam0} indicates that $\check{\mathcal{H}}_{\lambda} (\textbf{x} \boldsymbol{\mid} \textbf{y}) = 0$ when $x_j = y_{\mathfrak{m}_j}$ for each $j \in [1, \ell]$, if there exists some $k \in [1, n]$ such that $\mathfrak{l}_{N - k} > \mathfrak{m}_{\ell - k}$. 
	
\end{example}

\begin{rem} 
	
	\label{hfunctionn1} 
	
	For any partition $\lambda = (\lambda_1, \lambda_2, \ldots , \lambda_N)$, \Cref{h0n1} in fact enables us to identify $\check{\mathcal{H}}_{\lambda} (\textbf{x} \boldsymbol{\mid} \textbf{y})$ in the case $n = 1$ as a \emph{factorial Schur function} $s_{\lambda} (\textbf{x} \boldsymbol{\mid} \textbf{y})$. Introduced as equation (6.4) of \cite{FV} and equation (4) of \cite{TRF}, the latter is defined by
	\begin{flalign}
		\label{sxy} 
		s_{\lambda} (\textbf{x} \boldsymbol{\mid} \textbf{y}) = \displaystyle\prod_{1 \le i < j \le N} (x_i - x_j)^{-1} \det \Bigg[ \displaystyle\prod_{k = 1}^{\lambda_i + N - i} (x_j - y_k) \Bigg]_{1 \le i, j \le N}. 
	\end{flalign} 
	
	\noindent Indeed, by Theorem 2.1 of \cite{FF}, these functions vanish under the same specializations described by \Cref{h0n1}. Moreover, by Theorem 3.1 of \cite{IIP}, $s_{\lambda} (\textbf{x} \boldsymbol{\mid} \textbf{y})$ is (up to a constant factor) is the unique symmetric polynomial in $\textbf{x}$ of degree $|\lambda|$ satisfying this vanishing property. Thus, $\check{\mathcal{H}}_{\lambda} (\textbf{x} \boldsymbol{\mid} \textbf{y}) = (-1)^{|\lambda|} s_{\lambda} (\textbf{x} | \textbf{y})$, as it is quickly verified that the coefficient of $\prod_{i = 1}^n x_i^{\lambda_i}$ in both equals $(-1)^{|\lambda|}$.
	
\end{rem}

For general $n$ and $N$, it appears that a complete coordinate-wise description (as in \Cref{h0n1}) for the vanishing points described by \Cref{lambdam0} would be intricate to state. However, we can still identify some of these points, which are similar to those considered in \Cref{h0n1}. 

\begin{example} 
	
	\label{h0nkappa0} 
	
	Let us analyze when a marked sequence $(\mathfrak{m}, \kappa)$ satisfies the condition of \Cref{lambdam0} if $\kappa = 0^{\ell}$. In this case, $\mathfrak{m}$ admits the unique splitting $\mathfrak{M} = (\mathfrak{m}, \mathfrak{m}, \ldots , \mathfrak{m}) \in \SeqSign_n$, since each entry $\mathfrak{m}_j \in \mathfrak{m}$ must appear in $n - \kappa_j = n$ signatures of $\mathfrak{M}$. In particular, $(\mathfrak{m}, 0^{\ell})$ satisfies the condition from \Cref{lambdam0} if and only if there exists an index $h \in [1, n]$ for which $\mathfrak{l}^{(h)} \nprec \mathfrak{m}$. Thus, $\mathcal{H}_{\boldsymbol{\lambda}} (\textbf{x}; \infty \boldsymbol{\mid} \textbf{y}; \infty) = 0$ when $x_j = y_{\mathfrak{m}_j}$ for each $j \in [1, \ell]$, if there exists an index $k \in [1, N]$ such that $\max_{h \in [1, n]} \mathfrak{l}_{N - k}^{(h)} > \mathfrak{m}_{\ell - k}$.
	
\end{example}

The following lemma provides an explicit coordinate-wise description of the vanishing points described in \Cref{lambdam0} in the simplest ($N > 1$) case\footnote{We will also assume that $\ell = 2$, from which the $\ell = 1$ case can be recovered by setting $\mathfrak{m}_1 = \infty$.} not provided by the examples above, namely, $N = 2 = n$. 

\begin{lem} 
	
	\label{hlambda02}
	
	Adopt the notation of \Cref{lambdam0}, and suppose $N = n = \ell = 2$. In each of the following four cases, we have $\mathcal{H}_{\boldsymbol{\lambda}} (\textbf{\emph{x}}; \infty \boldsymbol{\mid} \textbf{\emph{y}}; \infty) = 0$.
	
	\begin{enumerate}
		\item We have that $(x_1, x_2) = (y_{\mathfrak{m}_1}, y_{\mathfrak{m}_2})$ and either
		\begin{flalign}
			\label{xy00} 
			\max \big\{ \mathfrak{l}_1^{(1)}, \mathfrak{l}_1^{(2)} \big\} > \mathfrak{m}_1, \quad \text{or} \quad \max \big\{ \mathfrak{l}_2^{(1)}, \mathfrak{l}_2^{(2)} \big\} > \mathfrak{m}_2.
		\end{flalign}
		
		\item We have that $(x_1, x_2) = (q^{-1} y_{\mathfrak{m}_1}, y_{\mathfrak{m}_2})$ and either
		\begin{flalign}
			\label{xy10}
			\min \big\{ \mathfrak{l}_1^{(1)}, \mathfrak{l}_1^{(2)} \big\} > \mathfrak{m}_1, \quad \text{or} \quad \max \big\{ \mathfrak{l}_2^{(1)}, \mathfrak{l}_2^{(2)} \big\} > \mathfrak{m}_2.
		\end{flalign}
		
		\item We have that $(x_1, x_2) = (y_{\mathfrak{m}_1}, q^{-1} y_{\mathfrak{m}_2})$ and both 
		\begin{flalign}
			\label{1kappa01}
			\mathfrak{l}_2^{(1)} > \mathfrak{m}_2, \quad \text{or} \quad \max \big\{ \mathfrak{l}_1^{(1)}, \mathfrak{l}_2^{(2)} \big\} > \mathfrak{m}_1.
		\end{flalign} 
		
		\noindent and 
		\begin{flalign}
			\label{2kappa01} 
			\mathfrak{l}_2^{(2)} > \mathfrak{m}_2, \quad \text{or} \quad \max \big\{ \mathfrak{l}_1^{(2)}, \mathfrak{l}_2^{(1)} \big\} > \mathfrak{m}_1.
		\end{flalign}
		
		\item We have that $(x_1, x_2) = (q^{-1} y_{\mathfrak{m}_1}, q^{-1} y_{\mathfrak{m}_2})$ and either
		\begin{flalign}
			\label{xy11} 
			\min \big\{ \mathfrak{l}_1^{(2)}, \mathfrak{l}_2^{(1)} \big\} > \mathfrak{m}_1, \quad \text{or} \quad \min \big\{ \mathfrak{l}_1^{(1)}, \mathfrak{l}_2^{(2)} \big\} > \mathfrak{m}_1, \quad \text{or} \quad \min \big\{ \mathfrak{l}_2^{(1)}, \mathfrak{l}_2^{(2)} \big\} > \mathfrak{m}_2.
		\end{flalign}
		
	\end{enumerate}
\end{lem} 

\begin{proof}
	
	The four parts of the lemma correpsond to the cases when $\kappa$ from \Cref{lambdam0} is equal to $(0, 0)$, $(1, 0)$, $(0, 1)$, and $(1, 1)$, respectively. Let us analyze when $\mathfrak{m}$ satisfies the vanishing condition indicated there in each of these four cases separately. 
	
	The first scenario $\kappa = (0, 0)$ was addressed by \Cref{h0nkappa0}, in which case $\mathfrak{m}$ satisfies the vanishing condition if and only if \eqref{xy00} holds. By \Cref{lambdam0}, this verifies the first part of the lemma.
	
	In the second scenario $\kappa = (1, 0)$, a signature sequence $\mathfrak{M} = \big( \mathfrak{m}^{(1)}, \mathfrak{m}^{(2)} \big)$ is a splitting of $\mathfrak{m}$ if and only if $\mathfrak{m}_1 \in \mathfrak{m}^{(1)} \cup \mathfrak{m}^{(2)}$ and $\mathfrak{m}_2 \in \mathfrak{m}^{(1)} \cap \mathfrak{m}^{(2)}$. Thus, we may assume that either $\mathfrak{M} = \big( (\mathfrak{m}_1, \mathfrak{m}_2), (\mathfrak{m}_2) \big)$ or $\mathfrak{M} = \big( (\mathfrak{m}_2), (\mathfrak{m}_1, \mathfrak{m}_2) \big)$. In the former case, there exists some $h \in \{ 1, 2 \}$ with $\mathfrak{l}^{(h)} \nprec \mathfrak{m}^{(h)}$ if and only if either $\max \big\{ \mathfrak{l}_2^{(1)}, \mathfrak{l}_2^{(2)} \big\} > \mathfrak{m}_2$ or $\mathfrak{l}_1^{(1)} > \mathfrak{m}_1$. In the latter case, such an $i$ exists if and only if either $\max \big\{ \mathfrak{l}_2^{(1)}, \mathfrak{l}_2^{(2)} \big\} > \mathfrak{m}_2$ or $\mathfrak{l}_1^{(2)} > \mathfrak{m}_1$. Thus, if $\kappa = (1, 0)$, all splittings $\mathfrak{M}$ of $(\mathfrak{m}, \kappa)$ satisfy the condition from \Cref{lambdam0} if and only if \eqref{xy10} holds. This verifies the second of the lemma.
	
	In the third scenario $\kappa = (0, 1)$, a signature sequence $\mathfrak{M} = \big( \mathfrak{m}^{(1)}, \mathfrak{m}^{(2)} \big)$ is a splitting of $\mathfrak{m}$ if and only if $\mathfrak{m}_1 \in \mathfrak{m}^{(1)} \cap \mathfrak{m}^{(2)}$ and $\mathfrak{m}_2 \in \mathfrak{m}^{(1)} \cup \mathfrak{m}^{(2)}$. Thus, we may assume that either $\mathfrak{M} = \big( (\mathfrak{m}_1, \mathfrak{m}_2), (\mathfrak{m}_1) \big)$ or $\mathfrak{M} = \big( (\mathfrak{m}_1), (\mathfrak{m}_1, \mathfrak{m}_2) \big)$. In the former case, there exists $h \in \{ 1, 2 \}$ with $\mathfrak{l}^{(h)} \nprec \mathfrak{m}^{(h)}$ if and only if \eqref{1kappa01} holds; in the latter case, such an $i$ exists if and only if \eqref{2kappa01} holds. Thus, if $\kappa = (0, 1)$, all splittings $\mathfrak{M}$ of $(\mathfrak{m}, \kappa)$ satisfy the condition from \Cref{lambdam0} if and only if both \eqref{1kappa01} and \eqref{2kappa01} hold. This verifies the third part of the lemma. 
	
	In the fourth scenario $\kappa = (1, 1)$, a signature sequence $\mathfrak{M} = \big( \mathfrak{m}^{(1)}, \mathfrak{m}^{(2)} \big)$ is a splitting of $\mathfrak{m}$ if and only if $\mathfrak{m}^{(1)} \cup \mathfrak{m}^{(2)} = \{ \mathfrak{m}_1, \mathfrak{m}_2 \}$. We may assume each $\mathfrak{m}_j$ appears in at most one signature of $\mathfrak{M}$, so $\mathfrak{M} = \big\{  \big( (\mathfrak{m}_1, \mathfrak{m}_2), \varnothing \big), \big( (\mathfrak{m}_1), (\mathfrak{m}_2) \big), \big( (\mathfrak{m}_2), (\mathfrak{m}_1) \big), \big( \varnothing, (\mathfrak{m}_1, \mathfrak{m}_2) \big) \big\}$. If $\mathfrak{M} = \big( (\mathfrak{m}_1, \mathfrak{m}_2), \varnothing \big)$, then there exists an $h \in \{ 1, 2 \}$ with $\mathfrak{l}^{(h)} \nprec \mathfrak{m}^{(h)}$ if and only if either $\mathfrak{l}_2^{(1)} > \mathfrak{m}_2$ or $\mathfrak{l}_1^{(1)} > \mathfrak{m}_1$. Similarly, if $\mathfrak{M} = \big( \varnothing, (\mathfrak{m}_1, \mathfrak{m}_2) \big)$, then such an $i$ exists if and only if either $\mathfrak{l}_2^{(2)} > \mathfrak{m}_2$ or $\mathfrak{l}_1^{(2)} > \mathfrak{m}_1$. If instead $\mathfrak{M} = \big( (\mathfrak{m}_1), (\mathfrak{m}_2) \big)$, then we must have $\mathfrak{l}_2^{(1)} > \mathfrak{m}_1$ or $\mathfrak{l}_2^{(2)} > \mathfrak{m}_2$, and if $\mathfrak{M} = \big( (\mathfrak{m}_2), (\mathfrak{m}_1) \big)$ then must have $\mathfrak{l}_2^{(1)} > \mathfrak{m}_2$ or $\mathfrak{l}_2^{(2)} > \mathfrak{m}_1$. Hence, if $\kappa = (1, 1)$, all splittings $\mathfrak{M}$ of $(\mathfrak{m}, \kappa)$ satisfy the condition from \Cref{lambdam0} if and only if \eqref{xy11} holds. This verifies the fourth part of the lemma.
\end{proof}

\section{Blocking Vertices}

\label{SetsVertex}

In this section we introduce the notion of a blocking  vertex with respect to a path ensemble, which will be useful for the proof of \Cref{lambdam0}. We will only require them here in the case $n = 1$ (that is, when there is only one fermionic color), so throughout this section we will assume this holds. 

Then, any sequence of $n = 1$ skew-shapes $\boldsymbol{\lambda} / \boldsymbol{\mu}$ consists of a single skew-shape $\lambda / \mu$. So, in what follows, we identify $\boldsymbol{\lambda} / \boldsymbol{\mu} = \lambda / \mu$; for instance, recalling \Cref{pgpfph} we write $\mathfrak{P}_H (\lambda / \mu) = \mathfrak{P}_H (\boldsymbol{\lambda} / \boldsymbol{\mu})$ (and $\mathfrak{P}_H (\lambda) = \mathfrak{P}_H (\boldsymbol{\lambda} / \boldsymbol{\varnothing})$). Moreover, since $n = 1$, any arrow configuration $(\textbf{A}, \textbf{B}; \textbf{C}, \textbf{D})$ is of the form $(\textbf{e}_a, \textbf{e}_b; \textbf{e}_c, \textbf{e}_d)$ for some $a, b, c, d \in \{ 0, 1 \}$. Hence, we will also abbreviate $(a, b; c, d) = (\textbf{e}_a, \textbf{e}_b; \textbf{e}_c, \textbf{e}_d)$, as well as the weight $\mathcal{W}_z (a, b; c, d) = \mathcal{W}_z (\textbf{e}_a, \textbf{e}_b; \textbf{e}_c, \textbf{e}_d \boldsymbol{\mid} \infty, \infty)$ from \eqref{limitw}, for any $a, b, c, d \in \{ 0, 1 \}$. Recalling \Cref{weightesum}, we also abbreviate the weight of any path ensemble $\mathcal{E}$ under $\mathcal{W}_z (a, b; c, d)$ by $W (\mathcal{E} \boldsymbol{\mid} \textbf{x} \boldsymbol{\mid} \textbf{y}) = W (\mathcal{E} \boldsymbol{\mid} \textbf{x}; \infty \boldsymbol{\mid} \textbf{y}; \infty)$.

Now recall from \eqref{gfhe} and \Cref{weightesum} that, for any $\lambda, \mu \in \Sign_n$, the quantity $\mathcal{H}_{\lambda / \mu} (\textbf{x}; \infty \boldsymbol{\mid} \textbf{y}; \infty)$ is the partition function under the $\mathcal{W}_{x / y} (a, b; c, d)$ weights for the vertex model $\mathfrak{P}_H (\lambda / \mu)$. Due to the factor of $(xy^{-1}; q)_{b - c}$ present in these weights, if $x = y$ then $\mathcal{W}_{x / y} (a, b; c, d) = 0$ whenever $(b, c) = (1, 0)$. In particular, if a path ensemble $\mathcal{E} \in \mathfrak{P}_H (\lambda / \mu)$ has nonzero weight $W (\mathcal{E} \boldsymbol{\mid} \textbf{x} \boldsymbol{\mid} \textbf{y}) \ne 0$, and there exists some vertex $(i, j)$ with $x_j = y_i$, then the arrow configuration $(a, b; c, d)$ at $(i, j)$ must satisfy $(b, c) \ne (1, 0)$. Thus, we introduce the following definition.

\begin{definition}
	
	\label{vertexe} 
	
	Let $\mathcal{E}$ denote a path ensemble on some domain $\mathcal{D} \subseteq \mathbb{Z}_{> 0}^2$; for each vertex $u \in \mathcal{D}$, let $\big( a(u), b(u); c(u), d(u) \big)$ denote the arrow configuration at $u$ under $\mathcal{E}$. We call any $v \in \mathcal{D}$ a \emph{blocking vertex} with respect to $\mathcal{E}$ if $\big( b(v), c(v) \big) \ne (1, 0)$. 
	
\end{definition}

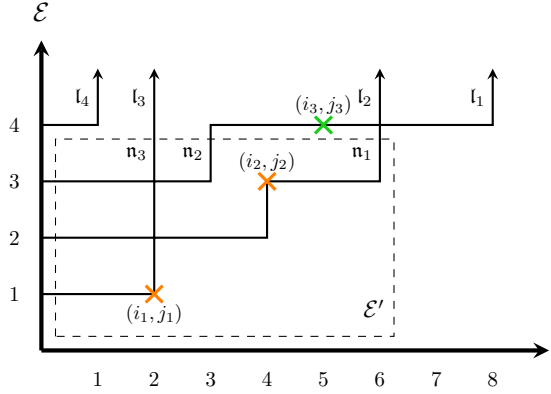
\begin{figure}

	\begin{center}

		\begin{tikzpicture}[
			>=stealth,
			scale = .75
			]
			\draw[ultra thick, ->] (0, 0) -- (9, 0);
			\draw[ultra thick, ->] (0, 0) -- (0, 5.5);
			
			\draw[thick, ->] (0, 1) -- (2, 1) -- (2, 2) -- (3, 2) -- (4, 2) -- (4, 3) -- (6, 3) -- (6, 4) -- (8, 4) -- (8, 5);
			\draw[thick, ->] (0, 2) -- (2, 2) -- (2, 3) -- (3, 3) -- (3, 4) -- (6, 4) -- (6, 5);
			\draw[thick, ->] (0, 3) -- (2, 3) -- (2, 5);
			\draw[thick, ->] (0, 4) -- (1, 4) -- (1, 5);
			
			\draw[] (8, 4.5) circle[radius = 0] node[left, scale = .8]{$\mathfrak{l}_1$};
			\draw[] (6, 4.5) circle[radius = 0] node[left, scale = .8]{$\mathfrak{l}_2$};
			\draw[] (2, 4.5) circle[radius = 0] node[left, scale = .8]{$\mathfrak{l}_3$};
			\draw[] (1, 4.5) circle[radius = 0] node[left, scale = .8]{$\mathfrak{l}_4$};
			
			\draw[] (6, 3.5) circle[radius = 0] node[left, scale = .8]{$\mathfrak{n}_1$};
			\draw[] (3, 3.5) circle[radius = 0] node[left, scale = .8]{$\mathfrak{n}_2$};
			\draw[] (2, 3.5) circle[radius = 0] node[left, scale = .8]{$\mathfrak{n}_3$};
			
			\draw[] (-.25, 1) circle[radius = 0] node[left, scale = .8]{$1$};
			\draw[] (-.25, 2) circle[radius = 0] node[left, scale = .8]{$2$};
			\draw[] (-.25, 3) circle[radius = 0] node[left, scale = .8]{$3$};
			\draw[] (-.25, 4) circle[radius = 0] node[left, scale = .8]{$4$};
			
			\draw[] (1, -.25) circle[radius = 0] node[below, scale = .8]{$1$};
			\draw[] (2, -.25) circle[radius = 0] node[below, scale = .8]{$2$};
			\draw[] (3, -.25) circle[radius = 0] node[below, scale = .8]{$3$};
			\draw[] (4, -.25) circle[radius = 0] node[below, scale = .8]{$4$};
			\draw[] (5, -.25) circle[radius = 0] node[below, scale = .8]{$5$};
			\draw[] (6, -.25) circle[radius = 0] node[below, scale = .8]{$6$};
			\draw[] (7, -.25) circle[radius = 0] node[below, scale = .8]{$7$};
			\draw[] (8, -.25) circle[radius = 0] node[below, scale = .8]{$8$};
			
			\draw[] (0, 6) circle[radius = 0] node[scale = 1]{$\mathcal{E}$};
			
			\draw[dashed] (.25, .25) -- (.25, 3.75) -- (6.25, 3.75) -- (6.25, .25) -- (.25, .25);
			
			\draw[] (5.9, .75) circle[radius = 0] node[scale = .9]{$\mathcal{E}'$};

			\draw[orange, very thick] (1.85, .85) -- (2.15, 1.15);
			\draw[orange, very thick] (2.15, .85) -- (1.85, 1.15);
			
			\draw[orange, very thick] (3.85, 2.85) -- (4.15, 3.15);
			\draw[orange, very thick] (4.15, 2.85) -- (3.85, 3.15);
			
			\draw[green, very thick] (4.85, 3.85) -- (5.15, 4.15);
			\draw[green, very thick] (5.15, 3.85) -- (4.85, 4.15);
			
			\draw[] (2, .65) circle[radius = 0] node[scale = .75]{$(i_1, j_1)$};
			\draw[] (4, 3.35) circle[radius = 0] node[scale = .75]{$(i_2, j_2)$};
			\draw[] (5, 4.35) circle[radius = 0] node[scale = .75]{$(i_3, j_3)$};

		\end{tikzpicture}
		
	\end{center}
	
	\caption{\label{lambdahvertexset} Depicted above is an attempt for an increasing blocking vertex set, given by the orange and green crosses. However, only the orange crosses are blocking with respect to $\mathcal{E}$; the green one is not. }
\end{figure}

Observe in particular that any path in $\mathcal{E}$ horizontally entering some blocking vertex $v \in \mathcal{D}$ must exit $v$ vertically. In this sense, $v$ ``blocks'' the horizontal trajectory of this path. We refer to \Cref{lambdahvertexset} for a depiction, where the orange crosses at $(2, 1)$ and $(4, 3)$ are blocking with respect to $\mathcal{E}$ but the green one at $(5, 4)$ is not. 

In this section we establish the following proposition, which indicates how blocking vertices may arrange themselves in a path ensemble $\mathcal{E} \in \mathfrak{P}_H (\lambda / \mu)$.

\begin{prop}
	
	\label{lambdabnc}
	
	Fix a signature $\lambda = (\lambda_1, \lambda_2, \ldots , \lambda_N) \in \Sign_N$, and denote $\mathfrak{l} = \mathfrak{T} (\lambda) \in \Sign_N$. For any fixed path ensemble $\mathcal{E} \in \mathfrak{P}_H (\lambda)$, there does not exist a set of $K \le N$ blocking vertices $(v_1, v_2, \ldots , v_K) \subset \mathcal{D}_N = \mathbb{Z}_{> 0} \times \{ 1, 2, \ldots , N \}$ with respect to $\mathcal{E}$ satisfying the following two properties. In the below we set $v_k = (i_k, j_k)$ for each $k \in [1, K]$. 
	
	\begin{enumerate}
		\item We have $1 \le j_1 < j_2 < \cdots < j_K \le N$ and $1 \le i_1 < i_2 < \cdots < i_K$.
		\item We have $i_K < \mathfrak{l}_{N - K + 1}$.
	\end{enumerate}
\end{prop}

Stated alternatively, \Cref{lambdabnc} indicates that, for any vertex sequence $(v_1, v_2, \ldots , v_K)$ satisfying its conditions and path ensemble $\mathcal{E} \in \mathfrak{P}_H (\lambda)$, there must exist some $j \in [1, K]$ for which $\big( b(v_j), c(v_j) \big) = (1, 0)$.

\begin{proof}[Proof of \Cref{lambdabnc}]
	
	In what follows, for any path ensemble $\mathcal{E}$ on some domain $\mathcal{D}$, we denote the arrow configuration at any vertex $v \in \mathcal{D}$ under $\mathcal{E}$ by $\big( a_{\mathcal{E}} (v), b_{\mathcal{E}} (v); c_{\mathcal{E}} (v), d_{\mathcal{E}} (v) \big)$. Throughout this proof, we call any vertex set $(v_1, v_2, \ldots , v_K) \in \mathcal{D}_N$ satisfying the two properties listed in the propositon an \emph{increasing blocking vertex set} with respect to $\mathcal{E}$. We must show that no path ensemble $\mathcal{E} \in \mathfrak{P}_H (\lambda)$ admits an increasing blocking vertex set.
	
	To that end, we induct on $N + K \ge 2$. If $N + K = 2$, then $N = 1 = K$, so there is a unique path ensemble $\mathcal{E} \in \mathfrak{P}_H (\lambda)$ consisting of one path that horizontally enters $\mathcal{D}_1$ through the vertex $(1, 1)$ and proceeds east until it reaches $(\mathfrak{l}_1, 1)$, where it vertically exits $\mathcal{D}_1$. Thus, $\big( b_{\mathcal{E}} (i, 1), c_{\mathcal{E}} (i, 1) \big) = (1, 0)$ holds for each $i \in [0, \mathfrak{l}_1 - 1]$, and hence there does not exist a blocking vertex satisfying the second condition in the proposition. In particular, no increasing blocking vertex set exists.
	
	So, let us suppose that the proposition holds whenever $N + K < m$ for some integer $m \ge 3$, and we will show it also holds if $N + K = m$. Suppose to the contrary that there exists a path ensemble $\mathcal{E} \in \mathfrak{P}_H (\lambda)$ that admits an increasing blocking vertex set $\mathcal{V} = (v_1, v_2, \ldots , v_K) \subset \mathcal{D}_N$; for each $k \in [1, K]$, set $v_k = (i_k, j_k) \in \mathcal{D}_N$. 
	
	Let $\mathcal{E}'$ denote the restriction of $\mathcal{E}$ to the subdomain $\mathcal{D}_{N - 1} = \mathbb{Z}_{> 0} \times \{ 1, 2, \ldots , N - 1 \} \subset \mathcal{D}_N$ (that is, to the bottommost $N - 1$ rows of $\mathcal{D}_N$). Then, there exists a signature $\nu \in \Sign_{N - 1}$ such that $\mathcal{E}' \in \mathfrak{P}_H (\nu)$. Let $\mathfrak{T} (\nu) = \mathfrak{n} = (\mathfrak{n}_1, \mathfrak{n}_2, \ldots , \mathfrak{n}_{N - 1}) \in \Sign_{N- 1}$, which denote the $x$-coordinates of the locations where paths in $\mathcal{E}'$ vertically exit the row $\mathbb{Z}_{> 0} \times \{ N - 1 \}$; see \Cref{lambdahvertexset} for a depiction. Since these paths (in $\mathcal{E}$) exit the row $\mathbb{Z}_{> 0} \times \{ N \}$ through the $x$-coordinates $(\mathfrak{l}_1, \mathfrak{l}_2, \ldots , \mathfrak{l}_{N - 1})$, we have the interlacing property 
	\begin{flalign}
		\label{nui1ji1j1ij} 
		\mathfrak{l}_{j + 1} \le \mathfrak{n}_j \le \mathfrak{l}_j, \qquad \text{for any index $j \in [1, N - 1]$},
	\end{flalign}
	
	\noindent which follows from the fact that no two paths may share an edge. 
	
	Now, let us show that $\mathcal{E}'$ admits an increasing blocking vertex set in each of the cases $j_K \ne N$ and $j_K = N$ separately. If $j_K \ne N$, then we claim that the vertex set $\mathcal{V} = (v_1, v_2, \ldots , v_K) \subset \mathcal{D}_{N - 1}$ is increasing blocking with respect to $\mathcal{E}'$. Indeed, the fact that each of its vertices is blocking with respect to $\mathcal{E}'$ follows from the fact that they are with respect to $\mathcal{E}$. Moreover, the facts that $1 \le j_1 < j_2 < \cdots < j_K \le N - 1$ and $1 \le i_1 < i_2 < \cdots < i_K$ follow from the facts that $\mathcal{V}$ is increasing blocking with respect to $\mathcal{E}$ and that $j_K \ne N$. Additionally, the bound $i_K < \mathfrak{n}_{N - K}$ follows from the fact that $i_K < \mathfrak{l}_{N - K + 1} \le \mathfrak{n}_{N - K}$, where the first inequality holds again since $\mathcal{V}$ is increasing blocking with respect to $\mathcal{E}$ and the second holds by \eqref{nui1ji1j1ij}. This confirms that the vertex set $\mathcal{V}$ is increasing blocking with respect to $\mathcal{E}'$ if $j_K \ne N$, which contradicts the inductive hypothesis that no path ensemble in $\mathfrak{P}_H (\nu)$ can admit such a set (since $\ell (\nu) + K = N + K - 1 = m - 1 < m$). 
	
	So, let us assume instead that $j_K = N$, in which case we claim that the vertex set $\mathcal{V}' = (v_1, v_2, \ldots , v_{K - 1}) \subset \mathcal{D}_{N - 1}$ is increasing blocking with respect to $\mathcal{E}'$. Again, the facts that each vertex $v_i \in \mathcal{V}'$ is blocking and that $1 \le j_1 < j_2 < \cdots < j_{K - 1} \le N - 1$ and $1 \le i_1 < i_2 < \cdots < i_{K - 1}$ follow from the fact that $\mathcal{V}$ is increasing blocking with respect to $\mathcal{E}$. Thus, it remains to verify the bound $i_{K - 1} < \mathfrak{n}_{N - K + 1}$. To that end, since $i_{K - 1} < i_K$, it suffices to show that $i_K \le \mathfrak{n}_{N - K + 1}$.  
	
	Assume to the contrary that this is false, so $\mathfrak{n}_{N - K + 1} < i_K < \mathfrak{l}_{N - K + 1}$, where the last inequality holds since $\mathcal{V}$ is increasing blocking with respect to $\mathcal{E}$. Since there exists a path in $\mathcal{E}$ that enters the row $\mathbb{Z}_{> 0} \times \{ N \}$ vertically through the vertex $(\mathfrak{n}_{N - K + 1}, N)$ and proceeds east until it reaches $(\mathfrak{l}_{N - K + 1}, N)$, we must have $\big( b_{\mathcal{E}} (i), c_{\mathcal{E}} (i) \big) = (1, 0)$ for any $i \in (\mathfrak{n}_{N - K + 1}, \mathfrak{l}_{N - K + 1})$; see \Cref{lambdahvertexset}, where $(N, K) = (4, 3)$ and $\big\{ (i_1, j_1), (i_2, j_2), (i_3, j_3) \big\} = \big\{ (2, 1), (4, 3), (5, 4) \big\}$. In particular, this holds for $i = i_K$, contradicting the fact that $v_K = (i_K, N)$ is a blocking vertex with respect to $\mathcal{E}$. 
	
	Hence, $i_K \le \mathfrak{n}_{N - K + 1}$, so $\mathcal{V}'$ is increasing blocking with respect to $\mathcal{E}'$. This again contradicts the inductive hypothesis that no path ensemble in $\mathfrak{P}_H (\nu)$ can admit such a set, thereby establishing the proposition.
\end{proof}

\section{Proof of \Cref{lambdam0}} 

\label{H0Proof}

In this section we establish \Cref{lambdam0}, which will largely follow from \Cref{lambdabnc}.

\begin{proof}[Proof of \Cref{lambdam0}]
	
	Throughout this proof, we abbreviate the vertex model $\mathfrak{P}_H (\boldsymbol{\lambda}) = \mathfrak{P}_H (\boldsymbol{\lambda} / \boldsymbol{\varnothing})$ from \Cref{pgpfph}; the vertex weights $\mathcal{W}_z (\textbf{A}, \textbf{B}; \textbf{C}, \textbf{D}) = \mathcal{W}_z (\textbf{A}, \textbf{B}; \textbf{C}, \textbf{D} \boldsymbol{\mid} \infty, \infty)$ from \eqref{limitw}; and the weight $W (\mathcal{E} \boldsymbol{\mid} \textbf{x} \boldsymbol{\mid} \textbf{y}) = W (\mathcal{E} \boldsymbol{\mid} \textbf{x}; \infty \boldsymbol{\mid} \textbf{y}; \infty)$ of any path ensemble $\mathcal{E}$ from \Cref{weightesum}. We moreover recall from \Cref{SetsVertex} the $n = 1$ notation $\mathfrak{P}_H (\lambda) = \mathfrak{P}_H (\boldsymbol{\lambda})$ if $\boldsymbol{\lambda} = (\lambda) \in \SeqSign_1$. Let us also assume throughout this proof that $x_{N - j + 1} = y_{\mathfrak{m}_j}$ (instead of $x_j = y_{\mathfrak{m}_j}$) for each $j \in [1, \ell]$, which we may do since \Cref{gxfxsigma} implies that $\mathcal{H}_{\boldsymbol{\lambda}}$ is symmetric in $\textbf{x}$ up to a factor. 
	
	Now, assume to the contrary that $\mathcal{H}_{\boldsymbol{\lambda}} (\textbf{x}; \infty \boldsymbol{\mid} \textbf{y}; \infty) \ne 0$. Since \eqref{gfhe} and \eqref{weightesum} together imply that $\mathcal{H}_{\boldsymbol{\lambda}} (\textbf{x}; \infty \boldsymbol{\mid} \textbf{y}; \infty)$ is the partition function, under the weights $\mathcal{W}_z$, for the vertex model $\mathfrak{P}_H (\boldsymbol{\lambda})$, there must exist some path ensemble $\mathcal{E} \in \mathfrak{P}_H (\boldsymbol{\lambda})$ on $\mathcal{D}_N = \mathbb{Z}_{> 0} \times \{ 1, 2, \ldots , N \}$ with nonzero weight $W (\mathcal{E} \boldsymbol{\mid} \textbf{x} \boldsymbol{\mid} \textbf{y}) \ne 0$. Let $\big( \textbf{A} (v), \textbf{B} (v); \textbf{C} (v), \textbf{D} (v) \big)$ denote the arrow configuration at any vertex $v \in \mathcal{D}_N$ under $\mathcal{E}$, and set $\textbf{X} (v) = \big( X_1 (v), X_2 (v), \ldots , X_n (v) \big) \in \{ 0, 1 \}^n$ for each index $X \in \{ A, B, C, D \}$.  
	
	For each $i \in [1, n]$, let $\mathcal{E}^{(i)} \in \mathfrak{P}_H \big( \lambda^{(i)} \big)$ denote the restriction to color $i$ of $\mathcal{E}$. Stated alternatively, the arrow configuration at any vertex $v \in \mathcal{D}_N$ under $\mathcal{E}^{(i)}$ is given by $\big( A_i (v), B_i (v); C_i (v), D_i (v) \big)$. In this way, each $\mathcal{E}^{(i)}$ is a path ensemble with one fermionic color, and so the framework of \Cref{SetsVertex} applies to it. We refer to \Cref{e1e2} for a depiction, where there red is color $1$ and blue is color $2$.
	
	Under this notation, for each $j \in [1, \ell]$, let $\mathcal{I}_j \subseteq [1, n]$ denote the set of colors $i \in [1, n]$ for which $(\mathfrak{m}_j, N - j + 1) \in \mathcal{D}_N$ is a blocking vertex with respect to $\mathcal{E}^{(i)}$. We claim that $|\mathcal{I}_j| \ge n - \kappa_j$. Indeed, assume to the contrary that there exists a subset $\mathcal{U} \subseteq [1, n]$ of $\kappa_j + 1$ distinct colors such that $(\mathfrak{m}_j, N - j + 1) \in \mathcal{D}_N$ is not a blocking vertex with respect to $\mathcal{E}^{(u)}$ for each $u \in \mathcal{U}$; so, $B_u (\mathfrak{m}_j, N - j + 1) - C_u (\mathfrak{m}_j, N - j + 1) = 1$ for each such $u$. Thus, abbreviating the arrow configuration under $\mathcal{E}$ at $(\mathfrak{m}_j, N - j + 1) \in \mathcal{D}_N$ by
	\begin{flalign*} 
		(\textbf{A}, \textbf{B}; \textbf{C}, \textbf{D}) = \big( \textbf{A} (\mathfrak{m}_j, N - j + 1), \textbf{B} (\mathfrak{m}_j, N - j + 1); \textbf{C} (\mathfrak{m}_j, N - j + 1), \textbf{D} (\mathfrak{m}_j, N - j + 1) \big),
	\end{flalign*} 
	
	\noindent we deduce either that $|\textbf{B}| - | \textbf{C}| \ge \kappa_j + 1$ or that $\textbf{B} \ge \textbf{C}$ does not hold. Recalling that $x_{N - j + 1} = q^{-\kappa_j} y_{\mathfrak{m}_j}$, it follows from the factor of 
	\begin{flalign*} 
		(x_{N - j + 1} y_{\mathfrak{m}_j}^{-1}; q)_{|\textbf{B}|  - |\textbf{C}|} \textbf{1}_{\textbf{B} \ge \textbf{C}} = (q^{-\kappa_j}; q)_{|\textbf{B}| - |\textbf{C}|} \textbf{1}_{\textbf{B} \ge \textbf{C}}, 
	\end{flalign*} 
	
	\noindent appearing in the vertex weight $\mathcal{W}_{x_{N - j + 1}; y_{\mathfrak{m}_j}} (\textbf{A}, \textbf{B}; \textbf{C}, \textbf{D} \boldsymbol{\mid} \infty; \infty)$ (from the first statement of \eqref{limitw}) that this weight equals $0$. This contradicts the fact that $\mathcal{E}$ has nonzero weight $W(\mathcal{E} \boldsymbol{\mid} \textbf{x} \boldsymbol{\mid} \textbf{y}) \ne 0$ under $\mathcal{W}_z$, thus verifying $|\mathcal{I}_j| \ge n - \kappa_j$.

	\begin{figure}
		
		\begin{center}

			\begin{tikzpicture}[
				>=stealth,
				scale = .5
				]
				
				\draw[->, thick, red] (0, .95) -- (1.95, .95) -- (1.95, 1.95) -- (3.95, 1.95) -- (3.95, 2.95) -- (4.95, 2.95) -- (4.95, 4);
				\draw[->, thick, red] (0, 1.95) -- (1.95, 1.95) -- (1.95, 2.95) -- (2.95, 2.95) -- (2.95, 4);
				\draw[->, thick, red] (0, 2.95) -- (.95, 2.95) -- (.95, 4);
				
				\draw[->, thick, blue] (0, 1.05) -- (3.05, 1.05) -- (3.05, 2.05) -- (5.05, 2.05) -- (5.05, 3.05) -- (6.05, 3.05) -- (6.05, 4);
				\draw[->, thick, blue] (0, 2.05) -- (2.05, 2.05) -- (2.05, 3.05) -- (5.05, 3.05) -- (5.05, 4);
				\draw[->, thick, blue] (0, 3.05) -- (2.05, 3.05) -- (2.05, 4);

				\draw[-] (1, 4.5) -- (1, 5) -- (7, 5) -- (7, 4.5);
				
				\draw[] (4, 5) circle[radius = 0] node[above, scale = .7]{$\boldsymbol{\lambda}$};
				\draw[] (0, 4.75) circle[radius = 0] node[above, scale = .8]{$\mathcal{E}$};

				\draw[->, thick, red] (10, 1) -- (12, 1) -- (12, 2) -- (14, 2) -- (14, 3) -- (15, 3) -- (15, 4);
				\draw[->, thick, red] (10, 2) -- (12, 2) -- (12, 3) -- (13, 3) -- (13, 4);
				\draw[->, thick, red] (10, 3) -- (11, 3) -- (11, 4);

				\draw[-] (11, 4.5) -- (11, 5) -- (16, 5) -- (16, 4.5);
				
				\draw[] (13.5, 5) circle[radius = 0] node[above, scale = .7]{$\lambda^{(1)}$};
				\draw[] (10, 4.75) circle[radius = 0] node[above, scale = .8]{$\mathcal{E}^{(1)}$};

				\draw[->, thick, blue] (20, 1) -- (23, 1) -- (23, 2) -- (25, 2) -- (25, 3) -- (26, 3) -- (26, 4);
				\draw[->, thick, blue] (20, 2) -- (22, 2) -- (22, 3) --(25, 3) -- (25, 4);
				\draw[->, thick, blue] (20, 3) -- (22, 3) -- (22, 4);
				
				\draw[-] (21, 4.5) -- (21, 5) -- (26, 5) -- (26, 4.5);
				
				\draw[] (23.5, 5) circle[radius = 0] node[above, scale = .7]{$\lambda^{(2)}$};
				\draw[] (20, 4.75) circle[radius = 0] node[above, scale = .8]{$\mathcal{E}^{(2)}$};

				\draw[->, very thick] (0, 0) -- (0, 4.5);
				\draw[->, very thick] (0, 0) -- (7.5, 0);
				
				\draw[->, very thick] (10, 0) -- (10, 4.5);
				\draw[->, very thick] (10, 0) -- (17.5, 0);
				
				\draw[->, very thick] (20, 0) -- (20, 4.5);
				\draw[->, very thick] (20, 0) -- (27.5, 0);
				
			\end{tikzpicture}
			
		\end{center}
		
		\caption{\label{e1e2} Shown to the left, middle, and right are examples of the path ensembles $\mathcal{E}$, $\mathcal{E}^{(1)}$, and $\mathcal{E}^{(2)}$, respectively, from the proof of \Cref{lambdam0}. } 
		
	\end{figure}
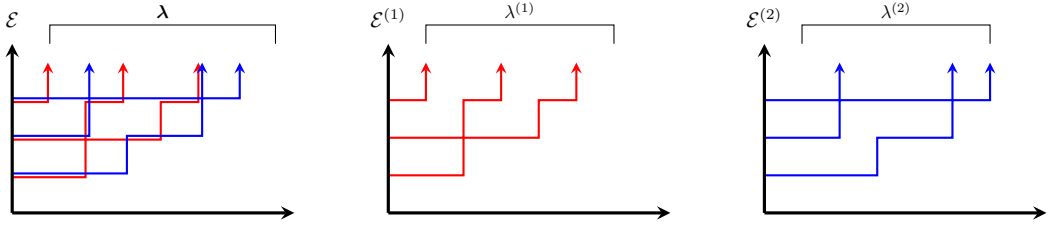
	
	Next, for each $i \in [1, n]$, let $\mathcal{V}^{(i)}$ denote the set of vertices of the form $(\mathfrak{m}_j, N - j + 1) \in \mathcal{D}_N$, for some $j \in [1, \ell]$, that are blocking with respect to $\mathcal{E}^{(i)}$. We order the vertices in $\mathcal{V}^{(i)}$ by
	\begin{flalign*}
		\mathcal{V}^{(i)} = \big( v_1^{(i)}, v_2^{(i)}, \ldots , v_{\ell_i}^{(i)} \big) = \Big\{ \big( \mathfrak{m}_{j_1^{(i)}}, N - j_1^{(i)} + 1 \big), \big( \mathfrak{m}_{j_2^{(i)}}, N - j_2^{(i)} + 1 \big), \ldots , \big( \mathfrak{m}_{j_{\ell_i}^{(i)}}, N - j_{\ell_i}^{(i)} + 1 \big) \Big\},
	\end{flalign*} 
	
	\noindent so that $j_1^{(i)} > j_2^{(i)} > \cdots > j_{\ell_i}^{(i)}$, which since $\mathfrak{m}$ is decreasing implies 
	\begin{flalign*}
		\mathfrak{m}_{j_1^{(i)}} < \mathfrak{m}_{j_2^{(i)}} < \cdots < \mathfrak{m}_{j_{\ell_i}^{(i)}}.
	\end{flalign*}
	
	\noindent Hence, the blocking vertex set $\mathcal{V}^{(i)} = \big( v_1^{(i)}, v_2^{(i)}, \ldots , v_{\ell_i}^{(i)} \big)$ (with respect to $\mathcal{E}^{(i)}$) satisfies the first property listed in \Cref{lambdabnc}, as does its truncation $\mathcal{V}_K^{(i)} = \big( v_1^{(i)}, v_2^{(i)}, \ldots , v_K^{(i)} \big)$, for any $K \in [1, \ell_i]$. Thus, that proposition implies none of the $\mathcal{V}_K^{(i)}$ can satisfy the second property listed there, meaning that $\mathfrak{m}_{j_K^{(i)}} \le \mathfrak{l}_{N - K + 1}^{(i)}$ holds for each $K \in [1, \ell_i]$. So, letting
	\begin{flalign*}
		\mathfrak{m}^{(i)} = \Big( \mathfrak{m}_{j_1^{(i)}}, \mathfrak{m}_{j_2^{(i)}}, \ldots , \mathfrak{m}_{j_{\ell^{(i)}}} \Big),
	\end{flalign*}
	
	\noindent we deduce that 
	\begin{flalign} 
		\label{msequence2} 
		\mathfrak{m}^{(i)} \nprec \mathfrak{l}^{(i)} \qquad \text{does not hold for any $i \in [1, n]$}.
	\end{flalign} 
	
	Now, observe that $\mathfrak{M} = \big( \mathfrak{m}^{(1)}, \mathfrak{m}^{(2)}, \ldots , \mathfrak{m}^{(n)} \big)$ is a splitting of $(\mathfrak{m}, \kappa)$. Indeed, it satisfies the first condition listed in \Cref{sequencem}, since each entry of every $\mathfrak{m}^{(i)}$ is also an entry of $\mathfrak{m}$. It also satisfies the second, since $|\mathcal{I}_j| \ge n - \kappa_j$ for each $j \in [1, \ell]$ implies that any $\mathfrak{m}_j$ is an element in at least $n - \kappa_j$ of the $\mathfrak{m}^{(i)}$. 
	
	Thus, since the marked sequence $(\mathfrak{m}, \kappa)$ satisfies the condition of the theorem, there must exist some $h \in [1, n]$ such that $\mathfrak{l}^{(h)} \nprec \mathfrak{m}^{(h)}$. This contradicts \eqref{msequence2}, and so $\mathcal{H}_{\boldsymbol{\lambda}} (\textbf{x}; \infty \boldsymbol{\mid} \textbf{y}; \infty) = 0$. 
\end{proof}

\chapter{Vertex Models for Nonsymmetric Macdonald Polynomials}

\label{Polynomials}

In this chapter we provide an expression, given by \Cref{thm:f-formula} below, for nonsymmetric Macdonald polynomais in terms of fermionic, colored vertex model partition functions. This result is similar to Theorem 4.2 of \cite{NPIVM} (which established an expression for these polynomials through a bosonic partition function) and will be established in \Cref{ProofPolynomial} below, largely following the framework of \cite{NPIVM}. Throughout this section, we fix complex numbers $q, t \in \mathbb{C}$.

\section{Cherednik--Dunkl Operators and Nonsymmetric Polynomials}

\label{Polynomialf}

In what follows, we use the same conventions as in \cite{NPIVM} (which are in turn based on reversing all indices in the notation of \cite{RKP}) in defining the nonsymmetric Macdonald polynomials $f_{\mu} (\textbf{x}) = f_{\mu}(x_1,\dots,x_n)$, for any composition $\mu = (\mu_1, \mu_2, \ldots , \mu_n)$ and sequence $\textbf{x} = (x_1, x_2, \ldots , x_n)$ of complex variables. % All operators in this section carry a tilde; this is to distinguish them from a modified form of the operators that we write down in \Cref{ssec:reverse}. 
We begin by recalling the definition of the {\it Hecke algebra} of type $A_{n-1}$. It is the algebra generated by the family of generators $T_1, T_2, \dots T_{n-1}$, which satisfy the relations
\begin{align}
\label{hecke1}
(T_i - t)(T_i + 1) = 0,
\quad 
\text{for $i \in [1, n - 1]$};
\qquad
T_i T_{i+1} T_i = T_{i+1} T_i 
T_{i+1},
\quad
\text{for $i \in [1, n - 2]$},
\end{align}
as well as the commutativity property
\begin{align}
\label{hecke2}
[T_i,T_j] = 0,
\quad \text{for all $i,j$ such that $|i-j| > 1$}.
\end{align}
\index{T@$T_i$; generators of Hecke algebra}
A well-known realization of this algebra is its {\it polynomial representation}. In this representation, one identifies the abstract generator $T_i$ and its inverse with explicit \emph{Demazure--Lusztig} operators on $\mathbb{C}[\textbf{x}] = \mathbb{C}[x_1, x_2, \dots,x_n]$, the ring of polynomials in $n$ variables, by setting for each $1 \le i \le n - 1$
\begin{align}
\label{Hecke}
T_i \mapsto t - \frac{x_i - t x_{i+1}}{x_i-x_{i+1}} (1-\s_i),
\qquad
T_i^{-1} \mapsto t^{-1}\left(1 - \frac{x_i - t x_{i+1}}{x_i-x_{i+1}} (1-\s_i) \right),
\end{align}
where, as previously, $\mathfrak{s}_i$ denotes the transposition operator on neighboring variables, namely
\begin{align*}
\mathfrak{s}_i \cdot h(x_1, x_2, \dots,x_n)
=
h(x_1,\dots, x_{i - 1}, x_{i+1},x_i, x_{i + 2}, \dots,x_n),
\end{align*}
for any polynomial $h \in \mathbb{C}[\textbf{x}]$. Let us introduce a further operator $\omega$, defined by setting
\begin{align}
\label{omega2} 
\omega \cdot
h(x_1,\dots,x_n)
=
h(x_2, x_3, \ldots , x_n, q x_1),
\end{align} \index{0@$\omega$}
for any  $h \in \mathbb{C}[\textbf{x}]$. Collectively, the operators $T_1,\dots,T_{n-1},\omega$ give a polynomial representation of the affine Hecke algebra of type $A_{n-1}$. The {\it Cherednik--Dunkl operators} $Y_i$\index{Y@$Y_i$; Cherednik--Dunkl operators} generate an Abelian subalgebra of the affine Hecke algebra. They are given for each $1 \le i \le n$ by
\begin{align}
\label{Yi}
Y_i = 
T_{i - 1} \cdots T_1
\cdot
\omega 
\cdot
T_{n - 1}^{-1}
\cdots 
T_{i}^{-1},
\end{align}
and for each $i, j \in [1, n]$ satisfy the commutation relation
\begin{align}
\label{Y-commute}
[Y_i,Y_j] = 0.
\end{align}

\noindent In view of their commutativity, one can seek to jointly diagonalize the operators \eqref{Yi}. This brings us to the definition of the nonsymmetric Macdonald polynomials; see Proposition 2.3 of \cite{NPIVM}. In the below, for any composition $\lambda = (\lambda_1, \lambda_2, \ldots , \lambda_n)$ and sequence of $n$ variables $\textbf{x} = (x_1, x_2, \ldots , x_n)$, we denote $\textbf{x}^{\lambda} = \prod_{i = 1}^N x_i^{\lambda_i}$.

\begin{definition}
	Introduce the following two orders on length $n$ compositions 
	$\mu = (\mu_1,\dots,\mu_n)$ and $\nu = (\nu_1,\dots,\nu_n)$. The first is the dominance order, denoted $<$, and given by
	\begin{align*}
	\nu < \mu 
	\iff 
	\nu \ne \mu, \qquad \text{and} \qquad \sum_{i=1}^{j} \nu_i \le \sum_{i=1}^{j} \mu_i, \qquad \text{for each $j \in [1, n]$}.
	\end{align*}
	The second order, denoted $\prec$, is given by
	\begin{align*}
	\nu \prec \mu 
	\iff
	\Big(
	\nu^{+} < \mu^{+}
	\quad
	\text{or}
	\quad
	\nu^{+} = \mu^{+},
	\
	\nu < \mu
	\Big),
	\end{align*}
	where $\lambda^{+}$ is the unique partition that can be obtained by permuting the parts of a composition $\lambda$. The nonsymmetric Macdonald polynomials $f_{\mu} (\textbf{x}) = f_{\mu} (x_1, x_2, \ldots,  x_n) = f_{\mu}(x_1, x_2, \dots,x_n;q,t)$ are the unique family of polynomials in $\mathbb{C}[\textbf{x}]$ satisfying the triangularity property
	\begin{align}
	\label{monic-f}
	f_{\mu} (\textbf{x}) &= 
	\textbf{x}^{\mu}
	+ 
	\sum_{\overleftarrow{\nu} \prec \overleftarrow{\mu}} c_{\mu,\nu}(q,t)
	\textbf{x}^{\nu},
	\qquad 
	\text{for some $c_{\mu,\nu}(q,t) \in \mathbb{Q}(q,t)$},
	\end{align} \index{F@$f_{\mu} (\textbf{x})$; nonsymmetric Macdonald polynomial}
	with the sum taken over all compositions $\nu$ such that $\overleftarrow{\nu} \prec \overleftarrow{\mu}$, as well as the eigenvalue equation
	\begin{align}
	\label{eig-Yi}
	Y_i\cdot
	f_{\mu}(\textbf{x})
	=
	y_i(\mu;q,t)
	f_{\mu}(\textbf{x}),
	\quad
	\text{for each $i \in [1, n]$}, 
	\end{align}
	with the eigenvalues on the right-hand side of \eqref{eig-Yi} given by
	\begin{align}
	\label{rho}
	y_i(\mu;q,t)
	=
	q^{\mu_i}
	t^{\eta_i(\mu)+i-1};
	\quad
	\eta_i(\mu)= -\#\{j < i : \mu_j > \mu_i\} - \#\{j > i : \mu_j \ge \mu_i\}.
	\end{align} \index{Y@$y_i (\mu; q, t)$; eigenvalues of Cherednik--Dunkl operators} \index{0@$\eta_i (\mu)$}

\end{definition}

\begin{rem}
	\label{rem:alt-def}
	In the current work, it will be convenient to make use of an alternative characterization of the nonsymmetric Macdonald polynomials. Let $\mathbb{C}^D [\textbf{x}] = \mathbb{C}^{D}[x_1,\dots,x_n]$ be the set of all polynomials in 
	$\textbf{x} = (x_1, x_2, \dots,x_n)$ of total degree $\le D$. Viewing the latter as a vector space, we can write $\mathbb{C}^{D}[\textbf{x}] = {\rm Span}\{\textbf{x}^{\mu}\}_{|\mu| \le D}$.
	
	Introduce a generating series $Y(w) = \sum_{i=1}^{n} Y_i w^{i-1}$ of the operators \eqref{Yi}, which can be viewed as a linear operator $Y(w) : \mathbb{C}^{D}[\textbf{x}] \rightarrow \mathbb{C}^{D}[\textbf{x}]$. The corresponding eigenvalues $y(\mu;q,t;w) = \sum_{i=1}^{n} y_i(\mu;q,t) w^{i-1}$ are pairwise distinct for all $|\mu| \le D$, with $q,t,w$ living inside some open set in $\mathbb{C}$. The simplicity of the spectrum of $Y(w)$ (this simplicity plays an important role in Cherednik's theory of double affine Hecke algebras, \cf\ \cite[Theorem 8.2(i)]{IDA}) uniquely determines each eigenvector $f_{\mu}(\textbf{x})$, with $|\mu| \le D$, up to a multiplicative constant. One can therefore define $f_{\mu}$ as the unique polynomial solution of the equations \eqref{eig-Yi} such that $f_{\mu} \in \mathbb{C}^{|\mu|}[\textbf{x}]$ and ${\rm Coeff}[f_{\mu}; \textbf{x}^{\mu}] = 1$, where $\text{Coeff} [h; \textbf{x}^{\mu}]$ denotes the coefficient of $\textbf{x}^{\mu}$ in any $h \in \mathbb{C}[\textbf{x}]$. We will tacitly assume this definition in what follows.
	
\end{rem}

\section{Fermionic $L$-Matrix and Row Operators}

Our subsequent partition functions will be constructed in terms of another object, which we term the {\it $L$-matrix}. For any $x \in \mathbb{C}$, its entries $L_x (\textbf{A}, b; \textbf{C}, d)$ will be obtained from the specializations of the weights $W_{x; q} (\textbf{A}, b; \textbf{C}, d \boldsymbol{\mid} q^{-1 / 2}, s)$ given in \Cref{rql1} (where here we sometimes emphasize their dependence on $q$ in the subscript), by setting 
\begin{flalign}
\label{wl} 
L_x (\textbf{A}, b; \textbf{C}, d) = \displaystyle\lim_{s \rightarrow 0} (-s)^{-\textbf{1}_{d > 0}} W_{x / s; t} (\textbf{A}, b; \textbf{C}, d \boldsymbol{\mid} t^{-1/2}, s).
\end{flalign}
\index{L@$L_x (\textbf{A}, b; \textbf{C}, d)$}

\noindent Observe here that we have replaced the quantization parameter $q$ used for the weights in previous chapters with the Macdonald parameter $t$ here; this is because the symbol $q$ historically plays a different role in the context of Macdonald polynomials.  Further observe that these $L_x$ weights are given by the $r = q^{-1 / 2}$ cases of the $\mathcal{W}_x (\textbf{A}, \textbf{B}; \textbf{C}, \textbf{D} \boldsymbol{\mid} r)$ weights from \eqref{limitw} and are more explicitly through the following definition.

\begin{definition} 
	
	\label{lz} 
	
	For any $\textbf{A}, \textbf{C} \in \mathbb{Z}_{\ge 0}^n$ and $b, d \in \{ 0, 1, \ldots ,n \}$, set $L_x (\textbf{A}, b; \textbf{C}, d) = 0$ unless $\textbf{A}, \textbf{C} \in \{ 0, 1 \}^n$. Letting $\textbf{A}_i^+, \textbf{A}_j^-, \textbf{A}_{ij}^{+-}$ be as in \eqref{aij}, for each $\textbf{A} \in \{ 0, 1 \}^n$ and $i \in [1, n]$ set
\begin{flalign}
\label{lz1} 
\begin{aligned}
L_x (\textbf{A}, 0; \textbf{A}, 0) & = 1; \qquad L_x (\textbf{A}, 0; \textbf{A}_i^-, i) = t^{A_{[i + 1, n]}} x (1 - t); \\
  L_x (\textbf{A}, i; \textbf{A}_i^+, 0) & = 1; \qquad L_x (\textbf{A}, i; \textbf{A}; i) = (-1)^{A_i} t^{A_{[i, n]}} x.
  \end{aligned} 
\end{flalign}

\noindent Moreover, for any $1 \le i < j \le n$, set 
\begin{flalign}
\label{lz2} 
L_x (\textbf{A}, i; \textbf{A}_{ij}^{+-}, j) = t^{A_{[j + 1, n]}} x (1 - t); \qquad L_x (\textbf{A}, j; \textbf{A}_{ji}^{+-}, i) = 0.
\end{flalign}

\noindent We further set $L_x (\textbf{A}, \textbf{B}; \textbf{C}, \textbf{D} \boldsymbol{\mid} r, s) = 0$ if $(\textbf{A}, \textbf{B}; \textbf{C}, \textbf{D})$ is not of the above form. 

\end{definition} 

As before, we interpret $L_x (\textbf{A}, b; \textbf{C}, d)$ as the weight associated with a vertex $v$ whose arrow configuration is $(\textbf{A}, \textbf{e}_b; \textbf{C}, \textbf{e}_d)$. We refer to \Cref{vertexfigurerql} for a depiction of these weights; it is quickly verified that \eqref{wl} holds for them.

\begin{figure}[t]
	
	\begin{center}
		
		\begin{tikzpicture}[
		>=stealth,
		scale = .85
		]

		\draw[-, black] (-7.5, 3.1) -- (10, 3.1);
		\draw[-, black] (-7.5, -2.1) -- (10, -2.1);
		\draw[-, black] (-7.5, -1.1) -- (10, -1.1);
		\draw[-, black] (-7.5, -.4) -- (10, -.4);
		\draw[-, black] (-7.5, 2.4) -- (10, 2.4);
		\draw[-, black] (-7.5, -2.1) -- (-7.5, 3.1);
		\draw[-, black] (10, -2.1) -- (10, 3.1);
		\draw[-, black] (7.5, -2.1) -- (7.5, 3.1);
		\draw[-, black] (-5, -2.1) -- (-5, 2.4);
		\draw[-, black] (5, -2.1) -- (5, 3.1);
		\draw[-, black] (-2.5, -2.1) -- (-2.5, 2.4);
		\draw[-, black] (2.5, -2.1) -- (2.5, 2.4);
		\draw[-, black] (0, -2.1) -- (0, 3.1);

		\draw[->, thick, blue] (-6.3, .1) -- (-6.3, 1.9); 
		\draw[->, thick, green] (-6.2, .1) -- (-6.2, 1.9); 
		
		\draw[->, thick, blue] (-3.8, .1) -- (-3.8, 1) -- (-2.85, 1);
		\draw[->, thick, green] (-3.7, .1) -- (-3.7, 1.9);
		
		\draw[->, thick, blue] (-1.35, .1) -- (-1.35, 1.9);
		\draw[->, thick, green] (-1.25, .1) -- (-1.25, 1.9);
		\draw[->, thick,  orange] (-2.15, 1.1) -- (-1.15, 1.1) -- (-1.15, 1.9);
		
		\draw[->, thick, red] (.35, 1) -- (1.15, 1) -- (1.15, 1.9);
		\draw[->, thick, blue] (1.25, .1) -- (1.25, 1.9);
		\draw[->, thick, green] (1.35, .1) -- (1.35, 1.1) -- (2.15, 1.1);
		
		\draw[->, thick, blue] (3.65, .1) -- (3.65, 1) -- (4.65, 1);
		\draw[->, thick, green] (3.75, .1) -- (3.75, 1.9);
		\draw[->, thick, orange] (2.85, 1.1) -- (3.85, 1.1) -- (3.85, 1.9); 
		
		\draw[->, thick, red] (5.35, 1) -- (7.15, 1); 
		\draw[->, thick, blue] (6.2, .1) -- (6.2, 1.9);
		\draw[->, thick, green] (6.3, .1) -- (6.3, 1.9); 
		
		\draw[->, thick, blue] (7.85, 1) -- (9.65, 1); 
		\draw[->, thick, blue] (8.7, .1) -- (8.7, 1.9);
		\draw[->, thick, green] (8.8, .1) -- (8.8, 1.9); 
		
		\filldraw[fill=gray!50!white, draw=black] (-2.85, 1) circle [radius=0] node [black, right = -1, scale = .75] {$i$};
		\filldraw[fill=gray!50!white, draw=black] (2.15, 1) circle [radius=0] node [black, right = -1, scale = .75] {$j$};
		\filldraw[fill=gray!50!white, draw=black] (4.65, 1) circle [radius=0] node [black, right = -1, scale = .75] {$i$};
		\filldraw[fill=gray!50!white, draw=black] (7.15, 1) circle [radius=0] node [black, right = -1, scale = .75] {$i$};
		\filldraw[fill=gray!50!white, draw=black] (9.65, 1) circle [radius=0] node [black, right = -1, scale = .75] {$i$};
		
		\filldraw[fill=gray!50!white, draw=black] (7.85, 1) circle [radius=0] node [black, left = -1, scale = .75] {$i$};
		\filldraw[fill=gray!50!white, draw=black] (5.35, 1) circle [radius=0] node [black, left = -1, scale = .75] {$i$};
		\filldraw[fill=gray!50!white, draw=black] (2.85, 1) circle [radius=0] node [black, left = -1, scale = .75] {$j$};
		\filldraw[fill=gray!50!white, draw=black] (.35, 1) circle [radius=0] node [black, left = -1, scale = .75] {$i$};
		\filldraw[fill=gray!50!white, draw=black] (-2.15, 1) circle [radius=0] node [black, left = -1, scale = .75] {$i$};
		
		\filldraw[fill=gray!50!white, draw=black] (-6.25, 1.9) circle [radius=0] node [black, above = -1, scale = .75] {$\textbf{A}$};
		\filldraw[fill=gray!50!white, draw=black] (-3.75, 1.9) circle [radius=0] node [black, above = -1, scale = .75] {$\textbf{A}_i^-$};
		\filldraw[fill=gray!50!white, draw=black] (-1.25, 1.9) circle [radius=0] node [black, above = -1, scale = .75] {$\textbf{A}_i^+$};
		\filldraw[fill=gray!50!white, draw=black] (1.25, 1.9) circle [radius=0] node [black, above = -1, scale = .75] {$\textbf{A}_{ij}^{+-}$};
		\filldraw[fill=gray!50!white, draw=black] (3.75, 1.9) circle [radius=0] node [black, above = -1, scale = .75] {$\textbf{A}_{ji}^{+-}$};
		\filldraw[fill=gray!50!white, draw=black] (6.25, 1.9) circle [radius=0] node [black, above = -1, scale = .75] {$\textbf{A}$};
		\filldraw[fill=gray!50!white, draw=black] (8.75, 1.9) circle [radius=0] node [black, above = -1, scale = .75] {$\textbf{A}$};

		\filldraw[fill=gray!50!white, draw=black] (-6.25, .1) circle [radius=0] node [black, below = -1, scale = .7] {$\textbf{A}$};
		\filldraw[fill=gray!50!white, draw=black] (-3.75, .1) circle [radius=0] node [black, below = -1, scale = .7] {$\textbf{A}$};
		\filldraw[fill=gray!50!white, draw=black] (-1.25, .1) circle [radius=0] node [black, below = -1, scale = .7] {$\textbf{A}$};
		\filldraw[fill=gray!50!white, draw=black] (1.25, .1) circle [radius=0] node [black, below = -1, scale = .7] {$\textbf{A}$};
		\filldraw[fill=gray!50!white, draw=black] (3.75, .1) circle [radius=0] node [black, below = -1, scale = .7] {$\textbf{A}$};
		\filldraw[fill=gray!50!white, draw=black] (6.25, .1) circle [radius=0] node [black, below = -1, scale = .7] {$\textbf{A}$};
		\filldraw[fill=gray!50!white, draw=black] (8.75, .1) circle [radius=0] node [black, below = -1, scale = .75] {$\textbf{A}$};

		\filldraw[fill=gray!50!white, draw=black] (-3.75, 2.75) circle [radius=0] node [black] {$1 \le i \le n$};
		\filldraw[fill=gray!50!white, draw=black] (2.5, 2.75) circle [radius=0] node [black] {$1 \le i < j \le n$}; 
		\filldraw[fill=gray!50!white, draw=black] (6.25, 2.75) circle [radius=0] node [black] {$A_i = 0$};
		\filldraw[fill=gray!50!white, draw=black] (8.75, 2.75) circle [radius=0] node [black] {$A_i = 1$};
		
		\filldraw[fill=gray!50!white, draw=black] (-6.25, -.75) circle [radius=0] node [black, scale = .9] {$(\textbf{A}, 0; \textbf{A}, 0)$};
		\filldraw[fill=gray!50!white, draw=black] (-3.75, -.75) circle [radius=0] node [black, scale = .9] {$\big( \textbf{A}, 0; \textbf{A}_i^-, i \big)$};
		\filldraw[fill=gray!50!white, draw=black] (-1.25, -.75) circle [radius=0] node [black, scale = .9] {$\big( \textbf{A}, i; \textbf{A}_i^+, 0 \big)$};
		\filldraw[fill=gray!50!white, draw=black] (1.25, -.75) circle [radius=0] node [black, scale = .9] {$\big( \textbf{A}, i; \textbf{A}_{ij}^{+-}, j \big)$};
		\filldraw[fill=gray!50!white, draw=black] (3.75, -.75) circle [radius=0] node [black, scale = .9] {$\big( \textbf{A}, j; \textbf{A}_{ji}^{+-}, i \big)$};
		\filldraw[fill=gray!50!white, draw=black] (6.25, -.75) circle [radius=0] node [black, scale = .9] {$(\textbf{A}, i; \textbf{A}, i)$};
		\filldraw[fill=gray!50!white, draw=black] (8.75, -.75) circle [radius=0] node [black, scale = .9] {$(\textbf{A}, i; \textbf{A}, i)$};
		
		\filldraw[fill=gray!50!white, draw=black] (-6.25, -1.6) circle [radius=0] node [black, scale = .8] {$1$};
		\filldraw[fill=gray!50!white, draw=black] (-3.75, -1.6) circle [radius=0] node [black, scale = .8] {$t^{A_{[i + 1, n]}} x (1 - t)$};
		\filldraw[fill=gray!50!white, draw=black] (-1.25, -1.6) circle [radius=0] node [black, scale = .8] {$1$};
		\filldraw[fill=gray!50!white, draw=black] (1.25, -1.6) circle [radius=0] node [black, scale = .8] {$t^{A_{[j + 1, n]}} x (1 - t)$};
		\filldraw[fill=gray!50!white, draw=black] (3.75, -1.6) circle [radius=0] node [black, scale = .8] {$0$};
		\filldraw[fill=gray!50!white, draw=black] (6.25, -1.6) circle [radius=0] node [black, scale = .8] {$t^{A_{[i, n]}} x$};
		\filldraw[fill=gray!50!white, draw=black] (8.75, -1.6) circle [radius=0] node [black, scale = .8] {$- t^{A_{[i, n]}} x$};

		\end{tikzpicture}
		
	\end{center}
	
	\caption{\label{vertexfigurerql} The $L_x$ weights, and their arrow configurations, are depicted above. Here red, blue, green, and orange denote the colors $1$, $2$, $3$, and $4$, respectively.}
\end{figure}

\begin{rem}
	\label{rem:x-depend}
	
	The weights from \Cref{lz} admit a simple dependence on the parameter $x$. If $d \ge 1$, then $L_x (\textbf{A}, b; \textbf{C}, d)$ is a linear, homogeneous polynomial in $x$. If instead $d = 0$, then 
	$L_x (\textbf{A}, b; \textbf{C}, d)$ does not depend on $x$. Put another way, we obtain a factor of $x$ every time a path horizontally exits through a vertex.
\end{rem}

Given these weights, we next define row transfer operators on a vector space, as in (the homogeneous specializations of) \Cref{Row12Functions} and \Cref{OperatorRow}. To that end, we first require the associated single row partition functions, which are the analogs of those appearing in \Cref{mwabcd}. In what follows, we will typically index finitary sequences $\mathscr{X} = (\textbf{X}_0, \textbf{X}_1, \ldots)$ starting from $0$ instead of from $1$.

\begin{definition}
	
	\label{lwabcd}
	
	Fix a complex number $x \in \mathbb{C}$. For any finitary sequences $\mathscr{A} = (\textbf{A}_0, \textbf{A}_1, \ldots)$ and $\mathscr{C} = (\textbf{C}_0, \textbf{C}_1, \ldots )$ of elements in $\{ 0, 1 \}^n$, and indices $b, d \in [0, n]$, define 
	\begin{flalign}
	\label{lxac} 
	L_x (\mathscr{A}, b; \mathscr{C}, d) = \displaystyle\sum_{\mathfrak{J}} \displaystyle\prod_{i = 0}^{\infty} L_x (\textbf{A}_i, j_i; \textbf{C}_i, j_{i + 1}),
	\end{flalign} \index{L@$L_x (\mathscr{A}; b; \mathscr{C}; d)$}
	
	\noindent where the sum is over all sequences $\mathfrak{J} = (j_0, j_1, \ldots )$ of indices in $[0, n]$ such that $j_0 = b$ and $j_i = d$ for sufficiently large $i$ (we assume that the infinite product on the right side of \eqref{lxac} converges). By arrow conservation, the sum in \eqref{lxac} is supported on at most one term. 
\end{definition}

As in \Cref{mwabcd}, $L_x (\mathscr{A}, b; \mathscr{C}, d)$ is the partition function under the weights $L_x$ for the vertex model on a row, whose vertical entrance data is given by $(\textbf{A}_0, \textbf{A}_1, \ldots )$; horizontal entrance data by $\textbf{e}_b$; vertical exit data by $(\textbf{C}_0, \textbf{C}_1, \ldots)$; and horizontal exit data by $\textbf{e}_d$. See the left side of \Cref{wxyabcd} for the depiction of a similar row.

Next, we define row operators. As in \Cref{OperatorRow}, let $\mathbb{V}$ denote the infinite-dimensional vector space spanned by basis vectors of the form $| \mathscr{A} \rangle$, where $\mathscr{A}$ ranges over all finitary sequences of elements in $\{ 0, 1 \}^n$. Similarly, let $\mathbb{V}^*$ denote the space spanned all dual basis vectors of the form $\langle \mathscr{C} |$, over finitary $\mathscr{C}$. We impose an inner product on $\mathbb{V}^* \times \mathbb{V}$ determined by setting $\langle \mathscr{C} | \mathscr{A} \rangle = \textbf{1}_{\mathscr{A} = \mathscr{C}}$ for any finitary $\mathscr{A}$ and $\mathscr{C}$. 

\begin{definition} 
	
\label{dci} 

For any complex number $x \in \mathbb{C}$ and index $i \in [1, n]$, define the operators $\mathsf{C}_i (x): \mathbb{V} \rightarrow \mathbb{V}$\index{C@$\mathsf{C}_i (x)$} and $\mathsf{D} (x): \mathbb{V} \rightarrow \mathbb{V}$\index{D@$\mathds{D} (x; r)$!$\mathsf{D} (x)$} by, for any finitary sequence $\mathscr{A}$, setting
\begin{flalign*}
\mathsf{C}_i (x) | \mathscr{A} \rangle = \displaystyle\sum_{\mathscr{C}} L_x (\mathscr{A}, i; \mathscr{C}, 0) | \mathscr{C} \rangle; \qquad \mathsf{D} (x) | \mathscr{A} \rangle = \displaystyle\sum_{\mathscr{C}} L_x (\mathscr{A}, 0; \mathscr{C}, 0) | \mathscr{C} \rangle,
\end{flalign*}

\noindent where both sums are over all finitary sequences $\mathscr{C}$ of elements in $\{ 0, 1 \}^n$. In view of the inner product between $\mathbb{V}^*$ and $\mathbb{V}$, the operators $\mathsf{C}_i (x)$ and $\mathsf{D} (x)$ admit dual actions on $\mathbb{V}^*$ given by
\begin{flalign*}
& \langle \mathscr{C} | \mathsf{C}_i (x) = \displaystyle\sum_{\mathscr{A}} L_x ( \mathscr{A}, i; \mathscr{C}, 0) \langle \mathscr{A} |; \qquad  \langle \mathscr{C} | \mathsf{D} (x) = \displaystyle\sum_{\mathscr{A}}  L_x ( \mathscr{A}, 0; \mathscr{C}, 0) \langle \mathscr{A}|.
\end{flalign*}

\end{definition}

In particular, under $\mathsf{C}_i$ one color $i$ arrow horizontally enters the vertex model and none horizontally exit, and under $\mathsf{D}$ no arrows horizontally enter or exit the vertex model. 

By virtue of the ($r = q^{-1 / 2} = s$ case of) Yang--Baxter equation, \Cref{wabcdproduct2}, one can now derive the following commutation relations of the row operators $\mathsf{C}_i$ through the following proposition. We omit its proof, which is very similar to that of Theorem 3.2.1 of \cite{SVMST} (see also the proof of \Cref{bd2limit}). 
\begin{prop}
	\label{prop:CC}
	
	For any complex numbers $x, y \in \mathbb{C}$, we have 
	\begin{flalign}
	\label{dl1d}
	\mathsf{D} (x) \mathsf{D} (y) = \mathsf{D} (y) \mathsf{D} (x).
	\end{flalign}
	
	\noindent Moreover, for any integers $0 \le i < j \le n$ (where here we set $\mathsf{C}_0 (z) = \mathsf{D} (z)$ for any $z \in \mathbb{C}$), we have
	\begin{align}
	% \label{CC=}
	% (y-tx)\mathsf{C}_i(x) \mathsf{C}_i(y) &= (x-ty)\mathsf{C}_i(y) \mathsf{C}_i(x),
	% \\
	% \label{CC<}
	% t\, \mathsf{C}_i(x) \mathsf{C}_j(y) &= 
	% \frac{x-ty}{x-y}\, \mathsf{C}_j(y) \mathsf{C}_i(x) 
	% - 
	% \frac{(1-t)x}{x-y}\, \mathsf{C}_j(x) \mathsf{C}_i(y),
	% \quad\quad
	% i<j,
	% \\
	\label{CC>}
	\mathsf{C}_i(x) \mathsf{C}_j(y) &= 
	\frac{x-ty}{x-y}\, \mathsf{C}_j(y) \mathsf{C}_i(x) 
	- 
	\frac{(1-t)y}{x-y}\, \mathsf{C}_j(x) \mathsf{C}_i(y). 
	\end{align}
\end{prop}

\section{Partition Function Formula for $f_{\mu}(\textbf{x})$}

\label{ssec:states}

In this section we provide an expression for $f_{\mu} (\textbf{x})$ in terms of the $\mathsf{C}_i$ operators from \Cref{dci}. Before doing so, we require some definitions regarding the vector space $\mathbb{V}$. For any nonnegative composition $\lambda \in \mathbb{Z}_{\ge 0}^p$ of length $p \in [1, n]$, define the sequence $\mathscr{I} (\lambda) = \big( \textbf{I}_0 (\lambda), \textbf{I}_1 (\lambda), \ldots \big)$ of elements in $\{ 0, 1 \}^n$ so that the $k$-th coordinate of $\textbf{I}_j (\lambda)$ is given by $\textbf{1}_{\lambda_k = j}$, for each $k \in [1, p]$ and $j \ge 0$.

%In this section we present our main result, namely an explicit partition function expression for the nonsymmetric Macdonald polynomial $f_{\mu}(x_1,\dots,x_n)$, for any composition $\mu \in \mathbb{N}^n$; this is the content of \Cref{thm:f-formula}. The proof of this theorem will be given progressively over the course of the section; the elements of the proof are showing that the proposed formula \eqref{f-ansatz} satisfies {\bf 1.} certain exchange relations under the action of Hecke generators \eqref{Hecke}, and {\bf 2.} a particular cyclic relation under the action of the operator \eqref{omega}. Combining {\bf 1} and {\bf 2} allows one to show that \eqref{f-ansatz} satisfies the eigenvalue equation \eqref{eig-Yi}, which uniquely characterizes the nonsymmetric Macdonald polynomials (up to normalization).
%
%\Cref{thm:f-formula} expresses the nonsymmetric Macdonald polynomials in terms of a certain linear form which acts on a product of the row operators \eqref{C-row} and returns a scalar quantity; this is often termed a ``matrix product formula'' in the literature. In particular, the formula 
%\eqref{f-ansatz} was inspired by the earlier work of \cite{CantiniGW}, where a matrix product expression was obtained for a different family of nonsymmetric polynomials, which also yield (symmetric) Macdonald polynomials under appropriate symmetrization. We will comment more on the difference between the present work and that of \cite{CantiniGW} in \Cref{rem:cdgw}.

\begin{definition}
	
	\label{xlinear} 
	
	Let $\mu$ be a nonnegative composition of length $n$, and let its largest part be $N = \max_{i \in [1, n]} \mu_i$. Fix a linear operator $\mathsf{X} \in {\rm End}(\mathbb{V})$ and a set of indeterminates $\textbf{v} = (v_{i,j})$, which we will refer to as \emph{twist parameters} below, where $(i, j)$ ranges over all integer pairs $(i, j) \in [1, n] \times [0, N]$. Define the linear form
	\begin{align}
	\label{form}
	\langle \mathsf{X} \big\rangle_{\mu}
	\equiv
	\langle \mathsf{X} \big\rangle_{\mu}(\textbf{v})
	:=
	\sum_{\mathscr{M}}
	\prod_{i=1}^{n}
	\prod_{j=0}^{N}
	v_{i,j}^{M_{i,j}}
	\big\langle \mathscr{M} + \mathscr{I} (\mu) \big| \mathsf{X} \big| \mathscr{M} \big\rangle,
	\end{align}
	\noindent where we sum over sequences $\mathscr{M} = (\textbf{M}_0, \textbf{M}_1, \ldots )$ of elements in $\{ 0, 1 \}^n$ such that $\textbf{M}_k = \textbf{e}_0$ for $k > N$, and we have denoted $\textbf{M}_j = (M_{1, j}, M_{2, j}, \ldots , M_{n, j}) \in \{ 0, 1 \}^n$ for each $j \ge 0$. Observe that this sum is supported on only finitely many nonzero terms. 
	% Here, have defined
	%
	% \begin{align}
	% \ket{\{m\} + \mu}
	% =
	% \bigotimes_{j=0}^{N}
	% \ket{\bm{M}(j) + \bm{\mu}(j)}_j,
	% \qquad
	% \text{with}\ \ 
	% \bm{M}(j) + \bm{\mu}(j) = 
	% \sum_{i=1}^{n} (m_{i,j} + \bm{1}(\mu_i = j)) \textbf{e}_i.
	% \end{align}
	%
	% The summation in \eqref{form} is taken over $0 \le m_{i,j} \le 1$, for each $1 \le i \le n$ and $0 \le j \le N$.
	%; the exception is that the $m_{i,\mu_i}$ for $1 \le i \le n$ are each set to $0$.
\end{definition}

Now we can state the following theorem providing an expression for the nonsymmetric Macdonald polynomial $f_{\lambda} (\textbf{x})$ through the fermionic row operators $\mathsf{C}_i$ from \Cref{dci}; its proof will appear in \Cref{ssec:recast} and \Cref{ProofPolynomial} below. Let us mention that a similar bosonic expression for $f_{\lambda}$ appeared as Theorem 4.2 of \cite{NPIVM}, and our results and proofs here are largely based on those there. However, one subtantial difference between that result and the following one is that the latter provides a finite summation identity for $f_{\lambda}$, while the one from \cite{NPIVM} involved an infinite sum.

\begin{thm}
	\label{thm:f-formula}
	Fix a nonnegative composition $\mu$ of length $n$, set $N = \max_{i \in [1, n]} \mu_i$. Then, $f_{\mu}(\textbf{\emph{x}})$ is given in terms of the linear form \eqref{form} by
	\begin{align}
	\label{f-ansatz}
	f_{\mu}(\textbf{\emph{x}})
	=
	\Omega_{\mu}(q,t)
	\big\langle \mathsf{C}_n(x_n) \cdots \mathsf{C}_1(x_1) \big\rangle_{\mu},
	\end{align}
	where for $i \in [1, n]$ and $j \in [1, N]$ the underlying parameters $v_{i,j}$ are chosen to be
	\begin{align}
	\label{v_ij}
	v_{i,j}
	=
	q^{\mu_i-j}
	t^{\gamma_{i,j}(\mu)}
	\bm{1}_{\mu_i > j},
	\end{align} \index{V@$v_{i, j} (\mu)$; twist parameters}
	with exponents $\gamma_{i,j}(\mu)$ given by
	\begin{align}
	\label{gamma}
	\gamma_{i,j}(\mu)
	=
	-\#\{k < i : \mu_k > \mu_i\} + \#\{k > i : j \le \mu_k < \mu_i\}.
	\end{align} \index{0@$\gamma_{i, j} (\mu)$}
	The normalization factor $\Omega_{\mu}(q,t)$ appearing on the right-hand side of \eqref{f-ansatz} is given by
	\begin{align}
	\label{Omega}
	\Omega_{\mu}(q,t)
	=
	\prod_{i=1}^{n}
	\prod_{j=0}^{\mu_i-1}
	\frac{1}{1-q^{\mu_i-j} t^{\alpha_{i,j}(\mu)+1}},
	\end{align} \index{0@$\Omega_{\mu} (q, t)$}
	where we have defined
	\begin{align}
	\label{alpha}
	\alpha_{i,j}(\mu)
	=
	\#\{k < i: \mu_k = \mu_i\}
	+
	\#\{k \not= i: j < \mu_k < \mu_i\}
	+
	\#\{k > i: j = \mu_k\}.
	\end{align} \index{0@$\alpha_{i, j} (\mu)$}
\end{thm}

\section{Proof of \Cref{thm:f-formula}}
\label{ssec:recast}

In this section we establish \Cref{thm:f-formula}, assuming the following two propositions. The first, which we will establish in \Cref{ssec:proof} below, is a \emph{cyclic relation} that will essentially amount to the eigenvalue equation \eqref{eig-Yi} for $f_{\mu}$. The second, which we will establish in \Cref{ssec:norm} below, fixes the normalization \eqref{monic-f} for $f_{\mu}$. Here, we recall the exponents $\gamma_{i, j}$ and the quantity $\Omega_{\mu} (q, t)$ from \eqref{gamma} and \eqref{Omega}, respectively. 

\begin{prop}
	
	\label{omegaif}
	
	For any integer $i \in [1, n]$, we have 
	\begin{align}
	\label{cyclic-rel}
	\begin{aligned}
	\langle \mathsf{C}_i(qx_i) \mathsf{C}_n(x_n) \cdots & \mathsf{C}_{i + 1} (x_{i + 1}) \mathsf{C}_{i - 1} (x_{i - 1}) \cdots \mathsf{C}_1(x_1) \big\rangle_{\mu} \\
	& =
	q^{\mu_i}
	t^{\gamma_{i,0}(\mu)}
	\langle \mathsf{C}_n(x_n) \cdots \mathsf{C}_{i + 1} (x_{i + 1}) \mathsf{C}_{i - 1} (x_{i - 1}) \cdots \mathsf{C}_1(x_1) \mathsf{C}_i(x_i) \big\rangle_{\mu},
	\end{aligned} 
	\end{align}
	
	\noindent where the twist parameters associated with the linear form \eqref{form} are given by \eqref{v_ij}--\eqref{gamma}.
\end{prop}

\begin{prop}
	
	\label{coefficientc} 
	
	We have 
	\begin{align}
	\label{Omega-compute}
	{\rm Coeff}\Big[ 
	\Omega_{\mu}(q,t)
	\big\langle \mathsf{C}_n(x_n) \cdots \mathsf{C}_1(x_1) \big\rangle_{\mu}; \textbf{\emph{x}}^{\mu}
	\Big]
	=
	1,
	\end{align}
	
	\noindent where the twist parameters associated with the linear form \eqref{form} are given by \eqref{v_ij}--\eqref{gamma}.
\end{prop}

Before proceeding to the proof of \Cref{thm:f-formula} (assuming the two results above), we first modify \eqref{eig-Yi} slightly, by transferring the string of Hecke generators $T_{i-1} \cdots T_1$ contained within the Cherednik--Dunkl operator $Y_i$ to the right-hand side of the relation. One obtains
\begin{align}
\label{eig-recast}
\omega \cdot T^{-1}_{n-1} \cdots T^{-1}_i
\cdot
f_{\mu}(\textbf{x})
=
y_i(\mu;q,t)
T^{-1}_1 \cdots T^{-1}_{i-1}
\cdot
f_{\mu}(\textbf{x}).
\end{align}
We wish to substitute the Ansatz \eqref{f-ansatz} into \eqref{eig-recast}, and check its validity. In order to do so, we require the following auxiliary result.
\begin{prop} 
	For any integers $i \in [1, n - 1]$ and $1 \le j<k \le n$, we have
	\begin{align}
	\label{T-CC}
	T^{-1}_i
	\cdot
	\mathsf{C}_k(x_{i+1}) \mathsf{C}_j(x_i)
	=
	t^{-1}
	\mathsf{C}_j(x_{i+1}) \mathsf{C}_k(x_i).
	\end{align}
\end{prop}

\begin{proof}
	This follows by direct calculation, using the explicit form \eqref{Hecke} of $T^{-1}_i$. We compute
	\begin{align*}
	T^{-1}_i
	\cdot
	\mathsf{C}_k(x_{i+1}) \mathsf{C}_j(x_i)
	&=
	t^{-1}\left(
	\frac{(t-1)x_i}{x_{i+1}-x_i}
	\mathsf{C}_k(x_{i+1}) \mathsf{C}_j(x_i)
	+
	\frac{x_{i+1}-tx_i}{x_{i+1}-x_i}
	\mathsf{C}_k(x_i) \mathsf{C}_j(x_{i+1})
	\right)
	\\
	&=
	t^{-1}
	\mathsf{C}_j(x_{i + 1}) \mathsf{C}_k(x_i),
	\end{align*}
	where the final equality follows from the exchange relation \eqref{CC>}.
\end{proof}

Now we can establish \Cref{thm:f-formula}, assuming \Cref{omegaif} and \Cref{coefficientc}. 

\begin{proof}[Proof of \Cref{thm:f-formula} Assuming \Cref{omegaif} and \Cref{coefficientc}]
	
	By \Cref{rem:alt-def}, it suffices to verify that our Ansatz \eqref{f-ansatz} satisfies both the eigenvalue equation \eqref{eig-Yi} and the normalization constraint \eqref{Omega-compute}, with $\Omega_{\mu}(q,t)$. The latter follows by \Cref{coefficientc}, so it remains to establish \eqref{eig-Yi} or, equivalently, \eqref{eig-recast}. 
	
	To that end, let us examine the left-hand side of \eqref{eig-recast}, assuming the substitution \eqref{f-ansatz}. It becomes
\begin{multline*}
\Omega_{\mu}(q,t) 
\cdot
\omega 
\cdot
T^{-1}_{n-1} \cdots T^{-1}_i
\cdot
\langle \mathsf{C}_n(x_n) \cdots \mathsf{C}_1(x_1) \big\rangle_{\mu} 
\\
=
\Omega_{\mu}(q,t) t^{i-n}
\cdot
\omega
\cdot
\langle \mathsf{C}_i (x_n) \mathsf{C}_n(x_{n - 1}) \cdots \mathsf{C}_{i+1}(x_i)
\mathsf{C}_{i-1}(x_{i - 1}) \cdots \mathsf{C}_1(x_1) \big\rangle_{\mu},
\end{multline*}
\noindent where to deduce the equality we used the linearity of the form \eqref{form}, together with $n-i$ applications of the identity \eqref{T-CC} to transfer the operator $\mathsf{C}_i$ towards the left of the operator product. Finally we act explicitly with $\omega$, which by \eqref{omega2} yields
\begin{multline}
\label{lhs}
\Omega_{\mu}(q,t) 
\cdot
\omega 
\cdot
T^{-1}_{n-1} \cdots T^{-1}_i
\cdot
\langle \mathsf{C}_n(x_n) \cdots \mathsf{C}_1(x_1) \big\rangle_{\mu} 
\\
=
\Omega_{\mu}(q,t) t^{i-n}
\cdot
\langle \mathsf{C}_i (qx_1) \mathsf{C}_n(x_n) \cdots \mathsf{C}_{i+1}(x_{i + 1})
\mathsf{C}_{i - 1}(x_i) \cdots \mathsf{C}_1(x_2) \big\rangle_{\mu}.
\end{multline}

Similarly, under the substitution \eqref{f-ansatz}, the right-hand side of \eqref{eig-recast} is given by
\begin{multline}
\label{rhs}
y_i(\mu;q,t)
\Omega_{\mu}(q,t)
\cdot
T^{-1}_1 \cdots T^{-1}_{i-1}
\cdot
\langle \mathsf{C}_n(x_n) \cdots \mathsf{C}_1(x_1) \big\rangle_{\mu},
\\
=
t^{1-i}
y_i(\mu;q,t)
\Omega_{\mu}(q,t)
\cdot
\langle 
\mathsf{C}_n(x_n) \cdots \mathsf{C}_{i+1}(x_{i + 1})
\mathsf{C}_{i-1}(x_i) \cdots \mathsf{C}_1(x_2)
\mathsf{C}_i(x_1)
\big\rangle_{\mu}.
\end{multline}
\noindent where to deduce the equality we again used the linearity of the form \eqref{form}, together with $i-1$ applications of \eqref{T-CC} to send the operator $\mathsf{C}_i$ towards the right of the operator product. 

Matching the right-hand sides of \eqref{lhs} and \eqref{rhs}, and applying $\mathfrak{s}_{i - 1} \cdots \mathfrak{s}_1$ to both sides of the resulting equation, we are left with the task of proving that
\begin{flalign}
\label{cii1}
\begin{aligned}
\big\langle \mathsf{C}_i (qx_i) \mathsf{C}_n (x_n) & \cdots \mathsf{C}_{i + 1} (x_{i + 1}) C_{i - 1} (x_{i - 1}) \cdots \mathsf{C}_1 (x_1) \big\rangle_{\mu} \\
&  =
t^{n-2i+1} y_i(\mu;q,t)
\big\langle \mathsf{C}_n (x_n) \cdots \mathsf{C}_{i + 1} (x_{i + 1}) \mathsf{C}_{i - 1} (x_{i - 1}) \cdots \mathsf{C}_1 (x_1) \mathsf{C}_i(x_i) \big\rangle_{\mu}.
\end{aligned} 
\end{flalign}
Using the explicit form of the eigenvalue $y_i (\mu; q, t)$ from \eqref{rho}, together with the fact that 
\begin{align*}
n-i+\eta_i(\mu)
=
-\#\{j < i : \mu_j > \mu_i\} + \#\{j > i : \mu_j < \mu_i\}
\equiv
\gamma_{i,0}(\mu),
\end{align*}

\noindent \eqref{cii1} follows from \Cref{omegaif}. This establishes \eqref{eig-recast} and therefore the theorem.
\end{proof}

\chapter{Proof of \Cref{omegaif} and \Cref{coefficientc}} 

\label{sec:matrix-prod}

\label{ProofPolynomial}

By the content in \Cref{ssec:recast}, to establish \Cref{thm:f-formula} it suffices to show the cyclic relation \Cref{omegaif} and the normalization condition \Cref{coefficientc}. In this chapter we perform these tasks, closely following the proofs of equation (4.14) and Section 4.10 of \cite{NPIVM}, respectively. 

\section{Column Operators and Diagrammatic Notation}

\label{Column}

Our main tool for proving \Cref{omegaif} will be another class of objects called column operators, constructed in terms of the following single-column partition functions. 

\begin{definition} 
	
\label{column} 

Fix two sequences $\mathfrak{B} = (b_1, b_2, \ldots , b_n)$ and $\mathfrak{D} = (d_1, d_2, \ldots , d_n)$ of indices in $[0, n]$, as well as two elements $\textbf{A}, \textbf{C} \in \{ 0, 1 \}^n$. For any set of complex numbers $\textbf{x} = (x_1, x_2, \ldots , x_n)$, let 
\begin{flalign}
\label{columnl}
L_{\textbf{x}} (\textbf{A}, \mathfrak{B}; \textbf{C}, \mathfrak{D}) = \displaystyle\sum_{\mathscr{K}} \displaystyle\prod_{i = 1}^n L_{x_i} (\textbf{K}_i, b_i; \textbf{K}_{i + 1}, d_i),
\end{flalign} \index{L@$L_{\textbf{x}} (\textbf{A}, \mathfrak{B}; \textbf{C}, \mathfrak{D})$}

\noindent where we sum over all sequences $\mathscr{K} = (\textbf{K}_1, \textbf{K}_2, \ldots , \textbf{K}_{n + 1})$ of elements in $\{ 0, 1 \}^n$ such that $\textbf{K}_1 = \textbf{A}$ and $\textbf{K}_{n + 1} = \textbf{C}$. By arrow conservation, the sum in \eqref{columnl} is supported on at most one term. We also set $L_{\textbf{x}} (\textbf{A}, \mathfrak{B}; \textbf{C}, \mathfrak{D}) = 0$ if either $\textbf{A}$ or $\textbf{C}$ is not in $\{ 0, 1 \}^n$. 

\end{definition} 

In particular, $L_{\textbf{x}} (\textbf{A}, \mathfrak{B}; \textbf{C}, \mathfrak{D})$ denotes the partition function for a column with vertical entrance data $\textbf{A}$; horizontal entrance data $(b_1, b_2, \ldots , b_n)$; vertical exit data $\textbf{C}$; and horizontal exit data $(d_1, d_2, \ldots , d_n)$. In what follows, instead of using the algebraic expressions \eqref{columnl} for partition functions, it will often be convenient to express them diagrammatically, for example by 
\begin{align}
\label{tower}
L_{\textbf{x}} (\textbf{A}, \mathfrak{B}; \textbf{C}, \mathfrak{D})
:=
\tikz{.75}{
	\foreach\y in {1,...,4}{
		\draw[lgray,line width=1.5pt] (1.5,0.5+\y) -- (2.5,0.5+\y);
	}
	\draw[lgray,line width=1.5pt] (1.5,0.5) -- (2.5,0.5) -- (2.5,5.5) -- (1.5,5.5) -- (1.5,0.5);
	\node at (2,1.5) {\footnotesize$\bullet$}; \draw (2,1.5) -- (4,1.5); \node[right] at (4,1.5) {\footnotesize$\textbf{K}_2$};
	\node at (2,2.5) {\footnotesize$\bullet$}; \draw (2,2.5) -- (4,2.5); \node[right] at (4,2.5) {\footnotesize$\textbf{K}_3$}; 
	\node at (2,4.5) {\footnotesize$\bullet$}; \draw (2,4.5) -- (4,4.5);
	\node[right] at (4,4.5) {\footnotesize$\textbf{K}_{n-1}$};
	%spectral parameters
	\node[text centered] at (2,1) {$x_1$};
	\node[text centered] at (2,2) {$x_2$};
	\node[text centered] at (2,3.1) {$\vdots$};
	\node[text centered] at (2,4.1) {$\vdots$};
	\node[text centered] at (2,5) {$x_n$};
	%bottom labels
	\node[below] at (2,0.5) {\footnotesize$\textbf{A}$};
	%top labels
	\node[above] at (2,5.5) {\footnotesize$\textbf{C}$};
	%right labels
	\node[right] at (2.5,1) {$d_1$};
	\node[right] at (2.5,2) {$d_2$};
	\node[right] at (2.5,3.1) {$\vdots$};
	\node[right] at (2.5,4.1) {$\vdots$};
	\node[right] at (2.5,5) {$d_n$};
	%left labels
	\node[left] at (1.5,1) {$b_1$};
	\node[left] at (1.5,2) {$b_2$};
	\node[left] at (1.5,3.1) {$\vdots$};
	\node[left] at (1.5,4.1) {$\vdots$};
	\node[left] at (1.5,5) {$b_n$};
}
\end{align}
where the indices in the diagram on the right side of \eqref{tower} are as follows. The $x_i$ denote the spectral parameters assigned to each vertex (represented by a face in \eqref{tower}); the indices $(b_1, b_2, \ldots , b_n)$ along the left boundary denote the colors of the arrows horizontally entering the column; the indices $(d_1, d_2, \ldots , d_n)$ along the right boundary denote the colors of the arrows horizontally exiting the colum; the $n$-tuples $\textbf{A}$ and $\textbf{C}$ denote the boundary data for vertical arrows along the bottom and top boundaries, respectively; and the $n$-tuples $\textbf{K}_i$ along the middle of the column denote the collection of vertical arrows proceeding from vertex $i - 1$ to vertex $i$. 

Next, to define column operators, let $W \cong \mathbb{C}^{n+1}$ be an 
$(n+1)$-dimensional vector space with basis $\{\ket{0},\ket{1},\dots,\ket{n}\}$.
The basis of the dual vector space $W^{*}$ will be written as $\{\bra{0},\bra{1},\dots,\bra{n}\}$ with $\langle i | j \rangle = \textbf{1}_{i = j}$. We wish to consider an $n$-fold tensor product of $W$, denoted by
\begin{align*}
\mathbb{W}
=
W_1 \otimes \cdots \otimes W_n
=
{\rm Span}\left\{
\bigotimes_{k=1}^{n}
\ket{i_k}_k
\right\}_{0 \le i_1,\dots,i_n \le n}
\equiv
{\rm Span} \big\{
\ket{\mathfrak{I}} \big\}_{\mathfrak{I} \in [0, n]^n},
\end{align*}

\noindent where we have set $\mathfrak{I} = (i_1, i_2, \ldots , i_n)$. We now define operators acting on $\mathbb{W}$.

\begin{definition}
	Fix sequences $\textbf{x} = (x_1, x_2, \ldots , x_n)$ and $\textbf{v} = (v_1,\dots,v_n)$ of complex numbers, and fix an element $\textbf{C} \in \{0,1\}^n$. We introduce the \emph{column operator} $\Psi_{\textbf{v}} (\textbf{C}) = \Psi_{\textbf{v}; \textbf{x}} (\textbf{C}) \in {\rm End}(\mathbb{W})$, defined by setting for any sequence $\mathfrak{D} = (d_1, d_2, \ldots , d_n)$ of indices in $[0, n]$ 
	\begin{flalign}
	\label{Q-def}
	\Psi_{\textbf{v}}(\textbf{C}) \ket{\mathfrak{D}} 
	= \sum_{\mathfrak{B}}
	\Bigg(
	\sum_{\textbf{M} \in \{0,1\}^n}
	\prod_{k=1}^{n}
	v_k^{M_k}
	L_{\textbf{x}} (\textbf{M}, \mathfrak{B}; \textbf{M} + \textbf{C}, \mathfrak{D})
	\Bigg)
	\ket{\mathfrak{B}},
	\end{flalign} \index{0@$\Psi_{\textbf{v}} (\textbf{C})$}

	\noindent and then extending its action to all of $\mathbb{W}$ by linearity. The first sum in \eqref{Q-def} is over all sequences $\mathfrak{B} = (b_1, b_2, \ldots , b_n)$ of indices in $[0, n]$, and the second over all elements $\textbf{M} = (M_1, M_2, \dots, M_n) \in \{0,1\}^n$.
\end{definition}

The reason for introducing the column operators \eqref{Q-def} is that they offer an alternative algebraic setup for operator products of the form \eqref{f-ansatz}, which will be convenient for studying the cyclic relation \eqref{cyclic-rel}. In particular, we note the following, where in the below we recall from \Cref{ssec:states} the sequence $\big( \textbf{I}_0 (\lambda), \textbf{I}_1 (\lambda), \ldots \big)$ of elements in $\{ 0, 1 \}^n$. 
\begin{prop}
	Fix a composition $\mu$ of length $n$, and let $N = \max_{i \in [1, n]} \mu_i$. Then, denoting $\textbf{\emph{v}}^{(j)} = (v_{1, j}, v_{2, j}, \ldots, v_{n, j})$ for $j \in [0, N]$, where the parameters $v_{i, j}$ are as in \eqref{v_ij}, we have 
	\begin{align}
	\label{column-decomp}
	\big\langle \mathsf{C}_n(x_n) \cdots \mathsf{C}_1(x_1) \big\rangle_{\mu}
	=
	\big\langle 1,\dots,n \big|
	\Psi_{\textbf{\emph{v}}^{(0)}} \big( \textbf{\emph{I}}_0 (\mu) \big) \cdots  \Psi_{\textbf{\emph{v}}^{(N)}} \big( \textbf{\emph{I}}_j (\mu) \big)
	\big| 0,\dots,0 \big\rangle.
	\end{align}
\end{prop}

\begin{proof}
	By the expression of the row operators from \Cref{dci} and the linear form \eqref{form}, we obtain the following partition function representation of $\big\langle\mathsf{C}_1(x_1) \cdots \mathsf{C}_n(x_n) \big\rangle_{\mu}$: 
	\begin{align}
	\label{lattice}
	\big\langle\mathsf{C}_n(x_n) \cdots \mathsf{C}_1(x_1) \big\rangle_{\mu}
	=
	\sum_{\mathscr{M}}
	\prod_{i=1}^{n}
	\prod_{j=0}^{N}
	v_{i,j}^{M_{i,j}}
	\times
	\tikz{.85}{
		\foreach\y in {0,...,5}{
			\draw[lgray,line width=1.5pt] (1.5,0.5+\y) -- (7.5,0.5+\y);
		}
		\foreach\x in {0,...,6}{
			\draw[lgray,line width=1.5pt] (1.5+\x,0.5) -- (1.5+\x,5.5);
		}
		%spectral parameters
		\node[text centered] at (2,1) {$x_1$};
		\node[text centered] at (4,1) {$\cdots$};
		\node[text centered] at (5,1) {$\cdots$};
		\node[text centered] at (7,1) {$x_1$};
		\node[text centered] at (2,2) {$x_2$};
		\node[text centered] at (4,2) {$\cdots$};
		\node[text centered] at (5,2) {$\cdots$};
		\node[text centered] at (7,2) {$x_2$};
		\node[text centered] at (2,5) {$x_n$};
		\node[text centered] at (4,5) {$\cdots$};
		\node[text centered] at (5,5) {$\cdots$};
		\node[text centered] at (7,5) {$x_n$};
		\node[text centered] at (2,3.1) {$\vdots$};
		\node[text centered] at (2,4.1) {$\vdots$};
		\node[text centered] at (7,3.1) {$\vdots$};
		\node[text centered] at (7,4.1) {$\vdots$};
		%bottom labels
		\node[below] at (2,0.5) {\scriptsize$\textbf{M}_0$};
		\node[above,text centered] at (4,0) {$\cdots$};
		\node[above,text centered] at (5,0) {$\cdots$};
		\node[below] at (7,0.5) {\scriptsize$\textbf{M}_N$};
		%top labels
		\node[above] at (2,5.5) {\scriptsize$\textbf{M}_0+ \textbf{I}_0 (\mu)$};
		\node[above,text centered] at (4,5.5) {$\cdots$};
		\node[above,text centered] at (5,5.5) {$\cdots$};
		\node[above] at (7,5.5) {\scriptsize$\textbf{M}_N+\textbf{I}_N (\mu)$};
		%right labels
		\node[right] at (7.5,1) {$0$};
		\node[right] at (7.5,2) {$0$};
		\node[right] at (7.5,3.1) {$\vdots$};
		\node[right] at (7.5,4.1) {$\vdots$};
		\node[right] at (7.5,5) {$0$};
		%left labels
		\node[left] at (1.5,1) {$1$};
		\node[left] at (1.5,2) {$2$};
		\node[left] at (1.5,3.1) {$\vdots$};
		\node[left] at (1.5,4.1) {$\vdots$};
		\node[left] at (1.5,5) {$n$};
	}
	\end{align}
	where we sum over all sequences $\mathscr{M} = (\textbf{M}_0, \textbf{M}_1, \ldots , \textbf{M}_N)$ of elements in $\{ 0, 1 \}^n$, with $\textbf{M}_j = (M_{1, j}, M_{2, j}, \ldots , M_{n, j})$ for each $j \in [0, N]$. Under this diagrammatic notation, summation over all possible states is implied at each internal lattice edge.
	
	Now we read the partition function \eqref{lattice} column-by-column, starting from the rightmost and working to the left. We distribute the product $\prod_{i=1}^{n} \prod_{j=0}^{N} v_{i,j}^{M_{i,j}}$ over the columns by assigning the factor $\prod_{i=1}^{n} v_{i,j}^{M_{i,j}}$ to column $j$, and for each $0 \le j \le N$ we compute separately the sums over $\textbf{M}_j = (M_{1,j},\dots,M_{n,j})$ (which only affects the boundary states of column $j$). Referring to the definition \eqref{Q-def} of 
	$\Psi_{\textbf{v}^{(j)}} \big( \textbf{I}_j (\mu) \big)$, equation \eqref{column-decomp} is nothing more than the algebraic realization of this column-by-column decomposition.
\end{proof}

The partition function \eqref{lattice} is useful for visualizing the structure of $\big\langle\mathsf{C}_1(x_1) \dots \mathsf{C}_n(x_n) \big\rangle_{\mu}$, but it is still not in the most convenient form. In particular, we would like to blend the summation over $\mathscr{M}$ into our graphical conventions. With that in mind, we set up the following notation.
\begin{definition}
	Let $\Psi_{\textbf{v}}(\textbf{C}) = \Psi_{\textbf{v}; \textbf{x}} (\textbf{C})$ be a column operator as defined in \eqref{Q-def}. We shall simplify the graphical representation of its components by for any $\mathfrak{B} = (b_1, b_2, \ldots , b_n)$ and $\mathfrak{D} = (d_1, d_2, \ldots , d_n)$ writing
	\begin{align*}
	\bra{\mathfrak{B}}
	\Psi_{\textbf{v}}(\textbf{C})
	\ket{\mathfrak{D}}
	=
	\sum_{\textbf{M} \in \{0,1\}^n}
	\prod_{k=1}^{n}
	v_k^{M_k} \times
	\tikz{.9}{
		\foreach\y in {1,...,4}{
			\draw[lgray,line width=1.5pt] (1.5,0.5+\y) -- (2.5,0.5+\y);
		}
		\draw[lgray,line width=1.5pt] (1.5,0.5) -- (2.5,0.5) -- (2.5,5.5) -- (1.5,5.5) -- (1.5,0.5);
		%spectral parameters
		\node[text centered] at (2,1) {$x_1$};
		\node[text centered] at (2,2) {$x_2$};
		\node[text centered] at (2,3.1) {$\vdots$};
		\node[text centered] at (2,4.1) {$\vdots$};
		\node[text centered] at (2,5) {$x_n$};
		%bottom labels
		\node[below] at (2,0.5) {\footnotesize$\textbf{M}$};
		%top labels
		\node[above] at (2,5.5) {\footnotesize$\textbf{M}+\textbf{C}$};
		%right labels
		\node[right] at (2.5,1) {$d_1$};
		\node[right] at (2.5,2) {$d_2$};
		\node[right] at (2.5,3.1) {$\vdots$};
		\node[right] at (2.5,4.1) {$\vdots$};
		\node[right] at (2.5,5) {$d_n$};
		%left labels
		\node[left] at (1.5,1) {$b_1$};
		\node[left] at (1.5,2) {$b_2$};
		\node[left] at (1.5,3.1) {$\vdots$};
		\node[left] at (1.5,4.1) {$\vdots$};
		\node[left] at (1.5,5) {$b_n$};
	}
	=
	\tikz{.9}{
		\foreach\y in {1,...,4}{
			\draw[lgray,line width=1.5pt] (1.5,0.5+\y) -- (2.5,0.5+\y);
		}
		\draw[lgray,line width=1.2pt] (2.5,0.5) -- (2.5,5.5);
		\draw[double,lgray,line width=1.7pt] (2.5,5.5) -- (1.5,5.5);
		\draw[lgray,line width=1.2pt] (1.5,5.5) -- (1.5,0.5);
		\draw[double,lgray,line width=1.7pt] (1.5,0.5) -- (2.5,0.5);
		%spectral parameters
		\node[text centered] at (2,1) {$x_1$};
		\node[text centered] at (2,2) {$x_2$};
		\node[text centered] at (2,3.1) {$\vdots$};
		\node[text centered] at (2,4.1) {$\vdots$};
		\node[text centered] at (2,5) {$x_n$};
		%bottom labels
		\node[below] at (2,0.5) {\phantom{\footnotesize$\textbf{M}$}};
		%top labels
		\node[above] at (2,5.5) {\footnotesize$\textbf{C}$};
		%right labels
		\node[right] at (2.5,1) {$d_1$};
		\node[right] at (2.5,2) {$d_2$};
		\node[right] at (2.5,3.1) {$\vdots$};
		\node[right] at (2.5,4.1) {$\vdots$};
		\node[right] at (2.5,5) {$d_n$};
		%left labels
		\node[left] at (1.5,1) {$b_1$};
		\node[left] at (1.5,2) {$b_2$};
		\node[left] at (1.5,3.1) {$\vdots$};
		\node[left] at (1.5,4.1) {$\vdots$};
		\node[left] at (1.5,5) {$b_n$};
	}
	\end{align*}
	where $\textbf{M} = (M_1, M_2, \ldots , M_n) \in \{ 0, 1 \}^n$, and the dependence on $\textbf{v}$ is kept implicit in the object appearing on the right-hand side. Observe here the double line along the top and bottom of the column on the right side.
\end{definition}

\section{Explicit Computation of Column Operator Components}

In this section we turn to the question of explicitly computing the matrix elements of the linear operators \eqref{Q-def}. We begin with some auxiliary definitions for admissible vectors, color data, and coordinates. 

\begin{definition}
	\label{definition:admiss}
	For any sequences $\mathfrak{B} = (b_1, b_2, \dots, b_n)$ and $\mathfrak{D} = (d_1, d_2, \dots,d_n)$ of indices in $[0, n]$, we say that the pair $(\mathfrak{B},\mathfrak{D})$ is {\it admissible} if for each $k \in [1, n]$ we have 
	\begin{align}
	\label{mult-admiss}
	0 \le m_k(\mathfrak{D}) \le m_k(\mathfrak{B}) \le 1,
	\end{align}
	\noindent where $m_k (\mathfrak{B})$ and $m_k(\mathfrak{D})$ denote the multiplicities of $k$ in $\mathfrak{B}$ and $\mathfrak{D}$, respectively.
\end{definition}

\begin{definition}
	\label{definition:data}
	Given an admissible pair of sequences $(\mathfrak{B},\mathfrak{D})$, as in \Cref{definition:admiss}, we introduce two disjoint sets $\mathcal{P} \subset \{1,\dots,n\}$ and $\mathcal{Q} \subset \{1,\dots,n\}$, where $p \in \mathcal{P}$ if and only if $m_p(\mathfrak{B}) = 1$ and
	$m_p(\mathfrak{D}) = 0$,  while $p \in \mathcal{Q}$ if and only if 
	$m_p(\mathfrak{B}) = m_p(\mathfrak{D}) = 1$. We refer to $(\mathcal{P},\mathcal{Q})$ as the {\it color data} associated to $(\mathfrak{B},\mathfrak{D})$.
\end{definition}

\begin{definition}
	\label{definition:coord}
	Let $\mathfrak{B} = (b_1, b_2, \dots, b_n)$ and $\mathfrak{D} = (d_1, d_2, \dots, d_n)$ be sequences such that the pair $(\mathfrak{B}, \mathfrak{D})$ is admissible, as in \Cref{definition:admiss}; also let $(\mathcal{P},\mathcal{Q})$ be their associated color data, as in \Cref{definition:data}. Introduce another two vectors
	\begin{align}
	\label{coord-def}
	\mathcal{I} = (i_p)_{p \in \mathcal{P} \cup \mathcal{Q}},
	\qquad 
	\mathcal{J} = (j_p)_{p \in \mathcal{Q}},
	\end{align} 
	such that
	\begin{align}
	\label{coord-def2}
	b_{i_p} = p, \quad \text{for each $p \in \mathcal{P} \cup \mathcal{Q}$};
	\qquad\qquad
	d_{j_p} = p, \quad \text{for each $p \in \mathcal{Q}$}.
	\end{align}
	We shall call $(\mathcal{I},\mathcal{J})$ the {\it coordinates} of 
	$(\mathfrak{B},\mathfrak{D})$.
\end{definition}

\begin{example}
	In order to illustrate the above definitions, let us choose $n=5$, and $\mathfrak{B} = (0,2,1,5,3)$, $\mathfrak{D} = (2,0,0,0,5)$. These vectors clearly satisfy the requirements \eqref{mult-admiss}, and are thus admissible. We find that $m_p (\mathfrak{B}) = 1$ and $m_p(\mathfrak{D}) = 0$ for $p \in \{1,3\}$, while $m_p(\mathfrak{B}) = m_p (\mathfrak{D}) = 1$ for $p \in \{2,5\}$. Hence, $\mathcal{P} = \{1,3\}$ and $\mathcal{Q} = \{2,5\}$ is the color data associated to $(\mathfrak{B},\mathfrak{D})$. Finally, one finds $\mathcal{I} = (i_1,i_2,i_3,i_5) = (3,2,5,4)$ by reading the positions of $\{1,2,3,5\}$ in $\mathfrak{B}$, and $\mathcal{J} = (j_2, j_5) = (1,5)$ by reading the positions of $\{2,5\}$ in $\mathfrak{D}$.
	
	To illustrate how we will make use of such quantities, consider the following example of the column \eqref{tower} for $n=5$, $\textbf{A} = (0,1,0,0,0)$, $\textbf{C} = (1,1,1,0,0)$, and $(\mathfrak{B}, \mathfrak{D})$ as above, given by
	\begin{align*}
	L(\textbf{A}, \mathfrak{B}; \textbf{C}, \mathfrak{D})
	=
	
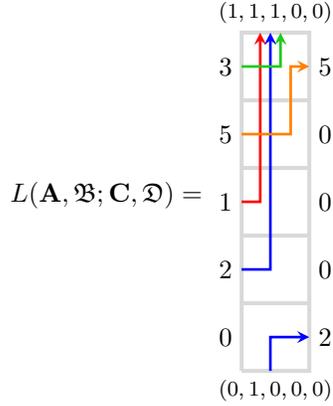
\tikz{.9}{
		\foreach\y in {1,...,4}{
			\draw[lgray,line width=1.5pt] (1.5,0.5+\y) -- (2.5,0.5+\y);
		}
		\draw[lgray,line width=1.5pt] (1.5,0.5) -- (2.5,0.5) -- (2.5,5.5) -- (1.5,5.5) -- (1.5,0.5);
		%paths
		\draw[red,line width=1pt,->] (1.5,3) -- (1.775,3) -- (1.775,5.5);
		\draw[blue,line width=1pt,->] (1.5,2) -- (1.925,2) -- (1.925,5.5);
		\draw[blue,line width=1pt,->] (1.925,0.5) -- (1.925,1) -- (2.5,1);
		\draw[green,line width=1pt,->] (1.5,5) -- (2.075,5) -- (2.075,5.5);
		\draw[orange,line width=1pt,->] (1.5,4) -- (2.225,4) -- (2.225,5) -- (2.5,5);
		%bottom labels
		\node[below] at (2,0.5) {\footnotesize$(0,1,0,0,0)$};
		%top labels
		\node[above] at (2,5.5) {\footnotesize$(1,1,1,0,0)$};
		%right labels
		\node[right] at (2.5,1) {$2$};
		\node[right] at (2.5,2) {$0$};
		\node[right] at (2.5,3) {$0$};
		\node[right] at (2.5,4) {$0$};
		\node[right] at (2.5,5) {$5$};
		%left labels
		\node[left] at (1.5,1) {$0$};
		\node[left] at (1.5,2) {$2$};
		\node[left] at (1.5,3) {$1$};
		\node[left] at (1.5,4) {$5$};
		\node[left] at (1.5,5) {$3$};
	}
	\end{align*}
	\noindent where we omit the spectral parameters $\textbf{x} = (x_1, x_2, \ldots , x_5)$ in both the weight $L = L_{\textbf{x}}$ and the column diagram. In this picture, $\mathfrak{B}$ and $\mathfrak{D}$ label the colors (read from bottom to top) of the arrows entering and exiting through the left and right boundaries of the column, respectively. The admissibility of $(\mathfrak{B},\mathfrak{D})$ translates into the fact that each color $\{1, 2, \dots,5\}$ appears at most once on the left or right of the tower, with the requirement that all colors that appear on the right must have also appeared on the left. Then $\mathcal{P} \cup \mathcal{Q}$ gives the set of all colors entering via the left edges of the tower, $\mathcal{Q}$ gives the set of all colors leaving via the right edges. Here red, blue, green, and orange are colors $1$, $2$, $3$, and $5$ respectively. 
\end{example}

Now the following proposition provides explicit forms for the matrix entries of column operators.

\begin{prop}
	\label{thm:comp}
	Let $\mathfrak{B} = (b_1, b_2, \dots, b_n)$, $\mathfrak{D} = (d_1 d_2, \dots,j_n)$ be sequences of indices in $[0, n]$ such that the pair $(\mathfrak{B}, \mathfrak{D})$ is admissible, as in \Cref{definition:admiss}. Associate to it the color data $(\mathcal{P},\mathcal{Q})$ in the same way as in \Cref{definition:data}. Let $(\mathcal{I},\mathcal{J})$ be the coordinates of $(\mathfrak{B}, \mathfrak{D})$, as in equations \eqref{coord-def}--\eqref{coord-def2}. Fix sequences $\textbf{\emph{x}} = (x_1, x_2, \ldots , x_n)$ and $\textbf{\emph{v}} = (v_1, v_2, \ldots , v_n)$ of $n$ complex numbers such that $v_r = 0$ for all $r \not\in \mathcal{P} \cup \mathcal{Q}$, and define the $\{ 0, 1 \}^n$ element
	\begin{align*}
	\textbf{\emph{e}}_{\mathcal{P}}
	=
	\sum_{p \in \mathcal{P}} 
	\textbf{\emph{e}}_p.
	\end{align*}
	Under the above set of assumptions, we have
	\begin{multline}
	\label{components}
	\big\langle \mathfrak{B} \big|
	\Psi_{\textbf{\emph{v}}; \textbf{\emph{x}}}
	\left( \textbf{\emph{e}}_{\mathcal{P}} \right)
	\big| \mathfrak{D} \big\rangle
	=
	\prod_{\substack{
			p > \ell
			\\
			p \in \mathcal{P} \cup \mathcal{Q}
			\\
			\ell \in \mathcal{Q}
	}}
	\bm{1}_{i_p \not= j_{\ell}}
	\prod_{p \in \mathcal{P}} t^{g(p)}
	\prod_{p \in \mathcal{Q}} x_{j_p}
	\prod_{\substack{
			p \in \mathcal{Q} \\ i_p=j_p
	}}
	(1-v_p t^{f(p)+1})
	\prod_{\substack{
			p\in \mathcal{Q} \\ i_p \not= j_p
	}}
	v_p^{\textbf{\emph{1}}_{i_p > j_p}} t^{h(p)} (1 - t),
	\end{multline}
	where we have defined the combinatorial exponents
	\begin{align}
	\label{exponents}
	\begin{split}
	f(p) 
	= 
	\# \{\ell \in & \mathcal{Q}: \ell < p\}; \qquad
	g(p)
	=
	\# \{\ell \in \mathcal{Q}: \ell < p, i_p < j_{\ell}\}; \\
	& h(p)
	=
	\# \big\{\ell \in \mathcal{Q}: \ell < p, j_{\ell} \in (i_p, j_p) \big\},
	\end{split}
	\end{align}
	with the interval $(i_p, j_p)$ appearing in $h(p)$ to be interpreted in a cyclic sense; namely, for all integers $1 \le a,b \le n$ we define
	\begin{align*}
	(a,b) 
	:=
	\left\{
	\begin{array}{ll}
	\{a+1,\dots,b-1\},
	&
	a < b,
	\\
	\{a+1,\dots,n\} \cup \{1,\dots,b-1\},
	&
	a > b,
	\\
	\varnothing,
	&
	a = b.
	\end{array}
	\right.
	\end{align*}
	
\end{prop} 

\begin{proof}
	We compute the components by calculating
	\begin{align}
	\label{comp-calc1}
	\big\langle \mathfrak{B} \big|
	\Psi_{\textbf{v}; \textbf{x}}
	( \textbf{e}_{\mathcal{P}} )
	\big| \mathfrak{D} \big\rangle
	=
	\sum_{\textbf{M} \in \{0,1\}^n}
	\prod_{k=1}^{n}
	v_k^{M_k}
	L_{\textbf{x}} (\textbf{M}, \mathfrak{B}; 
	\textbf{M}+\textbf{e}_{\mathcal{P}}, \mathfrak{D}),
	\end{align}
	where we set $\textbf{M} = (M_1, M_2, \ldots , M_n) \in \{ 0, 1 \}^n$, and the $L_{\textbf{x}}$ here is the weight of the column \eqref{tower} (with $\textbf{A} = \textbf{M}$ and $\textbf{C} = \textbf{M} + \textbf{e}_{\mathcal{P}}$). We begin by remarking that this column has weight zero if $b_k > d_k \ge 1$ for any $k \in [1, n]$, due to the vanishing of the fifth weight in \Cref{vertexfigurerql} (namely, $L_x (\textbf{A}, j; \textbf{A}_{ji}^{+-}, i) = 0$ for any $1 \le i < j \le n$). This gives rise to the product of indicator functions present in \eqref{components}. In the rest of the proof, we assume that the constraints imposed by these indicators are always obeyed.
	
	Noting that the vertex weights from \Cref{lz} are factorized across different colors, \ie\ over the $n$ components of $\textbf{A}$, we find that the $n$ sums in \eqref{comp-calc1}, over $M_1, M_2, \dots, M_n \ge 0$,  can be computed independently of each other (once $\mathfrak{B}$ and $\mathfrak{D}$ are fixed). This leads to the factorization
	\begin{align}
	\label{blambdad}
	\big\langle \mathfrak{B} \big|
	\Psi_{\textbf{v}; \textbf{x}}
	\left( \textbf{e}_{\mathcal{P}} \right)
	\big| \mathfrak{D} \big\rangle
	=
	\prod_{p \in \mathcal{P}}
	\phi_p(v_p)
	\prod_{\substack{p \in \mathcal{Q} \\ i_p = j_p}}
	\chi_p(v_p)
	\prod_{\substack{p \in \mathcal{Q} \\ i_p \not= j_p}}
	\psi_p(v_p),
	\end{align}
	where we have introduced the functions
	\begin{align}
	\label{phi-sum}
	\phi_p(v)
	=
	t^{g(p)}; \qquad
	\chi_p(v)
	=
	x_{j_p}
	(1-v t^{f(p)+1}); \qquad
	\psi_p(v) 
	=
	x_{j_p}
	v^{\bm{1}_{i_p > j_p}}
	t^{h(p)}(1-t),
	\end{align}
	\noindent and for notational convenience we have omitted the dependence of these functions on $\textbf{x}$.
	
	The quantities in \eqref{phi-sum} can be obtained by tracing paths of a fixed color; see \Cref{fig:phi-chi-psi}.
	\begin{figure}[t]
	\begin{center}
		\tikz{1}{
			\foreach\y in {1,...,4}{
				\draw[lgray,line width=1.5pt] (1.5,0.5+\y) -- (2.5,0.5+\y);
			}
			\draw[lgray,line width=1.5pt] (1.5,0.5) -- (2.5,0.5) -- (2.5,5.5) -- (1.5,5.5) -- (1.5,0.5);
			\draw[lgray,dashed,line width=2.2pt] (2.6,3.5) -- (2.6,5.5);
			%paths
			\draw[red,line width=1pt,->] (1.5,3) -- (2.225,3) -- (2.225,5.5);
			%\draw[red,line width=1pt,->] (1.925,0.5) -- (1.925,5.5);
			%\draw[red,line width=1pt,->] (2.075,0.5) -- (2.075,5.5);
			%
			\draw[blue,line width=1pt,->] (1.5,2) -- (1.775,2) -- (1.775,4) -- (2.5,4);
			%
			%\draw[green,line width=1pt,->] (1.5,1) -- (2.5,1);
			%bottom labels
			%\node[below] at (2,0.5) {\footnotesize$m$\ {\rm times}};
			%right labels
			\node[right] at (2.6,4) {$\ell_2$};
			%\node[right] at (2.6,1) {$\ell_1$};
			%left labels
			\node[left] at (1.5,3) {$p$};
		}
		\qquad\qquad\qquad
		\tikz{1}{
			\foreach\y in {1,...,4}{
				\draw[lgray,line width=1.5pt] (1.5,0.5+\y) -- (2.5,0.5+\y);
			}
			\draw[lgray,line width=1.5pt] (1.5,0.5) -- (2.5,0.5) -- (2.5,5.5) -- (1.5,5.5) -- (1.5,0.5);
			%paths
			\draw[red,line width=1pt,->] (1.5,3) -- (2.5,3);
			%\draw[red,line width=1pt,->] (1.925,0.5) -- (1.925,5.5);
			%\draw[red,line width=1pt,->] (2.075,0.5) -- (2.075,5.5);
			%
			\draw[green,line width=1pt,->] (1.5,1) -- (2.25,1) -- (2.25,2) -- (2.5,2);
			%bottom labels
			\node[below] at (2,0.5) {\footnotesize$0$\ {\rm times}};
			%right labels
			\node[right] at (2.6,2) {$\ell_1$};
			%left labels
			\node[left] at (1.5,3) {$p$};
		}
		\tikz{1}{
			\foreach\y in {1,...,4}{
				\draw[lgray,line width=1.5pt] (1.5,0.5+\y) -- (2.5,0.5+\y);
			}
			\draw[lgray,line width=1.5pt] (1.5,0.5) -- (2.5,0.5) -- (2.5,5.5) -- (1.5,5.5) -- (1.5,0.5);
			%paths
			\draw[red,line width=1pt,->] (1.5,3) -- (2.5,3);
			%\draw[red,line width=1pt,->] (1.925,0.5) -- (1.925,5.5);
			\draw[red,line width=1pt,->] (2.0,0.5) -- (2.0,5.5);
			\draw[green,line width=1pt,->] (1.5,1) -- (2.25,1) -- (2.25,2) -- (2.5,2);
			%bottom labels
			\node[below] at (2,0.5) {\footnotesize$1$\ {\rm time}};
			%right labels
			\node[right] at (2.6,2) {$\ell_1$};
			%left labels
			\node[left] at (1.5,3) {$p$};
		}
		\qquad\qquad\qquad
		\tikz{1}{
			\foreach\y in {1,...,4}{
				\draw[lgray,line width=1.5pt] (1.5,0.5+\y) -- (2.5,0.5+\y);
			}
			\draw[lgray,line width=1.5pt] (1.5,0.5) -- (2.5,0.5) -- (2.5,5.5) -- (1.5,5.5) -- (1.5,0.5);
			\draw[lgray,dashed,line width=2.2pt] (2.6,2.5) -- (2.6,4.5);
			%paths
			\draw[red,line width=1pt,->] (1.5,2) -- (2.225,2) -- (2.225,5) -- (2.5,5);
			%\draw[red,line width=1pt,->] (1.925,0.5) -- (1.925,5.5);
			%\draw[red,line width=1pt,->] (2.075,0.5) -- (2.075,5.5);
			%
			\draw[blue,line width=1pt,->] (1.5,3) -- (1.775,3) -- (1.775,4) -- (2.5,4);
			%
			%\draw[green,line width=1pt,->] (1.5,1) -- (2.5,1);
			%bottom labels
			%\node[below] at (2,0.5) {\footnotesize$m$\ {\rm times}};
			%right labels
			\node[right] at (2.6,4) {$\ell_2$};
			%\node[right] at (2.6,1) {$\ell_1$};
			%left labels
			\node[left] at (1.5,2) {$p$};
		}
		\caption{Left panel: a color $p$ path enters via the left boundary, and exits via the top. $g(p)$ counts the number of colors $\ell_2 < p$ that exit through one of the shaded right edges, which are situated strictly above the position where $p$ entered.
			\\
			{\color{white} .}
			\hspace{0.1em} Middle panel: a color $p$ path enters via the left boundary, and circles $0$ or $1$ times around it before exiting via the right edge directly opposite its entry point. $f(p)$ counts the number of colors $\ell_1 < p$ which enter from the left of the column and leave anywhere on its right.
			\\
			{\color{white} .}
			\hspace{0.1em} Right panel: a color $p$ path of enters via the left boundary, and then exits via a right edge in a different row to its starting point. $h(p)$ counts the number of colors $\ell_2 < p$ that exit through one of the shaded right edges, which are situated strictly between the initial and final locations of path $p$ (the dashed interval on the right boundary).}
		\label{fig:phi-chi-psi}
	
	\end{center}
	\end{figure}
	Pictorially, the first quantity $\phi_p (v)$ in \eqref{phi-sum} corresponds to a color $p$ path that entered through the left boundary of the column and then exited through its top. The second quantity $\chi_p (v)$ in \eqref{phi-sum} corresponds to a color $p$ path that entered through the left boundary, circled $0$ or $1$ times over the column, and then exited on its right through the {\it same} row as where it entered. The third quantity $\psi_p (v)$ in \eqref{phi-sum} corresponds to a color $p$ path that entered through the left boundary of the column, and then exited on its right but via a {\it different} row from the one where it entered (after wrapping around the column, if required).
	
	The result \eqref{components} now follows by inserting \eqref{phi-sum} into \eqref{blambdad}.
\end{proof}

\section{Column Rotation}
\label{ssec:rotate}

In this section we define a specific ``rotation'' operation which acts on the column operator components. The key result will be the fact that, up to a simple overall factor depending on $t$ and $\textbf{v} = (v_1, v_2, \dots,v_n)$, the components \eqref{components} remain invariant under this operation. 
\begin{prop}
	\label{prop:rotate}
	Assuming the setup in \Cref{thm:comp}, consider the effect of rotating the edge states of the column, which is equivalent to sending $(b_1, \dots,b_{n-1}, b_n) \mapsto (b_n, b_1,\dots,b_{n-1})$ and $(d_1,\dots, d_{n-1},d_n) \mapsto (d_n,d_1,\dots,d_{n-1})$ in \eqref{components}, together with the corresponding cyclic shift of the row parameters $(x_1,\dots, x_{n-1}, x_n) \mapsto (x_n,x_1,\dots, x_{n-1})$. We find that
	\begin{align}
	\label{rotation-const}
	\frac{
		\bra{b_1,\dots,b_n}
		\Psi_{\textbf{\emph{v}}; \textbf{\emph{x}}} (\textbf{\emph{e}}_{\mathcal{P}})
		\ket{d_1,\dots, d_n}
	}
	{
		\bra{b_n,b_1,\dots,b_{n-1}}
		\mathfrak{s}_{n-1} \cdots \mathfrak{s}_1 \cdot
		\Psi_{\textbf{\emph{v}}; \textbf{\emph{x}}}
		\left( \textbf{\emph{e}}_{\mathcal{P}} \right)
		\ket{d_n,d_1,\dots, d_{n-1}}
	}
	=
	\kappa_{\textbf{\emph{v}}; t}(\mathcal{I},\mathcal{J}),
	\end{align}
	where the right-hand side is given by
	\begin{align}
	\label{rotation-const2}
	\kappa_{\textbf{\emph{v}}; t} (\mathcal{I},\mathcal{J})
	=
	\prod_{\substack{(p,\ell) \in \mathcal{P} \times \mathcal{Q} \\ p > \ell}}
	t^{\bm{1}_{j_{\ell}=n}-\bm{1}_{i_p = n}}
	\prod_{p \in \mathcal{Q}}
	v_p^{\bm{1}_{i_p=n}-\bm{1}_{j_p = n}}.
	\end{align}
\end{prop}

\begin{proof}
	We can explicitly compute the ratio on the left-hand side of \eqref{rotation-const} by using the formula \eqref{components} for the column components. The main observation is that, under the proposed rotation, almost all of the factors present in \eqref{components} remain invariant.
	
	Let us recall that $\mathcal{P}\cup\mathcal{Q}$ is the set of colors entering on the $n$ left edges of the column, while $\mathcal{Q}$ is the set of colors exiting via the $n$ right edges of the column. The vector $(i_p)_{p \in \mathcal{P} \cup \mathcal{Q}}$ records the positions where colors traverse left edges, while $(j_p)_{p \in \mathcal{Q}}$ records the positions where colors traverse right edges.
	
	First, we examine the term $\prod_{p \in \mathcal{Q}} x_{j_p}$ in \eqref{components}, where the product ranges over all colors which exit via the right edges of the column. The required cyclic rotation of the exiting positions is achieved by the replacement $j_p \mapsto j_p+1\ (\text{mod}\ n)$, which is clearly negated by the cyclic permutation $\mathfrak{s}_{n-1} \cdots \mathfrak{s}_1$, which sends $x_i \mapsto x_{i-1\ (\text{mod}\ n)}$. Therefore this term cancels out in the ratio on the left-hand side of \eqref{rotation-const}. 
	
	Similarly, it is clear from their definition that both of the exponents $f(p)$ and $h(p)$ in \eqref{exponents} remain unchanged under the rotation $i_p \mapsto i_p+1\ (\text{mod}\ n)$ for $p \in \mathcal{P} \cup \mathcal{Q}$ and $j_p \mapsto j_p+1\ (\text{mod}\ n)$ for $p \in \mathcal{Q}$. This allows us to cancel all factors involving those exponents in the ratio \eqref{rotation-const}, and we then read
	\begin{align}
	\label{rotation-const3}
	\frac{
		\bra{b_1,\dots,b_n}
		\Psi_{\textbf{v}; \textbf{x}}(\textbf{e}_{\mathcal{P}})
		\ket{d_1,\dots,d_n}
	}
	{
		\bra{b_n,b_1,\dots,b_{n-1}}
		\mathfrak{s}_{n-1} \cdots \mathfrak{s}_1 \cdot
		\Psi_{\textbf{v}; \textbf{x}}
		\left( \textbf{e}_{\mathcal{P}} \right)
		\ket{d_n,d_1,\dots,d_{n-1}}
	}
	=
	\frac{
		\prod_{p \in \mathcal{P}} t^{g(p)} 
		\prod_{p \in \mathcal{Q}} v_p^{\bm{1}_{i_p > j_p}}
	}
	{
		\prod_{p \in \mathcal{P}} t^{\tilde{g}(p)}
		\prod_{p \in \mathcal{Q}} v_p^{\bm{1}_{\tilde{i}_p > \tilde{j}_p}}
	},
	\end{align}
	where we have defined $\tilde{i}_p = i_p+1\ (\text{mod}\ n)$, 
	$\tilde{j}_p = j_p+1\ (\text{mod}\ n)$ and 
	\begin{align}
	\tilde{g}(p)
	=
	\#\{\ell \in \mathcal{Q} : \ell < p, \tilde{i}_p < \tilde{j}_{\ell}\}.
	\end{align}
	We are able to simplify the right-hand side of \eqref{rotation-const3} yet further, by noticing that it is only in the case $i_p = n$ for some $p \in \mathcal{P} \cup \mathcal{Q}$, or $j_p = n$ for some $p \in \mathcal{Q}$, where we see a discrepancy between the numerator and denominator. In those cases, one has $\tilde{i}_p=1$ or $\tilde{j}_p=1$, which can cause inequalities of the form $i_p > j_{\ell}$ or $i_{\ell} < j_p$ that previously held to now be violated. Analyzing these cases yields the final result, given by 
	\begin{multline*}
	\frac{
		\bra{b_1,\dots,b_n}
		\Psi_{\textbf{v}; \textbf{x}}(\textbf{e}_{\mathcal{P}})
		\ket{d_1,\dots,d_n}
	}
	{
		\bra{b_n,b_1,\dots,b_{n-1}}
		\mathfrak{s}_{n-1} \cdots \mathfrak{s}_1 \cdot
		\Psi_{\textbf{v}; \textbf{x}}
		\left( \textbf{e}_{\mathcal{P}} \right)
		\ket{d_n,d_1,\dots,d_{n-1}}
	}
	\\
	=
	\prod_{\substack{(p,\ell) \in \mathcal{P} \times \mathcal{Q} \\ p > \ell}}
	t^{\bm{1}_{j_{\ell}=n}-\bm{1}_{i_p = n}}
	\prod_{p \in \mathcal{Q}}
	v_p^{\bm{1}_{n = i_p > j_p}-\bm{1}_{i_p < j_p = n}},
	\end{multline*}
	where we can simplify the exponent in the second product by noting that $\bm{1}_{n=i_p > j_p} -\bm{1}_{i_p < j_p = n} = \bm{1}_{i_p=n} - \bm{1}_{j_p=n}$. 	
\end{proof}

\begin{rem}
	In the calculations that follow, it is more useful to recast the right-hand side of \eqref{rotation-const} in terms of the color data $\mathcal{P}$, $\mathcal{Q}$ and the entries $b_n$, $d_n$, without making reference to the coordinates of \Cref{definition:coord}. We find that 
	\begin{align}
	\label{rot-const}
	\kappa_{\textbf{v}; t}
	=
	\frac{
		t^{\#\{ p \in \mathcal{P} : p > d_n \} \bm{1}_{d_n \ge 1}}
	}
	{
		t^{\#\{ p \in \mathcal{Q} : b_n > p\} \bm{1}_{b_n \in \mathcal{P}}}
	}
	\cdot
	\frac{
		(v_{b_n})^{\bm{1}_{b_n \in \mathcal{Q}}}
	}
	{
		(v_{d_n})^{\bm{1}_{d_n \ge 1}}
	},
	\end{align}
	with the same $\mathcal{P}$ and $\mathcal{Q}$ as in \Cref{definition:data}.
\end{rem}

\section{Proof of \Cref{omegaif}} 

\label{ssec:proof}

We are now ready to return to the proof of \Cref{omegaif}. Fix a composition $\mu$ of length $n$, and set $N = \max_{i \in [1, n]} \mu_i$. For each integer $a \in [1, N]$, further fix a sequence $\mathfrak{K}^{(a)} = \big( k^{(a)}_1,\dots,k^{(a)}_n \big)$ of indices in $[0, n]$. We introduce a function $Z_{\rm l}$ obtained by concatenating $N$ column operators, as
\begin{align}
\label{Z_l}
Z_{\rm l}\Big[ \mathfrak{K}^{(1)}, \mathfrak{K}^{(2)}, \ldots , \mathfrak{K}^{(N)} \Big]
=
\tikz{.9}{
	\foreach\y in {0}{
		\draw[lgray,line width=1.7pt,double] (1.5,0.5+\y) -- (7.5,0.5+\y);
	}
	\foreach\y in {5}{
		\draw[lgray,line width=1.7pt,double] (1.5,0.5+\y) -- (7.5,0.5+\y);
	}
	\foreach\y in {1,...,4}{
		\draw[lgray,line width=1.5pt] (1.5,0.5+\y) -- (7.5,0.5+\y);
	}
	\foreach\x in {0,...,6}{
		\draw[lgray,line width=1.5pt] (1.5+\x,0.5) -- (1.5+\x,5.5);
	}
	%spectral parameters
	\node[text centered] at (2,1) {$x_1$};
	\node[text centered] at (4,1) {$\cdots$};
	\node[text centered] at (5,1) {$\cdots$};
	\node[text centered] at (7,1) {$x_1$};
	\node[text centered] at (2,4) {$x_n$};
	\node[text centered] at (4,4) {$\cdots$};
	\node[text centered] at (5,4) {$\cdots$};
	\node[text centered] at (7,4) {$x_n$};
	\node[text centered] at (2,5) {$qx_i$};
	\node[text centered] at (4,5) {$\cdots$};
	\node[text centered] at (5,5) {$\cdots$};
	\node[text centered] at (7,5) {$qx_i$};
	\node[text centered] at (2,2.1) {$\vdots$};
	\node[text centered] at (2,3.1) {$\vdots$};
	\node[text centered] at (7,2.1) {$\vdots$};
	\node[text centered] at (7,3.1) {$\vdots$};
	%top labels
	\node[above] at (2,5.5) {\footnotesize$\textbf{I}_0 (\mu)$};
	\node[above,text centered] at (4,5.5) {$\cdots$};
	\node[above,text centered] at (5,5.5) {$\cdots$};
	\node[above] at (7,5.5) {\footnotesize$\textbf{I}_N (\mu)$};
	%right labels
	\node[right] at (7.5,1) {$0$};
	\node[right] at (7.5,2.1) {$\vdots$};
	\node[right] at (7.5,3.1) {$\vdots$};
	\node[right] at (7.5,4) {$0$};
	\node[right] at (7.5,5) {$0$};
	%left labels
	\node[left] at (1.5,1) {$1$};
	\node[left] at (1.5,2.1) {$\vdots$};
	\node[left] at (1.5,3.1) {$\vdots$};
	\node[left] at (1.5,4) {$n$};
	\node[left] at (1.5,5) {$i$};
	%internal labels
	\node[right, text centered] at (2.5,5) {\footnotesize $k^{(1)}_i$};
	\node at (2.5,5) {$\bullet$};
	\node[left, text centered] at (6.5,5) {\footnotesize $k^{(N)}_i$};
	\node at (6.5,5) {$\bullet$};
	\node[right, text centered] at (2.5,4) {\footnotesize $k^{(1)}_n$};
	\node at (2.5,4) {$\bullet$};
	\node[left, text centered] at (6.5,4) {\footnotesize $k^{(N)}_n$};
	\node at (6.5,4) {$\bullet$};
	\node[right, text centered] at (2.5,1) {\footnotesize $k^{(1)}_1$};
	\node at (2.5,1) {$\bullet$};
	\node[left, text centered] at (6.5,1) {\footnotesize $k^{(N)}_1$};
	\node at (6.5,1) {$\bullet$};
}
\end{align}
where the indices in this picture are as follows. The colors of the arrows horizontally entering the rows from bottom to top are $1, 2, \ldots , i - 1, i + 1, i + 2, \ldots , n, i$. In what follows, we will label each row according to the color of the arrow entering through its left edge; accordingly, the top row will continue to be called the $i$-th row. For each $b \in [1, n]$, the internal arrow colors in the $b$-th row are indexed from left to right by $\big( k^{(1)}_b,\dots,k^{(N)}_b \big)$. These colors are fixed, not summed over.

We introduce a similar function $Z_{\rm r}$, given by 
\begin{align}
\label{Z_r}
Z_{\rm r} \Big[ \mathfrak{K}^{(1)}, \mathfrak{K}^{(2)}, \ldots , \mathfrak{K}^{(N)} \Big]
=
\tikz{.9}{
	\foreach\y in {0}{
		\draw[lgray,line width=1.7pt,double] (1.5,0.5+\y) -- (7.5,0.5+\y);
	}
	\foreach\y in {5}{
		\draw[lgray,line width=1.7pt,double] (1.5,0.5+\y) -- (7.5,0.5+\y);
	}
	\foreach\y in {1,...,4}{
		\draw[lgray,line width=1.5pt] (1.5,0.5+\y) -- (7.5,0.5+\y);
	}
	\foreach\x in {0,...,6}{
		\draw[lgray,line width=1.5pt] (1.5+\x,0.5) -- (1.5+\x,5.5);
	}
	%spectral parameters
	\node[text centered] at (2,1) {$x_i$};
	\node[text centered] at (4,1) {$\cdots$};
	\node[text centered] at (5,1) {$\cdots$};
	\node[text centered] at (7,1) {$x_i$};
	\node[text centered] at (2,2) {$x_1$};
	\node[text centered] at (4,2) {$\cdots$};
	\node[text centered] at (5,2) {$\cdots$};
	\node[text centered] at (7,2) {$x_1$};
	\node[text centered] at (2,5) {$x_n$};
	\node[text centered] at (4,5) {$\cdots$};
	\node[text centered] at (5,5) {$\cdots$};
	\node[text centered] at (7,5) {$x_n$};
	\node[text centered] at (2,3.1) {$\vdots$};
	\node[text centered] at (2,4.1) {$\vdots$};
	\node[text centered] at (7,3.1) {$\vdots$};
	\node[text centered] at (7,4.1) {$\vdots$};
	%top labels
	\node[above] at (2,5.5) {\footnotesize$\textbf{I}_0 (\mu)$};
	\node[above,text centered] at (4,5.5) {$\cdots$};
	\node[above,text centered] at (5,5.5) {$\cdots$};
	\node[above] at (7,5.5) {\footnotesize$\textbf{I}_N (\mu)$};
	%right labels
	\node[right] at (7.5,1) {$0$};
	\node[right] at (7.5,2) {$0$};
	\node[right] at (7.5,3.1) {$\vdots$};
	\node[right] at (7.5,4.1) {$\vdots$};
	\node[right] at (7.5,5) {$0$};
	%left labels
	\node[left] at (1.5,1) {$i$};
	\node[left] at (1.5,2) {$1$};
	\node[left] at (1.5,3.1) {$\vdots$};
	\node[left] at (1.5,4.1) {$\vdots$};
	\node[left] at (1.5,5) {$n$};
	%internal labels
	\node[right, text centered] at (2.5,5) {\footnotesize $k^{(1)}_n$};
	\node at (2.5,5) {$\bullet$};
	\node[left, text centered] at (6.5,5) {\footnotesize $k^{(N)}_n$};
	\node at (6.5,5) {$\bullet$};
	\node[right, text centered] at (2.5,2) {\footnotesize $k^{(1)}_1$};
	\node at (2.5,2) {$\bullet$};
	\node[left, text centered] at (6.5,2) {\footnotesize $k^{(N)}_1$};
	\node at (6.5,2) {$\bullet$};
	\node[right, text centered] at (2.5,1) {\footnotesize $k^{(1)}_i$};
	\node at (2.5,1) {$\bullet$};
	\node[left, text centered] at (6.5,1) {\footnotesize $k^{(N)}_i$};
	\node at (6.5,1) {$\bullet$};
}
\end{align}
\noindent where again the left arrow colors are ordered sequentially from bottom to top, with the exception of state $i$; it is omitted from $(1,\dots,n)$ and transferred to the bottom row. This time, the bottom row is called the $i$-th row. As before, we assign the colors $\big(k_b^{(1)},\dots, k_b^{(N)} \big)$ to the internal arrows in the $b$-th row, for each $b \in [1, n]$. 

We note that $Z_{\rm l}$ and $Z_{\rm r}$ are equivalent under application of the rotation operation of \Cref{ssec:rotate} to each of the $N+1$ columns that are present (up to the variable shift $qx_i \mapsto x_i$ in the $i$-th row). $Z_{\rm l}$ and $Z_{\rm r}$ also serve as refinements of the left and right-hand sides of \eqref{cyclic-rel}, respectively. One sees that
\begin{align}
\label{left-refined}
\big\langle \mathsf{C}_i (qx_i) \mathsf{C}_n(x_n) \cdots \mathsf{C}_{i + 1} (x_{i + 1}) \mathsf{C}_{i - 1} (x_{i - 1}) \cdots \mathsf{C}_1 (x_1) \big\rangle_{\mu}
&=
\sum
Z_{\rm l}\Big[\mathfrak{K}^{(1)}, \mathfrak{K}^{(2)}, \ldots , \mathfrak{K}^{(N)} \Big],
\\
\label{right-refined}
\big\langle \mathsf{C}_n (x_n) \cdots \mathsf{C}_{i +1} (x_{i + 1}) \mathsf{C}_{i - 1} (x_{i - 1}) \cdots \mathsf{C}_1(x_1) \mathsf{C}_i (x_i) \big\rangle_{\mu}
&=
\sum
Z_{\rm r}\Big[ \mathfrak{K}^{(1)}, \mathfrak{K}^{(2)}, \ldots , \mathfrak{K}^{(N)} \Big],
\end{align}

\noindent where both sums are over all length $n$ sequences $\big\{ \mathfrak{K}^{(a)} \big\}_{a \in [1, N]}$ of indices in $[0, n]$.

\begin{example} Let us give an explicit example of the correspondence between \eqref{left-refined} and \eqref{right-refined}, in the case $n=5$, $i=3$ and $\mu =(0,4,4,1,5)$. In this case we have 
	\begin{align*}
	\textbf{I}_0 (\mu) = \textbf{e}_1,\qquad 
	\textbf{I}_1 (\mu) = \textbf{e}_4,\qquad  
	\textbf{I}_2 (\mu) = \textbf{e}_0 = \textbf{I}_3 (\mu), \qquad
	\textbf{I}_4 (\mu) = \textbf{e}_2 + \textbf{e}_3, \qquad 
	\textbf{I}_5 (\mu) = \textbf{e}_5.
	\end{align*}
	We extract a sample term from the sums \eqref{left-refined} and 
	\eqref{right-refined} given by 
	\begin{align*}
	\tikz{.9}{
		\foreach\y in {0}{
			\draw[lgray,line width=1.5pt,double] (1.5,0.5+\y) -- (7.5,0.5+\y);
		}
		\foreach\y in {5}{
			\draw[lgray,line width=1.5pt,double] (1.5,0.5+\y) -- (7.5,0.5+\y);
		}
		\foreach\y in {1,...,4}{
			\draw[lgray,line width=1.5pt] (1.5,0.5+\y) -- (7.5,0.5+\y);
		}
		\foreach\x in {0,...,6}{
			\draw[lgray,line width=1.5pt] (1.5+\x,0.5) -- (1.5+\x,5.5);
		}
		%top labels
		\node[above] at (2,5.5) {\footnotesize$\textbf{e}_1$};
		\node[above] at (3,5.5) {\footnotesize$\textbf{e}_4$};
		\node[above] at (4,5.5) {\footnotesize$\textbf{e}_0$};
		\node[above] at (5,5.5) {\footnotesize$\textbf{e}_0$};
		\node[above] at (6,5.5) {\footnotesize$\textbf{e}_2+\textbf{e}_3$};
		\node[above] at (7,5.5) {\footnotesize$\textbf{e}_5$};
		%right labels
		\node[right] at (7.5,1) {$0$};
		\node[right] at (7.5,2) {$0$};
		\node[right] at (7.5,3.1) {$0$};
		\node[right] at (7.5,4.1) {$0$};
		\node[right] at (7.5,5) {$0$};
		%spectral parameters
		\node[left] at (1,1) {$x_1$};
		\node[left] at (1,2) {$x_2$};
		\node[left] at (1,3) {$x_4$};
		\node[left] at (1,4) {$x_5$};
		\node[left] at (1,5) {$qx_3$};
		%left labels
		\node[left] at (1.5,1) {$1$};
		\node[left] at (1.5,2) {$2$};
		\node[left] at (1.5,3) {$4$};
		\node[left] at (1.5,4) {$5$};
		\node[left] at (1.5,5) {$3$};
		%paths
		\draw[ultra thick,yellow,->] (1.5,5) -- (2.925,5) -- (2.925,5.5);
		\draw[ultra thick,yellow,->] (2.925,0.5) -- (2.925,2) -- (4,2) -- (4,3) -- (5,3) -- (5,5.5);
		\draw[ultra thick,yellow,->] (5,0.5) -- (5,1) -- (6.075,1) -- (6.075,5.5);
		\draw[ultra thick,orange,->] (1.5,4) -- (4,4) -- (4,5) -- (7,5) -- (7,5.5);
		\draw[ultra thick,green,->] (1.5,3) -- (3.075,3) -- (3.075,5.5);
		\draw[ultra thick,blue,->] (1.5,2) -- (2.075,2) -- (2.075,5.5); 
		\draw[ultra thick,blue,->] (2.075,0.5) -- (2.075,1) -- (5,1) -- (5,2) -- (5.925,2) -- (5.925,5.5);
		\draw[ultra thick,red,->] (1.5,1) -- (1.925,1) -- (1.925,5.5);
	}
	\quad\quad
	\tikz{.9}{
		\foreach\y in {0}{
			\draw[lgray,line width=1.5pt,double] (1.5,0.5+\y) -- (7.5,0.5+\y);
		}
		\foreach\y in {5}{
			\draw[lgray,line width=1.5pt,double] (1.5,0.5+\y) -- (7.5,0.5+\y);
		}
		\foreach\y in {1,...,4}{
			\draw[lgray,line width=1.5pt] (1.5,0.5+\y) -- (7.5,0.5+\y);
		}
		\foreach\x in {0,...,6}{
			\draw[lgray,line width=1.5pt] (1.5+\x,0.5) -- (1.5+\x,5.5);
		}
		%top labels
		\node[above] at (2,5.5) {\footnotesize$\textbf{e}_1$};
		\node[above] at (3,5.5) {\footnotesize$\textbf{e}_4$};
		\node[above] at (4,5.5) {\footnotesize$\textbf{e}_0$};
		\node[above] at (5,5.5) {\footnotesize$\textbf{e}_0$};
		\node[above] at (6,5.5) {\footnotesize$\textbf{e}_2+\textbf{e}_3$};
		\node[above] at (7,5.5) {\footnotesize$\textbf{e}_5$};
		%right labels
		\node[right] at (7.5,1) {$0$};
		\node[right] at (7.5,2) {$0$};
		\node[right] at (7.5,3.1) {$0$};
		\node[right] at (7.5,4.1) {$0$};
		\node[right] at (7.5,5) {$0$};
		%spectral parameters
		\node[left] at (1,1) {$x_3$};
		\node[left] at (1,2) {$x_1$};
		\node[left] at (1,3) {$x_2$};
		\node[left] at (1,4) {$x_4$};
		\node[left] at (1,5) {$x_5$};
		%left labels
		\node[left] at (1.5,1) {$3$};
		\node[left] at (1.5,2) {$1$};
		\node[left] at (1.5,3) {$2$};
		\node[left] at (1.5,4) {$4$};
		\node[left] at (1.5,5) {$5$};
		%paths
		\draw[ultra thick,yellow,->] (1.5,1) -- (2.925,1) -- (2.925,3) -- (4,3) -- (4,4) -- (5,4) -- (5,5.5);
		\draw[ultra thick,yellow,->] (5,0.5) -- (5,2) -- (6.075,2) -- (6.075,5.5);
		\draw[ultra thick, orange,->] (1.5,5) -- (4,5) -- (4,5.5); 
		\draw[ultra thick,orange,->] (4,0.5)-- (4,1) -- (7,1) -- (7,5.5);
		\draw[ultra thick,green,->] (1.5,4) -- (3.075,4) -- (3.075,5.5);
		\draw[ultra thick,blue,->] (1.5,3) -- (2.075,3) -- (2.075,5.5); 
		\draw[ultra thick,blue,->] (2.075,0.5) -- (2.075,2) -- (5,2) -- (5,3) -- (5.925,3) -- (5.925,5.5);
		\draw[ultra thick, red, ->] (1.5,2) -- (1.925,2) -- (1.925,5.5);
	}
	\end{align*}
	where the configuration on the left is one term in the sum \eqref{left-refined}, while the configuration on the right is one term in the sum \eqref{right-refined}. The configuration on the right is obtained from the configuration on the left by rotating each of the $N+1$ columns; namely, by bumping each row upwards by one step, and transferring the top row to the bottom of the lattice (and then dividing the spectral parameter in the new bottommost row by $q$).
\end{example}

\begin{prop}
	\label{prop:refined-cyclic-rel}
	
	For each integer $a \in [1, N]$, let $\mathfrak{K}^{(a)} = \big( k_1^{(a)}, k_2^{(a)}, \ldots , k_n^{(a)} \big)$ denote a sequence of indices in $[0, n]$, such that $Z_{\rm r} \big[ \mathfrak{K}^{(1)}, \mathfrak{K}^{(2)}, \ldots , \mathfrak{K}^{(N)} \big] \ne 0$. Then $Z_{\rm l} \big[ \mathfrak{K}^{(1)}, \mathfrak{K}^{(2)}, \ldots , \mathfrak{K}^{(N)} \big]$ is also non-vanishing, and we have
	\begin{align}
	\label{refined-cyclic-rel}
	\frac{Z_{\rm l}\Big[ \mathfrak{K}^{(1)}, \mathfrak{K}^{(2)}, \ldots , \mathfrak{K}^{(N)} \Big]}
	{Z_{\rm r}\Big[ \mathfrak{K}^{(1)}, \mathfrak{K}^{(2)}, \ldots , \mathfrak{K}^{(N)} \Big]}
	=
	q^{\mu_i} t^{\gamma_{i,0}(\mu)},
	\end{align}
	
	\noindent where the underlying parameters $v_{i, j}$ are given by \eqref{v_ij}.
\end{prop}

\begin{proof}
	The proof is based on $N+1$ applications of \eqref{rotation-const} and \eqref{rot-const}, noting that any of the $N+1$ columns present in $Z_{\rm l}$ differs from the corresponding column in $Z_{\rm r}$ by the rotation operation of \Cref{prop:rotate}. The ratio $Z_{\rm l} / Z_{\rm r}$ will thus be of the form $q^a t^b$ for some $a,b \in \mathbb{Z}$; let us proceed to calculate these exponents.
	
	For the purpose of calculating the right-hand side of \eqref{refined-cyclic-rel}, \Cref{prop:rotate} tells us that it is sufficient to focus on the internal edge states of the $i$-th row; in particular we will be interested in those $k^{(j)}_i$ which take non-zero values, since these contribute non-trivially to the right-hand side of \eqref{rotation-const}. As we will not need to specify $k^{(j)}_b$ for $b \not= i$, we hereafter lighten the notation by writing $k^{(j)}_i = K_j$ for all $1 \le j \le N$. We also set $K_0 = i$ and $K_{N+1} = 0$.
	
	By iterating \eqref{rot-const} over the $N+1$ columns and assigning a factor of $q$ to each integer $1 \le j \le N$ for which $K_j \ge 1$,\footnote{We obtain a power of $q$ for every horizontal step by a path in the $i$-th row of $Z_{\rm l}$, due to the $q$-shifted argument of $\mathcal{C}_i(qx_i)$ and \Cref{rem:x-depend}.} we find that
	\begin{align}
	\label{N-columns}
	\frac{Z_{\rm l}}{Z_{\rm r}}
	=
	\prod_{j=0}^{N}
	\frac{
		(v_{K_j,j})^{\bm{1}_{K_j \in \mathcal{Q}_j}}
	}
	{
		(v_{K_{j+1},j})^{\bm{1}_{K_{j+1} \ge 1}}
	}
	\cdot
	\frac{
		t^{\#\{a \in \mathcal{P}_j: a>K_{j+1}\} \bm{1}_{K_{j+1} \ge 1}}
	}
	{
		t^{\#\{a \in \mathcal{Q}_j: a<K_{j}\} \bm{1}_{K_j \in \mathcal{P}_j}}
	}
	\cdot
	\prod_{j=1}^{N}
	q^{\bm{1}_{K_j \ge 1}},
	\end{align}
	where $\mathcal{P}_j$, $\mathcal{Q}_j$ are the color data associated to the $j$-th column of $Z_{\rm l}$. Since $\mathcal{P}_j$ is the set of colors which exit column $j$ via its top edge and $\mathcal{Q}_j$ is the set of colors exiting via columns $j+1,\dots,N$, we can write
	\begin{align}
	\label{PQ_j}
	\mathcal{P}_j = \{a: \mu_a = j\},
	\qquad
	\mathcal{Q}_j = \{a: \mu_a > j\}.
	\end{align}
	The right-hand side of \eqref{N-columns} becomes
	\begin{align}
	\label{N-columns-2}
	\frac{Z_{\rm l}}{Z_{\rm r}}
	&=
	\frac{
		(v_{i,0})^{\bm{1}_{i \in \mathcal{Q}_0}}
	}
	{t^{\#\{ a<i : \mu_a > 0 \} \bm{1}_{i \in \mathcal{P}_0}}}
	\prod_{j=1}^{N}
	\frac{
		(v_{K_j,j})^{\bm{1}_{K_j \in \mathcal{Q}_j}}
	}
	{
		(v_{K_j,j-1})^{\bm{1}_{K_j \ge 1}}
	}
	\cdot
	\frac{
		t^{\#\{a>K_{j} : \mu_a=j-1\} \bm{1}_{K_{j} \ge 1}}
	}
	{
		t^{\#\{a<K_{j} : \mu_a > j\} \bm{1}_{K_j \in \mathcal{P}_j}}
	}
	\cdot
	q^{\bm{1}_{K_j \ge 1}},
	\end{align}
	where we have used the fact that $K_0 = i$ and $K_{N+1} = 0$ to redistribute the factors in the product.  
	
	Now let us compute the $j$-th term in the product on the right-hand side of \eqref{N-columns-2}. It is clearly sufficient to restrict our attention to $j$ such that $K_j \ge 1$, since for each $j$ such that $K_j = 0$ the factors inside the product are all equal to $1$; we tacitly assume $K_j \ge 1$ in what follows. Invoking (for the first time) the explicit form of the parameters \eqref{v_ij}, we find that
	\begin{align}
	\label{j-ratio}
	\frac{
		(v_{K_j,j})^{\bm{1}_{K_j \in \mathcal{Q}_j}}
	}
	{
		(v_{K_j,j-1})^{\bm{1}_{K_j \ge 1}}
	}
	=
	\left\{
	\begin{array}{ll}
	q^{-1} t^{\gamma_{K_j,j}(\mu)-\gamma_{K_j,j-1}(\mu)}, \quad
	& \mu_{K_j} > j,
	\\
	\\
	q^{-1} t^{-\gamma_{K_j,j-1}(\mu)}, \quad
	& \mu_{K_j} = j,
	\end{array}
	\right.
	\end{align}
	where the case $\mu_{K_j} < j$ never appears (it would mean that the color $K_j$ has traversed beyond column $\mu_{K_j}$, which is forbidden). From \eqref{gamma}, we can write down the $t$ exponents appearing in \eqref{j-ratio}. In the case $\mu_{K_j} > j$ (equivalent to $K_j \in \mathcal{Q}_j$), one has
	\begin{align*}
	\gamma_{K_j,j}(\mu)
	-
	\gamma_{K_j,j-1}(\mu)
	&=
	\#\{a > K_j : j \le \mu_a < \mu_{K_j} \}
	-
	\#\{a > K_j : j-1 \le \mu_a < \mu_{K_j} \}
	\\
	&=
	-\#\{a > K_j : j-1 = \mu_a \}.
	\end{align*}
	Similarly, for $\mu_{K_j} = j$ (equivalent to $K_j \in \mathcal{P}_j$), we find that
	\begin{align*}
	-\gamma_{K_j,j-1}(\mu)
	&=
	\#\{a < K_j : \mu_a > \mu_{K_j} \}
	-
	\#\{a > K_j : j-1 \le \mu_a < \mu_{K_j} \}
	\\
	&=
	\#\{a < K_j : \mu_a > j \}
	-
	\#\{a > K_j : j-1 = \mu_a \}.
	\end{align*}
	Using these facts in \eqref{j-ratio}, we read
	\begin{align}
	\label{j-ratio-2}
	\frac{
		(v_{K_j,j})^{\bm{1}_{K_j \in \mathcal{Q}_j}}
	}
	{
		(v_{K_j,j-1})^{\bm{1}_{K_j \ge 1}}
	}
	=
	q^{-1}
	\cdot
	\left\{
	\begin{array}{ll}
	t^{-\#\{a > K_j : \mu_a = j-1 \}}, \quad\quad
	& K_j \in \mathcal{Q}_j,
	\\
	\\
	t^{\#\{a < K_j : \mu_a > j \}
		- \#\{a > K_j : \mu_a = j-1 \}}, \quad\quad
	& K_j \in \mathcal{P}_j.
	\end{array}
	\right.
	\end{align}
	We can now see that in either case, the $t$ exponents cancel perfectly with the remaining factors in the $j$-th term of the product \eqref{N-columns-2}, and hence
	\begin{align*}
	\frac{
		(v_{K_j,j})^{\bm{1}_{K_j \in \mathcal{Q}_j}}
	}
	{
		(v_{K_j,j-1})^{\bm{1}_{K_j \ge 1}}
	}
	\cdot
	\frac{
		t^{\#\{a>K_{j} : \mu_a=j-1\} \bm{1}_{K_{j} \ge 1}}
	}
	{
		t^{\#\{a<K_{j} : \mu_a > j\} \bm{1}_{K_j \in \mathcal{P}_j}}
	}
	=
	q^{-\bm{1}_{K_j \ge 1}}.
	\end{align*}
	Returning to the expression \eqref{N-columns-2}, we have shown that
	\begin{align}
	\frac{Z_{\rm l}}{Z_{\rm r}}
	&=
	\frac{
		(v_{i,0})^{\bm{1}_{i \in \mathcal{Q}_0}}
	}
	{t^{\#\{ a<i : \mu_a > 0 \} \bm{1}_{i \in \mathcal{P}_0}}},
	\end{align}
	which expresses the remarkable fact that the ratio $Z_{\rm l}/Z_{\rm r}$ does not depend on any of the values of the colors $k^{(j)}_b$, not even those for which $b=i$. Finally, we check that
	\begin{align}
	\frac{
		(v_{i,0})^{\bm{1}_{i \in \mathcal{Q}_0}}
	}
	{t^{\#\{ a<i : \mu_a > 0 \} \bm{1}_{i \in \mathcal{P}_0}}}
	=
	\left\{
	\begin{array}{ll}
	v_{i,0} = q^{\mu_i} t^{\gamma_{i,0}(\mu)}, \quad\quad
	& i \in \mathcal{Q}_0,
	\\
	\\
	t^{-\#\{ a<i : \mu_a > 0 \}}
	=
	q^{\mu_i} t^{\gamma_{i,0}(\mu)}, \quad\quad
	& i \in \mathcal{P}_0,
	\end{array}
	\right.
	\end{align}
	where in the latter case, we have noted that $i \in \mathcal{P}_0$ implies $\mu_i = 0$.
\end{proof}

Now we can quickly establish \Cref{omegaif}. 

\begin{proof}[Proof of \Cref{omegaif}]
	Since we have a one-to-one pairing of each non-vanishing term in the sum \eqref{left-refined} with a corresponding term in \eqref{right-refined}, and we have also shown that the two terms have the correct proportionality constant \eqref{refined-cyclic-rel} irrespective of the values of the internal states $k^{(j)}_b$, \eqref{cyclic-rel} follows immediately. 
	%Note that, by repeating the above arguments with ``permuted'' boundary conditions at the left edges of the lattice, we prove \Cref{prop2-intro} from the introduction.
\end{proof}

\section{Proof of \Cref{coefficientc}}
\label{ssec:norm}

In this section we establish \Cref{coefficientc}. To that end, we make the following claim.
\begin{prop}
	
	\label{coefficientxc}
	
	The quantity ${\rm Coeff}[\big\langle\mathsf{C}_n(x_n) \cdots \mathsf{C}_1(x_1) \big\rangle_{\mu}; \textbf{\emph{x}}^{\mu}]$ is generically non-vanishing and is given by the weight of the unique lattice configuration of the form
	\begin{align}
	\label{frozen-xmu}
	\sum_{\mathscr{M}}
	\prod_{i=1}^{n}
	\prod_{j=0}^{N}
	v_{i,j}^{M_{i,j}}
	\times
	\tikz{.9}{
		\foreach\y in {0,...,5}{
			\draw[lgray,line width=1.5pt] (1.5,0.5+\y) -- (7.5,0.5+\y);
		}
		\foreach\x in {0,...,6}{
			\draw[lgray,line width=1.5pt] (1.5+\x,0.5) -- (1.5+\x,5.5);
		}
		%bottom labels
		\node[below] at (2,0.5) {\footnotesize$\textbf{\emph{M}}_0$};
		\node[above,text centered] at (4,0) {$\cdots$};
		\node[above,text centered] at (5,0) {$\cdots$};
		\node[below] at (7,0.5) {\footnotesize$\textbf{\emph{M}}_N$};
		%top labels
		\node[above] at (2,5.5) {\footnotesize$\textbf{\emph{M}}_0+\textbf{\emph{I}}_0 (\mu)$};
		\node[above,text centered] at (4,5.5) {$\cdots$};
		\node[above,text centered] at (5,5.5) {$\cdots$};
		\node[above] at (7,5.5) {\footnotesize$\textbf{\emph{M}}_N+\textbf{\emph{I}}_N (\mu)$};
		%right labels
		\node[right] at (7.5,1) {$0$};
		\node[right] at (7.5,2) {$0$};
		\node[right] at (7.5,3.1) {$\vdots$};
		\node[right] at (7.5,4.1) {$\vdots$};
		\node[right] at (7.5,5) {$0$};
		%spectral parameters
		\node[left] at (1,1) {$x_1$};
		\node[left] at (1,2) {$x_2$};
		\node[left] at (1,3.1) {$\vdots$};
		\node[left] at (1,4.1) {$\vdots$};
		\node[left] at (1,5) {$x_n$};
		%left labels
		\node[left] at (1.5,1) {$1$};
		\node[left] at (1.5,2) {$2$};
		\node[left] at (1.5,3.1) {$\vdots$};
		\node[left] at (1.5,4.1) {$\vdots$};
		\node[left] at (1.5,5) {$n$};
		%paths
		\draw[ultra thick,red,->] (1.5,5) -- (2.5,5) -- (4,5) -- (4,5.5);
		\draw[ultra thick,orange,->] (1.5,4) -- (2,4) -- (2,5.5);
		\draw[ultra thick,yellow,->] (1.5,3) -- (2.5,3) -- (5.1,3) -- (5.1,5.5);
		\draw[ultra thick,green,->] (1.5,2) -- (7,2) -- (7,5.5); 
		\draw[ultra thick,blue,->] (1.5,1) -- (2.5,1) -- (4.9,1) -- (4.9,5.5);
	}
	\end{align}
	in which path $i$ travels straight for $\mu_i$ consecutive horizontal steps, before turning and exiting at the top of the $\mu_i$-th column, for all 
	$1 \le i \le n$.
\end{prop}

\begin{proof}
	In what follows let $N = \max_{i \in [1, n]} \mu_i$. Let $i$ be the largest integer such that $\mu_i = N$ (there may be other integers $j$ such that $\mu_j = N$, and we assume that $i>j$ for all such $j$). Then there is no other path which travels further, horizontally, than the path of color $i$ does in order to reach its final destination in column $N$. 
	
	Given that we wish to calculate ${\rm Coeff}[\big\langle\mathsf{C}_n(x_n) \cdots \mathsf{C}_1 (x_1) \big\rangle_{\mu}; \textbf{x}^{\mu}]$, we restrict our attention to lattice configurations which give rise to a factor of $x_i^{\mu_i} = x_i^N$. We claim that the only possible way in which this factor is obtained is when the path of color $i$ travels horizontally for $N$ consecutive steps before turning into the $N$-th column. In other words, isolating the $i$-th row of the lattice,
	
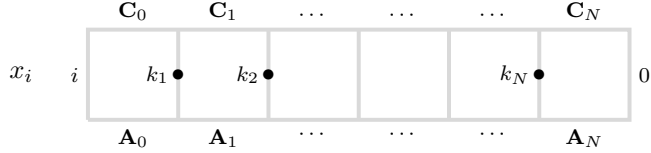
\begin{align*}
	\tikz{1.2}{
		\draw[lgray,line width=1.5pt] (0.5,-0.5) -- (6.5,-0.5) -- (6.5,0.5) -- (0.5,0.5) -- (0.5,-0.5);
		\foreach\x in {1,...,5}{
			\draw[lgray,line width=1.5pt] (0.5+\x,-0.5) -- (0.5+\x,0.5);
		}
		\node[left] at (0,0) {$x_i$};
		\node[left] at (1.5,0) {\fs $k_1$}; \node at (1.5,0) {$\bullet$};
		\node[left] at (2.5,0) {\fs $k_2$}; \node at (2.5,0) {$\bullet$};
		\node[left] at (5.5,0) {\fs $k_N$}; \node at (5.5,0) {$\bullet$};
		\node[left] at (0.5,0) {\fs $i$};\node[right] at (6.5,0) {\fs $0$};
		\node[below] at (6,-0.5) {\fs $\textbf{A}_N$};\node[above] at (6,0.5) {\fs $\textbf{C}_N$};
		\node[below] at (5,-0.5) {\fs $\cdots$};\node[above] at (5,0.5) {\fs $\cdots$};
		\node[below] at (4,-0.5) {\fs $\cdots$};\node[above] at (4,0.5) {\fs $\cdots$};
		\node[below] at (3,-0.5) {\fs $\cdots$};\node[above] at (3,0.5) {\fs $\cdots$};
		\node[below] at (2,-0.5) {\fs $\textbf{A}_1$};\node[above] at (2,0.5) {\fs $\textbf{C}_1$};
		\node[below] at (1,-0.5) {\fs $\textbf{A}_0$};\node[above] at (1,0.5) {\fs $\textbf{C}_0$};
	}
	\end{align*} 
	and labelling its internal edges as $k_1,\dots,k_N$ as shown above, we claim that the only possible contribution to $x_i^N$ comes in the case $k_1 = \cdots = k_N = i$. To prove this, we will assume that some other lattice configuration exists that gives rise to the factor $x_i^N$, in which not all $k_1,\dots,k_N$ are equal to $i$, and show that we obtain a contradiction.
	
	To begin, we note that all $k_1,\dots,k_N$ must be nonzero, or else we cannot recover the required degree in $x_i$. Assume that for some integer $p_1$ one has a lattice configuration in which $k_{p_1} = i$ and 
	$k_{p_1+1} = j_1$, where $i < j_1$ (note that one cannot have $i > j_1 \ge 1$, since the corresponding vertex weight in \Cref{lz} vanishes). Now the path of color $j_1$ must turn out of the $i$-th row somewhere before the $N$-th column, since by assumption $i$ is the largest color such that $\mu_i = N$. Let the column where this turning happens be labelled $p_2$. Then we must have $k_{p_2} = j_1$ and $k_{p_2+1} = j_2$, where $j_1 < j_2$. We again find that the path of color $j_2$ must turn out of the $i$-th row somewhere before the $N$-th column. One may now iterate this reasoning, ultimately arriving at the contradiction that $k_N$ will be forced to assume a value greater than $i$, which is impossible.
	
	We have thus proved the required statement about the path of color $i$, completely constraining its motion to the $i$-th row, and it plays no role in configurations of the remaining rows. Then let us define 
	$\hat{\mu} = (\hat\mu_1,\dots,\hat\mu_{n-1}) = (\mu_1,\dots,\mu_{i-1},\mu_{i+1},\dots,\mu_n)$, \ie\ we omit the $i$-th part from $\mu$. One may now apply a similar reasoning to the largest integer $j$ such that 
	$\mu_j = \max_{i \in [1, n - 1]} \hat\mu_i$, arriving at the same conclusion for the path of color $j$; namely, it is forced to take 
	$\max_{i \in [1, n - 1]} \hat\mu_i$ consecutive horizontal steps before turning into the $\mu_j$-th column, where it terminates. By iterating this argument over each of the colors, we arrive at the statement of the proposition.
\end{proof}

Now we can establish \Cref{coefficientc}. 

\begin{proof}[Proof of \Cref{coefficientc}]
	
Let us compute the weight of the configuration shown in \eqref{frozen-xmu}. To do this, we use the formula \eqref{components} to calculate the weight of each of the columns in \eqref{frozen-xmu}, and multiply them together. Each of the columns that appear in the frozen configuration \eqref{frozen-xmu} take a much simpler form than the generic column components \eqref{components}; we denote the weight of the $j$-th column by $X_j$, and remark that it is a function only of the colors $\{a: \mu_a = j\}$ which exit via its top and those $\{b: \mu_b > j\}$ which horizontally traverse it. We are able to write
\begin{align}
\label{something}
{\rm Coeff}\Big[ 
\big\langle\mathsf{C}_n(x_n) \cdots \mathsf{C}_1 (x_1) \big\rangle_{\mu};
\textbf{x}^{\mu}
\Big]
=
\prod_{j=0}^{N}
X_j(\mathcal{P}_j ; \mathcal{Q}_j),
\end{align}
with $\mathcal{P}_j$ and $\mathcal{Q}_j$ as given by \eqref{PQ_j}, and where we have defined the function
\begin{align}
\label{coeff-prod}
X_j(\mathcal{P} ; \mathcal{Q})
=
\prod_{p \in \mathcal{Q}}
(1-v_{p,j} t^{f(p)+1})
=
\prod_{p \in \mathcal{Q}}
(1-v_{p,j} t^{\#\{\ell \in \mathcal{Q} : \ell < p\}+1}).
\end{align}
Rewriting the product \eqref{coeff-prod} more explicitly, we have
\begin{flalign}
\label{coeff-prod2}
\begin{aligned}
{\rm Coeff}\Big[ 
\big\langle\mathsf{C}_n(x_n) \cdots \mathsf{C}_1 (x_1) \big\rangle_{\mu};
\textbf{x}^{\mu}
\Big] 
& =
\prod_{j=0}^{N}
\prod_{i: \mu_i > j}
(1-v_{i,j} t^{\#\{\ell < i : j < \mu_{\ell}\}+1}) \\
& = 
\prod_{i=1}^{n}
\prod_{j=0}^{\mu_i - 1}
(1-v_{i,j} t^{\#\{\ell < i : j < \mu_{\ell}\}+1}).
\end{aligned}
\end{flalign}
To conclude our calculation, we recall the explicit form \eqref{v_ij} of the parameters $v_{i,j}$. They contribute a factor of $t^{\gamma_{i,j}(\mu)}$ with $\gamma_{i,j}(\mu)$ given by \eqref{gamma}, which can be combined with the $t^{\#\{\ell < i : j < \mu_{\ell}\}}$ term appearing in \eqref{coeff-prod2}. Collecting the exponents proceeds as follows:
\begin{align*}
\#\{\ell < i : j < \mu_{\ell}\}
+
\gamma_{i,j}(\mu)
&=
\#\{\ell < i : j < \mu_\ell\}
-
\#\{k < i : \mu_k > \mu_i\}
+
\#\{k>i : j \le \mu_k < \mu_i\}
\\
&=
\#\{k < i : j < \mu_k \le \mu_i\}
+
\#\{k > i : j \le \mu_k < \mu_i\}
\\
&=
\#\{k < i : \mu_k = \mu_i\}
+
\#\{k \not=i : j < \mu_k < \mu_i\}
+
\#\{k > i : j = \mu_k \}
\\
&=
\alpha_{i,j}(\mu),
\end{align*}
with $\alpha_{i,j}(\mu)$ given by \eqref{alpha}. Putting everything together, we have shown that
\begin{align}
\label{coeff-prod3}
{\rm Coeff}\Big[ 
\big\langle\mathsf{C}_n(x_n) \cdots \mathsf{C}_1(x_1) \big\rangle_{\mu};
\textbf{x}^{\mu}
\Big]
=
\prod_{i=1}^{n}
\prod_{j=0}^{\mu_i-1}
(1-q^{\mu_i-j} t^{\alpha_{i,j}(\mu)+1})
=
\frac{1}{\Omega_{\mu}(q,t)},
\end{align}
with $\Omega_{\mu}(q,t)$ given by \eqref{Omega}. 
\end{proof}

\chapter{Vertex Models for Symmetric Macdonald Polynomials}

In this chapter we describe a symmetrization procedure which yields the symmetric Macdonald polynomials, in their integral form. Throughout this section we fix complex numbers $q, t \in \mathbb{C}$ and adopt the diagrammatic conventions introduced in \Cref{Column}.

\section{Partition Function for Integral Macdonald Polynomials} 

\label{SumJ}

In this section we first provide an expression for the integral, symmetric Macdonald polynomials in terms of the $\mathsf{D} (x)$ operators from \Cref{dci}, which implies a decomposition of these integral Macdonald polynomials into sums of the $G_{\boldsymbol{\lambda} / \boldsymbol{\mu}}$ functions from \Cref{fgdefinition}. To define the former, for any signature $\lambda = (\lambda_1, \lambda_2, \ldots , \lambda_{\ell}) \in \Sign$ and sequence $\textbf{x} = (x_1, x_2, \ldots , x_N)$ of complex numbers, recall the Macdonald polynomial $P_{\lambda} (\textbf{x}) = P_{\lambda} (\textbf{x}; q, t)$\index{P@$P_{\lambda} (\textbf{x})$; symmetric Macdonald polynomial} from (6.4.7) of \cite{SFP}. We further recall from equation (6.8.3) of \cite{SFP} the \emph{integral, symmetric Macdonald polynomial} $J_{\lambda} (\textbf{x}) = J_{\lambda} (\textbf{x}; q, t)$,\index{J@$J_{\lambda} (\textbf{x})$; integral Macdonald polynomial} defined by setting
\begin{flalign}
\label{cqt} 
J_{\lambda} (\textbf{x}) = c_{\lambda} (q, t) P_{\lambda} (x), \qquad \text{where} \quad c_{\lambda} (q, t) = \displaystyle\prod_{b \in \mathscr{Y} (\lambda)} (1 - q^{a(b)} t^{l(b) + 1}),
\end{flalign}\index{C@$c_{\lambda} (q, t)$}

\noindent and the product is over all boxes $b$ in the Young diagram $\mathscr{Y} (\lambda)$ for $\lambda$. Here, $a(b)$ and $l(b)$ denote the \emph{arm} and \emph{leg} lengths of $b \in \mathscr{Y} (\lambda)$, given by the number of boxes in $\mathscr{Y} (\lambda)$ east and south of $b$, respectively; see equation (6.6.14) of \cite{SFP}.  

Now let us define the following polynomial, which we will eventually show to coincide with the integral Macdonald polynomial, in terms of the $\mathsf{D} (x)$ operators. In what follows, we recall the sequence $\mathscr{I} (\nu) = \big( \textbf{I}_0 (\nu), \textbf{I}_1 (\nu), \ldots \big)$ of elements in $\{ 0, 1 \}^n$ from \Cref{ssec:states}, for any composition $\nu$ of length at most $n$. 
\begin{definition} 
	
	\label{pd} 
	
	Fix an integer $p \in [1, n]$; a sequence $\textbf{x} = (x_1, x_2, \ldots , x_n)$ of complex numbers; and a (positive) composition $\nu = (\nu_1, \nu_2, \ldots , \nu_p)$ of length $p$. Denoting $\nu_k = 0$ for $k > p$ and letting $N = \max_{i \in [1, p]} \nu_i$, define the set of twist parameters $\textbf{w} = (w_{i, j})$, where $(i, j)$ ranges over all pairs of integers $(i, j) \in [1, n] \times [1, N]$, by setting
	\begin{align}
	\label{w_ij}
	w_{i,j}
	=
	q^{\nu_i-j}
	\bm{1}_{\nu_i > j}.
	\end{align}
	
	\noindent Then, define
	\begin{flalign}
	\label{pnu}
	\mathcal{P}_{\nu} (\textbf{x}) = \displaystyle\sum_{\mathscr{M}} \displaystyle\prod_{i = 1}^n \displaystyle\prod_{j = 1}^N w_{i, j}^{M_{i, j}} \Big\langle \mathscr{M} + \mathscr{I} (\nu) \big| \mathsf{D} (x_n) \cdots \mathsf{D} (x_1) \big| \mathscr{M} + \big(\textbf{e}_{[1, p]}, \textbf{e}_0, \textbf{e}_0, \ldots \big) \Big\rangle,
	\end{flalign} 
	
	\noindent where we sum over all sequences $\mathscr{M} = (\textbf{M}_0, \textbf{M}_1, \ldots )$ of elements in $\{ 0, 1 \}^n$ such that $\textbf{M}_k = \textbf{e}_0$ for $k = 0$ and each $k > N$; we also set $\textbf{M}_j = (M_{1, j}, M_{2, j}, \ldots , M_{n, j}) \in \{ 0, 1 \}^n$ for each $j \in [1, N]$, and $\textbf{e}_{[1, p]} = (1^p, 0^{n - p}) \in \{ 0, 1 \}^n$ to be the vector whose first $p$ coordinates are equal to $1$ and whose last $n - p$ coordinates are equal to $0$. 
\end{definition} 

In particular, $\mathcal{P}_{\nu} (\textbf{x})$ is diagrammatically given by the partition function
\begin{align}
\label{sym-pf}
\mathcal{P}_{\nu}(\textbf{x})
=
\sum_{\mathscr{M}}
\prod_{i=1}^n
\prod_{j=1}^{N}
w_{i,j}^{M_{i,j}}
\times
\tikz{1}{
	\foreach\y in {0,...,5}{
		\draw[lgray,line width=1.5pt] (1.5,0.5+\y) -- (7.5,0.5+\y);
	}
	\foreach\x in {0,...,6}{
		\draw[lgray,line width=1.5pt] (1.5+\x,0.5) -- (1.5+\x,5.5);
	}
	%spectral parameters
	\node[text centered] at (2,1) {$x_1$};
	\node[text centered] at (4,1) {$\cdots$};
	\node[text centered] at (5,1) {$\cdots$};
	\node[text centered] at (7,1) {$x_1$};
	\node[text centered] at (2,2) {$x_2$};
	\node[text centered] at (4,2) {$\cdots$};
	\node[text centered] at (5,2) {$\cdots$};
	\node[text centered] at (7,2) {$x_2$};
	\node[text centered] at (2,5) {$x_n$};
	\node[text centered] at (4,5) {$\cdots$};
	\node[text centered] at (5,5) {$\cdots$};
	\node[text centered] at (7,5) {$x_n$};
	\node[text centered] at (2,3.1) {$\vdots$};
	\node[text centered] at (2,4.1) {$\vdots$};
	\node[text centered] at (7,3.1) {$\vdots$};
	\node[text centered] at (7,4.1) {$\vdots$};
	%bottom labels
	\node[below] at (2,0.5) {\scriptsize$\textbf{e}_{[1, p]}$};
	\node[below] at (3,0.5) {\scriptsize$\textbf{M}_1$};
	\node[above,text centered] at (4,0) {$\cdots$};
	\node[above,text centered] at (5,0) {$\cdots$};
	\node[below] at (7,0.5) {\scriptsize$\textbf{M}_N$};
	%top labels
	\node[left = 5, above] at (2,5.5) {\scriptsize $\textbf{e}_0$};
	\node[above] at (3,5.5) {\scriptsize$\textbf{M}_1+\textbf{I}_1$};
	\node[above,text centered] at (4,5.5) {$\cdots$};
	\node[above,text centered] at (5,5.5) {$\cdots$};
	\node[above] at (7,5.5) {\scriptsize$\textbf{M}_N+\textbf{I}_N$};
	%right labels
	\node[right] at (7.5,1) {$0$};
	\node[right] at (7.5,2) {$0$};
	\node[right] at (7.5,3.1) {$\vdots$};
	\node[right] at (7.5,4.1) {$\vdots$};
	\node[right] at (7.5,5) {$0$};
	%left labels
	\node[left] at (1.5,1) {$0$};
	\node[left] at (1.5,2) {$0$};
	\node[left] at (1.5,3.1) {$\vdots$};
	\node[left] at (1.5,4.1) {$\vdots$};
	\node[left] at (1.5,5) {$0$};
}
\end{align} \index{P@$\mathcal{P}_{\nu} (\textbf{x})$}
\noindent where we have abbreviated $\textbf{I}_j = \textbf{I}_j (\nu)$ for each $j \ge 1$. Here, colors $1, 2, \dots, p$ all enter at the base of the leftmost column, and this column does not have cylindrical boundary conditions (while the others do).

Now we can state the following theorem, equating the integral, symmetric Macdonald polynomial $J_{\lambda} (\textbf{x})$ with the partition function $\mathcal{P}_{\lambda} (\textbf{x})$ from \Cref{pd}. We will establish this result in \Cref{ProofJ} below. 

\begin{thm}
	
	\label{jpd} 
	
	Fix an integer $p \in [1, n]$; a composition $\nu = (\nu_1, \nu_2, \ldots , \nu_p)$ of length $p$ that is anti-dominant, namely, $1 \le \nu_1 \le \nu_2 \le \cdots \le \nu_p$; and a set $\textbf{\emph{x}} = (x_1, x_2, \ldots , x_n)$ of complex numbers. Then, $J_{\nu^+} (\textbf{\emph{x}}) = \mathcal{P}_{\nu} (\textbf{\emph{x}})$, where we recall that $\nu^+$ denotes the dominant ordering of $\nu$.
	
\end{thm}

\begin{rem}
	
	\label{psumg} 
	
Observe by \eqref{wl} that the $\big\langle \mathscr{M} + \mathscr{I} (\nu) \big| \mathsf{D} (x_n) \cdots \mathsf{D} (x_1) \big| \mathscr{M} + \big( \textbf{e}_{[1, p]}, \textbf{e}_0, \textbf{e}_0, \ldots \big) \big\rangle$ appearing in \eqref{pnu} are specializations of the functions $G_{\boldsymbol{\lambda} / \boldsymbol{\mu}}$ from\footnote{Here, the constituent signatures in $\boldsymbol{\lambda} \in \SeqSign_n$ and in $\boldsymbol{\mu} \in \SeqSign_n$ have different lengths (see \Cref{fgmn})} \Cref{fgdefinition}, whose $r$ and $s$ parameters are all equal to $t^{-1 / 2}$ and $0$, respectively (and whose quantization parameter is equal to $t$). In particular, \Cref{jpd} decomposes the integral Macdonald polynomial $J_{\nu^+} (\textbf{x})$ as a linear combination of such functions. 

\end{rem} 

Before proceeding to the proof of \Cref{jpd}, let us explain how it can be used to recover an expression from \cite{CP} for the modified Macdonald polynomials as a linear combination of skew LLT polynomials. To define the former, for any partition $\lambda = (\lambda_1, \lambda_2, \ldots , \lambda_M)$ and sequence of complex variables $\textbf{x} = (x_1, x_2, \ldots,  x_N)$, let the \emph{modified Macdonald polynomial} $\widetilde{J}_{\lambda} (\textbf{x})$\index{J@$J_{\lambda} (\textbf{x})$; integral Macdonald polynomial!$\widetilde{J}_{\lambda} (\textbf{x})$; modified Macdonald polynomial} be 
\begin{flalign}
	\label{jxjx}
	\widetilde{J}_{\lambda} (\textbf{x}) = J_{\lambda} \Bigg( \bigcup_{j = 0}^{\infty} t^j \textbf{x} \Bigg).
\end{flalign}

\noindent Equivalently, (see, for instance, equation (4) of \cite{CP}) it is given by the plethystic substitution $\widetilde{J}_{\lambda} [X] = J_{\lambda} \big[ (1 - t)^{-1} X \big]$, where $X$ denotes the formal sum $X = \sum_{i = 1}^N x_i$. In the below, we recall the sequence $\mathscr{S} (\boldsymbol{\lambda})$ of elements in $\{ 0, 1 \}^n$ for any $\boldsymbol{\lambda} \in \{ 0, 1 \}^n$ from \Cref{Symmetric}.

\begin{cor}[{\cite[Theorem 2.2, Equation (23), and Proposition 3.4]{CP}}]
	
	\label{psuml}
	 
	Fix integers $p, n, N \ge 1$ with $p \le n$; an anti-dominant composition $\nu$ of length $p$; and a sequence of complex numbers $\textbf{\emph{x}} = (x_1, x_2, \ldots,  x_N)$. For each $1 \le i \le n$ and $j \ge 1$, define
	\begin{flalign*}
		u_{i, j} = u_{i, j} (\nu) = q^{\nu_i - j} t^{\beta_{i, j} (\nu)} \textbf{\emph{1}}_{\nu_i > j},
	\end{flalign*}
	
	\noindent where $\nu_k = 0$ for $k > p$, and the exponents $\beta_{i, j} (\nu)$ are given by 
	\begin{align*}
		\beta_{i,j} (\nu) = \displaystyle\frac{1}{2} \Big( \#\{ k < i: \nu_k > j \} - p + i \Big).
	\end{align*}
	
	\noindent For any infinite sequence $\mathscr{K} = (\textbf{\emph{K}}_1, \textbf{\emph{K}}_2, \ldots )$ of elements in $\{ 0, 1 \}^n$, define $\boldsymbol{\mu} (\mathscr{K}) \in \SeqSign_n$ so that $\mathscr{S} \big( \boldsymbol{\mu} (\mathscr{K}) \big) = \mathscr{K}$. Then,
	\begin{flalign}
	\label{jnusum} 
		\widetilde{J}_{\nu^+} (\textbf{\emph{x}}) = \displaystyle\sum_{\mathscr{M}} \mathcal{L}_{\boldsymbol{\mu} (\mathscr{I}_0 (\nu) + \mathscr{M}_0) / \boldsymbol{\mu} (\mathscr{M}_p)} (\textbf{\emph{x}}; q) \displaystyle\prod_{i = 1}^n \displaystyle\prod_{j = 1}^{\infty} u_{i, j}^{M_{i, j}},
	\end{flalign}
	
	\noindent where the sum is over all infinite sequences $\mathscr{M} = (\textbf{\emph{M}}_1, \textbf{\emph{M}}_2, \ldots )$ of elements in $\{ 0, 1 \}^n$, with $\textbf{\emph{M}}_j = (M_{1, j}, M_{2, j}, \ldots , M_{n, j})$, such that $\textbf{\emph{I}}_k (\nu) + \textbf{\emph{M}}_k \in \{ 0, 1 \}^n$ for each $k \ge 1$. Here, we have denoted $\mathscr{M}_0 = (\textbf{\emph{e}}_0, \textbf{\emph{M}}_1, \textbf{\emph{M}}_2, \ldots )$; $\mathscr{I}_0 (\nu) = \big( \textbf{\emph{e}}_0, \textbf{\emph{I}}_1, \textbf{\emph{I}}_2, \ldots \big)$; and $\mathscr{M}_p = \big( \textbf{\emph{e}}_{[1, p]}, \textbf{\emph{M}}_1, \textbf{\emph{M}}_2, \ldots \big)$.
	
\end{cor}

\begin{proof} 
	
	By \eqref{pnu} and \Cref{jpd}, we have for any sequence of complex variables $\textbf{y} = (y_1, y_2, \ldots , y_K)$ that 
	\begin{flalign*}
	J_{\nu^+} (\textbf{y}) = \displaystyle\sum_{\mathscr{M}} \displaystyle\prod_{i = 1}^n \displaystyle\prod_{j = 1}^{\infty} w_{i, j}^{M_{i, j}} \big\langle \mathscr{M}_0 + \mathscr{I}_0 (\nu) \big| \mathsf{D} (y_K) \cdots \mathsf{D} (y_1) \big| \mathscr{M}_p \big\rangle,
	\end{flalign*} 

	\noindent where the sum is as in \eqref{jnusum}, and we have used the fact that $w_{i, j} = 0$ for $j > \nu_p$. Moreover, by \Cref{fgdefinition}, \eqref{wl}, \Cref{lwabcd}, \Cref{dci}, and the last statement of \eqref{limitg} we have
	\begin{flalign*}
	\mathcal{G}_{\boldsymbol{\mu} (\mathscr{M}_0 + \mathscr{I}_0 (\nu)) / \boldsymbol{\mu} (\mathscr{M}_p)} (\textbf{y}; q^{-1 / 2} \boldsymbol{\mid} 0; 0) = \big\langle \mathscr{M}_0 + \mathscr{I}_0 (\nu) \big| \mathsf{D} (y_K) \cdots \mathsf{D} (y_1) \big| \mathscr{M}_p \big\rangle,
	\end{flalign*}	
	
	\noindent and so
	\begin{flalign}
		\label{jnuy} 
		J_{\nu^+} (\textbf{y}) = \displaystyle\sum_{\mathscr{M}} \mathcal{G}_{\boldsymbol{\mu} (\mathscr{M}_0 + \mathscr{I}_0 (\nu)) / \boldsymbol{\mu} (\mathscr{M}_p)} (\textbf{y}; q^{-1 / 2} \boldsymbol{\mid} 0; 0) \displaystyle\prod_{i = 1}^n \displaystyle\prod_{j = 1}^{\infty} w_{i, j}^{M_{i, j}}.
	\end{flalign} 

	Taking $\textbf{y} = \bigcup_{j = 1}^N \{ x_j, q x_j, \ldots \}$ in \eqref{jnuy}, and then appyling \eqref{jxjx} and (the $L_i = \infty$ limiting case of) \Cref{gq}, yields 
	\begin{flalign}
		\label{sumj} 
		\widetilde{J}_{\nu^+} (\textbf{x}) = \displaystyle\sum_{\mathscr{M}} \mathcal{G}_{\boldsymbol{\mu} (\mathscr{M}_0 + \mathscr{I}_0 (\nu)) / \boldsymbol{\mu} (\mathscr{M}_p)} (\textbf{x}; \infty \boldsymbol{\mid} 0; 0) \displaystyle\prod_{i = 1}^n \displaystyle\prod_{j = 1}^{\infty} w_{i, j}^{M_{i, j}}.
	\end{flalign} 
	
	\noindent By the first statement \Cref{limitg0} (and recalling the function $\psi$ from \eqref{lambdamupsi}), it follows that 
	\begin{flalign}
		\label{jxnu}
		\widetilde{J}_{\nu^+} (\textbf{x}) = \displaystyle\sum_{\mathscr{M}} q^{\psi (\boldsymbol{\mu} (\mathscr{M}_0 + \mathscr{I}_0 (\nu))) - \psi (\boldsymbol{\mu} (\mathscr{M}_p))} \mathcal{L}_{\boldsymbol{\mu} (\mathscr{M}_0 + \mathscr{I}_0 (\nu)) / \boldsymbol{\mu} (\mathscr{M}_p)} (\textbf{x}; q) \displaystyle\prod_{i = 1}^n \displaystyle\prod_{j = 1}^{\infty} w_{i, j}^{M_{i, j}}.
	\end{flalign}
	
	\noindent Let us analyze the exponent of $q$ appearing on the right side of \eqref{jxnu}, assuming that $M_{i, j} = 0$ whenever $i > p$ (which we may do, since $w_{i, j} = 0$ whenever $i > p$). Since $\nu$ is anti-dominant, we have from \eqref{lambdamupsi} that 
	\begin{flalign*}
	\psi \Big( \boldsymbol{\mu} \big( \mathscr{I}_0 (\nu) + \mathscr{M}_0 \big) \Big) & = \displaystyle\frac{1}{2} \displaystyle\sum_{1 \le i < h \le n} \displaystyle\sum_{j = \nu_h + 1}^{\infty} M_{i, j} + \displaystyle\frac{1}{2} \displaystyle\sum_{1 \le h < i \le n} \displaystyle\sum_{j = 1}^{\nu_h - 1} M_{i, j} + \displaystyle\frac{1}{2} \displaystyle\sum_{1 \le h < i \le n} \displaystyle\sum_{1 \le j < k} M_{h, k} M_{i, j}; \\
	\psi \big( \boldsymbol{\mu} (\mathscr{M}_p) \big) & = \displaystyle\frac{1}{2} \displaystyle\sum_{1 \le h < i \le n} \displaystyle\sum_{1 \le j < k} M_{h, k} M_{i, j} + \displaystyle\frac{1}{2} \displaystyle\sum_{i = 1}^p \displaystyle\sum_{j = 1}^{\infty} (p - i) M_{i, j},
	\end{flalign*} 

	\noindent meaning that
	\begin{flalign}
		\label{sumpsi2}
		\begin{aligned}
		\psi \Big( \boldsymbol{\mu} \big( \mathscr{I}_0 (\nu) + \mathscr{M}_0 \big) \Big) - \psi \big( \boldsymbol{\mu} (\mathscr{M}_p) \big) & = \displaystyle\frac{1}{2} \displaystyle\sum_{i = 1}^p \displaystyle\sum_{j = 1}^{\infty} M_{i, j} \Big( \#\{ h > i: \nu_h < j \} + \#\{ h < i: \nu_h > j \} - p + i \Big). 
		\end{aligned} 
	\end{flalign}

	\noindent Since $\nu$ is anti-dominant, there exists some $k > i$ with $\nu_k < j$ only if $\nu_i < j$. This would imply $w_{i, j} = 0$, in which case the $(i, j)$ summand on the right side of \eqref{jxnu} only contributes if $M_{i, j} = 0$. Hence, we may assume in what follows that $M_{i, j} = 0$ whenever there exists some $k > i$ with $\nu_k < j$; inserting this into \eqref{sumpsi2} yields
	\begin{flalign*}
			\psi \Big( \boldsymbol{\mu} \big( \mathscr{I}_0 (\nu) + \mathscr{M}_0 \big) \Big) - \psi \big( \boldsymbol{\mu} (\mathscr{M}_p) \big) & = \displaystyle\frac{1}{2} \displaystyle\sum_{i = 1}^p \displaystyle\sum_{j = 1}^{\infty} M_{i, j} \Big( \#\{ k < i: \nu_k > j \} - p + i \Big) = \displaystyle\sum_{i = 1}^p \displaystyle\sum_{j = 1}^{\infty} \beta_{i, j},
	\end{flalign*}

	\noindent which yields the corollary upon insertion into \eqref{jxnu}, since $u_{i, j} = q^{\beta_{i, j}} w_{i, j}$. 
\end{proof}

\section{Identifying the Macdonald Polynomial} 

In this section we show that $\mathcal{P}_{\nu} (\textbf{x})$ coincides with a symmetric Macdonald polynomial, up to a constant factor that we determine in \Cref{ProofJ} below. We first quickly show that it is symmetric in $\textbf{x}$, through the following lemma.

\begin{lem}
	\label{prop-sym}
	The polynomial $\mathcal{P}_{\nu}(\textbf{\emph{x}})$ is symmetric in $\textbf{\emph{x}} = (x_1, x_2, \ldots , x_n)$.
\end{lem}

\begin{proof}
	
	By \eqref{pnu}, $\mathcal{P}_{\nu}(\textbf{x})$ is given by an appropriate linear form acting on $\mathsf{D}(x_1) \dots \mathsf{D}(x_n)$. Since the $\mathsf{D}(x_i)$ operators commute by \eqref{dl1d}, the symmetry is immediate.
\end{proof}

We next show that any $\mathcal{P}_{\nu} (\textbf{x})$ is in the span of all nonsymmetric Macdonald polynomials $f_{\mu} (\textbf{x})$ with $\mu^+ = \nu^+$. To that end, we begin by providing an alternative expression for the latter nonsymmetric polynomials involving both $\mathsf{C}_i$ and $\mathsf{D}$ operators. This will be given by \Cref{pdc} below, to establish which we first show the following proposition that essentially ``separates out'' the leftmost column from the remaining part of the $\mathsf{C}_n (x_n) \cdots \mathsf{C}_1 (x_1)$ type partition function for this nonsymmetric polynomial $f_{\mu}$ from \Cref{thm:f-formula}. In what follows, we let $\mu$ be a composition formed by concatenating $p$ non-zero parts and $n-p$ zeros, namely,
\begin{align}
\label{concat-comp}
\mu = (\nu_1,\dots,\nu_p,0^{n-p}),
\qquad
\text{with $\nu_i \ge 1$, for all $i \in [1, p]$}.
\end{align}
\begin{prop}
	
	\label{csumdc} 
	
	Let $\mu = (\mu_1,\dots,\mu_n)$ and $\nu = (\nu_1,\dots,\nu_p)$ be two compositions related as in \eqref{concat-comp}, and denote the largest part of $\mu$ by $N = \max_{i \in [1, n]} \mu_i = \max_{i \in [1, p]} \nu_i$. Fix a set $\textbf{\emph{v}} = (v_{i, j})$ of complex numbers, where $(i, j)$ ranges over all integer pairs $(i, j) \in [1, n] \times [0, N]$, such that $v_{i, j} = 0$ for any $(i, j) \in [p + 1, n] \times [0, N]$. Then, we have 
	\begin{align}
	\label{modified-match}
	\big\langle\mathsf{C}_n (x_n) \cdots \mathsf{C}_1 (x_1) \big\rangle_{\mu}(\textbf{\emph{v}})
	=
	\prod_{i=1}^{p} (1-v_{i,0} t^i) x_i
	\cdot
	\big\langle
	\mathsf{D} (x_n)
	\cdots
	\mathsf{D} (x_{p + 1})
	\mathsf{C}_p (x_p) \cdots \mathsf{C}_1 (x_1) \big\rangle_{\nu}'(\textbf{\emph{v}}),
	\end{align}
	with the two linear forms are given by
	\begin{align*}
	\big\langle \mathsf{X} \big\rangle_{\mu}(\textbf{\emph{v}})
	=
	\sum_{\mathscr{M}}
	\prod_{i=1}^{n}
	\prod_{j=0}^{N}
	v_{i,j}^{M_{i,j}}
	\big\langle \mathscr{M} + \mathscr{I} (\mu) \big| \mathsf{X} \big| \mathscr{M} \big\rangle; \qquad 
	\big\langle \mathsf{X} \big\rangle_{\nu}'(\textbf{\emph{v}})
	=
	\sum_{\mathscr{M}}
	\prod_{i=1}^n
	\prod_{j=1}^{N}
	v_{i,j}^{M_{i,j}}
	\big\langle \mathscr{M} + \mathscr{I} (\nu) \big|' \mathsf{X} \big| \mathscr{M} \big\rangle',
	\end{align*}

	\noindent where the first and second sums are over all sequences of $\{ 0, 1 \}^n$ elements $\mathscr{M} = (\textbf{\emph{M}}_0, \textbf{\emph{M}}_1, \ldots , \textbf{\emph{M}}_N)$ and $\mathscr{M} = (\textbf{\emph{M}}_1, \textbf{\emph{M}}_2, \ldots , \textbf{\emph{M}}_N)$, respectively, setting $\textbf{\emph{M}}_j = (M_{1, j}, M_{2, j}, \ldots , M_{n, j}) \in \{ 0, 1 \}^n$ for each $j$. Here, the primed states have the same definition as in \Cref{ssec:states}, but with all tensor products starting at $j=1$. 
\end{prop}

\begin{proof}
	This is a consequence of expanding the partition function on the left hand side of \eqref{modified-match} with respect to its zeroth column. This column is pictured as
	\begin{align}
	\label{zeroth-column}
	\sum_{\textbf{M}_0}
	\prod_{k=1}^{n}
	v_{k,0}^{M_{k, 0}}
	\tikz{.9}{
		\foreach\y in {1,...,4}{
			\draw[lgray,line width=1.5pt] (1.5,0.5+\y) -- (2.5,0.5+\y);
		}
		\draw[lgray,line width=1.5pt] (1.5,0.5) -- (2.5,0.5) -- (2.5,5.5) -- (1.5,5.5) -- (1.5,0.5);
		%spectral parameters
		\node[text centered] at (2,1) {$x_1$};
		\node[text centered] at (2,2) {$x_2$};
		\node[text centered] at (2,2.6) {$\vdots$};
		\node[text centered] at (2,3) {$x_p$};
		\node[text centered] at (2,4) {$x_{p+1}$};
		\node[text centered] at (2,4.6) {$\vdots$};
		\node[text centered] at (2,5) {$x_n$};
		%bottom labels
		\node[below] at (2,0.5) {\footnotesize$\textbf{M}_0$};
		%top labels
		\node[above] at (2,5.5) {\footnotesize$\textbf{M}_0+(0^p,1^{n-p})$};
		%right labels
		\node[right] at (2.5,1) {$d_1$};
		\node[right] at (2.5,2) {$d_2$};
		\node[right] at (2.5,2.6) {$\vdots$};
		\node[right] at (2.5,3) {$d_p$};
		\node[right] at (2.5,4) {$d_{p+1}$};
		\node[right] at (2.5,4.6) {$\vdots$};
		\node[right] at (2.5,5) {$d_n$};
		%left labels
		\node[left] at (1.5,1) {$1$};
		\node[left] at (1.5,2) {$2$};
		\node[left] at (1.5,2.6) {$\vdots$};
		\node[left] at (1.5,3) {$p$};
		\node[left] at (1.5,4) {$p+1$};
		\node[left] at (1.5,4.6) {$\vdots$};
		\node[left] at (1.5,5) {$n$};
	}
	\end{align}
	where $(d_1, d_2, \dots, d_n) \in \{0,1,\dots,p\}^n$. One can see that there is only one possible choice of $(d_1, d_2, \dots, d_n)$ which yields a non-vanishing result; namely, $(d_1, d_2, \dots,d_n) = (1,\dots,p,0^{n-p})$. Indeed, in view of the vanishing of the fifth weight in \Cref{vertexfigurerql} (that is, $L_x (\textbf{A}, j; \textbf{A}_{ji}^{+-}, i) = 0$ for $1 \le i < j \le n$) we see that necessarily $d_{p+1} = \cdots = d_n = 0$ in \eqref{zeroth-column}. This constrains $(d_1, d_2, \dots , d_p)$ to be a permutation of $(1, 2, \dots,p)$, and the only permutation which does not lead to vertices of the form $(\textbf{A}, j; \textbf{A}_{ji}^{+-}, i)$ for $1\le i < j \le n$ is precisely $(d_1, d_2, \dots , d_p) = (1, 2, \dots,p)$. 
	
	Using the third and last two weights in \Cref{vertexfigurerql}, the weight of the column is then readily computed as 
	\begin{flalign*}
 	\displaystyle\sum \displaystyle\prod_{i = 1}^p (-v_{i, 0} t^i)^{M_{i, 0}} x_i = \prod_{i=1}^{p} (1-v_{i,0} t^i) x_i,
	\end{flalign*}
	
	\noindent where the left side is summed over all $(M_{1, 0}, M_{2, 0}, \ldots , M_{p, 0}) \in \{ 0, 1 \}^p$. We now recognize this as the formula on the right hand side of \eqref{modified-match}.
\end{proof}

\begin{cor} 
	
	\label{pdc} 
	
	Let $\mu = (\mu_1,\dots,\mu_n)$ and $\nu = (\nu_1,\dots,\nu_p)$ be compositions which satisfy \eqref{concat-comp}, with $\nu$ anti-dominant. The nonsymmetric Macdonald polynomial $f_{\mu}(\textbf{\emph{x}})$ is given by
	\begin{align}
	\label{f-ansatz2}
	f_{\mu}(\textbf{\emph{x}})
	=
	\prod_{i=1}^p
	(1 - t^i w_{i, 0}) x_i
	\cdot
	\Omega_{\nu}(q,t)
	\big\langle
	\mathsf{D} (x_n)
	\cdots
	\mathsf{D} (x_{p + 1})
	\mathsf{C}_p (x_p) \cdots \mathsf{C}_1 (x_1) \big\rangle_{\nu}'(\textbf{\emph{w}})
	\end{align}
	where the parameters $\textbf{\emph{w}} = (w_{i,j})$ for $(i, j) \in [1, p] \times [1, N]$ are chosen to be as in \eqref{w_ij}. The normalization factor appearing on the right side of \eqref{f-ansatz2} is given by \eqref{Omega} and \eqref{alpha}. 
\end{cor}

\begin{proof}
	This follows from our previous formula \eqref{f-ansatz}, combined with \eqref{modified-match}, together with the fact that the exponent $\gamma_{i, j}$ from \eqref{gamma} is equal to $0$ since $\nu$ anti-dominant.
\end{proof}

Now we can establish the following proposition, indicating that $\mathcal{P}_{\nu} (\textbf{x})$ is in the linear span of certain nonsymmetric Macdonald polynomials. 

\begin{prop}
	\label{prop-V}
	For each partition $\lambda = \lambda_1 \ge \cdots \ge \lambda_n \ge 0$, let $\mathcal{V}_{\lambda} = {\rm Span}_{\mathbb{C}} \{ f_{\mu}(\textbf{\emph{x}}) \} _{\mu^{+} = \lambda}$. For any anti-dominant composition $\nu = (\nu_1,\dots,\nu_p)$ such that $\nu_i \ge 1$ for each $i \in [1, p]$, we have $\mathcal{P}_{\nu}(\textbf{\emph{x}}) \in \mathcal{V}_{\kappa}$, where $\kappa = (\nu^+, 0^{n - p})$ is the partition obtained by sorting parts into dominant order, and then appending $n-p$ zeros.
\end{prop}

\begin{proof}
	The proof proceeds by expanding over the configurations of the leftmost column in \eqref{sym-pf}. Performing that expansion, we obtain
	\begin{align}
	\label{sum}
	\mathcal{P}_{\nu}(x_1,\dots,x_n)
	=
	(1-t)^p
	\sum_{\mathfrak{D}}
	t^{{\rm inv}'(\overleftarrow{\mathfrak{D}})}
	\prod_{i=1}^{n}
	x_i^{\bm{1}_{d_i \ge 1}}
	\cdot
	\big\langle
	\mathsf{C}_{d_n}(x_n)
	\cdots
	\mathsf{C}_{d_1}(x_1)
	\big\rangle_{\nu}' (\textbf{w}),
	\end{align}
	\noindent where $\mathfrak{D} = (d_1,\dots,d_n)$ is summed over all permutations of the vector $(1,\dots,p,0^{n-p})$; we have let $\inv' (\overleftarrow{\mathfrak{D}})$ denote\footnote{Observe that this is slightly different from the definition of $\inv (\overleftarrow{\mathfrak{D}})$, which would count distinct indices $i, j$ with $d_i = 0$ and $d_j > 0$.} the number of pairs $1 \le i < j \le n$ such that $1 \le d_i < d_j$; and we have set $\mathsf{C}_0 (x) = \mathsf{D} (x)$. Indeed, the factor $(1 - t)^p \prod_{i = 1}^n x_i^{\textbf{1}_{d_i \ge 1}}$ arise from to the factors of $(1 - t) x$ in the second and fourth weights in \Cref{vertexfigurerql}, and the $t^{\inv' (\overleftarrow{\mathfrak{D}})}$ factor arises as the product of the $t^{A_{[d_i + 1, n]}}$ exponents in these weights.
	
	Examining the right side of \eqref{sum}, we see that the $\mathfrak{D} = (1,\dots,p,0^{n-p})$ term is known already; up to overall multiplicative constants, that term is a nonsymmetric Macdonald polynomial, by \eqref{f-ansatz2}. Our aim is to show that all terms in the sum \eqref{sum} are related to the leading one, $(d_1,\dots,d_n) = (1,\dots,p,0^{n-p})$, under the action of Hecke generators \eqref{Hecke}. To that end, we will use that for any polynomial $g \in \mathbb{C} [\textbf{x}]$ we have the identity
	\begin{align}
	\label{daha-relation}
	T^{-1}_i \cdot \Big( x_i g(\textbf{x}) \Big)
	=
	t^{-1}
	x_{i+1}
	T_i
	\cdot
	g(\textbf{x}),
	\end{align}

	\noindent which follows from \eqref{Hecke}. Using the commutativity of the Hecke generator $T^{-1}_i$ with the multiplication by $x_i x_{i+1}$ (which follows from \eqref{Hecke}, since $[\mathfrak{s}_i, x_i x_{i + 1}] = 0$) and the exchange relation \eqref{T-CC}, we see that for $1\le d_i < d_{i+1}$ we have
	\begin{multline*}
	T^{-1}_i
	\cdot
	t^{{\rm inv}' (\overleftarrow{\mathfrak{D}})}
	\prod_{k=1}^{n}
	x_k^{\bm{1}_{d_k \ge 1}}
	\cdot
	\big\langle
	\mathsf{C}_{d_n}(x_n)
	\cdots
	\mathsf{C}_{d_1}(x_1)
	\big\rangle_{\nu}'(\textbf{w})
	\\
	=
	t^{{\rm inv}' (\overleftarrow{\mathfrak{D}})}
	\prod_{k=1}^{n}
	x_k^{\bm{1}_{d_k \ge 1}}
	T^{-1}_i
	\cdot
	\big\langle
	\mathsf{C}_{d_n}(x_n)
	\cdots
	\mathsf{C}_{d_1}(x_1)
	\big\rangle_{\nu}'(\textbf{w})
	\\
	=
	t^{{\rm inv}' (\overleftarrow{\mathfrak{D}'})}
	\prod_{k=1}^{n}
	x_k^{\bm{1}_{d'_k \ge 1}}
	\big\langle
	\mathsf{C}_{d'_n}(x_n)
	\cdots
	\mathsf{C}_{d'_1}(x_1)
	\big\rangle_{\nu}'(\textbf{w}),
	\end{multline*}
	where $\mathfrak{D}' = (d'_1,\dots,d'_n) = \mathfrak{s}_i (\mathfrak{D}) = (d_1, \ldots , d_{i - 1}, d_{i + 1}, d_i, d_{i + 2}, \ldots , d_n)$. Similarly, using \eqref{daha-relation} and \eqref{T-CC} (in which we act on both sides by $T_i$), we find that for $d_i \ge 1$ and $d_{i+1}=0$ there holds
	\begin{multline*}
	T^{-1}_i
	\cdot
	t^{{\rm inv}' (\overleftarrow{\mathfrak{D}})}
	\prod_{k=1}^{n}
	x_k^{\bm{1}_{d_k \ge 1}}
	\cdot
	\big\langle
	\mathsf{C}_{d_n}(x_n)
	\cdots
	\mathsf{C}_{d_1}(x_1)
	\big\rangle_{\nu}'(\textbf{w})
	\\
	=
	t^{{\rm inv}' (\overleftarrow{\mathfrak{D}})-1}
	\prod_{k=1}^{n}
	x_k^{\bm{1}_{d'_k \ge 1}}
	T_i
	\cdot
	\big\langle
	\mathsf{C}_{d_n}(x_n)
	\cdots
	\mathsf{C}_{d_1}(x_1)
	\big\rangle_{\nu}'(\textbf{w})
	\\
	=
	t^{{\rm inv}' (\overleftarrow{\mathfrak{D}'})}
	\prod_{k=1}^{n}
	x_k^{\bm{1}_{d'_k \ge 1}}
	\big\langle
	\mathsf{C}_{d'_n}(x_n)
	\cdots
	\mathsf{C}_{d'_1}(x_1)
	\big\rangle_{\nu}'(\textbf{w}),
	\end{multline*}
	where again $\mathfrak{D}' = (d'_1,\dots,d'_n) = \mathfrak{s}_i (\mathfrak{D})$. Letting $\mu = (\nu, 0^{n - p})$, these two exchange relations and \Cref{pdc} together yield 
	\begin{align}
	\label{nuptsigma}
	\begin{aligned}
	\mathcal{P}_{\nu}(\textbf{x})
	&=
	(1-t)^p
	\left(
	\sum_{\sigma \in \mathfrak{S}_n}
	T^{-1}_{\sigma}
	\right)
	\cdot
	t^{\binom{p}{2}}
	\prod_{i=1}^p
	x_i
	\cdot
	\big\langle
	\mathsf{D} (x_n)
	\cdots
	\mathsf{D}(x_{p + 1})
	\mathsf{C}_p(x_p)
	\cdots
	\mathsf{C}_1(x_1)
	\big\rangle_{\nu}'(\textbf{w})
	\\
	&=
	\frac{(1-t)^p t^{\binom{p}{2}}}{\Omega_{\nu}(q,t)} \displaystyle\prod_{i = 1}^p (1 - t^i w_{i, 0})^{-1}
	\left(
	\sum_{\sigma \in \mathfrak{S}_n}
	T^{-1}_{\sigma}
	\right)
	\cdot
	f_{\mu}(\textbf{x}),
	\end{aligned} 
	\end{align}
	\noindent where we recall the $w_{i, j}$ from \eqref{w_ij}, and for any permutation $\sigma \in \mathfrak{S}_n$ with reduced word form $\sigma = \mathfrak{s}_{i_1} \cdots \mathfrak{s}_{i_{\ell}}$ we have set $T_{\sigma} = T_{i_1} \cdots T_{i_{\ell}}$ (which is well-defined by the last relation in \eqref{hecke1}). Here, the factor of $t^{\binom{p}{2}}$ follows from the fact that $\inv' (\overleftarrow{\mathfrak{D}'}) = \binom{p}{2}$ if $\mathfrak{D}' = (1, 2, \ldots , p, 0^{n - p})$.
	
	Now $f_{\mu}(\textbf{x}) \in \mathcal{V}_{\kappa}$, with $\kappa = (\nu^+, 0^{n - p})$. Moreover, it quickly follows from \eqref{hecke1} and the recursive properties of nonsymmetric Macdonald polynomials with respect to the Hecke generators $T_i^{-1}$ (see Theorem 4.2 of \cite{IVN}) that $T_{\sigma}^{-1} h \in \mathcal{V}_{\kappa}$ for any $h \in \mathcal{V}_{\kappa}$ and $\sigma \in \mathfrak{S}_n$. This, together with \eqref{nuptsigma}, implies the theorem. 
\end{proof}

\begin{cor}
	
	\label{pdp} 
	
	For any anti-dominant composition $\nu = (\nu_1,\dots,\nu_p)$ such that $\nu_i \ge 1$ for each $i \in [1, p]$, let $\kappa = (\nu^{+},0^{n-p})$. Up to an (at this stage) unspecified multiplicative constant, we have
	\begin{align*}
	\mathcal{P}_{\nu}(\textbf{\emph{x}})
	=
	d_{\nu}(q,t)
	P_{\kappa}(\textbf{\emph{x}};q,t),
	\end{align*}
	where $P_{\kappa}(\textbf{\emph{x}};q,t)$ denotes a symmetric Macdonald polynomial.
\end{cor}

\begin{proof}
	This follows from the fact that, up to an overall multiplicative constant, the only polynomial which satisfies the properties of \Cref{prop-sym} and \Cref{prop-V} is $P_{\kappa}(\textbf{x};q,t)$.
\end{proof}

\section{Identifying the Leading Coefficient}

\label{ProofJ}

In this section we establish \Cref{jpd}, which will follow from \Cref{pdp}, together with an analysis of the leading coefficient of $\mathcal{P}_{\nu} (\textbf{x})$.

\begin{proof}[Proof of \Cref{jpd}] 
	
	We must show that 
	\begin{align}
	\label{J}
	J_{\nu^{+}}(\textbf{x})
	=
	\sum_{\mathscr{M}}
	\prod_{i=1}^{p}
	\prod_{j=1}^{N}
	w_{i,j}^{M_{i,j}}
	\times
	\tikz{1}{
		\foreach\y in {0,...,5}{
			\draw[lgray,line width=1.5pt] (1.5,0.5+\y) -- (7.5,0.5+\y);
		}
		\foreach\x in {0,...,6}{
			\draw[lgray,line width=1.5pt] (1.5+\x,0.5) -- (1.5+\x,5.5);
		}
		%spectral parameters
		\node[text centered] at (2,1) {$x_1$};
		\node[text centered] at (4,1) {$\cdots$};
		\node[text centered] at (5,1) {$\cdots$};
		\node[text centered] at (7,1) {$x_1$};
		\node[text centered] at (2,2) {$x_2$};
		\node[text centered] at (4,2) {$\cdots$};
		\node[text centered] at (5,2) {$\cdots$};
		\node[text centered] at (7,2) {$x_2$};
		\node[text centered] at (2,5) {$x_n$};
		\node[text centered] at (4,5) {$\cdots$};
		\node[text centered] at (5,5) {$\cdots$};
		\node[text centered] at (7,5) {$x_n$};
		\node[text centered] at (2,3.1) {$\vdots$};
		\node[text centered] at (2,4.1) {$\vdots$};
		\node[text centered] at (7,3.1) {$\vdots$};
		\node[text centered] at (7,4.1) {$\vdots$};
		%bottom labels
		\node[below] at (2,0.5) {\scriptsize$\textbf{e}_{[1, p]}$};
		\node[below] at (3,0.5) {\scriptsize$\textbf{M}_1$};
		\node[above,text centered] at (4,0) {$\cdots$};
		\node[above,text centered] at (5,0) {$\cdots$};
		\node[below] at (7,0.5) {\scriptsize$\textbf{M}_N$};
		%top labels
		\node[left = 5, above] at (2,5.5) {\scriptsize $\textbf{e}_0$};
		\node[above] at (3,5.5) {\scriptsize$\textbf{M}_1+\textbf{I}_1$};
		\node[above,text centered] at (4,5.5) {$\cdots$};
		\node[above,text centered] at (5,5.5) {$\cdots$};
		\node[above] at (7,5.5) {\scriptsize$\textbf{M}_N+\textbf{I}_N$};
		%right labels
		\node[right] at (7.5,1) {$0$};
		\node[right] at (7.5,2) {$0$};
		\node[right] at (7.5,3.1) {$\vdots$};
		\node[right] at (7.5,4.1) {$\vdots$};
		\node[right] at (7.5,5) {$0$};
		%left labels
		\node[left] at (1.5,1) {$0$};
		\node[left] at (1.5,2) {$0$};
		\node[left] at (1.5,3.1) {$\vdots$};
		\node[left] at (1.5,4.1) {$\vdots$};
		\node[left] at (1.5,5) {$0$};
	}
	\end{align}
	where we abbreviate $\textbf{I}_j = \textbf{I}_j (\nu)$ for each $j \ge 1$, and for each $(i, j) \in [1, p] \times [1, N]$ the twist parameter $w_{i, j}$ is given by $w_{i,j} =\bm{1}_{\nu_i > j} q^{\nu_i - j}$ (recall \eqref{w_ij}).
	 
	  As \Cref{prop-V}, we proved that \eqref{J} gives a symmetric Macdonald polynomial, modulo normalization. So, it remains to verify that the coefficient of the $\prod_{i=1}^{p} x_i^{\nu^{+}_i} =  \prod_{i=1}^{p} x_i^{\nu_{p-i+1}}$ monomial is $c_{\nu^+} (q, t)$ from \eqref{cqt}, in order to verify the normalization. 
	  
	  Observe that any path ensemble in the above partition function contributing to the coefficient of this monomial $\prod_{i=1}^{p} x_i^{\nu_{p-i+1}}$ cannot have any path take a horizontal step in some row $k > p$. Otherwise, \Cref{rem:x-depend} implies that such a step would contribute a factor of $x_k$ to the weight of this ensemble, and so it cannot give rise to the monomial $\prod_{i=1}^{p} x_i^{\nu_{p-i+1}}$. In particular, in any such path ensemble, each path exits the leftmost column of the model in one of the rows $1, 2, \ldots , p$. 
	  
	Given this, the procedure for analyzing the coefficient of $\prod_{i=1}^{p} x_i^{\nu_{p-i+1}}$ closely parallels the argument in \Cref{ssec:norm}. In particular, one can show (we omit the proof, since it is entirely analogous to that of \Cref{coefficientxc}) that there is an {\it almost} unique configuration in the partition function \eqref{J} which gives rise to the monomial  $\prod_{i=1}^{p} x_i^{\nu_{p-i+1}}$; these occur when the path of color $i$ makes a right turn out of column zero into row $p-i'+1$ for some $i'$ with $\nu_{i'} = \nu_i$, and then takes $\nu_i$ consecutive horizontal steps in that row, before making a vertical step into its destination column. The ``almost'' part of the claim comes parts of equal size in $\nu$; if a part is repeated $k$ times in $\nu$, then there are $k!$ possibilities for the corresponding colors to permute themselves in any way, before making the right turn out of the zeroth column. We thus need to compute the weight of configurations of the form
	
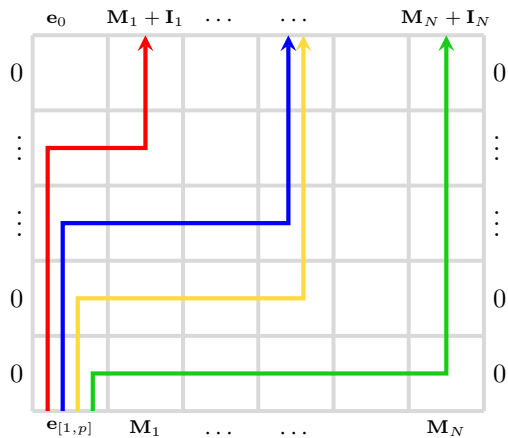
\begin{align*}
	\tikz{1}{
		\foreach\y in {0,...,5}{
			\draw[lgray,line width=1.5pt] (1.5,0.5+\y) -- (7.5,0.5+\y);
		}
		\foreach\x in {0,...,6}{
			\draw[lgray,line width=1.5pt] (1.5+\x,0.5) -- (1.5+\x,5.5);
		}
		%bottom labels
		\node[below] at (2,0.5) {\scriptsize$\textbf{e}_{[1, p]}$};
		\node[below] at (3,0.5) {\scriptsize$\textbf{M}_1$};
		\node[above,text centered] at (4,0) {\small$\cdots$};
		\node[above,text centered] at (5,0) {\small $\cdots$};
		\node[below] at (7,0.5) {\scriptsize$\textbf{M}_N$};
		%top labels
		\node[left = 5, above] at (2,5.5) {\scriptsize $\textbf{e}_0$};
		\node[above] at (3,5.5) {\scriptsize$\textbf{M}_1+\textbf{I}_1$};
		\node[above,text centered] at (4,5.5) {\small$\cdots$};
		\node[above,text centered] at (5,5.5) {\small$\cdots$};
		\node[above] at (7,5.5) {\scriptsize$\textbf{M}_N+\textbf{I}_N$};
		%right labels
		\node[right] at (7.5,1) {$0$};
		\node[right] at (7.5,2) {$0$};
		\node[right] at (7.5,3.1) {$\vdots$};
		\node[right] at (7.5,4.1) {$\vdots$};
		\node[right] at (7.5,5) {$0$};
		%left labels
		\node[left] at (1.5,1) {$0$};
		\node[left] at (1.5,2) {$0$};
		\node[left] at (1.5,3.1) {$\vdots$};
		\node[left] at (1.5,4.1) {$\vdots$};
		\node[left] at (1.5,5) {$0$};
		%paths
		\draw[ultra thick,red,->] (1.7,0.5) -- (1.7,4) -- (2.5,4) -- (3,4) -- (3,5.5);
		\draw[ultra thick,blue,->] (1.9,0.5) -- (1.9,3) -- (2.5,3) -- (4.9,3) -- (4.9,5.5);
		\draw[ultra thick,yellow,->] (2.1,0.5) -- (2.1,2) -- (2.5,2) -- (5.1,2) -- (5.1,5.5);
		\draw[ultra thick,green,->] (2.3,0.5) -- (2.3,1) -- (7,1) -- (7,5.5); 
	}
	\end{align*}
	and perform the aforementioned summation over permutations of colors in the zeroth column, for parts of repeated size. 
	
	Denoting the colors of the arrows exiting the leftmost column, from bottom to top, by $\mathfrak{D} = (d_1, d_2, \ldots , d_n)$, we find that the weight of this column is given by $(1 - t)^p t^{\inv' (\overleftarrow{\mathfrak{D}})} \prod_{i = 1}^p x_i$, where (as in the proof of \Cref{prop-V}) $\inv' (\overleftarrow{\mathfrak{D}})$ denotes the number of pairs $1 \le i < j \le n$ such that $1 \le d_i < d_j$. Indeed, the factor of $(1 - t)^p \prod_{i = 1}^p x_i$ arises from to the factors of $(1 - t) x$ in the second weight in \Cref{vertexfigurerql}, and the $t^{\inv' (\overleftarrow{\mathfrak{D}})}$ factor arises as the product of the $t^{A_{[d_i + 1, n]}}$ exponents in these weights.
	
	Moreover, consulting equation \eqref{components} for the column weights, the contribution to this configuration weight coming from the remaining $N$ columns is equal to
	\begin{align*}
	\prod_{i=1}^{p} x_i^{\nu_{p-i+1} - 1}
	\cdot
	\prod_{j=1}^{N}
	\prod_{i: \nu_i > j}
	(1-w_{i,j} t^{\#\{k<i : \nu_k > j\}+1}).
	\end{align*}
	
	\noindent In particular, this quantity is independent of $\mathfrak{D}$ satisfying $\nu_{d_1} \ge \nu_{d_2} \ge \cdots \ge \nu_{d_p}$ and $d_{p + 1} = d_{p + 2} = \cdots = d_n = 0$, such that $\nu_{d_{K_i + 1}} = \nu_{d_{K_i + 2}} = \cdots = \nu_{d_{K_{i + 1}}}$, where $K_i = \sum_{j = 1}^{i - 1} m_j (\nu)$, for each $i \in [1, p]$.
	
	Taking the product of the contributions from the leftmost column and the remaining columns; summing over $\mathfrak{D}$ satisfying the three properties listed above; and using the fact that $\sum_{\sigma \in \mathfrak{S}_k} t^{\inv(\sigma)} = (1 - t)^{-k} (t; t)_k$ for any $k \ge 0$ yields that the coefficient of $\prod_{i = 1}^p x_i^{\nu_{p - i + 1}}$ in $\mathcal{P}_{\nu} (\textbf{x})$ is equal to
	\begin{flalign}
	\label{coefficientp}
	\displaystyle\prod_{i = 1}^p (t; t)_{m_i (\nu)} \prod_{j=1}^{N}
	\prod_{i: \nu_i > j}
	(1-w_{i,j} t^{\#\{k<i : \nu_k > j\}+1}).
	\end{flalign}
	
	Further	observe that, for any partition $\lambda = (\lambda_1, \lambda_2, \ldots , \lambda_p)$, we may write the quantity $c_{\lambda}(q,t)$ from \eqref{cqt} as
	\begin{flalign}
	\label{clambdaqt2} 
	c_{\lambda}(q,t)
	=
	\prod_{i = 1}^p
	\prod_{j=1}^{\lambda_i-1}
	(1-q^{\lambda_i-j} t^{\#\{ k > i : \lambda_k \ge j \}+1})
	\cdot
	\prod_{i = 1}^p
	(1- t^{\#\{k <i : \lambda_k = \lambda_i\}+1}),
	\end{flalign}
	where the first product comes from analyzing all boxes in \eqref{cqt} which have positive arm length, while the second product deals with all boxes that have arm length equal to zero. Now the theorem follows from the fact that \eqref{coefficientp} and \eqref{clambdaqt2} coincide for $\lambda = \nu^+$.
\end{proof}

\chapter{Expansion of LLT Polynomials in the Modified Hall--Littlewood Basis} 

\label{Expansion}

In this chapter we provide a combinatorial formula, as partition functions for a fused $U_q \big( \widehat{\mathfrak{sl}} (2 | n) \big)$ vertex model, for the expansion coefficients of the LLT polynomials in the modified Hall--Littlewood basis. This result is given by \Cref{thm:comb}, and its proof will appear in \Cref{sec:proof} below. Throughout this chapter, we fix an integer $n \ge 1$. Moreover, for any (possibly infinite) set of variables $\textbf{x} = (x_1, x_2, \ldots , x_k)$; complex number $r \in \mathbb{C}$; signature $\lambda \in \Sign_{\ell}$; and signature sequences $\boldsymbol{\lambda}, \boldsymbol{\mu} \in \SeqSign_n$, we recall the LLT polynomials $\mathcal{L}_{\boldsymbol{\lambda} / \boldsymbol{\mu}} (\textbf{x})$ from \eqref{llambdamu} and \eqref{lambdamug}; the (standard) Hall--Littlewood polynomials $Q_{\lambda} (\textbf{x})$ from Section 3.2 of \cite{SFP}; the modified Hall--Littlewood polynomials $Q_{\lambda}' (\textbf{x})$ from \Cref{StabilityPolynomials}; and the functions $\mathcal{G}_{\boldsymbol{\lambda} / \boldsymbol{\mu}} (\textbf{x}; r \boldsymbol{\mid} 0; 0)$ from \eqref{limitg}. We additionally recall the sequences $\mathscr{S} (\boldsymbol{\lambda}) = \big( \textbf{S}_1 (\boldsymbol{\lambda}), \textbf{S}_2 (\boldsymbol{\lambda}), \ldots \big)$ and $\mathscr{I} (\lambda) = \big( \textbf{I}_0 (\lambda), \textbf{I}_1 (\lambda), \ldots \big)$ of elements in $\{ 0, 1 \}^n$ from \Cref{Symmetric} and \Cref{ssec:states}, respectively (assuming in the latter that $\ell (\lambda) \le n$).

\section{Combinatorial Formula for Expansion Coefficients}

\label{sec:comb}

As explained in \Cref{r0g}, the polynomials $\mathcal{G}_{\boldsymbol{\lambda} / \boldsymbol{\mu}} (\textbf{x}; \infty \boldsymbol{\mid} 0; 0)$ and $\mathcal{L}_{\boldsymbol{\lambda} / \boldsymbol{\mu}} (\textbf{x})$ (and, by \Cref{qlambdaqlambda}, also $Q_{\lambda}' (\textbf{x})$) satisfy the compatibility condition of \Cref{g0} that enables us to view them as elements in the ring $\Lambda (\textbf{x})$ of symmetric functions in an infinite set of variables $\textbf{x} = (x_1, x_2, \ldots )$. Since (by, for example, Chapter 3.2 of \cite{SFP}) the (modified) Hall--Littlewood functions form a basis of $\Lambda (\textbf{x})$, there exist coefficients $f_{\boldsymbol{\lambda}/\boldsymbol{\mu}}^{\nu}(q) \in \mathbb{C} (q)$ such that
\index{F@$f_{\boldsymbol{\lambda} / \boldsymbol{\mu}}^{\nu} (q)$}
\begin{align}
	\label{L-into-Q}
	\mathcal{G}_{\boldsymbol{\lambda}/\boldsymbol{\mu}} (\textbf{x}; \infty \boldsymbol{\mid} 0; 0)
	=
	\sum_{\nu}
	f_{\boldsymbol{\lambda}/\boldsymbol{\mu}}^{\nu}(q) 
	Q'_{\nu} (\textbf{x}),
\end{align}
with the sum taken over all partitions $\nu$. Recalling the function $\psi$ from \eqref{lambdamupsi}, \eqref{1gl} implies that $\mathcal{G}_{\boldsymbol{\lambda}/\boldsymbol{\mu}} (\textbf{x}; \infty \boldsymbol{\mid} 0; 0) = q^{\psi (\boldsymbol{\lambda}) - \psi (\boldsymbol{\mu})} \mathcal{L}_{\boldsymbol{\lambda} / \boldsymbol{\mu}} (\textbf{x})$, and so (up to an explicit power of $q$) the $f_{\boldsymbol{\lambda} / \boldsymbol{\mu}}^{\nu} (q)$ provide the expansion coefficients of the LLT polynomials in the modified Hall--Littlewood basis.

The main result to be stated in this section is an explicit combinatorial formula for these coefficients $f_{\boldsymbol{\lambda}/\boldsymbol{\mu}}^{\nu}(q)$ for all $\boldsymbol{\lambda}, \boldsymbol{\mu} \in \SeqSign_{n; M}$ and partitions $\nu$; see \Cref{thm:comb} below. In order to state our main formula, we introduce two types of vertex weights; their origin will become clear in \Cref{sec:proof} below.

\begin{definition}
	
\label{ssec:weights1}

Fix nonnegative integers $\mathfrak{a}, \mathfrak{b}, \mathfrak{c}, \mathfrak{d} \ge 0$ and $n$-tuples $\textbf{A}, \textbf{B}, \textbf{C}, \textbf{D} \in \{ 0, 1 \}^n$. Define $\textbf{V} = (V_1,\dots,V_n) \in \{ 0, 1 \}^n$ with $V_i = \min\{A_i,B_i,C_i,D_i\}$ for each $i \in [1, n]$. We introduce the lattice weights
\begin{align}
	\begin{aligned}
	\label{weights-light}
	\tikz{0.6}{
		\draw[lgray,line width=1.5pt] (-1,-1) -- (1,-1) -- (1,1) -- (-1,1) -- (-1,-1);
		\node[left] at (-1,0) {\tiny $(\mathfrak{b},\textbf{B})$};\node[right] at (1,0) {\tiny $(\mathfrak{d},\textbf{D})$};
		\node[below] at (0,-1) {\tiny $(\mathfrak{a},\textbf{A})$};\node[above] at (0,1) {\tiny $(\mathfrak{c},\textbf{C})$};
	}
	& =
	(-1)^{\mathfrak{c}+|\textbf{V}|} q^{\chi}
	\frac{(q;q)_{\mathfrak{b}-|\textbf{B}|}}{(q;q)_{\mathfrak{d}-|\textbf{D}|}} 
	\frac{(q^{\mathfrak{b}-|\textbf{B}|-\mathfrak{a}+1};q)_\mathfrak{c}}{(q;q)_\mathfrak{c}}
	{\bm 1}_{|\textbf{B}| \le \mathfrak{b}}
	{\bm 1}_{|\textbf{D}| \le \mathfrak{d}}
	{\bm 1}_{\mathfrak{a} + \mathfrak{d} = \mathfrak{b} +  \mathfrak{c}}
	\\
	& \quad \times
	{\bm 1}_{\textbf{A}+\textbf{B}=\textbf{C}+\textbf{D}}
	\prod_{j:B_j - D_j=1}
	\big( 1-q^{\mathfrak{d}-B_{[j+1,n]}-D_{[1, j - 1]}} \big),
	\end{aligned} 
\end{align}
where the exponent $\chi \equiv \chi(\mathfrak{a}, \mathfrak{b}, \mathfrak{c}, \mathfrak{d};\textbf{A},\textbf{B},\textbf{C},\textbf{D})$ is given by
\begin{align*}
	\chi
	=
	\binom{\mathfrak{d}-|\textbf{D}|}{2}
	+
	\binom{\mathfrak{c}+1}{2}
	-
	\big( \mathfrak{c}+|\textbf{C}| \big)\mathfrak{d}
	+
	|\textbf{V}| \big( \mathfrak{d}-|\textbf{D}|+1 \big)
	+
	\varphi(\textbf{D},\textbf{C})+\varphi(\textbf{V},\textbf{D}-\textbf{B}).
\end{align*}

\end{definition} 

Aside from the vanishing constraints imposed by the indicator functions in \eqref{weights-light}, these vertex weights have additional vanishing properties that result from the $q$-Pochhammer function $(q^{\mathfrak{b}-|\textbf{B}|-\mathfrak{a}+1};q)_\mathfrak{c}$ and the factor $\big( 1 - q^{\mathfrak{d}-B_{[j+1,n]}-D_{[1,j-1]}} \big)$ present in the final product. 

We will require a further set of weights, whose notation are distinguished from those in \eqref{weights-light} by means of shading.

\begin{definition}
\label{ssec:weights2}

As in \Cref{ssec:weights1}, let us fix nonnegative integers $\mathfrak{a}, \mathfrak{b}, \mathfrak{c}, \mathfrak{d} \ge 0$ and $n$-tuples $\textbf{A}, \textbf{B}, \textbf{C}, \textbf{D} \in \{ 0, 1 \}^n$ with coordinates indexed by $[1, n]$. We define the weights
\begin{align}
	\label{weights-dark}
	\tikz{0.6}{
		\filldraw[lgray,line width=1.5pt,fill=llgray] (-1,-1) -- (1,-1) -- (1,1) -- (-1,1) -- (-1,-1);
		\node[left] at (-1,0) {\tiny $(\mathfrak{b},\textbf{B})$};\node[right] at (1,0) {\tiny $(\mathfrak{d},\textbf{D})$};
		\node[below] at (0,-1) {\tiny $(\mathfrak{a},\textbf{A})$};\node[above] at (0,1) {\tiny $(\mathfrak{c},\textbf{C})$};
	}
	=
	{\bm 1}_{|\textbf{B}| = \mathfrak{b} \le 1}
	{\bm 1}_{|\textbf{D}| \le \mathfrak{d} \le 1}
	{\bm 1}_{\mathfrak{a} + \mathfrak{d} = \mathfrak{b} +  \mathfrak{c}}
	{\bm 1}_{\textbf{A}+\textbf{B}=\textbf{C}+\textbf{D}}
	\cdot
	W \big(\mathfrak{a},\mathfrak{b},\mathfrak{c},\mathfrak{d};\textbf{A},\textbf{B},\textbf{C},\textbf{D}\big),
\end{align}
where (recalling $\textbf{A}_i^+, \textbf{A}_j^-, \textbf{A}_{ij}^{+-}$ from \eqref{aij}) the final function appearing in \eqref{weights-dark} is given by
\begin{flalign}
	\begin{aligned}
	\label{W-def}
	W & \big( \mathfrak{a}, \mathfrak{b}, \mathfrak{c}, \mathfrak{d};\textbf{A},\textbf{B},\textbf{C},\textbf{D} \big) =
	\\
	& \left\{
	\begin{array}{ll}
		q^{|\textbf{A}|},
		&
		\quad {\text{if}} \quad
		(\mathfrak{b},\textbf{B}) = (0,\textbf{e}_0),
		\quad
		(\mathfrak{c},\textbf{C}) = (\mathfrak{a},\textbf{A}),
		\quad
		(\mathfrak{d},\textbf{D}) = (0,\textbf{e}_0),
		\\
		\\
		1,
		&
		\quad {\text{if}} \quad
		(\mathfrak{b},\textbf{B}) = (0,\textbf{e}_0),
		\quad
		(\mathfrak{c},\textbf{C}) = (\mathfrak{a}+1,\textbf{A}),
		\quad
		(\mathfrak{d},\textbf{D}) = (1,\textbf{e}_0),
		\\
		\\
		(1-q^{A_i})q^{A_{[i+1,n]}},
		&
		\quad {\text{if}} \quad
		(\mathfrak{b},\textbf{B}) = (0,\textbf{e}_0),
		\quad
		(\mathfrak{c},\textbf{C}) = (\mathfrak{a}+1,\textbf{A}^{-}_i),
		\quad
		(\mathfrak{d},\textbf{D}) = (1,\textbf{e}_i),
		\\
		\\
		(1 - q^{A_j}) q^{A_{[j + 1, n]}} \textbf{1}_{A_i = 0}
		&
		\quad {\text{if}} \quad
		(\mathfrak{b},\textbf{B}) = (1, \textbf{e}_i),
		\quad
		(\mathfrak{c},\textbf{C}) = (\mathfrak{a},\textbf{A}^{+-}_{ij}),
		\quad
		(\mathfrak{d},\textbf{D}) = (1,\textbf{e}_j),
		\\
		\\
		\bm{1}_{A_i=0},
		&
		\quad {\text{if}} \quad
		(\mathfrak{b},\textbf{B}) = (1,\textbf{e}_i),
		\quad
		(\mathfrak{c},\textbf{C}) = (\mathfrak{a},\textbf{A}^{+}_i),
		\quad
		\quad
		(\mathfrak{d},\textbf{D}) = (1,\textbf{e}_0),
		\\
		\\
		(-1)^{A_i}
		q^{A_{[i, n]}},
		&
		\quad {\text{if}} \quad
		(\mathfrak{b},\textbf{B}) = (1,\textbf{e}_i),
		\quad
		(\mathfrak{c},\textbf{C}) = (\mathfrak{a},\textbf{A}),
		\quad
		(\mathfrak{d},\textbf{D}) = (1,\textbf{e}_i),
	\end{array}
	\right.
	\end{aligned} 
\end{flalign}

\noindent where we assume that $i < j$ whenever $i$ and $j$ both appear. In all cases of $(\mathfrak{a},\mathfrak{b},\mathfrak{c},\mathfrak{d};\textbf{A},\textbf{B},\textbf{C},\textbf{D})$ not listed in \eqref{W-def}, we set $W \big(\mathfrak{a},\mathfrak{b},\mathfrak{c},\mathfrak{d};\textbf{A},\textbf{B},\textbf{C},\textbf{D}\big) = 0$.

\end{definition}

Observe that there are two quadruples we consider in \Cref{ssec:weights1} and \Cref{ssec:weights2}. The first is the quadruple $(\textbf{A}, \textbf{B}; \textbf{C}, \textbf{D})$ of elements in $\{ 0, 1 \}^n$, which is stipulated to satisfy $\textbf{A} + \textbf{B} = \textbf{C} + \textbf{D}$. We once again view $\textbf{A}$, $\textbf{B}$, $\textbf{C}$, and $\textbf{D}$ as indexing the up-right directed fermionic paths passing through the south, west, north, and east boundaries of a tile, respectively. The second is the quadruple $(\mathfrak{a}, \mathfrak{b}; \mathfrak{c}, \mathfrak{d})$ of integers, which is stipulated to satisfy $\mathfrak{a} + \mathfrak{d} = \mathfrak{b} +  \mathfrak{c}$. In this way, we may view $\mathfrak{a}$, $\mathfrak{b}$, $\mathfrak{c}$, and $\mathfrak{d}$ as counting the numbers of directed down-right paths (that is, they proceed in a different direction from the fermionic paths) of the same bosonic color that pass through the bottom, left, top, and right boundaries of a tile, respectively. We provide depictions of these tiles in \Cref{Models2n} below. 

Now we present the main result of this section, \Cref{thm:comb}; its proof will appear in \Cref{sec:proof} below.

\begin{thm}
	\label{thm:comb}
	
	Fix integers $M \ge 0$ and $m \in [1, n]$; signature sequences $\boldsymbol{\lambda}, \boldsymbol{\mu} \in \SeqSign_{n; M}$; and a partition $\nu \in \Sign_{m}$, such that $|\boldsymbol{\lambda}|-|\boldsymbol{\mu}| = |\nu|$. Fix any integer 
	\begin{align*}
		K \ge \max \Big\{ \max_{i \in [1, n]} \mathfrak{T} \big(\lambda^{(i)} \big),\nu_1 + 1 \Big\},
	\end{align*}
	and write $\bar{\nu}=(\bar{\nu}_1, \bar{\nu}_2, \cdots , \bar{\nu}_{m})$ for the signature obtained by complementing $\nu$ in a $(K - 1) \times m$ box; its parts are given by $\bar{\nu}_j = K - \nu_{m - j +1} - 1$, for each $j \in [1, m]$. Denote $\mathfrak{T} (\bar{\nu}) = (\mathfrak{n}_1, \mathfrak{n}_2, \ldots , \mathfrak{n}_{m})$ so, for each $j \in [1, m]$, 
	\begin{align}
		\label{maya-seq1}
		\mathfrak{n}_j = \bar{\nu}_j + m - j + 1 = K - \nu_{m - j + 1} + m - j.
	\end{align}
	The coefficients \eqref{L-into-Q} are given by the partition function
	\begin{align}
		\label{main-formula}
		f^{\nu}_{\boldsymbol{\lambda}/\boldsymbol{\mu}}(q)
		=
		\frac{(-1)^{|\bar{\nu}|}}{b_{\nu}(q)(q;q)_{m}}
		\times
		\tikz{1.1}{
			\filldraw[lgray,line width=1.5pt,fill=llgray] (1.5,3.5) -- (7.5,3.5) -- (7.5,4.5) -- (1.5,4.5) -- (1.5,3.5);
			\filldraw[lgray,line width=1.5pt,fill=llgray] (1.5,0.5) -- (7.5,0.5) -- (7.5,1.5) -- (1.5,1.5) -- (1.5,0.5);
			\foreach\y in {0,...,6}{
				\draw[lgray,line width=1.5pt] (1.5,0.5+\y) -- (7.5,0.5+\y);
			}
			\foreach\x in {0,...,6}{
				\draw[lgray,line width=1.5pt] (1.5+\x,0.5) -- (1.5+\x,6.5);
			}
			%bottom labels
			\node[below] at (2,0.5) {\footnotesize $\big( 0, \textbf{\emph{S}}_1 (\boldsymbol{\mu}) \big)$};
			\node[above,text centered] at (4,0) {$\cdots$};
			\node[above,text centered] at (5,0) {$\cdots$};
			\node[below] at (7,0.5) {\footnotesize $\big(0, \textbf{\emph{S}}_K (\boldsymbol{\mu}) \big)$};
			%top labels
			\node[above] at (2,6.5) {\footnotesize $\big(0,\textbf{\emph{S}}_1 (\boldsymbol{\lambda}) \big)$};
			\node[above,text centered] at (4,6.5) {$\cdots$};
			\node[above,text centered] at (5,6.5) {$\cdots$};
			\node[above] at (7,6.5) {\footnotesize $\big( 0, \textbf{\emph{S}}_K (\boldsymbol{\lambda}) \big)$};
			%right labels
			\node[right] at (7.5,1) {\footnotesize $(1,\textbf{\emph{e}}_0)$};
			\node[right] at (7.5,2) {\footnotesize $(0,\textbf{\emph{e}}_0)$};
			\node[right] at (7.5,3.1) {$\vdots$};
			\node[right] at (7.5,4) {\footnotesize $(1,\textbf{\emph{e}}_0)$};
			\node[right] at (7.5,5) {$\vdots$};
			\node[right] at (7.5,6) {\footnotesize $(0,\textbf{\emph{e}}_0)$};
			%left labels
			\node[left] at (1.5,1) {\footnotesize $(0,\textbf{\emph{e}}_0)$};
			\node[left] at (1.5,2) {\footnotesize $(0,\textbf{\emph{e}}_0)$};
			\node[left] at (1.5,3.1) {$\vdots$};
			\node[left] at (1.5,4) {\footnotesize $(0,\textbf{\emph{e}}_0)$};
			\node[left] at (1.5,5) {$\vdots$};
			\node[left] at (1.5,6) {\footnotesize $(m,\textbf{\emph{e}}_0)$};
		}
	\end{align}

	 \noindent consisting of $\mathfrak{n}_1 + 1$ rows, where the $i$-th row of the lattice (counted from the top to bottom, starting at $i=0$) takes the form
	\begin{align}
		\label{maya-seq2}
		\left\{
		\begin{array}{ll}
			\tikz{0.4}{
				\filldraw[lgray,line width=1.5pt,fill=llgray] (-1,-1) -- (7,-1) -- (7,1) -- (-1,1) -- (-1,-1);
				\draw[lgray,line width=1.5pt] (1,-1) -- (1,1);
				\draw[lgray,line width=1.5pt] (3,-1) -- (3,1);
				\draw[lgray,line width=1.5pt] (5,-1) -- (5,1);
				\node[left] at (-1,0) {\tiny $(0,\textbf{\emph{e}}_0)$};\node[right] at (7,0) {\tiny $(1,\textbf{\emph{e}}_0)$};
			}
			\quad\quad & i \in \mathfrak{T} (\bar{\nu}),
			\\ \\
			\tikz{0.4}{
				\draw[lgray,line width=1.5pt] (-1,-1) -- (7,-1) -- (7,1) -- (-1,1) -- (-1,-1);
				\draw[lgray,line width=1.5pt] (1,-1) -- (1,1);
				\draw[lgray,line width=1.5pt] (3,-1) -- (3,1);
				\draw[lgray,line width=1.5pt] (5,-1) -- (5,1);
				\node[left] at (-1,0) {\tiny $(0,\textbf{\emph{e}}_0)$};\node[right] at (7,0) {\tiny $(0,\textbf{\emph{e}}_0)$};
			}
			\quad\quad & i \not\in \mathfrak{T} (\bar{\nu}),\ i\not=0,
			\\ \\ 
			\tikz{0.4}{
				\draw[lgray,line width=1.5pt] (-1,-1) -- (7,-1) -- (7,1) -- (-1,1) -- (-1,-1);
				\draw[lgray,line width=1.5pt] (1,-1) -- (1,1);
				\draw[lgray,line width=1.5pt] (3,-1) -- (3,1);
				\draw[lgray,line width=1.5pt] (5,-1) -- (5,1);
				\node[left] at (-1,0) {\tiny $(m,\textbf{\emph{e}}_0)$};\node[right] at (7,0) {\tiny $(0,\textbf{\emph{e}}_0)$};
			}
			\quad\quad & i=0.
		\end{array} 
		\right.
	\end{align}
	The constant $b_{\nu}(q)$ appearing in \eqref{main-formula} is given by $b_{\nu}(q) = \prod_{j = 1}^{\infty} (q;q)_{m_j (\nu)}$.\index{B@$b_{\nu} (q)$}
\end{thm}

\section{Expansion Coefficients for Modified Macdonald Polynomials}

Fix complex numbers $q, t \in \mathbb{C}$; let $\textbf{x} = (x_1, x_2, \ldots , x_k)$ denote a (possibly infinite) set of variables; and let $\lambda \in \Sign$ denote a signature. Recalling the modified Macdonald polynomial $\widetilde{J}_{\lambda} (\textbf{x}; q, t)$ from \eqref{jxjx}, we next provide combinatorial formulas for the coefficients $g_{\lambda}^{\nu}(q,t)$\index{G@$g_{\lambda}^{\nu} (q, t)$} in the expansion
\begin{align}
	\label{mod-mac-exp}
	\widetilde{J}_{\lambda}(\textbf{x};q,t)
	=
	\sum_{\nu}
	g_{\lambda}^{\nu}(q,t)
	Q'_{\nu}(\textbf{x};t),
\end{align}

\noindent where $Q'_{\nu}(\textbf{x}; t)$ denotes a modified Hall--Littlewood polynomial (with its parameter $q$ replaced by $t$). Since \Cref{psuml} expands $\widetilde{J}_{\lambda}$ over the LLT polynomials, this follows as a quick corollary of \Cref{thm:comb}.

\begin{cor}
	
	\label{jsumg}
	
	Fix integers $p, m \in [1, n]$; an anti-dominant composition $\lambda$ of length $p$; a partition $\nu$ of length $m$; and an integer $K \ge \max \{ \lambda_p, \nu_1 + 1 \}$. For any integers $i \in [1, p]$ and $j \ge 1$, set $w_{i,j} = q^{\lambda_i-j} \bm{1}_{\lambda_i > j}$. Then, the coefficient $g_{\lambda^+}^{\nu}(q,t)$ is given by
	\begin{align}
		\label{gsumfunction} 
		g_{\lambda^+}^{\nu}(q,t)
		=
		\frac{(-1)^{|\bar{\nu}|}}{b_{\nu}(t)(t;t)_{m}}
		\sum_{\mathscr{M}}
		\prod_{i=1}^p
		\prod_{j = 1}^K 
		%\left(
		%q^{\lambda_i-j}
		%\cdot
		%\bm{1}_{\lambda_i > j}
		%\right)^{N_{i,j}}
		w_{i,j}^{M_{i,j}} \times
		\tikz{1}{
			\filldraw[lgray,line width=1.5pt,fill=llgray] (1.5,3.5) -- (7.5,3.5) -- (7.5,4.5) -- (1.5,4.5) -- (1.5,3.5);
			\filldraw[lgray,line width=1.5pt,fill=llgray] (1.5,0.5) -- (7.5,0.5) -- (7.5,1.5) -- (1.5,1.5) -- (1.5,0.5);
			\foreach\y in {0,...,6}{
				\draw[lgray,line width=1.5pt] (1.5,0.5+\y) -- (7.5,0.5+\y);
			}
			\foreach\x in {0,...,6}{
				\draw[lgray,line width=1.5pt] (1.5+\x,0.5) -- (1.5+\x,6.5);
			}
			%bottom labels
			\node[above] at (1.875,-.0625) {\footnotesize $(0,\textbf{\emph{e}}_{[1, p]})$};
			\node[above] at (3.125,0) {\footnotesize $(0,\textbf{\emph{M}}_1)$};
			\node[above,text centered] at (5,0) {$\cdots$};
			\node[above] at (7,0) {\footnotesize $(0,\textbf{\emph{M}}_K)$};
			%top labels
			\node[above] at (1.875,6.5) {\footnotesize $(0,\textbf{\emph{e}}_0)$};
			\node[above] at (3.25,6.5) {\footnotesize $(0,\textbf{\emph{M}}_1+\textbf{\emph{I}}_1)$};
			\node[above,text centered] at (5,6.5) {$\cdots$};
			\node[above] at (7,6.5) {\footnotesize $(0,\textbf{\emph{M}}_K + \textbf{\emph{I}}_K)$};
			%right labels
			\node[right] at (7.5,1) {\footnotesize $(1,\textbf{\emph{e}}_0)$};
			\node[right] at (7.5,2) {\footnotesize $(0,\textbf{\emph{e}}_0)$};
			\node[right] at (7.5,3.1) {$\vdots$};
			\node[right] at (7.5,4) {\footnotesize $(1,\textbf{\emph{e}}_0)$};
			\node[right] at (7.5,5) {$\vdots$};
			\node[right] at (7.5,6) {\footnotesize $(0,\textbf{\emph{e}}_0)$};
			%left labels
			\node[left] at (1.5,1) {\footnotesize $(0,\textbf{\emph{e}}_0)$};
			\node[left] at (1.5,2) {\footnotesize $(0,\textbf{\emph{e}}_0)$};
			\node[left] at (1.5,3.1) {$\vdots$};
			\node[left] at (1.5,4) {\footnotesize $(0,\textbf{\emph{e}}_0)$};
			\node[left] at (1.5,5) {$\vdots$};
			\node[left] at (1.5,6) {\footnotesize $(m,\textbf{\emph{e}}_0)$};
		}
	\end{align}

	\noindent where the rows of the partition function \eqref{gsumfunction} are specified by \eqref{maya-seq1} and \eqref{maya-seq2}, using the same weights as in \Cref{ssec:weights1} and \Cref{ssec:weights2} but with the $q$ there replaced by $t$ here. The sum is over all sequences $\mathscr{M} = (\textbf{\emph{M}}_1, \textbf{\emph{M}}_2, \ldots, \textbf{\emph{M}}_K)$ of elements in $\{ 0, 1 \}^n$, with $\textbf{\emph{M}}_j = (M_{1,j}, M_{2, j}, \ldots,M_{n,j})$, such that $M_{i, k} = 0$ for $i > p$ and $\textbf{\emph{I}}_k (\lambda) + \textbf{\emph{M}}_k \in \{ 0, 1 \}^n$ for each $k \in [1, K]$. Here, we abbreviated $\textbf{\emph{I}}_k = \textbf{\emph{I}}_k (\lambda)$ for each $k$, and we recalled $b_{\nu} (t)$ from \Cref{thm:comb}.
\end{cor}

%\begin{rem}
%The formula \eqref{mac-comb} provides two vastly different combinatorial expressions for the Kostka--Foulkes polynomials. On the one hand, we may choose $q$ generic, $t=0$, when \eqref{mod-mac-exp} becomes
%%
%\begin{align*}
%Q'_{\lambda}(x_1,\dots,x_n;q)
%=
%\sum_{\nu}
%g_{\lambda}^{\nu}(q,0)
%s_{\nu}(x_1,\dots,x_n).
%\end{align*}
%%
%Comparing with \eqref{kf-1}, we see that $g_{\lambda}^{\nu}(q,0) = K_{\lambda}^{\nu}(q)$. On the other hand, the pre-fused version of \eqref{mod-mac-exp} at $q=t$ yields
%%
%\begin{align*}
%c_{\lambda}(t,t)
%\cdot
%s_{\lambda}(x_1,\dots,x_n)
%=
%\sum_{\nu}
%g_{\lambda}^{\nu}(t,t)
%b_{\nu}(t)
%\cdot
%P_{\nu}(x_1,\dots,x_n;t),
%\end{align*}
%%
%from which it follows that $g_{\lambda}^{\nu}(t,t) \cdot \frac{b_{\nu}(t)}{c_{\lambda}(t,t)} = K^{\lambda}_{\nu}(t)$. A more detailed analysis of these two different expressions for the Kostka--Foulkes polynomials is interesting in its own right; either formula, for example, could be used in \eqref{c-into-f}.
%\end{rem}

\begin{proof}
	
	Together \eqref{limitw}, \eqref{gfhe}, and \Cref{weightesum} imply the diagrammatic equality 
		\begin{align}
		\label{g0sum}
		\mathcal{G}_{\boldsymbol{\lambda} / \boldsymbol{\mu}} (\textbf{x}; \infty \boldsymbol{\mid} 0; 0)
		=
		\tikz{1}{
			\foreach\y in {0,...,5}{
				\draw[lgray,line width=1.5pt] (1.5,0.5+\y) -- (7.5,0.5+\y);
			}
			\foreach\x in {0,...,6}{
				\draw[lgray,line width=1.5pt] (1.5+\x,0.5) -- (1.5+\x,5.5);
			}
			%spectral parameters
			\node[text centered] at (2,1) {$x_1$};
			\node[text centered] at (4,1) {$\cdots$};
			\node[text centered] at (5,1) {$\cdots$};
			\node[text centered] at (7,1) {$x_1$};
			\node[text centered] at (2,2) {$x_2$};
			\node[text centered] at (4,2) {$\cdots$};
			\node[text centered] at (5,2) {$\cdots$};
			\node[text centered] at (7,2) {$x_2$};
			\node[text centered] at (2,5) {$x_n$};
			\node[text centered] at (4,5) {$\cdots$};
			\node[text centered] at (5,5) {$\cdots$};
			\node[text centered] at (7,5) {$x_n$};
			\node[text centered] at (2,3.1) {$\vdots$};
			\node[text centered] at (2,4.1) {$\vdots$};
			\node[text centered] at (7,3.1) {$\vdots$};
			\node[text centered] at (7,4.1) {$\vdots$};
			%bottom labels
			\node[below] at (2,0.5) {\scriptsize$\textbf{e}_{[1, p]}$};
			\node[below] at (3,0.5) {\scriptsize$\textbf{S}_1 (\boldsymbol{\mu})$};
			\node[above,text centered] at (4,0) {$\cdots$};
			\node[above,text centered] at (5,0) {$\cdots$};
			\node[below] at (7,0.5) {\scriptsize$\textbf{S}_K (\boldsymbol{\mu})$};
			%top labels
			\node[left = 5, above] at (2,5.5) {\scriptsize $\textbf{e}_0$};
			\node[above] at (3,5.5) {\scriptsize$\textbf{S}_1 (\boldsymbol{\lambda})$};
			\node[above,text centered] at (4,5.5) {$\cdots$};
			\node[above,text centered] at (5,5.5) {$\cdots$};
			\node[above] at (7,5.5) {\scriptsize$\textbf{S}_K (\boldsymbol{\lambda})$};
			%right labels
			\node[right, scale = .8] at (7.5,1) {$\textbf{e}_0$};
			\node[right, scale = .8] at (7.5,2) {$\textbf{e}_0$};
			\node[right] at (7.5,3.1) {$\vdots$};
			\node[right] at (7.5,4.1) {$\vdots$};
			\node[right, scale = .8] at (7.5,5) {$\textbf{e}_0$};
			%left labels
			\node[left, scale = .8] at (1.5,1) {$\textbf{e}_0$};
			\node[left, scale = .8] at (1.5,2) {$\textbf{e}_0$};
			\node[left] at (1.5,3.1) {$\vdots$};
			\node[left] at (1.5,4.1) {$\vdots$};
			\node[left, scale = .8] at (1.5,5) {$\textbf{e}_0$};
		}
	\end{align}

	\noindent where the weights in the $j$-th row (from the bottom) are given by $\mathcal{W}_{x_j} (\textbf{A}, \textbf{B}; \textbf{C}, \textbf{D} \boldsymbol{\mid} \infty, 0) = x_j^{|\textbf{D}|} q^{\varphi (\textbf{D}, \textbf{C} + \textbf{D})} \textbf{1}_{|\textbf{V}| = 0} \textbf{1}_{\textbf{A} + \textbf{B} = \textbf{C} + \textbf{D}}$. Combining this with \eqref{sumj} gives
	\begin{align*}
		\widetilde{J}_{\lambda^+} (\textbf{x}; q, t)
		=
		\sum_{\mathscr{M}}
		\prod_{i=1}^p
		\prod_{j=1}^{N}
		w_{i,j}^{M_{i,j}}
		\times
		\tikz{1}{
			\foreach\y in {0,...,5}{
				\draw[lgray,line width=1.5pt] (1.5,0.5+\y) -- (7.5,0.5+\y);
			}
			\foreach\x in {0,...,6}{
				\draw[lgray,line width=1.5pt] (1.5+\x,0.5) -- (1.5+\x,5.5);
			}
			%spectral parameters
			\node[text centered] at (2,1) {$x_1$};
			\node[text centered] at (4,1) {$\cdots$};
			\node[text centered] at (5,1) {$\cdots$};
			\node[text centered] at (7,1) {$x_1$};
			\node[text centered] at (2,2) {$x_2$};
			\node[text centered] at (4,2) {$\cdots$};
			\node[text centered] at (5,2) {$\cdots$};
			\node[text centered] at (7,2) {$x_2$};
			\node[text centered] at (2,5) {$x_n$};
			\node[text centered] at (4,5) {$\cdots$};
			\node[text centered] at (5,5) {$\cdots$};
			\node[text centered] at (7,5) {$x_n$};
			\node[text centered] at (2,3.1) {$\vdots$};
			\node[text centered] at (2,4.1) {$\vdots$};
			\node[text centered] at (7,3.1) {$\vdots$};
			\node[text centered] at (7,4.1) {$\vdots$};
			%bottom labels
			\node[below] at (2,0.5) {\scriptsize$\textbf{e}_{[1, p]}$};
			\node[below] at (3,0.5) {\scriptsize$\textbf{M}_1$};
			\node[above,text centered] at (4,0) {$\cdots$};
			\node[above,text centered] at (5,0) {$\cdots$};
			\node[below] at (7,0.5) {\scriptsize$\textbf{M}_K$};
			%top labels
			\node[left = 5, above] at (2,5.5) {\scriptsize $\textbf{e}_0$};
			\node[above] at (3,5.5) {\scriptsize$\textbf{M}_1+\textbf{I}_1$};
			\node[above,text centered] at (4,5.5) {$\cdots$};
			\node[above,text centered] at (5,5.5) {$\cdots$};
			\node[above] at (7,5.5) {\scriptsize$\textbf{M}_K+\textbf{I}_K$};
			%right labels
			\node[right] at (7.5,1) {$0$};
			\node[right] at (7.5,2) {$0$};
			\node[right] at (7.5,3.1) {$\vdots$};
			\node[right] at (7.5,4.1) {$\vdots$};
			\node[right] at (7.5,5) {$0$};
			%left labels
			\node[left] at (1.5,1) {$0$};
			\node[left] at (1.5,2) {$0$};
			\node[left] at (1.5,3.1) {$\vdots$};
			\node[left] at (1.5,4.1) {$\vdots$};
			\node[left] at (1.5,5) {$0$};
		}
	\end{align*}

	\noindent where the vertices in the $j$-row are again taken with respect to the weights $\mathcal{W}_x (\textbf{A}, \textbf{B}; \textbf{C}, \textbf{D} \boldsymbol{\mid} \infty, 0)$. This, together with \eqref{g0sum} and \Cref{thm:comb}, implies the corollary. 
\end{proof}

\section{Examples}

\label{Models2n}

Before moving on to the proof of \Cref{thm:comb}, we present some explicit examples of the combinatorial formula \eqref{main-formula}. Naively one would expect \eqref{main-formula} to be a huge sum, but in practice it consists of just a handful of non-zero lattice configurations, as the following examples indicate. All examples in this section were generated by a \texttt{Mathematica} implementation of \eqref{main-formula}; the code is available from the authors upon request. 

In what follows, we draw bosonic paths in the color blue. Unlike in previous diagrams, to simplify the figures, we depict paths as sometimes making diagonal steps, which are equivalent to making one step up or down and one step right. 

\begin{example}
	Let $n=1$, $M = 3$, $\boldsymbol{\lambda} = \left( \lambda^{(1)} \right)$, and $\boldsymbol{\mu} = \left( \mu^{(1)} \right)$ with $\lambda^{(1)} = (2,1,1)$ and $\mu^{(1)} = (0,0,0)$.  Let $\nu = (3,1)$ and choose $K = 5$, so that $\bar{\nu} = (3, 1)$ and $\mathfrak{T} (\bar{\nu}) = (5,2)$. Given that $|\bar{\nu}| = 4$, $m = 2$, and $b_{\nu}(q) = (1-q)^2$, we have
	\begin{align}
		\label{factor1}
		\frac{(-1)^{|\bar{\nu}|}}{b_{\nu}(q)(q;q)_{m}}
		=
		\frac{1}{(1-q)^3(1-q^2)}.
	\end{align}
	The formula \eqref{main-formula} provides two non-vanishing lattice configurations, indicated below with their weights, where we have multiplied by the overall factor \eqref{factor1}. 
	\begin{align*}
		\begin{array}{cc}
			\tikz{0.7}{
				\filldraw[lgray,line width=1.5pt,fill=llgray] (2.5,3.5) -- (7.5,3.5) -- (7.5,4.5) -- (2.5,4.5) -- (2.5,3.5);
				\filldraw[lgray,line width=1.5pt,fill=llgray] (2.5,0.5) -- (7.5,0.5) -- (7.5,1.5) -- (2.5,1.5) -- (2.5,0.5);
				\foreach\y in {0,...,6}{
					\draw[lgray,line width=1.5pt] (2.5,0.5+\y) -- (7.5,0.5+\y);
				}
				\foreach\x in {1,...,6}{
					\draw[lgray,line width=1.5pt] (1.5+\x,0.5) -- (1.5+\x,6.5);
				}
				%%right labels
				%\node[right] at (7.5,1) {\footnotesize $(1,\textbf{e}_0)$};
				%\node[right] at (7.5,2) {\footnotesize $(0,\textbf{e}_0)$};
				%\node[right] at (7.5,3) {\footnotesize $(0,\textbf{e}_0)$};
				%\node[right] at (7.5,4) {\footnotesize $(1,\textbf{e}_0)$};
				%\node[right] at (7.5,5) {\footnotesize $(0,\textbf{e}_0)$};
				%\node[right] at (7.5,6) {\footnotesize $(0,\textbf{e}_0)$};
				%%left labels
				%\node[left] at (2.5,1) {\footnotesize $(0,\textbf{e}_0)$};
				%\node[left] at (2.5,2) {\footnotesize $(0,\textbf{e}_0)$};
				%\node[left] at (2.5,3) {\footnotesize $(0,\textbf{e}_0)$};
				%\node[left] at (2.5,4) {\footnotesize $(0,\textbf{e}_0)$};
				%\node[left] at (2.5,5) {\footnotesize $(0,\textbf{e}_0)$};
				%\node[left] at (2.5,6) {\footnotesize $(2,\textbf{e}_0)$};
				%boson paths
				\draw[blue,line width=1pt,<-] (7.5,1-0.1) -- (2.5,6-0.1);
				\draw[blue,line width=1pt,<-] (7.5,4) -- (5.5,4) -- (4.5,5) -- (3.5,5) -- (2.5,6);
				%fermion paths
				\draw[red,line width=1pt,->] (3,0.5) -- (3,4.5) -- (4,5.5) -- (4,6.5);
				\draw[red,line width=1pt,->] (4,0.5) -- (4,3.5) -- (5,4.5) -- (5,6.5);
				\draw[red,line width=1pt,->] (5-0.1,0.5) -- (5-0.1,3.5) -- (5.5,4+0.1) -- (6.5,4+0.1) -- (7-0.1,4.5) -- (7-0.1,6.5);
			}
			&
			\quad
			%%%%%%
			\tikz{0.7}{
				\filldraw[lgray,line width=1.5pt,fill=llgray] (2.5,3.5) -- (7.5,3.5) -- (7.5,4.5) -- (2.5,4.5) -- (2.5,3.5);
				\filldraw[lgray,line width=1.5pt,fill=llgray] (2.5,0.5) -- (7.5,0.5) -- (7.5,1.5) -- (2.5,1.5) -- (2.5,0.5);
				\foreach\y in {0,...,6}{
					\draw[lgray,line width=1.5pt] (2.5,0.5+\y) -- (7.5,0.5+\y);
				}
				\foreach\x in {1,...,6}{
					\draw[lgray,line width=1.5pt] (1.5+\x,0.5) -- (1.5+\x,6.5);
				}
				%%right labels
				%\node[right] at (7.5,1) {\footnotesize $(1,\textbf{e}_0)$};
				%\node[right] at (7.5,2) {\footnotesize $(0,\textbf{e}_0)$};
				%\node[right] at (7.5,3) {\footnotesize $(0,\textbf{e}_0)$};
				%\node[right] at (7.5,4) {\footnotesize $(1,\textbf{e}_0)$};
				%\node[right] at (7.5,5) {\footnotesize $(0,\textbf{e}_0)$};
				%\node[right] at (7.5,6) {\footnotesize $(0,\textbf{e}_0)$};
				%%left labels
				%\node[left] at (2.5,1) {\footnotesize $(0,\textbf{e}_0)$};
				%\node[left] at (2.5,2) {\footnotesize $(0,\textbf{e}_0)$};
				%\node[left] at (2.5,3) {\footnotesize $(0,\textbf{e}_0)$};
				%\node[left] at (2.5,4) {\footnotesize $(0,\textbf{e}_0)$};
				%\node[left] at (2.5,5) {\footnotesize $(0,\textbf{e}_0)$};
				%\node[left] at (2.5,6) {\footnotesize $(2,\textbf{e}_0)$};
				%boson paths
				\draw[blue,line width=1pt,<-] (7.5,1-0.1) -- (6.5,2-0.1) -- (5.5,2-0.1) -- (4-0.1,3.5) -- (4-0.1,4.5) -- (2.5,6-0.1);
				\draw[blue,line width=1pt,<-] (7.5,4) -- (6.5,4) -- (5.5,5) -- (3.5,5) -- (2.5,6);
				%fermion paths
				\draw[red,line width=1pt,->] (3,0.5) -- (3,4.5) -- (4,5.5) -- (4,6.5);
				\draw[red,line width=1pt,->] (4,0.5) -- (4,4.5) -- (5,5.5) -- (5,6.5);
				\draw[red,line width=1pt,->] (5,0.5) -- (5,1.5) -- (6,2.5) -- (6,3.5) -- (7,4.5) -- (7,6.5);
			}
			\\
			\\
			-q^2  & \quad -q(1-q)
		\end{array}
	\end{align*}
	We therefore find that $f^{\nu}_{\boldsymbol{\lambda}/\boldsymbol{\mu}}(q)
	=
	-q^2 - q(1-q)
	=
	-q$.
\end{example}

\begin{example}
	Let $n=1$, $M = 3$, $\boldsymbol{\lambda} = \left( \lambda^{(1)} \right)$, and $\boldsymbol{\mu} = \left( \mu^{(1)} \right)$ with $\lambda^{(1)} = (2, 2, 1)$ and $\mu^{(1)} = (0,0,0)$. Let $\nu = (3,2)$ and choose $K=5$, so that $\bar{\nu} = (2,1)$ and $\mathfrak{T} (\bar{\nu}) = (4,2)$. Given that $|\bar{\nu}| = 3$, $m = 2$, and $b_{\nu}(q) = (1-q)^2$, we have
	\begin{align}
		\label{factor2}
		\frac{(-1)^{|\bar{\nu}|}}{b_{\nu}(q)(q;q)_{m}}
		=
		\frac{-1}{(1-q)^3(1-q^2)}.
	\end{align}
	The formula \eqref{main-formula} provides three non-vanishing lattice configurations, indicated below with their weights, where we have multiplied by the overall factor \eqref{factor2}.
	\begin{align*}
		\begin{array}{ccc}
			\tikz{0.7}{
				\filldraw[lgray,line width=1.5pt,fill=llgray] (2.5,3.5) -- (7.5,3.5) -- (7.5,4.5) -- (2.5,4.5) -- (2.5,3.5);
				\filldraw[lgray,line width=1.5pt,fill=llgray] (2.5,1.5) -- (7.5,1.5) -- (7.5,2.5) -- (2.5,2.5) -- (2.5,1.5);
				\foreach\y in {1,...,6}{
					\draw[lgray,line width=1.5pt] (2.5,0.5+\y) -- (7.5,0.5+\y);
				}
				\foreach\x in {1,...,6}{
					\draw[lgray,line width=1.5pt] (1.5+\x,1.5) -- (1.5+\x,6.5);
				}
				%%right labels
				%\node[right] at (7.5,2) {\footnotesize $(1,\textbf{e}_0)$};
				%\node[right] at (7.5,3) {\footnotesize $(0,\textbf{e}_0)$};
				%\node[right] at (7.5,4) {\footnotesize $(1,\textbf{e}_0)$};
				%\node[right] at (7.5,5) {\footnotesize $(0,\textbf{e}_0)$};
				%\node[right] at (7.5,6) {\footnotesize $(0,\textbf{e}_0)$};
				%%left labels
				%\node[left] at (2.5,2) {\footnotesize $(0,\textbf{e}_0)$};
				%\node[left] at (2.5,3) {\footnotesize $(0,\textbf{e}_0)$};
				%\node[left] at (2.5,4) {\footnotesize $(0,\textbf{e}_0)$};
				%\node[left] at (2.5,5) {\footnotesize $(0,\textbf{e}_0)$};
				%\node[left] at (2.5,6) {\footnotesize $(2,\textbf{e}_0)$};
				%%boson paths
				\draw[blue,line width=1pt,<-] (7.5,2) -- (6.5,3) -- (5.5,3) -- (2.5,6);
				\draw[blue,line width=1pt,<-] (7.5,4+0.1) -- (5.5,4+0.1) -- (4.5,5+0.1) -- (3.5,5+0.1) -- (2.5,6+0.1);
				%fermion paths
				\draw[red,line width=1pt,->] (3,1.5) -- (3,4.5) -- (4,5.5) -- (4,6.5);
				\draw[red,line width=1pt,->] (4-0.1,1.5) -- (4-0.1,3.5) -- (4.5,4+0.1) -- (5.5,4+0.1) -- (6-0.1,4.5) -- (6-0.1,6.5);
				\draw[red,line width=1pt,->] (5,1.5) -- (5,2.5) -- (7,4.5) -- (7,6.5);
			}
			%%%
			&
			\quad
			
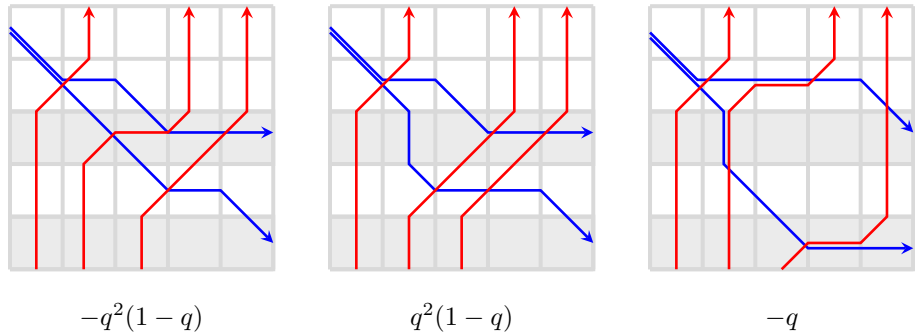
\tikz{0.7}{
				\filldraw[lgray,line width=1.5pt,fill=llgray] (2.5,3.5) -- (7.5,3.5) -- (7.5,4.5) -- (2.5,4.5) -- (2.5,3.5);
				\filldraw[lgray,line width=1.5pt,fill=llgray] (2.5,1.5) -- (7.5,1.5) -- (7.5,2.5) -- (2.5,2.5) -- (2.5,1.5);
				\foreach\y in {1,...,6}{
					\draw[lgray,line width=1.5pt] (2.5,0.5+\y) -- (7.5,0.5+\y);
				}
				\foreach\x in {1,...,6}{
					\draw[lgray,line width=1.5pt] (1.5+\x,1.5) -- (1.5+\x,6.5);
				}
				%%right labels
				%\node[right] at (7.5,2) {\footnotesize $(1,\textbf{e}_0)$};
				%\node[right] at (7.5,3) {\footnotesize $(0,\textbf{e}_0)$};
				%\node[right] at (7.5,4) {\footnotesize $(1,\textbf{e}_0)$};
				%\node[right] at (7.5,5) {\footnotesize $(0,\textbf{e}_0)$};
				%\node[right] at (7.5,6) {\footnotesize $(0,\textbf{e}_0)$};
				%%left labels
				%\node[left] at (2.5,2) {\footnotesize $(0,\textbf{e}_0)$};
				%\node[left] at (2.5,3) {\footnotesize $(0,\textbf{e}_0)$};
				%\node[left] at (2.5,4) {\footnotesize $(0,\textbf{e}_0)$};
				%\node[left] at (2.5,5) {\footnotesize $(0,\textbf{e}_0)$};
				%\node[left] at (2.5,6) {\footnotesize $(2,\textbf{e}_0)$};
				%%boson paths
				\draw[blue,line width=1pt,<-] (7.5,2) -- (6.5,3) -- (4.5,3) -- (4,3.5) -- (4,4.5) -- (2.5,6);
				\draw[blue,line width=1pt,<-] (7.5,4+0.1) -- (5.5,4+0.1) -- (4.5,5+0.1) -- (3.5,5+0.1) -- (2.5,6+0.1);
				%fermion paths
				\draw[red,line width=1pt,->] (3,1.5) -- (3,4.5) -- (4,5.5) -- (4,6.5);
				\draw[red,line width=1pt,->] (4,1.5) -- (4,2.5) -- (6,4.5) -- (6,6.5);
				\draw[red,line width=1pt,->] (5,1.5) -- (5,2.5) -- (7,4.5) -- (7,6.5);
			}
			%%%
			&
			\quad
			\tikz{0.7}{
				\filldraw[lgray,line width=1.5pt,fill=llgray] (2.5,3.5) -- (7.5,3.5) -- (7.5,4.5) -- (2.5,4.5) -- (2.5,3.5);
				\filldraw[lgray,line width=1.5pt,fill=llgray] (2.5,1.5) -- (7.5,1.5) -- (7.5,2.5) -- (2.5,2.5) -- (2.5,1.5);
				\foreach\y in {1,...,6}{
					\draw[lgray,line width=1.5pt] (2.5,0.5+\y) -- (7.5,0.5+\y);
				}
				\foreach\x in {1,...,6}{
					\draw[lgray,line width=1.5pt] (1.5+\x,1.5) -- (1.5+\x,6.5);
				}
				%%right labels
				%\node[right] at (7.5,2) {\footnotesize $(1,\textbf{e}_0)$};
				%\node[right] at (7.5,3) {\footnotesize $(0,\textbf{e}_0)$};
				%\node[right] at (7.5,4) {\footnotesize $(1,\textbf{e}_0)$};
				%\node[right] at (7.5,5) {\footnotesize $(0,\textbf{e}_0)$};
				%\node[right] at (7.5,6) {\footnotesize $(0,\textbf{e}_0)$};
				%%left labels
				%\node[left] at (2.5,2) {\footnotesize $(0,\textbf{e}_0)$};
				%\node[left] at (2.5,3) {\footnotesize $(0,\textbf{e}_0)$};
				%\node[left] at (2.5,4) {\footnotesize $(0,\textbf{e}_0)$};
				%\node[left] at (2.5,5) {\footnotesize $(0,\textbf{e}_0)$};
				%\node[left] at (2.5,6) {\footnotesize $(2,\textbf{e}_0)$};
				%%boson paths
				\draw[blue,line width=1pt,<-] (7.5,2-0.1) -- (5.5,2-0.1) -- (4-0.1,3.5) -- (4-0.1,4.5) -- (2.5,6-0.1);
				\draw[blue,line width=1pt,<-] (7.5,4+0.1) -- (6.5,5+0.1) -- (3.5-0.1,5+0.1) -- (2.5,6);
				%fermion paths
				\draw[red,line width=1pt,->] (3,1.5) -- (3,4.5) -- (4,5.5) -- (4,6.5);
				\draw[red,line width=1pt,->] (4,1.5) -- (4,4.5) -- (4.5,5) -- (5.5,5) -- (6,5.5) -- (6,6.5);
				\draw[red,line width=1pt,->] (5,1.5) -- (5.5,2) -- (6.5,2) -- (7,2.5) -- (7,6.5);
			}
			\\
			\\
			-q^2(1-q) &\quad q^2(1-q) &\quad -q
		\end{array}
	\end{align*}
	We therefore find that $f^{\nu}_{\boldsymbol{\lambda}/\boldsymbol{\mu}}(q)
	=
	-q^2(1-q) + q^2 (1-q) -q
	=
	-q$.
\end{example}

\begin{example}
	
	\label{lambdamun2} 
	
	Let $n=2$, $M = 2$, $\boldsymbol{\lambda}=\left(\lambda^{(1)}, \lambda^{(2)}\right)$, and $\boldsymbol{\mu} = \left(\mu^{(1)}, \mu^{(2)}\right)$ with $\lambda^{(1)} = (3,1)$, $\lambda^{(2)}=(2,2)$, $\mu^{(1)} = (2, 0)$, and $\mu^{(2)}=(1,1)$. Let $\nu = (2,1,1)$ and choose $K = 5$, so that $\bar{\nu} = (3, 3, 2)$ and $\mathfrak{T} (\bar{\nu}) = (6,5,3)$. Given that $|\bar{\nu}| = 8$, $m = 3$, and $b_{\nu}(q) = (1-q)^2(1-q^2)$, we have
	\begin{align}
		\label{factor3}
		\frac{(-1)^{|\bar{\nu}|}}{b_{\nu}(q)(q;q)_{m}}
		=
		\frac{1}{(1-q)^3(1-q^2)^2(1-q^3)}.
	\end{align}
	The formula \eqref{main-formula} provides five non-vanishing lattice configurations, indicated below with their weights, where we have multiplied by the overall factor \eqref{factor3}.
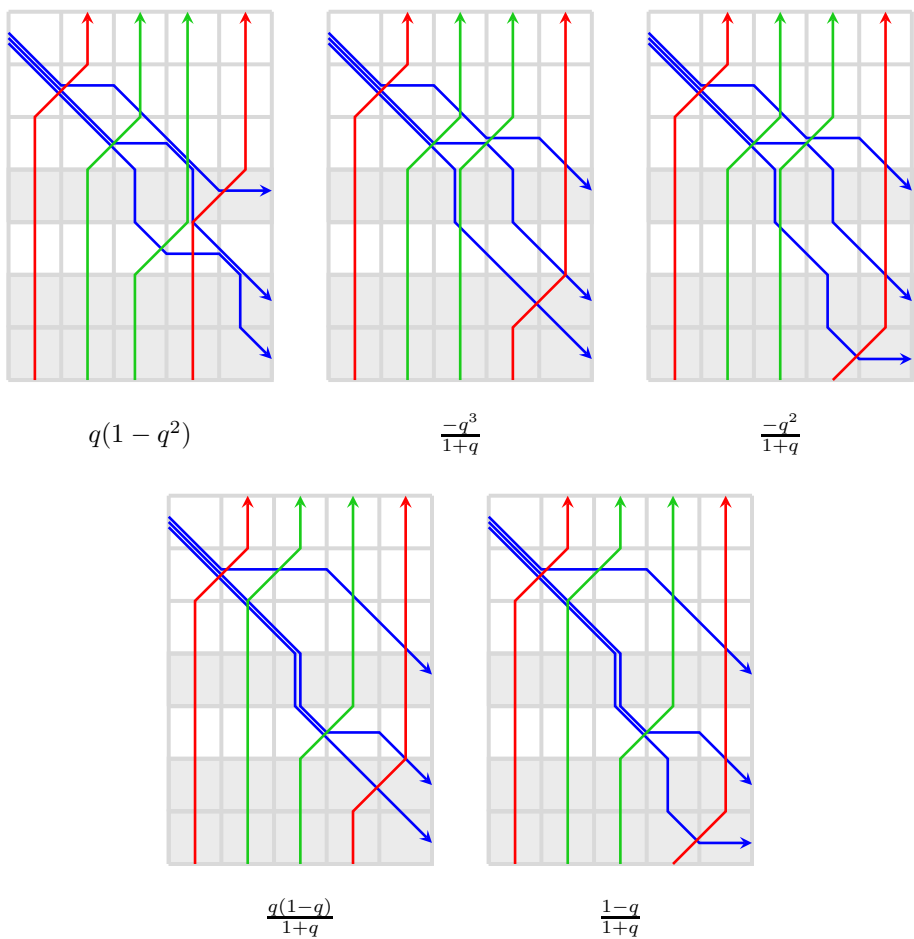
	\begin{align*}
		\begin{array}{ccc}
			\tikz{0.7}{
				\filldraw[lgray,line width=1.5pt,fill=llgray] (2.5,2.5) -- (7.5,2.5) -- (7.5,3.5) -- (2.5,3.5) -- (2.5,2.5);
				\filldraw[lgray,line width=1.5pt,fill=llgray] (2.5,0.5) -- (7.5,0.5) -- (7.5,1.5) -- (2.5,1.5) -- (2.5,0.5);
				\filldraw[lgray,line width=1.5pt,fill=llgray] (2.5,-0.5) -- (7.5,-0.5) -- (7.5,0.5) -- (2.5,0.5) -- (2.5,-0.5);
				\foreach\y in {-1,...,6}{
					\draw[lgray,line width=1.5pt] (2.5,0.5+\y) -- (7.5,0.5+\y);
				}
				\foreach\x in {1,...,6}{
					\draw[lgray,line width=1.5pt] (1.5+\x,-0.5) -- (1.5+\x,6.5);
				}
				%%right labels
				%\node[right] at (7.5,0) {\footnotesize $(1,\textbf{e}_0)$};
				%\node[right] at (7.5,1) {\footnotesize $(1,\textbf{e}_0)$};
				%\node[right] at (7.5,2) {\footnotesize $(0,\textbf{e}_0)$};
				%\node[right] at (7.5,3) {\footnotesize $(1,\textbf{e}_0)$};
				%\node[right] at (7.5,4) {\footnotesize $(0,\textbf{e}_0)$};
				%\node[right] at (7.5,5) {\footnotesize $(0,\textbf{e}_0)$};
				%\node[right] at (7.5,6) {\footnotesize $(0,\textbf{e}_0)$};
				%%left labels
				%\node[left] at (2.5,0) {\footnotesize $(0,\textbf{e}_0)$};
				%\node[left] at (2.5,1) {\footnotesize $(0,\textbf{e}_0)$};
				%\node[left] at (2.5,2) {\footnotesize $(0,\textbf{e}_0)$};
				%\node[left] at (2.5,3) {\footnotesize $(0,\textbf{e}_0)$};
				%\node[left] at (2.5,4) {\footnotesize $(0,\textbf{e}_0)$};
				%\node[left] at (2.5,5) {\footnotesize $(0,\textbf{e}_0)$};
				%\node[left] at (2.5,6) {\footnotesize $(3,\textbf{e}_0)$};
				%%boson paths
				\draw[blue,line width=1pt,<-] (7.5,0-0.1) -- (7-0.1,0.5) -- (7-0.1,1.5) -- (6.5,2-0.1) -- (5.5,2-0.1) -- (5-0.1,2.5) -- (5-0.1,3.5) -- (2.5,6-0.1);
				\draw[blue,line width=1pt,<-] (7.5,1) -- (6,2.5) -- (6,3.5) -- (5.5,4) -- (4.5,4) -- (2.5,6);
				\draw[blue,line width=1pt,<-] (7.5,3+0.1) -- (6.5,3+0.1) -- (4.5,5+0.1) -- (3.5,5+0.1) -- (2.5,6+0.1);
				%red fermion paths
				\draw[red,line width=1pt,->] (3,-0.5) -- (3,4.5) -- (4,5.5) -- (4,6.5);
				\draw[red,line width=1pt,->] (6,-0.5) -- (6,2.5) -- (7,3.5) -- (7,6.5);
				%green fermion paths
				\draw[green,line width=1pt,->] (4,-0.5) -- (4,3.5) -- (5,4.5) -- (5,6.5);
				\draw[green,line width=1pt,->] (5-0.1,-0.5) -- (5-0.1,1.5) -- (6-0.1,2.5) -- (6-0.1,6.5);
			}
			&
			\quad
			\tikz{0.7}{
				\filldraw[lgray,line width=1.5pt,fill=llgray] (2.5,2.5) -- (7.5,2.5) -- (7.5,3.5) -- (2.5,3.5) -- (2.5,2.5);
				\filldraw[lgray,line width=1.5pt,fill=llgray] (2.5,0.5) -- (7.5,0.5) -- (7.5,1.5) -- (2.5,1.5) -- (2.5,0.5);
				\filldraw[lgray,line width=1.5pt,fill=llgray] (2.5,-0.5) -- (7.5,-0.5) -- (7.5,0.5) -- (2.5,0.5) -- (2.5,-0.5);
				\foreach\y in {-1,...,6}{
					\draw[lgray,line width=1.5pt] (2.5,0.5+\y) -- (7.5,0.5+\y);
				}
				\foreach\x in {1,...,6}{
					\draw[lgray,line width=1.5pt] (1.5+\x,-0.5) -- (1.5+\x,6.5);
				}
				%%right labels
				%\node[right] at (7.5,0) {\footnotesize $(1,\textbf{e}_0)$};
				%\node[right] at (7.5,1) {\footnotesize $(1,\textbf{e}_0)$};
				%\node[right] at (7.5,2) {\footnotesize $(0,\textbf{e}_0)$};
				%\node[right] at (7.5,3) {\footnotesize $(1,\textbf{e}_0)$};
				%\node[right] at (7.5,4) {\footnotesize $(0,\textbf{e}_0)$};
				%\node[right] at (7.5,5) {\footnotesize $(0,\textbf{e}_0)$};
				%\node[right] at (7.5,6) {\footnotesize $(0,\textbf{e}_0)$};
				%%left labels
				%\node[left] at (2.5,0) {\footnotesize $(0,\textbf{e}_0)$};
				%\node[left] at (2.5,1) {\footnotesize $(0,\textbf{e}_0)$};
				%\node[left] at (2.5,2) {\footnotesize $(0,\textbf{e}_0)$};
				%\node[left] at (2.5,3) {\footnotesize $(0,\textbf{e}_0)$};
				%\node[left] at (2.5,4) {\footnotesize $(0,\textbf{e}_0)$};
				%\node[left] at (2.5,5) {\footnotesize $(0,\textbf{e}_0)$};
				%\node[left] at (2.5,6) {\footnotesize $(3,\textbf{e}_0)$};
				%%boson paths
				\draw[blue,line width=1pt,<-] (7.5,0-0.1) -- (5.5,2-0.1) -- (5-0.1,2.5) -- (5-0.1,3.5) -- (2.5,6-0.1);
				\draw[blue,line width=1pt,<-] (7.5,1) -- (6,2.5) -- (6,3.5) -- (5.5,4) -- (4.5,4) -- (2.5,6);
				\draw[blue,line width=1pt,<-] (7.5,3+0.1) -- (6.5,4+0.1) -- (5.5,4+0.1) -- (4.5,5+0.1) -- (3.5,5+0.1) -- (2.5,6+0.1);
				%red fermion paths
				\draw[red,line width=1pt,->] (3,-0.5) -- (3,4.5) -- (4,5.5) -- (4,6.5);
				\draw[red,line width=1pt,->] (6,-0.5) -- (6,0.5) -- (7,1.5) -- (7,6.5);
				%green fermion paths
				\draw[green,line width=1pt,->] (4,-0.5) -- (4,3.5) -- (5,4.5) -- (5,6.5);
				\draw[green,line width=1pt,->] (5,-0.5) -- (5,3.5) -- (6,4.5) -- (6,6.5);
			}
			&
			\quad
			\tikz{0.7}{
				\filldraw[lgray,line width=1.5pt,fill=llgray] (2.5,2.5) -- (7.5,2.5) -- (7.5,3.5) -- (2.5,3.5) -- (2.5,2.5);
				\filldraw[lgray,line width=1.5pt,fill=llgray] (2.5,0.5) -- (7.5,0.5) -- (7.5,1.5) -- (2.5,1.5) -- (2.5,0.5);
				\filldraw[lgray,line width=1.5pt,fill=llgray] (2.5,-0.5) -- (7.5,-0.5) -- (7.5,0.5) -- (2.5,0.5) -- (2.5,-0.5);
				\foreach\y in {-1,...,6}{
					\draw[lgray,line width=1.5pt] (2.5,0.5+\y) -- (7.5,0.5+\y);
				}
				\foreach\x in {1,...,6}{
					\draw[lgray,line width=1.5pt] (1.5+\x,-0.5) -- (1.5+\x,6.5);
				}
				%%right labels
				%\node[right] at (7.5,0) {\footnotesize $(1,\textbf{e}_0)$};
				%\node[right] at (7.5,1) {\footnotesize $(1,\textbf{e}_0)$};
				%\node[right] at (7.5,2) {\footnotesize $(0,\textbf{e}_0)$};
				%\node[right] at (7.5,3) {\footnotesize $(1,\textbf{e}_0)$};
				%\node[right] at (7.5,4) {\footnotesize $(0,\textbf{e}_0)$};
				%\node[right] at (7.5,5) {\footnotesize $(0,\textbf{e}_0)$};
				%\node[right] at (7.5,6) {\footnotesize $(0,\textbf{e}_0)$};
				%%left labels
				%\node[left] at (2.5,0) {\footnotesize $(0,\textbf{e}_0)$};
				%\node[left] at (2.5,1) {\footnotesize $(0,\textbf{e}_0)$};
				%\node[left] at (2.5,2) {\footnotesize $(0,\textbf{e}_0)$};
				%\node[left] at (2.5,3) {\footnotesize $(0,\textbf{e}_0)$};
				%\node[left] at (2.5,4) {\footnotesize $(0,\textbf{e}_0)$};
				%\node[left] at (2.5,5) {\footnotesize $(0,\textbf{e}_0)$};
				%\node[left] at (2.5,6) {\footnotesize $(3,\textbf{e}_0)$};
				%%boson paths
				\draw[blue,line width=1pt,<-] (7.5,0-0.1) -- (6.5,0-0.1) -- (6-0.1,0.5) -- (6-0.1,1.5) -- (5.5,2-0.1) -- (5-0.1,2.5) -- (5-0.1,3.5) -- (2.5,6-0.1);
				\draw[blue,line width=1pt,<-] (7.5,1) -- (6,2.5) -- (6,3.5) -- (5.5,4) -- (4.5,4) -- (2.5,6);
				\draw[blue,line width=1pt,<-] (7.5,3+0.1) -- (6.5,4+0.1) -- (5.5,4+0.1) -- (4.5,5+0.1) -- (3.5,5+0.1) -- (2.5,6+0.1);
				%red fermion paths
				\draw[red,line width=1pt,->] (3,-0.5) -- (3,4.5) -- (4,5.5) -- (4,6.5);
				\draw[red,line width=1pt,->] (6,-0.5) -- (7,0.5) -- (7,6.5);
				%green fermion paths
				\draw[green,line width=1pt,->] (4,-0.5) -- (4,3.5) -- (5,4.5) -- (5,6.5);
				\draw[green,line width=1pt,->] (5,-0.5) -- (5,3.5) -- (6,4.5) -- (6,6.5);
			}
			\\
			\\
			q (1-q^2) &\quad \frac{-q^3}{1+q} &\quad \frac{-q^2}{1+q}
		\end{array}
	\end{align*}
	
	\begin{align*}
		\begin{array}{cc}
			\tikz{0.7}{
				\filldraw[lgray,line width=1.5pt,fill=llgray] (2.5,2.5) -- (7.5,2.5) -- (7.5,3.5) -- (2.5,3.5) -- (2.5,2.5);
				\filldraw[lgray,line width=1.5pt,fill=llgray] (2.5,0.5) -- (7.5,0.5) -- (7.5,1.5) -- (2.5,1.5) -- (2.5,0.5);
				\filldraw[lgray,line width=1.5pt,fill=llgray] (2.5,-0.5) -- (7.5,-0.5) -- (7.5,0.5) -- (2.5,0.5) -- (2.5,-0.5);
				\foreach\y in {-1,...,6}{
					\draw[lgray,line width=1.5pt] (2.5,0.5+\y) -- (7.5,0.5+\y);
				}
				\foreach\x in {1,...,6}{
					\draw[lgray,line width=1.5pt] (1.5+\x,-0.5) -- (1.5+\x,6.5);
				}
				%%right labels
				%\node[right] at (7.5,0) {\footnotesize $(1,\textbf{e}_0)$};
				%\node[right] at (7.5,1) {\footnotesize $(1,\textbf{e}_0)$};
				%\node[right] at (7.5,2) {\footnotesize $(0,\textbf{e}_0)$};
				%\node[right] at (7.5,3) {\footnotesize $(1,\textbf{e}_0)$};
				%\node[right] at (7.5,4) {\footnotesize $(0,\textbf{e}_0)$};
				%\node[right] at (7.5,5) {\footnotesize $(0,\textbf{e}_0)$};
				%\node[right] at (7.5,6) {\footnotesize $(0,\textbf{e}_0)$};
				%%left labels
				%\node[left] at (2.5,0) {\footnotesize $(0,\textbf{e}_0)$};
				%\node[left] at (2.5,1) {\footnotesize $(0,\textbf{e}_0)$};
				%\node[left] at (2.5,2) {\footnotesize $(0,\textbf{e}_0)$};
				%\node[left] at (2.5,3) {\footnotesize $(0,\textbf{e}_0)$};
				%\node[left] at (2.5,4) {\footnotesize $(0,\textbf{e}_0)$};
				%\node[left] at (2.5,5) {\footnotesize $(0,\textbf{e}_0)$};
				%\node[left] at (2.5,6) {\footnotesize $(3,\textbf{e}_0)$};
				%%boson paths
				\draw[blue,line width=1pt,<-] (7.5,0-0.1) -- (5.5,2-0.1) -- (5-0.1,2.5) -- (5-0.1,3.5) -- (2.5,6-0.1);
				\draw[blue,line width=1pt,<-] (7.5,1) -- (6.5,2) -- (5.5,2) -- (5,2.5) -- (5,3.5) -- (2.5,6);
				\draw[blue,line width=1pt,<-] (7.5,3+0.1) --(5.5,5+0.1) -- (3.5,5+0.1) -- (2.5,6+0.1);
				%red fermion paths
				\draw[red,line width=1pt,->] (3,-0.5) -- (3,4.5) -- (4,5.5) -- (4,6.5);
				\draw[red,line width=1pt,->] (6,-0.5) -- (6,0.5) -- (7,1.5) -- (7,6.5);
				%green fermion paths
				\draw[green,line width=1pt,->] (4,-0.5) -- (4,4.5) -- (5,5.5) -- (5,6.5);
				\draw[green,line width=1pt,->] (5,-0.5) -- (5,1.5) -- (6,2.5) -- (6,6.5);
			}
			& \quad
			\tikz{0.7}{
				\filldraw[lgray,line width=1.5pt,fill=llgray] (2.5,2.5) -- (7.5,2.5) -- (7.5,3.5) -- (2.5,3.5) -- (2.5,2.5);
				\filldraw[lgray,line width=1.5pt,fill=llgray] (2.5,0.5) -- (7.5,0.5) -- (7.5,1.5) -- (2.5,1.5) -- (2.5,0.5);
				\filldraw[lgray,line width=1.5pt,fill=llgray] (2.5,-0.5) -- (7.5,-0.5) -- (7.5,0.5) -- (2.5,0.5) -- (2.5,-0.5);
				\foreach\y in {-1,...,6}{
					\draw[lgray,line width=1.5pt] (2.5,0.5+\y) -- (7.5,0.5+\y);
				}
				\foreach\x in {1,...,6}{
					\draw[lgray,line width=1.5pt] (1.5+\x,-0.5) -- (1.5+\x,6.5);
				}
				%%right labels
				%\node[right] at (7.5,0) {\footnotesize $(1,\textbf{e}_0)$};
				%\node[right] at (7.5,1) {\footnotesize $(1,\textbf{e}_0)$};
				%\node[right] at (7.5,2) {\footnotesize $(0,\textbf{e}_0)$};
				%\node[right] at (7.5,3) {\footnotesize $(1,\textbf{e}_0)$};
				%\node[right] at (7.5,4) {\footnotesize $(0,\textbf{e}_0)$};
				%\node[right] at (7.5,5) {\footnotesize $(0,\textbf{e}_0)$};
				%\node[right] at (7.5,6) {\footnotesize $(0,\textbf{e}_0)$};
				%%left labels
				%\node[left] at (2.5,0) {\footnotesize $(0,\textbf{e}_0)$};
				%\node[left] at (2.5,1) {\footnotesize $(0,\textbf{e}_0)$};
				%\node[left] at (2.5,2) {\footnotesize $(0,\textbf{e}_0)$};
				%\node[left] at (2.5,3) {\footnotesize $(0,\textbf{e}_0)$};
				%\node[left] at (2.5,4) {\footnotesize $(0,\textbf{e}_0)$};
				%\node[left] at (2.5,5) {\footnotesize $(0,\textbf{e}_0)$};
				%\node[left] at (2.5,6) {\footnotesize $(3,\textbf{e}_0)$};
				%%boson paths
				\draw[blue,line width=1pt,<-] (7.5,0-0.1) -- (6.5,0-0.1) -- (6-0.1,0.5) -- (6-0.1,1.5) -- (5.5,2-0.1) -- (5-0.1,2.5) -- (5-0.1,3.5) -- (2.5,6-0.1);
				\draw[blue,line width=1pt,<-] (7.5,1) -- (6.5,2) -- (5.5,2) -- (5,2.5) -- (5,3.5) -- (2.5,6);
				\draw[blue,line width=1pt,<-] (7.5,3+0.1) --(5.5,5+0.1) -- (3.5,5+0.1) -- (2.5,6+0.1);
				%red fermion paths
				\draw[red,line width=1pt,->] (3,-0.5) -- (3,4.5) -- (4,5.5) -- (4,6.5);
				\draw[red,line width=1pt,->] (6,-0.5) -- (7,0.5) -- (7,6.5);
				%green fermion paths
				\draw[green,line width=1pt,->] (4,-0.5) -- (4,4.5) -- (5,5.5) -- (5,6.5);
				\draw[green,line width=1pt,->] (5,-0.5) -- (5,1.5) -- (6,2.5) -- (6,6.5);
			}
			\\
			\\
			\frac{q(1-q)}{1+q} &\quad \frac{1-q}{1+q}
		\end{array}
	\end{align*}
	We therefore find that $f^{\nu}_{\boldsymbol{\lambda}/\boldsymbol{\mu}}(q)
	=
	q(1-q^2)-\frac{q^3}{1+q}-\frac{q^2}{1+q}+\frac{q(1-q)}{1+q}+\frac{1-q}{1+q}
	=
	1-q^2 -q^3$.
\end{example}

\section{Pre-fused LLT and Hall--Littlewood Polynomials}
\label{sec:prefused}

Before beginning with the proof of \Cref{thm:comb}, in this section we make a slight detour and give an alternative definition of the coefficients $f_{\boldsymbol{\lambda}/\boldsymbol{\mu}}^{\nu}(q)$, which will be useful in that proof. This is done via a ``pre-fused\footnote{Our reason for this terminology is that the associated vertex models under this transformation will only be partially fused, in that multiple arrows can exist along a vertical edge, but only one can exist along a horizontal edge.}  formulation'' \eqref{L-into-Q}, which at the level of the participating symmetric functions means applying the plethysm $X \mapsto (1-q) X$ to \eqref{L-into-Q}, where $X = \sum_{x \in \textbf{x}} x$ (recall \Cref{Polynomialsr} for our notation on plethysms). By \Cref{gl}, this converts $\mathcal{G}_{\boldsymbol{\lambda}/\boldsymbol{\mu}} (\textbf{x}; \infty \boldsymbol{\mid} 0; 0)$ into the \emph{pre-fused LLT function} $\mathcal{G}_{\boldsymbol{\lambda}/\boldsymbol{\mu}} (\textbf{x}; q^{-1 / 2} \boldsymbol{\mid} 0; 0)$, and $Q'_{\nu} (\textbf{x})$ into an ordinary Hall--Littlewood function $Q_{\nu} (\textbf{x})$. This plethysm does not affect the value of $f_{\boldsymbol{\lambda}/\boldsymbol{\mu}}^{\nu}(q)$.

\begin{lem}
	
	\label{fgq} 
	
	Fix an integer $M \ge 0$, and let $\textbf{\emph{x}} = (x_1, x_2, \ldots )$ denote an infinite set of variables. For any $\boldsymbol{\lambda}, \boldsymbol{\mu} \in \SeqSign_{n; M}$ and partition $\nu$, we have that $f_{\boldsymbol{\lambda}/\boldsymbol{\mu}}^{\nu}(q)$ is equal to the coefficient of $Q_{\nu}$ in the Hall--Littlewood expansion of a pre-fused LLT function $G_{\boldsymbol{\lambda}/\boldsymbol{\mu}} (\textbf{\emph{x}}; q^{-1 / 2} \boldsymbol{\mid} 0; 0)$, Explicitly,
	\begin{align}
		\label{pre-fused-exp}
		\mathcal{G}_{\boldsymbol{\lambda}/\boldsymbol{\mu}} (\textbf{\emph{x}}; q^{-1 / 2} \boldsymbol{\mid} 0; 0)
		=
		\sum_{\nu}
		f_{\boldsymbol{\lambda}/\boldsymbol{\mu}}^{\nu}(q)
		Q_{\nu} (\textbf{\emph{x}}),
	\end{align}
	with the sum taken over all partitions $\nu$.
\end{lem}

\begin{proof}
	
	This follows from applying the plethystic substitution $X \mapsto (1 - q) X$ to both sides of \eqref{L-into-Q}. Indeed, \eqref{lgr} and \eqref{1gl} together imply that this substitution maps $\mathcal{G}_{\boldsymbol{\lambda} / \boldsymbol{\mu}} (\textbf{x}; \infty \boldsymbol{\mid} 0; 0)$ to $\mathcal{G}_{\boldsymbol{\lambda} / \boldsymbol{\mu}} (\textbf{x}; q^{-1 / 2} \boldsymbol{\mid} 0; 0)$ and, by the definition (from \Cref{StabilityPolynomials}) of the modified Hall--Littlewood polynomials, it also maps $Q_{\nu}' (\textbf{x})$ to $Q_{\nu} (\textbf{x})$. 
\end{proof}

This pre-fused formulation of the $f_{\boldsymbol{\lambda} / \boldsymbol{\mu}}^{\nu}$ coefficients turns out to be valuable in the proof presented in \Cref{sec:proof} below, as the pre-fused functions $\mathcal{G}_{\boldsymbol{\lambda}/\boldsymbol{\mu}} (\textbf{x}; q^{-1 / 2} \boldsymbol{\mid} 0; 0)$ and $Q_{\nu} (\textbf{x})$ have a simpler algebraic structure than their fused counterparts. Indeed, \eqref{pre-fused-exp} will be our starting point in the proof of \Cref{thm:comb}; one of its main advantages is that it allows us to write an integral formula for the coefficients $f_{\boldsymbol{\lambda}/\boldsymbol{\mu}}^{\nu}(q)$, provided by the following proposition and corollary.

\begin{prop}
	
	\label{fc}
	
	Let $\textbf{\emph{x}} = (x_1, x_2, \ldots )$ denote an infinite set of variables, and let $F \in \Lambda (\textbf{\emph{x}})$ be a symmetric function with the expansion $F = \sum_{\nu} c_{\nu}(q) Q_{\nu}$ over the Hall--Littlewood basis. Write $F(x_1,\dots,x_N)$ for the symmetric polynomial obtained by setting $x_i=0$, for all $i > N$. For any partition $\nu = (\nu_1, \nu_2, \cdots , \nu_N)$, we have
	\begin{align}
		\label{coeff-int}
		c_{\nu}(q)
		=
		\frac{1}{b_{\nu}(q)}
		\displaystyle\frac{1}{(2 \pi \textbf{\emph{i}})^N}
		\oint_C
		\cdots
		\oint_C
		F(x_1^{-1},\dots,x_N^{-1})
		\prod_{1 \le i < j \le N}
		 \frac{x_j-x_i}{x_j-qx_i}
		\prod_{i=1}^N
		x_i^{\nu_i-1} dx_i,
	\end{align}
	%\begin{align}
	%\label{coeff-int}
	%c_{\nu}(q)
	%=
	%\frac{1}{b_{\nu}(q)}
	%\cdot
	%\oint_C \frac{dx_1}{2\pi \textbf{i}}
	%\cdots
	%\oint_C \frac{dx_n}{2\pi \textbf{i}}
	%\prod_{1 \le i < j \le n}
	%\left( \frac{x_i-x_j}{x_i-qx_j} \right)
	%\prod_{i=1}^{n}
	%x_i^{-\nu_i-1}
	%F(x_1,\dots,x_n),
	%\end{align}
	%
	where the contour $C$ is a positively oriented circle centered on the origin, and where we assume that $qC$ is contained in $C$ (equivalently, $|q| < 1$). Here, we recall $b_{\nu} (q) = \prod_{j = 1}^{\infty} (q; q)_{m_j (\nu)}$.
\end{prop}

\begin{proof}
	The Hall--Littlewood polynomials $P_{\lambda}$ (see Chapter 3.2 of \cite{SFP}) have the following orthogonality property, valid for any $\lambda, \nu \in \Sign_N$:
	\begin{align}
		\label{orthog}
		\displaystyle\frac{1}{(2 \pi \textbf{i})^N} \oint_C
		\cdots
		\oint_C 
		P_{\lambda}(x_1^{-1},\dots,x_N^{-1})
		\prod_{1 \le i < j \le N}
		\frac{x_j-x_i}{x_j-qx_i} 
		\prod_{i=1}^N
		x_i^{\nu_i-1} dx_i
		=
		\bm{1}_{\lambda=\nu}.
	\end{align}
	%\begin{align}
	%\label{orthog}
	%\oint_C \frac{dx_1}{2\pi \textbf{i}}
	%\cdots
	%\oint_C \frac{dx_n}{2\pi \textbf{i}}
	%\prod_{1 \le i < j \le n}
	%\left( \frac{x_i-x_j}{x_i-qx_j} \right)
	%\prod_{i=1}^{n}
	%x_i^{-\nu_i-1}
	%P_{\lambda}(x_1,\dots,x_n)
	%=
	%\bm{1}_{\lambda=\nu}.
	%\end{align}
	%
	The relation \eqref{orthog} can be deduced as a consequence of equation (3.2.15) of \cite{SFP}. Indeed, the latter states for any sequence of $M \ge N$ complex numbers $\textbf{y} = (y_1, y_2, \ldots , y_M)$ that
	\begin{align}
		\label{mac-eqn}
		Q_{\nu}(y_1,\dots,y_M)
		=
		{\rm Coeff}
		\left[
		\prod_{1 \le i<j \le N}
		\left( \frac{1-x_i x_j^{-1}}{1-qx_i x_j^{-1}} \right)
		\prod_{i=1}^N
		\prod_{j=1}^M
		\left( \frac{1-q y_j x_i^{-1}}{1-y_j x_i^{-1}} \right); 
		\prod_{i=1}^N
		x_i^{-\nu_i}
		\right],
	\end{align}
	where the right hand side calls for the coefficient of $\prod_{i=1}^N x_i^{-\nu_i}$ in the Laurent series expansion (about $x_1=\cdots=x_N=0$) of the indicated product, and where this Laurent series is computed under the assumption that 
	\begin{align*}
		|y_j x_i^{-1}| <1,
		\quad
		&
		\text{for all}\ \ 
		1 \le i \le N,\quad 
		1 \le j \le M,
		\\ 
		|q x_i x_j^{-1}| < 1;
		\quad
		&
		\text{for all}\ \ 
		1 \le i < j \le N.
	\end{align*}
	Casting \eqref{mac-eqn} in terms of integrals, it becomes
	\begin{align*}
		Q_{\nu}(y_1,\dots,y_M)
		=
		\displaystyle\frac{1}{(2 \pi \textbf{i})^N}
		\oint_C
		\cdots
		\oint_C
		\prod_{1 \le i<j \le N}
		 \frac{x_j-x_i}{x_j-qx_i} 
		\prod_{i=1}^N
		\prod_{j=1}^M
		\frac{1-qy_j x_i^{-1}}{1-y_j x_i^{-1}} 
		\prod_{i=1}^N
		x_i^{\nu_i-1} dx_i,
	\end{align*}
	where the contour $C$ is a positively oriented circle centered on the origin, such that $q C$ and $\{y_1,\dots,y_M \}$ are contained in $C$. Finally, we use the Cauchy identity for Hall--Littlewood polynomials (see equation (3.4.4) \cite{SFP}) in the preceding integral, which yields
	\begin{flalign*}
		Q_{\nu}(y_1,\dots,y_M)
		=
		\displaystyle\frac{1}{(2 \pi \textbf{i})^N} \oint_C
		\cdots
		\oint_C &
		\sum_{\lambda \in \Sign_N}
		P_{\lambda}(x_1^{-1},\dots,x_N^{-1})
		Q_{\lambda}(y_1,\dots,y_M)
		\\
		& \times
		\prod_{1 \le i<j \le N}
		 \frac{x_j-x_i}{x_j-qx_i} 
		\prod_{i=1}^N
		x_i^{\nu_i-1} dx_i,
	\end{flalign*}
	and \eqref{orthog} follows from the fact that for $\lambda \in \Sign_N$ the polynomials $Q_{\lambda}(y_1,\dots,y_M)$ are linearly independent as elements in the ring of symmetric polynomials in the $M \ge N$ variables $(y_1,\dots,y_M)$.
	
	Now starting from the symmetric function identity $F = \sum_{\lambda} c_{\lambda}(q) Q_{\lambda}$, we set $x_i=0$ for all $i>N$, which yields
	\begin{align}
		\label{finite-expansion}
		F(x_1,\dots,x_N)
		=
		\sum_{\lambda: \ell(\lambda) \le N}
		c_{\lambda}(q)
		Q_{\lambda}(x_1,\dots,x_N).
	\end{align}
	First reciprocating all variables in equation \eqref{finite-expansion}; then multiplying by 
	$\prod_{i=1}^N x_i^{\nu_i-1} \prod_{1 \le i < j \le N} (x_j-x_i) (x_j-qx_i)^{-1}$, and next performing $n$ integrations as in \eqref{orthog}, the result \eqref{coeff-int} then follows, recalling that $Q_{\lambda} = b_{\lambda}(q) P_{\lambda}$.
\end{proof}

\begin{cor}
	
	Fix integers $M \ge 0$ and $N \ge 1$; signature sequences $\boldsymbol{\lambda}, \boldsymbol{\nu} \in \SeqSign_{n; M}$; and a partition $\nu$. Recalling $b_{\nu} (q)$ from \Cref{fc}, the coefficient $f_{\boldsymbol{\lambda}/\boldsymbol{\mu}}^{\nu}(q)$ is given by
	\begin{align}
		\label{f-int}
		f_{\boldsymbol{\lambda}/\boldsymbol{\mu}}^{\nu}(q)
		=
		\frac{1}{b_{\nu}(q)}
		\displaystyle\frac{1}{(2 \pi \textbf{\emph{i}})^N} 
		\cdot
		\oint_C 
		\cdots
		\oint_C \mathcal{G}_{\boldsymbol{\lambda}/\boldsymbol{\mu}}(\textbf{\emph{x}}^{-1}; r^{-1 / 2} \boldsymbol{\mid} 0; 0)
		\prod_{1 \le i < j \le N}
		 \frac{x_j-x_i}{x_j-qx_i}
		\prod_{i=1}^N
		x_i^{\nu_i-1} dx_i,
	\end{align}
	where we have denoted $\textbf{\emph{x}} = (x_1, x_2, \ldots , x_N)$, and the contour $C$ is a positively oriented circle centered on the origin, and where we assume that $qC$ is contained in $C$ (equivalently, $|q| < 1$).
\end{cor}

\begin{proof}
	This follows from \Cref{fgq} and \Cref{fc}.
\end{proof}

\chapter{Proof of \Cref{thm:comb}}
\label{sec:proof}

In this chapter we provide the proof of \Cref{thm:comb}, which will proceed as follows. First, in \Cref{ssec:hybrid-model}, we write down vertex weights $w_x$ that are slightly more general than the $L_x$ weights from \eqref{wl}; see equations \eqref{hybrid1}--\eqref{hybrid2}. At the level of the underlying quantized affine Lie algebras, the $L_x$ weights are based on $U_q(\widehat{\mathfrak{sl}}(1|n))$, whereas the $w_x$ weights that we provide in \Cref{ssec:hybrid-model} are based on $U_q(\widehat{\mathfrak{sl}}(2|n))$. This is the first time in our text that we make use of a model of mixed boson-fermion nature for computational purposes. At a conceptual level, one may think of the $n$ fermionic species in \eqref{hybrid1}--\eqref{hybrid2} as the progenitors of a pre-fused LLT polynomial, while the single bosonic family present in \eqref{hybrid1}--\eqref{hybrid2} is the forerunner of a Hall--Littlewood polynomial; the model \eqref{hybrid1}--\eqref{hybrid2} thus provides us with the necessary algebraic framework for studying expansions of the form \eqref{pre-fused-exp}.

Second, in \Cref{ssec:hybrid-pf}, we construct a certain family of partition functions $Z^{\mathcal{A}}_{\boldsymbol{\lambda}/\boldsymbol{\mu}}$ using the weights of the model \eqref{hybrid1}--\eqref{hybrid2}. The Yang--Baxter equation provides us with exchange relations on the family $Z^{\mathcal{A}}_{\boldsymbol{\lambda}/\boldsymbol{\mu}}$ which, combined with an explicit initial condition, suffice to determine the family uniquely; see \Cref{ssec:properties}. In \Cref{ssec:evaluation} we solve these exchange relations and obtain an integral formula for $Z^{\mathcal{A}}_{\boldsymbol{\lambda}/\boldsymbol{\mu}}$.

The integral formula obtained in \Cref{ssec:evaluation} bears obvious resemblance to the right hand side of \eqref{f-int}, and can in fact be brought identically to that form, via a certain limiting procedure; the details of this procedure are the content of \Cref{ssec:procedure}. Finally, one needs to examine what happens to the partition function $Z^{\mathcal{A}}_{\boldsymbol{\lambda}/\boldsymbol{\mu}}$ under the same limiting procedure; this involves some delicate calculations at the level of the fully-fused $U_q(\widehat{\mathfrak{sl}}(2|n))$ model, and is the subject of \Cref{ssec:pf-procedure}. The final result of these computations is the combinatorial formula quoted in \Cref{thm:comb}, which is then matched with $f^{\nu}_{\boldsymbol{\lambda}/\boldsymbol{\mu}}$ courtesy of \eqref{f-int}; this completes the proof of \Cref{thm:comb}.

Here, for any $\boldsymbol{\lambda} \in \SeqSign_n$, we recall the sequence $\mathscr{S} (\boldsymbol{\lambda}) = \big( \textbf{S}_1 (\boldsymbol{\lambda}), \textbf{S}_2 (\boldsymbol{\lambda}), \ldots \big)$ of elements in $\{ 0, 1 \}^n$ from \Cref{Symmetric}.

\section{\texorpdfstring{$U_q \big( \widehat{\mathfrak{sl}}(2|n) \big)$}{} Model}
\label{ssec:hybrid-model}

We begin with the following definition for the weights of the vertex model we will use here. 

\begin{definition} 
	
\label{wxabcd}
 
Fix two $n$-tuples of integers $\textbf{A} = (A_1,A_2,\dots,A_{n+1}) \in \mathbb{Z}_{\ge 0}^{n + 1}$ and $\textbf{C} = (C_1,C_2,\dots,C_{n+1}) \in \mathbb{Z}_{\ge 0}^{n + 1}$ such that $A_k, B_k \in \{ 0, 1 \}$ for each $k \in [2, n + 1]$ and $A_1, C_1 \ge 0$. Further fix two integers $b, d \in [0, n + 1]$. For any complex number $x \in \mathbb{C}$, define the vertex weights
\begin{align}
	\label{hybrid1}
	\tikz{0.8}{
		\draw[lgray,line width=1.5pt] (-1,-1) -- (1,-1) -- (1,1) -- (-1,1) -- (-1,-1);
		\node[left] at (-1,0) {\footnotesize $b$};\node[right] at (1,0) {\footnotesize $d$};
		\node[below] at (0,-1) {\footnotesize $\textbf{A}$};\node[above] at (0,1) {\footnotesize $\textbf{C}$};
	}
	= w_x(\textbf{A},b;\textbf{C},d)
\end{align}

\noindent where the values of $w_x (\textbf{A}, b; \textbf{C}, d)$\index{W@$w_x (\textbf{A}, b; \textbf{C}, d)$} are tabulated below (recalling $\textbf{A}_i^+, \textbf{A}_i^-, \textbf{A}_{ij}^{+-}$ from \eqref{aij}), assuming that $1 \le i < j \le n + 1$:
\begin{align}
	\label{hybrid2}
	\begin{tabular}{|c|c|}
		\hline &
		\\
		\quad
		$w_x(\textbf{A},0;\textbf{A},0)$
		\quad
		&
		\quad
		$w_x(\textbf{A},i;\textbf{A},i)$
		\quad
		\\[0.5cm]
		\quad
		$=1$
		\quad
		& 
		\quad
		$=x(-q)^{\bm{1}_{i\ge 2} \cdot A_i}q^{\textbf{A}_{[i+1, n+1]}}$
		\quad
		\\[0.5cm]
		\hline &
		\\
		\quad
		$w_x(\textbf{A},0;\textbf{A}^{-}_i,i)$
		\quad
		& 
		\quad
		$w_x(\textbf{A},i;\textbf{A}^{+}_i,0)$
		\quad
		\\[0.5cm]
		\quad
		$=x(1-q^{A_i}) q^{\textbf{A}_{[i+1, n+1]}}$
		\quad
		&
		\quad
		$={\bm 1}_{i=1}+\bm{1}_{i \ge 2} {\bm 1}_{A_i=0}$
		\quad
		\\[0.5cm]
		\hline &
		\\
		\quad
		$w_x(\textbf{A},i;\textbf{A}^{+-}_{ij},j)$
		\quad
		&
		\quad
		$w_x(\textbf{A},j;\textbf{A}^{+-}_{ji},i)$
		\quad
		\\[0.5cm]  
		\quad
		$=x(1-q^{A_j}) q^{\textbf{A}_{[j+1, n+1]}}
		\Big( {\bm 1}_{i=1}+\bm{1}_{i \ge 2} {\bm 1}_{A_i=0} \Big)$
		\quad
		&
		\quad
		$=0$
		\quad
		\\[0.5cm]
		\hline
	\end{tabular} 
\end{align}

\noindent with $w_x (\textbf{A}, b; \textbf{C}, d) = 0$ for any $(\textbf{A}, b; \textbf{C}, d)$ not listed above. 

\end{definition} 

We once again view $\textbf{A}$, $b$, $\textbf{C}$, and $d$ as indexing the colors of the paths vertically entering, horizontally entering, vertically exiting, and horizontally exiting through a vertex, respectively. Observe here that, since we impose $A_k, C_k \in \{ 0, 1 \}$ for each $k \in [2, n + 1]$, the colors $\{ 2, 3, \ldots , n + 1 \}$ are all fermionic. However, since the integers $A_1$ and $C_1$ are permitted to be of unbounded size, the color $1$ is bosonic. Thus, the model \eqref{hybrid1}--\eqref{hybrid2} is of boson-fermion nature. In fact, the $w_x (\textbf{A}, b; \textbf{C}, d)$ are special limits of the $U_q \big( \widehat{\mathfrak{sl}} (2 | n) \big)$ weights $\mathcal{R}_{x, y}^{(2; n)} (\textbf{A}, \textbf{B}; \textbf{C}, \textbf{D})$ introduced and evaluated in \Cref{rijkh} and \Cref{rxyml}, respectively; we will see this more explicitly in \Cref{ssec:limits=0} below. 

\begin{rem} 
	
	\label{wxlx} 
	
	The $w_x (\textbf{A}, b; \textbf{C}, d)$ weights reduce to the $L_x (\textbf{A}, b; \textbf{C}, d)$ ones from \Cref{lz} by choosing $A_1 = C_1 = 0$; restricting $b,d$ to values in the set $\{2,\dots,n+1\}$; shifting all color labels downwards by one unit; and replacing $q$ with $t$. This can be verified directly, using the explicit forms for the $w_x$ and $L_x$ weights. It can also be viewed as a consequence of the facts that these $w_x$ and $L_x$ can both be obtained as the same limit of the $\mathcal{R}^{(2; n)}$ and $\mathcal{R}^{(1; n)}$ fused weights (from \Cref{rijkh}), respectively, and that the former degenerate to the latter in the absence of any color $1$ arrows.

\end{rem}

Next, similarly to as in \Cref{OperatorRow}, we will define a certain vector space and a family of operators on it. To that end, let $\mathfrak{A}(n)$\index{A@$\mathfrak{A} (n)$} denote the set of all $(n + 1)$-tuples $(A_1,A_2,\dots,A_{n+1}) \in \mathbb{Z}_{\ge 0}^{n + 1}$ such that $A_k \in \{ 0, 1 \}$, for each $k \in [2, n + 1]$. We associate with this a vector space $U (n)$\index{U@$U(n)$} (over $\mathbb{C}$) as follows:
	\begin{align*}
		U(n):= {\rm Span}_{\mathbb{C}}\big\{ | \textbf{A} \rangle \big\}_{\textbf{A} \in \mathfrak{A}(n)},
	\end{align*}

	\noindent and its dual 
	\begin{flalign*}
		U(n)^*:= {\rm Span}_{\mathbb{C}}\big\{ \langle \textbf{A} | \big\}_{\textbf{A} \in \mathfrak{A}(n)}.
	\end{flalign*}

	\noindent We introduce an inner product between $U(n)$ and $U(n)^*$ by imposing that $\langle \textbf{C} | \textbf{A} \rangle = \textbf{1}_{\textbf{A} = \textbf{C}}$. 

	Next, fix a nonnegative integer $K$ and consider the $(K+1)$-fold tensor products of such spaces, 
	\begin{align*}
		& \mathbb{U}^K = \mathbb{U}^K(n) := U_0 (n) \otimes U_1(n) \otimes 
		\cdots \otimes U_K (n); \\ 
		& (\mathbb{U}^K)^* = \mathbb{U}^K(n)^* := U_0 (n)^* \otimes U_1(n)^* \otimes 
		\cdots \otimes U_K (n)^*,
	\end{align*}\index{U@$\mathbb{U}^K (n)^*$} \index{U@$\mathbb{U}^K (n)^*$}
	where each $U_i (n)$ (and $U_i (n)^*$) denotes a copy of $U(n)$ (and $U(n)^*$, respectively). The inner product between $U(n)$ and $U(n)^*$ induces one between $\mathbb{U}^K$ and $(\mathbb{U}^K)^*$.
	
	For a $(K + 1)$-tuple $\mathscr{A} = (\textbf{A}_0,\textbf{A}_1,\dots,\textbf{A}_K)$ of basis vectors in $\mathfrak{A}(n)$, denote the elements
	\begin{flalign*} 
	| \mathscr{A} \rangle = | \textbf{A}_0 \rangle \otimes | \textbf{A}_1 \rangle \otimes \cdots \otimes | \textbf{A}_K \rangle \in \mathbb{U}^K; \qquad \langle \mathscr{A} | = \langle \textbf{A}_0 | \otimes \langle \textbf{A}_1 | \otimes \cdots \otimes \langle \textbf{A}_K | \in (\mathbb{U}^K)^*.
	\end{flalign*} \index{A@$\mid$$\mathscr{A} \rangle$}\index{C@$\langle \mathscr{C}$$\mid$}

	\noindent Then, $\mathbb{U}^K$ and $(\mathbb{U}^K)^*$ are spanned by such vectors of the form $| \mathscr{A} \rangle$ and $\langle  \mathscr{A} |$, respectively.
	
	We next define a family of transfer operators on $\mathbb{U}^K$ and $(\mathbb{U}^K)^*$, similarly to in \Cref{lwabcd} and \Cref{dci}.
	
	\begin{definition} 
		
	\label{ubk} 
	
	Fix a complex number $x \in \mathbb{C}$; two $(K + 1)$-tuples $\mathscr{A} = (\textbf{A}_0, \textbf{A}_1, \ldots , \textbf{A}_K)$ and $\mathscr{C} = (\textbf{C}_0, \textbf{C}_1, \ldots , \textbf{C}_K)$ of basis vectors in $\mathfrak{A} (n)$; and two indices $b, d \in [0, n + 1]$. Define
	\begin{flalign*}
	w_x (\mathscr{A}, b; \mathscr{C}, d) = \displaystyle\sum_{\mathfrak{J}} \displaystyle\prod_{i = 0}^K w_x (\textbf{A}_i, j_i; \textbf{C}_i, j_{i + 1}),
	\end{flalign*}

	\noindent where we sum over all sequences $\mathfrak{J} = (j_0, j_1, \ldots , j_K)$ of indices in $[0, n + 1]$ such that $j_0 = b$ and $j_{K + 1} = d$. By arrow conservation, this sum is supported on at most one term.
	
	Next, for any index $0 \in [1, n + 1]$, we introduce the \emph{transfer operator} $\mathfrak{B}_i (x): \mathbb{U}^K \rightarrow \mathbb{U}^K$ by, for any $(K + 1)$-tuple $\mathscr{A} = (\textbf{A}_0, \textbf{A}_1, \ldots , \textbf{A}_K)$ of basis vectors in $\mathfrak{A} (n)$, setting
	\begin{flalign*}
	\mathfrak{B}_i (x) | \mathscr{A} \rangle = \displaystyle\sum_{\mathscr{C}} w_x (\mathscr{A}, 0; \mathscr{C}, i) | \mathscr{C} \rangle,
	\end{flalign*} \index{B@$\mathfrak{B}_i (x)$}

	\noindent where we sum over all $(K + 1)$-tuples $\mathscr{C} = (\textbf{C}_0, \textbf{C}_1, \ldots , \textbf{C}_K)$ of basis vectors in $\mathfrak{A} (n)$. This induces a dual action of $\mathfrak{B}_i (x)$ on $(\mathbb{U}^K)^*$, defined by
	\begin{flalign*}
		\langle \mathscr{C} | \mathfrak{B}_i (x) = \displaystyle\sum_{\mathscr{A}} w_x (\mathscr{A}, 0; \mathscr{C}, i) \langle \mathscr{A} |,
	\end{flalign*}

	\noindent where now we sum over all $(K + 1)$-tuples $\mathscr{A} = (\textbf{A}_0, \textbf{A}_1, \ldots , \textbf{A}_K)$ of basis vectors in $\mathfrak{A} (n)$.

	\end{definition} 

	The operators $\mathfrak{B}_i (x)$ can also be expressed in diagrammatic notation, as follows:
	\begin{flalign}
		\label{B-row}
		\mathfrak{B}_i(x)
			| \mathscr{A} \rangle = \sum_{\mathscr{C}}
		\left(
		\tikz{1.2}{
			\draw[lgray,line width=1.5pt] (0.5,-0.5) -- (6.5,-0.5) -- (6.5,0.5) -- (0.5,0.5) -- (0.5,-0.5);
			\foreach\x in {1,...,5}{
				\draw[lgray,line width=1.5pt] (0.5+\x,-0.5) -- (0.5+\x,0.5);
			}
			\node[left] at (0.5,0) {$0$};\node[right] at (6.5,0) {$i$};
			\node[below] at (6,-0.5) {\footnotesize $\textbf{A}_K$};\node[above] at (6,0.5) {\footnotesize $\textbf{C}_K$};
			\node[below] at (5,-0.5) {\fs $\cdots$};\node[above] at (5,0.5) {\fs $\cdots$};
			\node[below] at (4,-0.5) {\fs $\cdots$};\node[above] at (4,0.5) {\fs $\cdots$};
			\node[below] at (3,-0.5) {\fs $\cdots$};\node[above] at (3,0.5) {\fs $\cdots$};
			\node[below] at (2,-0.5) {\footnotesize $\textbf{A}_1$};\node[above] at (2,0.5) {\footnotesize $\textbf{C}_1$};
			\node[below] at (1,-0.5) {\footnotesize $\textbf{A}_0$};\node[above] at (1,0.5) {\footnotesize $\textbf{C}_0$};
		}
		\right) 
		| \mathscr{C} \rangle.
	\end{flalign}

	\noindent Note that each face weight appearing on the right-hand side of \eqref{B-row} depends on the same parameter, $x$. The quantity
	
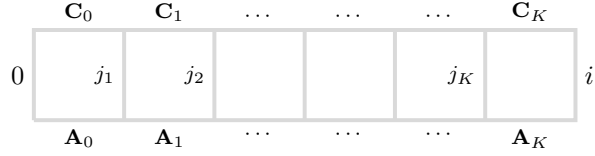
\begin{align*}
		\tikz{1.2}{
			\draw[lgray,line width=1.5pt] (0.5,-0.5) -- (6.5,-0.5) -- (6.5,0.5) -- (0.5,0.5) -- (0.5,-0.5);
			\foreach\x in {1,...,5}{
				\draw[lgray,line width=1.5pt] (0.5+\x,-0.5) -- (0.5+\x,0.5);
			}
			\node[left] at (0.5,0) {$0$};\node[right] at (6.5,0) {$i$};
			\node[below] at (6,-0.5) {\footnotesize $\textbf{A}_K$};\node[above] at (6,0.5) {\footnotesize $\textbf{C}_K$};
			\node[below] at (5,-0.5) {\fs $\cdots$};\node[above] at (5,0.5) {\fs $\cdots$};
			\node[below] at (4,-0.5) {\fs $\cdots$};\node[above] at (4,0.5) {\fs $\cdots$};
			\node[below] at (3,-0.5) {\fs $\cdots$};\node[above] at (3,0.5) {\fs $\cdots$};
			\node[below] at (2,-0.5) {\footnotesize $\textbf{A}_1$};\node[above] at (2,0.5) {\footnotesize $\textbf{C}_1$};
			\node[below] at (1,-0.5) {\footnotesize $\textbf{A}_0$};\node[above] at (1,0.5) {\footnotesize $\textbf{C}_0$};
			\node[left] at (1.5, 0) {\footnotesize $j_1$};
			\node[left] at (2.5, 0) {\footnotesize $j_2$};
			\node[left] at (5.5, 0) {\footnotesize $j_K$};
		}
	\end{align*}
	is a one-row partition function in the model \eqref{hybrid2}, and can be calculated by multiplying the weights of each face from left to right, noting that the integer values $j_i$ prescribed to all internal vertical edges are fixed by arrow conservation.

	The following exchange relation for the transfer operators $\mathfrak{B}_i (x)$ will be useful for us. For $K = 0$, it can either be verified directly or deduced from the Yang--Baxter equation \Cref{wabcdproduct} for the $\mathcal{R}_{x, 1}^{(2; n)}$ weights (since the content of \Cref{ssec:limits=0} below implies that the weights $w_x$ \eqref{hybrid2} are limits of these fused weights). For $K \ge 1$, it then follows from $K + 1$ applications of the Yang--Baxter equation. We omit further details of this proof, since it is very similar to that of \Cref{w2w}. 
	
\begin{prop}
	
	Fix complex numbers $x, y \in \mathbb{C}$; an integer $K \ge 0$; and indices $0 \le i < j \le n + 1$. The transfer operators \eqref{B-row} on $\mathbb{U}^K$ satisfy the exchange relations
	\begin{align}
		\label{b-com1}
		\mathfrak{B}_i(x) \mathfrak{B}_j(y) 
		= 
		\frac{(1-q)y}{y-qx} 
		\mathfrak{B}_i (y) \mathfrak{B}_j (x) 
		+ 
		\frac{y-x}{y-qx} 
		\mathfrak{B}_j (y) \mathfrak{B}_i(x).
	\end{align} 

\end{prop}

\section{\texorpdfstring{$U_q \big( \widehat{\mathfrak{sl}}(2|n) \big)$}{} Partition Functions}
\label{ssec:hybrid-pf}

In this section we define certain partition functions for the vertex model with weights $w_x$ from \Cref{wxabcd}, which will be useful for combinatorially interpreting the integral appearing in \eqref{f-int}. To that end, we first require some notation.

Fix an integer $M \ge 0$ and signature sequences $\boldsymbol{\lambda}, \boldsymbol{\mu} \in \SeqSign_{n; M}$. Further let $K$ be any integer such that $K \ge \max_{i \in [1, n]} \max \mathfrak{T} \big(\boldsymbol{\lambda}^{(i)} \big)$. Additionally fix two integers $m, N \ge 1$ as well as a set $\mathcal{A} = \{a_1 < \cdots < a_m\} \subset [1, N + m]$, and denote its complement by $\bar{\mathcal{A}} = \{\bar{a}_1 < \cdots < \bar{a}_N \} = [1,N+m] \setminus \mathcal{A}$. From these sets, for all $1 \le i \le N+m$ we define
\begin{align*}
	b_i(\mathcal{A})
	=
	\left\{
	\begin{array}{ll}
		0, \qquad & i = a_j,\ \ j \in [1,m],
		\\ \\
		1, \qquad & i = \bar{a}_j, \ \ j \in [1,N].
	\end{array}
	\right.
\end{align*} \index{B@$b_i (\mathcal{A})$}
We introduce a family\footnote{We use the word {\it family} in the sense of varying the set $\mathcal{A} \subset [1,N+m]$ over all $N+m \choose m$ possibilities.} of partition functions $Z_{\boldsymbol{\lambda}/\boldsymbol{\mu}}^{\mathcal{A}}(x_1,\dots,x_{N+m})$\index{Z@$Z_{\boldsymbol{\lambda} / \boldsymbol{\mu}}^{\mathcal{A}} (x_1, \ldots , x_{N + m})$} in the model \eqref{hybrid1}--\eqref{hybrid2} as follows:
\begin{align}
	\label{pf-definition2}
	Z_{\boldsymbol{\lambda}/\boldsymbol{\mu}}^{\mathcal{A}}\left(x_1,\dots,x_{N + m}\right)
	=
	\tikz{1.2}{
		\foreach\y in {0,...,5}{
			\draw[lgray,line width=1.5pt] (1.5,0.5+\y) -- (7.5,0.5+\y);
		}
		\foreach\x in {0,...,6}{
			\draw[lgray,line width=1.5pt] (1.5+\x,0.5) -- (1.5+\x,5.5);
		}
		%spectral parameters
		\node[text centered] at (2,1) {$x_1$};
		\node[text centered] at (4,1) {$\cdots$};
		\node[text centered] at (5,1) {$\cdots$};
		\node[text centered] at (7,1) {$x_1$};
		\node[text centered] at (2,2) {$x_2$};
		\node[text centered] at (4,2) {$\cdots$};
		\node[text centered] at (5,2) {$\cdots$};
		\node[text centered] at (7,2) {$x_2$};
		\node[text centered] at (2,5) {$x_{N+m}$};
		\node[text centered] at (4,5) {$\cdots$};
		\node[text centered] at (5,5) {$\cdots$};
		\node[text centered] at (7,5) {$x_{N+m}$};
		\node[text centered] at (2,3.1) {$\vdots$};
		\node[text centered] at (2,4.1) {$\vdots$};
		\node[text centered] at (7,3.1) {$\vdots$};
		\node[text centered] at (7,4.1) {$\vdots$};
		%bottom labels
		\node[below] at (2,0.5) {\footnotesize $(N,\textbf{e}_0)$};
		\node[below] at (3,0.5) {\footnotesize $\big( 0, \textbf{S}_1 (\boldsymbol{\mu}) \big)$};
		\node[above,text centered] at (4,0) {$\cdots$};
		\node[above,text centered] at (5,0) {$\cdots$};
		\node[below] at (7,0.5) {\footnotesize $\big( 0,\textbf{S}_K (\boldsymbol{\mu}) \big)$};
		%top labels
		\node[above] at (2,5.5) {\footnotesize $(0,\textbf{e}_0)$};
		\node[above] at (3,5.5) {\footnotesize $\big(0,\textbf{S}_1 (\boldsymbol{\lambda}) \big)$};
		\node[above,text centered] at (4,5.5) {$\cdots$};
		\node[above,text centered] at (5,5.5) {$\cdots$};
		\node[above] at (7,5.5) {\footnotesize $\big( 0, \textbf{S}_K (\boldsymbol{\lambda}) \big)$};
		%right labels
		\node[right] at (7.5,1) {$b_1(\mathcal{A})$};
		\node[right] at (7.5,2) {$b_2(\mathcal{A})$};
		\node[right] at (7.5,3.1) {$\vdots$};
		\node[right] at (7.5,4.1) {$\vdots$};
		\node[right] at (7.5,5) {$b_{N+m}(\mathcal{A})$};
		%left labels
		\node[left] at (1.5,1) {$0$};
		\node[left] at (1.5,2) {$0$};
		\node[left] at (1.5,3.1) {$\vdots$};
		\node[left] at (1.5,4.1) {$\vdots$};
		\node[left] at (1.5,5) {$0$};
	}
\end{align}
where the boundary conditions prescribed to the leftmost column inject $N$ bosonic arrows at the base of the lattice; all other all columns have analogous boundary conditions to those used to define the $G_{\boldsymbol{\lambda} / \boldsymbol{\mu}}$ functions from \Cref{pgpfph} (see also the left side of \Cref{fgpaths}), namely, they inject fermionic arrows at the base and eject fermionic arrows at the top of the lattice. Note that exactly $N$ of the integers $b_i(\mathcal{A})$ for $i \in [1, N + m]$ are equal to $1$; the remaining $m$ are equal to $0$. Hence the number of bosonic arrows injected into the lattice \eqref{pf-definition2} is equal to the number ejected from it, and the partition function \eqref{pf-definition2} is generically nonzero.

We may use the row operators \eqref{B-row} to write the partition function \eqref{pf-definition2} algebraically as
\begin{align}
	\label{Z-family2}
	Z^{\mathcal{A}}_{\boldsymbol{\lambda}/\boldsymbol{\mu}}
	(x_1,\dots,x_{N+m})
	=
	\Big\langle 0, \boldsymbol{\lambda} \Big| \mathfrak{B}_{b_{N + m} (\mathcal{A})} (x_{N + m}) \cdots \mathfrak{B}_{b_1 (\mathcal{A})}(x_1)
	\Big| N, \boldsymbol{\mu} \Big\rangle,
\end{align}

\noindent with row operators given by
\begin{align*}
	\mathfrak{B}_{b_i(\mathcal{A})}(x)
	=
	\left\{
	\begin{array}{ll}
		\mathfrak{B}_0(x), \qquad & i = a_j,\ \ j \in [1,m],
		\\ \\
		\mathfrak{B}_1(x), \qquad & i = \bar{a}_j, \ \ j \in [1,N].
	\end{array}
	\right.
\end{align*}
The vectors in \eqref{Z-family2} are defined as
\begin{align*}
	\Big\langle 0, \boldsymbol{\lambda} \Big|
	=
	\Big\langle 0,\textbf{e}_0\Big| 
	\otimes
	\bigotimes_{j=1}^K \Big\langle 0, \textbf{S}_j (\boldsymbol{\lambda}) \Big|
	\in 
	\mathbb{U}^K (n)^*; \qquad
	\Big| N,\boldsymbol{\mu} \Big\rangle 
	=
	\Big| N,\textbf{e}_0\Big\rangle 
	\otimes
	\bigotimes_{j=1}^K \Big| 0, \textbf{S}_j (\boldsymbol{\mu}) \Big\rangle 
	\in 
	\mathbb{U}^K (n).
\end{align*}

\section{Properties of \texorpdfstring{$Z^{\mathcal{A}}_{\boldsymbol{\lambda}/\boldsymbol{\mu}}$}{}}

\label{ssec:properties}

The partition functions \eqref{pf-definition2} satisfy two basic properties which are easily seen to characterize them. We derive these properties in the following two propositions. In what follows, only the $r = q^{-1 / 2}$ case of the function $\mathcal{G}_{\boldsymbol{\lambda} / \boldsymbol{\mu}} (\textbf{x}; r \boldsymbol{\mid} 0; 0)$ will be useful for us. So, to ease notation, we will abbreviate in the remainder of this chapter  
\begin{flalign}
	\label{gxlx} 
	\mathcal{G}_{\boldsymbol{\lambda} / \boldsymbol{\mu}} (\textbf{x}) = \mathcal{G}_{\boldsymbol{\lambda} / \boldsymbol{\mu}} (\textbf{x}; q^{-1 / 2} \boldsymbol{\mid} 0; 0).
\end{flalign}\index{G@$G_{\boldsymbol{\lambda} / \boldsymbol{\mu}} (\textbf{x}; \textbf{r} \boldsymbol{\mid} \textbf{y}; \textbf{s})$!$\mathcal{G}_{\boldsymbol{\lambda} / \boldsymbol{\mu}} (\textbf{x}; r \boldsymbol{\mid} 0; 0)$!$G_{\boldsymbol{\lambda} / \boldsymbol{\mu}} (\textbf{x})$}

\noindent As previously, we refer to $\mathcal{G}_{\boldsymbol{\lambda} / \boldsymbol{\mu}} (\textbf{x})$ as a pre-fused LLT function.

\begin{prop}
	\label{prop:init}
	
	Fix an integer $M \ge 1$ and signature sequences $\boldsymbol{\lambda}, \boldsymbol{\mu} \in \SeqSign_{n; M}$. When $\mathcal{A} = \{N+1,\dots,N+m\}$, we have
	\begin{align}
		\label{Z-init}
		Z^{\{N+1,\dots,N+m\}}_{\boldsymbol{\lambda}/\boldsymbol{\mu}}
		(x_1,\dots,x_{N+m})
		=
		q^{nMN} (q; q)_N 
		\prod_{j=1}^{N} x_j^{K+1} \cdot	\mathcal{G}_{\boldsymbol{\lambda}/\boldsymbol{\mu}}(x_{N+1},\dots,x_{N+m}),
	\end{align}
	where $\mathcal{G}_{\boldsymbol{\lambda}/\boldsymbol{\mu}}$ denotes a pre-fused LLT polynomial \eqref{gxlx}.
\end{prop}

\begin{proof}
	In the case $\mathcal{A} = \{N+1,\dots,N+m\}$, \eqref{Z-family2} becomes
	\begin{align*}
		Z^{\{N+1,\dots,N+m\}}_{\boldsymbol{\lambda}/\boldsymbol{\mu}}
		(x_1,\dots,x_{N+m})
		=
		\Big\langle 0, \boldsymbol{\lambda} \Big|
		\mathfrak{B}_0 (x_{N + m}) \cdots \mathfrak{B}_0(x_{N + 1})
		\mathfrak{B}_1 (x_N) \cdots \mathfrak{B}_1(x_1)
		\Big| N, \boldsymbol{\mu} \Big\rangle.
	\end{align*}
	Now consider the partition function
	
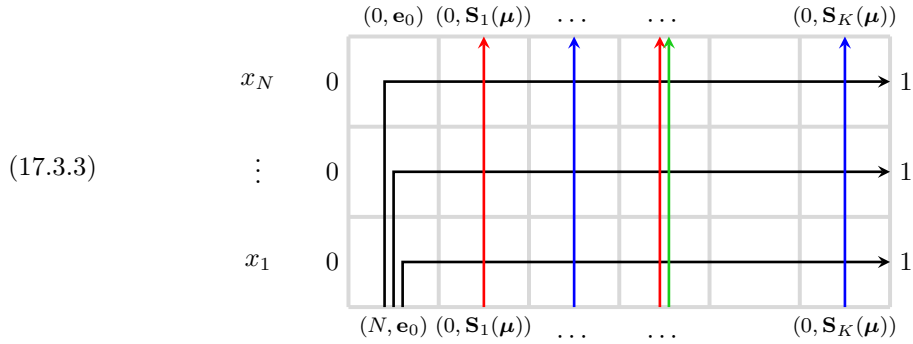
\begin{align}
		\label{frozen-config}
		\tikz{1.2}{
			\foreach\y in {0,...,3}{
				\draw[lgray,line width=1.5pt] (1.5,0.5+\y) -- (7.5,0.5+\y);
			}
			\foreach\x in {0,...,6}{
				\draw[lgray,line width=1.5pt] (1.5+\x,0.5) -- (1.5+\x,3.5);
			}
			%spectral parameters
			\node[text centered] at (0.5,1) {$x_1$};
			%\node[text centered] at (4,1) {$\cdots$};
			%\node[text centered] at (5,1) {$\cdots$};
			%\node[text centered] at (7,1) {$x_1$};
			%\node[text centered] at (2,2) {$x_2$};
			%\node[text centered] at (4,2) {$\cdots$};
			%\node[text centered] at (5,2) {$\cdots$};
			%\node[text centered] at (7,2) {$x_2$};
			\node[text centered] at (0.5,3) {$x_{N}$};
			%\node[text centered] at (4,3) {$\cdots$};
			%\node[text centered] at (5,3) {$\cdots$};
			%\node[text centered] at (7,3) {$x_{N}$};
			\node[text centered] at (0.5,2.1) {$\vdots$};
			%bottom labels
			\node[below] at (2,0.5) {\footnotesize $(N,\textbf{e}_0)$};
			\node[below] at (3,0.5) {\footnotesize $(0,\textbf{S}_1 (\boldsymbol{\mu}))$};
			\node[above,text centered] at (4,0) {$\cdots$};
			\node[above,text centered] at (5,0) {$\cdots$};
			\node[below] at (7,0.5) {\footnotesize $(0,\textbf{S}_K (\boldsymbol{\mu}))$};
			%top labels
			\node[above] at (2,3.5) {\footnotesize $(0,\textbf{e}_0)$};
			\node[above] at (3,3.5) {\footnotesize $(0,\textbf{S}_1 (\boldsymbol{\mu}))$};
			\node[above,text centered] at (4,3.5) {$\cdots$};
			\node[above,text centered] at (5,3.5) {$\cdots$};
			\node[above] at (7,3.5) {\footnotesize $(0,\textbf{S}_K (\boldsymbol{\mu}))$};
			%right labels
			\node[right] at (7.5,1) {$1$};
			\node[right] at (7.5,2) {$1$};
			\node[right] at (7.5,3) {$1$};
			%left labels
			\node[left] at (1.5,1) {$0$};
			\node[left] at (1.5,2) {$0$};
			\node[left] at (1.5,3) {$0$};
			%paths
			\draw[line width=1pt,->] (2.1,0.5) -- (2.1,1) -- (7.5,1);
			\draw[line width=1pt,->] (2,0.5) -- (2,2) -- (7.5,2);
			\draw[line width=1pt,->] (1.9,0.5) -- (1.9,3) -- (7.5,3);
			\draw[red,line width=1pt,->] (3,0.5) -- (3,3.5);
			\draw[blue,line width=1pt,->] (4,0.5) -- (4,3.5);
			\draw[red,line width=1pt,->] (4.95,0.5) -- (4.95,3.5);
			\draw[green,line width=1pt,->] (5.05,0.5) -- (5.05,3.5);
			\draw[blue,line width=1pt,->] (7,0.5) -- (7,3.5);
		}
	\end{align}
	which consists of a unique configuration with weight $q^{nMN} (q; q)_N \prod_{j=1}^{N} x_j^{K+1}$. Indeed, the $N$ bosonic paths which enter via the base of the leftmost column turn right into each of the $N$ rows, and completely saturate all edges in those rows. This in turn forces all $nM$ fermionic paths which enter at the base of the remaining columns to propagate in vertical straight lines to the top of the lattice. It is then clear that the partition function has a unique configuration and its weight is readily computed from \eqref{hybrid1}--\eqref{hybrid2}.
	
	Matching each of the $N$ rows in \eqref{frozen-config} with their algebraic equivalents \eqref{B-row}, we obtain
	\begin{align*}
		\mathfrak{B}_1 (x_N) \cdots \mathfrak{B}_1(x_1) \Big| N, \boldsymbol{\mu} \Big\rangle
		=
		q^{nMN}
		(q; q)_N
		\prod_{j=1}^{N} x_j^{K+1}
		\cdot
		\Big| 0, \boldsymbol{\mu} \Big\rangle.
	\end{align*}
	\noindent Thus, to establish \eqref{Z-init}, it suffices to verify
	\begin{align}
		\label{gb0} 
		\Big\langle 0, \boldsymbol{\lambda} \Big|
		\mathfrak{B}_0(x_{N+m}) \cdots \mathfrak{B}_0(x_{N+1})
		\Big| 0, \boldsymbol{\mu} \Big\rangle
		=
		\mathcal{G}_{\boldsymbol{\lambda}/\boldsymbol{\mu}}(x_{N+1},\dots,x_{N+m}).
	\end{align}
	
	\noindent To that end, recall the transfer operators $\mathsf{D} (x)$ from \Cref{dci} for the $L_x$ weights from \Cref{lz}. By \Cref{wxlx}, \Cref{ubk}, \Cref{lwabcd}, and \Cref{dci}, we have 
	\begin{flalign}
	\label{gb01}
		\Big\langle 0, \boldsymbol{\lambda} \Big|
		\mathfrak{B}_0(x_{N+m}) \cdots \mathfrak{B}_0(x_{N+1})
		\Big| 0, \boldsymbol{\mu} \Big\rangle
		=
		\big\langle \boldsymbol{\lambda} \big| \mathsf{D} (x_{N + m}) \cdots \mathsf{D} (x_{N + 1}) \big| \boldsymbol{\mu} \big\rangle.
	\end{flalign}

	\noindent Moreover, \Cref{dci}, \eqref{wl}, the last statement of \eqref{limitg}, and \Cref{fgdefinition} together imply
	\begin{flalign}
	\label{gb02} 
	\big\langle \boldsymbol{\lambda} \big| \mathsf{D} (x_{N + m}) \cdots \mathsf{D} (x_{N + 1}) \big| \boldsymbol{\mu} \big\rangle = \mathcal{G}_{\boldsymbol{\lambda}/\boldsymbol{\mu}}(x_{N+1},\dots,x_{N+m}).
	\end{flalign}

	\noindent Together, \eqref{gb01} and \eqref{gb02} establish \eqref{gb0} and thus the proposition.
\end{proof}

\begin{prop}
	\label{prop:rec}
	Fix a set $\mathcal{A} = \{a_1< \cdots < a_m\} \subset [1,N+m]$ and an integer $\alpha \in [1,N+m-1]$ such that $\alpha \in \mathcal{A}$ and $\alpha+1 \not\in \mathcal{A}$. Define the shifted set
	\begin{align*}
		\mathcal{A}(\alpha \rightarrow \alpha+1)
		=
		\{a_1 < \cdots < a_{\ell-1} < \alpha+1 < a_{\ell+1} < \cdots < a_m\},
	\end{align*}
	where $\ell \in [1,n]$ is the integer such that $a_{\ell} = \alpha$ (necessarily $a_{\ell+1} > \alpha+1$, since $\alpha+1 \not\in\mathcal{A}$). We have the recursion relation
	\begin{align}
		\label{Z-rec}
		Z^{\mathcal{A}}_{\boldsymbol{\lambda}/\boldsymbol{\mu}}
		=
		\frac{x_{\alpha+1}-qx_{\alpha}}{x_{\alpha+1}-x_{\alpha}}
		\cdot
		\mathfrak{s}_{\alpha}
		\left(
		Z^{\mathcal{A}(\alpha \rightarrow \alpha+1)}_{\boldsymbol{\lambda}/\boldsymbol{\mu}}
		\right)
		+
		\frac{(1-q)x_{\alpha+1}}{x_{\alpha}-x_{\alpha+1}}
		\cdot
		Z^{\mathcal{A}(\alpha \rightarrow \alpha+1)}_{\boldsymbol{\lambda}/\boldsymbol{\mu}},
		%=
		%\widetilde{T}_{\alpha} \left( Z^{\mathcal{A}}_{\boldsymbol{\lambda}/\boldsymbol{\mu}} \right).
	\end{align}
	where $\mathfrak{s}_{\alpha}$ denotes the transposition of the variables $(x_{\alpha},x_{\alpha+1})$.
\end{prop}

\begin{proof}
	Recall the commutation relation \eqref{b-com1}. After choosing $i=0$, $j=1$, $x = x_{\alpha}$, $y = x_{\alpha+1}$, and rearranging, \eqref{b-com1} becomes
	\begin{align}
		\label{B-relation}
		\mathfrak{B}_1(x_{\alpha + 1}) \mathfrak{B}_0(x_{\alpha})
		=
		\frac{x_{\alpha+1}-qx_{\alpha}}{x_{\alpha+1}-x_{\alpha}}
		\mathfrak{B}_0(x_{\alpha}) \mathfrak{B}_1(x_{\alpha + 1})
		+
		\frac{(1-q)x_{\alpha+1}}{x_{\alpha}-x_{\alpha+1}}
		\mathfrak{B}_0(x_{\alpha + 1}) \mathfrak{B}_1 (x_{\alpha}).
	\end{align}
	We use this relation in the algebraic formula \eqref{Z-family2} for $Z^{\mathcal{A}}_{\boldsymbol{\lambda}/\boldsymbol{\mu}}(x_1,\dots,x_{N+m})$; it affects the operators standing at positions $\alpha$ and $\alpha+1$ in the operator product, namely, $\mathfrak{B}_1(x_{\alpha + 1}) \mathfrak{B}_0(x_{\alpha})$. One sees that the two terms on the right hand side of \eqref{B-relation} correspond to the respective two terms on the right hand side of \eqref{Z-rec}, and the proof is completed.
\end{proof}

\begin{rem} 
	
	\label{zza}

	\Cref{prop:rec} allows us to incrementally increase the elements of the set
	$\mathcal{A}$ which indexes $Z^{\mathcal{A}}_{\boldsymbol{\lambda}/\boldsymbol{\mu}}$, at the expense of generating two terms for each application of the rule \eqref{Z-rec}. At the end of this ordering procedure, one obtains a sum over partition functions of the form $Z^{\{N+1,\dots,N+m\}}_{\boldsymbol{\lambda}/\boldsymbol{\mu}}$; these are known explicitly by \Cref{prop:init}. \Cref{prop:rec} and \Cref{prop:init} thus completely determine the family \eqref{Z-family2}.\footnote{In fact, these propositions overdetermine the family \eqref{Z-family2}, since there are multiple ways of performing the ordering $\mathcal{A} \rightarrow \{N+1,\dots,N+m\}$. The fact that all ways of doing this ordering lead to the same answer is a consequence of the Yang--Baxter equation.}
\end{rem}

\section{Evaluation of the Partition Function \texorpdfstring{$Z^{\mathcal{A}}_{\boldsymbol{\lambda}/\boldsymbol{\mu}}$}{}}
\label{ssec:evaluation}

In this section we solve for $Z^{\mathcal{A}}_{\boldsymbol{\lambda}/\boldsymbol{\mu}}$ that satisfies the two defining relations \eqref{Z-init} and \eqref{Z-rec}. Before we do so, it is convenient to define a version of the partition function \eqref{pf-definition2} with reciprocated variables:
\begin{align}
	\label{Z-recip}
	\bar{Z}^{\mathcal{A}}_{\boldsymbol{\lambda}/\boldsymbol{\mu}}
	(x_1,\dots,x_{N+m})
	=
	\prod_{i=1}^{N+m} x_i^{K+1}
	\cdot
	Z^{\mathcal{A}}_{\boldsymbol{\lambda}/\boldsymbol{\mu}}
	(x_1^{-1},\dots,x_{N+m}^{-1}).
\end{align}\index{Z@$\bar{Z}_{\boldsymbol{\lambda} / \boldsymbol{\mu}}^{\mathcal{A}} (x_1, \ldots , x_{N + m})$}
In terms of these new conventions, \eqref{Z-init} and \eqref{Z-rec} translate to
\begin{align}
	\label{Z-init2}
	\bar{Z}^{\{N+1,\dots,N+m\}}_{\boldsymbol{\lambda}/\boldsymbol{\mu}}
	(x_1,\dots,x_{N+m})
	=
	q^{nMN} (q; q)_N
	\prod_{i=1}^m (x_{N+i})^{K + 1} \cdot \mathcal{G}_{\boldsymbol{\lambda}/\boldsymbol{\mu}}(x_{N+1}^{-1},\dots,x_{N+m}^{-1}),
\end{align}
and
\begin{align}
	\label{Z-rec22}
	\bar{Z}^{\mathcal{A}}_{\boldsymbol{\lambda}/\boldsymbol{\mu}}
	=
	\frac{x_{\alpha}-qx_{\alpha+1}}{x_{\alpha}-x_{\alpha+1}}
	\cdot
	\mathfrak{s}_{\alpha}
	\left(
	\bar{Z}^{\mathcal{A}(\alpha\rightarrow \alpha+1)}_{\boldsymbol{\lambda}/\boldsymbol{\mu}}
	\right)
	+
	\frac{(1-q)x_{\alpha}}{x_{\alpha+1}-x_{\alpha}}
	\cdot
	\bar{Z}^{\mathcal{A}(\alpha\rightarrow \alpha+1)}_{\boldsymbol{\lambda}/\boldsymbol{\mu}}.
\end{align}
We may also write \eqref{Z-rec22} as
\begin{align}
	\label{Z-rec2}
	\bar{Z}^{\mathcal{A}}_{\boldsymbol{\lambda}/\boldsymbol{\mu}}
	=
	T_{\alpha}
	\left( \bar{Z}^{\mathcal{A}(\alpha\rightarrow \alpha+1)}_{\boldsymbol{\lambda}/\boldsymbol{\mu}} \right),
	\qquad
	T_{\alpha}
	=
	q-
	\frac{x_{\alpha}-qx_{\alpha+1}}{x_{\alpha}-x_{\alpha+1}}
	(1-\mathfrak{s}_{\alpha}),
\end{align}
where $T_{\alpha}$ denotes a Hecke algebra generator in its polynomial representation (which we recall from \Cref{Polynomialf}).

Now we can state the following contour integral formula for $\bar{Z}_{\boldsymbol{\lambda} / \boldsymbol{\mu}}^{\mathcal{A}}$. In the below, we recall the pre-fused LLT polynomial $\mathcal{G}_{\boldsymbol{\lambda} / \boldsymbol{\mu}}$ from \eqref{gxlx}. 

\begin{thm}
	The partition function \eqref{pf-definition2} is given by the following multiple integral formula:
	\begin{flalign}
		\label{Z-integral}
		\begin{aligned}
		\bar{Z}^{\mathcal{A}}_{\boldsymbol{\lambda}/\boldsymbol{\mu}}
		(x_1,\dots,x_{N+m})
		=
		\displaystyle\frac{q^{nMN}}{(2 \pi \textbf{\emph{i}})^m} (q; q)_N
		\prod_{i=1}^m x_{a_i}
		& \oint_C 
		\cdots
		\oint_C
		\mathcal{G}_{\boldsymbol{\lambda}/\boldsymbol{\mu}}(w_1^{-1},\dots,w_m^{-1})
		\prod_{1 \le i < j \le m}
		\frac{w_i-w_j}{w_i-qw_j} \\
		& \times 
		\prod_{i=1}^m
		\prod_{j=a_i+1}^{N+m}
		\frac{w_i-qx_j}{w_i-x_j}
		\cdot
		\prod_{i=1}^m
		\frac{w_i^K dw_i}{w_i-x_{a_i}},
		\end{aligned}
	\end{flalign}
	where $C$ is a positively oriented contour that surrounds all points in $qC$ and in the set $\{0\} \cup \{x_1,\dots,x_{N+m} \}$. 
\end{thm}

\begin{proof}
	
	By \Cref{zza}, it suffices to check that \eqref{Z-integral} satisfies the properties \eqref{Z-init2} and \eqref{Z-rec2}. 
	
	\medskip
	\noindent
	\underline{\it Property \eqref{Z-init2}.} Setting $a_i = N+i$ for all $i \in [1, m]$, formula \eqref{Z-integral} becomes
	\begin{flalign}
		\label{Z-init-int}
		\begin{aligned} 
		\bar{Z}^{\mathcal{A}}_{\boldsymbol{\lambda}/\boldsymbol{\mu}}
		(x_1,\dots,x_{N+m})
		=
		\displaystyle\frac{q^{nMN}}{(2 \pi \textbf{i})^m} (q; q)_N
		\prod_{i=1}^m x_{N+i}
		\oint_C 
		\cdots
		\oint_C 
		\mathcal{G}_{\boldsymbol{\lambda}/\boldsymbol{\mu}}(w_1^{-1},\dots,w_m^{-1}) 
		\prod_{1 \le i < j \le m}
		\frac{w_i-w_j}{w_i-qw_j} 
		\\
		\times
		\prod_{i=1}^m
		\prod_{j=N+i+1}^{N+m}
		\frac{w_i-qx_j}{w_i-x_j}
		\cdot
		\prod_{i=1}^m
		\frac{w_i^K dw_i}{w_i-x_{N+i}}.
	\end{aligned} 
	\end{flalign}
	The integrals in \eqref{Z-init-int} may now be performed sequentially, starting with the integral over $w_m$, and working through to $w_1$. We see that the integrand of \eqref{Z-init-int} picks up a simple pole at $w_m = x_{N+m}$, and has no other singularity in that variable. Indeed, by virtue of the fact that $C$ is contained in $q^{-1} C$, the points $q^{-1} w_i$ for $i \in [1, m - 1]$ are not enclosed by the $w_m$ integration contour; furthermore, the integrand is non-singular at $w_m=0$, since $\mathcal{G}_{\boldsymbol{\lambda}/\boldsymbol{\mu}}(w_1^{-1},\dots,w_m^{-1})$ is at most of degree $K$ in $w_m^{-1}$ (which can can be deduced directly from the vertex model interpretation for $\mathcal{G}_{\boldsymbol{\lambda} / \boldsymbol{\mu}} (w_1^{-1}, w_2^{-1}, \ldots , w_m^{-1})$). Evaluating the residue of the pole at $w_m = x_{N+m}$, we find that
	\begin{flalign*}
		\bar{Z}^{\mathcal{A}}_{\boldsymbol{\lambda}/\boldsymbol{\mu}}
		(x_1,\dots,x_{N+m})
		& =
		\displaystyle\frac{q^{nMN}}{(2 \pi \textbf{i})^{m - 1}} (q; q)_N
		\prod_{i=1}^{m-1} x_{N+i} \cdot (x_{N+m})^{K+1}
		\\
		& \qquad \times
		\oint_C 
		\cdots
		\oint_C 
		\mathcal{G}_{\boldsymbol{\lambda}/\boldsymbol{\mu}}(w_1^{-1},\dots,w_{n-1}^{-1},x_{N+m}^{-1})
		\prod_{1 \le i < j \le m-1}
		 \frac{w_i-w_j}{w_i-qw_j}
		\\
		& \qquad \qquad \qquad \times 
		\prod_{i=1}^{m-1}
		\prod_{j=N+i+1}^{N+m}
		\frac{w_i-qx_j}{w_i-x_j}
		\cdot
		\prod_{i=1}^{m-1}
		\frac{w_i-x_{N+m}}{w_i-qx_{N+m}} \frac{w_i^K dw_i}{w_i-x_{N+i}}.
	\end{flalign*}
	After some cancellation of factors within the integrand, this simplifies to
	\begin{flalign*}
		\bar{Z}^{\mathcal{A}}_{\boldsymbol{\lambda}/\boldsymbol{\mu}}
		(x_1,\dots,x_{N+m})
		& =
		\displaystyle\frac{q^{nMN}}{(2 \pi \textbf{i})^m} (q; q)_N
		\prod_{i=1}^{m - 1} x_{N+i} \cdot (x_{N+m})^{K+1} \\
		& \qquad \times 
		\oint_C 
		\cdots
		\oint_C 
		\mathcal{G}_{\boldsymbol{\lambda}/\boldsymbol{\mu}}(w_1^{-1},\dots,w_{n-1}^{-1},x_{N+m}^{-1})
		\prod_{1 \le i < j \le m-1}
		 \frac{w_i-w_j}{w_i-qw_j} \\
		 & \qquad \qquad \qquad \times 
		\prod_{i=1}^{m-1}
		\prod_{j=N+i+1}^{N+m-1}
		\frac{w_i-qx_j}{w_i-x_j}
		\cdot
		\prod_{i=1}^{m-1}
		\frac{w_i^K dw_i}{w_i-x_{N+i}}.
	\end{flalign*}
	We are left with an $(m-1)$-fold integral that has the same structure as that which we started with. Indeed, the $w_{m-1}$ variable develops a unique simple pole at $w_{n-1} = x_{N+m-1}$; this integration may be explicitly performed, similarly to above. Iterating this procedure all the way down to $w_1$, we arrive at the desired result, \eqref{Z-init2}.

	\medskip
	\noindent
	\underline{\it Property \eqref{Z-rec2}.}
	The proof of \eqref{Z-rec2} makes use of the elementary identity
	\begin{align}
		\label{hecke-elem}
		T_{\alpha} \left(  \frac{x_{\alpha+1}}{w-x_{\alpha+1}} \right)
		&=
		\frac{x_{\alpha}-qx_{\alpha+1}}{x_{\alpha}-x_{\alpha+1}}
		\cdot
		\frac{x_{\alpha}}{w-x_{\alpha}}
		+
		\frac{(1-q)x_{\alpha}}{x_{\alpha+1}-x_{\alpha}}
		\cdot
		\frac{x_{\alpha+1}}{w-x_{\alpha+1}}
		%\\
		%&=
		%\frac{x_i(x_i-qx_{i+1})(w-x_{i+1})-(1-q)x_i x_{i+1} (w-x_i)}
		%{(x_i-x_{i+1})(w-x_i)(w-x_{i+1})}
		%\\
		%&=
		%\frac{x_i(x_i-x_{i+1})(w-qx_{i+1})}{(x_i-x_{i+1})(w-x_i)(w-x_{i+1})}
		=
		\frac{x_{\alpha}(w-qx_{\alpha+1})}{(w-x_{\alpha})(w-x_{\alpha+1})},
	\end{align}
	as well as the fact that 
	\begin{align}
		\label{hecke-sym}
		T_{\alpha} \Big( f(x_{\alpha},x_{\alpha+1}) g(x_{\alpha},x_{\alpha+1}) \Big)
		=
		f(x_{\alpha},x_{\alpha+1})
		\cdot
		T_{\alpha} \Big( g(x_{\alpha},x_{\alpha+1}) \Big)
	\end{align}
	for functions $f(x_{\alpha},x_{\alpha+1})$ which are symmetric in $(x_{\alpha},x_{\alpha+1})$. Choose $\alpha \in [1,N+m-1]$ such that $\alpha \in \mathcal{A}$ and $\alpha+1 \not\in \mathcal{A}$. Further, fix $\ell \in [1, m]$ such that $a_{\ell}=\alpha$ (and note that $a_{\ell+1} > \alpha+1$). We act on $\bar{Z}^{\mathcal{A}(\alpha\rightarrow \alpha+1)}_{\boldsymbol{\lambda}/\boldsymbol{\mu}}$, as given by \eqref{Z-integral}, with the operator $T_{\alpha}$. Denoting $\tilde{a}_i = a_i$ for each $i \in [1, m] \setminus \{ \ell \}$ and $\tilde{a}_{\ell} = a_{\ell} + 1 = \alpha + 1$, the product
	\begin{align*}
		\prod_{i\not=\ell} x_{a_i}
		\cdot
		\prod_{i=1}^{m}
		\prod_{j = \tilde{a}_i+1}^{N+m}
		\frac{w_i-q x_j}{w_i-x_j}
		\cdot
		\prod_{i\not=\ell}
		\frac{1}{w_i-x_{a_i}}
	\end{align*}
	is symmetric in $(x_{\alpha},x_{\alpha+1})$; by virtue of \eqref{hecke-sym}, we therefore have
	\begin{flalign*}
		T_{\alpha} &
		\left(
		\bar{Z}^{\mathcal{A}(\alpha\rightarrow \alpha+1)}_{\boldsymbol{\lambda}/\boldsymbol{\mu}}
		(x_1,\dots,x_{N+m})
		\right) \\
		& =
		\displaystyle\frac{q^{nMN}}{(2 \pi \textbf{i})^m} (q; q)_N
		\prod_{i \ne \ell} x_{a_i}
		\oint_C 
		\cdots
		\oint_C 
		\mathcal{G}_{\boldsymbol{\lambda}/\boldsymbol{\mu}}(w_1^{-1},\dots,w_m^{-1})
		\prod_{1 \le i < j \le m}
		\frac{w_i-w_j}{w_i-qw_j}
		\\
		& \qquad \qquad \qquad \times 
		\prod_{i=1}^m
		\prod_{j = \tilde{a}_i+1}^{N+m}
		\frac{w_i-q x_j}{w_i-x_j}
		\cdot
		\prod_{i\not=\ell}
		\frac{1}{w_i-x_{a_i}}
		\cdot
		T_{\alpha} \left( \frac{x_{\alpha+1}}{w_{\ell}-x_{\alpha+1}} \right) \displaystyle\prod_{i = 1}^m w_i^K dw_i.
	\end{flalign*}
	Using the relation \eqref{hecke-elem}, we conclude that
	\begin{flalign*}
		T_{\alpha} &
		\left(
		\bar{Z}^{\mathcal{A}(\alpha\rightarrow \alpha+1)}_{\boldsymbol{\lambda}/\boldsymbol{\mu}}
		(x_1,\dots,x_{N+m})
		\right) \\
		& =
		\displaystyle\frac{q^{nMN}}{(2 \pi \textbf{i})^m} (q; q)_N
		\prod_{i\not=\ell} x_{a_i}
		\cdot
		x_{\alpha}
		\oint_C 
		\cdots
		\oint_C 
		\mathcal{G}_{\boldsymbol{\lambda}/\boldsymbol{\mu}}(w_1^{-1},\dots,w_m^{-1})
		\prod_{1 \le i < j \le m}
		\frac{w_i-w_j}{w_i-qw_j}
		\\
		& \qquad \qquad \qquad \times
		\prod_{i=1}^m
		\prod_{j = \tilde{a}_i+1}^{N+m}
		\frac{w_i-q x_j}{w_i-x_j}
		\cdot
		\frac{w_{\ell}-q x_{\alpha+1}}{w_{\ell}-x_{\alpha+1}}
		\cdot
		\prod_{i\not=\ell}
		\frac{1}{w_i-x_{a_i}}
		\cdot
		\frac{1}{w_{\ell}-x_{\alpha}}
		\displaystyle\prod_{i = 1}^m w_i^K dw_i \\
		& =
		\displaystyle\frac{q^{nMN}}{(2 \pi \textbf{i})^m} (q; q)_N
		\prod_{i = 1}^m x_{a_i}
		\oint_C  1
		\cdots
		\oint_C 
		\mathcal{G}_{\boldsymbol{\lambda}/\boldsymbol{\mu}}(w_1^{-1},\dots,w_m^{-1})
		\prod_{1 \le i < j \le m}
		\frac{w_i-w_j}{w_i-qw_j}
		\\
		& \qquad \qquad \qquad \qquad \qquad \qquad \quad \times
		\prod_{i=1}^m
		\prod_{j = a_i+1}^{N+m}
		\frac{w_i-q x_j}{w_i-x_j}
		\cdot
		\prod_{i = 1}^m
		\frac{1}{w_i-x_{a_i}}
		\displaystyle\prod_{i = 1}^m w_i^K dw_i,
	\end{flalign*}
	where the final expression is recognized as $\bar{Z}^{\mathcal{A}}_{\boldsymbol{\lambda}/\boldsymbol{\mu}} (x_1,\dots,x_{N+m})$. This validates \eqref{Z-rec2}.
\end{proof}

\section{Degree Counting}

In this section we pause to determine the degree of $\bar{Z}^{\mathcal{A}}_{\boldsymbol{\lambda}/\boldsymbol{\mu}}$, as a polynomial in $(x_1,\dots,x_{N+m})$.

\begin{prop}
	\label{prop:degree}
	$\bar{Z}^{\mathcal{A}}_{\boldsymbol{\lambda}/\boldsymbol{\mu}}(x_1,\dots,x_{N+m})$ is a homogeneous polynomial in $(x_1,\dots,x_{N+m})$ of degree $m (K+1)-|\boldsymbol{\lambda}|+|\boldsymbol{\mu}|$, and is divisible by $x_{a_i}$ for all $i \in [1,m]$.
\end{prop}

\begin{proof}
	Examining the weights \eqref{hybrid2}, we see that vertices (or faces, which we use to diagrammatically represent them) whose right vertical edge is occupied by a path (of any color) receive a weight of $x$; conversely, vertices in which this edge is unoccupied have degree zero in $x$. Working out the degree of $Z^{\mathcal{A}}_{\boldsymbol{\lambda}/\boldsymbol{\mu}}(x_1,\dots,x_{N+m})$ is then a matter of counting the total number of right vertical edges which are occupied by a path; this number is invariant across all lattice configurations, which immediately yields the homogeneity claim. The number of right vertical edges occupied by fermionic colors $\{2,\dots,n+1\}$ is equal to $|\boldsymbol{\lambda}| - |\boldsymbol{\mu}|$. The number of right vertical edges occupied by the bosonic color $1$ is equal to $N(K+1)$. It follows that
	\begin{align*}
		{\rm deg}
		Z^{\mathcal{A}}_{\boldsymbol{\lambda}/\boldsymbol{\mu}}(x_1,\dots,x_{N+m})
		=
		|\boldsymbol{\lambda}| - |\boldsymbol{\mu}| + N(K+1),
	\end{align*}
	and accordingly (by \eqref{Z-recip}),
	\begin{align*}
		{\rm deg}
		\bar{Z}^{\mathcal{A}}_{\boldsymbol{\lambda}/\boldsymbol{\mu}}(x_1,\dots,x_{N+m})
		=
		(N+m)(K+1)-|\boldsymbol{\lambda}| + |\boldsymbol{\mu}|-N(K+1)
		=
		m (K+1)-|\boldsymbol{\lambda}| + |\boldsymbol{\mu}|.
	\end{align*}
	For the divisibility claim, note that the effect of reciprocating all rapidities and multiplying the partition function by $\prod_{i=1}^{N+m} x_i^{K+1}$ (this is exactly one $x_i$ factor for every vertex in row $i$, for all $i \in [1, N + m]$) is to reverse the above rule; namely, vertices in \eqref{hybrid2} with an empty right vertical edge will receive a weight of $x$, while those with an occupied right edge have degree zero in $x$. As a result, any row of the partition function $\bar{Z}^{\mathcal{A}}_{\boldsymbol{\lambda}/\boldsymbol{\mu}} \left(x_1,\dots,x_{N+m}\right)$ that has an empty right vertical edge will be divisible by the rapidity variable associated to that row; the rows which have this property are precisely those encoded by $\mathcal{A} = \{a_1<\cdots<a_m\}$.
\end{proof}

\section{Limiting Procedure}
\label{ssec:procedure}

The integral formula \eqref{Z-integral} for $Z^{\mathcal{A}}_{\boldsymbol{\lambda}/\boldsymbol{\mu}}$ depends on a set of variables $(x_1,\dots,x_{N+m})$, which up to this point were indeterminates. We now provide a certain limiting procedure on 
$(x_1,\dots,x_{N+m})$ which transforms \eqref{Z-integral} to match directly 
with \eqref{f-int}.

Let us define
\begin{align}
	\label{Z-norm}
	\mathfrak{Z}^{\mathcal{A}}_{\boldsymbol{\lambda}/\boldsymbol{\mu}}
	(x_1,\dots,x_{N+m})
	=
	q^{-nMN} (q; q)_N^{-1} \prod_{i=1}^m x_{a_i}^{-1} \cdot \bar{Z}^{\mathcal{A}}_{\boldsymbol{\lambda}/\boldsymbol{\mu}}
		(x_1,\dots,x_{N+m}).
\end{align}\index{Z@$\mathfrak{Z}_{\boldsymbol{\lambda} / \boldsymbol{\mu}}^{\mathcal{A}} (x_1, \ldots , x_{N + m})$}

\subsection{Sending $x_{a_1},\dots,x_{a_m}$ to $0$}

We begin by taking each of the variables $x_{a_1},\dots,x_{a_m}$ in \eqref{Z-norm} to zero; using \eqref{Z-integral}, this limit is computed as
\begin{flalign}
	\label{xa-to-zero}
	\begin{aligned}
	\lim_{x_{a_1} \rightarrow 0}
	\cdots
	\lim_{x_{a_m} \rightarrow 0}
	\mathfrak{Z}^{\mathcal{A}}_{\boldsymbol{\lambda}/\boldsymbol{\mu}}
	(x_1,\dots,x_{N+m})
	=
	 \frac{1}{(2\pi \textbf{i})^m} \displaystyle\oint_C 
	\cdots
	\oint_C 	
	\mathcal{G}_{\boldsymbol{\lambda}/\boldsymbol{\mu}}(w_1^{-1},\dots,w_m^{-1})
	\prod_{1 \le i < j \le m}
	 \frac{w_i-w_j}{w_i-qw_j}  \\
	\times \prod_{i=1}^m
	\prod_{\substack{j=a_i+1 \\ j \not= a_{i+1},\dots,a_m}}^{N+m}
	\frac{w_i-qx_j}{w_i-x_j}
	\prod_{i=1}^m
	w_i^{K - 1} dw_i.
	\end{aligned}
\end{flalign}

\subsection{Principal Specializations}
For the remaining variables $x_b$ with $b \not\in \mathcal{A}$, we make the following choice. We let  $a_1=1$, and write $a_{m+1} = N+m+1$ by agreement. We further assume for each $i \in [1,  m]$ that $a_{i+1}-a_i-1=L g_i$, where $L \ge 1$ and $g_i \ge 0$ are fixed integers. Then for each $i \in [1,  m]$ we set
\begin{align}
	\label{cov-1}
	\left( x_{a_i+1},\dots,x_{a_{i+1}-1} \right)
	=
	\left( y^{(1)}_1,\dots,y^{(1)}_{L} \right)
	\cup
	\cdots
	\cup
	\left( y^{(g_i)}_1,\dots,y^{(g_i)}_{L} \right)
\end{align}
with 
\begin{align}
	\label{cov-2}
	y^{(j)}_c = y q^{L-c},
	\qquad \text{for all $j \in [1, g_i]$ and $c \in [1, L]$}.
\end{align}
This amounts to splitting the group of variables $x_b$ for $b \in (a_i,a_{i+1})$ into $g_i$ bundles each of cardinality $L$, and principally specializing within the bundles. Note that we use the same indeterminate, $y$, as base of every geometric progression. In the equations that follow, we denote change of variables \eqref{cov-1}--\eqref{cov-2} by the symbol $\dagger$. Applying the change of variables \eqref{cov-1}--\eqref{cov-2} to \eqref{xa-to-zero}, we obtain
\begin{flalign}
	\label{dagger}
	\begin{aligned}
	& \left[ 
	\lim_{x_{a_1} \rightarrow 0}
	\cdots
	\lim_{x_{a_m} \rightarrow 0}
	\mathfrak{Z}^{\mathcal{A}}_{\boldsymbol{\lambda}/\boldsymbol{\mu}}
	(x_1,\dots,x_{N+m})
	\right]^{\dagger}
	\\
	& \qquad \qquad \qquad =
	\displaystyle\frac{1}{(2 \pi \textbf{i})^m} \oint_C 
	\cdots
	\oint_C 
	\mathcal{G}_{\boldsymbol{\lambda}/\boldsymbol{\mu}}(w_1^{-1},\dots,w_m^{-1}) 
	\prod_{1 \le i < j \le m}
	 \frac{w_i-w_j}{w_i-qw_j} \\
	 & \qquad \qquad \qquad \qquad \qquad \qquad \times  \prod_{1 \le c\le d \le m}
	\left( \frac{w_c- y q^{L}}{w_c-y} \right)^{g_d}
	\prod_{i=1}^{m}
	w_i^{K-1} dw_i.
	\end{aligned} 
\end{flalign}
Due to \Cref{prop:degree} and the fact that we divide by $x_{a_1},\dots,x_{a_m}$ in the definition \eqref{Z-norm}, when \eqref{dagger} is viewed as a polynomial in $y$ it must take the form
\begin{align}
	\label{homog-y}
	\left[
	\lim_{x_{a_1} \rightarrow 0}
	\cdots
	\lim_{x_{a_m} \rightarrow 0}
	\mathfrak{Z}^{\mathcal{A}}_{\boldsymbol{\lambda}/\boldsymbol{\mu}}
	(x_1,\dots,x_{N+m})
	\right]^{\dagger}
	=
	y^{mK-|\boldsymbol{\lambda}|+|\boldsymbol{\mu}|}
	\cdot
	C_{\{g_1,\dots,g_m\}}(q^{L})
\end{align}
where $C_{\{g_1,\dots,g_m\}}$ is independent of $y$, but depends polynomially on $q^{L}$.

\subsection{Sending $q^{L} \rightarrow \infty$}

Using \eqref{dagger}, we compute
\begin{multline}
	\label{limit2}
	(-y)^{-\sum_{i=1}^m i g_i}
	\lim_{q^L \rightarrow \infty}
	q^{-L \sum_{i=1}^m i g_i}
	\left[
	\lim_{x_{a_1} \rightarrow 0}
	\cdots
	\lim_{x_{a_m} \rightarrow 0}
	\mathfrak{Z}^{\mathcal{A}}_{\boldsymbol{\lambda}/\boldsymbol{\mu}}
	(x_1,\dots,x_{N+m})
	\right]^{\dagger}
	\\
	= \displaystyle\frac{1}{(2 \pi \textbf{i})^m} \oint_C
	\cdots
	\oint_C 
	\mathcal{G}_{\boldsymbol{\lambda}/\boldsymbol{\mu}}(w_1^{-1},\dots,w_m^{-1})
	%\\
	%\times
	\prod_{1 \le i < j \le m}
	\frac{w_i-w_j}{w_i-qw_j} 
	\prod_{1 \le c\le d \le m}
	\frac{1}{(w_c-y)^{g_d}}
	\prod_{i=1}^{m}
	w_i^{K-1} dw_i.
\end{multline}
Now let us specify the nonnegative integers $\{g_1,\dots,g_m\}$. Fix a partition $\nu = (\nu_1, \nu_2, \cdots , \nu_m)$ such that $|\nu| = |\boldsymbol{\lambda}|-|\boldsymbol{\mu}|$ and let $\bar{\nu} = (\bar{\nu}_1, \bar{\nu}_2, \cdots , \bar{\nu}_m)$ denote the partition obtained from $\nu$ by complementation in an $m \times (K - 1)$ box; namely, $\bar{\nu}_i = K - \nu_{m-i+1} - 1$, for all $i \in [1, m]$. We define $g_i = \bar{\nu}_i-\bar{\nu}_{i+1}$, for all $i \in [1, m]$, where $\bar{\nu}_{m+1} = -1$ by agreement. This yields
\begin{align*}
	\sum_{i=1}^m i g_i
	=
	\sum_{i=1}^m i (\bar{\nu}_i - \bar{\nu}_{i+1})
	=
	\sum_{i=1}^m \bar{\nu}_i + m
	=
	mK-|\nu|
	=
	mK-|\boldsymbol{\lambda}|+|\boldsymbol{\mu}|.
\end{align*}
Comparing \eqref{homog-y} and \eqref{limit2}, we see that the integral on the right hand side of \eqref{limit2} is necessarily independent of $y$. Noting that the contour $C$ has been chosen such that $0$ and $y$ are contained in $C$, we are at liberty to take $y=0$ in the right hand side of \eqref{limit2}: 
\begin{multline}
	\label{after-limit}
	(-y)^{-\sum_{i=1}^m i g_i}
	\lim_{q^L \rightarrow \infty}
	q^{-L \sum_{i=1}^m i g_i}
	\left[
	\lim_{x_{a_1} \rightarrow 0}
	\cdots
	\lim_{x_{a_m} \rightarrow 0}
	\mathfrak{Z}^{\mathcal{A}}_{\boldsymbol{\lambda}/\boldsymbol{\mu}}
	(x_1,\dots,x_{N+m})
	\right]^{\dagger}
	\\
	=
	\displaystyle\frac{1}{(2 \pi \textbf{i})^m} \oint_C 
	\cdots
	\oint_C 
	\mathcal{G}_{\boldsymbol{\lambda}/\boldsymbol{\mu}}(w_1^{-1},\dots,w_m^{-1})
	\prod_{1 \le i < j \le m}
	\frac{w_i-w_j}{w_i-qw_j} 
	\prod_{i=1}^m
	w_i^{K-\bar{\nu}_i - 2} dw_i
	\\
	=
	\displaystyle\frac{1}{(2 \pi \textbf{i})^m} \oint_C 
	\cdots
	\oint_C 
	\mathcal{G}_{\boldsymbol{\lambda}/\boldsymbol{\mu}}(w_1^{-1},\dots,w_m^{-1})
	\prod_{1 \le i < j \le m}
	\frac{w_i-w_j}{w_i-qw_j} 
	\prod_{i=1}^m
	w_i^{\nu_{m-i+1}-1} dw_i,
\end{multline}
where we have used the fact that $\sum_{i \le d \le m} g_d = \bar{\nu}_i - \bar{\nu}_{m + 1} = \bar{\nu}_i + 1$.

\subsection{Match With $f_{\boldsymbol{\lambda}/\boldsymbol{\mu}}^{\nu}(q)$}

After relabelling integration variables in \eqref{after-limit} by replacing $w_i$ with $w_{m-i+1}$ for all $ i \in [1, m]$, and recalling that $\mathcal{G}_{\boldsymbol{\lambda}/\boldsymbol{\mu}}$ is symmetric in its arguments, we obtain
\begin{multline}
	\label{after-limit2}
	(-y)^{-\sum_{i=1}^m i g_i}
	\lim_{q^L \rightarrow \infty}
	q^{-L \sum_{i=1}^m i g_i}
	\left[
	\lim_{x_{a_1} \rightarrow 0}
	\cdots
	\lim_{x_{a_m} \rightarrow 0}
	\mathfrak{Z}^{\mathcal{A}}_{\boldsymbol{\lambda}/\boldsymbol{\mu}}
	(x_1,\dots,x_{N+m})
	\right]^{\dagger}
	\\
	=
	\displaystyle\frac{1}{(2 \pi \textbf{i})^m} \oint_C 
	\cdots
	\oint_C
	\mathcal{G}_{\boldsymbol{\lambda}/\boldsymbol{\mu}}(w_1^{-1},\dots,w_m^{-1})
	\prod_{1 \le i < j \le m}
	 \frac{w_j-w_i}{w_j-qw_i} 
	\prod_{i=1}^{m}
	w_i^{\nu_i-1} dw_i.
\end{multline}
Matching with \eqref{f-int}, we have proved that
\begin{align}
	\label{limit-match}
	(-y)^{-\sum_{i=1}^m i g_i}
	\lim_{q^L \rightarrow \infty}
	q^{-L \sum_{i=1}^m i g_i}
	\left[
	\lim_{x_{a_1} \rightarrow 0}
	\cdots
	\lim_{x_{a_m} \rightarrow 0}
	\mathfrak{Z}^{\mathcal{A}}_{\boldsymbol{\lambda}/\boldsymbol{\mu}}
	(x_1,\dots,x_{N+m})
	\right]^{\dagger}
	=
	b_{\nu}(q) f_{\boldsymbol{\lambda}/\boldsymbol{\mu}}^{\nu}(q);
\end{align}
the expansion coefficients of \eqref{pre-fused-exp} thus emerge as a direct limit of the partition function \eqref{pf-definition2}.

\section{Limiting Procedure Applied to \texorpdfstring{$Z^{\mathcal{A}}_{\boldsymbol{\lambda}/\boldsymbol{\mu}}$}{}}
\label{ssec:pf-procedure}

We begin with some preliminaries regarding vertex weights in Sections \ref{ssec:double-fuse}--\ref{ssec:one-sat}, before turning to the explicit computation of the limits \eqref{limit-match} applied to the the partition function \eqref{Z-norm} in Sections \ref{ssec:pf-limits1}--\ref{ssec:pf-limits2}.

The analysis in this section will make use of the $(m, n) = (2, n)$ cases of the $\mathcal{R}_{x, y}^{(m; n)}$ fused weights introduced in \Cref{rijkh} and provided explicitly by \Cref{rxyml}. It will be convenient for us to alter several aspects of the notation for these weights in this section. First, these fused weights implicitly depend on two additional parameters $(L, M)$; here, we will relabel them to $(\mathsf{L}, \mathsf{M})$, in order to disambiguate them from the parameters used in previous sections of this chapter. Second, the argument for $\mathcal{R}_{x, y}^{(m; n)}$ is a quadruple of $(n + 2)$-tuples $\textbf{A}', \textbf{B}', \textbf{C}', \textbf{D}' \in \mathbb{Z}_{\ge 0}^{n + 2}$, with coordinates indexed by $[0, n + 1]$, such that $|\textbf{A}'| = \mathsf{M} = |\textbf{C}'|$ and $|\textbf{B}'| = \mathsf{L} = |\textbf{D}'|$. Here, we omit the $0$-th coordinate of these $(n + 2)$-tuples from consideration. In particular, denoting $\textbf{X}' = \big( \mathsf{M} - |\textbf{X}|, \textbf{X} \big)$ for $X \in \{ A, C \}$ and $\textbf{X}' = \big( \mathsf{L} - |\textbf{X}|, \textbf{X} \big)$ for $X \in \{ B, D \}$, we write\footnote{Recall from \Cref{rxyml} that these weights only depend on $(x, y)$ through $xy^{-1}$.} 
\begin{flalign*} 
\mathcal{R}_{x / y} (\textbf{A}, \textbf{B}; \textbf{C}, \textbf{D}) = \mathcal{R}_{x, y}^{(2; n)} (\textbf{A}', \textbf{B}'; \textbf{C}', \textbf{D}'). 
\end{flalign*}

\subsection{Doubly-Fused \texorpdfstring{$U_q \big( \widehat{\mathfrak{sl}}(2|n) \big)$}{} Model}
\label{ssec:double-fuse}

Fix four $(n+1)$-tuples $\textbf{A}, \textbf{B}, \textbf{C}, \textbf{D} \in \mathbb{Z}_{\ge 0}^{n + 1}$, with coordinates indexed by $[1, n + 1]$, such that 
\begin{align*}
	A_1,B_1,C_1,D_1 \ge 0,
	\qquad 
	\text{and}
	\qquad 
	{\rm max}\{A_j,B_j,C_j,D_j\} \le 1, 
	\qquad 
	\text{for} 
	\qquad
	j \in [2,n+1].
\end{align*}
Construct the $n$-tuple $\textbf{V} = (V_1,V_2,\dots,V_{n+1})$, where $V_1 = 0$ and $V_j = {\rm min}\{A_j,B_j,C_j,D_j\}$ for $j \in [2,n+1]$, and write $v=\sum_{i=1}^{n+1} V_i$. By \Cref{rxyml}, the doubly-fused $U_q \big( \widehat{\mathfrak{sl}}(2|n) \big)$ weights read
\begin{align}
	\label{LM-weights}
	\tikz{0.7}{
		\draw[lgray,line width=4pt,->] (-1,0) -- (1,0);
		\draw[lgray,line width=4pt,->] (0,-1) -- (0,1);
		\node[left] at (-1,0) {\tiny $\textbf{B}$};\node[right] at (1,0) {\tiny $\textbf{D}$};
		\node[below] at (0,-1) {\tiny $\textbf{A}$};\node[above] at (0,1) {\tiny $\textbf{C}$};
		%spectral parameters and spins
		\node[left] at (-1.5,0) {$(\l)$};
		\node[below] at (0,-1.5) {$(\m)$};
	}
	&=
	\mathcal{R}_z(\textbf{A},\textbf{B};\textbf{C},\textbf{D})
	\\
	\nonumber
	&=
	\bm{1}_{\textbf{A}+\textbf{B}=\textbf{C}+\textbf{D}}
	\bm{1}_{|\textbf{A}| \le \m}
	\bm{1}_{|\textbf{C}| \le \m}
	\bm{1}_{|\textbf{B}| \le \l}
	\bm{1}_{|\textbf{D}| \le \l}
	\\
	& \quad \times
	(-1)^v
	q^{\varphi(\textbf{V},\textbf{A})-v\m}
	\frac{(q^{\l-v+1}z;q)_v}{(q^{\l-\m-v}z;q)_v}
	\omega^{(\l-v,\m)}_z(\textbf{A},\textbf{B}-\textbf{V};\textbf{C},\textbf{D}-\textbf{V}).
	\nonumber
\end{align}

\noindent Here, we used the fact that $\sum_{h: V_h = 1} A_{[h + 1, n + 1]} = \varphi (\textbf{V}, \textbf{A})$ and we recalled from \eqref{omegaxy} that 
\begin{flalign}
	\label{LM-weights2}
	\begin{aligned} 
	\omega^{(\l-v,\m)}_z(\textbf{A},\textbf{B}-\textbf{V};\textbf{C},\textbf{D}-\textbf{V})
	& =
	z^{|\textbf{D}|-|\textbf{B}|}
	q^{|\textbf{A}|\l-|\textbf{D}|\m+v\m-v|\textbf{A}|}
	\\
	& \qquad \times
	\sum_{\textbf{P}}
	\Phi(\textbf{C}-\textbf{P},\textbf{C}+\textbf{D}-\textbf{V}-\textbf{P};q^{\l-v-\m}z,q^{-\m}z) \\
	& \qquad \qquad \times
	\Phi(\textbf{P},\textbf{B}-\textbf{V};q^{-\l+v}z^{-1},q^{-\l+v}),
	\end{aligned} 
\end{flalign}
with the sum over all $\textbf{P} = (P_1,P_2,\dots,P_{n+1}) \in \mathbb{Z}_{\ge 0}^{n + 1}$ such that $P_i \le {\rm min}\{C_i,B_i-V_i\}$ for all $i \in [1,n+1]$, with the function $\Phi$ given by \eqref{lambdamuxyfunction}. 

Now fix complex numbers $s, x \in \mathbb{C}$. Inserting \eqref{LM-weights2} into \eqref{LM-weights}, analytically continuing in $q^{\mathsf{M}}$, and setting $q^{-\mathsf{M}} = s^2$ and $z = s^{-1} x$ gives
\begin{flalign}
	\label{srxs} 
	\begin{aligned} 
	(- & s)^{|\textbf{D}|} \mathcal{R}_{x / s} (\textbf{A}, \textbf{B}; \textbf{C}, \textbf{D}) \\
	& = \textbf{1}_{\textbf{A} + \textbf{B} = \textbf{C} + \textbf{D}} \textbf{1}_{|\textbf{B}| \le \mathsf{L}} \textbf{1}_{|\textbf{D}| \le \mathsf{L}} \cdot (-1)^{|\textbf{D}| + v} q^{\varphi (\textbf{V}, \textbf{A})} s^{|\textbf{B}|} x^{|\textbf{D}| - |\textbf{B}|} q^{(|\textbf{A}| - |\textbf{C}|)(\mathsf{L} - v)} \displaystyle\frac{(s^{-1} q^{\mathsf{L} - v + 1} x; q)_v}{(s q^{\mathsf{L}- v} x; q)_v} \\
	& \quad \times \displaystyle\sum_{\textbf{P}} \displaystyle\frac{(s q^{\mathsf{L} - v} x; q)_{|\textbf{C}| - |\textbf{P}|} (q^{v - \mathsf{L}}; q)_{|\textbf{D}| - v}}{(s q^{\mathsf{L} - v} x; q)_{|\textbf{C}| + |\textbf{D}| - |\textbf{P}| - v}} \displaystyle\frac{(s q^{v - \mathsf{L}} x^{-1}; q)_{|\textbf{P}|} (s^{-1} x; q)_{|\textbf{B}| - |\textbf{P}| - v}}{(q^{v - \mathsf{L}}; q)_{|\textbf{B}| - v}} (s^{-1} x)^{|\textbf{P}|} \\
	& \qquad \times q^{(\mathsf{L} - v) |\textbf{P}| + \varphi (\textbf{D} - \textbf{V}, \textbf{C} - \textbf{P}) + \varphi (\textbf{B} - \textbf{V} - \textbf{P}, \textbf{P})} \displaystyle\prod_{i = 1}^{n + 1} \displaystyle\frac{(q; q)_{C_i + D_i - V_i - P_i}}{(q; q)_{C_i - P_i} (q; q)_{D_i - V_i}} \displaystyle\frac{(q; q)_{B_i - V_i}}{(q; q)_{P_i} (q; q)_{B_i - P_i}}.
	\end{aligned} 
\end{flalign}

\subsection{Recovering the Model \eqref{hybrid1}--\eqref{hybrid2}}
\label{ssec:limits=0}

We consider the limits of the weights \eqref{srxs} as $s$ tends to $0$. Let us denote this limit by
\begin{align*}
	\mathbb{R}_x(\textbf{A},\textbf{B};\textbf{C},\textbf{D})
	=
	\lim_{s\rightarrow 0}
	\left(
	(-s)^{-|\textbf{D}|}
	\mathcal{R}_z(\textbf{A},\textbf{B};\textbf{C},\textbf{D}) \Big|_{q^{\m} \mapsto s^{-2},\ z \mapsto x/s}
	\right).
\end{align*}

\noindent Inserting the limits  
\begin{flalign*}
	&\displaystyle\lim_{s \rightarrow 0} s^v (s^{-1} q^{\mathsf{L} - v + 1} x; q)_v = (-x)^v q^{(\mathsf{L} - v + 1) v + \binom{v}{2}}; \\
	& \displaystyle\lim_{s \rightarrow 0} s^{|\textbf{B}| - |\textbf{P}| - v} (s^{-1} x; q)_{|\textbf{B}| - |\textbf{P}| - v} = (-x)^{|\textbf{B}| - |\textbf{P}| - v} q^{\binom{|\textbf{B}| - |\textbf{P}| - v}{2}},
\end{flalign*} 

\noindent into \eqref{srxs}, we obtain
\begin{flalign}
	\label{s=0weights}
	\begin{aligned}
	\mathbb{R}_x(\textbf{A},\textbf{B};\textbf{C},\textbf{D})
	& =
	\bm{1}_{\textbf{A}+\textbf{B}=\textbf{C}+\textbf{D}}
	\bm{1}_{|\textbf{B}| \le \l}
	\bm{1}_{|\textbf{D}| \le \l}
	\cdot
	(-1)^{|\textbf{D}|-|\textbf{B}|+v}
	q^{\varphi(\textbf{V},\textbf{A})}
	x^{|\textbf{D}|}
	q^{(|\textbf{A}|-|\textbf{C}|)(\l-v)}
	\\
	& \qquad \times
	q^{(\l-v+1)v + \binom{v}{2}}
	q^{\varphi(\textbf{D}-\textbf{V},\textbf{C})}
	\frac{(q^{v-\l};q)_{|\textbf{D}|-v}}{(q^{v-\l};q)_{|\textbf{B}|-v}} \\
	& \qquad \times
	\sum_{\textbf{P}}
	(-1)^{|\textbf{P}|} q^{(\l-v)|\textbf{P}|} q^{-\varphi(\textbf{D}-\textbf{V},\textbf{P})}
	q^{\binom{|\textbf{B}|-v-|\textbf{P}|}{2}} 
	q^{\varphi(\textbf{B}-\textbf{V}-\textbf{P},\textbf{P})} \\
	& \qquad \qquad \times
	\prod_{i=1}^{n+1}
	\displaystyle\frac{(q; q)_{C_i+D_i-V_i-P_i}}{(q; q)_{C_i-P_i} (q; q)_{D_i - V_i}}
	\displaystyle\frac{(q; q)_{B_i-V_i}}{(q; q)_{P_i} (q; q)_{B_i - V_i - P_i}}.
	\end{aligned}
\end{flalign}

Using the bilinearity of $\varphi$, we deduce
\begin{flalign*}
(\mathsf{L} & - v + 1) v + \binom{v}{2} + \binom{|\textbf{B}| - v - |\textbf{P}|}{2} + \varphi (\textbf{B} - \textbf{V} - \textbf{P}, \textbf{P}) \\
& = \big( \mathsf{L} - |\textbf{B}| + 1 \big)v + \binom{|\textbf{B}|}{2} + \binom{|\textbf{P}| + 1}{2} - \varphi (\textbf{P}, \textbf{P}) + |\textbf{P}| \big( v - |\textbf{B}| \big) + \varphi (\textbf{B} - \textbf{V}, \textbf{P}).
\end{flalign*} 

\noindent Inserting this, together with the fact that 
\begin{align*}
	\frac{|\textbf{P}|^2}{2}
	-
	\varphi(\textbf{P},\textbf{P})
	=
	\frac{1}{2}
	\sum_{i=1}^{n+1}
	\sum_{j=1}^{n+1}
	P_i P_j
	-
	\sum_{1 \le i<j \le n+1}
	P_i P_j
	=
	\frac{1}{2}
	\sum_{i=1}^{n+1}
	P_i^2.
\end{align*}

\noindent into \eqref{s=0weights}, it becomes
\begin{flalign}
	\label{rsxs0}
	\begin{aligned} 
	\mathbb{R}_x(\textbf{A},\textbf{B};\textbf{C},\textbf{D})
	& =
	\bm{1}_{\textbf{A}+\textbf{B}=\textbf{C}+\textbf{D}}
	\bm{1}_{|\textbf{B}| \le \l}
	\bm{1}_{|\textbf{D}| \le \l}
	\cdot
	(-1)^{|\textbf{D}|-|\textbf{B}|+v}
	q^{\varphi(\textbf{V},\textbf{A})}
	x^{|\textbf{D}|}
	q^{(|\textbf{A}|-|\textbf{C}|)(\l-v)}
	\\
	& \qquad \times
	q^{(\l- |\textbf{B}| + 1)v + \binom{|\textbf{B}|}{2}}
	q^{\varphi(\textbf{D}-\textbf{V},\textbf{C})}
	\frac{(q^{v-\l};q)_{|\textbf{D}|-v}}{(q^{v-\l};q)_{|\textbf{B}|-v}} \\
	& \qquad \times
	\sum_{\textbf{P}}
	(-1)^{|\textbf{P}|} q^{\l |\textbf{P}| +\varphi(\textbf{B}-\textbf{D},\textbf{P}) - |\textbf{P}| |\textbf{B}|} \\
	& \qquad \qquad \times
	\prod_{i=1}^{n+1} q^{\binom{P_i + 1}{2}}
	\displaystyle\frac{(q; q)_{C_i+D_i-V_i-P_i}}{(q; q)_{C_i-P_i} (q; q)_{D_i - V_i}}
	\displaystyle\frac{(q; q)_{B_i-V_i}}{(q; q)_{P_i} (q; q)_{B_i - V_i - P_i}}.
	\end{aligned}
\end{flalign} 

\noindent The sum on the right side of \eqref{rsxs0} factors, thereby yielding
\begin{flalign}
	\begin{aligned}
	\label{s=0weights2}
	\mathbb{R}_x(\textbf{A},\textbf{B};\textbf{C},\textbf{D})
	& =
	\bm{1}_{\textbf{A}+\textbf{B}=\textbf{C}+\textbf{D}}
	\bm{1}_{|\textbf{B}| \le \l}
	\bm{1}_{|\textbf{D}| \le \l} \cdot 
	(-1)^{|\textbf{D}|-|\textbf{B}|+v}
	q^{\varphi(\textbf{V},\textbf{A})}
	x^{|\textbf{D}|}
	q^{(|\textbf{A}|-|\textbf{C}|)(\l-v)} \\
	& \quad \times q^{v(\l-|\textbf{B}|+1)+\binom{|\textbf{B}|}{2}}
	q^{\varphi(\textbf{D}-\textbf{V},\textbf{C})}
	\frac{(q^{v-\l};q)_{|\textbf{D}|-v}}{(q^{v-\l};q)_{|\textbf{B}|-v}}
	\\
	& \quad \times
	\prod_{i=1}^{n+1}
	\sum_{P_i \le \min\{C_i,B_i-V_i\}}
	(-1)^{P_i}
	q^{\binom{P_i + 1}{2}}
	q^{P_i(\l-B_{[i,n+1]}-D_{[1,i-1]})} \\
	& \qquad \qquad \qquad \qquad \times
	\displaystyle\frac{(q; q)_{C_i+D_i-V_i-P_i}}{(q; q)_{C_i-P_i} (q; q)_{D_i - V_i}}
	\displaystyle\frac{(q; q)_{B_i-V_i}}{(q; q)_{P_i} (q; q)_{B_i - V_i - P_i}}.
	\end{aligned}
\end{flalign}
Of the $n+1$ sums in \eqref{s=0weights2}, only the first one is non-trivial; the sums over $(P_2,\dots,P_{n+1})$ each consist of either one or two nonzero summands. Explicitly computing those final $n$ sums, we arrive at the expression
\begin{flalign}
	\label{s=0weights3}
	\begin{aligned}
	\mathbb{R}_x(\textbf{A},\textbf{B};\textbf{C},\textbf{D})
	& =
	\bm{1}_{\textbf{A}+\textbf{B}=\textbf{C}+\textbf{D}}
	\bm{1}_{|\textbf{B}| \le \l}
	\bm{1}_{|\textbf{D}| \le \l}
	\cdot
	(-1)^{|\textbf{D}|-|\textbf{B}|+v}
	x^{|\textbf{D}|}
	q^{\varphi(\textbf{V},\textbf{A}) + (|\textbf{A}|-|\textbf{C}|)(\l-v)}
	\\
	& \quad \times
	q^{v(\l-|\textbf{B}|+1)+\binom{|\textbf{B}|}{2} + \varphi(\textbf{D}-\textbf{V},\textbf{C})}
	\frac{(q^{v-\l};q)_{|\textbf{D}|-v}}{(q^{v-\l};q)_{|\textbf{B}|-v}}
	\prod_{\substack{i\in [2,n+1] \\ B_i - D_i = 1}}
	\left(1-q^{\l-B_{[i+1,n+1]}-D_{[1,i-1]}}\right)
	\\
	& \quad \times
	\sum_{P_1 \le \min\{C_1,B_1\}}
	(-1)^{P_1}
	q^{\binom{P_1 + 1}{2}}
	q^{P_1(\l-|\textbf{B}|)}
	\displaystyle\frac{(q; q)_{C_1+D_1-P_1}}{(q; q)_{C_1-P_1} (q; q)_{D_1}}
	\displaystyle\frac{(q; q)_{B_1}}{(q; q)_{P_1} (q; q)_{B_1 - P_1}}.
	\end{aligned} 
\end{flalign}

\noindent Observe the similarity between these weights and those given by \Cref{ws0}, in that they both involve a product of terms of the form $1 - u q^{- B_{[j + 1, n]} - D_{[1, j - 1]}}$ for $u$ either equal to $q^{\mathsf{L}}$ or $r^{-2}$. Indeed, it can be shown that the former degenerate to the latter under the correspondence described in \Cref{wxlx}. 

\begin{prop}
	When $\l=1$, the model \eqref{s=0weights3} reduces to the model \eqref{hybrid1}--\eqref{hybrid2}; namely
	\begin{align*}
		\mathbb{R}_x(\textbf{\emph{A}},\bm{\emph{e}}_b;\textbf{\emph{C}},\bm{\emph{e}}_d)
		\Big|_{\l=1}
		=
		w_x(\textbf{\emph{A}},b;\textbf{\emph{C}},d),
	\end{align*}
	for all $(n+1)$-tuples $\textbf{\emph{A}}$, $\textbf{\emph{C}} \in \mathbb{Z}_{\ge 0}^{n + 1}$ and integers $b, d \in [0, n + 1]$. 
\end{prop}

\begin{proof}
	
	This follows from inserting $\mathsf{L} = 1$ in \eqref{s=0weights3}. 
\end{proof}

\subsection{Two Saturated Horizontal Edges}
\label{ssec:two-sat}

We now examine what happens to the formula \eqref{s=0weights3} when the states assigned to the horizontal edges of vertices, namely $\textbf{B}$ and $\textbf{D}$, are almost saturated by bosonic arrows (which have color index $1$). More precisely, we fix integers $\mathfrak{a},\mathfrak{b},\mathfrak{c},\mathfrak{d} \ge 0$ and choose four $(n + 1)$-tuples of integers $\widetilde{\textbf{A}} = (\mathfrak{a}, \textbf{A})$, $\widetilde{\textbf{B}} = (\l-\mathfrak{b},B_2,\dots,B_{n+1})$, $\widetilde{\textbf{C}} = (\mathfrak{c},\textbf{C})$, and $\widetilde{\textbf{D}} = (\l-\mathfrak{d}, \textbf{D})$, where for each $X \in \{ A, B, C, D \}$ we write $\textbf{X} = (X_2, X_3, \ldots , X_{n + 1}) \in \mathbb{Z}_{\ge 0}^{n + 1}$. Then,
\begin{align}
	\label{new-weights}
	|\textbf{A}| = \sum_{i=2}^{n+1} A_i,
	\qquad
	|\textbf{B}| = \sum_{i=2}^{n+1} B_i,
	\qquad
	|\textbf{C}| = \sum_{i=2}^{n+1} C_i,
	\qquad
	|\textbf{D}| = \sum_{i=2}^{n+1} D_i.
\end{align}
We continue to write $\textbf{V} = (0,V_2,\dots,V_{n+1})$ with $V_j = {\rm min}\{A_j,B_j,C_j,D_j\}$ for $j \in [2,n+1]$. By arrow conservation, we must have $\mathfrak{a}-\mathfrak{b} = \mathfrak{c}-\mathfrak{d}$. By the finite capacity $\l$ of the horizontal line, we also necessarily have $\mathfrak{b} \ge |\textbf{B}|$ and $\mathfrak{d} \ge |\textbf{D}|$.

For reasons that will become clear later, we let $x = y^{-1} q^{-\l+1}$ in \eqref{s=0weights3} and multiply our weights by $y^{\l} q^{\binom{\l}{2}}$. Since we have
\begin{flalign*} 
	\binom{\mathsf{L}}{2} & + (1 - \mathsf{L}) \big( \mathsf{L} - \d + |\textbf{D}| \big) + \big( |\textbf{A}| + \a - |\textbf{C}| - \c \big) (\mathsf{L} - v) + v (\b - \textbf{B} + 1) + \binom{L - \b + |\textbf{B}|}{2}\\
	& = \mathsf{L} \big( |\textbf{A}| + |\textbf{B}| - |\textbf{C}| - |\textbf{D}| + \a - \b - \c + \d \big) + \big( |\textbf{C}| + \c - |\textbf{A}| - \a - |\textbf{B}| + \b + 1 \big) v \\
	& \qquad + |\textbf{D}| - \d + \binom{\b - |\textbf{B}| + 1}{2} =|\textbf{D}| - \d + \binom{\b - |\textbf{B}| + 1}{2} + \big( \d - |\textbf{D}| + 1 \big) v,
\end{flalign*} 

\noindent if $\textbf{A} + \textbf{B} = \textbf{C} + \textbf{D}$ and $\a - \b = \c - \d$, and since we also have
\begin{flalign*} 
	\varphi ( \widetilde{\textbf{D}}, \widetilde{\textbf{C}}) = \varphi (\textbf{D}, \textbf{C}) + (\mathsf{L} - \d) |\textbf{C}|,
\end{flalign*} 

\noindent the result of all of these choices is the following expression:
\begin{flalign}
	\label{s=0weights4}
	\begin{aligned}
	y^{\l} & q^{\binom{\l}{2}}
	\mathbb{R}_{y^{-1}q^{-\l+1}} \big(\widetilde{\textbf{A}},\widetilde{\textbf{B}};\widetilde{\textbf{C}},\widetilde{\textbf{D}} \big)
	 \\
	 & =
	{\bm 1}_{\textbf{A}+\textbf{B}=\textbf{C}+\textbf{D}}
	{\bm 1}_{\mathfrak{a}-\mathfrak{b} = \mathfrak{c}-\mathfrak{d}}
	{\bm 1}_{|\textbf{B}| \le \mathfrak{b}}
	{\bm 1}_{|\textbf{D}| \le \mathfrak{d}}
	\\
	& \quad \times
	(-1)^{|\textbf{D}|-|\textbf{B}|-\mathfrak{d}+\mathfrak{b}+v}
	y^{\mathfrak{d}-|\textbf{D}|}
	q^{\varphi(\textbf{V},\textbf{A}) + |\textbf{D}|-\mathfrak{d}+\binom{\b-|\textbf{B}| + 1}{2}+(\mathfrak{d}-|\textbf{D}|+1)v}
	\\
	& \quad \times
	q^{(\l-\mathfrak{d})|\textbf{C}|+\varphi(\textbf{D}-\textbf{V},\textbf{C})}
	\frac{(q^{v-\l};q)_{\l-\d+|\textbf{D}|-v}}{(q^{v-\l};q)_{\l-\mathfrak{b}+|\textbf{B}|-v}}
	\prod_{\substack{i\in [2,n+1] \\ B_i - D_i=1}}
	\left(1-q^{\mathfrak{d}-B_{[i+1,n+1]}-D_{[2,i-1]}}\right)
	\\
	& \quad \times
	\sum_{p=0}^{\mathfrak{c}}
	(-1)^{p}
	q^{\binom{p + 1}{2}}
	q^{p(\mathfrak{b}-|\textbf{B}|)}
	\displaystyle\frac{(q; q)_{\mathfrak{c}-\mathfrak{d}+\l-p}}{(q; q)_{\mathfrak{c}-p} (q; q)_{\mathsf{L} - \mathfrak{d}}}
	\displaystyle\frac{(q; q)_{\l-\mathfrak{b}}}{(q; q)_p (q; q)_{\mathsf{L} - \mathfrak{b} - p}}.
	\end{aligned} 
\end{flalign}
Our goal is to bring \eqref{s=0weights4} into a form whereby all $\l$ dependence occurs in the combination $q^{\l}$. Manipulating the $q$-Pochhammer and $q$-binomial functions in \eqref{s=0weights4}, we have
\begin{align*}
	\frac{(q^{v-\l};q)_{\l-\mathfrak{d}+|\textbf{D}|-v}}{(q^{v-\l};q)_{\l-\mathfrak{b}+|\textbf{B}|-v}}
	=
	\frac{(q^{-1};q^{-1})_{\mathfrak{b}-|\textbf{B}|}}{(q^{-1};q^{-1})_{\mathfrak{d}-|\textbf{D}|}},
\end{align*}
\begin{align*}
	\displaystyle\frac{(q; q)_{\mathfrak{c}-\mathfrak{d}+\l-p}}{(q; q)_{\mathfrak{c}-p} (q; q)_{\mathsf{L} - \mathfrak{d}}}
	=
	\frac{\prod_{i=1}^{\mathfrak{c}-p} (1-q^{i+\l-\mathfrak{d}})}{(q;q)_{\mathfrak{c}-p}},
	\qquad
	\displaystyle\frac{(q; q)_{\l-\mathfrak{b}}}{(q; q)_p (q; q)_{\mathsf{L} - \mathfrak{b}}}
		=
	\frac{\prod_{i=1}^{p}(1-q^{i+\l-\mathfrak{b}-p})}{(q;q)_p}.
\end{align*}
Making these substitutions in \eqref{s=0weights4}, we obtain
\begin{flalign}
	\begin{aligned}
	\label{analytic-L}
	y^{\l} q^{\binom{\l}{2}}
	\mathbb{R}_{y^{-1}q^{-\l+1}} \big(\widetilde{\textbf{A}},\widetilde{\textbf{B}};\widetilde{\textbf{C}},\widetilde{\textbf{D}} \big)
	& =
	{\bm 1}_{\textbf{A} + \textbf{B} = \textbf{C} + \textbf{D}}
	{\bm 1}_{\mathfrak{a}-\mathfrak{b} = \mathfrak{c}-\mathfrak{d}}
	{\bm 1}_{|\textbf{B}| \le \mathfrak{b}}
	{\bm 1}_{|\textbf{D}| \le \mathfrak{d}}
	(-1)^{|\textbf{D}|-|\textbf{B}|-\mathfrak{d}+\mathfrak{b}+v}
	y^{\mathfrak{d}-|\textbf{D}|} \\
	& \quad \times
	q^{\l(\mathfrak{c}+|\textbf{C}|) + \varphi(\textbf{V},\textbf{A}) + |\textbf{D}|-\mathfrak{d}+\binom{\mathfrak{b}-|\textbf{B}|+1}{2}+(\mathfrak{d}-|\textbf{D}|+1)v -\mathfrak{d}|\textbf{C}|+\varphi(\textbf{D}-\textbf{V},\textbf{C})} \\
	& \quad \times 
	\frac{(q^{-1};q^{-1})_{\mathfrak{b}-|\textbf{B}|}}{(q^{-1};q^{-1})_{\mathfrak{d}-|\textbf{D}|}}
	\prod_{\substack{i\in [2,n+1] \\ B_i - D_i=1}}
	\left(1-q^{\mathfrak{d}-B_{[i+1,n+1]}-D_{[2,i-1]}}\right)
	\\
	& \quad \times
	\sum_{p=0}^{\mathfrak{c}}
	(-1)^{p}
	q^{\binom{p + 1}{2}}
	q^{p(\mathfrak{b}-|\textbf{B}|)}
	\frac{\prod_{i=1}^{\mathfrak{c}-p} (q^{-\l}-q^{i-\mathfrak{d}})}{(q;q)_{\mathfrak{c}-p}}
	\frac{\prod_{i=1}^{p}(q^{-\l}-q^{i-\mathfrak{b}-p})}{(q;q)_p},
	\end{aligned} 
\end{flalign}
and indeed our expression depends rationally on $q^{\l}$, with no other $\l$ dependence. Treating $q^{\l}$ as an arbitrary complex parameter (by analytic continuation) we conclude, by computing the limit $q^{\l} \rightarrow \infty$: 
\begin{flalign}
\begin{aligned}
	\label{s=0weights5}
	\lim_{q^{\l} \rightarrow \infty} 
	\Big[
	q^{-\l(\mathfrak{c}+|\textbf{C}|)} &
	y^{\l} q^{\binom{\l}{2}}
	\mathbb{R}_{y^{-1}q^{-\l+1}} \big( \widetilde{\textbf{A}},\widetilde{\textbf{B}};\widetilde{\textbf{C}},\widetilde{\textbf{D}})
	\Big] \\
	& 
	=
	{\bm 1}_{\textbf{A}+\textbf{B}=\textbf{C}+\textbf{D}}
	{\bm 1}_{\mathfrak{a}-\mathfrak{b} = \mathfrak{c}-\mathfrak{d}}
	{\bm 1}_{|\textbf{B}| \le \mathfrak{b}}
	{\bm 1}_{|\textbf{D}| \le \mathfrak{d}}
	\cdot
	(-1)^{|\textbf{D}|-|\textbf{B}|-\mathfrak{d}+\mathfrak{b}+v}
	y^{\mathfrak{d}-|\textbf{D}|} \\
	& \quad \times 
	q^{\varphi(\textbf{V},\textbf{A}) + |\textbf{D}|-\mathfrak{d}+\binom{\mathfrak{b}-|\textbf{B}|+1}{2}+(\mathfrak{d}-|\textbf{D}|+1)v -\mathfrak{d}|\textbf{C}|+\varphi(\textbf{D}-\textbf{V},\textbf{C})} \\
	& \quad \times
	\frac{(q^{-1};q^{-1})_{\mathfrak{b} - |\textbf{B}|}}{(q^{-1};q^{-1})_{\mathfrak{d}-|\textbf{D}|}}
	\prod_{\substack{i\in [2,n+1] \\ B_i - D_i=1}}
	\left(1-q^{\mathfrak{d}-B_{[i+1,n+1]}-D_{[2,i-1]}}\right)
	\\
	& \quad \times
	(-1)^\mathfrak{c}
	q^{\binom{\mathfrak{c} + 1}{2} - \mathfrak{c} \mathfrak{d}}
	\sum_{p=0}^{c}
	(-1)^{p}
	q^{\binom{p + 1}{2}}
	q^{p(\mathfrak{b}-|\textbf{B}|-\mathfrak{a})}
	\frac{1}{(q;q)_{\mathfrak{c}-p}(q;q)_p},
	\end{aligned}
	\end{flalign}

\noindent where we have used the fact that for $\a - \b = \c - \d$ we have
\begin{flalign*}
	\displaystyle\sum_{i = 1}^{\c - p} (i - \d) + \displaystyle\sum_{i = 1}^p (i - \b - p) & = \binom{\c - p + 1}{2} - \d (\c - p) + \binom{p + 1}{2} - p (\b + p) \\
	& = \binom{\c + 1}{2} - \c \d + p (\d - \b - \c) = \binom{\c + 1}{2} - \c \d - p \a.
\end{flalign*}
The sum on the final line of \eqref{s=0weights5} may now be computed explicitly, using the $q$-binomial theorem:
\begin{align*}
	\sum_{p=0}^{\mathfrak{c}}
	(-1)^{p}
	q^{\binom{p+1}{2}}
	q^{p(\mathfrak{b}-|\textbf{B}|-\mathfrak{a})}
	\frac{1}{(q;q)_{\mathfrak{c}-p}(q;q)_p}
	=
	\frac{(q^{\mathfrak{b}-|\textbf{B}|-\mathfrak{a}+1};q)_{\mathfrak{c}}}{(q;q)_{\mathfrak{c}}}.
\end{align*}
Using this in \eqref{s=0weights5}, together with the fact that 
\begin{flalign*} 
\displaystyle\frac{(q^{-1}; q^{-1})_{\b - |\textbf{B}|}}{(q^{-1}; q^{-1})_{\d - |\textbf{D}|}} = (-1)^{\b - |\textbf{B}| - \d + |\textbf{D}|} q^{\binom{\d - |\textbf{D}| + 1}{2} - \binom{\b - |\textbf{B}| + 1}{2}} \displaystyle\frac{(q; q)_{\b - |\textbf{B}|}}{(q; q)_{\d - |\textbf{D}|}},	
\end{flalign*}

\noindent we arrive at our final expression
\begin{flalign*}
	\lim_{q^{\l} \rightarrow \infty} 
	\Big[
	q^{-\l(\mathfrak{c}+|\textbf{C}|)}
	& 
	y^{\l} q^{\binom{\mathsf{L}}{2}}
	\mathbb{R}_{y^{-1}q^{-\l+1}} \big( \widetilde{\textbf{A}},\widetilde{\textbf{B}};\widetilde{\textbf{C}},\widetilde{\textbf{D}} \big)
	\Big] \\
	& =
	{\bm 1}_{\textbf{A}+\textbf{B}= \textbf{C}+ \textbf{D}}
	{\bm 1}_{\mathfrak{a}-\mathfrak{b} = \mathfrak{c}-\mathfrak{d}}
	{\bm 1}_{|\textbf{B}| \le \mathfrak{b}}
	{\bm 1}_{|\textbf{D}| \le \mathfrak{d}}
	\cdot
	(-1)^{\mathfrak{c}+v}
	y^{\mathfrak{d}-|\textbf{D}|}
	q^{\chi} \\
	& \qquad \times
	\frac{(q;q)_{\mathfrak{b}-|\textbf{B}|}}{(q;q)_{\mathfrak{d}-|\textbf{D}|}}
	\cdot
	\frac{(q^{\mathfrak{b}-|\textbf{B}|-\mathfrak{a}+1};q)_{\mathfrak{c}}}{(q;q)_{\mathfrak{c}}}
	\prod_{\substack{i\in [2,n+1] \\ B_i - D_i=1}}
	\left(1-q^{\mathfrak{d}-B_{[i+1,n+1]}-D_{[2,i-1]}}\right),
\end{flalign*}
where the exponent $\chi$ is given by
\begin{align*}
	\chi
	=
	\varphi(\textbf{V},\textbf{A})+\binom{\mathfrak{d}-|\textbf{D}|}{2}+v(\mathfrak{d}-|\textbf{D}|+1)+\binom{\mathfrak{c}+1}{2} - (\mathfrak{c}+|\textbf{C}|) \mathfrak{d} + 
	\varphi(\textbf{D}-\textbf{V},\textbf{C}).
\end{align*}

\subsection{One Saturated Horizontal Edge}
\label{ssec:one-sat}

Let us consider a further case of the formula \eqref{s=0weights3}, when one horizontal edge state ($\widetilde{\textbf{D}}$) is almost saturated by bosonic arrows, while the other horizontal edge state ($\widetilde{\textbf{B}}$) is empty. Specifically, we fix integers $\mathfrak{a},\mathfrak{c},\mathfrak{d} \ge 0$ and $N > \l$, and choose $(n + 1)$-tuples $\widetilde{\textbf{A}} = (N- \mathfrak{a},\textbf{A})$, $\widetilde{\textbf{B}} = \textbf{e}_0$, $\widetilde{\textbf{C}} = (N-\l-\mathfrak{c}, \textbf{C})$, $\textbf{D} = (\l-\mathfrak{d}, \textbf{D})$, with the coordinates of $\textbf{A}, \textbf{C}, \textbf{D} \in \{ 0, 1 \}^n$ indexed by $[2, n + 1]$.

In this case, given that $\textbf{B} = \textbf{e}_0$, we have $\textbf{V} = \textbf{e}_0$. By arrow conservation, we must have $\mathfrak{a}=\mathfrak{c}+\mathfrak{d}$. By the finite capacity $\l$ of the horizontal line, we also necessarily have $\mathfrak{d} \ge |\textbf{D}|$. As we did previously, we let $x = y^{-1} q^{-\l+1}$ in \eqref{s=0weights3} and multiply our weights by $y^{\l} q^{\binom{\l}{2}}$, which results in the following expression:
\begin{flalign}
	\begin{aligned}
	\label{one-edge}
	y^{\l} q^{\binom{\l}{2}}
	\mathbb{R}_{y^{-1}q^{-\l+1}} \big( \widetilde{\textbf{A}},\textbf{e}_0; \widetilde{\textbf{C}}, \widetilde{\textbf{D}} \big)
	& =
	{\bm 1}_{\textbf{A}=\textbf{C}+\textbf{D}}
	{\bm 1}_{\a = \c+\d}
	{\bm 1}_{|\textbf{D}| \le \d}
	(-1)^{\l-\d+|\textbf{D}|}
	y^{\d-|\textbf{D}|} \\
	& \qquad \times 
	q^{\binom{\l+1}{2}+\l |\textbf{C}| + |\textbf{D}|-\d-|\textbf{C}|\d + \varphi(\textbf{D},\textbf{C})}
	\displaystyle\frac{(q; q)_{N-\a} (q^{-\l};q)_{\l-\d+|\textbf{D}|}}{(q; q)_{N-\l-\c} (q; q)_{\mathsf{L} + \c - \a}}.
	\end{aligned}
\end{flalign}
Once again, we wish to rearrange so that all $\l$ dependence occurs via $q^{\l}$. We do this making use of the relations
\begin{align*}
	& (q^{-\l};q)_{\l-\d+|\textbf{D}|}
	=
	\frac{(q^{-\l};q)_{\l}}{(q^{|\textbf{D}|-\d};q)_{\d-|\textbf{D}|}}
	=
	q^{-\binom{\l+1}{2}} (-1)^{\l}
	\cdot
	\frac{(q;q)_{\l}}{(q^{|\textbf{D}|-\d};q)_{\d-|\textbf{D}|}}; \\
	& \frac{(q;q)_{N-\a}}{(q;q)_{N-\l-\c}(q;q)_{\l-\d}}
	=
	\frac{(q;q)_{N-\a}}{(q;q)_{N-\l-\c}}
	\cdot
	\frac{(q^{\l - \d+1};q)_{\d}}{(q;q)_{\l}},
\end{align*}
which when substituted into \eqref{one-edge} produce the formula
\begin{flalign}
	\begin{aligned}
	\label{one-edge2}
	y^{\l} q^{\binom{\l}{2}}
	\mathbb{R}_{y^{-1}q^{-\l+1}}(\widetilde{\textbf{A}},\textbf{e}_0;\widetilde{\textbf{C}},\widetilde{\textbf{D}})
	& =
	{\bm 1}_{\textbf{A}= \textbf{C} + \textbf{D}}
	{\bm 1}_{\a = \c+\d}
	{\bm 1}_{|\textbf{D}| \le \d}
	\cdot
	(-1)^{\d-|\textbf{D}|}
	y^{\d-|\textbf{D}|}
	\\ 
	& \quad \times q^{\l |\textbf{C}| + |\textbf{D}|-\d-|\textbf{C}| \d + \varphi(\textbf{D},\textbf{C})}
	\frac{(q^{\l - \d+1};q)_{\d}}{(q^{|\textbf{D}| - \d};q)_{\d-|\textbf{D}|}}
	\cdot
	\frac{(q;q)_{N-\a}}{(q;q)_{N-\l - \c}}	\\
	& =
	{\bm 1}_{\textbf{A}= \textbf{C}+\textbf{D}}
	{\bm 1}_{\a = \c+\d}
	{\bm 1}_{|\textbf{D}| \le \d}
	\cdot
	q^{\l(|\textbf{C}|+\d)}
	\frac{(q;q)_{N-\a}}{(q;q)_{N-\l-\c}}
	\\
	& \quad \times 
	(-1)^{\d-|\textbf{D}|}
	y^{\d-|\textbf{D}|}
	q^{|\textbf{D}|-\d-|\textbf{C}| \d+\varphi(\textbf{D},\textbf{C})}
	\frac{\prod_{i=1}^{\d} (q^{-\l} - q^{-\d+i})}{(q^{|\textbf{D}|-\d};q)_{\d-|\textbf{D}|}}.
	\end{aligned} 
\end{flalign}
We conclude by taking the limit $q^{\l} \rightarrow \infty$:
\begin{multline}
	\label{one-edge3}
	\lim_{q^{\l} \rightarrow \infty}
	\left[
	q^{-\l(|\textbf{C}| + \mathfrak{d})}
	\frac{(q;q)_{N-\l-\mathfrak{c}}}{(q;q)_{N-\mathfrak{a}}}
	y^{\l} q^{\binom{\l}{2}}
	\mathbb{R}_{y^{-1}q^{-\l+1}} \big( \widetilde{\textbf{A}}, \widetilde{\textbf{B}}; \widetilde{\textbf{C}}, \widetilde{\textbf{D}} \big)
	\right]
	\\
	=
	{\bm 1}_{\textbf{A}=\textbf{C}+\textbf{D}}
	{\bm 1}_{\mathfrak{a} = \mathfrak{c}+ \mathfrak{d}}
	{\bm 1}_{|\textbf{D}| \le \mathfrak{d}}
	\cdot
	(-1)^{|\textbf{D}|}
	y^{\mathfrak{d}-|\textbf{D}|}
	q^{|\textbf{D}|-|\textbf{C}|d-\binom{\mathfrak{d}+1}{2}+\varphi(\textbf{D},\textbf{C})}
	\frac{1}{(q^{|\textbf{D}|-\mathfrak{d}};q)_{\mathfrak{d}-|\textbf{D}|}}.
\end{multline}

\subsection{Partition Function Limits (1)}
\label{ssec:pf-limits1}

\begin{prop}
	\label{prop:limits}
	Computing the same limits and change of variables as in \eqref{dagger}, we have
	\begin{multline}
		\label{limit*}
		\left[
		\lim_{x_{a_1} \rightarrow 0}
		\cdots
		\lim_{x_{a_m} \rightarrow 0}
		\mathfrak{Z}^{\mathcal{A}}_{\boldsymbol{\lambda}/\boldsymbol{\mu}}
		(x_1,\dots,x_{N+m})
		\right]^{\dagger}
		=
		\\
		q^{-nMN} (q; q)_N^{-1}
		\times
		\tikz{1.15}{
			\filldraw[lgray,line width=1.5pt,fill=llgray] (1.5,5.5) -- (7.5,5.5) -- (7.5,6.5) -- (1.5,6.5) -- (1.5,5.5);
			\filldraw[lgray,line width=1.5pt,fill=llgray] (1.5,3.5) -- (7.5,3.5) -- (7.5,4.5) -- (1.5,4.5) -- (1.5,3.5);
			\filldraw[lgray,line width=1.5pt,fill=llgray] (1.5,0.5) -- (7.5,0.5) -- (7.5,1.5) -- (1.5,1.5) -- (1.5,0.5);
			\filldraw[lgreen,line width=1.5pt,fill=lgreen] (6.5,0.5) -- (7.5,0.5) -- (7.5,1.5) -- (6.5,1.5) -- (6.5,0.5);
			\filldraw[lgreen,line width=1.5pt,fill=lgreen] (6.5,3.5) -- (7.5,3.5) -- (7.5,4.5) -- (6.5,4.5) -- (6.5,3.5);
			\filldraw[lgreen,line width=1.5pt,fill=lgreen] (6.5,5.5) -- (7.5,5.5) -- (7.5,6.5) -- (6.5,6.5) -- (6.5,5.5);
			\foreach\y in {0,...,7}{
				\draw[lgray,line width=1.5pt] (1.5,0.5+\y) -- (7.5,0.5+\y);
			}
			\foreach\x in {0,...,6}{
				\draw[lgray,line width=1.5pt] (1.5+\x,0.5) -- (1.5+\x,7.5);
			}
			\draw[red,line width=1.5pt] (1.5,1.5) -- (2.5,1.5) -- (2.5,2.5) -- (1.5,2.5) -- (1.5,1.5);
			\draw[red,line width=1.5pt] (1.5,2.5) -- (2.5,2.5) -- (2.5,3.5) -- (1.5,3.5) -- (1.5,2.5);
			\draw[red,line width=1.5pt] (1.5,4.5) -- (2.5,4.5) -- (2.5,5.5) -- (1.5,5.5) -- (1.5,4.5);
			\draw[red,line width=1.5pt] (1.5,6.5) -- (2.5,6.5) -- (2.5,7.5) -- (1.5,7.5) -- (1.5,6.5);
			%bottom labels
			\node[below] at (2,0.5) {\footnotesize $(N,\textbf{\emph{e}}_0)$};
			\node[below] at (3,0.5) {\footnotesize $(0,\textbf{\emph{S}}_1 (\boldsymbol{\mu}))$};
			\node[above,text centered] at (4.5,0) {$\cdots$};
			\node[above,text centered] at (5.5,0) {$\cdots$};
			\node[below] at (7,0.5) {\footnotesize $(0,\textbf{\emph{S}}_K (\boldsymbol{\mu}))$};
			%top labels
			\node[above] at (2,7.5) {\footnotesize $(0,\textbf{\emph{e}}_0)$};
			\node[above] at (3,7.5) {\footnotesize $(0,\textbf{\emph{S}}_1 (\boldsymbol{\lambda}))$};
			\node[above,text centered] at (4.5,7.5) {$\cdots$};
			\node[above,text centered] at (5.5,7.5) {$\cdots$};
			\node[above] at (7,7.5) {\footnotesize $(0,\textbf{\emph{S}}_K (\boldsymbol{\lambda}))$};
			%right labels
			\node[right] at (7.5,1) {\footnotesize $0$};
			\node[right] at (7.5,2) {\footnotesize $(L,\textbf{\emph{e}}_0)$};
			\node[right] at (7.5,3) {\footnotesize $(L,\textbf{\emph{e}}_0)$};
			\node[right] at (7.5,4) {\footnotesize $0$};
			\node[right] at (7.5,5) {\footnotesize $(L,\textbf{\emph{e}}_0)$};
			%\node at (7.6,5.6) {\huge $\vdots$};
			\node[right] at (7.5,6) {\footnotesize $0$};
			\node[right] at (7.5,7) {\footnotesize $(L,\textbf{\emph{e}}_0)$};
			%left labels
			\node[left] at (1.5,1) {\footnotesize $0$};
			\node[left] at (1.5,2) {\footnotesize $(0,\textbf{\emph{e}}_0)$};
			\node[left] at (1.5,3) {\footnotesize $(0,\textbf{\emph{e}}_0)$};
			\node[left] at (1.5,4) {\footnotesize $0$};
			\node[left] at (1.5,5) {\footnotesize $(0,\textbf{\emph{e}}_0)$};
			\node at (0.4,6) {\huge $\vdots$};
			\node[left] at (1.5,6) {\footnotesize $0$};
			\node[left] at (1.5,7) {\footnotesize $(0,\textbf{\emph{e}}_0)$};
			\draw [decorate,decoration={brace,amplitude=5pt},xshift=4pt,yshift=0pt]
			(0.5,1.5) -- (0.5,3.5) node [black,midway,xshift=-0.6cm] {\footnotesize $g_1$};
			\draw [decorate,decoration={brace,amplitude=5pt},xshift=4pt,yshift=0pt]
			(0.5,4.5) -- (0.5,5.5) node [black,midway,xshift=-0.6cm] {\footnotesize $g_2$};
			\draw [decorate,decoration={brace,amplitude=5pt},xshift=4pt,yshift=0pt]
			(0.5,6.5) -- (0.5,7.5) node [black,midway,xshift=-0.6cm] {\footnotesize $g_m$};
		}
	\end{multline}
	where we use different vertex weights throughout the lattice, as indicated by different colorings/shadings:
	\begin{itemize} 
		\item
		\tikz{0.5}{\filldraw[lgray,line width=1.5pt,fill=llgray] 
			(1.5,1.5) -- (2.5,1.5) -- (2.5,2.5) -- (1.5,2.5) -- (1.5,1.5);
			\node[below] at (2,1.5) {\footnotesize $\textbf{\emph{A}}$}; \node[left] at (1.5,2) {\footnotesize $b$}; \node[right] at (2.5,2) {\footnotesize $d$}; \node[above] at (2,2.5) {\footnotesize $\textbf{\emph{C}}$};
		} is given by \eqref{hybrid1}--\eqref{hybrid2} with $w_x(\textbf{\emph{A}},b;\textbf{\emph{C}},d) \mapsto \displaystyle\lim_{x \rightarrow 0} x^{\textbf{\emph{1}}_{b > 0} - \textbf{\emph{1}}_{d > 0} + 1} w_{1/x}(\textbf{\emph{A}},b;\textbf{\emph{C}},d)$;
		\item
		\tikz{0.5}{\filldraw[lgray,line width=1.5pt,fill=lgreen] 
			(1.5,1.5) -- (2.5,1.5) -- (2.5,2.5) -- (1.5,2.5) -- (1.5,1.5);
			\node[below] at (2,1.5) {\footnotesize $\textbf{\emph{A}}$}; \node[left] at (1.5,2) {\footnotesize $b$}; \node[right] at (2.5,2) {\footnotesize $0$}; \node[above] at (2,2.5) {\footnotesize $\textbf{\emph{C}}$};
		} is given by \eqref{hybrid1}--\eqref{hybrid2};
		\item
		\tikz{0.5}{\draw[lgray,line width=1.5pt] 
			(1.5,1.5) -- (2.5,1.5) -- (2.5,2.5) -- (1.5,2.5) -- (1.5,1.5);
			\node[below] at (2,1.5) {\footnotesize $\widetilde{\textbf{\emph{A}}}$}; \node[left] at (1.5,2) {\footnotesize $\widetilde{\textbf{\emph{B}}}$}; \node[right] at (2.5,2) {\footnotesize $\widetilde{\textbf{\emph{D}}}$}; \node[above] at (2,2.5) {\footnotesize $\widetilde{\textbf{\emph{C}}}$};
		} is given by $y^L q^{\binom{L}{2}}\cdot \mathbb{R}_{y^{-1}q^{-L+1}} \big( \widetilde{\textbf{\emph{A}}}, \widetilde{\textbf{\emph{B}}}; \widetilde{\textbf{\emph{C}}}, \widetilde{\textbf{\emph{D}}} \big)$ as in \eqref{analytic-L};
		\item
		\tikz{0.5}{\draw[red,line width=1.5pt] 
			(1.5,1.5) -- (2.5,1.5) -- (2.5,2.5) -- (1.5,2.5) -- (1.5,1.5);
			\node[below] at (2,1.5) {\footnotesize $\widetilde{\textbf{\emph{A}}}$}; \node[left] at (1.5,2) {\footnotesize $\widetilde{\textbf{\emph{B}}}$}; \node[right] at (2.5,2) {\footnotesize $\widetilde{\textbf{\emph{D}}}$}; \node[above] at (2,2.5) {\footnotesize $\widetilde{\textbf{\emph{C}}}$};
		} is given by $y^L q^{\binom{L}{2}}\cdot \mathbb{R}_{y^{-1}q^{-L+1}} \big( \widetilde{\textbf{\emph{A}}}, \widetilde{\textbf{\emph{B}}}; \widetilde{\textbf{\emph{C}}}, \widetilde{\textbf{\emph{D}}} \big)$ as in \eqref{one-edge2}.
	\end{itemize}
\end{prop}

\begin{proof}
	Recall the definition \eqref{Z-norm} of $\mathfrak{Z}^{\mathcal{A}}_{\boldsymbol{\lambda}/\boldsymbol{\mu}}(x_1,\dots,x_{N+m})$. It may be rewritten as 
	\begin{align}
		\label{rewrite}
		\mathfrak{Z}^{\mathcal{A}}_{\boldsymbol{\lambda}/\boldsymbol{\mu}}
		(x_1,\dots,x_{N+m})
		=
		q^{-nMN} (q; q)_N^{-1} \displaystyle\prod_{i = 1}^m x_{a_i}^{-1} \prod_{i=1}^{N+m} x_i^{K+1} \cdot
			Z^{\mathcal{A}}_{\boldsymbol{\lambda}/\boldsymbol{\mu}}
			(x_1^{-1},\dots,x_{N+m}^{-1}) 
	\end{align}
	where the partition function $Z^{\mathcal{A}}_{\boldsymbol{\lambda}/\boldsymbol{\mu}}
	(x_1^{-1},\dots,x_{N+m}^{-1})$ is given by \eqref{pf-definition2}. The $q$-dependent factors within \eqref{rewrite} match with those on the right hand side of \eqref{limit*}, so we may neglect them in the remainder of the proof. 
	
	Taking the limits as $x_{a_1},\dots,x_{a_m}$ tend to $ 0$ affects only the weights of rows $a_1,\dots,a_m$ within $Z^{\mathcal{A}}_{\boldsymbol{\lambda}/\boldsymbol{\mu}}$. We compute these limits by distributing the factor $x_{a_i}^K$ across the weights of the leftmost $K$ vertices of row $a_i$, prior to taking $x_{a_i}$ to $0$; we do this for each $i \in [1,m]$. After performing the limits the effective weights in row $a_i$ are given by \eqref{hybrid1}--\eqref{hybrid2} with $w_x(\textbf{A},b;\textbf{C},d) \mapsto \lim_{x \rightarrow 0} \left( x^{\textbf{1}_{b > 0} - \textbf{1}_{d > 0} + 1} \cdot w_{1/x}(\textbf{A},b;\textbf{C},d) \right)$ for any vertex not in the rightmost column of the lattice, and by \eqref{hybrid1}--\eqref{hybrid2} for the rightmost vertex within the row. The resulting rows $a_1,\dots,a_m$ ultimately become the shaded rows of \eqref{limit*}.
	
	Next we study the effect of specializing variables, as in \eqref{cov-1}--\eqref{cov-2}, within $\mathfrak{Z}^{\mathcal{A}}_{\boldsymbol{\lambda}/\boldsymbol{\mu}}$. Bearing in mind that the alphabet $(x_1,\dots,x_{N+m})$ is reciprocated within 
	$Z^{\mathcal{A}}_{\boldsymbol{\lambda}/\boldsymbol{\mu}}(x_1^{-1},\dots,x_{N+m}^{-1})$, the change of variables \eqref{cov-2} instigates fusion (recall \Cref{WeightsR}, in particular \Cref{rijkh}, \Cref{zabcd2}, and \Cref{zxy1}) of the rows $a_i+1,\dots,a_{i+1}-1$ of the lattice, for all $i \in [1,m]$; these rows get replaced by $g_i$ rows within the model \eqref{s=0weights3}. The parameter $x$ appearing in \eqref{s=0weights3} is replaced by $y^{-1} q^{-L+1}$, which is the base of the geometric progression \eqref{cov-2} modulo the aforementioned reciprocation of variables. The weights \eqref{s=0weights3} also need to be multiplied by $y^{\l} q^{\binom{\l}{2}}$, which is an artifact of the factor $x_j^{K+1}$ appearing in \eqref{rewrite}, once that factor is distributed over the $K+1$ vertices within row $j$ of $Z^{\mathcal{A}}_{\boldsymbol{\lambda}/\boldsymbol{\mu}}$, for each $j \in (a_i,a_{i+1})$. These considerations lead to the unshaded rows of \eqref{limit*}.
	
	Finally, after performing fusion each of the unshaded rows in \eqref{limit*} has a right edge state $(L,\textbf{e}_0)$; that is, it consists of $L$ bosonic arrows and no fermionic arrows. It is therefore convenient to work in terms of saturated left/right edge states as in Sections \ref{ssec:two-sat} and \ref{ssec:one-sat}.
\end{proof}

\subsection{Partition Function Limits (2)}
\label{ssec:pf-limits2}

\begin{prop}
	Computing the same limits and change of variables as in \eqref{limit2}, we have
	\begin{multline}
		\label{star}
		(-y)^{-\sum_{i=1}^m i g_i}
		\lim_{q^\l \rightarrow \infty}
		q^{-L \sum_{i=1}^m i g_i}
		\left[
		\lim_{x_{a_1} \rightarrow 0}
		\cdots
		\lim_{x_{a_m} \rightarrow 0}
		\mathfrak{Z}^{\mathcal{A}}_{\boldsymbol{\lambda}/\boldsymbol{\mu}}
		(x_1,\dots,x_{N+m})
		\right]^{\dagger}
		=
		\\
		\frac{(-1)^{|\bar{\nu}|}}{(q;q)_m}
		\times
		\tikz{1.15}{
			\filldraw[lgray,line width=1.5pt,fill=llgray] (2.5,5.5) -- (7.5,5.5) -- (7.5,6.5) -- (2.5,6.5) -- (2.5,5.5);
			\filldraw[lgray,line width=1.5pt,fill=llgray] (2.5,3.5) -- (7.5,3.5) -- (7.5,4.5) -- (2.5,4.5) -- (2.5,3.5);
			\filldraw[lgray,line width=1.5pt,fill=llgray] (2.5,0.5) -- (7.5,0.5) -- (7.5,1.5) -- (2.5,1.5) -- (2.5,0.5);
			\filldraw[lgreen,line width=1.5pt,fill=lgreen] (6.5,0.5) -- (7.5,0.5) -- (7.5,1.5) -- (6.5,1.5) -- (6.5,0.5);
			\filldraw[lgreen,line width=1.5pt,fill=lgreen] (6.5,3.5) -- (7.5,3.5) -- (7.5,4.5) -- (6.5,4.5) -- (6.5,3.5);
			\filldraw[lgreen,line width=1.5pt,fill=lgreen] (6.5,5.5) -- (7.5,5.5) -- (7.5,6.5) -- (6.5,6.5) -- (6.5,5.5);
			\foreach\y in {0,...,7}{
				\draw[lgray,line width=1.5pt] (2.5,0.5+\y) -- (7.5,0.5+\y);
			}
			\foreach\x in {1,...,6}{
				\draw[lgray,line width=1.5pt] (1.5+\x,0.5) -- (1.5+\x,7.5);
			}
			%\draw[red,line width=1.5pt] (1.5,1.5) -- (2.5,1.5) -- (2.5,2.5) -- (1.5,2.5) -- (1.5,1.5);
			%\draw[red,line width=1.5pt] (1.5,2.5) -- (2.5,2.5) -- (2.5,3.5) -- (1.5,3.5) -- (1.5,2.5);
			%\draw[red,line width=1.5pt] (1.5,4.5) -- (2.5,4.5) -- (2.5,5.5) -- (1.5,5.5) -- (1.5,4.5);
			%\draw[red,line width=1.5pt] (1.5,6.5) -- (2.5,6.5) -- (2.5,7.5) -- (1.5,7.5) -- (1.5,6.5);
			%bottom labels
			\node[below] at (3,0.5) {\footnotesize $(0,\textbf{\emph{S}}_1 (\boldsymbol{\mu}))$};
			\node[above,text centered] at (4.5,0) {$\cdots$};
			\node[above,text centered] at (5.5,0) {$\cdots$};
			\node[below] at (7,0.5) {\footnotesize $(0,\textbf{\emph{S}}_K (\boldsymbol{\mu}))$};
			%top labels
			\node[above] at (3,7.5) {\footnotesize $(0,\textbf{\emph{S}}_1 (\boldsymbol{\lambda}))$};
			\node[above,text centered] at (4.5,7.5) {$\cdots$};
			\node[above,text centered] at (5.5,7.5) {$\cdots$};
			\node[above] at (7,7.5) {\footnotesize $(0,\textbf{\emph{S}}_K (\boldsymbol{\lambda}))$};
			%right labels
			\node[right] at (7.5,1) {\footnotesize $0$};
			\node[right] at (7.5,2) {\footnotesize $(L,\textbf{\emph{e}}_0)$};
			\node[right] at (7.5,3) {\footnotesize $(L,\textbf{\emph{e}}_0)$};
			\node[right] at (7.5,4) {\footnotesize $0$};
			\node[right] at (7.5,5) {\footnotesize $(L,\textbf{\emph{e}}_0)$};
			%\node at (7.6,5.6) {\huge $\vdots$};
			\node[right] at (7.5,6) {\footnotesize $0$};
			\node[right] at (7.5,7) {\footnotesize $(L,\textbf{\emph{e}}_0)$};
			%left labels
			\node[left] at (2.5,1) {\footnotesize $1$};
			\node[left] at (2.5,2) {\footnotesize $(L,\textbf{\emph{e}}_0)$};
			\node[left] at (2.5,3) {\footnotesize $(L,\textbf{\emph{e}}_0)$};
			\node[left] at (2.5,4) {\footnotesize $1$};
			\node[left] at (2.5,5) {\footnotesize $(L,\textbf{\emph{e}}_0)$};
			\node at (0.4,6) {\huge $\vdots$};
			\node[left] at (2.5,6) {\footnotesize $1$};
			\node[left] at (2.5,7) {\footnotesize $(L-m,\textbf{e}_0)$};
			%%paths
			%\draw[line width=1pt,->] (2,0.5)--(2,1)--(2.5,1);
			%\draw[line width=2pt,->] (2,0.5)--(2,2)--(2.5,2);
			%\draw[line width=2pt,->] (2,0.5)--(2,3)--(2.5,3);
			%\draw[line width=1pt,->] (2,0.5)--(2,4)--(2.5,4);
			%\draw[line width=2pt,->] (2,0.5)--(2,5)--(2.5,5);
			%\draw[line width=1pt,->] (2,0.5)--(2,6)--(2.5,6);
			%\draw[line width=2pt,->] (2,0.5)--(2,7)--(2.5,7);
			\draw [decorate,decoration={brace,amplitude=5pt},xshift=4pt,yshift=0pt]
			(0.5,1.5) -- (0.5,3.5) node [black,midway,xshift=-0.6cm] {\footnotesize $g_1$};
			\draw [decorate,decoration={brace,amplitude=5pt},xshift=4pt,yshift=0pt]
			(0.5,4.5) -- (0.5,5.5) node [black,midway,xshift=-0.6cm] {\footnotesize $g_2$};
			\draw [decorate,decoration={brace,amplitude=5pt},xshift=4pt,yshift=0pt]
			(0.5,6.5) -- (0.5,7.5) node [black,midway,xshift=-0.6cm] {\footnotesize $g_m$};
		}
	\end{multline}
	where the vertices in \eqref{star} are the same as those used in \Cref{prop:limits}, but with $y = 1$.
\end{prop}

\begin{proof}
	We study the limit as $q^{\l}$ tends to $\infty$ of the partition function \eqref{limit*}. The first thing to observe is the cancellation of the factors $q^{-nMN}$ and $(q; q)_N^{-1}$ multiplying \eqref{limit*} with certain terms in the weights of the vertices. 
	
	The factor $q^{-nMN}$ cancels perfectly with $q^{\l |\textbf{C}|}$ present in \eqref{analytic-L} and \eqref{one-edge2}. Indeed, the factor $q^{\l |\textbf{C}|}$ contributes one power of $q^{\l}$ for every fermionic path which passes through the top of a vertex \tikz{0.7}{\draw[lgray,line width=1.5pt] 
		(1.5,1.5) -- (2.5,1.5) -- (2.5,2.5) -- (1.5,2.5) -- (1.5,1.5);
		%\node[below] at (2,1.5) {\footnotesize $\textbf{A}$}; \node[left] at (1.5,2) {\footnotesize $\textbf{B}$}; \node[right] at (2.5,2) {\footnotesize $\textbf{D}$}; \node[above] at (2,2.5) {\footnotesize $\textbf{C}$};
	} or \tikz{0.7}{\draw[red,line width=1.5pt] 
		(1.5,1.5) -- (2.5,1.5) -- (2.5,2.5) -- (1.5,2.5) -- (1.5,1.5);
		%\node[below] at (2,1.5) {\footnotesize $\textbf{A}$}; \node[left] at (1.5,2) {\footnotesize $\textbf{B}$}; \node[right] at (2.5,2) {\footnotesize $\textbf{D}$}; \node[above] at (2,2.5) {\footnotesize $\textbf{C}$};
	}\ . This leads to $n M$ powers of $q^{\l}$ per unshaded row, and there are $\sum_{i=1}^{n} g_i$ such rows; the total contribution from the $q^{\l |\textbf{C}|}$ factors is thus $q^{n M \cdot \l \sum_{j=1}^{n} g_j} = q^{n MN}$ as claimed. The factor $(q;q)_N^{-1}$ cancels perfectly with $(q;q)_{N-\a} (q;q)_{N-\l-\c}^{-1}$ present in \eqref{one-edge2} and the factor $(1-q^j)$ arising from vertices \tikz{0.7}{\filldraw[lgray,line width=1.5pt,fill=llgray] 
		(1.5,1.5) -- (2.5,1.5) -- (2.5,2.5) -- (1.5,2.5) -- (1.5,1.5);
		\node[below] at (2,1.5) {\footnotesize $(j,\textbf{e}_0)$}; \node[left] at (1.5,2) {\footnotesize $0$}; \node[right] at (2.5,2) {\footnotesize $1$}; \node[above] at (2,2.5) {\footnotesize $(j-1,\textbf{e}_0)$};} in the leftmost column of \eqref{limit*}. Indeed, in every configuration of the leftmost column, these factors must always telescope to yield $(q;q)_N$, in view of the top/bottom boundary conditions of this column.

	After cancelling $q^{\l |\textbf{C}|}$ from \eqref{analytic-L}, the renormalized weights are polynomials in $q^{\l}$ of degree $\c$. Similarly, after cancelling $q^{\l |\textbf{C}|} (q;q)_{N-\a} (q;q)_{N-\l-\c}^{-1}$ from \eqref{one-edge2}, the renormalized weights are polynomial in $q^{\l}$ of degree $\d$. Since we wish to take the limit $q^{\l} \rightarrow \infty$, we now examine which lattice configurations contribute maximal degree in $q^{\l}$. 
	
	The leftmost column can contribute at most degree $m$ in $q^L$; this occurs when a bosonic path turns right into each of the $m$ shaded rows of \eqref{limit*}, leaving exactly $m$ ``vacancies'' distributed across the unshaded rows, and results from taking the top degree term of the renormalized weights \eqref{one-edge2} for every vertex
	\tikz{0.7}{\draw[red,line width=1.5pt] 
		(1.5,1.5) -- (2.5,1.5) -- (2.5,2.5) -- (1.5,2.5) -- (1.5,1.5);
		%\node[below] at (2,1.5) {\footnotesize $\textbf{A}$}; \node[left] at (1.5,2) {\footnotesize $\textbf{B}$}; \node[right] at (2.5,2) {\footnotesize $\textbf{D}$}; \node[above] at (2,2.5) {\footnotesize $\textbf{C}$};
	} in the column. If we choose the configuration with all $m$ vacancies occurring in the highest unshaded row, we obtain the following picture:
	
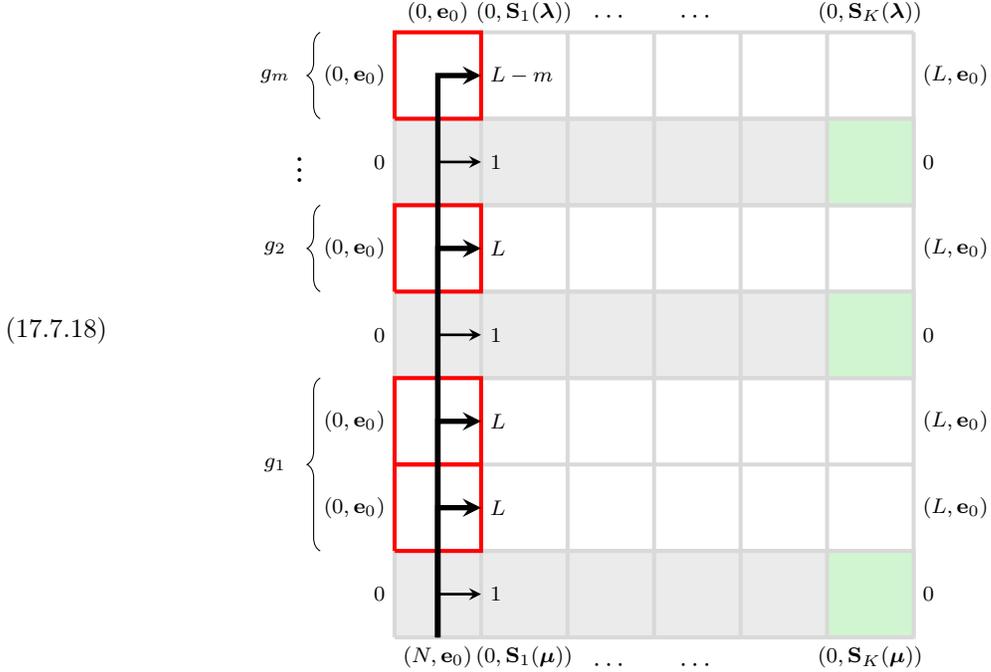
\begin{align}
		\label{dag}
		\tikz{1.15}{
			\filldraw[lgray,line width=1.5pt,fill=llgray] (1.5,5.5) -- (7.5,5.5) -- (7.5,6.5) -- (1.5,6.5) -- (1.5,5.5);
			\filldraw[lgray,line width=1.5pt,fill=llgray] (1.5,3.5) -- (7.5,3.5) -- (7.5,4.5) -- (1.5,4.5) -- (1.5,3.5);
			\filldraw[lgray,line width=1.5pt,fill=llgray] (1.5,0.5) -- (7.5,0.5) -- (7.5,1.5) -- (1.5,1.5) -- (1.5,0.5);
			\filldraw[lgreen,line width=1.5pt,fill=lgreen] (6.5,0.5) -- (7.5,0.5) -- (7.5,1.5) -- (6.5,1.5) -- (6.5,0.5);
			\filldraw[lgreen,line width=1.5pt,fill=lgreen] (6.5,3.5) -- (7.5,3.5) -- (7.5,4.5) -- (6.5,4.5) -- (6.5,3.5);
			\filldraw[lgreen,line width=1.5pt,fill=lgreen] (6.5,5.5) -- (7.5,5.5) -- (7.5,6.5) -- (6.5,6.5) -- (6.5,5.5);
			\foreach\y in {0,...,7}{
				\draw[lgray,line width=1.5pt] (1.5,0.5+\y) -- (7.5,0.5+\y);
			}
			\foreach\x in {0,...,6}{
				\draw[lgray,line width=1.5pt] (1.5+\x,0.5) -- (1.5+\x,7.5);
			}
			\draw[red,line width=1.5pt] (1.5,1.5) -- (2.5,1.5) -- (2.5,2.5) -- (1.5,2.5) -- (1.5,1.5);
			\draw[red,line width=1.5pt] (1.5,2.5) -- (2.5,2.5) -- (2.5,3.5) -- (1.5,3.5) -- (1.5,2.5);
			\draw[red,line width=1.5pt] (1.5,4.5) -- (2.5,4.5) -- (2.5,5.5) -- (1.5,5.5) -- (1.5,4.5);
			\draw[red,line width=1.5pt] (1.5,6.5) -- (2.5,6.5) -- (2.5,7.5) -- (1.5,7.5) -- (1.5,6.5);
			%bottom labels
			\node[below] at (2,0.5) {\footnotesize $(N,\textbf{e}_0)$};
			\node[below] at (3,0.5) {\footnotesize $(0,\textbf{S}_1 (\boldsymbol{\mu}))$};
			\node[above,text centered] at (4,0) {$\cdots$};
			\node[above,text centered] at (5,0) {$\cdots$};
			\node[below] at (7,0.5) {\footnotesize $(0,\textbf{S}_K (\boldsymbol{\mu}))$};
			%top labels
			\node[above] at (2,7.5) {\footnotesize $(0,\textbf{e}_0)$};
			\node[above] at (3,7.5) {\footnotesize $(0,\textbf{S}_1 (\boldsymbol{\lambda}))$};
			\node[above,text centered] at (4,7.5) {$\cdots$};
			\node[above,text centered] at (5,7.5) {$\cdots$};
			\node[above] at (7,7.5) {\footnotesize $(0,\textbf{S}_K (\boldsymbol{\lambda}))$};
			%right labels
			\node[right] at (7.5,1) {\footnotesize $0$};
			\node[right] at (7.5,2) {\footnotesize $(L,\textbf{e}_0)$};
			\node[right] at (7.5,3) {\footnotesize $(L,\textbf{e}_0)$};
			\node[right] at (7.5,4) {\footnotesize $0$};
			\node[right] at (7.5,5) {\footnotesize $(L,\textbf{e}_0)$};
			%\node at (7.6,5.6) {\huge $\vdots$};
			\node[right] at (7.5,6) {\footnotesize $0$};
			\node[right] at (7.5,7) {\footnotesize $(L,\textbf{e}_0)$};
			%left labels
			\node[left] at (1.5,1) {\footnotesize $0$};
			\node[left] at (1.5,2) {\footnotesize $(0,\textbf{e}_0)$};
			\node[left] at (1.5,3) {\footnotesize $(0,\textbf{e}_0)$};
			\node[left] at (1.5,4) {\footnotesize $0$};
			\node[left] at (1.5,5) {\footnotesize $(0,\textbf{e}_0)$};
			\node at (0.4,6) {\huge $\vdots$};
			\node[left] at (1.5,6) {\footnotesize $0$};
			\node[left] at (1.5,7) {\footnotesize $(0,\textbf{e}_0)$};
			\node[right] at (2.5,1) {\footnotesize $1$};
			\node[right] at (2.5,2) {\footnotesize $L$};
			\node[right] at (2.5,3) {\footnotesize $L$};
			\node[right] at (2.5,4) {\footnotesize $1$};
			\node[right] at (2.5,5) {\footnotesize $L$};
			%\node at (2.4,5.6) {\huge $\vdots$};
			\node[right] at (2.5,6) {\footnotesize $1$};
			\node[right] at (2.5,7) {\footnotesize $L-m$};
			%paths
			\draw[line width=1pt,->] (2,0.5)--(2,1)--(2.5,1);
			\draw[line width=2pt,->] (2,0.5)--(2,2)--(2.5,2);
			\draw[line width=2pt,->] (2,0.5)--(2,3)--(2.5,3);
			\draw[line width=1pt,->] (2,0.5)--(2,4)--(2.5,4);
			\draw[line width=2pt,->] (2,0.5)--(2,5)--(2.5,5);
			\draw[line width=1pt,->] (2,0.5)--(2,6)--(2.5,6);
			\draw[line width=2pt,->] (2,0.5)--(2,7)--(2.5,7);
			\draw [decorate,decoration={brace,amplitude=5pt},xshift=4pt,yshift=0pt]
			(0.5,1.5) -- (0.5,3.5) node [black,midway,xshift=-0.6cm] {\footnotesize $g_1$};
			\draw [decorate,decoration={brace,amplitude=5pt},xshift=4pt,yshift=0pt]
			(0.5,4.5) -- (0.5,5.5) node [black,midway,xshift=-0.6cm] {\footnotesize $g_2$};
			\draw [decorate,decoration={brace,amplitude=5pt},xshift=4pt,yshift=0pt]
			(0.5,6.5) -- (0.5,7.5) node [black,midway,xshift=-0.6cm] {\footnotesize $g_m$};
		}
	\end{align}
	It is straightforward to show that the remaining columns in \eqref{dag} contribute degree $\sum_{i=1}^m i g_i - m$ in $q^L$ (again, this results from taking the top degree term of the renormalized weights \eqref{analytic-L} for every vertex \tikz{0.7}{\draw[lgray,line width=1.5pt] 
		(1.5,1.5) -- (2.5,1.5) -- (2.5,2.5) -- (1.5,2.5) -- (1.5,1.5);
		%\node[below] at (2,1.5) {\footnotesize $\textbf{A}$}; \node[left] at (1.5,2) {\footnotesize $\textbf{B}$}; \node[right] at (2.5,2) {\footnotesize $\textbf{D}$}; \node[above] at (2,2.5) {\footnotesize $\textbf{C}$};
	} in the lattice). Furthermore, it easy to show that any other configuration of the leftmost column would lead to smaller overall degree in $q^L$. It follows that when taking the limit $q^L \rightarrow \infty$, only the configuration \eqref{dag} needs to be considered, and it has degree 
	$\left( \sum_{i=1}^m i g_i - m \right) + m = \sum_{i=1}^{n} i g_i$ in $q^L$.
	
	Multiplying \eqref{limit*} by $q^{-L \sum_{i=1}^m i g_i} (-y)^{-\sum_{i=1}^m i g_i}$ and letting $q^L$ tend to $\infty$, we thus obtain \eqref{star}. The factor of $(-1)^{|\bar{\nu}| + m}$ is due to $(-1)^{-\sum_{i=1}^m i g_i}$ that we multiply by, and the factor of $\frac{(-1)^{m}}{(q;q)_m}$ is the result of computing the weight of the top vertex in the leftmost column of \eqref{dag} (all other vertices in the leftmost column of \eqref{dag} can be considered to have weight $1$, in view of previously cancelled-off factors).
\end{proof}

\subsection{Final Matching}

Given the analysis in the previous sections, we can now quickly establish \Cref{thm:comb}. 

\begin{proof}[Proof of \Cref{thm:comb}]

Comparing \eqref{limit-match} and \eqref{star} now yields equation \eqref{main-formula}. In producing the match we perform a complementation by $L$ of bosonic states that live on left and right edges of vertices; namely, such a state $(L-\b,\textbf{B})$ in \eqref{star} gets replaced by $(\b,\textbf{B})$ in \eqref{main-formula}. We also implement an analogous complementation for labelling of left and right edges of shaded tiles in going from \eqref{star} to \eqref{main-formula} (namely, we replace an arrow configuration $(\textbf{A}, b; \textbf{C}, d)$ with $(\textbf{A}, 1 - b; \textbf{C}, 1 - d)$). It is quickly verified under this complementation that the white tiles of \Cref{prop:limits} transform at $y = 1$ into the unshaded ones of \Cref{ssec:weights1} and that both the green and gray tiles of \Cref{prop:limits} transform into the shaded oned ones of \Cref{ssec:weights2}.

Shaded rows in \eqref{star} occur at row $m-i+\sum_{j=i}^m g_j$ for each $i \in [1,m]$, where rows are counted top to bottom, starting at zero. This matches with the vector 
$\mathfrak{n}$ in \eqref{main-formula}, since $g_j = \bar{\nu}_j - \bar{\nu}_{j + 1}$ and $\bar{\nu}_{m + 1} = -1$ together imply $\sum_{j=i}^m g_j = \bar{\nu}_i + 1$.
\end{proof}

\printindex

\end{document}